\documentclass[12pt]{amsart}
\usepackage{amsmath, amsthm, amssymb}
\usepackage{fullpage}
\usepackage{xcolor}
\usepackage{hyperref}
\usepackage{soul}
\usepackage{enumitem}
\usepackage[normalem]{ulem}
\usepackage{mathtools}
\mathtoolsset{showonlyrefs}

\newtheorem{theorem}{Theorem}
\newtheorem{lemma}{Lemma}
\newtheorem{proposition}[lemma]{Proposition}
\newtheorem{corollary}[lemma]{Corollary}
\newtheorem{definition}[lemma]{Definition}
\newtheorem{remark}[lemma]{Remark}
\newtheorem{conjecture}{Conjecture}
\numberwithin{lemma}{section}

\numberwithin{equation}{section}

\newcommand{\R}{{\mathbb R}}

\renewcommand{\H}{{\mathcal H }}
\newcommand{\B}{{\mathcal B }}
\newcommand{\dB}{{\dot{\mathcal B }}}

\renewcommand{\R}{\mathbb R}

\newcommand{\calL}{\mathcal L}

\newcommand{\mfQ}{\mathfrak Q}

\renewcommand{\L}{\mathfrak L}
\newcommand{\M}{\mathfrak M}

\newcommand{\mfl}{\mathfrak l}
\newcommand{\mfm}{\mathfrak m}

\newcommand{\bR}{\mathbf R}

\newcommand{\bI}{\mathbf I}
\newcommand{\bJ}{\mathbf J}

\newcommand{\xis}{\xi^\sharp}
\newcommand{\calv}{\mathfrak{v}}

\renewcommand{\L}{\mathfrak L}
\newcommand{\mfu}{\mathfrak u}

\newcommand{\calN}{\mathfrak{N}}
\newcommand{\bal}{\text{bal}}
\newcommand{\mfq}{\mathfrak q}

\newcommand{\bu}{{\bar u}}
\newcommand{\bv}{{\bar v}}
\newcommand{\bfu}{{\mathbf u}}
\newcommand{\bfv}{{\mathbf v}}
\newcommand{\bfw}{{\mathbf w}}

\newcommand{\dH}{{\dot{\mathcal H} }}

\newcommand{\bw}{{\bar w}}

\renewcommand{\r}{{\mathbf r}}

\renewcommand{\b}{{\mathbf b}}

\newcommand{\la}{\langle}
\newcommand{\ra}{\rangle}

\newcommand{\ol}{\overline}
\newcommand{\ms}{\fm^\sharp}

\newcommand{\calR}{\mathcal{R}}

\newcommand{\gv}{\varphi}
\newcommand{\dn}{\mathfrak{g}}
\newcommand{\DN}{\mathfrak{G}}

\newcommand{\fP}{{\mathfrak P}}
\newcommand{\fM}{{\mathfrak M}}
\newcommand{\fp}{{\mathfrak p}}
\newcommand{\tfp}{{\tilde{\mathfrak p}}}
\newcommand{\fm}{{\mathfrak m}}
\newcommand{\ff}{{\mathfrak f}}
\newcommand{\fF}{{\mathfrak F}}
\newcommand{\fe}{{\mathfrak e}}
\newcommand{\fE}{{\mathfrak E}}

\newcommand{\LPs}{P^\sharp}
\newcommand{\lps}{p^\sharp}
\newcommand{\gs}{g^\sharp}
\newcommand{\tps}{{\tilde p}^\sharp}
\newcommand{\tPs}{{\tilde P}^\sharp}

\newcommand{\fsj}{\ff^{\sharp,j}}

\newcommand{\Fsj}{\fF^{\sharp,j}}

\newcommand{\Fsk}{\fF^{\sharp,k}}

\newcommand{\esj}{\fe^{\sharp,jk}}

\newcommand{\Esj}{\fE^{\sharp,jk}}

\newcommand{\step}{\mathfrak c}

\newcommand{\psj}{\fp^{\sharp,j}}
\newcommand{\Ms}{\fM^\sharp}
\newcommand{\Ps}{\fP^\sharp}
\newcommand{\Psj}{\fP^{\sharp,j}}
\newcommand{\Psk}{\fP^{\sharp,k}}
\newcommand{\Mss}{\fM^{\sharp\sharp}}
\newcommand{\Pss}{\fP^{\sharp\sharp}}
\newcommand{\Pssj}{\fP^{\sharp\sharp,j}}
\newcommand{\Msnl}{\fM^{\sharp,NL}}

\newcommand{\Psnlj}{\fP^{\sharp,j,NL}}
\newcommand{\Mslin}{\fM^{\sharp,lin}}
\newcommand{\Pslin}{\fP^{\sharp,lin}}

\newcommand{\paraequal}{\overset{bal}{=}}
\newcommand{\coeff}{\mathcal B}
\newcommand{\coeffu}{\mathcal B^{ub}}
\newcommand{\coeffs}{\mathcal B^\sharp}
\newcommand{\coeffuu}{\mathcal B^{\sharp,ub}}
\newcommand{\coeffb}{\mathcal B^{bal}}

\newcommand{\tcoeffu}{\tilde {\mathcal B}^{ub}}
\newcommand{\tcoeffs}{\tilde {\mathcal B}^\sharp}
\newcommand{\tcoeffb}{\tilde {\mathcal B}^{bal}}

\newcommand{\rT}{T^R_{D_t}}

\begin{document}

\title{Rough global solutions for 3D gravity waves with nonlocalized data}

\author{Mihaela Ifrim}
\address{Department of Mathematics, University of Wisconsin, Madison}
\email{ifrim@wisc.edu}

\author{Ben Pineau}
\address{Courant Institute for Mathematical Sciences\\
New York University
} \email{brp305@nyu.edu}

\author{ Daniel Tataru}
\address{Department of Mathematics, University of California at Berkeley}
\email{tataru@math.berkeley.edu}

\vspace*{-2cm}
\begin{abstract}
\vspace{1in}
In this article we study the problem of local and global well-posedness in the small data, low regularity regime for gravity waves in a fluid of infinite depth and infinite width, in spatial dimension $n \geq 3$, which corresponds to the fluid interface 
evolution in $d=n-1$ dimensions. We do this 
under the assumption that the initial data is small in a critical homogeneous Besov space $\dot B^{s_c}_{2,1} \times \dot B^{s_c+\frac12}_{2,1}$ in $\R^d$.

The main local well-posedness result provides  Hadamard well-posedness in Sobolev spaces at $1/12$ derivative above scaling in dimension $n = 3$, in any Sobolev space above scaling  in  dimension $n = 4$, and in the scale invariant Besov space in dimension $n \geq 5$. This represents the first time the critical threshold was reached for any genuinely quasilinear dispersive flow.

The main global well-posedness results asserts that the above local solutions are global in time, under
\emph{the sole} assumption of smallness for the initial data in the critical Besov space. Moreover, we show that they satisfy a large family of translation invariant
linear and bilinear bounds, globally in time. To put this result in the context of existing results, this represents the first global well-posedness result in water waves which does not simultaneously assume both high regularity and localization for the initial data.

Critical regularity is out of reach of existing methods. So is the case of nonlocalized data. To attain our results in both of these settings, we 
develop an entirely new  set of tools for quasilinear dispersive flows: (a) a \emph{layered paradifferential expansion} of the water wave equations, (b) a new class of bilinear estimates, which we call \emph{translation averaged bilinear bounds}, (c) a new class of \emph{paratransported symbols}, whose 
spatial dependence fits into an appropriate class of \emph{paracontrolled functions}, as well as (d) 
an efficient strategy to handle the presence of the transport term in the proof of the bilinear estimates, via \emph{advected interaction Morawetz} techniques.

One may also view  our global result as the proof, in the context of gravity waves in deep water, of    
the broad global well-posedness conjecture, proposed a few years ago by two of the authors in the general setting of cubic quasilinear dispersive equations. 
The conjecture, called \emph{The Nonlocalized Data Global Wellposedness Conjecture},  applied to problems in two and higher space dimensions, asserts that global well-posedness and scattering holds for small initial data in Sobolev spaces.  As a first validation in dimension $d\geq 2$, it was initially proved by the first and third authors for  quasilinear Schr\"odinger flows, and  then, jointly with the second author, for the case of quasilinear ultrahyperbolic Schr\"odinger flows. The $d=1$ version the conjecture was also validated  for some models with general dispersion relations by the first and third authors, but not in higher dimensions.  

By comparison, the water wave problem is far harder. The equations are fully nonlinear and
nonlocal, with a dispersion relation not of Schr\"odinger type, the nonlinearities are quadratic rather than cubic, and the
solutions are transported by a strong advection field; each of these alone
defeats the methods used for the Schr\"odinger models.

\bigskip

\end{abstract}


\maketitle

\newpage

\setcounter{tocdepth}{1}
\tableofcontents


\section{Introduction}
In this article, we consider the problem of local and global well-posedness in the small data, low regularity regime for gravity waves in a fluid of infinite depth and infinite width, in spatial dimension $n \geq 3$. The governing equations within the fluid domain, which will be described shortly, are the \emph{incompressible Euler equations}, with appropriate boundary conditions on the fluid surface, which
can be viewed as a free boundary.

Under the assumption that the fluid is also \emph{irrotational}, the
velocity of the fluid can be described as the gradient of a harmonic velocity potential,
which is uniquely determined by its trace on the surface. Thus, the original $n$-dimensional Euler flow can be reduced to a $d = n-1$ dimensional flow modeling the evolution of the free surface together
with the trace of the velocity potential on the free surface. This is commonly referred to as the
\emph{water wave problem}, which is a fully nonlinear, nonlocal dispersive evolution.

Our objective in this article is to study both the local well-posedness question at low regularity,
and the global well-posedness under the assumption that the initial data is small, but without any localization assumptions. Both of these problems have been intensely studied in the last 25 years or so, but with global well-posedness primarily considered in the restrictive context where the initial data is assumed to be both smooth, small and localized. 

Our results simultaneously provide a drastic improvement for both the local and the global theory for these equations, reaching  sharp Sobolev thresholds.
More precisely, our main results in this article are as follows:
\begin{enumerate}[label=\roman*)]
 \item We sharply improve the local well-posedness theory for the gravity wave problem, $1/12$ derivatives above scaling in dimension $n = 3$, up to scaling 
  in dimension $n = 4$, and at scaling in dimension $n \geq 5$. 
\item We prove global well-posedness
for initial data as above which is small
in a critical Besov space, without any localization assumption.
\end{enumerate}

The global result also provides a full family of global in time bounds for the solutions, including uniform energy bounds, Strichartz type estimates as well as bilinear $L^2_{t,x}$ bounds, which capture the dispersive character of the solutions in a translation invariant setting. In addition, for the local result we also obtain local in time (global if $n \geq 5$) difference bounds for the solutions of the same type but in a weaker topology. 

 This is the first work where a local well-posedness result was proved for a genuinely quasilinear dispersive problem at scaling regularity. To place the current results   in the context of previous work in water waves, this is also the first global well-posedness result for water waves without localization assumptions on the initial data, at any regularity.

These results should also be understood in the context of the
general global well-posedness conjecture for two dimensional quasilinear dispersive flows,  proposed in Ifrim-Tataru~\cite{IT-qnls2} and proved there for quasilinear Schr\"odinger flows, and then in Ifrim-Pineau-Tataru~\cite{IPT-ultra} for ultrahyperbolic Schr\"odinger flows (see also the one dimensional version of the conjecture, discussed below).
By comparison our problem is fully nonlinear and
nonlocal, with a dispersion relation not of Schr\"odinger type and singular at zero frequency, the nonlinearities are quadratic rather than cubic, and the
solutions are transported by a strong advection field; each of these alone
defeats the methods used for the Schr\"odinger models.

\subsection{An overview of the global well-posedness conjectures}
To better motivate the present work, we now give a brief overview of the  \emph{Nonlocalized Data Global Well-posedness Conjectures}, formulated by the first and third authors, beginning with the work \cite{IT-global}. In the one-dimensional setting (for both quasilinear and semilinear equations), the conjecture can be roughly phrased as follows:
\begin{conjecture}[Non-localized data defocusing GWP conjecture]\label{c:nld}
One dimensional dispersive flows on the real line with cubic defocusing nonlinearities and small initial data have global in time, scattering solutions.
\end{conjecture}
The main results of \cite{IT-global} and \cite{IT-qnls} provided the first validation for this conjecture in the semilinear, and quasilinear settings, respectively. These results represented the first global in time well-posedness results of this type at this level of generality. More precisely, in these articles it was shown that if the equation has \emph{phase rotation symmetry} and is \emph{conservative} and \emph{defocusing}, then small data in Sobolev spaces give rise to global, scattering solutions. Here we remark that scattering should be interpreted in a suitable weak sense, meaning that the solution satisfies a global $L^6_{t,x}$ Strichartz bound as well as bilinear $L^2_{t,x}$ bounds. Classical scattering is precluded due to the strong nonlinear effects for this problem. The case of semilinear flows with more general dispersion relations, not of Schr\"odinger type, was later considered in \cite{IT-general}; see also the survey \cite{IT-conjecture}.

We crucially note that the defocusing condition is imperative for the global result in one dimension. In the focusing case, the existence of small solitons serves as an obstruction to obtaining global, scattering solutions. Nevertheless, an analogue of the above conjecture, proposed by two of the authors in \cite{IT-focusing}, suggests a potentially optimal lifespan bound for solutions.

\begin{conjecture}[Non-localized data long time well-posedness conjecture]\label{c:nld-focusing}
One dimensional\\ dispersive flows on the real line with cubic conservative nonlinearities and initial data of size $\epsilon \ll 1$ have
long time solutions on the $\epsilon^{-8}$ time scale.
\end{conjecture}

The conservative condition (which is also relevant in the two-dimensional setting) is a natural assumption. Roughly, one can interpret this as preventing nonlinear ODE blow-up of solutions with wave packet localization.
\medskip

In contrast to the one-dimensional setting, in two and higher dimensions solutions exhibit stronger dispersive decay. In particular, the $L^4_{t,x}$ space-time norm of the solutions will play the leading role, as opposed to $L^6_{t,x}$ in one dimension. We remark that in the semilinear case,
the small data problem may be treated in a relatively standard fashion by Strichartz estimates, which yields global well-posedness and scattering. The focusing/defocusing assumption is not required in this setting, and only plays a role in the large data case, see e.g. Dodson's results in \cite{D,D-focus}.

Therefore, the most interesting case to consider is that of quasilinear flows, where Strichartz estimates are not typically directly applicable. With this in mind, the recent work of the first and third authors in \cite{IT-qnls2} posed the following higher dimensional analogue of the above conjectures.

\begin{conjecture}[$d=2$ non-localized data GWP conjecture]\label{c:nld2}
Cubic quasilinear dispersive flows with small data give rise to global in time, scattering solutions in two and higher dimension.
\end{conjecture}

In contrast to the one dimensional case, here we interpret scattering in the classical fashion, but with the caveat that, due to the nature of quasilinear
well-posedness, we only expect convergence to a linear evolution in a weaker topology. In the article \cite{IT-qnls2}, this conjecture was established for all cubic quasilinear Schr\"odinger flows with positive-definite metrics:

\begin{theorem}[\cite{IT-qnls2}]\label{t:qnls}
Consider a cubic quasilinear Schr\"odinger flow in $\R^d$, $d \geq 2$,
with a positive definite metric. Then small initial data in $H^s(\R^d)$
yields global, scattering solutions, provided that either (i) $d \geq 3$ and
$s > \frac{d+1}2$, or (ii) $d = 2$ and $s \geq \frac74$, or (iii) $d = 2$,
$s > \frac32$ and the problem is conservative.
\end{theorem}

This was subsequently extended by the authors in \cite{IPT-ultra} to merely non-degenerate metrics, i.e. to ultrahyperbolic Schr\"odinger flows, which was the first setting where the dispersion relation is not convex. In both cases the global result holds at the same Sobolev regularity as the corresponding local well-posedness result, only one half of a derivative above scaling.

Although the above conjecture is stated in all dimensions, the most interesting case by far is that of two dimensions, where the dispersive character is weakest. Prior to the work of the first and third author in \cite{IT-qnls2}, no small data global well-posedness results existed in this setting. The reason for this is morally due to the fact that in two space dimensions $L^4_{t,x}$ is a sharp Strichartz norm, whereas in higher
dimensions there is room to interpolate. In fact, in dimension three and higher, global well-posedness was known in some specific model cases, though only at higher regularity; see for instance \cite{HLT}.

The gravity water wave problem in $n = 3$ dimensions, which reduces to a two dimensional flow for the free surface, is a prime candidate for Conjecture~\ref{c:nld2}: it is a physically fundamental quasilinear dispersive flow which, as discussed below, is effectively cubic. At the same time, it lies well outside the class of problems considered in \cite{IT-qnls2, IPT-ultra}. One main aim in the present paper is to establish Conjecture~\ref{c:nld2} for this problem, in the physical dimension $n = 3$ as well as in higher dimension. A  conservative-type assumption turns out to be satisfied for our water wave model,
and plays a key role in dimension $n=3$, but only a minimal role  in dimension $n=4$.

\subsection{Gravity water waves in deep water}

We consider the motion of an incompressible, inviscid fluid occupying a time-dependent domain
\[
\Omega(t)=\{(x,z)\in\mathbb R^d\times\mathbb R:\ z<\eta(t,x)\}, \qquad d = n-1,
\]
where $\eta(t,x)$ describes the free surface
\[
\Gamma(t)=\{(x,z)\in\mathbb R^d\times\mathbb R:\ z = \eta(t,x)\}.
\]

This motion is described using the fluid velocity $u=u(t,x,z)$ and pressure $p=p(t,x,z)$, which satisfy the incompressible Euler equations
\begin{equation}\label{Euler}
\begin{cases}
\partial_t u + (u\cdot\nabla_{x,z})u = -\nabla_{x,z} p - g\,e_z, \\
\nabla_{x,z}\cdot u = 0,
\end{cases}
\qquad \text{in } \Omega(t),
\end{equation}
where $g>0$ is the gravitational constant and $e_z$ is the vertical unit vector.

Under the additional assumption that the flow is irrotational,
\[
\nabla_{x,z}\times u = 0,
\]
the velocity $u$ can be represented using a velocity potential $\phi$ in the form
\[
u = \nabla_{x,z}\phi,
\]
where
\[
\Delta_{x,z}\phi = 0 \quad \text{in } \Omega(t).
\]

The motion of the free surface is uniquely determined by two boundary conditions on $z=\eta(t,x)$, as follows:

\medskip

\noindent
\textbf{Kinematic condition:} the surface is transported by the flow,
\[
\partial_t \eta = \partial_z \phi - \nabla \eta \cdot \nabla \phi,
\qquad \text{ on } \Gamma(t).
\]

\medskip

\noindent
\textbf{Dynamic condition:} the pressure is constant (normalized to zero) on the surface. This can be seen as equivalent to Bernoulli's law,
\[
\partial_t \phi + \tfrac12 |\nabla_{x,z}\phi|^2 + g\,\eta = 0
\qquad \text{ on } \Gamma(t).
\]

\medskip

In the deep water regime, one imposes decay at infinity,
\[
\nabla_{x,z} \phi(t,x,z)\to 0 \quad \text{as } z\to -\infty.
\]

\medskip

As observed by Zakharov~\cite{Zakharov}, restricting the velocity potential to the free surface,
\[
\psi(t,x):=\phi(t,x,\eta(t,x)),
\]
the state of the system is uniquely determined by the functions $\eta$
and $\psi$, which in particular reduces the dimensionality of the problem by one.
Indeed, one obtains a closed system for $(\eta,\psi)$, see also Craig-Sulem~\cite{CraigSulem1993}, namely
\begin{equation}\label{WW}
\begin{cases}
\partial_t \eta = G(\eta)\psi, \\
\partial_t \psi = -g\,\eta - \dfrac12 |\nabla \psi|^2 + \dfrac12 \dfrac{(G(\eta)\psi + \nabla \eta \cdot \nabla \psi)^2}{1+|\nabla \eta|^2},
\end{cases}
\end{equation}
involving the Dirichlet--Neumann operator $G(\eta)$ associated to the domain $\Omega(t)$. This is defined by (see the extensive discussion in Section~\ref{sec:DN})
\[
G(\eta)\psi
=
\bigl(\partial_z\phi - \nabla\eta\cdot\nabla\phi\bigr)\big|_{z=\eta(x)}.
\]
The system \eqref{WW} is Hamiltonian, with conserved energy
\[
E(\eta,\psi) = \frac12 \int_{\R^d} g\, \eta^2 + \psi\, G(\eta) \psi \, dx,
\]
which is comparable to the square of the $L^2 \times \dot H^\frac12$ norm of $(\eta,\psi)$.

\subsection{Linear dispersion and resonances}
Linearizing around the equilibrium solution $(\eta,\psi)=(0,0)$ yields the system
\[
\begin{cases}
\partial_t \eta = |D|\psi, \\
\partial_t \psi = -g\,\eta,
\end{cases}
\]
and hence the linear equation
\[
\partial_t^2 \eta + g|D|\,\eta = 0,
\]
which shows the dispersive character of the problem, with the
dispersion relation
\[
\omega^{\pm}(\xi) = \pm \sqrt{g |\xi|}.
\]
This differs from the Schr\"odinger dispersion relation in several important ways:
\begin{itemize}
\item The group velocity $\nabla \omega^+(\xi) = \frac{\sqrt g}2 |\xi|^{-\frac32}\xi$ decreases with the frequency, and is unbounded as $\xi \to 0$. Thus high frequency waves travel slower than low frequency waves.
\item The Hessian of $\omega^+$ is nondegenerate away from $\xi = 0$, but of indefinite signature, with one negative eigenvalue in the radial direction and $d-1$ positive eigenvalues in the angular directions.
\end{itemize}
In particular, localized solutions to the linear equation have $t^{-\frac{d}2}$ dispersive decay, and the linear flow satisfies Strichartz estimates with the same exponents as the Schr\"odinger flow in $\R^d$, but with a different balance of derivatives. For instance, in the physical case $n = 3$ we have
\begin{equation}\label{Str-lin}
\| e^{-it\omega^+(D)} u_0 \|_{L^4_{t,x}(\R \times \R^2)} \lesssim \| u_0\|_{\dot H^{\frac38}(\R^2)},
\end{equation}
where, exactly as in the case of the two dimensional Schr\"odinger equation, $L^4_{t,x}$ is a sharp Strichartz norm. 

At the nonlinear level, the system \eqref{WW} contains quadratic interactions. However, the gravity wave dispersion relation admits no nontrivial three wave resonances: the relation
\[
\omega^{\pm}(\xi_1) + \omega^{\pm}(\xi_2) + \omega^{\pm}(\xi_3) = 0, \qquad \xi_1+\xi_2+\xi_3 = 0,
\]
can only hold if one of the three frequencies vanishes. Hence, at least formally, the quadratic terms can be removed by a normal form transformation, after which the leading nonlinear interactions are cubic; it is in this sense that \eqref{WW} should be viewed as a \emph{cubic} quasilinear dispersive flow, within the scope of Conjecture~\ref{c:nld2}. Two difficulties arise here: the normal form transformation is unbounded, which reflects the quasilinear character of the problem, and its symbol is singular at zero frequency. In addition, in contrast with the case $n = 2$, if $n \geq 3$ then the set of four wave resonances is large and the resonant cubic interactions are nontrivial.

\subsection{Scaling, criticality and the good variable}
The water wave equation is invariant with respect to the scaling
law
\begin{equation}\label{scaling}
(\eta(t,x),\psi(t,x)) \to (\lambda^{-1} \eta(\lambda^{\frac12} t, \lambda x), \lambda^{-\frac32} \psi(\lambda^{\frac12} t, \lambda x)).
\end{equation}
This suggests one should work with $(\eta,\psi)$ in Sobolev spaces of the form $\H^s := H^s \times H^{s+\frac12}$, which corresponds to $\phi \in H^{s+1}(\Omega(t))$. The scale invariant space
is
\[
\dot \H^{s_c} = \dot H^{s_c} \times \dot H^{s_c+\frac12}, \qquad s_c = \frac{n+1}2 = \frac{d}2 + 1,
\]
and will play a prominent role in our analysis. This corresponds to the Lipschitz norm of the free surface, which is also scale invariant; in particular the problem is energy supercritical in all dimensions. We will also use the associated $1$-Besov space,
\[
\B^s =B^s_{2,1} \times B^{s+\frac12}_{2,1},
\]
and its homogeneous version $\dB^s$, and in particular the scale invariant space $\dB^{s_c}$.

However, as is well understood by now, we cannot work directly with $\psi$ and $\eta$ in the above spaces. Even though the space $H^{s+1}$
would be the appropriate one for $\phi$, the trace theorem does not apply properly at this level for an $H^s$ boundary, and instead
$\psi$ inherits some of the lower, $H^s$ regularity of $\eta$.
To resolve this difficulty we will instead need to work with some well-chosen good variables, for which, at least initially, we employ the choices discovered by Alazard-M\'etivier~\cite{AM}.
To describe this we begin with the functions $V$ and $B$, which represent the traces of the horizontal and vertical components of the velocity field, respectively,
\begin{equation*}
V\coloneq (\nabla_x\phi)_{|z=\eta},\hspace{5mm}B\coloneq (\partial_z\phi)_{|z=\eta}.
\end{equation*}
These relate to $\psi$ and $G(\eta) \psi$ via the relations
\begin{equation}\label{BV}
\left\{
\begin{aligned}
V + B \nabla \eta = & \ \nabla \psi,
\\
B - V \cdot \nabla \eta = & \ G(\eta) \psi,
\end{aligned}
\right.
\end{equation}
which in turn can be used to rewrite the system \eqref{WW}
in the form
\begin{equation}\label{WW-BV0}
\begin{cases}
\partial_t \eta = B - V \cdot \nabla \eta, \\
\partial_t \psi = -g\eta+\frac{1}{2}B^2-\frac{1}{2}|V|^2 -BV\cdot\nabla \eta.
\end{cases}
\end{equation}

While the function $\psi$ is not a good variable,
$V$ and $B$ can be thought of as good variables,
though at the level of $\nabla \psi$.
Correspondingly, at the level of $\psi$
Alazard and M\'etivier introduced the good variable substitute
\begin{equation}\label{good-var}
\varphi = \psi - T_B \eta,
\end{equation}
as suggested by the first equation in \eqref{BV}.
Then the strategy is to recast the system \eqref{WW} as a system
for the independent variables $(\eta, \varphi)$, where it is now
appropriate to work in the $\H^s$ Sobolev spaces or their homogeneous counterparts $\dot \H^s$, as well as in the critical Besov space $\dB^{s_c}$.

\subsection{A brief history of the problem}
The mathematical study of gravity water waves has its modern analytical origins in the Hamiltonian formulation of Zakharov  and the reduction to the free surface via the Dirichlet--Neumann operator  developed in \cite{Zakharov}. The resulting equations are quasilinear and nonlocal, combining elliptic and dispersive features, and their well-posedness theory has been a central topic in nonlinear PDE. Other formulations of the water wave  equations are also available, including in Lagrangian coordinates, or the holomorphic coordinates 
in dimension $n = 2$.

\medskip
\noindent
\textbf{Local well-posedness.} After Zakharov~\cite{Zakharov} reduced the problem to the free surface, and
after early results for small or analytic data, due to
Nalimov~\cite{Nalimov}, and Yosihara~\cite{Yosihara},
local well-posedness for large data in Sobolev spaces $H^s$, $s \geq 4$, was
proved by Wu~\cite{Wu2d,Wu3D} in two and three space dimensions. 
 Regardless of the spatial dimension, the Dirichlet--Neumann operator is
always present in the equations. A systematic analysis of the
Dirichlet--Neumann operator, and of the full water wave system, was developed
by Lannes~\cite{Lannes}.  These results were complemented by many others, considering also rotational flows, finite depth or surface tension; without being exhaustive we mention \cite{CL, Lindblad, CoSh, ShZ, ZZ}, and refer the reader to the book \cite{Lannes-book} for further references.

Subsequent work using paradifferential calculus significantly lowered the regularity requirements. The paralinearization of the Dirichlet--Neumann operator, together with the good variable described above, was introduced by Alazard--M\'etivier \cite{AM}. Building on this, Alazard--Burq--Zuily \cite{ABZ} proved the following result, which we state in the infinite depth case and in the notations of this paper\footnote{The exponent $s$ in \cite{ABZ} corresponds to $s - \frac12$ here.}:

\begin{theorem}[\cite{ABZ}]\label{t:ABZ}
Let $d \geq 1$ and $s > s_c + \frac12$. Then the gravity wave system \eqref{WW} is locally well-posed for initial data
\[
(\eta_0,\psi_0) \in H^{s} \times H^{s}, \qquad (V_0,B_0) \in H^{s-\frac12}.
\]
\end{theorem}

In terms of Sobolev embeddings, this corresponds to a Lipschitz velocity field and a $C^{\frac32}$ free surface, and is essentially optimal from the quasilinear perspective, i.e. for a proof based exclusively on energy estimates, without taking dispersion into account; in two space dimensions it was slightly relaxed to BMO type control norms in \cite{HIT}. The counterpart of this result for data in uniformly local Sobolev spaces, without any decay at infinity, was proved in \cite{ABZ-nonloc}.

To go below this threshold one needs to take advantage of the dispersive character of the equations. This was first achieved by Alazard--Burq--Zuily~\cite{ABZ-disp} where, by proving Strichartz estimates with derivative losses for the paradifferential flow, the above threshold was improved to $s > s_c + \frac12 - \mu$, with $\mu = \frac1{24}$ for $d = 1$ and $\mu = \frac1{12}$ for $d \geq 2$. The Strichartz estimates, and correspondingly the range of $s$, were further improved by Ai \cite{Ai1,Ai2}, to $\mu = \frac1{8}$ for $d = 1$ and $\mu = \frac1{5}$ for $d \geq 2$.

A different route was opened in two space dimensions by Ai--Ifrim--Tataru~\cite{AIT}, where it was shown that the energy estimates themselves can be substantially improved, by a normal form analysis implemented in a paradifferential fashion at the level of modified energies. The resulting \emph{balanced cubic energy estimates} led to a gain over Theorem~\ref{t:ABZ} which exceeds what was previously obtained using Strichartz estimates by another $1/8$ derivative. The same principle has since proved effective in a variety of quasilinear problems, see e.g. \cite{AIT-minimal}, and was combined with bilinear $L^2_{t,x}$ bounds, obtained via interaction Morawetz estimates, in the Schr\"odinger setting in \cite{IT-qnls, IT-qnls2, IPT-ultra}. The local well-posedness result of the present paper brings these ideas to bear on the water wave problem in dimension $n \geq 3$.

\medskip
\noindent
\textbf{Global solutions for small localized data.} The question of global behavior is substantially more delicate, due to the weak dispersion and the presence of resonant interactions, and has so far been studied almost exclusively for initial data which is small, smooth and \emph{localized}. In dimension $n = 3$ (two horizontal dimensions), global regularity was established by Germain--Masmoudi--Shatah \cite{GMS} and independently by Wu \cite{Wu-global}:

\begin{theorem}[\cite{GMS, Wu-global}]\label{t:GMS-Wu}
Let $n = 3$. Assume that the initial data for \eqref{WW} is small, sufficiently regular and localized. Then the solution is global in time, and decays pointwise at the linear dispersive rate $t^{-1}$.
\end{theorem}

The proof in \cite{GMS} is carried out in Eulerian coordinates and is based on the space-time resonance method, with localization measured using weighted norms for the linear profile of the solution; it also shows that the solutions scatter. The proof in \cite{Wu-global} uses instead Lagrangian coordinates and a nonlinear change of unknowns which reveals the cubic structure of the equations, combined with the vector field method.  Further global results in three dimensions were obtained for capillary waves \cite{GMS-cap}, for gravity-capillary waves \cite{DIPP}, and for gravity waves in finite depth \cite{Wang-3D}.

In dimension $n = 2$ (one horizontal dimension) the dispersive decay rate is only $t^{-\frac12}$, and the solutions exhibit modified scattering. Almost global solutions were first obtained by Wu~\cite{Wu-almost}, and global solutions by Ionescu--Pusateri \cite{IP} and, independently, Alazard--Delort \cite{AD}. More recently, new approaches have emphasized the interaction between the quasilinear structure and dispersive effects. In particular, the use of holomorphic (or conformal) coordinates, of modified energies and of wave packet testing has led to streamlined proofs and sharper descriptions of the long-time behavior, including modified scattering and asymptotic dynamics; see \cite{HIT, IT-g}, as well as \cite{AIT2} for global solutions at low regularity, and \cite{IT-c, IP-cap, Wang-inf, IT-vort, HGIT} for related problems. These methods provide a flexible framework that captures both the geometric structure of the equations and their dispersive features. We also refer the reader to the survey \cite{IP-survey}.

A common feature of all of these results is that the proofs combine (i) energy estimates for high Sobolev norms, where the energy growth is controlled by a pointwise norm of the solution, with (ii) dispersive decay bounds for this pointwise norm, where localization of the data, propagated along the flow via weighted norms or vector fields, plays an essential role.

\medskip
\noindent
\textbf{Nonlocalized data.} If the initial data is merely small in a Sobolev space then no uniform pointwise decay can be expected, and step (ii) above breaks down. In this setting, the results available so far for water waves are primarily enhanced lifespan bounds. Precisely, the normal form structure yields a lifespan of at least $\epsilon^{-2}$ for data of size $\epsilon$, see \cite{HIT} for $n = 2$, which was further improved in two dimensions in \cite{Wu-quartic, DIP}. In three dimensions, the best result prior to our work is due to Zheng:

\begin{theorem}[\cite{Zheng}]\label{t:Zheng}
Let $n = 3$. For initial data which is small in an unweighted Sobolev space and sufficiently regular, the solutions to \eqref{WW} are almost global. If in addition the data is small in a weighted Sobolev space with weight $|x|^\alpha$ for some $\alpha > 0$, then the solutions are global.
\end{theorem}

More precisely, the nonlocalized result in
\cite{Zheng} gives a lifespan $T_\epsilon\gtrsim\exp(c/\epsilon)$
for initial data satisfying
$\|\eta_0\|_{H^{N+\frac12}}
+\||D|^{\frac12}\psi_0\|_{H^N}\leq\epsilon$ for $N\geq 11$.
The proof uses quadratic normal forms and quartic energy
estimates, together with Strichartz estimates (with a derivative loss) by treating the equation as a perturbation of the linear evolution $(\partial_t+i|D|^{\frac{1}{2}})u=0$. The low frequency part of cubic errors in the energy estimate is controlled 
using the $L_t^2L_x^\infty$ endpoint estimate, which in two dimensions has a logarithmic in time loss. This obstruction appears at the linear level and cannot be removed by strengthening the regularity assumptions on the data. This is fundamentally what leads to the exponential lifespan bound in Theorem \ref{t:Zheng} rather than global existence. 

If the initial data is additionally small in weighted Sobolev
norms with weight $(1+|x|)^\alpha$ for some $\alpha>0$, Zheng
obtains global existence. Additionally, the required regularity depends
on $\alpha$ (in particular, one requires $N\to\infty$ as $\alpha\downarrow 0$).

In contrast, our global result removes the localization assumption altogether, and at the same time lowers the regularity to the scaling critical level.

\subsection{The new results}

The first objective of this article is to study the question of local well-posedness in Sobolev spaces. To avoid a number of technical difficulties, we will consider the well-posedness question
only in the small data regime. However, we will only ask for smallness in a critical norm for the good variables (see \eqref{good-var}) $(\eta,\varphi)$, namely
\begin{equation}\label{small-data}
 \|(\eta_0,\varphi_0)\|_{\dB^{s_c}} \leq \epsilon \ll 1.
\end{equation}
Even in the small data regime, this comes with a multitude of difficulties due to the infinite speed of propagation (in sharp contrast to corresponding questions for nonlinear wave equations).

\begin{theorem}\label{t:local}
The water wave problem \eqref{WW} is locally well-posed for initial data $(\eta_0,\varphi_0)$  satisfying the smallness condition \eqref{small-data},
in Sobolev spaces as follows:
\begin{enumerate}[label=(\roman*)]
    \item $n = 3$: in $\H^s$ for $s \geq s_c+\frac1{12}$.
    \item $n = 4$: in $\H^s$ for $s > s_c$.
\item $n\geq 5$: in the critical Besov space $\B^{s_c}$.
\end{enumerate}
 
\end{theorem}
This result represents a major improvement over the best prior result in \cite{Ai1}, by almost $1/3$ derivatives in dimension $n=3$, and to the sharp scaling threshold in higher dimensions $n \geq 4$. By comparison, the energy threshold in Theorem~\ref{t:ABZ}  is one half of a derivative above scaling in all dimensions.

\medskip

As it is standard in quasilinear settings, we interpret well-posedness in an enhanced Hadamard sense, see \cite{IT-primer}, which in this article includes the following suite of properties:
\begin{itemize}
    \item existence of solutions in $C([0,T]; \H^s)$.
    \item uniqueness of regular solutions.
    \item uniqueness of rough solutions as
    uniform limits of regular solutions.
    \item continuous dependence on the initial data.
    \item Lipschitz dependence of the solutions on the initial data in a weaker topology  $\H^0$.
    \item higher regularity: more regular data yields more regular solutions.
\end{itemize}
\smallskip

We continue with our global result:
\begin{theorem}\label{t:global}
Let $s$ be as in the local well-posedness result. Then the 3D water wave problem \eqref{WW} is
globally well-posed in $\H^s$ for  initial data which is small in the critical Besov space $\dB^{s_c} := \dot B^{2}_{2,1} \times \dot B^{\frac52}_{2,1}$, as in \eqref{small-data}. The same result holds for $n \geq 4$, provided that the initial data is in the local well-posedness space 
in Theorem~\ref{t:local}.
\end{theorem}

We emphasize that the smallness condition \eqref{small-data} involves no weights and no vector fields. Thus Theorem~\ref{t:global} applies in particular to initial data which is a superposition of infinitely many widely separated wave packets, for which the solutions have no uniform pointwise decay rate. Instead, the dispersive character of our global solutions is captured by a family of space-time bounds, which include:
\begin{enumerate}[label=(\roman*)]
\item Uniform energy bounds.
\item Bilinear $L^2_{t,x}$ bounds, in both balanced and unbalanced form.
\item Strichartz type bounds.
\end{enumerate}
The reader is referred to Theorem~\ref{t:local-fe3} for a complete, frequency envelope based description of these bounds.

\medskip

We continue with several remarks concerning these results.

\begin{remark}
Compared with Theorem~\ref{t:GMS-Wu} and Theorem~\ref{t:Zheng}, our result removes all localization assumptions, and also lowers the regularity of the data to the scaling critical level. On the other hand, for localized data the results in \cite{GMS, Wu-global} provide pointwise decay for the solutions, which cannot hold in our setting.
\end{remark}

\begin{remark}[The case $n = 2$]
Our results are restricted to dimension $n \geq 3$, where the free surface
has dimension at least two. The two dimensional problem, where the free
surface is a curve, is of a different nature. There the dispersive decay is
weaker, the balanced cubic interactions are strongly resonant rather than merely
borderline, and the bilinear $L^2_{t,x}$ bounds must be supplemented by
$L^6_{t,x}$ Strichartz type bounds. The natural counterparts of
Theorem~\ref{t:global} are then Conjecture~\ref{c:nld} and
Conjecture~\ref{c:nld-focusing}. Since the cubic NLS approximation for
gravity waves in deep water is focusing, it is the latter conjecture which is
relevant here; thus for nonlocalized data one should expect long time, rather
than global, well-posedness. Global solutions for two dimensional gravity
waves are known only for small localized data \cite{AD}, \cite{IP}, \cite{HIT}-\cite{IT-g}, and exhibit modified scattering. Many of the tools
developed in this article apply equally in the case $n = 2$; this
problem will be considered in forthcoming work.
\end{remark}

\begin{remark}
 Unlike the authors' recent work on quasilinear Schr\"odinger flows \cite{IT-qnls2, IPT-ultra}, here scattering is 
 unlikely to hold for solutions at the low
 regularity in the above global result. The main culprit is the advection field $V$, which is unbounded at the regularity we consider, and may cause a substantial drift compared to linear dispersion. However, we  do expect scattering 
 to hold if we improve the regularity
 of the initial data both at low frequency
 and at high frequency, e.g. for data 
 $(\eta_0,\varphi_0) \in \dH^0 \cap \dH^{s_c+\frac12}$ in dimension $n=3$. This will be discussed in follow-up work. 
\end{remark}

We conclude with a brief comparison with the quasilinear Schr\"odinger flows studied in \cite{IT-qnls2, IPT-ultra}. While the overall strategy is similar, the water wave problem brings a number of new difficulties:
\begin{itemize}
\item The equations are fully nonlinear and nonlocal, and only become quasilinear after paralinearization, when expressed in terms of the good variables.
\item The quadratic terms cannot be removed by a bounded normal form transformation, and must instead be dealt with at the level of the density-flux identities.
\item The dispersion relation is not of Schr\"odinger type; its Hessian is indefinite, and it is singular at zero frequency, where the group velocities are unbounded.
\item In addition to the dispersive term, the equations also contain an advection term,
with a velocity $V$ which is in general unbounded, and strongly interferes with dispersion.

\item Similar to the Schr\"odinger case, the study of the balanced cubic interactions 
plays an essential role, particularly in dimension $n = 3$. This is also where the 
$\frac1{12}$ loss in the local well-posedness result for $n = 3$ originates. 
However, here working with pure
cubic expressions is no longer sufficient; instead, we need a full accounting of the corresponding low frequency paracoefficients. 

\end{itemize}

\subsection{An outline of the paper}

 The three dimensional case $n=3$ of our results is both the most interesting and the most difficult. Hence, in order to 
 keep our notations cleaner and our presentation focused, the bulk of the paper 
is devoted to the three dimensional case.
Only at end, in  the last section of the paper we return to the higher dimensional case and provide a discussion of the differences in the proof, which may be seen as a strict subset of the three dimensional 
proof.

 As in \cite{IT-qnls, IT-qnls2, IPT-ultra}, the linear and bilinear bounds for the solutions cannot be proved independently, and the entire proof is wrapped within an extensive bootstrap argument. Nevertheless, we are able to organize the proof in a modular fashion, where the main ideas and logical steps are as follows.

\medskip
\noindent
\emph{I. The Dirichlet--Neumann operator and 
the paradifferential expansion of the equations.} As a starting point for the study of the water wave equation, the first step is to recast the water wave system~\ref{WW} as an evolution for the good pair of variables $(\eta,\varphi)$. 
This is best done using Bony's paradifferential formalism, for which 
we follow the approach of Alazard-Metivier~\cite{AM}, Lannes~\cite{Lannes}
and Alazard-Burq-Zuily~\cite{ABZ}, writing the equations as a system of the form
\begin{equation}\label{Z-para-retro}
\begin{cases}
& \ (\partial_t + T_{V} \nabla) \eta - T_{G(\eta)} \varphi = N_\eta(\eta,\varphi)  
\\
& \ (\partial_t + T_{V} \nabla) \varphi +
T_{a} \eta = N_\varphi(\eta,\varphi),
\end{cases}
\end{equation}
which is first discussed in Section~\ref{sec:DN}.
This involves a paradifferential form of 
the Dirichlet to Neumann operator $G(\eta)$,
as well as the \emph{Taylor coefficient} $a$, given by 
\[
a = - \partial_z p_{|z = \eta(t,x)}
\]
i.e. minus the vertical derivative of the pressure on the free boundary. Since $G$
is a positive operator, it is heuristically obvious that well-posedness even for the 
linear part of the above system requires
the Taylor sign condition $a \geq 0$. This will certainly be satisfied for the class of solutions we study here, where we show that 
$a = g+O(\epsilon)$, with $g$ normalized by scaling to $g = 1$.

For the study of the above system it will be very useful to diagonalize it and use a single complex variable 
\[
u = \eta+ i T_{\mfq} \varphi, \qquad  
\mfq = (\dn/a)^{\frac12}.
\]
where $\dn$ represents the principal symbol of the Dirichlet to Neumann map, see \eqref{dn-symbol}.
The equation for $u$ then will be written in  the schematic form\footnote{Also with some further more subtle adjustments which we omit for now.} 
\begin{equation}
(\partial_t + T_V \cdot \nabla + \frac12 (\nabla \cdot V)+ i T_\omega) u = N_u,    
\end{equation}
with a dispersion relation 
\[
\omega =  (a \dn(\eta))^\frac12.
\]

 We will also need the corresponding linearized equation around a solution $(\eta,\varphi)$, which after choosing suitable good variables  can also be recast in a complex form
 \begin{equation}
(\partial_t + T_V \cdot \nabla + \frac12 (\nabla \cdot V)+ i T_\omega) v = N_u^{lin} v    .
\end{equation}
In understanding the properties of both the full system and the linearized one, 
it will also be useful to consider the associated linear paradifferential evolution
\begin{equation}
(\partial_t + T_V \cdot \nabla + \frac12 (\nabla \cdot V)+ i T_\omega) w = 0,    \end{equation}
as well as the corresponding inhomogeneous flow.

\medskip

\emph{II. Layered paradifferential expansions.} 
Once we have the above paradifferential 
representation of the water wave equation, 
the standard approach, see e.g. \cite{ABZ}, is  to estimate perturbatively the nonlinearity $N_u$. 
However, both at regularity near scaling 
and for the study of global in time dynamics
this is no longer sufficient. Instead, we develop a more refined expansion of the nonlinearity $N$, 
which we will call \emph{layered paradifferential expansion}, of the form
\[
N_u =  N_u^{[1]} + N_u^{[2]} + N_u^{[3]} + N_u^{[4]} + N_u^{err}.
\]
The layers that play a nonperturbative role and thus need to be separated and explicitly computed 
in the 3D setting are as follows:
\begin{enumerate}[label=(\roman*)]
    \item The paradifferential layer $N^{[1]}$, 
following Alazard-Metivier~\cite{AM}, with the 
additional proviso that  the (unbalanced part of the) subprincipal symbol is also necessary and explicitly computed.
\item The $high \times high$ interactions $N^{[2]}$,
referred to as ``hh" in the sequel.  
These are at heart bilinear interactions of two high frequencies with a lower frequency output, but which also carry an explicit  lower frequency paracoefficient.

\item The cubic balanced interactions $N^{[3]}$, further separated into a component which has phase rotation symmetry and carries the worst resonant interactions, and a remaining part that can be largely viewed as nonresonant. 

\item The quartic balanced interactions $N^{[4]}$,
which only play a mildly nonperturbative role and can also be largely viewed as nonresonant. 
\end{enumerate}
The situation simplifies in higher dimensions. The balanced quartic term $N^{[4]}$ goes in the perturbative box for $n \geq 4$, and the balanced cubic term $N^{[3]}$ goes in the perturbative box for $n \geq 5$.

In the context of the above expansion, the 
structure of the layers plays an important role. Precisely, we will take full advantage 
of the \emph{null structure} of the water wave equations, which manifests in two of the above layers:
\begin{itemize}
    \item A \emph{null condition} for the $N^{[2]}$ type terms, which heuristically guarantees that the associated normal form corrections are bounded in a suitable sense,
    and
    \item A double null condition for the 
    $N^{[3]}$ component with phase rotation symmetry, which roughly guarantees that (near) double resonant interactions do not significantly increase energy. This represents  a stronger form of the conservative condition of \cite{IT-qnls2}.  
\end{itemize}

The full layered decomposition discussed above is used at the level of the full equation. For the linearized equation, the null structure of $N^{[2]}$ type terms remains, but the  null structure of $N^{[3]}$ type terms disappears due to a loss of symmetry. This causes us to switch the $N^{[3]}$ and also $N^{[4]}$ terms to the perturbative box, which has two consequences for the linearized equation bounds in dimensions $n=3,4$: (i) they are 
only local in time, and (ii) they require 
some additional regularity, namely $1/2$ 
extra derivatives for $n = 3$, respectively $\epsilon$ extra derivatives for $n = 4$.

As an intermediate step, many of the bounds
in the article are first established at the level of the corresponding linear paradifferential flow.

For clarity, the layered paradifferential expansions are computed early on, in Section~\ref{sec:layered}, 
at the level of the $(\eta,\varphi)$ system,
first for the full equation and then for the linearized equation. However, they are only usable once the remaining error is 
estimated favourably; this is only done later, as part of our main bootstrap argument, and the proof of the error bounds is only completed in Section~\ref{s:u}
for the full equation, and in Section~\ref{subsec:wfinal} for the linearized equation. 

\medskip
\noindent
\emph{ III. Linear and bilinear estimates.} 
From a historical perspective, Strichartz 
estimates have initially played a leading role in semilinear dispersive equations, followed next by bilinear and multilinear estimates. In the quasilinear setting, 
however, establishing such estimates has been troublesome, at least in low regularity settings, and has generally been possible 
only with a loss of derivatives, and only locally in time. This is why,
in the authors' earlier work on QNLS \cite{IT-qnls2}, \cite{IPT-ultra}, the leading role has been played by bilinear 
$L^2_{t,x}$ estimates. Considering for simplicity the context of the linear paradifferential flow for the diagonal complex\footnote{ See the diagonalization discussion in Section~\ref{subsec:dispersive-cc}, also following \eqref{def-u}, and finally the full definitions in Section~\ref{sec:complexvar}.}  variable $w$, these estimates appear at the dyadic level in two 
forms, namely 
\begin{itemize}
    \item \emph{The unbalanced case,} for transversal interactions,
 \[
\| w_\lambda \bar w_\mu\|_{L^2_{t,x}} \lesssim \mu^{\frac{d}2-\frac14} \|w_\lambda(0)\|_{L^2_x} \|w_\mu(0)\|_{L^2_x}, \qquad \mu < \lambda;
 \]
 \item \emph{The balanced case,} for near parallel interactions,
 \[
\|P_\alpha (w_\lambda \bw_\mu)\|_{L^2_{t,x}} \lesssim  \alpha^\frac{d-3}2  \lambda^{\frac54} \|w_\lambda(0)\|_{L^2_x} \|w_\mu(0)\|_{L^2_x}, \qquad  \alpha \lesssim \mu \approx \lambda.
 \]
 \end{itemize}
We remark that in the classical Schr\"odinger case, with a convex 
dispersion relation, one also has square summability in $\alpha$ in the latter case,
providing a clean $\dot H^{\frac{3-n}2}$ bound for the product, see \cite{PV}, \cite{IT-qnls2}.
However, this summability appears to be lost in the nonconvex case, see \cite{IPT-ultra}.

Another idea introduced beginning with  the semilinear work \cite{IT-global}, and then 
continuing with the quasilinear results in \cite{IT-qnls} and \cite{IT-qnls2}, is that when working with multilinear forms it is very useful to have translated versions of 
these bounds, allowing for a loss in the translation parameter. In our case, these 
translated estimates have the form
 \[
\| w_\lambda \bar w_\mu^{x_0}\|_{L^2_{t,x}} \lesssim \mu^{\frac{d}2-\frac14} \|w_\lambda(0)\|_{L^2_x} \|w_\mu(0)\|_{L^2_x}(1+\mu |x_0|)^\frac12, \qquad \mu < \lambda,
 \]
respectively
 \[
\|P_\alpha (w_\lambda \bw_\mu^{x_0})\|_{L^2_{t,x}} \lesssim  \alpha^\frac{d-3}2  \lambda^{\frac54} \|w_\lambda(0)\|_{L^2_x} \|w_\mu(0)\|_{L^2_x}(1+\lambda |x_0|)^\frac12, \qquad  \alpha \lesssim \mu \approx \lambda ,
 \]
where $x_0$ represents the translation parameter for the second function.  Compared to the earlier quasilinear Schr\"odinger results \cite{IT-qnls}, \cite{IT-qnls2}, \cite{IPT-ultra}, in the present article we were able 
to lower the power of the translation factor in the balanced case from $1$ to $1/2$. On the other hand, no such factor was present there in the unbalanced case, heuristically 
because transversality is unaffected by translations. By contrast, in the present article we have the the unbounded advection velocity $V$, which at large translation 
parameters may defeat the group velocity difference coming from dispersion.

We now arrive to a key finding in this 
article, which is that in effect the 
above bounds are superseded by a new 
class of estimates, which we call 
\emph{translation averaged bilinear $L^2$
bounds}, which also come in an unbalanced 
and a balanced version:
\begin{itemize}
    \item \emph{The unbalanced translation averaged bounds:} 
 \[
\| w_\lambda(x) \bar w_\mu(y)\|_{L^2_{t,xy}(\{|x-y|<r\})} \lesssim \mu^{\frac14}r^\frac12 \|w_\lambda(0)\|_{L^2_{x}} \|w_\mu(0)\|_{L^2_x}, \qquad \mu < \lambda, \quad r >\mu^{-1}.
 \]
 \item \emph{The balanced translation averaged bounds:} 
 \[
\|(\partial_x+\partial_y)(w_\lambda(x) \bw_\mu(y))\|_{L^2_{t,xy}(\{|x-y|<r\})} \lesssim   \lambda^{\frac54} r^\frac12 \|w_\lambda(0)\|_{L^2_x} \|w_\mu(0)\|_{L^2_x}, \quad  \alpha \lesssim \mu \approx \lambda, \quad r >\mu^{-1},
 \]
 \[
\|(\partial_x-\partial_y)(w_\lambda(x) w_\mu(y))\|_{L^2_{t,xy}(\{|x-y|<r\})} \lesssim   \lambda^{\frac54} r^\frac12 \|w_\lambda(0)\|_{L^2_x} \|w_\mu(0)\|_{L^2_x}, \quad  \alpha \lesssim \mu \approx \lambda, \quad r >\mu^{-1}.
 \]
 \end{itemize}
Indeed, we show in Lemma~\ref{lem:boot}
that the translation averaged bounds above 
imply all the earlier bilinear $L^2_{t,x}$ bounds.

To conclude, we briefly discuss Strichartz type  bounds, which in this article play a secondary role. We begin with the dimension $d=2$, where $L^4_{t,x}$ would be a sharp Strichartz bound, which would allow us to treat cubic balanced terms perturbatively at scaling. We cannot prove 
that, so instead we contend ourselves 
with weaker bounds obtained from the balanced bilinear bounds with $w_\lambda = w_\mu$ via Sobolev embeddings. Directly 
this yields a scale invariant $L^4_tL^p_x$ bound
with $p > 8$, while interpolating with the energy norm we arrive at $L^{6,\infty}_t L^4_x$.
The latter bound is what allows us to treat 
balanced cubic nonlinearities perturbatively 
at a cost of $1/12$ derivatives, leading to the local well-posedness result in Theorem~\ref{t:local}.

By contrast, in dimension $d = 3$ the sharp Strichartz exponent would be $L^{\frac{10}3}_{t,x}$.
Instead, we obtain $L^4_t L^{4+}_x$, which almost suffices to treat cubic nonlinearities at scaling, and yields local well-posedness in all Sobolev spaces above 
scaling. 

Finally, in higher dimension $d \geq 4$ we directly obtain the scale invariant $L^4_{t,x}$ bound, which suffices to treat cubic nonlinearities at scaling, and yields local well-posedness in the scale invariant Besov space.

\medskip
\noindent
\emph{IV. Frequency envelopes and the bootstrap setup.} While stated in a short form in the introduction, a full form of our  main local and global well-posedness
should include a complete set of estimates
satisfied by solutions and by differences of solutions, including (i) energy bounds, 
(ii) bilinear $L^2_{t,x}$ bounds and Strichartz type bounds. Following the cue of our earlier
works, in Section~\ref{sec:bootstrap} we provide a complete frequency envelope version of our main results, which includes a full family of dyadic linear and bilinear bounds for the solutions to the full equation (Theorem~\ref{t:local-fe3}), the paradifferential flow (Theorem~\ref{t:local-fe3-lin-para}), and the linearized equation 
(Theorem~\ref{t:local-fe3-lin}). 
In brief, these theorems assert that if we place the initial data under a suitable envelope, then the solutions satisfy 
linear and bilinear frequency localized bounds controlled by the same frequency envelope. 

To prove all these bounds, it is very convenient to assume such bounds already hold with worse constants; this can indeed
be assumed by a now standard continuity 
argument.

This is exactly what we do in Section~\ref{sec:bootstrap}, initiating a grand bootstrap loop spanning 
most of the rest of the paper. The loop closes piecewise in several installments:
\begin{itemize}
    \item The Taylor coefficient bounds, in Section~\ref{sec:aexpansions}.
    \item The energy bounds, in Section~\ref{sec:density-flux}.
    \item The translation averaged bilinear 
$L^2$ bounds, in Section~\ref{sec:Morawetz}. 
\end{itemize}
The proof of the frequency envelope bounds
in Theorems~\ref{t:local-fe3}, \ref{t:local-fe3-lin} and \ref{t:local-fe3-lin-para} 
is concluded at the end of Section ~\ref{sec:Morawetz}.

\medskip
\noindent
\emph{V. Continuous Littlewood-Paley expansions and the key error bounds.} As mentioned above, a significant challenge in properly isolating the layers $N_u^{[j]}$ is in obtaining sufficiently precise estimates for the error terms. One (among many) of the chief difficulties stems from the fact that the elliptic operator $G(\eta)$, and the variables $V$, $B$, $a$ depend in a highly nonlinear and non-local way on the variable $u$. In order to algebraically to implement the normal form corrections for each layer $N_u^{[j]}$, and to reveal their null structures mentioned above, it is necessary to represent these terms explicitly as multi-linear pseudo-differential operators applied to the variables $u$, $\eta$ and $\varphi$, with explicit, bounded low-frequency para-coefficients. This analysis makes up the bulk of Sections \ref{sec:elliptic}--\ref{subsec:wfinal}, whose contents we now briefly describe.

We begin in Section~\ref{sec:elliptic} by collecting most of the necessary elliptic estimates for the Dirichlet-Neumann
operator and other related boundary operators. These
include $L^p$ based linear bounds, together with estimates
for the errors arising when we replace these operators
by their principal parts. To obtain the bilinear bounds
needed later, we must also identify the structure of some of the subprincipal
terms in these elliptic operators more explicitly. We will then use some of these estimates to transfer the bootstrap bounds satisfied by $u$ to the variables $\eta$ and $\varphi$. At the end of this section, we also introduce a useful
multi-linear functional framework that will allow us to efficiently apply the
translation-averaged bounds in our error analysis. This will be of particular importance for proving multi-linear estimates involving elliptic operators that are singular at low frequency (and thus have kernels with limited decay).

In Section~\ref{sec:VB}, we also transfer the bootstrap
bounds for $u$, $\eta$ and $\varphi$ to the
velocity traces $V$ and $B$; the main estimates are stated in  Proposition \ref{VBunbil}. This will allow us to identify the corresponding principal first order layers $B^{[1]}$ and $V^{[1]}$ of $B$ and $V$, respectively, with good remainder estimates for the errors. These serve as basic ingredients for obtaining the full layered expansions for $N_u$, where the variables $V$ and $B$ will appear frequently. 

In Sections~\ref{sec:Neta} and \ref{sec:Nphi}, we carry
out the layered expansions for the nonlinearities $N_\eta$ and $N_\varphi$, corresponding to $\eta$ and $\varphi$ (along with their counterparts for the linearized equation), starting directly from the two equations in the Zakharov system. The main tool here is a continuous
Littlewood-Paley decomposition, applied both to the good unknowns
and to the nonlinear expressions and elliptic operators obtained by regularizing
the free boundary. Using this tool, we can successively expand by differentiating with respect to the regularization parameter. These successive expansions (or linearizations)  will then allow us to extract
the different layers for $N_{\eta}$ and $N_{\varphi}$, while retaining explicitly all of the lower frequency paracoefficients. At each iteration of this procedure, the
preceding linear and multilinear estimates are used to
peel off the perturbative errors and to establish several auxiliary
bounds needed later to close the bootstrap argument. Carrying out this analysis at the critical regularity level is highly delicate (especially due to the need to make use of our bilinear estimates), and presents a substantial number of additional technical issues compared to what one would encounter in the high regularity regime.

In the remaining sections, our goal is to appropriately symmetrize the systems of layered paradifferential equations for $\eta$ and $\varphi$ (and their linearized counterparts) using the diagonal variable $u$ mentioned above. A key difficulty here is that the symmetrizing operator $T_{\mathfrak{q}}$ 
depends on the Taylor coefficient $a$, and commuting it
with the para-material derivative (which roughly speaking corresponds to suitable modification of the operator  $(\partial_t+T_V\cdot\nabla)$)  introduces errors involving both $a$ and its para-material derivative. To deal with these expressions, in Sections~\ref{sec:aexpansions} and
\ref{sec:Dtaexpansions}, we need to derive corresponding layered
expansions (only to first order) for these coefficients and prove the necessary linear and
multi-linear bounds needed to control their errors. In the latter of these two sections, the singular nature of many of the elliptic operators at low frequency manifests very strongly. Here, the translation-averaged bounds (which we implement in this section efficiently through the multi-linear framework established in Section \ref{sec:elliptic}) are of crucial importance. Finally, in Sections~\ref{s:u} and \ref{subsec:wfinal}, we combine
these estimates with the expansions for the
$(\eta,\varphi)$ system to obtain
the equations for the symmetrized variable $u$ (and its linearized counterpart $w$), completing the error
estimates for their layered paradifferential expansions.

\medskip
\noindent
\emph{VI. Paracontrolled functions and paratransported symbols.} The presence in the equations of the advection velocity $V$, which is possibly unbounded,  adds another difficulty to the mix, namely that 
even very nice multipliers, such as the classical Littlewood-Paley projectors, do not commute well with the material derivative. To rectify this, we introduce a class of so called \emph{paratransported symbols}, which do commute approximatively with the material derivative.
Then the challenge, resolved in Section~\ref{sec:para-functions}, is to construct such symbols (e.g. Littlewood-Paley projector substitutes).  For this, the price to pay is to expand the class of admissible paracoefficients, so that it is closed with respect to paraproducts. This is achieved using suitable classes of paracontrolled functions.

\medskip
\noindent
\emph{VII. Conservation laws in density-flux form.} In Section~\ref{sec:density-flux} we first establish frequency localized mass and momentum conservation laws for the linear paradifferential flow, in \emph{density-flux form}. One new feature here is the presence 
of the advection term, which is nonperturbative and is instead included 
in the principal part in a paradifferential fashion. Thus we work with \emph{advected density-flux} identities,
which for the frequency localized mass
density, for instance, have the form 
\[
(\partial_t + \nabla \cdot V_{<\lambda})
\Ms_\lambda(w) = \partial_j \Fsj_\lambda(w) + R_\lambda.
\]
Importantly, the symbols for both the mass and momentum densities belong to the class of paratransported functions;
this in particular avoids the necessity of infinitely\footnote{With only one exception,
which is fortunately very narrow in scope.}
many reiterated normal form type corrections at the level of the principal symbols.
The $\sharp$ superscript indicates that 
the corresponding densities have been further modified, akin to the \emph{balanced modified energies} of \cite{AIT} and \cite{AIT-minimal}, in order to improve accuracy and to guarantee that the remaining
source term $R_\lambda$ is balanced quartic and thus perturbative.

Once we have a good understanding of the 
paradifferential equation, we expand these results to both the full equation and the linearized equation, by adding normal form type corrections corresponding to their respective layered expansions. 
At the conclusion 
of Section~\ref{sec:density-flux}, we are in particular able to close the bootstrap loop for the uniform energy bounds.

\medskip
\noindent
\emph{VIII. Interaction Morawetz bounds.} These are our main tool in the proof of the bilinear $L^2_{t,x}$ bounds, and are developed in Section~\ref{sec:Morawetz}, building on the NLS ideas in \cite{CKSTT, PV} as implemented in \cite{IT-qnls2, IPT-ultra} in the quasilinear Schr\"odinger setting, as well as on the adaptation of these ideas to general dispersion relations in \cite{IT-general}; compare also with the Morawetz estimates for water waves in \cite{AIT-Morawetz}. 

The new difficulties here arise from the nonlocal character of the densities and fluxes, the presence of the advection component, as well as from the lack of convexity of the dispersion relation.
At the technical level, the use of paratransported symbols for the mass 
and momentum densities forces a loss of symmetry in the leading component of the time derivative of the interaction functional, compared to the setting of 
\cite{IT-qnls2, IPT-ultra}. This can only be rectified by introducing additional corrections to the interaction Morawetz functional itself.

The bilinear $L^2_{t,x}$ bounds proved in Section~\ref{sec:Morawetz} close the bootstrap loop for the bilinear $L^2_{t,x}$ bounds,
and thus complete the proof of 
Theorem~\ref{t:local-fe3-lin-para} for the linear paradifferential flow, 
Theorems~\ref{t:local-fe3} for the full nonlinear flow, and Theorem~\ref{t:local-fe3-lin} for the linearized flow.

\medskip
\noindent
\emph{IX. Rough solutions as limits of smooth solutions.} 
The frequency envelope bounds in 
Theorem~\ref{t:local-fe3} and Theorem~\ref{t:local-fe3-lin} apply to sufficiently regular solutions to the water wave system. This leaves us with the task of 
constructing rough solutions as unique limits of smooth solutions.
This is carried out in Section~\ref{sec:rough}, using the global frequency envelope bounds for the full equation (tame estimates) combined with the low regularity, local in time weak Lipschitz bounds for the linearized equation. This argument is relatively standard by now, 
following the general strategy described in \cite{IT-primer}. We also refer the reader to the abstract result of \cite{ABITZ}, which
unfortunately does not apply directly due to the intermediate role played by the good variable. To efficiently handle the good variable at this stage, we take our cue from 
\cite{AIT}, where a similar difficulty was resolved in lower dimension.

\subsection{Acknowledgments} 
The first  author gratefully acknowledges support from  the NSF grant DMS-2348908, and from a Vilas Associate Fellowship. The second author was supported
by a fellowship in the Simons Society of Fellows. The third author was supported by the NSF grants DMS-2054975 and DMS-2554866.

\subsection{AI disclosure:} All the research in the paper was developed prior to March 2026, when this result was first announced at the Oberwolfach Research Institute for Mathematics \cite{mfo}, Workshop 2610. AI  was used only in a minor way by the second author to proofread a few of the arguments, and only before 
September 2026
(GPT 5.6 Sol and 6.0 Astra). The authors take full responsibility for the contents of this article.


\section{Notation and Preliminaries}
\label{sec:notations}
 We work on $\mathbb{R}^2$, with spatial Fourier transform
\[
\hat{f}(\xi) =\int_{\mathbb{R}^2}e^{-ix\cdot \xi}f(x)\, dx, \quad f(x)=(2\pi)^{-2} \int_{\mathbb{R}^2}e^{ix\cdot \xi}\hat{f}(\xi)\, d\xi.
\]
We write $D=-i\nabla$, $\vert D\vert ^s=(-\Delta)^{s/2}$. Unless otherwise indicated, $L^p$ denotes a spatial norm, while $L^p_{t,x}(I\times \mathbb{R}^2)$ denotes a spacetime norm over the interval $I$. We use $A\lesssim B$ when $A\leq CB$ for a $C$ a positive constant and $A\sim B $ when $A\lesssim B\lesssim A$. Throughout, implicit constants may depend on the choice of Littlewood–Paley cutoffs. For estimates involving pseudodifferential or multilinear operators, they may also depend on finitely many seminorms of the corresponding symbols. Unless otherwise indicated, all estimates are uniform in the frequency parameters.  A spatial translation is denoted by $f^{y}:=f(x-y)$.
\subsection{Littlewood-Paley decomposition}\label{subsec:littlewoodpaley} 

Throughout our analysis, we localize functions to homogeneous frequency annuli in $\mathbb R^2$. We fix a constant
\begin{equation*}
1<\step
\end{equation*}
such that
\begin{equation*}
\step-1\ll1.
\end{equation*}
At certain points, we will need the sum of two frequency vectors pointing in the same direction and belonging to the same annulus $\{\lambda\leq |\xi|\leq \step \lambda\}$ to lie outside that annulus. Our choice of $c$ ensures this property. The precise value will not be important, and it will suffice to choose $\step=1+1/10$. However, to simplify the formulas, we retain $\step$ as a parameter.

Let $p_0$ be a smooth real-valued radial function on $\mathbb R^2$ which is equal to $1$ on $B_1(0)$ and supported on $B_\step(0)$. For every $\lambda>0$, define the low-frequency projection $P_{<\lambda}$ and the continuous Littlewood-Paley projection $S_\lambda$ by
\begin{equation*}
\widehat{P_{<\lambda}f}(\xi)
=p_0(\xi/\lambda)\widehat f(\xi),
\qquad
S_\lambda\coloneq\lambda\frac{d}{d\lambda}P_{<\lambda}.
\end{equation*}
The symbol of $S_\lambda$ is supported in the annulus
\begin{equation*}
\{\xi\in\mathbb R^2:\lambda\leq|\xi|\leq \step \lambda\}.
\end{equation*}
The fundamental theorem of calculus gives, in the homogeneous sense,
\begin{equation*}
f=\int_0^\infty S_\lambda f\frac{d\lambda}{\lambda},
\qquad
P_{<\lambda}f=\int_0^\lambda S_\nu f\frac{d\nu}{\nu}.
\end{equation*}
The first equality is understood, for example, in $L^2$ for $f\in L^2$, or modulo polynomials in the homogeneous tempered-distribution setting. The second identity also holds in $L^2$ for $f\in L^2$, but requires an appropriate low-frequency normalization in the homogeneous distributional setting.

From this, we can define an associated discrete decomposition with respect to $\step^{\mathbb{Z}}$. Let $\lambda\in \step^{\mathbb Z}$ and set
\begin{equation*}
P_\lambda
\coloneq P_{<\lambda}-P_{<\lambda/\step}
=\int_{\lambda/\step}^\lambda S_\nu\frac{d\nu}{\nu}.
\end{equation*}
Thus $P_\lambda$ has symbol
\begin{equation*}
p_\lambda(\xi)
=
p_0(\xi/\lambda)-p_0(\step \xi/\lambda),
\end{equation*}
and
\begin{equation*}
\operatorname{supp}p_\lambda
\subset
\{\xi\in\mathbb R^2:\step^{-1}\lambda\leq|\xi|\leq \step\lambda\},
\qquad
\sum_{\lambda\in \step^{\mathbb Z}}p_\lambda(\xi)=1,
\quad \xi\neq0.
\end{equation*}
We write $f_\lambda=S_\lambda f$ in continuous frequency integrals and $f_\lambda=P_\lambda f$ in discrete frequency sums. We refer to the elements of $\step^{\mathbb Z}$ as dyadic frequencies. 
\\

For positive frequency scales $\lambda,\mu>0$, we use the convention 
\[
\lambda\sim_\step \mu  \quad  \mbox{ when } \quad \step^{-4}\mu<\lambda< \step^4\mu
\]
and say $\mu\not\sim_\step \lambda$ otherwise. For discrete scales $\lambda,\mu\in \step^{\mathbb Z}$, the relation $\lambda\sim_\step \mu$ means that the scales are at most three lattice steps apart. We reserve the notation $a\approx b$ to mean $C^{-1}a\leq b\leq Ca$ for some fixed universal constant $C$ (possibly much larger than $c$ above).
The projections $P_{\lambda}$ are smooth multipliers. Their multiplier kernels satisfy, for every $N\geq 0$,
\[
\vert K_{\lambda} (x)\vert \lesssim _{N} \lambda^2 (1+\lambda \vert x\vert )^{-N}
\]
with 
\[
\Vert P_{\lambda}f \Vert_{L^{q}} \lesssim \Vert f\Vert_{L^q}, \qquad 1\leq q\leq \infty .
\]
We will use the Bernstein estimates
\[
\Vert \vert D\vert ^s  f_{\lambda}\Vert_{L^q}\lesssim_{s} \lambda^s\Vert f_{\lambda}\Vert_{L^q}, \qquad  \Vert \vert f_{\lambda}\Vert_{L^q}\lesssim_{s} \lambda^{2(1/p-1/q)}\Vert f_{\lambda}\Vert_{L^q}, \quad 1\leq p\leq q\leq \infty,
\]
with the analogous reverse bound for $\vert D\vert ^s$ on an annulus. In particular, the usual square-function
characterizations of Sobolev spaces hold, with constants permitted to depend on $\step$.

\subsection{Paraproducts and remainders.}\label{subsec:paraproductssec2}
We use the Littlewood-Paley decomposition introduced above
to define the paraproduct
\[
T_f g
:=
\sum_{\lambda\in c^{\mathbb Z}}
P_{<\lambda/\step^K}f\,P_\lambda g,
\]
where $K$ is a fixed sufficiently large integer.
Thus, $T_f g$ corresponds to the low-high frequency
interactions, with the frequency of $f$ well below
that of $g$. Alternatively we can define the continuous version 
\[
T_f g
:=
\int 
P_{<\lambda/\step^K}f\,P_\lambda g \frac{d\lambda}{\lambda}.
\]
Both of these versions will be used in the paper. While not equal,
the choice\footnote{They can also be seen as identical relative to different choices of the spectral projectors.}  makes essentially no difference at any point, so we will 
largely use them interchangeably. 

We also define the remainder
\[
R(f,g):=fg-T_f g.
\]
In particular, $R(f,g)$ contains both the high-low
and the balanced frequency interactions.

\subsection{Linear and multilinear operators}

Since our problem is invariant with respect to translations, translation invariant multilinear forms play an important role in the analysis. Below we review some of our notations in this regard, following 
\cite{IT-global}, \cite{IT-qnls} and \cite{IT-qnls2}.

We begin with multipliers $m(D)$, which correspond to multiplication in the Fourier space by their symbol $m(\xi)$, and which in the physical space  are interpreted as convolution operators,
\[
m(D) u (x) = \int K(y) u(x-y) \,dy, \quad \hat K = m.
\]
More generally, we will denote by $\calL$  any convolution operator 
\[
\calL u (x) = \int K(y) u^y(x)\, dy,
\]
where $K$ is an integrable kernel, or 
 a bounded measure, with a universal bound and where for translations we use the notation
\[
u^y(x) := u(x-y).
\]

\medskip

 Moving on to multilinear operators,  we denote by $\calL$ any multilinear form 
\[
\calL(u_1,\cdots, u_k) (x) = 
\int K(y_1, \cdots,y_k) u_1^{y_1}(x) \cdots u_k^{y_k}(x) \, dy,
\]
where, again, $K$ is assumed to have a kernel which is integrable, or more generally,  a bounded measure (in order to be able to include products here). Such multilinear forms satisfy the same bounds as 
corresponding products do, as long as we work in translation invariant Sobolev spaces. In many estimates we need to enforce some kernel localization. 
Given some frequency scale $\lambda$, define $\calL_\lambda$ to have kernels satisfying
\[
|K_\lambda(y)| \lesssim \frac{\lambda^d}{(
1+\lambda |x|)^N} .
\]

Such multilinear forms may also be described by their symbol $a(\xi^1, \cdots \xi^k)$, as
\[
\mathcal F \calL(u_1,\cdots, u_k)(\xi) = (2\pi)^{-\frac{n(k-1)}2}\int_{\xi^1+ \cdots \xi^k = \xi} 
a(\xi^1, \cdots \xi^k) \hat u_1(\xi^1) \cdots \hat u_k(\xi^k) d \xi^1 \cdots d \xi^{k-1},
\]
where $a$ is nothing but the  inverse Fourier transform of the kernel $K$.

\medskip

In the study of our two dimensional problem we also need to consider multilinear forms involving both $u$ and $\bar u$. Since there is no phase rotation symmetry, the conjugated and nonconjugated arguments need not alternate. For such forms we adapt the notations from our earlier paper \cite{IT-global}.

Precisely, for an integer $k \geq 2$ and a choice of signs $\epsilon=(\epsilon_1,\ldots,\epsilon_k)\in\{+1,-1\}^k$, we will use translation invariant $k$-linear forms
\[
(\mathcal D(\R^2))^{k} \ni (u_1, \cdots, u_{k}) \to Q_\epsilon(u_1^{\epsilon_1},\cdots,u_k^{\epsilon_k}) \in \mathcal D'(\R^2),
\]
where $u_j^+=u_j$ and $u_j^-=\bar u_j$. There is no restriction on the choice of signs.

Such a form is uniquely described by its symbol $q_\epsilon(\xi^1,\xi^2,\cdots,\xi^k)$ via
\[
\begin{aligned}
Q_\epsilon(u_1^{\epsilon_1},\cdots,u_k^{\epsilon_k})(x) = (2\pi)^{-2k} &
\int e^{i\epsilon_1(x-x_1)\cdot\xi^1}\cdots e^{i\epsilon_k(x-x_k)\cdot\xi^k}
q_\epsilon(\xi^1,\cdots,\xi^k)
\\ & \qquad
u_1^{\epsilon_1}(x_1)\cdots u_k^{\epsilon_k}(x_k)
\,dx_1\cdots dx_k\,d\xi^1\cdots d\xi^k,
\end{aligned}
\]
or equivalently on the Fourier side
\[
\mathcal F Q_\epsilon(u_1^{\epsilon_1},\cdots,u_k^{\epsilon_k})(\xi)
= (2\pi)^{-(k-1)} \int_{D_\epsilon}
q_\epsilon(\xi^1,\cdots,\xi^k)
\widehat{u_1^{\epsilon_1}}(\epsilon_1\xi^1)\cdots
\widehat{u_k^{\epsilon_k}}(\epsilon_k\xi^k)
\,d\xi^1\cdots d\xi^{k-1},
\]
where
\[
D_\epsilon := \{\xi=\epsilon_1\xi^1+\cdots+\epsilon_k\xi^k\}.
\]
Thus the Fourier factor corresponding to $\epsilon_j=+1$ is $\hat u_j(\xi^j)$, while that corresponding to $\epsilon_j=-1$ is $\overline{\hat u_j(\xi^j)}$.

They can also be described via their kernel
\[
Q_\epsilon(u_1^{\epsilon_1},\cdots,u_k^{\epsilon_k})(x) =
\int K_\epsilon(x-x_1,\cdots,x-x_k)
u_1^{\epsilon_1}(x_1)\cdots u_k^{\epsilon_k}(x_k)
\,dx_1\cdots dx_k,
\]
where $K_\epsilon$ is defined in terms of the Fourier transform of $q_\epsilon$
\[
K_\epsilon(x_1,\cdots,x_k) =
(2\pi)^{-k}\widehat{q_\epsilon}(-\epsilon_1x_1,\cdots,-\epsilon_kx_k).
\]

These notations allow us to treat all combinations of conjugated and nonconjugated arguments in the same way; our usage will be clear from the context. In the absence of phase rotation symmetry, the roles of multilinear forms as functions or densities are not determined by whether their degree is odd or even. We will use forms of any degree as functions in evolution equations, or as densities in density-flux pairs, as appropriate.

\medskip
In particular,
to each $k$-linear form $Q_\epsilon$ which is used as a density we will associate
a corresponding $k$-linear functional $\mathbf{Q}_\epsilon$ defined by
\[
\mathbf{Q}_\epsilon(u_1,\cdots,u_k) := \int_{\R^2} Q_\epsilon(u_1^{\epsilon_1},\cdots,u_k^{\epsilon_k})(x)\, dx,
\]
which takes real or complex values.
This may be alternatively expressed
on the Fourier side as
\[
\mathbf{Q}_\epsilon(u_1,\cdots,u_k) = (2\pi)^{2-k} \int_{D_{\epsilon,0}}
q_\epsilon(\xi^1,\cdots,\xi^k)
\widehat{u_1^{\epsilon_1}}(\epsilon_1\xi^1)\cdots
\widehat{u_k^{\epsilon_k}}(\epsilon_k\xi^k)
\,d\xi^1\cdots d\xi^{k-1},
\]
where the diagonal $D_{\epsilon,0}$ is given by
\[
D_{\epsilon,0} = \{0=\epsilon_1\xi^1+\cdots+\epsilon_k\xi^k\}.
\]
In order to define the multilinear functional $\mathbf{Q}_\epsilon$ it suffices to know the symbol $q_\epsilon$ on $D_{\epsilon,0}$.

We denote by $S^m_{1,0} $ the class of smooth symbols $a(x,\xi)$ for which 
\begin{equation}
\vert \partial_x^{\alpha}\partial_{\xi}^{\beta}a(x,\xi)\vert \leq C_{\alpha\beta}\langle \xi\rangle^{m-\vert\beta\vert},\quad \alpha, \beta \in \mathbb{N}^2_{0}.
\end{equation}
For both linear and multilinear forms,
we use $S^m_\lambda$ to denote symbols 
localized at frequency $\lambda$ with regularity on the dyadic scale,
\begin{equation}
\vert \partial_{\xi}^{\alpha}a(\xi)\vert \leq C_{\alpha\beta}\lambda^{m-\vert\alpha\vert},\quad \alpha \in \mathbb{N}^2_{0}.
\end{equation}

\subsection{Strichartz and bilinear estimates}\label{subsec:dispersive-cc}
Here we collect the Strichartz and bilinear $L^2_{t,x}$ estimates for the
linear water wave flow around the zero solution. This has the form 
\begin{equation}
\label{lin-sys}
\left\{\begin{array}{l}
h_t+|D|q= 0
\cr
q_t - h = 0,
\end{array}\right.
\end{equation}
for real valued linearized variables $(h,q)$. For the dispersive analysis it is 
more convenient to diagonalize and switch to a single complex valued function 
\[
u= h+ i|D|^\frac12 q.
\]
This solves an equation of the form
\begin{equation}
\label{lin}
\left\{\begin{array}{l}
(\partial_t + i \omega_0(D)u=f,\qquad u:\mathbb R\times\mathbb R^2\to\mathbb C,\\[2mm]
u(0,x)=u_0(x),
\end{array}\right.
\end{equation}
with the dispersion relation 
\begin{equation}\label{omega0-def}
\omega_0(\xi) =|\xi|^{\frac12},
\end{equation}
where we have also added a source term $f$ for the discussion below. A similar diagonalization will be carried out later 
for the full equation at the paradifferential level.

We write $u_\lambda=P_\lambda u$, $u_{\lambda,0}=P_\lambda u_0$,
and $f_\lambda=P_\lambda f$, where $\lambda\in c^{\mathbb Z}$.
All spacetime norms below are taken over $I\times\mathbb R^2$,
where $I$ is an interval containing $0$.  Constants are independent
of $I$ and of the frequency scales.

The frequency-localized Strichartz estimates reflect the curvature
of the dispersion relation $\varphi(\xi)=|\xi|^{\frac12}$.
We call a pair $(p,q)$ admissible if
\begin{equation*}
\frac{2}{p}+\frac{2}{q}=1,\qquad 2\leq p,q\leq\infty,
\qquad (p,q)\neq(2,\infty).
\end{equation*}
The endpoint $(2,\infty)$ is excluded in two space dimensions.

\begin{lemma}[Strichartz estimates]
Let $u$ solve \eqref{lin}. For every admissible pair $(p,q)$,
\begin{equation}
\label{eq:halfwave-strichartz}
\|u_\lambda\|_{L^p_tL^q_x}
\lesssim \lambda^{\frac{3}{2p}}
\bigl(\|u_{\lambda,0}\|_{L^2_x}
      +\|f_\lambda\|_{L^1_tL^2_x}\bigr).
\end{equation}
In particular, for a solution of the homogeneous equation,
\begin{equation}
\label{eq:halfwave-L4}
\|u_\lambda\|_{L^4_{t,x}}
\lesssim \lambda^{\frac38}\|u_{\lambda,0}\|_{L^2_x}.
\end{equation}
\end{lemma}

We next record the bilinear $L^2_{t,x}$ estimates. In the unbalanced case,
the different group velocities of the two factors yield a gain over
an application of the Strichartz inequality to each factor.

\begin{lemma}[Unbalanced bilinear estimate]
Let $u$ and $v$ solve \eqref{lin}, with initial data $u_0,v_0$
and forcing terms $f,g$, respectively. If $\mu\ll\lambda$, then
\begin{equation}
\label{eq:halfwave-bilinear-unbalanced}
\|u_\lambda v_\mu\|_{L^2_{t,x}}
\lesssim \mu^{\frac34}
\bigl(\|u_{\lambda,0}\|_{L^2_x}+\|f_\lambda\|_{L^1_tL^2_x}\bigr)
\bigl(\|v_{\mu,0}\|_{L^2_x}+\|g_\mu\|_{L^1_tL^2_x}\bigr).
\end{equation}
The same estimate holds if either factor is replaced by its
complex conjugate.
\end{lemma}

\noindent In the balanced case we use Strichartz and the localization of the
product to frequencies $\lesssim\lambda$.

\begin{lemma}[Balanced bilinear estimate]
With $u,v$ as above, at a common frequency $\lambda$ we have
\begin{equation}
\label{eq:halfwave-bilinear-balanced}
\bigl\||D|^{\frac12}(u_\lambda\overline{v_\lambda})\bigr\|_{L^2_{t,x}}
\lesssim \lambda^{\frac54}
\bigl(\|u_{\lambda,0}\|_{L^2_x}+\|f_\lambda\|_{L^1_tL^2_x}\bigr)
\bigl(\|v_{\lambda,0}\|_{L^2_x}+\|g_\lambda\|_{L^1_tL^2_x}\bigr).
\end{equation}
\end{lemma}

Indeed, the last bound follows from the $L^4_{t,x}$ estimate,
H\"older's inequality and Bernstein's inequality. The
unbalanced estimate instead exploits the separation of the group
velocities $\nabla\varphi(\xi)=\frac12|\xi|^{-\frac32}\xi$
at widely separated frequency scales. These bounds will serve as
reference estimates for the variable-coefficient flows considered later.

\subsection{Resonant interactions}
\label{subsec:resonant}

We discuss both quadratic and cubic nonlinear
interactions, relative to the constant coefficient dispersion relation $\omega_0$ in \eqref{omega0-def}; however, the discussion will carry over ad-litteram to the full dispersion relation $\omega$ in the nonlinear problem.

\medskip
We begin with the quadratic terms. For simplicity, consider an equation
of the form
\[
(\partial_t+i|D|^{\frac12})u=C^{[2]}(\mathbf u,\mathbf u),
\]
where $C^{[2]}$ is a translation invariant bilinear form with symbol $c^{[2]}(\xi^1,\xi^2)$ and
$\mathbf u=\{u,\bar u\}$ allows for either $u$ or $\bar u$ in each argument.
This includes $C^{[2]}(u,u)$, $C^{[2]}(u,\bar u)$ and $C^{[2]}(\bar u,\bar u)$.

Given two input frequencies $\xi^1,\xi^2$, the output frequency is
\begin{equation}\label{eq:hw-quadratic-frequencies}
\xi^3=\pm\xi^1\pm\xi^2,
\end{equation}
where the signs depend on the arguments of $C^{[2]}$. Such an interaction is
\emph{resonant} if the corresponding time frequencies satisfy the same
relation,
\begin{equation}\label{eq:hw-quadratic-time}
\omega_0(\xi^3)=\pm\omega_0(\xi^1)\pm\omega_0(\xi^2),
\end{equation}
with the same choice of signs as in \eqref{eq:hw-quadratic-frequencies}.

If both signs are minus, \eqref{eq:hw-quadratic-time} reads
$\omega_0(\xi^1)+\omega_0(\xi^2)+\omega_0(\xi^3)=0$ and all three frequencies
vanish. In the two remaining cases, after relabeling the frequencies,
\eqref{eq:hw-quadratic-frequencies} takes the form
\[
\Delta^3\xi=0,\qquad \Delta^3\xi:=\xi^1+\xi^2-\xi^3,
\]
and \eqref{eq:hw-quadratic-time} becomes
\[
\Delta^3\omega_0=0,\qquad
\Delta^3\omega_0:=\omega_0(\xi^1)+\omega_0(\xi^2)-\omega_0(\xi^3)
=|\xi^1|^{\frac12}+|\xi^2|^{\frac12}-|\xi^3|^{\frac12}.
\]
On $\{\Delta^3\xi=0\}$ we have
\[
\omega_0(\xi^3)=|\xi^1+\xi^2|^{\frac12}
\le\big(|\xi^1|+|\xi^2|\big)^{\frac12}
\le|\xi^1|^{\frac12}+|\xi^2|^{\frac12},
\]
and the second inequality is strict unless $\xi^1=0$ or $\xi^2=0$. Hence,
unlike the cubic case, which we will discuss next, there are no resonant quadratic interactions away
from the zero frequency.
To determine the size of $\Delta^3\omega_0$ 
we note that the two highest frequencies must be comparable, and compare the corresponding $\omega_0$ signs, distinguishing two cases:

\begin{enumerate}[label=(\alph*)]
    \item Strongly nonresonant, when the $\omega_0$ signs for the two highest frequencies are matched,
    \item Weakly nonresonant, when the  $\omega_0$ signs for the two highest frequencies are mismatched.
\end{enumerate}
For convenience, we place the case of three 
comparable frequencies in the strongly nonresonant box.  In the first case, the size of $\Delta^3 \omega_0$ is given by the highest frequency, while in the second case is governed by the smallest of the three frequencies. We summarize our observations in
the following lemma.

\begin{lemma}\label{lem:hw-quadratic-resonance}
Except for interactions in which one of the frequencies
$\xi^1,\xi^2,\xi^3$ vanishes, all quadratic interactions are nonresonant.
Quantitatively, on $\{\Delta^3\xi=0\}$, we have 

a) in the strongly resonant case,
\[
\max\{|\xi^1|,|\xi^2|,|\xi^3|\}^{\frac12}
\le |\Delta^3\omega_0|
\le2\max\{|\xi^1|,|\xi^2|,|\xi^3|\}^{\frac12}.
\]
b) in the weakly resonant case,
\[
\tfrac12\min\{|\xi^1|,|\xi^2|,|\xi^3|\}^{\frac12}
\le |\Delta^3\omega_0|
\le2\min\{|\xi^1|,|\xi^2|,|\xi^3|\}^{\frac12}.
\]
\end{lemma}

We thus conclude that, away from the zero frequency, all quadratic
interactions are nonresonant, and the quadratic terms can be removed by a
normal form transformation
\begin{equation}
 v = u + B^{[2]}(\bfu,\bfu).
\end{equation}
The normal form correction eliminates the quadratic terms in the $v$ equation but generates cubic terms,
\[
(\partial_t+i|D|^{\frac12})v=C^{[3]}(\bfu,\bfu,\bfu),
\]
where, assuming $B^{[2]}$ is symmetric,
we schematically have
\[
C^{[3],NF}(\bfu,\bfu,\bfu) = 2 B^{[2]}(\mathbf C^{[2]}(\bfu,\bfu),\bfu)
\]
with the bold $C$ indicating the possible conjugation of this entry.
Here $C^{[3],NF}$ plays the role of  a translation invariant trilinear form and, as before,  $\mathbf u=\{u,\bar u\}$ allows for either $u$ or $\bar u$ in each argument.
The $C^{[3],NF}$ above must be added to $C^{[3]}$ below, before the
resonant analysis. As it will become apparent next, a similar analysis as the one classifying the types of quadratic interactions will be needed in the cubic case as well.   In particular, the
symbol of $C^{res}$ (defined below) on the diagonal $\mathcal R_2$ (again, defined below) receives a contribution
from $C^{[3],NF}$.

A direct computation shows that the symbol $b^{[2]}$ of the normal form correction 
has the form, depending on the sign choices above,
\begin{equation}
\b^{[2]}(\xi^1,\xi^2) =     \frac{i c^{[2]}(\xi^1,\xi^2)}{\pm \omega_0(\xi^1) \pm \omega_0(\xi^2) - \omega_0(\xi^3)}. 
\end{equation}

To track the size of the correction it is natural to separate $C^{[2]}$ into
\[
C^{[2]}(\bfu,\bfu)=C^{[2],wnr}(\bfu,\bfu)+C^{[2],nr}(\bfu,\bfu),
\]
where
\begin{itemize}
\item $C^{[2],wnr}(\bfu,\bfu)$ contains the weakly nonresonant interactions;
\item $C^{[2],nr}(\bfu,\bfu)$ contains the strongly nonresonant
interactions.
\end{itemize}
For clarity, we note that the weakly nonresonant interactions include only 
quadratic terms of the form
\begin{equation}\label{wnr}
C^{[2]}(u^l,u^h), \qquad C^{[2]}(\bu^l,u^h),
\qquad C^{[2]}(\bu^h,u^h),
\end{equation}
where the ``h" and ``l" superscripts are used to indicate high, respectively low frequencies. By
Lemma~\ref{lem:hw-quadratic-resonance}, in this case the symbol of the normal form
carries a factor of size $(\Delta^3\omega_0)^{-1}\approx\min\{|\xi^1|,|\xi^2|,|\xi^3|\}^{-\frac12}$,
so the transformation is singular at low frequencies unless the symbol of
$C^{[2]}$ vanishes there; this is the reason for separating $C^{[2],wnr}$ -- and we broadly call this \emph{null structure}. 

The way this null structure will manifest itself in the analysis of the full nonlinear
equation depends on which terms in \eqref{wnr} we consider:
\begin{itemize}
    \item The two ``lh" interactions in \eqref{wnr} will be included in the paradifferential flow and will provide at leading order the transport term in the equation, which will be \emph{conservative}.
  \item The ``hh" interaction  in \eqref{wnr}
  will be separated, and its symbol $c^{[2],hh}(\xi^1,\xi^2)$ will vanish when $\xi^1= \xi^2$; this, we will directly refer to as the \emph{null condition}.
\end{itemize}

\bigskip

We now focus on the cubic case, where we 
\[
(i\partial_t+|D|^{\frac12})u=C^{[3]}(\mathbf u,\mathbf u,\mathbf u).
\]
Again, recall $C^{[3]}$ is a translation invariant trilinear form and, as before, 
$\mathbf u=\{u,\bar u\}$ allows for either $u$ or $\bar u$ in each argument.
This includes the following possible trilinear forms: $C(u,\bar u,u)$, $C(u,u,u)$, $C(u,\bar u,\bar u)$ and
$C(\bar u,\bar u,\bar u)$. The first expression has phase rotation symmetry, but overall we make no phase rotation symmetry assumption.

A linear wave at spatial frequency $\xi$ oscillates in time at frequency
$\omega_0(\xi)$, and its complex conjugate sits at spatial frequency $-\xi$
and time frequency $-\omega_0(\xi)$.

\smallskip

Given three input frequencies $\xi^1,\xi^2,\xi^3$, the output frequency is
\begin{equation}\label{eq:hw-cubic-frequencies}
\xi^4=\pm\xi^1\pm\xi^2\pm\xi^3,
\end{equation}
where the signs depend on the arguments of $C$. We call an interaction
\emph{balanced} if all the frequencies $\xi^j$ have comparable size, and
\emph{unbalanced} otherwise. Such an interaction is \emph{resonant} if the
corresponding time frequencies satisfy the same relation,
\begin{equation}\label{eq:hw-cubic-time}
\omega_0(\xi^4)=\pm\omega_0(\xi^1)\pm\omega_0(\xi^2)\pm\omega_0(\xi^3),
\qquad\text{i.e.}\qquad
|\xi^4|^{\frac12}=\pm|\xi^1|^{\frac12}\pm|\xi^2|^{\frac12}\pm|\xi^3|^{\frac12},
\end{equation}
with the same choice of signs as in \eqref{eq:hw-cubic-frequencies}.

From the viewpoint of bilinear estimates we single out the
\emph{transversal} interactions, namely those for which no more than two of
the four frequencies $\xi^1,\xi^2,\xi^3,\xi^4$ coincide. We now classify
the interactions from this perspective:

\medskip
\noindent
(i) \emph{The case of phase rotation symmetry}, with two-plus one-minus interactions. After relabeling the inputs
so that the conjugated argument sits in the second slot,
\eqref{eq:hw-cubic-frequencies} takes the form
\[
\Delta^4\xi=0,\qquad \Delta^4\xi:=\xi^1-\xi^2+\xi^3-\xi^4,
\]
and \eqref{eq:hw-cubic-time} becomes
\[
\Delta^4\omega_0=0,\qquad
\Delta^4\omega_0:=\omega_0(\xi^1)-\omega_0(\xi^2)+\omega_0(\xi^3)-\omega_0(\xi^4)
=|\xi^1|^{\frac12}-|\xi^2|^{\frac12}+|\xi^3|^{\frac12}-|\xi^4|^{\frac12}.
\]
Accordingly we define the resonant set for this sign configuration as
\[
\mathcal R:=\{\Delta^4\xi=0,\ \Delta^4\omega_0=0\}.
\]
For the quadratic dispersion relation $|\xi|_g^2:=g^{ij}\xi_i\xi_j$ of
the Schr\"odinger case in \cite{IPT-ultra}, where $g$ is a symmetric
nondegenerate matrix, restricting to $\{\Delta^4\xi=0\}$ one has
\[
\Delta^4|\xi|_g^2=2\,g^{ij}(\xi^1-\xi^2)_i(\xi^2-\xi^3)_j,
\]
so that $\mathcal R$ consists of the quadruples for which $\xi^1-\xi^2$ is
orthogonal to $\xi^2-\xi^3$ with respect to $g$, and in one dimension
$\mathcal R$ reduces to the trivial resonances
$\{\xi^1,\xi^3\}=\{\xi^2,\xi^4\}$. Neither statement survives for
$\omega_0(\xi)=|\xi|^{\frac12}$. The resonant set is strictly larger already
when $n=1$.

In one dimension we claim that
\[
\mathcal R=\big\{\{\xi^1,\xi^3\}=\{\xi^2,\xi^4\}\big\}\cup\mathcal R_{nt},
\]
where $\mathcal R_{nt}$ is the nontrivial branch obtained in (c) below. We
go through the sign configurations of the frequencies.

\emph{(a) $\xi^1,\xi^3$ have the same sign and $\xi^2,\xi^4$ have the same
sign.} Then $\Delta^4\xi=0$ forces all four signs to agree (or all
frequencies to vanish), and $|\xi^1|+|\xi^3|=|\xi^2|+|\xi^4|$. Squaring
$\Delta^4\omega_0=0$ gives $|\xi^1\xi^3|=|\xi^2\xi^4|$, hence
$\xi^1\xi^3=\xi^2\xi^4$; equal sums and equal products give
\[
\{\xi^1,\xi^3\}=\{\xi^2,\xi^4\}.
\]

\emph{(b) Each of the pairs $(\xi^1,\xi^3)$ and $(\xi^2,\xi^4)$ consists of
frequencies of opposite signs.} Then the square roots of the moduli have
equal sums by $\Delta^4\omega_0 =0$ and, dividing $\Delta^4\xi=0$ by this,
equal differences, so the pairs coincide and the resonance is trivial.

\emph{(c) Three frequencies share a sign and the fourth has the opposite
sign.} Up to the symmetries of $\mathcal R$ we may take
$\xi^1,\xi^3,\xi^4>0>\xi^2$, and with $\xi^1=\gamma^2$, $\xi^3=\delta^2$
the two relations are solved by
\[
\xi^2=-\frac{\gamma^2\delta^2}{(\gamma+\delta)^2},\qquad
\xi^4=\frac{(\gamma^2+\gamma\delta+\delta^2)^2}{(\gamma+\delta)^2},
\qquad \gamma,\delta>0 .
\]
This is the nontrivial branch $\mathcal R_{nt}$; on it at most two
frequencies coincide, which proves the claim.

However, in higher dimensions the resonant set $\mathcal R$ is larger
still, and admits no description comparable to the orthogonality relation
of the quadratic case.

Most resonant quadruples split into two pairs of mutually separated
frequencies; these are the \emph{transversal resonant} quadruples. The
non-transversal resonant interactions are given by the diagonal
\[
\mathcal R_2:=\{\xi^1=\xi^2=\xi^3=\xi^4\},
\]
since three coincident frequencies force the fourth by $\Delta^4\xi=0$.
We refer to these as the \emph{doubly resonant} interactions. Note that,
restricted to $\{\Delta^4\xi=0\}$, $\Delta^4\omega_0$ vanishes to second
order along $\mathcal R_2$.

\medskip
\noindent
(ii) \emph{The remaining sign configurations.} Here we describe the worst
case, namely interactions which are both resonant and non-transversal. The
first requirement is
\[
\xi^1\pm\xi^2\pm\xi^3\pm\xi^4=0,\qquad
\omega_0(\xi^1)\pm\omega_0(\xi^2)\pm\omega_0(\xi^3)\pm\omega_0(\xi^4)=0,
\]
with the same signs, and the second is that three of the frequencies are
equal. Without any loss in generality assume that $\xi^1=\xi^2=\xi^3:=\xi$.
Depending on the choice of signs, and excluding the two-plus one-minus
case, we obtain one of the systems
\[
\xi^4=\pm3\xi,\qquad \omega_0(\xi^4)=\pm3\,\omega_0(\xi),
\]
with matched signs, or
\[
\xi^4=-\xi,\qquad \omega_0(\xi^4)=-\omega_0(\xi).
\]
Since $\omega_0(\pm3\xi)=\sqrt3\,\omega_0(\xi)$, $\omega_0(-\xi)=\omega_0(\xi)$
and $\omega_0$ vanishes only at the origin, $\xi=0$ is the only solution in
each case.

This is the analogue of the positive definite Schr\"odinger case: for the quadratic
dispersion relation the same computation leads to $|\xi|_g^2=0$, which
has only the trivial solution if $g$ is positive definite, but is solved
by the null cone of $g$ in the ultrahyperbolic case, where it gives rise
to the weakly transversal interactions. Here the null set of $\omega_0 $ is
trivial, so no weakly transversal interactions occur. We summarize our
observations in the following lemma.

\begin{lemma}\label{lem:hw-cubic-resonance}
Except for the case $(0,0,0,0)$, all cubic interactions other than the
ones with phase rotation symmetry are either nonresonant or transversal. For the
cubic interactions with phase rotation symmetry, the non-transversal resonant set is the
diagonal $\mathcal R_2$.
\end{lemma}

We will use the above classification to separate the cubic interactions
using smooth frequency cut-offs as
\[
C^{[3]}(\mathbf u,\mathbf u,\mathbf u)
=C^{[3],res}(u,\bar u,u)
+C^{[3],nr}(\mathbf u,\mathbf u,\mathbf u)+C^{[3],tr}(\mathbf u,\mathbf u,\mathbf u),
\]
where, after relabeling the arguments, the conjugate in $C^{[3],res}$ is placed
in the middle slot, and
\begin{itemize}
\item $C^{[3],res}(u,\bar u,u)$ contains the balanced 
interactions with phase rotation symmetry localized near the diagonal $\mathcal R_2$;
\item $C^{[3],nr}(\mathbf u,\mathbf u,\mathbf u)$ contains only balanced
nonresonant interactions;
\item $C^{[3],tr}(\mathbf u,\mathbf u,\mathbf u)$ contains only transversal
interactions, balanced or unbalanced.
\end{itemize}
Unlike the ultrahyperbolic case, no weakly transversal term $C^{wt}$ is
needed. Heuristically, the nonresonant interactions can be corrected via normal form, while transversal interactions can be 
controlled using bilinear $L^2$ estimates. 
This is why the resonant term will play a leading role in our analysis later in the paper.

In $L^2$ energy estimates, the energy growth caused by the resonant term is by the integral 
\[
\bR^{[4]}(u,\bu,u,\bu) = 2 \int \Re C^{[3],res}(u,\bar u,u) \cdot \bu \, dx
\]
whose leading contribution is determined by the behavior of the (fully symmetrized) symbol $r^{[4]}(\xi^1,\xi^2,\xi^3,\xi^4)$ near the doubly resonant set $\calR_2$:
\begin{itemize}
\item  The simple vanishing of $r^{[4]}$ on $\calR_2$, referred to as the \emph{conservative} condition in \cite{IT-qnls2}, \cite{IPT-ultra}, will be called here
\emph{null condition}.

\item The double vanishing of $r^{[4]}$ on $\calR_2$ will be called here the
\emph{double null condition}, and plays an important role in this paper.
\end{itemize}

Finally, we note  that the correct interpretation of both the null condition and the double null condition requires combining the direct cubic terms with the 
secondary cubic terms obtained from the 
quadratic ones via the normal form analysis.


\section{Water wave variables and their relations}
\label{sec:wwvar}
This section serves as a starting point for determining a good formulation of the water waves equations for our purposes. We introduce all the main variables associated 
with the water wave system, and then 
describe first the static and then the dynamic relations between them, concluding with the derivation of  the associated linearized system. This is needed in order to unveil the underlying structure of the problem, which is reflected not only in the analysis we perform, but also in the results we end up proving. We do this in a modular fashion, to allow for an easy tracking of the main variables and their properties.

\subsection{The main variables} 
For completeness, we begin with the 
original variables $(\eta,\psi)$, as well as $\phi$, as follows:

\begin{itemize}
    \item $\eta$ denotes the defining of the free boundary of  $\Gamma_t$, i.e. 
    \[
    \Gamma_t := \{ z = \eta(t,x)\}.
    \]
    In our setting, $\eta$ will always be Lipschitz continuous in $x$, with a small Lipschitz constant, uniformly in $t$.
    \item $\phi$ denotes the velocity potential in the fluid domain 
    \[
    \Omega_t = \{ z < \eta(t,x)\},
    \]
    and  is a harmonic function in the fluid domain, whose gradient vanishes
    at infinite depth. It is only defined modulo constants.
    
    \item $\psi$ denotes the trace of the velocity potential on the free surface $\Gamma_t$,
    \[
    \psi(t,x) := \phi(t,x,\eta(t,x)).
    \]
\end{itemize}

\subsection{The  horizontal and vertical velocities \texorpdfstring{$(V,B)$}{}}
With the above  notations, the fluid velocity
is given by $\nabla_{x,z} \phi$; we split its  trace on the free surface,  into two components:

\begin{itemize}
    \item The horizontal velocity $V$ on the free surface,
    \[
    V(t,x) := \nabla_x \phi(t,x,\eta(t,x)).
    \]
    \item The vertical  velocity $B$ on the free surface,
    \[
    B(t,x) := \partial_z \phi(t,x,\eta(t,x)).
    \]
\end{itemize}
Currently, these variables have a physical meaning, though relative to the free moving surface, they do not track its evolution as closely as one would need. 
\subsection{The Dirichlet to Neuman map}\label{sec:DN}  The Dirichlet to Neumann operator $G(\eta)$ associated to the domain $\Omega(t)$ plays an essential role 
 in the water wave equations, and is defined by 
\[
G(\eta)\psi
:=
\bigl(\partial_z\phi - \nabla\eta\cdot\nabla\phi\bigr)\big|_{z=\eta(x)}.
\]
For fixed $\eta$, this can be seen as an elliptic pseudodifferential operator  of order 1
whose principal symbol is given by
\begin{equation}\label{dn-symbol}
\mathfrak{g}(x,\xi) = [(1+|\nabla \eta|^2) \vert \xi\vert ^2 - (\nabla \eta \cdot \xi)^2]^\frac12.
\end{equation}
On the other hand as a function of $\eta$ this is a nonlinear  operator, which makes it interesting to 
compute its linearization. Following \cite{Lannes}, this is given by
\begin{equation}
\label{G-lin}
DG(\eta)[h] \psi = -G(\eta)(hB) - \nabla_x(hV). 
\end{equation}

In later sections we use the notation $\mathfrak{g}(\eta)$ instead of $\mathfrak{g}(x, \xi)$, omitting the dependence on $(x,\xi)$; this would help clear the exposition in the sections where we compute various linearizations of our variables.

\subsection{ The Taylor coefficient}
The Taylor coefficient $a$ represents the normal derivative of the pressure on the free boundary,
\begin{equation*}
a\coloneq -\partial_z p_{|z=\eta}.
\end{equation*}
where, normalizing  gravity to  $g=1$, $p$ is 
the solution to the elliptic equation
\[
\left\{
\begin{aligned}
&\  \Delta p = -(\partial_i \partial_j \phi)^2\\
&\  p = 0 \qquad on \  \{z = \eta(x)\} \\
&\  p + z \to 0 \ \text{ as } z \to -\infty .
\end{aligned}
\right.
\]
Its role will become apparent
shortly, once we compute the  linearized equations.

\subsection{Static equations}
\label{ss:static eq}
A short computation connects the horizontal velocity $V$ and the vertical
velocity $B$ on the free surface with the trace $\psi$ of the velocity
potential there, through the Dirichlet-to-Neumann operator $G(\eta)$:
\begin{equation}\label{BV-re}
\left\{
\begin{aligned}
V + B \nabla \eta = & \  \nabla \psi
\\
B - V \cdot \nabla \eta = & \ G(\eta) \psi.
\end{aligned}
\right.
\end{equation}
In our analysis, the variables $(V,B)$ 
will be regarded as dependent variables (they depend on $(\eta, \psi)$).
Solving the system \eqref{BV} allows us 
to express them in terms of $\psi$ as 
\[
B = L \psi, \qquad V = M \psi,
\]
where the \textbf{commuting}  operators $L$ and $M$ are defined by 
\begin{equation}
\label{LM}
\begin{aligned}
    L(\eta) :=& \ \frac{1}{1+|\nabla \eta|^2} G(\eta)
+ \nabla\eta \cdot \frac{1}{1+\nabla \eta \otimes \nabla \eta} \nabla_x, \\
M(\eta) := & \ \frac{1}{1+\nabla \eta \otimes \nabla \eta} \nabla - \frac{1}{1+|\nabla \eta|^2} \nabla \eta \, G(\eta),
\end{aligned}
\end{equation}
and will also play an important role in what follows. We denote their principal symbols  
by $\mfm$ and $\mfl$, and the corresponding operators in the left quantization by
$\M$ and $\L$ respectively.

\begin{lemma}
 The differential operators $L(\eta)$, $M(\eta)$ defined in \eqref{LM}  commute, and are related via the dual identities
\begin{equation}\label{LvsM}
 L^*G - M^* \nabla = 0,  \qquad \nabla M + GL = 0.
\end{equation}
\end{lemma}
For later use we also observe that, as a corollary, we have the relation
\begin{equation}\label{GB+DV}
G(\eta) B = - \nabla \cdot V  .  
\end{equation}

\begin{proof}
 To prove this we apply the last expression to $\psi$, 
\[
\nabla M(\eta)\psi  + GL(\eta)\psi = \nabla V + GB .
\]
With $B = \phi_z$, $V= \phi_x$ and $\psi = \phi$
on the boundary we have 
\begin{equation}
\left\{
\begin{aligned}
&\ V + B \nabla \eta =  \  \nabla \psi
\\
& \ B - V \cdot \nabla \eta = \ G(\eta) \psi .
\end{aligned}
\right.
\end{equation}
Replacing $\phi$ with $\phi_z$, the second relation gives
\[
\phi_{zz} - \phi_{zx} \cdot \nabla \eta = G(\eta)B.
\]
Replacing $\phi$ with $\phi_x$, the second relation gives the matrix identity
\[
\phi_{xx} + \phi_{zx}  \nabla \eta = \nabla V,
\]
where taking the trace we get
\[
\Delta_x \phi + \phi_{zx} \cdot \nabla \eta = \nabla \cdot  V,
\]
and the identities \eqref{LvsM} follow.
\end{proof}

We next turn our attention to $p$, which we represent as
\[
p = - \frac12 |\nabla \phi|^2 - z + \tilde{p}
\]
with $\tilde{p}$ harmonic decaying at the infinite bottom. Then $\tilde{p}$ is the harmonic extension of 
\[
\tilde{p}_0:=\frac12 |\nabla \phi|^2 + z =  \frac12(|B|^2+|V|^2) + \eta
\]
and its $z$ derivative satisfies
\[
\partial_z \tilde{p}  = L(\eta) \tilde{p}_0.
\]
On the other hand 
\[
\phi_{zz} =  L(\eta) B, 
\qquad \phi_{zx} = L(\eta) V =M(\eta)B,
\]
so 
\[
\begin{aligned}
\phi_{z} \phi_{zz} + \phi_{x} \phi_{xz}
= & \ B L(\eta) B + V \cdot L(\eta) V .
\end{aligned}
\]

Hence we have
\[
\begin{aligned}
a = & \ \partial_z \tilde{p}_0 - \partial_z \tilde{p} 
\\
= & \ 1 +  \nabla \phi_z \cdot \nabla \phi 
- \frac{1}{1+|\nabla \eta|^2} \left( G(\eta)(\frac12 (|B|^2 +|V|^2)+\eta)
+ \frac12 \nabla \eta \nabla(B^2+V^2) + |\nabla \eta|^2\right) .
\end{aligned}
\]
Using the definition of the operator $L$,
we arrive at
\begin{equation}\label{a-formula}
a =     1+ B L(\eta) B + V \cdot L(\eta) V - L(\eta)
 [\frac12(|B|^2+|V|^2) + \eta].
\end{equation}

\subsection{Dynamic equations}

We first use the velocity fields $(B,V)$
and the associated system \eqref{BV}
to rewrite the system \eqref{WW}
in the form
\begin{equation}\label{WW-BV}
\begin{cases}
(\partial_t+ V \cdot \nabla) \eta =   B  \\
(\partial_t + V \cdot \nabla) \psi = -g\eta+\frac{1}{2}(B^2+|V|^2).
\end{cases}
\end{equation}

On the other hand, recalling that the original Euler equations \eqref{Euler} were written  in terms of the full velocity vector field, restricting these
equations to the free boundary and taking 
horizontal and vertical components
we obtain a set of transport equations for $(V,B)$, namely 
\begin{equation}\label{VBeqn}
\begin{cases}
&(\partial_t+V\cdot\nabla)B=a-1    
\\
&(\partial_t+V\cdot\nabla)V=-a\nabla \eta .
\end{cases}
\end{equation}
The reformulation \eqref{WW-BV} 
of the water wave system will be helpful 
in Section~\ref{sec:layered} when we we carry out the expansion of
the nonlinearity and the resonant analysis. The equations \eqref{VBeqn} 
will help us clarify the role of the Taylor coefficient $a$ in the expansions for the full equation, and even more directly in the derivation of the the linearized equation in Subsection~\ref{subsec:linearized}.

\subsection{Linearizing \texorpdfstring{$B$}{} and \texorpdfstring{$V$}{}}
\label{subsec:linearized}

This is a preliminary step in computing the linearized equations, and will also be needed when computing multilinear expansions of the variables we will be working with. $B$ and $V$ are linear in $\psi$, so we only need to compute their differential with respect to $\eta$.
\begin{lemma}\label{l:DBV}
We have
\begin{equation}\label{dVB}
\begin{aligned}
 D_\eta B[h] = & \ -L(hB) + h LB
 \\
 D_\eta V[h] = & \ - M(hB) + h M B.
\end{aligned}    
\end{equation}
\end{lemma}
As a corollary, since $B = L(\eta) \psi$
and $V = M(\eta) \psi$, this also gives the linearization of the operators $L$ and $M$,
\begin{equation}\label{dML}
\begin{aligned}
 D_\eta L[h] = & \ -LhL + h L L
 \\
 D_\eta M[h] = & \ - MhL + h M L.
\end{aligned}    
\end{equation}
\begin{proof}
We denote the linearized variables by
\[
b = D_\eta B[h], \qquad v = D_\eta V[h].
\]
Linearizing \eqref{BV}, we get
\[
\left\{
\begin{aligned}
 v + b \nabla \eta + B \nabla h = & \ 0 \\
 b - v \cdot \nabla \eta - V \cdot \nabla h 
= & \ - G(\eta)(B h) - \nabla (v h).
\end{aligned}
\right.
\]
Eliminating $v$ we get
\begin{align}
b - V \cdot \nabla h + b |\nabla \eta|^2 + B \nabla h \cdot \nabla \eta
= - G(\eta)(B h) - \nabla (V h),
\end{align}
so
\[
\begin{aligned}
b 
= & \   \frac{
- G(\eta)(B h) - \nabla \cdot (V h) - B \nabla h \cdot \nabla \eta + V \cdot \nabla h
}{
1 + |\nabla \eta|^2
}
\\
= & \  \frac{
- G(\eta)(B h) - \nabla \cdot V \, h - B \nabla h \cdot \nabla \eta
}{
1 + |\nabla \eta|^2
}.
\end{aligned}
\]
This can be rewritten  using the operator $L = L(\eta)$ as
\[
\begin{aligned}
b =& \  - L(Bh) + \frac{\nabla \eta \cdot \nabla (Bh) - \nabla \cdot V \, h - B \nabla h \cdot \nabla \eta}{1 + |\nabla \eta|^2}.
\\
= & \  
 - L(Bh) + \frac{h \big( \nabla B \cdot \nabla \eta - \nabla \cdot v \big)}{1 + |\nabla \eta|^2}.
\\
= & \ -L(Bh) + h LB. 
\end{aligned}
\]
\bigskip
Now for $v$:
\begin{align}
v + (\nabla \eta \otimes \nabla \eta)\, v + v \cdot \nabla h \, \nabla \eta + B \nabla h
= G(\eta)(Bh)\nabla \eta + \nabla (v h)\nabla \eta.
\end{align}
Thus
\[
\begin{aligned}
v 
= &\ \frac{1}{
1 + \nabla \eta \otimes \nabla \eta
}(
G(\eta)(Bh)\nabla \eta + \nabla \cdot (V h)\nabla \eta - V \cdot \nabla h \, \nabla \eta - B \nabla h
).
\\
= & \ \frac{1}{
1 + \nabla \eta \otimes \nabla \eta
}
(G(\eta)(Bh)\nabla \eta + (\nabla \cdot V)\, h \, \nabla \eta - B \nabla h
).
\end{aligned}
\]
\medskip

Rewriting using $M$,
\[
\begin{aligned}
v 
=& \  - M(Bh) + \frac{1}{1 + \nabla \eta \otimes \nabla \eta}
+ h (\nabla \cdot V \, \nabla \eta + \nabla B)
\\
= & \ - M(Bh) + h \cdot M B.
\end{aligned}
\]
\end{proof}

\subsection{The linearized water wave equations}

We denote the linearized variables associated to $(\eta,\psi)$ by $(h,q)$.
Using the differential of $G$ in \eqref{G-lin} as well as the relations \eqref{BV}
we compute the linearization 
\begin{equation}\label{WW-lin0}
\begin{cases}
\partial_t h = G(\eta)(q - hB) - \nabla_x(hV)  \\
\partial_t q = -g h  - \nabla \psi \cdot \nabla q - \nabla\eta \cdot \nabla h  \dfrac{(G(\eta)\psi + \nabla \eta \cdot \nabla \psi)^2}{(1+|\nabla \eta|^2)^2}
\\
\qquad +
\dfrac{(G(\eta)\psi + \nabla \eta \cdot \nabla \psi)
( G(\eta)(q - hB) - \nabla_x(hV) + \nabla h \nabla \psi + \nabla \eta \nabla q)
}{1+|\nabla \eta|^2}.
\end{cases}
\end{equation}
We simplify the second equation of the system above \eqref{WW-lin0} using the relation 
\[
B = \dfrac{G(\eta)\psi + \nabla \eta \cdot \nabla \psi}{1+|\nabla \eta|^2}
\] to arrive at
\begin{equation}\label{WW-lin1}
\begin{cases}
\partial_t h = \ G(\eta)(q - hB) - \nabla_x(hV)  \\
\partial_t q =\  -g h  - \nabla \psi \cdot \nabla q - B^2 \nabla\eta \cdot \nabla h  +
B( G(\eta)(q - hB) \!-\! \nabla_x(hV) ) + B(\nabla h\nabla \psi + \nabla \eta \nabla q).
\end{cases}
\end{equation}
Further eliminating $\psi$ gives the following more manageable system
\begin{equation}\label{WW1+}
\begin{cases}
\partial_t h = G(\eta)(q - hB) - \nabla_x(hV)  \\
\partial_t q = -g h  - V \cdot \nabla q 
 - B V \cdot \nabla h  +
B( G(\eta)(q - hB) - \nabla_x(hV)  ) ,
\end{cases}
\end{equation}
or  equivalently
\begin{equation}\label{WW1++}
\begin{cases}
\partial_t h = G(\eta)(q - hB) - \nabla_x(hV)  \\
\partial_t q = -g h - V \cdot \nabla B h - V \cdot \nabla (q-Bh) 
  +
B( G(\eta)(q - hB) - \nabla_x(hV)  ). 
\end{cases}
\end{equation}

Now we introduce the diagonal variable 
\begin{equation}
\label{def diag variable}
 r := q-Bh , 
\end{equation}
which is the outcome of diagonalizing the above system which is  hyperbolic, nondegenerate, quasilinear with a double speed.  Under this transformation the equations become
\begin{equation}\label{WW2}
\begin{cases}
\partial_t h = G(\eta)r - \nabla_x(hV) \\
\partial_t r = -g h  - V \cdot \nabla r  - (B_t + V\cdot \nabla B) h  .
\end{cases}
\end{equation}
The diagonal variable defined in \eqref{def diag variable} is closely related to the \emph{good variable of Alazard and Metivier} \cite{AM}, which is described in the next section.

Since 
\[
B_t = -V \nabla \cdot B +a-1
\]
using also the first relation in \eqref{BV}
we arrive at
the final linearized system
\begin{equation}\label{WW-lin}
\begin{cases}
&\partial_t h + \nabla_x(hV) - G(\eta) r = 0 \\
&\partial_t r + V\cdot \nabla r + a h = 0  .
\end{cases}
\end{equation}
Since $G$ is a positive operator, it is now clear that the well-posedness of this system hinges on sign of $a$, which corresponds to the \emph{Taylor sign condition}
\begin{equation}
a > 0.    
\end{equation}
In the context of the small data problem in this paper, $a$ ranges in a
small neighborhood of $1$. Indeed, by the formula above,
$a-1=-L(\eta)\eta+O(\text{quadratic})$, so the leading term in $a-1$ is essentially controlled controlled by one
derivative of by $\|\nabla \eta\|_{L^\infty}$; similarly the quadratic term is roughly at the level of  $\|D^\frac12 (B,V)\|_{L^\infty}$. Thus $a$ remains pointwise  small as long as
the above pointwise norms remain small. Two consequences follow. First, the Taylor sign condition $a>0$ holds automatically, and need not be imposed as an
assumption. Second, the dispersion relation $\omega(x,\xi) = \sqrt{a \dn(\xi)}$ of the
linearized problem is a small perturbation of $|\xi|^{\frac12}$.

Finally, we remark that while the above choice of good variables is algebraically convenient and perfectly adequate at higher 
regularity, for our low regularity results
it is necessary to switch to a paradifferential choice for $r$, akin to the definition of the good variable $\varphi$
in \eqref{good-var}. This is discussed in Section~\ref{sec:layered-lin}.

\section{The classical paradifferential expansion of the equation}

Our aim here is to provide a brief overview of the 
classical paradifferential expansion of the water wave equations,
roughly following ideas from Lannes~\cite{Lannes}, Alazard-Metivier~\cite{AM} and Alazard-Burq-Zuily~\cite{ABZ}.

\subsection{The good variable}
We begin with some motivation for the need for a good variable. A good starting place
is the energy space, which corresponds to $\eta \in L^2$ and $\nabla \phi \in L^2(\Omega_t)$.
Under reasonable assumption on $\eta$, the latter corresponds, by the trace theorem, to $\psi \in \dot H^\frac12$. 

Based on this, at higher regularity it is natural to seek to work with $\eta \in H^k$,
and $\nabla \phi \in H^k(\Omega_t)$. By the trace theorem applied to the components $(V,B)$ 
of $\nabla \phi$, we then expect $(V,B) \in H^{k-\frac12}$. However, the trace theorem no longer applies directly to $\phi$ in order to conclude that $\psi \in H^{k+\frac12}$, due 
to the fact that the free boundary only has $H^k$ regularity. This is most readily seen in the first equation in \eqref{BV}, where the term $B \nabla \eta$ can only have the  $H^{k-1}$   
regularity inherited from $\eta$. 

To rectify the above issue, one can write the first equation in \eqref{BV} in the form
\[
\nabla \psi - \nabla T_B  \eta = V + T_{\nabla B} \eta + R(B,\nabla \eta),
\]
where the right hand side can be placed into $H^{k-\frac12}$. It then follows that the 
function that can be placed in $H^{k+\frac12}$ is not $\phi$, but instead its substitute 
\begin{equation}\label{good-var-AM}
\varphi := \psi - T_B  \eta,
\end{equation}
which is usually referred to as the \emph{good variable}; it was first derived by Alazard and Metivier.

Given the good variable $\varphi$ of Alazard-Metivier (often referred to as Alinhac's good unknown, as the idea
originates in Alinhac's work on hyperbolic problems
\cite{Alinhac-interaction}; its nonlinear precursor $\nabla\psi-B\nabla\eta=V$ appears
in Lannes \cite{Lannes})
$V$ and $B$ should be thought of as functions of 
$\nabla \eta$ and $\nabla \varphi$, linear in the latter. 
This can be done perturbatively by rewriting the system \eqref{BV}, using also \eqref{G2+}, in the form
\[
\left\{
\begin{aligned}
V + R(B, \nabla \eta) = & \  \nabla \varphi + T_{\nabla B} \eta
\\
B - R(V, \nabla \eta) = & T_{G(\eta)} \varphi + N_{\eta},
\end{aligned} 
\right.
\]
which can be solved perturbatively for $V$ and $B$, provided 
one has good bounds for $N_\eta$.

\medskip

We next consider the Dirichlet to Neumann map $G(\eta)$. While for fixed $\eta$, $G(\eta)$ is a linear operator, as a function of $\eta$ it is fully nonlinear, with its differential given by \eqref{G-lin}.
Given this expression, Bony's paradifferential formalism 
 suggests that we should have a paradifferential expansion for 
 $G(\eta) \psi$ of the form
\begin{equation}\label{G2+}
\begin{aligned}
G(\eta)\psi  = &  T_{G(\eta)} \psi - T_{G(\eta)} (T_B \eta) -  T_V \nabla \eta + N_\eta
\\
 = & \  T_{G(\eta)} \varphi -  T_V \nabla \eta + N_\eta,
\end{aligned}
\end{equation}
where $T_{G(\eta)}$ is the paradifferential operator with symbol $g(\eta)$, $\varphi$ represents exactly the good variable, and the remainder $N_\eta$ is expected to be better, at least in higher regularity.

\subsection{The paradifferential expansion of the water wave system }

Using the form \eqref{WW-BV} of the water wave system, 
and some simple algebraic manipulations, 
the full system can be written in the paradifferential form 
\begin{equation}\label{Z-para}
\begin{cases}
(\partial_t + T_{V} \nabla) \eta - T_{G(\eta)} \varphi &= N_\eta(\eta,\varphi)  
\\
(\partial_t + T_{V} \nabla) \varphi +
T_{a} \eta &= N_\varphi(\eta,\varphi),
\end{cases}
\end{equation}
where the nonlinear expressions $N_\eta$ and $N_{\varphi}$ are at least quadratic, and are roughly as follows:
\begin{itemize}
    \item $N_\eta(\eta,\varphi)$ is a combination of order zero multilinear expressions in $\nabla \eta$ and $\nabla \varphi$, linear in the latter.
    \item $N_\varphi(\eta,\varphi)$ is a combination of order zero multilinear expressions in $\nabla \eta$ and $\nabla \varphi$, which is either independent of the latter, or exactly quadratic in it. 
\end{itemize}

This is exactly the paradifferential expansion of the $(\eta,\varphi)$ formulation of the water wave system in Bony's formalism.

The key point is that the nonlinearities $N_\eta$  and $N_\varphi$ are expected to be better behaved, at sufficient regularity, and can be treated perturbatively. This observation lies at the heart of all prior
low regularity work on the water wave problem \cite{ABZ,ABZ1,Ai1,Ai2},
where establishing it constituted a substantial part, if not the bulk, of
the analysis.

By contrast, for both the local and global well-posedness results in our paper, it no longer suffices to treat these nonlinearities perturbatively.  Instead, we will need to develop higher order expansions for the
nonlinearities, and to uncover the deeper
structure of the water wave system, which will then be exploited in order
to obtain our main result. This is the aim of the next section.

\section{Enhanced good variables and layered paradifferential expansions}
\label{sec:layered}

We begin this section with a discussion of the nonlinearities
$N_\eta$ and $N_{\varphi}$ in the $(\eta,\varphi)$ equation in \eqref{Z-para}. As noted earlier, there are components of these nonlinearities which we cannot treat perturbatively. The first objective in this section is to identify and classify these nonperturbative components. Then we provide an explicit heuristic computation for them.
This we can do here, in this section, without getting into the technical part of the paper. It will provide a road-map  to understanding the sections that would rigorously deal with the nonlinearities. The remaining task, namely showing that everything else
is perturbative in an appropriate sense, 
will be carried out in Sections~\ref{s:freq+boot}--\ref{s:u}.

We begin with a brief enumeration
of these ``layers" of the nonlinearity:

\begin{enumerate}[label=(\roman*)]
    \item The quadratic nonlinearities are always a problem for global well-posedness, and, in our case, also for the low regularity local well-posedness theory. We distinguish two distinct types:
    \begin{itemize} 
\item $lh$ interactions. A-priori these should be part of the paradifferential flow, so we would find them inside $N_{\eta}$ 
and $N_\varphi$. However, it is much more convenient to work with principal symbols
in the paradifferential flow, which leaves out lower order paradifferential terms,
which we will denote by $N^{[1]}_{\eta}$ and $N^{[1]}_\varphi$.

\item $hh$ interactions, which correspond to 
two high frequencies interacting to produce a lower frequency output.   We will denote their 
contributions by $N^{[2]}_{\eta}$ and $N^{[2]}_\varphi$. 

    \end{itemize}

    \item The cubic nonlinearities will generally fit inside the perturbative box, unless they are balanced. The balanced case 
    contains potentially parallel resonant interactions, which are 
    unfavourable exactly in low dimension $d = 3$, and to a much lesser extent in $d=4$.
    We will denote their contributions by $N^{[3]}_{\eta}$ and $N^{[3]}_\varphi$.

    \item The quartic nonlinearities will also generally fit inside the perturbative box, unless they are balanced.
    In the balanced case, they have to be either transversal or nonresonant; only the latter case is nonperturbative, 
    though it can be treated using normal form methods.
    We will denote their contributions by $N^{[4]}_{\eta}$ and $N^{[4]}_\varphi$. These are needed only in dimension $d = 3$. Furthermore, any lower frequency paracoefficient tilts the balance and moves the corresponding terms into the perturbative box.  
\end{enumerate}

There are two important provisos that apply to the above discussion:

\begin{description}
    \item[a) Lower order paracoefficients] While the above discussion starts from a naive multilinear expansion of the nonlinearity,
in each of the above layers we need to also
allow for lower frequency paracoefficients, 
where ``lower" generally means below the output frequency. There are two mild exceptions here. In $N^{[1]}$ terms, ``lower"
needs to be interpreted as low   relative to the  second highest frequency, which is in turn below the 
highest frequency, namely the output frequency. On the other hand in $N^{[4]}$ (nonconstant) lower frequency paracoefficients yield perturbative contributions, so they can be omitted.

\item[b) Principal part] In each of the the $N^{[j]}$ terms, with $j=\overline{1,4}$, only the principal part is significant, with all lower order terms being perturbative. There is also a mild exception here in the case of $N^{[1]}$, where terms are directly perturbative if either\footnote{So in particular terms where 
derivatives are favourably rebalanced only between the two highest frequencies  are still nonperturbative.}
the second and third highest frequencies are comparable or if we have additional derivatives on factors below the highest frequency.
\end{description}

Given the above discussion, we will seek multilinear expansions
for the nonlinear terms in terms of $\eta$ and $\varphi$, of the form
\begin{equation}\label{N-expansion}
\left\{
\begin{aligned}
N_\eta = & \ N_\eta^{[1]} + N_\eta^{[2]} + N_\eta^{[3]} + N_\eta^{[4]} + N_\eta^{err}
\\
N_\varphi = & \ N_\varphi^{[1]}+  N_\varphi^{[2]} + N_\varphi^{[3]} + N_\varphi^{[4]} +N_{\varphi}^{err},
\end{aligned}
\right.
\end{equation}
where we have added the terms $N_\eta^{[1]}$
and $N_\varphi^{[1]}$ to account for 
unbalanced subprincipal paradifferential terms.

Before we begin to compute the expansions in \eqref{N-expansion}, there is one more difficulty to consider, concerning the coefficients $V$, $B$ and $a$ appearing not only in the paradifferential equation \eqref{Z-para}, but also in the Alazard-Metivier choice of the good variable $\varphi$ in \eqref{good-var-AM}. 
Viewed as functions of $\eta$ and $\varphi$, these also admit expansions as above; the problem is that the quadratic and cubic terms in these expansions could contribute  to the cubic, respectively quartic terms in the expansion for 
$N_\eta$ and $N_{\varphi}$. This is only important in dimension $d = 3$, and marginally in dimension $d=4$. To address this issue, our strategy will be to refine the choice of both the good variable $\varphi$ and of the paracoefficients in the paradifferential equation. This will ensure 
that there are no cubic or quartic balanced interactions hidden inside the paradifferential flow.

In a nutshell, given a dyadic frequency $\lambda$ we would like 
to define the \emph{enhanced good variable} $\varphi_\lambda$ at frequency $\lambda$
to satisfy 
\begin{equation}
\label{def new good var}
\varphi_\lambda := \psi_\lambda - B_{[\lambda]} \eta_{\lambda},
\end{equation}
where the paracoefficient $B_{[<\lambda]}$ is defined as
\begin{equation}\label{Bl}
  B_{[\lambda]} := P_{<\lambda} B(\eta_{\ll \lambda}, \psi_{\ll \lambda})  .
\end{equation}
Here we emphasize that $B_{[\lambda]}$ has both entries at much lower frequencies than $\lambda$. However $B$ is fully nonlinear so this does not guarantee that $B(\eta_{\ll \lambda}, \psi_{\ll \lambda})$ is also at lower frequency, which is why the additional $P_{<\lambda}$ truncation is needed.
Similarly, we define the paracoefficients $V_{[\lambda]}$, 
$a_{[\lambda]}$ well as the  paradifferential operator $\DN_{[\lambda]}$ with symbol
\begin{equation}
  \dn_{[\lambda]}(x,\xi) := P_{<\lambda} \dn(\eta_{\ll \lambda},\xi) . 
\end{equation}
Note that  $\DN_{[\lambda]}$ is the principal part of the Dirichlet to Neumann operator $G(\eta_{\ll \lambda})$, which plays a major role in all that follows. 

As above, we also define the operators $L_{[\lambda]}$ and $M_{[\lambda]}$. Their principal symbols are denoted by $\mfl_{[\lambda]}$ and $\mfl_{[\lambda]}$, and the associated 
operators with the left quantization are denoted by $\L_{[\lambda]}$ and $\M_{[\lambda]}$.

With these notations, we will write the equations \eqref{WW}
in the paradifferential form 
\begin{equation}\label{Z-para+}
\begin{cases}
(\partial_t + \nabla \cdot V_{[\lambda]}) \eta_\lambda - \DN_{[\lambda]} \varphi_\lambda &= N_{\eta,\lambda}(\eta,\varphi)  
\\
(\partial_t + V_{[\lambda]} \cdot \nabla) \varphi_\lambda +
a_{[\lambda]} \eta_\lambda &= N_{\varphi,\lambda}(\eta,\varphi),
\end{cases}
\end{equation}
where the nonlinear expressions $N_{\eta,\lambda}$ and $N_{\varphi,\lambda}$ are at least quadratic, and are roughly as follows:
\begin{itemize}
    \item $N_\eta(\eta,\varphi)$ is a combination of order zero multilinear expressions in $\nabla \eta$ and $\nabla \varphi$, linear in the latter.
    \item $N_\varphi(\eta,\varphi)$ is a combination of order zero multilinear expressions in $\nabla \eta$ and $\nabla \varphi$, which is either independent of the latter, or exactly quadratic in it.
\end{itemize}
\begin{remark}
An important point concerning the equations \eqref{Z-para+} is that we do not use standard Littlewood-Paley projections to localize $\eta$ and $\varphi$ to frequency $\lambda$. Instead, they should be interpreted as identities when applying Littlewood-Paley projections to both sides of the equation. This is because 
the standard Littlewood-Paley projections do not paracommute\footnote{See section~\ref{sec:para-functions}.} with the material derivative, 
so commutator terms would be leading order errors (which we carefully avoid).

\end{remark}

At a later point we will switch from a system with two real valued unknowns $(\eta,\varphi)$ to a complex scalar formulation of the equations, by setting 
\begin{equation}\label{def-u}
u_\lambda = \eta_\lambda + i \mfq_{[\lambda]}(x,D) \varphi_\lambda,   
\end{equation}
where 
\begin{equation}\label{def-mfq}
\mfq_{[\lambda]}:= (\dn/a)^\frac12_{[\lambda]} .
\end{equation}
This will solve a paradifferential equation of the form\footnote{ This form  is given only for orientation here, and will be refined later.}
\begin{equation}\label{u-para-init}
   (\partial_t + V_{[\lambda]} \cdot \nabla +
i \omega_{[\lambda]}(x,D)) u_\lambda  = N_\lambda(u),
\end{equation}
with the dispersion relation
\begin{equation}
 \omega_{[\lambda]} = (a \dn)^\frac12_{[\lambda]}   ,\end{equation}
which is a pseudodifferential operator of order $1/2$.
Up to less significant commutator terms, the expression for the nonlinearity in \eqref{u-para-init} 
is primarily obtained from the nonlinearities in \eqref{Z-para+}.

Returning to our layered paradifferential expansions,  we want to identify the important components of the nonlinearities, modulo lower order terms which will play a perturbative role in our analysis. 

Fixing a dyadic output frequency $\lambda$, we will seek multilinear expansions for frequency $\lambda$ components of the nonlinear terms in terms of $\eta$ and $\varphi$, of the form
\begin{equation}\label{N-expansion+}
\left\{
\begin{aligned}
N_{\eta,\lambda} = & N_{\eta,\lambda}^{[1]}+ N_{\eta,\lambda}^{[2]} + N_{\eta,\lambda}^{[3]} + N_{\eta,\lambda}^{[4]} + R_{\eta,\lambda}
\\
N_{\varphi,\lambda} = & \ N_{\varphi,\lambda}^{[1]}+ N_{\varphi,\lambda}^{[2]} + N_{\varphi,\lambda}^{[3]} + N_{\varphi,\lambda}^{[4]} +R_{\varphi,\lambda}.
\end{aligned}
\right.
\end{equation}

We will in particular produce an expansion for the Dirichlet-to-Neumann map.
We will in particular need a similar expansion for the DN map,
\begin{equation}\label{DN-expansion}
G(\eta) \psi = (G(\eta)\psi)^{[1]} + (G(\eta)\psi)^{[2]} + (G(\eta)\psi)^{[3]} + \cdots, 
\end{equation}
where we have also included a linear paradifferential term, which in the first equation in  \eqref{Z-para+} appears on the left,
\begin{equation}
(G(\eta)\psi)^{[1]}  = - T_{V^l}\cdot \nabla \eta - T_{G(\eta^l)} \varphi .
\end{equation}
Here the superscript $l$, from "low", indicates the same localization as in \eqref{Bl}, but without a reference 
frequency $\lambda$. A more precise definition of this nonlinear quantization choice will be given later in Subsection~\ref{subsec:principal}.

Similar expansions will be needed 
for $B$, $V$ and $a$, of the form 
\begin{equation}\label{BV-expansion}
\begin{aligned}
B = & \ B^{[1]} + B^{[2]}  + \cdots
\\
V = & \ V^{[1]} + V^{[2]}  + \cdots
\end{aligned}
\end{equation}
\begin{equation}\label{a-expansion}
a = a^{[1]}  + \cdots, \qquad D_t a = (D_t a)^{[1]}  + \cdots
\end{equation}

The aim of the rest of this section is to algebraically compute  the terms in the above expansion, and to investigate their null structure. We set aside for now the fixed time elliptic analysis, needed in order to insure that the classical good variable and our improved one are suitably equivalent; see Section~\ref{sec:elliptic} later on. The final step in the study of these expansions, which is to show that the error terms are perturbative, also will happen
much later, in Sections~\ref{sec:Neta}-\ref{sec:Dtaexpansions}, due to the fact that (i) the error bounds are necessarily space-time bounds and (ii) they need to be proved within our main bootstrap argument.

We continue with the statement of the main result concerning the expansions.
Precisely, we have:

\begin{theorem}\label{t:expansions}
The expansions \eqref{N-expansion+} hold with:
\begin{equation}
\begin{aligned}
N^{[1]}_{\eta,\lambda} = & \ - \sum_{\mu \ll \lambda} [\DN_{[\mu]},\eta_\mu]^{2+} \mathfrak \L_{[<\mu]} \varphi_\lambda ,
\\
N^{[1]}_{\varphi,\lambda} = & \  0,
\end{aligned}
\end{equation}
\begin{equation}
\begin{aligned}
N^{[2]}_{\eta,\lambda} = & \  - \DN_{[\lambda]}(\eta_{\gtrsim \lambda}  \L_{[\lambda]} \varphi_{\gtrsim \lambda})
- \nabla \cdot ( \eta_{\gtrsim\lambda}  \M_{[\lambda]} \varphi_{\gtrsim\lambda}),
\\
N^{[2]}_{\varphi,\lambda} = & \ \frac12 \L_{[\lambda]} \varphi_{\gtrsim \lambda} \DN_{[\lambda]} \varphi_{\gtrsim \lambda} - \frac12 \M_{[\lambda]} \varphi_{\gtrsim \lambda} \cdot \nabla \varphi_{\gtrsim \lambda},
\end{aligned}
\end{equation}
respectively
\begin{equation}
\begin{aligned}
N^{[3]}_{\eta,\lambda} 
= & \ 
\L_{[\lambda]}^{*} (\eta_\lambda \DN_{[\lambda]}(\eta_\lambda \L_{[\lambda]}\varphi_\lambda))
- \frac12 \L_{[\lambda]}^{*}( \eta_\lambda^2 \DN_{[\lambda]}\L_{[\lambda]} \varphi_\lambda) - \frac12 \L_{[\lambda]}^{*} \DN_{[\lambda]} (\eta_\lambda^2 \L_{[\lambda]} \varphi_\lambda)
\\
= &\ \DN_{[\lambda]} (\eta_\lambda \L_{[\lambda]} ( \eta_\lambda \L_{[\lambda]} \varphi_\lambda ))
\! + \! \nabla ( \eta_\lambda \M_{[\lambda]}(\eta_\lambda \L_{[\lambda]} \varphi_\lambda))
 \! -\! \frac12 \DN_{[\lambda]}( \eta_\lambda^2 \L_{[\lambda]} \L_{[\lambda]} \varphi_\lambda) \!-\!\frac12\nabla (\eta_\lambda^2 \M_{[\lambda]} \L_{[\lambda]} \varphi_\lambda),
\\
N^{[3]}_{\varphi,\lambda} = & \ - \L_{[\lambda]}\varphi_\lambda \DN_{[\lambda]}(\eta_\lambda\L_{[\lambda]}\varphi_\lambda) + \L_{[\lambda]}\varphi_\lambda \eta_\lambda \DN_{[\lambda]} \L_{[\lambda]}\varphi_\lambda,
\end{aligned}
\end{equation}
and
\begin{equation}
\begin{aligned}
  N^{[4]}_{\eta,\lambda} = & \ \calL_\lambda(\nabla \eta_\lambda, \nabla \eta_\lambda,\nabla \varphi_\lambda)
  \\
  N^{[4]}_{\eta,\lambda} = & \ \calL_\lambda(\nabla \eta_\lambda,\nabla \eta_\lambda, \nabla \varphi_\lambda,\nabla \varphi_\lambda).
\end{aligned}
\end{equation}
We also have the expansions \eqref{BV-expansion} for $B$ and $V$ with 
\begin{equation}
\begin{aligned}
B^{[1]}_\lambda = & \ \L_{[\lambda]} \varphi_\lambda,
\qquad 
B^{[2]}_\lambda = - \L_{[\lambda]}(\eta_{\gtrsim \lambda} \L_{[\lambda]}\varphi_{\gtrsim \lambda} + \eta_{\gtrsim \lambda} \L_{[\lambda]} \L_{[\lambda]} \varphi_{\gtrsim \lambda} ,
\\
    V^{[1]}_\lambda = & \ \M_{[\lambda]}  \varphi_\lambda,
\qquad 
V^{[2]}_\lambda = - \M_{[\lambda]} (\eta_{\gtrsim \lambda} L_{[\lambda]} \varphi_{\gtrsim \lambda}) + \eta_{\gtrsim \lambda} \M_{[\lambda]} \L_{[\lambda]} \varphi_{\gtrsim \lambda},
\end{aligned}    
\end{equation}
as well as the $a$ and $D_t a$ expansions \eqref{a-expansion} with
\begin{equation}
   a^{[1]}=  - T_a \L_{[\lambda]} \eta_\lambda, \qquad  (T_{D_t} a)^{[1]}= - T_a \L_{[\lambda]} \DN_{[\lambda]} \varphi_\lambda.
\end{equation}
\end{theorem}

The expansions in the above theorem play a major role in this article, and will not be rigurously established until Section~\ref{s:u}.
One key difficulty is that such expansions are only useful if we have 
good perturbative bounds on the remaining errors.
However, since the full analysis is 
extensive, in this section we provide a heuristic derivation of these expansions; this will give us a clear 
objective later on.

\begin{remark} A pure multilinear expansion for the water wave equation
was previously computed in \cite{GMS} up to third order. This agrees with the 
our expressions for $N^{[2]}$ and $N^{[3]}$ above, and
corresponds to taking $\DN = \L = |D|$ and $\M = \nabla$ in the above formulas. 
\end{remark}

In what follows we successively compute the above expansion keeping in mind the following principles:

\begin{itemize}
    \item Only the leading order terms in each interaction mode are important, the lower order terms are all perturbative.
    \item All instances of $G$, $L$ and $M$
    can be replaced with the associated principal parts $\DN$, $\L$, $\M$.
    \item The exact form of the $N^{[4]}$ terms 
is not essential, which is why their computation is omitted. 
The reason for this is that we can exclude as perturbative (i) any lower terms with extra derivatives on lower frequency factors, (ii) any quartic expressions with transversal interactions (i.e. with at least two dyadically separated frequencies) and (iii) any quartic terms satisfying a null condition.
\end{itemize}


\subsection{The layered paradifferential expansion for the \texorpdfstring{$\eta$}{}-equation}

We start with the formal expansion for $G(\eta) \psi$. Along the way, this will also provide us with an expansion for 
$B$ and $V$. We recall the differential of $G$ from \eqref{G-lin},
\begin{equation*}
DG(\eta)[h] \psi = -G(\eta)(hB) - \nabla_x(hV) .
\end{equation*}
Following Bony, the linear paradifferential expansion of $G(\eta)\psi$
should be, at the  leading order
\[
(G(\eta)\psi)^{[1]} = -\nabla (T_{V^l}  \eta) -  G(\eta^l)(B^{l} \eta) + G(\eta^l) \psi = -\nabla(T_{V^l} \eta) + G(\eta^l) \varphi,
\]
where we have carefully restricted  
the coefficients to depend only on the 
lower frequencies in $\eta$ and $\psi$.

Substituting this in the $(V,B)$ system \eqref{BV} written paradifferentially,
we get
\begin{equation}\label{BV1}
\left\{
\begin{aligned}
V^{[1]} +  T_{\nabla \eta} B^{[1]}  \approx & \  \nabla \varphi
\\
B^{[1]} - T_{\nabla \eta} V^{[1]} \approx &  T_{\dn} \varphi,
\end{aligned}
\right.
\end{equation}
which solving linearly and neglecting commutators yields the expressions
\[
\left\{
\begin{aligned}
V^{[1]}  = & \ T_{\frac{1}{1+\nabla \eta \otimes \nabla \eta}} (\nabla \varphi - T_{\nabla \eta}   T_{\dn} \varphi)  
\\
B^{[1]} = & \  T_{\frac{1}{1+|\nabla \eta|^2}} ( T_{\dn} \varphi + T_{\nabla \eta} \nabla \varphi).
\end{aligned}
\right.
\]
We can rewrite this using the principal symbols $\mfl$ and $\mfl$ of the operators $L$ and $M$,
\begin{equation}
V^{[1]} = T_{\mfm}\varphi, \qquad B^{[1]} = T_{\mfl} \varphi .    
\end{equation}
where $T_{\mathfrak{m}}$ and $T_{\mathfrak{l}}$ will be defined precisely later in Subsection~\ref{subsec:principal}.
\begin{remark}
     Here one may also ask to consider the paradifferential 
component with the high frequency on $\eta$. However, this would be at the level of $T_{\nabla^2 \omega} \eta$, and thus it is lower order.
\end{remark}

\bigskip

We next consider the  second term in the expansion for $G(\eta) \psi$, which corresponds to having two balanced high frequency inputs, possibly with lower frequency coefficients. Here we will consider 
 the case when one high frequency comes from $\eta$ and the second one comes from $\psi$.  For this we  use again the linearization of $G$, denoting $h = \eta^h$ and assuming $\varphi$ is also at high frequency, arriving at
\[
\begin{aligned}
(G(\eta)\psi)^{[2]}  = & \ -T_{G(\eta)^l}(hB^{[1]}) - \nabla_x(hV^{[1]})
\\ 
= & \ -T_{G(\eta)^l}(h T_{L(\phi)^l} \varphi) - \nabla_x(h T_{M(\psi)^l}\psi).
\end{aligned}
\]
\begin{remark} Apriori this expansion should 
also contain a term with two high frequencies on $\eta$, but that should be  at the level of $T_{\nabla^2 \omega} \Pi(\eta,\nabla \eta)$ and thus lower order and directly perturbative.
\end{remark}

Now we compute the quadratic terms in $B$ and $V$ by substituting into the $(B,V)$ system, \eqref{BV} and discarding all lower order terms:
\[
\left\{
\begin{aligned}
V^{[2]} + T_{\nabla \eta}B^{[2]}  + B^{[1]} \nabla h  = & \  0
\\
B^{[2]} - V^{[2]} \cdot \nabla \eta -V^{[1]}\cdot \nabla h = & (G(\eta)\psi)^{[2]} ,
\end{aligned}
\right.
\]
where we get
\[
\left\{
\begin{aligned}
&(I +\nabla \eta \otimes \nabla \eta)V^{[2]} = -B^{[1]}\nabla h-\nabla\eta(V^{[1]}\cdot \nabla h)
-\nabla\eta (G(\eta)\psi)^{[2]} \\
&(1+ \vert \nabla \eta \vert ^2  )B^{[2]} =-\nabla \eta B^{[1]}\nabla h+V^{[1]}\cdot \nabla h +(G(\eta)\psi)^{[2]}.
\end{aligned}
\right.
\]
This further gives
\begin{equation}\label{BV2}
\left\{
\begin{aligned}
    V^{[2]} = & \ T_{\frac{1}{I +\nabla \eta \otimes \nabla \eta}}\left( -B^{[1]}\nabla h
    + \nabla \eta G(\eta)(hB^{[1]})
    + h \nabla\eta \nabla_x \cdot  V^{[1]}\right)
\\
B^{[2]} = & \ T_{\frac{1}{1+ \vert \nabla \eta \vert ^2  }}
\left(-\nabla \eta B^{[1]}\nabla h
-G(\eta)(h B^{[1]}) - h \nabla \cdot V^{[1]}\right),
\end{aligned}
\right.
\end{equation}
or using \eqref{LvsM},
\begin{equation}\label{BV2-LM}
\left\{
\begin{aligned}
    V^{[2]} = & \ - T_\mfm (h B^{[1]}) + h T_\mfm B^{[1]} 
\\
B^{[2]} = & \ - T_\mfl(h B^{[1]}) + h T_\mfl B^{[1]},
\end{aligned}
\right.
\end{equation}
which is unsurprisingly akin to the formulas
\eqref{dVB} for the linearizations of $B$ and $V$. Substituting $B^{[1]}$ we arrive at 
\begin{equation}\label{BV2-LM+}
\left\{
\begin{aligned}
    V^{[2]} = & \ - T_\mfm (\eta^h T_\mfl \varphi^h) + \eta^h T_\mfm T_\mfl \varphi^h 
\\
B^{[2]} = & \ - T_\mfl(\eta^h T_\mfl \varphi^h) + h T_\mfl T_\mfl \varphi^h.
\end{aligned}
\right.
\end{equation}

\bigskip

Next we consider the balanced cubic term in $G(\eta) \psi$, which involves two $\eta$ and one $\varphi$ factor. This requires the second differential of $G(\eta)$ with respect to $\eta$. Hence we differentiate again with respect to $\eta$,
\begin{equation}
\label{G-2lin}
\begin{aligned}
D^2G(\eta)[h_1,h_2] \psi = &\ -DG(\eta)[h_2](h_1 B) 
- G(\eta) ( h_1 DB[h_2])
- \nabla_x(h_1 DV[h_2]) 
\\ = & \ 
 G(\eta)(h_2 B_1)+ \nabla_x (h_2 V_1)  
- G(\eta) ( h_1 DB[h_2])
- \nabla_x(h_1 DV[h_2]), 
\end{aligned}
\end{equation}
where $B_1$ and $V_1$ are the $B$ and $V$ associated to $h_1 B$ instead of $\psi$,
i.e. 
\[
B_1 = L(h_1 B), \qquad  V_1 = M(h_1 B)
\]

Now we take the cubic part, thinking of $h_1$,
$h_2$ and $\varphi$ as the high frequencies, to arrive at
\[
2(G(\eta)\psi)^{[3]} = G(\eta)(h_2 B^{[1]}( \varphi_1))+ \nabla_x (h_2 V^{[1]}(\varphi_1)) 
- G(\eta) ( h_1 B^{[2]}(h_2,\varphi))
- \nabla_x(h_1 V^{[2]}(h_2,\varphi))
\]
where $\varphi_1$ is associated to $h_1 B^{[1]}$.
Here we have 
\[
\left\{
\ \begin{aligned}
V^{[1]}(\varphi_1) = & \ T_\mfm ( h_1 B^{[1]}), 
\\
\  B^{[1]}( \varphi_1)= & \ T_\mfl ( h_1 B^{[1]})
\end{aligned}
\right.
\]
We  substitute and observe the symmetry in $h_1$ and $h_2$,
\[
\begin{aligned}
2(G(\eta)\psi)^{[3]} = & \ 
T_{\dn}(h_2 T_\mfl ( h_1 B^{[1]}))+ \nabla_x (h_2  T_\mfm ( h_1 B^{[1]})) 
- T_{\dn} ( h_1 (- T_\mfl(h_2 B^{[1]}) + h_2 T_\mfl B^{[1]}))
\\ &\ - \nabla_x(h_1( - T_\mfm (h_2 B^{[1]}) + h_2 T_\mfm B^{[1]}))
\\ 
= & \ T_{\dn}(h_2 T_\mfl ( h_1 B^{[1]}))+T_{\dn}(h_1 T_\mfl ( h_2 B^{[1]})) + \nabla_x (h_2  T_\mfm ( h_1 B^{[1]})) + \nabla_x (h_1  T_\mfm ( h_2 B^{[1]}))
\\ & \ - T_\dn(h_1h_2 T_\mfl B^{[1]}) - \nabla_x(h_1h_2 T_\mfm B^{[1]}).
\end{aligned}
\]
Substituting also the expression for $B^{[1]}$ this gives at leading order
\begin{equation}
 N_\eta^{[3]}=   T_{\dn}(h T_\mfl ( h T_\mfl \varphi ))
 + \nabla_x( h T_\mfm (hT_\mfl \varphi))
 -\frac12 T_\dn( h^2 T_\mfl T_\mfl \varphi) -\frac12\nabla_x (h^2 T_\mfm T_\mfl \varphi)),
\end{equation}
where $h = \eta^h$ and $\varphi = \varphi^h$
are the three high frequencies, which we think of as comparable.
For later use it is also convenient to rewrite this expression in the equivalent form  
\[
N_\eta^{[3]}= T_{\mfl^{*}} [h T_\dn(h T_\mfl \varphi) 
- \frac12 h^2 T_\dn T_\mfl \varphi - \frac12 T_\dn (h^2 T_\mfl \varphi)]
\]
which is identical modulo commutator terms. 
The verification is purely algebraic, using the principal symbols $\dn,\mfl$ and $\mfm$,
and is left for the reader.


\subsection{The layered paradifferential expansion of the \texorpdfstring{$\varphi$}{}-equation}

We start as before with the $\psi$ equation
\[
\partial_t \psi = -\eta+\frac{1}{2}B^2-\frac{1}{2}|V|^2 -BV\cdot\nabla \eta
\]
On the other hand 
using the  $B$ equation
 \[
 (\partial_t + V \nabla)B = a-1.
 \]
we get
\[
\begin{aligned}
\partial_t T_{B^l} \eta  \approx & \ T_{\partial_t  B^l}
\eta + T_{B^{l}} \partial_t  \eta
\\
\approx  & \ T_{a^l-1 
- V^l \nabla B^l} \eta 
+
T_{B^l}  G(\eta) \psi 
\\
\approx  & \ T_{a^l-1} \eta 
- T_{V^l \nabla B^l} \eta
+ T_{B^l} (B-V \nabla \eta).
\end{aligned}
\]
So we arrive at 
\[
\begin{aligned}
(\partial_t + T_{V^l} \cdot \nabla)\varphi +T_{a^l} \eta
= & \ \frac12 B^2  -  \frac12 V^2 - BV \nabla \eta 
+ T_{V^l \nabla B^l} \eta
\\ & \ - T_{B^l} (B-V \nabla \eta)+ T_{V^l}( V+B\nabla \eta)- T_{V^l} \nabla T_{B^l} \eta 
\\
= &\ ( \frac12 B^2  - T_{B^l} B)
- ( \frac12 V^2  - T_{V^l} V)
\\
& \ - \left[ BV\nabla \eta - T_{B^l}(V \nabla \eta)
- T_{V^l}(B \nabla \eta)+ T_{V^l} T_{B^l} \nabla \eta\right].
\end{aligned}
\]

Here we first extract the quadratic part,
\begin{equation}
 \begin{aligned}   
N_\varphi^{[2]} = \frac12 B^{[1]}  B^{[1]}
- \frac12 
V^{[1]}  V^{[1]} - T_{\nabla\eta} V^{[1]} B^{[1]},
\end{aligned}
\end{equation}
which we can rewrite as 
\begin{equation}
 \begin{aligned}   
N_\varphi^{[2]} = \frac12 B^{[1]}  T_\dn \varphi
- \frac12 
V^{[1]} \nabla \varphi ,
\end{aligned}
\end{equation}
where substituting the expressions for $B^{[1]}$ and $V^{[1]}$ we arrive at the desired expression in the proposition.

Finally taking the cubic part we arrive at
\[
\begin{aligned}
N_\varphi^{[3]} = & \  B^{[1]} B^{[2]} - V^{[1]} V^{[2]} - B^{[1]} V^{[1]} \nabla h 
- \nabla \eta V^{[1]} B^{[2]} 
- B^{[1]} \nabla \eta V^{[2]}, 
\end{aligned}
\]
or using the system for $B^{[2]}$, $V^{[2]}$,
\[
\begin{aligned}
N_\varphi^{[3]} = & \  B^{[1]} B^{[2]} 
- B^{[1]} \nabla \eta V^{[2]} 
\\ = & \ B^{[1]}( V^{[1]} \nabla h  -T_\dn(h T_\mfl \varphi) - \nabla_x(h T_\mfm\varphi))
\\ = & \ B^{[1]}( hT_\dn B^{[1]}  -T_\dn(h T_\mfl \varphi)).
\end{aligned}
\]
This finally gives
\begin{equation}\label{Nphi3}
 N_\varphi^{[3]} =  - T_\mfl\varphi \, T_\dn(hT_\mfl\varphi) + T_\mfl\varphi \,  h \, T_\dn T_\mfl \varphi ,
\end{equation}
where again $h = \eta^h$ and $\varphi = \varphi^h$
are the three balanced high frequencies.

\subsection{The paradifferential \texorpdfstring{$G(\eta)$}{}
expansion}

In our analysis later on it is very convenient to replace the paradifferential form of the operators $G(\eta)$ and then also 
$L(\eta)$ and $M(\eta)$ with their principal symbols. This substitution is free in the $N^{[2]}$, $N^{[3]}$ and $N^{[4]}$ terms as it only produces lower order contributions which can be treated perturbatively. 
However, this is no longer the case in the leading paradifferential term in the $\eta$ 
equation. There, substituting $G(\eta^l)$ by $T_{\dn}$ yields an error which contains 
nonperturbative, possibly quadratic terms.
So the challenge is to identify the leading order, nonperturbative part. Heuristically,
the difference is a pseudodifferential operator of order zero, with paracoefficients depending nonlinearly on the low frequency part of $\eta$. Thinking of these paracoefficients in a multilinear fashion, we can further split them in terms of the into a \emph{balanced} and an \emph{unbalanced} part, depending on whether the two highest $\eta$ frequencies in the paracoefficient are 
comparable or not. The balanced part 
in not explicit but will be estimated   perturbatively later on, so we still need to compute the 
unbalanced part. At the quadratic level we can express the difference as
\[
(G(\eta^l)- T_{\dn})\varphi =\calL_{lh}(\nabla^2 \eta,\varphi),
\]
possibly with additional lower frequency paracoefficients at the level of $\nabla \eta$. These paracoefficients must be at frequency below 
the $\nabla^2 \eta$ entry, lest this term is balanced and thus perturbative.

Using the linearization of $G$, the unbalanced terms in $G(\eta_{<\lambda}) \varphi_\lambda$
can be expanded as\footnote{Here one might prefer 
to use a continuous Littlewood-Paley expansion instead.}
\[
 G(\eta_{\ll\lambda}) \varphi_\lambda = - \sum_{\mu \ll \lambda} G(\eta_{<\mu})(\eta_\mu B) + \nabla( \eta_\mu V) + balanced
 \]
with
\[
B = L(\eta_{<\mu}) \varphi_\lambda, \qquad V= M(\eta_{<\mu}) \varphi_\lambda .
\]
Since $GB+\nabla V = 0$,  the above expression can be rewritten as 
\[
G(\eta_{\ll \lambda}) \varphi_\lambda  = - \sum_{\mu \ll \lambda} [G(\eta_{<\mu}),\eta_\mu] L(\eta_{<\mu}) \varphi_\lambda + \nabla \eta_\mu M(\eta_{<\mu}) \varphi_\lambda + balanced.
\]
Here we still have the full $G,L,M$ on the right. However, the lower order terms 
in $G,L,M$ have derivatives on lower frequencies, so their contributions are balanced. Thus we can substitute $G$ with $\DN$ above, also $L$ and $M$ with their principal parts 
\[
 G(\eta_{\ll \lambda}) \varphi_\lambda = - \sum_{\mu \ll \lambda} [\DN(\eta_{<\mu}),\eta_\mu] \L(\eta_{<\mu}) \varphi_\lambda + \nabla \eta_\mu \M(\eta_{<\mu}) \varphi_\lambda + balanced,
\]
where everything is now explicit. The leading order term in the commutator (i.e. the Poisson bracket)
and the second term combined give exactly the corresponding term in the similar expansion of $\DN$. 
So we arrive at 
\[
(G(\eta_{\ll \lambda}) - T_{\dn}) \varphi_\lambda = - \sum_{\mu \ll \lambda} [\DN(\eta_{\ll \mu}),\eta_\mu]^{2+} \mathfrak L(\eta_{<\mu}) \varphi_\lambda + balanced,
\]
where the $2+$ superscript indicates that we select only the subprincipal terms in the commutator (i.e. which have at least two derivatives on $\eta_\mu$),
\begin{equation}
    [\DN(\eta_{\ll \mu}),\eta_\mu]^{2+}
    = [\DN(\eta_{\ll\mu}),\eta_\mu]
    - Op(i \{\dn(\eta_{\ll \mu}),\eta_\mu\}).
\end{equation}


\subsection{The expansion for \texorpdfstring{$a$}{}}
 The starting point in this computation is provided by the formula \eqref{a-formula} which we restate 
below,
\[
\begin{aligned}
a  = & \ 1+ B L(\eta) B + V \cdot L(\eta) V - L(\eta)
 [\frac12(|B|^2+|V|^2) + \eta]
\\
:= & \ 1+ a_1+a_2 ,
\end{aligned}
\]
where $a_1$ depends only on $\eta$,
\[
a_1 = - L(\eta) \eta
\]
and $a_2$ is quadratic in $\omega$,
\[
a_2 = - \frac12 L(\eta)(B^2 + V^2) + B L(\eta)B - V L
(\eta) V .
\]
Here $a_1$ can be thought of
a multilinear expression in $\nabla \eta$, and $a_2$
as a multilinear expression 
in $\nabla \eta$, $(B,V)$ 
and $\nabla (B,V)$, but with a leading order Leibnitz type cancellation when the 
$\nabla (B,V)$ term has the highest frequency. So when we seek a multilinear expansion for $a$, the main term should 
be the one involving $\nabla \eta$ with a paracoefficient.
We will also need the leading term in the expansion for $D_t a$, which will appear in the derivation of the $u$ equation:

\begin{lemma}\label{l:a-exp}
The leading terms in $a$, respectively $D_t a$ are given by
\begin{equation}\label{a-Dta}
 a^{[1]}  =   - T_L T_a \eta  + \bal, 
\qquad (D_t a)^{[1]} = - T_{a} T_L T_G \varphi  + \bal.
\end{equation}
\end{lemma}

We first identify the leading order term in $a$, namely $a^{[1]}$, as in \eqref{a-Dta}, using the linearization of $L$ in \eqref{dML}.  We begin with the $a_1$ component, where we compute
\[
\begin{aligned}
D a_1[h] = & \  - L(\eta) h - L(\eta) ( h L(\eta)\eta)+ h L(\eta) L(\eta) \eta,
\\
\approx & \  - L(\eta) (h (1+a_1)),
\end{aligned}
\]
where we have neglected the last term, which is lower order when $h$ is high frequency.
Thus, the paradifferential term in $a_1$ is 
\begin{equation}\label{a1[1]}
a_1^{[1]} = T_\mfl(1+ T_{a_1}) \eta^h.   
\end{equation}

Denoting by $\sigma$ and $\gamma$ the linearizations
of $\psi$ and $\varphi$, and $B^0 = L(\eta) \eta$, we similarly have
\[
\begin{aligned}
Da_2[h,\gamma] =  & \ \frac12 \left(L(\eta)(L(\eta)(B^2+V^2) h) - h  L(\eta)(L(\eta)(B^2+V^2)
\right)
\\
& - B (L(\eta)(h L(\eta)B) - h L(\eta)L(\eta)B)  - V (L(\eta)(h L(\eta)V) - h L(\eta)L(\eta)V) 
\\ & - L(\eta) (B DB[h,\gamma]+V DV[h,\gamma]) + B L(\eta) (DB[h,\gamma]+V DV[h,\gamma]) .
\end{aligned}
\]
The terms having undifferentiated $h$ as a coefficient are lower order.
Also the last line has a commutator structure,
\[
- [L(\eta),B] DB[h,\gamma] - [L(\eta),V] DV[h,\gamma],
\]
which is at the level of $DB[h,\gamma]$
and thus lower order. 
Placing these terms in the perturbative box, at leading order we are left with 
\begin{equation}\label{a2[1]}
a_2^{[1]} = - T_L T_{a_2} \eta^h.
\end{equation}
Combining the two relations \eqref{a1[1]}
and \eqref{a2[1]} we obtain the first relation of \eqref{a-Dta}.

\bigskip

It remains to justify the $D_t a$ expansion,
where it is necessary to compute the commutator 
$[D_t,L(\eta)]$. We first commute $\partial_t$,
with $B = L\psi = \partial_z \phi_{|top}$:
\begin{align}
\partial_t \partial_z \phi(t,x,\eta(t,x))
&= \partial_z \partial_t \phi + \partial_t \eta \cdot \partial_z^2 \phi \\
&= L \partial_t \phi + \partial_t \eta \cdot L B \\
&= L\bigl(\partial_t \psi - \partial_t \eta \cdot B\bigr)+ \partial_t \eta \cdot L B.
\end{align}
Hence,
\begin{align}
\partial_t L \psi
&= L\bigl(\partial_t \psi - B \partial_t \eta\bigr)
+ \partial_t \eta \cdot L B.
\end{align}
Similarly with $\partial_t$ replaced by $\partial_x$,
\begin{align}
V \cdot \nabla L \psi
&= V \cdot \left[
L\bigl(\partial_x \psi - B \partial_x \eta\bigr)
+ \partial_x \eta \cdot L B
\right].
\end{align}
This yields the commutator formula
\begin{align}
\left[D_t, L\right]\psi
&= -[V,L]\bigl(\partial_x \psi - L\psi\, \partial_x \eta\bigr)
- L\bigl(L\psi \cdot D_t \eta\bigr)
+ D_t \eta \cdot L^2 \psi.
\end{align}
where the last two terms also exhibit a commutator structure.
We use this to compute $D_t a$. First replacing $\psi$ with $\eta$ we have
\begin{align}
D_t a_1 & = \ - D_t L \eta
\\
&= - L D_t \eta
- [V,L]\bigl(\partial_x \eta (1 - L\eta)\bigr)  + L\bigl(L\eta \cdot D_t \eta\bigr)
- D_t \eta \cdot L^2 \eta 
\\
&= - L (T_{1+a_1} D_t \eta)
- [V,L]\bigl(\partial_x \eta (1 - L\eta)\bigr)  + L\bigl(\Pi(L\eta \cdot D_t \eta\bigr))
\\ & - \Pi(D_t \eta \cdot L^2 \eta) + [ L,T_{\partial_t \eta}] L\eta + T_{L^2\eta} D_t \eta
\\
&\approx \, - L  T_{1+a_1}  D_t \eta,
\end{align}
where all but the first term are balanced.

Similarly, we replace successively $\psi$ by $B$, $V$
and $B^2+V^2$, for instance
 \begin{align}
 D_t L B
&= L D_t B - [V,L]\bigl(\nabla B - L B \nabla \eta\bigr)
- L (L B D_t\eta) + D_t \eta L^2 B
\\
& \approx  L D_t B - L(T_{LB} D_t \eta).
\end{align}
This leads to 
\[
\begin{aligned}
D_t a_2 \approx & \  - \frac12 L D_t(B^2+V^2) + B L D_t B
+ V L D_t V 
\\& \
+ L(T_{\frac12L(B^2+V^2)} \partial_t \eta) - 
B L(T_{LB} \partial_t \eta)- V L(T_{LV} \partial_t \eta)
\\ & \
+ D_t B LB + D_t V LV
\\
\approx & \ 
-L(T_{a_2} \partial_t \eta) ,
\end{aligned}
\]
where the first and last lines are balanced due to a commutator structure, and the second line yields the leading term after another commutation.

\subsection{Double null condition verification} 

In this subsection we present a  heuristic verification of the double null condition satisfied by the water wave system.  This is just for clarity, in a simpler setting; the full computation will happen at the level of energy/flux corrections but it will emulate the computation below. This is most readily done at the level of the $u$ evolution \eqref{u-para-init}, where 
we use the quadratic and cubic terms $N_\lambda^{[2]}(u)$, $N_{\lambda}^{[3]}(u)$ in the expansion of $N_\lambda(u)$, based on the discussion in Section~\ref{subsec:resonant}.
We recall that $u$ and $\eta$ and $\varphi$ are connected by 
\begin{equation}\label{to-u}
2\eta = u+\bu, \qquad 2\varphi= -i T_\mfq (u- \bu), 
\end{equation}
and then 
\[
N_\lambda^{[2]}(u) = N_{\eta,\lambda}^{[2]}(u)+i T_{\mfq} N_{\varphi,\lambda}^{[2]}(u),
\qquad
N_\lambda^{[3]}(u) = N_{\eta,\lambda}^{[3]}(u)+i T_{\mfq}N_{\varphi,\lambda}^{[3]}(u).
\]

At the cubic level, we are only interested in the nonlinearity component which has phase rotations symmetry and balanced entries, else the cubic interaction is either transversal or nonresonant.

As discussed in Section~\ref{subsec:resonant}, the quadratic interactions are by themselves nonresonant. However, correcting them with a normal form type transformation, the quadratic 
term in the nonlinearity yields cubic terms, of which we need to select again those with phase rotation symmetry. 

Based on the above discussion, one would expect to compute two cubic expressions, one from the cubic nonlinearity and the other from the double use of the quadratic nonlinearity via normal forms,
and then show that their sum satisfies
the double null condition. Fortunately, it turns out that the two cubic expressions satisfy the double null condition separately. 
Even better, the cubic part also naturally splits into two components
which can be dealt with separately.
Finally, all the computations in this section are done assuming that paradifferential operators commute; the logic here is that any commutators would yield lower order terms, with additional derivatives on the paracoefficients, which are directly perturbative without any need for a null condition.

\subsubsection{ The \texorpdfstring{$N_\eta^{[3]}$}{} terms}
We expand $N_\eta^{[3]}$ using the relations \eqref{to-u}, and selecting only the terms with phase rotation symmetry, i.e. with two $u$ factors an a $\bu$ factor. We denote the three (close) input frequencies by $\xi^1,\xi^2,\xi^3$ where the middle one corresponds to the $\bu$ factor.
Then the output is at frequency $\xi^4:= \xi^1-\xi^2 + \xi^3$.

We have two scenarios for the two $\eta$ entries, either both yielding $u$ factors or $u$ and $\bu$, which we label  $++$, respectively $+-$. In the first case $\varphi$ will yield the $\bu$ entry and in the second it will yield the $u$ entry.
Using the first expression in Theorem~\ref{t:expansions} for $N_\eta^{[3]}$, 
\[
\begin{aligned}
N^{[3]}_\eta = & T_{\mfl^*} [\eta T_{\dn}(\eta T_{\mfl} \varphi)  -  \frac12 \eta^2 T_{\dn} T_{\mfl} \varphi - \frac12 T_{\dn} (\eta^2 T_{\mfl} \varphi)]
\end{aligned}
 \]
in the first case we obtain the 
symmetrized symbol
\[
16 N_{\eta,\lambda}^{[3],++}
= i \mfq^{-1}(\xi^2) \bar\mfl(\xi^4)
\bar\mfl(\xi^2)(\dn(\xi^1-\xi^2)+\dn(\xi^3-\xi^2) - \dn(\xi^2) - \dn(\xi^4)),
\]
while in the second case we obtain
\[
\begin{aligned}
16 N_{\eta,\lambda}^{[3],+-}
= & \ -i \mfq^{-1}(\xi^3)
\bar\mfl(\xi^4) \mfl(\xi^3) 
(\dn(\xi^1+\xi^3) + \dn(\xi^2-\xi^3)
- \dn(\xi^3)-\dn(\xi^4))
\\ 
\approx & \  -i \mfq^{-1}(\xi^3)
\bar\mfl(\xi^4) \mfl(\xi^3) 
(\dn(\xi^1) + \dn(\xi^2-\xi^3)
-\dn(\xi^4)),
\end{aligned}
\]
where "$\approx$" indicates equality 
modulo  expressions vanishing to second order on the diagonal and symmetrization.

\subsubsection{ The \texorpdfstring{$N_\varphi^{[3]}$}{} terms}
Recalling that $N_\varphi^{[3]}$ has the form
\[
\begin{aligned}
N_\varphi^{[3]} = & - T_{\mfl} \varphi T_{\dn}(\eta T_{\mfl}\varphi) + T_{\mfl} \varphi \eta T_{\dn} T_{\mfl} \varphi  , 
\end{aligned}
\]
we compute the symbols for the same two cases of resonant interactions.
We obtain
\[
\begin{aligned}
16 N^{[3],++}_\varphi = & \ - 2\mfq^{-1}(\xi^1)\mfq^{-1}(\xi^3) \mfl(\xi^1)
\mfl(\xi^3) ( - \dn(\xi^3-\xi^2) +\dn(\xi^1))
\\
\approx & \ - \mfq^{-1}(\xi^1)\mfq^{-1}(\xi^3) \mfl(\xi^1)
\mfl(\xi^3) ( - \dn(\xi^3-\xi^2) -\dn(\xi^1-\xi^2)+\dn(\xi^1)+\dn(\xi^3)),
\end{aligned}
\]
respectively 
\[
\begin{aligned}
16 N^{[3],+-}_\varphi = & \   \mfq^{-1}(\xi^1)\mfq^{-1}(\xi^2)
  \mfl(\xi^1)
\bar \mfl(\xi^2) ( -\dn(\xi^1+\xi^3)- \dn(\xi^3-\xi^2) + \dn(\xi^1)+\dn(\xi^2)) 
\\
= & \  \mfq^{-1}(\xi^1)\mfq^{-1}(\xi^2)
 \mfl(\xi^1)
\bar \mfl(\xi^2) ( -\dn(\xi^3)- \dn(\xi^3-\xi^2) +\dn(\xi^2)) .
\end{aligned}
\]

\subsubsection{The double null condition for $N^{[3]}_u$}
Here we combine the $N^{[3]}_\eta$
and the $N^{[3]}_\varphi$ contributions. It is convenient to keep separately the $++$ and the $+-$
terms. Noting that 
\[
\mfq(\xi^1) \mfq(\xi^3) \approx 
\mfq(\xi^2) \mfq(\xi^4), \qquad \mfq^{-1}(\xi^2) + \mfq^{-1}(\xi^4)
\approx 2 \mfq^{-1}(\frac{\xi^2+\xi^4}2) \approx 2 \mfq^{-1}(\frac{\xi^1+\xi^3}2), 
\]
in the $++$ case we 
obtain
\[
\begin{aligned}
16  N^{[3],++}_u = & \  16 N^{[3],++}_\eta
+ 16 i \mfq(\xi^4) N^{[3],++}_\varphi
\\
\approx & \ i \mfq^{-1}((\xi^2+\xi^4)/2) \left[\ \bar\mfl(\xi^4)
\bar\mfl(\xi^2)(\dn(\xi^1-\xi^2)+\dn(\xi^3-\xi^2) - \dn(\xi^2) - \dn(\xi^4)) \right.
\\ & \ \qquad \qquad \qquad \qquad  \left. -\mfl(\xi^1)
\mfl(\xi^3) ( - \dn(\xi^3-\xi^2) -\dn(\xi^1-\xi^3)+\dn(\xi^1)+\dn(\xi^3))\right]
\\ \approx & \ 0,
\end{aligned}
\]
where at the last step we have 
applied the $(\xi^1,\xi^3) \leftrightarrow (\xi^2,\xi^4)$ symmetry which also requires a complex conjugation.

Next we consider the $+-$ case,
using symmetries to simplify the 
expression,
\[
\begin{aligned}
16  N^{[3],+-}_u = & \   16 N^{[3],+-}_\eta
+ 16 i \mfq(\xi^4) N^{[3],+-}_\varphi
\\
\approx &   -i \mfq^{-1}(\xi^3)
\bar\mfl(\xi^4) \mfl(\xi^3) 
(\dn(\xi^1) + \dn(\xi^2-\xi^3)
-\dn(\xi^4)) 
\\ & \ 
+ i \mfq(\xi^4) \mfq^{-1}(\xi^1)\mfq^{-1}(\xi^2)
 \mfl(\xi^1)
\bar \mfl(\xi^2) ( -\dn(\xi^3)- \dn(\xi^3-\xi^2) +\dn(\xi^2))
\\
\approx & \ -i ( \mfq^{-1}(\xi^1) + \mfq(\xi^4) \mfq^{-1}(\xi^1)\mfq^{-1}(\xi^2))
 \mfl(\xi^1)
\bar \mfl(\xi^2) ( \dn(\xi^3)+ \dn(\xi^3-\xi^2) - \dn(\xi^2))
\\
=  &  \ -i ( \mfq(\xi^2) + \mfq(\xi^4)) 
 \frac{\mfl}{\mfq}(\xi^1)
\frac{\bar \mfl}{\mfq}(\xi^2) ( \dn(\xi^3)+ \dn(\xi^3-\xi^2) - \dn(\xi^2))
\\
\approx & \  \frac{i}2 \Delta^4 \mfq(\xi) \frac{\mfl}{\mfq}(\xi^1)
\frac{\bar \mfl}{\mfq}(\xi^2) (\dn(\xi^3-\xi^2)+\dn(\xi^1) +\dn(\xi^2) - \dn(\xi^3) - \dn(\xi^4))
\\ & \ + (\mfq(\xi^1) + \mfq(\xi^2)+
\mfq(\xi^3) +\mfq(\xi^4))\frac{\mfl}{\mfq}(\xi^1)
\frac{\bar \mfl}{\mfq}(\xi^2)
\Delta^4 \dn(\xi).
\end{aligned}
\]
Now it suffices to observe that 
both $\Delta^4 \mfq$ and $\Delta^4 \dn$ vanish to second order 
on $\Delta^4 \xi = 0$ at $\xi^1= \xi^2= \xi^3 = \xi^4$.

\bigskip

\subsubsection{The double null condition for the iterated \texorpdfstring{$N^{[2]}$}{}
contribution}

Here we consider the contributions arising from the quadratic terms via a normal form transformation. 
Our starting point is given by the quadratic contributions
\begin{equation}
\left\{
\begin{aligned}
N^{[2]}_{\eta} = & \  - G(\eta  L \varphi)
- \nabla \cdot ( \eta M \varphi)
\\
N^{[2]}_{\varphi} = & \ \frac12 L \varphi G \varphi - \frac12 M \varphi \cdot \nabla \varphi,
\end{aligned}
\right.
\end{equation}
which have symbols
\[
\left\{
\begin{aligned}
N^{[2]}_{\eta}(\xi^1,\xi^2) = & \ 
\dn(\xi^1+\xi^2) \mfl(\xi^2)
- i (\xi^1+\xi^2) \cdot \mfm(\xi^2)
\\
N^{[2]}_{\varphi}(\xi^1,\xi^2) = & \ 
\frac12( \mfl(\xi^1) \dn(\xi^2) - i \mfm(\xi^1) \cdot \xi^2 ),
\end{aligned}
\right.
\]
where the latter can also be symmetrized. These symbols satisfy a 
simple null condition, namely they vanish when $\xi^1+\xi^2 = 0$ (i.e. when the output is at zero frequency). A-priori there are two  potential balanced cubic contributions
to the transformed equation, arising from
\medskip

(i) $lh$ quadratic interactions.
These are contained within the paradifferential equation, and will be treated at the level of the energy estimates. Fortunately, since the paradifferential equation is conservative at leading order, any such errors 
(e.g. $\nabla \cdot V_{[<\lambda]}
|u_\lambda|^2$) have an extra derivative on the paracoefficient, which already yields a simple null condition. This is coupled with the fact the the quadratic contributions themselves satisfy the simple null condition.
Heuristically, these two facts yield the desired double null condition.

\medskip

(ii) $hh$ quadratic interactions.
If we denote by $\xi$ the cubic output frequency, here we only need to consider the nondegenerate case when the two $N^{[2]}$ input frequencies are close to $2\xi$ and $-\xi$, and then use the equation bifurcate the frequency $2\xi$ factor in the correction into two frequency $\xi$ terms. In brief, this is the computation we need to do here.

We begin by using the expressions above for $N^{[2]}_{\eta}$ and $N^{[2]}_{\varphi}$ in order to compute the principal symbol for $N^{[2]}_u$, which we separate into 
\[
N^{[2]}_u = N^{[2],++}(u,u) +N^{[2],+-}(u,\bu) +N^{[2],--}(\bu,\bu).
\]
Here a direct computation yields the principal symbols
\[
\begin{aligned}
4n^{[2],++}(\xi,\eta) = & \ i \mfq^{-1}(\eta) (\dn(\xi+\eta) \mfl(\eta) 
+ i(\xi+\eta) \mfm(\eta)) + \frac{i}2 \mfq(\xi+\eta) \mfq^{-1}(\xi) \mfq^{-1}(\eta)(\mfl(\xi) \dn(\eta) - i\mfm(\xi) \cdot \eta) \\
4n^{[2],--}(\xi,\eta) = & \ i \mfq^{-1}(\eta) (\dn(\xi+\eta) \bar\mfl(\eta) 
- i(\xi+\eta) \bar\mfm(\eta)) + \frac{i}2 \mfq(\xi+\eta) \mfq^{-1}(\xi) \mfq^{-1}(\eta)(\bar \mfl(\xi) \dn(\eta) + i\bar \mfm(\xi) \cdot \eta) \\
 4n^{[2],+-}(\xi,\eta)= & \  i \mfq^{-1}(\xi) (\dn(\xi-\eta) \mfl(\xi) + i(\xi-\eta) \mfm(\xi))  - i \mfq^{-1}(\eta) (\dn(\xi-\eta) \bar \mfl(\eta) + i(\xi-\eta) \bar\mfm(\eta)) 
\\
& \ 
+ \frac{i}{2}\mfq(\xi-\eta) \mfq^{-1}(\xi)\mfq^{-1}(\eta)( \mfl(\xi) \dn(\eta) + \bar \mfl(\eta) \dn(\xi) +i \mfm(\xi) \cdot \eta - i \bar \mfm(\eta) \cdot \xi), 
\end{aligned}
\]
where the first two relations should also be symmetrized in $\xi$, $\eta$. Formally each of these three
terms contribute to the corresponding normal form correction $B^{[2]}_u$,
\[
B^{[2]}_u = B^{[2],++}(u,u) +B^{[2],+-}(u,\bu) +B^{[2],--}(\bu,\bu),
\]
where these three terms have principal symbols
\[
b^{[2],ab}(\xi,\eta) = \frac{i \, n^{[2],ab}(\xi,\eta)}{a \omega(\xi)+b\omega(\eta) - \omega(a\xi+b\eta)}.
\]
Applying the time derivative to the 
above normal form corrections, the $N^{[2]}_u$ will yield cubic terms in the transformed equation. Our task is to show 
that these cubic terms satisfy the double null condition.
We are interested in normal form scenarios which yield cubic expressions with phase rotation symmetry, 
but which do not go through zero frequency, where the double null condition is automatic.
This leaves us with the following remaining cubic expressions:
\begin{equation}\label{nf-cubic}
N^{[3],NF,1} = B^{[2],+-}(N^{[2],++}(u,u),\bu), \qquad 
N^{[3],NF,2} =
2 B^{[2],--}(\ol{N^{[2],--}(\bu,\bu)},\bu). 
\end{equation}

Next we simplify the above expressions modulo 
terms that vanish quadratically. We will need such expansions near $\xi = \eta$ and near $\xi = -2\eta$ for $n^{[2],++}$
and $n^{[2],--}$, and near $\xi = 2\eta$ for $n^{[2],+-}$.
We begin with the expansion near $\xi = \eta$.
Starting from \eqref{LvsM}, we obtain the relations
\[
\dn(\xi) \mfl(\xi) + i \xi \cdot \mfm(\xi) = 0, \qquad \dn(\xi) \bar \mfl(-\xi) - i \xi \cdot \bar \mfm(-\xi) = 0.
\]

As a slight extension of these, it follows that
\[
\begin{aligned}
&\dn(\xi) \mfl(\eta) + i \xi \cdot \mfm(\eta) \approx 0 \qquad \text{ near }  \xi \in \R^+ \eta\\
&\dn(\xi) \bar \mfl(\eta) - i \xi \cdot \bar \mfm(\eta) \approx  0  \qquad \text{ near }  \xi \in \R^- \eta.
\end{aligned}
\]
This yields, after symmetrization,
\begin{equation}\label{n2ppmm0}
\left\{ 
\begin{aligned}
 8n^{[2],++}(\xi,\eta) \approx & \   i \mfq(\xi+\eta) \mfq^{-1}(\xi) \mfq^{-1}(\eta)(\mfl(\xi) \dn(\eta) 
 + \mfl(\eta) \dn(\xi))
 \\
  8n^{[2],--}(\xi,\eta)  \approx & \   i\mfq(\xi+\eta) \mfq^{-1}(\xi) \mfq^{-1}(\eta)(\bar \mfl(\xi) \dn(\eta)+ \bar \mfl(\eta) \dn(\xi))
\end{aligned}
\right.
\qquad \text{near } \xi = \eta,
\end{equation}
and similarly
\begin{equation}
\left\{ 
\begin{aligned}
 4n^{[2],++}(\xi,\eta) \approx & \   i \mfq^{-1}(\eta) \dn(\xi+\eta) \mfl(\eta) 
 \\
 4n^{[2],--}(\xi,\eta)  \approx & \  i \mfq^{-1}(\eta) \dn(\xi+\eta) \bar\mfl(\eta) 
\end{aligned}
\right.
\qquad \text{near } \xi = -2\eta,
\end{equation}
with the mixed sign representation as follows
\begin{equation}
\begin{aligned}
4n^{[2],+-}(\xi,\eta)\approx & \ - i \mfq^{-1}(\eta) \dn(\xi-\eta) \bar \mfl(\eta)
\end{aligned}
\qquad \text{near } \xi = 2\eta.
\end{equation}

We now compute the symbols of the expressions in \eqref{nf-cubic}. As usual we denote by $\xi^1$ and $\xi^3$ the frequencies of the $u$ factors, and by $\xi^2$ the frequency\footnote{before conjugation} of the $u$ that gets conjugated.
Finally we set $\xi^4 = \xi^1+\xi^2-\xi^3$ 
to be the output frequency.

For the first cubic form in \eqref{nf-cubic} we compute a succession of equivalent  expressions using various symmetrizations
\[
\begin{aligned}
64 N^{[3],NF,1}(\xi^1,\xi^2,\xi^3) = & \ 
b^{[2],\pm}(\xi^1+ \xi^3,\xi^2) n^{[2],++}(\xi^1,\xi^3)
\\ =& \
\frac{i n^{[2],\pm}(\xi^1+ \xi^3,\xi^2) n^{[2],++}(\xi^1,\xi^3)}{\omega(\xi^1+\xi^3)- \omega(\xi^2) - \omega(\xi^4)}
\\
\approx & \ \frac{i \mfq^{-1}(\xi^2) \dn(\xi^4) \bar \mfl(\xi^2) \mfq(\xi^1+\xi^3) \mfq^{-1}(\xi^1) \mfq^{-1}(\xi^3)(\mfl(\xi^1) \dn(\xi^3) 
 + \mfl(\xi^3) \dn(\xi^1))}{\omega(\xi^1+\xi^3)- \omega(\xi^2) - \omega(\xi^4)}
 \\
 = & \ \frac{\mfq(\xi^1+\xi^3) \prod\mfq^{-1}(\xi^j) }{\omega(\xi^1+\xi^3)- \omega(\xi^2) - \omega(\xi^4)} i\,\mfq(\xi^2) \dn(\xi^4) \bar \mfl(\xi^2) (\mfl(\xi^1) \dn(\xi^3) 
 + \mfl(\xi^3) \dn(\xi^1)).
\end{aligned}
\]
Here the fraction is symmetric to second order, so we examine the remaining expression, symmetrizing first in $(\xi^2,\xi^4)$,
\[
\begin{aligned}
2 \mfq(\xi^2) \dn(\xi^4) \bar \mfl(\xi^2)
\approx & \ 
(\mfq(\xi^2)\!+\!\mfq(\xi^4)) (\dn(\xi^4) \bar \mfl(\xi^2)\!+\!\dn(\xi^2) \bar \mfl(\xi^4))
+(\mfq(\xi^2)\!-\!\mfq(\xi^4)) (\dn(\xi^4) \bar \mfl(\xi^2)\!-\!\dn(\xi^2) \bar \mfl(\xi^4))
\\
\approx & \  (\mfq(\xi^2)+\mfq(\xi^4)) (\dn(\xi^4) \bar \mfl(\xi^2)+\dn(\xi^2) \bar \mfl(\xi^4)).
\end{aligned}
\]
Inserting this into the previous expression 
we obtain the fully symmetric factor
\[
 (\dn(\xi^4) \bar \mfl(\xi^2)+\dn(\xi^2) \bar \mfl(\xi^4))(\mfl(\xi^1) \dn(\xi^3) 
 + \mfl(\xi^3) \dn(\xi^1)),
\]
and the remaining factor $i (\mfq(\xi^2)+\mfq(\xi^4))$. Symmetrizing in $(\xi^1,\xi^3)$ and $(\xi^2,\xi^4)$,
\[
2i (\mfq(\xi^2)+\mfq(\xi^4)) \approx -i \Delta^4 \xi,
\]
which is fully symmetric.

\bigskip

Similarly, we have 
\[
\begin{aligned}
64 n^{[3],NF,2}(\xi^1,\xi^2,\xi^3) = & \
2 b^{[2],--}(-\xi^1,-\xi^3,\xi^2)
\ol{n^{[2],--}(\xi^1,\xi^3)}
\\
= & \ \frac{2i\, n^{[2],--}(-\xi^1-\xi^3,\xi^2)\ol{n^{[2],--}(\xi^1,\xi^3)}}{
-\omega(\xi^1+\xi^3) - \omega(\xi^2)-\omega(\xi^4)}
\\
\approx & \ \frac{2i\, \mfq^{-1}(\xi^2) \dn(\xi^4) \bar\mfl(\xi^2)\mfq(\xi^1+\xi^3) \mfq^{-1}(\xi^1) \mfq^{-1}(\xi^3)( \mfl(\xi^1) \dn(\xi^3)+  \mfl(\xi^1) \dn(\xi^3))}{
-\omega(\xi^1+\xi^3) - \omega(\xi^2)-\omega(\xi^4)}.
\end{aligned}
\]
From this point on, the symmetrizations proceed exactly as in the previous case.

\subsection{The layered paradifferential expansion for the linearized equation}
\label{sec:layered-lin}
Our starting point is the linearized equation \eqref{WW1++},
which we recall here in a slightly modified form:
\begin{equation}\label{WW-lin2}
\begin{cases}
\ \partial_t h = G(\eta)(q - hB) - \nabla_x(hV)  \\
\ \partial_t q = -g h  - V \cdot \nabla  q + B V \cdot \nabla h  
  +
B( G(\eta)(q - hB) - \nabla_x(hV)  ). 
\end{cases}
\end{equation}
This equation is not properly diagonalized, and in Section~\ref{sec:wwvar} we saw that 
algebraically it is very natural to introduce a good variable $r = q-Bh$.
This would suffice at high regularity. However, in our low regularity setting
it is much better to make a refined choice of the good variable, namely 
\begin{equation}
 r = q - T_{B^l} h,   
\end{equation}
which is akin to the good variable we use for the full equation.

Our task now is to compute the system for $(h,r)$ and then its 
layered paradifferential expansion.

We seek to rewrite this in the frequency localized form  
\begin{equation}\label{Z-para-lin}
\begin{cases}
&(\partial_t + \nabla \cdot V_{[\lambda]}) h_\lambda - \DN_{[\lambda]} r_\lambda = N_{\eta,\lambda}^{lin}(h,r)  
\\
&(\partial_t + V_{[\lambda]} \cdot \nabla) r_\lambda +
a_{[\lambda]} h_\lambda = N^{lin}_{\varphi,\lambda}(h,r),
\end{cases}
\end{equation}
where for the nonlinear terms on the right 
we want an expansion of the form
\begin{equation}\label{N-expansion-lin}
\begin{aligned}
N_{\eta,\lambda}^{lin} = & 
 N_{\eta,\lambda}^{lin,[0]} + N_{\eta,\lambda}^{lin,[1]}+ N_{\eta,\lambda}^{lin,[2]} + R^{lin}_{\eta,\lambda}
\\
N_{\varphi,\lambda}^{lin} = & \  N_{\varphi,\lambda}^{lin,[0]} +N_{\varphi,\lambda}^{lin,[1]}+ N_{\varphi,\lambda}^{lin,[2]} +R^{lin}_{\varphi,\lambda}.
\end{aligned}
\end{equation}
This is comparable to the expansion \eqref{N-expansion+} with the following differences:
\begin{itemize}
    \item The $N^{[3]}$ and $N^{[4]}$ terms are omitted, and instead go into the perturbative box. The key difference occurs in the $N^{[3]}$ term, where symmetry is lost in the linearization, which in turn destroys the double null condition appearing in the study of the full equation after full symmetrization.  The price to pay for this is the loss of $\frac1{12}$ derivatives in the 3D local well-posedness result.

    \item The new term $N^{[0]}$ contains the leading, quadratic paradifferential contributions where  
    the $(h,r)$ entries are at low frequency. This also appears due to the loss of symmetry in the linearisation, but it will turn out to be mostly harmless.
\end{itemize}

\begin{proposition}
The expansions \eqref{N-expansion-lin} hold with:
\begin{equation}
\begin{aligned}
& N^{lin,[0]}_{\eta,\lambda} =   - \nabla h_{\ll \lambda} \M_{[\lambda]} \varphi_{\lambda}   - 
[G_{[\lambda]},h_{\ll\lambda}] \L_{[\lambda]} \varphi_{\lambda}  - \nabla(\eta_{ \lambda} \M_{[\lambda]} r_{\ll \lambda}) - G_{[\lambda]} ( \eta_{ \lambda} \L_{[\lambda]} r_{\ll \lambda}),\\
& N^{lin,[0]}_{\varphi,\lambda} =  \ 
\M_{[\lambda]} \varphi_{ \lambda}  G_{[\lambda]} r_{\ll \lambda} - \L_{[\lambda]} \varphi_{ \lambda}  \nabla  r_{\ll \lambda},\\
 &N^{lin,[1]}_{\eta,\lambda} =  \ - \sum_{\mu \ll \lambda} [\DN_{[\mu]},\eta_\mu]^{2+} \mathfrak \L_{[<\mu]} r_\lambda , \\
& N^{[1]}_{\varphi,\lambda} =  \  0,\\
& N^{[2]}_{\eta,\lambda} =  \  - \DN_{[\lambda]}(h_{\gtrsim \lambda}  \L_{[\lambda]} \varphi_{\gtrsim \lambda})
- \nabla \cdot ( h_{\gtrsim\lambda}  \M_{[\lambda]} \varphi_{\gtrsim\lambda})
 - \DN_{[\lambda]}(\eta_{\gtrsim \lambda}  \L_{[\lambda]} r_{\gtrsim \lambda})
- \nabla \cdot ( \eta_{\gtrsim\lambda}  \M_{[\lambda]} r_{\gtrsim\lambda}), \\
& N^{[2]}_{\varphi,\lambda} =  \ \M_{[\lambda]} \varphi_{\gtrsim \lambda}  G_{[\lambda]} r_{\gtrsim \lambda} - \L_{[\lambda]} \varphi_{\gtrsim \lambda}  \nabla  r_{\gtrsim \lambda}.
\end{aligned}
\end{equation}
\end{proposition}

\begin{proof}
To begin with, we start from the system \eqref{WW-lin2} and derive a corresponding system for the good variables $(h,r)$. This is easily done 
in the first equation, 
which takes the form 
\begin{equation*}
  \partial_t h = -\nabla_x(hV)  
  + G(\eta)r  - G(\eta)((B -T_{B^l}) h)  .
\end{equation*}
The computation for the second equation 
is slightly more involved,
\[
\begin{aligned}
\partial_t r = & \  \partial_t q - T_{B^l} \partial_t h - T_{\partial_t B^l} h 
\\
= & \ -g h - V \cdot \nabla (r + T_{B^l} h) -BV \cdot \nabla h +(B-T_{B^l})( G(\eta)(r + T_{B^l} h - hB) - \nabla_x(hV)  )\\
&\ -T_{\partial_t B^l} h
\\
= & \ - V \cdot \nabla r - 
(g  + V \cdot T_{\nabla B^l} + T_{\partial_tB^l}) h
 - V \cdot ( T_{B^l}-B)\nabla h
 \\ 
 & \ +(B-T_{B^l})( G(\eta)(r + T_{B^l} h - hB) - \nabla_x(hV)  ).
 \\
\end{aligned}
\]
Reorganizing the two equations and using the $R$ notation, we arrive at
\begin{equation*}
\begin{cases}
& \partial_t h + \nabla (Vh)= G(\eta) r - G(\eta) R(B,h)    \\
& \partial_t r  + V \nabla r =-(g  + V \cdot T_{\nabla B^l} + T_{\partial_tB^l}) h
 - V \cdot R(B,\nabla h) + R(B, G(\eta)r) \\
 &\qquad \qquad  \qquad- R(B, G(\eta) R(B,h) - \nabla_x(hV)) .
\end{cases}
\end{equation*}

We begin with the first equation, where we have
\[
N_{\eta}^{lin}(h,r) = - \nabla R(V,h) - G(\eta) R(B,h) + R(G(\eta),r) + T_{G-\DN} r.
\]
Examining the terms on the right, the only paradifferential contribution
is the one arising from the last term, which is identical to the one for the full equation, as expected.

We next identify the ``hh" contribution, i.e. the $N^{lin,[2]}_\eta$ term. This arises from the $V^{[1]}$,
$B^{[1]}$ and $G^{[1]}$ type terms, computed earlier,
\[
N^{lin,[2]}_{\eta,\lambda} = - \nabla (\M_{[\lambda]} \varphi_{\gtrsim \lambda} h_{\gtrsim \lambda})  - 
G_{[\lambda]} (\L_{[\lambda]} \varphi_{\gtrsim \lambda} h_{\gtrsim \lambda}) - \nabla(\eta_{\gtrsim \lambda} \M_{[\lambda]} r_{\gtrsim \lambda}) - G_{[\lambda]} ( \eta_{\gtrsim \lambda} \L_{[\lambda]} r_{\gtrsim \lambda}),
\]
which represents the exact polarization of the $N^{[2]}_{\eta,\lambda}$ term in the full equation.
 
Finally, we need to identify the ``hl" contribution $N^{[0]}$, which corresponds to the case when $h$ and $r$ are at low frequency.  This also arises from the $V^{[1]}$,
$B^{[1]}$ and $G^{[1]}$ type terms, so it is akin to $N^{lin,[2]}_{\eta,\lambda}$ but with adjusted 
frequency localizations,
\[
N^{lin,[0]}_{\eta,\lambda} = - \nabla (\M_{[\lambda]} \varphi_{\lambda} h_{\ll \lambda})  - 
G_{[\lambda]} (\L_{[\lambda]} \varphi_{\lambda} h_{\ll\lambda}) - \nabla(\eta_{ \lambda} \M_{[\lambda]} r_{\ll \lambda}) - G_{[\lambda]} ( \eta_{ \lambda} \L_{[\lambda]} r_{\ll \lambda}).
\]
Here the last two terms are apparently better behaved since we have a derivative on the low frequency factor $r_{\ll \lambda}$. However, recalling the identity $\nabla \M_{[\lambda]} = G_{[\lambda]} \L_{[\lambda]}$, it turns out that the first two terms exhibit a commutator structure which make them the better terms,
\[
N^{lin,[0]}_{\eta,\lambda} = - \nabla h_{\ll \lambda} \M_{[\lambda]} \varphi_{\lambda}   - 
[ G_{[\lambda]},h_{\ll\lambda}] \L_{[\lambda]} \varphi_{\lambda}  - \nabla(\eta_{ \lambda} \M_{[\lambda]} r_{\ll \lambda}) - G_{[\lambda]} ( \eta_{ \lambda} \L_{[\lambda]} r_{\ll \lambda}).
\]
Finally, we remark that in the above expression one may further localize $G_{[\lambda]}$, $\M_{[\lambda]}$
and $\L_{[\lambda]}$ below the frequency of $(h,r)$, and also substitute $G_{[\lambda]}$ by $\DN_{[\lambda]}$ modulo perturbative terms. 

\medskip

We continue with the second equation, where we have 
\[\begin{aligned}
N^{lin}_{\varphi} = & \ 
 - R(V, \nabla B^l h)
 - V \cdot R(B,\nabla h) + R(B, G(\eta)r) - R(V,\nabla r) - R(B, G(\eta) R(B,h) 
 \\ &  \ - \nabla_x(hV))
 \ -(  g  +  V^l \cdot T_{\nabla B^l} + T_{\partial_tB^l}-a^l) h.
 \end{aligned}
\]
The expression on the second line is paradifferential, and, as noted earlier in the context of the full equation, vanishes at the leading order, so its contribution here will also belong to the perturbative box. It remains to consider the expression on the first line.  Just as in the case of the full equation, there are no remaining paradifferential contributions here, so we need to examine $N^{lin,[2]}_{\varphi}$
and $N^{lin,[0]}_{\varphi}$. We first dispense with some perturbative terms:
\begin{itemize}
    \item The first term has a high frequency on $V$ and a lower frequency on $B$, which is also differentiated. So it is cubic with a favorable frequency balance, and thus perturbative.
    \item The same can be said about the $V^l$ contributions in the second and the last term,
    \[
 V^l \cdot R(B,\nabla h) - R(B, V^l \nabla h) - R(B, h \nabla V^l ) ,
    \]
   where the first two terms exhibit a commutator structure. 
\end{itemize}
This leaves us with 
\[\begin{aligned}
N^{lin}_{\varphi} \approx
  - R(V, R(B,\nabla h)) + R(B, G(\eta)r) - R(V,\nabla r) - R(B, G(\eta) R(B,h) - \nabla_x R(V,h)) .
 \end{aligned}
\]
We are now ready to consider the $N^{lin,[2]}_{\varphi}$ contribution. There is none such in the first and the last term, which contain at least two frequencies higher than that of $h$, so we are left with the two middle terms, where the only interesting ``hh" contribution comes from $(B^{[1]}, V^{[1]})$ and $r$,
\[
N^{lin,[2]}_{\varphi,\lambda} = \M_{[\lambda]} \varphi_{\gtrsim \lambda}  G_{[\lambda]} r_{\gtrsim \lambda} - \L_{[\lambda]} \varphi_{\gtrsim \lambda}  \nabla  r_{\gtrsim \lambda},
\]
which is also a nonsymmetric polarization of $ N^{[2]}_{\varphi,\lambda}$.

Lastly, we consider the $N^{lin,[0]}_{\varphi}$ contribution, which corresponds to the case when the linearized variables $(h,r)$ represent the second highest frequency, much smaller than the highest frequency.
No such configurations can occur in the first and the last term, so the only contributions are from the middle term, 
exactly mirroring the previous case with contributions from $(B^{[1]}, V^{[1]})$ and $r$,
\[
N^{lin,[0]}_{\varphi,\lambda} = \M_{[\lambda]} \varphi_{ \lambda}  G_{[\lambda]} r_{\ll \lambda} - \L_{[\lambda]} \varphi_{ \lambda}  \nabla  r_{\ll \lambda}.
\]

\end{proof}


\section{The complex variable \texorpdfstring{$u$}{}}\label{sec:complexvar}

\newcommand{\bb}{\vartheta}

So far we have worked with the physical real valued variables $\eta$ and $\psi$, as well as 
the good variable $\varphi$. The advantage 
in doing that was that we were able to produce clean heuristic expansions of the nonlinearities in the corresponding system for $(\eta,\varphi)$. The same set-up will be used in order to formalize these expansions
later on. However, for the resonant analysis in our density-flux identities and in the interaction Morawetz bounds  it is more convenient to aggregate the two physical variables into a single complex valued function $u$, defined as follows:

\begin{definition}\label{complexvariable}
For each $\lambda>0$, define the symbols
\begin{equation*}
\mathfrak{q}_{<\lambda}(x,\xi)
\coloneq
\left(
\frac{\mathfrak{g}_{<\lambda}(x,\xi)}
{a_{<\lambda}(x)}
\right)^{\frac12},
\qquad
\mathfrak{\omega}_{<\lambda}(x,\xi)
\coloneq
\left(
a_{<\lambda}(x)\mathfrak{g}_{<\lambda}(x,\xi)
\right)^{\frac12},
\end{equation*}
where $a_{<\lambda}$ is the Taylor coefficient corresponding to the truncated graph and velocity trace $\eta_{<\lambda}$ and $\psi_{<\lambda}$ and $\mathfrak{g}_{<\lambda}$ is the principal symbol of $G_{<\lambda}\coloneq G(\eta_{<\lambda})$. We define the complex variable $u$ by
\begin{equation}\label{complexudef}
u\coloneq\eta+iT_{\mathfrak{q}}\varphi,
\qquad
T_{\mathfrak{q}}\varphi
\coloneq
\int_{0}^{\infty}
Op(\mathfrak{q}_{<\delta\lambda}^{\,l})
\varphi_{\lambda}\frac{d\lambda}{\lambda}.
\end{equation}
Here $\mathfrak{q}_{<\delta\lambda}^{\,l}$ denotes the spatially truncated symbol defined according to the conventions of Subsection \ref{subsec:principal}.
\end{definition}

We also define the auxiliary real valued symbol
\begin{equation}
\bb_{<\lambda} = \frac12\{\tau + V_{[\lambda]} \cdot \xi), \ln \mfq_{<\lambda}\}_{para},       
\end{equation}
where the para-Lie bracket\footnote{See also the definition of the paramaterial derivative of a symbol in Section~\ref{sec:para-functions}.} of the material derivative and $\ln \mfq$ above is defined 
exactly like the standard Lie bracket 
but using paraproducts instead,
\[
\{\tau + V_{[\lambda]} \cdot \xi), \ln \mfq_{<\lambda}\}_{para}: = 
(\partial_t + T_{V^l} \cdot \nabla)\ln \mfq(t,x,\xi)
- T_{\partial_\xi \ln \mfq(t,x,\xi)}
\partial_x V \cdot \xi.
\]
The above symbols are chosen so that $u$ 
will solve a scalar equation of the form 
\begin{equation}\label{u-para-re}
   (\partial_t + V_{[\lambda]} \cdot \nabla
  +\frac12 (\nabla \cdot V_{[\lambda]})
 - \bb_{[\lambda]}(x,D) ) + i
\omega_{[\lambda]}(x,D)) u_\lambda  = N_\lambda(u),
\end{equation}
where the linear component was carefully chosen so that the remaining paradifferential type terms in the nonlinearity 
$N_\lambda$ are lower order. By lower order, we mean either
\begin{itemize}
    \item paradifferential operators of order $-\frac12$ applied to $u_\lambda$, or 
    \item paradifferential operators of order\footnote{A higher order is allowed here due to the fact that such terms are 
    strongly nonresonant, see the discussion in Section~\ref{subsec:resonant}.}
    $0$ applied to $\bu_\lambda$.
\end{itemize}
In particular we note that the contribution of $\vartheta$ does not qualify as lower order.

The form of the full equation above suggests 
that it should be instructive to consider first the associated linear paradifferential 
equation associated to a solution $u$ (and by extension, $(\eta,\psi)$) for the full equation,
\begin{equation}\label{u-para-only}
  (\partial_t + V_{[\lambda]} \cdot \nabla
  +\frac12 (\nabla \cdot V_{[\lambda]})
 - \bb_{[\lambda]}(x,D)  +
i \omega_{[\lambda]}) v_\lambda  = f_\lambda,
\end{equation}
Indeed, in later sections we will first carry out the density-flux analysis on this equation, and then use this as a starting point for the study of the full equation.

In a similar fashion, corresponding to the linearized system \eqref{Z-para-lin} we introduce the complex variable $w$ by 
\begin{equation}\label{complexwdef}
w\coloneq h+iT_{\mathfrak{q}}r,
\end{equation}  
which in turn will solve a scalar equation of the form 
\begin{equation}\label{u-para-lin-re}
  (\partial_t + V_{[\lambda]} \cdot \nabla
  +\frac12 (\nabla \cdot V)
 - \bb(x,D) )  +
i \omega_{[\lambda]}) w_\lambda  = N^{lin}_\lambda(u) w.
\end{equation}
Estimates for the linearized equation will play an essential role in the proof of our local well-posedness result.

\section{Frequency envelopes and the bootstrap set-up}
\label{s:freq+boot}
\label{sec:bootstrap}

The aim of this section is threefold:
\begin{enumerate}[label=\roman*)]
\item to introduce the bilinear bounds for solutions to the water wave
equations and to their linearization;
\item to restate our main results in the frequency envelope formalism,
with a full set of linear and bilinear bounds included;
\item to set up the bootstrap argument used in the proofs.
\end{enumerate}
This section thus opens a bootstrap loop spanning most of the rest of the paper, and which is closed only at the end
of Section~\ref{sec:Morawetz}.

\subsection{Bounds for the full equation}
\label{subsec:boot}

Let $(\eta_0,\varphi_0)\in\mathcal H^2\cap\mathcal B^2$, where
\begin{equation*}
\mathcal H^2:=H^2\times\big(\dot H^{\frac12}\cap\dot H^{\frac52}\big),
\hspace{5mm}
\mathcal B^2:=\dot B^2_{2,1}\times\dot B^{\frac52}_{2,1}.
\end{equation*}
We assume the data is sufficiently small in the critical Besov norm
\begin{equation*}
\|(\eta_0,\varphi_0)\|_{\mathcal B^2}\leq\epsilon.
\end{equation*}
Equivalently (after possibly relabeling $\epsilon$), we have $u_0\in H^2$
with
\begin{equation*}
\|u_0\|_{\dot B^2_{2,1}}\leq\epsilon.
\end{equation*}
We note that, since $\dot B^2_{2,1}\subset\dot H^2$, the assumption
$u_0\in H^2$ only adds the finiteness of the $L^2$ norm of $u_0$, on which
no smallness is assumed; this norm is controlled by the conserved energy.

We will have to track the smallness of the corresponding solution $u$ in
the critical norm. To do this, we place the initial data $u_0$ under a
critical frequency envelope $c_\lambda$,
\begin{equation}
\label{eq:env-crit}
\lambda^2\|u_{0,\lambda}\|_{L^2}\leq\epsilon c_\lambda,
\end{equation}
where $c_\lambda$ is assumed to be approximately $\ell^1$ normalized. That
is,
\begin{equation*}
\sum_{\lambda\in2^{\mathbb Z}}c_\lambda\lesssim1.
\end{equation*}
Here we use an  unbalanced slowly varying frequency envelope $c_\lambda$
\begin{equation*}
\Big(\frac\mu\lambda\Big)^{-C}\leq\frac{c_\mu}{c_\lambda}
\leq\Big(\frac\mu\lambda\Big)^{\delta}
\qquad\text{for all }\mu\geq\lambda,
\end{equation*}
where $C$ is a fixed, possibly large, constant and 
with $\delta$ a small universal constant. Such an envelope always exists,
and can be chosen so that
\[
\|u_0\|_{\dot B^2_{2,1}} \approx \epsilon \sum_\lambda c_\lambda,
\]
as well as 
\[
\|u_0\|_{\dot H^s}^2 \approx \epsilon^2 \sum_\lambda c_\lambda^2 \lambda^{2(s-s_c)}, \qquad s_c \leq s < s_c+C.
\]

 The envelope may thus decay rapidly at high frequencies, and grow at a slow rate. All the results in this paper are
compatible with this choice of frequency envelope, which in turn allows us to directly propagate higher Sobolev regularity for the solutions.

\begin{remark}
For convenience, we work with homogeneous frequency envelopes, which
allows us to close the critical and higher-regularity estimates
simultaneously. In the local well-posedness results, we also require
control\footnote{ This is qualitative, not quantitative, simply in order to be able to use existing local well-posedness results
as a starting point; likely this constraint can be removed. It also does not affect any implicit constants in our quantitative bounds.}
of the full inhomogeneous norm. This convention is harmless
because the conserved Hamiltonian controls the $L^2$ component of
$u$, while the frequency envelope controls its homogeneous
$\dot H^s$ norm.
\end{remark}

Our goal will then be to show that the solutions satisfy global frequency envelope bounds at critical regularity, as follows:
\begin{enumerate}[label=(\roman*)]
\item Uniform frequency envelope bound:
\begin{equation}\label{uk-ee}
\| u_\lambda \|_{L^\infty_t L^2_x} \lesssim \epsilon c_\lambda \lambda^{-s_c};
\end{equation}
\item Unbalanced bilinear $L^2_{t,x}$ bounds: Let $u_{\mu}$ and $u_{\lambda}$ be unbalanced, with $\mu\ll \lambda$. Then 
\begin{equation} \label{uab-bi-unbal}
\| u_\lambda \bu_\mu^{x_0}  \|_{L^2_{t,x}} \lesssim \epsilon^2 c_\lambda c_\mu \lambda^{-s_c} \mu^{-s_c+\frac{3}{4}}(1+\mu |x_0|)^\frac12;
\end{equation}
\item Translation-averaged bilinear unbalanced $L^2_{t,x}$ bounds: Let $u_{\mu}$ and $u_{\lambda}$ be unbalanced, with $\mu\ll\lambda$. Then
\begin{equation}\label{uab-bi-unbal-trans}
\left\|
u_{\lambda}(x)\overline{u_{\mu}}(y)
\right\|_{L_t^2L_{x,y}^2(B_r)}\lesssim  \epsilon^2c_{\lambda}c_{\mu}\lambda^{-s_c}\mu^{-s_c} \mu^{\frac{1}{4}} r^{\frac{1}{2}},\qquad r\geq \mu^{-1},
\end{equation}
where $B_r\subset \mathbb{R}^2\times\mathbb{R}^2$ is the ball $|x-y|\leq r$;
\item Balanced bilinear $L^2_{t,x}$ bound: 
\begin{equation} \label{uab-bi-bal}
\|P_{\alpha} (u_\lambda \bu_\mu^{x_0})  \|_{L^2_{t,x}} \lesssim \epsilon^2 c_\lambda c_\mu \lambda^{-2s_c+\frac{5}{4}}  \alpha^{-\frac12} (1+ \alpha |x_0|)^\frac12
, \qquad \alpha \lesssim \lambda \approx \mu ;
\end{equation}
\item Translation-averaged bilinear balanced $L^2_{t,x}$ bounds: 
\begin{equation} \label{uab-bi-bal-trans}
  \| (\partial_x+\partial_y)(u_\lambda(x) \bu_\mu  (y))  \|_{L^2_{t,x,y}(B_r)} \lesssim \epsilon^2 c_\lambda c_\mu \lambda^{-2s_c+\frac{5}{4}}   r^{\frac12}, \qquad \lambda \approx \mu ,\quad  r\geq \mu^{-1},
\end{equation}
as well as 
\begin{equation} \label{uab-bi-bal-trans+}
 \| (\partial_x-\partial_y)(u_\lambda(x) u_\mu(y))  \|_{L^2_{t,x,y}(B_r)} \lesssim \epsilon^2 c_\lambda c_\mu \lambda^{-2s_c+\frac{5}{4}}   r^{\frac12}
, \qquad \lambda \approx \mu, \quad  r\geq \mu^{-1};
\end{equation}

\item $L_t^4L_x^p$ bound: For $8<p\leq\infty$, 
\begin{equation}\label{uL4Lp}
\|u_{\lambda}\|_{L_t^4L_x^p}\lesssim_p \epsilon c_{\lambda}\lambda^{-s_c+\frac{7}{8}-\frac{2}{p}}. 
\end{equation}
as well as\footnote{We only use one pair of exponents i.e., $(6,4)$ but the range of exponents one can obtain in the interpolation is larger.}
\begin{equation}\label{uL6L4}
\|u_{\lambda}\|_{L_t^{6,\infty} L_x^4}\lesssim  \epsilon c_{\lambda}\lambda^{\frac{5}{12}-s_c}.
\end{equation}
\end{enumerate}

\medskip

Our main result is then
\begin{theorem}\label{t:local-fe3}
Let $T>0$, and $0<\epsilon\leq\epsilon_0$, where
$\epsilon_0>0$ is sufficiently small. Consider the equation \eqref{u-para-re}, and let
\begin{equation*}
u\in C([0,T];H^{2}\cap \dot{B}^2_{2,1})
\end{equation*}
be a smooth solution with initial data
$u_0\in H^{s_c}\cap\dot B_{2,1}^2$ satisfying
\begin{equation*}
\|u_0\|_{\dot B_{2,1}^2}\leq\epsilon.
\end{equation*}
Let $c_\lambda$ be a critical frequency envelope for the initial
data as in \eqref{eq:env-crit}. Then the solution $u$ satisfies the following bounds on $[0,T]$, uniformly with respect to $x_0\in\mathbb R^2 :$  \eqref{uk-ee}, \eqref{uab-bi-unbal},
\eqref{uab-bi-unbal-trans}, \eqref{uab-bi-bal}, \eqref{uab-bi-bal-trans},  \eqref{uab-bi-bal-trans+},
\eqref{uL4Lp} and \eqref{uL6L4}.
\end{theorem}
To prove the theorem, we first fix a sufficiently large constant
$C\gg1$. We then make the following  bootstrap assumptions
\begin{itemize}
\item[A.] Bootstrapped energy bounds:
\begin{equation}\label{uk-ee-boot}
\| u_\lambda \|_{L^\infty_t L^2_x} \lesssim C\epsilon c_\lambda \lambda^{-s_c};
\end{equation} 
\item[B.] Bootstrapped translation-averaged bilinear unbalanced $L^2_{t,x}$ bounds:
\begin{equation}\label{uab-bi-unbal-trans-boot}
\left\|
u_{\lambda}(x)\overline{u_{\mu}}(y)
\right\|_{L_t^2L_{x,y}^2(B_r)}\lesssim  C^2\epsilon^2c_{\lambda}c_{\mu}\lambda^{-s_c}\mu^{-s_c} \mu^{\frac{1}{4}} r^{\frac{1}{2}},\qquad\mu\ll\lambda,\hspace{2mm} r\geq \mu^{-1},
\end{equation}
where $B_r\subset\mathbb{R}^2\times\mathbb{R}^2$ is the ball $|x-y|\leq r$;
 \item[C.]  Bootstrapped translation-averaged bilinear balanced $L^2_{t,x}$ bounds: 
\begin{equation} \label{uab-bi-bal-trans-boot}
  \| (\partial_x+\partial_y)(u_\lambda(x) \bu_\mu  (y))  \|_{L^2_{t,x,y}(B_r)} \lesssim C^2\epsilon^2 c_\lambda c_\mu \lambda^{-2s_c+\frac{5}{4}}   r^{\frac12}, \qquad \lambda \approx \mu ,\quad  r\geq \mu^{-1};
\end{equation}
\item[D.] Bootstrapped bound for the Taylor coefficient,
\begin{equation}\label{Taylorcoeffsmall}
\sup_{\rho\in2^{\mathbb Z}}
\left(
\|a_{<\rho}-1\|_{L_t^\infty L_x^\infty}
+
\|\nabla_xa_{<\rho}\|_{L_t^\infty\dot{B}^0_{2,1}}
\right)
\lesssim C\epsilon,
\end{equation}
\end{itemize}
and  prove the same bounds
with constant $C/2$ (for instance), provided that $\epsilon\leq\epsilon_0(C)$.
The result then follows by a standard continuity argument.

\medskip 

All these bootstrap assumptions will be closed in the following sections.
Explicitly, in Section~\ref{sec:aexpansions} the bootstrap hypothesis
\eqref{Taylorcoeffsmall} will be bounded in terms of the constants in
the bootstrap assumptions \eqref{uk-ee}--\eqref{uL4Lp} for $u$, once we
relate the principal part of the Taylor coefficient to $\eta$. The energy
bound \eqref{uk-ee-boot} is closed in Section~\ref{sec:energy-close}. The
translation-averaged unbalanced bilinear $L^2_{t,x}$ bounds
\eqref{uab-bi-unbal-trans-boot} are closed in
Section~\ref{sec:unbalanced-close}, and the translation-averaged balanced
bilinear $L^2_{t,x}$ bounds \eqref{uab-bi-bal-trans-boot} in
Section~\ref{sec:balanced-close}.

\bigskip
For the proof of the theorem it will be convenient to also have access to the rest of the bounds, respectively:
\begin{itemize}

\item Unbalanced bilinear $L^2_{t,x}$ bounds: Let $u_{\mu}$ and $u_{\lambda}$ be unbalanced, with $\mu\ll \lambda$. Then 
\begin{equation} \label{uab-bi-unbal-boot}
\| u_\lambda \bu_\mu^{x_0}  \|_{L^2_{t,x}} \lesssim C^2 \epsilon^2 c_\lambda c_\mu \lambda^{-s_c} \mu^{-s_c+\frac{3}{4}}(1+\mu |x_0|)^\frac12;
\end{equation}

\item Balanced bilinear $L^2_{t,x}$ bound: 
\begin{equation} \label{uab-bi-bal-boot}
\|P_{\alpha}|D|^{\frac12} (u_\lambda \bu_\mu^{x_0})  \|_{L^2_{t,x}} \lesssim C^2\epsilon^2 c_\lambda c_\mu \lambda^{-2s_c+\frac{5}{4}}(1+ \alpha |x_0|)^\frac12
, \qquad \alpha\lesssim\lambda \approx \mu ;
\end{equation}

\item   Translation-averaged bilinear balanced $L^2_{t,x}$ bounds:  \begin{equation} \label{uab-bi-bal-trans+-boot}
 \| (\partial_x-\partial_y)(u_\lambda(x) u_\mu(y))  \|_{L^2_{t,x,y}(B_r)} \leq C^2\epsilon^2 c_\lambda c_\mu \lambda^{-2s_c+\frac{5}{4}}   r^{\frac12}
, \qquad \lambda \approx \mu, \quad  r\geq \mu^{-1};
\end{equation}
\item $L_t^4L_x^p$ bound: For $8<p\leq\infty$, 
\begin{equation}\label{L4Lp-boot}
\|u_{\lambda}\|_{L_t^4L_x^p}\lesssim_p C\epsilon c_{\lambda}\lambda^{-s_c+\frac{7}{8}-\frac{2}{p}} 
\end{equation}
as well as the weak-$L^6_t$ bound
\begin{equation}\label{uL6L4-boot}
\|u_{\lambda}\|_{L_t^{6,\infty} L_x^4}\lesssim  C\epsilon c_{\lambda}\lambda^{\frac{5}{12}-s_c}.
\end{equation}
\end{itemize}
The next lemma shows that these bounds can be seen as consequences of \eqref{uab-bi-unbal-trans-boot} and 
\eqref{uab-bi-bal-trans-boot}. However, we choose to state them separately, both for ease of reference and for comparison with the bounds in earlier work.

\begin{lemma}\label{lem:boot}
a) Assume that the bound \eqref{uab-bi-unbal-trans-boot} holds. Then we also have \eqref{uab-bi-unbal-boot}.

b) Assume that the bound \eqref{uab-bi-bal-trans-boot} holds. Then we also have \eqref{uab-bi-bal-boot},\eqref{uab-bi-bal-trans+-boot}, \eqref{L4Lp-boot} and
\eqref{uL6L4-boot}.
\end{lemma}

\begin{proof}
a) Denote $r = |x_0|+\mu^{-1}$. Then the 
bound \eqref{uab-bi-unbal-trans-boot}
implies that 
\[
\|u_\lambda(x) \bu_\mu^{x_0} (y)\|_{L^2\{|x-y|< \mu^{-1}\}} 
\lesssim C^2 \epsilon^2 c_\lambda c_\mu
\lambda^{-s_c} \mu^{-s_c} r^\frac12 \mu^{\frac14},
\]
or equivalently 
\[
\|u_\lambda(x) \bu_\mu^{x_0} (x+z)\|_{L^2\{|z|< \mu^{-1}\}} 
\lesssim C^2 \epsilon^2 c_\lambda c_\mu
\lambda^{-s_c} \mu^{-s_c} r^\frac12 \mu^{\frac14}.
\]
Now we use the frequency localization of $u_\mu^{x_0}$ to write its value at $x$ as a local average on the scale $\mu^{-1}$.
By H\"older's inequality, this implies
\[
\|u_\lambda(x) \bu_\mu^{x_0} (x)\|_{L^2_{t,x}} 
\lesssim C^2 \epsilon^2 c_\lambda c_\mu
\lambda^{-s_c} \mu^{-s_c} r^\frac12 \mu^{\frac54}.
\]
This suffices for \eqref{uab-bi-unbal-boot} since 
\[
r^{\frac12} \mu^\frac54 \approx \mu^{\frac34}(1+\mu|x_0|)^\frac12.
\]

b) Here we argue as in the Schr\"odinger case 
via the identity\footnote{ First exploited in the NLS case in \cite{PV}.}
\begin{equation}\label{ident}
| v(x) \partial_j \bar w(y) + 
\partial_j v(x) \bar w(y)|^2 
= | v(x) \partial_j  w(y) -
\partial_j v(x) w(y)|^2 + 
\partial_j |v(x)|^2
 \partial_j |w(y)|^2,
\end{equation}
used with $v = u_\lambda$ and $w= u_\mu$.

To ensure a good sign for the contribution of the last term on the right in \eqref{uab-bi-bal-trans-boot}, we consider first the case when $w = v$ and then use a soft
cutoff, e.g. a Gaussian
\[
\chi_r = e^{-\frac{|x|^2}{r^2}}.
\]
Then the bound \eqref{uab-bi-bal-trans-boot} implies 
\begin{equation}\label{bi-v=w}
\iint \chi_r(x-y) \partial_j |v(x)|^2
 \partial_j |v(y)|^2 \, dx dy dt \lesssim C^4 \epsilon^4 c_\lambda^4 \lambda^{-4 s_c}
\cdot \lambda^{\frac52} r,
\end{equation}
where the left hand side is now positive definite. Returning to \eqref{uab-bi-bal-trans-boot}, without assuming $v = w$,
using the identity \eqref{ident}, we can employ \eqref{bi-v=w} for both $v$ and $w$ 
in order to estimate the contribution of 
the second term on the right. This yields
a similar bound for the first term on the right, namely
\begin{equation}\label{bi-u-v}
\iint \chi_r(x-y)  |\nabla v_\lambda (x) w_\lambda(y) - v_\lambda(x) \nabla  w_\lambda(y))|^2   \, dx dy dt \lesssim  C^4 \epsilon^4 c_\lambda^4 \lambda^{-4 s_c} \cdot \lambda^{\frac52} r,
\end{equation}
which in turn yields the bound \eqref{uab-bi-bal-trans-boot}.

Returning to \eqref{bi-v=w}, on  the Fourier side this yields
\[
\| e^{-\frac{r^2 \xi^2}{8}} |\xi| \mathcal F | v|^2(\xi) \|_{L^2}^2  \lesssim  C^4 \epsilon^4 c_\lambda^4 \lambda^{-4 s_c} \cdot \lambda^{\frac52} r^{-1}.
\]
Restricting $\xi$ to the range $|\xi| \approx r^{-1}$ associated to the $r$ scale yields
\begin{equation}\label{bi-bal-u}
\| P_{r^{-1}}  |v|^2 \|_{L^2_{t,x}}^2 \lesssim  C^4 \epsilon^4 c_\lambda^4 \lambda^{-4 s_c} \cdot r \lambda^\frac52
\end{equation}
as needed for the bilinear estimate \eqref{uab-bi-bal-boot}.
Here by polarization we can separate 
the arguments and also obtain 
\begin{equation}\label{bi-bal-uv}
\| P_{r^{-1}}  (v \bar w) \|_{L^2_{t,x}}^2 \lesssim  C^4 \epsilon^4 c_\lambda^4 \lambda^{-4 s_c}  r \lambda^\frac52 .
\end{equation}
This proves \eqref{uab-bi-bal-boot} in the case $x_0 = 0$. For a general $x_0$, we adapt the argument used in (a). Precisely, the estimate \eqref{uab-bi-bal-trans-boot}
implies that for $v = u_\lambda$, $w= u_\mu^{x_0}$ we have 
\[
\int \chi_r(x-y) | v(x) \partial_j \bar w(y) + 
\partial_j v(x) \bar w(y)|^2 \, dx dy dt
\lesssim  C^4 \epsilon^4 c_\lambda^4 \lambda^{-4 s_c}\cdot \lambda^{\frac52} r, \qquad r \geq r_0:=\lambda^{-1}+|x_0|.
\]
By the earlier argument, this gives the lossless bound \eqref{bi-bal-uv} for $r \geq r_0$. For smaller $r$, we have the above bound with $r_0$ on the right, so we obtain instead
\begin{equation}\label{bi-bal-uv-r0}
\| P_{r^{-1}}  (v \bar w) \|_{L^2_{t,x}}^2 \lesssim  C^4 \epsilon^4 c_\lambda^4 \lambda^{-4 s_c}  r_0 \lambda^\frac52,
\end{equation}
 exactly as needed for \eqref{uab-bi-bal-boot}.

Finally, the $L_t^4L_x^p$ bound \eqref{L4Lp-boot} follows from \eqref{uab-bi-bal-boot}, by 
\begin{equation*}
\|u_{\lambda}\|_{L_t^4L_x^p}=\||u_{\lambda}|^2\|_{L_t^2L_x^{\frac{p}{2}}}^{\frac{1}{2}}\lesssim_p \lambda^{\frac{1}{4}-\frac{2}{p}}\sup_{\alpha\lesssim\lambda}\|P_{\alpha}|D|^{\frac{1}{2}}|u_{\lambda}|^2\|^{\frac{1}{2}}_{L_{t,x}^2}   , 
\end{equation*}
where we used Bernstein in the last inequality.

Finally, for the bound \eqref{uL6L4} we interpolate linearly between the energy bound
for $u_\lambda$ which gives
\[
\||u_\lambda|^2\|_{L^\infty L^1} 
\lesssim C^2 \epsilon^2 c_\lambda^2 \lambda^{-2 s_c}, 
\]
and the $l^\infty$ Besov type bound 
\[
\|P_\alpha |D|^\frac12 |u_\lambda|^2\|_{L^2_{t,x}}\lesssim C^2 \epsilon^2 c_\lambda^2 \lambda^{-2 s_c+\frac54}, \qquad \alpha \lesssim\lambda,
\]
which follows from \eqref{uab-bi-bal-boot}.
\end{proof}

The expressions appearing in our translation averaged balanced bilinear $L^2$ bounds can be thought of as 
\emph{null forms}, and will be denoted by $\mfQ$, i.e. 
\[
\mfQ(u,v) = (\partial_x + \partial_y) (u(x) \bv(y)) \quad \text{or}
\quad \mfQ(u,v) = (\partial_x - \partial_y) (u(x) v(y))
\]
By a slight abuse of notation, we will also 
 call \emph{null forms}, still denoted by $\mfQ$, 
 the single variable expressions 
\[
\mfQ(u,v) = \partial_x  (u(x) \bv(x)) \quad \text{or}
\quad \mfQ(u,v) = (\partial_x u(x) v(x)- u(x) \partial_x v(x))
\]
These expressions will play a major role in our analysis.
In particular, we will use the following bounds:

\begin{corollary}
Assume that the bound \eqref{uab-bi-bal-trans-boot} holds. Then we also have the null form bound 
\begin{equation}
\|\mfQ(u_\lambda,v_\nu)\|_{L^2_{t,x}} \lesssim 
C^2 \epsilon^2 c_\lambda^2 \lambda^{-2s_c+\frac74},\hspace{5mm}\nu\approx\lambda
\end{equation}
\end{corollary}
\begin{proof}
  For the first choice of $\mfQ$ above this is a direct consequence of 
  \eqref{uab-bi-bal-boot}.
  For the second, we use in addition the relation \eqref{ident} with $x = y$.
\end{proof}
\subsection{Bounds for the paradifferential 
equation and for the linearized equation}
\label{subsec:para-boot}

Here we repeat the discussion in the previous subsection in the context of the paradifferential equation \eqref{u-para-only}
and the linearized equation \eqref{u-para-lin}. A-priori this should be considered in Sobolev spaces $\H^s$, but 
for convenience we will set $s = 0$ in what follows, and simply work in $\H^0 = L^2 \times \dot H^{\frac12}$. In both cases we denote by $w_\lambda$ the corresponding complex variable, which will be considered in $L^2$.

We denote by $d_\lambda$ a slowly varying  $\H^0$ frequency envelope for the initial data
$(h_0,r_0)$,
\begin{equation}
\| (h_{0,\lambda},r_{0,\lambda})\|_{\H^0} \leq d_\lambda,
\end{equation}
or equivalently an $L^2$ frequency envelope for $w_0$,
\begin{equation}
\| w_{0,\lambda}\|_{L^2}\lesssim d_\lambda.
\end{equation}

Our goal will then be to show that the solutions $w$ satisfy \emph{local} frequency envelope bounds as follows:

\begin{enumerate}[label=(\roman*)]
\item Uniform energy frequency envelope bound:
\begin{equation}\label{wk-ee}
\| w_\lambda \|_{L^\infty_t L^2_x} \lesssim d_\lambda ;
\end{equation}

\item Unbalanced bilinear $L^2_{t,x}$ bounds: Let $\mu\ll \lambda$. Then 
\begin{equation} \label{uw-bi-unbal}
\| u_\lambda \bw_\mu^{x_0}  \|_{L^2_{t,x}} \lesssim \epsilon  c_\lambda d_\mu\, \lambda^{-s_c} \mu^{\frac{3}{4}}(1+\mu |x_0|)^\frac12,
\end{equation}
\begin{equation} \label{wu-bi-unbal}
\| w_\lambda \bu_\mu^{x_0}  \|_{L^2_{t,x}} \lesssim \epsilon  d_\lambda c_\mu\, \mu^{-s_c+\frac{3}{4}}(1+\mu |x_0|)^\frac12,
\end{equation}
respectively 
\begin{equation} \label{ww-bi-unbal}
\| w_\lambda \bw_\mu^{x_0}  \|_{L^2_{t,x}} \lesssim  d_\lambda d_\mu\,  \mu^{\frac{3}{4}}(1+\mu |x_0|)^\frac12;
\end{equation}

\item Translation-averaged bilinear unbalanced $L^2_{t,x}$ bounds: Let
$\mu\ll\lambda$ and $r\geq\mu^{-1}$. Then
\begin{equation}\label{uw-bi-unbal-trans}
\left\| u_{\lambda}(x)\overline{w_{\mu}}(y) \right\|_{L_t^2L_{x,y}^2(B_r)}
\lesssim \epsilon c_{\lambda}d_{\mu}\,\lambda^{-s_c} \mu^{\frac{1}{4}} r^{\frac{1}{2}},
\end{equation}
\begin{equation}\label{wu-bi-unbal-trans}
\left\| w_{\lambda}(x)\overline{u_{\mu}}(y) \right\|_{L_t^2L_{x,y}^2(B_r)}
\lesssim \epsilon d_{\lambda}c_{\mu}\,\mu^{-s_c+\frac{1}{4}} r^{\frac{1}{2}},
\end{equation}
respectively
\begin{equation}\label{ww-bi-unbal-trans}
\left\| w_{\lambda}(x)\overline{w_{\mu}}(y) \right\|_{L_t^2L_{x,y}^2(B_r)}
\lesssim d_{\lambda}d_{\mu}\,\mu^{\frac{1}{4}} r^{\frac{1}{2}},
\end{equation}
where $B_r$ is the ball $|x-y|\leq r$;

\item Balanced bilinear $L^2_{t,x}$ bounds: Let $\lambda\approx\mu$ and $\alpha\lesssim\lambda$. Then
\begin{equation} \label{wuab-bi-bal}
\|P_{\alpha}|D|^{\frac12} (w_\lambda \bu_\mu^{x_0})  \|_{L^2_{t,x}} \lesssim \epsilon d_\lambda c_\mu\, \lambda^{-s_c+\frac{5}{4}} (1+ \alpha |x_0|)^\frac12,
\end{equation}
respectively
\begin{equation} \label{ww-bi-bal}
\|P_{\alpha}|D|^{\frac12} (w_\lambda \bw_\mu^{x_0})  \|_{L^2_{t,x}} \lesssim  d_\lambda d_\mu\, \lambda^{\frac{5}{4}} (1+ \alpha |x_0|)^\frac12 ;
\end{equation}

\item Translation-averaged bilinear balanced $L^2_{t,x}$ bounds: Let $\lambda\approx\mu$ and $r\geq\lambda^{-1}$. Then
\begin{equation} \label{wuab-bi-bal-trans}
\| (\partial_x+\partial_y)(w_\lambda(x) \bu_\mu(y))  \|_{L^2_{t,x,y}(B_r)} \lesssim \epsilon d_\lambda c_\mu\, \lambda^{-s_c+\frac{5}{4}}   r^{\frac12},
\end{equation}
respectively
\begin{equation} \label{ww-bi-bal-trans}
\| (\partial_x+\partial_y)(w_\lambda(x) \bw_\mu(y))  \|_{L^2_{t,x,y}(B_r)} \lesssim  d_\lambda d_\mu\, \lambda^{\frac{5}{4}}   r^{\frac12},
\end{equation}
as well as the corresponding bounds for the non-conjugated products,
\begin{equation} \label{wuab-bi-bal-trans+}
\| (\partial_x-\partial_y)(w_\lambda(x) u_\mu(y))  \|_{L^2_{t,x,y}(B_r)} \lesssim \epsilon d_\lambda c_\mu\, \lambda^{-s_c+\frac{5}{4}}   r^{\frac12},
\end{equation}
respectively
\begin{equation} \label{ww-bi-bal-trans++}
\| (\partial_x-\partial_y)(w_\lambda(x) w_\mu(y))  \|_{L^2_{t,x,y}(B_r)} \lesssim  d_\lambda d_\mu\, \lambda^{\frac{5}{4}}   r^{\frac12}.
\end{equation}
\item $L_t^4L_x^p$ bound: For $8<p\leq\infty$, 
\begin{equation}\label{wL4Lp}
\|w_{\lambda}\|_{L_t^4L_x^p}\lesssim d_{\lambda}\lambda^{\frac{7}{8}-\frac{2}{p}}, 
\end{equation}
as well as 
\begin{equation}\label{wL6L4}
\|w_{\lambda}\|_{L_t^6 L_x^4}\lesssim d_{\lambda}\lambda^{\frac{5}{12}}, 
\end{equation}
\end{enumerate}

\medskip

Our main result for the linearized equation is then
\begin{theorem}\label{t:local-fe3-lin} Let $0<T\leq 1$. Let $u$ be a solution to the full equation on $[0,T]$ satisfying the bounds \eqref{uk-ee-boot} -- \eqref{uL6L4-boot} with an $\ell^1$ frequency envelope $\epsilon c_\lambda$ which in addition satisfies 
$c_\lambda \lesssim \lambda^{-\frac1{12}}$.
Let $w$ be a solution to the linearized equation \eqref{u-para-lin} in the time interval $[0,T]$, with $L^2$ 
initial data $w(0)$ with frequency envelope $d_\lambda$.
Then the corresponding complex variable $w$ satisfies
\eqref{wk-ee}-\eqref{wL4Lp} on $[0,T]$,
uniformly with respect to $x_0\in\mathbb R^2$.
\end{theorem}

As discussed in the previous subsection~ \ref{subsec:boot}, to prove  Theorem ~\ref{t:local-fe3-lin} we follow the same strategy as before. We begin by fixing  a sufficiently large constant
$C\gg1$ and then making the following  bootstrap assumptions
\begin{enumerate}[label=(\roman*)]
\item[A'.] Bootstrapped Uniform energy frequency envelope bound:
\begin{equation}\label{wk-ee-boot}
\| w_\lambda \|_{L^\infty_t L^2_x} \lesssim C d_\lambda ;
\end{equation}
\item[B'.] Bootstrapped translation-averaged bilinear unbalanced $L^2_{t,x}$ bounds: Let
$\mu\ll\lambda$ and $r\geq\mu^{-1}$. Then
\begin{equation}\label{uw-bi-unbal-trans-boot}
\left\| u_{\lambda}(x)\overline{w_{\mu}}(y) \right\|_{L_t^2L_{x,y}^2(B_r)}
\lesssim  C^2\epsilon c_{\lambda}d_{\mu}\,\lambda^{-s_c} \mu^{\frac{1}{4}} r^{\frac{1}{2}},
\end{equation}
\begin{equation}\label{wu-bi-unbal-trans-boot}
\left\| w_{\lambda}(x)\overline{u_{\mu}}(y) \right\|_{L_t^2L_{x,y}^2(B_r)}
\lesssim C^2\epsilon d_{\lambda}c_{\mu}\,\mu^{-s_c+\frac{1}{4}} r^{\frac{1}{2}},
\end{equation}
respectively
\begin{equation}\label{ww-bi-unbal-trans-boot}
\left\| w_{\lambda}(x)\overline{w_{\mu}}(y) \right\|_{L_t^2L_{x,y}^2(B_r)}
\lesssim C^2 d_{\lambda}d_{\mu}\,\mu^{\frac{1}{4}} r^{\frac{1}{2}},
\end{equation}
where $B_r$ is the ball $|x-y|\leq r$;

\item[C'.] 
Bootstrapped translation-averaged bilinear balanced $L^2_{t,x}$ bounds: Let $\lambda\approx\mu$ and $r\geq\lambda^{-1}$. Then
\begin{equation} \label{wuab-bi-bal-trans-boot}
\| (\partial_x+\partial_y)(w_\lambda(x) \bu_\mu(y))  \|_{L^2_{t,x,y}(B_r)} \lesssim  C^2\epsilon d_\lambda c_\mu\, \lambda^{-s_c+\frac{5}{4}}   r^{\frac12},
\end{equation}
respectively
\begin{equation} \label{ww-bi-bal-trans-boot}
\| (\partial_x+\partial_y)(w_\lambda(x) \bw_\mu(y))  \|_{L^2_{t,x,y}(B_r)} \lesssim  C^2 d_\lambda d_\mu\, \lambda^{\frac{5}{4}}   r^{\frac12},
\end{equation}

\end{enumerate}
As before, we do not include the remaining bounds  in the bootstrap loop,
since they will can be recovered from the bootstrap bounds $A'$--$C'$ above exactly as in Lemma~\ref{lem:boot}. We state them below nonetheless:

\begin{enumerate}[label=(\roman*)]
    \item[D'.] Bootstrapped unbalanced bilinear $L^2_{t,x}$ bounds: Let $\mu\ll \lambda$. Then 
\begin{equation} \label{uw-bi-unbal-boot}
\| u_\lambda \bw_\mu^{x_0}  \|_{L^2_{t,x}} \lesssim C^2\epsilon  c_\lambda d_\mu\, \lambda^{-s_c} \mu^{\frac{3}{4}}(1+\mu |x_0|)^\frac12,
\end{equation}
\begin{equation} \label{wu-bi-unbal-boot}
\| w_\lambda \bu_\mu^{x_0}  \|_{L^2_{t,x}} \lesssim C^2\epsilon  d_\lambda c_\mu\, \mu^{-s_c+\frac{3}{4}}(1+\mu |x_0|)^\frac12,
\end{equation}
respectively 
\begin{equation} \label{ww-bi-unbal-boot}
\| w_\lambda \bw_\mu^{x_0}  \|_{L^2_{t,x}} \lesssim C^2 d_\lambda d_\mu\,  \mu^{\frac{3}{4}}(1+\mu |x_0|)^\frac12;
\end{equation}
\item[E'.] Bootstrapped balanced bilinear $L^2_{t,x}$ bounds: Let $\lambda\approx\mu$ and $\alpha\lesssim\lambda$. Then
\begin{equation} \label{wuab-bi-bal-boot}
\|P_{\alpha}|D|^{\frac12} (w_\lambda \bu_\mu^{x_0})  \|_{L^2_{t,x}} \lesssim C^2 \epsilon d_\lambda c_\mu\, \lambda^{-s_c+\frac{5}{4}} (1+ \alpha |x_0|)^\frac12,
\end{equation}
respectively
\begin{equation} \label{ww-bi-bal-boot}
\|P_{\alpha}|D|^{\frac12} (w_\lambda \bw_\mu^{x_0})  \|_{L^2_{t,x}} \lesssim  C^2 d_\lambda d_\mu\, \lambda^{\frac{5}{4}} (1+ \alpha |x_0|)^\frac12 ;
\end{equation}

\item[F'.] Bootstrapped additional translation-averaged bilinear balanced $L^2_{t,x}$ bounds: Let $\lambda\approx\mu$ and $r\geq\lambda^{-1}$. Then
\begin{equation} \label{wuab-bi-bal-trans+-boot}
\| (\partial_x-\partial_y)(w_\lambda(x) u_\mu(y))  \|_{L^2_{t,x,y}(B_r)} \lesssim C^2\epsilon d_\lambda c_\mu\, \lambda^{-s_c+\frac{5}{4}}   r^{\frac12},
\end{equation}
respectively
\begin{equation} \label{ww-bi-bal-trans+}
\| (\partial_x-\partial_y)(w_\lambda(x) w_\mu(y))  \|_{L^2_{t,x,y}(B_r)} \lesssim  C^2 d_\lambda d_\mu\, \lambda^{\frac{5}{4}}   r^{\frac12}.
\end{equation}
\item[G'.] Bootstrapped $L_t^4L_x^p$ bound: For $8<p\leq\infty$, 
\begin{equation}\label{wL4Lp-boot}
\|w_{\lambda}\|_{L_t^4L_x^p}\lesssim C d_{\lambda}\lambda^{\frac{7}{8}-\frac{2}{p}}, 
\end{equation}
as well as 
\begin{equation}\label{wL6L4-boot}
\|w_{\lambda}\|_{L_t^{6,\infty} L_x^4}\lesssim Cd_{\lambda}\lambda^{\frac{5}{12}}. 
\end{equation}
\end{enumerate}
\medskip 

As an intermediate step in the proof of  Theorem~\ref{t:local-fe3} and Theorem~\ref{t:local-fe3-lin}, we will prove a similar result for the linear paradifferential equation:
\begin{theorem}\label{t:local-fe3-lin-para} Let $u$ be a solution to the full equation satisfying the bounds \eqref{uk-ee-boot} -- \eqref{L4Lp-boot} with an $\ell^1$ frequency envelope $\epsilon c_\lambda$.
Let $w$ be a solution to the linearized equation \eqref{u-para-only}, with $L^2$ 
initial data $w(0)$ with frequency envelope $d_\lambda$.
Then the corresponding complex variable $w$ satisfies
\eqref{wk-ee}-\eqref{wL4Lp} on $[0,T]$,
uniformly with respect to $x_0\in\mathbb R^2$.
\end{theorem}

Here we can use again the bootstrap assumptions $A'$--$C'$.  One key difference, though, is that in the 
paradifferential result we no longer need to impose any 
constraint on the length of the time interval.

\section{Elliptic Estimates and Multi-linear Estimates and Expansions}
\label{sec:elliptic}

This section collects the elliptic and multilinear estimates needed for
the expansions of the water-waves equation proved in the following sections. We first set out our
paraproduct conventions, which are somewhat non-standard, and are defined in terms of not only frequency regularizations, but also "geometric" regularizations for operators and functions that depend nonlinearly on the good variables $\eta$ and $\varphi$. We also record the
linearization formulas for $G$, $L$, and $M$ with respect to these so-called geometric regularizations. We then establish
linear estimates for these operators, isolate their principal and
subprincipal parts, and control the high and low frequency tails
arising from their nonlinear dependence on $\eta$. We also derive the linear and bilinear bounds for the good variable
$\varphi$ from the bootstrap assumptions on $u$. 
\\

To handle multi-linear operators
whose kernels have limited spatial decay, we also develop a multi-linear functional framework based on the translation averaged bilinear bound in Section \ref{sec:bootstrap}. This framework is better adapted to these slowly decaying kernels that we will frequently encounter later. We also establish
the corresponding multi-linear elliptic and commutator
bounds that are compatible with this framework. These estimates will be crucial for retaining the bilinear structure
needed to control many of the error terms in the later expansions.
\subsection{Paradifferential Operators}\label{multpara}

In this subsection, we introduce a generalized continuous Little\-wood--Paley decomposition and the regularization conventions used in the layered paradifferential expansions in the following sections. In addition to the
standard Fourier-based truncations of the main variables $\eta$, $\psi$, and $\varphi$, we will need
nonlinear truncations of the operators $G$, $L$, and $M$, as well as
of the associated velocity traces $B$ and $V$. Here, we also record their linearization
formulas, define the corresponding generalized paraproducts, and
isolate their corresponding principal symbols that will be used in the expansions later in subsequent sections.
\\

\noindent We begin by describing a more general version of the Littlewood-Paley decomposition that was defined in Subsection \ref{subsec:littlewoodpaley}. For this, let $f$ be Schwartz function $f:\mathbb{R}^2\to\mathbb{C}$. Given a one parameter family of Schwartz functions $f_{<\lambda}:\mathbb{R}^2\to\mathbb{C}$ (not necessarily of the form $P_{<\lambda}f$) depending smoothly on $\lambda$ and $x$ and satisfying $\lim_{\lambda\to 0^+}f_{<0}=0$ and  $\lim_{\lambda\to\infty}f_{<\lambda}=f$ in $\mathcal{S}'$, we can write 
\[
f = \int_{0}^\infty \frac{d}{d\lambda} f_{<\lambda} d \lambda \coloneq \int_{0}^\infty  f_{\lambda} \frac{d \lambda}{\lambda}
\]
By taking
\begin{equation}\label{standardreg}
f_{<\lambda}\coloneq P_{<\lambda}f,\hspace{5mm}f_{\lambda}= \lambda \frac{d}{d\lambda} P_{<\lambda}f\coloneq S_{\lambda}f.
\end{equation}
we recover the standard (continuous) Littlewood-Paley decomposition defined in Subsection \ref{subsec:littlewoodpaley}. From here on, we make the the convention that the variables $\eta$, $\varphi$ and $\psi$ are regularized according to the standard decomposition above. That is, 
\begin{equation}
\eta_{<\lambda}\coloneq P_{<\lambda}\eta, \hspace{5mm}\varphi_{<\lambda}\coloneq P_{<\lambda}\varphi,\hspace{5mm} \psi_{<\lambda}\coloneq P_{<\lambda}\psi.
\end{equation}
We similarly define $\eta_{\lambda}$, $\varphi_{\lambda}$ and $\psi_{\lambda}$ as well as the regularizations for the linearized variables $h$ and $r$ according to this convention. 

On the other hand, many objects appearing in the water waves system (for instance, $V$ and $B$) depend nonlinearly (and non-locally) on these "standard" variables $\eta$, $\psi$ and $\varphi$, and we will not want to regularize them directly. Instead, we will need the more flexible definition of $f_{<\lambda}$ above to regularize these objects; in fact, they may not be strictly frequency localized like in \eqref{standardreg}. For instance, we can define nonlinearly truncated variants of the operators $G(\eta)$, $L(\eta)$ and $M(\eta)$ where the latter two operators are defined in \eqref{LM}
\begin{equation*}
\begin{split}
G_{<\lambda}\coloneq G(\eta_{<\lambda}),\hspace{5mm}L_{<\lambda}\coloneq L(\eta_{<\lambda}),\hspace{5mm} M_{<\lambda}\coloneq M(\eta_{<\lambda}).    
\end{split}    
\end{equation*}
We will frequently encounter functions of the form $f_{<\lambda}=A_{<\lambda}P_{<\lambda}g$ where $A\in\{G,L,M\}$ is one of the three operators above. In this case, $f_{<\lambda}$ is not strictly frequency localized because of high-tail leakage coming from the nonlinear dependence of $A_{<\lambda}$ on $\eta_{<\lambda}$. One important place this occurs is in our definition of the vertical and horizontal \textit{truncated} velocity traces,
\begin{equation}\label{BVop}
B_{<\lambda}=L_{<\lambda}\psi_{<\lambda},\hspace{5mm} V_{<\lambda}=M_{<\lambda}\psi_{<\lambda}\end{equation}
To compute the linearization
\begin{equation*}
\lambda\frac{d}{d\lambda}f_{<\lambda}   
\end{equation*}
for the above three operators (and more complicated examples later), we record the following linearization formulas for $G_{<\lambda}$, $L_{<\lambda}$ and $M_{<\lambda}$,
\begin{equation}\label{Blin}
\lambda\frac{d}{d\lambda}L_{<\lambda}g=-L_{<\lambda}(L_{<\lambda}g\eta_{\lambda})+\eta_{\lambda}L_{<\lambda}^2g    
\end{equation}
and
\begin{equation}\label{Vlin}
\lambda\frac{d}{d\lambda}M_{<\lambda}g=-M_{<\lambda}(L_{<\lambda}g\eta_{\lambda})+\eta_{\lambda}M_{<\lambda}L_{<\lambda}g    
\end{equation}
and the linearization formula for $G_{<\lambda}$,
\begin{equation}\label{Glin}
\lambda\frac{d}{d\lambda}G_{<\lambda}g=-G_{<\lambda}(\eta_{\lambda}L_{<\lambda}g)-\nabla\cdot (\eta_{\lambda}M_{<\lambda}g).
\end{equation}
The formula for $G_{<\lambda}$ follows from \eqref{G-lin}, whereas the formulas for $L$ and $M$ follow from Lemma \ref{l:DBV} by taking $g=\psi$ and substituting $B$ and $V$ for the equivalent expressions $L\psi$ and $M\psi$. In the above formulas, $g$ is fixed and independent of $\lambda$. If in addition, we replace $g$ by $g_{<\lambda}$, then to each of the formulas above, we add an  additional term $A_{<\lambda}g_{\lambda}$ where $A_{<\lambda}$ denotes one of the above three operators. These formulas will play a central role in our refined paradifferential expansion of the water waves equation later on.
\subsection{Paraproducts}\label{paraconventions}

With the above conventions, we can define the corresponding family of generalized paradifferential operators. First, from subsection \ref{subsec:paraproductssec2}, we recall the standard Bony decomposition, where we define
\begin{equation*}
T_fg\coloneq \int_{0}^{\infty}P_{<\delta\lambda}fS_{\lambda}g\frac{d\lambda}{\lambda},\hspace{5mm}\Pi(f,g)\coloneq fg-T_fg-T_gf   
\end{equation*}
and also
\begin{equation*}
R(f,g)\coloneq fg-T_fg    
\end{equation*}
to isolate the portion of $fg$ where $f$ is at higher or comparable frequency to $g$. Above, $0<\delta\ll 1$ is some fixed, small positive constant of the form $\mathfrak{c}^{-K}$ (to match the convention in subsection \ref{subsec:paraproductssec2}).
\\

Let $A$ denote either $G$, $L$ or $M$, as above. With $\delta$ as above, fix some large universal constant $M\gg 1$ such that $\delta M\ll 1$ (for instance $M=\delta^{-\frac{1}{2}})$. We define the operator-valued paradifferential operator $T_A$ by
\begin{equation*}
T_{A}f\coloneq \int_{0}^{\infty}P_{<M\lambda}(A_{<\delta\lambda}f_{\lambda})\frac{d\lambda}{\lambda}.    
\end{equation*}
\begin{remark}
The additional projection enforces localization to output frequencies
$\lesssim M\lambda$ and separates off the high-output tail arising
from the fact that $A_{<\delta\lambda}f_\lambda$ need not be exactly
frequency localized. The complementary high-output contribution will typically
be treated separately as a perturbative error, in practice.
\end{remark}

\subsection{Truncated variables}

In some of our paradifferential expansions later, the variables $B_{<\lambda}$ and $V_{<\lambda}$ will also naturally appear as para-coefficients. These can be thought of as nonlinear (and non-local) functions of the independent variables $\eta_{<\lambda}$ and $\psi_{<\lambda}$. These functions are not strictly localized to frequency $\lesssim\lambda$, but rather only up to high frequency tails. Therefore, we introduce the truncated variants of $V$ and $B$
\begin{equation*}
B^l_{<\lambda}\coloneq P_{< M\lambda}B_{<\lambda},\hspace{5mm}V^l_{<\lambda}\coloneq P_{<M\lambda} V_{<\lambda}   
\end{equation*}
Correspondingly, let $A_{<\lambda}$ denote either $G_{<\lambda}$, $L_{<\lambda}$ or $M_{<\lambda}$. For $\mu>0$ and $f_{<\mu}=P_{<\mu}f$, we also define the truncated expressions 
\begin{equation*}
A_{<\lambda}^lf_{<\mu}\coloneq P_{<M\max\{\lambda,\mu\}}A_{<\lambda}f_{<\mu}    
\end{equation*}
for some large, fixed constant $M>1$. This is defined to account for the fact that $G_{<\lambda}f_{<\mu}$, $L_{<\lambda}f_{<\mu}$ and $M_{<\lambda}f_{<\mu}$ are not strictly localized to frequencies $\lesssim \lambda+\mu$ (even though $\eta_{<\lambda}$ is itself frequency localized). 

Finally, it is also useful to define the "geometric" paraproducts associated to $B$ and $V$,
\begin{equation*}
T_{B^l}f=\int_{0}^{\infty}B^l_{<\delta\lambda}S_{\lambda}f\frac{d\lambda}{\lambda},\hspace{5mm}T_{V^l}f=\int_{0}^{\infty}V^l_{<\delta\lambda}S_{\lambda}f\frac{d\lambda}{\lambda}    
\end{equation*}
We can similarly define $T_{a^l}$ where $a$ is the Taylor coefficient, and will later also define other geometric paraproducts with this convention, as needed. We note above that $B_{<\delta \lambda}^l$ can be equivalently described using the notation $B_{[\delta M\lambda]}$ which is defined in \eqref{Bl}. A similar description holds with $V^l_{<\delta\lambda}$ and also $a^l_{<\delta\lambda}$. We will primarily adopt the $F^l_{<\delta\lambda}$ notation in the error analysis leading to the layered paradifferential equation for $u$, and then switch to the other notation afterwards. 
\subsection{Principal terms}\label{subsec:principal} 

As mentioned in previous sections, in our layered paradifferential expansion of the water-waves equation, whenever the operators $G$, $L$ and $M$ appear, we will try (as often as is possible) to reduce them to their principal parts for more manageable computation. We recall from Section \ref{sec:wwvar} that the operator $G_{<\lambda}$ has principal symbol given by
\begin{equation*}
\mathfrak{g}_{<\lambda}(x,\xi)\coloneq \sqrt{(1+|\nabla\eta_{<\lambda}|^2)|\xi|^2-(\nabla\eta_{<\lambda}\cdot\xi)^2}.    
\end{equation*}
Similarly, from the principal symbol for $G_{<\lambda}$, we can easily compute the principal symbols for $L_{<\lambda}$ and
$M_{<\lambda}$, which are given by
\begin{equation*}
\mathfrak{l}_{<\lambda}(x,\xi)
\coloneq
\frac{
\mathfrak{g}_{<\lambda}(x,\xi)
+i\nabla\eta_{<\lambda}\cdot\xi
}{
1+|\nabla\eta_{<\lambda}|^2
}
\end{equation*}
and
\begin{equation*}
\mathfrak{m}_{<\lambda}(x,\xi)
\coloneq
i\xi-\nabla\eta_{<\lambda}\,
\mathfrak{l}_{<\lambda}(x,\xi),
\end{equation*}
respectively. One can then define the truncated versions $\mathfrak{g}^l$, $\mathfrak{m}^l$ and $\mathfrak{l}^l$ as above by projecting the spatial part of the symbol to frequency $\lesssim M\lambda$. We can then define the principal paradifferential operators
\begin{equation*}
T_{\mathfrak{g}}f=\int_{0}^{\infty}Op(\mathfrak{g}^l_{<\delta\lambda})f_{\lambda}\frac{d\lambda}{\lambda}\hspace{5mm}T_{\mathfrak{l}}f=\int_{0}^{\infty}Op(\mathfrak{l}^l_{<\delta\lambda})f_{\lambda}\frac{d\lambda}{\lambda},\hspace{5mm}T_{\mathfrak{m}}f=\int_{0}^{\infty}Op(\mathfrak{m}^l_{<\delta\lambda})f_{\lambda}\frac{d\lambda}{\lambda},     
\end{equation*}
whose precise form we had postponed defining up until now. Above, and in the sequel, the pseudo-differential operator $Op(\mathfrak{a})$ is defined through the standard left-quantization. This way, the expression  $Op(\mathfrak{g}^l_{<\delta\lambda})$ is consistent with the notation $\mathfrak{G}_{[\delta M\lambda]}$ defined above the equation \eqref{Z-para+}, and likewise for the other two operators. Later, we will obtain subprincipal error bounds between the full operators $G$, $L$ and $M$, and their principal paradifferential operators above in various linear and multilinear norms.


\subsection{Paraproduct estimates} 

In this subsection, we record some standard $L^p$ based and homogeneous Besov estimates for the para-products introduced above. These bounds will be used repeatedly to estimate various errors in our paradifferential decomposition of the equation. We will also recall some standard Moser types estimates and as well as a useful bilinear Leibniz type identity involving commutators with Littlewood-Paley projections.
\\

We begin with the following standard estimates, which are essentially a consequence of Bony's decomposition and Coifman-Meyer estimates; see, for instance, \cite[Theorems 2.82 and 2.85]{bahouricheminDanchin} and \cite{grafakosTorres}.
\begin{proposition}\label{besovparabounds}
Let $\sigma\geq 0$. We have the following standard para-product estimates.
\begin{enumerate}
\item ($L^p$ estimates) Let $1\leq p_1<\infty$ and $1<p_2,p_3\leq\infty$ satisfy $\frac{1}{p_1}=\frac{1}{p_2}+\frac{1}{p_3}.$ Then the low-high para-product $T_fg$ satisfies
\begin{equation*}
\||D|^{\sigma}(T_fg)\|_{L^{p_1}}\lesssim \|f\|_{L^{p_2}}\||D|^{\sigma}g\|_{L^{p_3}},
\end{equation*}
and the high-high para-product $\Pi(f,g)$ satisfies
\begin{equation*}
\||D|^{\sigma}\Pi(f,g)\|_{L^{p_1}}\lesssim \|f\|_{L^{p_2}}\||D|^{\sigma}g\|_{L^{p_3}}.
\end{equation*}

\item (Besov estimates) Let $1\leq p_1,p_2,p_3\leq\infty$, $1\leq q\leq\infty$, and $\frac{1}{p_1}=\frac{1}{p_2}+\frac{1}{p_3}.$ The low-high para-product $T_fg$ satisfies
\begin{equation}\label{low-high}
\||D|^{\sigma}(T_fg)\|_{\dot{B}^0_{p_1,q}}\lesssim \|f\|_{\dot{B}^0_{p_2,1}}\||D|^{\sigma}g\|_{\dot{B}^0_{p_3,q}}.
\end{equation}
If $\sigma>0$, the high-high para-product $\Pi(f,g)$ satisfies
\begin{equation*}
\||D|^{\sigma}\Pi(f,g)\|_{\dot{B}^0_{p_1,1}}\lesssim \|f\|_{\dot{B}^0_{p_2,1}}\||D|^{\sigma}g\|_{\dot{B}^0_{p_3,1}}.
\end{equation*}
\end{enumerate}
Above, the implicit constants depend on the corresponding exponents and on $\sigma$.
\end{proposition}
We also need the following standard Moser estimates.
\begin{proposition}\label{MoserBesov}
Let $\sigma>0$, let $G$ be a smooth function with $G(0)=0$, and let $u\in L^\infty$. We have the following estimates:
\begin{enumerate}
\item ($L^p$ estimate) If $1<p<\infty$, then
\begin{equation*}
\||D|^{\sigma}G(u)\|_{L^p}
\lesssim_{p,\sigma}
C_G(\|u\|_{L^{\infty}})
\||D|^{\sigma}u\|_{L^p}.
\end{equation*}
\item (Besov estimate) If $1\leq p\leq\infty$, then
\begin{equation*}
\||D|^{\sigma}G(u)\|_{\dot{B}^0_{p,1}}
\lesssim_{p,\sigma}
C_G(\|u\|_{L^{\infty}})
\||D|^{\sigma}u\|_{\dot{B}^0_{p,1}}.
\end{equation*}
\end{enumerate}
Here, $C_G$ is a nondecreasing function depending on finitely many derivatives of $G$ on bounded sets.
\end{proposition}
The $L^p$ estimate is the homogeneous counterpart of the standard Moser estimate in \cite[Chapter~2, Section~4]{taylor2000tools}. The Besov estimate follows from \cite[Theorem~2.61]{bahouricheminDanchin}.
\\

Next, to efficiently deal with commutators of low frequency coefficients with Littlewood-Paley projections, we will frequently use the following Leibniz type rule from \cite{tao2001global}.
\begin{proposition}\label{Leibnizcommutator}
Let $S_{\lambda}$ denote either $P_{<\lambda}$ or $P_{\lambda}$. We have the following commutator identity
\begin{equation*}
[S_{\lambda},f]g=\lambda^{-1}L(\nabla_xf,g)    
\end{equation*}
where $L$ is a translation-invariant bilinear form with rapidly decaying kernel at scale $\lambda$.
\end{proposition}
\subsection{Flattened coordinates and elliptic preliminaries}\label{subsec:flattenedcoords}
Our next goal is to establish various linear elliptic estimates involving the operators $L$, $M$ and $G$. A recurring step in our analysis will be to express the harmonic extension operator in an appropriate flattened coordinate system. Throughout this section, let $u$ be a solution satisfying the bootstrap hypotheses of Section \ref{sec:bootstrap}, and let $\eta=\Re(u)$ denote the defining graph function for the fluid domain
\begin{equation*}
\Omega_t\coloneq \{(x,z)\in\mathbb{R}^2\times\mathbb{R}:z<\eta(t,x)\}.
\end{equation*}
In particular, the smallness built into the bootstrap assumptions will allow us to treat the variable coefficient elliptic estimates below as perturbations of their flat counterparts. This will drastically simplify a lot of our analysis compared to the large data scenario.
\\

For a function $f(t,\cdot)$ defined on the closure of $\Omega_t$, we define its flattened pullback by
\begin{equation*}
f_*(t,x,z)\coloneq f(t,x,z+\eta(t,x)),
\end{equation*}
so that $f_*$ is defined on the closed half-space
\begin{equation*}
\mathbb{H}\coloneq\{(x,z)\in\mathbb{R}^2\times\mathbb{R}:z\leq0\}.
\end{equation*}
The pullback map transforms ordinary derivatives into the variable coefficient operators
\begin{equation*}
D_i\coloneq\partial_i-\partial_i\eta\,\partial_z,\hspace{5mm}i=1,2,3,
\end{equation*}
where $\partial_3=\partial_z$ and we use the convention $\partial_3\eta=0$. We then have
\begin{equation*}
D_if_*=(\partial_if)_*.
\end{equation*}
Here and below, repeated indices are summed over $i=1,2,3$. With this notation, we have
\begin{equation}\label{variablelaplacian}
\Delta_{\eta}\coloneq D_iD_i
=(1+|\nabla\eta|^2)\partial_z^2+\Delta_x
-2\nabla\eta\cdot\nabla\partial_z-(\Delta_x\eta)\partial_z.
\end{equation}

We identify functions on $\partial\Omega_t$ with their pullbacks under the graph parametrization $x\mapsto(x,\eta(t,x))$. Thus, a boundary datum will always be written as a function $g=g(t,x)$ on $\mathbb{R}^2$. Let $\mathcal{H}_{\eta}g$ denote its harmonic extension to $\Omega_t$, and define
\begin{equation*}
w_g\coloneq(\mathcal{H}_{\eta}g)_*.
\end{equation*}
By \eqref{variablelaplacian}, we have
\begin{equation*}
\Delta_{\eta}w_g=0,\hspace{5mm}
\left.w_g\right|_{z=0}=g,\hspace{5mm}
\nabla_{x,z}w_g\to0\quad\text{as}\quad z\to-\infty.
\end{equation*}
By definition of $L(\eta)$, we also have the useful trace identity
\begin{equation*}
\partial_zw_g|_{z=0}=L(\eta)g.    
\end{equation*}
Moreover, the boundary operators $M(\eta)$ and $G(\eta)$ can be written in terms of $L(\eta)$ as follows
\begin{equation}\label{MtoL}
M(\eta)g
=
\left.(D_1w_g,D_2w_g)\right|_{z=0}
=
\nabla g-\nabla\eta\,L(\eta)g
\end{equation}
and
\begin{equation}\label{MtoG}
\begin{split}
G(\eta)g=
\left.\left(D_3w_g-\nabla\eta\cdot(D_1w_g,D_2w_g)\right)\right|_{z=0}=
(1+|\nabla\eta|^2)L(\eta)g-\nabla\eta\cdot\nabla g.
\end{split}
\end{equation}
Using these identities, we will later be able to easily convert most estimates for the operator $L(\eta)$ to estimates for $G(\eta)$ and $M(\eta)$ using the above identities. The benefit isolating $L(\eta)$ is its interpretation as the boundary trace for the vertical derivative of the flattened harmonic extension. This will have slightly nicer commutation properties with various vector fields.
\\

Finally, for the remainder of this section, for a dyadic number $\beta$, we will also write
\begin{equation*}
g_{\beta}\coloneq P_{\beta}g,\hspace{5mm}
g_{<\beta}\coloneq P_{<\beta}g.
\end{equation*}
The same conventions will be used with $\eta$ replaced by a truncated graph $\eta_{<\mu}$. We adopt this notation (for $\eta$ specifically) throughout the remainder of the paper. To conclude this subsection, we record a basic estimate for the flattened harmonic extension that we will use many times throughout this paper.
\begin{proposition}\label{flattenedbasicbounds}
Let $1<p\leq \infty$, and let $\epsilon>0$ be sufficiently small. Assume that
$\eta\in \dot{B}_{2,1}^2$ satisfies
\begin{equation*}
\|\eta\|_{\dot{B}_{2,1}^2}\leq \epsilon.
\end{equation*}
Then, for every boundary function $g$ such that
$\nabla_xg\in\dot{B}_{p,1}^0$, its flattened harmonic extension with respect to
$\eta$ satisfies
\begin{equation*}
\|\nabla_{z,x}w_g\|_{L_z^{\infty}\dot{B}_{p,1}^0}
\lesssim_p
\|\nabla_xg\|_{\dot{B}_{p,1}^0}.
\end{equation*}
\end{proposition}

\begin{proof}
The idea will be to treat $\Delta_{\eta}$ as a small perturbation of the flat
Laplacian $\Delta_{z,x}$ and to compare $w_g$ with the flat harmonic extension
$e^{z|D|}g$. It will be useful to write the variable coefficient part of
$\Delta_{\eta}$ in divergence form. Indeed, we have
\begin{equation}
\begin{split}
-|\nabla\eta|^2\partial_z^2h
+2\nabla\eta\cdot\nabla_x\partial_zh
+\Delta_x\eta\partial_zh
&=
\nabla_x\cdot(\nabla_x\eta\partial_zh)+
\partial_z\left(
\nabla_x\eta\cdot\nabla_xh
-|\nabla\eta|^2\partial_zh
\right)
\\
&\coloneq
\nabla_{z,x}\cdot(A_{\eta}\nabla_{z,x}h).
\end{split}
\end{equation}
Here, with the convention
$\nabla_{z,x}=(\partial_z,\nabla_x)$, we may take
\begin{equation*}
A_{\eta}
=
\begin{pmatrix}
-|\nabla\eta|^2 & (\nabla\eta)^{T}\\
\nabla\eta & 0
\end{pmatrix}.
\end{equation*}
In particular, the entries of $A_{\eta}$ are either linear or quadratic in
$\nabla\eta$. Using $\Delta_{\eta}w_g=0$, we find
\begin{equation*}
\Delta_{z,x}w_g
=
\nabla_{z,x}\cdot
(A_{\eta}\nabla_{z,x}w_g).
\end{equation*}

Let $\Delta_{z,x}^{-1}$ denote the inverse of the Dirichlet Laplacian on the
half-space $\mathbb{H}$, with zero trace at $z=0$ and vanishing gradient as
$z\to-\infty$, and define
\begin{equation*}
\mathcal{R}_0
\coloneq
\nabla_{z,x}\Delta_{z,x}^{-1}\nabla_{z,x}\cdot.
\end{equation*}
Since $w_g-e^{z|D|}g$ has zero boundary trace, we obtain
\begin{equation}\label{wgduhamel}
\nabla_{z,x}w_g
=
(|D|,\nabla_x)e^{z|D|}g
+
\mathcal{R}_0(A_{\eta}\nabla_{z,x}w_g).
\end{equation}

The standard Poisson kernel estimate gives
\begin{equation*}
\|(|D|,\nabla_x)e^{z|D|}g\|_{L_z^{\infty}\dot{B}_{p,1}^0}
\lesssim_p
\|\nabla_xg\|_{\dot{B}_{p,1}^0}.
\end{equation*}
On the other hand, by the reflection formula for the Dirichlet Laplacian on
the half-space, the kernel $K_{\alpha}(z,z',y)$ of
$P_{\alpha}\mathcal{R}_0$ satisfies
\begin{equation*}
\sup_{z\leq0}
\int_{-\infty}^{0}
|K_{\alpha}(z,z',y)|\,dz'
\lesssim_N
\alpha^2\langle\alpha y\rangle^{-N},
\qquad N>2,
\end{equation*}
where the integral is understood as a point evaluation at $z'=z$ for the "local"
$\delta(z-z')$ component of the kernel. Consequently,
\begin{equation*}
\|\mathcal{R}_0F\|_{L_z^{\infty}\dot{B}_{p,1}^0}
\lesssim_p
\|F\|_{L_z^{\infty}\dot{B}_{p,1}^0}.
\end{equation*}
Applying these estimates to \eqref{wgduhamel}, we obtain
\begin{equation}\label{wgpreabsorb}
\|\nabla_{z,x}w_g\|_{L_z^{\infty}\dot{B}_{p,1}^0}
\lesssim_p
\|\nabla_xg\|_{\dot{B}_{p,1}^0}
+
\|A_{\eta}\nabla_{z,x}w_g\|_{L_z^{\infty}\dot{B}_{p,1}^0}.
\end{equation}

It remains to control the last term by a small multiple of
$\|\nabla_{z,x}w_g\|_{L_z^{\infty}\dot{B}_{p,1}^0}$. Choose a finite
$r\geq2$ sufficiently large that $\frac{1}{p}+\frac{1}{r}\leq1,$ and define $q\geq1$ by $\frac{1}{q}=\frac{1}{p}+\frac{1}{r}$. Such a finite $r$ exists because $p>1$. By Proposition
\ref{besovparabounds}, for each fixed $z\leq0$ we have
\begin{equation*}
\begin{split}
\|A_{\eta}\nabla_{z,x}w_g\|_{\dot{B}_{p,1}^0}\lesssim
\|A_{\eta}\|_{\dot{B}_{\infty,1}^0}
\|\nabla_{z,x}w_g\|_{\dot{B}_{p,1}^0}+\|\Pi(A_{\eta},\nabla_{z,x}w_g)\|_{\dot{B}_{p,1}^0}.
\end{split}
\end{equation*}
For the high--high term, Sobolev embedding and the high--high estimate in
Proposition \ref{besovparabounds} give
\begin{equation*}
\begin{split}
\|\Pi(A_{\eta},\nabla_{z,x}w_g)\|_{\dot{B}_{p,1}^0}
&\lesssim
\||D|^{\frac{2}{r}}
\Pi(A_{\eta},\nabla_{z,x}w_g)\|_{\dot{B}_{q,1}^0}\lesssim
\||D|^{\frac{2}{r}}A_{\eta}\|_{\dot{B}_{r,1}^0}
\|\nabla_{z,x}w_g\|_{\dot{B}_{p,1}^0}.
\end{split}
\end{equation*}
By Sobolev embedding, Proposition \ref{MoserBesov}, and the fact that the
entries of $A_{\eta}$ vanish when $\nabla\eta=0$, we have
\begin{equation*}
\begin{split}
\|A_{\eta}\|_{\dot{B}_{\infty,1}^0}
+
\||D|^{\frac{2}{r}}A_{\eta}\|_{\dot{B}_{r,1}^0}
&\lesssim
\||D|A_{\eta}\|_{\dot{B}_{2,1}^0}\lesssim
C(\|\nabla\eta\|_{L^{\infty}})
\|\nabla_x^2\eta\|_{\dot{B}_{2,1}^0}\lesssim\epsilon.
\end{split}
\end{equation*}
Taking the supremum over $z\leq0$, it follows that
\begin{equation*}
\|A_{\eta}\nabla_{z,x}w_g\|_{L_z^{\infty}\dot{B}_{p,1}^0}
\lesssim_p
\epsilon
\|\nabla_{z,x}w_g\|_{L_z^{\infty}\dot{B}_{p,1}^0}.
\end{equation*}
For $\epsilon$ sufficiently small, this term can be absorbed into the
left-hand side of \eqref{wgpreabsorb}, which concludes the proof.
\end{proof}


\subsection{Elliptic estimate for \texorpdfstring{$G$, $L$ and $M$}{}}

In this subsection, we establish some basic mapping and frequency localization properties of the truncated elliptic operators $G_{<\mu}$, $L_{<\mu}$ and $M_{<\mu}$. We first establish some relevant elliptic estimates for these operators and error estimates involving their principal and subprincipal symbols. We then show that these operators preserve the bilinear bounds satisfied by the functions they act on and quantify in a sufficiently precise way the high and low frequency tail leakage from these operators. 
\begin{proposition}\label{basicopbounds}
Let $1< p\leq \infty$ and let $\lambda,\mu>0$. Let $A_{<\mu}$ denote either $L_{<\mu}$, $G_{<\mu}$ or $M_{<\mu}$. Let $\mathfrak{a}_{<\mu}$ denote its principal symbol as defined in Subsection \ref{subsec:principal}. The following estimates hold
\begin{enumerate}
\item (Subprincipal error estimate). Let $1\leq r<\infty$ and let $\sigma\geq 0$. There holds
\begin{equation}\label{DNprincerror}
\||D|^{\sigma}(A_{<\mu}-Op(\mathfrak{a}_{<\mu}))f_{\lambda}\|_{\dot{B}_{p,1}^0}\lesssim \max(\mu,\lambda)^{\sigma}\||D|^{2(1+\frac{1}{r})}\eta_{<\mu}\|_{\dot{B}^0_{r,1}}\|f_{\lambda}\|_{L^p}.   
\end{equation}
\item (Linear estimates). Let $\sigma\geq 0$. There holds
\begin{equation}\label{Alinbound}
\||D|^{\sigma}A_{<\mu}f_{\lambda}\|_{\dot{B}_{p,1}^0}\lesssim \max(\lambda,\mu)^{\sigma}\lambda\|f_{\lambda}\|_{L_x^p}. 
\end{equation}
\item (Principal decomposition). There exists a decomposition
\begin{equation}\label{principaldecomposition}
A_{<\mu}f_{\lambda}=\int_{\mathbb{R}^2}K(x,y)f_{\lambda}(y)dy+\mathcal{R}_{\mu}f_{\lambda},
\end{equation}
where the kernel $K$ satisfies
\begin{equation*}
|K(x,y)|\lesssim_N \lambda^3(1+\lambda |x-y|)^{-N},\hspace{5mm} N\in\mathbb{N},
\end{equation*}
and $\mathcal{R}_{\mu}f_{\lambda}$ is a remainder term satisfying
\begin{equation*}
\|\mathcal{R}_{\mu}f_{\lambda}\|_{\dot{B}_{p,1}^0}\lesssim \||D|^{2(1+\frac{1}{r})}\eta_{<\mu}\|_{\dot{B}_{r,1}^0}\|f_{\lambda}\|_{L^p},\hspace{5mm}1\leq r<\infty.     
\end{equation*}
\end{enumerate}
\end{proposition}

Of course, from the embedding $\dot{B}_{p,1}^0\hookrightarrow L^p$, all of the above estimates hold with $\dot{B}_{p,1}^0$ replaced by $L^p$ as well.
\begin{remark}
The last estimate will be particularly useful for proving bilinear estimates for many of the error terms later in our analysis. For instance, it ensures that we have the decomposition
\begin{equation*}
\|A_{<\mu}f_{\lambda_1}g_{\lambda_2}\|_{L_{t,x}^2}\lesssim \lambda_1 \|L(f_{\lambda_1},g_{\lambda_2})\|_{L_{t,x}^2}+\|\mathcal{R}_{\mu}f_{\lambda}g_{\lambda_2}\|_{L_{t,x}^2}    
\end{equation*}
where $f_{\lambda_1}$ and $g_{\lambda_2}$ are functions localized to frequency $\lambda_1$ and $\lambda_2$, respectively and $L$ is a translation invariant bilinear form with rapidly decaying kernel. The first term on the right hand side of the expansion above is amenable to bilinear estimates, while the latter term possesses more favorable decay and balance of derivatives, and will therefore be something that will typically be controlled using linear estimates. This fact will heavily used in the sequel.
\end{remark}

\subsubsection{Subprincipal error bounds} We begin with the subprincipal error bounds, which require the most work. We will first handle the case $A_{<\mu}=G_{<\mu}$, and then deduce the other two bounds as a simple consequence. For this, we follow a crude version of the paralinearization procedure in \cite{AM}. Let $w_f$ denote the flattened harmonic extension of the boundary datum $f_{\lambda}$ relative to the graph $\eta_{<\mu}$, as in Subsection \ref{subsec:flattenedcoords}, with $\eta$ replaced by $\eta_{<\mu}$. It will be convenient to normalize the flattened operator $\Delta_{\eta}$ by the coefficient
\begin{equation*}
\rho\coloneq (1+|\nabla_x\eta_{<\mu}|^2)^{-1}.
\end{equation*}
Then the normalized form of \eqref{variablelaplacian} is
\begin{equation*}
\partial_z^2w_f+\rho\Delta_x w_f-2\rho\nabla_x\eta_{<\mu}\cdot\nabla_x\partial_z w_f-\rho\partial_z w_f\Delta_x\eta_{<\mu}=0.
\end{equation*}
Similarly to \cite{AM}, we then have the factorization
\begin{equation*}
\begin{split}
(\partial_z-T_{a_-})(\partial_z-T_{a_+})w_f
&=(T_{a_-}T_{a_+}-T_{a_-a_+})w_f
+\rho\partial_zw_f\Delta_x\eta_{<\mu}
-R(\rho,\Delta_x w_f)
\\
&\quad
+2R(\rho\nabla_x\eta_{<\mu},\nabla_x \partial_z w_f)
\\
&\coloneq E_{f},
\end{split}
\end{equation*}
where the symbols $a_\pm(x,\xi)$ are given by
\begin{equation}\label{apmdef}
\begin{split}
a_\pm(x,\xi)=i\rho\nabla_x\eta_{<\mu}\cdot\xi
\pm\sqrt{\rho|\xi|^2-(\rho\nabla_x\eta_{<\mu}\cdot\xi)^2}.
\end{split}
\end{equation}
Now, define
\begin{equation*}
b_{\lambda}(z)\coloneq (\partial_z-T_{a_+})w_f.
\end{equation*}
Our first (and main) goal will be to control $|D|^{\sigma}b_{\lambda}(0)$ in $\dot{B}_{p,1}^0$ by the right-hand side of \eqref{DNprincerror}. Our starting point is the following elliptic estimate for $b_{\lambda}$.
\begin{lemma}\label{Besovelliptic}
We have the estimate
\begin{equation*}
\begin{aligned}
\||D|^{\sigma}b_{\lambda}\|_{L_z^{\infty}\dot{B}^0_{p,1}}
&\lesssim
\mu^{\sigma}
\||D|^{2(1+\frac{1}{r})}\eta_{<\mu}\|_{\dot{B}^0_{r,1}}
\int_{-\infty}^{0}
\|\nabla_{x,z}w_f(z)\|_{\dot{B}^0_{p,1}}dz
\\
&\quad+
\||D|^{2(1+\frac{1}{r})}\eta_{<\mu}\|_{\dot{B}^0_{r,1}}
\int_{-\infty}^{0}
\|\nabla_{x,z}|D|^{\sigma}w_f(z)\|_{\dot{B}^0_{p,1}}dz,
\end{aligned}
\end{equation*}
where $1<p\leq \infty$, $\sigma\geq 0$ and $1\leq r<\infty$ are as above.
\end{lemma}

\begin{proof}

For a Schwartz function $g:\mathbb{R}^3\to\mathbb{C}$, let $h\coloneq (\partial_z+|D|)^{-1}g$ denote the solution to the equation
\begin{equation*}
(\partial_z+|D|)h=g
\end{equation*}
with zero boundary data at $-\infty$. That is,
\begin{equation*}
h(z,x)=\int_{-\infty}^{z}(e^{-(z-s)|D|}g(s,\cdot))(x)ds.
\end{equation*}
Using the fact that $b_{\lambda}(z)\to0$ as $z\to-\infty$, we compute
\begin{equation*}
b_{\lambda}(z)
=(\partial_z+|D|)^{-1}(|D|+T_{a_-})b_{\lambda}
+(\partial_z+|D|)^{-1}E_{f}.
\end{equation*}
For $z<0$ and each dyadic number $\alpha\in 2^{\mathbb{Z}}$, the operator $e^{z|D|}P_{\alpha}$ is bounded from $L^p$ to $L^p$ with norm $\lesssim_p e^{c\alpha z}$ uniformly in $z$. Therefore, we have
\begin{equation*}
\begin{split}
\||D|^{\sigma}(\partial_z+|D|)^{-1}(|D|+T_{a_{-}})b_{\lambda}\|_{L_z^{\infty}\dot{B}^0_{p,1}}
&\lesssim
\||D|^{\sigma-1}(|D|+T_{a_{-}})b_{\lambda}\|_{L_z^{\infty}\dot{B}^0_{p,1}}
\\
&\lesssim
\epsilon \||D|^{\sigma}b_{\lambda}\|_{L_z^{\infty}\dot{B}^0_{p,1}},
\end{split}
\end{equation*}
where in the second line, we used
\begin{equation*}
\|\nabla_x\eta_{<\mu}\|_{L_x^{\infty}}\lesssim\epsilon
\end{equation*}
and the fact that $a_-+|\xi|$ is an order one symbol whose normalized order one symbol seminorms satisfy
\begin{equation*}
\sup_{x\in\mathbb{R}^2,\,\xi\neq0}
|\xi|^{-1}|a_-(x,\xi)+|\xi||\lesssim\epsilon,
\end{equation*}
together with the analogous estimates for the required $\xi$ derivatives.

By Minkowski's inequality and the $\dot{B}^0_{p,1}\to \dot{B}^0_{p,1}$ boundedness of $e^{-(z-s)|D|}$ (which holds for $1<p\leq\infty$), we also have the estimate
\begin{equation*}
\||D|^{\sigma}(\partial_z+|D|)^{-1}E_f\|_{L_z^{\infty}\dot{B}^0_{p,1}}
\lesssim
\int_{-\infty}^{0}
\||D|^{\sigma}E_f(z)\|_{\dot{B}^0_{p,1}}dz.
\end{equation*}
Combining the above two estimates, we obtain for $\epsilon$ small enough,
\begin{equation*}
\||D|^{\sigma}b_{\lambda}\|_{L_z^{\infty}\dot{B}^0_{p,1}}
\lesssim
\int_{-\infty}^{0}
\||D|^{\sigma}E_f(z)\|_{\dot{B}^0_{p,1}}dz.
\end{equation*}
To conclude the proof of the lemma, it suffices to show the bound
\begin{equation*}
\begin{split}
\||D|^{\sigma}E_f(z)\|_{\dot{B}^0_{p,1}}
&\lesssim
\mu^{\sigma}
\||D|^{2(1+\frac{1}{r})}\eta_{<\mu}\|_{\dot{B}^0_{r,1}}
\|\nabla_{z,x}w_f(z)\|_{\dot{B}^0_{p,1}}
\\
&\quad+
\||D|^{2(1+\frac{1}{r})}\eta_{<\mu}\|_{\dot{B}^0_{r,1}}
\|\nabla_{z,x}|D|^{\sigma}w_f(z)\|_{\dot{B}^0_{p,1}}.
\end{split}
\end{equation*}
To show this, we go term by term in the definition of $E_f$. By the standard symbolic calculus for paradifferential operators
\cite[Corollary 3.4.H]{taylor1991pseudodifferential}, together with the
Besov-space boundedness in \cite[Chapter I, Proposition 12.1]{taylor2000tools} and the definition of $a_+$ and $a_-$, the paradifferential operator $T_{a_-}T_{a_+}-T_{a_-a_+}$ maps $\dot{B}_{p,1}^{\sigma+1}$ to $\dot{B}_{p,1}^{\sigma}$ with norm controlled by
\begin{equation*}
C(\|\nabla_x\eta_{<\mu}\|_{\dot{B}_{\infty,1}^0})
\|\nabla_x^2\eta_{<\mu}\|_{\dot{B}_{\infty,1}^0}.
\end{equation*}
Therefore, we have
\begin{equation*}
\begin{split}
\||D|^{\sigma}(T_{a_-}T_{a_+}-T_{a_-a_+})w_f\|_{\dot{B}^0_{p,1}}
&\lesssim
\|\nabla^2_x\eta_{<\mu}\|_{\dot{B}_{\infty,1}^0}
\|\nabla_x |D|^{\sigma}w_f\|_{\dot{B}^0_{p,1}}
\\
&\lesssim
\||D|^{2(1+\frac{1}{r})}\eta_{<\mu}\|_{\dot{B}^0_{r,1}}
\|\nabla_x|D|^{\sigma}w_f\|_{\dot{B}^0_{p,1}},
\end{split}
\end{equation*}
where we used Sobolev embedding in the last inequality. Next, we estimate the other terms in $E_{f}$. For instance, for
\begin{equation*}
R(\rho,\Delta_x w_f)=T_{\Delta_x w_f}\rho+\Pi(\Delta_x w_f,\rho),
\end{equation*}
by Proposition \ref{besovparabounds}, Proposition \ref{MoserBesov}, and Bernstein, applied to $\rho-1$, we have
\begin{equation*}
\begin{split}
\||D|^{\sigma}T_{\Delta_x w_f}\rho\|_{\dot{B}^0_{p,1}}
&\lesssim
\||D|^{\sigma+1}(\rho-1)\|_{\dot{B}_{\infty,1}^0}
\|\nabla_x w_f\|_{\dot{B}^0_{p,1}}\lesssim
\mu^{\sigma}
\||D|^{2(1+\frac{1}{r})}\eta_{<\mu}\|_{\dot{B}^0_{r,1}}
\|\nabla_x w_f\|_{\dot{B}^0_{p,1}}.
\end{split}
\end{equation*}
To estimate the high--high para-product, fix a finite exponent
$\widetilde r>\max\{r,p'\}$ and define $q$ by $\frac{1}{q}=\frac{1}{p}+\frac{1}{\widetilde r}.$ By Sobolev embedding, we have
\begin{equation*}
\||D|^{2(1+\frac{1}{\widetilde r})}\eta_{<\mu}\|_{\dot{B}_{\widetilde r,1}^0}
\lesssim
\||D|^{2(1+\frac{1}{r})}\eta_{<\mu}\|_{\dot{B}_{r,1}^0}.
\end{equation*}
Therefore, by Sobolev embedding and again Proposition \ref{MoserBesov} and Proposition \ref{besovparabounds}, applied to $\rho-1$, we have
\begin{equation*}
\begin{split}
\||D|^{\sigma}\Pi(\Delta_x w_f,\rho)\|_{\dot{B}^0_{p,1}}
&\lesssim
\||D|^{\frac{2}{\widetilde r}+\sigma}\Pi(\Delta_x w_f,\rho)\|_{\dot{B}_{q,1}^0}\lesssim
\|\nabla_x w_f\|_{\dot{B}^0_{p,1}}
\||D|^{\sigma+1+\frac{2}{\widetilde r}}(\rho-1)\|_{\dot{B}^0_{\widetilde r,1}}
\\
&\lesssim
\mu^{\sigma}
\|\nabla_x w_f\|_{\dot{B}^0_{p,1}}
\||D|^{2(1+\frac{1}{r})}\eta_{<\mu}\|_{\dot{B}^0_{r,1}}.
\end{split}
\end{equation*}
The remaining terms satisfy, by the same reasoning,
\begin{equation*}
\begin{split}
&\||D|^{\sigma}
(\rho\partial_zw_f\Delta_x\eta_{<\mu})\|_{\dot{B}^0_{p,1}}
+
\||D|^{\sigma}
R(\rho\nabla_x\eta_{<\mu},\nabla_x\partial_z w_f)\|_{\dot{B}^0_{p,1}}
\\
&\lesssim
\mu^{\sigma}
\||D|^{2(1+\frac{1}{r})}\eta_{<\mu}\|_{\dot{B}^0_{r,1}}
\|\nabla_{z,x}w_f\|_{\dot{B}^0_{p,1}}+
\||D|^{2(1+\frac{1}{r})}\eta_{<\mu}\|_{\dot{B}^0_{r,1}}
\|\nabla_{z,x}|D|^{\sigma}w_f\|_{\dot{B}^0_{p,1}}.
\end{split}
\end{equation*}
Combining the above estimates concludes the proof.
\end{proof}

Now we return to the proof of \eqref{DNprincerror}. The last main ingredient is to establish the elliptic estimate
\begin{equation}\label{ellipticw}
\int_{-\infty}^{0}
\||D|^{\beta}\nabla_{x,z}w_f(z)\|_{\dot{B}^0_{p,1}}dz
\lesssim
(\mu+\lambda)^{\beta}\|f_{\lambda}\|_{L^p}
\end{equation}
for any $\beta\geq 0$, which we will then apply with $\beta=0$ and $\beta=\sigma$ on the right hand side of Lemma \ref{Besovelliptic}. To obtain this, we compare $w_f$ to the harmonic extension of $f_{\lambda}$ on the half space $\{z<0\}$. Using the notation from the proof of Proposition \ref{flattenedbasicbounds}, with $\eta$ replaced by $\eta_{<\mu}$, we have
\begin{equation}\label{flatcomparison}
\Delta_{z,x}(w_f-e^{z|D|}f_{\lambda})
=
\nabla_{z,x}\cdot
\left(A_{\eta_{<\mu}}\nabla_{z,x}w_f\right).
\end{equation}
Recall also from the proof of Proposition \ref{flattenedbasicbounds} the half-space operator
\begin{equation*}
\mathcal{R}_0\coloneq
\nabla_{z,x}\Delta_{z,x}^{-1}\nabla_{z,x}\cdot.
\end{equation*}
Since $w_f-e^{z|D|}f_{\lambda}$ has zero boundary value, \eqref{flatcomparison} yields
\begin{equation*}
\nabla_{z,x}w_f
=
(|D|,\nabla_x)e^{z|D|}f_{\lambda}
+
\mathcal{R}_0
\left(A_{\eta_{<\mu}}\nabla_{z,x}w_f\right).
\end{equation*}
The same dyadic kernel bounds used in the proof of Proposition \ref{flattenedbasicbounds}, now integrated in the output variable $z$, give
\begin{equation*}
\|\mathcal{R}_0H\|_{L_z^1\dot{B}_{p,1}^0}
\lesssim
\|H\|_{L_z^1\dot{B}_{p,1}^0}.
\end{equation*}
Therefore, defining
\begin{equation*}
X_{\beta}\coloneq
\int_{-\infty}^{0}
\||D|^{\beta}\nabla_{x,z}w_f(z)\|_{\dot{B}^0_{p,1}}dz,
\end{equation*}
we obtain
\begin{equation}\label{young}
X_{\beta}
\lesssim
\lambda^{\beta+1}\int_{-\infty}^{0}
\|e^{z|D|}f_{\lambda}\|_{L_x^p}dz
+
\int_{-\infty}^{0}
\||D|^{\beta}
(A_{\eta_{<\mu}}\nabla_{z,x}w_f)(z)\|_{\dot{B}^0_{p,1}}dz.
\end{equation}
Since $f_{\lambda}$ is frequency localized, the dyadic Poisson kernel estimate gives
\begin{equation*}
\lambda^{\beta+1}\int_{-\infty}^{0}
\|e^{z|D|}f_{\lambda}\|_{L_x^p}dz
\lesssim
\lambda^{\beta}\|f_{\lambda}\|_{L_x^p}.
\end{equation*}

It remains to estimate the second term on the right hand side of \eqref{young}. Recall that the entries of $A_{\eta_{<\mu}}$ are either linear or quadratic in $\nabla_x\eta_{<\mu}$. To handle the high--high interactions, choose a finite exponent $\widetilde r>\max\{2,p'\}$ and define $q$ by $\frac{1}{q}=\frac{1}{p}+\frac{1}{\widetilde r}.$ By Sobolev embedding, Bernstein's inequality, and Proposition \ref{MoserBesov}, we have
\begin{equation*}
\||D|^{\frac{2}{\widetilde r}+\beta}
A_{\eta_{<\mu}}\|_{\dot{B}_{\widetilde r,1}^0}
\lesssim
\mu^{\beta}
\|\nabla_x^2\eta_{<\mu}\|_{\dot{B}_{2,1}^0}
\lesssim
\epsilon\mu^{\beta}.
\end{equation*}
Using the embedding $\dot{B}^{\frac{2}{\widetilde r}}_{q,1}
\subset
\dot{B}^0_{p,1}$ for the high--high interactions, Proposition \ref{besovparabounds} gives
\begin{equation}\label{Fparabounds}
\begin{split}
\||D|^{\beta}
(A_{\eta_{<\mu}}\nabla_{z,x}w_f)\|_{\dot{B}^0_{p,1}}
&\lesssim
\epsilon\|\nabla_{z,x}|D|^{\beta}w_f\|_{\dot{B}^0_{p,1}}
+
\epsilon\mu^{\beta}
\|\nabla_{z,x}w_f\|_{\dot{B}^0_{p,1}}.
\end{split}
\end{equation}
Combining \eqref{young} and \eqref{Fparabounds}, we obtain
\begin{equation*}
X_{\beta}
\lesssim
\lambda^{\beta}\|f_{\lambda}\|_{L_x^p}
+\epsilon X_{\beta}
+\epsilon\mu^{\beta}X_0.
\end{equation*}
When $\beta=0$, we may absorb the last two terms to obtain
\begin{equation*}
X_0\lesssim\|f_{\lambda}\|_{L_x^p}.
\end{equation*}
Substituting this estimate back into the preceding inequality and absorbing the $X_{\beta}$ term gives
\begin{equation*}
X_{\beta}
\lesssim
(\lambda^{\beta}+\mu^{\beta})\|f_{\lambda}\|_{L_x^p}
\lesssim
(\lambda+\mu)^{\beta}\|f_{\lambda}\|_{L_x^p}.
\end{equation*}
This proves \eqref{ellipticw}. We now establish \eqref{DNprincerror}. First, Lemma \ref{Besovelliptic} and \eqref{ellipticw} give
\begin{equation*}
\||D|^{\sigma}b_{\lambda}(0)\|_{\dot{B}^0_{p,1}}
\lesssim
\max(\mu,\lambda)^{\sigma}
\||D|^{2(1+\frac{1}{r})}\eta_{<\mu}\|_{\dot{B}^0_{r,1}}
\|f_{\lambda}\|_{L^p}.
\end{equation*}
Using the boundary trace identity from Subsection \ref{subsec:flattenedcoords}, we have
\begin{equation*}
\begin{split}
G_{<\mu}f_{\lambda}
&=\rho^{-1}\left.\partial_z w_f\right|_{z=0}
-\nabla_x\eta_{<\mu}\cdot\nabla_x f_{\lambda}=\rho^{-1}b_{\lambda}(0)
+(\rho^{-1}T_{a_+}-\nabla_x\eta_{<\mu}\cdot\nabla_x)f_{\lambda}.
\end{split}
\end{equation*}
The definition of $a_+$ gives the exact symbol identity
\begin{equation*}
\rho^{-1}a_+(x,\xi)
-i\nabla_x\eta_{<\mu}\cdot\xi
=
\mathfrak{g}_{<\mu}(x,\xi).
\end{equation*}
Since we use left quantization, it follows that
\begin{equation*}
\rho^{-1}Op(a_+)f_{\lambda}
-\nabla_x\eta_{<\mu}\cdot\nabla_x f_{\lambda}
=
Op(\mathfrak{g}_{<\mu})f_{\lambda}.
\end{equation*}
Consequently,
\begin{equation*}
(\rho^{-1}T_{a_+}-\nabla_x\eta_{<\mu}\cdot\nabla_x)f_{\lambda}
=
Op(\mathfrak{g}_{<\mu})f_{\lambda}
+\mathcal{F}_f,
\end{equation*}
where
\begin{equation*}
\mathcal{F}_f
\coloneq
\rho^{-1}(T_{a_+}-Op(a_+))f_{\lambda}.
\end{equation*}
Moreover, for each of the required multi-indices $\alpha$, the definition of $a_+$, Sobolev embedding and Proposition \ref{MoserBesov} give
\begin{equation*}
\sup_{\xi\neq0}
|\xi|^{|\alpha|-1}
\|\nabla_x\partial_{\xi}^{\alpha}a_+(\cdot,\xi)\|_{\dot{B}_{\infty,1}^0}
\lesssim_{\alpha}
\||D|^{2(1+\frac{1}{r})}\eta_{<\mu}\|_{\dot{B}_{r,1}^0}.
\end{equation*}
The operator $T_{a_+}-Op(a_+)$ contains only the interactions in which the spatial frequency of the symbol is not much lower than the frequency of the input. Therefore, Bernstein’s inequality and the standard paraproduct estimates, with the high--high interactions treated by the same finite-exponent argument as above, yield
\begin{equation*}
\||D|^{\sigma}\mathcal{F}_f\|_{\dot{B}^0_{p,1}}
\lesssim
\max(\lambda,\mu)^{\sigma}
\||D|^{2(1+\frac{1}{r})}\eta_{<\mu}\|_{\dot{B}^0_{r,1}}
\|f_{\lambda}\|_{L^p}.
\end{equation*}
To estimate the other remainder, write
\begin{equation*}
\rho^{-1}-1=|\nabla_x\eta_{<\mu}|^2.
\end{equation*}
The same product estimate as above, including the  treatment of the high--high interactions, yields
\begin{equation*}
\begin{split}
\||D|^{\sigma}(\rho^{-1}b_{\lambda}(0))\|_{\dot{B}^0_{p,1}}
&\lesssim
\||D|^{\sigma}b_{\lambda}(0)\|_{\dot{B}^0_{p,1}}
+\epsilon\mu^{\sigma}
\|b_{\lambda}(0)\|_{\dot{B}^0_{p,1}}
\\
&\lesssim
\max(\mu,\lambda)^{\sigma}
\||D|^{2(1+\frac{1}{r})}\eta_{<\mu}\|_{\dot{B}^0_{r,1}}
\|f_{\lambda}\|_{L^p},
\end{split}
\end{equation*}
where we used the estimates for $b_{\lambda}(0)$ with exponents $\sigma$ and $0$, respectively. Combining these estimates, we have the Besov bound
\begin{equation*}
\||D|^{\sigma}
(G_{<\mu}-Op(\mathfrak{g}_{<\mu}))f_{\lambda}
\|_{\dot{B}^0_{p,1}}
\lesssim
\max(\mu,\lambda)^{\sigma}
\||D|^{2(1+\frac{1}{r})}\eta_{<\mu}\|_{\dot{B}^0_{r,1}}
\|f_{\lambda}\|_{L^p}.
\end{equation*}
In particular, this proves \eqref{DNprincerror} in the case $A_{<\mu}=G_{<\mu}$.

It remains to deduce the corresponding estimates for $L_{<\mu}$ and $M_{<\mu}$. From their definitions, we have
\begin{equation*}
L_{<\mu}f
=
\rho\left(G_{<\mu}f+\nabla_x\eta_{<\mu}\cdot\nabla_x f\right),
\qquad
M_{<\mu}f
=
\nabla_x f-\nabla_x\eta_{<\mu}L_{<\mu}f.
\end{equation*}
The corresponding principal symbols satisfy
\begin{equation*}
\mathfrak{l}_{<\mu}(x,\xi)
=
\rho\left(
\mathfrak{g}_{<\mu}(x,\xi)
+i\nabla_x\eta_{<\mu}\cdot\xi
\right),
\qquad
\mathfrak{m}_{<\mu}(x,\xi)
=
i\xi-\nabla_x\eta_{<\mu}\mathfrak{l}_{<\mu}(x,\xi).
\end{equation*}
Consequently, for the left quantization used here,
\begin{equation*}
\begin{split}
(L_{<\mu}-Op(\mathfrak{l}_{<\mu}))f_{\lambda}
&=
\rho(G_{<\mu}-Op(\mathfrak{g}_{<\mu}))f_{\lambda},
\\
(M_{<\mu}-Op(\mathfrak{m}_{<\mu}))f_{\lambda}
&=
-\nabla_x\eta_{<\mu}
(L_{<\mu}-Op(\mathfrak{l}_{<\mu}))f_{\lambda}.
\end{split}
\end{equation*}
Applying the same product estimate to these identities, using the estimate for $G_{<\mu}-Op(\mathfrak{g}_{<\mu})$, gives
\begin{equation*}
\begin{split}
&\||D|^{\sigma}
(L_{<\mu}-Op(\mathfrak{l}_{<\mu}))f_{\lambda}
\|_{\dot{B}^0_{p,1}}
+
\||D|^{\sigma}
(M_{<\mu}-Op(\mathfrak{m}_{<\mu}))f_{\lambda}
\|_{\dot{B}^0_{p,1}}
\\
&\lesssim
\max(\mu,\lambda)^{\sigma}
\||D|^{2(1+\frac{1}{r})}\eta_{<\mu}\|_{\dot{B}^0_{r,1}}
\|f_{\lambda}\|_{L^p}.
\end{split}
\end{equation*}
This completes the proof of \eqref{DNprincerror}.

\subsubsection{Linear estimates} Next, we establish \eqref{Alinbound}. We begin with the following lemma.
\begin{lemma}\label{lincprincerrorlemma}
Let $A'_{<\mu}$ denote either $G_{<\mu}-|D|$, $L_{<\mu}-|D|$ or $M_{<\mu}-\nabla$. For every $1<p\leq\infty$ and $\sigma\geq 0$, there holds
\begin{equation}\label{linprincerror}
\||D|^{\sigma}A'_{<\mu}f_{\lambda}\|_{\dot{B}^0_{p,1}}
\lesssim_{p,\sigma}
\epsilon\max(\lambda,\mu)^{\sigma}
\lambda\|f_{\lambda}\|_{L_{x}^{p}}.
\end{equation}
\end{lemma}
\begin{proof}
We first establish the bound for $A'_{<\mu}=G(\eta_{<\mu})-|D|$. Let $w_f$ denote the flattened harmonic extension of $f_{\lambda}$ relative to the graph $\eta_{<\mu}$. Recalling the notation from the proof of Proposition \ref{flattenedbasicbounds}, we recall again
\begin{equation}\label{linearR0representation}
\nabla_{z,x}w_f
=
(|D|,\nabla_x)e^{z|D|}f_{\lambda}
+
\mathcal{R}_0
\left(A_{\eta_{<\mu}}\nabla_{z,x}w_f\right),
\end{equation}
We first record the product estimate that will be used below. Arguing as before, we have the product estimate
\begin{equation}\label{Aetaproductlinear}
\begin{split}
\||D|^{\beta}
(A_{\eta_{<\mu}}H)\|_{\dot{B}^0_{p,1}}
&\lesssim_{p,\beta}
\epsilon\||D|^{\beta}H\|_{\dot{B}^0_{p,1}}
+
\epsilon\mu^{\beta}\|H\|_{\dot{B}^0_{p,1}}
\end{split}
\end{equation}
for every $\beta\geq0$. Now, for $\beta\geq0$, define
\begin{equation*}
Y_{p,\beta}(w_f)
\coloneq
\||D|^{\beta}\nabla_{z,x}w_f\|_{L_z^{\infty}\dot{B}^0_{p,1}}.
\end{equation*}
The kernel bounds for $\mathcal{R}_0$ used in the proof of Proposition \ref{flattenedbasicbounds}, together with \eqref{linearR0representation} and \eqref{Aetaproductlinear}, give
\begin{equation*}
Y_{p,\beta}(w_f)
\lesssim_{p,\beta}
\lambda^{\beta+1}\|f_{\lambda}\|_{L_x^p}
+
\epsilon Y_{p,\beta}(w_f)
+
\epsilon\mu^{\beta}Y_{p,0}(w_f).
\end{equation*}
Taking first $\beta=0$ and choosing $\epsilon$ sufficiently small, we obtain
\begin{equation*}
Y_{p,0}(w_f)\lesssim_p\lambda\|f_{\lambda}\|_{L_x^p}.
\end{equation*}
Returning to the estimate for general $\beta\geq0$ and absorbing the second term on the right, we conclude that
\begin{equation}\label{linearYbound}
Y_{p,\beta}(w_f)
\lesssim_{p,\beta}
\lambda\max(\lambda,\mu)^{\beta}
\|f_{\lambda}\|_{L_x^p}.
\end{equation}

The same kernel bounds apply to the boundary trace of the vertical component of $\mathcal{R}_0$. Taking the vertical component of \eqref{linearR0representation} at $z=0$, we obtain
\begin{equation*}
(\partial_zw_f)_{|z=0}-|D|f_{\lambda}
=
\left.
\left[
\mathcal{R}_0
\left(A_{\eta_{<\mu}}\nabla_{z,x}w_f\right)
\right]_z
\right|_{z=0}.
\end{equation*}
Therefore, by \eqref{Aetaproductlinear} and \eqref{linearYbound},
\begin{equation}\label{betaerror}
\begin{split}
\||D|^{\beta}
\left((\partial_zw_f)_{|z=0}-|D|f_{\lambda}\right)
\|_{\dot{B}^0_{p,1}}
&\lesssim_{p,\beta}
\epsilon Y_{p,\beta}(w_f)
+
\epsilon\mu^{\beta}Y_{p,0}(w_f)
\\
&\lesssim_{p,\beta}
\epsilon\lambda\max(\lambda,\mu)^{\beta}
\|f_{\lambda}\|_{L_x^p}.
\end{split}
\end{equation}
The boundary trace identity from Subsection \ref{subsec:flattenedcoords} gives
\begin{equation*}
G(\eta_{<\mu})f_{\lambda}
=
(1+|\nabla_x\eta_{<\mu}|^2)(\partial_zw_f)_{|z=0}
-\nabla_x\eta_{<\mu}\cdot\nabla_xf_{\lambda},
\end{equation*}
and hence
\begin{equation*}
\begin{split}
(G(\eta_{<\mu})-|D|)f_{\lambda}
&=
(1+|\nabla_x\eta_{<\mu}|^2)
\left((\partial_zw_f)_{|z=0}-|D|f_{\lambda}\right)
\\
&\quad
+|\nabla_x\eta_{<\mu}|^2|D|f_{\lambda}
-\nabla_x\eta_{<\mu}\cdot\nabla_xf_{\lambda}.
\end{split}
\end{equation*}
Applying \eqref{betaerror}, together with the same paraproduct estimates needed to establish \eqref{Aetaproductlinear}, we obtain
\begin{equation*}
\||D|^{\sigma}
(G(\eta_{<\mu})-|D|)f_{\lambda}
\|_{\dot{B}^0_{p,1}}
\lesssim_{p,\sigma}
\epsilon\max(\lambda,\mu)^{\sigma}
\lambda\|f_{\lambda}\|_{L_x^p}.
\end{equation*}
This proves \eqref{linprincerror} for $G_{<\mu}-|D|$. It remains to deduce the corresponding estimates for $L_{<\mu}-|D|$ and $M_{<\mu}-\nabla$. From the definitions of $L_{<\mu}$ and $M_{<\mu}$, we have the identities
\begin{equation*}
\begin{split}
(L_{<\mu}-|D|)f_{\lambda}
&=
\rho(G_{<\mu}-|D|)f_{\lambda}
+\rho \nabla_x\eta_{<\mu}\cdot\nabla_xf_{\lambda}
+(\rho-1)|D|f_{\lambda},
\\
(M_{<\mu}-\nabla_x)f_{\lambda}
&=
-\nabla_x\eta_{<\mu}L_{<\mu}f_{\lambda}.
\end{split}
\end{equation*}
Proposition \ref{MoserBesov}, and Proposition \ref{besovparabounds}, as applied above allow us to deduce the corresponding bounds for $L_{<\mu}-|D|$, and then for $M_{<\mu}-\nabla_x$. This concludes the proof of the lemma.
\end{proof}

Finally, since $f_{\lambda}$ is frequency localized, we have
\begin{equation*}
\||D|^{\sigma+1}f_{\lambda}\|_{\dot{B}^0_{p,1}}
+
\||D|^{\sigma}\nabla f_{\lambda}\|_{\dot{B}^0_{p,1}}
\lesssim_{p,\sigma}
\lambda^{\sigma+1}\|f_{\lambda}\|_{L_x^p}
\lesssim_{p,\sigma}
\max(\lambda,\mu)^{\sigma}
\lambda\|f_{\lambda}\|_{L_x^p}.
\end{equation*}
Combining this with Lemma \ref{lincprincerrorlemma} proves \eqref{Alinbound}.
\subsubsection{Principal decomposition}
Now, we turn to the principal decomposition in Proposition \ref{basicopbounds}. We first show the details for $A_{<\mu}=G(\eta_{<\mu})$. We write
\begin{equation*}
G(\eta_{<\mu})f_{\lambda}
=
Op(\mathfrak{g}_{<\mu})f_{\lambda}
+
(G(\eta_{<\mu})-Op(\mathfrak{g}_{<\mu}))f_{\lambda}.
\end{equation*}
We can write
\begin{equation*}
Op(\mathfrak{g}_{<\mu})f_{\lambda}
=
\int_{\mathbb{R}^2}
K_{\lambda,\mu}(x,y)f_{\lambda}(y)dy,
\end{equation*}
where the kernel $K_{\lambda,\mu}$ is given by
\begin{equation*}
K_{\lambda,\mu}(x,y)
=
\frac{1}{(2\pi)^2}
\int_{\mathbb{R}^2}
e^{i(x-y)\cdot\xi}
\mathfrak{g}_{<\mu}(x,\xi)
\chi_0(\lambda^{-1}\xi)d\xi.
\end{equation*}
Here, $\chi_0\in C_c^\infty(\mathbb{R}^2\setminus\{0\})$ is a fixed smooth cutoff chosen so that $\chi_0(\lambda^{-1}\xi)=1$ on the support of $\widehat{f_{\lambda}}$ and supported in a slightly larger annulus. By the degree one homogeneity and smoothness of $\mathfrak{g}_{<\mu}(x,\xi)$ for $\xi\neq0$, together with the bound
\begin{equation*}
\|\nabla_x\eta_{<\mu}\|_{L^\infty}\lesssim\epsilon,
\end{equation*}
we have, for every multi-index $\alpha$,
\begin{equation*}
|\partial_{\xi}^{\alpha}
(\mathfrak{g}_{<\mu}(x,\xi)\chi_0(\lambda^{-1}\xi))|
\lesssim_{\alpha}
\lambda^{1-|\alpha|}.
\end{equation*}
Since the integrand is supported where $|\xi|\approx\lambda$, integration by parts gives
\begin{equation}\label{pointwisekern}
|K_{\lambda,\mu}(x,y)|
\lesssim_N
\lambda^3(1+\lambda|x-y|)^{-N},
\hspace{5mm}
N\in\mathbb{N}.
\end{equation}
To control the subprincipal contribution, define
\begin{equation*}
\mathcal{R}_{\mu}f_{\lambda}
\coloneq
(G(\eta_{<\mu})-Op(\mathfrak{g}_{<\mu}))f_{\lambda}.
\end{equation*}
By \eqref{DNprincerror} with $\sigma=0$, we have
\begin{equation*}
\|\mathcal{R}_{\mu}f_{\lambda}\|_{\dot{B}^0_{p,1}}
\lesssim
\||D|^{2(1+\frac{1}{r})}\eta_{<\mu}\|_{\dot{B}_{r,1}^0}
\|f_{\lambda}\|_{L^p},
\hspace{5mm}
1\leq r<\infty.
\end{equation*}

The cases $A_{<\mu}=L_{<\mu}$ and $A_{<\mu}=M_{<\mu}$ follow in the same way. Indeed, the corresponding symbols $\mathfrak{l}_{<\mu}$ and $\mathfrak{m}_{<\mu}$ satisfy the same order one bounds in the $\xi$ variable, while their subprincipal remainders are controlled by \eqref{DNprincerror}. This completes the proof of \eqref{principaldecomposition} and Proposition \ref{basicopbounds}.

A useful corollary of the third part of Proposition \ref{basicopbounds} is that $A_{<\mu}u_{\lambda}$ inherits the unbalanced bilinear bounds satisfied by $|D|u_{\lambda}$.
\begin{corollary}\label{Abilbounds}
Let $A_{<\mu}$ denote either $G(\eta_{<\mu})$, $L_{<\mu}$ or $M_{<\mu}$. Assume $\mu\lesssim\lambda$ and that either $\alpha\ll\lambda$ or $\lambda\ll\alpha$. There holds
\begin{equation}\label{Abilbound}
\|A_{<\mu}u_{\lambda}\overline{u}^{x_0}_{\alpha}\|_{L^2_{t,x}}
\lesssim
c_{\lambda}c_{\alpha}\epsilon^2
\lambda^{1-s_c}\alpha^{-s_c}
\min(\lambda,\alpha)^{\frac{3}{4}}
\left(1+\min(\lambda,\alpha)|x_0|\right)^{\frac12}.
\end{equation}
\end{corollary}
\begin{proof}
We begin by estimating the contribution of $\mathcal{R}_{\mu}u_{\lambda}$ in \eqref{principaldecomposition}. From the bootstrap hypotheses in Section \ref{sec:bootstrap}, Bernstein and dyadic summation we have
\begin{equation*}
\||D|^{\frac{17}{8}}\eta_{<\mu}\|_{L_t^4\dot{B}_{16,1}^0}
\lesssim
\epsilon\mu^{\frac{7}{8}}.
\end{equation*}
Now continuing with \eqref{Abilbound}, first suppose that $\alpha\ll\lambda$. Using the remainder estimate in Proposition \ref{basicopbounds} with $r=16$ and $p=\frac{16}{7}$ as well as the bootstrap hypotheses in Section \ref{sec:bootstrap}, we obtain
\begin{equation*}
\begin{split}
\|\mathcal{R}_{\mu}u_{\lambda}\overline{u}^{x_0}_{\alpha}\|_{L^2_{t,x}}
&\lesssim
\|u_{\alpha}\|_{L_t^4L_x^{16}}
\|\mathcal{R}_{\mu}u_{\lambda}\|_{L_t^4\dot{B}^{0}_{\frac{16}{7},1}}\lesssim
\|u_{\alpha}\|_{L_t^4L_x^{16}}
\||D|^{\frac{17}{8}}\eta_{<\mu}\|_{L_t^4\dot{B}_{16,1}^0}
\|u_{\lambda}\|_{L_t^{\infty}L_x^{\frac{16}{7}}}
\\
&\lesssim
\epsilon^3c_{\lambda}c_{\alpha}
\alpha^{\frac{3}{4}-s_c}
\mu^{\frac{7}{8}}
\lambda^{\frac{1}{8}-s_c}
\\
&\lesssim
\epsilon^3c_{\lambda}c_{\alpha}
\alpha^{\frac{3}{4}-s_c}
\lambda^{1-s_c}.
\end{split}
\end{equation*}
On the other hand, suppose that $\lambda\ll\alpha$. Using the remainder estimate with $r=16$ and $p=\infty$, we have
\begin{equation*}
\begin{split}
\|\mathcal{R}_{\mu}u_{\lambda}\overline{u}^{x_0}_{\alpha}\|_{L^2_{t,x}}&\lesssim
\|u_{\alpha}\|_{L_t^{\infty}L_x^2}
\|\mathcal{R}_{\mu}u_{\lambda}\|_{L_t^2\dot{B}^{0}_{\infty,1}}\lesssim
\|u_{\alpha}\|_{L_t^{\infty}L_x^2}
\||D|^{\frac{17}{8}}\eta_{<\mu}\|_{L_t^4\dot{B}_{16,1}^0}
\|u_{\lambda}\|_{L_t^4L_x^{\infty}}
\\
&\lesssim
\epsilon^3c_{\lambda}c_{\alpha}
\alpha^{-s_c}
\mu^{\frac{7}{8}}
\lambda^{\frac{7}{8}-s_c}
\\
&\lesssim
\epsilon^3c_{\lambda}c_{\alpha}
\alpha^{-s_c}
\lambda^{\frac{7}{4}-s_c}.
\end{split}
\end{equation*}
Since $\epsilon\ll1$, these estimates give the desired bound for the remainder contribution. It remains to estimate the principal part. Using the kernel bound in \eqref{principaldecomposition}, we obtain
\begin{equation*}
\begin{split}
\|(A_{<\mu}-\mathcal{R}_{\mu})u_{\lambda}
\overline{u}_{\alpha}^{x_0}\|_{L_{t,x}^2}
&\lesssim_N
\lambda
\sup_{\theta,\tau\in\mathbb{R}^2}
\langle\theta\rangle^{-N}
\langle\tau\rangle^{-N}
\|u_{\lambda}(\cdot+\lambda^{-1}\theta)
\overline{u}_{\alpha}^{x_0}(\cdot+\alpha^{-1}\tau)\|_{L_{t,x}^2}
\\
&\lesssim
c_{\lambda}c_{\alpha}\epsilon^2
\lambda^{1-s_c}\alpha^{-s_c}
\min(\lambda,\alpha)^{\frac{3}{4}}
\left(1+\min(\lambda,\alpha)|x_0|\right)^{\frac12},
\end{split}
\end{equation*}
where in the second inequality we applied the appropriate unbalanced bilinear estimate, interchanging the two factors when $\lambda\ll\alpha$ and used the rapid decay of the kernel to absorb the additional translation factors. Combining the principal and remainder estimates concludes the proof.
\end{proof}
A slightly inconvenient technical issue is that the operator $A_{<\mu}$ does not preserve the frequency localization of the functions on which it acts. Fortunately, the resulting high-frequency tails can be controlled, under mild conditions, by arbitrarily large powers of $\nabla_x^2\eta_{<\mu}$. This will allow us to treat these tails perturbatively in virtually every application.
\begin{proposition}\label{Bernsteinelliptic}
Let $\sigma\geq 0$, and let $1<p,q\leq\infty$ and $1<r\leq\infty$ satisfy $\frac{1}{p}=\frac{1}{q}+\frac{1}{r}.$ Let $\lambda,\mu>0$ with $\mu\lesssim\lambda$, and let $A_{<\mu}$ denote either $G(\eta_{<\mu})$, $L_{<\mu}$ or $M_{<\mu}$. For every integer $N\geq 1$, there exist constants $M=M(N)>1$ and $0<\delta_0=\delta_0(N)\ll 1$ such that the following estimates hold:
\begin{enumerate}
\item For every $0\leq\sigma_1\leq\min\{\sigma,1\}$, there holds
\begin{equation}\label{bernelliptic}
\begin{split}
\||D|^{\sigma}P_{>M\lambda}A_{<\mu}f_{<\lambda}\|_{\dot{B}^0_{p,1}}
&\lesssim_{M,N}
\lambda^{\sigma-\sigma_1-N}
\|\nabla_x^2\eta_{<\mu}\|_{\dot{B}^{0}_{rN,1}}^N
\||D|^{\sigma_1}L_{<\mu}f_{<\lambda}\|_{\dot{B}^0_{q,1}}.
\end{split}
\end{equation}
\item For every $0<\delta\leq\delta_0$, there holds
\begin{equation}\label{bernelliptic2}
\begin{split}
\|P_{<\delta\lambda}A_{<\delta\lambda}f_{\lambda}\|_{\dot{B}_{p,1}^0}
&\lesssim_{\delta_0,N}
\lambda^{1-N}
\|\nabla_x^2\eta_{<\delta\lambda}\|_{\dot{B}_{rN,1}^0}^N
\|f_{\lambda}\|_{L^q}.
\end{split}
\end{equation}
\end{enumerate}
The implicit constants may also depend on $p,q,r$ and $\sigma$.
\end{proposition}
\begin{remark}
The fact that the term on the right hand side of the first bound involves $L_{<\mu}f_{<\lambda}$ instead of $|D|f_{<\lambda}$ will be of critical importance for many estimates. For instance, the truncated trace $B_{<\lambda}$ is of the form $L_{<\lambda}\psi_{<\lambda}$ and is more regular that $|D|\psi_{<\lambda}$. This will allow us to safely control high-tail leakage for the horizontal and vertical traces $V_{<\lambda}$ and $B_{<\lambda}$.    
\end{remark}
\begin{proof}
We first establish \eqref{bernelliptic}. We begin with the case $A_{<\mu}=L_{<\mu}$. As usual, let $w_f$ denote the flattened harmonic extension of the boundary datum $f_{<\lambda}$ relative to the graph $\eta_{<\mu}$, as in Subsection \ref{subsec:flattenedcoords}. Define
\begin{equation*}
W_f\coloneq\partial_zw_f,
\qquad
B_f\coloneq W_f(0)=L_{<\mu}f_{<\lambda}.
\end{equation*}
Since the coefficients of $\Delta_{\eta_{<\mu}}$ are independent of $z$, the function $W_f$ satisfies the same flattened elliptic equation as $w_f$. Recalling the matrix $A_{\eta_{<\mu}}$ and the operator $\mathcal{R}_0$ from the proof of Proposition \ref{flattenedbasicbounds}, we have
\begin{equation}\label{WfR0rep}
\Delta_{z,x}W_f
=
\nabla_{z,x}\cdot
\left(A_{\eta_{<\mu}}\nabla_{z,x}W_f\right)
\end{equation}
and
\begin{equation}\label{WfR0rep2}
\nabla_{z,x}W_f
=
(|D|,\nabla_x)e^{z|D|}B_f
+
\mathcal{R}_0
\left(A_{\eta_{<\mu}}\nabla_{z,x}W_f\right).
\end{equation}
For a vector field $F=(F_z,F_x)$ on the half-space, define
\begin{equation*}
\mathcal{S}_0F
\coloneq
\int_{-\infty}^{0}
e^{z|D|}
\left(
\nabla_x\cdot F_x-|D|F_z
\right)(z)dz.
\end{equation*}
The vertical component of the boundary trace of $\mathcal{R}_0$ is given by
\begin{equation}\label{R0boundarytrace}
\left.[\mathcal{R}_0F]_z\right|_{z=0}
=F_z(0)+\mathcal{S}_0F.
\end{equation}
Set
\begin{equation*}
F_f\coloneq A_{\eta_{<\mu}}\nabla_{z,x}W_f.
\end{equation*}
In particular, we have
\begin{equation*}
(F_f)_z
=
-|\nabla_x\eta_{<\mu}|^2\partial_zW_f
+
\nabla_x\eta_{<\mu}\cdot\nabla_xW_f.
\end{equation*}
On the other hand, evaluating the flattened equation for $w_f$ at $z=0$ gives
\begin{equation*}
\begin{split}
(1+|\nabla_x\eta_{<\mu}|^2)\partial_zW_f(0)
&=
-\Delta_xf_{<\lambda}
+2\nabla_x\eta_{<\mu}\cdot\nabla_xB_f
+(\Delta_x\eta_{<\mu})B_f.
\end{split}
\end{equation*}
Consequently,
\begin{equation*}
\partial_zW_f(0)-(F_f)_z(0)
=
-\Delta_xf_{<\lambda}
+\nabla_x\eta_{<\mu}\cdot\nabla_xB_f
+(\Delta_x\eta_{<\mu})B_f.
\end{equation*}
Taking the vertical component of \eqref{WfR0rep2} at $z=0$ and using \eqref{R0boundarytrace}, we obtain the boundary formula
\begin{equation}\label{DBrep}
\begin{split}
|D|B_f
&=
-\Delta_xf_{<\lambda}
+\nabla_x\eta_{<\mu}\cdot\nabla_xB_f
+(\Delta_x\eta_{<\mu})B_f
-\mathcal{S}_0F_f.
\end{split}
\end{equation}
Next, we may assume without loss of generality that $\sigma\geq1$. Indeed, if $0\leq\sigma<1$, then
\begin{equation*}
\||D|^{\sigma}P_{>M\lambda}B_f\|_{\dot{B}_{p,1}^0}
\lesssim_M
\lambda^{\sigma-1}
\||D|P_{>M\lambda}B_f\|_{\dot{B}_{p,1}^0},
\end{equation*}
and we can apply the estimate with $\sigma$ replaced by $1$ and with the same $\sigma_1$. Next, define intermediate exponents $p_j$, $0\leq j\leq N$, by $\frac{1}{p_j}
=
\frac{1}{q}+\frac{j}{rN}.$ Thus $p_0=q$, $p_N=p$, and $\frac{1}{p_j}
=
\frac{1}{p_{j-1}}+\frac{1}{rN}.$ Now, choose constants
\begin{equation*}
1\ll M_1\ll M_2\ll\cdots\ll M_N=M
\end{equation*}
with the separations sufficiently large depending on $N$. For each $j$, let $\widetilde P_{M_j\lambda}$ denote a smooth projection to a fixed annular neighborhood of the transition region of $P_{>M_j\lambda}$. We choose the constants so that, for $2\leq j\leq N$,
\begin{equation*}
\operatorname{supp}\widetilde P_{M_j\lambda}
\subset
\{\xi:|\xi|>M_{j-1}\lambda\}.
\end{equation*}
We allow the projections $\widetilde P_{M_j\lambda}$ to be harmlessly fattened from line to line.

We next record the key one-step estimate that will be then be iterated $N$ times. For $1\leq j\leq N$, there holds
\begin{equation}\label{bernellipticonestep}
\begin{split}
\||D|^{\sigma}P_{>M_j\lambda}B_f\|_{\dot{B}_{p_j,1}^0}
&\lesssim_{M_j}
\lambda^{-1}
\|\nabla_x^2\eta_{<\mu}\|_{\dot{B}_{rN,1}^0}
\||D|^{\sigma}\widetilde P_{M_j\lambda}B_f\|_{\dot{B}_{p_{j-1},1}^0}.
\end{split}
\end{equation}
To see this, we first note that the dyadic kernel bounds from Proposition \ref{flattenedbasicbounds}, applied to \eqref{WfR0rep2}, together with the integrated estimate \eqref{ellipticw}, give
\begin{equation}\label{localizedWfestimate}
\begin{split}
&\||D|^{\sigma-1}P_{>M_j\lambda}\nabla_{z,x}W_f
\|_{L_z^{\infty}\dot{B}_{p_j,1}^0}
+
\||D|^{\sigma-1}P_{>M_j\lambda}\partial_z\nabla_{z,x}W_f
\|_{L_z^1\dot{B}_{p_j,1}^0}
\\
&\lesssim
\||D|^{\sigma}P_{>M_j\lambda}B_f\|_{\dot{B}_{p_j,1}^0}
+
\lambda^{-1}
\|\nabla_x^2\eta_{<\mu}\|_{\dot{B}_{rN,1}^0}
\||D|^{\sigma}\widetilde P_{M_j\lambda}B_f\|_{\dot{B}_{p_{j-1},1}^0}.
\end{split}
\end{equation}
Here the $L_z^{\infty}$ estimate is the localized version of Proposition \ref{flattenedbasicbounds}, while the $L_z^1$ estimate follows from the corresponding integrated elliptic estimate used above, or equivalently from \eqref{ellipticw} applied to $W_f$. In commuting $P_{>M_j\lambda}$ through the coefficients, Proposition \ref{Leibnizcommutator} gives a factor $(M_j\lambda)^{-1}$ and places the remaining factor in the transition region $|\xi|\approx M_j\lambda$; we bounded this factor by $\lambda^{-1}$ in \eqref{localizedWfestimate}. The terms in which the projection remains on $W_f$ are absorbed using the smallness bootstrap
\begin{equation*}
\|\nabla_x\eta_{<\mu}\|_{\dot{B}_{\infty,1}^0}
\lesssim\epsilon.
\end{equation*}
Next, we observe that at each dyadic frequency, the kernel of $\mathcal{S}_0$ satisfies
\begin{equation*}
\||D|^{\beta}\mathcal{S}_0H\|_{\dot{B}_{a,1}^0}
\lesssim
\||D|^{\beta}H\|_{L_z^{\infty}\dot{B}_{a,1}^0},
\qquad \beta\geq0,
\end{equation*}
since the derivative in the definition of $\mathcal{S}_0$ is cancelled by integration of the Poisson kernel in $z$. Applying $|D|^{\sigma-1}P_{>M_j\lambda}$ to \eqref{DBrep}, using Proposition \ref{Leibnizcommutator}, Proposition \ref{besovparabounds}, and \eqref{localizedWfestimate}, we obtain
\begin{equation*}
\begin{split}
\||D|^{\sigma}P_{>M_j\lambda}B_f\|_{\dot{B}_{p_j,1}^0}
&\lesssim
\epsilon
\||D|^{\sigma}P_{>M_j\lambda}B_f\|_{\dot{B}_{p_j,1}^0}+
\lambda^{-1}
\|\nabla_x^2\eta_{<\mu}\|_{\dot{B}_{rN,1}^0}
\||D|^{\sigma}\widetilde P_{M_j\lambda}B_f\|_{\dot{B}_{p_{j-1},1}^0}.
\end{split}
\end{equation*}
Indeed, all the low-high terms contain an undifferentiated factor of $\nabla_x\eta_{<\mu}$ and are included in the first term on the right. Every transition term either comes from a commutator, or contains the factor $\Delta_x\eta_{<\mu}$; in both cases it is controlled by the second term on the right. Choosing $\epsilon$ sufficiently small and absorbing the first term proves \eqref{bernellipticonestep}.

We now iterate \eqref{bernellipticonestep}, starting with $j=N$. By the support properties of the transition projections, we obtain
\begin{equation*}
\begin{split}
\||D|^{\sigma}P_{>M\lambda}B_f\|_{\dot{B}_{p,1}^0}
&\lesssim_{M,N}
\lambda^{-N}
\|\nabla_x^2\eta_{<\mu}\|_{\dot{B}_{rN,1}^0}^{N}
\||D|^{\sigma}\widetilde P_{M_1\lambda}B_f\|_{\dot{B}_{q,1}^0}
\\
&\lesssim_{M,N}
\lambda^{\sigma-\sigma_1-N}
\|\nabla_x^2\eta_{<\mu}\|_{\dot{B}_{rN,1}^0}^{N}
\||D|^{\sigma_1}B_f\|_{\dot{B}_{q,1}^0},
\end{split}
\end{equation*}
where the last inequality follows from Bernstein on the annulus $|\xi|\approx M_1\lambda$. This proves \eqref{bernelliptic} for $A_{<\mu}=L_{<\mu}$.

The proof for $G_{<\mu}$ and $M_{<\mu}$ follows from the bound for $L_{<\mu}$ by using the following formulas
\begin{equation*}
G_{<\mu}f
=
(1+|\nabla_x\eta_{<\mu}|^2)L_{<\mu}f
-\nabla_x\eta_{<\mu}\cdot\nabla_xf,\hspace{5mm}M_{<\mu}f
=
\nabla_xf-\nabla_x\eta_{<\mu}L_{<\mu}f.
\end{equation*}
For instance, in the formula for $G_{<\mu}f$, applying the projection $P_{>M\lambda}$ for $M$ large enough, we obtain
\begin{equation}\label{highfreqG}
P_{>M\lambda}G_{<\mu}f=P_{>M\lambda}\left((1+|\nabla_x\eta_{<\mu}|^2)\tilde{P}_{>M\lambda}L_{<\mu}f\right)    
\end{equation}
where $\tilde{P}_{>M\lambda}$ is a slightly fattened projection. The latter term satisfies $P_{>M\lambda}(\nabla_x\eta_{<\mu}\cdot\nabla_xf)=0$ because of the localization of $\eta_{<\mu}$ and $f$. The estimate for $P_{>M\lambda}G_{<\mu}f$ can be reduced to the corresponding estimate for $L_{<\mu}f$ by using Bernstein and Proposition \ref{besovparabounds} (noting that only the low-high interaction remains due to the localization of $\eta_{<\mu}$). The analysis of $M_{<\mu}f$ is similar. 
\medskip

We now prove \eqref{bernelliptic2}. Again, we first consider the case $A_{<\delta\lambda}=L_{<\delta\lambda}$. We use the same notation as above, with $w_f$ now denoting the flattened harmonic extension of $f_{\lambda}$ relative to the graph $\eta_{<\delta\lambda}$ and with
\begin{equation*}
W_f\coloneq\partial_zw_f.
\end{equation*}
The $L_z^{\infty}$ estimate from Proposition \ref{flattenedbasicbounds}, the linear estimate \eqref{Alinbound}, and the integrated estimate \eqref{ellipticw}, applied to $W_f$, give
\begin{equation}\label{lowtailuniform}
\begin{split}
\|\nabla_{z,x}w_f\|_{L_z^{\infty}\dot{B}_{q,1}^0}
+
\|\partial_z\nabla_{z,x}w_f\|_{L_z^1\dot{B}_{q,1}^0}
\lesssim_q
\lambda\|f_{\lambda}\|_{L^q}.
\end{split}
\end{equation}
Indeed, the second term is the integrated gradient estimate for the flattened harmonic extension $W_f$, whose boundary value satisfies
\begin{equation*}
\|W_f(0)\|_{\dot{B}_{q,1}^0}
=
\|L_{<\delta\lambda}f_{\lambda}\|_{\dot{B}_{q,1}^0}
\lesssim
\lambda\|f_{\lambda}\|_{L^q}.
\end{equation*}
Recalling the notation from the proof of Proposition \ref{flattenedbasicbounds}, we have
\begin{equation}\label{lowtailR0rep}
\nabla_{z,x}w_f
=
(|D|,\nabla_x)e^{z|D|}f_{\lambda}
+
\mathcal{R}_0
\left(A_{\eta_{<\delta\lambda}}\nabla_{z,x}w_f\right).
\end{equation}
The entries of $A_{\eta_{<\delta\lambda}}$ are linear or quadratic functions of $\nabla_x\eta_{<\delta\lambda}$. Consequently,
\begin{equation}\label{Aetalowbounds}
\|A_{\eta_{<\delta\lambda}}\|_{\dot{B}_{\infty,1}^0}
\lesssim\epsilon,
\qquad
\|\nabla_xA_{\eta_{<\delta\lambda}}\|_{\dot{B}_{rN,1}^0}
\lesssim
\|\nabla_x^2\eta_{<\delta\lambda}\|_{\dot{B}_{rN,1}^0}.
\end{equation}
Now, choose fixed constants
\begin{equation*}
0<\delta_0<\delta_1<\cdots<\delta_N<\delta_{N+1}\ll1
\end{equation*}
so that $P_{<\delta_{j+1}\lambda}$ is identically one on a fixed neighborhood of the support of $P_{<\delta_j\lambda}$ and
\begin{equation*}
P_{<\delta_N\lambda}f_{\lambda}=0.
\end{equation*}
We choose all the constants depending only on $N$, with $\delta_{N+1}$ still sufficiently small relative to the frequency support of $f_{\lambda}$. We retain the intermediate exponents $p_j$ introduced in the high-frequency argument, so that
\begin{equation*}
\frac{1}{p_j}
=
\frac{1}{q}+\frac{j}{rN},
\qquad
p_0=q,
\qquad
p_N=p.
\end{equation*}
Applying $P_{<\delta_j\lambda}$ to \eqref{lowtailR0rep}, the flat term vanishes. We can then commute $A_{\eta_{<\delta\lambda}}$ through $P_{<\delta_j\lambda}$ and argue exactly as in the derivation of \eqref{localizedWfestimate}. The support properties of the nested projections, Proposition \ref{Leibnizcommutator}, and \eqref{Aetalowbounds} give
\begin{equation}\label{lowtailonestep}
\begin{split}
&\|P_{<\delta_j\lambda}\nabla_{z,x}w_f\|_{L_z^{\infty}\dot{B}_{p_{N-j},1}^0}
+
\|\partial_zP_{<\delta_j\lambda}\nabla_{z,x}w_f\|_{L_z^1\dot{B}_{p_{N-j},1}^0}
\\
&\lesssim_{\delta_0,N}
\lambda^{-1}
\|\nabla_x^2\eta_{<\delta\lambda}\|_{\dot{B}_{rN,1}^0}
\left(
\|P_{<\delta_{j+1}\lambda}\nabla_{z,x}w_f\|_{L_z^{\infty}\dot{B}_{p_{N-j-1},1}^0}+
\|\partial_zP_{<\delta_{j+1}\lambda}\nabla_{z,x}w_f\|_{L_z^1\dot{B}_{p_{N-j-1},1}^0}
\right).
\end{split}
\end{equation}
Iterating \eqref{lowtailonestep} for $j=1,\ldots,N$ and then using \eqref{lowtailuniform}, we obtain
\begin{equation*}
\begin{split}
\|P_{<\delta_1\lambda}\nabla_{z,x}w_f\|_{L_z^{\infty}\dot{B}_{p,1}^0}
+
\|\partial_zP_{<\delta_1\lambda}\nabla_{z,x}w_f\|_{L_z^1\dot{B}_{p,1}^0}\lesssim_{\delta_0,N}
\lambda^{1-N}
\|\nabla_x^2\eta_{<\delta\lambda}\|_{\dot{B}_{rN,1}^0}^N
\|f_{\lambda}\|_{L^q}.
\end{split}
\end{equation*}
Since $P_{<\delta\lambda}P_{<\delta_1\lambda}=P_{<\delta\lambda}$, the boundary trace identity gives
\begin{equation*}
\begin{split}
\|P_{<\delta\lambda}L_{<\delta\lambda}f_{\lambda}\|_{\dot{B}_{p,1}^0}=
\|P_{<\delta\lambda}W_f(0)\|_{\dot{B}_{p,1}^0}\lesssim_{\delta_0,N}
\lambda^{1-N}
\|\nabla_x^2\eta_{<\delta\lambda}\|_{\dot{B}_{rN,1}^0}^N
\|f_{\lambda}\|_{L^q}.
\end{split}
\end{equation*}
This proves \eqref{bernelliptic2} for $L_{<\delta\lambda}$.

The estimates for $G_{<\delta\lambda}$ and $M_{<\delta\lambda}$ now follow from the same identities used in the high-frequency argument. Indeed, after decreasing $\delta_0$ by a fixed universal factor, frequency separation gives
\begin{equation*}
P_{<\delta\lambda}\nabla_xf_{\lambda}=0,
\qquad
P_{<\delta\lambda}
\left(
\nabla_x\eta_{<\delta\lambda}\cdot\nabla_xf_{\lambda}
\right)=0.
\end{equation*}
Since $\nabla_x\eta_{<\delta\lambda}$ and $|\nabla_x\eta_{<\delta\lambda}|^2$ are localized to frequencies $\lesssim\delta\lambda$, the standard para-product estimates give, for a fixed universal constant $C>1$,
\begin{equation*}
\begin{split}
\|P_{<\delta\lambda}M_{<\delta\lambda}f_{\lambda}\|_{\dot{B}_{p,1}^0}
+
\|P_{<\delta\lambda}G_{<\delta\lambda}f_{\lambda}\|_{\dot{B}_{p,1}^0}\lesssim
\|P_{<C\delta\lambda}L_{<\delta\lambda}f_{\lambda}\|_{\dot{B}_{p,1}^0}.
\end{split}
\end{equation*}
The preceding estimate for $L_{<\delta\lambda}$ applies with $P_{<C\delta\lambda}$ in place of $P_{<\delta\lambda}$. By applying this, we conclude the proof.
\end{proof}
A relatively simple consequence of the above estimates is the following Bernstein type estimate which we will use repeatedly.
\begin{corollary}\label{Bestcor}
Let $\lambda,\mu\in 2^{\mathbb{Z}}$ with $\mu\lesssim\lambda$,
and let $1<p\leq\infty$, $1<q<p$, $\sigma>0$ and
$0<\sigma_1\leq\min(1,\sigma)$.
If $Z_{<\lambda}=f_{<\lambda}$ or
$Z_{<\lambda}=L_{<\lambda}f_{<\lambda}$, then
\begin{equation}
\||D|^{\sigma}A_{<\mu}Z_{<\lambda}\|_{\dot{B}^0_{p,1}}
\lesssim_{\sigma,p,q}
\lambda^{\sigma-\sigma_1+1+\frac{2}{q}-\frac{2}{p}}
\||D|^{\sigma_1}Z_{<\lambda}\|_{L^q}.
\end{equation}
If $Z_{<\lambda}=M_{<\lambda}f_{<\lambda}$ or
$Z_{<\lambda}=G_{<\lambda}f_{<\lambda}$, then
\begin{equation}
\begin{split}
\||D|^{\sigma}A_{<\mu}Z_{<\lambda}\|_{\dot{B}^0_{p,1}}
\lesssim_{\sigma,p,q}
\lambda^{\sigma-\sigma_1+1+\frac{2}{q}-\frac{2}{p}}
\left(
\||D|^{\sigma_1}Z_{<\lambda}\|_{L^q}
+
\||D|^{\sigma_1}L_{<\lambda}f_{<\lambda}\|_{L^q}
\right).
\end{split}
\end{equation}
\end{corollary}
\begin{proof}
In the case $Z_{<\lambda}=f_{<\lambda}$, Proposition
\ref{basicopbounds} and Bernstein's inequality give
\begin{equation*}
\begin{split}
\||D|^{\sigma}A_{<\mu}f_{<\lambda}\|_{\dot{B}^0_{p,1}}\lesssim
\lambda^{\sigma}
\||D|f_{<\lambda}\|_{\dot{B}^0_{p,1}}\lesssim
\lambda^{\sigma-\sigma_1+1+\frac{2}{q}-\frac{2}{p}}
\||D|^{\sigma_1}f_{<\lambda}\|_{L^q},
\end{split}
\end{equation*}
where we used $q<p$.

For the remaining cases, it suffices to prove the estimate
\begin{equation}\label{weakerbound}
\begin{split}
\||D|^{\sigma}A_{<\mu}Z_{<\lambda}\|_{L^p}
\lesssim_{\sigma,p,q}
\lambda^{\sigma-\sigma_1+1+\frac{2}{q}-\frac{2}{p}}
\left(
\||D|^{\sigma_1}Z_{<\lambda}\|_{L^q}
+
\||D|^{\sigma_1}L_{<\lambda}f_{<\lambda}\|_{L^q}
\right).
\end{split}
\end{equation}
Indeed, for any integer $N\geq1$, Bernstein's inequality gives
\begin{equation*}
\begin{split}
\||D|^{\sigma}A_{<\mu}Z_{<\lambda}\|_{\dot{B}^0_{p,1}}\lesssim_N
\lambda^{-N}
\||D|^{\sigma+N}A_{<\mu}Z_{<\lambda}\|_{L^p}+
\||D|^{\sigma}P_{<\lambda}
A_{<\mu}Z_{<\lambda}\|_{\dot{B}^0_{p,1}}.
\end{split}
\end{equation*}
The first term is controlled by \eqref{weakerbound} with
$\sigma$ replaced by $\sigma+N$. For the second term, choose
$q<r<p$ and use
\begin{equation*}
\||D|^{\sigma}P_{<\lambda}
A_{<\mu}Z_{<\lambda}\|_{\dot{B}^0_{p,1}}
\lesssim
\lambda^{\frac{2}{r}-\frac{2}{p}}
\||D|^{\sigma}A_{<\mu}Z_{<\lambda}\|_{L^r}.
\end{equation*}
Applying \eqref{weakerbound} with $p$ replaced by $r$ gives
the desired bound.

To prove \eqref{weakerbound}, we first establish a Bernstein estimate for
$L_{<\lambda}f_{<\lambda}$. Proposition
\ref{Bernsteinelliptic}, with $r=\infty$ and
$N=1$, and $\sigma=\sigma_1=1$ gives,
for some sufficiently large fixed $M>1$,
\begin{equation*}
\begin{split}
\||D|L_{<\lambda}f_{<\lambda}\|_{\dot{B}^0_{p,1}}\lesssim
\||D|P_{<M\lambda}
L_{<\lambda}f_{<\lambda}\|_{\dot{B}^0_{p,1}}+
\lambda^{-1}
\|\nabla_x^2\eta_{<\lambda}\|_{\dot{B}^0_{\infty,1}}
\||D|L_{<\lambda}f_{<\lambda}\|_{\dot{B}^0_{p,1}}.
\end{split}
\end{equation*}
By Bernstein and the smallness of $\eta$ in $\dot{B}_{2,1}^1$ we have
\begin{equation*}
\lambda^{-1}
\|\nabla_x^2\eta_{<\lambda}\|_{\dot{B}^0_{\infty,1}}
\lesssim\epsilon.
\end{equation*}
Moreover, since $q<p$ and $\sigma_1\leq1$,
\begin{equation*}
\||D|P_{<M\lambda}
L_{<\lambda}f_{<\lambda}\|_{\dot{B}^0_{p,1}}
\lesssim
\lambda^{1-\sigma_1+\frac{2}{q}-\frac{2}{p}}
\||D|^{\sigma_1}L_{<\lambda}f_{<\lambda}\|_{L^q}.
\end{equation*}
Thus, taking $\epsilon$ sufficiently small and absorbing the
high-frequency term, we obtain
\begin{equation*}
\||D|L_{<\lambda}f_{<\lambda}\|_{\dot{B}^0_{p,1}}
\lesssim
\lambda^{1-\sigma_1+\frac{2}{q}-\frac{2}{p}}
\||D|^{\sigma_1}L_{<\lambda}f_{<\lambda}\|_{L^q}.
\end{equation*}

For each of the three operator choices of $Z_{<\lambda}$,
Proposition \ref{Bernsteinelliptic} similarly gives
\begin{equation*}
\begin{split}
\||D|Z_{<\lambda}\|_{\dot{B}^0_{p,1}}\lesssim
\||D|P_{<M\lambda}Z_{<\lambda}\|_{\dot{B}^0_{p,1}}+
\lambda^{-1}
\|\nabla_x^2\eta_{<\lambda}\|_{\dot{B}^0_{\infty,1}}
\||D|L_{<\lambda}f_{<\lambda}\|_{\dot{B}^0_{p,1}}.
\end{split}
\end{equation*}
Combining Bernstein's inequality with the preceding estimate
therefore yields
\begin{equation}\label{Zbernstein}
\begin{split}
\||D|Z_{<\lambda}\|_{\dot{B}^0_{p,1}}
\lesssim
\lambda^{1-\sigma_1+\frac{2}{q}-\frac{2}{p}}
\left(
\||D|^{\sigma_1}Z_{<\lambda}\|_{L^q}
+
\||D|^{\sigma_1}L_{<\lambda}f_{<\lambda}\|_{L^q}
\right).
\end{split}
\end{equation}

It remains to establish
\begin{equation}\label{simplerDBbound}
\begin{split}
\||D|^{\sigma}A_{<\mu}Z_{<\lambda}\|_{L^p}
\lesssim_{\sigma,p}
\lambda^{\sigma}
\left(
\||D|Z_{<\lambda}\|_{\dot{B}^0_{p,1}}
+
\||D|L_{<\lambda}f_{<\lambda}\|_{\dot{B}^0_{p,1}}
\right).
\end{split}
\end{equation}
Decomposing $Z_{<\lambda}$ into dyadic pieces and applying
Proposition \ref{basicopbounds}, we have
\begin{equation*}
\begin{split}
\||D|^{\sigma}A_{<\mu}Z_{<\lambda}\|_{L^p}
&\lesssim
\sum_{\alpha\leq M\lambda}
\lambda^{\sigma}\||D|P_\alpha Z_{<\lambda}\|_{L^p}
+
\sum_{\alpha>M\lambda}
\alpha^{\sigma}\||D|P_\alpha Z_{<\lambda}\|_{L^p}
\\
&\lesssim
\lambda^{\sigma}\||D|Z_{<\lambda}\|_{\dot{B}^0_{p,1}}
+
\||D|^{\sigma+1}P_{>M\lambda}
Z_{<\lambda}\|_{\dot{B}^0_{p,1}},
\end{split}
\end{equation*}
Here we used
$\mu\lesssim\lambda$ for the low-frequency contribution.
Applying Proposition \ref{Bernsteinelliptic} to the operator
defining $Z_{<\lambda}$, with $r=\infty$,
$N=1$, and derivative parameters $\sigma+1$ and $\sigma_1=1$, gives
\begin{equation*}
\begin{split}
\||D|^{\sigma+1}P_{>M\lambda}
Z_{<\lambda}\|_{\dot{B}^0_{p,1}}
&\lesssim
\lambda^{\sigma-1}
\|\nabla_x^2\eta_{<\lambda}\|_{\dot{B}^0_{\infty,1}}
\||D|L_{<\lambda}f_{<\lambda}\|_{\dot{B}^0_{p,1}}
\\
&\lesssim
\epsilon\lambda^{\sigma}
\||D|L_{<\lambda}f_{<\lambda}\|_{\dot{B}^0_{p,1}}.
\end{split}
\end{equation*}
This proves \eqref{simplerDBbound}. Combining it with
\eqref{Zbernstein} and the preceding Bernstein estimate for
$L_{<\lambda}f_{<\lambda}$ proves \eqref{weakerbound}.
When $Z_{<\lambda}=L_{<\lambda}f_{<\lambda}$, the two norms
on its right-hand side coincide, giving the first claimed
estimate. This concludes the proof.
\end{proof}

\subsection{A more precise paradifferential expansion of \texorpdfstring{$T_A$}{}} Let $A$ denote either $G$, $L$ or $M$. In our precise paradifferential expansion of the water-waves equation later, the paradifferential operator $T_A$ will appear quite frequently in our analysis. However, this is algebraically quite inconvenient to work with, especially for later when we have to perform normal form type corrections to the density-flux identities for the equation. Whenever possible, one would like to replace $T_A$ with the corresponding principal paradifferential operator $T_{\mathfrak{a}}$. In many cases, for this to be viable, we must isolate enough of the precise structure of the subprincipal part to still enable us to prove bilinear estimates for the error. For this purpose we have the following proposition.
\begin{proposition}\label{paraexp}
For every $L,M,N\in\mathbb{N}$ with $L\geq1$, there is a positive constant
$\delta_0=\delta_0(L,M,N)\ll1$ such that, for every
$0<\delta\leq\delta_0$ and every $\lambda>0$, there
holds
\begin{equation}\label{Gexpprecise}
A_{<\delta\lambda}^lf_{\lambda}
=
Op(\mathfrak{a}^l_{<\delta\lambda})f_{\lambda}
+
A^{0}_{<\delta\lambda}f_{\lambda}
+\mathcal{E}_f.
\end{equation}
Here $A^{0}_{<\delta\lambda}$ is a zeroth order operator with kernel
$K^0_{\delta\lambda}$ satisfying the bound
\begin{equation}\label{K0bound}
\begin{split}
|K^0_{\delta\lambda}(x,y)|\lesssim_{M,N}
\lambda^2
(1+\lambda|x-y|)^{-N}
\int_{\mathbb{R}^2}
\langle\theta\rangle^{-M}
\left|
(\nabla_x^2\eta_{<\delta\lambda})
(y+\lambda^{-1}\theta)
\right|d\theta.
\end{split}
\end{equation}
Moreover, for every $1<p\leq\infty$, the error term satisfies
\begin{equation}\label{paraexperror}
\|\mathcal{E}_f\|_{L^p}
\lesssim_L
\lambda^{1-L}
\|\nabla_x^2\eta_{<\delta\lambda}\|_{\dot{B}_{\infty,1}^0}^{L}
\|f_\lambda\|_{L^p}.
\end{equation}
\end{proposition}
\begin{proof}
We will show the details for the case $A_{<\delta\lambda}^l=G_{<\delta\lambda}^l$. The other cases follow as a simple corollary, by the usual relation between $G$ and the other boundary operators.
\\

Fix $L,M,N\in\mathbb{N}$, $\lambda>0$ and
$0<\delta\leq\delta_0$. We also suppress the time variable from
the notation. For simplicity, we also set
\begin{equation*}
q\coloneq\nabla_x\eta_{<\delta\lambda},
\qquad
\rho\coloneq(1+|q|^2)^{-1}.
\end{equation*}
We also abbreviate the integrated factor in \eqref{K0bound} bound by
\begin{equation*}
\mathcal A_{\lambda,M}(y)
\coloneq
\int_{\mathbb R^2}
\langle\theta\rangle^{-M}
\left|
\nabla_x^2\eta_{<\delta\lambda}
(y+\lambda^{-1}\theta)
\right|d\theta.
\end{equation*}
Now, let $w_f$ denote the flattened harmonic extension of $f_\lambda$ 
relative to the graph $\eta_{<\delta\lambda}$. As in the factorization in the proof of Proposition \ref{basicopbounds}, we define the symbols 
\begin{equation*}
a_\pm(x,\xi)
\coloneq
i\rho q\cdot\xi
\pm
\sqrt{
\rho|\xi|^2-(\rho q\cdot\xi)^2
}
\end{equation*}
and write $A_\pm\coloneq Op(a_\pm)$ to denote the corresponding (left quantized) pseudo-differential operators. As in Proposition \ref{basicopbounds}, we have the factorization
\begin{equation}\label{exactflattenedfactorization}
\begin{split}
(\partial_z-A_-)(\partial_z-A_+)w_f=\left(
A_-A_+-Op(a_-a_+)
\right)w_f+
\rho(\Delta_x\eta_{<\delta\lambda})\partial_zw_f.
\end{split}
\end{equation}
Here, we simply use the full operators $Op(a_\pm)$, due to the fact that $\eta_{<\delta\lambda}$ is already localized to lower frequency than $f_{\lambda}$. 
\\

Next, we move to the kernel representation. We begin by collecting a few preliminary estimates and identities. First, we note that since $\eta_{<\delta\lambda}$ is localized to frequency $\lesssim\lambda$, we can write
\begin{equation*}
\nabla_x q(x)
=
\int_{\mathbb R^2}
\chi(\theta)\nabla_x q(x-\lambda^{-1}\theta)d\theta
\end{equation*}
for a Schwartz function $\chi$. A linear change of variables and using the rapid decay of $\chi$, we obtain, for $K$ sufficiently
large relative to $M$,
\begin{equation*}
\begin{split}
|\nabla_x q(x)|
&\lesssim_K
\int_{\mathbb R^2}
\langle\lambda(x-y)-\vartheta\rangle^{-K}
|\nabla_x q(y+\lambda^{-1}\vartheta)|d\vartheta\lesssim_M
\langle\lambda(x-y)\rangle^M
\mathcal A_{\lambda,M}(y).
\end{split}
\end{equation*}
More generally, differentiating gives
\begin{equation}\label{localizedgeometrybound}
|\partial_x^\beta q(x)|
\lesssim_{\beta,M}
\lambda^{|\beta|-1}
\langle\lambda(x-y)\rangle^M
\mathcal A_{\lambda,M}(y),
\qquad
|\beta|\geq1.
\end{equation}
The same estimate holds for a spatial derivative of a smooth function
of $q$, by the chain rule and the smallness bound for $q$ in $L^{\infty}_{t,x}$. We will also later occasionally need the simple bound
\begin{equation}\label{Aweighttransfer}
\mathcal A_{\lambda,M_1}(x)
\lesssim_{M,M_1}
\langle\lambda(x-y)\rangle^M
\mathcal A_{\lambda,M}(y),
\qquad M_1\geq M.
\end{equation}
which follows directly from a linear change of variables and the
elementary inequality
$\langle x\rangle\lesssim\langle x-y\rangle\langle y\rangle$. We now establish a resolvent decomposition for the flattened harmonic
extension that we will need below
\begin{lemma}\label{wdiagdecomp}
Let $\widetilde{P}_{\lambda}$ denote a fattened Littlewood-Paley projection adapted to the annulus $|\xi|\approx\lambda$. Define the operators
\begin{equation*}
H_{\lambda}f\coloneq (|D|,\nabla_x)e^{z|D|}f,\hspace{5mm}
\mathcal{T}F\coloneq\mathcal{R}_0\left(A_{\eta_{<\delta\lambda}}F\right),
\hspace{5mm}
\mathcal{T}_{\lambda}F\coloneq\widetilde{P}_{\lambda}\mathcal{T}F,
\end{equation*}
where $\mathcal{R}_0$ and $A_{\eta_{<\delta\lambda}}$ are defined as in the proof of Proposition \ref{flattenedbasicbounds}. For every integer $L\geq1$, if $\delta_0$ is sufficiently small, then $\widetilde{P}_{\lambda}$ can be chosen, with its width depending on $L$, so that
\begin{equation}\label{wdiagdecomposition}
\nabla_{z,x}w_f
=
W_f+\mathcal{E}^{\mathrm{har}}_f,
\qquad
W_f\coloneq
(I-\mathcal{T}_{\lambda})^{-1}H_{\lambda}f_{\lambda}.
\end{equation}
For every $1<p\leq\infty$, the error term satisfies
\begin{equation}\label{wdiagerror}
\begin{split}
\|\mathcal{E}^{\mathrm{har}}_f\|_{L_z^\infty\dot{B}_{p,1}^0}
+
\lambda
\|\mathcal{E}^{\mathrm{har}}_f\|_{L_z^1\dot{B}_{p,1}^0}
\lesssim_L
\lambda^{1-L}
\|\nabla_x^2\eta_{<\delta\lambda}\|_{\dot{B}_{\infty,1}^0}^{L}
\|f_{\lambda}\|_{L^p}.
\end{split}
\end{equation}
Moreover, the operator $f_{\lambda}\mapsto W_f$ has a kernel satisfying
\begin{equation}\label{paraexpharmonickernel}
\int_{-\infty}^0
|K_{W_f}(z;x,y)|dz
\lesssim_L
\lambda^2\langle\lambda(x-y)\rangle^{-L}.
\end{equation}
\end{lemma}
\begin{proof}
We start with the representation formula
\begin{equation}\label{paraexpR0identity}
\nabla_{z,x}w_f
=
H_\lambda f_\lambda
+
\mathcal{T}\nabla_{z,x}w_f.
\end{equation}
Consequently,
\begin{equation*}
\nabla_{z,x}w_f
=
(I-\mathcal{T})^{-1}H_\lambda f_\lambda.
\end{equation*}
Now set $Q_\lambda\coloneq I-\widetilde{P}_\lambda$, so that $\mathcal{T}-\mathcal{T}_{\lambda}=Q_{\lambda}\mathcal{T}$. The resolvent identity then gives
\begin{equation}\label{paraexpresolventidentity}
\begin{split}
\nabla_{z,x}w_f=
(I-\mathcal{T}_\lambda)^{-1}H_\lambda f_\lambda+
(I-\mathcal{T}_\lambda)^{-1}
Q_\lambda\mathcal{T}
(I-\mathcal{T})^{-1}H_\lambda f_\lambda.
\end{split}
\end{equation}
We define the second term on the right to be
$\mathcal{E}^{\mathrm{har}}_f$. To estimate $\mathcal{E}_{f}^{\mathrm{har}}$, it will be convenient to define the auxiliary norm
\begin{equation}\label{auxxnorm}
\|F\|_{X_{p,\lambda}}
\coloneq
\|F\|_{L_z^\infty\dot{B}_{p,1}^0}
+
\lambda\|F\|_{L_z^1\dot{B}_{p,1}^0}.
\end{equation}
Boundedness of  $\mathcal{R}_0$, Proposition
\ref{besovparabounds}, and the smallness of $\eta$ give
\begin{equation*}
\|\mathcal{T}F\|_{X_{p,\lambda}}
+
\|\mathcal{T}_\lambda F\|_{X_{p,\lambda}}
\lesssim
\epsilon\|F\|_{X_{p,\lambda}}.
\end{equation*}
In particular, both inverse operators in
\eqref{paraexpresolventidentity} are given by convergent Neumann
series. Now, choose nested annular projections
\begin{equation*}
P_\lambda^L\prec P_\lambda^{L-1}\prec\cdots
\prec P_\lambda^0=\widetilde{P}_\lambda
\end{equation*}
such that $P_\lambda^L$ is equal to one on the Fourier support of
$f_\lambda$, and set
\begin{equation*}
Q_\lambda^j\coloneq I-P_\lambda^j,
\qquad 0\leq j\leq L.
\end{equation*}
Thus $Q_\lambda^0=Q_\lambda$ and
$Q_\lambda^Lf_\lambda=0$. If $\delta_0$ is sufficiently small
relative to $L$, Proposition \ref{Leibnizcommutator} and Proposition
\ref{besovparabounds}, iterated through these nested projections,
give
\begin{equation}\label{nestedtailiteration}
\begin{split}
\|Q_\lambda\mathcal{T}^kH_\lambda f_\lambda\|_{X_{p,\lambda}}
&\lesssim_L
\mathbf{1}_{k\geq L}\binom{k}{L}
(C\epsilon)^{k-L}
\lambda^{1-L}
\|\nabla_x^2\eta_{<\delta\lambda}\|_{\dot{B}_{\infty,1}^0}^{L}
\|f_\lambda\|_{L^p}+
(C\epsilon)^k\lambda
\|Q_\lambda^Lf_\lambda\|_{L^p}.
\end{split}
\end{equation}
This follows by simply iterating the commutator bound
\begin{equation*}
\|Q_{\lambda}^j\mathcal{T}g\|_{X_{p,\lambda}}\lesssim \lambda^{-1}\|\nabla^2_x\eta_{<\delta\lambda}\|_{\dot{B}_{\infty,1}^0}\|Q_{\lambda}^{j+1}g\|_{X_{p,\lambda}}+\epsilon\|Q_{\lambda}^jg\|_{X_{p,\lambda}}    
\end{equation*}
$k$ times, where one stops the iteration for any terms involving $L$ commutations of different factors of $A_{\eta<\delta\lambda}$ with one of the projectors $Q_{\lambda}^j$, and continues otherwise, until the scheme terminates at the final input $H_{\lambda}f_{\lambda}$ (in which case we can use Bernstein to replace $Q_{\lambda}^j$ with $Q_{\lambda}^L$). The $k-L$ uncommuted factors of $A_{\eta}$ contribute 
$(C\epsilon)^{k-L}$ in the first term on the right-hand side, while counting all the possible locations of these commutations gives the binomial coefficient in \eqref{nestedtailiteration}. The final
term in \eqref{nestedtailiteration} vanishes by the localization of
$f_\lambda$. Combining this with the boundedness of
$(I-\mathcal{T}_\lambda)^{-1}$ proves \eqref{wdiagerror}.
\\

It remains to prove the kernel bound, which is now straightforward thanks to the localization in $\mathcal{T}_{\lambda}$. First, the flat harmonic extension
satisfies
\begin{equation*}
|K_{H_\lambda}(z;x,y)|
\lesssim_R
\lambda^3e^{c\lambda z}
\langle\lambda(x-y)\rangle^{-R}.
\end{equation*}
The kernel bounds for $\mathcal{R}_0$, together
with $\|A_{\eta_{<\delta\lambda}}\|_{\dot{B}_{\infty,1}^0}
\lesssim\epsilon$ then give
\begin{equation*}
\int_{-\infty}^0
|K_{\mathcal{T}_\lambda^kH_\lambda}(z;x,y)|dz
\lesssim_L
(1+k)^{L+2}(C\epsilon)^k
\lambda^2
\langle\lambda(x-y)\rangle^{-L}.
\end{equation*}
The factor $\lambda^2$ comes from integrating the flat kernel
$\lambda^3e^{c\lambda z}$ in $z$. Summing the Neumann series
\begin{equation*}
W_f
=
\sum_{k=0}^\infty
\mathcal{T}_\lambda^kH_\lambda f_\lambda
\end{equation*}
proves \eqref{paraexpharmonickernel} and concludes the proof of the
lemma.
\end{proof} 
We now return to \eqref{exactflattenedfactorization}. As in Proposition  \ref{basicopbounds}, we define
\begin{equation*}
b\coloneq(\partial_z-A_+)w_f,
\end{equation*}
so that
\begin{equation}\label{bequationparaexp}
(\partial_z-A_-)b
=
\left(
A_-A_+-Op(a_-a_+)
\right)w_f
+
\rho(\Delta_x\eta_{<\delta\lambda})\partial_zw_f.
\end{equation}
We abbreviate the right-hand side of \eqref{bequationparaexp} by
$F$. By definition of $H_\lambda f_\lambda$ and $\mathcal R_0$, the localized term $W_f$ from Lemma
\ref{wdiagdecomp} may be written as
\begin{equation}\label{compositiongradientform}
W_f=\nabla_{z,x}w_f^0,
\qquad
w_f^0\big|_{z=0}=f_\lambda.
\end{equation}
Set $\mathcal C\coloneq A_-A_+-Op(a_-a_+)$. and choose four fattened annular projections (with the same convention as above)
\begin{equation*}
\widetilde P_\lambda\prec P_\lambda^1\prec P_\lambda^2
\prec P_\lambda^3\prec P_\lambda^4
\end{equation*}
and define $\mathcal C_\lambda^0
\coloneq
A_-P_\lambda^1A_+-Op(a_-a_+).$ On the Fourier support of $w_f^0$, the (left-quantized) symbol of
$\mathcal C_\lambda^0$ is
\begin{equation}\label{compositionerrorsymbol}
\begin{split}
c_\lambda^0(x,\xi)
&=
\frac{1}{(2\pi)^2}
\int_{\mathbb R^2}\int_{\mathbb R^2}
e^{-i\theta\cdot\zeta}
a_-(x,\xi+\zeta)\chi_1(\lambda^{-1}(\xi+\zeta))\left(
a_+(x+\theta,\xi)-a_+(x,\xi)
\right)d\theta d\zeta,
\end{split}
\end{equation}
where $\chi_1$ is the unit-scale symbol defining $P_\lambda^1$. Using
\begin{equation*}
a_+(x+\theta,\xi)-a_+(x,\xi)
=
\theta_j\int_0^1
\partial_{x_j}a_+(x+s\theta,\xi)ds
\end{equation*}
and integrating by parts in $\zeta$, together with the symbol bounds for $a_\pm$ and
\eqref{localizedgeometrybound}, the kernel for $\mathcal{C}_{\lambda}^0\tilde{P}_{\lambda}$ satisfies
\begin{equation}\label{compositiongeometrykernel}
\begin{split}
|K_{\mathcal{C}_{\lambda}^0\tilde{P}_{\lambda}}(x,y)|\lesssim_{L,M,N}
\lambda^3
\langle\lambda(x-y)\rangle^{-N}\int_{\mathbb R^2}
\langle\theta\rangle^{-M}
\left|
\nabla_x^2\eta_{<\delta\lambda}
(y+\lambda^{-1}\theta)
\right|d\theta.
\end{split}
\end{equation}
We now define the localized source from the right hand side of \eqref{bequationparaexp} using the fattened projection $P_{\lambda}^2$.
\begin{equation}\label{F0paraexp}
F^0
\coloneq
P_\lambda^2
\left(
\mathcal C_\lambda^0w_f^0
+
\rho(\Delta_x\eta_{<\delta\lambda})(W_f)_z
\right).
\end{equation}
Combining \eqref{compositiongradientform},
\eqref{compositiongeometrykernel}, and
\eqref{paraexpharmonickernel}, and the inequality $\langle x\rangle\lesssim \langle x-y\rangle\langle y\rangle$ gives
\begin{equation}\label{F0kernelparaexp}
\int_{-\infty}^0
|K_{F^0}(z;x,y)|dz
\lesssim_{L,M,N}
\lambda^2
\langle\lambda(x-y)\rangle^{-N}
\mathcal A_{\lambda,M}(y)
\end{equation}
where by $K_{F^0}$, we mean the kernel associated to the map $f\mapsto F^0$.
The error between the localized source and the full source term $F$ on the right hand side of \eqref{bequationparaexp}, $F^{\mathrm{err}}\coloneq F-F^0$ has the decomposition
\begin{equation}\label{Ferrorparaexp}
\begin{split}
F^{\mathrm{err}}
&=
\mathcal C(w_f-w_f^0)
+
\rho(\Delta_x\eta_{<\delta\lambda})
\left(
\partial_zw_f-(W_f)_z
\right)
\\
&\quad+
(I-P_\lambda^2)
\left(
\mathcal C_\lambda^0w_f^0
+
\rho(\Delta_x\eta_{<\delta\lambda})(W_f)_z
\right)
\\
&\quad+
A_-(I-P_\lambda^1)A_+w_f^0.
\end{split}
\end{equation}
We will estimate this error similarly to $\mathcal{E}_f^{\mathrm{har}}$ in the previous lemma. We define the flat half operator $\mathcal{S}$ and auxilliary norm $Y_{B,p}$,
\begin{equation*}
\mathcal S\coloneq(\partial_z+|D|)^{-1},
\qquad
\|h\|_{Y_{B,p}}
\coloneq
\|h\|_{L_z^\infty\dot B^0_{p,1}}
+
\||D|h\|_{L_z^1\dot B^0_{p,1}}.
\end{equation*}
The same Besov source estimate used in the proof of Lemma
\ref{Besovelliptic}, together with \eqref{wdiagerror}, controls the
first line of \eqref{Ferrorparaexp}. More precisely,
\begin{equation*}
\left\|
\mathcal S\left(
\mathcal C(w_f-w_f^0)+
\rho(\Delta_x\eta_{<\delta\lambda})
\partial_z(w_f-w_f^0)
\right)
\right\|_{Y_{B,p}}
\lesssim
\epsilon
\|\nabla_{z,x}(w_f-w_f^0)\|_{X_{p,\lambda}}.
\end{equation*}
The term $\nabla_{z,x}(w_f-w_f^0)$ is precisely $\mathcal{E}_{f}^{\mathrm{har}}$ from the previous lemma and satisfies the estimate \eqref{wdiagerror}. For the remaining two lines, we use the same principle as in the argument in Lemma \ref{wdiagdecomp}. The basic idea is that the nontrivial parts of both of these terms (in view of the localization of $w_f^0$) must involve a frequency projection $P_{\gtrsim\lambda}$ applied to a smooth function of $\nabla\eta_{<\delta\lambda}$, which can be estimated in terms of a large power of $\nabla^2\eta_{<\delta\lambda}$. More precisely, for a smooth function $H$, repeated applications of Bernstein and then Proposition \ref{MoserBesov} and Proposition \ref{besovparabounds} gives (provided $\delta$ is small enough)
\begin{equation}\label{nonlineartail}
\|P_{\gtrsim\lambda}H(\nabla\eta_{<\delta\lambda})\|_{\dot{B}_{\infty,1}^0}\lesssim_{L,H} \lambda^{-L}\|\nabla^2\eta_{<\delta\lambda}\|_{\dot{B}_{\infty,1}^0}^L    
\end{equation}
Now, returning to \eqref{Ferrorparaexp}, we note that for the first term in the second line, the localization of $w_f^0$ and the localization the output force require the symbol for $c_{\lambda}^0$ to be localized to frequency $\gtrsim\lambda$. Using \eqref{nonlineartail}, it is easy to obtain the symbol bound
\begin{equation}\label{sourcesymboltails}
\begin{split}
&\sup_{|\xi|\approx\lambda}|\xi|^{|\alpha|-1}
\left(\left\|
\partial_\xi^\alpha
P^x_{\gtrsim c\lambda}c_{\lambda}^0(\cdot,\xi)
\right\|_{\dot B^0_{\infty,1}}+
\lambda^{-1}\left\|
\nabla_x\partial_\xi^\alpha
P^x_{\gtrsim c\lambda}c_{\lambda}^0(\cdot,\xi)
\right\|_{\dot B^0_{\infty,1}}\right)
\\
&\hspace{5mm}\lesssim_{\alpha,L} \lambda^{1-L}
\|\nabla_x^2\eta_{<\delta\lambda}\|_{
\dot B^0_{\infty,1}}^L.
\end{split}
\end{equation}
Similarly, for the second term in the second line, we have
\begin{equation*}
\left\|
P^x_{\gtrsim \lambda}
\left(
\rho\Delta_x\eta_{<\delta\lambda}
\right)
\right\|_{\dot B^0_{\infty,1}}
\lesssim_{L}
\lambda^{1-L}
\|\nabla_x^2\eta_{<\delta\lambda}\|_{
\dot B^0_{\infty,1}}^L
\end{equation*}
Using these bounds, we can control the $\mathcal{S}$ applied to the second line of \eqref{Ferrorparaexp} in $Y_{B,p}$ by the right-hand side of \eqref{wdiagerror}. For the third line, we apply the same estimate to the normalized
symbol $|\xi|^{-1}a_+$, together with its first spatial derivative,
and use the boundedness of $\mathcal SA_-$ on $Y_{B,p}$. To summarize, we have
\begin{equation}\label{Ferrorboundparaexp}
\|\mathcal S F^{\mathrm{err}}\|_{Y_{B,p}}
\lesssim_L
\lambda^{1-L}
\|\nabla_x^2\eta_{<\delta\lambda}\|_{
\dot B^0_{\infty,1}}^L
\|f_\lambda\|_{L^p}.
\end{equation}
We also record the bound 
\begin{equation}\label{SF0bound}
\|\mathcal SF^0\|_{Y_{B,p}}
\lesssim
\|\nabla_x^2\eta_{<\delta\lambda}\|_{
\dot B^0_{\infty,1}}
\|f_\lambda\|_{L^p}.
\end{equation}
which follows from by \eqref{F0kernelparaexp}, which we will need below.
\\

One last ingredient we need is an analogue of Lemma \ref{wdiagdecomp} but for the operator $(\partial_z-A_-)^{-1}$ (defined as below).
\begin{lemma}\label{bdiagdecomp}
Let $1<p\leq\infty$ and let $F=F^0+F^{\mathrm{err}}$ be as above. The solution $b$ to the equation
\begin{equation*}
(\partial_z-A_-)b=F,
\qquad
b(z)\longrightarrow0
\quad\text{as }z\longrightarrow-\infty.
\end{equation*}
has a decomposition $b=b^0+\mathcal E_f^b$ such that the map $f\mapsto b^0(0)$ has kernel bound
\begin{equation}\label{bzerokernel}
|K_{b^0(0)}(x,y)|
\lesssim_{L,M,N}
\lambda^2
\langle\lambda(x-y)\rangle^{-N}
\mathcal A_{\lambda,M}(y)
\end{equation}
and the error $\mathcal{E}_f^b$ satisfies
\begin{equation}\label{bdiagerror}
\|\mathcal E_f^b\|_{Y_{B,p}}
\lesssim_L
\lambda^{1-L}
\|\nabla_x^2\eta_{<\delta\lambda}\|_{
\dot B^0_{\infty,1}}^L
\|f_\lambda\|_{L^p}.
\end{equation}
\end{lemma}
\begin{proof}
The proof is very similar to Lemma \ref{wdiagdecomposition}. The main difference here is that the role of $\mathcal{T}$ in the resolvent expansion for $\nabla_{z,x}w_f$ will be played instead by the operator $\mathcal{K}$ defined by
\begin{equation*}
R_-\coloneq A_-+|D|,
\qquad
\mathcal K\coloneq\mathcal SR_-.
\end{equation*}
The standard Poisson kernel estimates for $\mathcal S$ and the
smallness of the normalized symbol of $R_-$ give
\begin{equation*}
\begin{split}
\|\mathcal SF\|_{Y_{B,p}}\lesssim
\|F\|_{L_z^1\dot B^0_{p,1}},
\qquad
\|R_-h\|_{L_z^1\dot B^0_{p,1}}
\lesssim
\epsilon\||D|h\|_{L_z^1\dot B^0_{p,1}},\qquad\|\mathcal Kh\|_{Y_{B,p}}\lesssim
\epsilon\|h\|_{Y_{B,p}}.
\end{split}
\end{equation*}
Hence, we have the initial resolvent representation $b=(I-\mathcal K)^{-1}\mathcal SF$ in $Y_{B,p}$. Now, using the two remaining wider projections $P_\lambda^3$ and $P_\lambda^4$
chosen above, define
\begin{equation*}
\mathcal K_\lambda
\coloneq
P_\lambda^3\mathcal K P_\lambda^4,
\qquad
b^0\coloneq(I-\mathcal K_\lambda)^{-1}\mathcal SF^0.
\end{equation*}
The kernel estimate for $b^0(0)$ follows by a similar Neumann-series expansion to the proof of Lemma \ref{wdiagdecomp}. For the remaining error term, we have 
\begin{equation*}
\begin{split}
\mathcal{E}_{f}^b\coloneq b-b^0=
(I-\mathcal K)^{-1}\mathcal SF^{\mathrm{err}}+
(I-\mathcal K)^{-1}
(\mathcal K-\mathcal K_\lambda)b^0.
\end{split}
\end{equation*}
Expanding
\begin{equation*}
b^0
=
\sum_{j=0}^{\infty}
\mathcal K_\lambda^j\mathcal SF^0,
\end{equation*}
using \eqref{SF0bound} and \eqref{Ferrorboundparaexp} and arguing as in Lemma \ref{wdiagdecomposition}, we see that $\mathcal{E}_{f}^{b}$ has the claimed bound.
\end{proof}
Applying Lemma \ref{bdiagdecomp} to \eqref{bequationparaexp}, we write
\begin{equation*}
b=b^0+\mathcal E_f^b.
\end{equation*}
The boundary trace identity gives
\begin{equation*}
G_{<\delta\lambda}f_\lambda
=
\rho^{-1}\partial_zw_f(0)
-
q\cdot\nabla_xf_\lambda.
\end{equation*}
On the other hand,
\begin{equation*}
\begin{split}
\rho^{-1}a_+(x,\xi)-iq\cdot\xi
=
\sqrt{
(1+|q|^2)|\xi|^2-(q\cdot\xi)^2
}=
\mathfrak g_{<\delta\lambda}(x,\xi)
\end{split}
\end{equation*}
so that
\begin{equation*}
\rho^{-1}A_+-q\cdot\nabla_x
=
Op(\mathfrak g_{<\delta\lambda}).
\end{equation*}
Recalling that $b(0)=\partial_zw_f(0)-A_+f_\lambda$, we obtain the identity
\begin{equation}\label{exactGprincipalidentity}
G_{<\delta\lambda}f_\lambda
-
Op(\mathfrak g_{<\delta\lambda})f_\lambda
=
\rho^{-1}b(0).
\end{equation}

Let $P_\lambda^l$ denote the output projection in the definition of
$G^l_{<\delta\lambda}$, and recall that
$\widetilde P_\lambda f_\lambda=f_\lambda$. We define
\begin{equation*}
\begin{split}
G^0_{<\delta\lambda}f_\lambda \coloneq
P_\lambda^l\left(\rho^{-1}b^0(0)\right)+
\left(
P_\lambda^lOp(\mathfrak g_{<\delta\lambda})
-Op(\mathfrak g^l_{<\delta\lambda})
\right)\widetilde P_\lambda f_\lambda
\end{split}
\end{equation*}
and $\mathcal E_f
\coloneq
P_\lambda^l\left(\rho^{-1}\mathcal E_f^b(0)\right)$. Multiplication by the bounded function $\rho^{-1}$ and convolution
with the kernel of $P_\lambda^l$ preserve \eqref{bzerokernel}. The
principal-symbol truncation error in the second term on the right above satisfies the same kernel
bound. Indeed, we choose the output and spatial-symbol cutoffs
compatibly, so that the spatial-symbol cutoff lies strictly inside
the region on which the output cutoff is equal to one in Fourier space. Its low spatial-frequency part then vanishes. Bernstein then gives a factor
\begin{equation*}
\lambda^{-1}\partial_x\mathfrak g_{<\delta\lambda},
\qquad
|\partial_x\partial_\xi^\alpha
\mathfrak g_{<\delta\lambda}(x,\xi)|
\lesssim_\alpha
\lambda^{1-|\alpha|}|\nabla_xq(x)|
\end{equation*}
on $|\xi|\approx\lambda$. Repeated integration by parts in $\xi$,
together with \eqref{localizedgeometrybound}, therefore gives
\eqref{K0bound} for $G^0_{<\delta\lambda}$. Finally, Proposition \ref{besovparabounds}, Proposition
\ref{MoserBesov}, and \eqref{bdiagerror} give
\begin{equation*}
\begin{split}
\|\mathcal E_f\|_{L^p}\lesssim
\|\rho^{-1}\mathcal E_f^b(0)\|_{\dot B^0_{p,1}}
\lesssim
\|\mathcal E_f^b\|_{Y_{B,p}}\lesssim_L
\lambda^{1-L}
\|\nabla_x^2\eta_{<\delta\lambda}\|_{
\dot B^0_{\infty,1}}^L
\|f_\lambda\|_{L^p}.
\end{split}
\end{equation*}
This proves \eqref{Gexpprecise} and concludes the proof.
\end{proof}


\subsection{Bilinear estimates for \texorpdfstring{$\varphi$}{}} Here, we deduce linear and bilinear estimates for the good variable $\varphi$ from the bootstrap bounds \eqref{uk-ee-boot}-\eqref{L4Lp-boot}. The corresponding $(u,\eta)$ bounds come for free simply from the fact that $\eta=\Re(u)$. We begin with the relevant linear and unbalanced bilinear estimates involving $\varphi$. Afterwards, we will prove a variation of the balanced estimate which has to be phrased more carefully. 
\begin{proposition}\label{phiu-bootstrap}
Suppose that $u$ satisfies the bootstrap bounds in Section \ref{sec:bootstrap}.
Then the following bounds hold:
\begin{enumerate}[label=(\roman*)]
\item Uniform frequency envelope bound:
\begin{equation}\label{phik-ee}
\|\varphi_\lambda\|_{L_t^\infty L_x^2}
\lesssim
\epsilon c_\lambda\lambda^{-s_c-\frac12}.
\end{equation}

\item $L_t^4L_x^p$ frequency envelope bound: For every $8<p\leq\infty$,
\begin{equation}\label{phik-L4Lp}
\|\varphi_\lambda\|_{L_t^4L_x^p}
\lesssim_p
\epsilon c_\lambda
\lambda^{-s_c+\frac38-\frac2p}.
\end{equation}

\item Unbalanced bilinear bounds: If
$\mu\ll\lambda$ or $\mu\gg\lambda$, then
\begin{equation}\label{phiu-bi-unbal}
\|u_{\mu}^{x_0}\varphi_{\lambda}\|_{L^2_{t,x}}
\lesssim
\epsilon^2c_\lambda c_\mu
\lambda^{-\frac12-s_c}\mu^{-s_c}\min(\lambda,\mu)^{\frac34}
\left(1+\min(\lambda,\mu)|x_0|\right)^{\frac12}.
\end{equation}

\item Translation-averaged unbalanced bilinear bounds: If
$\min(\lambda,\mu)\ll\max(\lambda,\mu)$, then
\begin{equation}\label{phiu-bi-unbal-trans}
\left\|
\varphi_{\lambda}(x)u_{\mu}(y)
\right\|_{L_t^2L_{x,y}^2(B_r)}
\lesssim
\epsilon^2c_\lambda c_\mu
\lambda^{-\frac12-s_c}\mu^{-s_c}\min(\lambda,\mu)^{\frac14}r^{\frac12},
\qquad r\geq\min(\lambda,\mu)^{-1}.
\end{equation}
where $B_r=\{(x,y)\in\mathbb R^2\times\mathbb R^2:|x-y|\leq r\}$.
\item Unbalanced $(\varphi,\varphi)$ bilinear bounds: If
$\mu\ll\lambda$ or $\mu\gg\lambda$, then
\begin{equation}\label{phiphi-bi-unbal}
\|\varphi_{\mu}^{x_0}\varphi_{\lambda}\|_{L^2_{t,x}}
\lesssim
\epsilon^2c_\lambda c_\mu
\lambda^{-\frac12-s_c}\mu^{-\frac12-s_c}\min(\lambda,\mu)^{\frac34}
\left(1+\min(\lambda,\mu)|x_0|\right)^{\frac12}.
\end{equation}

\item Translation-averaged unbalanced $(\varphi,\varphi)$ bilinear bounds: If
$\min(\lambda,\mu)\ll\max(\lambda,\mu)$, then
\begin{equation}\label{phiphi-bi-unbal-trans}
\left\|
\varphi_{\lambda}(x)\varphi_{\mu}(y)
\right\|_{L_t^2L_{x,y}^2(B_r)}
\lesssim
\epsilon^2c_\lambda c_\mu
\lambda^{-\frac12-s_c}\mu^{-\frac12-s_c}\min(\lambda,\mu)^{\frac14}r^{\frac12},
\qquad r\geq\min(\lambda,\mu)^{-1}.
\end{equation}
\end{enumerate}
The implicit constants in the bounds \eqref{phiu-bi-unbal} and \eqref{phiphi-bi-unbal} are uniform with respect to $x_0\in\mathbb R^2$.
\end{proposition}

\begin{proof}
Define the multiplicative inverse symbol
\begin{equation*}
\mathfrak{q}^{-1}_{<\lambda}
\coloneq
\left(
\frac{a_{<\lambda}}{\mathfrak{g}_{<\lambda}}
\right)^{\frac12}.
\end{equation*}
It follows from the definition of $u$ that $\Im(u)=T_{\mathfrak{q}}\varphi.$ For the remainder of the proof, we use the norm
\begin{equation*}
X\coloneq L_x^\infty\cap\dot{B}_{2,1}^1,
\qquad
\|b\|_X\coloneq
\|b\|_{L_x^\infty}
+
\|\nabla_xb\|_{\dot{B}_{2,1}^0}
\end{equation*}
to estimate some of the kernels arising from composition errors.

The smallness bounds for $\nabla\eta$, together with
\eqref{Taylorcoeffsmall} and Proposition \ref{MoserBesov}, show that the
normalized symbols
\begin{equation*}
|\xi|^{-\frac12}\mathfrak{q}_{<\rho}(x,\xi),
\qquad
|\xi|^{\frac12}\mathfrak{q}^{-1}_{<\rho}(x,\xi)
\end{equation*}
are uniformly bounded in $X$ and differ from $1$ by $O(\epsilon)$ in
this norm. The analogous bounds hold for the required normalized
$\xi$ derivatives.

By the standard symbolic calculus for paradifferential operators, we have
\begin{equation*}
T_{\mathfrak q^{-1}}T_{\mathfrak q}
=
I+\mathcal C_{\mathfrak q},
\qquad
\mathcal C_{\mathfrak q}
\coloneq
T_{\mathfrak q^{-1}}T_{\mathfrak q}
-
T_{\mathfrak q^{-1}\mathfrak q},
\end{equation*}
where we used
$\mathfrak q^{-1}\mathfrak q=1$ and
$T_{\mathfrak q^{-1}\mathfrak q}=I$. Here
$\mathcal C_{\mathfrak q}$ is an order zero paradifferential
composition error whose normalized symbol seminorms are
$O(\epsilon)$. Consequently,
\begin{equation}\label{phiparainverse}
\varphi
=
T_{\mathfrak{q}^{-1}}\Im(u)
-
\mathcal{C}_{\mathfrak{q}}\varphi.
\end{equation}
A direct computation of the dyadic kernels of
$T_{\mathfrak{q}^{-1}}P_{\lambda}$ and
$\mathcal{C}_{\mathfrak{q}}P_{\lambda}$ gives the representations
\begin{equation}\label{qinversekernels}
\begin{split}
T_{\mathfrak q^{-1}}P_{\lambda}F(t,x)
&=
\lambda^{-\frac12}
\int_{\mathbb R^2}
K_{\lambda}(t,x,y)F(t,x-y)dy,
\\
\mathcal C_{\mathfrak q}P_{\lambda}F(t,x)
&=
\int_{\mathbb R^2}
E_{\lambda}(t,x,y)F(t,x-y)dy.
\end{split}
\end{equation}
For every $N\geq0$, repeated integration by parts in the frequency
variable gives
\begin{equation}\label{qinversekernelbounds}
\begin{split}
\int_{\mathbb R^2}
\langle\lambda y\rangle^N
\sup_t\|K_{\lambda}(t,\cdot,y)\|_Xdy
&\lesssim_N1,
\\
\int_{\mathbb R^2}
\langle\lambda y\rangle^N
\sup_t\|E_{\lambda}(t,\cdot,y)\|_Xdy
&\lesssim_N\epsilon.
\end{split}
\end{equation}
Indeed, after rescaling to unit frequency, we may write the normalized
symbols in the form
\begin{equation*}
|\xi|^{-\frac12}\mathfrak q_{<\lambda}
=
1+\widetilde{\mathfrak q}_{<\lambda},
\qquad
|\xi|^{\frac12}\mathfrak q^{-1}_{<\lambda}
=
1+\widetilde{\mathfrak q}^{-1}_{<\lambda},
\end{equation*}
where the normalized symbol seminorms of
$\widetilde{\mathfrak q}_{<\lambda}$ and
$\widetilde{\mathfrak q}^{-1}_{<\lambda}$ are $O(\epsilon)$ in $X$.
The flat symbols compose exactly to the identity, so every term in the
composition error contains at least one smallness factor $\epsilon$. Together
with the spatial localization of the symbols to frequencies
$\ll\lambda$, this gives \eqref{qinversekernelbounds}.
\\

With these preliminaries, we are now ready to establish
\eqref{phik-ee}-\eqref{phiphi-bi-unbal-trans}. We begin with the first two linear
estimates. From \eqref{phiparainverse}--\eqref{qinversekernelbounds}, we
have, for $Y=L_t^\infty L_x^2$ or $Y=L_t^4L_x^p$,
\begin{equation*}
\|\varphi_{\lambda}\|_Y
\lesssim
\lambda^{-\frac12}
\sum_{\nu\approx\lambda}
\|u_{\nu}\|_Y
+
\epsilon
\sum_{\nu\approx\lambda}
\|\varphi_{\nu}\|_Y.
\end{equation*}
The corresponding bounds for $u$, a bootstrap on the corresponding
$Y$ bound for $\varphi$, and the slow variation of the frequency
envelope give \eqref{phik-ee} and \eqref{phik-L4Lp}.

We next move to the unbalanced bilinear estimates. To simplify
notation, we set
\begin{equation*}
\Lambda\coloneq\max(\lambda,\mu),
\qquad
\alpha\coloneq\min(\lambda,\mu).
\end{equation*}
Using the first kernel representation in \eqref{qinversekernels} and
Minkowski's inequality, we obtain
\begin{equation*}
\begin{split}
\left\|
T_{\mathfrak{q}^{-1}}\Im(u_{\lambda})
\,u_{\mu}^{x_0}
\right\|_{L_{t,x}^2}
\lesssim
\lambda^{-\frac12}
\int_{\mathbb R^2}
\sup_{t,x}|K_{\lambda}(t,x,y)|
\left\|
u_{\lambda}u_{\mu}^{x_0+y}
\right\|_{L_{t,x}^2}dy.
\end{split}
\end{equation*}
Thus, after a translation in $x$, we may apply the unbalanced bootstrap
estimate for $u$. Moreover,
\begin{equation}\label{translationkernelweight}
1+\alpha|x_0\pm y|
\lesssim
\left(1+\alpha|x_0|\right)
\langle\lambda y\rangle.
\end{equation}
It follows from \eqref{qinversekernelbounds} that
\begin{equation*}
\begin{split}
\left\|
P_{\lambda}T_{\mathfrak{q}^{-1}}\Im(u)
\,u_{\mu}^{x_0}
\right\|_{L_{t,x}^2}
&\lesssim
\epsilon^2c_\lambda c_\mu
\lambda^{-\frac12}
\Lambda^{-s_c}\alpha^{-s_c+\frac34}
\left(1+\alpha|x_0|\right)^{\frac12}
\\
&=
\epsilon^2c_\lambda c_\mu
\lambda^{-\frac12-s_c}\mu^{-s_c}\alpha^{\frac34}
\left(1+\alpha|x_0|\right)^{\frac12}.
\end{split}
\end{equation*}
The error term in \eqref{phiparainverse} is estimated in exactly the
same way and contributes, by a continuity argument, at most
$O(\epsilon)$ times the right-hand side of \eqref{phiu-bi-unbal}.
Absorbing this contribution proves \eqref{phiu-bi-unbal}. The
translation-averaged bound \eqref{phiu-bi-unbal-trans} is handled in a
similar fashion. We leave the straightforward modifications to the
reader. Finally, \eqref{phiphi-bi-unbal} and
\eqref{phiphi-bi-unbal-trans} follow by applying
\eqref{phiparainverse} to the remaining $\varphi_\mu$ factor
and repeating the same kernel argument and absorption, using
\eqref{phiu-bi-unbal} and \eqref{phiu-bi-unbal-trans}, respectively.
\end{proof}
Before concluding this subsection, we state a weaker variation of the balanced estimate \eqref{uab-bi-bal}. We note importantly that this bound is not invariant under interchanging $u$ with its complex conjugate, unlike the unbalanced estimates. Because of this, we cannot prove an analogous balanced bound for the full real variables $\eta$ and $\varphi$. However, it is convenient to observe the following weaker variant that does hold in this setting, as long as the output frequency is not too large. This will be a convenient tool in our layered paradifferential expansions in the forthcoming sections; it will allow us to efficiently estimate certain bilinear products of the real variables $\eta$ and $\varphi$ without expanding them in terms of $u$ and $\overline{u}$. 
\begin{lemma}[Balanced bilinear $L^2$ estimates with restricted output]\label{balancedrestrictedoutput}
Assume the bootstrap bo\-unds above, and the conventions in Subsection \ref{subsec:littlewoodpaley}. Let $\lambda,\mu>0$ satisfy $\lambda\approx\mu$, and let $x_0\in\mathbb R^2$. For any combination
\begin{equation*}
f_\lambda\in\{u_\lambda,\overline u_\lambda,\lambda^{\frac12}\varphi_\lambda\},
\qquad
g_\mu\in\{u_\mu,\overline u_\mu,\mu^{\frac12}\varphi_\mu\},
\end{equation*}
there holds
\begin{equation}\label{uab-bi-bal-restricted}
\begin{split}
\sup_{
\rho<\frac32\max(\lambda,\mu)}
\left\|
P_\rho|D|^{\frac12}
\left(f_\lambda g_\mu^{x_0}\right)
\right\|_{L_{t,x}^2}
\lesssim
\epsilon^2c_\lambda c_\mu\lambda^{-2s_c+\frac54}\left(1+\lambda|x_0|\right)^{\frac12}.
\end{split}
\end{equation}
\end{lemma}
This will be useful for efficiently estimating products of $\eta$ and $\varphi$ in situations where they have comparable input frequencies, but sufficiently small output frequency.
\begin{proof}
Set $\Lambda=\max(\lambda,\mu)$. We first consider $f_\lambda\in\{u_\lambda,\overline u_\lambda\}$ and $g_\mu\in\{u_\mu,\overline u_\mu\}$. The case $f_{\lambda}=u_{\lambda}$ and $g_{\mu}=\overline{u}_{\mu}$ follows from \eqref{uab-bi-bal}. Without loss of generality, it remains to consider $f_\lambda=u_\lambda$ and $g_\mu=u_\mu$.

The main observation that we will need is that the low output frequency restriction $\rho<\frac32\max(\lambda,\mu)$ ensures that $f_{\lambda}$ and $g_{\mu}$ are frequency separated. Indeed, let $\xi$ be the frequency variable belonging to the input at scale $\Lambda$, and let $\eta$ be the other input frequency. If this pair contributes to an output with $\rho<\frac32\Lambda$, then
\begin{equation*}
|\xi|\geq \step^{-1}\Lambda,
\qquad
|\xi+\eta|\leq \step\rho<\frac32\step\Lambda.
\end{equation*}
Above, and below $\step$ is the constant defined in the Littlewood-Paley conventions in subsection \ref{subsec:littlewoodpaley}. Consequently,
\begin{equation*}
|\xi-\eta|
=|2\xi-(\xi+\eta)|
\geq\left(\frac2\step-\frac32\step\right)\Lambda
\coloneq d\Lambda,
\qquad d>0.
\end{equation*}
Next, we decompose the two input functions into finitely many dyadic pieces of diameter at most $d\Lambda/8$, and write
\begin{equation*}
u_\lambda=\sum_a f_a,
\qquad
u_\mu=\sum_b g_b.
\end{equation*}
We retain only pairs whose sum of Fourier supports intersects the ball $|\zeta|\leq\frac32\step\Lambda$. Every other pair vanishes after applying $P_{\rho}$. By the above bound, each retained pair is $\approx\Lambda$ separated in teh sense that we have 
\begin{equation*}
\operatorname{supp}\mathcal F_x
\left(f_a\overline g_b^{x_0}\right)
\subset
\left\{\zeta:\frac d2\Lambda\leq|\zeta|\leq2c\Lambda\right\}.
\end{equation*}
To conclude, we now make use of the simple fact that the balanced bound \eqref{uab-bi-bal} is "stable" under slightly shrinking the frequency inputs. Indeed by frequency localization, we can write $f_a=K_a*u_\lambda$ and $g_b=J_b*u_\mu$, where the kernels satisfy
\begin{equation*}
\int_{\mathbb R^2}\langle\Lambda y\rangle
\left(|K_a(y)|+|J_b(y)|\right)dy\lesssim1.
\end{equation*}
Minkowski's inequality and translation invariance then gives
\begin{equation*}
\begin{split}
\left\|P_\tau|D|^{\frac12}
\left(f_a\overline g_b^{x_0}\right)\right\|_{L_{t,x}^2}&\leq
\int_{\mathbb R^2}\int_{\mathbb R^2}
|K_a(y)||J_b(z)|
\left\|P_\tau|D|^{\frac12}
\left(u_\lambda\overline u_\mu^{x_0+y-z}\right)
\right\|_{L_{t,x}^2}dydz
\\
&\lesssim
\epsilon^2c_\lambda c_\mu
\Lambda^{-2s_c+\frac54}
\left(1+\Lambda|x_0|\right)^{\frac12},
\end{split}
\end{equation*}
uniformly in $\tau$. Here we used \eqref{uab-bi-bal} and $1+\Lambda|x_0+y-z|
\lesssim
\left(1+\Lambda|x_0|\right)
\langle\Lambda y\rangle\langle\Lambda z\rangle.$ We then have by Bernstein (and the $\approx_d\Lambda$ frequency separation)
\begin{equation*}
\begin{split}
\left\|f_ag_b^{x_0}\right\|_{L_{t,x}^2}=\left\|f_a\overline g_b^{x_0}\right\|_{L_{t,x}^2}\lesssim
\epsilon^2c_\lambda c_\mu
\Lambda^{-2s_c+\frac34}\left(1+\Lambda|x_0|\right)^{\frac12},
\end{split}
\end{equation*}
which gives again by Bernstein, for each $\rho<\frac32\Lambda$,
\begin{equation*}
\begin{split}
\left\|P_\rho|D|^{\frac12}
\left(u_\lambda u_\mu^{x_0}\right)\right\|_{L_{t,x}^2}\lesssim
\rho^{\frac12}\sum_{a,b}
\left\|f_ag_b^{x_0}\right\|_{L_{t,x}^2}\lesssim
\epsilon^2c_\lambda c_\mu\Lambda^{-2s_c+\frac54}\left(1+\Lambda|x_0|\right)^{\frac12}.
\end{split}
\end{equation*}
Taking the supremum in $\rho$ proves the bound for the case $f_{\lambda}=u_{\lambda}$ and $g_{\mu}=u_{\mu}$, as desired.

For the terms containing $\varphi$, we can first use the Littlewood Paley trichotomy to obtain the general bound
\begin{equation}\label{balancedmultiplier}
\sup_\tau\left\|P_\tau|D|^{\frac12}(bF)\right\|_{L_{t,x}^2}
\lesssim
\|b\|_{L_t^\infty X}
\sup_\tau\left\|P_\tau|D|^{\frac12}F\right\|_{L_{t,x}^2},
\end{equation}
where $X$ is defined in the preceding proof. Applying this to the kernels in \eqref{qinversekernels}-\eqref{qinversekernelbounds}, and absorbing the small composition errors as in the proof of the previous proposition, transfers the bound \eqref{uab-bi-bal} to products in which either or both occurrences of $u_\lambda$ and $u_\mu$ are replaced by
\begin{equation*}
\lambda^{\frac12}P_\lambda(T_{\mathfrak q})^{-1}u,
\qquad
\mu^{\frac12}P_\mu(T_{\mathfrak q})^{-1}u,
\end{equation*}
respectively. These pieces are localized to the same dyadic frequencies as $u_\lambda$ and $u_\mu$, so the frequency separation argument above also gives the bound for all their conjugation choices. Since $T_{\mathfrak q}$ commutes with complex conjugation and $\varphi=\Im((T_{\mathfrak q})^{-1}u)$, expanding the imaginary parts proves \eqref{uab-bi-bal-restricted} for the remaining choices of $f_\lambda$ and $g_\mu$. This concludes the proof.
\end{proof}
\begin{remark}
Both Lemma \ref{balancedrestrictedoutput} and Proposition \ref{phiu-bootstrap} have analogues for the linearized variables, $w$, $r$ and $h$. The proofs in this case are identical to the ones above, except one uses the linearized bootstrap bounds in place of the nonlinear ones. Instead of enumerating all of the relevant estimates here, we will refer back to this subsection when necessary, and simply remark on the relevant modifications.    
\end{remark}
\subsection{Multi-linear norms and translation-averaged estimates}\label{subsec:multilinearnorms}

Here, we collect some useful multi-linear estimates. As usual, we assume throughout that the bootstrap hypotheses outlined in Section \ref{sec:bootstrap}. The main purpose of this subsection will be to setup a functional framework to more easily exploit the translation-averaged bound \eqref{uab-bi-unbal-trans} in the error estimates throughout this paper. This will mostly be relevant in scenarios when we are dealing with multi-linear kernels possessing relatively limited decay. A prototypical example would be when $|D|$ is applied to a product of two functions at balanced frequencies with a low frequency output. In this case, the kernel is not rapidly decaying at the input scale; in fact, it is not even integrable against the weight $|y|$. In such situations, we will use \eqref{uab-bi-unbal-trans}, which is consistent with these kernels possessing better decay when measured in $L^2_y$ instead of $L^1_y$ (against a suitable polynomial weight).
\\

We begin by defining the following weighted multi-linear "norm" that will be helpful for obtaining certain multi-linear estimates in situations where the corresponding multilinear kernel has relatively limited decay. 
\begin{definition}
Let $\sigma>0$, $1\leq p<\infty$, and let
$\lambda\in 2^{\mathbb Z}$. Define
\begin{equation}
M_{\lambda,p}^{\sigma}(G,F)
:=
\left(
\lambda^2\int_{\mathbb R^2}
\frac{\|G\,\tau_yF\|_{L_{t,x}^p}^p}
{\langle\lambda y\rangle^{p\sigma}}\,dy
\right)^{\frac{1}{p}}.
\end{equation}
\end{definition}
The motivation for this definition is provided by the following lemma that shows that $M_{\lambda,2}^{\sigma}(G,F)$ can be controlled in terms of $\lambda^{\frac{1}{2}}
\sup_{\substack{r\in 2^{\mathbb Z}\\r\geq\lambda^{-1}}}
r^{-\frac{1}{2}}
\left\|
\chi_r(y)G\tau_yF
\right\|_{L_y^2L_t^2L_x^2}$, which averages in $L^2$ over translations. As mentioned above, this is sometimes a bit more flexible than using the norm
\begin{equation*}
\sup_{y}\langle y\rangle^{-1}\|G\tau_yF\|_{L_{t,x}^2}    
\end{equation*}
which is more susceptible to endpoint losses when we encounter operators possessing kernels with only limited decay (i.e. those for which $\langle y\rangle K$ fails to be integrable).
\begin{lemma}[Annular control of the translation norm]\label{annularcontrol}
Let $\lambda\in 2^{\mathbb Z}$ and let $\sigma>\frac{1}{2}$. Let
$\{\chi_r\}_{r\geq\lambda^{-1}}$ be a dyadic partition
of unity in the $y$ variable such that $\chi_{\lambda^{-1}}$ is
supported on the ball $\{|y|\lesssim\lambda^{-1}\}$ and, for
$r>\lambda^{-1}$,
\begin{equation*}
\chi_r(y)\coloneq\chi(r^{-1}y),
\end{equation*}
where $\chi$ is supported on the annulus
$\{\frac{1}{2}\leq|y|\leq2\}$. Then
\begin{equation*}
M_{\lambda,2}^{\sigma}(G,F)
\lesssim_{\sigma}
\lambda^{\frac{1}{2}}
\sup_{\substack{r\in 2^{\mathbb Z}\\r\geq\lambda^{-1}}}
r^{-\frac{1}{2}}
\left\|
\chi_r(y)G\tau_yF
\right\|_{L_y^2L_t^2L_x^2},
\end{equation*}
where the implicit constant is uniform in $\lambda$.
\end{lemma}

\begin{proof}
We set
\begin{equation*}
\mathcal N\coloneq
\sup_{\substack{r\in 2^{\mathbb Z}\\r\geq\lambda^{-1}}}
r^{-\frac{1}{2}}
\left\|
\chi_r(y)G\tau_yF
\right\|_{L_y^2L_t^2L_x^2}.
\end{equation*}
On the support of $\chi_r$, one has
$\langle\lambda y\rangle\approx\lambda r$. Therefore,
\begin{equation*}
\begin{split}
\left(M_{\lambda,2}^{\sigma}(G,F)\right)^2
&\lesssim
\lambda^2
\sum_{\substack{r\in 2^{\mathbb Z}\\r\geq\lambda^{-1}}}
(\lambda r)^{-2\sigma}
\left\|
\chi_r(y)G\tau_yF
\right\|_{L_y^2L_t^2L_x^2}^2\\
&\leq
\lambda^2\mathcal N^2
\sum_{\substack{r\in 2^{\mathbb Z}\\r\geq\lambda^{-1}}}
(\lambda r)^{-2\sigma}r\\
&\lesssim_\sigma
\lambda\mathcal N^2.
\end{split}
\end{equation*}
In the last inequality, we used $\sigma>\frac{1}{2}$. Taking the
square root concludes the proof.
\end{proof}
Next, we collect some basic properties satisfied by these multi-linear norms.
\begin{proposition}[Basic properties of the $M_{\lambda,p}^{\sigma}$ bilinear norm]\label{basicMproperties}
The following estimates hold.
\begin{enumerate}
\item  Let $\lambda\in 2^{\mathbb Z}$ and $\sigma_1,\sigma_2,\sigma\geq 0$ with $\sigma=\sigma_1+\sigma_2$. Moreover, let $1\leq p,p_1,p_2<\infty$ with $\frac{1}{p}=\frac{1}{p_1}+\frac{1}{p_2}$. Then there holds
\begin{equation}\label{M-holder-product}
M_{\lambda,p}^{\sigma}
\left(G_1G_2,F_1F_2\right)
\leq
M_{\lambda,p_1}^{\sigma_1}(G_1,F_1)
M_{\lambda,p_2}^{\sigma_2}(G_2,F_2).
\end{equation}
\item Let $\lambda\in 2^{\mathbb Z}$, let $1\leq p<\infty$, and let $0<\sigma<\frac{2}{p}$. Suppose that $\frac{1}{a}+\frac{1}{a'}=\frac{1}{p}$ and $\frac{1}{r}+\frac{1}{s}
>
\frac{2}{p}-\frac{\sigma}{2}$ where $p\leq a,a',r,s\leq\infty$. Then there holds
\begin{equation}\label{M-translation-free}
M_{\lambda,p}^{\sigma}(G,F)
\lesssim_{\sigma,p,r,s}
\lambda^{\frac{2}{r}+\frac{2}{s}-\frac{2}{p}}
\|G\|_{L_t^aL_x^r}
\|F\|_{L_t^{a'}L_x^s}.
\end{equation}
\item Let $\sigma\geq 0$ and let $\alpha,\beta\in 2^{\mathbb Z}$ with $\alpha\leq\beta$, then there holds
\begin{equation}\label{M-scale-eq}
M_{\beta,2}^{\sigma}(G,F)
\gtrsim_{\sigma}
\left(\frac{\beta}{\alpha}\right)^{1-\sigma}
M_{\alpha,2}^{\sigma}(G,F).
\end{equation}
Moreover, if $0\leq\sigma<1$ and if $\widehat{F}$ is supported at frequency $|\xi|\lesssim\alpha$, then the reverse inequality is also true.
\end{enumerate}
\end{proposition}

\begin{proof}
The first estimate \eqref{M-holder-product} is a straightforward application of H\"older's inequality and the definition of $M_{\lambda,p}^{\sigma}$. For the second estimate, set
\begin{equation*}
K_{\lambda,p,\sigma}(y)
\coloneq
\lambda^2\langle\lambda y\rangle^{-p\sigma}.
\end{equation*}
Fubini's theorem gives
\begin{equation}
\left(M_{\lambda,p}^{\sigma}(G,F)\right)^p
=
\int_{\mathbb R}\int_{\mathbb R^2}
|G(t,x)|^p
\left(K_{\lambda,p,\sigma}*|F(t,\cdot)|^p\right)(x)
\,dx\,dt.
\end{equation}
Define $\ell$ by $\frac{1}{\ell}
=
2-\frac{p}{r}-\frac{p}{s}$ so that (by hypothesis) $\ell>\frac{2}{p\sigma}$. Hence,
\begin{equation*}
\begin{split}
\|K_{\lambda,p,\sigma}\|_{L^\ell}
&=
\lambda^{2-\frac{2}{\ell}}
\left\|\langle y\rangle^{-p\sigma}\right\|_{L^\ell}
\lesssim_{\sigma,p,r,s}
\lambda^{\frac{2p}{r}+\frac{2p}{s}-2}.
\end{split}
\end{equation*}
Moreover, $1+\frac{1}{(r/p)'}
=
\frac{1}{\ell}+\frac{1}{s/p}.$ Hence, Young's inequality gives
\begin{equation*}
\begin{split}
\int_{\mathbb R^2}
|G(t,x)|^p
\left(K_{\lambda,p,\sigma}*|F(t,\cdot)|^p\right)(x)\,dx
\lesssim_{\sigma,p,r,s}
\lambda^{\frac{2p}{r}+\frac{2p}{s}-2}
\|G(t)\|_{L_x^r}^p
\|F(t)\|_{L_x^s}^p.
\end{split}
\end{equation*}
H\"older's inequality in $t$, followed by taking the $p$-th root, proves \eqref{M-translation-free}. Finally, we move to \eqref{M-scale-eq}. First, we handle the $\gtrsim_{\sigma}$ part of the inequality. This follows from the elementary estimate
\begin{equation*}
\alpha^2\langle\alpha y\rangle^{-2\sigma}
\leq
\left(\frac{\alpha}{\beta}\right)^{2-2\sigma}
\beta^2\langle\beta y\rangle^{-2\sigma}
\end{equation*}
and then multiplying by
$\|G\,\tau_yF\|_{L_{t,x}^2}^2$, integrating in $y$, and taking the
square root. Next, we prove the reverse inequality under the hypotheses $0\leq\sigma<1$ and with $\widehat{F}$ supported at frequency $\lesssim\alpha$. In this setting, there exists a nonnegative rapidly decreasing kernel $\Phi_\alpha(y)
=
\alpha^2\Phi(\alpha y)$ such that
\begin{equation*}
|F(t,x)|^2
\lesssim
\left(\Phi_\alpha*|F(t,\cdot)|^2\right)(x).
\end{equation*}
Thus, setting $H(y)
\coloneq
\left\|G\,\tau_yF\right\|_{L_{t,x}^2}^2$, we have $H\lesssim\Phi_\alpha*H$. Since $0\leq\sigma<1$, we also have the bound
\begin{equation*}
\beta^2\langle\beta\cdot\rangle^{-2\sigma}
*\Phi_\alpha
\lesssim_\sigma
\left(\frac{\beta}{\alpha}\right)^{2-2\sigma}
\alpha^2\langle\alpha\cdot\rangle^{-2\sigma}
\end{equation*}
which follows by scaling, splitting into the regions $|v|\leq 2$, $|y-v|\leq |v|/2$ and $|y-v|>|v|/2$. Therefore,
\begin{equation*}
\begin{split}
\left(M_{\beta,2}^{\sigma}(G,F)\right)
&\lesssim_\sigma
\left(\frac{\beta}{\alpha}\right)^{1-\sigma}
\left(M_{\alpha,2}^{\sigma}(G,F)\right).
\end{split}
\end{equation*}
This concludes the proof.
\end{proof}
Next, we need the following proposition which will allow us to prove bilinear estimates for the nonlocal elliptic operator $\mathcal{P}f\coloneq \partial_z\Delta_{\Omega}^{-1}f|_{z=\eta(t,x)}$. This operator arises naturally in several places. For instance, it appears in formulas involving commutations with the operator $L$, and also in the transport equation for the Taylor coefficient.
\begin{proposition}[Weighted bilinear elliptic estimate]\label{weightedbilinear}
Define the operator $D_i\coloneq (\partial_i-\partial_i\eta\partial_z)$. Let $1<\sigma<2$. Given a smooth function $f$ with sufficient decay, let $f_*(z,x)\coloneq f(x,z+\eta(x))$ denote the corresponding flattened map. Let $\beta\in 2^{\mathbb{Z}}$. Let $h:\mathbb{R}^2\to\mathbb{R}$ be a Schwartz function. The following two bounds are true.
\begin{enumerate}
\item (divergence bound) Assume that $f_*=D\cdot g$ or $f_*=\nabla_{z,x}\cdot g$ for some function $g:\{z<0\}\to\mathbb{R}^3$. Define $f_{<\beta}$ such that $(f_{<\beta})_*=P_{<\beta}f_*$. Then there holds
\begin{equation}\label{weighted-bilinear-interior-derivative}
\|h \mathcal{P}f_{<\beta}\|_{L_t^1L_x^1}\lesssim_{\sigma}\sup_{z\leq 0}M_{\beta,1}^{\sigma}(h,g(z)).
\end{equation}
\item (non-divergence bound) Let $W=\nabla_{z,x}(\Delta_D^{-1}f)_*$. If $h=h_1h_2$, where $h_1$ and $h_2$, there holds
\begin{equation}\label{weighted-bilinear-exterior-derivative}
\begin{split}
\left\|
h\nabla_xP_{<\beta}\mathcal Pf
\right\|_{L_t^1L_x^1}
\lesssim_\sigma
\sup_{z\leq0}
M_{\beta,1}^{\sigma}
\left(h,f_*(z)\right)+
M_{\beta,2}^{\frac\sigma2}
\left(h_1,\nabla^2_x\eta\right)
\|h_2\|_{L_t^\infty L_x^2}
\|W\|_{L_z^\infty L_t^2\dot B_{\infty,1}^0}.
\end{split}
\end{equation}
\end{enumerate}
\end{proposition}

\begin{proof}
We begin with the first bound \eqref{weighted-bilinear-interior-derivative}. We will show the details in the case $f_*=D\cdot g$. The other case is essentially identical. Now set
\begin{equation*}
v_\beta
\coloneq
\left(
\Delta_D^{-1}f_{<\beta}
\right)_*,
\qquad
W_\beta
\coloneq
\nabla_{z,x}v_\beta.
\end{equation*}
From the definition of $D$, we have
\begin{equation*}
D\cdot g
=
\nabla_{z,x}\cdot\widetilde g,
\qquad
\widetilde g
\coloneq
\left(
g_1,g_2,g_3-\partial_i\eta\,g_i
\right).
\end{equation*}
Since $v_{\beta}$ vanishes at $z=0$ and $\Delta_{\eta}v_{\beta}\coloneq D_iD_iv_{\beta}=P_{<\beta}f_*$, we also have
\begin{equation*}
P_{<\beta}f_*
=
\nabla_{z,x}\cdot P_{<\beta}\widetilde g,
\end{equation*}
which gives the formula
\begin{equation}\label{flattened-weighted-elliptic-equation}
W_\beta
=
\mathcal R_0P_{<\beta}\widetilde g
-
\mathcal R_0(A_\eta W_\beta).
\end{equation}
where we write $\mathcal{R}_0\coloneq\nabla_{z,x}\Delta_{z,x}^{-1}\nabla_{z,x}\cdot$ Above, $A_\eta$ is a matrix whose entries are linear or quadratic in
$\nabla\eta$, so that in particular,
\begin{equation}\label{weighted-elliptic-coefficient-smallness}
\|A_\eta\|_{L_t^\infty\dot B_{\infty,1}^0}
+
\|A_\eta\|_{L_t^\infty\dot B_{2,1}^1}
\lesssim
\epsilon,
\end{equation}
As in the other elliptic estimates in this section, the strategy will be to treat (using the small data hypothesis) $W_{\beta}+\mathcal{R}_0(A_{\eta}W_{\beta})$ as a perturbation of $W_{\beta}$ in a suitable bilinear norm, and then estimate the remaining source term directly. For this purpose, we introduce the following bilinear norm,
\begin{equation*}
\|F\|_{X_h^\sigma}
\coloneq
\sum_\alpha
\sup_{z\leq0}
M_{\alpha,1}^\sigma
\left(
h,P_\alpha F(z)
\right).
\end{equation*}
We first record the flat estimate
\begin{equation}\label{flat-weighted-elliptic-estimate}
\|\mathcal R_0H\|_{X_h^\sigma}
\lesssim_\sigma
\sum_\alpha
\sup_{z\leq0}
M_{\alpha,1}^\sigma
\left(
h,\widetilde P_\alpha H(z)
\right),
\end{equation}
where $\widetilde P_\alpha$ is a sufficiently fattened annular
projection. Indeed, the kernel of
$P_\alpha\mathcal R_0\widetilde P_\alpha$ satisfies (with the natural modification for the $\delta(z-z')$ component of the kernel)
\begin{equation*}
\begin{split}
\sup_{z\leq 0}\int_{-\infty}^{0}|K_\alpha(z,z',y)|dz'
\lesssim_N \alpha^2\langle\alpha y\rangle^{-N}.
\end{split}
\end{equation*}
Indeed, using Minkowski's inequality and the bound
\begin{equation*}
\langle\alpha(y+u)\rangle^{-\sigma}
\lesssim_\sigma
\langle\alpha u\rangle^\sigma
\langle\alpha y\rangle^{-\sigma},
\end{equation*}
we obtain
\begin{equation*}
\begin{split}
M_{\alpha,1}^\sigma
\left(h,P_\alpha\mathcal R_0H(z)\right)
&\lesssim_\sigma
\int_{-\infty}^0\int_{\mathbb R^2}
|K_\alpha(z,z',u)|
\langle\alpha u\rangle^\sigma
M_{\alpha,1}^\sigma
\left(h,\widetilde P_\alpha H(z')\right)
dudz'
\\
&\lesssim_\sigma
\sup_{z'\leq0}
M_{\alpha,1}^\sigma
\left(h,\widetilde P_\alpha H(z')\right),
\end{split}
\end{equation*}
where in the last line we used $N>\sigma+2$. Taking the supremum in
$z$ and summing in $\alpha$ proves
\eqref{flat-weighted-elliptic-estimate}. Before proceeding, we also record that for $1<\sigma<2$, and uniformly in $\alpha,\beta$, there holds 
\begin{equation}\label{M-LP-boundedness}
M_{\beta,2}^{\frac{\sigma}{2}}
\left(h,P_\alpha F\right)
\lesssim_\sigma
M_{\beta,2}^{\frac{\sigma}{2}}(h,F).
\end{equation}
That is, for this range of $\sigma$, $P_{\alpha}$ is bounded in the second argument with respect to the $M$ norm. Indeed, setting $\omega_\beta(x)
\coloneq
\beta^2\langle\beta x\rangle^{-\sigma}$ and $W_{\beta,h}
\coloneq
\omega_\beta*|h|^2,$ a translation change of variables gives
\begin{equation*}
\left(
M_{\beta,2}^{\frac{\sigma}{2}}(h,F)
\right)^2
=
\int_{\mathbb R}\int_{\mathbb R^2}
W_{\beta,h}(t,x)|F(t,x)|^2dxdt.
\end{equation*}
Furthermore, by separating the cases
$\beta R\leq\frac12\langle\beta x\rangle$ and
$\beta R>\frac12\langle\beta x\rangle$, and using polar coordinates
in the latter case, we have, uniformly in $R,\beta$ and $x$,
\begin{equation*}
R^{-2}\int_{|y|\leq R}\omega_\beta(x-y)dy
\lesssim_\sigma
\omega_\beta(x)
\end{equation*}
for $\sigma<2$. Therefore, if $K_{\alpha}$ is the kernel for $P_{\alpha}$, its rapid decay $|K_{\alpha}(y)|\lesssim_N \alpha^2\langle\alpha y\rangle^{-N}$, a dyadic decomposition in space, along with the above bound gives $|K_\alpha|*\omega_\beta
\lesssim_\sigma
\omega_\beta.$ Consequently, associativity of convolution, Cauchy--Schwarz and Fubini's theorem give
\begin{equation*}
\begin{split}
\left(
M_{\beta,2}^{\frac{\sigma}{2}}
\left(h,P_\alpha F\right)
\right)^2
&\lesssim
\int_{\mathbb R}\int_{\mathbb R^2}
|F(t,x)|^2
\left(|K_\alpha|*W_{\beta,h}\right)(t,x)
dxdt\lesssim_\sigma
\left(
M_{\beta,2}^{\frac{\sigma}{2}}(h,F)
\right)^2.
\end{split}
\end{equation*}
The same estimate holds for the fattened projections
$\widetilde P_\alpha$. Now, returning to the proof of \eqref{weighted-bilinear-interior-derivative}, we next claim that under the small data bootstrap, we have
\begin{equation}\label{weighted-elliptic-absorption}
\sum_\alpha
\sup_{z\leq0}
M_{\alpha,1}^\sigma
\left(
h,P_\alpha(A_\eta F)(z)
\right)
\lesssim_\sigma
\epsilon \|F\|_{X_h^\sigma}.
\end{equation}
To see this, decompose $A_\eta F$ using the Littlewood-Paley trichotomy. The low-high
interactions satisfy
\begin{equation*}
M_{\kappa,1}^\sigma
\left(
h,
P_\kappa(A_{\eta,<\kappa}F_\kappa)
\right)
\lesssim
\|A_{\eta,<\kappa}\|_{L_{t,x}^\infty}
M_{\kappa,1}^\sigma(h,F_\kappa).
\end{equation*}
For the high--low, we claim that when $\nu\ll\kappa$, we have
\begin{equation}\label{Mhilo}
M_{\kappa,1}^\sigma
\left(
h,
P_\kappa(A_{\eta,\kappa}F_\nu)
\right)
\lesssim_\sigma
\kappa
\|A_{\eta,\kappa}\|_{L_t^\infty L_x^2}
M_{\nu,1}^\sigma(h,F_\nu).
\end{equation}
By boundedness of $P_{\kappa}$ on $M_{\kappa,1}^{\sigma}$, it suffices to prove this bound for the unprojected $A_{\eta,\kappa}F_{\nu}$ in the second argument above. Since $F_\nu$ is localized to frequency $\nu$, a reproducing kernel gives
\begin{equation}\label{reproducingbound}
|F_\nu(t,x+y)|
\lesssim_N
\nu^2
\int_{\mathbb R^2}
\langle\nu(y-y')\rangle^{-N}
|F_\nu(t,x+y')|\,dy'.
\end{equation}
After inserting this into the definition of
$M_{\kappa,1}^\sigma$ and applying Cauchy--Schwarz in $y$, we use
\begin{equation*}
\begin{split}
&\left(
\int_{\mathbb R^2}
\langle\kappa y\rangle^{-2\sigma}
\langle\nu(y-y')\rangle^{-2N}\,dy
\right)^{\frac12}
\lesssim_{\sigma,N}
\kappa^{-1}
\langle\nu y'\rangle^{-\sigma},
\qquad
\kappa\gtrsim\nu.
\end{split}
\end{equation*}
Here $N$ is chosen sufficiently large, and the estimate is uniform
in $\nu$ and $\kappa$ because $\sigma>1$. Consequently,
\begin{equation*}
\begin{split}
M_{\kappa,1}^\sigma
\left(
h,A_{\eta,\kappa}F_\nu
\right)
&\lesssim_\sigma
\kappa^2\nu^2\kappa^{-1}
\|A_{\eta,\kappa}\|_{L_t^\infty L_x^2}
\int_{\mathbb R^2}
\langle\nu y'\rangle^{-\sigma}
\|h\,\tau_{y'}F_\nu\|_{L_{t,x}^1}\,dy'
\\
&=
\kappa
\|A_{\eta,\kappa}\|_{L_t^\infty L_x^2}
M_{\nu,1}^\sigma(h,F_\nu).
\end{split}
\end{equation*}
which proves \eqref{Mhilo}. In particular, the high-frequency coefficient is summable using
$\|A_\eta\|_{L_t^\infty\dot B_{2,1}^1}$. Finally, if
$\kappa'\approx\kappa$ and $\alpha\lesssim\kappa$, then boundedness of
$P_\alpha$ with respect to $M_{\alpha,1}^{\sigma}$ (which follows from the reproducing kernel bound \eqref{reproducingbound} for $\nu=\alpha$) gives
\begin{equation*}
\begin{split}
M_{\alpha,1}^\sigma
\left(
h,
P_\alpha(A_{\eta,\kappa}F_{\kappa'})
\right)
\lesssim
M_{\alpha,1}^\sigma
\left(
h,
A_{\eta,\kappa}F_{\kappa'}
\right)\lesssim
\left(
\frac{\alpha}{\kappa}
\right)^{2-\sigma}
\|A_{\eta,\kappa}\|_{L_{t,x}^\infty}
M_{\kappa,1}^\sigma
\left(
h,F_{\kappa'}
\right).
\end{split}
\end{equation*}
Here we used
\begin{equation*}
\alpha^2\langle\alpha y\rangle^{-\sigma}
\lesssim
\left(
\frac{\alpha}{\kappa}
\right)^{2-\sigma}
\kappa^2\langle\kappa y\rangle^{-\sigma}.
\end{equation*}
The sum over $\alpha\lesssim\kappa$, followed by the sum over $\kappa$, converges because $\sigma<2$ and because $\|A_\eta\|_{\dot{B}_{\infty,1}^0}\lesssim \epsilon$.
Summing the three interaction regimes proves
\eqref{weighted-elliptic-absorption}. Applying \eqref{flat-weighted-elliptic-estimate} and
\eqref{weighted-elliptic-absorption} to the identity \eqref{flattened-weighted-elliptic-equation}, we obtain
\begin{equation*}
\|W_\beta\|_{X_h^\sigma}
\lesssim_\sigma
\epsilon\|W_\beta\|_{X_h^\sigma}
+
\sum_{\alpha\lesssim\beta}
\sup_{z\leq0}
M_{\alpha,1}^\sigma
\left(
h,\widetilde P_\alpha\widetilde g(z)
\right).
\end{equation*}
For $\alpha\lesssim\beta$, boundedness of
$\widetilde P_\alpha$ and the same weight comparison as above gives
\begin{equation*}
M_{\alpha,1}^\sigma
\left(
h,\widetilde P_\alpha\widetilde g
\right)
\lesssim_\sigma
\left(
\frac{\alpha}{\beta}
\right)^{2-\sigma}
M_{\beta,1}^\sigma
\left(
h,\widetilde g
\right).
\end{equation*}
Since $\sigma<2$, the sum in $\alpha$ converges. Moreover,
$\|\nabla\eta\|_{L^\infty}\ll1$ implies
\begin{equation*}
M_{\beta,1}^\sigma
\left(
h,\widetilde g(z)
\right)
\lesssim
M_{\beta,1}^\sigma
\left(
h,g(z)
\right).
\end{equation*}
Absorbing the first term therefore yields
\begin{equation}\label{weighted-elliptic-interior-estimate}
\|W_\beta\|_{X_h^\sigma}
\lesssim_\sigma
\sup_{z\leq0}
M_{\beta,1}^\sigma
\left(
h,g(z)
\right).
\end{equation}
Finally, the reproducing-kernel estimate \eqref{reproducingbound} applied to $P_{\alpha}W_{\beta}$ at frequency $\alpha$
gives
\begin{equation*}
\left\|
hP_\alpha W_\beta(0)
\right\|_{L_t^1L_x^1}
\lesssim_\sigma
M_{\alpha,1}^\sigma
\left(
h,P_\alpha W_\beta(0)
\right).
\end{equation*}
The flattening preserves the vertical derivative (i.e. $\partial_z$ commutes with $f\mapsto f_*$), so
\begin{equation*}
\mathcal P f_{<\beta}
=
\left.\partial_zv_\beta\right|_{z=0}.
\end{equation*}
Summing in $\alpha$ and using
\eqref{weighted-elliptic-interior-estimate}, we conclude that
\begin{equation*}
\left\|
h\mathcal P f_{<\beta}
\right\|_{L_t^1L_x^1}
\lesssim_\sigma
\sup_{z\leq0}
M_{\beta,1}^\sigma
\left(
h,g(z)
\right),
\end{equation*}
as desired. Next, we prove \eqref{weighted-bilinear-exterior-derivative}, which is similar, but requires some minor adjustments. In this case, set
\begin{equation*}
v\coloneq\left(\Delta_D^{-1}f\right)_*,
\qquad
W\coloneq\nabla_{z,x}v.
\end{equation*}
Then we have
\begin{equation}\label{flattened-weighted-scalar-equation}
W
=
\mathcal S_0f_*
-
\mathcal R_0(A_\eta W),
\qquad
\mathcal S_0\coloneq
\nabla_{z,x}\Delta_{z,x}^{-1}.
\end{equation}
For each $\alpha\lesssim\beta$, we repeat the Littlewood--Paley
decomposition used in the proof of
\eqref{weighted-elliptic-absorption}, keeping the exterior derivative
on the output. The low--high interactions
$A_{\eta,<c\alpha}\widetilde P_\alpha W$ can be absorbed by smallness of $\nabla\eta$, as in the other bound. Every remaining interaction has $A_{\eta}$ at frequency $\kappa\gtrsim\alpha$. The resulting
frequency-$\alpha$ kernels are rapidly decreasing at scale
$\alpha^{-1}$, and hence
\begin{equation}\label{weighted-exterior-dyadic-bound}
\begin{split}
\left\|
h\nabla_xP_\alpha W(0)
\right\|_{L_t^1L_x^1}&\lesssim_\sigma
\sup_{z\leq0}
M_{\alpha,1}^{\sigma}
\left(h,\widetilde P_\alpha f_*(z)\right)
\\
&+
\|h_2\|_{L_t^\infty L_x^2}
\|W\|_{L_z^\infty L_t^2\dot B_{\infty,1}^0}
\sum_{\kappa\gtrsim\alpha}
\frac{\alpha}{\kappa}
M_{\alpha,2}^{\frac\sigma2}
\left(h_1,\widetilde P_\kappa\nabla_xA_\eta\right).
\end{split}
\end{equation}
Indeed, for the second term the $z'$-integrated kernel of
$\nabla_xP_\alpha\mathcal R_0$ is bounded by
$C_N\alpha^3\langle\alpha y\rangle^{-N}$. Cauchy--Schwarz in
$(t,x)$ and in the translation variable, followed by Bernstein's
inequality on $A_{\eta,\kappa}$, gives the factor
\begin{equation*}
\frac{\alpha}{\kappa}
M_{\alpha,2}^{\frac\sigma2}
\left(h_1,\widetilde P_\kappa\nabla_xA_\eta\right)
\|h_2\|_{L_t^\infty L_x^2}.
\end{equation*}
The remaining $W$-factor is bounded by its
$L_z^\infty L_t^2\dot B_{\infty,1}^0$ norm.  The preceding low--high iterations preserve the kernel bound and sum
geometrically by \eqref{weighted-elliptic-coefficient-smallness}. Using the same weight comparison as above for the first term, and
Proposition \ref{basicMproperties} together with
\eqref{M-LP-boundedness} for the second term, we obtain
\begin{equation*}
\begin{split}
\sum_{\kappa\gtrsim\alpha}
\frac{\alpha}{\kappa}
M_{\alpha,2}^{\frac\sigma2}
\left(h_1,\widetilde P_\kappa\nabla_xA_\eta\right)\lesssim_\sigma
\left(\frac{\alpha}{\beta}\right)^{1-\frac\sigma2}
M_{\beta,2}^{\frac\sigma2}
\left(h_1,\nabla_xA_\eta\right).
\end{split}
\end{equation*}
Both this expression and the first term in
\eqref{weighted-exterior-dyadic-bound} are summable over
$\alpha\lesssim\beta$, since
$1-\frac\sigma2>0$ and $2-\sigma>0$. By definition of $A_{\eta}$ and smallness of $\nabla\eta$, we have
\begin{equation*}
M_{\beta,2}^{\frac\sigma2}
\left(h_1,\nabla_xA_\eta\right)
\lesssim
M_{\beta,2}^{\frac\sigma2}
\left(h_1,\nabla^2\eta\right).
\end{equation*}
Using $\mathcal Pf=W_3(0)$ now proves
\eqref{weighted-bilinear-exterior-derivative}.
\end{proof}
Our next goal will be to show that these multi-linear norms are somewhat well-behaved when one of the entries is a suitably localized harmonic extension of a boundary function. We begin with the following proposition that handles the "diagonal" frequency case.
\begin{proposition}[Localized bilinear harmonic extension bound]\label{bilinearext}
Let $0<\sigma<1$, and let
$\beta\in \mathfrak{c}^{\mathbb{Z}}$. For $g_{\beta}=P_{\beta}g$, we have
\begin{equation}
\sup_{z\leq0}
M_{\beta,2}^{\sigma}
\left(f,\widetilde{P}_\beta w_{g_\beta}(z)\right)
\lesssim_\sigma
M_{\beta,2}^{\sigma}
\left(f,g_\beta \right).
\end{equation}
Moreover, there holds
\begin{equation}\label{bilinear-gradient-extension}
\begin{split}
\sup_{z\leq0}
M_{\beta,2}^{\sigma}
\left(f,\widetilde P_\beta\nabla_{z,x}w_{g_{\beta}}\right)
\lesssim_\sigma{}&
M_{\beta,2}^{\sigma}
\left(f,\nabla_xg_\beta\right)+\sup_{z\leq0}
M_{\beta,2}^{\sigma}
\left(f,P_{\gg\beta}\nabla_{z,x}w_{g_{\beta}}\right).
\end{split}
\end{equation}
The implicit constants are independent of $\beta$.
\end{proposition}

\begin{proof}
We begin with the first bound. Set
\begin{equation*}
W_{\beta,\sigma}(y)
\coloneq
\frac{\beta^2}{\langle\beta y\rangle^{2\sigma}},
\qquad
\omega_f(t,x)
\coloneq
\left(W_{\beta,\sigma}*|f(t,\cdot)|^2\right)(x).
\end{equation*}
By Fubini and a change of variables, for every function $F$ we have
\begin{equation}\label{weighted-M-identity}
\left(M_{\beta,2}^{\sigma}(f,F)\right)^2
=
\int_{\mathbb R}\int_{\mathbb R^2}
|F(t,x)|^2\omega_f(t,x)\,dx\,dt.
\end{equation}
Since $0<2\sigma<2$, the weight $W_{\beta,\sigma}$ belongs to
$A_1(\mathbb R^2)$, with a constant independent of $\beta$. In
particular,
\begin{equation*}
\mathcal M W_{\beta,\sigma}
\lesssim_\sigma
W_{\beta,\sigma},
\end{equation*}
where $\mathcal M$ is the Hardy--Littlewood maximal operator.
Consequently, using positivity,
\begin{equation*}
\mathcal M\omega_f(t,\cdot)
\leq
\left(\mathcal M W_{\beta,\sigma}\right)
*|f(t,\cdot)|^2
\lesssim_\sigma
\omega_f(t,\cdot).
\end{equation*}
Thus $\omega_f(t,\cdot)\in A_1\subset A_2$, uniformly in $t$ and
$\beta$.

The flattened harmonic extension satisfies
\begin{equation*}
w_{g_\beta}
=
e^{z|D|}g_\beta
+
\Delta_{z,x}^{-1}\nabla_{z,x}\cdot
\left(A_\eta\nabla_{z,x}w_{g_\beta}\right),
\end{equation*}
where $\Delta_{z,x}^{-1}$ denotes the Dirichlet inverse on the lower
half space. We extend $\omega_f(t,x)$ constantly in $z$. This extension
belongs to $A_2(\mathbb R^3)$ uniformly in $t$ and $\beta$. Applying
the weighted estimates pointwise in $t$, the flat half-space operators
satisfy
\begin{equation*}
\begin{split}
\sup_{z\leq0}
\left\|
\widetilde P_\beta e^{z|D|}g_\beta
\right\|_{L_{t,x}^2(\omega_f)}
&\lesssim_\sigma
\left\|g_\beta\right\|_{L_{t,x}^2(\omega_f)},
\\
\left\|
\nabla_{z,x}e^{z|D|}g_\beta
\right\|_{L_{z,t,x}^2(\omega_f)}
&\lesssim_\sigma
\beta^{\frac12}
\left\|g_\beta\right\|_{L_{t,x}^2(\omega_f)},
\\
\left\|
\nabla_{z,x}\Delta_{z,x}^{-1}
\nabla_{z,x}\cdot H
\right\|_{L_{z,t,x}^2(\omega_f)}
&\lesssim_\sigma
\left\|H\right\|_{L_{z,t,x}^2(\omega_f)},
\\
\sup_{z\leq0}
\left\|
\widetilde P_\beta\Delta_{z,x}^{-1}
\nabla_{z,x}\cdot H(z)
\right\|_{L_{t,x}^2(\omega_f)}
&\lesssim_\sigma
\beta^{-\frac12}
\left\|H\right\|_{L_{z,t,x}^2(\omega_f)}.
\end{split}
\end{equation*}
The third estimate follows from the reflection formula for the
Dirichlet kernel on the half space and the weighted $L^2$ bounds for
the Riesz transforms (using that $\omega_f\in A_2$). The remaining estimates follow from the
corresponding frequency-localized Poisson and Dirichlet kernel
estimates. The $\beta$-dependent factors arise from integration in
$z$.

Applying the second and third estimates to the equation for
$w_{g_\beta}$ gives
\begin{equation*}
\begin{split}
\left\|
\nabla_{z,x}w_{g_\beta}
\right\|_{L_{z,t,x}^2(\omega_f)}
\lesssim_\sigma{}&
\beta^{\frac12}
\left\|g_\beta\right\|_{L_{t,x}^2(\omega_f)}+
\|A_\eta\|_{L_{t,z,x}^\infty}
\left\|
\nabla_{z,x}w_{g_\beta}
\right\|_{L_{z,t,x}^2(\omega_f)}.
\end{split}
\end{equation*}
By the small-data assumption, the last term can be absorbed into the
left-hand side, and thus
\begin{equation}\label{weighted-gradient-bound}
\left\|
\nabla_{z,x}w_{g_\beta}
\right\|_{L_{z,t,x}^2(\omega_f)}
\lesssim_\sigma
\beta^{\frac12}
\left\|g_\beta\right\|_{L_{t,x}^2(\omega_f)}.
\end{equation}
Using the first and fourth flat estimates together with
\eqref{weighted-gradient-bound}, we obtain
\begin{equation*}
\begin{split}
\sup_{z\leq0}
\left\|
\widetilde P_\beta w_{g_\beta}(z)
\right\|_{L_{t,x}^2(\omega_f)}
&\lesssim_\sigma
\left\|g_\beta\right\|_{L_{t,x}^2(\omega_f)}+
\beta^{-\frac12}
\|A_\eta\|_{L_{t,z,x}^\infty}
\left\|
\nabla_{z,x}w_{g_\beta}
\right\|_{L_{z,t,x}^2(\omega_f)}\\
&\lesssim_\sigma
\left\|g_\beta\right\|_{L_{t,x}^2(\omega_f)}.
\end{split}
\end{equation*}
The conclusion now follows from \eqref{weighted-M-identity}. We next prove \eqref{bilinear-gradient-extension}. Let
$S_\beta=P_{\lesssim\beta}$ be a fattened projection equal to the identity on the Fourier support of
$\widetilde P_\beta$, and define $P_{\gg\beta}\coloneq I-S_\beta.$ Set for convenience
\begin{equation*}
U\coloneq\nabla_{z,x}w_{g_\beta}.
\end{equation*}
Recall from \eqref{weighted-M-identity} that
\begin{equation*}
\left(M_{\beta,2}^{\sigma}(f,F)\right)^2
=
\left\|F\right\|_{L_{t,x}^2(\omega_f)}^2
\end{equation*}
and that $\omega_f(t,\cdot)\in A_2$ uniformly in $t$ and $\beta$.

For $H:\{z<0\}\to\mathbb R^3$, the Dirichlet Green formula gives
\begin{equation*}
\nabla_{z,x}\Delta_{z,x}^{-1}\nabla_{z,x}\cdot H
=
\mathbf e_zH_3+\mathcal NH,
\end{equation*}
up to an irrelevant sign in front of the local term. Every kernel
occurring in $S_\beta\mathcal N$ is a finite linear combination of
operators of the form
\begin{equation*}
S_\beta R^a|D|e^{-\tau|D|},
\qquad
\tau=|z-z'|
\quad\text{or}\quad
\tau=|z|+|z'|,
\end{equation*}
where $R^a$ denotes a finite product Riesz
transforms. The weighted Poisson and Riesz-transform estimates give
\begin{equation*}
\left\|
S_\beta R^a|D|e^{-\tau|D|}
\right\|_{L_x^2(\omega_f(t))\to L_x^2(\omega_f(t))}
\lesssim_\sigma
\min\left\{\beta,\tau^{-1}\right\}.
\end{equation*}
Since
\begin{equation*}
\int_0^\infty
\min\left\{\beta,\tau^{-1}\right\}^2\,d\tau
\lesssim\beta,
\end{equation*}
Minkowski's inequality and Cauchy--Schwarz in $z'$ imply
\begin{equation}\label{localized-nonlocal-gradient-bound}
\sup_{z\leq0}
\left\|S_\beta\mathcal NH(z)\right\|_{L_{t,x}^2(\omega_f)}
\lesssim_\sigma
\beta^{\frac12}
\left\|H\right\|_{L_{z,t,x}^2(\omega_f)}.
\end{equation}
The weighted estimates from above give
\begin{equation}
\left\|U\right\|_{L_{z,t,x}^2(\omega_f)}
\lesssim_\sigma
\beta^{\frac12}
\left\|g_\beta\right\|_{L_{t,x}^2(\omega_f)}
\lesssim_\sigma
\beta^{-\frac12}
\left\|\nabla_xg_\beta\right\|_{L_{t,x}^2(\omega_f)}.
\label{weighted-bulk-gradient-derivative}
\end{equation}
On the other hand, the equation for $w_{g_\beta}$ gives
\begin{equation*}
U
=
\nabla_{z,x}e^{z|D|}g_\beta
+
\mathbf e_z(A_\eta U)_3
+
\mathcal N(A_\eta U).
\end{equation*}
Applying $S_\beta$, using
\eqref{localized-nonlocal-gradient-bound} and
\eqref{weighted-bulk-gradient-derivative}, and decomposing
\begin{equation*}
U=S_\beta U+P_{\gg\beta}U
\end{equation*}
in the local term, we obtain
\begin{equation*}
\begin{split}
\sup_{z\leq0}
\left\|S_\beta U(z)\right\|_{L_{t,x}^2(\omega_f)}
\lesssim_\sigma{}&
\left\|\nabla_xg_\beta\right\|_{L_{t,x}^2(\omega_f)}
+
\|A_\eta\|_{L_{t,z,x}^\infty}
\sup_{z\leq0}
\left\|S_\beta U(z)\right\|_{L_{t,x}^2(\omega_f)}
\\
&+
\|A_\eta\|_{L_{t,z,x}^\infty}
\sup_{z\leq0}
\left\|P_{\gg\beta}U(z)\right\|_{L_{t,x}^2(\omega_f)}.
\end{split}
\end{equation*}
By the small-data assumption, the second term on the right-hand side
can be absorbed. Since
$\widetilde P_\beta=\widetilde P_\beta S_\beta$, it follows from
\eqref{weighted-M-identity} that
\begin{equation}
\begin{split}
\sup_{z\leq0}
M_{\beta,2}^{\sigma}
\left(f,\widetilde P_\beta U(z)\right)
\lesssim_\sigma{}&
M_{\beta,2}^{\sigma}
\left(f,\nabla_xg_\beta\right)+
\|A_\eta\|_{L_{t,z,x}^\infty}
\sup_{z\leq0}
M_{\beta,2}^{\sigma}
\left(f,P_{\gg\beta}U(z)\right).
\end{split}
\end{equation}
Using $|A_{\eta}|\lesssim 1$, this proves
\eqref{bilinear-gradient-extension}.
\end{proof}
To compliment the "matched" frequency bound above we also need the following off-diagonal bound.
\begin{proposition}[Weighted off-diagonal harmonic-extension bound]
\label{weightedoffdiag}
Assume the bootstrap hypotheses for $u$. Let
$\beta,\lambda,\mu\in\mathfrak{c}^{\mathbb Z}$ satisfy $\lambda\ll\mu$.
Suppose also that $\eta=\eta_{<\lambda}$ and that
$g=g(t,x)$ is a Schwartz function in $x$, localized away from
dyadic scale $\beta$,
\begin{equation*}
g=
\left(
P_{\ll\beta}+P_{\gg\beta}
\right)g.
\end{equation*}
Let $2\leq p\leq\infty$ and $1<\sigma<2$. Then, for every
$x_0\in\mathbb R^2$, there holds
\begin{equation}\label{weighted-offdiag-main}
\begin{split}
\sup_{z\leq0}
M_{\beta,2}^{\frac{\sigma}{2}}
\left(
u_\mu^{x_0},
\nabla_{z,x}P_\beta w_g(z)
\right)\lesssim_{\sigma,p}
\epsilon^2
\beta^{-\frac14+\frac2p}
\mu^{-s_c}c_\mu c_\beta
\left(1+\mu|x_0|\right)^{\frac12}
\left\|
\nabla_xg
\right\|_{L_t^\infty\dot B_{p,1}^0}.
\end{split}
\end{equation}
\end{proposition}

\begin{proof}
Set
\begin{equation*}
W_g\coloneq\nabla_{z,x}w_g.
\end{equation*}
The flattened harmonic-extension equation takes the form
\begin{equation*}
W_g
=
e^{z|D|}
\left(
|D|g,\nabla_xg
\right)
+
\mathcal R_0
\left(
A_\eta W_g
\right).
\end{equation*}
Since the boundary datum is supported away from frequency $\beta$,
the flat term vanishes after applying $P_\beta$. Using sufficiently
fattened projections and fixing a sufficiently small $\delta>0$,
we therefore write
\begin{equation*}
\begin{split}
P_\beta W_g
&=
\widetilde P_\beta\mathcal R_0
\left(
P_{<\delta\beta}A_\eta P_\beta W_g
\right)
+
\widetilde P_\beta\mathcal R_0
\left(
[P_\beta,P_{<\delta\beta}A_\eta]W_g
\right)
\\
&\quad+
\widetilde P_\beta\mathcal R_0P_\beta
\left(
P_{\geq\delta\beta}A_\eta W_g
\right)
\\
&\coloneq
\mathcal T_1+\mathcal T_2+\mathcal T_3.
\end{split}
\end{equation*}
Smallness of $A_\eta$ in
$L_t^\infty\dot B_{\infty,1}^0$, the frequency-$\beta$ kernel
bound for $\mathcal R_0$ proved as in
\eqref{flat-weighted-elliptic-estimate}, Proposition
\ref{basicMproperties}, and H\"older's inequality give
\begin{equation}\label{weighted-offdiag-T1}
\begin{split}
\sup_{z\leq0}
M_{\beta,2}^{\frac{\sigma}{2}}
\left(
u_\mu^{x_0},
\mathcal T_1(z)
\right)
\lesssim
\epsilon
\sup_{z\leq0}
M_{\beta,2}^{\frac{\sigma}{2}}
\left(
u_\mu^{x_0},
P_\beta W_g(z)
\right),
\end{split}
\end{equation}
which can be absorbed.

We next estimate the commutator term. The frequency-localized
kernel bounds for $P_\beta$ and $\mathcal R_0$, together with the
fundamental theorem of calculus, give, for $a=a(t,x)$,
\begin{equation}\label{weighted-commutator-product}
\begin{split}
&\sup_{z\leq0}
M_{\beta,2}^{\frac{\sigma}{2}}
\left(
u,
\widetilde P_\beta\mathcal R_0
\left(
[P_\beta,a]F
\right)(z)
\right)\lesssim_{\sigma,p}
\beta^{-1+\frac2p}
M_{\beta,2}^{\frac{\sigma}{2}}
\left(
u,
\nabla_xa
\right)
\|F\|_{L_z^\infty L_t^\infty L_x^p}.
\end{split}
\end{equation}
Indeed, if $K_\beta$ denotes the kernel of $P_\beta$, then
\begin{equation*}
\begin{split}
[P_\beta,a]F(x)
&=
\int_{\mathbb R^2}
K_\beta(y)
\left(
a(x-y)-a(x)
\right)
F(x-y)\,dy
\\
&=
-\int_0^1\int_{\mathbb R^2}
K_\beta(y)y\cdot\nabla_xa(x-\vartheta y)
F(x-y)\,dy\,d\vartheta.
\end{split}
\end{equation*}
Let $L_\beta$ be an $L^1$-normalized, rapidly decreasing kernel at
scale $\beta^{-1}$ such that
\begin{equation*}
|y||K_\beta(y)|
\lesssim
\beta^{-1}L_\beta(y).
\end{equation*}
After Cauchy--Schwarz in $y$ and then in $\vartheta$, we obtain
\begin{equation*}
\begin{split}
\left|
[P_\beta,a]F(x)
\right|^2
\lesssim
\beta^{-2}
\int_0^1
&\left(
\int_{\mathbb R^2}
L_\beta(y)
|\nabla_xa(x-\vartheta y)|^2\,dy
\right)
\left(
\int_{\mathbb R^2}
L_\beta(y)
|F(x-y)|^2\,dy
\right)
d\vartheta.
\end{split}
\end{equation*}
For $2\leq p\leq\infty$, Young's inequality gives
\begin{equation*}
\left\|
L_\beta*|F|^2
\right\|_{L_x^\infty}
\lesssim_p
\beta^{\frac4p}
\|F\|_{L_x^p}^2.
\end{equation*}
Multiplying by $|u(t,x)|^2$, integrating against the weight in the
definition of $M_{\beta,2}^{\frac\sigma2}$, and using the weight
comparison from Proposition \ref{basicMproperties} and the rapid
decay of $L_\beta$ gives
\eqref{weighted-commutator-product}. The same argument applies after
composition with
$\widetilde P_\beta\mathcal R_0$, by its corresponding
frequency-$\beta$ kernel bound. Applying \eqref{weighted-commutator-product} with
$a=P_{<\delta\beta}A_\eta$ and $F=W_g$, we obtain
\begin{equation}\label{remainingMfactor}
\begin{split}
\sup_{z\leq0}
M_{\beta,2}^{\frac{\sigma}{2}}
\left(
u_\mu^{x_0},
\mathcal T_2(z)
\right)\lesssim_{\sigma,p}
\beta^{-1+\frac2p}
M_{\beta,2}^{\frac{\sigma}{2}}
\left(
u_\mu^{x_0},
\nabla_xP_{<\delta\beta}A_\eta
\right)
\|W_g\|_{L_z^\infty L_t^\infty L_x^p}.
\end{split}
\end{equation}
To estimate the bilinear factor on the right-hand side, we need the following simple lemma.
\begin{lemma}\label{Aetabilinear}
Assume the bootstrap bounds in Section \ref{sec:bootstrap},
and let $A_\eta$ be as above. For $1<\sigma<2$, $\kappa\ll\mu$ and
$x_0\in\mathbb R^2$, there holds
\begin{equation*}
M_{\kappa,2}^{\frac{\sigma}{2}}
(u_\mu^{x_0},\widetilde P_\kappa A_\eta)
\lesssim_\sigma
\epsilon^2c_\mu c_\kappa\mu^{-s_c}\kappa^{-\frac14}
(1+\mu|x_0|)^{\frac12},
\end{equation*}
\end{lemma}
where $\tilde{P}_{\kappa}$ is a fattened projection at frequency $\approx\kappa$.
\begin{proof}
Clearly this estimate holds for the factors of $A_{\eta}$ that are linear in $\nabla\eta$. To estimate the remaining factor $|\nabla\eta|^2$, we can use the Littlewood-Paley trichotomy to write
\begin{equation*}
\tilde{P}_{\kappa}|\nabla\eta|^2=2\tilde{P}_{\kappa}(\nabla\eta_{<\delta\kappa}\cdot\nabla\widehat{P}_{\kappa}\eta)+\sum_{\alpha\approx\gamma\gtrsim\kappa}\tilde{P}_{\kappa}(\nabla\eta_\alpha\cdot\nabla\eta_{\gamma}) 
\end{equation*}
where $\widehat{P}_{\kappa}$ is a fattening of $\tilde{P}_{\kappa}$. Using \eqref{M-LP-boundedness} and Lemma \ref{annularcontrol} and the bootstrap in Section \ref{sec:bootstrap}, we have
\begin{equation*}
M^{\frac{\sigma}{2}}_{\kappa,2}(u_{\mu}^{x_0},2\tilde{P}_{\kappa}(\nabla\eta_{<\delta\kappa}\cdot\nabla\widehat{P}_{\kappa}\eta))\lesssim \epsilon M^{\frac{\sigma}{2}}_{\kappa,2}(u_{\mu}^{x_0},\nabla\widehat{P}_{\kappa}\eta)\lesssim \epsilon^2 c_{\mu}c_{\kappa}\mu^{-s_c}\kappa^{-\frac{1}{4}}(1+\mu|x_0|)^{\frac12}.    
\end{equation*}
For the high-high interactions, choose $p$ with 
$4<p<4/(2-\sigma)$. Proposition \ref{basicMproperties},
Bernstein and the bootstrap bounds give
\begin{equation*}
\begin{split}
\sum_{\alpha\approx\gamma\gtrsim\kappa}
M_{\kappa,2}^{\frac{\sigma}{2}}
\left(u_\mu^{x_0},
\widetilde P_\kappa(\partial_i\eta_\alpha\partial_j\eta_\gamma)
\right)&\lesssim
\kappa^{\frac2p}\|u_\mu\|_{L_t^\infty L_x^2}
\sum_{\alpha\approx\gamma\gtrsim\kappa}
\|\nabla\eta_\alpha\|_{L_t^4L_x^{2p}}
\|\nabla\eta_\gamma\|_{L_t^4L_x^{2p}}
\\
&\lesssim
\epsilon^3c_\mu\mu^{-s_c}\kappa^{\frac2p}
\sum_{\alpha\gtrsim\kappa}(c_\alpha)^2
\alpha^{-\frac14-\frac2p}
\\
&\lesssim_\sigma
\epsilon^3c_\mu c_\kappa\mu^{-s_c}\kappa^{-\frac14}.
\end{split}
\end{equation*}
\end{proof}
Returning to estimating \eqref{remainingMfactor}, we can use Proposition \ref{basicMproperties} and Lemma \ref{Aetabilinear} to obtain
\begin{equation*}
\begin{split}
\beta^{-1}
M_{\beta,2}^{\frac{\sigma}{2}}
\left(
u_\mu^{x_0},
\nabla_xP_{<\delta\beta}A_\eta
\right)&\lesssim
\beta^{-1}
\sum_{\kappa<\delta\beta}
\left(
\frac{\beta}{\kappa}
\right)^{1-\frac{\sigma}{2}}
\kappa
M_{\kappa,2}^{\frac{\sigma}{2}}
\left(
u_\mu^{x_0},
\widetilde P_\kappa A_\eta
\right)
\\
&\quad\lesssim
\epsilon^2
\beta^{-\frac14}
\mu^{-s_c}c_\mu
\left(1+\mu|x_0|\right)^{\frac12}
\sum_{\kappa<\delta\beta}
\left(
\frac{\kappa}{\beta}
\right)^{\frac{\sigma}{2}-\frac14}
c_\kappa
\\
&\quad\lesssim_\sigma
\epsilon^2
\beta^{-\frac14}
c_\beta
\mu^{-s_c}c_\mu
\left(1+\mu|x_0|\right)^{\frac12}.
\end{split}
\end{equation*}
Here the last sum converges because $\sigma>1$ and by the slowly
varying property of the frequency envelope. Minkowski's inequality and the flattened harmonic
extension estimate Proposition \ref{flattenedbasicbounds} gives
\begin{equation*}
\left\|
W_g
\right\|_{L_z^\infty L_t^\infty L_x^p}
\lesssim
\left\|
W_g
\right\|_{L_z^\infty L_t^\infty\dot B_{p,1}^0}
\lesssim
\left\|
\nabla_xg
\right\|_{L_t^\infty\dot B_{p,1}^0}.
\end{equation*}
Consequently,
\begin{equation}\label{weighted-offdiag-T2}
\begin{split}
\sup_{z\leq0}
M_{\beta,2}^{\frac{\sigma}{2}}
\left(
u_\mu^{x_0},
\mathcal T_2(z)
\right)\lesssim_{\sigma,p}
\epsilon^2
\beta^{-\frac14+\frac2p}
c_\beta
\mu^{-s_c}c_\mu
\left(1+\mu|x_0|\right)^{\frac12}
\left\|
\nabla_xg
\right\|_{L_t^\infty\dot B_{p,1}^0}.
\end{split}
\end{equation}

We finally estimate $\mathcal T_3$. From the localization assumption
on $\eta$, this term vanishes unless
$\beta\lesssim\lambda\ll\mu$. We use the output-product
estimate
\begin{equation}\label{weighted-output-product}
M_{\beta,2}^{\frac{\sigma}{2}}
\left(
u,
P_\beta(FG)
\right)
\lesssim_{\sigma,p}
\beta^{\frac2p}
M_{\beta,2}^{\frac{\sigma}{2}}
\left(
u,F
\right)
\|G\|_{L_t^\infty L_x^p}.
\end{equation}
Indeed, Cauchy--Schwarz in the convolution variable gives
\begin{equation*}
\begin{split}
\left|
P_\beta(FG)(t,x+y)
\right|^2
&\lesssim
\left(
|K_\beta|*|F|^2
\right)(t,x+y)
\left(
|K_\beta|*|G|^2
\right)(t,x+y),
\end{split}
\end{equation*}
while Young's inequality gives
\begin{equation*}
\left\|
|K_\beta|*|G|^2
\right\|_{L_x^\infty}
\lesssim_p
\beta^{\frac4p}
\|G\|_{L_x^p}^2.
\end{equation*}
Multiplying by $|u(t,x)|^2$, integrating in $t,x,y$ against the
weight defining $M_{\beta,2}^{\frac{\sigma}{2}}$, and using the
rapid decay of $K_\beta$ proves
\eqref{weighted-output-product}. The same estimate holds with
$\widetilde P_\beta\mathcal R_0P_\beta$ in place of $P_\beta$,
by the corresponding frequency-$\beta$ kernel bound for
$\mathcal R_0$. Therefore,
\begin{equation*}
\begin{split}
\sup_{z\leq0}
M_{\beta,2}^{\frac{\sigma}{2}}
\left(
u_\mu^{x_0},
\mathcal T_3(z)
\right)\lesssim_p
\beta^{\frac2p}
\left\|
W_g
\right\|_{L_z^\infty L_t^\infty\dot B_{p,1}^0}
\sum_{\kappa\gtrsim\beta}
M_{\beta,2}^{\frac{\sigma}{2}}
\left(
u_\mu^{x_0},
\widetilde P_\kappa A_\eta
\right).
\end{split}
\end{equation*}
By Proposition \ref{basicMproperties} and Lemma \ref{Aetabilinear} we have
\begin{equation*}
\begin{split}
M_{\beta,2}^{\frac{\sigma}{2}}
\left(
u_\mu^{x_0},
\widetilde P_\kappa A_\eta
\right)&\lesssim
\left(
\frac{\beta}{\kappa}
\right)^{1-\frac{\sigma}{2}}
M_{\kappa,2}^{\frac{\sigma}{2}}
\left(
u_\mu^{x_0},
\widetilde P_\kappa A_\eta
\right)
\\
&\lesssim
\epsilon^2
\left(
\frac{\beta}{\kappa}
\right)^{1-\frac{\sigma}{2}}
\kappa^{-\frac14}
c_\kappa
\mu^{-s_c}c_\mu
\left(1+\mu|x_0|\right)^{\frac12}.
\end{split}
\end{equation*}
Here $\mu\gg\lambda$ and $\eta=\eta_{<\lambda}$ ensure that the
bilinear estimate is applied in the unbalanced frequency setting.
Since $\sigma<2$, the slowly varying property of the frequency
envelope yields
\begin{equation*}
\begin{split}
\sum_{\kappa\gtrsim\beta}
\left(
\frac{\beta}{\kappa}
\right)^{1-\frac{\sigma}{2}}
\kappa^{-\frac14}
c_\kappa=
\beta^{-\frac14}
\sum_{\kappa\gtrsim\beta}
\left(
\frac{\beta}{\kappa}
\right)^{\frac54-\frac{\sigma}{2}}
c_\kappa\lesssim_\sigma
\beta^{-\frac14}c_\beta.
\end{split}
\end{equation*}
Consequently,
\begin{equation}\label{weighted-offdiag-T3}
\begin{split}
\sup_{z\leq0}
M_{\beta,2}^{\frac{\sigma}{2}}
\left(
u_\mu^{x_0},
\mathcal T_3(z)
\right)\lesssim_{\sigma,p}
\epsilon^2
\beta^{-\frac14+\frac2p}
c_\beta
\mu^{-s_c}c_\mu
\left(1+\mu|x_0|\right)^{\frac12}
\left\|
\nabla_xg
\right\|_{L_t^\infty\dot B_{p,1}^0}.
\end{split}
\end{equation}
\end{proof}

\subsection{Multi-linear Commutator Estimates}
Here, we collect several commutator estimates we will need for our analysis. These are a bit more delicate than standard linear commutator estimates because we have to preserve in a sufficiently accurate way, the multi-linear structure in these commutators to use the bilinear estimates in Section \ref{sec:bootstrap}. 
\\

We begin with the following analogue of Proposition \ref{Bernsteinelliptic} which will enable us to control any high-frequency tails perturbatively.
\begin{proposition}\label{highfreqtailcomm}
Let $\lambda,\mu>0$ satisfy $\mu\lesssim\lambda$.
Let $\sigma\geq0$, $\theta\in[0,1]$, and $N\in\mathbb N$, $N\geq1$,
and suppose that
\begin{equation*}
1<p,q_1,q_2\leq\infty,
\hspace{5mm}
1<r\leq\infty,
\hspace{5mm}
\frac1p=\frac1r+\frac1{q_1}+\frac1{q_2}.
\end{equation*}
Then, if $M$ is sufficiently large relative to $N$ and the fixed
frequency cutoffs, there holds
\begin{equation}\label{commutatorsummation}
\begin{split}
\sum_{\alpha<\lambda}
\left\|
|D|^\sigma P_{>M\lambda}
[A_{<\mu},f_{<\delta\alpha}]g_\alpha
\right\|_{L^p}
\lesssim{}&
\lambda^{\sigma-N}
\|\nabla_x^2\eta_{<\mu}\|_{\dot B^0_{rN,1}}^N
\||D|^{1-\theta}f_{<\lambda}\|_{\dot B^0_{q_1,1}}
\||D|^\theta g_{<\lambda}\|_{\dot B^0_{q_2,1}},
\end{split}
\end{equation}
where $A_{<\mu}$ denotes either $L_{<\mu}$, $M_{<\mu}$, or
$G(\eta_{<\mu})$.
\end{proposition}

The idea is that the output frequency $>M\lambda$ is much larger than
the input frequencies of $g_\alpha$, $f_{<\delta\alpha}$, and
$\eta_{<\mu}$. Thus, the only contribution can come from the nonlinear
dependence of $A_{<\mu}$ on $\nabla_x\eta_{<\mu}$, which produces,
with sufficiently large frequency separation, powers of
$\nabla_x^2\eta_{<\mu}$.

\begin{proof}
We abbreviate below $\eta_{<\mu}$ by just $\eta$. We show the details for $L_{<\mu}$. The other cases will follow from
the formulas relating $G(\eta_{<\mu})$ and $M_{<\mu}$ to
$L_{<\mu}$.  Now, we proceed with the proof. First, we note that by Bernstein, it suffices to consider $\sigma\geq1$. Now fix $\alpha<\lambda$, and let $w_{g_\alpha}$ and
$w_{f_{<\delta\alpha}g_\alpha}$ denote the flattened harmonic extensions
relative to $\eta_{<\mu}$. Set
\begin{equation*}
u_{\mathrm{com},\alpha}
\coloneq
w_{f_{<\delta\alpha}g_\alpha}
-f_{<\delta\alpha}w_{g_\alpha}.
\end{equation*}
Since
\begin{equation*}
u_{\mathrm{com},\alpha}\big|_{z=0}=0,
\qquad
\partial_zu_{\mathrm{com},\alpha}\big|_{z=0}
=
[L_{<\mu},f_{<\delta\alpha}]g_\alpha,
\end{equation*}
we therefore have
\begin{equation}\label{comm-tail1}
\left\|
|D|^\sigma P_{>M\lambda}
[L_{<\mu},f_{<\delta\alpha}]g_\alpha
\right\|_{L^p}
\lesssim
X_{p,\sigma}
\left(P_{>M\lambda}u_{\mathrm{com},\alpha}\right).
\end{equation}
where we define
\begin{equation*}
X_{p,\beta}(v)
\coloneq
\||D|^{\beta}\nabla_{z,x}v\|_{L_z^\infty\dot B^0_{p,1}}
+
\||D|^{\beta}\partial_z\nabla_{z,x}v\|_{L_z^1\dot B^0_{p,1}}.
\end{equation*}
We used this norm implicitly in the proof of Proposition \ref{Bernsteinelliptic} to efficiently perturb many of our estimates from those satisfied by the flat Laplacian $\Delta_{z,x}$. The equation satisfied by $\nabla_{z,x}u_{\mathrm{com},\alpha}$ is then given by
\begin{equation}\label{commtailequation}
\begin{split}
\nabla_{z,x}u_{\mathrm{com},\alpha}
&=\mathcal{R}_0 (A_{\eta}\nabla_{z,x}u_{\mathrm{com},\alpha})+\mathcal{S}_0F_{\mathrm{com,\alpha}}
\end{split}
\end{equation}
where $\mathcal{R}_0$ is as in Proposition \ref{flattenedbasicbounds} and $\mathcal{S}_0\coloneq \nabla_{z,x}\Delta_{z,x}^{-1}$. Above, we have
\begin{equation}\label{comm-tail-source}
\begin{split}
F_{\mathrm{com},\alpha}
\coloneq
2(\nabla_x\eta_{<\mu}\cdot\nabla_xf_{<\delta\alpha})
\partial_zw_{g_\alpha}
-2\nabla_xf_{<\delta\alpha}\cdot\nabla_xw_{g_\alpha}
-(\Delta_xf_{<\delta\alpha})w_{g_\alpha}.
\end{split}
\end{equation}
For $0\leq j\leq N$, define $p_j$ and $m_j$ by
\begin{equation*}
\frac1{p_j}
\coloneq
\frac{N-j}{rN}+\frac1{q_1}+\frac1{q_2},
\qquad
\frac1{m_j}
\coloneq
\frac{N-j}{rN}+\frac1{q_2}.
\end{equation*}
Thus $p_0=p$, and
\begin{equation*}
\frac1{p_j}=\frac1{rN}+\frac1{p_{j+1}},
\qquad
\frac1{p_j}=\frac1{q_1}+\frac1{m_j},
\qquad
\frac1{p_N}=\frac1{q_1}+\frac1{q_2}.
\end{equation*}
For $\beta\geq1$, write
\begin{equation*}
\mathcal E_{p,\beta}(F)
\coloneq
\||D|^{\beta-1}F\|_{L_z^\infty\dot B^0_{p,1}}
+
\||D|^{\beta-1}\partial_zF\|_{L_z^1\dot B^0_{\beta,1}}
\end{equation*}
which we will use to measure the error coming from the source terms in \eqref{commtailequation}. As in the proof of Proposition \ref{Bernsteinelliptic}, choose
\begin{equation*}
1\ll M_1\ll M_2\ll\cdots\ll M_N=M,
\end{equation*}
with enough separation to accommodate all the harmlessly fattened
projections below. Let $\widetilde P_{M_j\lambda}$ denote a smooth
projection to the transition annulus of $P_{>M_j\lambda}$, chosen so
that
\begin{equation*}
\operatorname{supp}\widetilde P_{M_j\lambda}
\subset
\{\xi:|\xi|>M_{j-1}\lambda\},
\qquad 2\leq j\leq N.
\end{equation*}
The localized estimate and commutator calculation used in the
proof of \eqref{localizedWfestimate}, now applied to \eqref{commtailequation}, give, for $0\leq j<N$,
\begin{equation}\label{commtailfirstbound}
\begin{split}
X_{p_j,\sigma}
\left(P_{>M_{N-j}\lambda}u_{\mathrm{com},\alpha}\right)
\lesssim{}&
\lambda^{-1}
\|\nabla_x^2\eta_{<\mu}\|_{\dot B^0_{rN,1}}
X_{p_{j+1},\sigma}
\left(\widetilde P_{M_{N-j}\lambda}u_{\mathrm{com},\alpha}\right)
\\
&+
\mathcal E_{p_j,\sigma}
\left(\widetilde P_{>M_{N-j}\lambda}F_{\mathrm{com},\alpha}\right).
\end{split}
\end{equation}
Here and below, a tilde on a high-frequency projection denotes a
fixed fattening of its cutoff. The first term in
\eqref{commtailfirstbound} comes from commuting the projection through
the coefficients; the terms in which the projection remains on
$u_{\mathrm{com},\alpha}$ are absorbed using the usual argument using the smallness
assumption on $\nabla_x\eta_{<\mu}$.

Because $\eta_{<\mu}$ and $f_{<\delta\alpha}$ have frequency
$\lesssim\lambda$, the projection in the second line of
\eqref{commtailfirstbound} forces the $w_{g_\alpha}$ factor into the
corresponding high-frequency region. H\"older's inequality,
Proposition \ref{besovparabounds}, and Bernstein's inequality therefore
give
\begin{equation*}
\begin{split}
\mathcal E_{p_j,\sigma}
\left(\widetilde P_{>M_{N-j}\lambda}F_{\mathrm{com},\alpha}\right)
\lesssim
\|\nabla_xf_{<\delta\alpha}\|_{\dot B^0_{q_1,1}}
X_{m_j,\sigma-1}
\left(\widehat P_{>M_{N-j}\lambda}w_{g_\alpha}\right),
\end{split}
\end{equation*}
where $\widehat P_{>M_{N-j}\lambda}$ is another fixed fattening. Its
cutoff is a fixed fraction of $M_{N-j}\lambda$, chosen still above
$M_{N-j-1}\lambda$ when $N-j\geq2$. Provided $M$ is large enough, the argument used in the proof of Proposition \ref{Bernsteinelliptic} gives
\begin{equation*}
\begin{split}
X_{m_j,\sigma-1}
\left(\widehat P_{>M_{N-j}\lambda}w_{g_\alpha}\right)
\lesssim{}&
\lambda^{-(N-j)}
\|\nabla_x^2\eta_{<\mu}\|_{\dot B^0_{rN,1}}^{N-j}
X_{q_2,\sigma-1}
\left(\widetilde P_{M_1\lambda}w_{g_\alpha}\right).
\end{split}
\end{equation*}
The terminal factor is localized to the annulus
$|\xi|\approx M_1\lambda$. Hence Bernstein and the standard flattened
harmonic-extension estimates imply
\begin{equation*}
\begin{split}
X_{q_2,\sigma-1}
\left(\widetilde P_{M_1\lambda}w_{g_\alpha}\right)
&\lesssim
\lambda^{\sigma-1}X_{q_2,0}(w_{g_\alpha})
\lesssim
\alpha\lambda^{\sigma-1}
\|g_\alpha\|_{\dot B^0_{q_2,1}}
\lesssim
\lambda^\sigma
\|g_\alpha\|_{\dot B^0_{q_2,1}}.
\end{split}
\end{equation*}
Consequently,
\begin{equation}\label{Fcommbound}
\begin{split}
\mathcal E_{p_j,\sigma}
\left(\widetilde P_{>M_{N-j}\lambda}F_{\mathrm{com},\alpha}\right)
\lesssim{}&
\lambda^{\sigma-(N-j)}
\|\nabla_x^2\eta_{<\mu}\|_{\dot B^0_{rN,1}}^{N-j}
\|\nabla_xf_{<\delta\alpha}\|_{\dot B^0_{q_1,1}}
\|g_\alpha\|_{\dot B^0_{q_2,1}}.
\end{split}
\end{equation}

Iterating \eqref{commtailfirstbound} and using the support properties
of the transition projections gives
\begin{equation}\label{penultimateudif}
\begin{split}
X_{p,\sigma}
\left(P_{>M\lambda}u_{\mathrm{com},\alpha}\right)
&\lesssim
\lambda^{-N}
\|\nabla_x^2\eta_{<\mu}\|_{\dot B^0_{rN,1}}^N
X_{p_N,\sigma}
\left(\widetilde P_{M_1\lambda}u_{\mathrm{com},\alpha}\right)
\\
&+
\sum_{j=0}^{N-1}
\lambda^{-j}
\|\nabla_x^2\eta_{<\mu}\|_{\dot B^0_{rN,1}}^j
\mathcal E_{p_j,\sigma}
\left(\widetilde P_{>M_{N-j}\lambda}F_{\mathrm{com},\alpha}\right).
\end{split}
\end{equation}
Using \eqref{Fcommbound} for the second line, it remains only to control the "terminal" factor in the first line. To
this end, we first write part of $F_{\mathrm{com},\alpha}$ as a divergence,
\begin{equation*}
F_{\mathrm{com},\alpha}
=
\nabla_{z,x}\cdot G_{\mathrm{com},\alpha}
+H_{\mathrm{com},\alpha},
\end{equation*}
where, 
\begin{equation*}
\begin{split}
G_{\mathrm{com},\alpha}
&\coloneq
\left(
2(\nabla_x\eta_{<\mu}\cdot\nabla_xf_{<\delta\alpha})w_{g_\alpha},
-w_{g_\alpha}\nabla_xf_{<\delta\alpha}
\right),
\\
H_{\mathrm{com},\alpha}
&\coloneq
-\nabla_xf_{<\delta\alpha}\cdot\nabla_xw_{g_\alpha}.
\end{split}
\end{equation*}
The usual small data absorption argument (used several times in this section) followed by the standard estimates for $\mathcal{S}_0\nabla_{z,x}\cdot $ and $\mathcal{S}_0$, we obtain
\begin{equation}\label{udifelliptic}
\begin{split}
X_{p_N,0}(u_{\mathrm{com},\alpha})
\lesssim{}&
\|G_{\mathrm{com},\alpha}\|_{L_z^\infty\dot B^0_{p_N,1}}
+
\|\partial_zG_{\mathrm{com},\alpha}\|_{L_z^1\dot B^0_{p_N,1}}
+
\|H_{\mathrm{com},\alpha}\|_{L_z^1\dot B^0_{p_N,1}}.
\end{split}
\end{equation}
In this final estimate, none of the factors are necessarily localized to frequency $\approx M\lambda$, so we have to handle low frequencies as well. Let us fix $\delta_0\ll c\ll1$ and set
$Q_\alpha=P_{<c\alpha}$ and $V_\alpha=Q_\alpha w_{g_\alpha}$.
Choose $k>2$ sufficiently large and define
\begin{equation*}
a=\frac2k,\qquad
\frac1b=\frac1{q_2}+\frac1k,\qquad
\frac1{p_N}+\frac1k<1.
\end{equation*}
Bernstein, Proposition \ref{MoserBesov} and the smallness of $\nabla\eta$ give $\|A_{\eta_{<\mu}}\|_{\dot B_{k,1}^{a}}\lesssim\epsilon$. Moreover, Proposition \ref{Leibnizcommutator} and Proposition
\ref{besovparabounds} give, uniformly in $\alpha$,
\begin{equation*}
\left\|
\widetilde Q_\alpha[Q_\alpha,A_{\eta_{<\mu}}]U
\right\|_{\dot B_{b,1}^{a}}
\lesssim
\|A_{\eta_{<\mu}}\|_{\dot B_{2,1}^{1}}
\|U\|_{\dot B_{q_2,1}^{0}}.
\end{equation*}
Here $\widetilde Q_\alpha$ is a fattened low-frequency projection.
Commuting $Q_\alpha$ through the flattened equation, and using that
the flat extension vanishes after applying $Q_\alpha$, gives
\begin{equation*}
\begin{split}
\nabla_{z,x}V_\alpha=-\mathcal R_0\widetilde Q_\alpha
\left(A_{\eta_{<\mu}}\nabla_{z,x}V_\alpha\right)-\mathcal R_0\widetilde Q_\alpha
\left([Q_\alpha,A_{\eta_{<\mu}}]
\nabla_{z,x}w_{g_\alpha}\right).
\end{split}
\end{equation*}
Absorbing the first term and using the above commutator bound, we
obtain
\begin{equation*}
\begin{split}
\|\nabla_{z,x}V_\alpha\|_{L_z^1\dot B_{b,1}^{a}}
\lesssim
\epsilon
\|\nabla_{z,x}w_{g_\alpha}\|_{L_z^1\dot B_{q_2,1}^{0}}
\lesssim
\epsilon\|g_\alpha\|_{\dot B_{q_2,1}^{0}}.
\end{split}
\end{equation*}
Since $V_\alpha(0)=0$, the fundamental theorem of calculus in $z$
implies the same bound for
$\|V_\alpha\|_{L_z^\infty\dot B_{b,1}^{a}}$. The same commutator
bound with $U=\nabla_{z,x}w_{g_\alpha}$ and
$U=\partial_z\nabla_{z,x}w_{g_\alpha}$, along with the standard Poisson bound
\begin{equation*}
X_{q_2,0}(w_{g_\alpha})
\lesssim
\alpha\|g_\alpha\|_{\dot B_{q_2,1}^{0}},
\end{equation*}
we get the same estimate as above with $V_{\alpha}$ replaced by $\alpha^{-1}\partial_z V_{\alpha}$, and directly by Bernstein, the same estimate for $\alpha^{-1}\nabla_x V_{\alpha}$. Absorbing the coefficient terms as before, we therefore obtain
\begin{equation}\label{comm-tail-low-extension}
\begin{split}
\|V_\alpha\|_{L_z^\infty\dot B_{b,1}^{a}}
+\|\nabla_{z,x}V_\alpha\|_{L_z^1\dot B_{b,1}^{a}}+
\alpha^{-1}\left(
\|\nabla_{z,x}V_\alpha\|_{L_z^\infty\dot B_{b,1}^{a}}
+\|\partial_z\nabla_{z,x}V_\alpha\|_{L_z^1\dot B_{b,1}^{a}}
\right)
\lesssim
\epsilon\|g_\alpha\|_{\dot B_{q_2,1}^{0}}.
\end{split}
\end{equation}
Now, using the product bound
\begin{equation*}
\|h_1h_2\|_{\dot B_{p_N,1}^{0}}
\lesssim
\|h_1\|_{\dot B_{q_1,1}^{0}}
\|h_2\|_{\dot B_{b,1}^{a}},
\end{equation*}
which follows from Proposition \ref{besovparabounds}, we estimate
the products involving $V_\alpha$ and its first derivatives.
For $(I-Q_\alpha)w_{g_\alpha}$, the frequency separation from
$f_{<\delta\alpha}$ gives the corresponding product bounds directly,
using
\begin{equation*}
\begin{split}
\|w_{g_\alpha}\|_{L_z^\infty\dot B_{q_2,1}^{0}}
+\|\nabla_{z,x}w_{g_\alpha}\|_{L_z^1\dot B_{q_2,1}^{0}}
\lesssim
\|g_\alpha\|_{\dot B_{q_2,1}^{0}}.
\end{split}
\end{equation*}
Multiplication by $\nabla_x\eta_{<\mu}$ is bounded on
$\dot B_{p_N,1}^{0}$ by its $\dot B_{2,1}^{1}$ norm.
Using these observations in \eqref{udifelliptic}, we conclude that
\begin{equation}\label{similarestimatescomm}
\begin{split}
X_{p_N,0}(u_{\mathrm{com},\alpha})
\lesssim
\|\nabla_xf_{<\delta\alpha}\|_{\dot B^0_{q_1,1}}
\|g_\alpha\|_{\dot B^0_{q_2,1}}.
\end{split}
\end{equation}
This proves \eqref{commutatorsummation} for $L_{<\mu}$. For the other two operators, we use
\begin{equation*}
\begin{split}
[G(\eta_{<\mu}),f]g
&=
(1+|\nabla_x\eta_{<\mu}|^2)[L_{<\mu},f]g
-(\nabla_x\eta_{<\mu}\cdot\nabla_xf)g,
\\
[M_{<\mu},f]g
&=
(\nabla_xf)g
-\nabla_x\eta_{<\mu}[L_{<\mu},f]g.
\end{split}
\end{equation*}
For $M$ sufficiently large, the explicit product terms vanish after
applying $P_{>M\lambda}$. The remaining terms reduce, after harmlessly
fattening the high-frequency projection, to the estimate already
proved for $L_{<\mu}$. This completes the proof.
\end{proof}
We will also need the following simpler commutator estimate when the
frequency of $\eta_{<\mu}$ is smaller than the input
frequency.

\begin{proposition}\label{highfreqcommutator}
Let $\lambda>0$ and $\mu\lesssim\lambda$. Let
$1<p,q_1,q_2\leq\infty$ satisfy
$\frac1p=\frac1{q_1}+\frac1{q_2}$. For every $\sigma\geq0$, there
holds, with a sufficiently small, fixed $\delta$,
\begin{equation*}
\left\|
|D|^\sigma[A_{<\mu},f_{<\delta\lambda}]g_{\geq\lambda}
\right\|_{L^p}
\lesssim_\sigma
\|\nabla_xf_{<\lambda}\|_{\dot B^0_{q_1,1}}
\||D|^\sigma g_{\geq\lambda}\|_{\dot B^0_{q_2,1}},
\end{equation*}
where $A_{<\mu}$ denotes either $L_{<\mu}$, $M_{<\mu}$, or
$G(\eta_{<\mu})$.
\end{proposition}

\begin{proof}
By summation, it suffices to prove the bound for $g_\alpha$ in place of $g_{\geq\lambda}$, where
$\alpha\geq\lambda$. We first consider $A_{<\mu}=L_{<\mu}$ and use
the notation from the preceding proof. The argument leading to
\eqref{similarestimatescomm} gives
\begin{equation}\label{highfreqcommutator-zero}
\begin{split}
X_{p,0}(u_{\mathrm{com},\alpha})
\lesssim
\|\nabla_xf_{<\delta\lambda}\|_{\dot B^0_{q_1,1}}
\|g_\alpha\|_{\dot B^0_{q_2,1}}.
\end{split}
\end{equation}
For $\sigma\geq1$, the same elliptic estimate as above
gives
\begin{equation*}
\begin{split}
X_{p,\sigma}(u_{\mathrm{com},\alpha})
\lesssim{}&
\mathcal E_{p,\sigma}(F_{\mathrm{com},\alpha})
+\epsilon X_{p,\sigma}(u_{\mathrm{com},\alpha})
+\epsilon\lambda^\sigma X_{p,0}(u_{\mathrm{com},\alpha}).
\end{split}
\end{equation*}
Commuting horizontal derivatives through the flattened equation and
using $\mu\lesssim\lambda\leq\alpha$ gives
\begin{equation*}
\begin{split}
X_{q_2,\sigma-1}(w_{g_\alpha})
\lesssim_\sigma
\alpha^\sigma\|g_\alpha\|_{\dot B^0_{q_2,1}}.
\end{split}
\end{equation*}
The coefficient terms are handled as above, with derivatives of
order $\sigma-1$ costing at most $\alpha^{\sigma-1}$.
For the low-frequency part of $w_{g_\alpha}$ we use
\eqref{comm-tail-low-extension}; for the remaining part we use the
separation from $f_{<\delta\lambda}$ and the preceding estimate.
Proposition \ref{besovparabounds} therefore gives
\begin{equation*}
\begin{split}
\mathcal E_{p,\sigma}(F_{\mathrm{com},\alpha})
\lesssim
\alpha^\sigma
\|\nabla_xf_{<\delta\lambda}\|_{\dot B^0_{q_1,1}}
\|g_\alpha\|_{\dot B^0_{q_2,1}}.
\end{split}
\end{equation*}
After absorption, \eqref{highfreqcommutator-zero} and
$\lambda\leq\alpha$ therefore imply
\begin{equation*}
\begin{split}
X_{p,\sigma}(u_{\mathrm{com},\alpha})
\lesssim
\alpha^\sigma
\|\nabla_xf_{<\delta\lambda}\|_{\dot B^0_{q_1,1}}
\|g_\alpha\|_{\dot B^0_{q_2,1}}.
\end{split}
\end{equation*}
Taking the boundary trace proves the corresponding estimate
for $[L_{<\mu},f_{<\delta\lambda}]g_\alpha$. The range
$0<\sigma<1$ follows by interpolation between the cases $\sigma=0$
and $\sigma=1$, and summation in $\alpha\geq\lambda$ proves the
stated estimate for $L_{<\mu}$. The estimates for
$G(\eta_{<\mu})$ and $M_{<\mu}$ follow from the same identities used
at the end of the preceding proof and Proposition
\ref{besovparabounds}. Indeed, the explicit product terms satisfy the
same dyadic estimate since $\mu\lesssim\lambda\leq\alpha$.
For the terms in which a coefficient depending on
$\eta_{<\mu}$ is multiplied by $[L_{<\mu},f_{<\delta\lambda}]g_\alpha$, we use
the order-$\sigma$ estimate above when the derivatives fall on the
commutator, and \eqref{highfreqcommutator-zero} together with
$\mu^\sigma\lesssim\alpha^\sigma$ when they fall on the coefficient.
Summing in $\alpha\geq\lambda$ completes the proof.
\end{proof}
Next, we will proceed with the establishing the multi-linear commutator estimates we will need in our analysis. First, we need the following lemma that gives an off-diagonal bound for the flattened harmonic extension.
\begin{lemma}\label{wgtailbound}
Let $g_\beta=P_\beta g$, and let $w_{g_\beta}$ denote its flattened
harmonic extension. Assume that $\|\eta\|_{L_t^\infty\dot B_{2,1}^2}\leq\epsilon$ where $\epsilon>0$ is sufficiently small. Let
$2\leq p<r\leq\infty$. If $c_0>0$ is sufficiently small and
$C_0>0$ is sufficiently large, then
\begin{equation*}
w_{g_\beta}^{\mathrm{off}}
\coloneq
\left(P_{<c_0\beta}+P_{>C_0\beta}\right)w_{g_\beta},
\end{equation*}
satisfies
\begin{equation}\label{wg-off-diagonal}
\begin{split}
\beta^{\frac{2}{r}}
\left(
\|\nabla_xw_{g_\beta}^{\mathrm{off}}\|_{L_t^2L_z^\infty L_x^r}
+\beta
\|w_{g_\beta}^{\mathrm{off}}\|_{L_t^2L_z^\infty L_x^r}
\right)\lesssim_{p,r}
\|\nabla^2\eta\|_{L_t^4\dot B_{\infty,1}^0}
\||D|^{\frac{2}{p}}g_\beta\|_{L_t^4L_x^p}.
\end{split}
\end{equation}
\end{lemma}

\begin{proof}
We first prove the bound for the low frequency part. By Bernstein it suffices to control just the $\|P_{<c_0\beta} w_{g_{\beta}}\|_{L_t^2L_z^{\infty}L_x^r}$ part of the bound. In proving this, we will need as the main ingredient, the following integrated bound
\begin{equation}\label{tailintegrated}
\beta^{1+\frac{2}{r}}\|P_{<c_0\beta} \nabla_{z,x}w_{g_{\beta}}^{\mathrm{off}}\|_{L_t^2L_z^{1}L_x^r}\lesssim_{p,r}
\|\nabla^2\eta\|_{L_t^4\dot B_{\infty,1}^0}
\||D|^{\frac{2}{p}}g_\beta\|_{L_t^4L_x^p}.
\end{equation}
Indeed, if $c_0$ is small enough, then
$P_{<c_0\beta}w_{g_\beta}$ vanishes at $z=0$. Therefore, from the equation for $w_{g_\beta}$ we have
\begin{equation*}
\begin{split}
\nabla_{z,x}P_{<c_0\beta}w_{g_\beta}
={}&
\nabla_{z,x}\Delta_{z,x}^{-1}\nabla_{z,x}\cdot
\widetilde P_{<c_0\beta}
\left(A_\eta P_{<c_0\beta}\nabla_{z,x}w_{g_\beta}\right)
\\
&+
\nabla_{z,x}\Delta_{z,x}^{-1}\nabla_{z,x}\cdot
\widetilde P_{<c_0\beta}
\left([P_{<c_0\beta},A_\eta]
\nabla_{z,x}w_{g_\beta}\right).
\end{split}
\end{equation*}
Absorbing the first term using smallness of $\nabla\eta$ and using Proposition
\ref{Leibnizcommutator} and Proposition \ref{besovparabounds}, we have
\begin{equation*}
\begin{split}
\beta^{\frac{2}{r}}
\|P_{<c_0\beta}\nabla_{z,x}w_{g_\beta}\|_{L_t^2L_z^1L_x^r}&\lesssim
\beta^{\frac{2}{\widetilde r}}
\|\widetilde P_{<c_0\beta}
([P_{<c_0\beta},A_\eta]\nabla_{z,x}w_{g_\beta})
\|_{L_t^2L_z^1L_x^{\widetilde r}}
\\
&\lesssim
\beta^{-1+\frac{2}{\widetilde r}}
\|\nabla_xA_\eta\|_{L_t^4\dot B_{\infty,1}^0}
\|\nabla_{z,x}w_{g_\beta}\|_{L_t^4L_z^1L_x^{\widetilde r}}
\\
&\lesssim
\beta^{-1+\frac{2}{\widetilde r}}
\|\nabla^2\eta\|_{L_t^4\dot B_{\infty,1}^0}
\||D|^{\frac{2}{p}-\frac{2}{\widetilde r}}
\nabla_{z,x}w_{g_\beta}\|_{L_t^4L_z^1L_x^p},
\end{split}
\end{equation*}
where $\widetilde r$ is an intermediate index such that
$p<\widetilde r<r$. We next note that, for
\begin{equation*}
\sigma=\frac{2}{p}-\frac{2}{\widetilde r}>0,
\end{equation*}
there holds
\begin{equation}\label{Dsigmawg}
\begin{split}
&\||D|^\sigma w_{g_\beta}\|_{L_t^4L_z^\infty L_x^p}
+\||D|^\sigma\nabla_{z,x}w_{g_\beta}\|_{L_t^4L_z^1L_x^p}
\lesssim
\||D|^\sigma g_\beta\|_{L_t^4L_x^p}.
\end{split}
\end{equation}
Indeed, comparing against the flat harmonic extension gives
\begin{equation*}
\begin{split}
|D|^\sigma\nabla_{z,x}w_{g_\beta}
=
e^{z|D|}\left(|D|,\nabla_x\right)|D|^\sigma g_\beta+
\nabla_{z,x}\Delta_{z,x}^{-1}\nabla_{z,x}\cdot
|D|^\sigma\left(A_\eta\nabla_{z,x}w_{g_\beta}\right).
\end{split}
\end{equation*}
The flat term is bounded in $L_t^4L_z^1\dot B_{p,1}^0$ by
$\||D|^\sigma g_\beta\|_{L_t^4L_x^p}$, by the standard Poisson kernel
bounds. For the variable-coefficient term, Proposition \ref{besovparabounds} and the fact that $\sigma>0$ give
\begin{equation*}
\begin{split}
\left\|
\nabla_{z,x}\Delta_{z,x}^{-1}\nabla_{z,x}\cdot
|D|^{\sigma}\left(A_\eta\nabla_{z,x}w_{g_\beta}\right)
\right\|_{L_t^4L_z^1\dot B_{p,1}^0}&\lesssim
\||D|^{\sigma}A_\eta\|_{L_t^\infty
\dot B_{\frac{p\widetilde r}{\widetilde r-p},1}^0}
\|\nabla_{z,x}w_{g_\beta}\|_{L_t^4L_z^1
\dot B_{\widetilde r,1}^0}
\\
&+
\|A_\eta\|_{L_t^\infty\dot B_{\infty,1}^0}
\||D|^{\sigma}\nabla_{z,x}w_{g_\beta}\|_{L_t^4L_z^1
\dot B_{p,1}^0}
\\
&\lesssim
\epsilon
\||D|^{\sigma}\nabla_{z,x}w_{g_\beta}\|_{L_t^4L_z^1
\dot B_{p,1}^0}.
\end{split}
\end{equation*}
Here we used Bernstein's inequality and Proposition \ref{MoserBesov}
\begin{equation*}
\|A_\eta\|_{L_t^\infty\dot B_{\infty,1}^0}
+\||D|^{\sigma}A_\eta\|_{L_t^\infty
\dot B_{\frac{p\widetilde r}{\widetilde r-p},1}^0}
\lesssim
\|\nabla^2\eta\|_{L_t^\infty\dot B_{2,1}^0}
\lesssim\epsilon.
\end{equation*}
as well as the embeddings
\begin{equation*}
\dot B_{p,1}^\sigma\hookrightarrow\dot B_{\widetilde r,1}^0,
\qquad
\dot B_{2,1}^1\hookrightarrow
\dot B_{\frac{p\widetilde r}{\widetilde r-p},1}^\sigma.
\end{equation*}
After absorbing the last term, the $L_z^\infty$ estimate follows
from
\begin{equation*}
|D|^{\sigma}w_{g_\beta}(z)
=
|D|^{\sigma}g_\beta-
\int_z^0|D|^{\sigma}\partial_zw_{g_\beta}(z')\,dz'.
\end{equation*}
This proves \eqref{Dsigmawg}.

Since $g_\beta$ is localized at frequency $\beta$, we have
\begin{equation*}
\begin{split}
\beta^{-1+\frac{2}{\widetilde r}}
\||D|^\sigma g_\beta\|_{L_t^4L_x^p}
\lesssim
\beta^{-1}
\||D|^{\frac{2}{p}}g_\beta\|_{L_t^4L_x^p}.
\end{split}
\end{equation*}
Consequently,
\begin{equation*}
\begin{split}
\beta^{\frac{2}{r}}
\|P_{<c_0\beta}\nabla_{z,x}w_{g_\beta}\|_{L_t^2L_z^1L_x^r}
\lesssim
\beta^{-1}
\|\nabla^2\eta\|_{L_t^4\dot B_{\infty,1}^0}
\||D|^{\frac{2}{p}}g_\beta\|_{L_t^4L_x^p}
\end{split}
\end{equation*}
which is \eqref{tailintegrated}. Moreover, since $P_{<c_0\beta}w_{g_\beta}$ vanishes at $z=0$,
integration in $z$ gives
\begin{equation*}
\begin{split}
\|P_{<c_0\beta}w_{g_\beta}\|_{L_t^2L_z^\infty L_x^r}
\lesssim
\|P_{<c_0\beta}\partial_zw_{g_\beta}\|_{L_t^2L_z^1L_x^r},
\end{split}
\end{equation*}
This proves \eqref{wg-off-diagonal} for the low-frequency part of $w_{g_{\beta}}^{\mathrm{off}}$.

It remains to estimate the high-frequency tail. By Bernstein, it suffices to estimate only the gradient component this time,
\begin{equation*}
W_g\coloneq\nabla_{z,x}w_{g_\beta}.
\end{equation*}
The flattened harmonic-extension equation is
\begin{equation*}
W_g
=
e^{z|D|}\left(|D|g_\beta,\nabla_xg_\beta\right)
-\mathcal R_0(A_\eta W_g),
\qquad
\mathcal R_0
\coloneq
\nabla_{z,x}\Delta_{z,x}^{-1}\nabla_{z,x}\cdot,
\end{equation*}
Again, the flat term vanishes
after applying $Q_\beta\coloneq P_{>C_0\beta}$, provided $C_0$ is
sufficiently large. Let $\widetilde Q_\beta$ be a slightly fattened
high-frequency projection such that
$\widetilde Q_\beta Q_\beta=Q_\beta$. We have
\begin{equation*}
\begin{split}
Q_\beta W_g=-\mathcal R_0\widetilde Q_\beta
\left(A_\eta Q_\beta W_g\right)
-\mathcal R_0\widetilde Q_\beta
\left([Q_\beta,A_\eta]W_g\right).
\end{split}
\end{equation*}
The first term can be absorbed by the same small-coefficient argument as above. For the second term, Proposition
\ref{Leibnizcommutator} and paradifferential calculus gives
\begin{equation*}
\begin{split}
\left\|
\widetilde Q_\beta[Q_\beta,A_\eta]W_g
\right\|_{L_t^2L_z^\infty\dot B_{r,1}^0}\lesssim
\beta^{-1}
\|\nabla_xA_\eta\|_{L_t^4\dot B_{\infty,1}^0}
\|W_g\|_{L_t^4L_z^\infty\dot B_{r,1}^0}.
\end{split}
\end{equation*}
Here the fattened projection is retained so that only output
frequencies $\gtrsim\beta$ occur. Consequently,
\begin{equation*}
\begin{split}
\|Q_\beta W_g\|_{L_t^2L_z^\infty\dot B_{r,1}^0}
\lesssim
\epsilon
\|Q_\beta W_g\|_{L_t^2L_z^\infty\dot B_{r,1}^0}+
\beta^{-1}
\|\nabla_xA_\eta\|_{L_t^4\dot B_{\infty,1}^0}
\|W_g\|_{L_t^4L_z^\infty\dot B_{r,1}^0}.
\end{split}
\end{equation*}
After absorbing the first term, Proposition
\ref{flattenedbasicbounds} and Bernstein's inequality give
\begin{equation*}
\begin{split}
\|W_g\|_{L_t^4L_z^\infty\dot B_{r,1}^0}\lesssim
\|\nabla_xg_\beta\|_{L_t^4\dot B_{r,1}^0}\lesssim
\beta^{1-\frac{2}{r}}
\||D|^{\frac{2}{p}}g_\beta\|_{L_t^4L_x^p}.
\end{split}
\end{equation*}
Consequently,
\begin{equation*}
\begin{split}
\beta^{-1+\frac{2}{r}}
\|P_{>C_0\beta}\nabla_{z,x}w_{g_\beta}\|_{L_t^2L_z^\infty L_x^r}&\lesssim
\beta^{-2+\frac{2}{r}}
\|\nabla_xA_\eta\|_{L_t^4\dot B_{\infty,1}^0}
\|\nabla_{z,x}w_{g_\beta}\|_{L_t^4L_z^\infty\dot B_{r,1}^0}
\\
&\lesssim
\beta^{-1}
\|\nabla^2\eta\|_{L_t^4\dot B_{\infty,1}^0}
\||D|^{\frac{2}{p}}g_\beta\|_{L_t^4L_x^p}.
\end{split}
\end{equation*}
Finally, inverse Bernstein on frequencies $>C_0\beta$ yields
\begin{equation*}
\begin{split}
\beta
\|P_{>C_0\beta}w_{g_\beta}\|_{L_t^2L_z^\infty L_x^r}
\lesssim
\|P_{>C_0\beta}\nabla_xw_{g_\beta}\|_{L_t^2L_z^\infty L_x^r}.
\end{split}
\end{equation*}
This proves the bound in \eqref{wg-off-diagonal} for the high-frequency part as well.
\end{proof}
The main multi-linear commutator estimate we need is the following.
\begin{proposition}[Weighted commutator bound]\label{weighted-commutator-bound}
Let $1<\sigma<2$, and let $\gamma,\beta\in 2^{\mathbb Z}$ satisfy
$\gamma\leq\beta$. Let $4<p<r<\frac{4}{2-\sigma}$ and assume that $\|\eta\|_{L_t^\infty\dot B_{2,1}^2}\leq\epsilon$ for some sufficiently small $\epsilon$. Let $g_\gamma=P_\gamma g$
and $f=f_{<\gamma}$. Then, for any Schwartz functions $h_1$ and
$h_2$, there holds
\begin{equation}\label{weighted-commutator-final}
\begin{split}
&\left\|
h_1h_2\nabla_xP_{<\beta}[L,f]g_\gamma
\right\|_{L_t^1L_x^1}\lesssim_{\sigma,p,r}
\beta
M_{\beta,2}^{\frac\sigma2}
\left(h_1,\nabla_xf\right)
M_{\beta,2}^{\frac\sigma2}
\left(h_2,g_\gamma\right)
\\
&\quad+
\left(\frac{\beta}{\gamma}\right)^{1+\frac2r}
M_{\beta,2}^{\frac\sigma2}
\left(h_1,\nabla_xf\right)
\|h_2\|_{L_t^\infty L_x^2}
\|\nabla_x^2\eta\|_{L_t^4\dot B_{\infty,1}^0}
\||D|^{\frac2p}g_\gamma\|_{L_t^4L_x^p}
\\
&\quad+
M_{\beta,2}^{\frac\sigma2}
\left(h_1,\nabla_x^2\eta\right)
\|h_2\|_{L_t^\infty L_x^2}
\|\nabla_xf\|_{L_t^4\dot B_{\infty,1}^0}
\||D|^{\frac2p}g_\gamma\|_{L_t^4L_x^p}.
\end{split}
\end{equation}
\end{proposition}
\begin{proof}
Let $w_{fg_\gamma}$ denote the flattened harmonic extension of
$fg_\gamma$, and set
\begin{equation*}
u_{\mathrm{com}}
\coloneq
w_{fg_\gamma}-fw_{g_\gamma},
\qquad
W_{\mathrm{com}}
\coloneq
\nabla_{z,x}u_{\mathrm{com}}.
\end{equation*}
Then $u_{\mathrm{com}}$ vanishes at $z=0$ and, by definition of $L$,
\begin{equation*}
[L,f]g_\gamma
=
\left.\partial_zu_{\mathrm{com}}\right|_{z=0}.
\end{equation*}
We have $\Delta_{\eta}u_{\mathrm{com}}=F_{\mathrm{com}}$, where a
simple computation gives
\begin{equation}\label{commutator-source-scalar}
F_{\mathrm{com}}
=
-2D_xw_{g_\gamma}\cdot\nabla_xf
-w_{g_\gamma}\Delta_xf,
\qquad
D_xw\coloneq
\nabla_xw-\nabla_x\eta\,\partial_zw.
\end{equation}
Equivalently,
\begin{equation*}
u_{\mathrm{com}}
=
\left(\Delta_D^{-1}F_{\mathrm{com}}\right)_*,
\qquad
[L,f]g_\gamma
=
\mathcal PF_{\mathrm{com}}.
\end{equation*}
It will be convenient to decompose part of $F_{\mathrm{com}}$ in
divergence form. More precisely, a simple algebraic manipulation gives
\begin{equation}\label{commutator-source-mixed}
F_{\mathrm{com}}
=
\nabla_{z,x}\cdot G_{\mathrm{com}}+H_{\mathrm{com}},
\end{equation}
where we define
\begin{equation*}
G_{\mathrm{com},x}
=
-w_{g_\gamma}\nabla_xf,
\qquad
G_{\mathrm{com},z}
=
2(\nabla_x\eta\cdot\nabla_xf)w_{g_\gamma},
\qquad
H_{\mathrm{com}}
=
-\nabla_xf\cdot\nabla_xw_{g_\gamma}.
\end{equation*}

Set $h\coloneq h_1h_2$. We apply the kernel argument used to prove
\eqref{weighted-exterior-dyadic-bound} directly to the decomposition
\eqref{commutator-source-mixed}. Namely, we first write
\begin{equation*}
W_{\mathrm{com}}
=
\mathcal R_0G_{\mathrm{com}}
+\mathcal S_0H_{\mathrm{com}}
-\mathcal R_0(A_\eta W_{\mathrm{com}}),
\qquad
\mathcal S_0
=
\nabla_{z,x}\Delta_{z,x}^{-1}.
\end{equation*}
At horizontal frequency $\alpha$, the bulk kernels contributing to the expressions
$\nabla_xP_\alpha(\mathcal R_0G_{\mathrm{com}})_3(0)$ and
$\nabla_xP_\alpha(\mathcal S_0H_{\mathrm{com}})_3(0)$ have weighted
$L^1_{z',y}$ sizes $O(\alpha)$ and $O(1)$, respectively. Integrating
the vertical-divergence term by parts in $z'$ also produces the
boundary term $\nabla_xP_\alpha G_{\mathrm{com},z}(0)$, whose weighted
spatial $L_y^1$ kernel has size $O(\alpha)$. Consequently, the same
argument as in the proof of \eqref{weighted-exterior-dyadic-bound}
gives
\begin{equation*}
\begin{split}
\left\|
h\nabla_xP_\alpha[L,f]g_\gamma
\right\|_{L_t^1L_x^1}
&\lesssim_\sigma
\alpha\sup_{z\leq0}
M_{\alpha,1}^{\sigma}
\left(h,\widetilde P_\alpha G_{\mathrm{com}}(z)\right)
+
\sup_{z\leq0}
M_{\alpha,1}^{\sigma}
\left(h,\widetilde P_\alpha H_{\mathrm{com}}(z)\right)
\\
&\quad+
\|h_2\|_{L_t^\infty L_x^2}
\|W_{\mathrm{com}}\|_{L_z^\infty L_t^2
\dot B_{\infty,1}^0}
\sum_{\kappa\gtrsim\alpha}
\frac{\alpha}{\kappa}
M_{\alpha,2}^{\frac\sigma2}
\left(h_1,\widetilde P_\kappa\nabla_xA_\eta\right).
\end{split}
\end{equation*}
Using Proposition \ref{basicMproperties}, the definitions of
$G_{\mathrm{com}}$ and $H_{\mathrm{com}}$, and the smallness of
$\nabla_x\eta$, we obtain
\begin{equation*}
\begin{split}
\left\|
h_1h_2\nabla_xP_\alpha[L,f]g_\gamma
\right\|_{L_t^1L_x^1}
&\lesssim_\sigma
M_{\alpha,2}^{\frac\sigma2}
\left(h_1,\nabla_xf\right)
\sup_{z\leq0}
\left[
\alpha
M_{\alpha,2}^{\frac\sigma2}
\left(h_2,w_{g_\gamma}(z)\right)
+
M_{\alpha,2}^{\frac\sigma2}
\left(h_2,\nabla_xw_{g_\gamma}(z)\right)
\right]
\\
&\quad+
\|h_2\|_{L_t^\infty L_x^2}
\|W_{\mathrm{com}}\|_{L_z^\infty L_t^2
\dot B_{\infty,1}^0}
\sum_{\kappa\gtrsim\alpha}
\frac{\alpha}{\kappa}
M_{\alpha,2}^{\frac\sigma2}
\left(h_1,\widetilde P_\kappa\nabla_xA_\eta\right).
\end{split}
\end{equation*}
For $\alpha<\beta$, Proposition \ref{basicMproperties} gives
\begin{equation*}
M_{\alpha,2}^{\frac\sigma2}(h_i,q)
\lesssim_\sigma
\left(\frac{\alpha}{\beta}\right)^{1-\frac\sigma2}
M_{\beta,2}^{\frac\sigma2}(h_i,q).
\end{equation*}
Summing over $\kappa\gtrsim\alpha$ for the second line and then both
lines over $\alpha<\beta$, using that $\sigma<2$, and using the
definition of $A_\eta$ and the smallness of $\nabla_x\eta$, we obtain
\begin{equation}\label{commutator-first-reduction}
\begin{split}
\left\|
h_1h_2\nabla_xP_{<\beta}[L,f]g_\gamma
\right\|_{L_t^1L_x^1}
&\lesssim_\sigma
M_{\beta,2}^{\frac\sigma2}
\left(h_1,\nabla_xf\right)
\sup_{z\leq0}
\left[
\beta
M_{\beta,2}^{\frac\sigma2}
\left(h_2,w_{g_\gamma}(z)\right)
+
M_{\beta,2}^{\frac\sigma2}
\left(h_2,\nabla_xw_{g_\gamma}(z)\right)
\right]
\\
&\quad+
M_{\beta,2}^{\frac\sigma2}
\left(h_1,\nabla_x^2\eta\right)
\|h_2\|_{L_t^\infty L_x^2}
\|W_{\mathrm{com}}\|_{L_z^\infty L_t^2
\dot B_{\infty,1}^0}.
\end{split}
\end{equation}

We next estimate the first term on the right. Choose
$\widetilde P_\gamma$ compatibly with the cutoffs in Lemma
\ref{wgtailbound}, so that the diagonal and off-diagonal parts satisfy
\begin{equation*}
w_{g_\gamma}
=
\widetilde P_\gamma w_{g_\gamma}
+w_{g_\gamma}^{\mathrm{off}}.
\end{equation*}
We use Proposition \ref{bilinearext} at frequency scale $\gamma$,
together with the scale comparison in Proposition
\ref{basicMproperties}. Indeed, since $\gamma\leq\beta$ and
$\widetilde P_\gamma w_{g_\gamma}$ is supported at frequencies
$\lesssim\gamma$, we have
\begin{equation*}
\begin{split}
\sup_{z\leq0}
M_{\beta,2}^{\frac\sigma2}
\left(h_2,\widetilde P_\gamma w_{g_\gamma}(z)\right)
&\lesssim_\sigma
\left(\frac{\beta}{\gamma}\right)^{1-\frac\sigma2}
\sup_{z\leq0}
M_{\gamma,2}^{\frac\sigma2}
\left(h_2,\widetilde P_\gamma w_{g_\gamma}(z)\right)
\\
&\lesssim_\sigma
\left(\frac{\beta}{\gamma}\right)^{1-\frac\sigma2}
M_{\gamma,2}^{\frac\sigma2}
\left(h_2,g_\gamma\right)
\\
&\lesssim_\sigma
M_{\beta,2}^{\frac\sigma2}
\left(h_2,g_\gamma\right).
\end{split}
\end{equation*}
Here the factors arising from the scale changes cancel. The same
estimate holds with $\widetilde P_\gamma$ replaced by a slightly
fatter annular projection $\widehat P_\gamma$. Moreover, using
\eqref{M-LP-boundedness} and harmlessly fattening the projection, we
obtain
\begin{equation*}
\begin{split}
\sup_{z\leq0}
M_{\beta,2}^{\frac\sigma2}
\left(h_2,\widetilde P_\gamma
\nabla_xw_{g_\gamma}(z)\right)
&\lesssim_\sigma
\gamma
\sup_{z\leq0}
M_{\beta,2}^{\frac\sigma2}
\left(h_2,\widehat P_\gamma w_{g_\gamma}(z)\right)
\\
&\lesssim_\sigma
\gamma
M_{\beta,2}^{\frac\sigma2}
\left(h_2,g_\gamma\right).
\end{split}
\end{equation*}
Since $\gamma\leq\beta$, the preceding estimates imply
\begin{equation}\label{commutator-on-diagonal}
\begin{split}
&\beta\sup_{z\leq0}
M_{\beta,2}^{\frac\sigma2}
\left(h_2,\widetilde P_\gamma w_{g_\gamma}(z)\right)
+
\sup_{z\leq0}
M_{\beta,2}^{\frac\sigma2}
\left(h_2,\widetilde P_\gamma
\nabla_xw_{g_\gamma}(z)\right)
\\
&\qquad\lesssim_\sigma
\beta
M_{\beta,2}^{\frac\sigma2}
\left(h_2,g_\gamma\right).
\end{split}
\end{equation}

To estimate the off-diagonal contribution, we use the following
consequence of Proposition \ref{basicMproperties}: since
$r<\frac{4}{2-\sigma}$, there holds
\begin{equation}\label{M-translation-free-bound}
\sup_{z\leq0}
M_{\beta,2}^{\frac\sigma2}
\left(h_2,q(z)\right)
\lesssim_{\sigma,r}
\beta^{\frac2r}
\|h_2\|_{L_t^\infty L_x^2}
\|q\|_{L_t^2L_z^\infty L_x^r}.
\end{equation}
Applying \eqref{M-translation-free-bound} and Lemma
\ref{wgtailbound} at scale $\gamma$ gives
\begin{equation}\label{commutator-off-diagonal}
\begin{split}
&\sup_{z\leq0}
\left[
\beta
M_{\beta,2}^{\frac\sigma2}
\left(h_2,w_{g_\gamma}^{\mathrm{off}}(z)\right)
+
M_{\beta,2}^{\frac\sigma2}
\left(h_2,\nabla_xw_{g_\gamma}^{\mathrm{off}}(z)\right)
\right]
\\
&\qquad\lesssim_{\sigma,p,r}
\left(\frac{\beta}{\gamma}\right)^{1+\frac2r}
\|h_2\|_{L_t^\infty L_x^2}
\|\nabla_x^2\eta\|_{L_t^4\dot B_{\infty,1}^0}
\||D|^{\frac2p}g_\gamma\|_{L_t^4L_x^p}.
\end{split}
\end{equation}
Indeed, Lemma \ref{wgtailbound} gives
\begin{equation*}
\begin{split}
\gamma^{\frac2r}
\left(
\|\nabla_xw_{g_\gamma}^{\mathrm{off}}\|_
{L_t^2L_z^\infty L_x^r}
+
\gamma
\|w_{g_\gamma}^{\mathrm{off}}\|_
{L_t^2L_z^\infty L_x^r}
\right)
\lesssim_{p,r}
\|\nabla_x^2\eta\|_{L_t^4\dot B_{\infty,1}^0}
\||D|^{\frac2p}g_\gamma\|_{L_t^4L_x^p}.
\end{split}
\end{equation*}
Thus, the term involving $w_{g_\gamma}^{\mathrm{off}}$ costs $\left(\frac{\beta}{\gamma}\right)^{1+\frac2r},$
while the term involving $\nabla_xw_{g_\gamma}^{\mathrm{off}}$ costs
only $\left(\frac{\beta}{\gamma}\right)^{\frac2r}.$ Since $\gamma\leq\beta$, this proves
\eqref{commutator-off-diagonal}. Combining
\eqref{commutator-on-diagonal} and
\eqref{commutator-off-diagonal}, we obtain
\begin{equation}\label{commutator-source-bound}
\begin{split}
&M_{\beta,2}^{\frac\sigma2}
\left(h_1,\nabla_xf\right)
\sup_{z\leq0}
\left[
\beta
M_{\beta,2}^{\frac\sigma2}
\left(h_2,w_{g_\gamma}(z)\right)
+
M_{\beta,2}^{\frac\sigma2}
\left(h_2,\nabla_xw_{g_\gamma}(z)\right)
\right]
\\
&\qquad\lesssim_{\sigma,p,r}
\beta
M_{\beta,2}^{\frac\sigma2}
\left(h_1,\nabla_xf\right)
M_{\beta,2}^{\frac\sigma2}
\left(h_2,g_\gamma\right)
\\
&\qquad+
\left(\frac{\beta}{\gamma}\right)^{1+\frac2r}
M_{\beta,2}^{\frac\sigma2}
\left(h_1,\nabla_xf\right)
\|h_2\|_{L_t^\infty L_x^2}
\|\nabla_x^2\eta\|_{L_t^4\dot B_{\infty,1}^0}
\||D|^{\frac2p}g_\gamma\|_{L_t^4L_x^p}.
\end{split}
\end{equation}

It remains to estimate $W_{\mathrm{com}}$ in the second line of
\eqref{commutator-first-reduction}. From
\eqref{commutator-source-mixed}, the flattened equation gives
\begin{equation*}
W_{\mathrm{com}}
=
\mathcal R_0G_{\mathrm{com}}
+\mathcal S_0H_{\mathrm{com}}
-\mathcal R_0(A_\eta W_{\mathrm{com}}),
\qquad
\mathcal S_0
=
\nabla_{z,x}\Delta_{z,x}^{-1}.
\end{equation*}
The flat half-space kernel bounds, Proposition
\ref{besovparabounds}, and the usual small-data absorption argument
give
\begin{equation}\label{Wcom-flat-bound}
\begin{split}
\|W_{\mathrm{com}}\|_{L_z^\infty L_t^2
\dot B_{\infty,1}^0}
\lesssim{}&
\|G_{\mathrm{com}}\|_{L_z^\infty L_t^2
\dot B_{\infty,1}^0}
+
\|H_{\mathrm{com}}\|_{L_z^1L_t^2
\dot B_{\infty,1}^0}.
\end{split}
\end{equation}
Repeating the argument used to prove \eqref{Dsigmawg}, with the
$L_t^4$ norm taken before the $z$ integration in the second term,
gives
\begin{equation}\label{commutator-natural-extension}
\begin{split}
&\|w_{g_\gamma}\|_{L_t^4L_z^\infty
\dot B_{p,1}^{\frac2p}}
+
\|\nabla_{z,x}w_{g_\gamma}\|_{L_z^1L_t^4
\dot B_{p,1}^{\frac2p}}
\lesssim_p
\||D|^{\frac2p}g_\gamma\|_{L_t^4L_x^p}.
\end{split}
\end{equation}
Using Proposition \ref{besovparabounds}, H\"older's inequality in time,
and
\begin{equation*}
\|\nabla_x\eta\|_{L_t^\infty\dot B_{2,1}^1}
+
\|\nabla_x^2\eta\|_{L_t^\infty\dot B_{2,1}^0}
\lesssim\epsilon,
\end{equation*}
we can estimate $G_{\mathrm{com}}$ and $H_{\mathrm{com}}$ in
\eqref{Wcom-flat-bound}. Consequently,
\begin{equation}\label{Wcom-final-bound}
\begin{split}
\|W_{\mathrm{com}}\|_{L_z^\infty L_t^2
\dot B_{\infty,1}^0}
\lesssim_p
\|\nabla_xf\|_{L_t^4\dot B_{\infty,1}^0}
\||D|^{\frac2p}g_\gamma\|_{L_t^4L_x^p}.
\end{split}
\end{equation}
This estimate holds for any $2<p<\infty$, and hence in particular for
the $p$ under consideration. Substituting
\eqref{commutator-source-bound} and \eqref{Wcom-final-bound} into
\eqref{commutator-first-reduction} proves
\eqref{weighted-commutator-final}.
\end{proof}


\section{Linear and Bilinear Estimates for \texorpdfstring{$V$}{} and \texorpdfstring{$B$}{}}\label{sec:VB}
Here we record some useful linear and bilinear estimates satisfied by the horizontal and vertical velocity traces $V$ and $B$. These will be used very liberally in our analysis below. At the end of this section, we also isolate the principal paradifferential parts of $V$ and $B$ in terms of the good variable $\varphi$, with suitable error bounds. We have the following estimates. 
\begin{proposition}\label{VBunbil}
Assume that $u$ satisfies the bootstrap bounds in Section \ref{sec:bootstrap}, as well as \eqref{Taylorcoeffsmall}. Fix $x_0\in\mathbb{R}^2$. Then $B_{\mu}=P_{\mu}B$ and $V_{\mu}=P_{\mu}V$ satisfy the bounds
\begin{enumerate}
\item ($L_t^{\infty}L_x^2$ bound)
\begin{equation}\label{VBL2}
\|B_{\mu}\|_{L_t^{\infty}L_x^2}+\|V_{\mu}\|_{L_t^{\infty}L_x^2}\lesssim \epsilon c_{\mu}\mu^{\frac{1}{2}-s_c}.
\end{equation}
\item ($L_t^{4}L_x^p$ bound)
\begin{equation}\label{BV48}
\|B_{\mu}\|_{L_t^4L_x^{p}}+\|V_{\mu}\|_{L_t^4L_x^{p}}\lesssim \epsilon c_{\mu}\mu^{\frac{9}{8}+2(\frac{1}{8}-\frac{1}{p})-s_c},\hspace{5mm}8<p\leq \infty.
\end{equation}
\item (Bilinear low-high bound)
\begin{equation}\label{BVulh}
\|B_{\mu}u_{\lambda}^{x_0}\|_{L_{t,x}^2}+\|V_{\mu}u_{\lambda}^{x_0}\|_{L_{t,x}^2}\lesssim c_{\lambda}c_{\mu}\epsilon^2\lambda^{-s_c}\mu^{\frac{5}{4}-s_c}(1+|x_0|\mu)^{\frac12},\hspace{5mm}\mu\ll\lambda.
\end{equation}
\item (Bilinear high-low bound)
\begin{equation}\label{BVuhl}
\|B_{\mu}u_{\lambda}^{x_0}\|_{L_{t,x}^2}+\|V_{\mu}u_{\lambda}^{x_0}\|_{L_{t,x}^2}\lesssim c_{\lambda}c_{\mu}\epsilon^2\lambda^{\frac{3}{4}-s_c}\mu^{\frac{1}{2}-s_c}(1+|x_0|\lambda)^{\frac12},\hspace{5mm}\mu\gg\lambda.
\end{equation}
\item Translation-averaged bilinear $L^2$ bound: For $\frac12<\sigma<1$, there holds
\begin{equation}\label{BVulhtrans}
M_{\mu,2}^{\sigma}(u_\lambda^{x_0},B_\mu)
+M_{\mu,2}^{\sigma}(u_\lambda^{x_0},V_\mu)
\lesssim_\sigma
\epsilon^2c_\lambda c_\mu\lambda^{-s_c}\mu^{\frac54-s_c}(1+\mu|x_0|)^{\frac12},
\hspace{5mm}\mu\ll\lambda.
\end{equation}
\end{enumerate}
Moreover, all of the above estimates hold with $B_\mu$ and $V_\mu$ replaced by $P_\mu B_{<\alpha}$ and $P_\mu V_{<\alpha}$, respectively, for any dyadic number $\alpha\in2^{\mathbb Z}$, with implicit constants uniform in $\alpha$.
\end{proposition}
\begin{remark}
We only prove the above translation-averaged bound in the high-low interaction regime and in terms of the weaker multi-linear norm defined in subsection \ref{subsec:multilinearnorms} as this is all we will need later.     
\end{remark}
\begin{remark}\label{linearizedVB}
We also remark that the natural linearized variants of the bilinear bounds \eqref{BVulh}, \eqref{BVuhl} and \eqref{BVulhtrans} also hold. Here, one replaces each instance of $u_{\lambda}^{x_0}$ in these bounds by the linearized variable $w_{\lambda}^{x_0}$, and also exchanges a factor of $\epsilon c_{\lambda}\lambda^{-s_c}$ with $d_{\lambda}$ on the right-hand side of each bound. The proof of this follows verbatim from the argument below by applying the linearized bootstrap hypotheses in Section \ref{sec:bootstrap} in place of the nonlinear ones. We omit the details for this and only focus on the nonlinear bounds in Proposition \ref{VBunbil} stated above.   
\end{remark}
\begin{proof}
We prove all of the estimates for $B$. At the end of the proof, we explain briefly how essentially the same argument gives all of the corresponding estimates for $V$. In our proof, it will be convenient to bootstrap the following bound for $B_{<\mu}\coloneq B(\eta_{<\mu},\psi_{<\mu})$ for each dyadic number $\mu$,
\begin{equation}\label{B48boot}
\|B_{<\mu}\|_{X^p}\coloneq \||D|^{\frac{7}{8}-2(\frac{1}{8}-\frac{1}{p})}B_{<\mu}\|_{L_t^4L_x^{p}}+\||D|^{\frac{3}{2}}B_{<\mu}\|_{L_t^{\infty}L_x^2}\leq K\epsilon,
\quad \mu\in \mathfrak{c}^{\mathbb{Z}},\quad 8<p\leq\infty.
\end{equation}
We note carefully that $B_{<\mu}$ above is defined by $B(\eta_{<\mu},\psi_{<\mu})$, and not by $P_{<\mu}B$. A priori, this bootstrap is made on some finite small time interval, which will be extended to include all of $[0,\infty)$ by improving the bootstrap constants. In our estimates below, we will show that $B_{<\mu}$ satisfies the same bound with an improved constant.
\subsubsection{Remainder bounds for $B$.} Our first goal is to isolate the principal linear contribution to $B_\mu$. We define the $\mu$-localized remainder term
\begin{equation}\label{Brdecomposition}
B_{r,\mu}\coloneq B_{\mu}-P_{\mu}\int_{0}^{\infty}P_{<M\lambda_0}(L_{<\lambda_0}\varphi_{\lambda_0})\frac{d\lambda_0}{\lambda_0}
\end{equation}
and also the full remainder $B_r\coloneq\sum_{\mu\in2^{\mathbb Z}}B_{r,\mu}$. We first fix an integer $N>s_c$ and then choose $M\gg_N 1$ sufficiently large for the high-frequency tail estimates below. We will show
\begin{equation}\label{Brestimates}
\begin{split}
\|B_{r,\mu}\|_{L_t^4L_x^{\frac{16}{7}}}\lesssim_K \epsilon^2c_{\mu}\mu^{\frac{1}{2}-s_c},
\hspace{2mm}
\|B_{r,\mu}\|_{L_t^2L_x^{\infty}}\lesssim_K \epsilon^2c_{\mu}\mu^{\frac{5}{4}-s_c},\hspace{2mm}\|B_{r,\mu}\|_{L_t^{\infty}L_x^2}&\lesssim_K \epsilon^2 c_{\mu}\mu^{\frac{1}{2}-s_c}.
\end{split}
\end{equation}
The first estimate gives, for $8<p\leq\infty$,
\begin{equation}\label{Br48}
\|B_{r,\mu}\|_{L_t^4L_x^{p}}\lesssim_K \mu^{\frac{5}{8}+2(\frac{1}{8}-\frac{1}{p})}\|B_{r,\mu}\|_{L_t^4L_x^{\frac{16}{7}}}\lesssim_K \epsilon^2 c_{\mu}\mu^{\frac{9}{8}+2(\frac{1}{8}-\frac{1}{p})-s_c}.
\end{equation}
It also gives
\begin{equation}\label{Brlh}
\|B_{r,\mu}u_{\lambda}^{x_0}\|_{L_{t,x}^2}\lesssim_K \|B_{r,\mu}\|_{L_t^4L_x^{\frac{16}{7}}}\|u_{\lambda}\|_{L_t^4L_x^{16}}\lesssim_K \epsilon^3c_{\lambda}c_{\mu}\lambda^{\frac{3}{4}-s_c}\mu^{\frac{1}{2}-s_c},\qquad\mu\gg\lambda.
\end{equation}
On the other hand, the $L_t^2L_x^{\infty}$ bound gives
\begin{equation}\label{Brhl}
\|B_{r,\mu}u_{\lambda}^{x_0}\|_{L_{t,x}^2}\lesssim_K \|B_{r,\mu}\|_{L_t^2L_x^{\infty}}\|u_{\lambda}\|_{L_t^{\infty}L_x^2}\lesssim_K c_{\lambda}c_{\mu}\epsilon^3\lambda^{-s_c}\mu^{\frac{5}{4}-s_c},\qquad\mu\ll\lambda.
\end{equation}
Thus these terms are perturbative as long as $\epsilon$ is sufficiently small relative to $K$. The third estimate in \eqref{Brestimates} will be needed for the $L_t^\infty L_x^2$ estimate for $B_\mu$.

To prove \eqref{Brestimates}, we will explicitly compute the remainder term $B_{r,\mu}$ by writing its multilinear paradifferential expansion up to a certain order. Indeed, we have from the linearization formula \eqref{Blin} and the fundamental theorem of calculus,
\begin{equation}\label{Bintegralexp}
\begin{split}
B
&=\int_{0}^{\infty}\lambda_0\frac{d}{d\lambda_0}P_{<M\lambda_0}B_{<\lambda_0}\frac{d\lambda_0}{\lambda_0}
\\
&=\int_{0}^{\infty}P_{<M\lambda_0}\left(L_{<\lambda_0}\Phi_{<\lambda_0}(\psi_{\lambda_0},\eta_{\lambda_0})+\eta_{\lambda_0}L_{<\lambda_0}B_{<\lambda_0}\right)\frac{d\lambda_0}{\lambda_0}
+\int_{0}^{\infty}\widetilde P_{M\lambda_0}B_{<\lambda_0}\frac{d\lambda_0}{\lambda_0},
\end{split}
\end{equation}
where we define
\begin{equation*}
\widetilde P_{M\lambda_0}\coloneq\lambda_0\frac{d}{d\lambda_0}P_{<M\lambda_0},
\qquad
\Phi_{<\lambda_0}(f,g)\coloneq f-B_{<\lambda_0}g.
\end{equation*}
In the sequel, we will continually adopt this notation for $\Phi_{<\lambda_0}(f,g)$. We will also use the abbreviation
\begin{equation}\label{Phi0definition}
\Phi^0_{<\lambda_0}\coloneq\Phi_{<\lambda_0}(\psi_{\lambda_0},\eta_{\lambda_0})
=\psi_{\lambda_0}-B_{<\lambda_0}\eta_{\lambda_0}.
\end{equation}
\begin{remark} If $\mu$ is a dyadic number and $\lambda_0=\delta\mu\ll\mu $ for some positive number $0<\delta\ll 1$, then for $(f,g)=(\psi_{\mu},\eta_{\mu})$, $\Phi_{<\lambda_0}(f,g)$ approximately corresponds to the good variable $\varphi$ localized at frequency $\mu$, at least, up to some mild commutator terms.
\end{remark}
Define
\begin{equation*}
B_r^1(\lambda_0)\coloneq P_{<M\lambda_0}\left(\eta_{\lambda_0}L_{<\lambda_0}B_{<\lambda_0}\right),
\qquad
B_{r,\mu}^1(\lambda_0)\coloneq P_\mu B_r^1(\lambda_0).
\end{equation*}
$B_r^1(\lambda_0)$ is the second term in the integrand of the first term on the right hand side of \eqref{Bintegralexp}. This term is directly perturbative. Indeed, Bernstein's inequality, Corollary \ref{Bestcor}, and \eqref{B48boot} give
\begin{equation*}
\begin{split}
\|B_{r,\mu}^1(\lambda_0)\|_{L_t^2L_x^{\infty}}\lesssim \|\eta_{\lambda_0}\|_{L_t^4L_x^{\infty}}\|L_{<\lambda_0}B_{<\lambda_0}\|_{L_t^4L_x^{\infty}}
\lesssim_K \epsilon^2 c_{\lambda_0}\lambda_0^{\frac{5}{4}-s_c}.
\end{split}
\end{equation*}
Similarly, we obtain the $L_t^4L_x^{\frac{16}{7}}$ bound
\begin{equation*}
\|B_{r,\mu}^1(\lambda_0)\|_{L_t^4L_x^{\frac{16}{7}}}
\lesssim_K \epsilon^2c_{\lambda_0}\lambda_0^{\frac{1}{2}-s_c}, 
\end{equation*}
as well as the $L_t^{\infty}L_x^2$ bound
\begin{equation*}
\begin{split}
\|B_{r,\mu}^1(\lambda_0)\|_{L_t^{\infty}L_x^2}&\lesssim \|\eta_{\lambda_0}\|_{L_t^{\infty}L_x^2}\|L_{<\lambda_0}B_{<\lambda_0}\|_{L_t^{\infty}L_x^{\infty}}
\\
&\lesssim \epsilon c_{\lambda_0}\lambda_0^{-s_c}\lambda_0^{\frac12}\||D|^{\frac32}B_{<\lambda_0}\|_{L_t^\infty L_x^2}
\\
&\lesssim_K\epsilon^2c_{\lambda_0}\lambda_0^{\frac12-s_c}.
\end{split}
\end{equation*}
Since $P_\mu P_{<M\lambda_0}=0$ when $\lambda_0\ll_M\mu$, Minkowski's inequality and the slowly varying property give
\begin{equation*}
\begin{split}
\int_0^\infty\|B_{r,\mu}^1(\lambda_0)\|_{L_t^4L_x^{\frac{16}{7}}}\frac{d\lambda_0}{\lambda_0}
+\int_0^\infty\|B_{r,\mu}^1(\lambda_0)\|_{L_t^\infty L_x^2}\frac{d\lambda_0}{\lambda_0}
&\lesssim_K\epsilon^2c_\mu\mu^{\frac12-s_c},
\end{split}
\end{equation*}
and also
\begin{equation*}
\begin{split}
\int_0^\infty\|B_{r,\mu}^1(\lambda_0)\|_{L_t^2L_x^\infty}\frac{d\lambda_0}{\lambda_0}
&\lesssim_K\epsilon^2c_\mu\mu^{\frac54-s_c}.
\end{split}
\end{equation*}
We also record the critical integrated bound which we will use later to close the bootstrap \eqref{B48boot}.
\begin{equation}\label{B1rsummed}
\begin{split}
\int_0^\infty\|B_r^1(\lambda_0)\|_{X^p}\frac{d\lambda_0}{\lambda_0}
&\lesssim
\sum_{\mu\in2^{\mathbb Z}}\mu^{\frac32}
\int_0^\infty\left(
\|B_{r,\mu}^1(\lambda_0)\|_{L_t^4L_x^{\frac{16}{7}}}
+\|B_{r,\mu}^1(\lambda_0)\|_{L_t^\infty L_x^2}
\right)\frac{d\lambda_0}{\lambda_0}
\lesssim_K\epsilon^2.
\end{split}
\end{equation}
Here and in the analogous $X^p$ summations below, including when closing \eqref{B48boot}, we apply the preceding estimates at the critical regularity $s=2$ and used the summability of the $\dot{B}_{2,1}^2$ envelope. Next, we turn to the first term in the second line on the right-hand side of \eqref{Bintegralexp}, which we can write as
\begin{equation}\label{Br23decomp}
\begin{split}
\int_0^\infty P_{<M\lambda_0}L_{<\lambda_0}\Phi^0_{<\lambda_0}\frac{d\lambda_0}{\lambda_0}
&=\int_0^\infty P_{<M\lambda_0}L_{<\lambda_0}(\Phi^0_{<\lambda_0}-\varphi_{\lambda_0})\frac{d\lambda_0}{\lambda_0}
\\
&+\int_0^\infty P_{<M\lambda_0}(L_{<\lambda_0}\varphi_{\lambda_0})\frac{d\lambda_0}{\lambda_0}.
\end{split}
\end{equation}
The second term is the principal contribution in \eqref{Brdecomposition}. The first is controlled by the following lemma.

\begin{lemma}\label{Phiest}
Under the hypotheses of Proposition \ref{VBunbil} and the bootstrap hypothesis \eqref{B48boot}, let $A_{<\lambda_0}$ denote either $L_{<\lambda_0}$ or $M_{<\lambda_0}$. Then, for every $8<p\leq\infty$,
\begin{equation}\label{Phi4inf}
\|A_{<\lambda_0}\Phi^0_{<\lambda_0}\|_{L_t^4\dot B^0_{p,1}}
\lesssim_p\epsilon c_{\lambda_0}\lambda_0^{\frac{11}{8}-\frac2p-s_c},
\qquad
\|A_{<\lambda_0}\Phi^0_{<\lambda_0}\|_{L_t^\infty\dot B^0_{2,1}}
\lesssim\epsilon c_{\lambda_0}\lambda_0^{\frac12-s_c}.
\end{equation}
Moreover, for $4<p\leq\infty$, the error term $E^A_{\lambda_0}\coloneq A_{<\lambda_0}(\Phi^0_{<\lambda_0}-\varphi_{\lambda_0})$, satisfies
\begin{equation}\label{B2errorbounds}
\begin{split}
\|E^A_{\lambda_0}\|_{L_t^4\dot B^0_{\frac{16}{7},1}}\lesssim_K\epsilon^2c_{\lambda_0}\lambda_0^{\frac12-s_c},
\hspace{2mm}
\|E^A_{\lambda_0}\|_{L_t^2\dot B^0_{p,1}}
\lesssim_{K,p}\epsilon^2c_{\lambda_0}\lambda_0^{\frac54-\frac{2}{p}-s_c},\hspace{2mm}\|E^A_{\lambda_0}\|_{L_t^\infty\dot B^0_{2,1}}\lesssim_K\epsilon^2c_{\lambda_0}\lambda_0^{\frac12-s_c}.
\end{split}
\end{equation}
\end{lemma}

\begin{proof}
It will be convenient to make the bootstrap hypothesis that \eqref{Phi4inf} holds, simultaneously for $A=L$ and $A=M$, with right-hand side multiplied by some large constant $K$. By \eqref{Phi0definition}, we have
\begin{equation}\label{Phithreeterms}
\Phi^0_{<\lambda_0}
=\varphi_{\lambda_0}+(T_{B^l}-B_{<\lambda_0})\eta_{\lambda_0}+[S_{\lambda_0},T_{B^l}]\eta,
\end{equation}
where $S_{\lambda_0}=\lambda_0\frac{d}{d\lambda_0}P_{<\lambda_0}$. The contribution of the first term can be handled by applying \ref{basicopbounds} and Proposition \ref{phiu-bootstrap}, to obtain
\begin{equation}\label{Phiprincipalbound}
\|A_{<\lambda_0}\varphi_{\lambda_0}\|_{L_t^4\dot B^0_{p,1}}
\lesssim_p\epsilon c_{\lambda_0}\lambda_0^{\frac{11}{8}-\frac2p-s_c},
\hspace{5mm}
\|A_{<\lambda_0}\varphi_{\lambda_0}\|_{L_t^\infty\dot B^0_{2,1}}
\lesssim\epsilon c_{\lambda_0}\lambda_0^{\frac12-s_c}.
\end{equation}
To deal with the remaining two terms (which define the contribution of the error term $E_{\lambda_0}^A$), we first record a few useful estimates. First, we will need for $8<p\leq\infty$, the bounds
\begin{equation}\label{BauxiliaryPhi}
\begin{split}
\||D|B_{<\lambda_0}\|_{L_t^4\dot B^0_{p,1}}\lesssim_{K,p}\epsilon\lambda_0^{\frac38-\frac2p},\hspace{5mm}
\||D|B_{<\lambda_0}\|_{L_t^\infty\dot B^0_{\infty,1}}\lesssim_K\epsilon\lambda_0^{\frac12}.
\end{split}
\end{equation}
For finite $p>8$, the first estimate follows by choosing $8<q<p$, applying Corollary \ref{Bestcor}, and using the $X^q$ bootstrap; the endpoint $p=\infty$ follows by choosing any finite $q>8$ and applying Bernstein. For the second estimate, we split at frequency $\lambda_0$, estimate the low-frequency part by Bernstein and \eqref{B48boot},
\begin{equation*}
\||D|P_{<C\lambda_0}B_{<\lambda_0}\|_{\dot B^0_{\infty,1}}
\lesssim\lambda_0^{\frac12}\||D|^{\frac32}B_{<\lambda_0}\|_{L_x^2},
\end{equation*}
and control the high-frequency tail using Proposition \ref{Bernsteinelliptic}.

Next, we record two auxiliary estimates. If $\delta\lambda_0\lesssim\nu\lesssim\lambda_0$, then the fundamental theorem of calculus and \eqref{Bintegralexp} give
\begin{equation}\label{Bberndecomp}
B_{<\nu}-B_{<\lambda_0}
=\int_{\lambda_0}^{\nu}\left(L_{<\alpha}\Phi^0_{<\alpha}+\eta_\alpha L_{<\alpha}B_{<\alpha}\right)\frac{d\alpha}{\alpha}.
\end{equation}
Moreover, \eqref{BauxiliaryPhi} and the paraproduct estimates give
\begin{equation*}
\begin{split}
\|\eta_\alpha L_{<\alpha}B_{<\alpha}\|_{L_t^4\dot B^0_{p,1}}\lesssim_{K,p}\epsilon^2c_\alpha\alpha^{\frac{11}{8}-\frac2p-s_c},\hspace{3mm}\|\eta_\alpha L_{<\alpha}B_{<\alpha}\|_{L_t^\infty\dot B^0_{2,1}}\lesssim_K\epsilon^2c_\alpha\alpha^{\frac12-s_c}.
\end{split}
\end{equation*}
Consequently, the bootstrap bounds imply the Bernstein type estimates
\begin{equation}\label{BBern}
\begin{split}
\|B_{<\nu}-B_{<\lambda_0}\|_{L_t^4\dot B^0_{p,1}}\lesssim_{K,p}\epsilon c_{\lambda_0}\lambda_0^{\frac{11}{8}-\frac2p-s_c},\hspace{3mm}\|B_{<\nu}-B_{<\lambda_0}\|_{L_t^\infty\dot B^0_{2,1}}\lesssim_K\epsilon c_{\lambda_0}\lambda_0^{\frac12-s_c}.
\end{split}
\end{equation}
We now estimate the second term in \eqref{Phithreeterms}. With $M\delta\ll1$, the paraproduct convention gives
\begin{equation}\label{phiparaerror}
\begin{split}
(T_{B^l}-B_{<\lambda_0})\eta_{\lambda_0}=-\int_{\alpha\approx\lambda_0}P_{\geq M\delta\alpha}B_{<\lambda_0}\,\widetilde\eta_\alpha\frac{d\alpha}{\alpha}+\int_{\alpha\approx\lambda_0}P_{<M\delta\alpha}\left(B_{<\delta\alpha}-B_{<\lambda_0}\right)\widetilde\eta_\alpha\frac{d\alpha}{\alpha},
\end{split}
\end{equation}
where $B^l_{<\delta\alpha}=P_{<M\delta\alpha}B_{<\delta\alpha}$ and $\widetilde\eta_\alpha$ is a harmless fattened localization of $\eta$ at frequency $\alpha$. In the first term on the right-hand side of \eqref{phiparaerror}, we estimate $B$ by \eqref{BauxiliaryPhi} and retain the frequency envelope factor from $\widetilde\eta_\alpha$. In the second term, \eqref{BBern} instead supplies the frequency envelope. Using Proposition \ref{basicopbounds} and Proposition \ref{besovparabounds}, and the slowly varying property, we obtain, for $A=L,M$ (with $8<p\leq\infty$ in the first line and $4<p\leq\infty$ in the third line),
\begin{equation}\label{PhiParaErrorGeneral}
\begin{split}
\|A_{<\lambda_0}(T_{B^l}-B_{<\lambda_0})\eta_{\lambda_0}\|_{L_t^4\dot B^0_{p,1}}
&\lesssim_{K,p}\epsilon^2c_{\lambda_0}\lambda_0^{\frac{11}{8}-\frac2p-s_c},
\\
\|A_{<\lambda_0}(T_{B^l}-B_{<\lambda_0})\eta_{\lambda_0}\|_{L_t^4\dot B^0_{\frac{16}{7},1}}
&\lesssim_K\epsilon^2c_{\lambda_0}\lambda_0^{\frac12-s_c},
\\
\|A_{<\lambda_0}(T_{B^l}-B_{<\lambda_0})\eta_{\lambda_0}\|_{L_t^2\dot B^0_{p,1}}
&\lesssim_{K,p}\epsilon^2c_{\lambda_0}\lambda_0^{\frac54-\frac2p-s_c},
\\
\|A_{<\lambda_0}(T_{B^l}-B_{<\lambda_0})\eta_{\lambda_0}\|_{L_t^\infty\dot B^0_{2,1}}
&\lesssim_K\epsilon^2c_{\lambda_0}\lambda_0^{\frac12-s_c}.
\end{split}
\end{equation}
For example, to estimate $A_{<\lambda_0}$ applied to the second term on the right hand side in \eqref{phiparaerror}, we use the critical bounds
\begin{equation*}
\lambda_0^{\frac{15}{8}}\|\eta_{\lambda_0}\|_{L_t^\infty L_x^{\frac{16}{7}}}
+\lambda_0^{\frac98+\frac1p}\|\eta_{\lambda_0}\|_{L_t^4L_x^{2p}}
+\lambda_0\|\eta_{\lambda_0}\|_{L_t^\infty L_x^\infty}
\lesssim\epsilon
\end{equation*}
along with \eqref{BBern}. For the middle estimate, when $4<p<\infty$, we use Hölder with both $L_t^4$ factors in $L_x^{2p}$ and apply \eqref{BauxiliaryPhi} and \eqref{BBern} with exponent $2p>8$; the case $p=\infty$ is unchanged. The contribution of the first term in \eqref{phiparaerror} is similar. Finally, Proposition \ref{Leibnizcommutator} gives the schematic representation
\begin{equation}\label{PhiCommRep}
[S_{\lambda_0},T_{B^l}]\eta
=\lambda_0^{-1}\int_{\alpha\approx\lambda_0}L\left(\nabla_xB^l_{<\delta\alpha},\widetilde\eta_\alpha\right)\frac{d\alpha}{\alpha},
\end{equation}
where $L$ is a translation-invariant bilinear form with rapidly decaying kernel. Thus the same estimates give
\begin{equation}\label{PhiCommGeneral}
\begin{split}
\|A_{<\lambda_0}[S_{\lambda_0},T_{B^l}]\eta\|_{L_t^4\dot B^0_{p,1}}
&\lesssim_{K,p}\epsilon^2c_{\lambda_0}\lambda_0^{\frac{11}{8}-\frac2p-s_c},
\\
\|A_{<\lambda_0}[S_{\lambda_0},T_{B^l}]\eta\|_{L_t^4\dot B^0_{\frac{16}{7},1}}
&\lesssim_K\epsilon^2c_{\lambda_0}\lambda_0^{\frac12-s_c},
\\
\|A_{<\lambda_0}[S_{\lambda_0},T_{B^l}]\eta\|_{L_t^2\dot B^0_{p,1}}
&\lesssim_{K,p}\epsilon^2c_{\lambda_0}\lambda_0^{\frac54-\frac2p-s_c},
\\
\|A_{<\lambda_0}[S_{\lambda_0},T_{B^l}]\eta\|_{L_t^\infty\dot B^0_{2,1}}
&\lesssim_K\epsilon^2c_{\lambda_0}\lambda_0^{\frac12-s_c}.
\end{split}
\end{equation}
Combining \eqref{Phiprincipalbound}, \eqref{PhiParaErrorGeneral}, and \eqref{PhiCommGeneral}, and then taking $\epsilon$ sufficiently small relative to $K$, improves the bootstrap and proves the lemma.
\end{proof}
Now, we return to estimating the errors in \eqref{Br23decomp}. Beginning with the first term on the right-hand side, we define
\begin{equation*}
B_r^2(\lambda_0)\coloneq P_{<M\lambda_0}E^L_{\lambda_0},
\qquad
B_{r,\mu}^2(\lambda_0)\coloneq P_\mu B_r^2(\lambda_0)
\end{equation*}
where $E^L_{\lambda_0}$ is defined as in Lemma \ref{Phiest}. Using Lemma \ref{Phiest}, Minkowski's inequality, and the same frequency-envelope summation and integration as above give
\begin{equation*}
\begin{split}
\int_0^\infty\left(
\|B_{r,\mu}^2(\lambda_0)\|_{L_t^4L_x^{\frac{16}{7}}}
+\|B_{r,\mu}^2(\lambda_0)\|_{L_t^\infty L_x^2}
\right)\frac{d\lambda_0}{\lambda_0}\lesssim_K\epsilon^2c_\mu\mu^{\frac12-s_c}
\end{split}
\end{equation*}
and
\begin{equation*}
\int_0^\infty\|B_{r,\mu}^2(\lambda_0)\|_{L_t^2L_x^\infty}\frac{d\lambda_0}{\lambda_0}\lesssim_K\epsilon^2c_\mu\mu^{\frac54-s_c}    
\end{equation*}
as well as the critical estimate
\begin{equation}\label{B2rsummed}
\int_0^\infty\|B_r^2(\lambda_0)\|_{X^p}\frac{d\lambda_0}{\lambda_0}\lesssim_K\epsilon^2.
\end{equation}
It remains to estimate
\begin{equation*}
B_r^3(\lambda_0)\coloneq\widetilde P_{M\lambda_0}B_{<\lambda_0},
\qquad
B_{r,\mu}^3(\lambda_0)\coloneq P_\mu B_r^3(\lambda_0).
\end{equation*}
It will be convenient to temporarily define
\begin{equation*}
H_{\lambda_0}\coloneq\|\nabla_x^2\eta_{<\lambda_0}\|_{L_t^\infty\dot B^0_{\infty,1}}.
\end{equation*}
The bootstrap bounds and slow variation give
$H_{\lambda_0}\lesssim\epsilon c_{\lambda_0}\lambda_0$.
Therefore,
\begin{equation}\label{Htamepower}
H_{\lambda_0}^{N-1}
\lesssim\epsilon^{N-1}c_{\lambda_0}\lambda_0^{N-1}.
\end{equation}
Proposition \ref{Bernsteinelliptic}, applied with the high-frequency projection $P_{>cM\lambda_0}$ for a fixed $0<c<1$, now yields
\begin{equation*}
\begin{split}
\|B_r^3(\lambda_0)\|_{L_t^2L_x^\infty}\lesssim_N\lambda_0^{-1-N}
\|\nabla_x^2\eta_{<\lambda_0}\|_{L_t^4\dot B^0_{\infty,1}}
H_{\lambda_0}^{N-1}
\||D|B_{<\lambda_0}\|_{L_t^4\dot B^0_{\infty,1}}\lesssim_K\epsilon^2c_{\lambda_0}\lambda_0^{\frac54-s_c}.
\end{split}
\end{equation*}
Here we used
\begin{equation*}
\|\nabla_x^2\eta_{<\lambda_0}\|_{L_t^4\dot B^0_{\infty,1}}\lesssim\epsilon\lambda_0^{\frac78},
\qquad
\||D|B_{<\lambda_0}\|_{L_t^4\dot B^0_{\infty,1}}\lesssim_K\epsilon\lambda_0^{\frac38}.
\end{equation*}
For the $L_t^4L_x^{16/7}$ estimate, we choose $q>8$ and $r<\infty$ with $\frac7{16}=\frac1q+\frac1r$ in Proposition \ref{Bernsteinelliptic}, use the $X^q$ bound for $B_{<\lambda_0}$, and apply the corresponding analogue of \eqref{Htamepower}. For the $L_t^\infty L_x^2$ estimate, we take $p=q=2$ and $r=\infty$. This gives
\begin{equation*}
\begin{split}
\int_0^\infty\left(
\|B_{r,\mu}^3(\lambda_0)\|_{L_t^4L_x^{\frac{16}{7}}}
+\|B_{r,\mu}^3(\lambda_0)\|_{L_t^\infty L_x^2}
\right)\frac{d\lambda_0}{\lambda_0}\lesssim_K\epsilon^2c_\mu\mu^{\frac12-s_c}.
\end{split}
\end{equation*}
and 
\begin{equation*}
\int_0^\infty\|B_{r,\mu}^3(\lambda_0)\|_{L_t^2L_x^\infty}\frac{d\lambda_0}{\lambda_0}\lesssim_K\epsilon^2c_\mu\mu^{\frac54-s_c},    
\end{equation*}
Moreover, we have
\begin{equation}\label{B3rsummed}
\int_0^\infty\|B_r^3(\lambda_0)\|_{X^p}\frac{d\lambda_0}{\lambda_0}\lesssim_K\epsilon^2,
\qquad 8<p\leq\infty.
\end{equation}
Collecting the bounds for $B_r^1$, $B_r^2$, and $B_r^3$ above
proves \eqref{Brestimates}. Although these estimates and
\eqref{Brdecomposition} are sufficient for the linear bounds for
$B$ and $V$, the term $L_{<\lambda_0}\varphi_{\lambda_0}$ may still
have output frequency much smaller than $\lambda_0$. We therefore
extract its principal paradifferential part.

\begin{lemma}\label{princextraction}
Define
\begin{equation*}
B_r^4(\lambda_0)
\coloneq
P_{<M\lambda_0}
\big(L_{<\lambda_0}\varphi_{\lambda_0}\big)
-
Op(\mathfrak l^l_{<\delta\lambda_0})
\varphi_{\lambda_0}.
\end{equation*}
Then for $4<p\leq\infty$ there holds
\begin{equation}\label{Br4Lp}
\int_0^\infty
\|P_\mu B_r^4(\lambda_0)\|_{L_t^2L^p}
\frac{d\lambda_0}{\lambda_0}
\lesssim
\epsilon^2c_\mu\mu^{\frac54-s_c}.
\end{equation}
\end{lemma}
\begin{proof}
Let $\mathfrak l_{<\lambda_0}$ denote the untruncated principal
symbol of $L_{<\lambda_0}$. We write
\begin{equation}\label{Br4threeerrors}
B_r^4(\lambda_0)
=E_{\lambda_0}^1+E_{\lambda_0}^2+E_{\lambda_0}^3,
\end{equation}
where
\begin{equation*}
\begin{split}
E_{\lambda_0}^1\coloneq
P_{<M\lambda_0}
\big(L_{<\lambda_0}-Op(\mathfrak l_{<\lambda_0})\big)
\varphi_{\lambda_0},\hspace{3mm}E_{\lambda_0}^2\coloneq
P_{<M\lambda_0}
Op\big(\mathfrak l_{<\lambda_0}
-\mathfrak l_{<\delta\lambda_0}\big)
\varphi_{\lambda_0},
\end{split}
\end{equation*}
and
\begin{equation*}
E_{\lambda_0}^3\coloneq
\big(
P_{<M\lambda_0}Op(\mathfrak l_{<\delta\lambda_0})
-Op(\mathfrak l^l_{<\delta\lambda_0})
\big)\varphi_{\lambda_0}.    
\end{equation*}
Fix $4<p<\infty$. Since $2p>8$, the bootstrap in Section \ref{sec:bootstrap} and Proposition \ref{phiu-bootstrap} give
\begin{equation}\label{Br4inputbounds}
\begin{split}
&\|\nabla_x^2\eta_{<\lambda_0}\|_
{L_t^4\dot B^0_{2p,1}}\lesssim
\epsilon\lambda_0^{\frac78-\frac1p},\hspace{3mm}\|\nabla_x(\eta_{<\lambda_0}
-\eta_{<\delta\lambda_0})\|_
{L_t^4\dot B^0_{2p,1}}\lesssim_\delta
\epsilon\lambda_0^{-\frac18-\frac1p},\hspace{3mm}\|\varphi_{\lambda_0}\|_{L_t^4L_x^{2p}}\lesssim
\epsilon c_{\lambda_0}
\lambda_0^{\frac38-\frac1p-s_c}.
\end{split}
\end{equation}
The subprincipal estimate in Proposition \ref{basicopbounds} then gives
\begin{equation*}
\|E_{\lambda_0}^1\|_{L_t^2L_x^p}
\lesssim
\|\nabla_x^2\eta_{<\lambda_0}\|_
{L_t^4\dot B^0_{2p,1}}
\|\varphi_{\lambda_0}\|_{L_t^4L_x^{2p}}
\lesssim
\epsilon^2c_{\lambda_0}
\lambda_0^{\frac54-\frac2p-s_c}.
\end{equation*}

For the second term, the symbol $\mathfrak l(q,\xi)$ is smooth as a function of
$\nabla_x\eta$ and homogeneous of degree one in $\xi$. Thus the
fundamental theorem of calculus shows that
$\mathfrak l_{<\lambda_0}-\mathfrak l_{<\delta\lambda_0}$ contains
one factor of
$\nabla_x(\eta_{<\lambda_0}-\eta_{<\delta\lambda_0})$, while the remaining $\nabla\eta$-dependent factors are harmless and can be placed in $L^{\infty}_{t,x}$. Consequently,
\begin{equation*}
\begin{split}
\|E_{\lambda_0}^2\|_{L_t^2L_x^p}\lesssim
\|\nabla_x(\eta_{<\lambda_0}
-\eta_{<\delta\lambda_0})\|_
{L_t^4\dot B^0_{2p,1}}
\|\nabla_x\varphi_{\lambda_0}\|_{L_t^4L_x^{2p}}\lesssim
\epsilon^2c_{\lambda_0}
\lambda_0^{\frac54-\frac2p-s_c}.
\end{split}
\end{equation*}

Finally, the output and spatial-symbol cutoffs defining
$\mathfrak l^l_{<\delta\lambda_0}$ are comparable, so the low
spatial-frequency part of the symbol cancels in $E_{\lambda_0}^3$.
The remaining symbol is supported at spatial frequencies
$\gtrsim\lambda_0$. Since
\begin{equation*}
|\nabla_x\mathfrak l_{<\delta\lambda_0}(x,\xi)|
\lesssim
|\xi|\,|\nabla_x^2\eta_{<\delta\lambda_0}(x)|,
\qquad |\xi|\approx\lambda_0,
\end{equation*}
Bernstein gives
\begin{equation*}
\|E_{\lambda_0}^3\|_{L_t^2L_x^p}
\lesssim
\|\nabla_x^2\eta_{<\delta\lambda_0}\|_
{L_t^4\dot B^0_{2p,1}}
\|\varphi_{\lambda_0}\|_{L_t^4L_x^{2p}}
\lesssim
\epsilon^2c_{\lambda_0}
\lambda_0^{\frac54-\frac2p-s_c}.
\end{equation*}
All three terms in \eqref{Br4threeerrors} have output
frequency $\lesssim_M\lambda_0$, so
$P_\mu B_r^4(\lambda_0)=0$ unless
$\lambda_0\gtrsim_M\mu$. Therefore, similarly to the estimates for $B_r^1$, $B_r^2$ and $B_r^3$, we obtain
\begin{equation*}
\begin{split}
\int_0^\infty
\|P_\mu B_r^4(\lambda_0)\|_{L_t^2\dot B^0_{p,1}}
\frac{d\lambda_0}{\lambda_0}\lesssim
\epsilon^2
\int_{\lambda_0\gtrsim_M\mu}
c_{\lambda_0}
\lambda_0^{\frac54-\frac2p-s_c}
\frac{d\lambda_0}{\lambda_0}\lesssim
\epsilon^2c_\mu\mu^{\frac54-\frac2p-s_c},
\end{split}
\end{equation*}
which proves \eqref{Br4Lp}.
\end{proof}
Now, we are ready to prove the linear and bilinear estimates in Proposition \ref{VBunbil} and close the bootstrap \eqref{B48boot}.
\subsubsection{Linear estimates for $B$}
Proposition \ref{basicopbounds} and Proposition \ref{phiu-bootstrap} imply
\begin{equation*}
\|L_{<\lambda_0}\varphi_{\lambda_0}\|_{L_t^4\dot B^0_{p,1}}
\lesssim_p\epsilon c_{\lambda_0}\lambda_0^{\frac98+2(\frac18-\frac1p)-s_c},
\qquad
\|L_{<\lambda_0}\varphi_{\lambda_0}\|_{L_t^\infty\dot B^0_{2,1}}
\lesssim\epsilon c_{\lambda_0}\lambda_0^{\frac12-s_c}.
\end{equation*}
Using \eqref{Brdecomposition} and \eqref{Brestimates}, we conclude
\begin{equation}\label{closedB}
\|B_\mu\|_{L_t^4L_x^p}\lesssim_p\epsilon c_\mu\mu^{\frac98+2(\frac18-\frac1p)-s_c},
\qquad
\|B_\mu\|_{L_t^\infty L_x^2}\lesssim\epsilon c_\mu\mu^{\frac12-s_c}.
\end{equation}
\subsubsection{The $X_p$ bootstrap for $B$} To close \eqref{B48boot}, we write
\begin{equation*}
\begin{split}
B_{<\mu}
&=\int_0^\mu\left[
P_{<M\lambda_0}\left(L_{<\lambda_0}\Phi^0_{<\lambda_0}+\eta_{\lambda_0}L_{<\lambda_0}B_{<\lambda_0}\right)
+\widetilde P_{M\lambda_0}B_{<\lambda_0}
\right]\frac{d\lambda_0}{\lambda_0}
+P_{\geq M\mu}B_{<\mu}.
\end{split}
\end{equation*}
The preceding estimates give
\begin{equation*}
\|B_{<\mu}\|_{X^p}\lesssim\epsilon+\|P_{\geq M\mu}B_{<\mu}\|_{X^p}.
\end{equation*}
Proposition \ref{Bernsteinelliptic} gives
\begin{equation*}
\|P_{\geq M\mu}B_{<\mu}\|_{X^p}\lesssim\epsilon\|B_{<\mu}\|_{X^p};
\end{equation*}
with an implicit constant independent of $K$ (as long as $\epsilon$ is small enough).
\subsubsection{Bilinear $(B,u)$ bounds}
If $\nu\ll\lambda$ or $\nu\gg\lambda$, the kernel argument used in Corollary \ref{Abilbounds}, together with Proposition \ref{phiu-bootstrap}, gives
\begin{equation}\label{Lphiu-bound}
\|L_{<\nu}\varphi_\nu\,u_\lambda^{x_0}\|_{L_{t,x}^2}
\lesssim\epsilon^2c_\nu c_\lambda
\nu^{\frac12-s_c}\lambda^{-s_c}\min(\nu,\lambda)^{\frac34}
\left(1+\min(\nu,\lambda)|x_0|\right)^{\frac12}.
\end{equation}
To obtain the full $(B,u)$ estimate, suppose first that $\lambda\ll\mu$. Since the principal integral in \eqref{Brdecomposition} is supported where $\nu\gtrsim\mu$, \eqref{Lphiu-bound} and \eqref{Brlh} yield
\begin{equation*}
\|B_\mu u_\lambda^{x_0}\|_{L_{t,x}^2}
\lesssim\epsilon^2c_\lambda c_\mu\lambda^{\frac34-s_c}\mu^{\frac12-s_c}(1+\lambda|x_0|)^{\frac12}.
\end{equation*}

Now suppose that $\mu\ll\lambda$. We split the principal integral into $\nu\approx\mu$ and $\nu\gg\mu$. The first range follows from \eqref{Lphiu-bound}. For $\nu\gg\mu$, the frequency support of $Op(\mathfrak l^l_{<\delta\nu})\varphi_\nu$ is contained in the region $|\xi|\approx\nu$. Hence
\begin{equation*}
P_\mu Op(\mathfrak l^l_{<\delta\nu})\varphi_\nu=0.
\end{equation*}
It follows directly from Lemma \ref{princextraction} that
\begin{equation}\label{Bprincipaloffdiag}
\int_{\nu\gg\mu}
\|P_\mu P_{<M\nu}L_{<\nu}\varphi_\nu\|_{L_t^2L_x^\infty}
\frac{d\nu}{\nu}
\lesssim\epsilon^2c_\mu\mu^{\frac54-s_c}.
\end{equation}
The bound \eqref{Lphiu-bound} remains valid with $L_{<\nu}$ replaced by $M_{<\nu}$ by the same kernel argument. Moreover, the proof of Lemma \ref{princextraction} applies without change with $L_{<\nu}$ and $\mathfrak l$ replaced by $M_{<\nu}$ and $\mathfrak m$, respectively. Thus \eqref{Bprincipaloffdiag} also remains valid with $L_{<\nu}$ replaced by $M_{<\nu}$. Pairing \eqref{Bprincipaloffdiag} with $u_\lambda$ in $L_t^\infty L_x^2$, and using \eqref{Brhl} for the remainder, proves
\begin{equation*}
\|B_\mu u_\lambda^{x_0}\|_{L_{t,x}^2}
\lesssim\epsilon^2c_\lambda c_\mu\lambda^{-s_c}\mu^{\frac54-s_c}(1+\mu|x_0|)^{\frac12}.
\end{equation*}

It remains to justify that the same bounds hold for the truncated $B_{<\alpha}$ uniformly in $\alpha$. The only substantive difference is that $B_{<\alpha}$ is linearized by integrating only up to the parameter $\alpha$,
\begin{equation}\label{Btruncatedexpansion}
\begin{split}
B_{<\alpha}
&=\int_0^\alpha\left[
P_{<M\nu}\left(L_{<\nu}\Phi^0_{<\nu}+\eta_\nu L_{<\nu}B_{<\nu}\right)
+\widetilde P_{M\nu}B_{<\nu}
\right]\frac{d\nu}{\nu}
+P_{\geq M\alpha}B_{<\alpha}.
\end{split}
\end{equation}
The integral terms are estimated exactly as in the proof of the corresponding bounds for $B$, since restricting the integration to $\nu\leq\alpha$ only shortens the frequency summations. The only additional term is the additional high-frequency tail on the right. Since
\begin{equation*}
P_{\geq M\alpha}B_{<\alpha}
=P_{\geq M\alpha}L_{<\alpha}\psi_{<\alpha},
\end{equation*}
Proposition \ref{Bernsteinelliptic}, with $N>s_c$ for the bilinear estimates, controls this term in the same norms used above. The preceding linear and bilinear arguments are therefore uniform in $\alpha$. Applying $P_\mu$ to \eqref{Btruncatedexpansion} proves all four estimates for $P_\mu B_{<\alpha}$.

\subsubsection{Translation-averaged $(B,u)$ bound} Finally, we prove \eqref{BVulhtrans}. By \eqref{Brdecomposition} and Lemma \ref{princextraction}, write
\begin{equation*}
B_\mu=P_\mu T_{\mathfrak l}\varphi+E_\mu,
\qquad
E_\mu\coloneq B_{r,\mu}
+\int_0^\infty P_\mu B_r^4(\nu)\frac{d\nu}{\nu}.
\end{equation*}
The preceding remainder calculations, keeping the two $L_t^4$ factors in $L_x^{2p}$, give
\begin{equation*}
\|E_\mu\|_{L_t^2L_x^p}
\lesssim_p\epsilon^2c_\mu\mu^{\frac54-\frac2p-s_c},
\qquad 4<p<\infty.
\end{equation*}
Fix $\frac12<\sigma<1$ and choose $4<p<2/(1-\sigma)$. Proposition \ref{basicMproperties} then gives
\begin{equation*}
\begin{split}
M_{\mu,2}^{\sigma}(u_\lambda^{x_0},E_\mu)
&\lesssim_\sigma
\mu^{\frac2p}
\|u_\lambda\|_{L_t^\infty L_x^2}
\|E_\mu\|_{L_t^2L_x^p}
\lesssim_\sigma
\epsilon^3c_\lambda c_\mu\lambda^{-s_c}\mu^{\frac54-s_c}.
\end{split}
\end{equation*}
For the remaining $P_\mu T_{\mathfrak l}\varphi$ component, only input frequencies for $\varphi$ comparable to $\mu$ contribute. The translation-averaged estimate for $\varphi$ in Proposition \ref{phiu-bootstrap}, Lemma \ref{annularcontrol}, and Minkowski's inequality applied to the rapidly decreasing principal-symbol kernels for $T_{\mathfrak{l}}$ give the desired bound, with the factor $\mu$ supplied by the fact that $T_{\mathfrak{l}}$ has order one. This proves \eqref{BVulhtrans} for $B$, and thus concludes the proof of the all the estimates involving $B$ in Proposition \ref{VBunbil}. Now, we remark on the minor changes needed to obtain the estimates for $V$. These simply follow by repeating the preceding proof using the linearization identity for $V$ instead of $B$,
\begin{equation*}
\nu\frac{d}{d\nu}V_{<\nu}
=M_{<\nu}\Phi^0_{<\nu}+\eta_\nu M_{<\nu}B_{<\nu}.
\end{equation*}
One then decomposes the resulting expression $V_r^1$, $V_r^2$, $V_r^3$, and $V_r^4$ corresponding, respectively, to $\eta_\nu M_{<\nu}B_{<\nu}$, $M_{<\nu}(\Phi^0_{<\nu}-\varphi_\nu)$, the high-frequency tail $\widetilde P_{M\nu}V_{<\nu}$, and the principal-symbol error
\begin{equation*}
P_{<M\nu}(M_{<\nu}\varphi_\nu)
-Op(\mathfrak m^l_{<\delta\nu})\varphi_\nu.
\end{equation*}
The preceding estimates apply without change with $L$ and $\mathfrak l$ replaced by $M$ and $\mathfrak m$. Indeed, Proposition \ref{basicopbounds} and Proposition \ref{Bernsteinelliptic} apply for both operators, Lemma \ref{Phiest} includes $A=M$, and the proof of Lemma \ref{princextraction} applies identically to the principal symbol $\mathfrak m$. The term containing $\eta_\nu M_{<\nu}B_{<\nu}$ is controlled using the estimates for $B$ already proved above. Repeating the preceding linear, bilinear, and truncated-variable arguments proves all the claimed estimates for $V$ as well.
\end{proof}
We finally collect some useful facts that are essentially corollaries of the the estimates used in the proof of above. First, from the estimates for $B_r^j$ for $j=1,2,3,4$ and the analogous estimates for the $V$ remainder, we have the principal paradifferential error bound
\begin{corollary}\label{principalVBexp}
Let $\mu>0$. Under the assumptions of Proposition \ref{VBunbil}, there holds
\begin{equation*}
\|P_\mu(B-T_{\mathfrak l}\varphi)\|_{L_t^2L_x^p}
+\|P_\mu(V-T_{\mathfrak m}\varphi)\|_{L_t^2L_x^p}
\lesssim_{\delta,p}
\epsilon^2c_\mu\mu^{\frac54-\frac2p-s_c},\hspace{5mm}4<p\leq\infty
\end{equation*}
The same estimate holds for $(B,V)$ and $(T_{\mathfrak{l}}\varphi,T_{\mathfrak{m}}\varphi)$ replaced by $(B_{<\alpha},V_{<\alpha})$ and $(T_{\mathfrak{l}}\varphi_{<\alpha},T_{\mathfrak{m}}\varphi_{<\alpha})$, respectively.
\end{corollary}
We will also need the following $L_{t,x}^2$ estimate for the error between (a truncation of) $\Phi^0_{<\delta\lambda_0}$ and the good variable $\varphi_{\lambda_0}$. The main observation is that the coefficient $B_{<\delta\lambda_0}$ in $\Phi^0_{<\delta\lambda_0}$ is localized to frequency much smaller than $\lambda_0$ (modulo tails), so that we have access to the unbalanced bilinear estimate. 
\begin{corollary}\label{Phil0varphibound}
Let $\lambda_0>0$ and define $\Phi^{l,0}_{<\delta\lambda_0}
\coloneq\psi_{\lambda_0}-B^l_{<\delta\lambda_0}\eta_{\lambda_0}$, which is obtained from $\Phi_{<\delta\lambda_0}$ by truncating $B_{<\delta\lambda_0}$ to output frequency less than $M\delta\lambda_0$ with $\delta\lambda_0\ll M\delta\lambda_0\ll\lambda_0$. Under the assumptions of Proposition \ref{VBunbil}, there holds
\begin{equation*}
\|\Phi^{l,0}_{<\delta\lambda_0}-\varphi_{\lambda_0}\|_{L_{t,x}^2}
\lesssim\epsilon^2c_{\lambda_0}\lambda_0^{-\frac34-s_c}.
\end{equation*}
\end{corollary}
\begin{proof}
Since $\varphi=\psi-T_{B^l}\eta$, we have
\begin{equation*}
\Phi^{l,0}_{<\delta\lambda_0}-\varphi_{\lambda_0}
=(T_{B^l}-B^l_{<\delta\lambda_0})\eta_{\lambda_0}
+[S_{\lambda_0},T_{B^l}]\eta.
\end{equation*}
Using Proposition \ref{VBunbil}, Proposition \ref{basicopbounds} and Proposition \ref{phiu-bi-unbal}, we have for $\mu\ll\lambda_0$,
\begin{equation}\label{Bunbalancedeta}
\|S_\mu B_{<\rho}\,\widetilde\eta_{\lambda_0}\|_{L_{t,x}^2}
+\|L_{<\mu}\varphi_\mu\,\widetilde\eta_{\lambda_0}\|_{L_{t,x}^2}
\lesssim\epsilon^2c_\mu c_{\lambda_0}
\lambda_0^{-s_c}\mu^{-\frac34},
\end{equation}
uniformly in the truncation parameter $\rho$. Above, we write $S_{\mu}=\mu\frac{d}{d\mu}P_{<\mu}$. For the commutator term, \eqref{PhiCommRep} and \eqref{Bunbalancedeta} give
\begin{equation*}
\begin{split}
\|[S_{\lambda_0},T_{B^l}]\eta\|_{L_{t,x}^2}
&\lesssim\lambda_0^{-1}
\int_{\alpha\approx\lambda_0}
\int_{\mu\lesssim\delta M \alpha}
\mu\,
\|S_{\mu} B_{<\delta\alpha}\,
\widetilde\eta_\alpha\|_{L_{t,x}^2}
\frac{d\mu}{\mu}\frac{d\alpha}{\alpha}
\\
&\lesssim\epsilon^2c_{\lambda_0}\lambda_0^{-1-s_c}
\int_{\mu\lesssim\delta M\lambda_0}
c_\mu\mu^{\frac14}\frac{d\mu}{\mu}
\\
&\lesssim\epsilon^2c_{\lambda_0}\lambda_0^{-\frac34-s_c}
\end{split}
\end{equation*}
By our paraproduct conventions, the other term is an integral over $\alpha\approx\lambda_0$ (with respect to the measure $\frac{d\alpha}{\alpha}$) involving $\left(B^l_{<\delta\alpha}-B^l_{<\delta\lambda_0}\right)
\widetilde\eta_\alpha.$ We decompose
\begin{equation*}
\begin{split}
B^l_{<\delta\alpha}-B^l_{<\delta\lambda_0}=P_{<M\delta\alpha}
\left(B_{<\delta\alpha}-B_{<\delta\lambda_0}\right)+
\left(P_{<M\delta\alpha}-P_{<M\delta\lambda_0}\right)
B_{<\delta\lambda_0}.
\end{split}
\end{equation*}
The second term is localized at frequencies comparable to $\delta\lambda_0$, so \eqref{Bunbalancedeta} gives the desired bound directly. For the first term, we use the linearization formula
\begin{equation*}
B_{<\delta\alpha}-B_{<\delta\lambda_0}
=\int_{\delta\lambda_0}^{\delta\alpha}
\left(
L_{<\beta}\Phi^0_{<\beta}
+\eta_\beta L_{<\beta}B_{<\beta}
\right)\frac{d\beta}{\beta}.
\end{equation*}
Here $\beta\approx\delta\lambda_0$. We decompose $L_{<\beta}\Phi^0_{<\beta}
=L_{<\beta}\varphi_\beta+E^L_\beta,$ where $E_{\beta}^L$ is as in Lemma \ref{Phiest}. The first term can be estimated using \eqref{Bunbalancedeta}. On the other hand, Lemma \ref{Phiest} and the estimate for $B_r^1$ in the proof of Proposition \ref{VBunbil} give
\begin{equation*}
\|E^L_\beta\|_{L_t^2L_x^\infty}
+\|P_{<M\beta}(\eta_\beta L_{<\beta}B_{<\beta})\|_{L_t^2L_x^\infty}
\lesssim\epsilon^2c_\beta\beta^{-\frac34}.
\end{equation*}
Pairing these terms with $\widetilde\eta_\alpha$ in $L_t^\infty L_x^2$ makes them perturbative. Since $\beta\approx\delta\lambda_0$, integration in $\alpha$ and $\beta$ therefore gives
\begin{equation*}
\|(T_{B^l}-B^l_{<\delta\lambda_0})\eta_{\lambda_0}\|_{L_{t,x}^2}
\lesssim\epsilon^2c_{\lambda_0}\lambda_0^{-\frac34-s_c}.
\end{equation*}
Combining the two estimates proves the corollary.
\end{proof}
We conclude this section by recording a bilinear estimate for the subprincipal errors $B-T_{\mathfrak{l}}\varphi$ and $V-T_{\mathfrak{m}}\varphi$. Later in our multi-linear expansions of the nonlinearity in the water wave equations, we will need a slightly more precise formula for these error terms, but we postpone this until it is needed.
\begin{proposition}\label{VBsecondord}
Under the assumptions of Proposition \ref{VBunbil}, define
\begin{equation*}
B^{\mathrm{sub}}\coloneq B-T_{\mathfrak l}\varphi,
\qquad
V^{\mathrm{sub}}\coloneq V-T_{\mathfrak m}\varphi,
\end{equation*}
and define their geometrically truncated analogues as in
Corollary \ref{principalVBexp}. Let $\lambda\ll\mu$, let
$x_0\in\mathbb R^2$, and let $L$ be a translation invariant
bilinear form with rapidly decaying kernel at the
corresponding input scales. For $F=B,V$, there holds, uniformly
in $\alpha>0$,
\begin{equation}\label{L1BVsubhl}
\begin{split}
\|\overline u_\mu^{x_0}P_\mu
L(u_\lambda,F^{\mathrm{sub}})\|_{L^1_{t,x}}
+\|\overline u_\mu^{x_0}P_\mu
L(u_\lambda,F^{\mathrm{sub}}_{<\alpha})\|_{L^1_{t,x}}\lesssim
\epsilon^4c_\mu^2\mu^{\frac14-2s_c}\lambda^{-\frac54}
(1+\mu|x_0|)^{\frac12}.
\end{split}
\end{equation}
The same estimate holds with $u_\lambda$ replaced by
$\lambda^{\frac12}\varphi_\lambda$ or
$\lambda^{-\frac12}P_\lambda(B,V)$.
\end{proposition}
The hypothesis $\lambda\ll\mu$ is to exclude the possibility that $L(u_{\mu},F^{\mathrm{sub}})$ contains a cubic term with three factors at the output frequency $\approx \mu$. Such terms are not expected to satisfy the same estimate.
\begin{proof}
We prove the estimate for $B^{\mathrm{sub}}$ with low-frequency
factor $u_\lambda$. The other cases follow by the same argument
with minor modifications. By \eqref{Bintegralexp} and
Lemma \ref{princextraction}, we have
\begin{equation}\label{B4terms}
B^{\mathrm{sub}}
=\sum_{j=1}^4\int_0^\infty B_r^j(\lambda_0)
\frac{d\lambda_0}{\lambda_0}.
\end{equation}
with each $B_r^j$ defined as above. The output projection and $\lambda\ll\mu$ restrict $F^{\mathrm{sub}}$, the second
factor in $L$ to frequencies comparable to $\mu$. Thus only
$\lambda_0\gtrsim_M\mu$ contributes in the integrals above. We retain these localizations  throughout the proof. For $j=1$, we recall
\begin{equation*}
B_r^1(\lambda_0)
=P_{<M\lambda_0}
(\eta_{\lambda_0}L_{<\lambda_0}B_{<\lambda_0}).
\end{equation*}
We begin by peeling off its principal quadratic part. Recalling the linearization formula for $B$,
\begin{equation*}
\rho\frac{d}{d\rho}(L_{<\rho}B_{<\rho})
=L_{<\rho}^2\Phi^0_{<\rho}
+\eta_\rho L_{<\rho}^2B_{<\rho}
\end{equation*}
and defining the error term
\begin{equation*}
\mathcal E_\rho\coloneq
L_{<\rho}^2(\Phi^0_{<\rho}-\varphi_\rho)
+(L_{<\rho}^2-T_{\mathfrak l}^2)\varphi_\rho
+\eta_\rho L_{<\rho}^2B_{<\rho},
\end{equation*}
we have
\begin{equation*}
L_{<\lambda_0}B_{<\lambda_0}
=\int_0^{\lambda_0}
\left(T_{\mathfrak l}^2\varphi_\rho+\mathcal E_\rho\right)
\frac{d\rho}{\rho},
\qquad
\|\mathcal E_\rho\|_{L_t^2L_x^\infty}
\lesssim\epsilon^2\rho^{\frac14},
\end{equation*}
The $L_t^2L_x^{\infty}$ error bound follows simply from Proposition \ref{VBunbil} and by essentially repeating the proofs of Lemmas
\ref{Phiest} and \ref{princextraction}. With this formula, we now estimate the contribution of $B_r^1$ in \eqref{L1BVsubhl}. We first estimate the principal contribution
\begin{equation*}
\overline u_\mu^{x_0}P_\mu
L\left(u_\lambda,
P_{<M\lambda_0}
(\eta_{\lambda_0}T_{\mathfrak l}^2\varphi_\rho)\right).
\end{equation*}
When $\rho\ll\lambda_0$, its output localization
implies $\mu\approx\lambda_0$, and we use the
unbalanced estimates for
$(u_\lambda,\eta_{\lambda_0})$ and
$(\overline u_\mu^{x_0},\varphi_\rho)$. On the other hand, when $\rho\approx\lambda_0$, we subdivide the frequency
annuli sufficiently narrowly. For output $\mu\sim_c\lambda_0$, Lemma
\ref{balancedrestrictedoutput} applies, we pair
$(\eta_{\lambda_0},\varphi_\rho)$ and
$(\overline u_\mu^{x_0},u_\lambda)$. Otherwise, we can use the unbalanced pairs
$(\overline u_\mu^{x_0},\eta_{\lambda_0})$ and
$(u_\lambda,\varphi_\rho)$. Combining this with the kernel bounds for $T_{\mathfrak{l}}^2$ then gives
\begin{equation*}
\begin{split}
\|\overline u_\mu^{x_0}P_\mu
L\left(u_\lambda,
P_{<M\lambda_0}
(\eta_{\lambda_0}T_{\mathfrak l}^2\varphi_\rho)\right)
\|_{L^1_{t,x}}\lesssim
\epsilon^4c_\mu c_{\lambda_0}\mu^{-s_c}
\lambda_0^{\frac14-s_c}\lambda^{-\frac54}
\left(\frac{\rho}{\lambda_0}\right)^{\frac14}
(1+\mu|x_0|)^{\frac12}.
\end{split}
\end{equation*}
For the contribution containing $\mathcal E_\rho$,
we use one unbalanced estimate with $(u_\lambda,\eta_{\lambda_0})$ in $L^2_{t,x}$ and then place
$\overline u_\mu^{x_0}$ in $L_t^\infty L_x^2$, and
$\mathcal E_\rho$ in $L_t^2L_x^\infty$. This gives
\begin{equation*}
\begin{split}
\|\overline u_\mu^{x_0}P_\mu
L\left(u_\lambda,
P_{<M\lambda_0}(\eta_{\lambda_0}\mathcal E_\rho)\right)
\|_{L^1_{t,x}}\lesssim
\epsilon^5c_\mu c_{\lambda_0}\mu^{-s_c}
\lambda_0^{\frac14-s_c}\lambda^{-\frac54}
\left(\frac{\rho}{\lambda_0}\right)^{\frac14}.
\end{split}
\end{equation*}
Integrating both estimates in $\rho$ and $\lambda_0$ shows $B_r^1$ contributes the required bound in \eqref{L1BVsubhl}.
\\

For $B_r^2$, formulas \eqref{phiparaerror},
\eqref{Bberndecomp}, and \eqref{PhiCommRep} give
\begin{equation}\label{Br2exp}
\begin{split}
B_r^2(\lambda_0)
={}&-P_{<M\lambda_0}L_{<\lambda_0}
\int_{\alpha\approx\lambda_0}
\widetilde\eta_\alpha
P_{\geq M\delta\alpha}B_{<\lambda_0}
\frac{d\alpha}{\alpha}
\\
&+P_{<M\lambda_0}L_{<\lambda_0}
\int_{\alpha\approx\lambda_0}
\int_{\lambda_0}^{\delta\alpha}
\widetilde\eta_\alpha P_{<M\delta\alpha}
\left(L_{<\rho}\Phi^0_{<\rho}
+\eta_\rho L_{<\rho}B_{<\rho}\right)
\frac{d\rho}{\rho}\frac{d\alpha}{\alpha}
\\
&+\lambda_0^{-1}P_{<M\lambda_0}L_{<\lambda_0}
\int_{\alpha\approx\lambda_0}
\mathcal L(\nabla_xB^l_{<\delta\alpha},
\widetilde\eta_\alpha)\frac{d\alpha}{\alpha}.
\end{split}
\end{equation}
Here $\mathcal L$ denotes the translation invariant
forms supplied by the kernel formulas.

In the first integral, decompose the $B$ factor into
frequencies $\beta\gtrsim M\delta\alpha$ and write
\begin{equation*}
P_\beta B_{<\lambda_0}
=P_\beta T_{\mathfrak l}\varphi_{<\lambda_0}
+P_\beta B^{\mathrm{sub}}_{<\lambda_0},
\qquad
\|P_\beta B^{\mathrm{sub}}_{<\lambda_0}
\|_{L_t^2L_x^\infty}
\lesssim\epsilon^2\beta^{-\frac34}
\end{equation*}
where the latter estimate follows from Corollary \ref{principalVBexp}. The contribution of the subprincipal term in \eqref{L1BVsubhl} is easy. We simply use one bilinear estimate for $(u_{\mu}^{x_0},u_{\lambda})$, estimate $\tilde{\eta}_{\alpha}$ in $L_t^{\infty}L_x^2$ and $P_{\beta}B_{<\lambda_0}^{\mathrm{sub}}$ in $L_t^2L_x^{\infty}$ using the bound above, and then integrate. For the remaining principal term term, we have
$M\delta\alpha\lesssim\beta\lesssim_M\lambda_0$. To estimates its contribution, we further decompose the product inside $L_{<\lambda_0}$
into frequencies $\gamma\gtrsim\lambda_0$ and $\gamma\ll\lambda_0$. For the latter case, using Proposition \ref{basicopbounds} and (the proof of) Lemma \ref{balancedrestrictedoutput}, we can estimate
\begin{equation*}
\|L_{<\lambda_0}P_{\ll\lambda_0}(\tilde{\eta}_{\alpha}P_{\beta}T_{\mathfrak{l}}\varphi_{<\lambda_0})\|_{L_{t,x}^2}\lesssim \lambda_0^{\frac{1}{2}}\sup_{\gamma\ll\lambda_0}\|P_{\gamma}|D|^{\frac{1}{2}}(\tilde{\eta}_{\alpha}P_{\beta}T_{\mathfrak{l}}\varphi_{<\lambda_0})\|_{L_{t,x}^2}\lesssim \epsilon^2\lambda^{\frac{1}{4}-s_c}_0c_{\lambda_0}    
\end{equation*}
Then using one bilinear estimate with $(u^{x_0}_{\mu},u_{\lambda})$, and integrating in $\lambda_0$, we see that the contribution of the above term is consistent with \eqref{L1BVsubhl}. On the other hand, for when the product inside $L_{<\lambda_0}$ is localized to frequency $\gamma\gtrsim\lambda_0$, we first replace $L_{<\lambda_0}$ by $T_{\mathfrak l}$. The kernel argument in the proof of Lemma \ref{princextraction} gives
\begin{equation*}
\|P_{<M\lambda_0}(L_{<\lambda_0}-T_{\mathfrak l})
P_{\gtrsim\lambda_0}
(\widetilde\eta_\alpha P_\beta T_{\mathfrak l}\varphi_{<\lambda_0})
\|_{L_{t,x}^2}
\lesssim_{M,\delta}
\epsilon^3c_{\lambda_0}\lambda_0^{\frac14-s_c}.
\end{equation*}
Here we put $\widetilde\eta_\alpha$ in $L_t^\infty L_x^2$ and the differentiated $\varphi$ factor and the differentiated coefficient from the operator difference in $L_t^4\dot{B}_{\infty,1}^0$. Pairing this error with $(\overline u_\mu^{x_0},u_\lambda)$ gives a contribution consistent with \eqref{L1BVsubhl}. For the remaining expression
\begin{equation*}
\overline u_\mu^{x_0}P_\mu L\left(u_\lambda,
P_{<M\lambda_0}T_{\mathfrak l}P_{\gtrsim\lambda_0}
(\widetilde\eta_\alpha P_\beta T_{\mathfrak l}\varphi_{<\lambda_0})\right),
\end{equation*}
we use the same frequency decomposition and bilinear pairings as for $B_r^1$, retaining the input and output localizations, and then integrate.
\\

Now, we turn to the second integral in \eqref{Br2exp}. First, we note the error bound
\begin{equation*}
\mathcal E_\rho^B\coloneq
L_{<\rho}\Phi^0_{<\rho}
+\eta_\rho L_{<\rho}B_{<\rho}
-T_{\mathfrak l}\varphi_\rho,
\qquad
\|\mathcal E_\rho^B\|_{L_t^2L_x^\infty}
\lesssim\epsilon^2\rho^{-\frac34}
\end{equation*}
which follows from Lemma \ref{Phiest} and the
proofs of the preceding remainder estimates. As before, the contribution of this error term in \eqref{L1BVsubhl} is perturbative. The remaining principal part is
$\widetilde\eta_\alpha
P_{<M\delta\alpha}T_{\mathfrak l}\varphi_\rho$.
The projection implies 
$\rho\lesssim M\delta\lambda_0\ll\lambda_0$,
so we pair
$(u_\lambda,\widetilde\eta_\alpha)$ and
$(\overline u_\mu^{x_0},\varphi_\rho)$ and then integrate to estimate this term as well. 
\\

For the last integral in \eqref{Br2exp}, we simply decompose
$B^l_{<\delta\alpha}$ into frequencies
$\rho\lesssim M\delta\alpha$ and pair 
$(u_\lambda,\widetilde\eta_\alpha)$ and
$(\overline u_\mu^{x_0},P_\rho B^l_{<\delta\alpha})$ using Proposition \ref{VBunbil}. Integrating and combining the above estimates, we find that the full contribution of $B_r^2$ in \eqref{L1BVsubhl} is perturbative.
\\

Next, we deal with
$B_r^3(\lambda_0)=\widetilde P_{M\lambda_0}B_{<\lambda_0}$. This is just a high-frequency tail error. By arguing as in the estimate for $B_r^3$ in Proposition \ref{VBunbil}, we have
\begin{equation*}
\|B_r^3(\lambda_0)\|_{L_t^1L_x^2}
\lesssim\epsilon^4c_{\lambda_0}\lambda_0^{-s_c},
\end{equation*}
and hence
\begin{equation*}
\begin{split}
\|\overline u_\mu^{x_0}P_\mu
L(u_\lambda,B_r^3(\lambda_0))\|_{L^1_{t,x}}\lesssim
\|u_\mu\|_{L_t^\infty L_x^2}
\|u_\lambda\|_{L^\infty_{t,x}}
\|B_r^3(\lambda_0)\|_{L_t^1L_x^2}\lesssim
\epsilon^6c_\mu c_{\lambda_0}
\mu^{-s_c}\lambda^{-1}\lambda_0^{-s_c}.
\end{split}
\end{equation*}
which suffices upon integration in $\lambda_0$.
\\
For $B_r^4$, we use the decomposition
\begin{equation}\label{Br4exp}
\begin{split}
B_r^4(\lambda_0)
=P_{<M\lambda_0}
(L_{<\lambda_0}-L_{<\delta\lambda_0})
\varphi_{\lambda_0}+\left(P_{<M\lambda_0}L_{<\delta\lambda_0}
-Op(\mathfrak l^l_{<\delta\lambda_0})\right)
\varphi_{\lambda_0}.
\end{split}
\end{equation}
For the first term, we apply the linearization formula for $L$ \eqref{Blin} and the fundamental theorem of calculus and Proposition \ref{basicopbounds}. After peeling off subprincipal errors as in $B_r^1$ and $B_r^2$, for $\rho\ll\lambda_0$, the contributions containing
$\nabla_x^2\eta_\rho$ and $\varphi_{\lambda_0}$
are estimated by pairing
$(\overline u_\mu^{x_0},\eta_\rho)$ and
$(u_\lambda,\varphi_{\lambda_0})$:
\begin{equation*}
\begin{split}
\rho^2
\|\overline u_\mu^{x_0}P_\mu
L\left(u_\lambda,P_{<M\lambda_0}
\mathcal L(\eta_\rho,\varphi_{\lambda_0})\right)
\|_{L^1_{t,x}}\lesssim
\epsilon^4c_\mu c_{\lambda_0}\mu^{-s_c}
\lambda_0^{\frac14-s_c}\lambda^{-\frac54}
\left(\frac{\rho}{\lambda_0}\right)^{\frac34}
(1+\mu|x_0|)^{\frac12}.
\end{split}
\end{equation*}
For $\rho\approx_\delta\lambda_0$, we again can use the frequency
division already described for $B_r^1$ (after peeling of subprincipal errors), now with
internal pair $(\eta_\rho,\varphi_{\lambda_0})$. The second term in \eqref{Br4exp} can be estimated in a similar way by using Proposition \ref{Bernsteinelliptic} and Proposition \ref{paraexp}. This handles the contribution of $B_r^4$. Combining the four contributions gives
\begin{equation*}
\begin{split}
&\sum_{j=1}^4\int_{\lambda_0\gtrsim_M\mu}
\|\overline u_\mu^{x_0}P_\mu
L(u_\lambda,B_r^j(\lambda_0))\|_{L^1_{t,x}}
\frac{d\lambda_0}{\lambda_0}\lesssim
\epsilon^4c_\mu^2\mu^{\frac14-2s_c}
\lambda^{-\frac54}(1+\mu|x_0|)^{\frac12}.
\end{split}
\end{equation*}
This concludes the proof.
\end{proof}


\section{Paradifferential expansion of \texorpdfstring{$N_{\eta}$}{}}\label{sec:Neta} 
Here, we use the methodology described in Subsection \ref{multpara} and the detailed estimates developed in Section \ref{multpara} and Section \ref{sec:VB} to obtain a sufficiently precise multilinear paradifferential expansion for the nonlinearity $N_{\eta}=G(\eta)\psi$,
as heuristically computed earlier in Theorem~\ref{t:expansions}. Our starting point is to use the linearization formula \eqref{Glin} for the Dirichlet--Neumann map,
\begin{equation}\label{Giter1}
\begin{split}
G(\eta)\psi&=\int_{0}^{\infty}\frac{d}{d\lambda_0}(G_{<\lambda_0}\psi_{<\lambda_0})d\lambda_0
\\
&=\int_{0}^{\infty}\left(G_{<\lambda_0}\psi_{\lambda_0}-G_{<\lambda_0}(\eta_{\lambda_0}B_{<\lambda_0})-\nabla\cdot (\eta_{\lambda_0}V_{<\lambda_0})\right)\frac{d\lambda_0}{\lambda_0}.
\end{split}
\end{equation}
To compactify notation somewhat, we define the operator
\begin{equation*}
\mathcal{T}_{1,<\lambda}(f,g)\coloneq G_{<\lambda}f-G_{<\lambda}(gB_{<\lambda})-\nabla\cdot (gV_{<\lambda}),
\end{equation*}
where $f,g$ are Schwartz functions and the coefficients $B_{<\lambda}$ and $V_{<\lambda}$ are determined by the background $(\eta_{<\lambda},\psi_{<\lambda})$. Note that if $f=\psi_{\lambda}$ and $g=\eta_{\lambda}$, then $\mathcal{T}_{1,<\lambda}(f,g)$ is just the integrand in \eqref{Giter1}. We can write $\mathcal{T}_{1,<\lambda}(f,g)$ more compactly as
\begin{equation*}
\mathcal{T}_{1,<\lambda}(f,g)=G_{<\lambda}\Phi_{<\lambda}(f,g)-\nabla\cdot (gV_{<\lambda}),
\end{equation*}
so that
\begin{equation}\label{Grep1}
G(\eta)\psi=\int_{0}^{\infty}\mathcal{T}_{1,<\lambda_0}(\psi_{\lambda_0},\eta_{\lambda_0})\frac{d\lambda_0}{\lambda_0}.
\end{equation}


\subsection{Extraction of \texorpdfstring{$N_{\eta}^{[1]}$}{}}\label{subsec:Neta1}
Our first main goal will be to extract the principal paradifferential part of \eqref{Grep1}, which we denote by $N_{\eta}^{[1]}$ or equivalently $(G(\eta)\psi)^{[1]}$. The remainder will consist of a perturbative error and also a term which will, upon reiteration, contribute further perturbative errors and terms of the form $N_{\eta}^{[k]}$ (or equivalently $(G(\eta)\psi)^{[k]}$) for $k\in\{2,3,4\}$, as described in \eqref{N-expansion+}. Below, we will frequently adopt the notation
\begin{equation*}
B^l_{<\alpha}\coloneq P_{<M\alpha}B_{<\alpha},\qquad
V^l_{<\alpha}\coloneq P_{<M\alpha}V_{<\alpha},\qquad
\Phi^l_{<\alpha}(f,g)\coloneq f-B^l_{<\alpha}g.
\end{equation*}
The first two definitions have already been used in the previous sections. On the other hand, $\Phi^l_{<\alpha}(f,g)$ is obtained from $\Phi_{<\alpha}(f,g)$ by truncating the $B_{<\alpha}$ coefficient to frequency $\lesssim M\alpha$. The difference between $\Phi_{<\alpha}(f,g)$ and $\Phi^l_{<\alpha}(f,g)$ contains a high-frequency tail of $B_{<\alpha}$ and will be perturbative in the applications below by Proposition \ref{Bernsteinelliptic}, together with the uniform estimates for the truncated variables in Proposition \ref{VBunbil}. For each fixed outer scale $\lambda_0$, we also use the shorthand notation
\begin{equation*}
\Phi^{l,0}_{<\delta\lambda_0}\coloneq
\Phi^l_{<\delta\lambda_0}(\psi_{\lambda_0},\eta_{\lambda_0})
=\psi_{\lambda_0}-B^l_{<\delta\lambda_0}\eta_{\lambda_0}.
\end{equation*}

We now have the following proposition for extracting $(G(\eta)\psi)^{[1]}$ from $G(\eta)\psi$.
\begin{proposition}\label{firstdecomposition}
Let $M$ and $\delta$ be fixed as in Subsection \ref{paraconventions}, with $M\gg1$ and $M\delta\ll1$, and suppose that the assumptions of Proposition \ref{VBunbil} hold. Define $(G(\eta)\psi)^{[1]}$ by the formula
\begin{equation}\label{G1def}
(G(\eta)\psi)^{[1]}\coloneq T_G\varphi-\nabla\cdot (T_{V^l}\eta).
\end{equation}
Then, for every $\mu>0$ and $x_0\in\mathbb R^2$, there holds
\begin{equation}\label{iter1decomp}
\begin{split}
(G(\eta)\psi)_{\mu}=\bigl((G(\eta)\psi)^{[1]}\bigr)_{\mu}
+\int_{0}^{\infty}\int_{\delta\lambda_0}^{\lambda_0}
\left(\mathcal{T}_{2,<\lambda_1}^{\lambda_0}\right)_{\mu}
\frac{d\lambda_1}{\lambda_1}\frac{d\lambda_0}{\lambda_0}
+\mathcal{R}_{\mu},
\end{split}
\end{equation}
where
\begin{equation*}
\mathcal{T}_{2,<\lambda_1}^{\lambda_0}\coloneq
\lambda_1\frac{d}{d\lambda_1}\left(
P_{<M\lambda_0}G_{<\lambda_1}
\Phi^l_{<\lambda_1}(\psi_{\lambda_0},\eta_{\lambda_0})
-\nabla\cdot (\eta_{\lambda_0}V^l_{<\lambda_1})
\right),
\end{equation*}
We also have the following four space-time bounds.
\begin{enumerate}
\item For every $\mu>0$, there holds 
\begin{equation}\label{Neta1mainpert}
\|\mathcal{R}_{\mu}\overline{u}_{\mu}^{x_0}\|_{L_{t,x}^1}
\lesssim \epsilon^4c_{\mu}^2\mu^{-2s_c}(1+\mu|x_0|)^{\frac12}.
\end{equation}
\item For every $0<\lambda\ll\mu$, the remainder also satisfies
\begin{equation}\label{firstdecompauxR}
\|\mathcal R_\lambda\overline u_\mu^{x_0}u_\mu\|_{L_{t,x}^1}
\lesssim\epsilon^5c_\mu^2\lambda^{-1}
\mu^{-2s_c}(1+\mu|x_0|)^{\frac12}.
\end{equation}
\item For every $\mu>0$, the remainder also satisfies
\begin{equation}\label{firstdecompL2R}
\|\mathcal R_\mu\|_{L_{t,x}^2}
\lesssim\epsilon^3c_\mu\mu^{\frac14-s_c}.
\end{equation}
\item Define $N_\eta^{\mathrm{err}}\coloneq
G(\eta)\psi-T_G\varphi+\nabla\cdot(T_{V^l}\eta)$. For every $\mu>0$, there holds
\begin{equation}\label{firstdecompaux}
\|(N_\eta^{\mathrm{err}})_\mu\|_{L_t^2L_x^\infty}
\lesssim\epsilon^2\mu^{-\frac34}.
\end{equation}
\end{enumerate}
\end{proposition}
The first bound \eqref{Neta1mainpert} will be the main perturbative error bound at this stage of the expansion for $N_{\eta}$, while the other three will be needed to estimate errors incurred from when the nonlinearity appears as a coefficient in our normal form analysis later. 
\begin{proof}
We begin by recording a simple algebraic identity that we will use later, 
\begin{equation}\label{paravanish}
\begin{split}
\int_0^\infty T_{\mathfrak{g}}(\Phi^{l,0}_{<\delta\lambda_0}-\varphi_{\lambda_0})\frac{d\lambda_0}{\lambda_0}=0.
\end{split}
\end{equation}
Indeed, since $\varphi=\psi-T_{B^l}\eta$, our paraproduct conventions and the definition of $\varphi$ and $\Phi^{l,0}_{<\delta\lambda_0}$ gives
\begin{equation*}
\begin{split}
\int_0^\infty(\Phi^{l,0}_{<\delta\lambda_0}-\varphi_{\lambda_0})\frac{d\lambda_0}{\lambda_0}&=\int_{0}^{\infty}S_{\lambda_0}(T_{B^l}\eta)
-B^l_{<\delta\lambda_0}\eta_{\lambda_0}\frac{d\lambda_0}{\lambda_0}=T_{B^l}\eta-T_{B^l}\eta=0.
\end{split}
\end{equation*}
Next, we show that the high-frequency tails coming from any nonlinear coefficients in $\mathcal{T}_{1,<\lambda_0}(\psi_{\lambda_0},\eta_{\lambda_0})$ contribute perturbative errors. More precisely, we show that
\begin{equation*}
\begin{split}
\mathcal{R}_{\eta}^1\coloneq \int_{0}^{\infty}\Big[
\mathcal{T}_{1,<\lambda_0}(\psi_{\lambda_0},\eta_{\lambda_0})
-P_{<M\lambda_0}G_{<\lambda_0}\Phi^l_{<\lambda_0}(\psi_{\lambda_0},\eta_{\lambda_0})+\nabla\cdot (V^l_{<\lambda_0}\eta_{\lambda_0})
\Big]\frac{d\lambda_0}{\lambda_0}
\end{split}
\end{equation*}
satisfies the three bounds
\begin{equation}\label{firstdecomperror}
\|(\mathcal{R}_{\eta}^1)_{\mu}\|_{L_t^1L_x^2}\lesssim\epsilon^4c_{\mu}\mu^{-s_c},\hspace{5mm}\|(\mathcal{R}_{\eta}^1)_{\mu}\|_{L_t^2L_x^\infty}\lesssim\epsilon^4\mu^{-\frac{3}{4}},\hspace{5mm}\|(\mathcal R_\eta^1)_\mu\|_{L_{t,x}^2}
\lesssim\epsilon^5c_\mu\mu^{\frac14-s_c}.
\end{equation}
uniformly for $\mu>0$.  
Here $M$ is the large, fixed, positive constant from Subsection \ref{paraconventions}. Below, we will show the details for the $L_t^1L_x^2$ bound. The other case is very similar and follows by modifying the H\"older exponents in the below estimates appropriately and using Bernstein and Sobolev inequalities. We first expand $\mathcal{R}_{\eta}^1$ as
\begin{equation*}
\begin{split}
\mathcal{R}_{\eta}^1= &-\int_{0}^{\infty}G_{<\lambda_0}(P_{\geq M\lambda_0}B_{<\lambda_0}\eta_{\lambda_0})\frac{d\lambda_0}{\lambda_0}
\\
&+\int_{0}^{\infty}P_{\geq M\lambda_0}G_{<\lambda_0}(\psi_{\lambda_0}-B^l_{<\lambda_0}\eta_{\lambda_0})\frac{d\lambda_0}{\lambda_0}
\\
&-\int_{0}^{\infty}\nabla\cdot (P_{\geq M\lambda_0}V_{<\lambda_0}\eta_{\lambda_0})\frac{d\lambda_0}{\lambda_0}
\\
&\coloneq \mathcal{R}^{1,1}+\mathcal{R}^{1,2}+\mathcal{R}^{1,3}.
\end{split}
\end{equation*}
To estimate $\mathcal{R}^{1,1}$, Proposition \ref{basicopbounds}, and the frequency separation between $P_{\geq M\lambda_0}B_{<\lambda_0}$ and $\eta_{\lambda_0}$ give, for every $\sigma\geq0$,
\begin{equation*}
\begin{split}
\||D|^{\sigma}G_{<\lambda_0}(P_{\geq M\lambda_0}B_{<\lambda_0}\eta_{\lambda_0})\|_{L_t^1L_x^2}
&\lesssim_{\sigma}\||D|^{\sigma+1}(P_{\geq M\lambda_0}B_{<\lambda_0}\eta_{\lambda_0})\|_{L_t^1\dot B^0_{2,1}}
\\
&\lesssim_{\sigma}\||D|^{\sigma+1}P_{\geq M\lambda_0}B_{<\lambda_0}\|_{L_t^1\dot B^0_{\infty,1}}\|\eta_{\lambda_0}\|_{L_t^\infty L_x^2}.
\end{split}
\end{equation*}
Proposition \ref{Bernsteinelliptic}, applied with $p=q=r=\infty$, $N=3$, and $\sigma_1=1$, together with \eqref{BauxiliaryPhi}, gives
\begin{equation*}
\begin{split}
\||D|^{\sigma+1}P_{\geq M\lambda_0}B_{<\lambda_0}\|_{L_t^1\dot B^0_{\infty,1}}
&\lesssim \lambda_0^{\sigma-3}
\|\nabla_x^2\eta_{<\lambda_0}\|_{L_t^4\dot B^0_{\infty,1}}^3
\||D|B_{<\lambda_0}\|_{L_t^4\dot B^0_{\infty,1}}\lesssim\epsilon^4\lambda_0^\sigma.
\end{split}
\end{equation*}
Therefore,
\begin{equation*}
\||D|^{\sigma}G_{<\lambda_0}(P_{\geq M\lambda_0}B_{<\lambda_0}\eta_{\lambda_0})\|_{L_t^1L_x^2}
\lesssim\epsilon^5c_{\lambda_0}\lambda_0^{\sigma-s_c}.
\end{equation*}
Hence, for $\sigma>s_c$ sufficiently large,
\begin{equation*}
\begin{split}
\|\mathcal{R}^{1,1}_{\mu}\|_{L_t^1L_x^2}
&\lesssim \mu^{-\sigma}\int_{0}^{\mu}\||D|^{\sigma}G_{<\lambda_0}(P_{\geq M\lambda_0}B_{<\lambda_0}\eta_{\lambda_0})\|_{L_t^1L_x^2}\frac{d\lambda_0}{\lambda_0}
\\
&\qquad+\int_{\mu}^{\infty}\|G_{<\lambda_0}(P_{\geq M\lambda_0}B_{<\lambda_0}\eta_{\lambda_0})\|_{L_t^1L_x^2}\frac{d\lambda_0}{\lambda_0}
\\
&\lesssim \mu^{-\sigma}\epsilon^5\int_{0}^{\mu}c_{\lambda_0}\lambda_0^{\sigma-s_c}\frac{d\lambda_0}{\lambda_0}
+\epsilon^5\int_{\mu}^{\infty}c_{\lambda_0}\lambda_0^{-s_c}\frac{d\lambda_0}{\lambda_0}
\\
&\lesssim \epsilon^5c_{\mu}\mu^{-s_c}.
\end{split}
\end{equation*}
To estimate $\mathcal{R}^{1,2}$, we fix an auxiliary constant $M_0\gg1$, with the constant $M$ in Subsection \ref{paraconventions} chosen at the outset so that $M\gg M_0$ with sufficient separation for the high-frequency tail estimate below. We split
\begin{equation}\label{lohiPhisplit}
\begin{split}
P_{\geq M\lambda_0}G_{<\lambda_0}\Phi^l_{<\lambda_0}(\psi_{\lambda_0},\eta_{\lambda_0})
={}&P_{\geq M\lambda_0}G_{<\lambda_0}P_{<M_0\lambda_0}\Phi^l_{<\lambda_0}(\psi_{\lambda_0},\eta_{\lambda_0})
\\
&+P_{\geq M\lambda_0}G_{<\lambda_0}P_{\geq M_0\lambda_0}\Phi^l_{<\lambda_0}(\psi_{\lambda_0},\eta_{\lambda_0}).
\end{split}
\end{equation}
As in the proof of Lemma \ref{Phiest}, we use the decomposition
\begin{equation*}
\Phi^l_{<\lambda_0}(\psi_{\lambda_0},\eta_{\lambda_0})
=\varphi_{\lambda_0}
+(T_{B^l}-B^l_{<\lambda_0})\eta_{\lambda_0}
+[S_{\lambda_0},T_{B^l}]\eta.
\end{equation*}
The estimates used there before applying $L_{<\lambda_0}$ or $M_{<\lambda_0}$, including \eqref{BauxiliaryPhi} for the coefficient terms, together with Proposition \ref{phiu-bootstrap}, give
\begin{equation*}
\|\Phi^l_{<\lambda_0}(\psi_{\lambda_0},\eta_{\lambda_0})\|_{L_t^\infty\dot B^0_{2,1}}
\lesssim \epsilon c_{\lambda_0}\lambda_0^{-\frac12-s_c}.
\end{equation*}
Consequently, Proposition \ref{basicopbounds} and Bernstein gives
\begin{equation*}
\|G_{<\lambda_0}P_{<M_0\lambda_0}\Phi^l_{<\lambda_0}(\psi_{\lambda_0},\eta_{\lambda_0})\|_{L_t^\infty\dot B^0_{2,1}}
\lesssim \epsilon c_{\lambda_0}\lambda_0^{\frac12-s_c}.
\end{equation*}
Applying Proposition \ref{Bernsteinelliptic} with $N=4$ to the first term in \eqref{lohiPhisplit}, we obtain
\begin{equation*}
\begin{split}
&\||D|^\sigma P_{\geq M\lambda_0}G_{<\lambda_0}P_{<M_0\lambda_0}\Phi^l_{<\lambda_0}(\psi_{\lambda_0},\eta_{\lambda_0})\|_{L_t^1L_x^2}
\\
&\qquad\lesssim \lambda_0^{\sigma-4}
\|\nabla_x^2\eta_{<\lambda_0}\|_{L_t^4\dot B^0_{\infty,1}}^4
\|L_{<\lambda_0}P_{<M_0\lambda_0}\Phi^l_{<\lambda_0}(\psi_{\lambda_0},\eta_{\lambda_0})\|_{L_t^\infty\dot B^0_{2,1}}
\\
&\qquad\lesssim\epsilon^5c_{\lambda_0}\lambda_0^{\sigma-s_c}.
\end{split}
\end{equation*}
On the other hand, since $\psi_{\lambda_0}$ is localized at frequency $\lambda_0$, the second term is, up to harmless fattening of the projections, of the form
\begin{equation*}
P_{\geq M\lambda_0}G_{<\lambda_0}P_{\geq M_0\lambda_0}
\left(P_{\gtrsim M_0\lambda_0}B^l_{<\lambda_0}\,\eta_{\lambda_0}\right).
\end{equation*}
After choosing $M_0$ sufficiently large, this is controlled by the same high-$B$ tail estimate as $\mathcal{R}^{1,1}$; our fixed choice of $M$ is sufficiently large relative to $M_0$. We conclude that
\begin{equation*}
\||D|^\sigma P_{\geq M\lambda_0}G_{<\lambda_0}(\psi_{\lambda_0}-B^l_{<\lambda_0}\eta_{\lambda_0})\|_{L_t^1L_x^2}
\lesssim \epsilon^5c_{\lambda_0}\lambda_0^{\sigma-s_c}.
\end{equation*}
The same argument as in the bound for $\mathcal{R}_{\mu}^{1,1}$ then gives 
\begin{equation*}
\|\mathcal{R}^{1,2}_{\mu}\|_{L_t^1L_x^2}\lesssim\epsilon^5c_{\mu}\mu^{-s_c}.
\end{equation*}
Now, we turn to $\mathcal{R}^{1,3}_{\mu}$. The analogue for $V$ of the uniform $X^p$ estimate obtained in the proof of Proposition \ref{VBunbil} gives
\begin{equation*}
\||D|V_{<\lambda_0}\|_{L_t^4\dot B^0_{p,1}}
\lesssim\epsilon\lambda_0^{\frac38-\frac2p},\hspace{5mm}8<p\leq\infty.
\end{equation*}
Bernstein's inequality and Proposition \ref{Bernsteinelliptic}, applied with $A=M$ and $N=3$, therefore give
\begin{equation*}
\begin{split}
\||D|^{\sigma+1}P_{\geq M\lambda_0}V_{<\lambda_0}\|_{L_t^1\dot B^0_{\infty,1}}\lesssim \lambda_0^{\sigma-3}
\|\nabla_x^2\eta_{<\lambda_0}\|_{L_t^4\dot B^0_{\infty,1}}^3
\||D|B_{<\lambda_0}\|_{L_t^4\dot B^0_{\infty,1}}
\lesssim\epsilon^4\lambda_0^\sigma.
\end{split}
\end{equation*}
Consequently,
\begin{equation*}
\begin{split}
\||D|^\sigma\nabla\cdot(P_{\geq M\lambda_0}V_{<\lambda_0}\eta_{\lambda_0})\|_{L_t^1L_x^2}
&\lesssim_\sigma\||D|^{\sigma+1}P_{\geq M\lambda_0}V_{<\lambda_0}\|_{L_t^1\dot B^0_{\infty,1}}
\|\eta_{\lambda_0}\|_{L_t^\infty L_x^2}
\\
&\lesssim\epsilon^5c_{\lambda_0}\lambda_0^{\sigma-s_c}.
\end{split}
\end{equation*}
Thus, as above, we have
\begin{equation*}
\|\mathcal{R}^{1,3}_{\mu}\|_{L_t^1L_x^2}\lesssim\epsilon^5c_{\mu}\mu^{-s_c}.
\end{equation*}
Combining the three estimates, and taking $\epsilon$ sufficiently small, proves the $L_t^1L_x^2$ component of \eqref{firstdecomperror}. As mentioned above, the $L_t^2L_x^{\infty}$ estimate is very similar, so we omit the routine modifications. To conclude the proof of the proposition, we use the fundamental theorem of calculus to obtain
\begin{equation*}
\begin{split}
&\int_{0}^{\infty}\left[
P_{<M\lambda_0}G_{<\lambda_0}\Phi^l_{<\lambda_0}(\psi_{\lambda_0},\eta_{\lambda_0})
-\nabla\cdot (V^l_{<\lambda_0}\eta_{\lambda_0})
\right]\frac{d\lambda_0}{\lambda_0}
\\
=&\int_{0}^{\infty}\left[
P_{<M\lambda_0}G_{<\delta\lambda_0}\Phi^{l,0}_{<\delta\lambda_0}
-\nabla\cdot (V^l_{<\delta\lambda_0}\eta_{\lambda_0})
\right]\frac{d\lambda_0}{\lambda_0}
\\
&+\int_{0}^{\infty}\int_{\delta\lambda_0}^{\lambda_0}
\mathcal{T}^{\lambda_0}_{2,<\lambda_1}
\frac{d\lambda_1}{\lambda_1}\frac{d\lambda_0}{\lambda_0}.
\end{split}
\end{equation*}
Using the conventions laid out in Subsection \ref{paraconventions}, we can write the term in the second line as
\begin{equation*}
(G(\eta)\psi)^{[1]}
+\int_{0}^{\infty}
\left(P_{<M\lambda_0}G_{<\delta\lambda_0}
-T_{\mathfrak g}\right)f_{\lambda_0}
\frac{d\lambda_0}{\lambda_0}
\end{equation*}
where we write $f_{\lambda_0}\coloneq \Phi^{l,0}_{<\delta\lambda_0}-\varphi_{\lambda_0}.$ We denote the second term on the right-hand side by $\mathcal{R}^{1,4}$
and set $\mathcal R=\mathcal R_\eta^1+\mathcal R^{1,4}$.
Above, we used \eqref{paravanish}. Therefore, \eqref{bernelliptic2}, applied with $N=2$, and Corollary \ref{Phil0varphibound} give
\begin{equation*}
P_\mu\mathcal{R}^{1,4}
=
\int_{\lambda_0\approx\mu}
P_\mu\left(
P_{<M\lambda_0}G_{<\delta\lambda_0}
-T_{\mathfrak g}
\right)f_{\lambda_0}
\frac{d\lambda_0}{\lambda_0}
+O_{L_t^1L_x^2}(\epsilon^4c_\mu\mu^{-s_c}).
\end{equation*}
Next, we can use Proposition \ref{paraexp} to expand 
\begin{equation*}
\left(P_{<M\lambda_0}G_{<\delta\lambda_0}-T_{\mathfrak g}\right)f_{\lambda_0}
=\left(G^0_{<\delta\lambda_0}
+Op(\mathfrak g^l_{<\delta\lambda_0})-T_{\mathfrak g}\right)f_{\lambda_0}
+\mathcal{E}_{f_{\lambda_0}}.
\end{equation*}
where $G_{<\delta\lambda_0}^0$ and $\mathcal{E}_{f_{\lambda_0}}$ are the main subprincipal, and error terms defined in Proposition \ref{paraexp}. By \eqref{K0bound}, \eqref{paraexperror} with integer parameter
$2$, the kernel bound for the symbol difference established
below, and Corollary \ref{Phil0varphibound}, Bernstein gives
\begin{equation*}
\begin{split}
\left\|
\left(P_{<M\lambda_0}G_{<\delta\lambda_0}
-T_{\mathfrak g}\right)f_{\lambda_0}
\right\|_{L_{t,x}^2}
&\lesssim\epsilon\lambda_0
\|f_{\lambda_0}\|_{L_{t,x}^2}
\lesssim\epsilon^3c_{\lambda_0}\lambda_0^{\frac14-s_c}.
\end{split}
\end{equation*}
Integrating over $\lambda_0\gtrsim_M\mu$, we obtain in particular,
\begin{equation*}
\|\mathcal R^{1,4}_\mu\|_{L_{t,x}^2}
\lesssim\epsilon^3
\int_{\lambda_0\gtrsim_M\mu}
c_{\lambda_0}\lambda_0^{\frac14-s_c}
\frac{d\lambda_0}{\lambda_0}
\lesssim\epsilon^3c_\mu\mu^{\frac14-s_c}.
\end{equation*}
Combining this with \eqref{firstdecomperror} proves
\eqref{firstdecompL2R}. It remains to prove the multi linear bounds
\eqref{Neta1mainpert}, \eqref{firstdecompauxR}, and stronger $L_t^2L_x^{\infty}$ error bound
\eqref{firstdecompaux}. To this end, we write
\begin{equation*}
\begin{split}
\mathcal{R}^{1,5}_{\mu}\coloneq
\int_{\lambda_0\approx\mu}
P_\mu\left(G^0_{<\delta\lambda_0}
+Op(\mathfrak g^l_{<\delta\lambda_0})-T_{\mathfrak g}\right)f_{\lambda_0}
\frac{d\lambda_0}{\lambda_0},
\\
\mathcal{R}^{1,6}_{\mu}
\coloneq
\int_{\lambda_0\approx\mu}
P_\mu\mathcal{E}_{f_{\lambda_0}}
\frac{d\lambda_0}{\lambda_0}.
\end{split}
\end{equation*}
We first note that the additional operator difference $Op(\mathfrak g^l_{<\delta\lambda_0})-T_{\mathfrak g}$ in
$\mathcal R^{1,5}$ obeys the same type of kernel bound as
$G^0_{<\delta\lambda_0}$. Indeed, by definition of $T_{\mathfrak g}$,
\begin{equation*}
\left(Op(\mathfrak g^l_{<\delta\lambda_0})-T_{\mathfrak g}\right)f_{\lambda_0}
=\int_{\nu\approx\lambda_0}
Op\left(\mathfrak g^l_{<\delta\lambda_0}
-\mathfrak g^l_{<\delta\nu}\right)S_\nu f_{\lambda_0}
\frac{d\nu}{\nu}.
\end{equation*}
We can then split the symbol difference as
\begin{equation*}
\begin{split}
\mathfrak g^l_{<\delta\lambda_0}-\mathfrak g^l_{<\delta\nu}
=P_{<M\delta\lambda_0}
\left(\mathfrak g_{<\delta\lambda_0}
-\mathfrak g_{<\delta\nu}\right)+\left(P_{<M\delta\lambda_0}-P_{<M\delta\nu}\right)g_{<\delta\nu},
\end{split}
\end{equation*}
From this, Bernstein and the fundamental theorem of calculus give,
for $\nu\approx\lambda_0$, $|\xi|\approx\lambda_0$, the symbol bound
\begin{equation*}
\begin{split}
\left|\partial_\xi^\gamma
\left(\mathfrak g^l_{<\delta\lambda_0}
-\mathfrak g^l_{<\delta\nu}\right)(x,\xi)\right|\lesssim_{\delta,N,\gamma}
\epsilon\lambda_0^{-|\gamma|}
\int_{\alpha\lesssim\delta\lambda_0}
\int_{\mathbb R^2}
\alpha^2\langle\alpha(x-z)\rangle^{-N}
\left|\nabla^2_x\eta_\alpha(z)\right|
\,dz\,\frac{d\alpha}{\alpha}
\end{split}
\end{equation*}
which by integration by parts in $\xi$, shows that the corresponding kernel $K_{\lambda_0,\nu}$ of the operator $Op(\mathfrak g^l_{<\delta\lambda_0}
-\mathfrak g^l_{<\delta\nu})S_\nu$
satisfies
\begin{equation*}
\begin{split}
|K_{\lambda_0,\nu}(x,y)|\lesssim_{\delta,N}
\epsilon\lambda_0^{2}
\langle\lambda_0(x-y)\rangle^{-N}\int_{\alpha\lesssim\delta\lambda_0}
\int_{\mathbb R^2}
\alpha^2\langle\alpha(x-z)\rangle^{-N}
\left|\nabla^2\eta_\alpha(z)\right|
\,dz\,\frac{d\alpha}{\alpha}.
\end{split}
\end{equation*}
In particular, for $\lambda_0\approx\mu$, we obtain by Bernstein,
\begin{equation*}
\begin{split}
\left\|
P_\mu\left(Op(\mathfrak g^l_{<\delta\lambda_0})-T_{\mathfrak g}\right)
f_{\lambda_0}\,\overline u_\mu^{x_0}
\right\|_{L_{t,x}^1}\lesssim_\delta
\epsilon\lambda_0^{\frac34}
\left(\int_{\alpha\lesssim\delta\lambda_0}
\|L(|D|^{\frac54}\eta_\alpha,\overline u_\mu^{x_0})\|_{L_{t,x}^2}
\frac{d\alpha}{\alpha}\right)
\|f_{\lambda_0}\|_{L_{t,x}^2}.
\end{split}
\end{equation*}
Here and below, $L$ denotes a translation-invariant bilinear form
with rapidly decaying kernels. 

Using \eqref{K0bound}, the preceding estimate,
Corollary \ref{Phil0varphibound}, and one unbalanced bilinear
estimate, we obtain
\begin{equation*}
\begin{split}
\|\mathcal{R}_{\mu}^{1,5}\overline{u}_{\mu}^{x_0}\|_{L_{t,x}^1}
&\lesssim
\int_{\lambda_0\approx\mu}
\lambda_0^{\frac34}
\left(
\int_{\alpha\lesssim\delta\lambda_0}
\|L(|D|^{\frac54}\eta_\alpha,\overline{u}_{\mu}^{x_0})\|_{L_{t,x}^2}
\frac{d\alpha}{\alpha}
\right)
\|f_{\lambda_0}\|_{L_{t,x}^2}
\frac{d\lambda_0}{\lambda_0}
\\
&\lesssim
\epsilon^4c_\mu^2\mu^{-2s_c}(1+\mu|x_0|)^{\frac12}.
\end{split}
\end{equation*}
Above, we used $\alpha$ integrability of the frequency envelope for $\eta_{\alpha}$. It remains to estimate $\mathcal{R}^{1,6}_{\mu}$. Taking the integer parameter in \eqref{paraexperror} equal to $2$ and using Corollary \ref{Phil0varphibound}, we have
\begin{equation*}
\begin{split}
\|\mathcal{E}_{f_{\lambda_0}}\|_{L_t^1L_x^2}\lesssim
\lambda_0^{-1}
\|\nabla_x^2\eta_{<\delta\lambda_0}\|_{L_t^4\dot B^0_{\infty,1}}^2
\|f_{\lambda_0}\|_{L_t^2L_x^2}\lesssim
\epsilon^4c_{\lambda_0}\lambda_0^{-s_c}.
\end{split}
\end{equation*}
Consequently,
\begin{equation*}
\|\mathcal{R}^{1,6}_{\mu}\overline{u}_{\mu}^{x_0}\|_{L_{t,x}^1}
\lesssim
\epsilon^5c_\mu^2\mu^{-2s_c}.
\end{equation*}
Combining \eqref{firstdecomperror} with the estimates above and the bound $\|u_\mu\|_{L_t^\infty L_x^2}
\lesssim\epsilon c_\mu\mu^{-s_c}$ gives
\begin{equation*}
\|\mathcal{R}_{\mu}\overline{u}_{\mu}^{x_0}\|_{L_{t,x}^1}
\lesssim
\epsilon^4c_\mu^2\mu^{-2s_c}(1+\mu|x_0|)^{\frac12}.
\end{equation*}
We next prove \eqref{firstdecompauxR}. The parts of $\mathcal R$
other than $\mathcal R^{1,5}$ have already been estimated in
$L_t^1L_x^2$. Using these estimates at the critical regularity $s_c=2$,
Bernstein, and the $L_t^{\infty}L_x^2$ bound for each factor of $u_\mu$, we get
\begin{equation*}
\begin{split}
\|(\mathcal R_\lambda-\mathcal R^{1,5}_\lambda)
\overline u_\mu^{x_0}u_\mu\|_{L_{t,x}^1}\lesssim
\lambda\|\mathcal R_\lambda-\mathcal R^{1,5}_\lambda
\|_{L_t^1L_x^2}\|u_\mu\|_{L_t^\infty L_x^2}^2\lesssim\epsilon^6c_\mu^2\lambda^{-1}\mu^{-2s_c}
\end{split}
\end{equation*}
which suffices. For $\mathcal R^{1,5}_\lambda$, we repeat the estimate above,
pairing $|D|^{5/4}\eta_\alpha$ with $\overline u_\mu^{x_0}$
in $L_{t,x}^2$, putting $f_{\lambda_0}$ in
$L_t^2L_x^\infty$, and the remaining $u_\mu$ in
$L_t^\infty L_x^2$. Indeed Corollary \ref{Phil0varphibound} (with $s=s_c=2$) and Bernstein give
\begin{equation*}
\|f_{\lambda_0}\|_{L_t^2L_x^\infty}
\lesssim\lambda_0\|f_{\lambda_0}\|_{L_{t,x}^2}
\lesssim\epsilon^2\lambda_0^{-\frac74}.
\end{equation*}
The kernel bound \eqref{K0bound}, the preceding bound for the
additional operator difference in $\mathcal R^{1,5}$, and one
unbalanced bilinear estimate therefore give
\begin{equation*}
\begin{split}
\|\mathcal R^{1,5}_\lambda
\overline u_\mu^{x_0}u_\mu\|_{L_{t,x}^1}
&\lesssim
\|u_\mu\|_{L_t^\infty L_x^2}
\int_{\lambda_0\approx\lambda}\lambda_0^{\frac34}
\left(
\int_{\alpha\lesssim\delta\lambda_0}
\|L(|D|^{\frac54}\eta_\alpha,\overline u_\mu^{x_0})
\|_{L_{t,x}^2}\frac{d\alpha}{\alpha}
\right)
\|f_{\lambda_0}\|_{L_t^2L_x^\infty}
\frac{d\lambda_0}{\lambda_0}
\\
&\lesssim
\epsilon^5c_\mu^2\mu^{-2s_c}(1+\mu|x_0|)^{\frac12}
\int_{\lambda_0\approx\lambda}\lambda_0^{-1}
\frac{d\lambda_0}{\lambda_0}
\\
&\lesssim
\epsilon^5c_\mu^2\lambda^{-1}\mu^{-2s_c}(1+\mu|x_0|)^{\frac12}.
\end{split}
\end{equation*}
Here $\alpha\lesssim\delta\lambda_0\approx\delta\lambda\ll\mu$,
and we used the integrability of the frequency envelope
in $\alpha$, as above. Combining the two estimates proves \eqref{firstdecompauxR}. It remains to prove \eqref{firstdecompaux}. Using the above identities, we have the initial decomposition
\begin{equation*}
N_\eta^{\mathrm{err}}
=\mathcal R_\eta^1
+\int_0^\infty P_{<M\lambda_0}G_{<\delta\lambda_0}
f_{\lambda_0}\frac{d\lambda_0}{\lambda_0}
+\int_0^\infty\int_{\delta\lambda_0}^{\lambda_0}
\mathcal T^{\lambda_0}_{2,<\lambda_1}
\frac{d\lambda_1}{\lambda_1}\frac{d\lambda_0}{\lambda_0}.
\end{equation*}
The contribution of $\mathcal{R}_{\eta}^1$ is handled by the second bound in \eqref{firstdecomperror}. For the endpoint integral term above, Proposition \ref{basicopbounds}, Corollary
\ref{Phil0varphibound}, and Bernstein give
\begin{equation*}
\|P_\mu P_{<M\lambda_0}G_{<\delta\lambda_0}
f_{\lambda_0}\|_{L_t^2L_x^\infty}
\lesssim\mu\lambda_0\|f_{\lambda_0}\|_{L_{t,x}^2}
\lesssim\epsilon^2\mu\lambda_0^{-\frac74}.
\end{equation*}
The output projection $P_{\mu}$ restricts the integral to
$\lambda_0\gtrsim_M\mu$, so this contribution is bounded by
$\epsilon^2\mu^{-3/4}$.

It remains to estimate the double integral term. To do this, we have to compute $\mathcal{T}_{2,<\lambda_1}^{\lambda_0}$ and collect a few auxiliary estimates. First, for
$\delta\lambda_0\leq\lambda_1\leq\lambda_0$, differentiating both
the coefficients and their output projections $P_{<M\lambda_1}$ gives
\begin{equation*}
\begin{split}
\lambda_1\frac{d}{d\lambda_1}B^l_{<\lambda_1}
&=P_{<M\lambda_1}\bigl(
L_{<\lambda_1}\Phi^0_{<\lambda_1}
+\eta_{\lambda_1}L_{<\lambda_1}B_{<\lambda_1}
\bigr)
+\left(\lambda_1\frac{d}{d\lambda_1}P_{<M\lambda_1}\right)
B_{<\lambda_1},
\\
\lambda_1\frac{d}{d\lambda_1}V^l_{<\lambda_1}
&=P_{<M\lambda_1}\bigl(
M_{<\lambda_1}\Phi^0_{<\lambda_1}
+\eta_{\lambda_1}M_{<\lambda_1}B_{<\lambda_1}
\bigr)
+\left(\lambda_1\frac{d}{d\lambda_1}P_{<M\lambda_1}\right)
V_{<\lambda_1}.
\end{split}
\end{equation*}
Lemma \ref{Phiest}, \eqref{BauxiliaryPhi}, Proposition \ref{basicopbounds}, Proposition \ref{Bernsteinelliptic} and Proposition \ref{VBunbil} therefore imply
\begin{equation*}
\left\|\lambda_1\frac{d}{d\lambda_1}B^l_{<\lambda_1}
\right\|_{L_t^4L_x^\infty}
+
\left\|\lambda_1\frac{d}{d\lambda_1}V^l_{<\lambda_1}
\right\|_{L_t^4L_x^\infty}
\lesssim\epsilon\lambda_1^{-\frac58}.
\end{equation*}
Integrating this estimate over $\lambda_1\approx_{\delta}\lambda_0$ and using the decomposition
\begin{equation*}
\Phi^l_{<\lambda_1}(\psi_{\lambda_0},\eta_{\lambda_0})
=\varphi_{\lambda_0}
+(T_{B^l}-B^l_{<\lambda_1})\eta_{\lambda_0}
+[S_{\lambda_0},T_{B^l}]\eta,
\end{equation*}
gives
\begin{equation}\label{PhiCoarseL4}
\begin{split}
\|\Phi^l_{<\lambda_1}(\psi_{\lambda_0},\eta_{\lambda_0})
\|_{L_t^4L_x^\infty}\lesssim_\delta\epsilon\lambda_0^{-\frac{13}{8}},\hspace{5mm}\|\Phi^l_{<\lambda_1}(\psi_{\lambda_0},\eta_{\lambda_0})
\|_{L_t^4\dot B^1_{\infty,1}}\lesssim_\delta\epsilon\lambda_0^{-\frac58}.
\end{split}
\end{equation}
Here the commutator is estimated using \eqref{PhiCommRep} and
\eqref{BauxiliaryPhi}, and the second bound follows from Bernstein and the first bound.
For $A=L,M$, Proposition \ref{basicopbounds} and Bernstein, with $\sigma\geq 0$ gives
\begin{equation*}
\begin{split}
\||D|^{\sigma}A_{<\lambda_1}
\Phi^l_{<\lambda_1}(\psi_{\lambda_0},\eta_{\lambda_0})
\|_{L_t^4L_x^\infty}
\lesssim_{\delta,\sigma}\epsilon\lambda_0^{-\frac58+\sigma}.
\end{split}
\end{equation*}
The linearization formula for $G$ now yields the identity
\begin{equation*}
\begin{split}
\mathcal T^{\lambda_0}_{2,<\lambda_1}
&=-P_{<M\lambda_0}G_{<\lambda_1}
\left(
\eta_{\lambda_1}L_{<\lambda_1}
\Phi^l_{<\lambda_1}(\psi_{\lambda_0},\eta_{\lambda_0})
\right)
\\
&-P_{<M\lambda_0}G_{<\lambda_1}
\left(
\eta_{\lambda_0}
\lambda_1\frac{d}{d\lambda_1}B^l_{<\lambda_1}
\right)
\\
&-P_{<M\lambda_0}\nabla\cdot
\left(
\eta_{\lambda_1}M_{<\lambda_1}
\Phi^l_{<\lambda_1}(\psi_{\lambda_0},\eta_{\lambda_0})
\right)
\\
&-\nabla\cdot
\left(
\eta_{\lambda_0}
\lambda_1\frac{d}{d\lambda_1}V^l_{<\lambda_1}
\right).
\end{split}
\end{equation*}
Using Bernstein for the frequency-truncated coefficients,
the product estimate in $\dot B^1_{\infty,1}$,
$\|\eta_{\lambda_j}\|_{L_t^4L_x^\infty}
\lesssim\epsilon\lambda_j^{-9/8}$ for $j=0,1$,
and the preceding bounds, we obtain
\begin{equation*}
\|\mathcal T^{\lambda_0}_{2,<\lambda_1}\|_{L_t^2L_x^\infty}
\lesssim_\delta
\lambda_0
(\epsilon\lambda_0^{-\frac98})
(\epsilon\lambda_0^{-\frac58})
=\epsilon^2\lambda_0^{-\frac34}.
\end{equation*}
Since the output is supported in $|\xi|\lesssim M\lambda_0$,
this gives
\begin{equation*}
\left\|P_\mu\int_0^\infty\int_{\delta\lambda_0}^{\lambda_0}
\mathcal T^{\lambda_0}_{2,<\lambda_1}
\frac{d\lambda_1}{\lambda_1}\frac{d\lambda_0}{\lambda_0}
\right\|_{L_t^2L_x^\infty}
\lesssim_\delta\epsilon^2\mu^{-\frac34}.
\end{equation*}
Combining these estimates proves the bound for
$N_\eta^{\mathrm{err}}$ in \eqref{firstdecompaux}. This concludes the proof.
\end{proof}


\subsection{Extraction of \texorpdfstring{$N_{\eta}^{[2]}$}{}}\label{subsec:Neta2}
Now, we move to the second iteration, which is by far the most delicate in this section. To compute $\mathcal{T}_{2,<\lambda_1}$ (dropping the $\lambda_0$ superscript for simplicity) we will need the following linearization formulas.
\begin{lemma}\label{vphlin}
For fixed Schwartz functions $f$ and $g$, there holds
\begin{equation*}
\begin{split}
\lambda\frac{d}{d\lambda}\Phi^l_{<\lambda}(f,g)
&=-gP_{<M\lambda}L_{<\lambda}\Phi^0_{<\lambda}
-gP_{<M\lambda}\left(\eta_{\lambda}L_{<\lambda}B_{<\lambda}\right)
-g\widetilde P_{M\lambda}B_{<\lambda},
\end{split}
\end{equation*}
and
\begin{equation*}
\begin{split}
\lambda\frac{d}{d\lambda}\nabla\cdot (gV^l_{<\lambda})
&=\nabla\cdot\left(gP_{<M\lambda}M_{<\lambda}\Phi^0_{<\lambda}\right)
+\nabla\cdot\left(gP_{<M\lambda}\left(\eta_{\lambda}M_{<\lambda}B_{<\lambda}\right)\right)
\\
&\quad+\nabla\cdot\left(g\widetilde P_{M\lambda}V_{<\lambda}\right),
\end{split}
\end{equation*}
where $\Phi^0_{<\lambda}
\coloneq\Phi_{<\lambda}(\psi_{\lambda},\eta_{\lambda})
=\psi_{\lambda}-B_{<\lambda}\eta_{\lambda}$, and $\widetilde P_{M\lambda}\coloneq\lambda\frac{d}{d\lambda}P_{<M\lambda}$.
\end{lemma}
\begin{proof}
Since $f$ and $g$ are fixed with respect to $\lambda$, the linearization identity for $B_{<\lambda}$ gives
\begin{equation*}
\begin{split}
\lambda\frac{d}{d\lambda}\Phi^l_{<\lambda}(f,g)
&=-g\lambda\frac{d}{d\lambda}\left(P_{<M\lambda}B_{<\lambda}\right)
\\
&=-gP_{<M\lambda}\left(
L_{<\lambda}\Phi^0_{<\lambda}
+\eta_{\lambda}L_{<\lambda}B_{<\lambda}
\right)
-g\widetilde P_{M\lambda}B_{<\lambda},
\end{split}
\end{equation*}
which proves the first formula. The second follows in exactly the same way from
\begin{equation*}
\lambda\frac{d}{d\lambda}V_{<\lambda}
=M_{<\lambda}\Phi^0_{<\lambda}
+\eta_{\lambda}M_{<\lambda}B_{<\lambda}.
\end{equation*}
\end{proof}
We will later use the following corollary, which extracts the leading part of the linearization of $\Phi^l_{<\lambda}$ above and allows us to discard high-frequency tails. We use the notation
\begin{equation}\label{optrunc}
A^{l,\beta}_{<\lambda}h\coloneq P_{<M\beta}A_{<\lambda}h,\hspace{5mm}A\in\{G,M,L\}.
\end{equation}
If $\beta=\lambda$, we simply write $A^l_{<\lambda}h$.
\begin{corollary}\label{Philintrunc}
Let $\lambda_j,\lambda_k>0$ be frequency parameters with comparable size $\lambda_j\approx\lambda_k$. For $m\in\{j,k\}$, write
\begin{equation}\label{Phitruncnot}
\Phi^{l,m}_{<\lambda_j}
\coloneq\Phi^l_{<\lambda_j}(\psi_{\lambda_m},\eta_{\lambda_m}),\hspace{5mm}\Phi^{m}_{<\lambda_j}\coloneq \Phi_{<\lambda_j}(\psi_{\lambda_m},\eta_{\lambda_m}).
\end{equation}
Then, there holds
\begin{equation*}
\lambda_j\frac{d}{d\lambda_j}\Phi^{l,k}_{<\lambda_j}
=-\eta_{\lambda_k}L^l_{<\lambda_j}\Phi^{l,j}_{<\lambda_j}
-\eta_{\lambda_k}\eta_{\lambda_j}L^l_{<\lambda_j}B^l_{<\lambda_j}
+O_{L_t^1L_x^2}\left(\epsilon^3c_{\lambda_k}\lambda_k^{-1-s_c}\right).
\end{equation*}
\end{corollary}
\begin{proof}
Applying Lemma \ref{vphlin} with $(f,g)=(\psi_{\lambda_k},\eta_{\lambda_k})$, we obtain
\begin{equation*}
\begin{split}
\lambda_j\frac{d}{d\lambda_j}\Phi^{l,k}_{<\lambda_j}=-\eta_{\lambda_k}P_{<M\lambda_j}L_{<\lambda_j}\Phi^j_{<\lambda_j}
-\eta_{\lambda_k}P_{<M\lambda_j}\left(\eta_{\lambda_j}L_{<\lambda_j}B_{<\lambda_j}\right)-\eta_{\lambda_k}\widetilde P_{M\lambda_j}B_{<\lambda_j}.
\end{split}
\end{equation*}
Set $P_j\coloneq P_{<M\lambda_j}$ and $Q_j\coloneq I-P_j$. Since $B^l_{<\lambda_j}=P_jB_{<\lambda_j}$, we have
\begin{equation*}
\Phi^j_{<\lambda_j}
=\Phi^{l,j}_{<\lambda_j}
-Q_jB_{<\lambda_j}\eta_{\lambda_j},
\end{equation*}
and thus
\begin{equation}\label{Philinerrorterms1}
-\eta_{\lambda_k}P_jL_{<\lambda_j}\Phi^j_{<\lambda_j}
=-\eta_{\lambda_k}L^l_{<\lambda_j}\Phi^{l,j}_{<\lambda_j}
+\eta_{\lambda_k}L^l_{<\lambda_j}
\left(Q_jB_{<\lambda_j}\eta_{\lambda_j}\right).
\end{equation}
Moreover,
\begin{equation}\label{Philinerrorterms2}
\begin{split}
\eta_{\lambda_k}P_j\left(\eta_{\lambda_j}L_{<\lambda_j}B_{<\lambda_j}\right)&=\eta_{\lambda_k}\eta_{\lambda_j}L^l_{<\lambda_j}B^l_{<\lambda_j}
+\eta_{\lambda_k}P_j\left(\eta_{\lambda_j}L_{<\lambda_j}Q_jB_{<\lambda_j}\right)
\\
&+\eta_{\lambda_k}[P_j,\eta_{\lambda_j}]L_{<\lambda_j}B^l_{<\lambda_j}.
\end{split}
\end{equation}
It remains to estimate the error terms in these identities and the term containing $\widetilde P_{M\lambda_j}B_{<\lambda_j}$. We will show the details for the error term in \eqref{Philinerrorterms1} and $\tilde{P}_{M\lambda_j}B_{<\lambda_j}$ along the way. The errors in \eqref{Philinerrorterms2} are similar. We first record the high-frequency bounds used below. First, by \ref{VBunbil} and the bootstrap in Section \ref{sec:bootstrap}
\begin{equation*}
\|\nabla_x^2\eta_{<\lambda_j}\|_{L_t^4\dot B^0_{\infty,1}}
\lesssim\epsilon\lambda_j^{\frac78},
\qquad
\||D|B_{<\lambda_j}\|_{L_t^4\dot B^0_{\infty,1}}
\lesssim\epsilon\lambda_j^{\frac38}.
\end{equation*}
Therefore, Proposition \ref{Bernsteinelliptic}, applied with $N=3$, $p=q=r=\infty$ and $\sigma=\sigma_1=1$, gives, for a fixed $0<c<1$,
\begin{equation*}
\begin{split}
\|\widetilde P_{M\lambda_j}B_{<\lambda_j}\|_{L_t^1L_x^\infty}
&\lesssim\lambda_j^{-1}
\||D|P_{\geq cM\lambda_j}B_{<\lambda_j}\|_{L_t^1\dot B^0_{\infty,1}}
\\
&\lesssim\lambda_j^{-4}
\|\nabla_x^2\eta_{<\lambda_j}\|_{L_t^4\dot B^0_{\infty,1}}^3
\||D|B_{<\lambda_j}\|_{L_t^4\dot B^0_{\infty,1}}
\\
&\lesssim\epsilon^4\lambda_j^{-1}.
\end{split}
\end{equation*}
Consequently, $\|\eta_{\lambda_k}\widetilde P_{M\lambda_j}B_{<\lambda_j}\|_{L_t^1L_x^2}
\lesssim\epsilon^5c_{\lambda_k}\lambda_k^{-1-s_c}.$ The same argument also gives
\begin{equation*}
\||D|P_{\geq cM\lambda_j}B_{<\lambda_j}\|_{L_t^1\dot B^0_{\infty,1}}
\lesssim\epsilon^4.
\end{equation*}
Proposition \ref{basicopbounds}, the preceding estimate, and the $L_t^{\infty}\dot{B}_{\infty,1}^0$ bound for $\eta_{\lambda_j}$ therefore give
\begin{equation*}
\left\|\eta_{\lambda_k}L^l_{<\lambda_j}
\left(Q_jB_{<\lambda_j}\eta_{\lambda_j}\right)\right\|_{L_t^1L_x^2}
\lesssim\epsilon^6c_{\lambda_k}\lambda_k^{-1-s_c}
\end{equation*}
which handles the error term in \eqref{Philinerrorterms1}. The error terms in \eqref{Philinerrorterms2} are handled similarly, observing that each term possesses at least one high-frequency tail factor. We omit the straightforward modifications. This proves the corollary.
\end{proof}
Similarly to Proposition \ref{firstdecomposition}, we now seek to extract the second term $(G(\eta)\psi)^{[2]}$ in the expansion of $G(\eta)\psi$. This will be the most technically demanding part of the expansion.
\begin{proposition}\label{seconddecomposition}
Let $M$ and $\delta$ be fixed as in Subsection \ref{paraconventions}, with $M\gg1$ and $M\delta\ll1$. Assume the bootstrap hypotheses in Section \ref{sec:bootstrap}. Let $(G(\eta)\psi)^{[1]}$ be as in \eqref{G1def}. Define
\begin{equation}\label{G2def}
(G(\eta)\psi)^{[2]}\coloneq -T_{\mathfrak{g}}\bigl(\Pi(\eta,T_{\mathfrak{l}}\varphi)\bigr)-\nabla\cdot \Pi(\eta,T_{\mathfrak{m}}\varphi).
\end{equation}
Then, for every $\mu\in2^{\mathbb Z}$ and $x_0\in\mathbb R^2$, there holds
\begin{equation}\label{iter2decomp}
\begin{split}
(G(\eta)\psi)_{\mu}&=(G(\eta)\psi)^{[1]}_{\mu}+(G(\eta)\psi)^{[2]}_{\mu}+\int_{0}^{\infty}\int_{\delta\lambda_0}^{\lambda_0}\int_{\delta\lambda_1}^{\lambda_1}(\mathcal{T}_{3,<\lambda_2})_{\mu}\frac{d\lambda_2}{\lambda_2}\frac{d\lambda_1}{\lambda_1}\frac{d\lambda_0}{\lambda_0}+(\mathcal{R}_{\eta}^{2})_{\mu},
\end{split}
\end{equation}
where we have the following error bounds
\begin{enumerate}
\item For every $\mu>0$ there holds
\begin{equation}\label{R2mainbound}
\|(\mathcal{R}_{\eta}^2)_{\mu}\overline{u}^{x_0}_{\mu}\|_{L_{t,x}^1}\lesssim \epsilon^4 c_{\mu}^2\mu^{-2s_c}(1+\mu|x_0|)^{\frac12}.
\end{equation}
\item For every $0<\lambda\ll\mu$, the remainder also satisfies
\begin{equation}\label{seconddecompaux}
\|(\mathcal R_\eta^2)_\lambda
\overline u_\mu^{x_0}u_\mu\|_{L_{t,x}^1}
\lesssim\epsilon^4c_\mu^2\lambda^{-1}
\mu^{-2s_c}(1+\mu|x_0|)^{\frac12}.
\end{equation}
\item For every $\mu>0$, there holds
\begin{equation}\label{seconddecompL2}
\|(\mathcal{R}^2_{\eta})_\mu\|_{L_{t,x}^2}
\lesssim_\delta\epsilon^3c_\mu\mu^{\frac14-s_c}.
\end{equation}
\end{enumerate}
Moreover, $\mathcal T_{3,<\lambda_2}$ is defined as follows: For fixed $\lambda_0$ and $\lambda_1$ we set
\begin{equation}\label{T2deriv}
\mathcal{T}_{3,<\lambda_2}\coloneq P_{<M\lambda_0}\lambda_2\frac{d}{d\lambda_2}\mathcal{T}^{l}_{2,<\lambda_2}\coloneq P_{<M\lambda_0}\lambda_2\frac{d}{d\lambda_2}(\mathcal{T}_{2,<\lambda_2}^{main}+\mathcal{T}_{2,<\lambda_2}^{err}) 
\end{equation}
where we define
\begin{equation*}
\begin{split}
\mathcal{T}_{2,<\lambda_2}^{main}\coloneq &-G^{l,\lambda_0}_{<\lambda_2}(\eta_{\lambda_1}L^{l,\lambda_0}_{<\lambda_2}\Phi^{l,0}_{<\lambda_2})-\nabla\cdot (\eta_{\lambda_1}M^{l,\lambda_0}_{<\lambda_2}\Phi^{l,0}_{<\lambda_2})
\\
&-G^{l,\lambda_0}_{<\lambda_2}(\eta_{\lambda_0}L^{l,\lambda_1}_{<\lambda_2}\Phi^{l,1}_{<\lambda_2})-\nabla\cdot (\eta_{\lambda_0}M^{l,\lambda_1}_{<\lambda_2}\Phi^{l,1}_{<\lambda_2}),
\end{split}
\end{equation*}
and
\begin{equation*}
\mathcal{T}_{2,<\lambda_2}^{err}\coloneq-G^{l,\lambda_0}_{<\lambda_2}(\eta_{\lambda_0}\eta_{\lambda_1}L^{l,\lambda_2}_{<\lambda_2}B^l_{<\lambda_2})-\nabla\cdot (\eta_{\lambda_0}\eta_{\lambda_1}M^{l,\lambda_2}_{<\lambda_2}B^l_{<\lambda_2}).
\end{equation*}
\end{proposition}
Above we have suppressed the dependence on $\lambda_0$ and $\lambda_1$ in $\mathcal T_{3,<\lambda_2}$ and $\mathcal T^l_{2,<\lambda_2}$ and have adopted the notation \eqref{optrunc} and \eqref{Phitruncnot}.
\begin{proof}
By Lemma \ref{vphlin} and the linearization formula \eqref{Glin}, we obtain the decomposition
\begin{equation*}
\int_{0}^{\infty}\int_{\delta\lambda_0}^{\lambda_0}\mathcal{T}_{2,<\lambda_1}^{\lambda_0}\frac{d\lambda_1}{\lambda_1}\frac{d\lambda_0}{\lambda_0}=\int_{0}^{\infty}\int_{\delta\lambda_0}^{\lambda_0}P_{<M\lambda_0}\mathcal{T}^{l}_{2,<\lambda_1}\frac{d\lambda_1}{\lambda_1}\frac{d\lambda_0}{\lambda_0}+\mathcal{R}^{2,1},
\end{equation*}
where $\mathcal{T}^{l}_{2,<\lambda_1}$ is defined as in \eqref{T2deriv}, with $\lambda_2=\lambda_1$, and $\mathcal{R}^{2,1}$ consists of high-frequency tails. Corollary \ref{Philintrunc}, Proposition \ref{Bernsteinelliptic}, and the same argument used to prove \eqref{firstdecomperror} give
\begin{equation}\label{R21bound}
\|\mathcal{R}^{2,1}_{\mu}\|_{L_t^1L_x^2}\lesssim \epsilon^4c_{\mu}\mu^{-s_c}.
\end{equation}
We omit the details, as this is a straightforward application of the high-frequency tail estimates already used in the proof of Proposition \ref{firstdecomposition}. By the fundamental theorem of calculus and the definition of $\mathcal{T}_{3,<\lambda_2}$, we have
\begin{equation*}
P_{<M\lambda_0}\mathcal{T}^{l}_{2,<\lambda_1}=P_{<M\lambda_0}\mathcal{T}^{l}_{2,<\delta\lambda_1}+\int_{\delta\lambda_1}^{\lambda_1}\mathcal{T}_{3,<\lambda_2}\frac{d\lambda_2}{\lambda_2}.
\end{equation*}
In view of Proposition \ref{firstdecomposition} and \eqref{R21bound}, it then suffices to show that
\begin{equation}\label{T2main}
\int_{0}^{\infty}\int_{\delta\lambda_0}^{\lambda_0}P_{<M\lambda_0}\mathcal{T}_{2,<\delta\lambda_1}^{main}\frac{d\lambda_1}{\lambda_1}\frac{d\lambda_0}{\lambda_0}=(G(\eta)\psi)^{[2]}+\int_{0}^{\infty}\int_{\delta\lambda_0}^{\lambda_0}\mathcal{R}\frac{d\lambda_1}{\lambda_1}\frac{d\lambda_0}{\lambda_0}
\end{equation}
and
\begin{equation}\label{T2err}
P_{<M\lambda_0}\mathcal{T}_{2,<\delta\lambda_1}^{err}=\mathcal{R},
\end{equation}
where, for the remainder of this proof, $\mathcal{R}$ denotes a perturbative error term satisfying
\begin{equation*}
\int_{0}^{\infty}\int_{\delta\lambda_0}^{\lambda_0}\|\mathcal{R}_{\mu}\overline{u}^{x_0}_{\mu}\|_{L_{t,x}^1}\frac{d\lambda_1}{\lambda_1}\frac{d\lambda_0}{\lambda_0}\lesssim \epsilon^4c_{\mu}^2\mu^{-2s_c}(1+\mu|x_0|)^{\frac12}
\end{equation*}
that is allowed to change from line to line. 
We will later verify that for $0<\lambda\ll\mu$, we also have the bound
\begin{equation}\label{seconddecompauxR}
\int_0^\infty\int_{\delta\lambda_0}^{\lambda_0}
\|\mathcal R_\lambda\overline u_\mu^{x_0}u_\mu\|_{L_{t,x}^1}
\frac{d\lambda_1}{\lambda_1}\frac{d\lambda_0}{\lambda_0}
\lesssim\epsilon^4c_\mu^2\lambda^{-1}
\mu^{-2s_c}(1+\mu|x_0|)^{\frac12}.
\end{equation}
We postpone this additional estimate until the end of the proof.

Below, after applying $P_\mu$, we may assume that $\lambda_0\gtrsim_M\mu$, since the integrands in \eqref{T2main} and \eqref{T2err} have output frequency $\lesssim M\lambda_0$. We begin with the most difficult part, which is to establish \eqref{T2main}. We will show in detail that the first term in the definition of $\mathcal{T}_{2,<\delta\lambda_1}^{main}$, given by
\begin{equation*}
-G^{l,\lambda_0}_{<\delta\lambda_1}(\eta_{\lambda_1}L^{l,\lambda_0}_{<\delta\lambda_1}\Phi^{l,0}_{<\delta\lambda_1}),
\end{equation*}
satisfies
\begin{equation}\label{firsttermT2main}
-G^{l,\lambda_0}_{<\delta\lambda_1}(\eta_{\lambda_1}L^{l,\lambda_0}_{<\delta\lambda_1}\Phi^{l,0}_{<\delta\lambda_1})=-T_{\mathfrak{g}}\bigl(\eta_{\lambda_1}T_{\mathfrak{l}}\varphi_{\lambda_0}\bigr)+\mathcal{R}.
\end{equation}
We will then remark at the end on the similar decompositions satisfied by the other terms in $\mathcal{T}_{2,<\delta\lambda_1}^{main}$, which are handled in an almost identical fashion, and use them to assemble \eqref{T2main}.

Our first goal is to replace $\Phi_{<\delta\lambda_1}^{l,0}$ by $\varphi_{\lambda_0}$ modulo a perturbative error. That is, we claim that
\begin{equation}\label{replacebyphi}
G^{l,\lambda_0}_{<\delta\lambda_1}(\eta_{\lambda_1}L^{l,\lambda_0}_{<\delta\lambda_1}(\Phi^{l,0}_{<\delta\lambda_1}-\varphi_{\lambda_0}))=\mathcal{R}.
\end{equation}
It is tempting to use Corollary \ref{Phil0varphibound}, but this will unfortunately not work in the balanced scenario $\mu\approx\lambda_0\approx\lambda_1$ because $\eta_{\lambda_1}$ and $u_{\mu}^{x_0}$ might interact to produce low output frequency. Instead, we have to expand $\Phi_{<\delta\lambda_1}^{l,0}-\varphi_{\lambda_0}$ more precisely. Let us write for shorthand
\begin{equation*}
\mathcal{F}_{\lambda_0,\lambda_1}\coloneq\Phi^{l,0}_{<\delta\lambda_1}-\varphi_{\lambda_0}.
\end{equation*}
By our paraproduct convention and Proposition \ref{Leibnizcommutator}, and the definition of $\Phi_{<\delta\lambda_1}^{l,0}$ and $\varphi_{\lambda_0}$, we have the representation formula
\begin{equation}\label{Phidifferror}
\begin{split}
\mathcal{F}_{\lambda_0,\lambda_1}
&=\int_{\lambda_2\approx\lambda_0}
L\left(B^l_{<\delta\lambda_2}-B^l_{<\delta\lambda_1},\widetilde\eta_{\lambda_2}\right)\frac{d\lambda_2}{\lambda_2}+\lambda_0^{-1}\int_{\lambda_2\approx\lambda_0}
L\left(\nabla B^l_{<\delta\lambda_2},\widetilde\eta_{\lambda_2}\right)\frac{d\lambda_2}{\lambda_2},
\end{split}
\end{equation}
where, as usual, $L$ denotes a translation-invariant multilinear form with rapidly decaying kernel (both $L$ above are possibly different bilinear forms). We further decompose the $B$ difference
\begin{equation*}
\begin{split}
B^l_{<\delta\lambda_2}-B^l_{<\delta\lambda_1}=P_{<M\delta\lambda_2}\left(B_{<\delta\lambda_2}-B_{<\delta\lambda_1}\right)+\left(P_{<M\delta\lambda_2}-P_{<M\delta\lambda_1}\right)B_{<\delta\lambda_1}.
\end{split}
\end{equation*}
The first difference factor can be represented using the linearization formula for $B$ as
\begin{equation}\label{Bfundcalc}
B_{<\delta\lambda_2}-B_{<\delta\lambda_1}
=\int_{\delta\lambda_1}^{\delta\lambda_2}\left(
L_{<\beta}\varphi_\beta+E^L_\beta+\eta_\beta L_{<\beta}B_{<\beta}
\right)\frac{d\beta}{\beta},
\end{equation}
where $E^L_\beta=L_{<\beta}(\Phi^0_{<\beta}-\varphi_\beta)$ is as in Lemma \ref{Phiest}. Notice that every frequency parameter $\beta$ in \eqref{Bfundcalc} satisfies $\beta\approx_\delta\lambda_0$ and $\beta\ll\lambda_0$. For $\lambda_2\approx\lambda_0$, define the principal low-frequency coefficient
\begin{equation*}
\begin{split}
\mathcal{C}_{\lambda_2,\lambda_1}\coloneq{}&
\int_{\delta\lambda_1}^{\delta\lambda_2}P_{<M\delta\lambda_2}L_{<\beta}\varphi_\beta\frac{d\beta}{\beta}
+\left(P_{<M\delta\lambda_2}-P_{<M\delta\lambda_1}\right)B_{<\delta\lambda_1}+\lambda_0^{-1}\nabla B^l_{<\delta\lambda_2},
\end{split}
\end{equation*}
and the remainder coefficient
\begin{equation*}
\mathcal{E}_{\lambda_2,\lambda_1}\coloneq
\int_{\delta\lambda_1}^{\delta\lambda_2}P_{<M\delta\lambda_2}\left(E^L_\beta+\eta_\beta L_{<\beta}B_{<\beta}\right)\frac{d\beta}{\beta}.
\end{equation*}
This gives the representation
\begin{equation*}
\mathcal{F}_{\lambda_0,\lambda_1}
=\int_{\lambda_2\approx\lambda_0}L\left(\mathcal{C}_{\lambda_2,\lambda_1}+\mathcal{E}_{\lambda_2,\lambda_1},\widetilde\eta_{\lambda_2}\right)\frac{d\lambda_2}{\lambda_2}
\end{equation*}
where, by a slight abuse of notation, all expressions involving $\mathcal C_{\lambda_2,\lambda_1}$ are understood componentwise, with a separate translation-invariant form $L$ for each scalar or vector component.
We first record the estimates for the principal coefficient. A slight variation of the bounds \eqref{Lphiu-bound}, \eqref{Bunbalancedeta}, and the truncated $(B,u)$ bounds in Proposition \ref{VBunbil}, followed by summation of the frequency envelope in the low-frequency slot, give the unbalanced bilinear estimates
\begin{equation}\label{Phidifflowpairings}
\begin{split}
\|L(\mathcal{C}_{\lambda_2,\lambda_1},\overline{u}^{x_0}_\mu)\|_{L_{t,x}^2}
&\lesssim\epsilon^2c_\mu\lambda_0^{-\frac34}\mu^{-s_c}(1+\mu|x_0|)^{\frac12},
\qquad \mu\approx\lambda_0,
\\
\|L(\mathcal{C}_{\lambda_2,\lambda_1},\widetilde\eta_{\lambda_2})\|_{L_{t,x}^2}
&\lesssim\epsilon^2c_{\lambda_0}\lambda_0^{-\frac34-s_c}.
\end{split}
\end{equation}
Lemma \ref{Phiest}, the $L_t^2L_x^\infty$ estimate for $B_r^1$, and the compatible high-frequency tail estimate used in the proof of \eqref{Brestimates}, applied at the critical regularity $s=2$ (along with summation in the output frequency parameter), also give
\begin{equation}\label{Phidifflowerror}
\|\mathcal{E}_{\lambda_2,\lambda_1}\|_{L_t^2L_x^\infty}
\lesssim_K\epsilon^2\lambda_0^{-\frac34}.
\end{equation}
Now, to estimate the contribution of \eqref{replacebyphi}, we separate into two cases. Let us write $\Lambda\coloneq\max(\lambda_1,\lambda_2)$. We also take the fattened projection defining $\widetilde\eta_{\lambda_2}$ to be supported in $c^{-2}\lambda_2\leq|\xi|\leq c^2\lambda_2$ where $c$ is as in Subsection \ref{subsec:littlewoodpaley}. The first case consider is $\mu\sim_c\Lambda$. To account for shifts from low frequency para-coefficients, we let $\widetilde P_\mu$ be a slightly fattened projection with support contained in an $O(M\delta\Lambda)$ enlargement of the $P_{\mu}$ annulus $c^{-1}\mu\leq|\xi|\leq c\mu$. The argument used in Lemma \ref{balancedrestrictedoutput} applied by expanding $\eta=\Re u$, gives
\begin{equation*}
\|\widetilde P_\mu|D|^{\frac12}L(\eta_{\lambda_1},\widetilde\eta_{\lambda_2})\|_{L_{t,x}^2}
\lesssim_\delta\epsilon^2c_{\lambda_0}\lambda_0^{-\frac34-s_c}.
\end{equation*}
Consequently, by Bernstein and $\mu\sim_c\Lambda$ and the fact that $\Lambda\approx\lambda_0$, we have
\begin{equation}\label{eta1esta2sep}
\|\widetilde P_\mu L(\eta_{\lambda_1},\widetilde\eta_{\lambda_2})\|_{L_{t,x}^2}
\lesssim_\delta\epsilon^2c_{\lambda_0}\lambda_0^{-\frac54-s_c}.
\end{equation}
Next, by applying the first estimate in Proposition \ref{basicopbounds}, together with its immediate analogue for $L$, dyadically, we separate each operator $G_{<\delta\lambda_1}^{l,\lambda_0}$ and $L_{<\delta\lambda_1}^{l,\lambda_0}$ in \eqref{replacebyphi} into its principal paradifferential part and a subprincipal error. To describe the setup, we write $G_{\mathrm{pr}}$ and $L_{\mathrm{pr}}$ for the corresponding principal operators, with the output localization retained. Suppressing for a moment the subscript and superscript parameters in each operator, we organize the decomposition according to the identity
\begin{equation}\label{GLdecomposition}
G(\eta_{\lambda_1}LF)=G_{\mathrm{pr}}(\eta_{\lambda_1}L_{\mathrm{pr}}F)+(G-G_{\mathrm{pr}})(\eta_{\lambda_1}LF)+G_{\mathrm{pr}}(\eta_{\lambda_1}(L-L_{\mathrm{pr}})F).
\end{equation}
Thus the last two terms include all contributions containing at least one subprincipal error. Let $\mathcal{P}_{\mathrm{pp}}^{\mathcal C}$ denote the first term with $F=L(\mathcal C_{\lambda_2,\lambda_1},\widetilde\eta_{\lambda_2})$. Passing to the principal operators retains the output localization $\widetilde{P}_{\mu}$ described above between $\eta_{\lambda_1}$ and $\tilde{\eta}_{\lambda_2}$. Hence, using the first bound in \eqref{Phidifflowpairings} (which applies since $\mathcal C_{\lambda_2,\lambda_1}$ has frequency much lower than $\mu$) and also using \eqref{eta1esta2sep}, we get
\begin{equation*}
\begin{split}
\|\overline{u}^{x_0}_\mu P_\mu\mathcal{P}_{\mathrm{pp}}^{\mathcal C}\|_{L_{t,x}^1}&\lesssim\lambda_0^2
\|\widetilde P_\mu L(\eta_{\lambda_1},\widetilde\eta_{\lambda_2})\|_{L_{t,x}^2}
\|L(\mathcal{C}_{\lambda_2,\lambda_1},\overline{u}^{x_0}_\mu)\|_{L_{t,x}^2}
\\
&\lesssim_\delta\epsilon^4c_{\lambda_0}c_\mu\lambda_0^{-s_c}\mu^{-s_c}(1+\mu|x_0|)^{\frac12}.
\end{split}
\end{equation*}
The principal-principal contribution associated with $\mathcal E_{\lambda_2,\lambda_1}$ is handled similarly by pairing this coefficient with $u_\mu$ using \eqref{Phidifflowerror} and the $L_t^\infty L_x^2$ bound for $u_\mu$.
\\

Next, we estimate the last two terms in the decomposition above. We apply the subprincipal error estimate in Proposition \ref{basicopbounds} with $\sigma=0$ and $r=p$, for some sufficiently large $8<p<\infty$. Since $\lambda_1\approx_\delta\lambda_0$, the bootstrap in Section \ref{sec:bootstrap} gives
\begin{equation*}
\begin{split}
\|\nabla_x^2\eta_{<\delta\lambda_1}\|_{L_t^4\dot B^0_{p,1}}&\lesssim\epsilon\lambda_0^{\frac78-\frac2p},
\qquad
\|\eta_{\lambda_1}\|_{L_t^4\dot B^0_{\infty,1}}\lesssim\epsilon\lambda_0^{-\frac98},
\qquad
\|\eta_{\lambda_1}\|_{L_t^4\dot B^0_{p,1}}\lesssim\epsilon\lambda_0^{-\frac98-\frac2p}.
\end{split}
\end{equation*}
For the terms involving $\mathcal C_{\lambda_2,\lambda_1}$, we always retain the bilinear estimate for $(\mathcal C_{\lambda_2,\lambda_1},\widetilde\eta_{\lambda_2})$. If one of the two operators in \eqref{GLdecomposition} is estimated using the subprincipal error bound in Proposition \ref{basicopbounds}, the other is estimated by its linear bound and contributes a factor of $\lambda_0$. The additional factor $\lambda_0^{\frac2p}$ below comes from the additional Bernstein estimate
\begin{equation*}
\||D|^{2(1+\frac1p)}\eta_{<\delta\lambda_1}\|_{\dot B^0_{p,1}}
\lesssim\lambda_0^{\frac2p}\|\nabla_x^2\eta_{<\delta\lambda_1}\|_{\dot B^0_{p,1}}.
\end{equation*}
Thus, denoting the sum of these terms by $\mathcal R_{\mathrm{sub}}^{\mathcal C}$, we obtain
\begin{equation*}
\begin{split}
\|\overline u^{x_0}_\mu P_\mu\mathcal R_{\mathrm{sub}}^{\mathcal C}\|_{L_{t,x}^1}&\lesssim\lambda_0^{1+\frac2p}
\|\nabla_x^2\eta_{<\delta\lambda_1}\|_{L_t^4\dot B^0_{p,1}}
\|\eta_{\lambda_1}\|_{L_t^4\dot B^0_{\infty,1}}
\|L(\mathcal C_{\lambda_2,\lambda_1},\widetilde\eta_{\lambda_2})\|_{L_{t,x}^2}
\|u_\mu\|_{L_t^\infty L_x^2}
\\
&\lesssim\epsilon^5c_\mu c_{\lambda_0}\mu^{-s_c}\lambda_0^{-s_c}.
\end{split}
\end{equation*}
For the terms involving $\mathcal E_{\lambda_2,\lambda_1}$, no bilinear estimate is used. We instead use one $L_t^2L_x^\infty$ estimate for $\mathcal{E}_{\lambda_2,\lambda_1}$, two $L_t^4\dot{B}_{p,1}^0$-based Besov estimates (as above), and the $L_t^\infty L_x^2$ estimates for $\widetilde\eta_{\lambda_2}$ and $u_\mu$. This gives for the corresponding error $\mathcal{R}_{\mathrm{sub}}^{\mathcal{E}}$,
\begin{equation*}
\begin{split}
\|u^{x_0}_\mu P_\mu\mathcal R_{\mathrm{sub}}^{\mathcal E}\|_{L_{t,x}^1}&\lesssim\lambda_0^{1+\frac4p}
\|\nabla_x^2\eta_{<\delta\lambda_1}\|_{L_t^4\dot B^0_{p,1}}
\|\eta_{\lambda_1}\|_{L_t^4\dot B^0_{p,1}}
\|\mathcal E_{\lambda_2,\lambda_1}\|_{L_t^2L_x^\infty}
\|\widetilde\eta_{\lambda_2}\|_{L_t^\infty L_x^2}
\|u_\mu\|_{L_t^\infty L_x^2}
\\
&\lesssim_K\epsilon^6c_\mu c_{\lambda_0}\mu^{-s_c}\lambda_0^{-s_c}.
\end{split}
\end{equation*}
Thus, all the subprincipal contributions are perturbative after integration in $\lambda_1$ and $\lambda_0$. This concludes the estimate in the region $\mu\sim_c\Lambda$. 
\\

Now suppose that $\mu\not\sim_c\Lambda$. Assume for simplicity that $\lambda_1=\Lambda$ (the other case is similar by interchanging the pairings below). We again use \eqref{GLdecomposition}. The above estimates for the subprincipal terms apply verbatim in this case. For the principal-principal part we now have separation between $\overline{u}_{\mu}^{x_0}$ and $\eta_{\lambda_1}$, we can use the bounds
\begin{equation*}
\begin{split}
\|L(\eta_{\lambda_1},\overline{u}^{x_0}_\mu)\|_{L_{t,x}^2}\lesssim\epsilon^2c_\mu\lambda_0^{-2}\mu^{\frac34-s_c}(1+\mu|x_0|)^{\frac12},\hspace{5mm}\|L(\mathcal{C}_{\lambda_2,\lambda_1},\widetilde\eta_{\lambda_2})\|_{L_{t,x}^2}\lesssim\epsilon^2c_{\lambda_0}\lambda_0^{-\frac34-s_c}
\end{split}
\end{equation*}
which gives
\begin{equation*}
\begin{split}
\|\overline{u}^{x_0}_\mu P_\mu\mathcal{P}_{\mathrm{pp}}^{\mathcal C}\|_{L_{t,x}^1}\lesssim\epsilon^4c_\mu c_{\lambda_0}\mu^{\frac34-s_c}\lambda_0^{-\frac34-s_c}(1+\mu|x_0|)^{\frac12}.
\end{split}
\end{equation*}
The principal-principal contribution associated with $\mathcal E_{\lambda_2,\lambda_1}$ is handled the same way, by pairing $\mathcal E_{\lambda_2,\lambda_1}$ with $\widetilde\eta_{\lambda_2}$ using \eqref{Phidifflowerror} and the $L_t^\infty L_x^2$ bound for $\widetilde\eta_{\lambda_2}$, while $(\eta_{\lambda_1},\overline u^{x_0}_\mu)$ is estimated as above. 
\\

Combining the above estimates, and integrating first in $\lambda_1$ and then in $\lambda_0$, and using the slowly varying property of the frequency envelope in the second case, gives
\begin{equation*}
\int_{\lambda_0\gtrsim_M\mu}\int_{\delta\lambda_0}^{\lambda_0}\|\overline{u}^{x_0}_{\mu}P_\mu G^{l,\lambda_0}_{<\delta\lambda_1}(\eta_{\lambda_1}L^{l,\lambda_0}_{<\delta\lambda_1}(\Phi^{l,0}_{<\delta\lambda_1}-\varphi_{\lambda_0}))\|_{L_{t,x}^1}\frac{d\lambda_1}{\lambda_1}\frac{d\lambda_0}{\lambda_0}
\lesssim\epsilon^4c_\mu^2\mu^{-2s_c}(1+\mu|x_0|)^{\frac12},
\end{equation*}
which proves \eqref{replacebyphi}. Our next goal is to show
\begin{equation}\label{paratruncation}
G^{l,\lambda_0}_{<\delta\lambda_1}(\eta_{\lambda_1}L^{l,\lambda_0}_{<\delta\lambda_1}\varphi_{\lambda_0})=P_{<M\lambda_0}T_G(\eta_{\lambda_1}T_L\varphi_{\lambda_0})+\mathcal{R}.
\end{equation}
For this, we set
\begin{equation*}
f_{\lambda_2}\coloneq S_{\lambda_2}(\eta_{\lambda_1}L^{l,\lambda_0}_{<\delta\lambda_1}\varphi_{\lambda_0}),
\qquad S_{\lambda_2}=\lambda_2\frac{d}{d\lambda_2}P_{<\lambda_2}.
\end{equation*}
The function inside $S_{\lambda_2}$ is localized to frequencies $\leq C_\delta\lambda_1$. Hence, using the linearization formula \eqref{Glin}, integrating from $\lambda_2=0$ to $\lambda_2=C_{\delta}\lambda_1$, and, in the range $\lambda_2\leq\lambda_1$, integrating in $\lambda_3$ from $\lambda_2$ to $\lambda_1$, we obtain
\begin{equation}\label{G2trunc}
\begin{split}
G^{l,\lambda_0}_{<\delta\lambda_1}(\eta_{\lambda_1}L^{l,\lambda_0}_{<\delta\lambda_1}\varphi_{\lambda_0})&=P_{<M\lambda_0}T_G(\eta_{\lambda_1}L^{l,\lambda_0}_{<\delta\lambda_1}\varphi_{\lambda_0})
\\
&-P_{<M\lambda_0}\int_{0}^{\lambda_1}\int_{\lambda_2}^{\lambda_1}\left[G_{<\delta\lambda_3}(\eta_{\delta\lambda_3}L_{<\delta\lambda_3}f_{\lambda_2})+\nabla\cdot (\eta_{\delta\lambda_3}M_{<\delta\lambda_3}f_{\lambda_2})\right]\frac{d\lambda_3}{\lambda_3}\frac{d\lambda_2}{\lambda_2}
\\
&+P_{<M\lambda_0}\int_{0}^{C_\delta\lambda_1}P_{\geq M\lambda_2}G_{<\delta\lambda_2}f_{\lambda_2}\frac{d\lambda_2}{\lambda_2}+\mathcal{R}_{cmp}.
\end{split}
\end{equation}
Here $\mathcal R_{cmp}$ corresponds to a contribution of the same form as the double integral in the second line above, but for the comparable range $\lambda_1<\lambda_2\leq C_{\delta}\lambda_1$. Proposition \ref{Bernsteinelliptic}, Corollary \ref{Bestcor}, and the estimates in Proposition \ref{phiu-bootstrap} show that the high-frequency tail in the last line is perturbative in $L_t^1L_x^2$ (and placing $u_{\mu}^{x_0}$ in $L^{\infty}_tL_x^2$). The term $\mathcal R_{cmp}$ is not in $L_t^1L_x^2$, but its contribution is controlled after pairing
with $\overline u_\mu^{x_0}$ by the principal/subprincipal decomposition
and frequency separation argument used above for \eqref{replacebyphi}. Here, 
$\lambda_1\leq\lambda_3\leq\lambda_2\leq C_\delta\lambda_1$,
so the additional frequency integrations are over comparable ranges. 
\\

Now, moving to the double integral, where $\lambda_2$ may not be comparable to $\lambda_1$. We will show in detail how to handle the more difficult first term; as the divergence term can be handled by identical reasoning (in fact, it is substantially easier because there is no non-local $G_{<\delta\lambda_3}$ operator present on the output). 
\\

Our first step is to observe that Proposition \ref{Bernsteinelliptic} and Proposition \ref{basicopbounds} allow us to replace the operators in the double integral by their output truncated versions $G^{l,\lambda_3}_{<\delta\lambda_3}$ and $L^{l,\lambda_3}_{<\delta\lambda_3}$. More precisely, after fixing an integer $N\geq3$, and choosing $M\gg_N 1$, the resulting error satisfies
\begin{equation*}
\int_{\lambda_0\gtrsim_M\mu}\int_{\delta\lambda_0}^{\lambda_0}\|\mathcal{R}_{\mu}\|_{L_t^1L_x^2}\frac{d\lambda_1}{\lambda_1}\frac{d\lambda_0}{\lambda_0}\lesssim_N \epsilon^Nc_{\mu}\mu^{-s_c}.
\end{equation*}
Consequently,
\begin{equation*}
\int_{\lambda_0\gtrsim_M\mu}\int_{\delta\lambda_0}^{\lambda_0}\|\overline{u}^{x_0}_{\mu}\mathcal{R}_{\mu}\|_{L_{t,x}^1}\frac{d\lambda_1}{\lambda_1}\frac{d\lambda_0}{\lambda_0}\lesssim \epsilon^4c_{\mu}^2\mu^{-2s_c}.
\end{equation*}
It therefore remains to estimate
\begin{equation*}
G^{l,\lambda_3}_{<\delta\lambda_3}(\eta_{\delta\lambda_3}L^{l,\lambda_3}_{<\delta\lambda_3}f_{\lambda_2}).
\end{equation*}
We begin by expanding the inner factor. The linearization formula for $L$ \eqref{Blin} and the fundamental theorem of calculus gives
\begin{equation*}
\begin{split}
L^{l,\lambda_3}_{<\delta\lambda_3}f_{\lambda_2}&=L^{l,\lambda_2}_{<\delta\lambda_2}f_{\lambda_2}+\int_{\lambda_2}^{\lambda_3}F_{\lambda_2,\lambda_4}\frac{d\lambda_4}{\lambda_4}+\mathcal R,
\end{split}
\end{equation*}
where
\begin{equation*}
F_{\lambda_2,\lambda_4}\coloneq P_{<M\lambda_4}\left[-L_{<\delta\lambda_4}(\eta_{\delta\lambda_4}L_{<\delta\lambda_4}f_{\lambda_2})+\eta_{\delta\lambda_4}L^2_{<\delta\lambda_4}f_{\lambda_2}\right].
\end{equation*}
This identity is obtained first with the fixed output projection $P_{<M\lambda_3}$. Replacing the projection at the lower endpoint by $P_{<M\lambda_2}$ and the projection in the integrand by $P_{<M\lambda_4}$ produces only high-frequency tails, since $\delta\lambda_2\leq\lambda_2$ and $\max(\lambda_2,\delta\lambda_4)\leq\lambda_4$. The resulting error $\mathcal R$ satisfies the preceding perturbative bound by Proposition \ref{Bernsteinelliptic}. Next, by Bernstein's inequality, Proposition \ref{basicopbounds}, and \eqref{principaldecomposition}, the contribution of $F_{\lambda_2,\lambda_4}$ can be estimated, for some $8<p<\infty$, by
\begin{equation}\label{delicatef2}
\begin{split}
\|G^{l,\lambda_3}_{<\delta\lambda_3}(\eta_{\delta\lambda_3}F_{\lambda_2,\lambda_4})\|_{L_t^1L_x^2}&\lesssim \lambda_3\|\eta_{\delta\lambda_3}\|_{L_t^4L_x^\infty}\|F_{\lambda_2,\lambda_4}\|_{L_t^{\frac43}L_x^2}
\\
&\lesssim \epsilon\lambda_3^{-\frac18}\lambda_4\lambda_2\|L(\eta_{\delta\lambda_4},f_{\lambda_2})\|_{L_t^{\frac43}L_x^2}
\\
&+\lambda_3\lambda_4^{\frac{15}{8}}\|\eta_{\delta\lambda_4}\|_{L_t^4L_x^\infty}\|\eta_{\delta\lambda_3}\|_{L_t^4L_x^\infty}\||D|^{\frac{11}{8}-2(\frac18-\frac1p)}\eta_{<\delta\lambda_4}\|_{L_t^4L_x^p}\|f_{\lambda_2}\|_{L_t^4L_x^2}.
\end{split}
\end{equation}
The bilinear form defining $L(\eta_{\delta\lambda_4},f_{\lambda_2})$ only has rapidly decaying kernel at scale $\lambda_2$ (which can be very small). Therefore, for the first term on the right-hand side, we keep $f_{\lambda_2}$ intact rather than pairing its internal factors with other functions. Using Proposition \ref{paraexp}, and the fact that $\lambda_2\leq\lambda_0\approx\lambda_1$, we can estimate
\begin{equation*}
\lambda_2^{\frac{1}{2}}\|f_{\lambda_2}\|_{L_{t,x}^2}\lesssim \lambda_1\|S_{\lambda_2}|D|^{\frac{1}{2}}L(\eta_{\lambda_1},\varphi_{\lambda_0})\|_{L_{t,x}^2}+\lambda_2^{\frac{1}{2}}\|\eta_{\lambda_1}\|_{L_t^{\infty}L_x^2}\|\varphi_{\lambda_0}\|_{L_t^4L_x^{\infty}}\|\nabla^2\eta_{<\delta\lambda_1}\|_{L_t^4\dot{B}_{\infty,1}^0}   
\end{equation*}
where the $L$ on the right-hand side here has rapidly decaying kernel at scale $\approx\lambda_1$ (which is allowed). The second term on the right-hand side accounts for the zeroth-order contribution in Proposition \ref{paraexp}; the higher-order remainder in \eqref{paraexperror} is absorbed into the same bound by taking the expansion order at least two, placing one factor of $\nabla^2\eta$ in $L_t^4\dot B^0_{\infty,1}$ and the remaining factors in $L_t^\infty\dot B^0_{\infty,1}$. Lemma \ref{balancedrestrictedoutput}, applicable
since $\lambda_2\leq\lambda_0$, alongside the
$L_t^4\dot{B}_{\infty,1}^0$ bounds, give
\begin{equation}\label{fbalancedbound}
\lambda_2^{\frac12}\|f_{\lambda_2}\|_{L_{t,x}^2}
\lesssim\epsilon^2c_{\lambda_1}\lambda_0^{-\frac14}\lambda_1^{-s_c}.
\end{equation}
Hence, the bilinear term in \eqref{delicatef2} can be bounded by H\"older's inequality and the $L_t^4L_x^{\infty}$ bound for $\eta_{\delta\lambda_4}$ to give
\begin{equation*}
\begin{split}
\|L(\eta_{\delta\lambda_4},f_{\lambda_2})\|_{L_t^{\frac43}L_x^2}\lesssim\|\eta_{\delta\lambda_4}\|_{L_t^4L_x^\infty}\|f_{\lambda_2}\|_{L_{t,x}^2}\lesssim\epsilon^3c_{\lambda_1}\lambda_4^{-\frac98}\lambda_2^{-\frac12}\lambda_0^{-\frac14}\lambda_1^{-s_c}.
\end{split}
\end{equation*}
Consequently, the first term in \eqref{delicatef2} is bounded by
\begin{equation*}
\epsilon^4c_{\lambda_1}\lambda_3^{-\frac18}\lambda_4^{-\frac18}\lambda_2^{\frac12}\lambda_0^{-\frac14}\lambda_1^{-s_c}.
\end{equation*}
For the second term, set $\theta=\frac2p\in(0,\frac14)$. Bernstein gives
\begin{equation*}
\begin{split}
\|f_{\lambda_2}\|_{L_t^4L_x^2}&\lesssim\lambda_2^\theta\|\eta_{\lambda_1}\|_{L_t^4L_x^p}\|L^{l,\lambda_0}_{<\delta\lambda_1}\varphi_{\lambda_0}\|_{L_t^\infty L_x^2}\lesssim\epsilon^2c_{\lambda_0}\lambda_2^\theta\lambda_1^{-\frac98-\theta}\lambda_0^{\frac12-s_c}.
\end{split}
\end{equation*}
Consequently, the second term above is bounded by
\begin{equation*}
\epsilon^5c_{\lambda_0}\lambda_3^{-\frac18}\lambda_4^{\frac34}\lambda_2^\theta\lambda_1^{-\frac58-\theta-s_c}.
\end{equation*}
We then obtain,
\begin{equation*}
\int_0^{\lambda_1}\int_{\lambda_2}^{\lambda_1}\int_{\lambda_2}^{\lambda_3}\lambda_3^{-\frac18}\lambda_4^{-\frac18}\lambda_2^{\frac12}\frac{d\lambda_4}{\lambda_4}\frac{d\lambda_3}{\lambda_3}\frac{d\lambda_2}{\lambda_2}\lesssim\lambda_1^{\frac14}
\end{equation*}
and
\begin{equation*}
\int_{\lambda_2}^{\lambda_3}\lambda_4^{\frac34}\frac{d\lambda_4}{\lambda_4}\lesssim\lambda_3^{\frac34},\qquad
\int_{\lambda_2}^{\lambda_1}\lambda_3^{\frac58}\frac{d\lambda_3}{\lambda_3}\lesssim\lambda_1^{\frac58},\qquad
\int_0^{\lambda_1}\lambda_2^\theta\frac{d\lambda_2}{\lambda_2}\lesssim\lambda_1^\theta.
\end{equation*}
Thus the first bound contributes $\epsilon^4c_{\lambda_1}\lambda_0^{-\frac14}\lambda_1^{\frac14-s_c}$ and the second contributes $\epsilon^5c_{\lambda_0}\lambda_1^{-s_c}$. Since $\lambda_1\approx_\delta\lambda_0$, both have the required summability. Defining
\begin{equation*}
\mathcal{R}_F\coloneq\int_{0}^{\lambda_1}\int_{\lambda_2}^{\lambda_1}\int_{\lambda_2}^{\lambda_3}G^{l,\lambda_3}_{<\delta\lambda_3}(\eta_{\delta\lambda_3}F_{\lambda_2,\lambda_4})\frac{d\lambda_4}{\lambda_4}\frac{d\lambda_3}{\lambda_3}\frac{d\lambda_2}{\lambda_2},
\end{equation*}
we obtain
\begin{equation*}
\int_{\lambda_0\gtrsim_M\mu}\int_{\delta\lambda_0}^{\lambda_0}\|\overline{u}_{\mu}^{x_0}(\mathcal{R}_F)_{\mu}\|_{L_{t,x}^1}\frac{d\lambda_1}{\lambda_1}\frac{d\lambda_0}{\lambda_0}\lesssim \epsilon^5 c_{\mu}^2\mu^{-2s_c}(1+\mu|x_0|)^{\frac12}.
\end{equation*}
Therefore, $F_{\lambda_2,\lambda_4}$ contributes a perturbative error.
\\

\noindent It remains to understand the contribution of
$G^{l,\lambda_3}_{<\delta\lambda_3}(\eta_{\delta\lambda_3}
L^{l,\lambda_2}_{<\delta\lambda_2}f_{\lambda_2})$. First, applying the fundamental theorem of calculus to the
Littlewood--Paley cutoff in the argument of $G$, splitting the
resulting $\lambda_4$-integration at $\lambda_4=\lambda_3$, and adding
and subtracting $G^{l,\lambda_4}_{<\delta\lambda_4}$ in the range
$\lambda_4\leq\lambda_3$, we obtain the following identity (where the
range $\lambda_4>\lambda_3$ is included in $\mathcal R$ as justified
below):
\begin{equation}\label{GLdecompbal}
\begin{split}
&P_{<M\lambda_0}\int_0^{\lambda_1}\int_{\lambda_2}^{\lambda_1}
G^{l,\lambda_3}_{<\delta\lambda_3}
\bigl(\eta_{\delta\lambda_3}
L^{l,\lambda_2}_{<\delta\lambda_2}f_{\lambda_2}\bigr)
\frac{d\lambda_3}{\lambda_3}\frac{d\lambda_2}{\lambda_2}
\\
&\quad=P_{<M\lambda_0}\int_0^{\lambda_1}\int_{\lambda_2}^{\lambda_1}
\int_0^{\lambda_3}
G^{l,\lambda_4}_{<\delta\lambda_4}
S_{\lambda_4}\bigl(\eta_{\delta\lambda_3}
L^{l,\lambda_2}_{<\delta\lambda_2}f_{\lambda_2}\bigr)
\frac{d\lambda_4}{\lambda_4}
\frac{d\lambda_3}{\lambda_3}\frac{d\lambda_2}{\lambda_2}
\\
&\quad+P_{<M\lambda_0}\int_0^{\lambda_1}\int_{\lambda_2}^{\lambda_1}
\int_0^{\lambda_3}
\bigl(G^{l,\lambda_3}_{<\delta\lambda_3}
-G^{l,\lambda_4}_{<\delta\lambda_4}\bigr)
S_{\lambda_4}\bigl(\eta_{\delta\lambda_3}
L^{l,\lambda_2}_{<\delta\lambda_2}f_{\lambda_2}\bigr)
\frac{d\lambda_4}{\lambda_4}
\frac{d\lambda_3}{\lambda_3}\frac{d\lambda_2}{\lambda_2}
+\mathcal R.
\end{split}
\end{equation}
We begin by briefly discussing the contribution of $\mathcal{R}$ which comes from the contribution of the output range $\lambda_4>\lambda_3$. The portion
$\lambda_4\gg_M\lambda_3$ contributes high-frequency tails and is straightforward to estimate (directly in $L_t^1L_x^2$). For the region
\begin{equation*}
\lambda_3<\lambda_4\lesssim_M\lambda_3,
\end{equation*}
after applying Proposition \ref{paraexp}, the subprincipal part of
$G^{l,\lambda_3}_{<\delta\lambda_3}$ is perturbative (also in $L_t^1L_x^2$) by the
corresponding remainder estimate and \eqref{fbalancedbound} and placing $\eta_{\delta\lambda_3}$ in $L_t^4L_x^{\infty}$. For its
principal part, the projection $P_\mu$ forces
$\mu\approx_{M,\delta}\lambda_4$, while frequency support forces
$\lambda_2\approx_{M,\delta}\lambda_4$. Since $M\delta\ll1$, the
factors $\overline u_\mu^{x_0}$ and $\eta_{\delta\lambda_3}$ are
frequency separated. We may therefore apply the unbalanced bilinear
estimate to $(\overline u_\mu^{x_0},\eta_{\delta\lambda_3})$, together
with \eqref{fbalancedbound}. Thus the entire range
$\lambda_4>\lambda_3$ is perturbative and may be included in
$\mathcal R$.
\\

Using a similar argument, the second triple integral in \eqref{GLdecompbal} is also
perturbative. We sketch the relevant details. First, modulo a high-frequency tail, we may replace
$G_{<\delta\lambda_4}^{l,\lambda_4}$ by
$G_{<\delta\lambda_4}^{l,\lambda_3}$. More precisely, when applied to
\begin{equation*}
S_{\lambda_4}\bigl(
\eta_{\delta\lambda_3}
L^{l,\lambda_2}_{<\delta\lambda_2}f_{\lambda_2}
\bigr),
\end{equation*}
the difference between these two output truncations is a 
high-frequency tail that can be dispensed with using Proposition \ref{Bernsteinelliptic}. Keeping the output projection
$P_{<M\lambda_3}$ fixed, we then apply the fundamental theorem of
calculus and \eqref{Glin} to
\begin{equation*}
G^{l,\lambda_3}_{<\delta\lambda_3}
-G^{l,\lambda_3}_{<\delta\lambda_4}
\end{equation*}
by integrating in a parameter
$\lambda_4\leq\lambda_5\leq\lambda_3$. By the linearization formula
\eqref{Glin}, every resulting term contains a factor
$\eta_{\delta\lambda_5}$. No additional high-frequency tail is
produced in this $\lambda_5$-integration.

We place
$\eta_{\delta\lambda_3}
L_{<\delta\lambda_2}^{l,\lambda_2}f_{\lambda_2}$ in
$L_t^{\frac43}L_x^2$, with $f_{\lambda_2}$ in $L_{t,x}^2$ and
$\eta_{\delta\lambda_3}$ in $L_t^4L_x^\infty$, and place
$\eta_{\delta\lambda_5}$ in $L_t^4L_x^\infty$ and
$u_\mu^{x_0}$ in $L_t^\infty L_x^2$. Proposition
\ref{basicopbounds} bounds the resulting integrands by
\begin{equation*}
\begin{split}
&\lambda_5\lambda_4\lambda_2
\|u_\mu^{x_0}\|_{L_t^\infty L_x^2}
\|\eta_{\delta\lambda_5}\|_{L_t^4L_x^\infty}
\|\eta_{\delta\lambda_3}\|_{L_t^4L_x^\infty}
\|f_{\lambda_2}\|_{L_{t,x}^2}
\\
&\qquad\lesssim
\epsilon^5c_\mu c_{\lambda_1}\mu^{-s_c}
\lambda_5^{-\frac18}\lambda_4\lambda_3^{-\frac98}
\lambda_2^{\frac12}\lambda_0^{-\frac14}\lambda_1^{-s_c}
(1+\mu|x_0|)^{\frac12}.
\end{split}
\end{equation*}
The successive frequency integrations are summable, since
\begin{equation*}
\int_{\lambda_4}^{\lambda_3}
\lambda_5^{-\frac18}\frac{d\lambda_5}{\lambda_5}
\lesssim\lambda_4^{-\frac18},
\qquad
\int_0^{\lambda_3}
\lambda_4^{\frac78}\frac{d\lambda_4}{\lambda_4}
\lesssim\lambda_3^{\frac78},
\end{equation*}
and
\begin{equation*}
\int_{\lambda_2}^{\lambda_1}
\lambda_3^{-\frac14}\frac{d\lambda_3}{\lambda_3}
\lesssim\lambda_2^{-\frac14},
\qquad
\int_0^{\lambda_1}
\lambda_2^{\frac14}\frac{d\lambda_2}{\lambda_2}
\lesssim\lambda_1^{\frac14}.
\end{equation*}
Thus this term is bounded by
\begin{equation*}
\epsilon^5c_\mu c_{\lambda_1}\mu^{-s_c}
\lambda_0^{-\frac14}\lambda_1^{\frac14-s_c}
(1+\mu|x_0|)^{\frac12},
\end{equation*}
and is perturbative after integrating in $\lambda_0$ and $\lambda_1$. In summary, we have
\begin{equation*}
\begin{split}
&P_{<M\lambda_0}\int_0^{\lambda_1}\int_{\lambda_2}^{\lambda_1}
\int_0^{\lambda_3}
\bigl(G^{l,\lambda_3}_{<\delta\lambda_3}
-G^{l,\lambda_4}_{<\delta\lambda_4}\bigr)
S_{\lambda_4}\bigl(\eta_{\delta\lambda_3}
L^{l,\lambda_2}_{<\delta\lambda_2}f_{\lambda_2}\bigr)
\frac{d\lambda_4}{\lambda_4}
\frac{d\lambda_3}{\lambda_3}\frac{d\lambda_2}{\lambda_2}
=\mathcal R.
\end{split}
\end{equation*}
We now estimate the main term in the second line of \eqref{GLdecompbal}. Put
$\nu=\delta\lambda_3$. After extracting the principal part of the
outer operator, the projection $P_\mu$ forces
$\lambda_4\approx\mu$, modulo high-frequency tails. If
$\nu$ and $\mu$ are separated, we can pair, as above,
$(\overline u_\mu^{x_0},\eta_{\delta\lambda_3})$ and use the
unbalanced estimate
\begin{equation*}
\|L(\eta_\nu,\overline u_\mu^{x_0})\|_{L_{t,x}^2}
\lesssim\epsilon^2c_\mu(1+\mu|x_0|)^{\frac12}
\begin{cases}
\mu^{-s_c}\nu^{-\frac54},&\nu\ll\mu,\\
\mu^{\frac34-s_c}\nu^{-2},&\mu\ll\nu.
\end{cases}
\end{equation*}
In the first case frequency support forces $\lambda_2\approx_{M,\delta}\mu$,
while in the second it forces $\lambda_2\approx_{M,\delta}\nu$. Together with
\eqref{fbalancedbound}, this shows that the contribution is
perturbative. This covers also, in particular, the case
$\mu\approx\lambda_0$, since
\begin{equation*}
\nu=\delta\lambda_3\leq\delta\lambda_1
\leq\delta\lambda_0\ll\mu.
\end{equation*}

It remains to consider $\mu\ll\lambda_0$ in the regime
$\nu\approx\mu$. In this region, we also have
$\lambda_2\lesssim_\delta\mu$, by frequency support. Similarly to \eqref{GLdecomposition}, we first peel off the principal part by writing
\begin{equation*}
\begin{split}
&G(\eta_\nu Lf)
-G_{\mathrm{pr}}(\eta_\nu L_{\mathrm{pr}}f^{\mathrm{pr}})
\\
&\quad=
(G-G_{\mathrm{pr}})(\eta_\nu Lf)
+G_{\mathrm{pr}}\bigl(\eta_\nu(L-L_{\mathrm{pr}})f\bigr)
+G_{\mathrm{pr}}\bigl(
\eta_\nu L_{\mathrm{pr}}(f-f^{\mathrm{pr}})
\bigr),
\end{split}
\end{equation*}
where the subscripts, superscripts and the outer
$S_{\lambda_4}$ localization have been temporarily suppressed. Here, we define $f_{\lambda_2}^{\mathrm{pr}}
\coloneq
S_{\lambda_2}\Bigl(
\eta_{\lambda_1}
\operatorname{Op}(\mathfrak l^l_{<\delta\lambda_1})
\varphi_{\lambda_0}
\Bigr)$. The first two terms are perturbative by Proposition
\ref{basicopbounds} and \eqref{fbalancedbound}, using two
$L_t^4\dot B_{\infty,1}^0$ bounds and placing $u_\mu^{x_0}$ in
$L_t^\infty L_x^2$. To handle the remaining term, we write
\begin{equation*}
\mathcal R_{L,0}
\coloneq
L_{<\delta\lambda_1}^{l,\lambda_0}
-\operatorname{Op}(\mathfrak l^l_{<\delta\lambda_1})
\end{equation*}
so that
\begin{equation*}
f_{\lambda_2}-f_{\lambda_2}^{\mathrm{pr}}
=
S_{\lambda_2}\bigl(
\eta_{\lambda_1}\mathcal R_{L,0}\varphi_{\lambda_0}
\bigr).
\end{equation*}
Since $\lambda_1\leq\lambda_0$, we can use Proposition \ref{paraexp} and Proposition \ref{phiu-bootstrap} to estimate
\begin{equation*}
\begin{split}
\|\mathcal R_{L,0}\varphi_{\lambda_0}\|_{L_t^2L_x^2}\lesssim\epsilon^2\lambda_0^{\frac{1}{4}-s_c}c_{\lambda_0} .
\end{split}
\end{equation*}
Consequently, by Bernstein and the $L_t^4L_x^{\infty}$ estimate for $\eta_{\lambda_1}$, we have
\begin{equation*}
\begin{split}
\|\overline{u}_{\mu}^{x_0}P_{\mu}G_{\mathrm{pr}}S_{\lambda_4}\left(\eta_{\nu}L_{\mathrm{pr}}S_{\lambda_2}\bigl(
\eta_{\lambda_1}\mathcal R_{L,0}\varphi_{\lambda_0}
\bigr)\right)\|_{L_{t,x}^1}
&\lesssim
\epsilon^3c_{\lambda_0}\lambda_2\lambda_4\lambda_0^{\frac{1}{4}-s_c}\lambda_1^{-\frac{9}{8}}\|u_{\mu}^{x_0}\|_{L_t^{\infty}L_x^2}\|\eta_{\nu}\|_{L_t^4L_x^{\infty}}
\\
&\lesssim \epsilon^5 c_{\lambda_0}c_{\mu}\lambda_2\lambda_4\lambda_0^{\frac{1}{4}-s_c}\lambda_1^{-\frac{9}{8}}\mu^{-s_c}\nu^{-\frac{9}{8}}.
\end{split}
\end{equation*}
After integration in $\lambda_0,\lambda_1,\lambda_2,\lambda_3,\lambda_4$, this is bounded by
\begin{equation*}
\epsilon^5c_\mu^2\mu^{-2s_c}\lesssim \epsilon^5c_\mu^2\mu^{-2s_c}(1+|x_0|\mu)^{\frac12}.
\end{equation*}

It remains to handle the contribution of the principal term, where we now reintroduce all the frequency parameters,
\begin{equation*}
\operatorname{Op}(\mathfrak g^l_{<\delta\lambda_4})S_{\lambda_4}
\bigl(
\eta_\nu
\operatorname{Op}(\mathfrak l^l_{<\delta\lambda_2})
f_{\lambda_2}^{\mathrm{pr}}
\bigr).    
\end{equation*}
Since $\mu\approx\nu$ (by assumption), the input factors in $f_{\lambda_2}$, namely
$\eta_{\lambda_1}$ and $\varphi_{\lambda_0}$, are frequency separated
from $\eta_\nu$ and $\overline u_\mu^{x_0}$, respectively. Ideally,
we would pair these factors and apply two unbalanced bilinear
estimates. In the low-output range $\lambda_2\ll\mu$, however, this strategy
encounters a serious obstruction; namely the kernel for
$\operatorname{Op}(\mathfrak l^l_{<\delta\lambda_2})
S_{\lambda_2}$ has limited decay. This issue occurs even in the
flat case, where this operator is simply $|D|S_{\lambda_2}$. More
precisely, its kernel is rapidly decreasing only at the low output
scale $\lambda_2^{-1}$, but
\begin{equation*}
\int_0^{c\mu}|D|S_{\lambda_2}\frac{d\lambda_2}{\lambda_2}
=|D|P_{<c\mu}.
\end{equation*}
The kernel of $|D|P_{<c\mu}$ is not integrable against the weight
$\langle\mu y\rangle$ which is a factor that can appear in certain cases when using two pairs of the ordinary unbalanced bilinear estimates at once (since each pair provides a translation factor of size $\langle\mu y\rangle^{\frac{1}{2}}$). To handle this and similar
endpoint cases encountered later, we use the translation-averaged
bilinear estimate \eqref{uab-bi-unbal-trans}, together with the
machinery developed in Subsection \ref{subsec:multilinearnorms}. The
main point is that this replaces the borderline weighted $L_y^1$
kernel estimate by weighted $L_y^2$-averaged control over
translations. To carry out this argument precisely, let
$k_\rho(q)$ be nonnegative kernel majorants satisfying
\begin{equation*}
k_\rho(q)\lesssim_N
\rho^2\langle\rho q\rangle^{-N},
\qquad
\rho\in\{\mu,\lambda_2,\lambda_0\}.
\end{equation*}
The three order-one principal operators contribute the factors
$\mu$, $\lambda_2$, and $\lambda_0$, respectively. Writing
\begin{equation*}
F^a\coloneq\tau_aF,\qquad \tau_aF(x)=F(x-a),
\end{equation*}
we expand their kernels to obtain
\begin{equation*}
\begin{split}
&I_{\lambda_2}
\coloneq
\|\overline u_\mu^{x_0}P_\mu
\operatorname{Op}(\mathfrak g^l_{<\delta\lambda_4})S_{\lambda_4}
\bigl(
\eta_\nu
\operatorname{Op}(\mathfrak l^l_{<\delta\lambda_2})
f_{\lambda_2}^{\mathrm{pr}}
\bigr)\|_{L_{t,x}^1}
\\
&\lesssim
\mu\lambda_2\lambda_0
\iiint
k_\mu(z)k_{\lambda_2}(y)k_{\lambda_0}(w)
\bigl\|
\overline u_\mu^{x_0}(x)\eta_\nu^{z}(x)
\eta_{\lambda_1}^{z+y}(x)
\varphi_{\lambda_0}^{z+y+w}(x)
\bigr\|_{L_{t,x}^1}
\,dz\,dy\,dw.
\end{split}
\end{equation*}
Changing variables, we obtain
\begin{equation*}
\begin{split}
I_{\lambda_2}\lesssim
\mu\lambda_2\lambda_0
\iint k_\mu(z)k_{\lambda_0}(w)M_{\lambda_2,1}^{\sigma}
\bigl(
\eta_\nu\tau_{-z}\overline u_\mu^{x_0},
\eta_{\lambda_1}\tau_w\varphi_{\lambda_0}
\bigr)\,dz\,dw
\end{split}
\end{equation*}
where we take $1<\sigma<2$. Here we used $k_{\lambda_2}(y)
\lesssim
\lambda_2^2\langle\lambda_2y\rangle^{-\sigma}$. Put $\theta=\sigma/2$. By \eqref{M-holder-product},
\begin{equation*}
\begin{split}
M_{\lambda_2,1}^{\sigma}
\bigl(
\eta_\nu\tau_{-z}\overline u_\mu^{x_0},
\eta_{\lambda_1}\tau_w\varphi_{\lambda_0}
\bigr)\leq
M_{\lambda_2,2}^{\theta}(\eta_\nu,\eta_{\lambda_1})
M_{\lambda_2,2}^{\theta}
(\tau_{-z}\overline u_\mu^{x_0},
\tau_w\varphi_{\lambda_0}).
\end{split}
\end{equation*}
Moreover, directly from the definition of $M_{\lambda_2,2}^\theta$,
\begin{equation*}
\begin{split}
M_{\lambda_2,2}^{\theta}
(\tau_{-z}\overline u_\mu^{x_0},
\tau_w\varphi_{\lambda_0})^2&=
\lambda_2^2\int_{\mathbb R^2}
\frac{
\|\overline u_\mu^{x_0}
\tau_{y+z+w}\varphi_{\lambda_0}\|_{L_{t,x}^2}^2}
{\langle\lambda_2y\rangle^{2\theta}}\,dy
\\
&\lesssim
\langle\lambda_2(z+w)\rangle^{2\theta}
M_{\lambda_2,2}^{\theta}
(\overline u_\mu^{x_0},\varphi_{\lambda_0})^2.
\end{split}
\end{equation*}
Since $\lambda_2\lesssim\mu,\lambda_0$, we have
\begin{equation*}
\iint k_\mu(z)k_{\lambda_0}(w)
\langle\lambda_2(z+w)\rangle^\theta\,dz\,dw
\lesssim1.
\end{equation*}
Therefore,
\begin{equation}\label{Ilambda2}
I_{\lambda_2}
\lesssim
\mu\lambda_2\lambda_0
M_{\lambda_2,2}^{\frac\sigma2}(\eta_\nu,\eta_{\lambda_1})
M_{\lambda_2,2}^{\frac\sigma2}
(\overline u_\mu^{x_0},\varphi_{\lambda_0}).
\end{equation}
Applying \eqref{M-scale-eq} to both
$M_{\lambda_2,2}^{\sigma/2}$ factors gives
\begin{equation}\label{GsecondMholder}
\begin{split}
I_{\lambda_2}\lesssim_\sigma
\mu\lambda_2\lambda_0
\left(\frac{\lambda_2}{\mu}\right)^{2-\sigma}
M_{\mu,2}^{\frac\sigma2}(\eta_\nu,\eta_{\lambda_1})
M_{\mu,2}^{\frac\sigma2}
(\overline u_\mu^{x_0},\varphi_{\lambda_0}).
\end{split}
\end{equation}
By Lemma \ref{annularcontrol}, together with the
translation-averaged unbalanced bounds in Section
\ref{sec:bootstrap} and Proposition \ref{phiu-bootstrap}, we obtain
\begin{equation}\label{Gsecondtranslatedpairs}
\begin{split}
M_{\mu,2}^{\frac\sigma2}(\eta_\nu,\eta_{\lambda_1})\lesssim\epsilon^2\nu^{-\frac54}\lambda_1^{-2},\hspace{3mm}M_{\mu,2}^{\frac\sigma2}
(\overline u_\mu^{x_0},\varphi_{\lambda_0})\lesssim\epsilon^2c_\mu c_{\lambda_0}
\mu^{\frac34-s_c}\lambda_0^{-\frac12-s_c}(1+\mu|x_0|)^{\frac12}.
\end{split}
\end{equation}
Since $\nu\approx\mu$ and
$\lambda_1\approx \lambda_0$, substitution into
\eqref{GsecondMholder} gives, before integration in $\lambda_2$,
\begin{equation*}
I_{\lambda_2}
\lesssim
\epsilon^4c_\mu c_{\lambda_0}
\mu^{\frac12-s_c}\lambda_0^{-\frac32-s_c}
\lambda_2
\left(\frac{\lambda_2}{\mu}\right)^{2-\sigma}
(1+\mu|x_0|)^{\frac12}.
\end{equation*}
Consequently, integration in $\lambda_2$ gives the contribution
\begin{equation}\label{Gsecondtransitionbound}
\epsilon^4c_\mu c_{\lambda_0}
\mu^{\frac32-s_c}\lambda_0^{-\frac32-s_c}(1+\mu|x_0|)^{\frac12}.
\end{equation}
This in turn integrates over the region $\lambda_0\gg\mu$ to the
usual perturbative bound
\begin{equation*}
\epsilon^4c_\mu^2\mu^{-2s_c}(1+\mu|x_0|)^{\frac12}.
\end{equation*}
Combining the above cases shows that both triple integrals on the
right-hand side of \eqref{GLdecompbal} are perturbative. Hence the
first term inside the brackets in the double integral in
\eqref{G2trunc} is perturbative. The divergence term in
\eqref{G2trunc} is treated by the same argument. This gives
\begin{equation*}
G^{l,\lambda_0}_{<\delta\lambda_1}
\bigl(\eta_{\lambda_1}
L^{l,\lambda_0}_{<\delta\lambda_1}\varphi_{\lambda_0}\bigr)
=
P_{<M\lambda_0}T_G
\bigl(\eta_{\lambda_1}
L^{l,\lambda_0}_{<\delta\lambda_1}\varphi_{\lambda_0}\bigr)
+\mathcal R.
\end{equation*}
By the same reasoning, using instead the linearization formula for
$L_{<\delta\lambda_1}$, we obtain
\begin{equation*}
P_{<M\lambda_0}T_G
\bigl(\eta_{\lambda_1}
L^{l,\lambda_0}_{<\delta\lambda_1}\varphi_{\lambda_0}\bigr)
=
P_{<M\lambda_0}T_G
\bigl(\eta_{\lambda_1}T_L\varphi_{\lambda_0}\bigr)
+\mathcal R,
\end{equation*}
which proves \eqref{paratruncation}. Proposition \ref{paraexp}, applied to the $\lambda_0$ component of
$T_L\varphi$, gives, modulo an error that is perturbative after
integration in $\lambda_0$ and $\lambda_1$,
\begin{equation*}
P_{<M\lambda_0}T_G
\bigl(\eta_{\lambda_1}T_L\varphi_{\lambda_0}\bigr)
=
P_{<M\lambda_0}T_G
\bigl(\eta_{\lambda_1}T_{\mathfrak l}\varphi_{\lambda_0}\bigr)
+\mathcal R.
\end{equation*}
Here the error terms generated are estimated using the same bilinear
pairings and $L_t^4L_x^p$ estimates as above, with the choice of
pairing depending on the frequency regime. Similarly, we may replace
the outer operator $T_G$ by $T_{\mathfrak g}$ modulo a perturbative
error. Thus,
\begin{equation}\label{finalT2main}
P_{<M\lambda_0}T_G
\bigl(\eta_{\lambda_1}T_L\varphi_{\lambda_0}\bigr)
=
T_{\mathfrak g}
\bigl(\eta_{\lambda_1}T_{\mathfrak l}\varphi_{\lambda_0}\bigr)
+\mathcal R.
\end{equation}
Here we used that the principal term on the right-hand side is
localized to frequencies $\ll M\lambda_0$ to remove the outer projection. Following the same steps for the other three terms in
$\mathcal T_{2,<\delta\lambda_1}^{main}$, we obtain
\begin{equation*}
\begin{split}
P_{<M\lambda_0}\mathcal T_{2,<\delta\lambda_1}^{main}
={}&-T_{\mathfrak g}
\bigl(\eta_{\lambda_1}T_{\mathfrak l}\varphi_{\lambda_0}\bigr)
-T_{\mathfrak g}
\bigl(\eta_{\lambda_0}T_{\mathfrak l}\varphi_{\lambda_1}\bigr)
\\
&-\nabla\cdot
\bigl(\eta_{\lambda_0}T_{\mathfrak m}\varphi_{\lambda_1}\bigr)
-\nabla\cdot
\bigl(\eta_{\lambda_1}T_{\mathfrak m}\varphi_{\lambda_0}\bigr)
+\mathcal R.
\end{split}
\end{equation*}
To obtain \eqref{T2main} with the definition of the balanced paraproduct $\Pi(f,g)$ in section \ref{multpara}, we just need to commute localization through the inner principal
operators. For $\mathfrak a\in\{\mathfrak l,\mathfrak m\}$, we have
\begin{equation*}
T_{\mathfrak a}\varphi_\lambda
=
T_{\mathfrak a}S_\lambda\varphi
=
S_\lambda(T_{\mathfrak a}\varphi)
+[T_{\mathfrak a},S_\lambda]\varphi.
\end{equation*}
The commutator contains a factor of
$\nabla^2\eta$ at low frequency. Proposition \ref{Leibnizcommutator}, together with
the same bilinear pairings as in the argument for \eqref{replacebyphi}, shows that all terms containing these commutators are perturbative.
We may therefore replace
$T_{\mathfrak a}\varphi_{\lambda_j}$ by
$S_{\lambda_j}(T_{\mathfrak a}\varphi)$ modulo $\mathcal R$. By the balanced paraproduct convention in Section \ref{multpara}, we have
\begin{equation*}
\begin{split}
&\int_{0}^{\infty}\int_{\delta\lambda_0}^{\lambda_0}
\left[
\eta_{\lambda_1}S_{\lambda_0}(T_{\mathfrak l}\varphi)
+\eta_{\lambda_0}S_{\lambda_1}(T_{\mathfrak l}\varphi)
\right]
\frac{d\lambda_1}{\lambda_1}\frac{d\lambda_0}{\lambda_0}
=
\Pi(\eta,T_{\mathfrak l}\varphi),
\\
&\int_{0}^{\infty}\int_{\delta\lambda_0}^{\lambda_0}
\left[
\eta_{\lambda_1}S_{\lambda_0}(T_{\mathfrak m}\varphi)
+\eta_{\lambda_0}S_{\lambda_1}(T_{\mathfrak m}\varphi)
\right]
\frac{d\lambda_1}{\lambda_1}\frac{d\lambda_0}{\lambda_0}
=
\Pi(\eta,T_{\mathfrak m}\varphi)
\end{split}
\end{equation*}
from which we finally conclude
\begin{equation*}
\begin{split}
\int_{0}^{\infty}\int_{\delta\lambda_0}^{\lambda_0}
P_{<M\lambda_0}\mathcal T_{2,<\delta\lambda_1}^{main}
\frac{d\lambda_1}{\lambda_1}\frac{d\lambda_0}{\lambda_0}&=
-T_{\mathfrak g}\bigl(\Pi(\eta,T_{\mathfrak l}\varphi)\bigr)
-\nabla\cdot\Pi(\eta,T_{\mathfrak m}\varphi)
\\
&\qquad+\int_{0}^{\infty}\int_{\delta\lambda_0}^{\lambda_0}
\mathcal R
\frac{d\lambda_1}{\lambda_1}\frac{d\lambda_0}{\lambda_0}.
\end{split}
\end{equation*}
This proves \eqref{T2main}. It remains to show that $\mathcal{T}_{2,<\delta\lambda_1}^{err}$ is perturbative; that is, we have \eqref{T2err}. This is far easier to handle than $\mathcal{T}_{2,<\delta\lambda_1}^{main}$. We give the relevant details for the first term
\begin{equation*}
-G^{l,\lambda_0}_{<\delta\lambda_1}(\eta_{\lambda_0}\eta_{\lambda_1}L^{l,\delta\lambda_1}_{<\delta\lambda_1}B^l_{<\delta\lambda_1}),
\end{equation*}
as the other term is handled in the same way. The same arguments as above above, using \ref{basicopbounds}, Proposition \ref{paraexp}, and Proposition \ref{Bernsteinelliptic} and Proposition \ref{VBunbil}, gives
\begin{equation*}
-G^{l,\lambda_0}_{<\delta\lambda_1}(\eta_{\lambda_0}\eta_{\lambda_1}L^{l,\delta\lambda_1}_{<\delta\lambda_1}B^l_{<\delta\lambda_1})=-T_{\mathfrak{g}}\bigl(\eta_{\lambda_0}\eta_{\lambda_1}T_{\mathfrak{l}}B^l_{<\delta\lambda_1}\bigr)+\mathcal{R}.
\end{equation*}
From here, we can argue as in the proof of \eqref{replacebyphi}. Decomposing $B^l_{<\delta\lambda_1}$ into its frequency
pieces $S_\rho B^l_{<\delta\lambda_1}$, where
$\rho\lesssim M\delta\lambda_1\ll\lambda_1$, we proceed as follows. When $\mu\sim_c\max(\lambda_0,\lambda_1)$, we use the pairings
$(\eta_{\lambda_0},\eta_{\lambda_1})$ and
$(\overline u_\mu^{x_0},S_\rho B^l_{<\delta\lambda_1})$ along with Proposition \ref{VBunbil}. When $\mu\not\sim_c\max(\lambda_0,\lambda_1)$, we instead use
the unbalanced pairings
$(\overline u_\mu^{x_0},\eta_{\lambda_0})$ and
$(S_\rho B^l_{<\delta\lambda_1},\eta_{\lambda_1})$,
recalling that $\lambda_1\leq\lambda_0$.
The kernel bounds for $\mathfrak l$ and $\mathfrak g$, together
with Proposition \ref{VBunbil}, then give, after integration in
$\rho$, $\lambda_1$ and $\lambda_0$,
\begin{equation*}
\begin{split}
&\int_{\lambda_0\gtrsim_M\mu}\int_{\delta\lambda_0}^{\lambda_0}
\|\overline{u}_{\mu}^{x_0}
(T_{\mathfrak{g}}(\eta_{\lambda_0}\eta_{\lambda_1}
T_{\mathfrak{l}}B^l_{<\delta\lambda_1}))_\mu\|_{L_{t,x}^1}
\frac{d\lambda_1}{\lambda_1}\frac{d\lambda_0}{\lambda_0}
\lesssim\epsilon^4c_\mu^2\mu^{-2s_c}(1+\mu|x_0|)^{\frac12}.
\end{split}
\end{equation*}
The second term in $\mathcal{T}_{2,<\delta\lambda_1}^{err}$
is handled in the same way. This proves \eqref{T2err}, and hence \eqref{iter2decomp} and \eqref{R2mainbound}.
\\

We now verify \eqref{seconddecompauxR} and \eqref{seconddecompaux}.
Instead of repeating all of the calculations in detail, we remark
on the small adjustments needed compared to the harder bound
\eqref{R2mainbound} that was proved above. The errors above consist,
in order of appearance, of the high-frequency tails, the
$\mathcal C_{\lambda_2,\lambda_1}$ and
$\mathcal E_{\lambda_2,\lambda_1}$ contributions in
\eqref{replacebyphi}, the term $\mathcal R_{cmp}$
in \eqref{G2trunc}, the contribution of $\mathcal R_F$, the terms
in \eqref{GLdecompbal} (including the output range
$\lambda_4>\lambda_3$), the subsequent operator replacements and
localization commutators, and finally the terms in \eqref{T2err}.

First, we can easily dispense with the parts of $\mathcal R$
that can be directly estimated in $L_t^1L_x^\infty$. By using
Bernstein, this includes all parts of $\mathcal{R}$ that satisfy a perturbative bound directly in $L_t^1L_x^2$.
More precisely, denoting the sum of these terms at fixed
$\lambda_0,\lambda_1$ by $\mathcal R^1$, using the $L_t^1L_x^2$ estimates from above at the critical regularity $s_c=2$ give
\begin{equation*}
\begin{split}
\int_0^\infty\int_{\delta\lambda_0}^{\lambda_0}
\|\mathcal R^1_\lambda
\overline u_\mu^{x_0}u_\mu\|_{L_{t,x}^1}
\frac{d\lambda_1}{\lambda_1}\frac{d\lambda_0}{\lambda_0}&\lesssim
\|u_\mu\|_{L_t^\infty L_x^2}^2
\int_0^\infty\int_{\delta\lambda_0}^{\lambda_0}
\|\mathcal R^1_\lambda\|_{L_t^1L_x^\infty}
\frac{d\lambda_1}{\lambda_1}\frac{d\lambda_0}{\lambda_0}
\\
&\lesssim_\delta
\epsilon^6c_\mu^2\lambda^{-1}\mu^{-2s_c}
\\
&\lesssim
\epsilon^4c_\mu^2\lambda^{-1}\mu^{-2s_c},
\end{split}
\end{equation*}
where we used $\epsilon\ll1$ in the last inequality. The terms of this form include the high-frequency tails that
appeared in the above analysis (which were estimated using Proposition \ref{Bernsteinelliptic}), as well as $\mathcal R_F$ and
the second triple integral in \eqref{GLdecompbal}. In addition
to these errors, the other main term in the above analysis that
belongs to the $\mathcal R^1_\lambda$ box comes from the term
$\mathcal E_{\lambda_2,\lambda_1}$ in \eqref{replacebyphi}.
This is because \eqref{Phidifflowerror} and two
$L_t^4L_x^\infty$ estimates give
\begin{equation*}
\lambda_0^2
\|\mathcal E_{\lambda_2,\lambda_1}\|_{L_t^2L_x^\infty}
\|\eta_{\lambda_1}\|_{L_t^4L_x^\infty}
\|\widetilde\eta_{\lambda_2}\|_{L_t^4L_x^\infty}
\lesssim_\delta\epsilon^4\lambda_0^{-1}.
\end{equation*}
Integration in $\lambda_2\approx\lambda_0$,
$\delta\lambda_0\leq\lambda_1\leq\lambda_0$, and
$\lambda_0\gtrsim_M\lambda$ therefore gives the required
$L_t^1L_x^\infty$ bound for this part.

For the remaining cubic part of $\mathcal R$, which we denote
by $\mathcal R^2$, we use the same operator decompositions as in the estimates for \eqref{R2mainbound}, just with different H\"older allocations. Estimating the
contribution of this term in \eqref{seconddecompauxR} is
substantially easier than \eqref{R2mainbound} because the
resulting expression we have to estimate in $L_{t,x}^1$ is
quintic rather than quartic. Nevertheless, we still need to
use at least one bilinear estimate here, so some care is still required. We give the details for the most delicate
triple integral term first, namely the first triple integral in
\eqref{GLdecompbal}.

First, using Proposition \ref{paraexp} and Proposition
\ref{Bernsteinelliptic}, we can replace both the outer operator
$G_{<\delta\lambda_4}^{l,\lambda_4}$ and
$L_{<\delta\lambda_2}^{l,\lambda_2}$ by their principal parts;
the subprincipal error terms can simply be estimated in
$L_t^1L_x^\infty$ using Bernstein, the bound
\eqref{fbalancedbound} and two $L_t^4L_x^\infty$-based bounds,
as above (in which case we can place these terms in the $\mathcal{R}^1$ box). What remains is to estimate is the principal contribution
\begin{equation*}
P_{<M\lambda_0}P_\lambda
\int_0^{\lambda_1}\int_{\lambda_2}^{\lambda_1}
\int_0^{\lambda_3}
T_{\mathfrak g}S_{\lambda_4}
\bigl(\eta_{\delta\lambda_3}T_{\mathfrak l}f_{\lambda_2}\bigr)
\frac{d\lambda_4}{\lambda_4}
\frac{d\lambda_3}{\lambda_3}\frac{d\lambda_2}{\lambda_2}.
\end{equation*}
Here the projection $P_\lambda$ restricts
$\lambda_4\approx\lambda$, so the $\lambda_4$ integration is
over a comparable range to the output frequency. Moreover, the kernel of the outer paradifferential operator $T_{\mathfrak{g}}S_{\lambda_4}$ is rapidly decaying at
scale $\lambda^{-1}$.

To estimate this term, we separate into two cases.
If $\lambda_0\ll\mu$, we pair $\eta_{\delta\lambda_3}$ with
$\overline u_\mu^{x_0}$ in $L_{t,x}^2$ and use one unbalanced
bilinear $L^2$ bound, put the remaining $u_\mu$ in
$L_t^\infty L_x^2$, and estimate using Bernstein,
\eqref{fbalancedbound} (with $s=s_c=2$), and
$\lambda_2\lesssim\lambda_0$,
\begin{equation*}
\begin{split}
\|f_{\lambda_2}\|_{L_t^2L_x^\infty}\lesssim
\lambda_2\|f_{\lambda_2}\|_{L_{t,x}^2}\lesssim_\delta
\epsilon^2\lambda_2^{\frac12}\lambda_0^{-\frac94}
\lesssim_\delta\epsilon^2\lambda_0^{-\frac74}.
\end{split}
\end{equation*}
Combining everything above and using the kernel bound for
$T_{\mathfrak g}$, we have
\begin{equation*}
\begin{split}
\|u_\mu\overline u_\mu^{x_0}P_\lambda T_{\mathfrak g}
S_{\lambda_4}
\bigl(\eta_{\delta\lambda_3}T_{\mathfrak l}f_{\lambda_2}\bigr)
\|_{L_{t,x}^1}&\lesssim
\lambda_4\lambda_2
\|L(\overline u_\mu^{x_0},\eta_{\delta\lambda_3})\|_{L_{t,x}^2}
\|f_{\lambda_2}\|_{L_t^2L_x^\infty}
\|u_\mu\|_{L_t^\infty L_x^2}
\\
&\qquad\lesssim_\delta
\epsilon^5c_\mu^2\mu^{-2s_c}(1+\mu|x_0|)^{\frac12}
\lambda_0^{-\frac74}
\lambda_4\lambda_2\lambda_3^{-\frac54}
\left(1+\frac{\delta\lambda_3}{\lambda_4}\right)^{\frac12}
\\
&\qquad\lesssim_\delta
\epsilon^5c_\mu^2\mu^{-2s_c}(1+\mu|x_0|)^{\frac12}
\lambda_0^{-\frac74}\lambda_3^{-\frac14}\lambda_2.
\end{split}
\end{equation*}
Here, we note carefully that the bilinear form $L$ includes the translations from the outer kernel,
which is only rapidly decaying at scale $\lambda_4^{-1}$ since
$\lambda_4\approx\lambda$. This, along with the translation factor in the unbalanced bilinear estimate is the reason for the last factor in the second line above. The last inequality follows from
\begin{equation*}
\lambda_4\left(1+\frac{\delta\lambda_3}{\lambda_4}\right)^{\frac12}
\leq\lambda_4+\delta\lambda_3\lesssim\lambda_3.
\end{equation*}
Consequently, for fixed $\lambda_0$ and $\lambda_1$, integration
first in $\lambda_4\approx\lambda$, and then in $\lambda_3$
and $\lambda_2$, gives
\begin{equation*}
\begin{split}
\epsilon^5c_\mu^2\mu^{-2s_c}(1+\mu|x_0|)^{\frac12}
\lambda_0^{-\frac74}
\int_0^{\lambda_1}\int_{\lambda_2}^{\lambda_1}
\lambda_3^{-\frac14}\lambda_2
\frac{d\lambda_3}{\lambda_3}\frac{d\lambda_2}{\lambda_2}\lesssim_\delta
\epsilon^5c_\mu^2\mu^{-2s_c}(1+\mu|x_0|)^{\frac12}
\lambda_0^{-1}.
\end{split}
\end{equation*}
This suffices upon integration over $\lambda_0\gtrsim\lambda$. If $\lambda_0\gtrsim\mu$, we instead keep
$f_{\lambda_2}$ in $L_{t,x}^2$ using \eqref{fbalancedbound},
and put $\eta_{\delta\lambda_3}$ in $L_t^4L_x^\infty$.
For fixed $\lambda_0$ and $\lambda_1$, integrating first in
$\lambda_4\approx\lambda$, and then in $\lambda_3$ and
$\lambda_2$, gives
\begin{equation*}
\begin{split}
&\left\|
P_{<M\lambda_0}P_\lambda
\int_0^{\lambda_1}\int_{\lambda_2}^{\lambda_1}
\int_0^{\lambda_3}
T_{\mathfrak g}S_{\lambda_4}
\bigl(\eta_{\delta\lambda_3}T_{\mathfrak l}f_{\lambda_2}\bigr)
\frac{d\lambda_4}{\lambda_4}
\frac{d\lambda_3}{\lambda_3}\frac{d\lambda_2}{\lambda_2}
\right\|_{L_t^{4/3}L_x^2}
\\
&\qquad\lesssim_\delta
\epsilon^3\lambda_0^{-\frac94}
\int_0^{\lambda_1}\int_{\lambda_2}^{\lambda_1}
\lambda_3^{-\frac18}\lambda_2^{\frac12}
\frac{d\lambda_3}{\lambda_3}\frac{d\lambda_2}{\lambda_2}
\\
&\qquad\lesssim_\delta\epsilon^3\lambda_0^{-\frac{15}{8}}.
\end{split}
\end{equation*}
Here we used $\lambda_4\leq\lambda_3$ to bound the contribution
from the outer operator by $\lambda_3$. We can then estimate $u_\mu^{x_0}$ and $u_\mu$ in
$L_t^4L_x^\infty$ and $L_t^\infty L_x^2$, respectively,
which bounds the product of this triple integral with the
two external factors $u_{\mu}$ and $\overline{u}_{\mu}^{x_0}$ in $L_{t,x}^1$ by
$\epsilon^5c_\mu^2\lambda_0^{-15/8}\mu^{7/8-2s_c}$.
Finally, integration in
$\delta\lambda_0\leq\lambda_1\leq\lambda_0$ and
$\lambda_0\gtrsim\mu$ gives
\begin{equation*}
\epsilon^5c_\mu^2\mu^{\frac78-2s_c}
\int_{\lambda_0\gtrsim\mu}
\int_{\delta\lambda_0}^{\lambda_0}
\lambda_0^{-\frac{15}{8}}
\frac{d\lambda_1}{\lambda_1}\frac{d\lambda_0}{\lambda_0}
\lesssim_\delta
\epsilon^5c_\mu^2\mu^{-1-2s_c},
\end{equation*}
which suffices in this region. This handles the contribution of the first triple integral in
\eqref{GLdecompbal}. The term $\mathcal R_{cmp}$ and the contribution of
$\lambda_3<\lambda_4\lesssim_M\lambda_3$ in
\eqref{GLdecompbal} are treated by the same pairing of
$\eta_{\delta\lambda_3}$ with $\overline u_\mu^{x_0}$
in $L_{t,x}^2$, since $\delta\lambda_3\ll\lambda\ll\mu$
for their principal contributions.
\medskip

The other cubic errors listed above also obey the same two bounds.
For the $\mathcal C_{\lambda_2,\lambda_1}$ contribution in
\eqref{replacebyphi}, once again, we can replace $G$ and $L$ by their
principal operators $T_{\mathfrak g}$ and $T_{\mathfrak l}$
modulo terms of the type $\mathcal R^1$.

If $\lambda_0\ll\mu$, we use the pairing
$(\mathcal C_{\lambda_2,\lambda_1},\overline u_\mu^{x_0})$
in $L_{t,x}^2$, put both $\eta_{\lambda_1}$ and
$\widetilde\eta_{\lambda_2}$ in $L_t^4L_x^\infty$,
and put the remaining $u_\mu$ in $L_t^\infty L_x^2$.
If $\lambda_0\gtrsim\mu$, we instead use the pairing
$(\mathcal C_{\lambda_2,\lambda_1},\widetilde\eta_{\lambda_2})$
in \eqref{Phidifflowpairings}, put $\eta_{\lambda_1}$ in
$L_t^4L_x^\infty$, and estimate $\overline u_\mu^{x_0}$
and $u_\mu$ in $L_t^4L_x^\infty$ and
$L_t^\infty L_x^2$, respectively.
The kernel estimates are handled as above. Integration in
$\lambda_2\approx\lambda_0$ then gives the same two bounds. The later cubic terms coming from replacements of the operators by their principal parts and commutators of these principal operators with Littlewood-Paley projections are treated in the same way, with
the low-frequency $\nabla^2\eta$ factor in place of
$\mathcal C_{\lambda_2,\lambda_1}$.
\medskip

It remains to handle \eqref{T2err}. Once again, modulo an error of the type $\mathcal{R}^1$, we can replace each operator by its principal part. Here, when $\lambda_0\ll\mu$, we pair one
$\eta_{\lambda_j}$ with $u_{\mu}^{x_0}$ and place the other factor of $u_{\mu}$ in $L_t^{\infty}L_x^2$, and use two $L_t^4L_x^{\infty}$ estimates for the remaining two factors. When $\lambda_0\gtrsim\mu$, use the
internal unbalanced $L_{t,x}^2$ pairing between a frequency piece of
$B^l_{<\delta\lambda_1}$ and one $\eta_{\lambda_j}$ using Proposition \ref{VBunbil}.
The remaining internal factors are placed in $L_t^4L_x^\infty$. Combining all of the cubic error bounds above, integration over $\lambda_0\gtrsim_M\lambda$ gives
\begin{equation*}
\begin{split}
\int_0^\infty\int_{\delta\lambda_0}^{\lambda_0}
\|\mathcal R^2_\lambda
\overline u_\mu^{x_0}u_\mu\|_{L_{t,x}^1}
\frac{d\lambda_1}{\lambda_1}\frac{d\lambda_0}{\lambda_0}\lesssim_\delta
\epsilon^5c_\mu^2\lambda^{-1}
\mu^{-2s_c}(1+\mu|x_0|)^{\frac12}.
\end{split}
\end{equation*}
Together with the preceding $L_t^1L_x^\infty$ bounds, this
proves \eqref{seconddecompauxR}.

Combining \eqref{seconddecompauxR}, \eqref{firstdecompauxR},
and \eqref{R21bound} (using Bernstein and two energy bounds
for the last term) proves \eqref{seconddecompaux}.
\\

Finally, we prove the $L_{t,x}^2$ error bound \eqref{seconddecompL2}. We begin with the contribution of $\mathcal{R}^2_{\eta}$. This term is essentially cubic. Consequently, the $L_{t,x}^2$ bound for this term is much easier than the multilinear bound proved above, and we do not need to rely too heavily on the precise multilinear structure of this term here. We describe the estimates for each of the relevant error term that appeared above. First, by Corollary \ref{Phil0varphibound} we obtain
\begin{equation*}
\|\mathcal F_{\lambda_0,\lambda_1}\|_{L_{t,x}^2}
\lesssim_\delta\epsilon^2c_{\lambda_0}
\lambda_0^{-\frac34-s_c}.
\end{equation*}
Proposition \ref{basicopbounds} and
$\|\eta_{\lambda_1}\|_{L_{t,x}^\infty}
\lesssim_\delta\epsilon\lambda_0^{-1}$ therefore bound the
left-hand side of \eqref{replacebyphi} in $L_{t,x}^2$ by
$\epsilon^3c_{\lambda_0}\lambda_0^{\frac14-s_c}$.
For the double integral in \eqref{G2trunc}, we can use
\eqref{fbalancedbound} to similarly bound
\begin{equation*}
\begin{split}
\|G^{l,\lambda_3}_{<\delta\lambda_3}
(\eta_{\delta\lambda_3}L^{l,\lambda_3}_{<\delta\lambda_3}
f_{\lambda_2})\|_{L_{t,x}^2}\lesssim_\delta
\lambda_3\lambda_2
\|\eta_{\delta\lambda_3}\|_{L_{t,x}^\infty}
\|f_{\lambda_2}\|_{L_{t,x}^2}
\lesssim_\delta\epsilon^3c_{\lambda_0}
\lambda_0^{-\frac14-s_c}\lambda_2^{\frac12}.
\end{split}
\end{equation*}
The divergence term has the same bound, and thus, upon integrating in $\lambda_2$ and $\lambda_3$, the full contribution of the double integral term also has $L_{t,x}^2$ size
$\epsilon^3c_{\lambda_0}\lambda_0^{\frac14-s_c}$.
The term $\mathcal R_{cmp}$ is handled similarly, while the first term in the third line of \eqref{G2trunc} is a high-frequency tail error that can be dispensed with by using Proposition \ref{Bernsteinelliptic}. The remaining "sub-principal" replacement errors and the contribution of \eqref{T2err} are at least cubic (with the lowest frequency factor differentiated enough times), and we can always use two $L_t^4L_x^{\infty}$ bounds alongside one $L_t^{\infty}L_x^2$ bound for such terms, giving (by H\"older) the same desired $O_{L_{t,x}^2}(\epsilon^3c_{\lambda_0}\lambda_0^{\frac{1}{4}-s_c})$ bound. Upon integrating in $\lambda_0$ and $\lambda_1$, we obtain
\begin{equation*}
\|(\mathcal R_\eta^2)_\mu\|_{L_{t,x}^2}
\lesssim_\delta\epsilon^3c_\mu\mu^{\frac14-s_c}.
\end{equation*}
This proves \eqref{seconddecompL2} and finally completes the proof of Proposition \ref{seconddecomposition}.
\end{proof}


\subsection{Extraction of \texorpdfstring{$N_{\eta}^{[3]}$}{} and \texorpdfstring{$N_{\eta}^{[4]}$}{}}\label{subsec:Neta34}
Finally, we extract the third and fourth terms, $(G(\eta)\psi)^{[3]}$ and $N_{\eta}^{[4]}$, in the expansion of $G(\eta)\psi$. We need to begin with a few notational preliminaries.
\\

While we have to identify $(G(\eta)\psi)^{[3]}$ precisely, fortunately, we will not need to track the form of $N_{\eta}^{[4]}$ as carefully as the other terms in the expansion, but rather only identify its basic structural form. The following definition will suffice for our purposes.

\begin{definition}\label{N4def}
For each $\lambda\in 2^{\mathbb Z}$, we use $N_{\eta,\lambda}^{[4]}$ to denote any finite linear combination of balanced quartic paradifferential expressions with the following properties:
\begin{itemize}
\item There are exactly four explicit inputs, consisting of three factors of $\eta$ and one factor of $\varphi$, all localized at frequencies comparable to $\lambda$, and the output is also localized at frequency comparable to $\lambda$.
\item The expression is scalar-valued and is obtained from pointwise products and contractions by composing the principal operators $T_{\mathfrak g}$, $T_{\mathfrak l}$, and $T_{\mathfrak m}$ with translation-invariant multilinear forms. It has total differential order four, where an exterior first-order derivative, such as the divergence appearing below, also counts toward the total order.
\end{itemize}
The low-frequency dependence of the symbols on $\eta$ is regarded as a paracoefficient and is not counted among the four explicit inputs. We also allow bounded integrations over comparable frequency parameters. For example, if $\lambda_j\approx\lambda$ for $1\leq j\leq4$, then
\begin{equation*}
P_{\lambda}T_{\mathfrak g}\left(
\eta_{\lambda_1}T_{\mathfrak l}\left(
\eta_{\lambda_2}\eta_{\lambda_3}T_{\mathfrak l}^2\varphi_{\lambda_4}
\right)\right)
\end{equation*}
is a term of the form $N_{\eta,\lambda}^{[4]}$. We write
\begin{equation*}
N_{\eta}^{[4]}\coloneq\sum_{\lambda}N_{\eta,\lambda}^{[4]}.
\end{equation*}
\end{definition}
Next, we proceed with the final iteration in the expansion of $G(\eta)\psi$. The principal cubic term $(G(\eta)\psi)^{[3]}$ in this expansion will have three comparable frequency inputs, and will therefore have to be precisely defined in terms of certain trilinear paradifferential operators. To define these cubic terms precisely, let
\begin{equation*}
\mathcal D\coloneq\left\{(\lambda_0,\lambda_1,\lambda_2):
0<\lambda_0<\infty,\quad
\delta\lambda_0\leq\lambda_1\leq\lambda_0,\quad
\delta\lambda_1\leq\lambda_2\leq\lambda_1\right\}
\end{equation*}
and set
\begin{equation*}
d\Omega\coloneq
\frac{d\lambda_2}{\lambda_2}
\frac{d\lambda_1}{\lambda_1}
\frac{d\lambda_0}{\lambda_0}.
\end{equation*}
Consistently with the balanced paraproduct convention for bilinear paraproducts in Section \ref{multpara}, in this subsection we define the trilinear balanced paraproduct
\begin{equation}\label{trilinearPiconvention}
\Pi(f,g,h)
\coloneq
\int_{\mathcal D}\sum_{\sigma\in\mathfrak S_{3}}
f_{\lambda_{\sigma(0)}}g_{\lambda_{\sigma(1)}}h_{\lambda_{\sigma(2)}}\,d\Omega.
\end{equation}
For $\mathfrak a\in\{\mathfrak l,\mathfrak m\}$, we also use the convention
\begin{equation}\label{nestedPiconvention}
\Pi\bigl(\eta,T_{\mathfrak a}(\eta,T_{\mathfrak l}\varphi)\bigr)
\coloneq
\int_{\mathcal D}\sum_{\sigma\in\mathfrak S_{3}}
\eta_{\lambda_{\sigma(0)}}
T_{\mathfrak a}\left(
\eta_{\lambda_{\sigma(1)}}
(T_{\mathfrak l}\varphi)_{\lambda_{\sigma(2)}}
\right)d\Omega.
\end{equation}
Above $\mathfrak{S}_3$ just denotes the group of permutations of $\{0,1,2\}$. We use $\mathfrak{S}_3$ instead of the more standard convention $S_3$ to avoid conflict with the $S_{\lambda}=\lambda\frac{d}{d\lambda}P_{<\lambda}$ projection notation. Here $(T_{\mathfrak l}\varphi)_\lambda=S_\lambda(T_{\mathfrak l}\varphi)$, and the right-hand side contains all six permutations of the three comparable-frequency inputs. The next proposition gives the final form of the expansion for $N_{\eta}=G(\eta)\psi$ that we will need in our analysis later.
\begin{proposition}\label{thirddecomposition}
Under the assumptions of Proposition \ref{seconddecomposition}, let $(G(\eta)\psi)^{[1]}$ and $(G(\eta)\psi)^{[2]}$ be as in \eqref{G1def} and \eqref{G2def}, respectively. Define
\begin{equation}\label{G3def}
\begin{split}
(G(\eta)\psi)^{[3]}\coloneq{}&
T_{\mathfrak g}\Pi\bigl(\eta,T_{\mathfrak l}(\eta,T_{\mathfrak l}\varphi)\bigr)
+\nabla\cdot\Pi\bigl(\eta,T_{\mathfrak m}(\eta,T_{\mathfrak l}\varphi)\bigr)
\\
&-\frac{1}{2}T_{\mathfrak g}\Pi\bigl(\eta,\eta,T_{\mathfrak l}^2\varphi\bigr)
-\frac{1}{2}\nabla\cdot\Pi\bigl(\eta,\eta,T_{\mathfrak m}T_{\mathfrak l}\varphi\bigr).
\end{split}
\end{equation}
Then there holds
\begin{equation}\label{iter3decomp}
G(\eta)\psi
=\sum_{k=1}^{3}(G(\eta)\psi)^{[k]}+N_{\eta}^{[4]}+\mathcal R_{\eta},
\end{equation}
where $N_{\eta}^{[4]}$ is as in Definition \ref{N4def}, and we have the following bounds
\begin{enumerate}
\item For every $\mu\in2^{\mathbb Z}$ and $x_0\in\mathbb R^2$, there holds
\begin{equation*}
\|(\mathcal R_{\eta})_{\mu}\overline u^{x_0}_{\mu}\|_{L_t^1L_x^1}
\lesssim \epsilon^4c_{\mu}^2\mu^{-2s_c}(1+\mu|x_0|)^{\frac12}.
\end{equation*}
\item For every $0<\lambda\ll\mu$ and $x_0\in\mathbb R^2$, there holds
\begin{equation}\label{thirddecompaux}
\begin{split}
&\|(\mathcal R_\eta)_\lambda
\overline u_\mu^{x_0}u_\mu\|_{L_{t,x}^1}
+\|((G(\eta)\psi)^{[3]})_\lambda
\overline u_\mu^{x_0}u_\mu\|_{L_{t,x}^1}+\|(N_\eta^{[4]})_\lambda
\overline u_\mu^{x_0}u_\mu\|_{L_{t,x}^1}
\\
&\qquad\lesssim\epsilon^4c_\mu^2\lambda^{-1}
\mu^{-2s_c}(1+\mu|x_0|)^{\frac12}.
\end{split}
\end{equation}
\item For every $\mu>0$, there holds
\begin{equation}\label{thirddecompL2}
\left\|\left(G(\eta)\psi-(G(\eta)\psi)^{[1]}
-(G(\eta)\psi)^{[2]}\right)_\mu\right\|_{L_{t,x}^2}
\lesssim_\delta\epsilon^3c_\mu\mu^{\frac14-s_c}.
\end{equation}
\end{enumerate}
\end{proposition}

\begin{proof}
We first compute the $\lambda_2$ derivative of each component of $\mathcal T^l_{2,<\lambda_2}$. Throughout this subsection, we will again adopt (as in Proposition \ref{seconddecomposition}) the notation $A_{<\lambda}^{l,\beta}=P_{<M\beta}A_{<\lambda}$ and again the notation \eqref{Phitruncnot}. To ensure simplicity of notation, we also adopt the following shorthand for the remainder of this subsection
\begin{equation*}
\mathbf P_0\coloneq P_{<M\lambda_0},\hspace{5mm}
\eta_j\coloneq\eta_{\lambda_j},\hspace{5mm}B_j\coloneq B^l_{<\lambda_j},\hspace{5mm}
A_j^i\coloneq A_{<\lambda_j}^{l,\lambda_i},
\end{equation*}
where $A\in\{G,L,M\}$ and $i,j\in\{0,1,2\}.$ We also write
\begin{equation*}
B^{l,\beta}_{<\rho}\coloneq P_{<M\beta}B_{<\rho}.
\end{equation*}
We will abbreviate $\Phi_{j}\coloneq \Phi_{<\lambda_2}^{l,j}$ (only the $<\lambda_2$ coefficient has to be isolated in our formulas below). When the lower index in $A_j^i$ is $2$, we write $A_{2}^i=A_i$.
\\

Now, we are finally ready to proceed with the proof of Proposition \ref{thirddecomposition}. Starting from Proposition \ref{seconddecomposition}, we begin by computing the first part of \eqref{T2deriv}. Using the linearization formulas for $G$, $L$, $M$, and $B$, we obtain
\begin{equation*}
\mathbf P_0\lambda_2\frac{d}{d\lambda_2}\mathcal T_{2,<\lambda_2}^{main}
=C_{3,<\lambda_2}^1+C_{3,<\lambda_2}^2+Q_{3,<\lambda_2}^1+Q_{3,<\lambda_2}^2+\mathcal E_{3,<\lambda_2}^{main},
\end{equation*}
where
\begin{equation*}
\begin{split}
C_{3,<\lambda_2}^1\coloneq\mathbf P_0\Bigl\{&
G_0\bigl(\eta_2L_0(\eta_1L_0\Phi_0)\bigr)
+\nabla\cdot\bigl(\eta_2M_0(\eta_1L_0\Phi_0)\bigr)
\\
&+G_0\bigl(\eta_1L_0(\eta_2L_0\Phi_0)\bigr)
-G_0\bigl(\eta_1\eta_2(L_0)^2\Phi_0\bigr)
+G_0\bigl(\eta_1L_0(\eta_0L_2\Phi_2)\bigr)
\\
&+\nabla\cdot\Bigl(
\eta_1M_0(\eta_2L_0\Phi_0)
-\eta_1\eta_2M_0L_0\Phi_0
+\eta_1M_0(\eta_0L_2\Phi_2)
\Bigr)\Bigr\}.
\end{split}
\end{equation*}
The second cubic term is
\begin{equation*}
\begin{split}
C_{3,<\lambda_2}^2\coloneq\mathbf P_0\Bigl\{&
G_0\bigl(\eta_2L_0(\eta_0L_1\Phi_1)\bigr)
+\nabla\cdot\bigl(\eta_2M_0(\eta_0L_1\Phi_1)\bigr)
\\
&+G_0\bigl(\eta_0L_1(\eta_2L_1\Phi_1)\bigr)
-G_0\bigl(\eta_0\eta_2(L_1)^2\Phi_1\bigr)
+G_0\bigl(\eta_0L_1(\eta_1L_2\Phi_2)\bigr)
\\
&+\nabla\cdot\Bigl(
\eta_0M_1(\eta_2L_1\Phi_1)
-\eta_0\eta_2M_1L_1\Phi_1
+\eta_0M_1(\eta_1L_2\Phi_2)
\Bigr)\Bigr\}.
\end{split}
\end{equation*}
The corresponding quartic terms are
\begin{equation*}
Q_{3,<\lambda_2}^1\coloneq\mathbf P_0\Bigl\{
G_0\bigl(\eta_1L_0(\eta_0\eta_2L_2B_2)\bigr)
+\nabla\cdot\bigl(\eta_1M_0(\eta_0\eta_2L_2B_2)\bigr)
\Bigr\},
\end{equation*}
and
\begin{equation*}
Q_{3,<\lambda_2}^2\coloneq\mathbf P_0\Bigl\{
G_0\bigl(\eta_0L_1(\eta_1\eta_2L_2B_2)\bigr)
+\nabla\cdot\bigl(\eta_0M_1(\eta_1\eta_2L_2B_2)\bigr)
\Bigr\}.
\end{equation*}
Similarly,
\begin{equation*}
\mathbf P_0\lambda_2\frac{d}{d\lambda_2}\mathcal T_{2,<\lambda_2}^{err}
=C_{3,<\lambda_2}^3+Q_{3,<\lambda_2}^3+\mathcal E_{3,<\lambda_2}^{err},
\end{equation*}
where $(L_2)^2$ denotes the square of $L_2$:
\begin{equation*}
C_{3,<\lambda_2}^3\coloneq
-\mathbf P_0\Bigl\{
G_0\bigl(\eta_0\eta_1(L_2)^2\Phi_2\bigr)
+\nabla\cdot\bigl(\eta_0\eta_1M_2L_2\Phi_2\bigr)
\Bigr\},
\end{equation*}
and
\begin{equation*}
\begin{split}
Q_{3,<\lambda_2}^3\coloneq\mathbf P_0\Bigl\{&
G_0\bigl(\eta_2L_0(\eta_0\eta_1L_2B_2)\bigr)
+\nabla\cdot\bigl(\eta_2M_0(\eta_0\eta_1L_2B_2)\bigr)
\\
&-G_0\bigl(\eta_0\eta_1\eta_2(L_2)^2B_2\bigr)
-\nabla\cdot\bigl(\eta_0\eta_1\eta_2M_2L_2B_2\bigr)
\Bigr\}.
\end{split}
\end{equation*}
In the last two equations, we used
\begin{equation*}
\begin{split}
\left.\rho\frac{d}{d\rho}
\bigl(L^{l,\lambda_2}_{<\rho}B^{l,\lambda_2}_{<\rho}\bigr)
\right|_{\rho=\lambda_2}
&=(L_2)^2\Phi_2+\eta_2(L_2)^2B_2,\\
\left.\rho\frac{d}{d\rho}
\bigl(M^{l,\lambda_2}_{<\rho}B^{l,\lambda_2}_{<\rho}\bigr)
\right|_{\rho=\lambda_2}
&=M_2L_2\Phi_2+\eta_2M_2L_2B_2,
\end{split}
\end{equation*}
modulo terms produced by the output projections (to be included in the $\mathcal E_{3,<\lambda_2}$ term).
\\

Finally, the expressions
$\mathcal E_{3,<\lambda_2}^{main}$ and
$\mathcal E_{3,<\lambda_2}^{err}$ are defined by collecting two types of terms:
the cutoff-shell terms produced by differentiating the matched
projections $P_{<M\rho}$, and the internal high-frequency tails
produced when fixed output projections are inserted between successive
operators. Schematically, these include terms containing
\begin{equation*}
\left(\rho\frac{d}{d\rho}P_{<M\rho}\right)A_{<\rho}F
\qquad\text{or}\qquad
P_{<M\beta}A_{<\rho}
\bigl(f\,P_{\gtrsim M\beta}A'_{<\rho}h\bigr),
\end{equation*}
where $A,A'\in\{G,L,M\}$ and $\beta$ is one of the fixed frequency
scales. In all the terms arising here, the frequencies entering $F$ and $h$
are bounded by $C_\delta\rho$ and $C_\delta\beta$, respectively. These terms all involve high-frequency tails and will be directly perturbative in $L_t^1L_x^2$ by Proposition \ref{Bernsteinelliptic}. We now set
\begin{equation*}
C_{3,<\lambda_2}\coloneq\sum_{i=1}^3C_{3,<\lambda_2}^i,\qquad
Q_{3,<\lambda_2}\coloneq\sum_{i=1}^3Q_{3,<\lambda_2}^i,
\end{equation*}
and
\begin{equation*}
\mathcal E_{3,<\lambda_2}
\coloneq\mathcal E_{3,<\lambda_2}^{main}
+\mathcal E_{3,<\lambda_2}^{err},
\end{equation*}
so that $\mathcal{T}_{3,<\lambda_2}$ (from \eqref{T2deriv}) can be written as
\begin{equation}\label{T3CQdecomp}
\mathcal T_{3,<\lambda_2}
=C_{3,<\lambda_2}+Q_{3,<\lambda_2}+\mathcal E_{3,<\lambda_2}.
\end{equation}
For the remainder of the proof, $\mathcal R$ denotes any term whose contribution, after the necessary integrations and pairing with $\overline u_\mu^{x_0}$, is bounded in $L_{t,x}^1$ by
\begin{equation*}
\epsilon^4c_\mu^2\mu^{-2s_c}(1+\mu|x_0|)^{\frac12}.
\end{equation*}
\subsubsection{The quartic term $Q_{3,<\lambda_2}$.}
We do not need to iterate the quartic term again. Since
\begin{equation*}
\lambda_0\approx_\delta\lambda_1\approx_\delta\lambda_2,
\end{equation*}
we separate it directly into balanced quartic terms, which contribute to $N_\eta^{[4]}$, and perturbative terms. The latter consist either of quartic interactions with a separated input or output frequency, or of quintic and higher-order terms. Isolating this will produce two kinds of errors. The first kind consists of quartic terms with at least one frequency separated pair which can be dispensed with using one unbalanced bilinear estimate, two $L_t^4L_x^\infty$ bounds, and one $L_t^{\infty}L_x^2$ bound. The second type of error consists of certain quintic terms, where we can place four factors in $L_t^4L_x^\infty$ and the remaining factors in $L_t^\infty L_x^2$.
In particular, $Q_{3,<\lambda_2}$ makes no contribution to $(G(\eta)\psi)^{[3]}$.
\\

We illustrate the argument for the representative term
\begin{equation*}
Q'_{<\lambda_2}\coloneq
\mathbf P_0G_0\bigl(\eta_1L_0(\eta_0\eta_2L_2B_2)\bigr),
\end{equation*}
as the remaining terms have the same structure. First, we show that this term is perturbative if $B_2$ is replaced by $P_{<M\lambda_2}B_{<\delta\lambda_2}$. More precisely, set
\begin{equation*}
\mathcal R_{\mathrm{low}}\coloneq
\mathbf P_0G_0\bigl(\eta_1L_0(\eta_0\eta_2L_2P_{<M\lambda_2}B_{<\delta\lambda_2})\bigr).
\end{equation*}
We claim that $\mathcal R_{\mathrm{low}}$ is perturbative. After the same principal-part extraction used above, Proposition \ref{basicopbounds}, Bernstein's inequality, and H\"older's inequality give
\begin{equation}\label{Qlowfirst}
\begin{split}
\|\mathcal R_{\mathrm{low}}\|_{L_t^1L_x^2}
&\lesssim \lambda_0^2
\|\eta_{\lambda_0}\|_{L_t^4L_x^\infty}
\|\eta_{\lambda_1}\|_{L_t^4L_x^\infty}
\|\mathcal L(\eta_{\lambda_2},L_2P_{<M\lambda_2}B_{<\delta\lambda_2})\|_{L_{t,x}^2}
\\
&\lesssim \epsilon^2c_{\lambda_0}\lambda_0^{\frac74-s_c}
\|\mathcal L(\eta_{\lambda_2},L_2P_{<M\lambda_2}B_{<\delta\lambda_2})\|_{L_{t,x}^2},
\end{split}
\end{equation}
where $\mathcal L$ is a translation-invariant bilinear form. By Proposition \ref{Bernsteinelliptic}, it suffices to estimate the remaining term with $B^l_{<\delta\lambda_2}$ in place of $P_{<M\lambda_2}B_{<\delta\lambda_2}$. We then have
\begin{equation}\label{Qlowinner}
\begin{split}
\|\mathcal L(\eta_{\lambda_2},L_2B^l_{<\delta\lambda_2})\|_{L_{t,x}^2}\lesssim
\|\eta_{\lambda_2}\|_{L_t^\infty L_x^2}
\|(L_2-T_L)B^l_{<\delta\lambda_2}\|_{L_t^2L_x^\infty}
+\|\mathcal L(\eta_{\lambda_2},T_LB^l_{<\delta\lambda_2})\|_{L_{t,x}^2}.
\end{split}
\end{equation}
By Propositions \ref{paraexp} and \ref{VBunbil} and the definition of $T_L$, the latter term can be easily estimated by
\begin{equation*}
\begin{split}
\|\mathcal L(\eta_{\lambda_2},T_LB^l_{<\delta\lambda_2})\|_{L_{t,x}^2}
&\lesssim \epsilon^2\lambda_2^{-\frac74}.
\end{split}
\end{equation*}
For the subprincipal contribution, we use the decomposition
\begin{equation*}
\begin{split}
(L_2-T_L)P_\alpha
=\bigl(L^{l,\lambda_2}_{<\lambda_2}-L^{l,\lambda_2}_{<\delta\alpha}\bigr)P_\alpha+(P_{<M\lambda_2}-P_{<M\alpha})L_{<\delta\alpha}P_\alpha
+(L^{l,\alpha}_{<\delta\alpha}-T_L)P_\alpha.
\end{split}
\end{equation*}
The middle term when applied to $B_{<\delta\lambda_2}$ contributes high-frequency tails that are easily dispensed with using Propositions \ref{Bernsteinelliptic} and \ref{VBunbil}. For the final term, we can estimate its contribution using Proposition \ref{paraexp} as above. For the first term, we apply the linearization formula \eqref{Blin} and the fundamental theorem of calculus to the difference $L_{<\lambda_2}-L_{<\delta\alpha}$ while keeping reference output projection fixed. We then obtain for $\alpha\lesssim\delta\lambda_2$,
\begin{equation*}
\begin{split}
&\|(L_2-T_L)P_\alpha B^l_{<\delta\lambda_2}\|_{L_t^2L_x^\infty}
\\
&\quad\lesssim
\int_{\delta\alpha}^{\lambda_2}
\bigl\|
-L^{l,\lambda_2}_{<\lambda_3}\bigl(
\eta_{\lambda_3}L^{l,\lambda_2}_{<\lambda_3}P_\alpha B^l_{<\delta\lambda_2}
\bigr)
+\eta_{\lambda_3}(L^{l,\lambda_2}_{<\lambda_3})^2P_\alpha B^l_{<\delta\lambda_2}
\bigr\|_{L_t^2L_x^\infty}
\frac{d\lambda_3}{\lambda_3}
+\epsilon^2\alpha^{\frac14}.
\end{split}
\end{equation*}
This term can be dispensed with using purely $L_t^4L_x^p$ estimates. Indeed, for some fixed $8<p<\infty$, Bernstein's inequality, Proposition \ref{basicopbounds}, and the bound
\begin{equation*}
\bigl\||D|^{\frac58+\frac2p}B^l_{<\delta\lambda_2}\bigr\|_{L_t^4L_x^p}
\lesssim\epsilon
\end{equation*}
along with Proposition \ref{besovparabounds} give
\begin{equation*}
\begin{split}
&\bigl\|
-L^{l,\lambda_2}_{<\lambda_3}\bigl(
\eta_{\lambda_3}L^{l,\lambda_2}_{<\lambda_3}P_\alpha B^l_{<\delta\lambda_2}
\bigr)
+\eta_{\lambda_3}(L^{l,\lambda_2}_{<\lambda_3})^2P_\alpha B^l_{<\delta\lambda_2}
\bigr\|_{L_t^2L_x^\infty}
\\
&\qquad\lesssim
\lambda_3\alpha\|P_\alpha B^l_{<\delta\lambda_2}\|_{L_t^4L_x^\infty}
\|\eta_{\lambda_3}\|_{L_t^4L_x^\infty}
\\
&\qquad\lesssim
\epsilon^2c'_{\lambda_3}\alpha^{\frac38}\lambda_3^{-\frac18}.
\end{split}
\end{equation*}
Consequently,
\begin{equation*}
\int_{\delta\alpha}^{\lambda_2}
c'_{\lambda_3}\alpha^{\frac38}\lambda_3^{-\frac18}
\frac{d\lambda_3}{\lambda_3}
\lesssim_\delta\alpha^{\frac14}.
\end{equation*}
Summing in $\alpha$ gives
\begin{equation}\label{L2TL}
\|(L_2-T_L)B^l_{<\delta\lambda_2}\|_{L_t^2L_x^\infty}
\lesssim\epsilon^2\lambda_2^{\frac14}.
\end{equation}
Using $\|\eta_{\lambda_2}\|_{L_t^\infty L_x^2}\lesssim\epsilon\lambda_2^{-2}$, we then conclude from \eqref{Qlowinner} that
\begin{equation*}
\|\mathcal L(\eta_{\lambda_2},L_2B^l_{<\delta\lambda_2})\|_{L_{t,x}^2}
\lesssim\epsilon^2\lambda_2^{-\frac74}.
\end{equation*}
Returning to \eqref{Qlowfirst}, we conclude that
\begin{equation}\label{thirditererror}
\|\mathcal R_{\mathrm{low}}\|_{L_t^1L_x^2}
\lesssim\epsilon^4c_{\lambda_0}\lambda_0^{-s_c}.
\end{equation}
After applying $P_\mu$ and using
$\|\overline u_\mu^{x_0}\|_{L_t^\infty L_x^2}\lesssim\epsilon c_\mu\mu^{-s_c}$,
the contribution of \eqref{thirditererror} is bounded by
$\epsilon^5c_\mu^2\mu^{-2s_c}$, which suffices.

We now return to identifying the leading part of $Q'_{<\lambda_2}$. Applying the fundamental theorem of calculus to $B$ with its output projection fixed gives
\begin{equation*}
\begin{split}
B_2=P_{<M\lambda_2}B_{<\delta\lambda_2}+\int_{\delta\lambda_2}^{\lambda_2}
P_{<M\lambda_2}\left(
L_{<\lambda_3}\Phi^3_{<\lambda_3}
+\eta_{\lambda_3}L_{<\lambda_3}B_{<\lambda_3}
\right)\frac{d\lambda_3}{\lambda_3}.
\end{split}
\end{equation*}
We can harmlessly replace $\Phi_{<\lambda_3}^3$ with $\Phi^{l,3}_{<\lambda_3}$ above, as the resulting error produces high-frequency tails, and can be safely included in $\mathcal R$. The term containing $\eta_{\lambda_3}L_{<\lambda_3}B_{<\lambda_3}$
is quintic. We put one high-frequency $\eta$ factor in
$L_t^\infty L_x^2$, the other three in $L_t^4L_x^\infty$, and
the differentiated $B$ factor in $L_t^4L_x^\infty$ using
\eqref{BauxiliaryPhi}. The remaining three differential operators in $Q'_{<\lambda_2}$ then show that this contribution can be estimated in $L_t^1L_x^2$ by
\begin{equation*}
\lambda_0^3
(\epsilon\lambda_0^{-\frac98})^3
(\epsilon\lambda_0^{\frac38})
(\epsilon c_{\lambda_0}\lambda_0^{-s_c})
=\epsilon^5c_{\lambda_0}\lambda_0^{-s_c},
\end{equation*}
which suffices. Therefore, we obtain
\begin{equation*}
\begin{split}
Q'_{<\lambda_2}
={}&\int_{\delta\lambda_2}^{\lambda_2}
\mathbf P_0G_0\left(
\eta_1L_0\left(
\eta_0\eta_2L_2L^{l,\lambda_2}_{<\lambda_3}\Phi^{l,3}_{<\lambda_3}
\right)\right)\frac{d\lambda_3}{\lambda_3}+\mathcal R
\\
={}&\int_{\delta\lambda_2}^{\lambda_2}
T_{\mathfrak g}\left(
\eta_{\lambda_1}T_{\mathfrak l}\left(
\eta_{\lambda_0}\eta_{\lambda_2}T_{\mathfrak l}^2\varphi_{\lambda_3}
\right)\right)\frac{d\lambda_3}{\lambda_3}+\mathcal R.
\end{split}
\end{equation*}
Here $\Phi^{l,3}_{<\lambda_3}-\varphi_{\lambda_3}$ is at least quadratic, while every error produced by replacing one of the operators by its principal part contains either a (high-frequency) factor of $\eta$ or a high-frequency tail. In making these replacements, we retain $\mathbf P_0$ until the full principal expression has been obtained. That expression is supported at frequencies $\ll M\lambda_0$, so the projection can then be removed as above.

It remains to separate the output frequency of the principal quartic term. The range $\mu\gg_M\lambda_0$ is absent from the principal expression by definition of the paradifferential operators involved and the fact that $\lambda_1\approx\lambda_2\approx\lambda_3\approx\lambda_0$. If $\mu\approx_{M,\delta}\lambda_0$, all four explicit inputs and the output are comparable, and this contribution is included in $N_{\eta,\mu}^{[4],Q}$. Finally, suppose that $\mu\ll_M\lambda_0$. Here, we can pair $\overline u_{\mu}^{x_0}$ with (for instance) $\eta_{\lambda_1}$ and estimate this in $L_{t,x}^2$. The remaining three factors can be estimated using two $L_t^4L_x^{\infty}$ bounds and one $L_t^{\infty}L_x^2$ bound. More precisely, we have
\begin{equation*}
\begin{split}
&\left\|\overline u_\mu^{x_0}P_\mu
\int_{\delta\lambda_2}^{\lambda_2}
T_{\mathfrak g}\left(
\eta_{\lambda_1}T_{\mathfrak l}\left(
\eta_{\lambda_0}\eta_{\lambda_2}T_{\mathfrak l}^2\varphi_{\lambda_3}
\right)\right)\frac{d\lambda_3}{\lambda_3}\right\|_{L_{t,x}^1}
\\
&\quad\lesssim
\mu\lambda_0^3
\int_{\delta\lambda_2}^{\lambda_2}
\|\mathcal L(\overline u_\mu^{x_0},\eta_{\lambda_0})\|_{L_{t,x}^2}
\|\eta_{\lambda_1}\|_{L_t^4L_x^\infty}
\|\eta_{\lambda_2}\|_{L_t^4L_x^\infty}
\|\varphi_{\lambda_3}\|_{L_t^\infty L_x^2}
\frac{d\lambda_3}{\lambda_3}
\\
&\quad\lesssim
\epsilon^5c_\mu c_{\lambda_0}
\mu^{\frac74-s_c}\lambda_0^{-\frac74-s_c}(1+\mu|x_0|)^{\frac12}.
\end{split}
\end{equation*}
Carrying out similar arguments for the other quartic terms in $Q_{3,<\lambda_2}$, we obtain the decomposition
\begin{equation}\label{N4Qiter}
\int_{\mathcal D}(Q_{3,<\lambda_2})_\mu\,d\Omega
=N_{\eta,\mu}^{[4],Q}+(\mathcal R_\eta^Q)_\mu,
\end{equation}
where 
\begin{equation*}
\|(\mathcal R_\eta^Q)_\mu\overline u_\mu^{x_0}\|_{L_t^1L_x^1}
\lesssim\epsilon^4c_\mu^2\mu^{-2s_c}(1+\mu|x_0|)^{\frac12}.
\end{equation*}
\subsubsection{The cubic term $C_{3,<\lambda_2}$.}
We now treat the remaining cubic term. We introduce another frequency parameter $\lambda_3$. We let
$C_{3,<\lambda_3}$ denote the expression obtained from $C_{3,<\lambda_2}$ by replacing each $<\lambda_2$ by $<\lambda_3$, while (as usual) keeping the $\eta_{\lambda_0}$, $\eta_{\lambda_1}$, and $\eta_{\lambda_2}$. Then we have by FTC,
\begin{equation}\label{C3FTC}
C_{3,<\lambda_2}
=C_{3,<\delta\lambda_2}
+\int_{\delta\lambda_2}^{\lambda_2}
\lambda_3\frac{d}{d\lambda_3}C_{3,<\lambda_3}
\frac{d\lambda_3}{\lambda_3}.
\end{equation}
After integration in $\lambda_0$, $\lambda_1$, and $\lambda_2$, the first term on the right-hand side gives precisely the terms in \eqref{G3def}, modulo perturbative errors. The differentiated term is a sum of quartic forms with comparable inputs and perturbative quintic errors. To extract the principal cubic part, we have the following lemma.
\begin{lemma}\label{quarticerrors}
We have
\begin{equation}\label{C3decomp}
\int_{\mathcal D}
\bigl(C_{3,<\delta\lambda_2}\bigr)_\mu\,d\Omega
=(G(\eta)\psi)^{[3]}_\mu
+(\mathcal R_\eta^{C,3})_\mu,
\end{equation}
where 
\begin{equation*}
\|(\mathcal R_\eta^{C,3})_\mu\overline u_\mu^{x_0}\|_{L_t^1L_x^1}
\lesssim\epsilon^4c_\mu^2\mu^{-2s_c}(1+\mu|x_0|)^{\frac12}.
\end{equation*}
\end{lemma}

\begin{proof}
Like the previous analysis, we retain $\mathbf P_0$ while extracting the principal parts of the operators in $C_{3,<\delta\lambda_2}$. Once the complete principal expressions have been obtained, their output is supported at frequencies $\ll M\lambda_0$, and $\mathbf P_0$ can be removed. Modulo errors which contribute to $(\mathcal R_\eta^{C,3})_\mu$, we have
\begin{equation*}
\begin{split}
C^1_{3,<\delta\lambda_2}&=
T_{\mathfrak g}\bigl(\eta_{\lambda_2}T_{\mathfrak l}(\eta_{\lambda_1}T_{\mathfrak l}\varphi_{\lambda_0})\bigr)
+\nabla\cdot\bigl(\eta_{\lambda_2}T_{\mathfrak m}(\eta_{\lambda_1}T_{\mathfrak l}\varphi_{\lambda_0})\bigr)
\\
&+T_{\mathfrak g}\Bigl(
\eta_{\lambda_1}T_{\mathfrak l}(\eta_{\lambda_2}T_{\mathfrak l}\varphi_{\lambda_0})
-\eta_{\lambda_1}\eta_{\lambda_2}T_{\mathfrak l}^2\varphi_{\lambda_0}
+\eta_{\lambda_1}T_{\mathfrak l}(\eta_{\lambda_0}T_{\mathfrak l}\varphi_{\lambda_2})
\Bigr)
\\
&+\nabla\cdot\Bigl(
\eta_{\lambda_1}T_{\mathfrak m}(\eta_{\lambda_2}T_{\mathfrak l}\varphi_{\lambda_0})
-\eta_{\lambda_1}\eta_{\lambda_2}T_{\mathfrak m}T_{\mathfrak l}\varphi_{\lambda_0}
+\eta_{\lambda_1}T_{\mathfrak m}(\eta_{\lambda_0}T_{\mathfrak l}\varphi_{\lambda_2})
\Bigr)
\end{split}
\end{equation*}
while $C^2_{3,<\delta\lambda_2}$ is obtained from this expression by interchanging $(\eta_{\lambda_0},\varphi_{\lambda_0})$ and $(\eta_{\lambda_1},\varphi_{\lambda_1})$. Finally (again modulo perturbative errors),
\begin{equation*}
\begin{split}
C^3_{3,<\delta\lambda_2}={}&
-T_{\mathfrak g}\bigl(\eta_{\lambda_0}\eta_{\lambda_1}T_{\mathfrak l}^2\varphi_{\lambda_2}\bigr)
-\nabla\cdot\bigl(\eta_{\lambda_0}\eta_{\lambda_1}T_{\mathfrak m}T_{\mathfrak l}\varphi_{\lambda_2}\bigr).
\end{split}
\end{equation*}
The errors produced in these replacements contain either a frequency-separated nonprincipal factor, in which case an unbalanced bilinear estimate is available, or a high-frequency tail, which is controlled by Proposition \ref{Bernsteinelliptic}. In each case, we can estimate the error directly in $L_t^1L_x^2$ and estimate $\overline u_{\mu}^{x_0}$ separately in $L_t^{\infty}L_x^2$. This gives the paired $L_{t,x}^1$ remainder estimate in the lemma.

Before extracting the principal cubic part, we commute the frequency projections through the principal operators. Thus,
\begin{equation*}
\begin{split}
T_{\mathfrak l}\varphi_\lambda
&=S_\lambda(T_{\mathfrak l}\varphi)+[T_{\mathfrak l},S_\lambda]\varphi,\\
T_{\mathfrak l}^2\varphi_\lambda
&=S_\lambda(T_{\mathfrak l}^2\varphi)+[T_{\mathfrak l}^2,S_\lambda]\varphi,\\
T_{\mathfrak m}T_{\mathfrak l}\varphi_\lambda
&=S_\lambda(T_{\mathfrak m}T_{\mathfrak l}\varphi)+[T_{\mathfrak m}T_{\mathfrak l},S_\lambda]\varphi.
\end{split}
\end{equation*}
Moreover,
\begin{equation*}
\begin{split}
[T_{\mathfrak l}^2,S_\lambda]
&=T_{\mathfrak l}[T_{\mathfrak l},S_\lambda]
+[T_{\mathfrak l},S_\lambda]T_{\mathfrak l},\\
[T_{\mathfrak m}T_{\mathfrak l},S_\lambda]
&=T_{\mathfrak m}[T_{\mathfrak l},S_\lambda]
+[T_{\mathfrak m},S_\lambda]T_{\mathfrak l}.
\end{split}
\end{equation*}
Each commutator contains a differentiated low-frequency
paracoefficient. Pairing this paracoefficient with the corresponding
higher-frequency input using an unbalanced bilinear estimate controls
the commutator in $L_{t,x}^2$. Multiplying by one of the remaining
factors in $L_t^4L_x^\infty$ gives an
$L_t^{\frac43}L_x^2$ block, and the last factor is estimated in
$L_t^4L_x^\infty$. Thus the full error is controlled directly in
$L_t^1L_x^2$, after which $\overline u_\mu^{x_0}$ is estimated separately in
$L_t^\infty L_x^2$. These terms therefore satisfy the paired remainder
estimate in the lemma.
\\

After these commutations, the six nested terms associated with each of $T_{\mathfrak g}$ and $\nabla\cdot$ in $C^1_{3,<\delta\lambda_2}$ and $C^2_{3,<\delta\lambda_2}$ combine to give, by \eqref{nestedPiconvention},
\begin{equation*}
T_{\mathfrak g}\Pi\bigl(\eta,T_{\mathfrak l}(\eta,T_{\mathfrak l}\varphi)\bigr)
+\nabla\cdot\Pi\bigl(\eta,T_{\mathfrak m}(\eta,T_{\mathfrak l}\varphi)\bigr).
\end{equation*}
For the remaining terms, the two $\eta$ entries are identical, so the symmetric trilinear paraproduct counts each placement of the last input twice. Consequently,
\begin{equation*}
\begin{split}
&\int_{\mathcal D}\Bigl(
\eta_{\lambda_1}\eta_{\lambda_2}F_{\lambda_0}
+\eta_{\lambda_0}\eta_{\lambda_2}F_{\lambda_1}
+\eta_{\lambda_0}\eta_{\lambda_1}F_{\lambda_2}
\Bigr)d\Omega
=\frac12\Pi(\eta,\eta,F).
\end{split}
\end{equation*}
Applying this identity first with $F=T_{\mathfrak l}^2\varphi$ and then with $F=T_{\mathfrak m}T_{\mathfrak l}\varphi$ gives the last two terms in \eqref{G3def}. This proves \eqref{C3decomp}.
\end{proof}

It remains to treat the differentiated term in \eqref{C3FTC}. The linearization formulas used above show that
\begin{equation*}
\lambda_3\frac{d}{d\lambda_3}C_{3,<\lambda_3}
\end{equation*}
is a sum of quartic terms and quintic errors. Indeed, differentiating an operator produces one additional factor $\eta_{\lambda_3}$, while differentiating $\Phi^{l,j}_{<\lambda_3}$ produces either a term containing
$\eta_{\lambda_j}L_{<\lambda_3}\Phi^l_{<\lambda_3}(\psi_{\lambda_3},\eta_{\lambda_3})$
or a term containing the two additional factors
$\eta_{\lambda_j}\eta_{\lambda_3}L_{<\lambda_3}B_{<\lambda_3}$. The latter is quintic. In the former terms, we replace $\Phi^l_{<\lambda_3}(\psi_{\lambda_3},\eta_{\lambda_3})$ by $\varphi_{\lambda_3}$ and extract the principal parts of the operators, retaining $\mathbf P_0$ until the output support of the full expression has been identified. Every surviving quartic term then has three explicit factors of $\eta$, one explicit factor of $\varphi$, and total differential order four. Moreover, by the integration ranges, we have
\begin{equation*}
\lambda_0\approx_\delta\lambda_1\approx_\delta\lambda_2\approx_\delta\lambda_3.
\end{equation*}
As in the treatment of $Q_{3,<\lambda_2}$, the separated output ranges are perturbative, and the only non-perturbative range is the range with output frequency $\mu\approx_{M,\delta}\lambda_0$, whose contribution we isolate as $N_{\eta,\mu}^{[4],C}$. The quintic terms are estimated using four $L_t^4L_x^\infty$ bounds and one $L_t^\infty L_x^2$ energy bound as in the treatment of $Q_{3,<\lambda_2}$. Thus,
\begin{equation}\label{N4Citer}
\begin{split}
\int_{\mathcal D}\int_{\delta\lambda_2}^{\lambda_2}
\left(\lambda_3\frac{d}{d\lambda_3}C_{3,<\lambda_3}\right)_\mu
\frac{d\lambda_3}{\lambda_3}\,d\Omega=N_{\eta,\mu}^{[4],C}+(\mathcal R_\eta^C)_\mu,
\end{split}
\end{equation}
where 
\begin{equation*}
\|(\mathcal R_\eta^C)_\mu\overline u_\mu^{x_0}\|_{L_t^1L_x^1}
\lesssim\epsilon^4c_\mu^2\mu^{-2s_c}(1+\mu|x_0|)^{\frac12}.
\end{equation*}
We then define $N_{\eta,\mu}^{[4]}\coloneq N_{\eta,\mu}^{[4],Q}+N_{\eta,\mu}^{[4],C}$ and set
\begin{equation*}
N_\eta^{[4]}\coloneq\sum_\mu N_{\eta,\mu}^{[4]}.
\end{equation*}
By construction, this is of the form described in Definition \ref{N4def}. Now combining \eqref{iter2decomp}, \eqref{T3CQdecomp}, \eqref{N4Qiter}, \eqref{C3FTC}, \eqref{C3decomp}, and \eqref{N4Citer}, together with the $L_t^1L_x^2$ estimates for $\mathcal E_{3,<\lambda_2}$ proves \eqref{iter3decomp} and the first remainder bound. We next verify \eqref{thirddecompaux}. The contribution of the remainder term $\mathcal{R}_{\eta}$ consists of the errors in 
\eqref{seconddecompaux} and the new errors above which are estimated directly in
$L_t^1L_x^2$. Arguing as in 
\ref{subsec:Neta2}, using Bernstein and two
$L_t^\infty L_x^2$ bounds for the factors of $u_{\mu}$ and $\overline{u}_{\mu}^{x_0}$, we have the bound
\begin{equation*}
\|(\mathcal{R}_{\eta})_{\lambda}\overline{u}_{\mu}^{x_0}u_{\mu}\|_{L_{t,x}^1}\lesssim \epsilon^4 c_{\mu}^2\lambda^{-1}\mu^{-2s_c}(1+\mu|x_0|)^{\frac12}.    
\end{equation*}
The quartic terms in
$N_\eta^{[4]}$ can also be
estimated in $L_t^1L_x^\infty$ by four
$L_t^4L_x^\infty$ bounds, and so, their contribution in \eqref{thirddecompaux} can be handled in the same way. Thus, it remains to estimate the principal cubic terms involved in the second term on the left hand side of \eqref{thirddecompaux}.
We use the expressions in the proof of \eqref{C3decomp},
and give the details for
\begin{equation*}
P_\lambda T_{\mathfrak g}
\bigl(\eta_{\lambda_2}T_{\mathfrak l}
(\eta_{\lambda_1}T_{\mathfrak l}\varphi_{\lambda_0})\bigr),
\end{equation*}
where $\delta\lambda_0\leq\lambda_1\leq\lambda_0$ and
$\delta\lambda_1\leq\lambda_2\leq\lambda_1$. The other cubic terms are handled in the same way. First, suppose that $\lambda_2\ll\mu$. In this case, we simply pair
$\eta_{\lambda_2}$ with $\overline u_\mu^{x_0}$ in
$L_{t,x}^2$, put $\eta_{\lambda_1}$ and
$\varphi_{\lambda_0}$ in $L_t^4L_x^\infty$, and put the
remaining $u_\mu$ in $L_t^\infty L_x^2$. Indeed, from the frequency supports of $\eta_{\lambda_1}$ and $\varphi_{\lambda_0}$ and the kernel bounds for $T_{\mathfrak{l}}$, we have
\begin{equation*}
\|T_{\mathfrak l}
(\eta_{\lambda_1}T_{\mathfrak l}\varphi_{\lambda_0})
\|_{L_t^2L_x^\infty}
\lesssim_\delta
\lambda_0^2
(\epsilon\lambda_0^{-\frac98})
(\epsilon\lambda_0^{-\frac{13}{8}})
=\epsilon^2\lambda_0^{-\frac34}.
\end{equation*}
The outer operator $T_{\mathfrak{g}}$ contributes a factor of $\lambda$, and its kernel is rapidly decaying at scale $\lambda^{-1}$. As in the preceding subsection,
the translations of $\eta_{\lambda_2}$ therefore cost a factor of
$(1+\lambda_2/\lambda)^{\frac12}$. Using the unbalanced bilinear
estimate and $\lambda_2\approx_\delta\lambda_0$, we obtain
\begin{equation*}
\begin{split}
\|u_\mu\overline u_\mu^{x_0}P_\lambda T_{\mathfrak g}
\bigl(\eta_{\lambda_2}T_{\mathfrak l}
(\eta_{\lambda_1}T_{\mathfrak l}\varphi_{\lambda_0})\bigr)
\|_{L_{t,x}^1}&\lesssim_\delta
\epsilon^5c_\mu^2\mu^{-2s_c}(1+\mu|x_0|)^{\frac12}
\lambda\left(1+\frac{\lambda_2}{\lambda}\right)^{\frac12}
\lambda_0^{-2}
\\
&\lesssim_\delta
\epsilon^5c_\mu^2\mu^{-2s_c}(1+\mu|x_0|)^{\frac12}
\lambda_0^{-1}.
\end{split}
\end{equation*}
Here we used $\lambda\lesssim_M\lambda_0$ and
$\lambda_2\leq\lambda_0$ in the last inequality.
Integrating in $\lambda_2$ and $\lambda_1$ over their
respective ranges, restricted to $\lambda_2\ll\mu$, and then
in $\lambda_0\gtrsim_M\lambda$, gives
\begin{equation*}
\epsilon^5c_\mu^2\mu^{-2s_c}(1+\mu|x_0|)^{\frac12}
\int_{\lambda_0\gtrsim_M\lambda}
\lambda_0^{-1}\frac{d\lambda_0}{\lambda_0}
\lesssim_\delta
\epsilon^5c_\mu^2\lambda^{-1}\mu^{-2s_c}(1+\mu|x_0|)^{\frac12}
\end{equation*}
which suffices. Now suppose that $\lambda_2\gtrsim\mu$. Since
$\lambda\ll\mu$, the output of
$T_{\mathfrak l}(\eta_{\lambda_1}
T_{\mathfrak l}\varphi_{\lambda_0})$
must have frequency comparable to $\lambda_2$.
Consequently, all three input frequencies and the output
of this quadratic factor are comparable to $\lambda_0$.
We can therefore use one internal bilinear estimate without
a large translation loss. There are now two cases to consider: either there is an unbalanced pair,
or, after subdividing the frequency supports, the low output frequency of the full cubic term
gives a pair of internal factors with frequency 
$|\xi_i\pm\xi_j|\gtrsim\lambda_0$, which allows us to apply Lemma \ref{balancedrestrictedoutput}. In either case, we can put the remaining unpaired input factor in $L_t^4L_x^\infty$. Using this reasoning, we obtain the bound
\begin{equation*}
\left\|P_\lambda T_{\mathfrak g}
\bigl(\eta_{\lambda_2}T_{\mathfrak l}
(\eta_{\lambda_1}T_{\mathfrak l}\varphi_{\lambda_0})\bigr)
\right\|_{L_t^{4/3}L_x^2}
\lesssim_\delta
\epsilon^3\lambda\lambda_0^{-\frac{23}{8}}.
\end{equation*}
We then put $\overline u_\mu^{x_0}$ in $L_t^4L_x^\infty$
and $u_\mu$ in $L_t^\infty L_x^2$.
Integration in $\lambda_2$ and $\lambda_1$ over their
defining ranges, restricted to $\lambda_2\gtrsim\mu$,
and then in $\lambda_0\gtrsim\mu$, gives
\begin{equation*}
\epsilon^5c_\mu^2\mu^{\frac78-2s_c}
\int_{\lambda_0\gtrsim\mu}
\lambda\lambda_0^{-\frac{23}{8}}
\frac{d\lambda_0}{\lambda_0}
\lesssim_\delta
\epsilon^5c_\mu^2\lambda\mu^{-2-2s_c}
\end{equation*}
which suffices.
The other principal cubic are handled similarly. Combining all these estimates proves \eqref{thirddecompaux}.

Finally, we prove \eqref{thirddecompL2}, which is quite straightforward. Since
$\lambda_0\approx_\delta\lambda_1\approx_\delta\lambda_2$,
we estimate the cubic terms in \eqref{T3CQdecomp} using two
$L_t^4L_x^\infty$ bounds for the $\eta$ factors and the
$L_t^\infty L_x^2$ bound for $\Phi_j$ along with Proposition \ref{basicopbounds}. For the quartic terms,
we put one $\eta$ factor in $L_t^\infty L_x^2$, another in
$L_{t,x}^\infty$, and the remaining $\eta$ factor and the
differentiated $B$ factor in $L_t^4L_x^\infty$, using Proposition \ref{VBunbil}. The high-frequency tails in
$\mathcal E_{3,<\lambda_2}$ satisfy the same bound by
Proposition \ref{Bernsteinelliptic}, as in \eqref{firstdecomperror}.
Thus,
\begin{equation*}
\|\mathcal T_{3,<\lambda_2}\|_{L_{t,x}^2}
\lesssim_\delta\epsilon^3c_{\lambda_0}
\lambda_0^{\frac14-s_c}.
\end{equation*}
The output projection restricts the integrand to $\lambda_0\gtrsim_M\mu$.
Integrating over the relevant frequency ranges and combining
\eqref{iter2decomp} with \eqref{seconddecompL2} proves
\eqref{thirddecompL2} and completes the proof.

\end{proof}


\subsection{The refined paradifferential equation for \texorpdfstring{$\eta$}{}}
\label{sec:N-eta}
We can now combine all of the expansions in this section to obtain the full multilinear expansion of $N_{\eta}$ that we will use in our analysis. The result below provides,
in a rigorous form, the layered paradifferential expansion of the $\eta$ equation, as derived in Theorem~\ref{t:expansions}:

\begin{theorem}\label{Netafinalcorollary}
Under the assumptions of Proposition \ref{thirddecomposition}, define
the quadratic balanced part
\begin{equation}\label{Neta2def}
N_{\eta}^{[2]}
\coloneq
(G(\eta)\psi)^{[2]}
=-T_{\mathfrak g}\bigl(\Pi(\eta,T_{\mathfrak l}\varphi)\bigr)
-\nabla\cdot\Pi(\eta,T_{\mathfrak m}\varphi),
\end{equation}
and the cubic balanced part
\begin{equation}\label{Neta3def}
\begin{split}
N_{\eta}^{[3]}
\coloneq
(G(\eta)\psi)^{[3]}
={}&
T_{\mathfrak g}\Pi\bigl(\eta,T_{\mathfrak l}(\eta,T_{\mathfrak l}\varphi)\bigr)
+\nabla\cdot\Pi\bigl(\eta,T_{\mathfrak m}(\eta,T_{\mathfrak l}\varphi)\bigr)
\\
&-\frac12T_{\mathfrak g}\Pi\bigl(\eta,\eta,T_{\mathfrak l}^2\varphi\bigr)
-\frac12\nabla\cdot\Pi\bigl(\eta,\eta,T_{\mathfrak m}T_{\mathfrak l}\varphi\bigr).
\end{split}
\end{equation}
Then the equation for $\eta$ takes the form
\begin{equation}\label{etafinalequation}
\partial_t\eta
-T_G\varphi+\nabla\cdot(T_{V^l}\eta)=N_{\eta},
\end{equation}
where $N_{\eta}$ has the multilinear paradifferential expansion
\begin{equation}\label{Netafinalexpansion}
N_{\eta}
=N_{\eta}^{[2]}+N_{\eta}^{[3]}
+N_{\eta}^{[4]}+\mathcal R_{\eta},
\end{equation}
for some $N_{\eta}^{[4]}$ as in Definition \ref{N4def}, and, for every
$\mu\in2^{\mathbb Z}$ and $x_0\in\mathbb R^2$, we have
\begin{equation*}
\|(\mathcal R_{\eta})_\mu
\overline u_\mu^{x_0}\|_{L_{t,x}^1}
\lesssim
\epsilon^4c_\mu^2\mu^{-2s_c}
(1+\mu|x_0|)^{\frac12}.
\end{equation*}
\end{theorem}
We also collect some auxiliary bounds for $N_{\eta}$ that will be needed for our density-flux corrections later.
\begin{corollary}\label{Netaauxbounds}
Under the assumptions and notation of Theorem
\ref{Netafinalcorollary}, we have the following bounds.
\begin{enumerate}
\item For every $\lambda>0$, there holds
\begin{equation}\label{NetaL2Linf}
\|(N_\eta)_\lambda\|_{L_t^2L_x^\infty}
\lesssim \epsilon^2\lambda^{-\frac34}.
\end{equation}
\item For every $0<\lambda\ll\mu$ and
$x_0\in\mathbb R^2$, we have a decomposition
\begin{equation}\label{Netaauxdecomp}
(N_\eta)_\lambda\overline u_\mu^{x_0}u_\mu
=\mathcal R_\eta+\nabla\cdot Q_\eta^{\mathrm{flux}},
\end{equation}
where for $\frac{4}{3}< p\leq 2$, we have
\begin{equation}\label{Netaauxtrilinear}
\|\mathcal R_\eta\|_{L_{t,x}^1}
+\lambda^{\frac{2}{p}-1}\|Q_\eta^{\mathrm{flux}}\|_{L_t^1L_x^p}
\lesssim
\epsilon^4c_\mu^2\lambda^{-1}
\mu^{-2s_c}(1+\mu|x_0|)^{\frac12}.
\end{equation}
\item For every $\lambda>0$, there holds
\begin{equation}\label{NetaL2cubic}
\|(N_\eta-N_\eta^{[2]})_\lambda\|_{L_{t,x}^2}
\lesssim_\delta\epsilon^3c_\lambda\lambda^{\frac14-s_c}.
\end{equation}
\end{enumerate}
\end{corollary}

\begin{proof}
The bound \eqref{NetaL2Linf} follows immediately from
\eqref{firstdecompaux} and \eqref{etafinalequation}.
For the second bound, let us write for simplicity
$H_\mu=\overline u_\mu^{x_0}u_\mu$.
By \eqref{thirddecompaux}, we already have
\begin{equation*}
(N_{\eta})_{\lambda}H_{\mu}
=(N_{\eta}^{[2]})_{\lambda}H_{\mu}
+O_{L_{t,x}^1}\bigl(
\epsilon^4c_\mu^2\lambda^{-1}
\mu^{-2s_c}(1+\mu|x_0|)^{\frac12}
\bigr).
\end{equation*}
It remains to consider the contribution of $N_\eta^{[2]}$.
This term is not directly perturbative, but is modulo a flux error.
We begin by splitting $(N_{\eta}^{[2]})_{\lambda}$ according to
whether its internal balanced frequencies are much smaller than
$\mu$ or comparable to or larger than $\mu$. That is, we write
\begin{equation*}
(N_\eta^{[2]})_\lambda
=F_\lambda^{\mathrm{lo}}+F_\lambda^{\mathrm{hi}},
\end{equation*}
according to whether the frequency $\lambda_0$ of one of the factors
in $N_{\eta,\lambda}^{[2]}$ satisfies $\lambda_0\ll\mu$ or
$\lambda_0\gtrsim\mu$, respectively. In both cases, the frequency
$\lambda_1$ of the second factor satisfies
$\delta\lambda_0\leq\lambda_1\leq\lambda_0$.

For $F_\lambda^{\mathrm{lo}}$, we can simply use two unbalanced
bilinear estimates; one where we pair $\overline u_{\mu}^{x_0}$
with the frequency $\lambda_0$ factor and one where we pair
$u_{\mu}$ with the frequency $\lambda_1$ factor. By
\eqref{uab-bi-unbal}, \eqref{phiu-bi-unbal}, and the kernel
bounds for $P_\lambda T_{\mathfrak g}$ and $T_{\mathfrak l}$,
we have
\begin{equation*}
\begin{split}
\left\|H_\mu P_\lambda T_{\mathfrak g}
\left(\eta_{\lambda_0}T_{\mathfrak l}\varphi_{\lambda_1}\right)
\right\|_{L_{t,x}^1}
&\lesssim_\delta
\epsilon^4c_\mu^2\mu^{-2s_c}(1+\mu|x_0|)^{\frac12}
\lambda\lambda_0^{-2}
\left(1+\frac{\lambda_0}{\lambda}\right).
\end{split}
\end{equation*}
Here the translations from the outer kernel have scale
$\lambda^{-1}$, and each of the two bilinear estimates costs
$\left(1+\frac{\lambda_0}{\lambda}\right)^{\frac12}$, since
$\lambda_1\approx\lambda_0$. This contributes the factor $(1+\frac{\lambda_0}{\lambda})$ on the right. The same estimate holds with $T_{\mathfrak g}$ replaced
by $\nabla\cdot$ and $\mathfrak l$ replaced by $\mathfrak m$ above,
and with the two input frequencies interchanged. Integrating gives
\begin{equation*}
\begin{split}
\|F_\lambda^{\mathrm{lo}}H_\mu\|_{L_{t,x}^1}
&\lesssim_\delta
\epsilon^4c_\mu^2\mu^{-2s_c}(1+\mu|x_0|)^{\frac12}
\int_{\lambda\lesssim_M\lambda_0\ll\mu}
\lambda\lambda_0^{-2}
\left(1+\frac{\lambda_0}{\lambda}\right)
\frac{d\lambda_0}{\lambda_0}
\\
&\lesssim_\delta
\epsilon^4c_\mu^2\lambda^{-1}\mu^{-2s_c}
(1+\mu|x_0|)^{\frac12},
\end{split}
\end{equation*}
which is the required bound.

For $F_\lambda^{\mathrm{hi}}$, we do not necessarily have two
separated pairs. Instead, we rely on two balanced bilinear estimates
to control this term, modulo a flux correction. First, since
$\lambda\ll\lambda_0$, we can argue as in the proof of Lemma
\ref{balancedrestrictedoutput}, with $s=s_c=2$, to obtain
\begin{equation*}
\begin{split}
\|F_\lambda^{\mathrm{hi}}\|_{L_{t,x}^2}\lesssim_\delta
\epsilon^2\lambda^{\frac12}
\int_{\lambda_0\gtrsim\mu}
\int_{\delta\lambda_0}^{\lambda_0}
\lambda_0^{-\frac94}
\frac{d\lambda_1}{\lambda_1}
\frac{d\lambda_0}{\lambda_0}\lesssim_\delta
\epsilon^2\lambda^{\frac12}\mu^{-\frac94}.
\end{split}
\end{equation*}
For every $\rho>0$, we also have the balanced estimate
\begin{equation*}
\|P_\rho H_\mu\|_{L_{t,x}^2}
\lesssim
\epsilon^2c_\mu^2\rho^{-\frac12}
\mu^{\frac54-2s_c}(1+\mu|x_0|)^{\frac12}.
\end{equation*}
Letting $\delta_0>0$ be a sufficiently small fixed constant,
frequency summation and Bernstein give
\begin{equation*}
\begin{split}
\|P_{\geq\delta_0\lambda}H_\mu\|_{L_{t,x}^2}
&\lesssim
\epsilon^2c_\mu^2\lambda^{-\frac12}
\mu^{\frac54-2s_c}(1+\mu|x_0|)^{\frac12},
\\
\lambda^{\frac{2}{p}-1}\|P_{<\delta_0\lambda}H_\mu\|_{L_t^2L_x^{\frac{2p}{2-p}}}
&\lesssim
\epsilon^2c_\mu^2\lambda^{\frac12}
\mu^{\frac54-2s_c}(1+\mu|x_0|)^{\frac12}.
\end{split}
\end{equation*}
In the second estimate, we used the hypothesis $p>\frac{4}{3}$. We then define
\begin{equation*}
Q_\eta^{\mathrm{flux}}
=\nabla\Delta^{-1}
\bigl(F_\lambda^{\mathrm{hi}}
P_{<\delta_0\lambda}H_\mu\bigr),
\end{equation*}
which is localized at frequency $\approx\lambda$.
Using the above estimates, we obtain
\begin{equation*}
\begin{split}
\|F_\lambda^{\mathrm{hi}}
P_{\geq\delta_0\lambda}H_\mu\|_{L_{t,x}^1}
+\lambda^{\frac{2}{p}-1}\|Q_\eta^{\mathrm{flux}}\|_{L_t^1L_x^p}&\lesssim
\|F_\lambda^{\mathrm{hi}}\|_{L_{t,x}^2}
\left(
\|P_{\geq\delta_0\lambda}H_\mu\|_{L_{t,x}^2}
+\lambda^{\frac{2}{p}-2}
\|P_{<\delta_0\lambda}H_\mu\|_{L_t^2L_x^{\frac{2p}{2-p}}}
\right)
\\
&\lesssim_\delta
\epsilon^4c_\mu^2\mu^{-1-2s_c}(1+\mu|x_0|)^{\frac12}.
\end{split}
\end{equation*}
Finally, setting
\begin{equation*}
\mathcal R_\eta
=(N_\eta-N_\eta^{[2]})_\lambda H_\mu
+F_\lambda^{\mathrm{lo}}H_\mu
+F_\lambda^{\mathrm{hi}}P_{\geq\delta_0\lambda}H_\mu
\end{equation*}
and combining the above estimates, we obtain
\eqref{Netaauxtrilinear}. Finally, \eqref{NetaL2cubic} follows immediately from
\eqref{thirddecompL2}, since
$N_\eta-N_\eta^{[2]}
=G(\eta)\psi-(G(\eta)\psi)^{[1]}-(G(\eta)\psi)^{[2]}$.
\end{proof}


\subsection{The refined expansion for the linearized interface}
\label{subsec:Neta-lin}
We now obtain the corresponding expansion for the linearized
equation for $\eta$ given by the first equation in \eqref{WW-lin2}. We recall that this is given by
\begin{equation}\label{etalinearized}
\partial_t h = G(\eta)(q - hB) - \nabla_x\cdot (hV).   
\end{equation}
The goal of the following proposition is to make precise the expansion of
$N_{\eta}^{\mathrm{lin}}$ in \eqref{N-expansion-lin}. Unlike in the nonlinear expansions above, we will estimate the cubic and higher terms perturbatively, using the additional
$\frac1{12}$ derivative of regularity for the background solution $u$. 

\begin{theorem}\label{Netalinearized}
Let $M>0$ and $\delta>0$ be fixed as in Subsection
\ref{paraconventions}, with $M\gg1$ and $M\delta\ll1$. Let $u$ and $w$ be as in Section \ref{sec:bootstrap}, satisfying the bootstrap bounds on $[0,T]$, where $0<T\leq1$, and with $c_\lambda\lesssim\lambda^{-\frac1{12}}$. Then the equation for $h$ takes the form
\begin{equation}\label{etalinearizedequation}
\partial_t h
-T_{\mathfrak g}r+\nabla\cdot(T_{V^l}h)
=N_{\eta}^{\mathrm{lin}},
\end{equation}
where $N_{\eta}^{\mathrm{lin}}$ has the paradifferential expansion
\begin{equation}\label{Netalinearizedexpansion}
N_{\eta}^{\mathrm{lin}}
=N_{\eta}^{\mathrm{lin},[0]}
+N_{\eta}^{\mathrm{lin},[1]}
+N_{\eta}^{\mathrm{lin},[2]}
+\mathcal R_{\eta}^{\mathrm{lin}}.
\end{equation} 
Above, we define the principal quadratic contribution when the linearized variable is at low frequency,
\begin{equation}\label{Netalin0def}
\begin{split}
N_{\eta}^{\mathrm{lin},[0]}
&\coloneq
-T_{\nabla h}\cdot T_{\mathfrak m}\varphi
-[T_{\mathfrak g},T_h]T_{\mathfrak l}\varphi
\\
&\hspace{5mm}-\nabla\cdot\bigl(T_{T_{\mathfrak m}r}\eta\bigr)
-T_{\mathfrak g}\bigl(T_{T_{\mathfrak l}r}\eta\bigr),
\end{split}
\end{equation}
The subprincipal
contribution with the linearized variable at high frequency is
\begin{equation}\label{Netalin1def}
N_{\eta}^{\mathrm{lin},[1]}
\coloneq
(T_G-T_{\mathfrak g})r,
\end{equation}
and the quadratic balanced part is defined by
\begin{equation}\label{Netalin2def}
\begin{split}
N_{\eta}^{\mathrm{lin},[2]}\coloneq
&-T_{\mathfrak g}\bigl(\Pi(h,T_{\mathfrak l}\varphi)\bigr)
-\nabla\cdot\Pi(h,T_{\mathfrak m}\varphi)
\\
&-T_{\mathfrak g}\bigl(\Pi(\eta,T_{\mathfrak l}r)\bigr)
-\nabla\cdot\Pi(\eta,T_{\mathfrak m}r).
\end{split}
\end{equation}
For every $\mu>0$ and $x_0\in\mathbb R^2$, the
remainder satisfies
\begin{equation}\label{Netalinearizedremainder}
\|(\mathcal R_{\eta}^{\mathrm{lin}})_\mu
\overline w_\mu^{x_0}\|_{L_{t,x}^1([0,T]\times\mathbb R^2)}
\lesssim
\epsilon^2d_\mu^2(1+\mu|x_0|)^{\frac12}.
\end{equation}
In addition, we have the following auxiliary bounds.
\begin{enumerate}
\item For every $\lambda>0$, there holds
\begin{equation}\label{NetalinL2remainder}
\|(N_{\eta}^{\mathrm{lin}}-N_{\eta}^{\mathrm{lin},[2]})_\lambda
\|_{L_{t,x}^2}
\lesssim_\delta
\epsilon d_\lambda\lambda^{\frac14}.
\end{equation}
\item For every $0<\nu\ll\mu$ and $x_0\in\mathbb R^2$, there holds
\begin{equation}\label{Netalinauxtrilinear}
\|(N_{\eta}^{\mathrm{lin}})_\nu
\overline w_\mu^{x_0}u_\mu\|_{L_{t,x}^1}
\lesssim_{\delta}
\epsilon^2d_\mu d_{\nu}\mu^{-\frac{5}{4}}\nu^{\frac{1}{4}}
(1+\mu|x_0|)^{\frac12}.
\end{equation}
\item For every $0<\lambda\ll\mu$ and $x_0\in\mathbb R^2$, there is a decomposition
\begin{equation}\label{Netawwauxdecomp}
(N_\eta)_\lambda w_\mu\overline w_\mu^{x_0}
=\mathcal R_{\eta,w}^{\mathrm{aux}}
+\nabla\cdot Q_{\eta,w}^{\mathrm{flux}},
\end{equation}
where, for $\frac43<p\leq2$, we have
\begin{equation}\label{Netawwauxtrilinear}
\|\mathcal R_{\eta,w}^{\mathrm{aux}}\|_{L_{t,x}^1}
+\lambda^{\frac2p-1}\|Q_{\eta,w}^{\mathrm{flux}}\|_{L_t^1L_x^p}
\lesssim_{\delta,p}
\epsilon^2d_\mu^2\lambda^{-1}(1+\mu|x_0|)^{\frac12}.
\end{equation}
\end{enumerate}
\end{theorem}
\begin{proof}
We begin by substituting $q=r+T_{B^l}h$ into \eqref{etalinearized}. We obtain
\begin{equation}\label{etalinearize2}
\partial_t h
=G(\eta)r-G(\eta)(Bh-T_{B^l}h)-\nabla\cdot(hV).
\end{equation}
We also record some basic estimates for the linearized variables $h$ and $r$,
\begin{equation}\label{hboot}
\begin{split}
\|h_\alpha\|_{L_t^\infty L_x^2}
+\alpha^{\frac12}\|r_\alpha\|_{L_t^\infty L_x^2}
&\lesssim d_\alpha,
\\
\|h_\alpha\|_{L_t^4L_x^p}
+\alpha^{\frac12}\|r_\alpha\|_{L_t^4L_x^p}
&\lesssim_p d_\alpha\alpha^{\frac78-\frac2p},
\qquad 8<p\leq\infty,
\\
\|h_\alpha\|_{L_t^{6,\infty}L_x^4}
+\alpha^{\frac12}\|r_\alpha\|_{L_t^{6,\infty}L_x^4}
&\lesssim d_\alpha\alpha^{\frac5{12}}.
\end{split}
\end{equation}
This follows from the fact that $h=\Re w$ and $T_{\mathfrak q}r=\Im w$, and the same
kernel estimates as in Proposition \ref{phiu-bootstrap}. The same argument also transfers the bilinear bounds in
Subsection \ref{subsec:para-boot} to $h$ and
$|D|^{\frac12}r$.
\\

Next, we expand the first term $G(\eta)r$ in the above equation. We apply the same successive linearizations as in the preceding
nonlinear expansions for $G(\eta)\psi$. The linearization formula for
\eqref{Glin} gives, for each $\alpha>0$, and $r_\alpha\coloneq S_{\alpha}r$,
\begin{equation}\label{Grexp}
\begin{split}
G(\eta)r_\alpha
=G_{<\delta\alpha}r_\alpha-\int_{\delta\alpha}^{\infty}
\left(
G_{<\beta}(\eta_\beta L_{<\beta}r_\alpha)
+\nabla\cdot(\eta_\beta M_{<\beta}r_\alpha)
\right)\frac{d\beta}{\beta}.
\end{split}
\end{equation}
After integration in $\alpha$ and estimating the high-frequency tails as in Proposition \ref{firstdecomposition}, the first term gives
$T_G r$, with the tail errors included in the
remainder below. We next apply the successive linearizations
from Subsections \ref{subsec:Neta2} and \ref{subsec:Neta34}
to the integral term in \eqref{Grexp}. As in
\eqref{GLdecomposition}, \eqref{G2trunc}, and
\eqref{T3CQdecomp}, we separate the errors from replacing
the elliptic operators by their principal parts into
balanced cubic and quartic contributions and the remaining
errors. This gives
\begin{equation}\label{linearizationetaexplanation}
\begin{split}
(G(\eta)r-T_Gr)_\mu\overline w_\mu^{x_0}
&=-\overline w_\mu^{x_0}P_\mu
\int_0^\infty\int_{\delta\alpha}^{\infty}
\Bigl(
T_{\mathfrak g}\bigl(\eta_\beta S_\alpha(T_{\mathfrak l}r)\bigr)
+\nabla\cdot\bigl(\eta_\beta S_\alpha(T_{\mathfrak m}r)\bigr)
\Bigr)
\frac{d\beta}{\beta}\frac{d\alpha}{\alpha}
\\
&\quad+
\bigl(N_{\eta}^{\mathrm{lin},[3]}
+N_{\eta}^{\mathrm{lin},[4]}\bigr)_\mu
\overline w_\mu^{x_0}
+O_{L_{t,x}^1}
\bigl(\epsilon^2d_\mu^2(1+\mu|x_0|)^{\frac12}\bigr).
\end{split}
\end{equation}
Above, and for the remainder of this subsection, $N_{\eta}^{\mathrm{lin},[3]}$ and
$N_{\eta}^{\mathrm{lin},[4]}$ will denote any balanced cubic and quartic contributions in the expansion of $N_{\eta}^{\mathrm{lin}}$. These terms serve as the counterparts of $N_{\eta}^{[3]}$ and $N_{\eta}^{[4]}$ in the nonlinear expansions in subsection \ref{subsec:Neta34}. We will not write out these terms explicitly here since it will not be necessary, but they arise from the same successive linearizations carried out in subsection \ref{subsec:Neta34}. For instance, in the expansion of $G(\eta)r$, $N^{[3],\mathrm{lin}}_{\eta}$ contains one factor of $r$ and two factors of $\eta$ all at comparable frequency, along with some bounded para-coefficients. Similarly, $N^{[4],\mathrm{lin}}_{\eta}$ contains terms with one factor of $r$ and three factors of $\eta$ all at comparable frequency.

The remaining perturbative $L_{t,x}^1$ errors in \eqref{linearizationetaexplanation} are
estimated similarly to the errors in subsections \ref{subsec:Neta2} and \ref{subsec:Neta34}, using instead \eqref{hboot} and the corresponding bilinear bounds for the linearized variables. We forego providing the detailed estimates as they are very similar to the analogous nonlinear bounds.

We now estimate the balanced cubic and quartic contributions
directly. This is in sharp contrast to the nonlinear estimates, where we had to identify the precise form of these terms. Here, we use the additional
$\frac1{12}$ derivative of regularity for the background solution $u$. To describe the treatment of these terms, we show the details for a representative cubic term, which is given by
\begin{equation*}
P_\mu T_{\mathfrak g}\bigl(
\eta_\alpha T_{\mathfrak l}
(\eta_\beta T_{\mathfrak l}r_\gamma)
\bigr),
\hspace{5mm}
\alpha\approx_\delta\beta\approx_\delta\gamma,
\quad \mu\lesssim_\delta\alpha.
\end{equation*}
Here, we can sacrifice a power of $T$ and use four $L_t^{6,\infty} L_x^4$ bounds to estimate this term directly. Indeed, the kernel bounds for $T_{\mathfrak{g}}$ and $T_{\mathfrak{l}}$ and H\"older gives
\begin{equation}\label{Netalincubicestimate}
\begin{split}
\left\|
P_\mu T_{\mathfrak g}\bigl(
\eta_\alpha T_{\mathfrak l}
(\eta_\beta T_{\mathfrak l}r_\gamma)
\bigr)\overline w_\mu^{x_0}
\right\|_{L_{t,x}^1}&\lesssim_\delta
T^{\frac13}\alpha^3
\|\eta_\alpha\|_{L_t^{6,\infty}L_x^4}
\|\eta_\beta\|_{L_t^{6,\infty}L_x^4}
\|r_\gamma\|_{L_t^{6,\infty}L_x^4}
\|w_\mu\|_{L_t^{6,\infty}L_x^4}
\\
&\lesssim_\delta
T^{\frac13}\alpha^3
\left(\epsilon c_\alpha\alpha^{-\frac{19}{12}}\right)^2
\left(d_\alpha\alpha^{-\frac1{12}}\right)
\left(d_\mu\mu^{\frac5{12}}\right)
\\
&=
T^{\frac13}\epsilon^2
\left(c_\alpha\alpha^{\frac1{12}}\right)^2
d_\alpha d_\mu
\left(\frac{\mu}{\alpha}\right)^{\frac5{12}}.
\end{split}
\end{equation}
Upon integration, using $T\leq1$,
$c_\alpha\lesssim\alpha^{-\frac1{12}}$, and slow variation
of $d_\alpha$, this gives
\begin{equation*}
\begin{split}
&\int_{\alpha\gtrsim_\delta\mu}
\int_{\beta\approx_\delta\alpha}
\int_{\gamma\approx_\delta\alpha}
\left\|
P_\mu T_{\mathfrak g}\bigl(
\eta_\alpha T_{\mathfrak l}
(\eta_\beta T_{\mathfrak l}r_\gamma)
\bigr)\overline w_\mu^{x_0}
\right\|_{L_{t,x}^1}
\frac{d\gamma}{\gamma}\frac{d\beta}{\beta}
\frac{d\alpha}{\alpha}
\\
&\qquad\lesssim_\delta
\epsilon^2d_\mu
\int_{\alpha\gtrsim_\delta\mu}
d_\alpha
\left(\frac{\mu}{\alpha}\right)^{\frac5{12}}
\frac{d\alpha}{\alpha}
\lesssim_\delta\epsilon^2d_\mu^2.
\end{split}
\end{equation*}
The other terms in $N_{\eta}^{\mathrm{lin},[3]}$ have a similar structure, and are estimated in the same way.

For the quartic terms, we place the additional factor $\eta$
factor in $L_{t,x}^\infty$, and then use the same cubic estimate above for the remaining balanced components. For example, when
$\alpha\approx_\delta\beta\approx_\delta\gamma
\approx_\delta\nu$, we have the representative bound
\begin{equation}\label{Netalinquarticestimate}
\begin{split}
\left\|
P_\mu T_{\mathfrak g}\bigl(
\eta_\alpha T_{\mathfrak l}
(\eta_\beta\eta_\gamma T_{\mathfrak l}^2r_\nu)
\bigr)\overline w_\mu^{x_0}
\right\|_{L_{t,x}^1}&\lesssim_\delta
T^{\frac13}\alpha^4
\|\eta_\alpha\|_{L_t^{6,\infty}L_x^4}
\|\eta_\beta\|_{L_t^{6,\infty}L_x^4}
\|\eta_\gamma\|_{L_{t,x}^\infty}
\|r_\nu\|_{L_t^{6,\infty}L_x^4}
\|w_\mu\|_{L_t^{6,\infty}L_x^4}
\\
&\lesssim_\delta
T^{\frac13}\epsilon^3c_\alpha
\left(c_\alpha\alpha^{\frac1{12}}\right)^2
d_\alpha d_\mu
\left(\frac{\mu}{\alpha}\right)^{\frac5{12}}.
\end{split}
\end{equation}
Using $T\leq1$, boundedness of $c_\alpha$, and the
hypothesis $c_\alpha\lesssim\alpha^{-\frac1{12}}$,
we can now integrate these estimates over the comparable
input frequencies parameters. Slow variation of $d_\alpha$ gives
\begin{equation*}
\begin{split}
\|(N_{\eta}^{\mathrm{lin},[3]})_\mu
\overline w_\mu^{x_0}\|_{L_{t,x}^1}
+\|(N_{\eta}^{\mathrm{lin},[4]})_\mu
\overline w_\mu^{x_0}\|_{L_{t,x}^1}\lesssim_\delta
\epsilon^2d_\mu
\int_{\alpha\gtrsim_\delta\mu}
d_\alpha
\left(\frac{\mu}{\alpha}\right)^{\frac5{12}}
\frac{d\alpha}{\alpha}
\lesssim_\delta\epsilon^2d_\mu^2.
\end{split}
\end{equation*}
Consequently, the balanced cubic and quartic terms can
also be estimated as a perturbative $O_{L_{t,x}^1}(\epsilon^2 d_{\mu}^2(1+\mu|x_0|)^{\frac12})$ error.

We can now split the remaining quadratic integral in \eqref{linearizationetaexplanation} 
to the relative sizes of $\alpha$ and $\beta$. The region
$\alpha<\delta\beta$, where the linearized input is at
low frequency, gives
\begin{equation*}
-T_{\mathfrak g}\bigl(T_{T_{\mathfrak l}r}\eta\bigr)
-\nabla\cdot\bigl(T_{T_{\mathfrak m}r}\eta\bigr),
\end{equation*}
while the balanced region 
$\delta\beta\leq\alpha\leq\delta^{-1}\beta$ gives
\begin{equation*}
-T_{\mathfrak g}\bigl(\Pi(\eta,T_{\mathfrak l}r)\bigr)
-\nabla\cdot\Pi(\eta,T_{\mathfrak m}r).
\end{equation*}
These are precisely the contributions containing
$r$ in \eqref{Netalin0def} and \eqref{Netalin2def}.

We next consider the terms containing $h$ in \eqref{etalinearize2}. There is an
additional cancellation to retain when $h$ is at low frequency. Using the basic identity
\begin{equation*}
G(\eta)B+\nabla\cdot V=0,
\end{equation*}
we obtain
\begin{equation*}
-G(\eta)(T_hB)-\nabla\cdot(T_hV)
=-[G(\eta),T_h]B-T_{\nabla h}\cdot V.
\end{equation*}
The commutator structure above ensures that $h$ is "differentiated" when at low frequency. Hence, we can write
\begin{equation*}
\begin{split}
&-G(\eta)(Bh-T_{B^l}h)-\nabla\cdot (hV)+\nabla\cdot (T_{V^l}h) 
\\
&=-G(\eta)(T_B-T_{B^l})h-G(\eta)(T_hB)-G(\eta)\Pi(h,B)-\nabla\cdot(hV)+\nabla\cdot (T_{V^l}h)  
\\
&=-G(\eta)(T_B-T_{B^l})h-[G(\eta),T_h]B-T_{\nabla h}\cdot V-G(\eta)\Pi(h,B)
\\
&\hspace{4mm}-\nabla\cdot\bigl((T_V-T_{V^l})h\bigr)-\nabla\cdot\Pi(V,h)
\end{split}
\end{equation*}
Using the above formula, and carrying out successive linearizations as in subsections \ref{subsec:Neta1}, \ref{subsec:Neta2} and \ref{subsec:Neta34} gives the decomposition
\begin{equation*}
\begin{split}
&\left(
-G(\eta)(Bh-T_{B^l}h)-\nabla\cdot(hV)
+\nabla\cdot(T_{V^l}h)
\right)_\mu\overline w_\mu^{x_0}
\\
&=\Bigl(
-[T_{\mathfrak g},T_h]T_{\mathfrak l}\varphi
-T_{\nabla h}\cdot T_{\mathfrak m}\varphi
-T_{\mathfrak g}\bigl(\Pi(h,T_{\mathfrak l}\varphi)\bigr)
-\nabla\cdot\Pi(h,T_{\mathfrak m}\varphi)+N_{\eta}^{\mathrm{lin},[3]}
+N_{\eta}^{\mathrm{lin},[4]}
\Bigr)_\mu\overline w_\mu^{x_0}
\\
&+O_{L_{t,x}^1}
\bigl(\epsilon^2d_\mu^2(1+\mu|x_0|)^{\frac12}\bigr).
\end{split}
\end{equation*}
where $N^{\mathrm{lin},[3]}_{\eta}$ and $N_{\eta}^{\mathrm{lin},[4]}$ have the same form as described above. The four quadratic terms are precisely the contributions
containing $h$ in \eqref{Netalin0def} and \eqref{Netalin2def}.

The errors are estimated similarly to those in Subsections
\ref{subsec:Neta1}, \ref{subsec:Neta2}, and
\ref{subsec:Neta34}, using the linearized bootstrap bounds
\eqref{hboot} and the corresponding mixed bilinear estimates in place of the nonlinear bounds.
The balanced cubic and quartic contributions are estimated
as in \eqref{Netalincubicestimate} and
\eqref{Netalinquarticestimate}. We omit the lengthy
calculations, since the same operator decompositions and
kernel estimates apply as with these substitutions.

Combining this with the expansion of $G(\eta)r$ above then gives
\begin{equation*}
\partial_t h-T_G r+\nabla\cdot(T_{V^l}h)
=N_\eta^{\mathrm{lin},[0]}
+N_\eta^{\mathrm{lin},[2]}
+\mathcal R_\eta^{\mathrm{lin}},
\end{equation*}
where
\begin{equation*}
\|(\mathcal R_\eta^{\mathrm{lin}})_\mu
\overline w_\mu^{x_0}\|_{L_{t,x}^1}
\lesssim
\epsilon^2d_\mu^2(1+\mu|x_0|)^{\frac12}.
\end{equation*}
Finally, writing
$T_G r=T_{\mathfrak g}r+(T_G-T_{\mathfrak g})r$
gives \eqref{etalinearizedequation},
\eqref{Netalinearizedexpansion}, and
\eqref{Netalinearizedremainder}.

It remains to prove the auxiliary bounds. We begin with
\eqref{NetalinL2remainder}. We first dispense with the contributions of $(N_{\eta}^{\mathrm{lin},[0]})_{\nu}$ and $(N_{\eta}^{\mathrm{lin},[1]})_{\nu}$. For $(N_{\eta}^{\mathrm{lin},[0]})_{\nu}$, we can simply one unbalanced bilinear $L_{t,x}^2$ estimate pairing the differentiated low-frequency linearized factor with the high frequency nonlinear factor (using Proposition \ref{Leibnizcommutator} for the commutator). For the contribution of $(N^{\mathrm{lin},[1]}_{\eta})_{\nu}$ we can use Proposition \ref{paraexp} with one unbalanced bilinear estimate for the main subprincipal term. In summary, we have obtain the bounds
\begin{equation*}
\begin{split}
\|(N_{\eta}^{\mathrm{lin},[0]})_\nu\|_{L_{t,x}^2}
&\lesssim_\delta
\epsilon c_\nu
\int_{\beta\ll\nu}
d_\beta\left(
\nu^{-1}\beta^{\frac54}
+\nu^{-\frac32}\beta^{\frac74}
\right)\frac{d\beta}{\beta}
\lesssim_\delta\epsilon d_\nu\nu^{\frac14},
\\
\|(N_{\eta}^{\mathrm{lin},[1]})_\nu\|_{L_{t,x}^2}
&\lesssim_\delta
\epsilon d_\nu\nu^{-\frac12}
\int_{\beta\ll\nu}
c_\beta\beta^{\frac34}\frac{d\beta}{\beta}
\lesssim_\delta\epsilon d_\nu\nu^{\frac14}.
\end{split}
\end{equation*}
For the cubic and higher terms, we use the decompositions
in the proofs of \eqref{firstdecompL2R},
\eqref{seconddecompL2}, and \eqref{thirddecompL2},
with the linearized bounds \eqref{hboot}. Here it is
important to retain the output frequency $\nu$ when
the input frequencies are much larger than $\nu$.
For example, the cubic term in
\eqref{Netalincubicestimate} satisfies
\begin{equation*}
\begin{split}
\left\|
P_\nu T_{\mathfrak g}\bigl(
\eta_\alpha T_{\mathfrak l}
(\eta_\beta T_{\mathfrak l}r_\gamma)
\bigr)
\right\|_{L_{t,x}^2}\lesssim_\delta
\nu\alpha^2
\|\eta_\alpha\|_{L_t^4L_x^\infty}
\|\eta_\beta\|_{L_t^4L_x^\infty}
\|r_\gamma\|_{L_t^\infty L_x^2}
\lesssim_\delta
\epsilon^2d_\alpha\nu\alpha^{-\frac34}
\end{split}
\end{equation*}
where $\alpha\approx_\delta\beta\approx_\delta\gamma$.
The corresponding terms containing $h$ have the same
bound, placing $h$ in $L_t^\infty L_x^2$ and the
background $\eta$ and $\varphi$ factors in
$L_t^4L_x^\infty$. The quartic terms are estimated
by placing the additional $\eta$ factor
in $L_{t,x}^\infty$.

The remaining cubic and quartic errors coming replacing elliptic operators by their principal parts are
estimated using an unbalanced bilinear estimate and placing
the remaining factors in $L_{t,x}^\infty$.
The quintic and higher order errors and high-frequency tails are
treated as in the proofs of \eqref{seconddecompL2} and
\eqref{thirddecompL2}, using the linearized bootstrap
bounds \eqref{hboot} in place of the nonlinear ones.
Consequently, after integration of the above bounds in the relevant input frequency variables, we obtain
\begin{equation*}
\begin{split}
\|(N_{\eta}^{\mathrm{lin}}-N_{\eta}^{\mathrm{lin},[2]})_\nu
\|_{L_{t,x}^2}
\lesssim_\delta
\epsilon d_\nu\nu^{\frac14}
+\epsilon^2\nu
\int_{\alpha\gtrsim_\delta\nu}
d_\alpha\alpha^{-\frac34}\frac{d\alpha}{\alpha}\lesssim_\delta
\epsilon d_\nu\nu^{\frac14}.
\end{split}
\end{equation*}
This proves \eqref{NetalinL2remainder}.

It remains to prove \eqref{Netalinauxtrilinear} and \eqref{Netawwauxtrilinear}; we begin with the former. Thanks to the additional $\frac{1}{12}$ regularity of $u$, this is very straightforward compared to the nonlinear counterpart. We also note that here, no flux correction is needed. Indeed, for the
contribution just estimated above, H\"older and
\eqref{NetalinL2remainder} immediately gives
\begin{equation*}
\begin{split}
\|(N_{\eta}^{\mathrm{lin}}-N_{\eta}^{\mathrm{lin},[2]})_\nu
\overline w_\mu^{x_0}u_\mu\|_{L_{t,x}^1}&\lesssim_\delta
T^{\frac16}\epsilon d_\nu\nu^{\frac14}
\|w_\mu\|_{L_t^{6,\infty}L_x^4}
\|u_\mu\|_{L_t^{6,\infty}L_x^4}\lesssim_\delta
\epsilon^2d_\nu d_\mu\nu^{\frac14}\mu^{-\frac54},
\end{split}
\end{equation*}
where we used $c_\mu\lesssim\mu^{-\frac1{12}}$.

For $N_{\eta}^{\mathrm{lin},[2]}$, we can simply use four
$L_t^6L_x^4$ bounds. For example, when
$\alpha\approx_\delta\beta\gtrsim_\delta\nu$, the
first term in \eqref{Netalin2def} satisfies
\begin{equation*}
\begin{split}
\left\|
P_\nu T_{\mathfrak g}
\bigl(h_\alpha T_{\mathfrak l}\varphi_\beta\bigr)
\overline w_\mu^{x_0}u_\mu
\right\|_{L_{t,x}^1}&\lesssim_\delta
T^{\frac13}\nu\alpha
\|h_\alpha\|_{L_t^{6,\infty}L_x^4}
\|\varphi_\beta\|_{L_t^{6,\infty}L_x^4}
\|w_\mu\|_{L_t^{6,\infty}L_x^4}
\|u_\mu\|_{L_t^{6,\infty}L_x^4}
\\
&\lesssim_\delta
T^{\frac13}\epsilon^2c_\alpha c_\mu
d_\alpha d_\mu
\nu\alpha^{-\frac23}\mu^{-\frac76}
\\
&\lesssim_\delta
\epsilon^2d_\alpha d_\mu
\nu\alpha^{-\frac34}\mu^{-\frac54}.
\end{split}
\end{equation*}
Here we used the additional $\frac1{12}$ derivative
of regularity on both background factors. The term in $N^{\mathrm{lin},[2]}_{\eta}$
containing $\eta_\alpha T_{\mathfrak l}r_\beta$
has the same bound. The two divergence terms are also estimated identically. Integrating in $\alpha$
and $\beta$, we obtain
\begin{equation*}
\begin{split}
\|(N_{\eta}^{\mathrm{lin},[2]})_\nu
\overline w_\mu^{x_0}u_\mu\|_{L_{t,x}^1}
&\lesssim_\delta
\epsilon^2d_\mu\nu\mu^{-\frac54}
\int_{\alpha\gtrsim_\delta\nu}
d_\alpha\alpha^{-\frac34}\frac{d\alpha}{\alpha}
\\
&\lesssim_\delta
\epsilon^2d_\nu d_\mu\nu^{\frac14}\mu^{-\frac54}
\end{split}
\end{equation*}
which proves \eqref{Netalinauxtrilinear}.

Finally, we prove \eqref{Netawwauxdecomp} and \eqref{Netawwauxtrilinear}.
We use the splitting
$(N_\eta^{[2]})_\lambda=F_\lambda^{\mathrm{lo}}+F_\lambda^{\mathrm{hi}}$
from the proof of \eqref{Netaauxtrilinear}, now with
$H_\mu=w_\mu\overline w_\mu^{x_0}$. For the contribution of
$N_\eta-N_\eta^{[2]}$, we repeat the proof of \eqref{thirddecompaux},
using the $L_t^{\infty}L_x^2$, $L_t^4L_x^\infty$, and mixed bilinear bounds
from Subsection \ref{subsec:para-boot} for the two $w_\mu$ factors.
For $F_\lambda^{\mathrm{lo}}$, we pair each $w_\mu$ factor with
one of the two lower-frequency background inputs using
\eqref{wu-bi-unbal-trans-boot} and the kernel bounds for
$T_{\mathfrak l}$ and $T_{\mathfrak m}$, as in the proof of
\eqref{Netaauxtrilinear}. Thus, we have
\begin{equation*}
\|(N_\eta-N_\eta^{[2]})_\lambda H_\mu\|_{L_{t,x}^1}
+\|F_\lambda^{\mathrm{lo}}H_\mu\|_{L_{t,x}^1}
\lesssim_\delta \epsilon^2d_\mu^2\lambda^{-1}(1+\mu|x_0|)^{\frac12}.
\end{equation*}
For $F_\lambda^{\mathrm{hi}}$, we reuse the bound
$\|F_\lambda^{\mathrm{hi}}\|_{L_{t,x}^2}
\lesssim_\delta\epsilon^2\lambda^{\frac12}\mu^{-\frac94}$
from the same proof. By \eqref{ww-bi-bal-boot}, we also have for each $\rho>0$,
\begin{equation*}
\|P_\rho H_\mu\|_{L_{t,x}^2}
\lesssim d_\mu^2\rho^{-\frac12}\mu^{\frac54}(1+\mu|x_0|)^{\frac12}.
\end{equation*}
Taking $Q_{\eta,w}^{\mathrm{flux}}
=\nabla\Delta^{-1}(F_\lambda^{\mathrm{hi}}P_{<\delta_0\lambda}H_\mu)$
with $\delta_0$ as in the proof of \eqref{Netaauxtrilinear},
the same argument gives
\begin{equation*}
\|F_\lambda^{\mathrm{hi}}P_{\geq\delta_0\lambda}H_\mu\|_{L_{t,x}^1}
+\lambda^{\frac2p-1}\|Q_{\eta,w}^{\mathrm{flux}}\|_{L_t^1L_x^p}
\lesssim_{\delta,p}\epsilon^2d_\mu^2\mu^{-1}(1+\mu|x_0|)^{\frac12}.
\end{equation*}
Since $\lambda\ll\mu$, combining this with the above bounds for $F^{\mathrm{lo}}_{\lambda}H_{\mu}$ and $(N_{\eta}-N_{\eta}^{[2]})_{\lambda}H_{\mu}$ proves \eqref{Netawwauxtrilinear}.
\end{proof}


\section{Paradifferential expansion of \texorpdfstring{$N_{\varphi}$}{}}
\label{sec:Nphi}

\subsection{Initial paradifferential equation for \texorpdfstring{$\varphi$}{}}\label{subsec:initialphi}

The purpose of this subsection is to write the equation for the good variable $\varphi\coloneq\psi-T_{B^l}\eta$ in a suitably precise preliminary paradifferential form. In the next subsection, we will then carry out a multilinear expansion of the "main" remaining nonlinearity as previously done in the equation for $\eta$. Unlike in the $\eta$ equation, where we simply started with $\partial_t\eta=G(\eta)\psi$, the initial setup here is substantially more involved. Our preliminary expansion is summarized by the following proposition.
\begin{proposition}
\label{paravarphiinitial}
Fix $0<\delta\ll1$, with $M\delta\ll1$, as in the paraproduct conventions outlined in Subsection \ref{multpara}. For each $\lambda>0$, let $\phi_{<\delta\lambda}$ be the harmonic extension of $\psi_{<\delta\lambda}$ into the domain below the graph $z=\eta_{<\delta\lambda}(t,x)$. Let $p_{<\delta\lambda}$ be the corresponding pressure, namely
\begin{equation*}
-\Delta_{x,z}p_{<\delta\lambda}
=
\left|\nabla^2_{x,z}\phi_{<\delta\lambda}\right|^2,
\qquad
p_{<\delta\lambda}\big|_{z=\eta_{<\delta\lambda}}=0,
\end{equation*}
with the normalization
\begin{equation*}
p_{<\delta\lambda}(t,x,z)+z\longrightarrow0
\qquad\text{as }z\longrightarrow-\infty.
\end{equation*}
Define the truncated Taylor coefficient $a_{<\delta\lambda}$ by
\begin{equation*}
a_{<\delta\lambda}
=
-\left.\partial_zp_{<\delta\lambda}\right|_{z=\eta_{<\delta\lambda}}.
\end{equation*}
Then the good variable $\varphi\coloneq \psi-T_{B^l}\eta$ satisfies the equation
\begin{equation}
\label{eq:geometric-truncated-good-unknown}
\left(\partial_t+T_{V^l}\cdot\nabla\right)\varphi
+T_{a^l}\eta
=N_{\varphi}+N_{\varphi}^{\mathrm{flux}}+\mathcal R,
\end{equation}
where
\begin{equation}\label{Nphidef}
\begin{split}
N_{\varphi}
&=
\left(\frac12 B^2-T_{B^l}B\right)
-
\left(\frac12 |V|^2-T_{V^l}\cdot V\right)
\\
&\quad-
\left(
BV\cdot\nabla\eta
-T_{B^l}(V\cdot\nabla\eta)
-T_{V^l}\cdot(B\nabla\eta)
+T_{V^l}\cdot T_{B^l}\nabla\eta
\right).
\end{split}
\end{equation}
Moreover, $\mathcal R$ is perturbative in the sense that
\begin{equation}\label{paraphiremainder}
\|\mathcal R_\mu\overline u_\mu^{x_0}\|_{L_{t,x}^1}
\lesssim
c_\mu^2\epsilon^4\mu^{-\frac12-2s_c}(1+\mu|x_0|)^{\frac12},
\end{equation}
and $N_{\varphi,\mu}^{\mathrm{flux}}$ contributes a flux term. That is,
\begin{equation*}
\overline u_\mu^{x_0}N_{\varphi,\mu}^{\mathrm{flux}}
=\nabla_x\cdot Q_{\varphi,\mu},
\end{equation*}
where for $\frac{4}{3}<p\leq 2$, there holds
\begin{equation}\label{QPhibound}
\mu^{\frac{2}{p}-1}\|Q_{\varphi,\mu}\|_{L_t^1L_x^p}
\lesssim_p
\epsilon^4c_\mu^2\mu^{-\frac12-2s_c}(1+\mu|x_0|)^{\frac12}.
\end{equation}
Moreover, the full nonlinear error term
\begin{equation}\label{Nphierrinitialdef}
N_\varphi^{\mathrm{err}}
\coloneq
\left(\partial_t+T_{V^l}\cdot\nabla\right)\varphi
+T_{a^l}\eta
=N_\varphi+N_\varphi^{\mathrm{flux}}+\mathcal R.
\end{equation}
satisfies the following three auxiliary bounds
\begin{enumerate}
\item For every $\lambda>0$, we have
\begin{equation}\label{NphiinitialL2Linf}
\|(N_\varphi^{\mathrm{err}}-N_{\varphi})_\lambda\|_{L_t^2L_x^\infty}
\lesssim\epsilon^2\lambda^{-\frac54}.
\end{equation}
\item For every $0<\lambda\ll\mu$ and $x_0\in\mathbb R^2$,
\begin{equation}\label{NphiinitialauxR}
\bigl\|(N_\varphi^{\mathrm{err}}-N_\varphi)_\lambda
\overline u_\mu^{x_0}u_\mu\bigr\|_{L_{t,x}^1}
\lesssim
\epsilon^5c_\mu^2\lambda^{-\frac32}
\mu^{-2s_c}(1+\mu|x_0|)^{\frac12}.
\end{equation}
\item For every $\lambda>0$, we have
\begin{equation}\label{NphiinitialL2}
\|(N_\varphi^{\mathrm{err}}-N_\varphi)_\lambda\|_{L_{t,x}^2}
\lesssim_\delta\epsilon^3c_\lambda\lambda^{-\frac14-s_c}.
\end{equation}
\end{enumerate}
\end{proposition}
The nonlinear term $N_{\varphi}$ above will be the basis for our multilinear expansions in the next subsection. The four error bounds \eqref{paraphiremainder}, \eqref{NphiinitialL2Linf}, \eqref{NphiinitialauxR} and \eqref{NphiinitialL2} should be thought of as analogues of the corresponding four error bounds in Proposition \ref{firstdecomposition}. Moreover, the flux term $N_{\varphi,\mu}^{\mathrm{flux}}$ will arise in a similar fashion to the one in Corollary \ref{Netaauxbounds} for the paradifferential expansion of $\eta$. The main difference here is that it arises in the "principal" nonlinear error, rather than the auxiliary bounds. Nevertheless, $N_{\varphi,\mu}^{\mathrm{flux}}$ will still play a perturbative role in the interaction Morawetz estimates established later.
\begin{proof}
Below, we will write $\mathcal R$ to generically denote any remainder of the form described above. As in earlier sections, the exact form of $\mathcal R$ will be allowed to change from line to line. 
\\

It will be convenient to temporarily introduce the "dynamic" truncated Taylor coefficient
\begin{equation}\label{tildeadef}
\widetilde a_{<\delta\lambda}
\coloneq
1+(\partial_t+V_{<\delta\lambda}\cdot\nabla)B_{<\delta\lambda}.
\end{equation}
This definition is motivated by the exact non-truncated equation for $B$,
\begin{equation*}
a=1+(\partial_t+V\cdot\nabla)B.
\end{equation*}
We will first write an equation for $\varphi$ in terms of $\widetilde a^l_{<\delta\lambda}$ (defined by applying $P_{<M\delta\lambda}$ to $\widetilde{a}_{<\delta\lambda}$), where the algebra is simpler. Afterwards, we will estimate the error between $\widetilde a^l_{<\delta\lambda}$ and the true truncated Taylor coefficient $a^l_{<\delta\lambda}$. First, we compute
\begin{equation*}
\partial_t(T_{B^l}\eta)
=T_{\partial_tB^l}\eta+T_{B^l}\partial_t\eta,
\qquad
\nabla(T_{B^l}\eta)
=T_{\nabla B^l}\eta+T_{B^l}\nabla\eta.
\end{equation*}
Then, using \eqref{tildeadef}, we obtain
\begin{equation*}
\begin{split}
\left(\partial_t+T_{V^l}\cdot\nabla\right)
\left(\psi-T_{B^l}\eta\right)
+T_{\widetilde a^l}\eta
={}&
\partial_t\psi+\eta
+T_{V^l}\cdot\nabla\psi
-T_{B^l}\partial_t\eta
\\
&+T_{(V\cdot\nabla B)^l}\eta
-T_{V^l}\cdot T_{\nabla B^l}\eta
-T_{V^l}\cdot T_{B^l}\nabla\eta.
\end{split}
\end{equation*}
Here, the paradifferential operator $T_{(V\cdot\nabla B)^l}$ is defined similarly to $T_{B^l}$ by
\begin{equation*}
\int_{0}^{\infty}P_{<M\delta\lambda}
\left(V_{<\delta\lambda}\cdot\nabla B_{<\delta\lambda}\right)S_{\lambda}\frac{d\lambda}{\lambda}.
\end{equation*}
We next use
\begin{equation*}
\partial_t\eta
=G(\eta)\psi
=B-V\cdot\nabla\eta,
\end{equation*}
and
\begin{equation*}
\partial_t\psi
=-\eta+\frac12B^2-\frac12|V|^2-BV\cdot\nabla\eta.
\end{equation*}
Substituting these, alongside the identity $\nabla\psi=V+B\nabla\eta$ gives
\begin{equation}\label{initialphipara}
\left(\partial_t+T_{V^l}\cdot\nabla\right)\varphi
+T_{a^l}\eta
=N_{\varphi}
+\left(
T_{(V\cdot\nabla B)^l}
-T_{V^l}\cdot T_{\nabla B^l}
\right)\eta
+T_{a^l-\widetilde a^l}\eta,
\end{equation}
where $N_{\varphi}$ is given by \eqref{Nphidef}. Our aim now is to show that the latter two terms on the right are perturbative, apart from the flux term in the proposition. We begin with the paracommutator term.

\begin{lemma}\label{paracomperturbative}
We have
\begin{equation*}
\left(
T_{(V\cdot\nabla B)^l}
-T_{V^l}\cdot T_{\nabla B^l}
\right)\eta
=\mathcal R,
\end{equation*}
where $\mathcal R$ satisfies the estimate \eqref{paraphiremainder}.
\end{lemma}

\begin{proof}
We have
\begin{equation*}
\begin{split}
\left(
T_{(V\cdot\nabla B)^l}
-T_{V^l}\cdot T_{\nabla B^l}
\right)\eta=
\int_0^\infty
\Bigl[
P_{<M\delta\lambda}
\left(V_{<\delta\lambda}\cdot\nabla B_{<\delta\lambda}\right)
\eta_\lambda
-T_{V^l}\cdot
\left(P_{<M\delta\lambda}\nabla B_{<\delta\lambda}\,\eta_\lambda\right)
\Bigr]
\frac{d\lambda}{\lambda}.
\end{split}
\end{equation*}
By Proposition \ref{Bernsteinelliptic}, we can replace $V_{<\delta\lambda}$ in the first term of the integrand by
\begin{equation*}
V^l_{<\delta\lambda}
=P_{<M\delta\lambda}V_{<\delta\lambda},
\end{equation*}
up to an error of the form $\mathcal R$. We can then expand the remaining part of the integrand as
\begin{equation*}
\begin{split}
[P_{<M\delta\lambda},V^l_{<\delta\lambda}]
\cdot\nabla B_{<\delta\lambda}\,\eta_\lambda+
\left(V^l_{<\delta\lambda}-T_{V^l}\right)\cdot
\left(P_{<M\delta\lambda}\nabla B_{<\delta\lambda}\,\eta_\lambda\right).
\end{split}
\end{equation*}
Proposition \ref{Leibnizcommutator} gives
\begin{equation*}
[P_{<M\delta\lambda},V^l_{<\delta\lambda}]
\cdot\nabla B_{<\delta\lambda}\,\eta_\lambda
=\lambda^{-1}\mathcal L
\left(
\nabla_x B_{<\delta\lambda},
\nabla_x V^l_{<\delta\lambda},
\eta_\lambda
\right).
\end{equation*}
By Propositions \ref{Bernsteinelliptic} and \ref{VBunbil},
\begin{equation*}
\begin{split}
\lambda^{-1}\mathcal L
\left(
\nabla_x B_{<\delta\lambda},
\nabla_x V^l_{<\delta\lambda},
\eta_\lambda
\right)=
\lambda^{-1}\mathcal L
\left(
\nabla_x B^l_{<\delta\lambda},
\nabla_x V^l_{<\delta\lambda},
\eta_\lambda
\right)
+O_{L_t^1L_x^2}
\left(\epsilon^3c_\lambda\lambda^{-\frac12-s_c}\right).
\end{split}
\end{equation*}
After integration in $\lambda$, the principal term on the right contributes a term of the form $\mathcal R$ by using two unbalanced bilinear estimates; one where we pair $(\eta_{\lambda},\nabla V^l_{<\delta\lambda})$ and one where we pair $(u_{\mu}^{x_0},\nabla_x B^l_{<\delta\lambda})$. Since $P_{<M\delta\lambda}\nabla B_{<\delta\lambda}\,\eta_\lambda$ is supported at frequency $\approx\lambda$, provided $M\delta\ll1$, we can write
\begin{equation*}
\begin{split}
&\left(V^l_{<\delta\lambda}-T_{V^l}\right)\cdot
\left(P_{<M\delta\lambda}\nabla B_{<\delta\lambda}\,\eta_\lambda\right)
\\
&\quad=
\int_{\lambda_1\approx\lambda}
\left(V^l_{<\delta\lambda}-V^l_{<\delta\lambda_1}\right)\cdot
S_{\lambda_1}
\left(P_{<M\delta\lambda}\nabla B_{<\delta\lambda}\,\eta_\lambda\right)
\frac{d\lambda_1}{\lambda_1}
\\
&\quad=
\int_{\lambda_1\approx\lambda}
P_{<M\delta\lambda}
\left(V_{<\delta\lambda}-V_{<\delta\lambda_1}\right)\cdot
S_{\lambda_1}
\left(P_{<M\delta\lambda}\nabla B_{<\delta\lambda}\,\eta_\lambda\right)
\frac{d\lambda_1}{\lambda_1}
\\
&\qquad+
O_{L_t^1L_x^2}
\left(\epsilon^3c_\lambda\lambda^{-\frac12-s_c}\right).
\end{split}
\end{equation*}
Here, the second equality follows from Proposition \ref{Bernsteinelliptic}.
By the fundamental theorem of calculus and the linearization formula \eqref{Vlin}, we have
\begin{equation*}
V_{<\delta\lambda}-V_{<\delta\lambda_1}
=
\int_{\delta\lambda_1}^{\delta\lambda}
\left(
M_{<\lambda_2}\Phi^2_{<\lambda_2}
+\eta_{\lambda_2}M_{<\lambda_2}B_{<\lambda_2}
\right)
\frac{d\lambda_2}{\lambda_2}.
\end{equation*}
Using Lemma \ref{Phiest}, with $A=M$, together with the $L_t^4L_x^\infty$ estimates for $\eta_{\lambda_2}$ and $M_{<\lambda_2}B_{<\lambda_2}$, we obtain
\begin{equation*}
V_{<\delta\lambda}-V_{<\delta\lambda_1}
=
\int_{\delta\lambda_1}^{\delta\lambda}
M_{<\lambda_2}\varphi_{\lambda_2}
\frac{d\lambda_2}{\lambda_2}
+O_{L_t^2L_x^\infty}
\left(\epsilon^2\lambda^{-\frac34}\right).
\end{equation*}
Using the unbalanced estimate
\begin{equation*}
\left\|
P_{<M\delta\lambda}\nabla B_{<\delta\lambda}\,\eta_\lambda
\right\|_{L_{t,x}^2}
\lesssim
\epsilon^2c_\lambda\lambda^{\frac14-s_c},
\end{equation*}
we conclude that
\begin{equation*}
\begin{split}
&\left(V^l_{<\delta\lambda}-T_{V^l}\right)\cdot
\left(P_{<M\delta\lambda}\nabla B_{<\delta\lambda}\,\eta_\lambda\right)
\\
&\quad=
\int_{\lambda_1\approx\lambda}
\int_{\delta\lambda_1}^{\delta\lambda}
\mathcal L
\left(
P_{<M\delta\lambda}M_{<\lambda_2}\varphi_{\lambda_2},
\nabla B^l_{<\delta\lambda},
\eta_\lambda
\right)
\frac{d\lambda_2}{\lambda_2}
\frac{d\lambda_1}{\lambda_1}
\\
&\qquad+
O_{L_t^1L_x^2}
\left(\epsilon^3c_\lambda\lambda^{-\frac12-s_c}\right).
\end{split}
\end{equation*}
After projecting the first term above to frequency $\mu$ and multiplying by $\overline u_\mu^{x_0}$, two unbalanced bilinear estimates show that it contributes a perturbative error of the form $\mathcal R$. This concludes the proof.
\end{proof}

\noindent Next, we identify the difference $\widetilde a^l-a^l$, which occurs (with a minus sign) in the last term in \eqref{initialphipara} as a para-coefficient. To do this, we need a sufficiently precise expansion of $a_{<\delta\lambda}-\widetilde a_{<\delta\lambda}$ for each $\lambda$. We begin by defining the truncated Bernoulli "source term"
\begin{equation*}
\mathcal F_{<\delta\lambda}\coloneq 
\partial_t\phi_{<\delta\lambda}
+\frac12|\nabla_{x,z}\phi_{<\delta\lambda}|^2+z.
\end{equation*}
Since $\phi_{<\delta\lambda}$ is harmonic, we have $\Delta_{x,z}\mathcal F_{<\delta\lambda}
=|\nabla^2_{x,z}\phi_{<\delta\lambda}|^2$, whereas $\Delta_{x,z}p_{<\delta\lambda}
=-|\nabla^2_{x,z}\phi_{<\delta\lambda}|^2.$ Consequently, $\mathcal F_{<\delta\lambda}+p_{<\delta\lambda}$ is harmonic in the truncated domain
\begin{equation*}
\Omega_{<\delta\lambda}
\coloneq
\{(x,z)\in\mathbb R^3:z<\eta_{<\delta\lambda}(t,x)\}
\end{equation*}
and has boundary trace
\begin{equation*}
\rho_{<\delta\lambda}
\coloneq
\left.
\left(
\partial_t\phi_{<\delta\lambda}
+\frac12|\nabla_{x,z}\phi_{<\delta\lambda}|^2+z
\right)
\right|_{z=\eta_{<\delta\lambda}}.
\end{equation*}
Therefore,
\begin{equation*}
\left.\partial_z\mathcal F_{<\delta\lambda}\right|_{z=\eta_{<\delta\lambda}}
=a_{<\delta\lambda}+L_{<\delta\lambda}\rho_{<\delta\lambda},
\end{equation*}
where we used the fact that, if $f$ is the trace of a harmonic function $g$ on $\Omega_{<\delta\lambda}$, then
\begin{equation*}
\left.\partial_zg\right|_{z=\eta_{<\delta\lambda}(x)}
=L_{<\delta\lambda}f.
\end{equation*}
On the other hand, differentiating
\begin{equation*}
B_{<\delta\lambda}
=
\left.\partial_z\phi_{<\delta\lambda}\right|_{z=\eta_{<\delta\lambda}}
\end{equation*}
along $\partial_t+V_{<\delta\lambda}\cdot\nabla$ gives
\begin{equation*}
\begin{split}
(\partial_t+V_{<\delta\lambda}\cdot\nabla)B_{<\delta\lambda}&=
\left.
\partial_z
\left(
\partial_t\phi_{<\delta\lambda}
+\frac12|\nabla_{x,z}\phi_{<\delta\lambda}|^2
\right)
\right|_{z=\eta_{<\delta\lambda}}
\\
&+
\left(
(\partial_t+V_{<\delta\lambda}\cdot\nabla)\eta_{<\delta\lambda}
-B_{<\delta\lambda}
\right)
\left.\partial_z^2\phi_{<\delta\lambda}\right|_{z=\eta_{<\delta\lambda}}.
\end{split}
\end{equation*}
Because $\partial_z\phi_{<\delta\lambda}$ is harmonic and has boundary trace $B_{<\delta\lambda}$, we have
\begin{equation*}
\left.\partial_z^2\phi_{<\delta\lambda}\right|_{z=\eta_{<\delta\lambda}}
=L_{<\delta\lambda}B_{<\delta\lambda}.
\end{equation*}
It follows that
\begin{equation*}
\begin{split}
(\partial_t+V_{<\delta\lambda}\cdot\nabla)B_{<\delta\lambda}
=a_{<\delta\lambda}-1
+L_{<\delta\lambda}\rho_{<\delta\lambda}+\kappa_{<\delta\lambda}L_{<\delta\lambda}B_{<\delta\lambda},
\end{split}
\end{equation*}
where
\begin{equation}\label{kappadef}
\kappa_{<\delta\lambda}
\coloneq
(\partial_t+V_{<\delta\lambda}\cdot\nabla)\eta_{<\delta\lambda}
-B_{<\delta\lambda}.
\end{equation}
Hence
\begin{equation*}
\widetilde a_{<\delta\lambda}-a_{<\delta\lambda}
=\mathcal E_{<\delta\lambda}
\coloneq
L_{<\delta\lambda}\rho_{<\delta\lambda}
+\kappa_{<\delta\lambda}L_{<\delta\lambda}B_{<\delta\lambda},
\end{equation*}
and therefore
\begin{equation*}
\left(\partial_t+T_{V^l}\cdot\nabla\right)\varphi
+T_{a^l}\eta
=N_{\varphi}-T_{\mathcal E^l}\eta+\mathcal R.
\end{equation*}
We are now in a position to decompose $T_{\mathcal E^l}\eta$ into a perturbative remainder and the flux term appearing in Proposition \ref{paravarphiinitial}. For this, we use the following lemma.

\begin{lemma}\label{aerrorinitial}
The error term $T_{\mathcal E^l}\eta$ contributes a perturbative error of the form $\mathcal R$ plus a flux term. That is, for each $\mu\in2^{\mathbb Z}$, there holds
\begin{equation}
\overline u_\mu^{x_0}P_\mu(T_{\mathcal E^l}\eta)
=\nabla_x\cdot Q_{\varphi,\mu}
+O_{L_t^1L_x^1}
\left(
\epsilon^4c_\mu^2\mu^{-\frac12-2s_c}(1+\mu|x_0|)^{\frac12}
\right),
\end{equation}
where $Q_{\varphi,\mu}$ satisfies \eqref{QPhibound}.
\end{lemma}

\begin{proof}
To avoid getting bogged down in notation, it will be convenient to define the deviation functions
\begin{equation*}
d_B\coloneq B_{<\delta\lambda}-P_{<\delta\lambda}B,
\qquad
d_V\coloneq V_{<\delta\lambda}-P_{<\delta\lambda}V.
\end{equation*}
We also define the following bilinear form that will appear quite frequently below
\begin{equation*}
\mathcal B_{<\delta\lambda}(f,g)
\coloneq
(P_{<\delta\lambda}f)(P_{<\delta\lambda}g)
-P_{<\delta\lambda}(fg).
\end{equation*}
The symbol of this bilinear form is
\begin{equation*}
b(\xi_1,\xi_2)
\coloneq
\chi(\Lambda^{-1}\xi_1)\chi(\Lambda^{-1}\xi_2)
-\chi(\Lambda^{-1}(\xi_1+\xi_2)),
\qquad
\Lambda\coloneq\delta\lambda.
\end{equation*}
Above, $\chi$ is a smooth cutoff supported on the ball $B_2(0)$ and equal to $1$ on $B_1(0)$. The symbol $b$ vanishes when either $\xi_1=0$ or $\xi_2=0$. Therefore, separating the low--low and high--high-to-low regions, we have the representation
\begin{equation}\label{bilformrep}
\mathcal B_{<\delta\lambda}(f,g)
=
\lambda^{-2}\mathcal L^l
\bigl(\nabla P_{<\lambda}f,\nabla P_{<\lambda}g\bigr)
+\mathcal L^l
\bigl(P_{\gtrsim\lambda}f,P_{\gtrsim\lambda}g\bigr),
\end{equation}
where the first term is a translation-invariant bilinear form with output frequency $\ll\lambda$. The second term (necessarily) denotes a dyadic sum of such forms over comparable high-frequency inputs with low output.

We next observe the identities
\begin{equation*}
V=\nabla\psi-B\nabla\eta,
\qquad
V_{<\delta\lambda}
=\nabla\psi_{<\delta\lambda}
-B_{<\delta\lambda}\nabla\eta_{<\delta\lambda},
\end{equation*}
which give
\begin{equation}\label{eq:dV-dB-relation}
d_V
=-d_B\nabla\eta_{<\delta\lambda}
-\mathcal B_{<\delta\lambda}(B,\nabla\eta).
\end{equation}
We will also need quantitative estimates for the deviation functions. We first write
\begin{equation}\label{dBVdecomp}
\begin{split}
d_B
={}&
\bigl(B_{<\delta\lambda}-T_{\mathfrak l}\varphi_{<\delta\lambda}\bigr)
-P_{<\delta\lambda}\bigl(B-T_{\mathfrak l}\varphi\bigr)
+[T_{\mathfrak l},P_{<\delta\lambda}]\varphi,
\\
d_V
={}&
\bigl(V_{<\delta\lambda}-T_{\mathfrak m}\varphi_{<\delta\lambda}\bigr)
-P_{<\delta\lambda}\bigl(V-T_{\mathfrak m}\varphi\bigr)
+[T_{\mathfrak m},P_{<\delta\lambda}]\varphi.
\end{split}
\end{equation}
We compare the full and truncated expansions used in the proof of Corollary \ref{principalVBexp}. In this comparison, the contributions whose "geometric" truncation scale is $\ll\delta\lambda$ cancel. The remaining contributions come from scales $\gtrsim\delta\lambda$, the high frequency tails in those expansions, and the two commutators in \eqref{dBVdecomp}. Corollary \ref{principalVBexp}, Proposition \ref{Leibnizcommutator}, and Proposition \ref{Bernsteinelliptic} therefore give
\begin{equation}\label{dvdbL2Linf}
\|d_B\|_{L_t^2L_x^\infty}
+\|d_V\|_{L_t^2L_x^\infty}
\lesssim_\delta
\epsilon^2\lambda^{-\frac34}.
\end{equation}
Consequently, if $\lambda\approx\mu$, the energy bound for $u_\mu^{x_0}$ gives
\begin{equation}\label{dvdBbil}
\begin{split}
\|\mathcal L(\overline u_\mu^{x_0},d_B)\|_{L_{t,x}^2}
+\|\mathcal L(\overline u_\mu^{x_0},d_V)\|_{L_{t,x}^2}\lesssim
\epsilon^2c_\lambda\lambda^{-\frac34-s_c}(1+\mu|x_0|)^{\frac12}.
\end{split}
\end{equation}
The analogous estimate with $\overline u_\mu^{x_0}$ replaced by $\eta_\lambda$ also holds (without the translation parameter). The same comparison, now using the $L_t^4L_x^p$ estimates in Proposition \ref{VBunbil}, gives
\begin{equation}\label{dBdVL4p}
\|d_B\|_{L_t^4L_x^p}
+\|d_V\|_{L_t^4L_x^p}
\lesssim
\epsilon\lambda^{-\frac78+2(\frac18-\frac1p)},
\qquad 8<p\leq\infty.
\end{equation}

We now estimate the contribution of
\begin{equation*}
\mathcal E_{\kappa,<\delta\lambda}
\coloneq
\kappa_{<\delta\lambda}L_{<\delta\lambda}B_{<\delta\lambda}
\end{equation*}
in $T_{\mathcal E^l}\eta$. We denote the corresponding term by $T_{\mathcal E_\kappa^l}\eta$. Using the kinematic boundary condition $\partial_t\eta=B-V\cdot\nabla\eta$, we have
\begin{equation*}
\begin{split}
\kappa_{<\delta\lambda}
&=
P_{<\delta\lambda}(B-V\cdot\nabla\eta)
+V_{<\delta\lambda}\cdot\nabla\eta_{<\delta\lambda}
-B_{<\delta\lambda}
\\
&=-d_B+d_V\cdot\nabla\eta_{<\delta\lambda}
+\mathcal B_{<\delta\lambda}(V,\nabla\eta).
\end{split}
\end{equation*}
Therefore, \eqref{eq:dV-dB-relation} gives the exact identity
\begin{equation}\label{kappaform}
\begin{split}
\kappa_{<\delta\lambda}
=
-\left(1+|\nabla\eta_{<\delta\lambda}|^2\right)d_B-\nabla\eta_{<\delta\lambda}\cdot
\mathcal B_{<\delta\lambda}(B,\nabla\eta)
+\mathcal B_{<\delta\lambda}(V,\nabla\eta).
\end{split}
\end{equation}
Using the representation \eqref{bilformrep}, the $L_t^4L_x^\infty$ bounds, and Proposition \ref{VBunbil}, we obtain
\begin{equation*}
\begin{split}
&\|\mathcal B_{<\delta\lambda}(B,\nabla\eta)\|_{L_t^2L_x^\infty}
+\|\mathcal B_{<\delta\lambda}(V,\nabla\eta)\|_{L_t^2L_x^\infty}
\\
&\qquad\lesssim
\lambda^{-2}
\left(
\|\nabla P_{<\lambda}B\|_{L_t^4L_x^\infty}
+\|\nabla P_{<\lambda}V\|_{L_t^4L_x^\infty}
\right)
\|\nabla^2P_{<\lambda}\eta\|_{L_t^4L_x^\infty}
\\
&\qquad+
\sum_{\nu\gtrsim\lambda}
\left(
\|B_\nu\|_{L_t^4L_x^\infty}
+\|V_\nu\|_{L_t^4L_x^\infty}
\right)
\|\nabla\eta_{\approx\nu}\|_{L_t^4L_x^\infty}
\\
&\qquad\lesssim
\epsilon^2\lambda^{-\frac34}.
\end{split}
\end{equation*}
Combining this with \eqref{kappaform}, the bound for $d_B$ above, and the bound $\|\nabla\eta_{<\delta\lambda}\|_{L_t^\infty L_x^\infty}
\lesssim\epsilon$, then gives
\begin{equation*}
\|\kappa_{<\delta\lambda}\|_{L_t^2L_x^\infty}
\lesssim
\epsilon^2\lambda^{-\frac34}.
\end{equation*}
It follows that, whenever $\lambda\approx\mu$,
\begin{equation}\label{kappaubil}
\|\mathcal L(\overline u_\mu^{x_0},\kappa_{<\delta\lambda})\|_{L_{t,x}^2}
\lesssim
\epsilon^2c_\lambda\lambda^{-\frac34-s_c}(1+\mu|x_0|)^{\frac12}.
\end{equation}

Using Propositions \ref{VBunbil} and \ref{Bernsteinelliptic}, we may also replace the term $L_{<\delta\lambda}B_{<\delta\lambda}$ in $\mathcal{E}_{\kappa}$ by its low-output truncated variant,
\begin{equation*}
\begin{split}
P_{<M\delta\lambda}
\left(
\kappa_{<\delta\lambda}L_{<\delta\lambda}B_{<\delta\lambda}
\right)\eta_\lambda=
P_{<M\delta\lambda}
\left(
\kappa_{<\delta\lambda}
L_{<\delta\lambda}^{l,\delta\lambda}B^l_{<\delta\lambda}
\right)\eta_\lambda
+O_{L_t^1L_x^2}
\left(
\epsilon^3c_\lambda\lambda^{-\frac12-s_c}
\right).
\end{split}
\end{equation*}
We also record the bound
\begin{equation*}
\|(L_{<\delta\lambda}^{l,\delta\lambda}-T_L)B^l_{<\delta\lambda}\|_{L_t^2L_x^\infty}
\lesssim
\epsilon^2\lambda^{\frac14}
\end{equation*}
which follows in exactly the same way as the bound \eqref{L2TL} proved previously. Together with Corollary \ref{Abilbounds} and the energy bound for $\eta_\lambda$, this yields
\begin{equation*}
\|\mathcal L(\eta_\lambda,L_{<\delta\lambda}^{l,\delta\lambda}B^l_{<\delta\lambda})\|_{L_{t,x}^2}
\lesssim
\epsilon^2c_\lambda\lambda^{\frac14-s_c}.
\end{equation*}
Therefore,
\begin{equation*}
\begin{split}
\|\overline u_\mu^{x_0}P_\mu(T_{\mathcal E_\kappa^l}\eta)\|_{L_t^1L_x^1}&\lesssim
\epsilon^4c_\mu^2\mu^{-\frac12-2s_c}(1+\mu|x_0|)^{\frac12}
\\
&+
\int_{\lambda\approx\mu}
\|\mathcal L(\overline u_\mu^{x_0},\kappa_{<\delta\lambda})\|_{L_{t,x}^2}
\|\mathcal L(\eta_\lambda,L_{<\delta\lambda}^{l,\delta\lambda}B^l_{<\delta\lambda})\|_{L_{t,x}^2}
\frac{d\lambda}{\lambda}
\\
&\lesssim
\epsilon^4c_\mu^2\mu^{-\frac12-2s_c}(1+\mu|x_0|)^{\frac12}.
\end{split}
\end{equation*}
Thus $T_{\mathcal E_\kappa^l}\eta$ is perturbative. Next, we deal with the contribution of the boundary trace term
\begin{equation*}
\mathcal E_{\rho,<\delta\lambda}
\coloneq
L_{<\delta\lambda}\rho_{<\delta\lambda}.
\end{equation*}
We denote its contribution in $T_{\mathcal E^l}\eta$ by $T_{\mathcal E_\rho^l}\eta$. We begin by writing $\rho_{<\delta\lambda}$ in a form that is easier to work with. First, the chain rule gives
\begin{equation*}
\left.\partial_t\phi_{<\delta\lambda}\right|_{z=\eta_{<\delta\lambda}}
=
\partial_t\psi_{<\delta\lambda}
-B_{<\delta\lambda}\partial_t\eta_{<\delta\lambda}.
\end{equation*}
On the other hand, the dynamic and kinematic boundary conditions give
\begin{equation*}
\partial_t\psi
=
B\partial_t\eta-\eta-\frac12|V|^2-\frac12B^2.
\end{equation*}
Using also $\partial_t\eta_{<\delta\lambda}
=P_{<\delta\lambda}(\partial_t\eta)$, we obtain
\begin{equation}
\label{eq:rho-first-form}
\begin{aligned}
\rho_{<\delta\lambda}
={}&
P_{<\delta\lambda}(B\partial_t\eta)
-B_{<\delta\lambda}P_{<\delta\lambda}(\partial_t\eta)
\\
&+\frac12\left(
|V_{<\delta\lambda}|^2-P_{<\delta\lambda}|V|^2
\right)
+\frac12\left(
(B_{<\delta\lambda})^2-P_{<\delta\lambda}(B^2)
\right).
\end{aligned}
\end{equation}
Before further manipulating this expression, it will be convenient to define the trilinear forms
\begin{equation*}
\mathcal C_{<\delta\lambda}(B,V,\nabla\eta)
=
(P_{<\delta\lambda}B)(P_{<\delta\lambda}V)
\cdot\nabla\eta_{<\delta\lambda}
-P_{<\delta\lambda}(BV\cdot\nabla\eta)
\end{equation*}
and
\begin{equation*}
\begin{split}
\mathcal S_{<\delta\lambda}(B,V,\nabla\eta)
=
\mathcal C_{<\delta\lambda}(B,V,\nabla\eta)-(P_{<\delta\lambda}B)
\mathcal B_{<\delta\lambda}(V,\nabla\eta)-(P_{<\delta\lambda}V)\cdot
\mathcal B_{<\delta\lambda}(B,\nabla\eta).
\end{split}
\end{equation*}
The symbol of the trilinear form $\mathcal S_{<\delta\lambda}(B,V,\nabla\eta)$ is (again taking $\Lambda\coloneq\delta\lambda$),
\begin{equation*}
\begin{split}
s(\xi_1,\xi_2,\xi_3)
={}&
-\chi(\Lambda^{-1}(\xi_1+\xi_2+\xi_3))
-\chi(\Lambda^{-1}\xi_1)
 \chi(\Lambda^{-1}\xi_2)
 \chi(\Lambda^{-1}\xi_3)
\\
&+\chi(\Lambda^{-1}\xi_1)
 \chi(\Lambda^{-1}(\xi_2+\xi_3))
+\chi(\Lambda^{-1}\xi_2)
 \chi(\Lambda^{-1}(\xi_1+\xi_3)).
\end{split}
\end{equation*}
This symbol vanishes when $\xi_1=0$ or $\xi_2=0$:
\begin{equation*}
s(0,\xi_2,\xi_3)=0,
\qquad
s(\xi_1,0,\xi_3)=0.
\end{equation*}
Therefore, separating the cases according to which of the first two inputs have frequency $\gtrsim\lambda$, we have the representation
\begin{equation}
\label{Srep}
\begin{split}
\mathcal S_{<\delta\lambda}(B,V,\nabla\eta)
&=
\lambda^{-2}\mathcal L^l
(\nabla P_{<\lambda}B,\nabla P_{<\lambda}V,\nabla\eta)+\lambda^{-1}\mathcal L^l
(\nabla P_{<\lambda}B,P_{\gtrsim\lambda}V,
P_{\gtrsim\lambda}\nabla\eta)
\\
&+\lambda^{-1}\mathcal L^l
(P_{\gtrsim\lambda}B,\nabla P_{<\lambda}V,
P_{\gtrsim\lambda}\nabla\eta)+\mathcal L^l
(P_{\gtrsim\lambda}B,P_{\gtrsim\lambda}V,\nabla\eta),
\end{split}
\end{equation}
where $\mathcal L^l$ denotes a translation-invariant trilinear form with output frequency $\ll\lambda$ (that in the sequel may change from line to line). Using the kinematic boundary condition, \eqref{eq:dV-dB-relation}, and \eqref{eq:rho-first-form}, we then obtain the useful form of $\rho_{<\delta\lambda}$,
\begin{equation}
\label{eq:rho-good-form}
\begin{split}
\rho_{<\delta\lambda}
={}&
\frac12\mathcal B_{<\delta\lambda}(V,V)
-\frac12\mathcal B_{<\delta\lambda}(B,B)
+\mathcal S_{<\delta\lambda}(B,V,\nabla\eta)
\\
&-d_B\mathcal B_{<\delta\lambda}(V,\nabla\eta)
+\frac12|d_V|^2
+\frac12d_B^2.
\end{split}
\end{equation}
We denote the first and second lines on the right by $\rho^1_{<\delta\lambda}$ and $\rho^2_{<\delta\lambda}$, respectively. Using \eqref{bilformrep} and \eqref{Srep}, together with \eqref{dBdVL4p}, Proposition \ref{VBunbil}, Proposition \ref{Bernsteinelliptic}, and Bernstein's inequality, we record the useful bound
\begin{equation}
\label{rhobasicbound}
\lambda^{-j}\|\nabla_x^j\rho_{<\delta\lambda}\|_{L_t^2L_x^{\infty}}
\lesssim
\epsilon^2\lambda^{-\frac54},\hspace{5mm}j=0,1,2.
\end{equation}
We will use the undifferentiated estimate (i.e. $j=0$) in \eqref{rhobasicbound} almost immediately, while the differentiated estimate will be needed later when we compute the material derivative of the Taylor coefficient, $D_ta$. Another useful consequence of this bound that we will need for the auxiliary bounds \eqref{NphiinitialL2Linf} and \eqref{NphiinitialauxR} is the elliptic bound
\begin{equation}\label{LrhoauxL2Linf}
\|L_{<\delta\lambda}\rho_{<\delta\lambda}
\|_{L_t^2L_x^\infty}
\lesssim\epsilon^2\lambda^{-\frac14}.
\end{equation}
Indeed, by Proposition \ref{basicopbounds} and interpolation, and then \eqref{rhobasicbound} we have 
\begin{equation*}
\begin{split}
\|L_{<\delta\lambda}\rho_{<\delta\lambda}\|_{L_t^2L_x^{\infty}}&\lesssim \||D|\rho_{<\delta\lambda}\|_{L_{t}^2\dot{B}_{\infty,1}^0}\lesssim \lambda \|\rho_{<\delta\lambda}\|_{L_t^2L_x^{\infty}}+\lambda^{-1}\|\nabla^2_x\rho_{<\delta\lambda}\|_{L_t^2L_x^{\infty}}\lesssim \epsilon^2\lambda^{-\frac{1}{4}}.
\end{split}    
\end{equation*}
Now, we turn to the estimate for $T_{\mathcal E_\rho^l}\eta$. We begin by controlling the error from replacing $L^{l}_{<\delta\lambda}\rho_{<\delta\lambda}$ by the principal term $T_L\rho_{<\delta\lambda}$. That is, we estimate the contribution of
\begin{equation*}
(L_{<\delta\lambda}^{l}-T_L)
\rho_{<\delta\lambda}.
\end{equation*}
Using Proposition \ref{Bernsteinelliptic}, it is easy to handle the contribution of $P_{\beta}\rho_{<\delta\lambda}$ when $\beta\gg\delta\lambda$, so we consider pieces at frequency $\lesssim\delta\lambda$, 
\begin{equation*}
\rho_\alpha=S_\alpha\rho_{<\delta\lambda},
\qquad
\alpha\lesssim\delta\lambda.
\end{equation*}
As in the estimate for $(L_2-T_L)B^l_{<\delta\lambda_2}$ in \eqref{L2TL}, we can write for $\alpha\lesssim\delta\lambda,$
\begin{equation}\label{rhoathreeterms}
\begin{split}
(L_{<\delta\lambda}^{l}-T_L)\rho_{\alpha}&=
P_{<M\delta\lambda}(L_{<\delta\lambda}
-L_{<\delta\alpha})\rho_{\alpha}+(P_{<M\delta\lambda}-P_{<M\alpha})
L_{<\delta\alpha}\rho_{\alpha}+(P_{<M\alpha}L_{<\delta\alpha}-T_L)\rho_{\alpha}.
\end{split}
\end{equation}
The middle term has high-frequency tails and its contribution is easily shown to be perturbative by using Proposition \ref{Bernsteinelliptic}. For the first term, write
\begin{equation*}
P_{<M\delta\lambda}(L_{<\delta\lambda}
-L_{<\delta\alpha})\rho_{\alpha}
=
P_{<M\delta\lambda}\int_{\delta\alpha}^{\delta\lambda}
\gamma\frac{d}{d\gamma}
L_{<\gamma}\rho_{\alpha}
\frac{d\gamma}{\gamma},
\end{equation*}
where we keep the output projection fixed. More precisely, the linearization formula \eqref{Blin} gives
\begin{equation*}
\begin{split}
\gamma\frac{d}{d\gamma}
L_{<\gamma}^{l,\delta\lambda}\rho_\alpha
=
P_{<M\delta\lambda}\Bigl[
-L_{<\gamma}\bigl(
\eta_\gamma L_{<\gamma}\rho_\alpha
\bigr)
+\eta_\gamma L_{<\gamma}^2\rho_\alpha
\Bigr].
\end{split}
\end{equation*}
Using \eqref{rhobasicbound}, Proposition \ref{paraexp}, and two
$L_t^4L_x^p$ bounds for $\eta$, we may extract the
principal parts from the linearized expression above. After carrying out
the $\gamma$ integration in the resulting error terms, we obtain
\begin{equation*}
\begin{split}
P_{<M\delta\lambda}
(L_{<\delta\lambda}-L_{<\delta\alpha})\rho_\alpha&=
P_{<M\delta\lambda}
\int_{\delta\alpha}^{\delta\lambda}
\left[
-T_{\mathfrak l}
\bigl(\eta_\gamma T_{\mathfrak l}\rho_\alpha\bigr)
+\eta_\gamma T_{\mathfrak l}^2\rho_\alpha
\right]\frac{d\gamma}{\gamma}
\\
&+O_{L_t^1L_x^\infty}
\left(
\epsilon^4\lambda^{-\frac54}\alpha^{\frac34}
\right).
\end{split}
\end{equation*}
Since we have
\begin{equation*}
\int_{\alpha\lesssim\delta\lambda}
\lambda^{-\frac54}\alpha^{\frac34}
\frac{d\alpha}{\alpha}
\lesssim_\delta\lambda^{-\frac12},
\end{equation*}
this gives a perturbative contribution in $T_{\mathcal E_\rho^l}\eta$ after
integration in $\alpha$ and placing the high-frequency factor of $\eta$ in
$L_t^\infty L_x^2$. Hence, it remains to estimate
\begin{equation*}
\begin{split}
\mathcal I\coloneq
\int_0^\infty
\int_{\alpha\lesssim\delta\lambda}
\int_{\delta\alpha}^{\delta\lambda}
\eta_\lambda P_{<M\delta\lambda}
\left[
-T_{\mathfrak l}
\bigl(\eta_\gamma T_{\mathfrak l}\rho_\alpha\bigr)
+\eta_\gamma T_{\mathfrak l}^2\rho_\alpha
\right]
\frac{d\gamma}{\gamma}
\frac{d\alpha}{\alpha}
\frac{d\lambda}{\lambda}.
\end{split}
\end{equation*}
This term is perturbative by a direct estimate in $L_t^1L_x^2$ (and placing $u_{\mu}^{x_0}$ in $L_t^{\infty}L_x^2$). Indeed,  we can place
$\rho_\alpha$ in $L_t^2L_x^\infty$ and use an unbalanced bilinear estimate for
$(\eta_\lambda,\eta_\gamma)$. For instance, for the first term, we can decompose the outer operator $P_{<M\delta\lambda}T_{\mathfrak{l}}=\int_{\beta\lesssim\gamma}P_{<M\delta\lambda}T_{\mathfrak{l}}S_{\beta}\frac{d\beta}{\beta}$.
The kernel $k_\beta(t,x,y)$ of $T_{\mathfrak l}S_\beta$ satisfies
\begin{equation*}
|k_\beta(t,x,y)|\lesssim_N
\beta^3\langle\beta y\rangle^{-N}.
\end{equation*}
Consequently, by \eqref{uab-bi-unbal},
\begin{equation*}
\begin{split}
\left\|
P_\mu\left(
\eta_\lambda P_{<M\delta\lambda}
T_{\mathfrak l}S_\beta
\bigl(\eta_\gamma T_{\mathfrak l}\rho_\alpha\bigr)
\right)
\right\|_{L_t^1L_x^2}&\lesssim\alpha\|\rho_\alpha\|_{L_t^2L_x^\infty}
\int_{\mathbb R^2}
\sup_{t,x}|k_\beta(t,x,y)|
\|\eta_\lambda\tau_y\eta_\gamma\|_{L_{t,x}^2}\,dy
\\
&\lesssim
\epsilon^2c_\lambda\lambda^{-s_c}\gamma^{-\frac54}
\alpha\|\rho_\alpha\|_{L_t^2L_x^\infty}
\beta\left(1+\frac{\gamma}{\beta}\right)^{\frac12}.
\end{split}
\end{equation*}
The square root translation factor is then mild enough to give the integrable bound
\begin{equation*}
\int_{\beta\lesssim_\delta\gamma}
\beta\left(1+\frac{\gamma}{\beta}\right)^{\frac12}
\frac{d\beta}{\beta}
\lesssim_\delta\gamma.
\end{equation*}
Together with \eqref{rhobasicbound}, the contribution at fixed
$\alpha$ and $\gamma$, after integration in $\beta$, is therefore bounded
by
\begin{equation*}
\epsilon^4c_\lambda
\lambda^{-s_c-\frac54}
\alpha\gamma^{-\frac14}.
\end{equation*}
For the term containing $\eta_\gamma T_{\mathfrak l}^2\rho_\alpha$,
a similar analysis gives the same bound (in fact the situation is cleaner here because the outer operator is just $P_{<M\delta\lambda}$ which has rapidly decaying kernel at scale $\lambda$), since
$T_{\mathfrak l}^2$ contributes $\alpha^2$ and
$\alpha\lesssim_\delta\gamma$. It follows that
\begin{equation*}
\|P_\mu\mathcal I\|_{L_t^1L_x^2}
\lesssim
\epsilon^4c_\mu\mu^{-\frac12-s_c},
\end{equation*}
where we used that the output localization forces
$\lambda\approx\mu$. Finally, placing $u_\mu^{x_0}$ in
$L_t^\infty L_x^2$, we obtain
\begin{equation*}
\|\overline u_\mu^{x_0}P_\mu\mathcal I\|_{L_t^1L_x^1}
\lesssim
\epsilon^4c_\mu^2\mu^{-\frac12-2s_c}
(1+\mu|x_0|)^{\frac12}.
\end{equation*}
This handles the first term in \eqref{rhoathreeterms}. The third term in \eqref{rhoathreeterms} is treated in essentially the same way. It therefore remains to treat the principal term
\begin{equation*}
\overline u_\mu^{x_0}P_\mu
\left(
P_{<M\delta\lambda}(T_L\rho_{<\delta\lambda})
\eta_\lambda
\right).
\end{equation*}
We first estimate the contribution of $\rho^2_{<\delta\lambda}$.
This requires no bilinear estimates. Indeed, using Proposition \ref{basicopbounds}, we have
\begin{equation*}
\begin{split}
\left\|
P_{<M\delta\lambda}
\bigl(T_L\rho^2_{<\delta\lambda}\bigr)
\right\|_{L_t^1L_x^\infty}\lesssim
\lambda\Bigl(
\|d_B\|_{L_t^2L_x^\infty}
\|\mathcal B_{<\delta\lambda}(V,\nabla\eta)\|_{L_t^2L_x^\infty}
+\|d_V\|_{L_t^2L_x^\infty}^2
+\|d_B\|_{L_t^2L_x^\infty}^2
\Bigr).
\end{split}
\end{equation*}
From \eqref{bilformrep} and the $L_t^4L_x^\infty$ estimates for
$V$ and $\eta$, we have
\begin{equation*}
\|\mathcal B_{<\delta\lambda}(V,\nabla\eta)\|_{L_t^2L_x^\infty}
\lesssim
\epsilon^2\lambda^{-\frac34}.
\end{equation*}
Combining this with \eqref{dvdbL2Linf}, we obtain
\begin{equation*}
\left\|
P_{<M\delta\lambda}
\bigl(T_L\rho^2_{<\delta\lambda}\bigr)
\right\|_{L_t^1L_x^\infty}
\lesssim
\epsilon^4\lambda^{-\frac12}.
\end{equation*}
Placing both $u_\mu^{x_0}$ and $\eta_\lambda$ in
$L_t^\infty L_x^2$, and using that the output localization forces
$\lambda\approx\mu$, we conclude that
\begin{equation*}
\begin{split}
\bigl\|
\overline u_\mu^{x_0}P_\mu
\left(
P_{<M\delta\lambda}(T_L\rho^2_{<\delta\lambda})
\eta_\lambda
\right)
\bigr\|_{L_{t,x}^1}\lesssim
\epsilon^4c_\mu^2\mu^{-\frac12-2s_c}
(1+\mu|x_0|)^{\frac12}.
\end{split}
\end{equation*}
Next, we handle the contribution of $\rho^1_{<\delta\lambda}$. This term is more difficult because the interactions between $B$ and $V$ can be balanced at frequency comparable to $\mu$, making it impossible to use only unbalanced bilinear estimates. The saving grace is that the principal part of $\rho^1_{<\delta\lambda}$ possesses a null structure that allows us to use balanced bilinear estimates instead. It is also this term that will contribute the flux component in Proposition \ref{paravarphiinitial}, which will play a role in the interaction Morawetz estimates proved later. We isolate this structure with the following lemma.

\begin{lemma}\label{nullprincipal}
Let $M>1$ with $\delta\ll M\delta\ll1$, so that
$P_{\geq M\delta\lambda}P_{<\delta\lambda}=0$. For a function $f$, define
\begin{equation*}
f^{\mathrm{hi}}\coloneq P_{\geq M\delta\lambda}f,
\qquad
f^{\mathrm{lo}}\coloneq P_{<M\delta\lambda}f,
\end{equation*}
so that $f=f^{\mathrm{lo}}+f^{\mathrm{hi}}$. For the dyadic components below, we use the convention $f_\alpha^{\mathrm{hi}}\coloneq(f^{\mathrm{hi}})_\alpha$; in particular, $\varphi_\alpha^{\mathrm{hi}}=(\varphi^{\mathrm{hi}})_\alpha$. Then
\begin{equation}
\rho^1_{<\delta\lambda}
=\mathcal N_{\rho,\lambda}^{\mathrm{null}}+H_{\rho,\lambda},
\end{equation}
where
\begin{equation}\label{nullformdef}
\begin{split}
\mathcal N_{\rho,\lambda}^{\mathrm{null}}
&\coloneq P_{<\delta\lambda}
\int_0^\infty\int_{\delta\alpha}^{\alpha}
\left[
T_{\mathfrak l}\bigl(\varphi_\alpha^{\mathrm{hi}}\bigr)
T_{\mathfrak l}\bigl(\varphi_\beta^{\mathrm{hi}}\bigr)
-
T_{\mathfrak m}\bigl(\varphi_\alpha^{\mathrm{hi}}\bigr)
\cdot T_{\mathfrak m}\bigl(\varphi_\beta^{\mathrm{hi}}\bigr)
\right]
\frac{d\beta}{\beta}\frac{d\alpha}{\alpha}
\\
&-
P_{<\delta\lambda}
\int_0^\infty\int_{\delta\alpha}^{\alpha}
\nabla\eta_{<\delta\lambda}\cdot
\left[
T_{\mathfrak l}\bigl(\varphi_\alpha^{\mathrm{hi}}\bigr)
T_{\mathfrak m}\bigl(\varphi_\beta^{\mathrm{hi}}\bigr)
+
T_{\mathfrak l}\bigl(\varphi_\beta^{\mathrm{hi}}\bigr)
T_{\mathfrak m}\bigl(\varphi_\alpha^{\mathrm{hi}}\bigr)
\right]
\frac{d\beta}{\beta}\frac{d\alpha}{\alpha}.
\end{split}
\end{equation}
Moreover, after integration in $\lambda$, the contribution of
$H_{\rho,\lambda}$ in $T_{\mathcal{E}^l_{\rho}}\eta$ satisfies the flux identity
\begin{equation}\label{rhoremainder}
\begin{split}
\overline u_{\mu}^{x_0}P_\mu\left(
\int_{\lambda\approx\mu}
P_{<M\delta\lambda}(T_LH_{\rho,\lambda})\eta_\lambda
\frac{d\lambda}{\lambda}
\right)=
\nabla_x\cdot Q_{\varphi,\mu}
+O_{L_{t,x}^1}\left(
\epsilon^4c_\mu^2\mu^{-\frac12-2s_c}
(1+\mu|x_0|)^{\frac12}
\right),
\end{split}
\end{equation}
where $Q_{\varphi,\mu}$ satisfies \eqref{QPhibound}.
\end{lemma}
The integrals defining $\mathcal{N}_{\rho,\lambda}^{\mathrm{null}}$ are a form of high--high paraproducts, but we prefer to write them in the expanded form above so as not to conflict with the established convention in Subsection \ref{multpara}.
\begin{proof}
To isolate the required null structure, we begin by rewriting $\rho_{<\delta\lambda}^1$ as follows:
\begin{equation}\label{rho1decomp1}
\rho^1_{<\delta\lambda}
=P_{<\delta\lambda}\left(
\frac12(B^{\mathrm{hi}})^2
-\frac12|V^{\mathrm{hi}}|^2
-B^{\mathrm{hi}}V^{\mathrm{hi}}\cdot\nabla\eta
\right)+H^1_{\rho,\lambda},
\end{equation}
where
\begin{equation}\label{initialrho1err}
\begin{split}
H^1_{\rho,\lambda}
={}&
\frac12\mathcal B_{<\delta\lambda}
(V^{\mathrm{lo}},V^{\mathrm{lo}})
-\frac12\mathcal B_{<\delta\lambda}
(B^{\mathrm{lo}},B^{\mathrm{lo}})
+\mathcal S_{<\delta\lambda}
(B^{\mathrm{lo}},V^{\mathrm{lo}},\nabla\eta)
\\
&-P_{<\delta\lambda}
(V^{\mathrm{lo}}\cdot V^{\mathrm{hi}})
+P_{<\delta\lambda}
(B^{\mathrm{lo}}B^{\mathrm{hi}})
\\
&+(P_{<\delta\lambda}B^{\mathrm{lo}})
P_{<\delta\lambda}(V^{\mathrm{hi}}\cdot\nabla\eta)
-P_{<\delta\lambda}
(B^{\mathrm{lo}}V^{\mathrm{hi}}\cdot\nabla\eta)
\\
&+(P_{<\delta\lambda}V^{\mathrm{lo}})\cdot
P_{<\delta\lambda}(B^{\mathrm{hi}}\nabla\eta)
-P_{<\delta\lambda}
(B^{\mathrm{hi}}V^{\mathrm{lo}}\cdot\nabla\eta).
\end{split}
\end{equation}
$H^1_{\rho,\lambda}$ contributes to the purely perturbative $L_{t,x}^1$ part of \eqref{rhoremainder}. We briefly explain this below.
\begin{itemize}
\item In the first row of \eqref{initialrho1err}, the $B$ and $V$ factors are localized to frequencies $\ll\lambda$. Combining \eqref{Srep} and \eqref{bilformrep}, we use two unbalanced bilinear estimates, pairing the exterior high-frequency factors $\eta$ and $u_{\mu}^{x_0}$ each with one low-frequency $B$ or $V$ factor. This contributes a perturbative error term.

\item Since
$P_{<\delta\lambda}V^{\mathrm{hi}}
=P_{<\delta\lambda}B^{\mathrm{hi}}=0$, the second row consists of commutators
\begin{equation*}
-[P_{<\delta\lambda},V^{\mathrm{lo}}]\cdot V^{\mathrm{hi}}
+[P_{<\delta\lambda},B^{\mathrm{lo}}]B^{\mathrm{hi}}.
\end{equation*}
The low output in each term forces $V^{\mathrm{hi}}$ and $B^{\mathrm{hi}}$ to frequency comparable to the respective low-frequency factors, so both inputs are in fact at frequency $\ll\lambda$. Proposition \ref{Leibnizcommutator} therefore gives a representation of the form
\begin{equation*}
\lambda^{-2}\mathcal L^l
(\nabla P_{\ll\lambda}V,\nabla P_{\ll\lambda}V)
+\lambda^{-2}\mathcal L^l
(\nabla P_{\ll\lambda}B,\nabla P_{\ll\lambda}B),
\end{equation*}
whose contribution is again controlled by two unbalanced bilinear estimates.

\item The last two rows of \eqref{initialrho1err} are
\begin{equation*}
\begin{split}
&\mathcal B_{<\delta\lambda}
(B^{\mathrm{lo}},V^{\mathrm{hi}}\cdot\nabla\eta)
+\mathcal B_{<\delta\lambda}
(V^{\mathrm{lo}},B^{\mathrm{hi}}\nabla\eta)
\\
&\quad=
\lambda^{-2}\mathcal L^l
(\nabla B^{\mathrm{lo}},\nabla V^{\mathrm{hi}},
\nabla\eta^{\mathrm{lo}})
+\lambda^{-2}\mathcal L^l
(\nabla V^{\mathrm{lo}},\nabla B^{\mathrm{hi}},
\nabla\eta^{\mathrm{lo}})
\\
&\qquad+
\lambda^{-1}\mathcal L^l
(\nabla B^{\mathrm{lo}},V^{\mathrm{hi}},
\nabla\eta^{\mathrm{hi}})
+\lambda^{-1}\mathcal L^l
(\nabla V^{\mathrm{lo}},B^{\mathrm{hi}},
\nabla\eta^{\mathrm{hi}}).
\end{split}
\end{equation*}
In the first two terms on the right, the low output forces the differentiated $B$ and $V$ factors to have frequency $\ll\lambda$, and two unbalanced bilinear estimates apply. In the last two terms, we place the two high-frequency factors in $L_t^4L_x^\infty$, pair the differentiated low-frequency factor with $\eta_\lambda$ in $L_{t,x}^2$, and estimate $\overline u_\mu^{x_0}$ in $L_t^\infty L_x^2$.
\end{itemize}
Thus all the terms in $H^1_{\rho,\lambda}$ give a perturbative contribution, as claimed.

Next, we replace $B^{\mathrm{hi}}$ and $V^{\mathrm{hi}}$ in the main term in \eqref{rho1decomp1} by their principal components. First, we separate the high-frequency part of $\nabla\eta$ in the last factor and write this whole expression as
\begin{equation}\label{lowhinullsplit}
\begin{split}
P_{<\delta\lambda}\left(
\frac12(B^{\mathrm{hi}})^2
-\frac12|V^{\mathrm{hi}}|^2
-B^{\mathrm{hi}}V^{\mathrm{hi}}\cdot\nabla\eta_{<\delta\lambda}
\right)-P_{<\delta\lambda}
(B^{\mathrm{hi}}V^{\mathrm{hi}}\cdot\nabla\eta_{\geq\delta\lambda}).
\end{split}
\end{equation}
In the first term, each of the quadratic factors involving $V$ and $B$ must be at comparable frequency (i.e. they are high--high paraproducts). We will use this fact shortly. First, we observe that for every $\alpha\gtrsim_{M,\delta}\lambda$, Corollary \ref{principalVBexp} gives
\begin{equation*}
\begin{split}
B_\alpha^{\mathrm{hi}}=T_{\mathfrak l}\varphi_\alpha^{\mathrm{hi}}
+O_{L_t^2L_x^\infty}
(\epsilon^2\alpha^{-\frac34}),\hspace{5mm}V_\alpha^{\mathrm{hi}}=T_{\mathfrak m}\varphi_\alpha^{\mathrm{hi}}
+O_{L_t^2L_x^\infty}
(\epsilon^2\alpha^{-\frac34}).
\end{split}
\end{equation*}
Here we have commuted both the projections through the principal operators using
\begin{equation*}
\begin{split}
P_\alpha P_{\geq M\delta\lambda}(T_{\mathfrak l}\varphi)
&=T_{\mathfrak l}\varphi_\alpha^{\mathrm{hi}}
+[P_\alpha P_{\geq M\delta\lambda},T_{\mathfrak l}]\varphi,\\
P_\alpha P_{\geq M\delta\lambda}(T_{\mathfrak m}\varphi)
&=T_{\mathfrak m}\varphi_\alpha^{\mathrm{hi}}
+[P_\alpha P_{\geq M\delta\lambda},T_{\mathfrak m}]\varphi.
\end{split}
\end{equation*}
The commutators contain a differentiated paracoefficient and satisfy the same perturbative bound by Proposition \ref{Leibnizcommutator} (after expanding each paradifferential operator). The portions coming from the high-frequency projection are confined to the transition range $\alpha\approx M\delta\lambda$, modulo coefficient tails, which are controlled by Proposition \ref{Bernsteinelliptic}. From Proposition \ref{phiu-bootstrap} and Proposition \ref{VBunbil}, we have
\begin{equation*}
\begin{split}
&\|B_\alpha^{\mathrm{hi}}\|_{L_t^\infty L_x^2}
+\|T_{\mathfrak l}\varphi_\alpha^{\mathrm{hi}}\|_{L_t^\infty L_x^2}
+\|V_\alpha^{\mathrm{hi}}\|_{L_t^\infty L_x^2}
+\|T_{\mathfrak m}\varphi_\alpha^{\mathrm{hi}}\|_{L_t^\infty L_x^2}
\lesssim\epsilon\alpha^{-\frac32}.
\end{split}
\end{equation*}
Combining this with
$\|\nabla\eta_{<\delta\lambda}\|_{L_{t,x}^\infty}\lesssim\epsilon$ and the high-high paraproduct structure in \eqref{lowhinullsplit},
this shows that the total replacement error in the first term in \eqref{lowhinullsplit} is bounded in $L_{t,x}^2$ by
\begin{equation*}
\epsilon^3
\sum_{\alpha\gtrsim_{M,\delta}\lambda}
\alpha^{-\frac94}
\lesssim_{M,\delta}
\epsilon^3\lambda^{-\frac94}.
\end{equation*}
The retained terms are precisely $\mathcal N_{\rho,\lambda}^{\mathrm{null}}$ as defined in \eqref{nullformdef}.
Consequently, we have the decomposition
\begin{equation*}
\begin{split}
P_{<\delta\lambda}\left(
\frac12(B^{\mathrm{hi}})^2
-\frac12|V^{\mathrm{hi}}|^2
-B^{\mathrm{hi}}V^{\mathrm{hi}}\cdot\nabla\eta
\right)=
\mathcal N_{\rho,\lambda}^{\mathrm{null}}
+\mathcal N_{\rho,\lambda}^{\mathrm{rep}}
-P_{<\delta\lambda}
(B^{\mathrm{hi}}V^{\mathrm{hi}}
\cdot\nabla\eta_{\geq\delta\lambda}),
\end{split}
\end{equation*}
where $\mathcal N_{\rho,\lambda}^{\mathrm{rep}}$ denotes the replacement errors above. We next note that the last term on the right satisfies the same $L_{t,x}^2$ bound as the replacement error. Indeed,
\begin{equation*}
\begin{split}
\left\|
P_{<\delta\lambda}
(B^{\mathrm{hi}}V^{\mathrm{hi}}
\cdot\nabla\eta_{\geq\delta\lambda})
\right\|_{L_{t,x}^2}\lesssim
\|B^{\mathrm{hi}}\|_{L_t^4L_x^\infty}
\|V^{\mathrm{hi}}\|_{L_t^4L_x^\infty}
\|\nabla\eta_{\geq\delta\lambda}\|_{L_t^\infty L_x^2}
\lesssim
\epsilon^3\lambda^{-\frac94}.
\end{split}
\end{equation*}
Thus, we can write $\mathcal N_{\rho,\lambda}^{\mathrm{err}}
\coloneq
\mathcal N_{\rho,\lambda}^{\mathrm{rep}}
-P_{<\delta\lambda}
(B^{\mathrm{hi}}V^{\mathrm{hi}}
\cdot\nabla\eta_{\geq\delta\lambda})$, where $\mathcal N_{\rho,\lambda}^{\mathrm{err}}$ satisfies the bound
\begin{equation}\label{rho1rem2}
\|\mathcal N_{\rho,\lambda}^{\mathrm{err}}\|_{L_{t,x}^2}
\lesssim
\epsilon^2\lambda^{-\frac94}.
\end{equation}
We finally verify that the remaining error term gives the required flux contribution (plus a perturbative error). We can write
\begin{equation*}
\begin{split}
\overline u_\mu^{x_0}P_\mu\left(
P_{<M\delta\lambda}
(T_L\mathcal N_{\rho,\lambda}^{\mathrm{err}})
\eta_\lambda
\right)&=
\nabla_x\cdot
\mathcal L\left(
\overline u_\mu^{x_0},\eta_\lambda,
|D|^{-1}T_L\mathcal N_{\rho,\lambda}^{\mathrm{err}}
\right)
\\
&+
\mathcal L\left(
|D|P_{<\mu/4}\mathcal L(\overline u_\mu^{x_0},\eta_\lambda),
|D|^{-1}T_L\mathcal N_{\rho,\lambda}^{\mathrm{err}}
\right).
\end{split}
\end{equation*}
Here the high frequency output of the balanced pair has been absorbed
into the divergence term above. The second term is controlled by Lemma \ref{balancedrestrictedoutput}
and \eqref{rho1rem2}:
\begin{equation*}
\begin{split}
&\mu^{\frac12}
\sup_{\nu\in2^{\mathbb Z}}
\left\|
P_\nu|D|^{\frac12}P_{<\mu/4}
\mathcal L(\overline u_\mu^{x_0},\eta_\lambda)
\right\|_{L_{t,x}^2}
\|\mathcal N_{\rho,\lambda}^{\mathrm{err}}\|_{L_{t,x}^2}\lesssim
\epsilon^4c_\mu^2\mu^{-\frac12-2s_c}
(1+\mu|x_0|)^{\frac12},
\end{split}
\end{equation*}
where we used that $\lambda\approx\mu$. The divergence term contributes
to $Q_{\varphi,\mu}$ and satisfies \eqref{QPhibound} by H\"older's
inequality, the $L_t^4L_x^q$ bounds (with $q>8$ and using the hypothesis $p>\frac{4}{3}$) for $u_\mu^{x_0}$ and
$\eta_\lambda$, and \eqref{rho1rem2}. Taking
\begin{equation*}
H_{\rho,\lambda}
\coloneq H^1_{\rho,\lambda}
+\mathcal N_{\rho,\lambda}^{\mathrm{err}}
\end{equation*}
and integrating in $\lambda$ proves \eqref{rhoremainder} and completes the proof.
\end{proof}

It remains to show that the contribution of the null term
$\mathcal N_{\rho,\lambda}^{\mathrm{null}}$ is either perturbative or produces a flux term of the same form as above. Arguing as in the reduction to \eqref{rho1rem2}, it suffices to prove
\begin{equation}\label{nullL2bound}
\|\mathcal N_{\rho,\lambda}^{\mathrm{null}}\|_{L_{t,x}^2}
\lesssim
\epsilon^2\lambda^{-\frac94}.
\end{equation}
For this, we have to study its symbol more closely. Consider a pair of
comparable input frequencies $\beta\leq\alpha$, with
$\delta\alpha\leq\beta\leq\alpha$. The paradifferential operators in
\eqref{nullformdef} initially give coefficients whose symbols depend on
$\eta_{<\delta\alpha}$ or $\eta_{<\delta\beta}$. We first replace these by
coefficients depending on $\eta_{<\delta\lambda}$, which has the fixed
$\lambda$ scale. By the fundamental theorem of calculus, each difference
introduces an integral over a frequency parameter $\gamma$ lying between
$\delta\alpha$ (or $\delta\beta$) and $\delta\lambda$, together with the linearization
of the symbol, which produces a factor $\nabla\eta_\gamma$. The two principal
operators present in each quadratic factor in
$\mathcal N_{\rho,\lambda}^{\mathrm{null}}$ each have order one. Therefore,
the contribution of the error for each quadratic term is bounded in
$L_{t,x}^2$ by
\begin{equation*}
\begin{split}
\alpha^2
\|\nabla\eta_\gamma\|_{L_t^\infty L_x^2}
\|\varphi_\alpha^{\mathrm{hi}}\|_{L_t^4L_x^\infty}
\|\varphi_\beta^{\mathrm{hi}}\|_{L_t^4L_x^\infty}
\lesssim_\delta
\epsilon^3\alpha^{-\frac54}\gamma^{-1}.
\end{split}
\end{equation*}
Since $\alpha\gtrsim_{M,\delta}\lambda$, integration in $\gamma$ gives
\begin{equation*}
\epsilon^3\alpha^{-\frac54}
\left|
\int_{\delta\alpha}^{\delta\lambda}
\gamma^{-1}\frac{d\gamma}{\gamma}
\right|
\lesssim_{M,\delta}
\epsilon^3\lambda^{-1}\alpha^{-\frac54}.
\end{equation*}
The same bound holds with $\alpha$ replaced by $\beta$, and integrating over the
comparable input frequencies $\alpha$ and $\beta\gtrsim_{M,\delta}\lambda$ gives an
$O_{L_{t,x}^2}(\epsilon^3\lambda^{-\frac94})$ error. To estimate the remaining
part, we can therefore introduce the frozen principal symbols
\begin{equation}
\label{eq:reduced-lm-symbols}
\mathfrak l_{<\delta\lambda}(x,\xi)
\coloneq
\frac{
\mathfrak g_{<\delta\lambda}(x,\xi)
+i\nabla\eta_{<\delta\lambda}(x)\cdot\xi
}{
1+|\nabla\eta_{<\delta\lambda}(x)|^2
},
\qquad
\mathfrak m_{<\delta\lambda}(x,\xi)
\coloneq
i\xi-
\nabla\eta_{<\delta\lambda}(x)
\mathfrak l_{<\delta\lambda}(x,\xi).
\end{equation}
If $f$ and $g$ are localized to comparable annuli at frequencies $\alpha$ and
$\beta$, where $\delta\alpha\leq\beta\leq\alpha$ and
$\alpha\gtrsim_{M,\delta}\lambda$, define the reduced bilinear form
\begin{equation}
\label{eq:reduced-null-expression}
\begin{aligned}
\mathcal N^{\mathrm{red}}_{\rho,\lambda}(f,g)
\coloneq{}&
\frac12
\bigl(T_{\mathfrak l_{<\delta\lambda}}f\bigr)
\bigl(T_{\mathfrak l_{<\delta\lambda}}g\bigr)
-\frac12
\bigl(T_{\mathfrak m_{<\delta\lambda}}f\bigr)\cdot
\bigl(T_{\mathfrak m_{<\delta\lambda}}g\bigr)
\\
&-\frac12\nabla\eta_{<\delta\lambda}\cdot
\Bigl[
\bigl(T_{\mathfrak l_{<\delta\lambda}}f\bigr)
\bigl(T_{\mathfrak m_{<\delta\lambda}}g\bigr)
+\bigl(T_{\mathfrak l_{<\delta\lambda}}g\bigr)
\bigl(T_{\mathfrak m_{<\delta\lambda}}f\bigr)
\Bigr].
\end{aligned}
\end{equation}
Its bilinear symbol is
\begin{equation*}
\begin{split}
\mathfrak n_{\rho,\lambda}(x,\xi_1,\xi_2)
\coloneq{}&
\frac12
\mathfrak l_{<\delta\lambda}(x,\xi_1)
\mathfrak l_{<\delta\lambda}(x,\xi_2)
\\
&-
\frac12
\mathfrak m_{<\delta\lambda}(x,\xi_1)
\cdot
\mathfrak m_{<\delta\lambda}(x,\xi_2)
\\
&-
\frac12
\mathfrak l_{<\delta\lambda}(x,\xi_1)
\nabla\eta_{<\delta\lambda}(x)\cdot
\mathfrak m_{<\delta\lambda}(x,\xi_2)
\\
&-
\frac12
\mathfrak l_{<\delta\lambda}(x,\xi_2)
\nabla\eta_{<\delta\lambda}(x)\cdot
\mathfrak m_{<\delta\lambda}(x,\xi_1).
\end{split}
\end{equation*}
A direct calculation using \eqref{eq:reduced-lm-symbols} gives
\begin{equation*}
\mathfrak n_{\rho,\lambda}(x,\xi_1,\xi_2)
=
\frac12\left[
\xi_1\cdot\xi_2
+\bigl(1+|\nabla\eta_{<\delta\lambda}|^2\bigr)
\mathfrak l_{<\delta\lambda}(x,\xi_1)
\mathfrak l_{<\delta\lambda}(x,\xi_2)
\right].
\end{equation*}
Moreover,
\begin{equation*}
\left(\mathfrak g_{<\delta\lambda}(x,\xi)\right)^2
=
\bigl(1+|\nabla\eta_{<\delta\lambda}|^2\bigr)|\xi|^2
-\bigl(\nabla\eta_{<\delta\lambda}\cdot\xi\bigr)^2,
\end{equation*}
and hence
\begin{equation*}
\bigl(1+|\nabla\eta_{<\delta\lambda}|^2\bigr)
\mathfrak l_{<\delta\lambda}(x,\xi)
\mathfrak l_{<\delta\lambda}(x,-\xi)
=|\xi|^2.
\end{equation*}
In particular, $\mathfrak{n}_{\rho,\lambda}$ vanishes at opposite input frequencies,
\begin{equation*}
\mathfrak n_{\rho,\lambda}(x,\xi,-\xi)=0.
\end{equation*}

For input functions at frequency $\alpha,\beta\gtrsim M\delta\lambda$, the output projection $P_{<\delta\lambda}$ then forces $|\xi_1+\xi_2|\lesssim M^{-1}\alpha$. Consequently, on the high--high-to-low region
\begin{equation*}
|\xi_1|\approx|\xi_2|\approx\alpha,
\qquad
|\xi_1+\xi_2|\leq c_1\alpha,
\end{equation*}
where $c_1>0$ is sufficiently small, the fundamental theorem of calculus gives (in this region) the factorization
\begin{equation}
\label{nullfactorization}
\mathfrak n_{\rho,\lambda}(x,\xi_1,\xi_2)
=i\sum_j
(\xi_1+\xi_2)_j
\mathfrak q_{j,\alpha}(x,\xi_1,\xi_2),
\end{equation}
where we have inserted a factor of $i$ so that the corresponding operator is a divergence. The symbols $\mathfrak q_{j,\alpha}$ are of order one and satisfy
\begin{equation*}
\left|
\partial_{\xi_1}^{\gamma_1}
\partial_{\xi_2}^{\gamma_2}
\mathfrak q_{j,\alpha}(x,\xi_1,\xi_2)
\right|
\lesssim_{\gamma_1,\gamma_2}
C(\|\nabla\eta\|_{L_{t,x}^\infty})
\alpha^{1-|\gamma_1|-|\gamma_2|},
\end{equation*}
and
\begin{equation*}
\left|
\partial_{\xi_1}^{\gamma_1}
\partial_{\xi_2}^{\gamma_2}
\nabla_x\mathfrak q_{j,\alpha}(x,\xi_1,\xi_2)
\right|
\lesssim_{\gamma_1,\gamma_2}
C(\|\nabla\eta\|_{L_{t,x}^\infty})
|\nabla^2\eta_{<\delta\lambda}(x)|
\alpha^{1-|\gamma_1|-|\gamma_2|}.
\end{equation*}
We denote the corresponding bilinear operator by
$\mathcal Q_{j,\alpha}$.

We now decompose the null form into its comparable frequency pairs. The reduction to  $\lambda$-scale para-coefficients above gives
\begin{equation}\label{nullfixedgeometry}
\begin{split}
\mathcal N_{\rho,\lambda}^{\mathrm{null}}&=2\int_0^\infty\int_{\delta\alpha}^{\alpha}
P_{<\delta\lambda}
\mathcal N^{\mathrm{red}}_{\rho,\lambda}
\bigl(\varphi_\alpha^{\mathrm{hi}},
\varphi_\beta^{\mathrm{hi}}\bigr)
\frac{d\beta}{\beta}\frac{d\alpha}{\alpha}+O_{L_{t,x}^2}
(\epsilon^3\lambda^{-\frac94}).
\end{split}
\end{equation}
The factor $2$ comes from the normalization in
\eqref{eq:reduced-null-expression}, and the integrand vanishes unless
$\alpha,\beta\gtrsim_{M,\delta}\lambda$. For every dyadic $\kappa$ and
$|\xi|\approx\alpha$, we have the normalized symbol bound
\begin{equation*}
\begin{split}
\alpha^{-1}\Bigl(
&\|\partial_x\mathfrak l_{<\kappa}(\cdot,\xi)\|_{L_t^\infty L_x^2}
+\|\partial_x\mathfrak m_{<\kappa}(\cdot,\xi)\|_{L_t^\infty L_x^2}
\Bigr)
\lesssim
\|\nabla^2\eta\|_{L_t^\infty L_x^2}
\lesssim\epsilon.
\end{split}
\end{equation*}
Since $\mathfrak q_{j,\alpha}$ depends on $x$,
\eqref{nullfactorization} gives, uniformly for
$\delta\alpha\leq\beta\leq\alpha$,
\begin{equation*}
\begin{split}
P_{<\delta\lambda}
\mathcal N^{\mathrm{red}}_{\rho,\lambda}
\bigl(\varphi_\alpha^{\mathrm{hi}},
\varphi_\beta^{\mathrm{hi}}\bigr)=
\sum_j\partial_jP_{<\delta\lambda}
\mathcal Q_{j,\alpha}
\bigl(\varphi_\alpha^{\mathrm{hi}},
\varphi_\beta^{\mathrm{hi}}\bigr)
+E_{\alpha,\beta,\lambda},
\end{split}
\end{equation*}
where $E_{\alpha,\beta,\lambda}$ consists of the terms in which the derivative
falls on $\mathfrak q_{j,\alpha}$. The preceding symbol bound, two
$L_t^4L_x^\infty$ estimates, and one energy estimate give
\begin{equation*}
\|E_{\alpha,\beta,\lambda}\|_{L_{t,x}^2}
\lesssim_\delta
\epsilon^3\alpha^{-\frac94}.
\end{equation*}
We can then write the (normalized) bilinear operator $\mathcal{Q}_{j,\alpha}$ as
\begin{equation*}
\bigl[\alpha^{-1}\mathcal Q_{j,\alpha}\bigr](f,g)(x)
=\int_{\mathbb R^2}\int_{\mathbb R^2}
k_{j,\alpha}(t,x;y,z)\tau_yf(x)\tau_zg(x)\,dy\,dz.
\end{equation*}
The symbol bounds above imply the rapid decay bounds
\begin{equation*}
\begin{split}
|k_{j,\alpha}(t,x;y,z)|
&\lesssim_{N,\delta}
\alpha^4\langle\alpha y\rangle^{-N}
\langle\alpha z\rangle^{-N},
\\
\|\nabla_xk_{j,\alpha}(\cdot,\cdot;y,z)\|_{L_t^\infty L_x^2}
&\lesssim_{N,\delta}
\epsilon\alpha^4\langle\alpha y\rangle^{-N}
\langle\alpha z\rangle^{-N}.
\end{split}
\end{equation*}
Our goal is now to apply the balanced bilinear estimate to $\partial_j P_{<\delta\lambda}\mathcal{Q}_{j,\alpha}(\varphi_{\alpha}^{hi},\varphi_{\beta}^{hi})$. Since this kernel depends on $x$, we first have to commute the full output operator through it before applying the balanced estimate. Thus,
\begin{equation*}
\begin{split}
\partial_jP_{<\delta\lambda}
\mathcal Q_{j,\alpha}(f,g)&=\alpha\int_{\mathbb R^2}\int_{\mathbb R^2}
k_{j,\alpha}(t,x;y,z)
\partial_jP_{<\delta\lambda}
\bigl(\tau_yf\,\tau_zg\bigr)(x)\,dy\,dz
\\
&+\alpha\int_{\mathbb R^2}\int_{\mathbb R^2}
\bigl[\partial_jP_{<\delta\lambda},
k_{j,\alpha}(t,\cdot;y,z)\bigr]
\bigl(\tau_yf\,\tau_zg\bigr)(x)\,dy\,dz.
\end{split}
\end{equation*}
Applying Lemma \ref{balancedrestrictedoutput} under the kernel integral, and then using Bernstein's inequality and summation over the output frequencies $\nu\lesssim\delta\lambda$, gives, for $f=\varphi_\alpha^{\mathrm{hi}}$ and $g=\varphi_\beta^{\mathrm{hi}}$,
\begin{equation*}
\begin{split}
\left\|
\alpha\int_{\mathbb R^2}\int_{\mathbb R^2}
k_{j,\alpha}(t,x;y,z)
\partial_jP_{<\delta\lambda}
\bigl(\tau_yf\,\tau_zg\bigr)(x)\,dy\,dz
\right\|_{L_{t,x}^2}\lesssim_\delta
\epsilon^2\lambda^{\frac12}\alpha^{-\frac{11}{4}}.
\end{split}
\end{equation*}
For the second term, Proposition \ref{Leibnizcommutator} gives
\begin{equation*}
\bigl\|[\partial_jP_{<\delta\lambda},k]F\bigr\|_{L_{t,x}^2}
\lesssim
\|\nabla_xk\|_{L_t^\infty L_x^2}
\|F\|_{L_t^2L_x^\infty}.
\end{equation*}
Therefore, using two $L_t^4L_x^\infty$ bounds for $f$ and $g$ and then integrating in $y$ and $z$, and using the rapid decay in $y$ and $z$ of $|\nabla_x k_{j,\alpha}(\cdot,\cdot;y,z)\|_{L_{t}^{\infty}L_x^2}$, we obtain
\begin{equation*}
\begin{split}
&\left\|
\alpha\int_{\mathbb R^2}\int_{\mathbb R^2}
\bigl[\partial_jP_{<\delta\lambda},
k_{j,\alpha}(t,\cdot;y,z)\bigr]
\bigl(\tau_yf\,\tau_zg\bigr)(x)\,dy\,dz
\right\|_{L_{t,x}^2}
\lesssim_\delta
\epsilon^3\alpha^{-\frac94}.
\end{split}
\end{equation*}
This commutator contribution is perturbative and is estimated together with $E_{\alpha,\beta,\lambda}$ below.
Combining this with \eqref{nullfixedgeometry}, integrating first in $\beta$
over $\delta\alpha\leq\beta\leq\alpha$, and then in
$\alpha\gtrsim_{M,\delta}\lambda$, we conclude that
\begin{equation*}
\begin{split}
\|\mathcal N_{\rho,\lambda}^{\mathrm{null}}\|_{L_{t,x}^2}
&\lesssim_{M,\delta}
\epsilon^2\lambda^{\frac12}
\int_{\alpha\gtrsim_{M,\delta}\lambda}
\alpha^{-\frac{11}{4}}\frac{d\alpha}{\alpha}
+\epsilon^3
\int_{\alpha\gtrsim_{M,\delta}\lambda}
\alpha^{-\frac94}\frac{d\alpha}{\alpha}
+\epsilon^3\lambda^{-\frac94}
\\
&\lesssim_{M,\delta}
\epsilon^2\lambda^{-\frac94},
\end{split}
\end{equation*}
which proves \eqref{nullL2bound}. The same flux reduction used for
$\mathcal N_{\rho,\lambda}^{\mathrm{err}}$ now shows that the contribution of
$\mathcal N_{\rho,\lambda}^{\mathrm{null}}$ is a flux term satisfying
\eqref{QPhibound}, modulo a remainder satisfying \eqref{paraphiremainder}.
Together with Lemma \ref{nullprincipal} and the estimate for
$\rho^2_{<\delta\lambda}$, this completes the proof of Lemma \ref{aerrorinitial}.
\end{proof}
To conclude the proof of Proposition \ref{paravarphiinitial}, it remains to prove the three auxiliary bounds \eqref{NphiinitialL2Linf}, \eqref{NphiinitialauxR} and \eqref{NphiinitialL2}. Here, it is not necessary to work with the above flux decomposition. Instead, we can work directly with the full error
\begin{equation}\label{Nphiinitialexacterror}
\begin{split}
N_\varphi^{\mathrm{err}}-N_\varphi
=\left(T_{(V\cdot\nabla B)^l}
-T_{V^l}\cdot T_{\nabla B^l}\right)\eta
-T_{\mathcal E^l}\eta.
\end{split}
\end{equation}
We begin by recording a few basic bounds for the para-coefficients in this error term. First, from Proposition \ref{basicopbounds} and Proposition \ref{VBunbil} and Proposition \ref{Bernsteinelliptic} and dyadic summation, we have the bounds 
\begin{equation*}
\begin{split}
\|\nabla B_{<\delta\lambda_0}\|_{L_t^4L_x^\infty}
+\|\nabla V_{<\delta\lambda_0}\|_{L_t^4L_x^\infty}
&\lesssim\epsilon\lambda_0^{\frac38},\hspace{5mm}\|L_{<\delta\lambda_0}B_{<\delta\lambda_0}\|_{L_{t,x}^\infty}\lesssim\epsilon\lambda_0^{\frac12}.
\end{split}
\end{equation*}
Together with the bound for
$\kappa_{<\delta\lambda_0}$ above and \eqref{LrhoauxL2Linf}, this gives
\begin{equation}\label{Ephiinitialcoefficient}
\|\mathcal E^l_{<\delta\lambda_0}\|_{L_t^2L_x^\infty}
\lesssim
\epsilon^2\lambda_0^{-\frac14}
+\epsilon^3\lambda_0^{-\frac14}
\lesssim\epsilon^2\lambda_0^{-\frac14}.
\end{equation}

The paracommutator in \eqref{Nphiinitialexacterror} contributes similar coefficient bounds. Indeed, its integrand, when expanding each of the paradifferential operators above is given by
\begin{equation}\label{phiparacomauxexact}
\begin{split}
&P_{<M\delta\lambda_0}\left(
(V_{<\delta\lambda_0}-V^l_{<\delta\lambda_0})
\cdot\nabla B_{<\delta\lambda_0}\right)\eta_{\lambda_0}
+[P_{<M\delta\lambda_0},V^l_{<\delta\lambda_0}]
\cdot\nabla B_{<\delta\lambda_0}\,\eta_{\lambda_0}
\\
&+
\int_{\lambda_1\approx\lambda_0}
(V^l_{<\delta\lambda_0}-V^l_{<\delta\lambda_1})
\cdot S_{\lambda_1}
\left(\nabla B^l_{<\delta\lambda_0}\,\eta_{\lambda_0}\right)
\frac{d\lambda_1}{\lambda_1}.
\end{split}
\end{equation}
Proposition \ref{Bernsteinelliptic} and Proposition \ref{VBunbil} gives the high frequency tail bound
\begin{equation*}
\|V_{<\delta\lambda_0}-V^l_{<\delta\lambda_0}\|_{L_t^4L_x^\infty}
\lesssim\epsilon\lambda_0^{-\frac58}.
\end{equation*}
Also, \eqref{dBdVL4p} and Bernstein and Proposition \ref{VBunbil} gives
\begin{equation*}
\|V^l_{<\delta\lambda_0}-V^l_{<\delta\lambda_1}
\|_{L_t^4L_x^\infty}
\lesssim\epsilon\lambda_0^{-\frac58},
\qquad\lambda_1\approx\lambda_0.
\end{equation*}
Combining this with Proposition \ref{Leibnizcommutator} and Proposition \ref{VBunbil} and Bernstein, we can estimate the full expression \eqref{phiparacomauxexact} in $L_t^2L_x^{\infty}$ by
\begin{equation*}
\epsilon^2\lambda_0^{-\frac{1}{4}}\|\eta_{\lambda_0}\|_{L_{t,x}^{\infty}}\lesssim \epsilon^3\lambda_0^{-\frac{5}{4}}.    
\end{equation*}
Consequently, by integrating over the region $\lambda_0\approx\lambda$, and using the above estimates, we obtain the bound \eqref{NphiinitialL2Linf} for $N_\varphi^{\mathrm{err}}-N_\varphi$.
\begin{equation}\label{NphiinitialerrorL2Linf}
\|(N_\varphi^{\mathrm{err}}-N_\varphi)_\lambda
\|_{L_t^2L_x^\infty}
\lesssim\epsilon^3\lambda^{-\frac54}.
\end{equation}
The same argument also proves \eqref{NphiinitialL2}. Indeed,
using \eqref{Ephiinitialcoefficient} and the coefficient bounds
for \eqref{phiparacomauxexact}, and placing $\eta_{\lambda_0}$
in $L_t^\infty L_x^2$ instead of $L_{t,x}^{\infty}$, we obtain
\begin{equation*}
\begin{split}
\|(N_\varphi^{\mathrm{err}}-N_\varphi)_\lambda\|_{L_{t,x}^2}\lesssim_\delta\epsilon^2
\int_{\lambda_0\approx\lambda}\lambda_0^{-\frac14}
\|\eta_{\lambda_0}\|_{L_t^\infty L_x^2}
\frac{d\lambda_0}{\lambda_0}\lesssim_\delta\epsilon^3c_\lambda\lambda^{-\frac14-s_c}.
\end{split}
\end{equation*}
Here only $\lambda_0\approx\lambda$ contributes, since the
coefficients have frequency $\ll\lambda_0$.
Next, we turn to obtaining the multilinear estimate \eqref{NphiinitialauxR}. This is a relatively simple consequence of the coefficient bounds above. First, from the representation \eqref{phiparacomauxexact}, and placing the coefficients involving $V$ and $B$ in $L_t^4L_x^{\infty}$ (using the above bounds), estimating one factor of $u_{\mu}^{x_0}$ in $L_t^{\infty}L_x^2$ and using one unbalanced $L^2$ bound for the pair $(u_{\mu},\eta_{\lambda_0})$, we can estimate (for $\mu\gg\lambda$)
\begin{equation*}
\begin{split}
\|P_{\lambda}\left((T_{(V\cdot\nabla B)^l}-T_{V^l}\cdot T_{\nabla B^l})\eta\right)\overline{u}_{\mu}^{x_0}u_{\mu}\|_{L_{t,x}^1}&\lesssim \epsilon^2\lambda^{-\frac{1}{4}}\|u_{\mu}^{x_0}\|_{L_{t}^{\infty}L_x^2}\|L(u_{\mu},\tilde{P}_{\lambda}\eta)\|_{L_{t,x}^2}
\\
&\lesssim \epsilon^5\lambda^{-\frac{3}{2}}\mu^{-2s_c}c_{\mu}^2   
\end{split}
\end{equation*}
where $\tilde{P}_{\lambda}$ is a fattened projection at scale $\lambda$. On the other hand, using \eqref{Ephiinitialcoefficient}, we have 
\begin{equation*}
\|\overline{u}_{\mu}^{x_0}u_{\mu}P_{\lambda}T_{\mathcal{E}^l}\eta\|_{L_{t,x}^1}\lesssim \epsilon^2\lambda^{-\frac{1}{4}}\|u_{\mu}^{x_0}\|_{L_{t}^{\infty}L_x^2}\|L(u_{\mu},\tilde{P}_{\lambda}\eta)\|_{L_{t,x}^2}
\end{equation*}
which gives the same bound. Combining the above gives
\begin{equation*}
\begin{split}
\bigl\|(N_\varphi^{\mathrm{err}}-N_\varphi)_\lambda
\overline u_\mu^{x_0}u_\mu\bigr\|_{L_{t,x}^1}\lesssim
\epsilon^5c_\mu^2\lambda^{-\frac32}\mu^{-2s_c}.
\end{split}
\end{equation*}
This proves the second bound \eqref{NphiinitialauxR}.
Combining Lemmas \ref{paracomperturbative} and \ref{aerrorinitial} with
\eqref{initialphipara}, and absorbing the the flux terms into
$Q_{\varphi,\mu}$, gives \eqref{eq:geometric-truncated-good-unknown},
\eqref{paraphiremainder}, and \eqref{QPhibound}. This completes the proof of
Proposition \ref{paravarphiinitial}.
\end{proof}


\subsection{Multilinear expansions for \texorpdfstring{$N_{\varphi}$}{}}\label{subsec:phimultilinear}
As in the expansion for $G(\eta)\psi$, we now aim to obtain a precise multilinear expansion for the remaining term $N_{\varphi}$ defined in
\eqref{Nphidef}. We use throughout the same definition of $\varphi$ and
the same continuous paraproduct convention as above.

We first introduce the class of quartic terms that will arise in the
expansion. As in Definition \ref{N4def}, we will not need to keep track
of their exact form, but just describe their generic structure.

\begin{definition}\label{Nphi4def}
For each $\lambda>0$, we use
$N_{\varphi,\lambda}^{[4]}$ to denote any finite linear combination of
balanced quartic paradifferential expressions with the following
properties:
\begin{itemize}
\item There are exactly four explicit inputs, consisting of two factors
of $\eta$ and two factors of $\varphi$, all localized at frequencies
comparable to $\lambda$, and the output is also localized at frequency
comparable to $\lambda$.
\item The expression is scalar-valued and is obtained from pointwise
products and contractions by composing $T_{\mathfrak l}$ and
$T_{\mathfrak m}$ with translation-invariant multilinear forms. It has
total differential order four.
\end{itemize}
The low-frequency dependence of the symbols on $\eta$ is regarded as a
paracoefficient and is not counted among the four explicit inputs. We
also allow bounded integrations over comparable frequency parameters.
For example, if $\lambda_j\approx\lambda$ for $1\leq j\leq4$, then
\begin{equation*}
P_\lambda\left(
\eta_{\lambda_1}\eta_{\lambda_2}
(T_{\mathfrak l}^2\varphi_{\lambda_3})
(T_{\mathfrak l}^2\varphi_{\lambda_4})
\right)
\end{equation*}
is a term of the form $N_{\varphi,\lambda}^{[4]}$. We use
$N_{\varphi}^{[4]}$ to denote a sum over $\lambda$ of finitely many
terms of this type. 
\end{definition}

Compared with Definition \ref{N4def}, one explicit factor of $\eta$ is
replaced by a second factor of $\varphi$. This accounts for the
additional half derivative in the corresponding estimates.
We use $N_{\varphi,B}^{[4]}$ below to denote the sum of the terms of
this form arising from the expansion of the $B$ component.

To precisely identify the cubic $N_{\varphi}^{[3]}$ components, we will adopt the symmetric trilinear high--high--high paraproduct convention from
\eqref{trilinearPiconvention}. For the nested terms, we use
the analogous convention for labelling each input: for
$\mathfrak a\in\{\mathfrak l,\mathfrak m\}$,
\begin{equation}\label{NphinestedPiconvention}
\Pi\bigl(F,T_{\mathfrak a}(G,H)\bigr)
\coloneq
\int_{\mathcal D}
\sum_{\sigma\in \mathfrak{S}_3}
F_{\lambda_{\sigma(0)}}
T_{\mathfrak a}\left(
G_{\lambda_{\sigma(1)}}H_{\lambda_{\sigma(2)}}
\right)d\Omega,
\end{equation}
which contains all six permutations of the three inputs, as in
\eqref{nestedPiconvention}. 

\subsubsection{Expansion of $\frac{1}{2}B^2-T_{B^l}B$.}
We begin with the expansion of
$\frac{1}{2}B^2-T_{B^l}B$. We claim that the quadratic balanced and cubic balanced principal parts are respectively given by
\begin{equation}\label{NphiB23def}
\begin{split}
N_{\varphi,B}^{[2]}\coloneq
\frac12\Pi(T_{\mathfrak l}\varphi,T_{\mathfrak l}\varphi),\hspace{5mm}N_{\varphi,B}^{[3]}\coloneq
\Pi(T_{\mathfrak l}\varphi,\eta,T_{\mathfrak l}^2\varphi)
-\Pi\bigl(
T_{\mathfrak l}\varphi,
T_{\mathfrak l}(\eta,T_{\mathfrak l}\varphi)
\bigr).
\end{split}
\end{equation}
This is made precise by the following proposition.
\begin{proposition}\label{B-TBexp}
Set $\mathcal B\coloneq\frac12B^2-T_{B^l}B.$ We have
\begin{equation*}
\mathcal B
=N_{\varphi,B}^{[2]}+N_{\varphi,B}^{[3]}
+N_{\varphi,B}^{[4]}+\mathcal R,
\end{equation*}
where $N_{\varphi,B}^{[4]}$ is of the form described in
Definition \ref{Nphi4def}. The following bounds hold:
\begin{enumerate}
\item For every $\mu>0$ and $x_0\in\mathbb R^2$, the remainder
satisfies
\begin{equation*}
\|\mathcal R_\mu\overline u_\mu^{x_0}\|_{L_{t,x}^1}
\lesssim
\epsilon^4c_\mu^2\mu^{-2s_c-\frac12}
(1+\mu|x_0|)^{\frac12}.
\end{equation*}

\item 
For every $\lambda>0$, we have
\begin{equation}\label{BauxL2Linf}
\|\mathcal B_\lambda\|_{L_t^2L_x^\infty}
\lesssim_\delta\epsilon^2\lambda^{-\frac54}.
\end{equation}

\item 
For every $0<\lambda\ll\mu$ and $x_0\in\mathbb R^2$, we have
a decomposition
\begin{equation}\label{Bauxdecomp}
\mathcal B_\lambda\overline u_\mu^{x_0}u_\mu
=\mathcal R_B^{\mathrm{aux}}+\nabla\cdot Q_B^{\mathrm{aux}},
\end{equation}
where, for $\frac43<p\leq2$,
\begin{equation}\label{Bauxtrilinear}
\|\mathcal R_B^{\mathrm{aux}}\|_{L_{t,x}^1}
+\lambda^{\frac2p-1}
\|Q_B^{\mathrm{aux}}\|_{L_t^1L_x^p}
\lesssim_{\delta,p}
\epsilon^4c_\mu^2\lambda^{-\frac32}
\mu^{-2s_c}(1+\mu|x_0|)^{\frac12}.
\end{equation}
\item For every $\lambda>0$, we have
\begin{equation}\label{BauxL2cubic}
\|(\mathcal B-N_{\varphi,B}^{[2]})_\lambda\|_{L_{t,x}^2}
\lesssim_\delta\epsilon^3c_\lambda\lambda^{-\frac14-s_c}.
\end{equation}
\end{enumerate}
\end{proposition}
As in the expansions in earlier subsections, the latter three bounds will be needed to measure error terms coming from the time derivative of certain normal form type corrections to the density-flux identities in our analysis later.
\begin{proof}
By FTC, we have
\begin{equation*}
\frac{1}{2}B^2
=\int_0^\infty B^{\mathrm{lin}}_{\lambda_0}B_{<\lambda_0}
\frac{d\lambda_0}{\lambda_0}
\end{equation*}
and
\begin{equation*}
T_{B^l}B
=\int_0^\infty
B^l_{<\delta\lambda_0}B^{\mathrm{lin}}_{\lambda_0}
\frac{d\lambda_0}{\lambda_0},
\end{equation*}
where, by the linearization formula for $B$, we define
\begin{equation*}
B^{\mathrm{lin}}_{\lambda_j}
\coloneq
\lambda_j\frac{d}{d\lambda_j}B_{<\lambda_j}
=L_{<\lambda_j}\Phi^j_{<\lambda_j}
+\eta_{\lambda_j}L_{<\lambda_j}B_{<\lambda_j}.
\end{equation*}
Hence, we have
\begin{equation}\label{Bminuspara}
\begin{split}
\frac{1}{2}B^2-T_{B^l}B=\int_0^\infty\int_{\delta\lambda_0}^{\lambda_0}
B^{\mathrm{lin}}_{\lambda_1}B^{\mathrm{lin}}_{\lambda_0}
\frac{d\lambda_1}{\lambda_1}
\frac{d\lambda_0}{\lambda_0}+\int_0^\infty
B_{\lambda_0}^{\mathrm{lin}}
P_{\geq M\delta\lambda_0}B_{<\delta\lambda_0}
\frac{d\lambda_0}{\lambda_0}.
\end{split}
\end{equation}
Similarly to the proof of Proposition \ref{firstdecomposition}, we
first localize each factor $B_{\lambda_j}^{\mathrm{lin}}$ to output
frequency $\lesssim_M\lambda_j$, and estimate any high-frequency tail errors perturbatively. More precisely, write
\begin{equation*}
B_{\lambda_j}^{\mathrm{lin}}
=B_{\lambda_j}^{\mathrm{loc}}+B_{\lambda_j}^{\mathrm{err}},
\end{equation*}
where
\begin{equation*}
B_{\lambda_j}^{\mathrm{loc}}
\coloneq
L^l_{<\lambda_j}\Phi^{l,j}_{<\lambda_j}
+\eta_{\lambda_j}L^l_{<\lambda_j}B^l_{<\lambda_j}
\end{equation*}
and
\begin{equation*}
\begin{split}
B_{\lambda_j}^{\mathrm{err}}
\coloneq{}&
P_{\geq M\lambda_j}L_{<\lambda_j}
\Phi^{l,j}_{<\lambda_j}
-L_{<\lambda_j}\left(
(P_{\geq M\lambda_j}B_{<\lambda_j})\eta_{\lambda_j}
\right)
\\
&+\eta_{\lambda_j}P_{\geq M\lambda_j}
L_{<\lambda_j}B^l_{<\lambda_j}
+\eta_{\lambda_j}L_{<\lambda_j}
P_{\geq M\lambda_j}B_{<\lambda_j}.
\end{split}
\end{equation*}
Thus $B_{\lambda_j}^{\mathrm{loc}}$ is localized to frequency
$\lesssim_M\lambda_j$, while every term in
$B_{\lambda_j}^{\mathrm{err}}$ is a
high-frequency tail.
It follows from \eqref{Bminuspara} that
\begin{equation}\label{principalpartBpara}
\frac{1}{2}B^2-T_{B^l}B
=\int_0^\infty\int_{\delta\lambda_0}^{\lambda_0}
B_{\lambda_1}^{\mathrm{loc}}B_{\lambda_0}^{\mathrm{loc}}
\frac{d\lambda_1}{\lambda_1}
\frac{d\lambda_0}{\lambda_0}
+\mathcal R,
\end{equation}
where
\begin{equation*}
\begin{split}
\mathcal R\coloneq
\int_0^\infty
B_{\lambda_0}^{\mathrm{lin}}
P_{\geq M\delta\lambda_0}B_{<\delta\lambda_0}
\frac{d\lambda_0}{\lambda_0}+\int_0^\infty\int_{\delta\lambda_0}^{\lambda_0}
\left(
B^{\mathrm{lin}}_{\lambda_1}B^{\mathrm{lin}}_{\lambda_0}
-B^{\mathrm{loc}}_{\lambda_1}B^{\mathrm{loc}}_{\lambda_0}
\right)
\frac{d\lambda_1}{\lambda_1}
\frac{d\lambda_0}{\lambda_0}.
\end{split}
\end{equation*}
We claim that this error term is perturbative, and even lies directly in $L_t^1L_x^2$ with the bound
\begin{equation}\label{BparaRbound1}
\|\mathcal R_\mu\|_{L_t^1L_x^2}
\lesssim c_\mu\epsilon^4\mu^{-\frac12-s_c}.
\end{equation}
To see this, fix an integer $N>\max\{s_c,3\}$ and then choose
$M\gg_N 1$. Arguing as in Proposition \ref{firstdecomposition}, and
using Bernstein's inequality, Proposition \ref{basicopbounds}, Proposition \ref{Bernsteinelliptic},
Lemma \ref{Phiest}, and Proposition \ref{VBunbil}, we obtain, for
every $\sigma\geq0$,
\begin{equation*}
\begin{split}
\||D|^\sigma B_\lambda^{\mathrm{err}}\|_{L_t^{\frac43}\dot B^0_{2,1}}
&\lesssim_{N,\sigma}
c_\lambda\epsilon^N\lambda^{\frac18+\sigma-s_c},\\
\||D|^\sigma B_\lambda^{\mathrm{err}}\|_{L_t^{\frac43}\dot B^0_{2,1}}
&\lesssim_{N,\sigma}
\epsilon^N\lambda^{\sigma-\frac{15}{8}},\\
\||D|^\sigma P_{\geq M\lambda}B_{<\lambda}\|_{L_t^{\frac43}\dot B^0_{2,1}}
&\lesssim_{N,\sigma}
\epsilon^N\lambda^{\sigma-\frac{15}{8}}.
\end{split}
\end{equation*}
We also have the estimates
\begin{equation}\label{Bloclinbounds}
\begin{split}
\||D|^\sigma B_\lambda^{\mathrm{lin}}\|_{L_t^4\dot B^0_{\infty,1}}
+\||D|^\sigma B_\lambda^{\mathrm{loc}}\|_{L_t^4\dot B^0_{\infty,1}}
&\lesssim_\sigma
\epsilon c_\lambda\lambda^{\frac{11}{8}+\sigma-s_c},\\
\||D|^\sigma B_\lambda^{\mathrm{lin}}\|_{L_t^4\dot B^0_{\infty,1}}
+\||D|^\sigma B_\lambda^{\mathrm{loc}}\|_{L_t^4\dot B^0_{\infty,1}}
&\lesssim_\sigma
\epsilon\lambda^{\sigma-\frac58}.
\end{split}
\end{equation}
By Bernstein's inequality, Proposition \ref{besovparabounds}, and the preceding bounds,
\begin{equation*}
\begin{split}
\|P_\mu|D|^\sigma(
B_{\lambda_0}^{\mathrm{lin}}
P_{\geq M\delta\lambda_0}B_{<\delta\lambda_0})
\|_{L_t^1L_x^2}\lesssim_N
\epsilon^N c_{\lambda_0}
\lambda_0^{-\frac12+\sigma-s_c},
\end{split}
\end{equation*}
and, since $\lambda_1\approx_\delta\lambda_0$,
\begin{equation*}
\begin{split}
\|P_\mu|D|^\sigma(
B^{\mathrm{lin}}_{\lambda_1}B^{\mathrm{lin}}_{\lambda_0}
-B^{\mathrm{loc}}_{\lambda_1}B^{\mathrm{loc}}_{\lambda_0})
\|_{L_t^1L_x^2}\lesssim_{N,\sigma}
\epsilon^N c_{\lambda_0}
\lambda_0^{-\frac12+\sigma-s_c}.
\end{split}
\end{equation*}
Minkowski's inequality, Bernstein's inequality, and slow variation of
the frequency envelope now give \eqref{BparaRbound1}. 
\\

Before expanding the principal term in \eqref{principalpartBpara},
it will be convenient at this point to prove the easier auxiliary bounds \eqref{BauxL2Linf}, \eqref{BauxL2cubic} and
\eqref{Bauxtrilinear}. For simplicity, write
$H_\mu=\overline u_\mu^{x_0}u_\mu$. We first consider the high frequency tail errors in
\eqref{principalpartBpara}, denoting their sum temporarily by
$\mathcal E_B$. A similar argument to the $L_t^1L_x^2$ bound
above gives the two bounds
\begin{equation}\label{EBtailbounds}
\|\mathcal E_{B,\lambda}\|_{L_t^1L_x^2}
\lesssim_N\epsilon^N\lambda^{-\frac52},
\hspace{5mm}
\|\mathcal E_{B,\lambda}\|_{L_t^2L_x^\infty}
\lesssim\epsilon^2\lambda^{-\frac54}.
\end{equation}
The $L_t^2L_x^\infty$ bound handles the contribution of
$\mathcal E_{B,\lambda}$ in \eqref{BauxL2Linf}, while the
$L_t^1L_x^2$ bound gives
\begin{equation*}
\begin{split}
\|\mathcal E_{B,\lambda}H_\mu\|_{L_{t,x}^1}
&\lesssim
\lambda\|\mathcal E_{B,\lambda}\|_{L_t^1L_x^2}
\|u_\mu\|_{L_t^\infty L_x^2}^2\lesssim_N
\epsilon^{N+2}c_\mu^2\lambda^{-\frac32}\mu^{-2s_c},
\end{split}
\end{equation*}
which is more than sufficient. We next consider the main double
integral in \eqref{principalpartBpara}. Its $L_t^2L_x^\infty$
bound follows immediately from using two
$L_t^4\dot B_{\infty,1}^0\subset L_t^4L_x^\infty$ bounds in
\eqref{Bloclinbounds} at $s=s_c=2$,
\begin{equation*}
\begin{split}
\left\|P_\lambda
\int_0^\infty\int_{\delta\lambda_0}^{\lambda_0}
B_{\lambda_1}^{\mathrm{loc}}B_{\lambda_0}^{\mathrm{loc}}
\frac{d\lambda_1}{\lambda_1}
\frac{d\lambda_0}{\lambda_0}\right\|_{L_t^2L_x^\infty}\lesssim_\delta
\epsilon^2
\int_{\lambda\lesssim_M\lambda_0}
\lambda_0^{-\frac54}\frac{d\lambda_0}{\lambda_0}\lesssim_\delta\epsilon^2\lambda^{-\frac54}.
\end{split}
\end{equation*}
Together with the second bound in \eqref{EBtailbounds}, this
proves \eqref{BauxL2Linf}. For \eqref{Bauxtrilinear}, the argument is similar to the proof
of \eqref{Netaauxdecomp}. We split the integral in the first
term on the right hand side of \eqref{principalpartBpara} into
a low and high frequency part,
$F_\lambda^{\mathrm{lo}}+F_\lambda^{\mathrm{hi}}$, according to
whether $\lambda_0\ll\mu$ or $\lambda_0\gtrsim\mu$.
For the low-frequency contribution, we pair each factor of
$B_{\lambda_j}^{\mathrm{loc}}$ with either $u_\mu$ or
$\overline u_\mu^{x_0}$ and use two unbalanced bilinear estimates.
Indeed, using Lemma \ref{Phiest}, Proposition \ref{paraexp},
Proposition \ref{VBunbil} and Proposition \ref{basicopbounds},
it is easy to show
\begin{equation}\label{Blambdalocbound}
\begin{split}
B_{\lambda_j}^{\mathrm{loc}}
=T_{\mathfrak l}\varphi_{\lambda_j}
+O_{L_t^2L_x^\infty}
(\epsilon^2\lambda_j^{-\frac34}),\hspace{5mm}B_{\lambda_j}^{\mathrm{loc}}
=T_{\mathfrak l}\varphi_{\lambda_j}
+O_{L_t^\infty L_x^2}
(\epsilon^2\lambda_j^{-\frac32}).
\end{split}
\end{equation}
For the two principal factors of $T_{\mathfrak{l}}\varphi_{\lambda_0}$ and $T_{\mathfrak{l}}\varphi_{\lambda_1}$, we use
\eqref{phiu-bi-unbal} and the kernel bounds for $P_\lambda$
and $T_{\mathfrak l}$.
If one factor is an error $B_{\lambda_j}^{\mathrm{loc}}-T_{\mathfrak{l}}\varphi_{\lambda_j}$, we place it in $L_t^2L_x^\infty$,
use one unbalanced bilinear estimate for the remaining
principal factor, and place the unpaired $u_\mu$ factor in
$L_t^\infty L_x^2$. If both factors are errors, we use two
$L_t^2L_x^\infty$ bounds and place both $u_\mu$ factors in
$L_t^\infty L_x^2$. Combining Proposition \ref{phiu-bootstrap} with H\"older's
inequality and the kernel bounds for $P_\lambda$, we then obtain
\begin{equation*}
\begin{split}
\|F_\lambda^{\mathrm{lo}}H_\mu\|_{L_{t,x}^1}
&\lesssim_\delta
\epsilon^4c_\mu^2\mu^{-2s_c}(1+\mu|x_0|)^{\frac12}
\int_{\lambda\lesssim_M\lambda_0\ll\mu}
\lambda_0^{-\frac32}
\left[
\left(1+\frac{\lambda_0}{\lambda}\right)
+\epsilon\left(1+\frac{\lambda_0}{\lambda}\right)^{\frac12}
+\epsilon^2
\right]
\frac{d\lambda_0}{\lambda_0}
\\
&\lesssim_\delta
\epsilon^4c_\mu^2\lambda^{-\frac32}
\mu^{-2s_c}(1+\mu|x_0|)^{\frac12}.
\end{split}
\end{equation*}
Here the translations from the outer kernel have scale
$\lambda^{-1}$, so each principal pairing costs
$\left(1+\lambda_0/\lambda\right)^{\frac12}$ (
$\lambda_1\approx \lambda_0$). The three terms in brackets correspond to zero, one,
and two $B_{\lambda_j}^{\mathrm{loc}}-T_{\mathfrak{l}}\varphi_{\lambda_j}$ errors respectively.

For the high-frequency contribution, we claim that
\begin{equation}\label{BauxhighL2}
\|F_\lambda^{\mathrm{hi}}\|_{L_{t,x}^2}
\lesssim_\delta
\epsilon^2\lambda^{-\frac12}\mu^{-\frac74}.
\end{equation}
First, we observe that by \eqref{Blambdalocbound} and
Proposition \ref{phiu-bootstrap}, we have
\begin{equation*}
\begin{split}
F_\lambda^{\mathrm{hi}}
-P_\lambda
\int_{\lambda_0\gtrsim\mu}
\int_{\delta\lambda_0}^{\lambda_0}
T_{\mathfrak l}\varphi_{\lambda_0}
T_{\mathfrak l}\varphi_{\lambda_1}
\frac{d\lambda_1}{\lambda_1}
\frac{d\lambda_0}{\lambda_0}
=O_{L_{t,x}^2}(\epsilon^3\mu^{-\frac94}),
\end{split}
\end{equation*}
which is controlled by the right-hand side of
\eqref{BauxhighL2}. Applying Bernstein in $\lambda$ and Lemma
\ref{balancedrestrictedoutput} to the product
$P_\lambda(T_{\mathfrak l}\varphi_{\lambda_0}
T_{\mathfrak l}\varphi_{\lambda_1})$, and using the error
bound above then gives \eqref{BauxhighL2}.

The proof now proceeds in essentially the same way as
\eqref{Netaauxdecomp}. Indeed, we first recall the balanced
estimate for $H_\mu$,
\begin{equation*}
\|P_\rho H_\mu\|_{L_{t,x}^2}
\lesssim
\epsilon^2c_\mu^2\rho^{-\frac12}
\mu^{\frac54-2s_c}(1+\mu|x_0|)^{\frac12}.
\end{equation*}
As in the proof of \eqref{Netaauxdecomp}, let $\delta_0>0$ be a
sufficiently small fixed constant. Frequency summation and
Bernstein give
\begin{equation*}
\begin{split}
\|P_{\geq\delta_0\lambda}H_\mu\|_{L_{t,x}^2}
&\lesssim
\epsilon^2c_\mu^2\lambda^{-\frac12}
\mu^{\frac54-2s_c}(1+\mu|x_0|)^{\frac12},
\\
\lambda^{\frac2p-1}
\|P_{<\delta_0\lambda}H_\mu\|_{L_t^2L_x^{\frac{2p}{2-p}}}
&\lesssim_p
\epsilon^2c_\mu^2\lambda^{\frac12}
\mu^{\frac54-2s_c}(1+\mu|x_0|)^{\frac12}.
\end{split}
\end{equation*}
and we used $p>\frac43$ to ensure the low-frequency sum
converges. We then simply define
\begin{equation*}
Q_B^{\mathrm{aux}}
=\nabla\Delta^{-1}
\bigl(F_\lambda^{\mathrm{hi}}
P_{<\delta_0\lambda}H_\mu\bigr).
\end{equation*}
The expression inside the parentheses is localized at frequency
comparable to $\lambda$. Consequently, \eqref{BauxhighL2} and the
preceding estimates imply
\begin{equation*}
\begin{split}
\|F_\lambda^{\mathrm{hi}}
P_{\geq\delta_0\lambda}H_\mu\|_{L_{t,x}^1}
+\lambda^{\frac2p-1}
\|Q_B^{\mathrm{aux}}\|_{L_t^1L_x^p}
&\lesssim_{\delta,p}
\epsilon^4c_\mu^2\lambda^{-1}
\mu^{-\frac12-2s_c}(1+\mu|x_0|)^{\frac12},
\end{split}
\end{equation*}
which is even stronger than \eqref{Bauxtrilinear}, since
$\lambda\ll\mu$. Combining the low-frequency contribution and
the preceding bounds proves \eqref{Bauxtrilinear}. 
\\

Finally, for \eqref{BauxL2cubic}, we begin again by using \eqref{principalpartBpara}. Similarly to \eqref{BparaRbound1}, we can estimate the remainder term in \eqref{principalpartBpara} in $L_{t,x}^2$. More precisely, we have from \eqref{Bernsteinelliptic}, and arguing as in the estimates following \eqref{principalpartBpara}, the inital expansion
\begin{equation}\label{principalpartBparaaux}
\mathcal{B}_{\lambda}
=P_{\lambda}\int_0^\infty\int_{\delta\lambda_0}^{\lambda_0}
B_{\lambda_1}^{\mathrm{loc}}B_{\lambda_0}^{\mathrm{loc}}
\frac{d\lambda_1}{\lambda_1}
\frac{d\lambda_0}{\lambda_0}
+O_{L_{t,x}^2}(\epsilon^4\lambda^{-\frac{1}{4}-s_c}).
\end{equation}
By using \eqref{Blambdalocbound} and then Proposition \ref{Leibnizcommutator}, we can peel off the principal part of each $B_{\lambda_j}^{\mathrm{loc}}$ factor to write
\begin{equation}\label{principalpartBparaaux2}
\mathcal{B}_{\lambda}
=P_{\lambda}\int_0^\infty\int_{\delta\lambda_0}^{\lambda_0}
S_{\lambda_1}(T_{\mathfrak{l}}\varphi)S_{\lambda_0}(T_{\mathfrak{l}}\varphi)
\frac{d\lambda_1}{\lambda_1}
\frac{d\lambda_0}{\lambda_0}
+O_{L_{t,x}^2}(\epsilon^4\lambda^{-\frac{1}{4}-s_c})
\end{equation}
where the errors from the commutator can be estimated using two
$L_t^4L_x^\infty$ bounds and placing the other factor in
$L_t^\infty L_x^2$. The remaining principal term on the right-hand side is precisely
$N_{\varphi,B}^{[2]}$. This proves \eqref{BauxL2cubic}.

Now, we turn to establish the main multilinear expansion in the statement of Proposition \ref{B-TBexp}. For the remainder of the proof, $\mathcal R$ denotes any term whose
contribution satisfies the paired $L_{t,x}^1$ estimate in Proposition
\ref{B-TBexp}, and it may change from line to line. Likewise,
$N_{\varphi,B}^{[4]}$ may denote different accumulated sums of terms
of the form described in Definition \ref{Nphi4def} until the final
collection.

We now expand the principal term on the right-hand side of
\eqref{principalpartBpara}. For positive parameters $\rho_0$ and
$\rho_1$, define
\begin{equation}\label{BBexp}
\begin{split}
\mathcal T^B_{<\rho_0,\rho_1}
&\coloneq
\bigl(L^{l,\lambda_0}_{<\rho_0}\Phi^{l,0}_{<\rho_0}\bigr)
\bigl(L^{l,\lambda_1}_{<\rho_1}\Phi^{l,1}_{<\rho_1}\bigr)
\\
&+\eta_{\lambda_0}
L^{l,\lambda_0}_{<\rho_0}B^l_{<\rho_0}
L^{l,\lambda_1}_{<\rho_1}\Phi^{l,1}_{<\rho_1}
\\
&+\eta_{\lambda_1}
L^{l,\lambda_1}_{<\rho_1}B^l_{<\rho_1}
L^{l,\lambda_0}_{<\rho_0}\Phi^{l,0}_{<\rho_0}
\\
&+\eta_{\lambda_0}\eta_{\lambda_1}
L^{l,\lambda_0}_{<\rho_0}B^l_{<\rho_0}
L^{l,\lambda_1}_{<\rho_1}B^l_{<\rho_1}
\\
\coloneq{}&
\sum_{j=1}^4\mathcal T^{B,j}_{<\rho_0,\rho_1}.
\end{split}
\end{equation}
Thus, we have $B_{\lambda_1}^{\mathrm{loc}}B_{\lambda_0}^{\mathrm{loc}}
=\mathcal T^B_{<\lambda_0,\lambda_1}.$
\\

Above, and for the remainder of this subsection, we use the notation
$A_{<\lambda}^{l,\beta}=P_{<M\beta}A_{<\lambda}$ to fix the output
projection independently of the parameter $\lambda$. As
in the preceding expansions in the $\eta$ equation, derivatives of the
"matched" projections occurring in $B^l_{<\rho}$, as well as projection
differences which are high-frequency tails in view of the relevant
frequency relations, are controlled by Proposition
\ref{Bernsteinelliptic}.
\\

Now, we proceed with expanding each term $\mathcal{T}^{B,j}_{<\lambda_0,\lambda_1}$. We begin with the last term
$\mathcal T^{B,4}_{<\lambda_0,\lambda_1}$, which is the simplest and does not contribute any non-perturbative quadratic or cubic balanced terms. 
First, we have for $j=0,1$,
\begin{equation}\label{LBetabil}
\begin{split}
\|\eta_{\lambda_j}
L^{l,\lambda_j}_{<\delta\lambda_j}B^l_{<\delta\lambda_j}\|_{L_{t,x}^2}&\lesssim
\|\eta_{\lambda_j}T_LB^l_{<\delta\lambda_j}\|_{L_{t,x}^2}+
\|\eta_{\lambda_j}
(L^{l,\lambda_j}_{<\delta\lambda_j}-T_L)
B^l_{<\delta\lambda_j}\|_{L_{t,x}^2}
\\
&\lesssim
\epsilon^2c_{\lambda_j}
\lambda_j^{\frac14-s_c}.
\end{split}
\end{equation}
For the first term, we used Proposition \ref{VBunbil} after extracting
the principal part with Proposition \ref{paraexp} and using $L_t^4\dot{B}_{\infty,1}^0$ estimates for the subprincipal error term. The second term is
estimated by the same argument as for \eqref{L2TL}, also using Proposition \ref{Bernsteinelliptic} to control any high-frequency tails. By the fundamental theorem of calculus, \eqref{LBetabil}, and the
bounds
\begin{equation*}
\|L_{<\rho}B^l_{<\rho}\|_{L_t^4L_x^\infty}
\lesssim\epsilon\rho^{\frac38},
\qquad
\|\eta_\rho\|_{L_t^4L_x^\infty}
\lesssim\epsilon\rho^{-\frac98},
\end{equation*}
we obtain
\begin{equation}\label{TB4int}
\begin{split}
\mathcal T^{B,4}_{<\lambda_0,\lambda_1}
&=\eta_{\lambda_0}\eta_{\lambda_1}
\int_{\delta\lambda_0}^{\lambda_0}
\int_{\delta\lambda_1}^{\lambda_1}
\left(
\lambda_2\frac{d}{d\lambda_2}
L^{l,\lambda_0}_{<\lambda_2}B^l_{<\lambda_2}
\right)\left(
\lambda_3\frac{d}{d\lambda_3}
L^{l,\lambda_1}_{<\lambda_3}B^l_{<\lambda_3}
\right)
\frac{d\lambda_3}{\lambda_3}
\frac{d\lambda_2}{\lambda_2}
\\
&+O_{L_t^1L_x^2}
(\epsilon^4c_{\lambda_0}\lambda_0^{-\frac12-s_c}).
\end{split}
\end{equation}
Here the $L_t^1L_x^2$ error contains $\mathcal{T}_{<\rho_0,\rho_1}^{B,4}$ evaluated at either $\rho_0=\delta\lambda_0$, $\rho_1=\delta\lambda_1$ or both, together with the
high-frequency tail terms described above.

Modulo high-frequency tail terms, the linearization formulas for $L$
and $B$ give
\begin{equation}\label{LBlinearizationNphi}
\begin{split}
\rho\frac{d}{d\rho}
\bigl(L^{l,\lambda_j}_{<\rho}B^l_{<\rho}\bigr)
=L^{l,\lambda_j}_{<\rho}L^l_{<\rho}
\Phi^l_{<\rho}(\psi_\rho,\eta_\rho)+\eta_\rho L^{l,\lambda_j}_{<\rho}
L^l_{<\rho}B^l_{<\rho}.
\end{split}
\end{equation}
More precisely, by Proposition \ref{Bernsteinelliptic}, the
contribution of the omitted tails to the preceding double integral is
controlled by
\begin{equation*}
O_{L_t^1L_x^2}
\left(
\epsilon^4c_{\lambda_0}\lambda_0^{-\frac12-s_c}
\right).
\end{equation*}
The contribution obtained by taking the first term in
\eqref{LBlinearizationNphi} for both differentiated factors can be
decomposed into a term of the form $N_{\varphi,B}^{[4]}$, terms that
are perturbative directly in $L_t^1L_x^2$, and terms satisfying the
usual $L_{t,x}^1$ paired remainder estimate. The procedure is very similar to
the one used to isolate the quartic terms in the expansion of
$N_\eta$ in Subsection \ref{subsec:Neta34}. We will provide the essential steps in the argument below for $\mathcal{T}^{B,4}$ and refer back to this set of steps when extracting other similar quartic terms that appear
in this section, but without providing all the details each time. We also refer to Subsection \ref{subsec:Neta34} for even more
comprehensive computations.
\\

We now give the relevant sequence of steps used to extract the
$N_{\varphi,B}^{[4]}$ term, which also applies (with minor modifications) to the remaining
$\mathcal T^{B,j}$ terms. 
First, any contribution containing at least one copy of the second
term in \eqref{LBlinearizationNphi},
\begin{equation*}
\eta_\rho L_{<\rho}^{l,\lambda_j}
L_{<\rho}^lB^l_{<\rho},
\end{equation*}
can be estimated in $L_t^2L_x^\infty$ using Proposition
\ref{basicopbounds} and Proposition \ref{VBunbil}. After removing this
block, the product of the exterior factors
$\eta_{\lambda_0}\eta_{\lambda_1}$ with the other differentiated
factor is at least cubic. Its leading cubic part is estimated using
two $L_t^4L_x^\infty$ bounds and one
$L_t^\infty L_x^2$ energy bound, with any additional factors placed
in $L_{t,x}^\infty$. Here we use Lemma \ref{Phiest} for the components
containing $\Phi^l_{<\rho}$.

It remains to consider the contribution obtained by taking two copies
of the first term in \eqref{LBlinearizationNphi} in \eqref{TB4int},
namely
\begin{equation}\label{TB4remaining}
\begin{split}
\mathcal I^{B,4}_{<\lambda_0,\lambda_1}
\coloneq{}&
\eta_{\lambda_0}\eta_{\lambda_1}
\int_{\delta\lambda_0}^{\lambda_0}
\int_{\delta\lambda_1}^{\lambda_1}
\left(
L^{l,\lambda_0}_{<\lambda_2}
L^l_{<\lambda_2}
\Phi^l_{<\lambda_2}
(\psi_{\lambda_2},\eta_{\lambda_2})
\right)
\\
&\qquad\times\left(
L^{l,\lambda_1}_{<\lambda_3}
L^l_{<\lambda_3}
\Phi^l_{<\lambda_3}
(\psi_{\lambda_3},\eta_{\lambda_3})
\right)
\frac{d\lambda_3}{\lambda_3}
\frac{d\lambda_2}{\lambda_2}.
\end{split}
\end{equation}
For $\rho=\lambda_2,\lambda_3$, we first replace each occurrence of
$\Phi^l_{<\rho}(\psi_\rho,\eta_\rho)$ in
\eqref{TB4remaining} by $\varphi_\rho$ using Lemma \ref{Phiest}. The
errors in these replacements are directly perturbative: we use the
two $L_t^4L_x^\infty$ bounds, place the other differentiated
factor in energy, and apply the bound in Lemma \ref{Phiest} to the
$\Phi-\varphi$ factor. We then extract the principal
$T_{\mathfrak l}$ part of each $L$ operator using Proposition
\ref{paraexp}. The nonprincipal errors either contain an additional
factor of $\eta$ or are high-frequency tails. In the former case the
complete expression is at least quintic and is estimated directly in
$L_t^1L_x^2$ using four $L_t^4L_x^\infty$ bounds and one
$L_t^\infty L_x^2$ energy bound, with any further factors placed in
$L_{t,x}^\infty$. We then estimate $u_\mu^{x_0}$ separately in
$L_t^\infty L_x^2$. The high-frequency tails are controlled by
Proposition \ref{Bernsteinelliptic}. Hence all the
errors produced by these replacements are included in $\mathcal R$.
The retained principal contribution in \eqref{TB4remaining} is therefore
\begin{equation}\label{T4Bsurviving}
\begin{split}
\mathcal J^{B,4}_{<\lambda_0,\lambda_1}
\coloneq{}&
\eta_{\lambda_0}\eta_{\lambda_1}
\int_{\delta\lambda_0}^{\lambda_0}
\int_{\delta\lambda_1}^{\lambda_1}
\bigl(T_{\mathfrak l}^2\varphi_{\lambda_2}\bigr)
\bigl(T_{\mathfrak l}^2\varphi_{\lambda_3}\bigr)
\frac{d\lambda_3}{\lambda_3}
\frac{d\lambda_2}{\lambda_2}.
\end{split}
\end{equation}
We now consider its output frequency $\mu$ relative to the inputs. Since
\begin{equation*}
\lambda_0\approx_\delta\lambda_1
\approx_\delta\lambda_2\approx_\delta\lambda_3,
\end{equation*}
all four inputs in \eqref{T4Bsurviving} have comparable frequency, and
the expression has total differential order four. The range of output frequencies
$\mu\gg_M\lambda_0$ is absent (by frequency support considerations), while the range
$\mu\approx_{M,\delta}\lambda_0$ is included in
$N_{\varphi,B}^{[4]}$. It remains to consider the range $\mu\ll_M\lambda_0$. Using the kernel
bounds for $T_{\mathfrak{l}}$, estimating the paracoefficients
in $L^\infty$, and pairing
$(\overline u_\mu^{x_0},\eta_{\lambda_0})$ with the unbalanced
estimate, we obtain
\begin{equation*}
\begin{split}
\|\overline u_\mu^{x_0}P_\mu
\mathcal J^{B,4}_{<\lambda_0,\lambda_1}\|_{L_{t,x}^1}&\lesssim
\int_{\delta\lambda_0}^{\lambda_0}
\int_{\delta\lambda_1}^{\lambda_1}
\lambda_2^2\lambda_3^2
\|\mathcal L(\overline u_\mu^{x_0},
\eta_{\lambda_0})\|_{L_{t,x}^2}\|\eta_{\lambda_1}\|_{L_t^4L_x^\infty}
\|\varphi_{\lambda_2}\|_{L_t^4L_x^\infty}
\|\varphi_{\lambda_3}\|_{L_t^\infty L_x^2}
\frac{d\lambda_3}{\lambda_3}
\frac{d\lambda_2}{\lambda_2}
\\
&\lesssim_\delta
\epsilon^5c_\mu c_{\lambda_0}
\mu^{\frac34-s_c}\lambda_0^{-\frac54-s_c}
(1+\mu|x_0|)^{\frac12}.
\end{split}
\end{equation*}
Here we used the estimates in Section \ref{sec:bootstrap} and Proposition \ref{phiu-bootstrap}. Finally, slow variation also gives
\begin{equation*}
\int_{\lambda_0\gg_M\mu}
c_{\lambda_0}\lambda_0^{-\frac54-s_c}
\frac{d\lambda_0}{\lambda_0}
\lesssim
c_\mu\mu^{-\frac54-s_c}.
\end{equation*}
Consequently,
\begin{equation*}
\begin{split}
\int_{\lambda_0\gg_M\mu}
\int_{\delta\lambda_0}^{\lambda_0}
\|\overline u_\mu^{x_0}P_\mu
\mathcal J^{B,4}_{<\lambda_0,\lambda_1}\|_{L_{t,x}^1}
\frac{d\lambda_1}{\lambda_1}
\frac{d\lambda_0}{\lambda_0}\lesssim
\epsilon^5c_\mu^2\mu^{-2s_c-\frac12}
(1+\mu|x_0|)^{\frac12},
\end{split}
\end{equation*}
which is controlled by the required remainder bound. Combining the
three output-frequency regimes with the errors estimated above, we
conclude that the contribution of
$\mathcal T^{B,4}_{<\lambda_0,\lambda_1}$ is included in
$N_{\varphi,B}^{[4]}+\mathcal R$. Slight variations of the above argument
will be used to extract the other principal quartic terms below.

We next expand $\mathcal T^{B,j}_{<\lambda_0,\lambda_1}$ for
$1\leq j\leq3$. By the fundamental theorem of calculus,
\begin{equation}\label{TBexp}
\begin{split}
\mathcal T^{B,j}_{<\lambda_0,\lambda_1}
=\mathcal T^{B,j}_{<\delta\lambda_0,\delta\lambda_1}
&+\int_{\delta\lambda_0}^{\lambda_0}
\lambda_2\frac{d}{d\lambda_2}
\mathcal T^{B,j}_{<\lambda_2,\delta\lambda_1}
\frac{d\lambda_2}{\lambda_2}+\int_{\delta\lambda_1}^{\lambda_1}
\lambda_3\frac{d}{d\lambda_3}
\mathcal T^{B,j}_{<\delta\lambda_0,\lambda_3}
\frac{d\lambda_3}{\lambda_3}
\\
&+\int_{\delta\lambda_0}^{\lambda_0}
\int_{\delta\lambda_1}^{\lambda_1}
\lambda_2\lambda_3
\frac{d}{d\lambda_2}\frac{d}{d\lambda_3}
\mathcal T^{B,j}_{<\lambda_2,\lambda_3}
\frac{d\lambda_3}{\lambda_3}
\frac{d\lambda_2}{\lambda_2}.
\end{split}
\end{equation}
As with the analysis of $\mathcal T_{<\lambda_0,\lambda_1}^{B,4}$, after isolating the principal quartic part, the expressions in the
double-derivative term have exactly two factors of $\eta$ and two
factors of $\varphi$ at comparable frequencies and total differential
order four. As in the previous term, we place the part of this term with output frequency $\mu\approx_{M,\delta}\lambda_0$ in the $N_{\varphi,B}^{[4]}$ group and the remaining part
in $\mathcal R$. It remains to analyze the lower endpoint and the two single-derivative
terms in \eqref{TBexp}. We claim that the following decomposition holds
\begin{equation}\label{Bquadprincipal}
\begin{split}
&\int_0^\infty\int_{\delta\lambda_0}^{\lambda_0}
\left(
\sum_{j=1}^3
\mathcal T^{B,j}_{<\delta\lambda_0,\delta\lambda_1}
-T_{\mathfrak l}\varphi_{\lambda_0}
T_{\mathfrak l}\varphi_{\lambda_1}
\right)
\frac{d\lambda_1}{\lambda_1}
\frac{d\lambda_0}{\lambda_0}
=\mathcal R.
\end{split}
\end{equation}
Indeed, Corollary \ref{Phil0varphibound} gives the bounds
\begin{equation}\label{L2phierr}
\begin{split}
\|\varphi_{\lambda_j}
-\Phi^{l,j}_{<\delta\lambda_j}\|_{L_{t,x}^2}
\lesssim
\epsilon^2c_{\lambda_j}
\lambda_j^{-\frac34-s_c},\hspace{5mm}\|\varphi_{\lambda_j}
-\Phi^{l,j}_{<\delta\lambda_j}\|_{L_{t,x}^2}
\lesssim
\epsilon^2\lambda_j^{-\frac{11}{4}}.
\end{split}
\end{equation}
Consequently, Proposition \ref{basicopbounds} and Bernstein's
inequality give
\begin{equation*}
\begin{split}
\|L^{l,\lambda_j}_{<\delta\lambda_j}
(\varphi_{\lambda_j}
-\Phi^{l,j}_{<\delta\lambda_j})\|_{L_{t,x}^2}\lesssim
\epsilon^2c_{\lambda_j}
\lambda_j^{\frac14-s_c},\hspace{5mm}
\|L^{l,\lambda_j}_{<\delta\lambda_j}
(\varphi_{\lambda_j}
-\Phi^{l,j}_{<\delta\lambda_j})\|_{L_{t,x}^2}\lesssim
\epsilon^2\lambda_j^{-\frac74}.
\end{split}
\end{equation*}
It follows that
\begin{equation*}
\begin{split}
\bigl(L^{l,\lambda_0}_{<\delta\lambda_0}
\Phi^{l,0}_{<\delta\lambda_0}\bigr)
\bigl(L^{l,\lambda_1}_{<\delta\lambda_1}
\Phi^{l,1}_{<\delta\lambda_1}\bigr)&=
L^{l,\lambda_0}_{<\delta\lambda_0}\varphi_{\lambda_0}
L^{l,\lambda_1}_{<\delta\lambda_1}\varphi_{\lambda_1}
\\
&+
L^{l,\lambda_0}_{<\delta\lambda_0}\varphi_{\lambda_0}
L^{l,\lambda_1}_{<\delta\lambda_1}
(\Phi^{l,1}_{<\delta\lambda_1}-\varphi_{\lambda_1})
\\
&+
L^{l,\lambda_0}_{<\delta\lambda_0}
(\Phi^{l,0}_{<\delta\lambda_0}-\varphi_{\lambda_0})
L^{l,\lambda_1}_{<\delta\lambda_1}\varphi_{\lambda_1}
\\
&+
O_{L_t^1L_x^2}
(\epsilon^4c_{\lambda_0}\lambda_0^{-\frac12-s_c}).
\end{split}
\end{equation*}
Using the representation \eqref{Phidifferror} and Proposition
\ref{paraexp}, we estimate the second and third terms by the
principal/subprincipal decomposition and bilinear pairings used
in the proof of \eqref{replacebyphi}, retaining the inherited
input and output localizations when applying Lemma
\ref{balancedrestrictedoutput}. Both terms therefore contribute
to $\mathcal R$.
\\

Proposition \ref{paraexp}, followed by the same pairings for the
subprincipal terms (as in the second and third terms above), gives
\begin{equation*}
\begin{split}
\int_0^\infty\int_{\delta\lambda_0}^{\lambda_0}
L^{l,\lambda_0}_{<\delta\lambda_0}\varphi_{\lambda_0}
L^{l,\lambda_1}_{<\delta\lambda_1}\varphi_{\lambda_1}
\frac{d\lambda_1}{\lambda_1}
\frac{d\lambda_0}{\lambda_0}=
\int_0^\infty\int_{\delta\lambda_0}^{\lambda_0}
T_{\mathfrak l}\varphi_{\lambda_0}
T_{\mathfrak l}\varphi_{\lambda_1}
\frac{d\lambda_1}{\lambda_1}
\frac{d\lambda_0}{\lambda_0}
+\mathcal R.
\end{split}
\end{equation*}
This handles the contribution of
$\mathcal T^{B,1}_{<\delta\lambda_0,\delta\lambda_1}$ in \eqref{Bquadprincipal}.

The other two lower-endpoint terms are perturbative. For
example, \eqref{L2phierr} and two $L_t^4L_x^{\infty}$ estimates give
\begin{equation*}
\begin{split}
\mathcal T^{B,2}_{<\delta\lambda_0,\delta\lambda_1}=\eta_{\lambda_0}
L^{l,\lambda_0}_{<\delta\lambda_0}B^l_{<\delta\lambda_0}
L^{l,\lambda_1}_{<\delta\lambda_1}\varphi_{\lambda_1}+O_{L_t^1L_x^2}
(\epsilon^4c_{\lambda_0}\lambda_0^{-s_c-\frac12}).
\end{split}
\end{equation*}
For the remaining term, we apply Proposition \ref{paraexp} and
use the frequency decomposition and bilinear pairings from the
proof of \eqref{replacebyphi}, and applying Lemma
\ref{balancedrestrictedoutput}. The subprincipal terms are
treated as above. After integration in $\lambda_0$ and $\lambda_1$, this gives
\begin{equation*}
\int_{\lambda_0\gtrsim_M\mu}
\int_{\delta\lambda_0}^{\lambda_0}
\|\overline u_\mu^{x_0}
(\mathcal T^{B,2}_{<\delta\lambda_0,\delta\lambda_1})_\mu
\|_{L_{t,x}^1}
\frac{d\lambda_1}{\lambda_1}
\frac{d\lambda_0}{\lambda_0}
\lesssim
\epsilon^4c_\mu^2\mu^{-2s_c-\frac12}
(1+\mu|x_0|)^{\frac12}.
\end{equation*}
The term
$\mathcal T^{B,3}_{<\delta\lambda_0,\delta\lambda_1}$ is treated in
the same way. This proves \eqref{Bquadprincipal}.

We next consider the two single-derivative terms in \eqref{TBexp}.
For the first component $\mathcal T^{B,1}$, the linearization formulas
give, modulo high-frequency cutoff tails,
\begin{equation}\label{LPhiNphilinearization}
\begin{split}
\rho\frac{d}{d\rho}
\bigl(L^{l,\lambda_j}_{<\rho}\Phi^{l,j}_{<\rho}\bigr)
={}&-L^{l,\lambda_j}_{<\rho}\bigl(
\eta_\rho L^{l,\lambda_j}_{<\rho}\Phi^{l,j}_{<\rho}
\bigr)
+\eta_\rho
\bigl(L^{l,\lambda_j}_{<\rho}\bigr)^2
\Phi^{l,j}_{<\rho}
\\
&-L^{l,\lambda_j}_{<\rho}\bigl(
\eta_{\lambda_j}L^l_{<\rho}
\Phi^l_{<\rho}(\psi_\rho,\eta_\rho)
\bigr)
\\
&-L^{l,\lambda_j}_{<\rho}\bigl(
\eta_{\lambda_j}\eta_\rho
L^l_{<\rho}B^l_{<\rho}
\bigr).
\end{split}
\end{equation}
After multiplication by the corresponding undifferentiated factor in
$\mathcal T^{B,1}$ and integration in $\rho$, the high-frequency tails
contribute
\begin{equation*}
O_{L_t^1L_x^2}
\left(
\epsilon^4c_{\lambda_0}\lambda_0^{-\frac12-s_c}
\right),
\end{equation*}
uniformly for $\lambda_1\approx_\delta\lambda_0$. Here the
undifferentiated factor is
$L^{l,\lambda_0}_{<\delta\lambda_0}
\Phi^{l,0}_{<\delta\lambda_0}$ when the derivative falls on the
$j=1$ factor, and
$L^{l,\lambda_1}_{<\delta\lambda_1}
\Phi^{l,1}_{<\delta\lambda_1}$ when it falls on the $j=0$ factor. 
\\

The last term in \eqref{LPhiNphilinearization} is cubic and therefore
becomes quartic after multiplication by this undifferentiated factor.
We decompose $B^l_{<\rho}$ dyadically. For its part at frequencies
comparable to $\rho$, Corollary \ref{principalVBexp} allows us to
replace $B$ by $T_{\mathfrak l}\varphi$. The resulting balanced
quartic term is treated as in the expansion of $\mathcal T^{B,4}$.
For the remaining low-frequency part, we use one $L_{t,x}^2$
unbalanced estimate, two $L_t^4L_x^\infty$ estimates, and one
$L_t^\infty L_x^2$ estimate. The $L$ operators acting on single (comparable or higher) frequency input can be replaced by $T_{\mathfrak l}$ using Proposition
\ref{paraexp}. When an outer $L_{<\rho}$ acts on a product, we first
decompose that product dyadically at output frequency $\beta$. For
$\beta\approx_\delta\rho$, Proposition \ref{paraexp} again applies directly.
For $\beta\ll\rho$, we use the fundamental theorem of calculus as in earlier computations in Subsection \ref{subsec:Neta2} to reduce to the case where the outer operator is $L_{<\beta}$ (where Proposition \ref{paraexp} applies again). The difference contains an additional factor
$\eta_\gamma$, with $\beta\lesssim\gamma\lesssim\rho$, and is therefore
quintic in the present expression; the subprincipal terms produced by
Proposition \ref{paraexp} have the same structure. These terms are
estimated directly using four
$L_t^4L_x^\infty$ bounds and one energy bound, after which
$u_\mu^{x_0}$ is estimated separately in energy. The accompanying
high-frequency tails are controlled by Proposition
\ref{Bernsteinelliptic}. Thus the last term contributes to
$N_{\varphi,B}^{[4]}+\mathcal R$.
\\

For the first three terms in \eqref{LPhiNphilinearization}, we first
isolate their cubic principal parts. For each occurrence of
$\Phi_{<\rho}^{l,j}$, we lower its $\rho$-scale coefficient to the
corresponding low-frequency coefficient at scale
$\delta\lambda_j$. For example, we write
\begin{equation*}
\begin{split}
\Phi^{l,j}_{<\rho}-\varphi_{\lambda_j}
={}&
\bigl(\Phi^{l,j}_{<\rho}
-\Phi^{l,j}_{<\delta\lambda_j}\bigr)
+\bigl(\Phi^{l,j}_{<\delta\lambda_j}
-\varphi_{\lambda_j}\bigr).
\end{split}
\end{equation*}
The latter term can be estimated in $L_{t,x}^2$ using Corollary \ref{Phil0varphibound}, which makes its contribution perturbative (by further using two $L_t^4L_x^{\infty}$ bounds and one $L_t^{\infty}L_x^2$ bound for $\overline{u}_{\mu}^{x_0}$). The
first difference is expanded using the fundamental theorem of calculus
and the linearization formula for $B$, introducing either a high-frequency tail or an additional
factor at a frequency comparable to the other inputs. After
multiplication by the undifferentiated factor, these errors are at
least quartic. Their balanced principal parts contribute to
$N_{\varphi,B}^{[4]}$, while the frequency-separated terms,
higher-order terms, and cutoff tails contribute to $\mathcal R$, as in
the expansion for $\mathcal T^{B,4}$. In the third term of
\eqref{LPhiNphilinearization}, we also replace
$\Phi^l_{<\rho}(\psi_\rho,\eta_\rho)$ by $\varphi_\rho$ using Lemma
\ref{Phiest}; after multiplication by the remaining factors, the error
is quartic and contributes to $N_{\varphi,B}^{[4]}+\mathcal R$ by the
same argument. Finally, we replace the inner $L$ operators by their
principal parts using Proposition \ref{paraexp}. For the outer $L$
operators acting on products, we use the dyadic output decomposition described above. The resulting errors are quartic: their balanced principal
parts contribute to $N_{\varphi,B}^{[4]}$, while their
frequency-separated and higher-order quintic parts contribute to $\mathcal R$.
Isolating the principal cubic terms from these replacements, the two
single-derivative contributions of $\mathcal T^{B,1}$ satisfy
\begin{equation}\label{TB1cubic}
\begin{split}
&\int_0^\infty\int_{\delta\lambda_0}^{\lambda_0}
\int_{\delta\lambda_0}^{\lambda_0}
\lambda_2\frac{d}{d\lambda_2}
\mathcal T^{B,1}_{<\lambda_2,\delta\lambda_1}
\frac{d\lambda_2}{\lambda_2}
\frac{d\lambda_1}{\lambda_1}
\frac{d\lambda_0}{\lambda_0}
\\
&\quad+
\int_0^\infty\int_{\delta\lambda_0}^{\lambda_0}
\int_{\delta\lambda_1}^{\lambda_1}
\lambda_3\frac{d}{d\lambda_3}
\mathcal T^{B,1}_{<\delta\lambda_0,\lambda_3}
\frac{d\lambda_3}{\lambda_3}
\frac{d\lambda_1}{\lambda_1}
\frac{d\lambda_0}{\lambda_0}
\\
&=
\int_0^\infty\int_{\delta\lambda_0}^{\lambda_0}
\int_{\delta\lambda_0}^{\lambda_0}
T_{\mathfrak l}\varphi_{\lambda_1}
\Bigl[
-T_{\mathfrak l}
(\eta_{\lambda_2}T_{\mathfrak l}\varphi_{\lambda_0})
+\eta_{\lambda_2}T_{\mathfrak l}^2\varphi_{\lambda_0}
-T_{\mathfrak l}
(\eta_{\lambda_0}T_{\mathfrak l}\varphi_{\lambda_2})
\Bigr]
\frac{d\lambda_2}{\lambda_2}
\frac{d\lambda_1}{\lambda_1}
\frac{d\lambda_0}{\lambda_0}
\\
&+
\int_0^\infty\int_{\delta\lambda_0}^{\lambda_0}
\int_{\delta\lambda_1}^{\lambda_1}
T_{\mathfrak l}\varphi_{\lambda_0}
\Bigl[
-T_{\mathfrak l}
(\eta_{\lambda_3}T_{\mathfrak l}\varphi_{\lambda_1})
+\eta_{\lambda_3}T_{\mathfrak l}^2\varphi_{\lambda_1}
-T_{\mathfrak l}
(\eta_{\lambda_1}T_{\mathfrak l}\varphi_{\lambda_3})
\Bigr]
\frac{d\lambda_3}{\lambda_3}
\frac{d\lambda_1}{\lambda_1}
\frac{d\lambda_0}{\lambda_0}
\\
&+N_{\varphi,B}^{[4]}+\mathcal R.
\end{split}
\end{equation}
We will shortly see that the two integrals on the right-hand side contribute (modulo perturbative errors) part of the cubic balanced term in \eqref{NphiB23def}. 
\\

Next, we isolate the contributions of $\mathcal T^{B,2}$ and
$\mathcal T^{B,3}$. The two single-derivative contributions associated
with $\mathcal T^{B,2}$ are
\begin{equation*}
\begin{split}
&\int_{\delta\lambda_0}^{\lambda_0}
\eta_{\lambda_0}
\left[
\lambda_2\frac{d}{d\lambda_2}
\bigl(
L^{l,\lambda_0}_{<\lambda_2}B^l_{<\lambda_2}
\bigr)
\right]
L^{l,\lambda_1}_{<\delta\lambda_1}
\Phi^{l,1}_{<\delta\lambda_1}
\frac{d\lambda_2}{\lambda_2}
\\
&\quad+
\int_{\delta\lambda_1}^{\lambda_1}
\eta_{\lambda_0}
L^{l,\lambda_0}_{<\delta\lambda_0}B^l_{<\delta\lambda_0}
\left[
\lambda_3\frac{d}{d\lambda_3}
\bigl(
L^{l,\lambda_1}_{<\lambda_3}
\Phi^{l,1}_{<\lambda_3}
\bigr)
\right]
\frac{d\lambda_3}{\lambda_3}.
\end{split}
\end{equation*}
In the first integral, the derivative falls on
$L^{l,\lambda_0}_{<\lambda_2}B^l_{<\lambda_2}$. The first term in
\eqref{LBlinearizationNphi}, after replacing the $\Phi$-factor by
$\varphi$ and extracting the principal parts of the $L$ operators,
produces the cubic term
\begin{equation*}
\eta_{\lambda_0}
T_{\mathfrak l}\varphi_{\lambda_1}
T_{\mathfrak l}^2\varphi_{\lambda_2}.
\end{equation*}
The second term in \eqref{LBlinearizationNphi} contains the additional
factor $\eta_{\lambda_2}$ and is therefore quartic. Together with the
errors produced in the $\Phi$- and $L$-replacements, it contributes to
$N_{\varphi,B}^{[4]}+\mathcal R$.

In the second integral, the undifferentiated exterior block is
\begin{equation*}
\eta_{\lambda_0}
L^{l,\lambda_0}_{<\delta\lambda_0}B^l_{<\delta\lambda_0},
\end{equation*}
while the differentiated factor is given by the right-hand side of
\eqref{LPhiNphilinearization} with $j=1$ and $\rho=\lambda_3$. This
contribution is directly perturbative: we place the undifferentiated
block in $L_{t,x}^2$ using \eqref{LBetabil}, and control the
differentiated factor in $L_t^2L_x^\infty$ using two
$L_t^4L_x^\infty$ bounds, with any remaining factors placed in
$L_{t,x}^\infty$. Similarly, the two single-derivative contributions associated with
$\mathcal T^{B,3}$ are
\begin{equation*}
\begin{split}
&\int_{\delta\lambda_0}^{\lambda_0}
\eta_{\lambda_1}
L^{l,\lambda_1}_{<\delta\lambda_1}B^l_{<\delta\lambda_1}
\left[
\lambda_2\frac{d}{d\lambda_2}
\bigl(
L^{l,\lambda_0}_{<\lambda_2}
\Phi^{l,0}_{<\lambda_2}
\bigr)
\right]
\frac{d\lambda_2}{\lambda_2}
\\
&\quad+
\int_{\delta\lambda_1}^{\lambda_1}
\eta_{\lambda_1}
\left[
\lambda_3\frac{d}{d\lambda_3}
\bigl(
L^{l,\lambda_1}_{<\lambda_3}B^l_{<\lambda_3}
\bigr)
\right]
L^{l,\lambda_0}_{<\delta\lambda_0}
\Phi^{l,0}_{<\delta\lambda_0}
\frac{d\lambda_3}{\lambda_3}.
\end{split}
\end{equation*}
As above, the first integral is directly perturbative, while the principal part
of the second gives
\begin{equation*}
\eta_{\lambda_1}
T_{\mathfrak l}\varphi_{\lambda_0}
T_{\mathfrak l}^2\varphi_{\lambda_3}.
\end{equation*}
The remaining terms in the second integral contribute to
$N_{\varphi,B}^{[4]}+\mathcal R$.

To summarize more precisely, the two directly perturbative contributions satisfy
\begin{equation*}
\begin{split}
\int_{\delta\lambda_1}^{\lambda_1}
\lambda_3\frac{d}{d\lambda_3}
\mathcal T^{B,2}_{<\delta\lambda_0,\lambda_3}
\frac{d\lambda_3}{\lambda_3}
+\int_{\delta\lambda_0}^{\lambda_0}
\lambda_2\frac{d}{d\lambda_2}
\mathcal T^{B,3}_{<\lambda_2,\delta\lambda_1}
\frac{d\lambda_2}{\lambda_2}=
O_{L_t^1L_x^2}
\left(
\epsilon^4c_{\lambda_0}\lambda_0^{-\frac12-s_c}
\right).
\end{split}
\end{equation*}
Combining the above expansions, we see that
$\mathcal T^{B,2}$ and $\mathcal T^{B,3}$ satisfy
\begin{equation}\label{T23cont}
\begin{split}
&\int_0^\infty\int_{\delta\lambda_0}^{\lambda_0}
\int_{\delta\lambda_0}^{\lambda_0}
\lambda_2\frac{d}{d\lambda_2}
\mathcal T^{B,2}_{<\lambda_2,\delta\lambda_1}
\frac{d\lambda_2}{\lambda_2}
\frac{d\lambda_1}{\lambda_1}
\frac{d\lambda_0}{\lambda_0}
\\
&\quad+
\int_0^\infty\int_{\delta\lambda_0}^{\lambda_0}
\int_{\delta\lambda_1}^{\lambda_1}
\lambda_3\frac{d}{d\lambda_3}
\mathcal T^{B,3}_{<\delta\lambda_0,\lambda_3}
\frac{d\lambda_3}{\lambda_3}
\frac{d\lambda_1}{\lambda_1}
\frac{d\lambda_0}{\lambda_0}
\\
&\quad=
\int_0^\infty\int_{\delta\lambda_0}^{\lambda_0}
\int_{\delta\lambda_0}^{\lambda_0}
\eta_{\lambda_0}
T_{\mathfrak l}\varphi_{\lambda_1}
T_{\mathfrak l}^2\varphi_{\lambda_2}
\frac{d\lambda_2}{\lambda_2}
\frac{d\lambda_1}{\lambda_1}
\frac{d\lambda_0}{\lambda_0}
\\
&\quad+
\int_0^\infty\int_{\delta\lambda_0}^{\lambda_0}
\int_{\delta\lambda_1}^{\lambda_1}
\eta_{\lambda_1}
T_{\mathfrak l}\varphi_{\lambda_0}
T_{\mathfrak l}^2\varphi_{\lambda_3}
\frac{d\lambda_3}{\lambda_3}
\frac{d\lambda_1}{\lambda_1}
\frac{d\lambda_0}{\lambda_0}
\\
&\quad+N_{\varphi,B}^{[4]}+\mathcal R.
\end{split}
\end{equation}
Finally, we isolate the quadratic term $N_{\varphi,B}^{[2]}$ from above, and will also show that the cubic principal terms in \eqref{T23cont} and \eqref{TB1cubic} add up, modulo perturbative errors, to $N_{\varphi,B}^{[3]}$ as defined in \eqref{NphiB23def}. For this, we have to commute the internal projections in the integrals above through the principal operators. Thus,
\begin{equation*}
\begin{split}
T_{\mathfrak l}\varphi_\lambda=S_\lambda(T_{\mathfrak l}\varphi)
+[T_{\mathfrak l},S_\lambda]\varphi,\hspace{5mm}T_{\mathfrak l}^2\varphi_\lambda
=S_\lambda(T_{\mathfrak l}^2\varphi)
+[T_{\mathfrak l}^2,S_\lambda]\varphi,
\end{split}
\end{equation*}
where
\begin{equation*}
[T_{\mathfrak l}^2,S_\lambda]
=T_{\mathfrak l}[T_{\mathfrak l},S_\lambda]
+[T_{\mathfrak l},S_\lambda]T_{\mathfrak l}.
\end{equation*}
Each commutator contains a differentiated paracoefficient and, for the
cubic terms, contributes directly to $\mathcal R$ in $L_t^1L_x^2$,
exactly as in the extraction of the cubic term $N_{\eta}^{[3]}$ in
Subsection \ref{subsec:Neta34}. The commutators in the quadratic
principal term coming from the lower endpoint in
\eqref{Bquadprincipal} also contribute to $\mathcal R$, but only with
the usual weaker paired bound
\begin{equation*}
\|\mathcal{R}_{\mu}\overline{u}_{\mu}^{x_0}\|_{L_{t,x}^1}\lesssim \epsilon^4 c_{\mu}^2\mu^{-2s_c-\frac{1}{2}}(1+\mu|x_0|)^{\frac12}.
\end{equation*}
Here we use the frequency decomposition and bilinear pairings
from the proof of \eqref{replacebyphi}, and applying Lemma
\ref{balancedrestrictedoutput}.
\\

The remaining part of the lower-endpoint term in \eqref{Bquadprincipal} now gives
\begin{equation*}
\frac12\Pi(T_{\mathfrak l}\varphi,T_{\mathfrak l}\varphi)
=N_{\varphi,B}^{[2]}.
\end{equation*}
It remains to identify the cubic principal terms. After splitting the
square integration regions in \eqref{TB1cubic} and
\eqref{T23cont} into their two ordered triangles and relabeling the
scale parameters, we enlarge the resulting ordered regions to
$\mathcal D$. The added portions have their smallest frequency
separated from their largest frequency and therefore contribute to
$\mathcal R$ by the same bilinear pairings used above. On the full
domain $\mathcal D$, the positive terms in these two identities give
the six ordered placements in
\begin{equation*}
\Pi(T_{\mathfrak l}\varphi,\eta,
T_{\mathfrak l}^2\varphi).
\end{equation*}
Similarly, modulo $\mathcal R$, the negative nested terms in
\eqref{TB1cubic} give the six ordered placements in
\begin{equation*}
-\Pi\bigl(
T_{\mathfrak l}\varphi,
T_{\mathfrak l}(\eta,T_{\mathfrak l}\varphi)
\bigr).
\end{equation*}
Their sum is $N_{\varphi,B}^{[3]}+\mathcal R$. Combining this identity with the
quadratic, quartic, and remainder estimates above proves Proposition
\ref{B-TBexp}.
\end{proof}

\subsubsection{Expansion of
$\frac{1}{2}|V|^2-T_{V^l}\cdot V$.}
We next consider the expansion of
$\frac{1}{2}|V|^2-T_{V^l}\cdot V$. Define
\begin{equation}\label{NphiV23def}
\begin{split}
N_{\varphi,V}^{[2]}\coloneq
\frac12\Pi(T_{\mathfrak m}\varphi,T_{\mathfrak m}\varphi),\hspace{5mm}N_{\varphi,V}^{[3]}\coloneq
\Pi(T_{\mathfrak m}\varphi,\eta,
T_{\mathfrak m}T_{\mathfrak l}\varphi)
-\Pi\bigl(
T_{\mathfrak m}\varphi,
T_{\mathfrak m}(\eta,T_{\mathfrak l}\varphi)
\bigr).
\end{split}
\end{equation}
\begin{proposition}\label{V-TVexp}
Set $\mathcal V\coloneq\frac12|V|^2-T_{V^l}\cdot V.$ We have
\begin{equation*}
\mathcal V=N_{\varphi,V}^{[2]}+N_{\varphi,V}^{[3]}
+N_{\varphi,V}^{[4]}+\mathcal R,
\end{equation*}
where $N_{\varphi,V}^{[4]}$ is of the form described in
Definition \ref{Nphi4def}. The following bounds hold:
\begin{enumerate}
\item For every $\mu>0$ and $x_0\in\mathbb R^2$, the remainder
satisfies
\begin{equation*}
\|\mathcal R_\mu\overline u_\mu^{x_0}\|_{L_{t,x}^1}
\lesssim
\epsilon^4c_\mu^2\mu^{-2s_c-\frac12}
(1+\mu|x_0|)^{\frac12}.
\end{equation*}

\item 
For every $\lambda>0$, we have
\begin{equation}\label{VauxL2Linf}
\|\mathcal V_\lambda\|_{L_t^2L_x^\infty}
\lesssim_\delta\epsilon^2\lambda^{-\frac54}.
\end{equation}

\item 
For every $0<\lambda\ll\mu$ and $x_0\in\mathbb R^2$, we have
a decomposition
\begin{equation}\label{Vauxdecomp}
\mathcal V_\lambda\overline u_\mu^{x_0}u_\mu
=\mathcal R_V^{\mathrm{aux}}+\nabla\cdot Q_V^{\mathrm{aux}},
\end{equation}
where, for $\frac43<p\leq2$,
\begin{equation}\label{Vauxtrilinear}
\|\mathcal R_V^{\mathrm{aux}}\|_{L_{t,x}^1}
+\lambda^{\frac2p-1}
\|Q_V^{\mathrm{aux}}\|_{L_t^1L_x^p}
\lesssim_{\delta,p}
\epsilon^4c_\mu^2\lambda^{-\frac32}
\mu^{-2s_c}(1+\mu|x_0|)^{\frac12}.
\end{equation}

\item For every $\lambda>0$, we have
\begin{equation}\label{VauxL2cubic}
\|(\mathcal V-N_{\varphi,V}^{[2]})_\lambda\|_{L_{t,x}^2}
\lesssim_\delta\epsilon^3c_\lambda\lambda^{-\frac14-s_c}.
\end{equation}
\end{enumerate}
\end{proposition}

\begin{proof}
The proof is the same as that of Proposition \ref{B-TBexp},
except some of the linearization formulas must use $M$ in
place of $L$ and $V$ in place of $B$, respectively. Otherwise,
the proof is virtually identical. We omit the details, which
are straightforward modifications.
\end{proof}
\subsubsection{Expansion of the remaining cubic terms.}
Finally, we expand the remaining cubic term
\begin{equation*}
F\coloneq-\left[
BV\cdot\nabla\eta-T_{B^l}(V\cdot\nabla\eta)
-T_{V^l}\cdot(B\nabla\eta)
+T_{V^l}\cdot T_{B^l}\nabla\eta
\right]
\end{equation*}
in the definition of $N_{\varphi}$. As in the preceding expansions, we will decompose $F$ into
quadratic, cubic, and quartic balanced terms, together with a
perturbative remainder. We use the same convention as in
\eqref{NphinestedPiconvention}, allowing also for
$\mathfrak a=\mathfrak g$. Set
\begin{equation}\label{NphiF2def}
\begin{split}
N_{\varphi,F}^{[2]}
\coloneq -T_{\nabla\eta}\cdot
\Pi\bigl(T_{\mathfrak l}\varphi,T_{\mathfrak m}\varphi\bigr),
\end{split}
\end{equation}
and
\begin{equation}\label{NphiF3def}
\begin{split}
N_{\varphi,F}^{[3]}
\coloneq -
\Pi\bigl(
T_{\mathfrak l}\varphi,
T_{\mathfrak g}(\eta,T_{\mathfrak l}\varphi)
\bigr)
+\Pi\bigl(
T_{\mathfrak l}\varphi,
\eta,
T_{\mathfrak g}T_{\mathfrak l}\varphi
\bigr)
-N_{\varphi,B}^{[3]}+N_{\varphi,V}^{[3]}.
\end{split}
\end{equation}
\begin{proposition}\label{cubicvarphi}
There exists a term $N_{\varphi,F}^{[4]}$ of the form described in
Definition \ref{Nphi4def} such that
\begin{equation}\label{Ffinalexpansion}
F=N_{\varphi,F}^{[2]}+N_{\varphi,F}^{[3]}
+N_{\varphi,F}^{[4]}+\mathcal R.
\end{equation}
where
\begin{equation}\label{Fcubicremainder}
\|\mathcal R_\mu\overline u_\mu^{x_0}\|_{L_{t,x}^1}
\lesssim
\epsilon^4c_\mu^2\mu^{-2s_c-\frac12}
(1+\mu|x_0|)^{\frac12}.
\end{equation}
\end{proposition}

\begin{proof}
As usual, $\mathcal R$ and $N_{\varphi,F}^{[4]}$ denote generic
expressions with the above definition that may change from line
to line. We begin by expanding the first two terms in the
definition of $F$. We have
\begin{equation}\label{cubictermexp1}
\begin{split}
BV\cdot\nabla\eta-T_{B^l}(V\cdot\nabla\eta)
={}&T_{V\cdot\nabla\eta}B+\int_0^\infty
P_{<\delta^{-1}\lambda_0}
(B-B_{<\delta\lambda_0})
(V\cdot\nabla\eta)_{\lambda_0}
\frac{d\lambda_0}{\lambda_0}
\\
&+\int_0^\infty
\bigl(P_{<\delta^{-1}\lambda_0}
-P_{<M\delta\lambda_0}\bigr)
B_{<\delta\lambda_0}
(V\cdot\nabla\eta)_{\lambda_0}
\frac{d\lambda_0}{\lambda_0}.
\end{split}
\end{equation}
Similarly,
\begin{equation}\label{cubictermexp2}
\begin{split}
T_{V^l}\cdot
\bigl(B\nabla\eta-T_{B^l}\nabla\eta\bigr)
={}&T_{V^l}\cdot T_{\nabla\eta}B+\int_0^\infty
T_{V^l}\cdot\left[
P_{<\delta^{-1}\lambda_0}
(B-B_{<\delta\lambda_0})
\nabla\eta_{\lambda_0}
\right]
\frac{d\lambda_0}{\lambda_0}
\\
&+\int_0^\infty
T_{V^l}\cdot\left[
\bigl(P_{<\delta^{-1}\lambda_0}
-P_{<M\delta\lambda_0}\bigr)
B_{<\delta\lambda_0}
\nabla\eta_{\lambda_0}
\right]
\frac{d\lambda_0}{\lambda_0}.
\end{split}
\end{equation}
Next, we analyze the difference between these two terms. We begin by
analyzing the difference between the first integrated term in
\eqref{cubictermexp1} and the corresponding term in
\eqref{cubictermexp2}. Here, there is a useful commutator structure
between the term
$P_{<\delta^{-1}\lambda_0}(B-B_{<\delta\lambda_0})$
and the $V$ paracoefficient above that will let us simplify this
component. To isolate this, we temporarily write
\begin{equation*}
b_{\lambda_0}
\coloneq
P_{<\delta^{-1}\lambda_0}
(B-B_{<\delta\lambda_0}).
\end{equation*}
Then we have
\begin{equation*}
\begin{split}
T_{V^l}\cdot
\bigl(b_{\lambda_0}\nabla\eta_{\lambda_0}\bigr)
-b_{\lambda_0}(V\cdot\nabla\eta)_{\lambda_0}
={}&
b_{\lambda_0}
\left[
T_{V^l}\cdot\nabla\eta_{\lambda_0}
-(V\cdot\nabla\eta)_{\lambda_0}
\right]
+\mathcal C_{\lambda_0},
\end{split}
\end{equation*}
where
\begin{equation*}
\mathcal C_{\lambda_0}
\coloneq
T_{V^l}\cdot
\bigl(b_{\lambda_0}\nabla\eta_{\lambda_0}\bigr)
-b_{\lambda_0}T_{V^l}\cdot\nabla\eta_{\lambda_0}.
\end{equation*}
From the definition of $T_{V^l}$, we can write
\begin{equation*}
\begin{split}
\mathcal C_{\lambda_0}
=
\int_0^\infty
V^l_{<\delta\lambda_1}\cdot
[S_{\lambda_1},b_{\lambda_0}]
\nabla\eta_{\lambda_0}
\frac{d\lambda_1}{\lambda_1}=
\int_0^\infty
(V^l_{<\delta\lambda_1}-V^l_{<\delta\lambda_0})\cdot
[S_{\lambda_1},b_{\lambda_0}]
\nabla\eta_{\lambda_0}
\frac{d\lambda_1}{\lambda_1}.
\end{split}
\end{equation*}
Thus the commutator contains both a difference of the $V$
paracoefficients and a commutator with the $B$-factor. Applying the
fundamental theorem of calculus to
$V^l_{<\delta\lambda_1}-V^l_{<\delta\lambda_0}$ and using the
linearization formula for $V$ expresses this difference as an
integral of linearized $V$-factors. In particular, no unrestricted
low-frequency $V$ paracoefficient remains. On the other hand,
Proposition \ref{Leibnizcommutator}, applied to
$[S_{\lambda_1},b_{\lambda_0}]$, produces a differentiated $B$
factor. From this, one can see that $\mathcal C_{\lambda_0}$ is
included in the perturbative $\mathcal R$ term. Indeed, after
projecting to output frequency $\mu$ and pairing with
$\overline u_\mu^{x_0}$, we use the principal/subprincipal
decomposition and bilinear pairings from the proof of
\eqref{replacebyphi} and Lemma \ref{balancedrestrictedoutput}.

Next, we analyze the difference between the second line in
\eqref{cubictermexp1} and the second line of
\eqref{cubictermexp2}. Set
\begin{equation*}
e_{\lambda_0}
\coloneq
\bigl(P_{<\delta^{-1}\lambda_0}
-P_{<M\delta\lambda_0}\bigr)
B_{<\delta\lambda_0}.
\end{equation*}
The integrand can be written as
\begin{equation*}
\begin{split}
&T_{V^l}\cdot
\bigl(e_{\lambda_0}\nabla\eta_{\lambda_0}\bigr)
-e_{\lambda_0}(V\cdot\nabla\eta)_{\lambda_0}
\\
&\quad=
e_{\lambda_0}
\left[
T_{V^l}\cdot\nabla\eta_{\lambda_0}
-(V\cdot\nabla\eta)_{\lambda_0}
\right]
+[T_{V^l},e_{\lambda_0}]
\cdot\nabla\eta_{\lambda_0}
\\
&\quad=
e_{\lambda_0}
\left[
(T_{V^l}-V)\cdot\nabla\eta_{\lambda_0}
\right]
-e_{\lambda_0}[S_{\lambda_0},V]\cdot\nabla\eta
+[T_{V^l},e_{\lambda_0}]
\cdot\nabla\eta_{\lambda_0}.
\end{split}
\end{equation*}
Since $e_{\lambda_0}$ is a high-frequency tail error and, arguing as
above for $[T_{V^l},e_{\lambda_0}]$, the $V$-factor is always either
``differentiated'' by a commutator, appears as a high-frequency
paracoefficient, or appears as an approximately
$\lambda_0$-localized linearized factor, this contribution is
included in $\mathcal R$ by Propositions
\ref{Bernsteinelliptic} and \ref{Leibnizcommutator}. Combining the
above cancellations, we obtain
\begin{equation}\label{Fprocessing}
F=F^1+F^2+\mathcal R,
\end{equation}
where we define
\begin{equation*}
\begin{split}
F^1\coloneq
T_{V^l}\cdot T_{\nabla\eta}B-T_{V\cdot\nabla\eta}B,\hspace{5mm}F^2\coloneq
\int_0^\infty
\left[
T_{V^l}\cdot\nabla\eta_{\lambda_0}
-(V\cdot\nabla\eta)_{\lambda_0}
\right]
b_{\lambda_0}
\frac{d\lambda_0}{\lambda_0}.
\end{split}
\end{equation*}
Next, we expand $F^1$. To keep our frequency indexing more compact,
we set for the remainder of this expansion,
\begin{equation*}
V^{\mathrm{lo}}\coloneq P_{<M\delta\lambda_0}V,
\qquad
V^{\mathrm{hi}}\coloneq V-V^{\mathrm{lo}}.
\end{equation*}
We use the same superscripts for $\eta$ and for other functions,
always relative to the current scale $\lambda_0$. Then
\begin{equation}\label{F1cubicexp}
\begin{split}
F^1&=\int_0^\infty
\left[
T_{V^l}\cdot
(B_{\lambda_0}\nabla\eta_{<\delta\lambda_0})
-B_{\lambda_0}
P_{<\delta\lambda_0}
(V^{\mathrm{lo}}\cdot\nabla\eta)
\right]
\frac{d\lambda_0}{\lambda_0}
\\
&\quad-\int_0^\infty
B_{\lambda_0}P_{<\delta\lambda_0}
(V^{\mathrm{hi}}\cdot\nabla\eta)
\frac{d\lambda_0}{\lambda_0}
\\
&\coloneq F^{1,1}+F^{1,2}.
\end{split}
\end{equation}
Above, we write $B_{\lambda_0}\coloneq S_{\lambda_0}B$.
\\

\paragraph{\textbf{The term $F^{1,1}$}:}
We begin with (by far) the most complicated term. Using
\eqref{bilformrep} (after inserting sufficiently fattened
projections), we obtain
\begin{equation}\label{F11firstsplit}
\begin{split}
F^{1,1}={}&\int_0^\infty
\bigl(T_{V^l}-P_{<\delta\lambda_0}V\bigr)\cdot
\bigl(B_{\lambda_0}\nabla\eta_{<\delta\lambda_0}\bigr)
\frac{d\lambda_0}{\lambda_0}+\int_0^\infty
\lambda_0^{-2}B_{\lambda_0}
\mathcal L\left(
\nabla V^{\mathrm{lo}},
\nabla^2\widetilde P_{<M\delta\lambda_0}\eta
\right)
\frac{d\lambda_0}{\lambda_0}.
\end{split}
\end{equation}
The second integral contributes to $\mathcal R$. Indeed, the
output support forces $\mu\approx\lambda_0$, and we use two
unbalanced pairings, between $\overline u_\mu^{x_0}$ and
$\nabla V^{\mathrm{lo}}$, and between $B_{\lambda_0}$ and
$\nabla^2\widetilde P_{<M\delta\lambda_0}\eta$. The resulting
contribution satisfies \eqref{Fcubicremainder}.
\\

For the first integral, we write
\begin{equation*}
\begin{split}
\bigl(T_{V^l}-P_{<\delta\lambda_0}V\bigr)\cdot
\bigl(B_{\lambda_0}\nabla\eta_{<\delta\lambda_0}\bigr)&=
\bigl(T_{V^l}-P_{<M\delta\lambda_0}V\bigr)\cdot
\bigl(B_{\lambda_0}\nabla\eta_{<\delta\lambda_0}\bigr)
+\mathcal A_{\lambda_0},
\end{split}
\end{equation*}
where we write as shorthand
\begin{equation*}
\mathcal A_{\lambda_0}\coloneq
B_{\lambda_0}\nabla\eta_{<\delta\lambda_0}\cdot
\bigl(P_{<M\delta\lambda_0}-P_{<\delta\lambda_0}\bigr)V.  
\end{equation*}
The principal part of $\mathcal{A}_{\lambda_0}$ (obtained by replacing $B$ and $V$ with $T_{\mathfrak{l}}\varphi$ and $T_{\mathfrak{m}}\varphi$, respectively) will cancel another term appearing a bit later in our analysis.
For this purpose, we define
\begin{equation}\label{Fannularsplit}
\begin{split}
\mathcal A^{\mathrm{prin}}
&\coloneq
\int_0^\infty
T_{\mathfrak l}\varphi_{\lambda_0}
\nabla\eta_{<\delta\lambda_0}\cdot
\bigl(P_{<M\delta\lambda_0}-P_{<\delta\lambda_0}\bigr)
T_{\mathfrak m}\varphi
\frac{d\lambda_0}{\lambda_0},
\end{split}
\end{equation}
and write
\begin{equation*}
\mathcal{A}^{\mathrm{sub}}\coloneq \int_{0}^{\infty}\mathcal{A}_{\lambda_0}\frac{d\lambda_0}{\lambda_0}-\mathcal{A}^{\mathrm{prin}}    
\end{equation*}
for the remainder. The principal term will cancel modulo $\mathcal R$ with similar
principal contributions below. The term $\mathcal A^{\mathrm{sub}}$
will be treated later on; it will contribute a term of the form $N_{\varphi}^{[4]}$ plus a perturbative remainder.

Using the definition of $T_{V^l}$, we can then write for the remaining integrand

\begin{equation}\label{TVminusVlo}
\begin{split}
\bigl(T_{V^l}-P_{<M\delta\lambda_0}V\bigr)\cdot
\bigl(B_{\lambda_0}\nabla\eta_{<\delta\lambda_0}\bigr)&=
\int_0^\infty
P_{<M\delta\lambda_1}
(V_{<\delta\lambda_1}-V)
\cdot
S_{\lambda_1}\left(
B_{\lambda_0}\nabla\eta_{<\delta\lambda_0}
\right)
\frac{d\lambda_1}{\lambda_1}
\\
&+
\int_0^\infty
\bigl(P_{<M\delta\lambda_1}
-P_{<M\delta\lambda_0}\bigr)V
\cdot
S_{\lambda_1}\left(
B_{\lambda_0}\nabla\eta_{<\delta\lambda_0}
\right)
\frac{d\lambda_1}{\lambda_1}.
\end{split}
\end{equation}
Both integrands are supported where
$\lambda_1\approx\lambda_0$. We first consider the second integral.
We begin by writing
\begin{equation}\label{Vprincipalfactor}
\begin{split}
\bigl(P_{<M\delta\lambda_1}
-P_{<M\delta\lambda_0}\bigr)V
={}&
\bigl(P_{<M\delta\lambda_1}
-P_{<M\delta\lambda_0}\bigr)
T_{\mathfrak m}\varphi
+\mathcal E^V_{\lambda_0,\lambda_1},
\end{split}
\end{equation}
where
\begin{equation*}
\begin{split}
V^{\mathrm{sub}}\coloneq V-T_{\mathfrak m}\varphi,\hspace{5mm}\mathcal E^V_{\lambda_0,\lambda_1}\coloneq
\bigl(P_{<M\delta\lambda_1}
-P_{<M\delta\lambda_0}\bigr)V^{\mathrm{sub}}.
\end{split}
\end{equation*}
To estimate the contribution of
$\mathcal E^V_{\lambda_0,\lambda_1}$, we commute
$\nabla\eta_{<\delta\lambda_0}$ through $S_{\lambda_1}$.
Since $\mathcal E^V_{\lambda_0,\lambda_1}$ is antisymmetric
in $\lambda_0,\lambda_1$, while $S_{\lambda_1}S_{\lambda_0}$
is symmetric, the term without the commutator contains
(after integrating in $\lambda_0$ and $\lambda_1$)
\begin{equation*}
\frac12\mathcal E^V_{\lambda_0,\lambda_1}\cdot
\bigl(\nabla\eta_{<\delta\lambda_0}
-\nabla\eta_{<\delta\lambda_1}\bigr)
S_{\lambda_1}S_{\lambda_0}B.
\end{equation*}
The remaining term contains
$\mathcal E^V_{\lambda_0,\lambda_1}\cdot
[S_{\lambda_1},\nabla\eta_{<\delta\lambda_0}]B_{\lambda_0}$.

The difference of the $\eta$ factors is localized at frequency
$\approx\delta\lambda_0$, while the commutator contributes
$\lambda_0^{-1}\nabla^2\eta_{<\delta\lambda_0}$.
We pair these $\eta$ factors with $\overline u_\mu^{x_0}$,
estimate $B_{\lambda_0}$ in $L_t^\infty L_x^2$, and use
Corollary \ref{principalVBexp} to estimate
$\mathcal E^V_{\lambda_0,\lambda_1}$ in $L_t^2L_x^\infty$.
Since $\lambda_1\approx\lambda_0\approx\mu$ and the latter
factor has frequency $\approx M\delta\lambda_0$, this gives
\begin{equation*}
\begin{split}
&\left\|\overline u_\mu^{x_0}P_\mu
\int_0^\infty\int_0^\infty
\mathcal E^V_{\lambda_0,\lambda_1}\cdot
S_{\lambda_1}\left(
B_{\lambda_0}\nabla\eta_{<\delta\lambda_0}\right)
\frac{d\lambda_1}{\lambda_1}
\frac{d\lambda_0}{\lambda_0}\right\|_{L_{t,x}^1}
\\
&\qquad\lesssim_{M,\delta}
\epsilon^5c_\mu^2\mu^{-\frac12-2s_c}(1+\mu|x_0|)^{\frac12}.
\end{split}
\end{equation*}
Thus this contribution is perturbative.

We can also replace $B_{\lambda_0}$ by its principal part
$T_{\mathfrak l}\varphi_{\lambda_0}$. Indeed, using the notation of
Proposition \ref{VBsecondord}, we have
\begin{equation*}
B_{\lambda_0}-T_{\mathfrak l}\varphi_{\lambda_0}
=
S_{\lambda_0}B^{\mathrm{sub}}
+[S_{\lambda_0},T_{\mathfrak l}]\varphi.
\end{equation*}
We estimate this error after pairing with
$\overline u_\mu^{x_0}$. Decomposing the principal $V$-factor
$T_{\mathfrak m}\varphi$ in \eqref{Vprincipalfactor} into a sum of
components localized at frequency $\beta\approx M\delta\lambda_0$,
we use the kernel bounds for $T_{\mathfrak l}$, apply Proposition
\ref{VBsecondord}, and estimate
$\nabla\eta_{<\delta\lambda_0}$ in $L_{t,x}^\infty$ to obtain
\begin{equation}\label{Blambda0principalcomparison}
\begin{split}
&\left\|\overline u_\mu^{x_0}P_\mu\left[
\bigl(P_{<M\delta\lambda_1}-P_{<M\delta\lambda_0}\bigr)
T_{\mathfrak m}\varphi\cdot
S_{\lambda_1}\left(
S_{\lambda_0}B^{\mathrm{sub}}
\nabla\eta_{<\delta\lambda_0}\right)\right]\right\|_{L_{t,x}^1}
\\
&\qquad\lesssim
\epsilon^5c_\mu^2\mu^{\frac14-2s_c}
\sum_{\beta\approx M\delta\lambda_0}\beta^{-\frac34}
(1+\mu|x_0|)^{\frac12}
\\
&\qquad\lesssim_{M,\delta}
\epsilon^5c_\mu^2\mu^{-\frac12-2s_c}(1+\mu|x_0|)^{\frac12}.
\end{split}
\end{equation}
Using Proposition \ref{Leibnizcommutator} and two unbalanced
bilinear estimates, the commutator term contributes the same bound. The remaining principal contribution of the second integral in
\eqref{TVminusVlo}, after integration in $\lambda_0$, is therefore
\begin{equation}\label{secondintppart}
\begin{split}
\int_0^\infty\int_0^\infty
\bigl(P_{<M\delta\lambda_1}
-P_{<M\delta\lambda_0}\bigr)
T_{\mathfrak m}\varphi\cdot
S_{\lambda_1}\left(
T_{\mathfrak l}\varphi_{\lambda_0}
\nabla\eta_{<\delta\lambda_0}
\right)
\frac{d\lambda_1}{\lambda_1}
\frac{d\lambda_0}{\lambda_0}.
\end{split}
\end{equation}
Next, we expand the first integral in \eqref{TVminusVlo}. The
linearization formula for $V$ gives
\begin{equation}\label{VtailFTCstart}
\begin{split}
P_{<M\delta\lambda_1}(V_{<\delta\lambda_1}-V)
={}&-P_{<M\delta\lambda_1}
\int_{\delta\lambda_1}^\infty
\mathcal T^V_{<\lambda_2}
\frac{d\lambda_2}{\lambda_2}
+\mathcal E^{V,\mathrm{tail}}_{\lambda_1},
\end{split}
\end{equation}
where, for each fixed $\lambda_2$, we suppress this frequency
projection scale from the notation and, for an independent
``geometric'' truncation scale $\rho>0$, set
\begin{equation}\label{TVfixeddef}
\mathcal T^V_{<\rho}
\coloneq
M^{l,\lambda_2}_{<\rho}\Phi^{l,2}_{<\rho}
+\eta_{\lambda_2}M^{l,\lambda_2}_{<\rho}B^l_{<\rho}.
\end{equation}
The term $\mathcal E^{V,\mathrm{tail}}_{\lambda_1}$ contains the
high-frequency tails produced by the fixed projections in passing
from the full linearization formula for $V$ to the fixed-projection
form \eqref{TVfixeddef}. Since
$\lambda_1\approx\lambda_0$, Proposition
\ref{Bernsteinelliptic}, together with the usual
$L_t^4L_x^\infty$ and energy bounds, gives
\begin{equation*}
\left\|
\int_0^\infty
\mathcal E^{V,\mathrm{tail}}_{\lambda_1}\cdot
S_{\lambda_1}\left(
B_{\lambda_0}\nabla\eta_{<\delta\lambda_0}
\right)
\frac{d\lambda_1}{\lambda_1}
\right\|_{L_t^1L_x^2}
\lesssim
\epsilon^4c_{\lambda_0}
\lambda_0^{-\frac12-s_c}.
\end{equation*}
After integration in $\lambda_0$ and estimating
$\overline u_\mu^{x_0}$ separately in $L_t^\infty L_x^2$, this
contribution satisfies \eqref{Fcubicremainder}.

To treat the main term in \eqref{VtailFTCstart}, we expand
$\mathcal T^V_{<\lambda_2}$ once more in the geometric truncation
parameter. The fundamental theorem of calculus gives
\begin{equation}\label{TVlinearize}
\mathcal T^V_{<\lambda_2}
=
\mathcal T^V_{<\delta\lambda_2}
+\int_{\delta\lambda_2}^{\lambda_2}
\lambda_3\frac{d}{d\lambda_3}
\mathcal T^V_{<\lambda_3}
\frac{d\lambda_3}{\lambda_3}.
\end{equation}
Here the suppressed output projection scale $\lambda_2$ remains
fixed when differentiating with respect to $\lambda_3$.

We first consider the contribution of the lower endpoint in
\eqref{TVlinearize}. Modulo high-frequency tails, controlled as above
using Proposition \ref{Bernsteinelliptic}, the exterior projection in
\eqref{VtailFTCstart} restricts it to
$\lambda_2\approx_M\delta\lambda_1$. We decompose the first term in
$\mathcal T^V_{<\delta\lambda_2}$ as
\begin{equation}\label{survivingprincipal}
M^{l,\lambda_2}_{<\delta\lambda_2}
\Phi^{l,2}_{<\delta\lambda_2}
=
T_{\mathfrak m}\varphi_{\lambda_2}
+\mathcal E^M_{\lambda_2},
\end{equation}
where
\begin{equation*}
\mathcal E^M_{\lambda_2}
=M^{l,\lambda_2}_{<\delta\lambda_2}
\bigl(
\Phi^{l,2}_{<\delta\lambda_2}
-\varphi_{\lambda_2}
\bigr)+
\bigl(
M^{l,\lambda_2}_{<\delta\lambda_2}
-T_{\mathfrak m}
\bigr)\varphi_{\lambda_2}.
\end{equation*}
For the first term, the representation \eqref{Phidifferror} gives
\begin{equation*}
\Phi^{l,2}_{<\delta\lambda_2}
-\varphi_{\lambda_2}
=
\int_{\beta\approx\lambda_2}
\mathcal L\left(
B^l_{<\delta\beta}-B^l_{<\delta\lambda_2},
\widetilde\eta_\beta
\right)
\frac{d\beta}{\beta}+\lambda_2^{-1}
\int_{\beta\approx\lambda_2}
\mathcal L\left(
\nabla B^l_{<\delta\beta},
\widetilde\eta_\beta
\right)
\frac{d\beta}{\beta}.
\end{equation*}
Arguing as in the estimates following \eqref{Phidifferror}, we
estimate the contribution of these two terms using the unbalanced
pairings
\begin{equation*}
(\overline u_\mu^{x_0},\widetilde\eta_\beta),
\qquad
\left(
B_{\lambda_0},
B^l_{<\delta\beta}-B^l_{<\delta\lambda_2}
\right),
\end{equation*}
and
\begin{equation*}
(\overline u_\mu^{x_0},\widetilde\eta_\beta),
\qquad
\left(
B_{\lambda_0},
\lambda_2^{-1}\nabla B^l_{<\delta\beta}
\right),
\end{equation*}
respectively. The first $B$-coefficient difference contains only a
bounded range of scales comparable to $\delta\lambda_2$, while a
dyadic component at frequency $\gamma\lesssim\delta\lambda_2$ in the
second $B$-coefficient carries the summable factor
$\gamma/\lambda_2$. Thus these terms give a perturbative contribution
in the region $\lambda_2\approx_M\delta\lambda_1$.

For the second term in $\mathcal E^M_{\lambda_2}$, Proposition
\ref{paraexp} gives a zeroth-order term containing
$\varphi_{\lambda_2}$ and an additional factor
$\nabla^2\eta_\gamma$, with
$\gamma\lesssim\delta\lambda_2$, together with higher-order geometric
errors and high-frequency tails. In the zeroth-order term, we
extract the factor $\gamma^2$, use the two pairings
\begin{equation*}
(\overline u_\mu^{x_0},\varphi_{\lambda_2}),
\qquad
(B_{\lambda_0},\eta_\gamma),
\end{equation*}
and then sum over $\gamma$. The higher-order geometric errors are
estimated using \eqref{paraexperror} and the same pairings, while
the high-frequency tails are controlled by Proposition
\ref{Bernsteinelliptic}. All remaining low-frequency factors,
including $\nabla\eta_{<\delta\lambda_0}$, are estimated in
$L_{t,x}^\infty$. Since
\begin{equation*}
\mu\approx\lambda_0,
\qquad
\beta\approx\lambda_2\approx_M\delta\lambda_0,
\end{equation*}
all the displayed pairings are frequency separated. Combining these
estimates, we conclude that the contribution of
$\mathcal E^M_{\lambda_2}$ satisfies \eqref{Fcubicremainder}.

The other term in $\mathcal T^V_{<\delta\lambda_2}$ is
\begin{equation*}
\eta_{\lambda_2}
M^{l,\lambda_2}_{<\delta\lambda_2}
B^l_{<\delta\lambda_2}.
\end{equation*}
This is directly perturbative by applying the two unbalanced
estimates to
\begin{equation*}
(\overline u_\mu^{x_0},\eta_{\lambda_2}),
\qquad
\left(
B_{\lambda_0},
M^{l,\lambda_2}_{<\delta\lambda_2}
B^l_{<\delta\lambda_2}
\right),
\end{equation*}
and estimating $\nabla\eta_{<\delta\lambda_0}$ in
$L_{t,x}^\infty$. Finally, we replace $B_{\lambda_0}$ by
$T_{\mathfrak l}\varphi_{\lambda_0}$ using
Proposition \ref{VBsecondord}, as in
\eqref{Blambda0principalcomparison}, now with the low-frequency
factor
$P_{<M\delta\lambda_1}T_{\mathfrak m}\varphi_{\lambda_2}$. Thus, modulo $\mathcal R$ satisfying \eqref{Fcubicremainder}, the contribution of the lower endpoint
in \eqref{TVlinearize} to the first integral in
\eqref{TVminusVlo}, after integration in $\lambda_0$, is
\begin{equation}\label{FVendpoint}
\begin{split}
-&\int_0^\infty\int_0^\infty
\int_{\delta\lambda_1}^\infty
P_{<M\delta\lambda_1}
\bigl(T_{\mathfrak m}\varphi_{\lambda_2}\bigr)
\cdot
S_{\lambda_1}\left(
T_{\mathfrak l}\varphi_{\lambda_0}
\nabla\eta_{<\delta\lambda_0}
\right)
\frac{d\lambda_2}{\lambda_2}
\frac{d\lambda_1}{\lambda_1}
\frac{d\lambda_0}{\lambda_0}.
\end{split}
\end{equation}
We now combine this principal contribution with \eqref{secondintppart} and with
$\mathcal A^{\mathrm{prin}}$. Indeed, using the basic identity
\begin{equation*}
\int_{\delta\lambda_1}^\infty
\varphi_{\lambda_2}\frac{d\lambda_2}{\lambda_2}
=
\varphi-P_{<\delta\lambda_1}\varphi,
\end{equation*}
we find that the sum of these three terms is
\begin{equation}\label{FVendpointcombined}
\begin{split}
\mathcal I\coloneq
\int_0^\infty\int_0^\infty
\Bigl(
P_{<M\delta\lambda_1}T_{\mathfrak m}
P_{<\delta\lambda_1}\varphi
-P_{<\delta\lambda_0}T_{\mathfrak m}\varphi
\Bigr)\cdot
S_{\lambda_1}\left(
T_{\mathfrak l}\varphi_{\lambda_0}
\nabla\eta_{<\delta\lambda_0}
\right)
\frac{d\lambda_1}{\lambda_1}
\frac{d\lambda_0}{\lambda_0}.
\end{split}
\end{equation}
We claim that this entire expression is also a perturbative error of the form $\mathcal{R}$; that is, it satisfies the $L_{t,x}^1$ bilinear bound \eqref{Fcubicremainder}. First, using
$P_{<M\delta\lambda_1}P_{<\delta\lambda_1}
=P_{<\delta\lambda_1}$, the coefficient
in \eqref{FVendpointcombined} becomes
\begin{equation*}
\begin{split}
P_{<M\delta\lambda_1}T_{\mathfrak m}
P_{<\delta\lambda_1}\varphi
-P_{<\delta\lambda_0}T_{\mathfrak m}\varphi=
\bigl(P_{<\delta\lambda_1}-P_{<\delta\lambda_0}\bigr)
T_{\mathfrak m}\varphi
+P_{<M\delta\lambda_1}
[T_{\mathfrak m},P_{<\delta\lambda_1}]\varphi.
\end{split}
\end{equation*}
We also have
\begin{equation*}
T_{\mathfrak l}\varphi_{\lambda_0}
=S_{\lambda_0}(T_{\mathfrak l}\varphi)
+[T_{\mathfrak l},S_{\lambda_0}]\varphi.
\end{equation*}

The terms containing either of these commutators contribute to
$\mathcal R$ since there are always to unbalanced pairs. Indeed, Proposition \ref{Leibnizcommutator} gives,
for $[T_{\mathfrak m},P_{<\delta\lambda_1}]\varphi$, terms of the form
$\mathcal L(\nabla^2\eta_\kappa,\varphi_\beta)$, where
$\beta\approx\delta\lambda_0$ and
$\kappa\lesssim_M\delta\beta$. We can then use the two unbalanced pairings
\begin{equation*}
(\overline u_\mu^{x_0},\eta_\kappa),
\qquad
(T_{\mathfrak l}\varphi_{\lambda_0},\varphi_\beta),
\end{equation*}
and estimate $\nabla\eta_{<\delta\lambda_0}$ in $L_{t,x}^\infty$. On the other hand, for the term involving $[T_{\mathfrak l},S_{\lambda_0}]\varphi$, we instead have terms of the
form $\mathcal L(\nabla^2\eta_\kappa,\varphi_{\lambda_0})$,
with $\kappa\lesssim_M\delta\lambda_0$. We pair
$\overline u_\mu^{x_0}$ with $\eta_\kappa$, and
$\varphi_{\lambda_0}$ with the component of $T_{\mathfrak m}\varphi$ at
frequency $\beta\approx\delta\lambda_0$. It remains to consider
\begin{equation*}
\int_0^\infty\int_0^\infty
\bigl(P_{<\delta\lambda_1}-P_{<\delta\lambda_0}\bigr)
T_{\mathfrak m}\varphi\cdot
S_{\lambda_1}\left(
\nabla\eta_{<\delta\lambda_0}
S_{\lambda_0}(T_{\mathfrak l}\varphi)
\right)
\frac{d\lambda_1}{\lambda_1}
\frac{d\lambda_0}{\lambda_0}.
\end{equation*}
Commuting $\nabla\eta_{<\delta\lambda_0}$ through $S_{\lambda_1}$ and using the antisymmetry of
$P_{<\delta\lambda_1}-P_{<\delta\lambda_0}$, this is
\begin{equation*}
\begin{split}
&\frac12\int_0^\infty\int_0^\infty
\bigl(P_{<\delta\lambda_1}-P_{<\delta\lambda_0}\bigr)
T_{\mathfrak m}\varphi\cdot
\bigl(\nabla\eta_{<\delta\lambda_0}
-\nabla\eta_{<\delta\lambda_1}\bigr)
S_{\lambda_1}S_{\lambda_0}(T_{\mathfrak l}\varphi)
\frac{d\lambda_1}{\lambda_1}
\frac{d\lambda_0}{\lambda_0}
\\
&+
\int_0^\infty\int_0^\infty
\bigl(P_{<\delta\lambda_1}-P_{<\delta\lambda_0}\bigr)
T_{\mathfrak m}\varphi\cdot
[S_{\lambda_1},\nabla\eta_{<\delta\lambda_0}]
S_{\lambda_0}(T_{\mathfrak l}\varphi)
\frac{d\lambda_1}{\lambda_1}
\frac{d\lambda_0}{\lambda_0}.
\end{split}
\end{equation*}
In the first integral, both the $T_{\mathfrak{m}}\varphi$ and $\nabla\eta$ are localized to frequencies comparable to $\delta\lambda_0$.
We pair one with $\overline u_\mu^{x_0}$ and the other with
$S_{\lambda_1}S_{\lambda_0}(T_{\mathfrak l}\varphi)$, using two unbalanced estimates. In the second integral, the commutator contributes a factor
$\lambda_0^{-1}\nabla^2\eta_\kappa$, with
$\kappa\lesssim\delta\lambda_0$. We pair $\eta_\kappa$ with
$\overline u_\mu^{x_0}$, and the localized
$T_{\mathfrak m}\varphi$ factor with
$S_{\lambda_0}(T_{\mathfrak l}\varphi)$ and then use two unbalanced bilinear estimates. Combining everything gives
\begin{equation}\label{FVendpointperturbative}
\|\mathcal I_\mu\overline u_\mu^{x_0}\|_{L_{t,x}^1}
\lesssim_{M,\delta}
\epsilon^4c_\mu^2\mu^{-2s_c-\frac12}
(1+\mu|x_0|)^{\frac12}.
\end{equation}
Hence, $\mathcal{I}$ is a perturbative error of the form $\mathcal{R}$, as claimed.

We next analyze the contribution of the differentiated term in \eqref{TVlinearize}. To keep the notation compact, set
\begin{equation*}
\mathcal T^{V,2}_{<\lambda_3}
\coloneq
\lambda_3\frac{d}{d\lambda_3}
\mathcal T^V_{<\lambda_3},
\end{equation*}
where the suppressed parameter $\lambda_2$ and all reference output
scales are held fixed. The linearization formula \eqref{Vlin} for $M$, Corollary
\ref{Philintrunc}, and the linearization formula \eqref{Blin} for $B$
give
\begin{equation}\label{TVsecondlinearization}
\begin{split}
\mathcal T^{V,2}_{<\lambda_3}
&=-M^{l,\lambda_2}_{<\lambda_3}
\left(
\eta_{\lambda_3}
L^{l,\lambda_2}_{<\lambda_3}
\Phi^{l,2}_{<\lambda_3}
\right)
-M^{l,\lambda_2}_{<\lambda_3}
\left(
\eta_{\lambda_2}L^{l,\lambda_3}_{<\lambda_3}
\Phi^{l,3}_{<\lambda_3}
\right)
\\
&+\eta_{\lambda_3}
M^{l,\lambda_2}_{<\lambda_3}
L^{l,\lambda_2}_{<\lambda_3}
\Phi^{l,2}_{<\lambda_3}
+\eta_{\lambda_2}
M^{l,\lambda_2}_{<\lambda_3}
L^{l,\lambda_3}_{<\lambda_3}
\Phi^{l,3}_{<\lambda_3}
\\
&-M^{l,\lambda_2}_{<\lambda_3}
\left(
\eta_{\lambda_2}\eta_{\lambda_3}
L^l_{<\lambda_3}B^l_{<\lambda_3}
\right)
+\eta_{\lambda_2}\eta_{\lambda_3}
M^{l,\lambda_2}_{<\lambda_3}
L^l_{<\lambda_3}B^l_{<\lambda_3}
+\mathcal E^{V,2}_{<\lambda_3}.
\end{split}
\end{equation}
Here $\mathcal E^{V,2}_{<\lambda_3}$ collects the high-frequency
tails produced by differentiating the matched output projections and
by comparing them with the fixed output projections in the
linearization formulas for $M$, $\Phi$, and $B$. More precisely,
Proposition \ref{Bernsteinelliptic}, Proposition \ref{basicopbounds}, Lemma
\ref{Phiest}, Proposition \ref{VBunbil}, and the usual
$L_t^4L_x^\infty$ and energy bounds give
\begin{equation*}
\begin{split}
&\left\|
\int_0^\infty
\left[
P_{<M\delta\lambda_1}
\int_{\delta\lambda_1}^\infty
\int_{\delta\lambda_2}^{\lambda_2}
\mathcal E^{V,2}_{<\lambda_3}
\frac{d\lambda_3}{\lambda_3}
\frac{d\lambda_2}{\lambda_2}
\right]\cdot
S_{\lambda_1}\left(
B_{\lambda_0}\nabla\eta_{<\delta\lambda_0}
\right)
\frac{d\lambda_1}{\lambda_1}
\right\|_{L_t^1L_x^2}
\\
&\qquad\lesssim
\epsilon^4c_{\lambda_0}
\lambda_0^{-\frac12-s_c}.
\end{split}
\end{equation*}
After integration in $\lambda_0$ and estimating
$\overline u_\mu^{x_0}$ separately in $L_t^\infty L_x^2$, this
contribution satisfies \eqref{Fcubicremainder}.

For fixed $\lambda_2$ and $\lambda_3$, let
$\mathcal T^{V,3}_{<\rho}$ denote the expression obtained from the
first six terms on the right-hand side of
\eqref{TVsecondlinearization}, excluding
$\mathcal E^{V,2}_{<\lambda_3}$, by replacing every geometric
truncation $<\lambda_3$ by $<\rho$, while keeping every other
parameter fixed. Applying the fundamental theorem of calculus once more gives
\begin{equation}\label{TVsecondFTC}
\mathcal T^{V,2}_{<\lambda_3}
=\mathcal T^{V,3}_{<\delta\lambda_3}
+\int_{\delta\lambda_3}^{\lambda_3}
\lambda_4\frac{d}{d\lambda_4}
\mathcal T^{V,3}_{<\lambda_4}
\frac{d\lambda_4}{\lambda_4}
+\mathcal E^{V,2}_{<\lambda_3}.
\end{equation}
The contribution of $\mathcal E^{V,2}_{<\lambda_3}$ was estimated
above, so we first consider the lower endpoint in
\eqref{TVsecondFTC}, which contributes terms in $\mathcal R$ together
with part of the cubic principal term $N_{\varphi,F}^{[3]}$ that we isolate below.
\\

First, the contribution of the last two terms in
\eqref{TVsecondlinearization} at this $<\delta\lambda_3$ endpoint belongs entirely to
$\mathcal R$. Indeed, after passing to the usual principal kernel
representations for $M$ and $L$, we can pair in $L_{t,x}^2$ the factor
$\eta_{\lambda_3}$ with the lower-frequency block
$L^l_{<\delta\lambda_3}B^l_{<\delta\lambda_3}$ in the fifth term and
with
$M^{l,\lambda_2}_{<\delta\lambda_3}L^l_{<\delta\lambda_3}
B^l_{<\delta\lambda_3}$ in the sixth term, using an unbalanced
bilinear estimate. We then estimate $\eta_{\lambda_2}$ and
$B_{\lambda_0}$ in $L_t^4L_x^\infty$ and
$\nabla\eta_{<\delta\lambda_0}$ in $L_{t,x}^\infty$. The factor
$\overline u_\mu^{x_0}$ is estimated separately in
$L_t^\infty L_x^2$.

For the first four terms, Corollary \ref{Phil0varphibound} allows us
to replace the truncated $\Phi$-factors by the corresponding
$\varphi$-factors, while estimating the remaining factor
$\eta_{\lambda_2}$ or $\eta_{\lambda_3}$ and $B_{\lambda_0}$ in
$L_t^4L_x^\infty$. Similarly, we replace the operators $L$ and $M$ by their
principal paradifferential operators $T_{\mathfrak l}$ and
$T_{\mathfrak m}$, respectively. The errors produced by these
replacements belong to $\mathcal R$. We estimate the replacement
of the exterior factor $B_{\lambda_0}$ below. The retained
principal expression inside
$\mathcal T^{V,3}_{<\delta\lambda_3}$ is therefore
\begin{equation}\label{Vquadraticprincipal}
\begin{split}
\mathcal{T}^{V,\mathrm{main}}_{<\delta\lambda_3}\coloneq -T_{\mathfrak m}\left(
\eta_{\lambda_3}T_{\mathfrak l}\varphi_{\lambda_2}
+\eta_{\lambda_2}T_{\mathfrak l}\varphi_{\lambda_3}
\right)+\eta_{\lambda_3}
T_{\mathfrak m}T_{\mathfrak l}\varphi_{\lambda_2}
+\eta_{\lambda_2}
T_{\mathfrak m}T_{\mathfrak l}\varphi_{\lambda_3}.
\end{split}
\end{equation}
The contribution of
$\mathcal T^{V,\mathrm{main}}_{<\delta\lambda_3}$
from \eqref{Vquadraticprincipal} to the first integral in
\eqref{TVminusVlo}, after integration in $\lambda_0$, is
\begin{equation}\label{TV3maincont}
\begin{split}
-&\int_0^\infty\int_0^\infty
\int_{\delta\lambda_1}^\infty
\int_{\delta\lambda_2}^{\lambda_2}
P_{<M\delta\lambda_1}
\bigl(\mathcal T^{V,\mathrm{main}}_{<\delta\lambda_3}\bigr)
\cdot
S_{\lambda_1}\left(
B_{\lambda_0}\nabla\eta_{<\delta\lambda_0}
\right)
\frac{d\lambda_3}{\lambda_3}
\frac{d\lambda_2}{\lambda_2}
\frac{d\lambda_1}{\lambda_1}
\frac{d\lambda_0}{\lambda_0}.
\end{split}
\end{equation}
We retain this term and will combine it below with the principal
contribution from $\mathcal A^{\mathrm{sub}}$ in \eqref{Fannularsplit}. Before doing this, we first consider the contribution of the integral term in \eqref{TVsecondFTC}. Its
treatment is the same as that of
$\mathcal T^{B,4}_{<\lambda_0,\lambda_1}$ above and gives, when
inserted into the first integral in \eqref{TVminusVlo}, a contribution of the form
\begin{equation*}
N_{\varphi,F}^{[4]}+\mathcal R.
\end{equation*}
Indeed, the $\lambda_4$ derivative in \eqref{TVsecondFTC}
produces at the principal level a cubic expression, with bounded
paracoefficients, in $\varphi$ and $\eta$ with input frequencies
$\lambda_2\approx_\delta\lambda_3\approx_\delta\lambda_4$
and $\lambda_2\gtrsim_\delta\lambda_0$. The remainder is quartic
and can be estimated in $L_t^1L_x^\infty$ using four
$L_t^4L_x^\infty$ bounds. We estimate $B_{\lambda_0}$ in
$L_t^\infty L_x^2$, which makes this contribution perturbative
in $L_t^1L_x^2$.

To replace $B_{\lambda_0}$ by
$T_{\mathfrak l}\varphi_{\lambda_0}$ in the principal cubic
contribution, we use Corollary \ref{principalVBexp} and the
corresponding commutator bound:
\begin{equation*}
\|B_{\lambda_0}-T_{\mathfrak l}\varphi_{\lambda_0}
\|_{L_t^2L_x^\infty}
\lesssim
\epsilon^2c_{\lambda_0}\lambda_0^{\frac54-s_c}.
\end{equation*}
We estimate the principal cubic part of
$\lambda_4\frac{d}{d\lambda_4}
\mathcal T^{V,3}_{<\lambda_4}$
in $L_{t,x}^2$ by placing its $\varphi$ input in
$L_t^\infty L_x^2$ and its two $\eta$ inputs in
$L_t^4L_x^\infty$. This gives
$O_\delta(\epsilon^3\lambda_2^{-7/4})$.
H\"older's inequality and integration in the internal
frequencies therefore make the replacement error
perturbative in $L_t^1L_x^2$.

When $\lambda_2\approx_\delta\lambda_0$, the remaining
expression has four comparable inputs of total differential
order four, (possibly with $\nabla\eta_{<\delta\lambda_0}$ as a bounded
coefficient). This gives a term of the form
$N_{\varphi,F}^{[4]}$. When
$\lambda_2\gg_\delta\lambda_0\approx\mu$, we use two
unbalanced estimates, pairing two of the internal inputs
with $T_{\mathfrak l}\varphi_{\lambda_0}$ and
$\overline u_\mu^{x_0}$, respectively, and estimating the
remaining factors in $L_{t,x}^\infty$.
This contribution belongs to $\mathcal R$.

At this point, it is now natural to return to identifying the contribution of $\mathcal A^{\mathrm{sub}}$ from
\eqref{Fannularsplit}, which is very similar to the above (but with slightly different frequency cutoffs). First, by applying the same expansion to its
$V$-factor, we obtain
\begin{equation}\label{FVannularquadratic}
\begin{split}
\mathcal A^{\mathrm{sub}}
&=
\int_0^\infty
\int_{\delta\lambda_0}^\infty
\int_{\delta\lambda_2}^{\lambda_2}
\bigl(P_{<M\delta\lambda_0}-P_{<\delta\lambda_0}\bigr)
\bigl(\mathcal T^{V,\mathrm{main}}_{<\delta\lambda_3}\bigr)
\cdot
\bigl(B_{\lambda_0}\nabla\eta_{<\delta\lambda_0}\bigr)
\frac{d\lambda_3}{\lambda_3}
\frac{d\lambda_2}{\lambda_2}
\frac{d\lambda_0}{\lambda_0}
\\
&+N_{\varphi,F}^{[4]}+\mathcal R.
\end{split}
\end{equation}
The errors at the lower endpoint are estimated as above,
using \eqref{Blambda0principalcomparison} for the replacement
of $B_{\lambda_0}$ by $T_{\mathfrak l}\varphi_{\lambda_0}$.
The higher-order contributions are treated by the preceding
argument, with
$P_{<M\delta\lambda_0}-P_{<\delta\lambda_0}$ in place of
$P_{<M\delta\lambda_1}$, and give a term of the form
$N_{\varphi,F}^{[4]}+\mathcal R$. We omit the details since they are similar. We now combine the first line of \eqref{FVannularquadratic}
with the contribution retained in \eqref{TV3maincont}. First, we replace the lower limit
$\delta\lambda_1$ in \eqref{TV3maincont} by $\delta\lambda_0$.
Since $\lambda_1\approx\lambda_0$, the difference has
$\lambda_2\approx\delta\lambda_0$ and
$\delta\lambda_2\leq\lambda_3\leq\lambda_2$.
Both inputs are therefore separated from
$\mu\approx\lambda_0$, and we estimate this error using two
unbalanced pairings with $B_{\lambda_0}$ and
$\overline u_\mu^{x_0}$, using Proposition \ref{VBunbil}.

We next replace $P_{<M\delta\lambda_1}$ by
$P_{<M\delta\lambda_0}$ modulo $\mathcal R$, using the same
commutation and antisymmetrization argument as for
$\mathcal E^V_{\lambda_0,\lambda_1}$ above, with the changes
in the lower integration limit estimated as above.
Here the quadratic expression satisfies the same
$L_t^2L_x^\infty$ bound used in that argument, by two
$L_t^4L_x^\infty$ estimates.
After integration in $\lambda_1$, the two terms involving
$P_{<M\delta\lambda_0}$ cancel. Thus, modulo $\mathcal R$,
the combined contribution of the first line of
\eqref{FVannularquadratic} and \eqref{TV3maincont} is
\begin{equation}\label{F11principal}
\begin{split}
-&\int_0^\infty
\int_{\delta\lambda_0}^\infty
\int_{\delta\lambda_2}^{\lambda_2}
P_{<\delta\lambda_0}
\bigl(\mathcal T^{V,\mathrm{main}}_{<\delta\lambda_3}\bigr)
\cdot
\bigl(B_{\lambda_0}\nabla\eta_{<\delta\lambda_0}\bigr)
\frac{d\lambda_3}{\lambda_3}
\frac{d\lambda_2}{\lambda_2}
\frac{d\lambda_0}{\lambda_0}.
\end{split}
\end{equation}

We now replace $B_{\lambda_0}$ by
$T_{\mathfrak l}\varphi_{\lambda_0}$ in \eqref{F11principal}.
When $\lambda_2\not\sim_c\mu$, we pair
$\overline u_\mu^{x_0}$ with the $\lambda_2$ input of
$\mathcal T^{V,\mathrm{main}}_{<\delta\lambda_3}$,
put its other input in $L_t^\infty L_x^2$, and estimate
$B_{\lambda_0}-T_{\mathfrak l}\varphi_{\lambda_0}$ in
$L_t^2L_x^\infty$ using Corollary \ref{principalVBexp}
and Proposition \ref{Leibnizcommutator}.

For $\lambda_2\sim_c\mu$, we first subdivide the input and
output frequencies contained in the principal quadratic part of $B_{\lambda_0}-T_{\mathfrak{l}}\varphi_{\lambda_0}$ sufficiently narrowly. The linearization
formula \eqref{Bintegralexp}, Lemma \ref{Phiest}, and
Proposition \ref{paraexp} give the decomposition
\begin{equation}\label{QEdecomp}
B_{\lambda_0}-T_{\mathfrak l}\varphi_{\lambda_0}
=Q_{\lambda_0}+E_{\lambda_0},
\qquad
\|E_{\lambda_0}\|_{L_{t,x}^2}
\lesssim\epsilon^2c_{\lambda_0}
\lambda_0^{\frac14-s_c}.
\end{equation}
Here $Q_{\lambda_0}$ collects the quadratic terms
$S_{\lambda_0}\mathcal L(\eta_\alpha,\varphi_\beta)$
whose comparable frequency inputs have output frequency above the range of
Lemma \ref{balancedrestrictedoutput}. The forms
$\mathcal L$ carry two derivatives and bounded
paracoefficients. The remaining quadratic components (which we group into $E_{\lambda_0}$) satisfy the $L_{t,x}^2$ error bound above by using Lemma \ref{balancedrestrictedoutput}. The contribution of $E_{\lambda_0}$ in \eqref{F11principal} can then be estimated by placing $u_{\mu}^{x_0}$ in $L_t^\infty L_x^2$ and using two  
$L_t^4L_x^\infty$ bounds for $\mathcal{T}^{V,\mathrm{main}}_{<\delta\lambda_3}$. For $Q_{\lambda_0}$, by its definition, one of its input frequencies must be separated from $\lambda_2\sim_c\mu$. We therefore pair one of its inputs with the
$\lambda_2$ input of
$\mathcal T^{V,\mathrm{main}}_{<\delta\lambda_3}$ and use an unbalanced $L_{t,x}^2$ bound. Two of the remaining factors can be put in $L_t^4L_x^\infty$,
and we can estimate $\overline u_\mu^{x_0}$ in
$L_t^\infty L_x^2$. Collecting everything above and integrating gives
\eqref{Fcubicremainder}. Thus, modulo $\mathcal R$,
the integrand becomes
\begin{equation}\label{peeledintegrand}
-P_{<\delta\lambda_0}
(\mathcal{T}^{V,\mathrm{main}}_{<\delta\lambda_3})\cdot
\left(
T_{\mathfrak l}\varphi_{\lambda_0}
\nabla\eta_{<\delta\lambda_0}
\right).
\end{equation}
This finally concludes our analysis of $F^{1,1}$.
We will later recall \eqref{peeledintegrand} in our final
assembly of $N_{\varphi,F}^{[3]}$.

\paragraph{\textbf{The term $F^{1,2}$.}}
We next consider the remaining component $F^{1,2}$ in
\eqref{F1cubicexp}, which is given by
\begin{equation*}
F^{1,2}=-\int_0^\infty
B_{\lambda_0}P_{<\delta\lambda_0}
(V^{\mathrm{hi}}\cdot\nabla\eta)
\frac{d\lambda_0}{\lambda_0}.
\end{equation*}
This term is far easier to handle. Indeed, the low output projection
$P_{<\delta\lambda_0}$ forces the $V$ and $\nabla\eta$ inputs to have
comparable frequencies $\alpha\gtrsim M\delta\lambda_0$. If $\alpha\gg_{M,\delta}\lambda_0$, then both the output frequency
$\mu\approx\lambda_0$ and the $B$ input at frequency $\lambda_0$ are
separated from $\alpha$. In this case, we can use two unbalanced
bilinear estimates, one in which we pair
$(B_{\lambda_0},V_\alpha)$ and the other in which we pair
$(\overline u_\mu^{x_0},\nabla\eta_\alpha)$, and then sum over
$\alpha$ to treat this contribution perturbatively. Therefore, it
remains to consider the range in which $B_{\lambda_0}$, $V$, and
$\eta$ are all localized at comparable frequencies. We first replace the relevant frequency part of $V^{\mathrm{hi}}$
by its principal part, using the $Q+E$ decomposition above and
its analogue for $V$. We estimate the $E$ contribution in
$L_{t,x}^2$, with $B_{\lambda_0}$ and $\nabla\eta_\alpha$
in $L_t^4L_x^\infty$ and $\overline u_\mu^{x_0}$ in
$L_t^\infty L_x^2$. For the $Q$ contribution, its inputs
are separated from its output frequency $\alpha$. We therefore
pair its $\varphi$ input with $\eta_\alpha$, put its other
input and $B_{\lambda_0}$ in $L_t^4L_x^\infty$, and again
estimate $\overline u_\mu^{x_0}$ in $L_t^\infty L_x^2$.
We next replace $B_{\lambda_0}$ by its principal part using
the frequency decomposition and pairings leading to
\eqref{peeledintegrand}, now with the quadratic factor
$P_{<\delta\lambda_0}
((T_{\mathfrak m}\varphi)_\alpha\cdot\nabla\eta_\alpha)$
in place of
$P_{<M\delta\lambda_1}
\mathcal T^{V,\mathrm{main}}_{<\delta\lambda_3}$.
Retaining the inherited input and output localizations and
integrating gives \eqref{Fcubicremainder} for all these errors. Combining this with the definition of $V^{hi}$ and the frequency support forced by the output projection $P_{<\delta\lambda_0}$, we can write
\begin{equation}\label{F12principal}
\begin{split}
F^{1,2}=-\int_0^\infty
\int_{M\delta\lambda_0}^{\delta^{-1}\lambda_0}
\int_{\delta\lambda_1}^{\delta^{-1}\lambda_1}
(T_{\mathfrak l}\varphi)_{\lambda_0}
P_{<\delta\lambda_0}
\left(
S_{\lambda_1}(T_{\mathfrak m}\varphi)
\cdot\nabla\eta_{\lambda_2}
\right)
\frac{d\lambda_2}{\lambda_2}
\frac{d\lambda_1}{\lambda_1}
\frac{d\lambda_0}{\lambda_0}
+\mathcal R.
\end{split}
\end{equation}
Although the low-output projection
actually forces $\lambda_2\approx\lambda_1$, we have harmlessly enlarged this
to the $[\delta\lambda_1,\delta^{-1}\lambda_1]$ range above. The motivation for writing the integration bounds this way is so that it is later transparent how this term assembles with the other components of the cubic balanced term $N_{\varphi,F}^{[3]}$ later. At this stage, we remark that \eqref{F12principal} gives only the portion selected
by the internal low-output projection $P_{<\delta\lambda_0}$. We retain this projection
until this term is combined with the corresponding part of
$F^{2,1}$ below. That term will supply the complementary output projection, and the remaining components to assemble the corresponding part of $N_{\varphi,F}^{[3]}$. This concludes our analysis for $F^{1,2}$, and thus the full term $F^1$. We now turn to $F^2$. Let us set $\widetilde S_{\lambda_0}
\coloneq P_{<\delta^{-1}\lambda_0}-P_{<\delta\lambda_0}.$ Simple algebraic manipulation, and using Proposition \ref{Leibnizcommutator}, we can expand $F^2$ into four terms
\begin{equation}\label{F2fourparts}
\begin{split}
F^2&=
\int_0^\infty
\left[
T_{V^l}\cdot\nabla\eta_{\lambda_0}
-(V\cdot\nabla\eta)_{\lambda_0}
\right]
\widetilde S_{\lambda_0}
(B-B_{<\delta\lambda_0})
\frac{d\lambda_0}{\lambda_0}
\\
&-
\int_0^\infty
\nabla\eta_{<\delta\lambda_0}\cdot V_{\lambda_0}
P_{<\delta\lambda_0}
(B-B_{<\delta\lambda_0})
\frac{d\lambda_0}{\lambda_0}
\\
&+\int_0^\infty
\lambda_0^{-1}\mathcal L
\left(
\widetilde P_{\lambda_0}V,
\nabla^2\eta^{\mathrm{lo}}
\right)
P_{<\delta\lambda_0}
(B-B_{<\delta\lambda_0})
\frac{d\lambda_0}{\lambda_0}
\\
&+\int_0^\infty
\left[
T_{V^l}\cdot\nabla\eta_{\lambda_0}
-(V\cdot\nabla\eta^{\mathrm{hi}})_{\lambda_0}
\right]
P_{<\delta\lambda_0}
(B-B_{<\delta\lambda_0})
\frac{d\lambda_0}{\lambda_0}
\\
&\coloneq F^{2,1}+F^{2,2}+F^{2,3}+F^{2,4}.
\end{split}
\end{equation}
Here the superscripts $\mathrm{lo}$ and $\mathrm{hi}$ are always
taken relative to the scale $\lambda_0$ (as above). Now, we analyze each of the four terms above (in order of convenience).
\\

\paragraph{\textbf{The term $F^{2,2}$.}}
The argument for $F^{2,2}$ follows the successive linearizations
\eqref{VtailFTCstart}, \eqref{TVlinearize}, and
\eqref{TVsecondFTC}, now applied to the $B$-difference
$P_{<\delta\lambda_0}(B-B_{<\delta\lambda_0})$.
We record the contributions of the first two endpoints, which
will be combined with the corresponding terms from
$F^{2,1}$ below. At the lower endpoint corresponding to \eqref{TVlinearize},
we obtain the analogue of \eqref{FVendpoint},
\begin{equation}\label{FBendpoint}
\begin{split}
\mathcal J_B\coloneq
-\int_0^\infty
S_{\lambda_0}(T_{\mathfrak m}\varphi)
\cdot\nabla\eta_{<\delta\lambda_0}
P_{<\delta\lambda_0}T_{\mathfrak l}(\varphi-\varphi_{<\delta\lambda_0})
\frac{d\lambda_0}{\lambda_0}.
\end{split}
\end{equation}
For the linearized $B$-factor, we use the analogue of
\eqref{survivingprincipal}, with $L$ and $T_{\mathfrak l}$
in place of $M$ and $T_{\mathfrak m}$, and estimate the errors
by the same pairings used there. The replacement of the $V_{\lambda_0}$ factor is justified as in
\eqref{Blambda0principalcomparison}, with the roles of $B$
and $V$ interchanged. We will combine $\mathcal J_B$ with a complimentary term from $F^{2,1}$ that we will identify below.

At the lower endpoint corresponding to \eqref{TVsecondFTC},
we obtain the quadratic principal analogue of
\eqref{Vquadraticprincipal},

\begin{equation}\label{TBmain2ord}
\begin{split}
\mathcal T^{B,\mathrm{main}}_{<\delta\lambda_2}
\coloneq
-T_{\mathfrak l}\left(
\eta_{\lambda_2}T_{\mathfrak l}\varphi_{\lambda_1}
+\eta_{\lambda_1}T_{\mathfrak l}\varphi_{\lambda_2}
\right)+\eta_{\lambda_2}T_{\mathfrak l}^2\varphi_{\lambda_1}
+\eta_{\lambda_1}T_{\mathfrak l}^2\varphi_{\lambda_2}.
\end{split}
\end{equation}
Arguing as in the passage from \eqref{F11principal} to
\eqref{peeledintegrand}, using the $V$ analogue of
\eqref{QEdecomp}, we obtain the
principal cubic contribution

\begin{equation*}
\begin{split}
&-\int_0^\infty
\int_{\delta\lambda_0}^\infty
\int_{\delta\lambda_1}^{\lambda_1}
T_{\mathfrak m}\varphi_{\lambda_0}\cdot
\nabla\eta_{<\delta\lambda_0}
P_{<\delta\lambda_0}
\left(
\mathcal T^{B,\mathrm{main}}_{<\delta\lambda_2}
\right)
\frac{d\lambda_2}{\lambda_2}
\frac{d\lambda_1}{\lambda_1}
\frac{d\lambda_0}{\lambda_0}.
\end{split}
\end{equation*}
The remaining contributions
are treated as in the argument for the integral term in
\eqref{TVsecondFTC} and are of the form
$N_{\varphi,F}^{[4]}+\mathcal R$.
\\

\paragraph{\textbf{The term $F^{2,3}$.}}
The contribution of $F^{2,3}$ is entirely perturbative. First, the unbalanced estimate in Proposition \ref{VBunbil} along with dyadic summation gives
\begin{equation}\label{Vetabil}
\lambda_0^{-1}
\left\|
\mathcal L(
\widetilde P_{\lambda_0}V,
\nabla^2\eta^{\mathrm{lo}})
\right\|_{L_{t,x}^2}
\lesssim
\epsilon^2c_{\lambda_0}
\lambda_0^{\frac14-s_c}.
\end{equation}
The linearization formula for $B$, Corollary
\ref{principalVBexp}, and Proposition \ref{basicopbounds} give
\begin{equation}\label{F23bexp}
\begin{split}
P_{<\delta\lambda_0}
(B-B_{<\delta\lambda_0})
=P_{<\delta\lambda_0}
T_{\mathfrak l}\widetilde\varphi_{\delta\lambda_0}
+O_{L_t^2L_x^{\infty}}(\epsilon^2\lambda_0^{-\frac{3}{4}}).
\end{split}
\end{equation}
where $\widetilde\varphi_{\delta\lambda_0}$ is localized to a fixed
fattening of the annulus $|\xi|\approx\delta\lambda_0$. Combining
\eqref{Vetabil} and \eqref{F23bexp}, and using one additional
unbalanced estimate with the pair
$(\overline u_\mu^{x_0},
P_{<\delta\lambda_0}
T_{\mathfrak l}\widetilde\varphi_{\delta\lambda_0})$, gives
\begin{equation*}
\|F^{2,3}_\mu\overline u_\mu^{x_0}\|_{L_{t,x}^1}
\lesssim
\epsilon^4c_\mu^2\mu^{-\frac12-2s_c}
(1+\mu|x_0|)^{\frac12}.
\end{equation*}

\paragraph{\textbf{The term $F^{2,4}$.}}
This term is also perturbative. Using \eqref{bilformrep}, we write
\begin{equation}\label{VetaF24}
-(V\cdot\nabla\eta^{\mathrm{hi}})_{\lambda_0}
=-V^{\mathrm{lo}}\cdot\nabla\eta_{\lambda_0}
+\mathcal E_{\lambda_0}.
\end{equation}
We estimate the contribution of $\mathcal E_{\lambda_0}$ after
multiplication by
$P_{<\delta\lambda_0}(B-B_{<\delta\lambda_0})$ and projecting the whole product to output frequency $\mu$, and pairing with $\overline u_\mu^{x_0}$.
Here, we have by frequency support, $\mu\approx\lambda_0$. We treat the contribution of $\mathcal{E}_{\lambda_0}$ as in the proof of \eqref{replacebyphi}, by peeling off the principal part of $V$, and applying Lemma
\ref{balancedrestrictedoutput} (and handling the subprincipal errors as in \eqref{replacebyphi}).

The remaining subprincipal errors are handled by the error bounds in
Corollary \ref{principalVBexp} and \eqref{F23bexp}, the $L_t^{\infty}L_x^2$ bound, and the same bilinear pairings. Consequently,
\begin{equation}\label{F24reduced}
F^{2,4}
=\int_0^\infty
(T_{V^l}-V^{\mathrm{lo}})\cdot\nabla\eta_{\lambda_0}
P_{<\delta\lambda_0}
(B-B_{<\delta\lambda_0})
\frac{d\lambda_0}{\lambda_0}
+\mathcal R.
\end{equation}
The first integral is perturbative by two unbalanced estimates. For
the second, Corollary \ref{principalVBexp} gives
\begin{equation*}
\begin{split}
P_{<M\delta\lambda_1}
(V_{<\delta\lambda_1}-V)
&=-P_{<M\delta\lambda_1}
T_{\mathfrak m}\widetilde\varphi_{\delta\lambda_1}
+O_{L_t^2L_x^{\infty}}(\epsilon^2\lambda_0^{-\frac34}).
\end{split}
\end{equation*}
Together with \eqref{F23bexp}, this reduces (modulo perturbative errors) the integrand
to a trilinear form with two inputs at frequencies
$\approx_M\delta\lambda_0$ and the factor
$\nabla\eta_{\lambda_0}$. Pairing each low input with one of
$\overline u_\mu^{x_0}$ and $\nabla\eta_{\lambda_0}$ and using two unbalanced bilinear estimates gives
\begin{equation*}
\|F^{2,4}_\mu\overline u_\mu^{x_0}\|_{L_{t,x}^1}
\lesssim
\epsilon^4c_\mu^2\mu^{-\frac12-2s_c}
(1+\mu|x_0|)^{\frac12}.
\end{equation*}

\paragraph{\textbf{The term $F^{2,1}$ and assembly of $N_{\varphi,F}$}} Here, our goal is to expand the last term
$F^{2,1}$ and then isolate $N_{\varphi,F}^{[2]}$ as well as combine the cubic balanced component of $F^{2,1}$ below with the components of $N_{\varphi,F}^{[3]}$ isolated above. First, Proposition \ref{Leibnizcommutator} gives
\begin{equation}\label{VetaLeibnizF21}
\begin{split}
(V\cdot\nabla\eta)_{\lambda_0}
={}&\lambda_0^{-1}\mathcal L
(\widetilde P_{\lambda_0}V,
\nabla^2\eta^{\mathrm{lo}})
+\lambda_0^{-1}\mathcal L
(\nabla V^{\mathrm{lo}},
\nabla\widetilde P_{\lambda_0}\eta)
\\
&+(V^{\mathrm{hi}}\cdot\nabla\eta^{\mathrm{hi}})_{\lambda_0}
+\nabla\eta^{\mathrm{lo}}\cdot V_{\lambda_0}
+V^{\mathrm{lo}}\cdot\nabla\eta_{\lambda_0}.
\end{split}
\end{equation}
The first two terms contribute to $\mathcal R$. After projecting to output frequency $\mu$ and pairing with
$\overline u_\mu^{x_0}$, we use the frequency decomposition and
bilinear pairings from the proof of \eqref{replacebyphi},
and applying Lemma \ref{balancedrestrictedoutput}, and then treating
the subprincipal errors as above. Consequently, we have
\begin{equation}\label{F21reduced}
\begin{split}
F^{2,1}={}&-\int_0^\infty
(V^{\mathrm{hi}}\cdot\nabla\eta^{\mathrm{hi}})_{\lambda_0}
\widetilde S_{\lambda_0}(B-B_{<\delta\lambda_0})
\frac{d\lambda_0}{\lambda_0}
\\
&-\int_0^\infty
\nabla\eta^{\mathrm{lo}}\cdot V_{\lambda_0}
\widetilde S_{\lambda_0}(B-B_{<\delta\lambda_0})
\frac{d\lambda_0}{\lambda_0}
\\
&+\int_0^\infty
(T_{V^l}-V^{\mathrm{lo}})\cdot\nabla\eta_{\lambda_0}
\widetilde S_{\lambda_0}(B-B_{<\delta\lambda_0})
\frac{d\lambda_0}{\lambda_0}
\,+\mathcal R.
\end{split}
\end{equation}
For the first term in \eqref{F21reduced}, the contribution
containing $\widetilde S_{\lambda_0}B_{<\delta\lambda_0}$
is perturbative by the same frequency separation argument as used in \eqref{replacebyphi}, applying Proposition
\ref{VBunbil} to its low-frequency part and Proposition
\ref{Bernsteinelliptic} to its high-frequency tails.

We then replace the remaining $B$ and $V$ factors by their
principal parts using the $Q+E$ decomposition in \eqref{QEdecomp} and its $V$ analogue.
For the $E$ contributions, we use the $L_{t,x}^2$ bound in \eqref{QEdecomp},
put the other two factors each in $L_t^4L_x^\infty$, and
estimate $\overline u_\mu^{x_0}$ in $L_t^\infty L_x^2$.
Terms containing both a $B-T_{\mathfrak{l}}\varphi$ error and a $V-T_{\mathfrak{m}}\varphi$ error can be handled by estimating individually each of these terms in $L_t^2L_x^\infty$ using Corollary \ref{principalVBexp}, and bounding both
$\nabla\eta$ and $\overline u_\mu^{x_0}$ in
$L_t^\infty L_x^2$. 

For a contribution containing a single $Q$, there are four
$\approx\lambda_0$ inputs. After subdividing the frequency annuli,
if $\mu$ is separated from the largest input frequency (call it $\Lambda$),
we pair $\overline u_\mu^{x_0}$ with that input, put one
remaining factor in $L_t^\infty L_x^2$, and the other two
in $L_t^4L_x^\infty$. Similarly if any of the four input frequencies are separated, we can estimate at least one pair of them using an unbalanced bilinear estimate. Otherwise, we have $|\xi_1+\xi_2+\xi_3+\xi_4|\sim_c\Lambda$ and $|\xi_j|\sim_c\Lambda$ for every $j$. In this case, at least one pair of the six possible pairs of input frequencies $\xi_i$ and $\xi_j$ satisfies $\delta\Lambda< |\xi_i+\xi_j|<\frac{3}{2}\Lambda$ for some $\delta>0$ which we can apply Lemma \ref{balancedrestrictedoutput} to; we put
the other two inputs in $L_t^4L_x^\infty$ and
$\overline u_\mu^{x_0}$ in $L_t^\infty L_x^2$. Integrating gives
\eqref{Fcubicremainder} for all these errors. Consequently,
\begin{equation}\label{Bannularprincipal}
\begin{split}
-\int_0^\infty
S_{\lambda_0}\left(
V^{\mathrm{hi}}\cdot\nabla\eta^{\mathrm{hi}}\right)
\widetilde S_{\lambda_0}(B-B_{<\delta\lambda_0})
\frac{d\lambda_0}{\lambda_0}=
-\int_0^\infty
S_{\lambda_0}\left(
(T_{\mathfrak m}\varphi)^{\mathrm{hi}}
\cdot\nabla\eta^{\mathrm{hi}}\right)
\widetilde S_{\lambda_0}(T_{\mathfrak l}\varphi)
\frac{d\lambda_0}{\lambda_0}
+\mathcal R.
\end{split}
\end{equation}
By Fubini, we can interchange $\tilde{S}_{\lambda_0}$ and $S_{\lambda_0}$, to write this as
\begin{equation*}
-\int_0^\infty
\widetilde S_{\lambda_0}\left(
(T_{\mathfrak m}\varphi)^{\mathrm{hi}}
\cdot\nabla\eta^{\mathrm{hi}}
\right) S_{\lambda_0}(T_{\mathfrak l}\varphi)
\frac{d\lambda_0}{\lambda_0}.
\end{equation*}
As in the analysis immediately preceding \eqref{F12principal}, we can further write this term as
\begin{equation}\label{F21tripleprincipal}
\begin{split}
-&\int_0^\infty
\int_{M\delta\lambda_0}^{\delta^{-1}\lambda_0}
\int_{\delta\lambda_1}^{\delta^{-1}\lambda_1}
S_{\lambda_0}(T_{\mathfrak l}\varphi)
\widetilde S_{\lambda_0}
\left(
S_{\lambda_1}(T_{\mathfrak m}\varphi)
\cdot\nabla\eta_{\lambda_2}
\right)
\frac{d\lambda_2}{\lambda_2}
\frac{d\lambda_1}{\lambda_1}
\frac{d\lambda_0}{\lambda_0}+\mathcal{R}.
\end{split}
\end{equation}
Since
\begin{equation*}
P_{<\delta\lambda_0}+\widetilde S_{\lambda_0}
=P_{<\delta^{-1}\lambda_0},
\end{equation*}
the terms in \eqref{F12principal} and
\eqref{F21tripleprincipal} combine to give
\begin{equation}\label{F21tripleprincipal+}
\begin{split}
-&\int_0^\infty
\int_{M\delta\lambda_0}^{\delta^{-1}\lambda_0}
\int_{\delta\lambda_1}^{\delta^{-1}\lambda_1}
S_{\lambda_0}(T_{\mathfrak l}\varphi)
P_{<\delta^{-1}\lambda_0}
\left(
S_{\lambda_1}(T_{\mathfrak m}\varphi)
\cdot\nabla\eta_{\lambda_2}
\right)
\frac{d\lambda_2}{\lambda_2}
\frac{d\lambda_1}{\lambda_1}
\frac{d\lambda_0}{\lambda_0}+\mathcal{R}.
\end{split}
\end{equation}
We extend the $\lambda_1$ integration down to $\delta\lambda_0$
modulo $\mathcal R$, using two bilinear estimates. Combining the above and using our
trilinear paraproduct convention gives
\begin{equation}\label{Fhhh}
\begin{split}
F^{1,2}-\int_0^\infty
(V^{\mathrm{hi}}\cdot\nabla\eta^{\mathrm{hi}})_{\lambda_0}
\widetilde S_{\lambda_0}(B-B_{<\delta\lambda_0})
\frac{d\lambda_0}{\lambda_0}=-\Pi\bigl(
\nabla\eta,
T_{\mathfrak l}\varphi,
T_{\mathfrak m}\varphi
\bigr)+\mathcal{R}
\end{split}
\end{equation}
We next consider the second term in \eqref{F21reduced}. First, by
two bilinear estimates, we can replace this term by
\begin{equation}\label{F21secondterm}
-T_{\nabla\eta}\cdot
\int_0^\infty
V_{\lambda_0}
\widetilde S_{\lambda_0}(B-B_{<\delta\lambda_0})
\frac{d\lambda_0}{\lambda_0}.
\end{equation}
Modulo terms whose contributions to \eqref{F21secondterm} are of the
form $N_{\varphi,F}^{[4]}+\mathcal R$, the formulas
\eqref{Vquadraticprincipal} and \eqref{TBmain2ord} give
\begin{equation}\label{F21VBmainexpansions}
\begin{split}
V_{\lambda_0}
&=S_{\lambda_0}(T_{\mathfrak m}\varphi)
+\int_{\delta\lambda_0}^{\infty}
\int_{\delta\lambda_1}^{\lambda_1}
S_{\lambda_0}\left(
\mathcal T^{V,\mathrm{main}}_{<\delta\lambda_2}
\right)
\frac{d\lambda_2}{\lambda_2}
\frac{d\lambda_1}{\lambda_1},
\\
\widetilde S_{\lambda_0}(B-B_{<\delta\lambda_0})
&=
\widetilde S_{\lambda_0}T_{\mathfrak l}
(\varphi-\varphi_{<\delta\lambda_0})+\int_{\delta\lambda_0}^{\infty}
\int_{\delta\lambda_1}^{\lambda_1}
\widetilde S_{\lambda_0}\left(
\mathcal T^{B,\mathrm{main}}_{<\delta\lambda_2}
\right)
\frac{d\lambda_2}{\lambda_2}
\frac{d\lambda_1}{\lambda_1}.
\end{split}
\end{equation}
Here the lower bound $\lambda_1\geq\delta\lambda_0$ follows from output
support.
\\

We now insert \eqref{F21VBmainexpansions} into
\eqref{F21secondterm}. We first combine the contribution of the product of the two
linear terms in \eqref{F21VBmainexpansions} with $\mathcal J_B$ from \eqref{FBendpoint}.
In $\mathcal J_B$, the factor
$P_{<\delta\lambda_0}T_{\mathfrak l}
(\varphi-\varphi_{<\delta\lambda_0})$ has frequency
$\beta\approx\delta\lambda_0$, so the output frequency satisfies
$\mu\approx\lambda_0$. As in the calculation leading up to
\eqref{F21secondterm}, we can replace multiplication
by $\nabla\eta_{<\delta\lambda_0}$ in $\mathcal{J}_{B}$ by the corresponding
paraproduct $T_{\nabla\eta}$ modulo a perturbative term of the form $\mathcal R$. Indeed, the
difference of the error produces a factor of $\nabla\eta$ localized at frequency
$\approx\delta\lambda_0$. As in the estimate
in the proof of \eqref{FVendpointperturbative}, where a similar term appeared, we pair this
factor and the other low-frequency factor in $\mathcal{J}_B$ with
$S_{\lambda_0}(T_{\mathfrak m}\varphi)$ and
$\overline u_\mu^{x_0}$, respectively, using two unbalanced estimates. Next, since $P_{<\delta\lambda_0}+\widetilde S_{\lambda_0}
=P_{<\delta^{-1}\lambda_0}$, the remaining part of $\mathcal{J}_B$ and the contribution of the linear principal parts of \eqref{F21VBmainexpansions} in \eqref{F21secondterm} add up to give
\begin{equation*}
\begin{split}
-T_{\nabla\eta}\cdot
\int_0^\infty
S_{\lambda_0}(T_{\mathfrak m}\varphi)
P_{<\delta^{-1}\lambda_0}
T_{\mathfrak l}(\varphi-\varphi_{<\delta\lambda_0})
\frac{d\lambda_0}{\lambda_0}
+\mathcal R.
\end{split}
\end{equation*}
We can then write (provided $\delta$ is small enough)
\begin{equation*}
\begin{split}
P_{<\delta^{-1}\lambda_0}
T_{\mathfrak l}(\varphi-\varphi_{<\delta\lambda_0})=
\widetilde S_{\lambda_0}T_{\mathfrak l}\varphi
-P_{<\delta^{-1}\lambda_0}
[T_{\mathfrak l},P_{<\delta\lambda_0}]\varphi.
\end{split}
\end{equation*}
The first term gives the balanced quadratic contribution
\begin{equation*}
\begin{split}
-T_{\nabla\eta}\cdot
\int_0^\infty
S_{\lambda_0}(T_{\mathfrak m}\varphi)
\widetilde S_{\lambda_0}T_{\mathfrak l}\varphi
\frac{d\lambda_0}{\lambda_0}=-T_{\nabla\eta}\cdot
\Pi(T_{\mathfrak l}\varphi,T_{\mathfrak m}\varphi)=N_{\varphi,F}^{[2]}.
\end{split}
\end{equation*}
On the other hand, the commutator contributes
\begin{equation*}
\begin{split}
\mathcal C_B\coloneq
T_{\nabla\eta}\cdot
\int_0^\infty
S_{\lambda_0}(T_{\mathfrak m}\varphi)
P_{<\delta^{-1}\lambda_0}
[T_{\mathfrak l},P_{<\delta\lambda_0}]\varphi
\frac{d\lambda_0}{\lambda_0}.
\end{split}
\end{equation*}

This term is estimated in the same way as the commutator
$[T_{\mathfrak m},P_{<\delta\lambda_1}]\varphi$ in the proof
of \eqref{FVendpointperturbative} using two unbalanced bilinear estimates. We omit the details since they are similar, but the resulting estimate is
\begin{equation*}
\|(\mathcal C_B)_\mu\overline u_\mu^{x_0}\|_{L_{t,x}^1}
\lesssim_{M,\delta}
\epsilon^5c_\mu^2\mu^{-2s_c-\frac12}(1+\mu|x_0|)^{\frac12}.
\end{equation*}
We have therefore shown that
\begin{equation}\label{F21quadraticprincipal}
\begin{split}
\mathcal J_B
&-T_{\nabla\eta}\cdot
\int_0^\infty
S_{\lambda_0}(T_{\mathfrak m}\varphi)
\widetilde S_{\lambda_0}T_{\mathfrak l}
(\varphi-\varphi_{<\delta\lambda_0})
\frac{d\lambda_0}{\lambda_0}
\\
&=-T_{\nabla\eta}\cdot
\Pi(T_{\mathfrak l}\varphi,T_{\mathfrak m}\varphi)
+\mathcal R
=N_{\varphi,F}^{[2]}+\mathcal R.
\end{split}
\end{equation}

Before assembling the cubic terms, we also account for the
correction to the linear $B$ factor when it is multiplied by
the quadratic $V$ term in \eqref{F21VBmainexpansions}. This
correction contains
$\widetilde S_{\lambda_0}T_{\mathfrak l}
\varphi_{<\delta\lambda_0}$, which is localized at frequency
$\beta\approx\delta\lambda_0$, and is perturbative. To see this,
we use the $V$ analogue of the decomposition \eqref{QEdecomp}.
For the $L_{t,x}^2$ component, we pair
$\widetilde S_{\lambda_0}T_{\mathfrak l}
\varphi_{<\delta\lambda_0}$ with $\overline u_\mu^{x_0}$
and use the bound
\begin{equation*}
\epsilon^4c_\mu^2
\mu^{\frac14-2s_c}\beta^{-\frac34}
(1+\mu|x_0|)^{\frac12}
\lesssim_\delta
\epsilon^4c_\mu^2\mu^{-\frac12-2s_c}
(1+\mu|x_0|)^{\frac12}.
\end{equation*}
For the remaining quadratic component, we use the narrow
frequency decomposition from \eqref{QEdecomp}. We pair an
input separated from $\mu\approx\lambda_0$ with
$\overline u_\mu^{x_0}$, and the other input with
$\widetilde S_{\lambda_0}T_{\mathfrak l}
\varphi_{<\delta\lambda_0}$, using two unbalanced estimates.
The same argument justifies the principal replacement of
the $V$ factor in \eqref{FBendpoint}.  The product of the linear $B$ term with the quadratic $V$ term
combines with \eqref{peeledintegrand} to give the first
contribution to $N_{\varphi,F}^{[3]}$. Indeed, after the same
commutations and frequency rearrangements used in deriving
\eqref{Fhhh}, their respective output projections on the quadratic
$V$ factor are $\widetilde S_{\lambda_0}$ and
$P_{<\delta\lambda_0}$. Using
\begin{equation*}
P_{<\delta\lambda_0}+\widetilde S_{\lambda_0}
=P_{<\delta^{-1}\lambda_0},
\end{equation*}
and extending the integration regions to be symmetric, as in the argument in obtaining
\eqref{Fhhh}, we obtain
\begin{equation}\label{FVnestedassembly}
-T_{\nabla\eta}\cdot
\Pi\bigl(
T_{\mathfrak l}\varphi,
(\eta,T_{\mathfrak m}T_{\mathfrak l}\varphi)
-T_{\mathfrak m}(\eta,T_{\mathfrak l}\varphi)
\bigr)+\mathcal R.
\end{equation}
The product of the linear $V$-term with the quadratic $B$-term,
together with the corresponding contribution from $F^{2,2}$, gives
the other nested contribution,
\begin{equation}\label{FBnestedassembly}
-T_{\nabla\eta}\cdot
\Pi\bigl(
T_{\mathfrak m}\varphi,
(\eta,T_{\mathfrak l}^2\varphi)
-T_{\mathfrak l}(\eta,T_{\mathfrak l}\varphi)
\bigr)+\mathcal R.
\end{equation}
Finally, the
product of the two quadratic terms in \eqref{F21VBmainexpansions}, together with the terms left
aside in \eqref{F21VBmainexpansions}, contributes
$N_{\varphi,F}^{[4]}+\mathcal R$.
\\

Finally, the last term in \eqref{F21reduced} is perturbative. Indeed,
it contains
\begin{equation*}
(T_{V^l}-V^{\mathrm{lo}})\cdot\nabla\eta_{\lambda_0}
\widetilde S_{\lambda_0}(B-B_{<\delta\lambda_0}).
\end{equation*}
Expanding $T_{V^l}-V^{\mathrm{lo}}$ as in
\eqref{TVminusVlo}, and peeling off the subprincipal errors as in many of the estimates above, one finds that the remaining terms when we apply $P_{\mu}$ and pair with $\overline{U}_{\mu}^{x_0}$ can always (as seen several times above) be controlled using two bilinear estimates.

\subsubsection{Assembly of $N_{\varphi,F}^{[3]}$}
Having identified the quadratic term $N_{\varphi,F}^{[2]}$ in
\eqref{F21quadraticprincipal}, it remains to collect the cubic
terms into $N_{\varphi,F}^{[3]}$. All the paraproducts below are understood with the symmetric
convention in \eqref{NphinestedPiconvention}; thus the three
explicit input scales are restricted to $\mathcal D$. The
principal-symbol identities
\begin{equation}\label{glmidentities}
\mathfrak g_{<\lambda}
=\mathfrak l_{<\lambda}-\nabla\eta_{<\lambda}\cdot\mathfrak m_{<\lambda},
\qquad
\mathfrak m_{<\lambda}=i\xi-\nabla\eta_{<\lambda}\mathfrak l_{<\lambda}
\end{equation}
are used at the paradifferential level modulo paraproduct and
projection commutators, whose contributions belong to $\mathcal R$.
Moreover, we recall the identities from the expansion of $\frac{1}{2}B^2-T_{B^l}B$ and $\frac{1}{2}|V|^2-T_{V^l}V$ earlier 
\begin{equation*}
\begin{split}
&N_{\varphi,B}^{[3]}=\Pi\bigl(
T_{\mathfrak l}\varphi,
\eta,
T_{\mathfrak l}^2\varphi
\bigr)-\Pi\bigl(
T_{\mathfrak l}\varphi,
T_{\mathfrak l}(\eta,T_{\mathfrak l}\varphi)
\bigr),
\\
&N_{\varphi,V}^{[3]}
=\Pi\bigl(
T_{\mathfrak m}\varphi,
\eta,
T_{\mathfrak m}T_{\mathfrak l}\varphi
\bigr)-\Pi\bigl(
T_{\mathfrak m}\varphi,
T_{\mathfrak m}(\eta,T_{\mathfrak l}\varphi)
\bigr).
\end{split}
\end{equation*}
Equations \eqref{Fhhh}, \eqref{FVnestedassembly}, and
\eqref{FBnestedassembly} show that, modulo $\mathcal R$, the full
cubic contribution to $F$ is
\begin{equation*}
\begin{split}
&-\Pi\bigl(
\nabla\eta,
T_{\mathfrak l}\varphi,
T_{\mathfrak m}\varphi
\bigr)-T_{\nabla\eta}\cdot
\Bigl(
\Pi\bigl(
T_{\mathfrak l}\varphi,
\eta,
T_{\mathfrak m}T_{\mathfrak l}\varphi
\bigr)-\Pi\bigl(
T_{\mathfrak l}\varphi,
T_{\mathfrak m}(\eta,T_{\mathfrak l}\varphi)
\bigr)
\Bigr)
\\
&\quad-T_{\nabla\eta}\cdot
\Bigl(
\Pi\bigl(
T_{\mathfrak m}\varphi,
\eta,
T_{\mathfrak l}^2\varphi
\bigr)-\Pi\bigl(
T_{\mathfrak m}\varphi,
T_{\mathfrak l}(\eta,T_{\mathfrak l}\varphi)
\bigr)
\Bigr).
\end{split}
\end{equation*}
By the second identity in \eqref{glmidentities}, the first and third
terms combine to give
\begin{equation*}
\begin{split}
-\Pi\bigl(
\nabla\eta,
T_{\mathfrak l}\varphi,
T_{\mathfrak m}\varphi
\bigr)-T_{\nabla\eta}\cdot
\Bigl(
\Pi\bigl(
T_{\mathfrak m}\varphi,
\eta,
T_{\mathfrak l}^2\varphi
\bigr)-\Pi\bigl(
T_{\mathfrak m}\varphi,
T_{\mathfrak l}(\eta,T_{\mathfrak l}\varphi)
\bigr)
\Bigr)
=N_{\varphi,V}^{[3]}+\mathcal R.
\end{split}
\end{equation*}
The first identity in \eqref{glmidentities} similarly gives
\begin{equation*}
\begin{split}
-T_{\nabla\eta}\cdot
\Bigl(
\Pi\bigl(
T_{\mathfrak l}\varphi,
\eta,
T_{\mathfrak m}T_{\mathfrak l}\varphi
\bigr)-\Pi\bigl(
T_{\mathfrak l}\varphi,
T_{\mathfrak m}(\eta,T_{\mathfrak l}\varphi)
\bigr)
\Bigr)&=-N_{\varphi,B}^{[3]}
-\Pi\bigl(
T_{\mathfrak l}\varphi,
T_{\mathfrak g}(\eta,T_{\mathfrak l}\varphi)
\bigr)
\\
&\quad+\Pi\bigl(
T_{\mathfrak l}\varphi,
\eta,
T_{\mathfrak g}T_{\mathfrak l}\varphi
\bigr)+\mathcal R.
\end{split}
\end{equation*}
Therefore the full cubic contribution to $F$ is of the form
$N_{\varphi,F}^{[3]}+\mathcal R$, where
$N_{\varphi,F}^{[3]}$ is defined by \eqref{NphiF3def}.
Together with the quadratic and quartic contributions
identified above, this proves \eqref{Ffinalexpansion}
and Proposition \ref{cubicvarphi}.
\end{proof}

We finally record the three auxiliary estimates for the full expression
$F$.
\begin{proposition}\label{Fauxbounds}
Let $F$ be as in Proposition \ref{cubicvarphi}. Then the following
bounds hold:
\begin{enumerate}
\item For every $\lambda>0$, we have
\begin{equation}\label{FauxL2Linf}
\|F_\lambda\|_{L_t^2L_x^\infty}
\lesssim\epsilon^2\lambda^{-\frac54}.
\end{equation}

\item For every $0<\lambda\ll\mu$ and $x_0\in\mathbb R^2$,
we have a decomposition
\begin{equation}\label{Fauxdecomp}
F_\lambda\overline u_\mu^{x_0}u_\mu
=\mathcal R_F^{\mathrm{aux}}+\nabla\cdot Q_F^{\mathrm{aux}},
\end{equation}
where, for $\frac43<p\leq2$,
\begin{equation}\label{Fauxtrilinear}
\|\mathcal R_F^{\mathrm{aux}}\|_{L_{t,x}^1}
+\lambda^{\frac2p-1}\|Q_F^{\mathrm{aux}}\|_{L_t^1L_x^p}
\lesssim_p
\epsilon^4c_\mu^2\lambda^{-\frac32}
\mu^{-2s_c}(1+\mu|x_0|)^{\frac12}.
\end{equation}
\item For every $\lambda>0$, we have
\begin{equation}\label{FauxL2cubic}
\|(F-N_{\varphi,F}^{[2]})_\lambda\|_{L_{t,x}^2}
\lesssim_\delta\epsilon^3c_\lambda\lambda^{-\frac14-s_c}.
\end{equation}
\end{enumerate}
\end{proposition}

\begin{proof}
The proof is simpler than that of the previous proposition
and does not require the full expansion obtained there. It will be convenient to define the operators
\begin{equation*}
T_B^{\mathrm{err}}=T_B-T_{B^l},\qquad
T_V^{\mathrm{err}}=T_V-T_{V^l}
\end{equation*}
which encode the error from replacing the "geometric" para products $T_{B^l}$ and $T_{V^l}$ with the standard ones $T_B$ and $T_V$, respectively. Correspondingly, we define
\begin{equation*}
F_{\mathrm{main}}
=-BV\cdot\nabla\eta+T_B(V\cdot\nabla\eta)
+T_V\cdot(B\nabla\eta)-T_V\cdot T_B\nabla\eta.
\end{equation*}
Then $F=F_{\mathrm{main}}+F_{\mathrm{err}}$ with
\begin{equation}\label{Fauxspatialdecomp}
\begin{split}
F_{\mathrm{err}}=-T_B^{\mathrm{err}}(V\cdot\nabla\eta)
+T_{V^l}\cdot T_B^{\mathrm{err}}\nabla\eta
-T_V^{\mathrm{err}}\cdot(B\nabla\eta-T_B\nabla\eta).
\end{split}
\end{equation}
We now prove the first bound \eqref{FauxL2Linf}, beginning with the contribution of
$F_{\mathrm{main}}$. Expanding $T_B$ gives
\begin{equation*}
\begin{split}
F_{\mathrm{main}}
=\int_0^\infty
\Bigl[
T_V\cdot\left(
(B-P_{<M\delta\nu}B)\nabla\eta_\nu
\right)
-(B-P_{<M\delta\nu}B)(V\cdot\nabla\eta)_\nu
\Bigr]
\frac{d\nu}{\nu}.
\end{split}
\end{equation*}
We now decompose $B$ and $V$ into their
Littlewood-Paley pieces. Before estimating the resulting
integrals, we identify the cancellations between the four
terms defining $F_{\mathrm{main}}$.

First, any contribution in which one input has frequency
much higher than the other two vanishes. For example, if
$V$ has frequency much higher than both $B$ and $\eta$,
then, for these frequency components,
\begin{equation*}
T_B(V\cdot\nabla\eta)=BV\cdot\nabla\eta,\qquad
T_V\cdot(B\nabla\eta)=T_V\cdot T_B\nabla\eta=0.
\end{equation*}
Thus the first and second terms cancel. If $B$ has the
highest frequency, the first and third terms cancel.
If $\eta$ has the highest frequency, both paraproducts
agree with multiplication and all four terms cancel.
Consequently, two of the largest input frequencies must
be comparable (with implicit constants depending on $M$ and $\delta$).

There is one further cancellation needed when one of
the $B,V$ frequencies is much smaller than the output
frequency $\lambda$. Consider first the contribution
of $V_{\lambda_1}$ with $\lambda_1\ll\lambda$, keeping
$B$ and $\eta$ unexpanded for this calculation.
After projection to frequency $\lambda$,
$T_{V_{\lambda_1}}$ agrees with multiplication by
$V_{\lambda_1}$. Therefore,
\begin{equation*}
P_\lambda\left(
-BV_{\lambda_1}\cdot\nabla\eta
+T_{V_{\lambda_1}}\cdot(B\nabla\eta)
\right)=0.
\end{equation*}
The other two terms give a commutator structure
\begin{equation*}
\begin{split}
P_\lambda\left(
T_B(V_{\lambda_1}\cdot\nabla\eta)
-T_{V_{\lambda_1}}\cdot T_B\nabla\eta
\right)&=
P_\lambda\left(
T_B(V_{\lambda_1}\cdot\nabla\eta)
-V_{\lambda_1}\cdot T_B\nabla\eta
\right)
\\
&=P_\lambda[T_B,V_{\lambda_1}]\cdot\nabla\eta.
\end{split}
\end{equation*}
Using the definition of $T_B$, this commutator is
\begin{equation*}
P_\lambda[T_B,V_{\lambda_1}]\cdot\nabla\eta
=
P_\lambda\int_{\nu\approx\lambda}
(P_{<M\delta\nu}B)
[S_\nu,V_{\lambda_1}]\cdot\nabla\eta
\frac{d\nu}{\nu}.
\end{equation*}
The restriction $\nu\approx\lambda$ follows from the
output support and $\lambda_1\ll\lambda$.
By \eqref{bilformrep}, we then hae
\begin{equation*}
[S_\nu,V_{\lambda_1}]\cdot\nabla\eta
=\frac{\lambda_1}{\lambda}\mathcal L
(V_{\lambda_1},\nabla\eta).
\end{equation*}
If instead $B_{\lambda_1}$ is the low-frequency factor,
with $\lambda_1\ll\lambda$, the contribution vanishes.
Indeed, $T_{B_{\lambda_1}}$ agrees with multiplication
on the $\eta$ input and on the frequency $\lambda$
output of $V\cdot\nabla\eta$, so the first two and last
two terms cancel in pairs. 
\\

We can now write the resulting frequency integrals.
We use $\lambda_0,\lambda_1$ for the frequencies of the
$B,V$ factors. When these frequencies are comparable,
the contributions have the form
\begin{equation}\label{Fauxprincipalterms1}
\int_{\lambda\lesssim\lambda_0}
\int_{\lambda_1\approx\lambda_0}
P_\lambda\mathcal L
\left(B_{\lambda_0},
V_{\lambda_1},
\nabla\eta_{<\delta^{-1}\lambda_0}\right)
\frac{d\lambda_1}{\lambda_1}
\frac{d\lambda_0}{\lambda_0}.
\end{equation}
When $\eta$ and one of $B,V$ have comparable frequencies,
the contributions have the form
\begin{equation}\label{Fauxprincipalterms2}
\int_{\lambda\lesssim\lambda_0}
\int_{\lambda_1\lesssim\lambda_0}
\min\left\{1,\frac{\lambda_1}{\lambda}\right\}
P_\lambda\mathcal L
\left(B_{\lambda_0},
V_{\lambda_1},
\nabla\widetilde P_{\lambda_0}\eta\right)
\frac{d\lambda_1}{\lambda_1}
\frac{d\lambda_0}{\lambda_0},
\end{equation}
with $B,V$ interchanged as needed.
Here $\widetilde P_{\lambda_0}$ is a fixed fattening
of the annulus at frequency $\lambda_0$. The integration
over the frequencies of $\eta$ comparable to $\lambda_0$
is included in $\mathcal L$. These forms involve only
smooth Littlewood-Paley projections and have rapidly decaying
kernels at the corresponding frequency scales. Using Proposition \ref{VBunbil} with $s=s_c=2$, we have the bounds
\begin{equation*}
\|B_{\lambda_j}\|_{L_t^4L_x^\infty}
+\|V_{\lambda_j}\|_{L_t^4L_x^\infty}
\lesssim\epsilon\lambda_j^{-\frac58}.
\end{equation*}
Hence, putting each of the two $B,V$ factors in $L_t^4L_x^\infty$,
the contributions of \eqref{Fauxprincipalterms1} and
\eqref{Fauxprincipalterms2} are bounded in $L_t^2L_x^{\infty}$ by
\begin{equation*}
\begin{split}
&\epsilon^2\int_{\lambda_0\gtrsim\lambda}
\lambda_0^{-\frac54}\frac{d\lambda_0}{\lambda_0}
\lesssim\epsilon^2\lambda^{-\frac54},
\\
&\epsilon^2\int_{\lambda_0\gtrsim\lambda}
\int_{\lambda_1\lesssim\lambda_0}
\min\left\{1,\frac{\lambda_1}{\lambda}\right\}
\lambda_0^{-\frac58}\lambda_1^{-\frac58}
\frac{d\lambda_1}{\lambda_1}\frac{d\lambda_0}{\lambda_0}
\lesssim\epsilon^2\lambda^{-\frac54}.
\end{split}
\end{equation*}
Thus, we find
\begin{equation*}
\|(F_{\mathrm{main}})_\lambda\|_{L_t^2L_x^\infty}
\lesssim\epsilon^2\lambda^{-\frac54}.
\end{equation*}
We next estimate the contribution of $F_{\mathrm{err}}$. We first group
the terms containing $T_B^{\mathrm{err}}$ as
\begin{equation}\label{Fauxerrcomm}
\begin{split}
&-T_B^{\mathrm{err}}(V\cdot\nabla\eta)
+T_{V^l}\cdot T_B^{\mathrm{err}}\nabla\eta
\\
&=[T_V,T_B^{\mathrm{err}}]\cdot\nabla\eta
-T_B^{\mathrm{err}}(V\cdot\nabla\eta-T_V\cdot\nabla\eta)
-T_V^{\mathrm{err}}\cdot T_B^{\mathrm{err}}\nabla\eta.
\end{split}
\end{equation}
The coefficients of $T_B^{\mathrm{err}}$ and
$T_V^{\mathrm{err}}$ can be expanded as in
\eqref{VtailFTCstart}, \eqref{survivingprincipal}, and
\eqref{F23bexp}. For each $\lambda_0>0$, this gives
\begin{equation}\label{Fauxcoefficientbounds}
\begin{split}
P_{<M\delta\lambda_0}(B-B_{<\delta\lambda_0})
&=P_{<M\delta\lambda_0}T_{\mathfrak l}
(\varphi-\varphi_{<\delta\lambda_0})
+O_{L_t^2L_x^\infty}(\epsilon^2\lambda_0^{-\frac34}),
\\
P_{<M\delta\lambda_0}(V-V_{<\delta\lambda_0})
&=P_{<M\delta\lambda_0}T_{\mathfrak m}
(\varphi-\varphi_{<\delta\lambda_0})
+O_{L_t^2L_x^\infty}(\epsilon^2\lambda_0^{-\frac34}).
\end{split}
\end{equation}
The principal linear terms have frequencies comparable
to $\delta\lambda_0$, up to constants depending on $M$,
and satisfy the $L_t^4L_x^\infty$ bound
$O(\epsilon\lambda_0^{-5/8})$.
The errors are also
$O_{L_{t,x}^\infty}(\epsilon\lambda_0^{-1/2})$,
so interpolation gives
$O_{L_t^4L_x^\infty}(\epsilon^{3/2}\lambda_0^{-5/8})$.
Consequently, the full coefficients satisfy
\begin{equation*}
\begin{split}
&\|P_{<M\delta\lambda_0}(B-B_{<\delta\lambda_0})
\|_{L_t^4L_x^\infty}
+\|P_{<M\delta\lambda_0}(V-V_{<\delta\lambda_0})
\|_{L_t^4L_x^\infty}
\lesssim\epsilon\lambda_0^{-\frac58}.
\end{split}
\end{equation*}

We also use the corresponding bounds after differentiating
in the frequency parameter. For example,
\begin{equation*}
\begin{split}
\lambda_0\frac{d}{d\lambda_0}
\left[P_{<M\delta\lambda_0}(B-B_{<\delta\lambda_0})\right]=
(\lambda_0\frac{d}{d\lambda_0}P_{<M\delta\lambda_0})
(B-B_{<\delta\lambda_0})
-P_{<M\delta\lambda_0}B^{\mathrm{lin}}_{\delta\lambda_0}.
\end{split}
\end{equation*}
The bounds for $B_\nu^{\mathrm{lin}}$ used in the proof
of \eqref{BparaRbound1}, together with FTC, give the same
$L_t^4L_x^\infty$ estimate for this expression.
The argument for $V$ is identical.

Using Proposition \ref{besovparabounds} and Proposition \ref{VBunbil}, we also have
\begin{equation*}
\begin{split}
&\|(B\nabla\eta-T_B\nabla\eta)_\lambda\|_{L_t^4L_x^\infty}
+\|(V\cdot\nabla\eta-T_V\cdot\nabla\eta)_\lambda
\|_{L_t^4L_x^\infty}
\lesssim\epsilon^2\lambda^{-\frac58}.
\end{split}
\end{equation*}
Combining the above estimates, we can control the last term in \eqref{Fauxspatialdecomp}
and the second term on the right-hand side of
\eqref{Fauxerrcomm} by
\begin{equation*}
O_{L_t^2L_x^{\infty}}(\epsilon\lambda^{-\frac58}
\epsilon^2\lambda^{-\frac58})
=O_{L_t^2L_x^{\infty}}(\epsilon^3\lambda^{-\frac54}).
\end{equation*}
For the last term in \eqref{Fauxerrcomm}, we use two
$L_t^4L_x^\infty$ bounds and
$\|\nabla\eta\|_{L_{t,x}^\infty}\lesssim\epsilon$,
giving the same estimate. It remains to estimate $[T_V,T_B^{\mathrm{err}}]\cdot\nabla\eta$.
We study its contribution $V_{\lambda_1}$ for each $\lambda_1>0$. First, the projection $P_{\lambda}$ restricts $V_{\lambda_1}$ to $\lambda_1\lesssim\lambda$.
When $\lambda_1\ll\lambda$, we can freely replace $T_{V_{\lambda_1}}$ by $V_{\lambda_1}$, so that
\begin{equation*}
\begin{split}
P_\lambda[T_{V_{\lambda_1}},T_B^{\mathrm{err}}]
\cdot\nabla\eta
&=-P_\lambda[T_B^{\mathrm{err}},V_{\lambda_1}]
\cdot\nabla\eta
\\
&=
-P_\lambda\int_{\lambda_0\approx\lambda}
P_{<M\delta\lambda_0}(B-B_{<\delta\lambda_0})
[S_{\lambda_0},V_{\lambda_1}]\cdot\nabla\eta
\frac{d\lambda_0}{\lambda_0}.
\end{split}
\end{equation*}
Then by Proposition \ref{Leibnizcommutator}, we have
\begin{equation*}
[S_{\lambda_0},V_{\lambda_1}]\cdot\nabla\eta
=\lambda_0^{-1}\mathcal L
(\nabla V_{\lambda_1},\nabla\eta).
\end{equation*}
Here, we put the coefficient
$P_{<M\delta\lambda_0}(B-B_{<\delta\lambda_0})$ and
$V_{\lambda_1}$ in $L_t^4L_x^\infty$ using the $L_t^4L_x^{\infty}$ bound following \eqref{Fauxcoefficientbounds} as well as Proposition \ref{VBunbil}, respectively. We place the factor
$\nabla\eta$ in $L_{t,x}^\infty$.
For $\lambda_1\approx\lambda$, the same bound follows
directly from the two $L_t^4L_x^\infty$ 
bounds, by directly expanding the commutator. In summary, combining these two estimates and Bernstein gives
\begin{equation*}
\begin{split}
\|P_\lambda[T_V,T_B^{\mathrm{err}}]\cdot\nabla\eta
\|_{L_t^2L_x^\infty}\lesssim
\epsilon^3\lambda^{-\frac58}
\int_{\lambda_1\lesssim\lambda}
\frac{\lambda_1}{\lambda}\lambda_1^{-\frac58}
\frac{d\lambda_1}{\lambda_1}\lesssim\epsilon^3\lambda^{-\frac54}.
\end{split}
\end{equation*}
Thus
\begin{equation*}
\|(F_{\mathrm{err}})_\lambda\|_{L_t^2L_x^\infty}
\lesssim\epsilon^3\lambda^{-\frac54}.
\end{equation*}
Combining the estimates for $F_{\mathrm{main}}$ and
$F_{\mathrm{err}}$ proves \eqref{FauxL2Linf}. For the remaining two auxiliary bounds, we will
also need the principal error bounds
\begin{equation}\label{Fauxfactorbounds}
\begin{split}
B_\alpha=T_{\mathfrak l}\varphi_\alpha
+O_{L_t^2L_x^\infty}(\epsilon^2\alpha^{-\frac34}),\hspace{5mm}V_\alpha=T_{\mathfrak m}\varphi_\alpha
+O_{L_t^2L_x^\infty}(\epsilon^2\alpha^{-\frac34})
\end{split}
\end{equation}
which follow from Corollary \ref{principalVBexp} and
Proposition \ref{Leibnizcommutator}.
The errors have frequency comparable to $\alpha$ and
are also bounded by
$O_{L_{t,x}^\infty}(\epsilon^2\alpha^{-1/2})$
and
$O_{L_t^\infty L_x^2}(\epsilon^2\alpha^{-\frac32})$,
using Proposition \ref{VBunbil} and the bootstrap bounds
for $\varphi$.
\\

We next prove \eqref{FauxL2cubic}, first estimating
$F_{\mathrm{err}}$ directly in $L_{t,x}^2$. 
Substituting \eqref{Fauxerrcomm} into
\eqref{Fauxspatialdecomp}, we have
\begin{equation}\label{Ferrwrittenout}
\begin{split}
F_{\mathrm{err}}&=[T_V,T_B^{\mathrm{err}}]\cdot\nabla\eta-T_B^{\mathrm{err}}
(V\cdot\nabla\eta-T_V\cdot\nabla\eta)
\\
&-T_V^{\mathrm{err}}\cdot T_B^{\mathrm{err}}\nabla\eta-T_V^{\mathrm{err}}\cdot
(B\nabla\eta-T_B\nabla\eta).
\end{split}
\end{equation}
We write $\kappa$ for the frequency of the $\nabla\eta$ factor.
For $\kappa\gtrsim_\delta\lambda$, we argue as in the
proof of \eqref{FauxL2Linf}, placing $\nabla\eta_{\kappa}$
in $L_t^\infty L_x^2$ and using two $L_t^4L_x^\infty$
bounds from Proposition \ref{VBunbil} and the estimates
following \eqref{Fauxcoefficientbounds}.
This applies to the commutator in \eqref{Ferrwrittenout}
and $T_V^{\mathrm{err}}\cdot T_B^{\mathrm{err}}\nabla\eta$
(these terms vanish when $\kappa\ll\lambda$ after applying
$P_{\lambda}$), and also to the remaining two terms for
this range of $\kappa$. Upon integration over the relevant frequency parameters,
these terms satisfy the required
$O_{L_{t,x}^2}(\epsilon^3c_{\lambda}
\lambda^{-\frac14-s_c})$ bound.

It remains to estimate
$T_B^{\mathrm{err}}(V\cdot\nabla\eta-T_V\cdot\nabla\eta)$
and
$T_V^{\mathrm{err}}\cdot(B\nabla\eta-T_B\nabla\eta)$
when $\kappa\ll_\delta\lambda$. We sketch the details
for the first term, as the latter is virtually identical.
In this case, $V$ must have frequency comparable to
$\lambda$ after applying $P_{\lambda}$.
By writing out the definitions of $T_B^{\mathrm{err}}$ and $T_V$, and using \eqref{Fauxcoefficientbounds} and
\eqref{Fauxfactorbounds}, it is easy to see that we have
\begin{equation*}
\begin{split}
P_\lambda T_B^{\mathrm{err}}
(V\cdot\nabla\eta_\kappa-T_V\cdot\nabla\eta_\kappa)&=
P_\lambda\int_{\lambda_0\approx\lambda}
\Bigl[
P_{<M\delta\lambda_0}T_{\mathfrak l}
(\varphi-\varphi_{<\delta\lambda_0})
\Bigr]
S_{\lambda_0}\bigl(
T_{\mathfrak m}\varphi\cdot\nabla\eta_\kappa
\bigr)
\frac{d\lambda_0}{\lambda_0}
\\
&+O_{L_{t,x}^2}\bigl(
\epsilon^4c_\lambda\lambda^{-\frac14-s_c}
\bigr).
\end{split}
\end{equation*}
Here $T_V\cdot\nabla\eta_\kappa$ has frequency
comparable to $\kappa$, so its contribution vanishes
after applying $P_\lambda T_B^{\mathrm{err}}$.
The errors in replacing $B$ and $V$ by their principal parts above are estimated by placing the
error from \eqref{Fauxcoefficientbounds} or
\eqref{Fauxfactorbounds} in $L_t^2L_x^\infty$,
the other factor in $L_t^\infty L_x^2$, and
$\nabla\eta$ in $L_{t,x}^\infty$. In the remaining integral above, the coefficient
$P_{<M\delta\lambda_0}T_{\mathfrak l}
(\varphi-\varphi_{<\delta\lambda_0})$
has frequency $\beta\approx_M\delta\lambda$, while
the $T_{\mathfrak m}\varphi$ factor has frequency
comparable to $\lambda$. One unbalanced bilinear
estimate, with the $\nabla\eta$ factor in $L_{t,x}^\infty$, bounds
this contribution by
\begin{equation*}
\epsilon^3c_\lambda\lambda^{\frac12-s_c}
\int_{\beta\approx_M\delta\lambda}
\beta^{-\frac34}\frac{d\beta}{\beta}
\lesssim_\delta
\epsilon^3c_\lambda\lambda^{-\frac14-s_c}.
\end{equation*}
The other term $T_V^{\mathrm{err}}(B\nabla\eta-T_B\nabla\eta)$ is estimated identically, with
$\mathfrak l$ and $\mathfrak m$ interchanged.
Combining the estimates for all four terms in
\eqref{Ferrwrittenout}, we obtain
\begin{equation*}
\|(F_{\mathrm{err}})_\lambda\|_{L_{t,x}^2}
\lesssim_\delta
\epsilon^3c_\lambda\lambda^{-\frac14-s_c}.
\end{equation*}
We now consider $F_{\mathrm{main}}$ and its decomposition
into \eqref{Fauxprincipalterms1} and
\eqref{Fauxprincipalterms2}. We can estimate the
contribution of \eqref{Fauxprincipalterms2} directly in
$L_{t,x}^2$ by placing
$\nabla\widetilde P_{\lambda_0}\eta$ in
$L_t^\infty L_x^2$ and $V_{\lambda_1}$ and
$B_{\lambda_0}$ in $L_t^4L_x^\infty$. By Proposition
\ref{VBunbil} and the bootstrap bounds, integrating
in $\lambda_0$ and $\lambda_1$ gives
$O_{L_{t,x}^2}(\epsilon^3c_\lambda
\lambda^{-\frac14-s_c})$, as required.
For \eqref{Fauxprincipalterms1}, we compare its principal
part with $N_{\varphi,F}^{[2]}$ using the algebra leading
to \eqref{F21quadraticprincipal}. The errors are of the
following types:
\begin{itemize}
\item Expressions where the $\nabla\eta$ frequency $\kappa$ satisfies
$\lambda\lesssim_\delta\kappa
\lesssim_\delta\lambda_0$.
We place $V_{\lambda_1}$ and $B_{\lambda_0}$ in
$L_t^4L_x^\infty$ and $\nabla\eta_\kappa$ in
$L_t^\infty L_x^2$. Since
$\lambda_1\approx\lambda_0$, the resulting bound is
$\epsilon^3c_\kappa\kappa^{1-s_c}\lambda_0^{-\frac54}$.
Integrating first in $\lambda_0,\lambda_1$ and then
in $\kappa$ gives
$O_{L_{t,x}^2}(\epsilon^3c_\lambda
\lambda^{-\frac14-s_c})$. The same estimate applies to these $\nabla\eta$
frequencies in $(N_{\varphi,F}^{[2]})_\lambda$.
\item Expressions where one or both velocity trace factors
are replaced by
$B_{\lambda_0}-T_{\mathfrak{l}}\varphi_{\lambda_0}$
or
$V_{\lambda_1}-T_{\mathfrak{m}}\varphi_{\lambda_1}$.
We use \eqref{Fauxfactorbounds} to estimate one error
in $L_t^2L_x^\infty$, and Proposition \ref{VBunbil}
and the bootstrap bounds to estimate the other factor
in $L_t^\infty L_x^2$, with $\nabla\eta$ placed in
$L_{t,x}^\infty$. If the other factor is also an error, it can also be placed in $L_t^{\infty}L_x^2$ by the triangle inequality and estimating each part separately.
The commutators needed to write the principal factors
as $S_{\lambda_j}(T_{\mathfrak l}\varphi)$ and
$S_{\lambda_j}(T_{\mathfrak m}\varphi)$ satisfy the
same bounds by Proposition \ref{Leibnizcommutator},
as in the proof of \eqref{BauxL2cubic}.
These contributions are bounded by
\begin{equation*}
\epsilon^4\int_{\lambda_0\gtrsim\lambda}
c_{\lambda_0}\lambda_0^{-\frac14-s_c}
\frac{d\lambda_0}{\lambda_0}
\lesssim_\delta
\epsilon^4c_\lambda\lambda^{-\frac14-s_c}.
\end{equation*}

\item Commutators with Littlewood-Paley projections. For $\kappa\ll_\delta\lambda$, commuting the
$\nabla\eta_{\kappa}$ factor through the projections gives a factor
$\kappa/\lambda$ by Proposition \ref{Leibnizcommutator}. In this case, the two velocity traces in
\eqref{Fauxprincipalterms1} can be estimated by putting one velocity factor
in $L_t^\infty L_x^2$, and the other velocity
factor and the differentiated $\nabla\eta$ factor from the commutator in
$L_t^4L_x^\infty$. Integrating in $\kappa$ and $\lambda_0$ gives
\begin{equation*}
\epsilon^3\lambda^{-\frac18}
\int_{\lambda_0\gtrsim\lambda}
c_{\lambda_0}\lambda_0^{-\frac18-s_c}
\frac{d\lambda_0}{\lambda_0}
\lesssim_\delta
\epsilon^3c_\lambda\lambda^{-\frac14-s_c}.
\end{equation*}
\end{itemize}
Adding up the above contributions, we obtain
\begin{equation*}
\|(F_{\mathrm{main}}-N_{\varphi,F}^{[2]})_\lambda
\|_{L_{t,x}^2}
\lesssim_\delta
\epsilon^3c_\lambda\lambda^{-\frac14-s_c}.
\end{equation*}
Together with the estimate for $F_{\mathrm{err}}$,
this proves \eqref{FauxL2cubic}.
\\

We now prove \eqref{Fauxtrilinear}. Again, we set
$H_\mu=\overline u_\mu^{x_0}u_\mu$. Like the $L_t^2L_x^{\infty}$ bound above, we will first consider the contribution of $F_{\mathrm{main}}$. As above, this can be written as a sum of terms of the form \eqref{Fauxprincipalterms1} or \eqref{Fauxprincipalterms2}.  We first use \eqref{Fauxfactorbounds} and \eqref{Fauxcoefficientbounds} to replace the $B,V$
factors by their principal parts $T_{\mathfrak{l}}\varphi$ and $T_{\mathfrak{m}}\varphi$ and estimate the resulting
errors. For $\lambda_0\ll\mu$, these errors contribute
directly to the remainder $\mathcal{R}_F^{\mathrm{aux}}$ while for $\lambda_0\gtrsim\mu$,
we estimate them in $L_{t,x}^2$ and include them in the
flux decomposition.We only sketch the necessary steps since they are similar to much of the analysis carried out earlier.  
\\

We first consider the case $\lambda_0\ll\mu$. The part of the replacement that is quadratic in the subprincipal errors can be estimated in $L_t^1L_x^{\infty}$ using two $L_t^2L_x^{\infty}$ error bounds from \eqref{Fauxfactorbounds} and \eqref{Fauxcoefficientbounds}. $u_{\mu}$ and $\overline{u}_{\mu}$ can then both be estimated in $L_t^{\infty}L_x^2$. Hence, this contribution belongs to $\mathcal{R}^{\mathrm{aux}}_F$. For the remaining part that consists of terms contatining one subprincipal error factor and one linear principal factor of either $T_{\mathfrak{l}}\varphi$ or $T_{\mathfrak{m}}\varphi$, we put the error in $L_t^2L_x^\infty$, and pair the remaining $T_{\mathfrak l}\varphi$ or
$T_{\mathfrak m}\varphi$ factor with one of $u_{\mu}$ or $\overline{u}_{\mu}^{x_0}$ using one unbalanced bilinear estimate. We place the other (unpaired) factor in $L_t^{\infty}L_x^2$. In each case, we always place $\nabla\eta$ in $L_{t,x}^{\infty}$. Inserting these errors into \eqref{Fauxprincipalterms1} and \eqref{Fauxprincipalterms2} shows that that the $\lambda_0\ll\mu$ contribution of $F_{\mathrm{main}}$ when multiplied by $H_{\mu}$ produces a term of the form $\mathcal{R}_F^{\mathrm{aux}}$ satisfying the $L_{t,x}^1$ bound \eqref{Fauxtrilinear}.
\\

For $\lambda_0\gtrsim\mu$, we directly estimate the replacement error in $L_{t,x}^2$.
If the two factors coming from $B$ and $V$ have comparable
frequencies, we put the error in $L_t^2L_x^\infty$ and the
other factor in $L_t^\infty L_x^2$, giving
$\epsilon^3\lambda_0^{-9/4}$.
If their frequencies are $\lambda_0$ and
$\lambda_1\ll\lambda_0$, then $\nabla\eta$ has frequency
comparable to $\lambda_0$ in view of \eqref{Fauxprincipalterms1} and \eqref{Fauxprincipalterms2}.
When the error has frequency $\lambda_0$, we use its
$L_t^2L_x^\infty$ bound, the $L_t^{\infty}L_x^2$ bound 
$\nabla\widetilde P_{\lambda_0}\eta$, and the
$L_{t,x}^\infty$ bound for the factor at frequency
$\lambda_1$. When the error has frequency $\lambda_1$,
we instead use the $L_t^{\infty}L_x^2$ for the factor at frequency
$\lambda_0$. The subprincipal errors when inserted into the terms \eqref{Fauxprincipalterms1} and \eqref{Fauxprincipalterms2} are then bounded directly in $L_{t,x}^2$ by
\begin{equation*}
\begin{split}
\epsilon^3\int_{\lambda_0\gtrsim\mu}\lambda_0^{-\frac94}
\frac{d\lambda_0}{\lambda_0}
&+\epsilon^3\int_{\lambda_0\gtrsim\mu}
\int_{\lambda_1\lesssim\lambda_0}
\min\left\{1,\frac{\lambda_1}{\lambda}\right\}
\left(
\lambda_0^{-\frac74}\lambda_1^{-\frac12}
+\lambda_0^{-\frac32}\lambda_1^{-\frac34}
\right)
\frac{d\lambda_1}{\lambda_1}\frac{d\lambda_0}{\lambda_0}
\\
&\lesssim\epsilon^3\left(
\mu^{-\frac94}
+\lambda^{-\frac12}\mu^{-\frac74}
+\lambda^{-\frac34}\mu^{-\frac32}
\right)
\lesssim\epsilon^3\lambda^{-1}\mu^{-\frac54}
\end{split}
\end{equation*}
We now estimate the main principal contribution in \eqref{Fauxprincipalterms1} and \eqref{Fauxprincipalterms2}. Again, we first
consider $\lambda_0\ll\mu$.
For \eqref{Fauxprincipalterms1}, we use the two unbalanced
pairings
\begin{equation*}
(\overline u_\mu^{x_0},T_{\mathfrak l}\varphi_{\lambda_0}),
\qquad
(u_\mu,T_{\mathfrak m}\varphi_{\lambda_1}),
\end{equation*}
and estimate $\nabla\eta_{<\delta^{-1}\lambda_0}$ in
$L_{t,x}^\infty$. Denoting this contribution by
$F_\lambda^{\mathrm{lo},1}$, we then obtain
\begin{equation*}
\begin{split}
\|F_\lambda^{\mathrm{lo},1}H_\mu\|_{L_{t,x}^1}
&\lesssim\epsilon^4c_\mu^2\mu^{-2s_c}(1+\mu|x_0|)^{\frac12}
\int_{\lambda\lesssim\lambda_0\ll\mu}
\lambda_0^{-\frac32}\left(1+\frac{\lambda_0}{\lambda}\right)
\frac{d\lambda_0}{\lambda_0}
\\
&\lesssim\epsilon^4c_\mu^2\lambda^{-\frac32}
\mu^{-2s_c}(1+\mu|x_0|)^{\frac12}.
\end{split}
\end{equation*}
Here we used \eqref{phiu-bi-unbal} and the kernel bounds for
$P_\lambda$, $T_{\mathfrak l}$, and $T_{\mathfrak m}$.
Since $\lambda_0\approx\lambda_1$, each of the two bilinear pairings costs a factor of
$\left(1+\lambda_0/\lambda\right)^{\frac12}$ for the translations coming from the output kernel. 

For \eqref{Fauxprincipalterms2}, it remains to consider
$\lambda_1\ll\lambda_0$. Here, we use one unbalanced pairing
\begin{equation*}
(u_\mu,T_{\mathfrak m}\varphi_{\lambda_1}),
\end{equation*}
place $T_{\mathfrak l}\varphi_{\lambda_0}$ and
$\nabla\widetilde P_{\lambda_0}\eta$ in $L_t^4L_x^\infty$,
and estimate $\overline u_\mu^{x_0}$ in $L_t^\infty L_x^2$.
Denoting this contribution by $F_\lambda^{\mathrm{lo},2}$,
these bounds and the kernel bounds for $P_{\lambda}$ give
\begin{equation*}
\begin{split}
\|F_\lambda^{\mathrm{lo},2}H_\mu\|_{L_{t,x}^1}
&\lesssim\epsilon^4c_\mu^2\mu^{-2s_c}(1+\mu|x_0|)^{\frac12}
\int_{\lambda\lesssim\lambda_0\ll\mu}
\int_{\lambda_1\lesssim\lambda_0}
\min\left\{1,\frac{\lambda_1}{\lambda}\right\}
\lambda_0^{-\frac34}\lambda_1^{-\frac34}
\left(1+\frac{\lambda_1}{\lambda}\right)^{\frac12}
\frac{d\lambda_1}{\lambda_1}\frac{d\lambda_0}{\lambda_0}
\\
&\lesssim\epsilon^4c_\mu^2\mu^{-2s_c}(1+\mu|x_0|)^{\frac12}
\lambda^{-1}
\int_{\lambda\lesssim\lambda_0\ll\mu}
\lambda_0^{-\frac12}\frac{d\lambda_0}{\lambda_0}
\\
&\lesssim\epsilon^4c_\mu^2\lambda^{-\frac32}
\mu^{-2s_c}(1+\mu|x_0|)^{\frac12}.
\end{split}
\end{equation*}
Next, we consider the principal contributions in
\eqref{Fauxprincipalterms1} and \eqref{Fauxprincipalterms2}
when $\lambda_0\gtrsim\mu$. Here, we estimate this term directly in $L_{t,x}^2$. In \eqref{Fauxprincipalterms1}, if $\nabla\eta$
and all the paracoefficients have frequencies much smaller
than $\lambda$, Lemma \ref{balancedrestrictedoutput} gives
\begin{equation*}
\left\|P_\lambda\mathcal L
\left(T_{\mathfrak l}\varphi_{\lambda_0},
T_{\mathfrak m}\varphi_{\lambda_1},\nabla\eta_{<\delta\lambda}
\right)\right\|_{L_{t,x}^2}
\lesssim\epsilon^3\lambda^{-\frac12}\lambda_0^{-\frac74}.
\end{equation*}
The contribution of $\nabla\eta_\kappa$ for
$\lambda\lesssim\kappa\lesssim\lambda_0$ satisfies the same bound. Indeed, integrating and using H\"older and Bernstien gives
\begin{equation*}
\begin{split}
\int_{\lambda\lesssim\kappa\lesssim\lambda_0}\left\|P_\lambda\mathcal L
\left(T_{\mathfrak l}\varphi_{\lambda_0},
T_{\mathfrak m}\varphi_{\lambda_1},\nabla\eta_\kappa
\right)\right\|_{L_{t,x}^2}\frac{d\kappa}{\kappa}&\lesssim
\int_{\lambda\lesssim\kappa\lesssim\lambda_0}\|T_{\mathfrak l}\varphi_{\lambda_0}\|_{L_t^\infty L_x^2}
\|T_{\mathfrak m}\varphi_{\lambda_1}\|_{L_t^4L_x^\infty}
\|\nabla\eta_\kappa\|_{L_t^4L_x^\infty}\frac{d\kappa}{\kappa}
\\
&\lesssim\epsilon^3\lambda^{-\frac{1}{8}}\lambda_0^{-\frac{17}{8}}
\\
&\lesssim \epsilon^2\lambda^{-\frac{1}{2}}\lambda_0^{-\frac{7}{4}}.
\end{split}
\end{equation*}
For \eqref{Fauxprincipalterms2}, we use instead
\begin{equation*}
\begin{split}
\|\min\left\{1,\frac{\lambda_1}{\lambda}\right\}
P_\lambda\mathcal L
\left(B_{\lambda_0},
V_{\lambda_1},
\nabla\widetilde P_{\lambda_0}\eta\right)\|_{L_{t,x}^2}&\lesssim\min\left\{1,\frac{\lambda_1}{\lambda}\right\}
\|T_{\mathfrak l}\varphi_{\lambda_0}\|_{L_t^\infty L_x^2}
\|T_{\mathfrak m}\varphi_{\lambda_1}\|_{L_t^4L_x^\infty}
\|\nabla\widetilde P_{\lambda_0}\eta\|_{L_t^4L_x^\infty}
\\
&\lesssim
\epsilon^3\min\left\{1,\frac{\lambda_1}{\lambda}\right\}
\lambda_1^{-\frac58}\lambda_0^{-\frac{13}{8}}.
\end{split}
\end{equation*}
Integrating the above bounds in $\lambda_0$ and $\lambda_1$ over the respective regions in \eqref{Fauxprincipalterms1} and \eqref{Fauxprincipalterms2}, we can estimate \eqref{Fauxprincipalterms1} and \eqref{Fauxprincipalterms2} directly in $L_{t,x}^2$ by (for instance) \begin{equation}\label{FauxprincipalhighL2}
\epsilon^2\left(
\lambda^{-\frac12}\mu^{-\frac74}+\lambda^{-\frac58}\mu^{-\frac{13}{8}}
\right)
\lesssim\epsilon^2\lambda^{-1}\mu^{-\frac54}.
\end{equation}

We next consider $F_{\mathrm{err}}$, using
formulas \eqref{Fauxspatialdecomp} and \eqref{Fauxerrcomm}. After inserting the principal terms from
\eqref{Fauxcoefficientbounds}, we apply
\eqref{Fauxfactorbounds} to the remaining $B,V$ factors.
The subprincipal errors are estimated similarly to above, while the
principal contributions have the forms
\eqref{Fauxprincipalterms1} and \eqref{Fauxprincipalterms2}
and satisfy the same estimates as above.
\\

In summary, by assembling all of the above estimates, we have the decomposition
\begin{equation*}
F_\lambda=F_\lambda^{\mathrm{lo}}+F_\lambda^{\mathrm{hi}},
\end{equation*}
where $F_\lambda^{\mathrm{lo}}$ collects the principal
contributions with $\lambda_0\ll\mu$ and the errors estimated
directly after multiplication by $H_\mu$. These terms contribute directly to $\mathcal{R}^{\mathrm{aux}}_F$. $F_\lambda^{\mathrm{hi}}$ collects the principal contributions
with $\lambda_0\gtrsim\mu$ and the subprincipal errors estimated directly
in $L_{t,x}^2$. Both parts include contributions from
$F_{\mathrm{main}}$ and $F_{\mathrm{err}}$.
By \eqref{FauxprincipalhighL2} and the preceding error bounds,
\begin{equation}\label{FauxhighL2}
\|F_\lambda^{\mathrm{hi}}\|_{L_{t,x}^2}
\lesssim\epsilon^2\lambda^{-1}\mu^{-\frac54}.
\end{equation}
We now argue exactly as in the proof of
\eqref{Bauxdecomp}--\eqref{Bauxtrilinear}.
The product $F_\lambda^{\mathrm{lo}}H_\mu$ contributes directly
to $\mathcal R_F^{\mathrm{aux}}$, while the product of
$F_\lambda^{\mathrm{hi}}$ with the low-frequency part of $H_\mu$
is incorporated into the flux correction $Q_F^{\mathrm{aux}}$ exactly as in the decomposition in Proposition \ref{B-TBexp}.
Its product with the remaining part of $H_\mu$ contributes
to $\mathcal R_F^{\mathrm{aux}}$. This finally gives \eqref{Fauxtrilinear} and completes the proof of Proposition \ref{Fauxbounds}.
\end{proof}
\subsection{The refined paradifferential equation for \texorpdfstring{$\varphi$}{}}
\label{sec:N-phi}
We can now combine all of the expansions in this section to obtain the full layered paradifferential expansion of the equation for $\varphi$, with the layers as computed in Theorem~\ref{t:expansions}. We have

\begin{theorem}\label{Nphifinalcorollary}
Define the quadratic balanced part,
\begin{equation}\label{Nphi2def}
N_{\varphi}^{[2]}
\coloneq
N_{\varphi,B}^{[2]}-N_{\varphi,V}^{[2]}
+N_{\varphi,F}^{[2]}=\frac12\Pi(T_{\mathfrak l}\varphi,T_{\mathfrak l}\varphi)-\frac12\Pi(T_{\mathfrak m}\varphi,T_{\mathfrak m}\varphi)-T_{\nabla\eta}\cdot
\Pi\bigl(
T_{\mathfrak l}\varphi,
T_{\mathfrak m}\varphi
\bigr),
\end{equation}
the cubic balanced part,
\begin{equation}\label{Nphi3def}
\begin{split}
N_{\varphi}^{[3]}
&\coloneq
N_{\varphi,B}^{[3]}-N_{\varphi,V}^{[3]}
+N_{\varphi,F}^{[3]}=-
\Pi\bigl(
T_{\mathfrak l}\varphi,
T_{\mathfrak g}(\eta,T_{\mathfrak l}\varphi)
\bigr)
+\Pi\bigl(
T_{\mathfrak l}\varphi,
\eta,
T_{\mathfrak g}T_{\mathfrak l}\varphi
\bigr),
\end{split}
\end{equation}
and the quartic balanced part,
\begin{equation}\label{Nphi4finaldef}
N_{\varphi}^{[4]}
\coloneq
N_{\varphi,B}^{[4]}-N_{\varphi,V}^{[4]}
+N_{\varphi,F}^{[4]}.
\end{equation}
Then the good variable $\varphi=\psi-T_{B^l}\eta$ satisfies
\begin{equation}\label{varphifinalequation}
\left(\partial_t+T_{V^l}\cdot\nabla\right)\varphi
+T_{a^l}\eta
=N_{\varphi}+N_{\varphi}^{\mathrm{flux}}+\mathcal R,
\end{equation}
where $N_{\varphi}$ has the paradifferential decomposition
\begin{equation}\label{Nphifinalexpansion}
N_{\varphi}
=N_{\varphi}^{[2]}+N_{\varphi}^{[3]}
+N_{\varphi}^{[4]}+\mathcal R_{\varphi}.
\end{equation}
For every $\mu\in2^{\mathbb Z}$ and $x_0\in\mathbb R^2$, the
two remainders satisfy
\begin{equation*}
\|\mathcal R_\mu\overline u_\mu^{x_0}\|_{L_{t,x}^1}
+\|(\mathcal R_{\varphi})_\mu
\overline u_\mu^{x_0}\|_{L_{t,x}^1}
\lesssim
\epsilon^4c_\mu^2\mu^{-2s_c-\frac12}
(1+\mu|x_0|)^{\frac12}.
\end{equation*}
The flux term is as in Proposition \ref{paravarphiinitial};
in particular,
\begin{equation*}
\overline u_\mu^{x_0}N_{\varphi,\mu}^{\mathrm{flux}}
=\nabla_x\cdot Q_{\varphi,\mu},
\end{equation*}
where $Q_{\varphi,\mu}$ satisfies \eqref{QPhibound}.

In addition, the full nonlinearity $N_\varphi$ satisfies the
following auxiliary bounds:
\begin{enumerate}
\item For every $\lambda>0$,
\begin{equation}\label{NphiL2Linf}
\|(N_\varphi)_\lambda\|_{L_t^2L_x^\infty}
\lesssim\epsilon^2\lambda^{-\frac54}.
\end{equation}
\item For every $0<\lambda\ll\mu$ and $x_0\in\mathbb R^2$, there is a
decomposition
\begin{equation}\label{Nphiauxdecomp}
(N_\varphi)_\lambda\overline u_\mu^{x_0}u_\mu
=\mathcal R_\varphi^{\mathrm{aux}}
+\nabla\cdot Q_\varphi^{\mathrm{aux}},
\end{equation}
where, for $\frac43<p\leq2$,
\begin{equation}\label{Nphiauxtrilinear}
\|\mathcal R_\varphi^{\mathrm{aux}}\|_{L_{t,x}^1}
+\lambda^{\frac2p-1}
\|Q_\varphi^{\mathrm{aux}}\|_{L_t^1L_x^p}
\lesssim_p
\epsilon^4c_\mu^2\lambda^{-\frac32}
\mu^{-2s_c}(1+\mu|x_0|)^{\frac12}.
\end{equation}
Here $\mathcal R_\varphi^{\mathrm{aux}}$ is the remainder in the
mixed-frequency product above, and $Q_\varphi^{\mathrm{aux}}$ is
its flux; these are distinct from $\mathcal R_\varphi$ in
\eqref{Nphifinalexpansion} and $Q_{\varphi,\mu}$ in
\eqref{QPhibound}.
\item 
For every $\lambda>0$,
\begin{equation}\label{NphiL2cubic}
\|(N_\varphi-N_\varphi^{[2]})_\lambda\|_{L_{t,x}^2}
\lesssim_\delta
\epsilon^3c_\lambda\lambda^{-\frac14-s_c}.
\end{equation}
\end{enumerate}
\end{theorem}

\begin{proof}
The equation and the multilinear expansion follow by combining
Propositions \ref{paravarphiinitial}, \ref{B-TBexp},
\ref{V-TVexp}, and \ref{cubicvarphi}.
For the auxiliary bounds, we use the identity
\begin{equation*}
N_\varphi=\mathcal B-\mathcal V+F.
\end{equation*}
The first bound follows from \eqref{BauxL2Linf}, its $V$ analogue,
and \eqref{FauxL2Linf}. For the second, take
\begin{equation*}
\mathcal R_\varphi^{\mathrm{aux}}
=\mathcal R_B^{\mathrm{aux}}
-\mathcal R_V^{\mathrm{aux}}+\mathcal R_F^{\mathrm{aux}},
\qquad
Q_\varphi^{\mathrm{aux}}
=Q_B^{\mathrm{aux}}-Q_V^{\mathrm{aux}}+Q_F^{\mathrm{aux}},
\end{equation*}
and apply \eqref{Bauxtrilinear}, its $V$ analogue, and
\eqref{Fauxtrilinear}. Since these estimates were proved for the
full expressions before collecting their multilinear principal
parts, they include the contribution of $\mathcal R_\varphi$ in
\eqref{Nphifinalexpansion}. Finally, subtracting \eqref{Nphi2def} from the identity above and
applying \eqref{BauxL2cubic}, \eqref{VauxL2cubic}, and
\eqref{FauxL2cubic} gives \eqref{NphiL2cubic}.
\end{proof}
\subsection{The refined paradifferential expansion for the linearized good variable}
\label{subsec:Nphi-lin}
We now obtain the corresponding expansion for the 
equation for the linearized good variable $r\coloneq q-T_{B^l}h$. Our starting point is to recall the second equation in
\eqref{WW-lin2} for the linearized potential trace $q$. Using the linearized equation for $h$, we can rewrite the equation for $q$ in the form
\begin{equation}\label{philinearized}
\partial_t q+V\cdot\nabla q+h
=B(\partial_t h+V\cdot\nabla h).
\end{equation}
This equivalent form involving the time derivative of $h$ will be slightly more convenient for isolating certain cancellations in the forthcoming expansions involving the linearization of the kinematic boundary condition. The main expansion for the linearized good variable $r$ is given by the following proposition. 

\begin{theorem}\label{Nphilinearized}
Using the hypotheses and notation of Theorem
\ref{Netalinearized}, the equation for $r=q-T_{B^l}h$ takes
the form
\begin{equation}\label{philinearizedequation}
\left(\partial_t+T_{V^l}\cdot\nabla\right)r
+T_{a^l}h=N_{\varphi}^{\mathrm{lin}},
\end{equation}
where $N_{\varphi}^{\mathrm{lin}}$ has the paradifferential expansion
\begin{equation}\label{Nphilinearizedexpansion}
N_{\varphi}^{\mathrm{lin}}
=N_{\varphi}^{\mathrm{lin},[0]}
+N_{\varphi}^{\mathrm{lin},[1]}
+N_{\varphi}^{\mathrm{lin},[2]}
+\mathcal R_{\varphi}^{\mathrm{lin}}.
\end{equation}
Above, we define the principal quadratic contribution when
the linearized variable is at low frequency,
\begin{equation}\label{Nphilin0def}
N_{\varphi}^{\mathrm{lin},[0]}
\coloneq
T_{T_{\mathfrak g}r}(T_{\mathfrak l}\varphi)
-T_{\nabla r}\cdot T_{\mathfrak m}\varphi.
\end{equation}
The subprincipal contribution is
$N_{\varphi}^{\mathrm{lin},[1]}=0$, and the quadratic
balanced part is defined by
\begin{equation}\label{Nphilin2def}
N_{\varphi}^{\mathrm{lin},[2]}
\coloneq
\Pi(T_{\mathfrak l}\varphi,T_{\mathfrak g}r)
-\Pi(T_{\mathfrak m}\varphi,\nabla r).
\end{equation}
For every $\mu>0$ and $x_0\in\mathbb R^2$, the remainder satisfies
\begin{equation}\label{Nphilinearizedremainder}
\|(\mathcal R_{\varphi}^{\mathrm{lin}})_\mu
\overline w_\mu^{x_0}\|_{L_{t,x}^1([0,T]\times\mathbb R^2)}
\lesssim
\epsilon^2d_\mu^2\mu^{-\frac12}(1+\mu|x_0|)^{\frac12}.
\end{equation}
In addition, we have the following auxiliary bounds.
\begin{enumerate}
\item For every $\lambda>0$, there holds
\begin{equation}\label{NphilinL2remainder}
\|(N_{\varphi}^{\mathrm{lin}}-N_{\varphi}^{\mathrm{lin},[2]})_\lambda
\|_{L_{t,x}^2}
\lesssim_\delta
\epsilon d_\lambda\lambda^{-\frac14}.
\end{equation}
\item For every $0<\nu\ll\mu$ and $x_0\in\mathbb R^2$, there holds
\begin{equation}\label{Nphilinauxtrilinear}
\|(N_{\varphi}^{\mathrm{lin}})_\nu
\overline w_\mu^{x_0}u_\mu\|_{L_{t,x}^1}
\lesssim_\delta
\epsilon^2d_\mu d_\nu\mu^{-\frac54}\nu^{-\frac14}
(1+\mu|x_0|)^{\frac12}.
\end{equation}
\item For every $0<\lambda\ll\mu$ and $x_0\in\mathbb R^2$,
the full nonlinear term $N_\varphi^{\mathrm{err}}$ in
\eqref{Nphierrinitialdef} admits a decomposition
\begin{equation}\label{Nphiwwauxdecomp}
(N_\varphi^{\mathrm{err}})_\lambda w_\mu\overline w_\mu^{x_0}
=\mathcal R_{\varphi,w}^{\mathrm{aux}}
+\nabla\cdot Q_{\varphi,w}^{\mathrm{flux}},
\end{equation}
where, for $\frac43<p\leq2$, we have
\begin{equation}\label{Nphiwwauxtrilinear}
\|\mathcal R_{\varphi,w}^{\mathrm{aux}}\|_{L_{t,x}^1}
+\lambda^{\frac2p-1}\|Q_{\varphi,w}^{\mathrm{flux}}\|_{L_t^1L_x^p}
\lesssim_{\delta,p}
\epsilon^2d_\mu^2\lambda^{-\frac32}(1+\mu|x_0|)^{\frac12}.
\end{equation}
\end{enumerate}
\end{theorem}

We divide the proof into two lemmas which mirror the expansions in subsection \ref{subsec:initialphi} and subsection \ref{subsec:phimultilinear}, respectively. For the remainder of the proof, $\mathcal R$ denotes a remainder satisfying
\eqref{Nphilinearizedremainder}, whose exact form may change
from line to line.

\begin{lemma}\label{philininitiallemma}
The variable $r$ satisfies
\begin{equation}\label{philininitial}
\left(\partial_t+T_{V^l}\cdot\nabla\right)r+T_{a^l}h
=BG(\eta)r-T_{B^l}G(\eta)r
-V\cdot\nabla r+T_{V^l}\cdot\nabla r+\mathcal R.
\end{equation}
\end{lemma}
\begin{proof}
Substituting $q=r+T_{B^l}h$ into \eqref{etalinearized} gives
\begin{equation*}
(\partial_t+V\cdot\nabla)h
=G(\eta)r-G(\eta)(Bh-T_{B^l}h)-h\nabla\cdot V.
\end{equation*}
We now apply $\partial_t+V\cdot\nabla$ to $r=q-T_{B^l}h$ and use
\eqref{philinearized}. Recalling the coefficient
$\widetilde a$ defined in \eqref{tildeadef}, and carrying out similar manipulations to those in subsection \ref{subsec:initialphi}, this gives
\eqref{philininitial} with
\begin{equation}\label{linphiRemainder}
\begin{split}
\mathcal R
&=T_{a^l-\widetilde a^l}h
+\left(T_{(V\cdot\nabla B)^l}
-T_{V^l}\cdot T_{\nabla B^l}\right)h
\\
&-(V-T_{V^l})\cdot T_{\nabla B^l}h
+[T_{B^l},V]\cdot\nabla h
\\
&-(B-T_{B^l})
\left(G(\eta)(Bh-T_{B^l}h)+h\nabla\cdot V\right).
\end{split}
\end{equation}
The term $T_{a^l-\widetilde a^l}h$ is the analogue of
$T_{a^l-\widetilde a^l}\eta$ in \eqref{initialphipara},
equivalently $-T_{\mathcal E^l}\eta$ in
\eqref{Nphiinitialexacterror}.  The term
$\left(T_{(V\cdot\nabla B)^l}
-T_{V^l}\cdot T_{\nabla B^l}\right)h$
is the analogue of
\begin{equation*}
\left(T_{(V\cdot\nabla B)^l}
-T_{V^l}\cdot T_{\nabla B^l}\right)\eta
\end{equation*}
in \eqref{initialphipara}, whose integrand is expanded in
\eqref{phiparacomauxexact}. Arguing almost exactly as in Lemma
\ref{paracomperturbative}, the second term in the first line of
\eqref{linphiRemainder} is perturbative by using \eqref{hboot}
and the corresponding bilinear bounds.

By successive linearization as in Subsection
\ref{subsec:Neta-lin}, the last two lines in
\eqref{linphiRemainder} give cubic and quartic contributions
analogous to $N_{\eta}^{\mathrm{lin},[3]}$ and
$N_{\eta}^{\mathrm{lin},[4]}$, together with high-frequency tails
and quintic or higher order errors. Where two frequency separated pairs are
available, we use two unbalanced bilinear estimates.
The remaining cubic and quartic contributions are estimated
as in \eqref{Netalincubicestimate} and
\eqref{Netalinquarticestimate}, using four
$L_t^{6,\infty}L_x^4$ bounds and placing the additional bounded
factors in $L_{t,x}^\infty$. This also includes terms
where the two background frequencies are comparable
but $h$ is at lower frequency. For the commutator in the second line of
\eqref{linphiRemainder}, we can expand the paraproduct and then apply Proposition \ref{Leibnizcommutator} and the
corresponding bilinear estimates. For the last line of
\eqref{linphiRemainder}, the identity
$G(\eta)B+\nabla\cdot V=0$ gives
\begin{equation*}
G(\eta)(Bh-T_{B^l}h)+h\nabla\cdot V
=[G(\eta),h]B-G(\eta)(T_{B^l}h).
\end{equation*}
We use this cancellation when $h$ is at low frequency,
as in the treatment of the terms containing $h$ in
Subsection \ref{subsec:Neta-lin}.
The high-frequency tails of the geometrically truncated
velocity traces $B_{<\alpha}$ and $V_{<\alpha}$ are estimated
using Proposition \ref{Bernsteinelliptic}.
Together with the preceding cubic and quartic estimates,
this gives
\begin{equation*}
\left\|
\left(\mathcal R-T_{a^l-\widetilde a^l}h\right)_\mu
\overline w_\mu^{x_0}
\right\|_{L_{t,x}^1}
\lesssim
\epsilon^2d_\mu^2\mu^{-\frac12}(1+\mu|x_0|)^{\frac12}.
\end{equation*}

It remains to estimate $T_{a^l-\widetilde a^l}h$. In the nonlinear expansion in Proposition \ref{paravarphiinitial}, this term was not directly perturbative and contributed a flux correction. Such a correction will not be necessary here thanks to the additional $\frac{1}{12}$ derivatives of regularity of the background solution $u$. Indeed, we can use the decompositions in the proof of Lemma
\ref{aerrorinitial} (where we replace the high frequency $\eta$ with the linearized variable $h$), including \eqref{rhoathreeterms}
and \eqref{nullformdef}. The quadratic coefficient
terms previously treated using a flux decomposition
can now be estimated directly using four $L_t^{6,\infty} L_x^4$ bounds. For instance, for $\lambda\approx\mu$ and
$\alpha\approx_\delta\beta\gtrsim_\delta\mu$,
a representative term in these estimates satisfies
\begin{equation}\label{Nphilintaylorerror}
\begin{split}
&\left\|
P_\mu\left(
\left[P_{<M\delta\lambda}T_{\mathfrak l}P_{<\delta\lambda}
\bigl((T_{\mathfrak l}\varphi_\alpha)
(T_{\mathfrak l}\varphi_\beta)\bigr)\right]h_\lambda
\right)\overline w_\mu^{x_0}
\right\|_{L_{t,x}^1}
\\
&\qquad\lesssim_\delta
T^{\frac13}\mu\alpha^2
\|\varphi_\alpha\|_{L_t^{6,\infty}L_x^4}
\|\varphi_\beta\|_{L_t^{6,\infty}L_x^4}
\|h_\lambda\|_{L_t^{6,\infty}L_x^4}
\|w_\mu\|_{L_t^{6,\infty}L_x^4}
\\
&\qquad\lesssim_\delta
T^{\frac13}\epsilon^2
\left(c_\alpha\alpha^{\frac1{12}}\right)^2
d_\mu^2\mu^{-\frac12}
\left(\frac{\mu}{\alpha}\right)^{\frac73}.
\end{split}
\end{equation}
Here the factor $\mu$ comes from the outer
$T_{\mathfrak l}$ acting at frequency $\lesssim\lambda$,
while the two inner operators contribute $\alpha^2$.
The bound is integrable in $\alpha$ and $\beta$ by the
additional regularity of $u$ and the hypothesis $T\leq1$.
The other terms in \eqref{nullformdef} satisfy the
same estimate.

We also need to estimate directly the para-coefficient error
$\mathcal N_{\rho,\lambda}^{\mathrm{err}}$ in
\eqref{rho1rem2}, whose contribution was also treated using
a flux correction there. By retaining the frequency
envelope in the estimates leading to \eqref{rho1rem2} (which was unnecessary at the time) gives the improved decay bound
\begin{equation*}
\|\mathcal N_{\rho,\lambda}^{\mathrm{err}}\|_{L_{t,x}^2}
\lesssim_\delta
\epsilon^3c_\lambda\lambda^{-\frac94}.
\end{equation*}
For $\lambda\approx\mu$, H\"older and \eqref{hboot}
therefore give
\begin{equation*}
\begin{split}
\left\|
P_\mu\left(
\left[P_{<M\delta\lambda}T_L
\mathcal N_{\rho,\lambda}^{\mathrm{err}}\right]h_\lambda
\right)\overline w_\mu^{x_0}
\right\|_{L_{t,x}^1}&\lesssim_\delta
T^{\frac16}\lambda
\|\mathcal N_{\rho,\lambda}^{\mathrm{err}}\|_{L_{t,x}^2}
\|h_\lambda\|_{L_t^{6,\infty}L_x^4}
\|w_\mu\|_{L_t^{6,\infty}L_x^4}
\\
&\lesssim_\delta
T^{\frac16}\epsilon^3c_\mu d_\mu^2\mu^{-\frac5{12}}
\lesssim_\delta
\epsilon^3d_\mu^2\mu^{-\frac12},
\end{split}
\end{equation*}
where the last inequality follows from
$c_\mu\mu^{\frac1{12}}\lesssim1$ and $T\leq1$.

The remaining errors in $T_{a^l-\tilde{a}^l}h$ can be estimated as in the proof of Lemma \ref{aerrorinitial}, but using \eqref{hboot} and the linearized
bilinear bounds in place of the nonlinear ones. Consequently,
\begin{equation*}
\left\|
(T_{a^l-\widetilde a^l}h)_\mu
\overline w_\mu^{x_0}
\right\|_{L_{t,x}^1}
\lesssim
\epsilon^2d_\mu^2\mu^{-\frac12}(1+\mu|x_0|)^{\frac12}.
\end{equation*}
Combining the preceding estimates gives
\begin{equation*}
\|\mathcal R_\mu\overline w_\mu^{x_0}\|_{L_{t,x}^1}
\lesssim
\epsilon^2d_\mu^2\mu^{-\frac12}(1+\mu|x_0|)^{\frac12},
\end{equation*}
which completes the proof of the lemma. 
\end{proof}
To conclude the proof of Theorem \ref{Nphilinearized}, it remains to expand the right hand side of \eqref{philininitial}. This is done in the following lemma.
\begin{lemma}\label{Nphilinexpansionlemma}
The nonlinearity in \eqref{philininitial} has the expansion
\begin{equation*}
BG(\eta)r-T_{B^l}G(\eta)r
-V\cdot\nabla r+T_{V^l}\cdot\nabla r
=N_{\varphi}^{\mathrm{lin},[0]}
+N_{\varphi}^{\mathrm{lin},[2]}+\mathcal R.
\end{equation*}
\end{lemma}
\begin{proof}
We begin with the paraproduct decompositions
\begin{equation*}
\begin{split}
BG(\eta)r-T_{B^l}G(\eta)r
&=T_{G(\eta)r}B+\Pi(B,G(\eta)r)
+(T_B-T_{B^l})G(\eta)r,
\\
V\cdot\nabla r-T_{V^l}\cdot\nabla r
&=T_{\nabla r}\cdot V+\Pi(V,\nabla r)
+(T_V-T_{V^l})\cdot\nabla r.
\end{split}
\end{equation*}
We then carry out the same successive linearizations
used in Propositions \ref{B-TBexp}, \ref{V-TVexp},
and Subsection \ref{subsec:Neta-lin}. This gives
\begin{equation*}
\begin{split}
BG(\eta)r-T_{B^l}G(\eta)r
&=T_{T_{\mathfrak g}r}(T_{\mathfrak l}\varphi)
+\Pi(T_{\mathfrak l}\varphi,T_{\mathfrak g}r)+N_{\varphi}^{\mathrm{lin},[3]}
+N_{\varphi}^{\mathrm{lin},[4]}+\mathcal R,
\\
V\cdot\nabla r-T_{V^l}\cdot\nabla r
&=T_{\nabla r}\cdot T_{\mathfrak m}\varphi
+\Pi(T_{\mathfrak m}\varphi,\nabla r)+N_{\varphi}^{\mathrm{lin},[3]}
+N_{\varphi}^{\mathrm{lin},[4]}+\mathcal R.
\end{split}
\end{equation*}
Here $N_{\varphi}^{\mathrm{lin},[3]}$,
$N_{\varphi}^{\mathrm{lin},[4]}$, and $\mathcal R$ denote
generic contributions whose exact form may change
between occurrences. In particular, we
use \eqref{Bfundcalc} and its $V$ analogue to extract
$T_{\mathfrak l}\varphi$ and $T_{\mathfrak m}\varphi$,
and the expansion leading to
\eqref{linearizationetaexplanation} to extract
$T_{\mathfrak g}r$ from $G(\eta)r$.

The resulting quadratic terms after subtracting
the second expansion from the first are precisely
\eqref{Nphilin0def} and \eqref{Nphilin2def}.
As in Subsection \ref{subsec:Neta-lin}, we denote the
remaining balanced cubic and quartic contributions
generically by $N_{\varphi}^{\mathrm{lin},[3]}$ and
$N_{\varphi}^{\mathrm{lin},[4]}$.
Here these contain one $r$ input, one $\varphi$ input,
and respectively one or two $\eta$ inputs at comparable
frequencies. The other errors are estimated as in the
preceding nonlinear expansions, using \eqref{hboot}
and the mixed bilinear bounds in place of the nonlinear ones.

For completeness, a representative term of the form $N_{\varphi}^{\mathrm{lin},[3]}$  for
$\alpha\approx_\delta\beta\approx_\delta\gamma$
and $\mu\lesssim_\delta\alpha$, can be estimated using four $L_t^6L_x^4$ bounds and H\"older in time by
\begin{equation}\label{Nphilincubicestimate}
\begin{split}
\left\|
P_\mu\left(
(T_{\mathfrak l}\varphi_\alpha)
T_{\mathfrak g}(\eta_\beta T_{\mathfrak l}r_\gamma)
\right)\overline w_\mu^{x_0}
\right\|_{L_{t,x}^1}&\lesssim_\delta
T^{\frac13}\alpha^3
\|\varphi_\alpha\|_{L_t^{6,\infty}L_x^4}
\|\eta_\beta\|_{L_t^{6,\infty}L_x^4}
\|r_\gamma\|_{L_t^{6,\infty}L_x^4}
\|w_\mu\|_{L_t^{6,\infty}L_x^4}
\\
&\lesssim_\delta
T^{\frac13}\epsilon^2
\left(c_\alpha\alpha^{\frac1{12}}\right)^2
d_\alpha d_\mu\mu^{-\frac12}
\left(\frac{\mu}{\alpha}\right)^{\frac{11}{12}}.
\end{split}
\end{equation}
The cubic terms with comparable background frequencies
and a lower-frequency $r$ input are also included in
$\mathcal R$. The same four $L_t^{6,\infty} L_x^4$ bounds work in thie case too, and give
the additional factor $(\gamma/\alpha)^{\frac{11}{12}}$
when $\gamma\ll\alpha\approx_\delta\beta$, so the
integration in $\gamma$ also converges by slow variation
of $d_\gamma$. The quartic terms are estimated by placing the
additional high frequency $\eta$ factor in $L_{t,x}^\infty$, as in
\eqref{Netalinquarticestimate}. Integrating and using
slow variation of $d_\alpha$, we then obtain
\begin{equation*}
\begin{split}
\|(N_{\varphi}^{\mathrm{lin},[3]})_\mu
\overline w_\mu^{x_0}\|_{L_{t,x}^1}
+\|(N_{\varphi}^{\mathrm{lin},[4]})_\mu
\overline w_\mu^{x_0}\|_{L_{t,x}^1}&\lesssim_\delta
\epsilon^2d_\mu\mu^{-\frac12}
\int_{\alpha\gtrsim_\delta\mu}
d_\alpha\left(\frac{\mu}{\alpha}\right)^{\frac{11}{12}}
\frac{d\alpha}{\alpha}
\\
&\lesssim_\delta
\epsilon^2d_\mu^2\mu^{-\frac12}.
\end{split}
\end{equation*}
Thus these contributions can also be included in
$\mathcal R$, which proves the lemma.
\end{proof}

\begin{proof}[Proof of Theorem \ref{Nphilinearized}]
Combining Lemmas \ref{philininitiallemma} and
\ref{Nphilinexpansionlemma}, and collecting their
remainders into $\mathcal R_{\varphi}^{\mathrm{lin}}$,
gives \eqref{philinearizedequation},
\eqref{Nphilinearizedexpansion}, and
\eqref{Nphilinearizedremainder}.

We now prove \eqref{NphilinL2remainder}. As in subsection \ref{subsec:Neta-lin}, we can estimate
$N_{\varphi}^{\mathrm{lin},[0]}$ directly by using an unbalanced
bilinear estimate for each of the two terms, alongside the basic kernel estimates for the paradifferential operators
$T_{\mathfrak l}$, $T_{\mathfrak m}$ and
$T_{\mathfrak g}$. This gives
\begin{equation*}
\begin{split}
\|(N_{\varphi}^{\mathrm{lin},[0]})_\nu\|_{L_{t,x}^2}
&\lesssim_\delta
\epsilon c_\nu\nu^{-\frac32}
\int_{\beta\ll\nu}
d_\beta\beta^{\frac54}\frac{d\beta}{\beta}
\lesssim_\delta
\epsilon d_\nu\nu^{-\frac14}.
\end{split}
\end{equation*}

For the first line of \eqref{linphiRemainder}, we use
\eqref{Ephiinitialcoefficient} and the paracoefficient
bounds following \eqref{phiparacomauxexact}, exactly as in
the proof of \eqref{NphiinitialL2}. Placing
$h_\lambda$ in $L_t^\infty L_x^2$ gives
\begin{equation*}
\begin{split}
\left\|
\left(
T_{a^l-\widetilde a^l}h
+\left(T_{(V\cdot\nabla B)^l}
-T_{V^l}\cdot T_{\nabla B^l}\right)h
\right)_\nu
\right\|_{L_{t,x}^2}\lesssim_\delta
\epsilon^2
\int_{\lambda\approx\nu}
\lambda^{-\frac14}\|h_\lambda\|_{L_t^\infty L_x^2}
\frac{d\lambda}{\lambda}
\lesssim_\delta
\epsilon^2d_\nu\nu^{-\frac14}.
\end{split}
\end{equation*}
Next, we estimate the contribution of the second and third line of. \eqref{linphiRemainder}. We first estimate
$-(V-T_{V^l})\cdot T_{\nabla B^l}h$, writing
\begin{equation*}
V-T_{V^l}=(V-T_V)+(T_V-T_{V^l}).
\end{equation*}
For the first contribution, we use Proposition \ref{VBunbil}
and \eqref{hboot}, placing the two velocity trace factors in
$L_t^4L_x^\infty$ and $h_\gamma$ in $L_t^\infty L_x^2$.
This gives
\begin{equation*}
\begin{split}
\|P_\nu((V-T_V)\cdot T_{\nabla B^l}h)\|_{L_{t,x}^2}
&\lesssim_\delta
\epsilon^2
\int_{\alpha\gtrsim_\delta\nu}
\int_{\gamma\lesssim_\delta\alpha}
d_\gamma\alpha^{-\frac58}\gamma^{\frac38}
\frac{d\gamma}{\gamma}\frac{d\alpha}{\alpha}
\\
&\lesssim_\delta
\epsilon^2d_\nu\nu^{-\frac14}.
\end{split}
\end{equation*}
The contribution containing $T_V-T_{V^l}$ satisfies the
same bound by using the paracoefficient estimates following
\eqref{Fauxcoefficientbounds}.

We next estimate $[T_{B^l},V]\cdot\nabla h$. Expanding $V$, $h$ in in their frequency localizations, when $V_{\beta}$ and $\nabla h_{\alpha}$ have comparable frequencies $\alpha\approx\beta$  and $B$
is at lower frequency, we combine it with the corresponding
contribution of $-(B-T_{B^l})(h\nabla\cdot V)$, using
\begin{equation*}
\begin{split}
[T_{B^l},V]\cdot\nabla h
-(B-T_{B^l})(h\nabla\cdot V)
=-(B-T_{B^l})\nabla\cdot(Vh)
+V\cdot(B-T_{B^l})\nabla h.
\end{split}
\end{equation*}
By expanding the $B$ para product, we can similarly estimate this contribution using the same two $L_t^4L_x^{\infty}$ bounds and the $L_t^{\infty}L_x^2$ bound for $h$.

When $V_\beta$ is at lower frequency than $h_\alpha$,
we expand $[T_{B^l},V_\beta]\cdot\nabla h_\alpha$
and apply Proposition \ref{Leibnizcommutator} to
$[S_\lambda,V_\beta]\cdot\nabla h_\alpha$, where
$\lambda\approx\alpha$. This gives a factor
$\beta/\alpha$, which cancels the derivative on
$h_\alpha$ and leaves a derivative on $V_\beta$. The same two
$L_t^4L_x^\infty$ bounds and the $L_t^{\infty}L_x^2$ bound give, for
$\beta\ll\alpha$,
\begin{equation*}
\begin{split}
\|P_\nu[T_{B^l},V_\beta]\cdot\nabla h_\alpha\|_{L_{t,x}^2}
+\|P_\nu[T_{B^l},V_\alpha]\cdot\nabla h_\beta\|_{L_{t,x}^2}\lesssim_\delta
\epsilon^2(d_\alpha+d_\beta)
\alpha^{-\frac58}\beta^{\frac38}
\end{split}
\end{equation*}
which upon integration in $\alpha$ and $\beta$ gives the required $L_{t,x}^2$ bound. A similar bound again applies to $-(B-T_{B^l})(h_{\alpha}\nabla\cdot V_{\beta})$ since the low frequency factor of $V_{\beta}$ is differentiated here. When $\alpha\ll\beta$, we can can estimate the commutator directly using the same three norm allocations by expanding
\begin{equation*}
[T_{B^l},V_{\beta}]\cdot\nabla h_{\alpha}=[T_{B^l}-B,V_{\beta}]\cdot\nabla h_{\alpha}.    
\end{equation*}
The important point here is that $h$ is differentiated when at low frequency. Unfortunately for this reason, we cannot directly estimate the contribution of 
\begin{equation*}
-(B-T_{B^l})(h_{\alpha}\nabla\cdot V_{\beta}) 
\end{equation*}
when $\alpha\ll\beta$. We therefore combine the remaining contribution
$-(B-T_{B^l})T_h\nabla\cdot V$ with the term containing
$G(\eta)$ in \eqref{linphiRemainder}. Indeed, by using
\begin{equation*}
Bh-T_{B^l}h=T_hB+\Pi(B,h)+(T_B-T_{B^l})h
\end{equation*}
and $G(\eta)B+\nabla\cdot V=0$, we obtain
\begin{equation*}
\begin{split}
G(\eta)(Bh-T_{B^l}h)+T_h\nabla\cdot V=[G(\eta),T_h]B
+G(\eta)\Pi(B,h)
+G(\eta)((T_B-T_{B^l})h).
\end{split}
\end{equation*}
Thus the three terms on the right, each acted on by
$-(B-T_{B^l})$, are the remaining contributions we need to estimate.

We first consider $G(\eta)\Pi(B,h)$. We keep $\alpha$
as the frequency of $h$, and write $\beta$ for the
frequency of the inner $B$ factor and $\gamma$ for that
of the outer $B$ factor. Here $\beta\approx_\delta\alpha$. By duality, Proposition \ref{basicopbounds} and Bernstein, we have for each $\rho>0$,
\begin{equation*}
\|P_{<\rho}G(\eta)f\|_{L_x^2}
\lesssim\rho\|f\|_{L_x^2}.
\end{equation*}
Consequently, by Proposition \ref{besovparabounds} we have
\begin{equation*}
\begin{split}
\left\|
P_\nu\left(
(B_\gamma-T_{B_\gamma})G(\eta)(B_\beta h_\alpha)
\right)
\right\|_{L_{t,x}^2}&\lesssim_\delta
\min\{\alpha,\gamma\}
\|B_\gamma\|_{L_t^4L_x^\infty}
\|B_\beta\|_{L_t^4L_x^\infty}
\|h_\alpha\|_{L_t^\infty L_x^2}
\\
&\lesssim_\delta
\epsilon^2d_\alpha
\alpha^{-\frac58}\gamma^{-\frac58}
\min\{\alpha,\gamma\}.
\end{split}
\end{equation*}
The term $G(\eta)((T_B-T_{B^l})h)$ has the same estimate
using the paracoefficient bound following
\eqref{Fauxcoefficientbounds}. The difference $T_B-T_{B^l}$ in the outer factor is estimated using the same paracoefficient bounds.

For $[G(\eta),T_h]B$, we use the same commutator
expansion as in the treatment of the low-frequency
$h$ contribution in Subsection \ref{subsec:Neta-lin} and then use the same three norm allocations as above again.  Together with the estimates above for the other terms in \eqref{linphiRemainder}, this finally gives
\begin{equation*}
\begin{split}
\left\|
P_\nu
\mathcal R
\right\|_{L_{t,x}^2}\lesssim_\delta
\epsilon^2d_\nu\nu^{-\frac14}.
\end{split}
\end{equation*}
The remaining cubic and quartic contributions from
both lemmas are estimated as in the proofs of
\eqref{BauxL2cubic} and \eqref{VauxL2cubic}, using
\eqref{hboot} in the linearized slot. For example,
when $\alpha\approx_\delta\beta\approx_\delta\gamma$,
the cubic term in \eqref{Nphilincubicestimate} satisfies
\begin{equation*}
\begin{split}
\left\|
P_\nu\left(
(T_{\mathfrak l}\varphi_\alpha)
T_{\mathfrak g}(\eta_\beta T_{\mathfrak l}r_\gamma)
\right)
\right\|_{L_{t,x}^2}\lesssim_\delta
\alpha^3
\|\varphi_\alpha\|_{L_t^4L_x^\infty}
\|\eta_\beta\|_{L_t^4L_x^\infty}
\|r_\gamma\|_{L_t^\infty L_x^2}
\lesssim_\delta
\epsilon^2d_\alpha\alpha^{-\frac14}.
\end{split}
\end{equation*}
The cubic terms containing $h$ have the same bound,
placing $h$ in $L_t^\infty L_x^2$ and the two
background $\varphi$ factors in $L_t^4L_x^\infty$. For $\gamma\ll\alpha\approx_\delta\beta$, the same norm allocations give
\begin{equation*}
\left\|
P_\nu\left(
(T_{\mathfrak l}\varphi_\alpha)
T_{\mathfrak g}(\eta_\beta T_{\mathfrak l}r_\gamma)
\right)
\right\|_{L_{t,x}^2}
\lesssim_\delta
\epsilon^2d_\gamma\alpha^{-\frac34}\gamma^{\frac12}.
\end{equation*}
Integrating in $\gamma$ gives the same bound as above by slow variation of $d_\gamma$ (the pertinent point here is that the linearized variable is differentiated in this frequency range to provide the needed integrability in $\gamma$). The quartic terms are estimated as usual by placing the
additional differentiated $\eta$ factor in
$L_{t,x}^\infty$. The remaining replacement errors from replacing elliptic operators by their principal part and the higher order errors are estimated using the same operator decompositions as in the proofs of \eqref{BauxL2cubic} and \eqref{VauxL2cubic}, with \eqref{hboot} and the linearized bilinear bounds in place of the nonlinear ones. The high-frequency tails are estimated using Proposition \ref{Bernsteinelliptic}. Combining everything, we obtain
\begin{equation*}
\begin{split}
\|(N_{\varphi}^{\mathrm{lin}}-N_{\varphi}^{\mathrm{lin},[2]})_\nu
\|_{L_{t,x}^2}
&\lesssim_\delta
\epsilon d_\nu\nu^{-\frac14}
+\epsilon^2
\int_{\alpha\gtrsim_\delta\nu}
d_\alpha\alpha^{-\frac14}\frac{d\alpha}{\alpha}
\\
&\lesssim_\delta
\epsilon d_\nu\nu^{-\frac14}.
\end{split}
\end{equation*}
This proves \eqref{NphilinL2remainder}.

It remains to prove \eqref{Nphilinauxtrilinear} and \eqref{Nphiwwauxtrilinear}; we begin with the former.
As in the proof of \eqref{Netalinauxtrilinear},
no flux correction is needed. First, H\"older and
\eqref{NphilinL2remainder} give
\begin{equation*}
\begin{split}
\|(N_{\varphi}^{\mathrm{lin}}-N_{\varphi}^{\mathrm{lin},[2]})_\nu
\overline w_\mu^{x_0}u_\mu\|_{L_{t,x}^1}\lesssim_\delta
T^{\frac16}\epsilon d_\nu\nu^{-\frac14}
\|w_\mu\|_{L_t^{6,\infty}L_x^4}
\|u_\mu\|_{L_t^{6,\infty}L_x^4}
\lesssim_\delta
\epsilon^2d_\nu d_\mu\nu^{-\frac14}\mu^{-\frac54},
\end{split}
\end{equation*}
where we used $T\leq1$ and
$c_\mu\lesssim\mu^{-\frac1{12}}$.

For $N_{\varphi}^{\mathrm{lin},[2]}$, we use four
$L_t^6L_x^4$ bounds. For example, when
$\alpha\approx_\delta\beta\gtrsim_\delta\nu$, the
first term in \eqref{Nphilin2def} satisfies
\begin{equation*}
\begin{split}
\left\|
P_\nu\left(
(T_{\mathfrak l}\varphi_\alpha)
(T_{\mathfrak g}r_\beta)
\right)\overline w_\mu^{x_0}u_\mu
\right\|_{L_{t,x}^1}&\lesssim_\delta
T^{\frac13}\alpha^2
\|\varphi_\alpha\|_{L_t^{6,\infty}L_x^4}
\|r_\beta\|_{L_t^{6,\infty}L_x^4}
\|w_\mu\|_{L_t^{6,\infty}L_x^4}
\|u_\mu\|_{L_t^{6,\infty}L_x^4}
\\
&\lesssim_\delta
T^{\frac13}\epsilon^2c_\alpha c_\mu
d_\alpha d_\mu\alpha^{-\frac16}\mu^{-\frac76}
\\
&\lesssim_\delta
\epsilon^2d_\alpha d_\mu
\alpha^{-\frac14}\mu^{-\frac54}.
\end{split}
\end{equation*}
The term containing
$T_{\mathfrak m}\varphi_\alpha\cdot\nabla r_\beta$
has the same bound. Integrating in $\alpha$ and
$\beta$, and using slow variation of $d_\alpha$, gives
\begin{equation*}
\begin{split}
\|(N_{\varphi}^{\mathrm{lin},[2]})_\nu
\overline w_\mu^{x_0}u_\mu\|_{L_{t,x}^1}
&\lesssim_\delta
\epsilon^2d_\mu\mu^{-\frac54}
\int_{\alpha\gtrsim_\delta\nu}
d_\alpha\alpha^{-\frac14}\frac{d\alpha}{\alpha}
\\
&\lesssim_\delta
\epsilon^2d_\nu d_\mu\nu^{-\frac14}\mu^{-\frac54}.
\end{split}
\end{equation*}
Combining the preceding estimates proves
\eqref{Nphilinauxtrilinear}.

Finally, we prove \eqref{Nphiwwauxdecomp} and \eqref{Nphiwwauxtrilinear}.
Again, we write $H_\mu=w_\mu\overline w_\mu^{x_0}$. For the contribution of
$N_\varphi^{\mathrm{err}}-N_\varphi$, we argue as in the proof of
\eqref{NphiinitialauxR}, pairing $w_\mu$ with the $\eta$ input
at frequency $\approx\lambda$ using \eqref{wu-bi-unbal-boot} and then placing the remaining $\overline w_\mu^{x_0}$ factor in $L_t^\infty L_x^2$.
This gives
\begin{equation*}
\|(N_\varphi^{\mathrm{err}}-N_\varphi)_\lambda H_\mu\|_{L_{t,x}^1}
\lesssim_\delta
\epsilon^3d_\mu^2\lambda^{-\frac32}(1+\mu|x_0|)^{\frac12}.
\end{equation*}

For $N_\varphi=\mathcal B-\mathcal V+F$, we use the splittings
in the proofs of \eqref{Bauxtrilinear}, \eqref{Vauxtrilinear},
and \eqref{Fauxtrilinear}. Write $F_\lambda^{\mathrm{hi}}$ for
the sum of their high-frequency contributions, with the same signs.
For the remaining contributions, we use
\eqref{wu-bi-unbal-boot} whenever a lower-frequency
background input is paired with a $w_\mu$ factor, and place any
unpaired $w_\mu$ factor in $L_t^\infty L_x^2$. As in the proof
of \eqref{Netawwauxtrilinear}, this replaces each factor
$\epsilon c_\mu\mu^{-s_c}$ by $d_\mu$ and gives the bound
\begin{equation*}
\|\bigl((N_\varphi)_\lambda-F_\lambda^{\mathrm{hi}}\bigr)H_\mu\|_{L_{t,x}^1}
\lesssim_\delta
\epsilon^2d_\mu^2\lambda^{-\frac32}(1+\mu|x_0|)^{\frac12}.
\end{equation*}

By \eqref{BauxhighL2} (and its $V$ analogue), and \eqref{FauxhighL2},
we have
$\|F_\lambda^{\mathrm{hi}}\|_{L_{t,x}^2}
\lesssim_\delta\epsilon^2\lambda^{-1}\mu^{-\frac54}$.
We therefore take
$Q_{\varphi,w}^{\mathrm{flux}}
=\nabla\Delta^{-1}(F_\lambda^{\mathrm{hi}}P_{<\delta_0\lambda}H_\mu)$,
with $\delta_0$ as in the proof of \eqref{Netawwauxtrilinear}.
Arguing as in that proof, we obtain
\begin{equation*}
\begin{split}
\|F_\lambda^{\mathrm{hi}}P_{\geq\delta_0\lambda}H_\mu\|_{L_{t,x}^1}
+\lambda^{\frac2p-1}\|Q_{\varphi,w}^{\mathrm{flux}}\|_{L_t^1L_x^p}
&\lesssim_{\delta,p}
\|F_\lambda^{\mathrm{hi}}\|_{L_{t,x}^2}
d_\mu^2\lambda^{-\frac12}\mu^{\frac54}(1+\mu|x_0|)^{\frac12}.
\end{split}
\end{equation*}
Combining the three estimates proves \eqref{Nphiwwauxtrilinear}.
\end{proof}


\section{Expansions for \texorpdfstring{$a$}{}}\label{sec:aexpansions}
In this section, our main goal will be to prove the necessary linear and bilinear estimates as well as error estimates for the Taylor coefficient $a$ and its geometrically truncated variant $a_{<\lambda}$ defined below. 
\subsection{Preliminary Identities and Linearization}
We begin by recalling our conventions for the pressure set out in Proposition \ref{paravarphiinitial}. We also retain the bootstrap assumptions from Section \ref{sec:bootstrap},
including \eqref{Taylorcoeffsmall} throughout this section. Let $p$ denote the pressure in
$\Omega_\eta=\{z<\eta(x)\}$, characterized by
\begin{equation*}
\Delta_{x,z}p=-|\nabla_{x,z}^2\phi|^2,
\qquad p|_{z=\eta}=0,\qquad p_{<\delta\lambda}(t,x,z)+z\longrightarrow0
\qquad\text{as }z\longrightarrow-\infty.
\end{equation*}
and the corresponding Taylor coefficient, $a\coloneq-\left.\partial_zp\right|_{z=\eta}.$ As in Proposition \ref{paravarphiinitial}, for each $\lambda$, we denote by $p_{<\lambda}$, the pressure obtained by the
same Laplace equation as above, but in the regularized domain
$\Omega_{<\lambda}=\{z<\eta_{<\lambda}(x)\}$. We define
\begin{equation*}
a_{<\lambda}\coloneq
-\left.\partial_zp_{<\lambda}\right|_{z=\eta_{<\lambda}}
\end{equation*}
to denote the corresponding Taylor coefficient for the regularized domain. We begin by recalling the following basic formula for $a$ from \eqref{a-formula} that algebraically relates the Taylor coefficient to the variables $\eta$, $B$ and $V$. We will use this identity below along with the estimates in Proposition \ref{VBunbil} to obtain several useful estimates for $a$.
\begin{lemma}\label{asimpleidentity} Define $Q_a(f,g)\coloneq \frac{1}{2}(fLg+gLf-L(fg))$. The Taylor coefficient satisfies the identity
\begin{equation*}
a = 1-L(\eta)\eta +Q_a(B,B) +Q_a(V,V).
\end{equation*}
The same identity holds for the geometrically truncated quantities.
For these quantities, $Q_a$ is defined using $L_{<\lambda}$.
For vector arguments, $Q_a(V,V)\coloneq\sum_jQ_a(V_j,V_j)$; products
and commutators involving vectors are contracted componentwise, and
$|MV|^2=\sum_{i,j}|M_iV_j|^2$.
\end{lemma}
\begin{proof}
This is trivially seen to be equivalent to the identity \eqref{a-formula}.    
\end{proof}
Our next aim is to use the above identity to linearize $a$ with respect to its geometric regularization. This is a relatively simple consequence of the above lemma.
\begin{lemma}[Linearization formula for $a$]
The truncated Taylor coefficient $a_{<\lambda}$ satisfies
\begin{equation}\label{linearizedaformula}
\begin{split}
\lambda\frac{d}{d\lambda}a_{<\lambda}&=-L_{<\lambda}(a_{<\lambda}\eta_{\lambda})+\frac{\Gamma_{<\lambda}}{1+|\nabla\eta_{<\lambda}|^2}\eta_{\lambda}
\\
&+L_{<\lambda}\Phi_{<\lambda}L_{<\lambda}B_{<\lambda}-[L_{<\lambda},B_{<\lambda}](L_{<\lambda}\Phi_{<\lambda})
\\
&+M_{<\lambda}\Phi_{<\lambda} L_{<\lambda}V_{<\lambda}-[L_{<\lambda},V_{<\lambda}](M_{<\lambda}\Phi_{<\lambda})
\end{split}
\end{equation}
where
\begin{equation*}
\begin{split}
\Gamma_{<\lambda}&\coloneq |M_{<\lambda} V_{<\lambda}|^2+|M_{<\lambda} B_{<\lambda}|^2+|L_{<\lambda} V_{<\lambda}|^2+|L_{<\lambda} B_{<\lambda}|^2
\\
&+a_{<\lambda}\Delta\eta_{<\lambda}+2\nabla\eta_{<\lambda}\cdot\nabla a_{<\lambda}
\end{split}
\end{equation*}
\end{lemma}
\begin{proof}
From the linearization formulas for $L$ and $V$, we find
\begin{equation*}
\begin{split}
\lambda\frac{d}{d\lambda}Q_{a}(V_{<\lambda},V_{<\lambda})&=-V_{<\lambda}L_{<\lambda}(\eta_{\lambda}L_{<\lambda}V_{<\lambda})+\frac{1}{2}L_{<\lambda}(\eta_{\lambda}L_{<\lambda}|V_{<\lambda}|^2)
\\
&+\eta_{\lambda}V_{<\lambda}L^2_{<\lambda}V_{<\lambda}-\frac{1}{2}\eta_{\lambda}L^2_{<\lambda}|V_{<\lambda}|^2
\\
&+(M_{<\lambda}\Phi_{<\lambda}+\eta_{\lambda} L_{<\lambda} V_{<\lambda})L_{<\lambda}V_{<\lambda}
\\
&-[L_{<\lambda},V_{<\lambda}](M_{<\lambda}\Phi_{<\lambda}+\eta_{\lambda}L_{<\lambda}V_{<\lambda})
\end{split}
\end{equation*}
The first two top terms can be rearranged to give
\begin{equation*}
\begin{split}
\lambda\frac{d}{d\lambda}Q_a(V_{<\lambda},V_{<\lambda})&=M_{<\lambda}\Phi_{<\lambda}L_{<\lambda}V_{<\lambda}-[L_{<\lambda},V_{<\lambda}](M_{<\lambda}\Phi_{<\lambda})
\\
&+\eta_{\lambda}\left(V_{<\lambda}L^2_{<\lambda}V_{<\lambda}-\frac{1}{2}L^2_{<\lambda}|V_{<\lambda}|^2+ |L_{<\lambda} V_{<\lambda}|^2\right)
\\
&=-L_{<\lambda}(\eta_{\lambda}Q_a(V_{<\lambda},V_{<\lambda}))
\end{split}
\end{equation*}
Similarly, we find
\begin{equation*}
\begin{split}
\lambda\frac{d}{d\lambda}Q_a(B_{<\lambda},B_{<\lambda})&=L_{<\lambda}\Phi_{<\lambda}L_{<\lambda}B_{<\lambda}-[L_{<\lambda},B_{<\lambda}]L_{<\lambda}\Phi_{<\lambda}
\\
&+\eta_{\lambda}\left(B_{<\lambda}L^2_{<\lambda}B_{<\lambda}-\frac{1}{2}L^2_{<\lambda}B_{<\lambda}^2+(L_{<\lambda}B_{<\lambda})^2\right)
\\
&-L_{<\lambda}(\eta_{\lambda}Q_a(B_{<\lambda},B_{<\lambda}))
\end{split}
\end{equation*}
Let us now write as shorthand
\begin{equation*}
\mathcal S_{<\lambda} := Q_{a}(B_{<\lambda},B_{<\lambda})+Q_a(V_{<\lambda},V_{<\lambda}).
\end{equation*}
From Lemma \ref{asimpleidentity}, we have
\begin{equation}\label{Sarelation}
\mathcal S_{<\lambda}=a_{<\lambda}-1+L_{<\lambda}\eta_{<\lambda}.
\end{equation}
Using the above formulas and the linearization formula for $L_{<\lambda}$, we then obtain
\begin{equation}\label{dlambdaa}
\begin{split}
\lambda\frac{d}{d\lambda}a_{<\lambda}&=-L_{<\lambda}\eta_{\lambda}+L_{<\lambda}(\eta_{\lambda}L_{<\lambda}\eta_{<\lambda})-L_{<\lambda}(\eta_{\lambda}\mathcal{S}_{<\lambda})
\\
&+\eta_{\lambda}\left(B_{<\lambda}L^2_{<\lambda}B_{<\lambda}-\frac{1}{2}L^2_{<\lambda}B^2_{<\lambda}+(L_{<\lambda}B_{<\lambda})^2\right)
\\
&+\eta_{\lambda}\left(V_{<\lambda}L^2_{<\lambda}V_{<\lambda}-\frac{1}{2}L^2_{<\lambda}|V_{<\lambda}|^2+|L_{<\lambda} V_{<\lambda}|^2\right)-\eta_{\lambda}L_{<\lambda}^2\eta_{<\lambda}
\\
&+M_{<\lambda}\Phi_{<\lambda}L_{<\lambda}V_{<\lambda}-[L_{<\lambda},V_{<\lambda}](M_{<\lambda}\Phi_{<\lambda})
\\
&+L_{<\lambda}\Phi_{<\lambda}L_{<\lambda}B_{<\lambda}-[L_{<\lambda},B_{<\lambda}](L_{<\lambda}\Phi_{<\lambda})
\end{split}
\end{equation}
The last two lines above combine to give precisely the last two lines in \eqref{linearizedaformula}. Using Lemma \ref{asimpleidentity}, the first line simply becomes
\begin{equation}\label{firstlinea}
\begin{split}
&-L_{<\lambda}\eta_\lambda
+L_{<\lambda}(\eta_\lambda L_{<\lambda}\eta_{<\lambda})
-L_{<\lambda}(\eta_\lambda\mathcal S_{<\lambda})
\\
&\qquad=-L_{<\lambda}\eta_\lambda
+L_{<\lambda}(\eta_\lambda L_{<\lambda}\eta_{<\lambda})
-L_{<\lambda}\bigl(\eta_\lambda
(a_{<\lambda}-1+L_{<\lambda}\eta_{<\lambda})\bigr)
\\
&\qquad=-L_{<\lambda}(a_{<\lambda}\eta_\lambda)
\end{split}
\end{equation}
which is the first term in \eqref{linearizedaformula} Next, we show that the second and third line in \eqref{dlambdaa} add up to precisely the term
\begin{equation*}
\frac{\Gamma_{<\lambda}}{1+|\nabla\eta_{<\lambda}|^2}\eta_{\lambda}.
\end{equation*}
We first record a useful coordinate formula for $L_{<\lambda}^2f$. Let $\phi_f$ denote the harmonic extension of $f$ with respect to $\eta_{<\lambda}$. Since $L_{<\lambda}f = \left.\partial_z \phi_f\right|_{z=\eta_{<\lambda}}$ and $\partial_z \phi_f$ is again harmonic, its boundary trace is $L_{<\lambda}f$, and therefore $L_{<\lambda}^2f = \left.\partial_z^2\phi_f\right|_{z=\eta_{<\lambda}}.$ The chain rule gives
\begin{equation*}
\nabla_x(L_{<\lambda}f) = \left. \left( \nabla_x\partial_z \phi_f +\partial_z^2\phi_f\,\nabla\eta_{<\lambda} \right) \right|_{z=\eta_{<\lambda}},
\end{equation*}
as well as
\begin{equation*}
\Delta f = \left. \left( \Delta_x\phi_f +2\nabla\eta_{<\lambda}\cdot\nabla_x\partial_z \phi_f +|\nabla\eta_{<\lambda}|^2\partial_z^2\phi_f +(\partial_z \phi_f)\Delta\eta_{<\lambda} \right) \right|_{z=\eta_{<\lambda}}.
\end{equation*}
Then using that $\Delta_x\phi_f+\partial_z^2\phi_f=0$ and eliminating $\nabla_x\partial_z \phi_f$, we obtain
\begin{equation}\label{L2formula}
L_{<\lambda}^2f = \frac{ 2\nabla\eta_{<\lambda}\cdot\nabla(L_{<\lambda}f) -\Delta f +(L_{<\lambda}f)\Delta\eta_{<\lambda} }{1+|\nabla\eta_{<\lambda}|^2}.
\end{equation}
Applying \eqref{L2formula} to $f$ and to $f^2$ gives
\begin{equation*}
\begin{split}
(1+|\nabla\eta_{<\lambda}|^2)\left( fL_{<\lambda}^2f-\frac12L_{<\lambda}^2(f^2) \right) &= 2f\nabla\eta_{<\lambda}\cdot\nabla(L_{<\lambda}f) -\nabla\eta_{<\lambda}\cdot\nabla L_{<\lambda}(f^2)
\\
&-f\Delta f+\frac12\Delta(f^2) +\left( fL_{<\lambda}f-\frac12L_{<\lambda}(f^2) \right)\Delta\eta_{<\lambda}.
\end{split}
\end{equation*}
Moreover, since
\begin{equation*}
\frac12\Delta(f^2) = f\Delta f+|\nabla f|^2
\end{equation*}
and
\begin{equation*}
\begin{split}
2\nabla\eta_{<\lambda}\cdot\nabla Q_a(f,f) &= 2(L_{<\lambda}f)\nabla\eta_{<\lambda}\cdot\nabla f +2f\nabla\eta_{<\lambda}\cdot\nabla(L_{<\lambda}f) 
\\
&-\nabla\eta_{<\lambda}\cdot\nabla L_{<\lambda}(f^2),
\end{split}
\end{equation*}
we conclude that
\begin{equation*}
\begin{split}
&fL_{<\lambda}^2f-\frac12L_{<\lambda}^2(f^2)
\\
&= \frac{1}{1+|\nabla\eta_{<\lambda}|^2}\Big( 2\nabla\eta_{<\lambda}\cdot\nabla Q_a(f,f) -2(L_{<\lambda}f)\nabla\eta_{<\lambda}\cdot\nabla f +|\nabla f|^2 +Q_a(f,f)\Delta\eta_{<\lambda} \Big).
\end{split}
\end{equation*}
Adding $(L_{<\lambda}f)^2$, we have
\begin{equation*}
\begin{split}
&|\nabla f|^2 -2(L_{<\lambda}f)\nabla\eta_{<\lambda}\cdot\nabla f +(1+|\nabla\eta_{<\lambda}|^2)(L_{<\lambda}f)^2
\\
&\qquad= \left|\nabla f-(L_{<\lambda}f)\nabla\eta_{<\lambda}\right|^2 +(L_{<\lambda}f)^2
\\
&\qquad= |M_{<\lambda}f|^2+(L_{<\lambda}f)^2.
\end{split}
\end{equation*}
Consequently,
\begin{equation}\label{Qffrel}
\begin{split}
&(1+|\nabla\eta_{<\lambda}|^2)\left( fL_{<\lambda}^2f-\frac12L_{<\lambda}^2(f^2)+(L_{<\lambda}f)^2 \right)
\\
&= 2\nabla\eta_{<\lambda}\cdot\nabla Q_a(f,f) +Q_a(f,f)\Delta\eta_{<\lambda} +(L_{<\lambda}f)^2+|M_{<\lambda}f|^2.
\end{split}
\end{equation}
On the other hand, applying \eqref{L2formula} with $f=\eta_{<\lambda}$ yields
\begin{equation}\label{L2etarelation}
-(1+|\nabla\eta_{<\lambda}|^2)L_{<\lambda}^2\eta_{<\lambda} = -2\nabla\eta_{<\lambda}\cdot\nabla(L_{<\lambda}\eta_{<\lambda}) +\Delta\eta_{<\lambda} -(L_{<\lambda}\eta_{<\lambda})\Delta\eta_{<\lambda}.
\end{equation}
Using $\mathcal{S}_{<\lambda}=Q_a(V_{<\lambda},V_{<\lambda})+Q_a(B_{<\lambda},B_{<\lambda})$, the identity $a=1-L\eta+\mathcal S$ and also \eqref{Qffrel} and \eqref{L2etarelation}, we find that the second and third line of \eqref{dlambdaa} collapses to
\begin{equation*}
\frac{\Gamma_{<\lambda}}{1+|\nabla\eta_{<\lambda}|^2}\eta_{\lambda}
\end{equation*}
as desired.
\end{proof}
Next, we show that the high frequency tails of $a_{<\lambda}$ satisfy good decay estimates in various norms. The mechanism for this is similar to Proposition \ref{Bernsteinelliptic}.
\begin{lemma}\label{hightailabound} Let $\lambda>0$, $\sigma\geq0$ and $1\leq p\leq\infty$. Under the standing bootstrap assumptions in Section \ref{sec:bootstrap} and \eqref{Taylorcoeffsmall}, for $M$ sufficiently large depending on $s_c$, there holds
\begin{equation}\label{atailbound}
\||D|^{\sigma}P_{>M\lambda}a_{<\lambda}\|_{L_t^pL_x^2}\lesssim \lambda^{\sigma+1-s_c-\frac{1}{2p}}\epsilon^5c_\lambda.
\end{equation}
\end{lemma}
\begin{proof}
By interpolation, we only need to prove this for $p=1,\infty$. We give the details for $p=1$ and leave the straightforward modifications to the reader. By Bernstein, it suffices to consider the case $\sigma\geq1$. Our strategy is to simply estimate individually each term on the right hand side in Lemma \ref{asimpleidentity}.
Throughout the proof, we will use the following consequences of the bootstrap bounds, Proposition \ref{basicopbounds} and
Proposition \ref{VBunbil}:
\begin{equation}\label{a-tail-bootstrap}
\begin{split}
\|\nabla_x^2\eta_{<\lambda}\|_{L_t^4\dot B_{r,1}^0}
&\lesssim_r\epsilon\lambda^{\frac78-\frac2r},
\qquad 8<r\leq\infty,
\\
\||D|L_{<\lambda}\eta_{<\lambda}\|_{L_t^\infty\dot B_{\infty,1}^0}
&\lesssim\epsilon\lambda,
\\
\||D|^{\frac12}V_{<\lambda}\|_{L_t^\infty\dot B_{\infty,1}^0}
&\lesssim\epsilon,
\qquad
\||D|^{\frac34}V_{<\lambda}\|_{L_t^4\dot B_{\infty,1}^0}
\lesssim\epsilon\lambda^{\frac18}.
\end{split}
\end{equation}
For the last bound, the sum over
frequencies $\lesssim\lambda$ is controlled by
\ref{VBunbil}; the remaining frequencies are absorbed using
Proposition \ref{Bernsteinelliptic} with $N=1$ and
$\lambda^{-1}\|\nabla_x^2\eta_{<\lambda}\|_{L_t^\infty\dot B_{\infty,1}^0}
\lesssim\epsilon$. Now, fix an integer $N>s_c+4$ and set
\begin{equation*}
H_\lambda\coloneq\|\nabla_x^2\eta_{<\lambda}\|_{L_t^\infty\dot B_{2N,1}^0}.
\end{equation*}
By Bernstein and interpolation, we have
\begin{equation}\label{a-tail-tame-curvature}
H_\lambda^m\lesssim\epsilon^m\sum_{\kappa\lesssim\lambda}c_\kappa\kappa^{m(1-\frac1N)+2-s_c}\lesssim\epsilon^mc_\lambda\lambda^{m(1-\frac1N)+2-s_c}
\end{equation}
Our choice of $N$ ensures that the sum converges. Now, by Proposition \ref{Bernsteinelliptic}, \eqref{a-tail-bootstrap} and \eqref{a-tail-tame-curvature}, we have
\begin{equation*}
\begin{split}
\|P_{>M\lambda}|D|^{\sigma}L_{<\lambda}\eta_{<\lambda}\|_{L_t^1L_x^2}
&\lesssim \lambda^{\sigma-1-N}\|\nabla^2\eta_{<\lambda}\|^4_{L_t^4\dot{B}_{2N,1}^0}H_\lambda^{N-4}\||D|L_{<\lambda}\eta_{<\lambda}\|_{L_t^{\infty}\dot{B}_{\infty,1}^0}
\\
&\lesssim \lambda^{\sigma+\frac{1}{2}-s_c}\epsilon^5c_\lambda.
\end{split}
\end{equation*}
Next, we study the terms $Q_a(V_{<\lambda},V_{<\lambda})$ and $Q_a(B_{<\lambda},B_{<\lambda})$. We show the details for the first term, as the latter is handled in an identical fashion. We begin by recording some technical estimates that we will use later. To compactify notation, we write $v=V_{<\lambda}$ and $v_\alpha=P_\alpha v$ for the remainder of
the proof. By Proposition \ref{Bernsteinelliptic}, with $N$ as above,
\begin{equation}\label{a-tail-v}
\begin{split}
\||D|^\beta P_{>K_0\lambda}v\|_{L_t^1\dot B_{2,1}^0}
&\lesssim
\lambda^{\beta-\frac34-N}
\|\nabla_x^2\eta_{<\lambda}\|_{L_t^4\dot B_{2N,1}^0}^{3}H_\lambda^{N-3}
\||D|^{\frac34}B_{<\lambda}\|_{L_t^4\dot B_{\infty,1}^0}
\\
&\lesssim\epsilon^4c_\lambda\lambda^{\beta-s_c}.
\end{split}
\end{equation}
We first apply the proposition for $\beta\geq\frac34$; the same
bound for $0\leq\beta<\frac34$ follows by Bernstein.
Here $K_0$ is a sufficiently large fixed constant. We also have
\begin{equation}\label{a-tail-v-basic}
\|\nabla_xv\|_{L_t^\infty\dot B_{\infty,1}^0}
\lesssim\epsilon\lambda^{\frac12},
\qquad
\|v_\alpha\|_{L_t^\infty L_x^\infty}
\lesssim\epsilon\alpha^{-\frac12}.
\end{equation}
The first estimate follows by splitting at $K_0\lambda$, using
\eqref{a-tail-bootstrap} on the low frequencies and Proposition
\ref{Bernsteinelliptic} on the remaining frequencies.

Now fix a sufficiently small parameter $\delta>0$. A Littlewood-Paley decomposition in the variable $v$ gives the identity
\begin{equation}\label{a-tail-Q-decomposition}
\begin{split}
Q_a(v,v)
&=\sum_\alpha\left(
v_\alpha L_{<\lambda}v_{<\delta\alpha}
-[L_{<\lambda},v_{<\delta\alpha}]v_\alpha
+R_\alpha\right),
\\
R_\alpha
&\coloneq
\sum_{\delta\alpha\leq\nu\leq\delta^{-1}\alpha}
\left(v_\nu L_{<\lambda}v_\alpha
-\frac12L_{<\lambda}(v_\nu v_\alpha)\right).
\end{split}
\end{equation}
That is, we only isolate the commutator structure for the dangerous low-high interaction, and expand it otherwise. Here we write $v_{<\delta\alpha}=\sum_{\nu<\delta\alpha}v_\nu$. Now, we choose fixed constants
$K$ and $M$ so that $K\gg_\delta K_0$ and $M\gg_\delta K$.
\\

We first estimate the commutator terms. For $\alpha<K\lambda$, we
apply Proposition \ref{highfreqtailcomm} with $N$ as above,
$\theta=\frac12$, $p=r=2$ and $q_1=q_2=\infty$ at
$K\lambda$. We obtain
\begin{equation}\label{a-tail-low-commutator}
\begin{split}
&\sum_{\alpha<K\lambda}
\||D|^\sigma P_{>M\lambda}
[L_{<\lambda},v_{<\delta\alpha}]v_\alpha
\|_{L_t^1L_x^2}
\\
&\qquad\lesssim
\lambda^{\sigma-N}
\|\nabla_x^2\eta_{<\lambda}\|_{L_t^4\dot B_{2N,1}^0}^{4}H_\lambda^{N-4}
\||D|^{\frac12}v\|_{L_t^\infty\dot B_{\infty,1}^0}^{2}
\lesssim\epsilon^6c_\lambda\lambda^{\sigma+\frac12-s_c}.
\end{split}
\end{equation}
For $\alpha\geq K\lambda$, Proposition
\ref{highfreqcommutator}, applied at scale $\alpha$, together with
\eqref{a-tail-v} and \eqref{a-tail-v-basic}, gives
\begin{equation}\label{a-tail-high-commutator}
\begin{split}
\sum_{\alpha\geq K\lambda}
\||D|^\sigma
[L_{<\lambda},v_{<\delta\alpha}]v_\alpha
\|_{L_t^1L_x^2}\lesssim
\|\nabla_xv\|_{L_t^\infty\dot B_{\infty,1}^0}
\sum_{\alpha\geq K\lambda}
\alpha^\sigma\|v_\alpha\|_{L_t^1L_x^2}
\lesssim\epsilon^5c_\lambda\lambda^{\sigma+\frac12-s_c}.
\end{split}
\end{equation}

It remains to estimate the other two terms in
\eqref{a-tail-Q-decomposition}. We first consider $\alpha<K\lambda$.
At each fixed time, Proposition \ref{Bernsteinelliptic} with $N$ as above
and $\sigma_1=0$ gives, for $\nu<\delta\alpha$,
\begin{equation*}
\begin{split}
\||D|^\sigma P_{>M\lambda}
(v_\alpha L_{<\lambda}v_\nu)\|_{L_x^2}\lesssim
\lambda^{\sigma-N}
\|\nabla_x^2\eta_{<\lambda}\|_{\dot B_{2N,1}^0}^{N}
\nu\|v_\nu\|_{L_x^\infty}\|v_\alpha\|_{L_x^\infty}.
\end{split}
\end{equation*}
Indeed, since $v_\alpha$ has frequency below $K\lambda$, only
frequencies $\gtrsim M\lambda$ of $L_{<\lambda}v_\nu$ can
contribute. We apply the elliptic tail bound to this factor in $L_x^2$ and use
$\|L_{<\lambda}v_\nu\|_{\dot B_{\infty,1}^0}
\lesssim\nu\|v_\nu\|_{L_x^\infty}$.
For $\delta\alpha\leq\nu\leq\delta^{-1}\alpha$, the same argument
for $v_\nu L_{<\lambda}v_\alpha$, and a direct application of the
elliptic tail bound to $L_{<\lambda}(v_\nu v_\alpha)$, give
\begin{equation*}
\begin{split}
&\||D|^\sigma P_{>M\lambda}
(v_\nu L_{<\lambda}v_\alpha)\|_{L_x^2}
+\||D|^\sigma P_{>M\lambda}
L_{<\lambda}(v_\nu v_\alpha)\|_{L_x^2}
\\
&\qquad\lesssim_\delta
\lambda^{\sigma-N}
\|\nabla_x^2\eta_{<\lambda}\|_{\dot B_{2N,1}^0}^{N}
\alpha\|v_\nu\|_{L_x^\infty}\|v_\alpha\|_{L_x^\infty}.
\end{split}
\end{equation*}
Here we used
\begin{equation*}
\||D|(v_\nu v_\alpha)\|_{\dot B_{\infty,1}^0}
\lesssim_\delta
\alpha\|v_\nu\|_{L_x^\infty}\|v_\alpha\|_{L_x^\infty},
\end{equation*}
which follows by summing the output frequencies $\lesssim\alpha$, keeping the derivative on the product. We then obtain
\begin{equation}\label{a-tail-low-products}
\begin{split}
&\sum_{\alpha<K\lambda}
\left(
\||D|^\sigma P_{>M\lambda}
(v_\alpha L_{<\lambda}v_{<\delta\alpha})\|_{L_t^1L_x^2}
+\||D|^\sigma P_{>M\lambda}R_\alpha\|_{L_t^1L_x^2}
\right)
\\
&\qquad\lesssim_\delta
\lambda^{\sigma-N}
\|\nabla_x^2\eta_{<\lambda}\|_{L_t^4\dot B_{2N,1}^0}^{4}H_\lambda^{N-4}
\||D|^{\frac12}v\|_{L_t^\infty\dot B_{\infty,1}^0}^{2}
\lesssim\epsilon^6c_\lambda\lambda^{\sigma+\frac12-s_c}.
\end{split}
\end{equation}
For $\alpha\geq K\lambda$, Proposition \ref{basicopbounds} and
\eqref{a-tail-v-basic} give
\begin{equation*}
\||D|^\beta L_{<\lambda}v_{<\delta\alpha}\|_{L_t^\infty\dot B_{\infty,1}^0}
\lesssim_{\beta}
\alpha^{\beta}\|\nabla_xv\|_{L_t^\infty\dot B_{\infty,1}^0}
\lesssim_\beta\epsilon\alpha^{\beta}\lambda^{\frac12},
\qquad \beta=0,\sigma.
\end{equation*}
Hence, by Proposition \ref{besovparabounds}, we have
\begin{equation}\label{QV1}
\begin{split}
\sum_{\alpha\geq K\lambda}
\||D|^\sigma(v_\alpha L_{<\lambda}v_{<\delta\alpha})
\|_{L_t^1L_x^2}\lesssim
\epsilon\lambda^{\frac12}
\sum_{\alpha\geq K\lambda}
\alpha^\sigma\|v_\alpha\|_{L_t^1L_x^2}\lesssim\epsilon^5c_\lambda\lambda^{\sigma+\frac12-s_c}.
\end{split}
\end{equation}
Finally, for the comparable-frequency terms, Proposition
\ref{basicopbounds} gives
\begin{equation*}
\begin{split}
\||D|^\sigma(v_\nu L_{<\lambda}v_\alpha)\|_{L_x^2}
+\||D|^\sigma L_{<\lambda}(v_\nu v_\alpha)\|_{L_x^2}\lesssim_{\sigma,\delta}
\alpha^{\sigma+1}
\|v_\nu\|_{L_x^\infty}\|v_\alpha\|_{L_x^2}.
\end{split}
\end{equation*}
For the second term, we again sum the output frequencies of the
product with a positive derivative before applying the operator
bound. Since both input frequencies are above $K_0\lambda$, placing
one factor in $L_t^1L_x^2$ and the other in $L_t^\infty L_x^\infty$ gives
\begin{equation}\label{QV2}
\begin{split}
\sum_{\alpha\geq K\lambda}
\||D|^\sigma R_\alpha\|_{L_t^1L_x^2}\lesssim_{\sigma,\delta}
\epsilon
\sum_{\nu\gtrsim_\delta K\lambda}
\nu^{\sigma+\frac12}\|v_\nu\|_{L_t^1L_x^2}\lesssim\epsilon^5c_\lambda\lambda^{\sigma+\frac12-s_c}.
\end{split}
\end{equation}
Combining \eqref{a-tail-low-commutator},
\eqref{a-tail-high-commutator}, \eqref{a-tail-low-products},
\eqref{QV1}, and \eqref{QV2}, we obtain
\begin{equation*}
\||D|^\sigma P_{>M\lambda}Q_a(V_{<\lambda},V_{<\lambda})
\|_{L_t^1L_x^2}
\lesssim\epsilon^5c_\lambda\lambda^{\sigma+\frac12-s_c}.
\end{equation*}
Doing the analogous estimate for $Q_a(B_{<\lambda},B_{<\lambda})$ proves the $p=1$ bound.
The $L_t^\infty L_x^2$ bound is similar and easier. We leave the details to the reader.
\end{proof}
Next, we establish linear and bilinear estimates for the Taylor coefficient $a$. The proof of this is analogous to Proposition \ref{VBunbil} which established the corresponding estimates for $V$ and $B$.
\begin{proposition}[Bilinear estimates for $a$]\label{abilinearbounds}
Assume that $u$ satisfies the bootstrap bounds in Section \ref{sec:bootstrap}, as well as \eqref{Taylorcoeffsmall}. Fix $x_0\in\mathbb{R}^2$. Then $a_{\mu}=P_{\mu}a$ satisfies the bounds
\begin{enumerate}
\item ($L_t^{\infty}L_x^2$ bound)
\begin{equation}\label{aL2}
\|a_{\mu}\|_{L_t^{\infty}L_x^2}\lesssim \epsilon c_{\mu}\mu^{1-s_c}.
\end{equation}
\item ($L_t^{4}L_x^p$ bound)
\begin{equation}\label{a48}
\|a_{\mu}\|_{L_t^4L_x^{p}}\lesssim \epsilon c_{\mu}\mu^{\frac{13}{8}+2(\frac{1}{8}-\frac{1}{p})-s_c},\hspace{5mm}8<p\leq \infty.
\end{equation}
\item (Bilinear low-high bound)
\begin{equation}\label{aulh}
\|a_{\mu}u_{\lambda}^{x_0}\|_{L_{t,x}^2}\lesssim c_{\lambda}c_{\mu}\epsilon^2\lambda^{-s_c}\mu^{\frac{7}{4}-s_c}(1+|x_0|\mu)^{\frac12},\hspace{5mm}\mu\ll\lambda.
\end{equation}
\item (Bilinear high-low bound)
\begin{equation}\label{auhl}
\|a_{\mu}u_{\lambda}^{x_0}\|_{L_{t,x}^2}\lesssim c_{\lambda}c_{\mu}\epsilon^2\lambda^{\frac{3}{4}-s_c}\mu^{1-s_c}(1+|x_0|\lambda)^{\frac12},\hspace{5mm}\mu\gg\lambda.
\end{equation}
\item (Translation-averaged bilinear $L^2$ bound)
For $\frac12<\sigma<1$, there holds
\begin{equation}\label{aulhtrans}
M_{\mu,2}^{\sigma}(u_\lambda^{x_0},a_\mu)
\lesssim_\sigma
\epsilon^2c_\lambda c_\mu\lambda^{-s_c}\mu^{\frac74-s_c}(1+\mu|x_0|)^{\frac12},
\hspace{5mm}\mu\ll\lambda.
\end{equation}
\end{enumerate}
Moreover, all of the above estimates hold with $a_\mu$ replaced by $P_\mu a_{<\alpha}$, for any $\alpha>0$, with implicit constants uniform in $\alpha$.
\end{proposition}

\begin{proof}
The proof strategy here is similar to \ref{VBunbil}, except one uses the
linearization formula for $a$ in place of the one for $B$. We use the notation $\Phi^0_{<\nu}$ from
\eqref{Phi0definition}.

It will be convenient to  bootstrap, uniformly in $\alpha>0$, the bound
\begin{equation}\label{a-auxiliary-bootstrap}
\|a_{<\alpha}-1\|_Y
\coloneq
\||D|a_{<\alpha}\|_{L_t^\infty\dot B_{2,1}^0}
+\||D|^{\frac14}(a_{<\alpha}-1)\|_{L_t^4\dot B_{16,1}^0}
\leq K\epsilon.
\end{equation}
In particular, $\|a_{<\alpha}-1\|_{L_{t,x}^\infty}\lesssim K\epsilon$.
We take $K$ large and subsequently take $\epsilon$ sufficiently
small relative to $K$. At the end of the proof, we show that the bound here holds with an improved constant. This in turn will close the bootstrap \eqref{Taylorcoeffsmall}. As in the proof of Proposition
\ref{VBunbil}, the bootstrap is apriori initially made on a short time
interval and is closed below.

We first record the versions of the high-frequency tail estimate
needed in this proof. For $M$ sufficiently large and $j\geq0$,
\begin{equation}\label{a-tail-remainder-norms}
\begin{split}
\||D|^jP_{>M\nu}a_{<\nu}\|_{L_t^2\dot B_{\infty,1}^0}
&\lesssim \epsilon^2c_\nu\nu^{j+\frac74-s_c},
\\
\||D|^jP_{>M\nu}a_{<\nu}\|_{L_t^4\dot B_{\frac{16}{7},1}^0}
+\||D|^jP_{>M\nu}a_{<\nu}\|_{L_t^\infty\dot B_{2,1}^0}
&\lesssim \epsilon^2c_\nu\nu^{j+1-s_c}.
\end{split}
\end{equation}
These follow from Lemma \ref{hightailabound} by interpolation and
Bernstein's inequality to slightly enlarge the number of derivatives applied (to enable summation over the output frequencies). Combining these tail bounds applied at $s=2$ with
\eqref{a-auxiliary-bootstrap} and Bernstein gives
\begin{equation}\label{a-coefficient-bounds}
\begin{split}
\|\nabla_xa_{<\nu}\|_{L_t^4\dot B_{16,1}^0}&\lesssim_K\epsilon\nu^{\frac34},\hspace{5mm}\|\nabla_xa_{<\nu}\|_{L_t^4\dot B_{\infty,1}^0}
\lesssim_K\epsilon\nu^{\frac78},
\\
\|\nabla_xa_{<\nu}\|_{L_t^\infty\dot B_{\frac83,1}^0}
&\lesssim_K\epsilon\nu^{\frac14},\hspace{5mm}\|\nabla_xa_{<\nu}\|_{L_t^\infty\dot B_{\infty,1}^0}
\lesssim_K\epsilon\nu.
\end{split}
\end{equation}
For the low frequencies these are direct Bernstein estimates.
For the high frequencies, the second line of
\eqref{a-tail-remainder-norms} along with Bernstein suffices. 
\\

Our next goal will be to  estimate the terms in the linearization \eqref{linearizedaformula}, except for the "principally linear" part
$-L_{<\nu}(a_{<\nu}\eta_\nu)$. Let $\mathcal E_\nu$ denote their
sum, so that
\begin{equation*}
\begin{split}
\mathcal E_\nu
={}&\frac{\Gamma_{<\nu}}{1+|\nabla_x\eta_{<\nu}|^2}\eta_\nu
+(L_{<\nu}\Phi^0_{<\nu})L_{<\nu}B_{<\nu}
-[L_{<\nu},B_{<\nu}](L_{<\nu}\Phi^0_{<\nu})
\\
&+(M_{<\nu}\Phi^0_{<\nu})\cdot L_{<\nu}V_{<\nu}
-[L_{<\nu},V_{<\nu}]\cdot(M_{<\nu}\Phi^0_{<\nu}).
\end{split}
\end{equation*}
We claim that each of these terms satisfies the low frequency bounds
\begin{equation}\label{a-source-remainder-bounds}
\begin{split}
\|P_{<M\nu}\mathcal E_\nu\|_{L_t^\infty L_x^2}
+\|P_{<M\nu}\mathcal E_\nu\|_{L_t^4L_x^{\frac{16}{7}}}
&\lesssim_K\epsilon^2c_\nu\nu^{1-s_c},
\\
\|P_{<M\nu}\mathcal E_\nu\|_{L_t^2L_x^\infty}
&\lesssim_K\epsilon^2c_\nu\nu^{\frac74-s_c}.
\end{split}
\end{equation}
Indeed, Proposition \ref{VBunbil}, Corollary \ref{Bestcor}, and
\eqref{a-coefficient-bounds} give
\begin{equation*}
\|\Gamma_{<\nu}\|_{L_t^\infty L_x^\infty}
\lesssim_K\epsilon\nu,
\qquad
\|\Gamma_{<\nu}\|_{L_t^4L_x^\infty}
\lesssim_K\epsilon\nu^{\frac78}.
\end{equation*}
Pairing these with $\eta_\nu$ in $L_t^\infty L_x^2$,
$L_t^\infty L_x^{16/7}$, and $L_t^4L_x^\infty$, respectively,
proves the claim for the first term in $\mathcal{E}_{\nu}$. 
\\

For the non-commutator product terms in $\mathcal{E}_\nu$, write $U_\nu=L_{<\nu}\Phi^0_{<\nu}$
or $M_{<\nu}\Phi^0_{<\nu}$ and $v=B_{<\nu}$ or $V_{<\nu}$.
Lemma \ref{Phiest} gives the energy and $L_t^2L_x^\infty$ bounds
directly. For the remaining norm, Bernstein and Proposition \ref{basicopbounds}, Proposition \ref{VBunbil} and Proposition \ref{Bernsteinelliptic} gives
\begin{equation*}
\begin{split}
\|P_{<M\nu}(U_\nu L_{<\nu}v)\|_{L_t^4L_x^{\frac{16}{7}}}
&\lesssim
\nu^{\frac14}\|U_\nu\|_{L_t^\infty L_x^2}
\|L_{<\nu}v\|_{L_t^4L_x^{16}}
\\
&\lesssim
\nu^{\frac14}
(\epsilon c_\nu\nu^{\frac12-s_c})
(\epsilon\nu^{\frac14})
=\epsilon^2c_\nu\nu^{1-s_c}.
\end{split}
\end{equation*}
For the commutator terms, we need the estimate
\begin{equation}\label{a-zero-order-commutator}
\|[L_{<\nu},v]U\|_{L_x^q}
\lesssim
\|\nabla_xv\|_{\dot B_{q_1,1}^0}
\|U\|_{\dot B_{q_2,1}^0},
\qquad \frac1q=\frac1{q_1}+\frac1{q_2},
\end{equation}
where $1<q<\infty$ and $1<q_1,q_2\leq\infty$.
To see this, the sum over the low-high $(v,U)$ interactions can be estimated using
\eqref{similarestimatescomm}. For the other interactions, where the undifferentiated term is at higher frequency, we directly expand the commutator and estimate
$v_\kappa L_{<\nu}U_\alpha$  directly and use
\begin{equation}\label{L(vU)}
\|L_{<\nu}(v_\kappa U_\alpha)\|_{L_x^q}
\lesssim
\|v_\kappa U_\alpha\|_{\dot B_{q,1}^1}
\lesssim
\kappa\|v_\kappa\|_{L_x^{q_1}}
\|U_\alpha\|_{L_x^{q_2}}
\end{equation}
when $\kappa\gg\alpha$. For high-high to low interactions, we use the operator $L_{<\nu}$ in \eqref{L(vU)} and Bernstein to sum over the low output frequency. Dyadic summation then proves
\eqref{a-zero-order-commutator}. Next, we apply this bound with
\begin{equation*}
(q,q_1,q_2)=(8,16,16),\quad(16/9,16,2),\quad(2,\infty,2)
\end{equation*}
and Bernstein and Proposition \ref{basicopbounds} for the contribution of the commutator terms in the three estimates in \eqref{a-source-remainder-bounds}. This proves \eqref{a-source-remainder-bounds}.

We now isolate the principal contribution. Fix a sufficiently
small $\delta>0$, with $M\delta\ll1$. We then have the identity
\begin{equation}\label{Laetaexp}
\begin{split}
L_{<\nu}(a_{<\nu}\eta_\nu)
&=(P_{<\delta\nu}a_{<\nu}) L_{<\nu}\eta_\nu
+[L_{<\nu},P_{<\delta\nu}a_{<\nu}]\eta_\nu
\\
&+L_{<\nu}((P_{\geq\delta\nu}a_{<\nu})\eta_\nu).
\end{split}
\end{equation}
Arguing as above, the latter two terms can be grouped with $\mathcal{E}_{\nu}$ and satisfy the same estimates as
\eqref{a-source-remainder-bounds}. For the commutator, this follows by using
\eqref{a-zero-order-commutator} and
\eqref{a-coefficient-bounds}. For the other term, we use Proposition
\ref{basicopbounds}, proposition \ref{besovparabounds},
and
\begin{equation*}
\|P_{\geq\delta\nu}a_{<\nu}\|_{\dot B_{q,1}^0}
\lesssim_\delta \nu^{-1}
\|\nabla_xa_{<\nu}\|_{\dot B_{q,1}^0}.
\end{equation*}
More explicitly, the support of $P_{\geq\delta\nu}a_{<\nu}$ gives
\begin{equation*}
\|(P_{\geq\delta\nu}a_{<\nu})\eta_\nu\|_{\dot B_{q,1}^1}
\lesssim_\delta
\|\nabla_xa_{<\nu}\|_{\dot B_{q_1,1}^0}
\|\eta_\nu\|_{L_x^{q_2}},
\qquad \frac1q=\frac1{q_1}+\frac1{q_2}
\end{equation*}
for a suitable set of exponents (depending on which of the three norms is being estimated), along with Bernstein.

The proof of Lemma \ref{princextraction}, with $\eta_\nu$ in
place of $\varphi_\nu$, controls the remaining difference
between $P_{<M\nu}((P_{<\delta\nu}a_{<\nu}) L_{<\nu}\eta_\nu)$ and
$(P_{<\delta\nu}a_{<\nu}) Op(\mathfrak l^l_{<\delta\nu})\eta_\nu$ in the same
three norms. More explicitly, the analogues of the first and third error terms in this lemma contain
a factor $\nabla_x^2\eta_{<\nu}$ and a factor $\eta_\nu$;
the second error contains
$\nabla_x(\eta_{<\nu}-\eta_{<\delta\nu})$ and
$\nabla_x\eta_\nu$. For the $L_t^4L_x^{16/7}$ estimate use
\begin{equation*}
\begin{split}
\|\nabla_x^2\eta_{<\nu}\|_{L_t^\infty\dot B_{8/3,1}^0}
&\lesssim\epsilon\nu^{\frac14},
\\
\|\nabla_x(\eta_{<\nu}-\eta_{<\delta\nu})\|_{L_t^\infty\dot B_{8/3,1}^0}
&\lesssim_\delta\epsilon\nu^{-\frac34},
\\
\|\eta_\nu\|_{L_t^4L_x^{16}}
&\lesssim\epsilon c_\nu\nu^{\frac34-s_c}.
\end{split}
\end{equation*}
These gives the required  $\epsilon^2c_\nu\nu^{1-s_c}$ factor. For the $L_t^{\infty}L_x^2$ bound, we place the extra coefficient in $L_{t,x}^\infty$ and $\eta_\nu$ in
$L_t^\infty L_x^2$. For the $L_t^2L_x^\infty$ bound, we put both factors
in $L_t^4L_x^{16}$ after applying Bernstein to the product; the result is
$\epsilon^2c_\nu\nu^{7/4-s_c}$. For the remaining term on the right hand side of \eqref{Laetaexp}, we can use Proposition \ref{basicopbounds} and \eqref{a-auxiliary-bootstrap} and Bernstein to replace $(P_{<\delta\nu}a_{<\nu})L_{<\nu}\eta_{\nu}$ by $(P_{<\delta\nu}a_{<\nu})Op(\mathfrak{l}^l_{<\delta\nu})\eta_{\nu}$ again modulo errors satisfying the same remainder bounds as in \eqref{a-source-remainder-bounds}.

With this in hand, we are now ready to close \eqref{a-auxiliary-bootstrap} and establish the bounds in Proposition \ref{abilinearbounds}. Indeed, the fundamental theorem of calculus and \eqref{linearizedaformula} now gives
\begin{equation}\label{a-integral-principal}
\begin{split}
a-1
&=\int_0^\infty
\left[P_{<M\nu}
\bigl(-L_{<\nu}(a_{<\nu}\eta_\nu)+\mathcal E_\nu\bigr)
+\widetilde P_{M\nu}(a_{<\nu}-1)\right]\frac{d\nu}{\nu}
\\
&=\int_0^\infty\left(-(P_{<\delta\nu}a_{<\nu})Op(\mathfrak l^l_{<\delta\nu})\eta_\nu+\mathcal R_\nu\right)\frac{d\nu}{\nu}.
\end{split}
\end{equation}
Here $\widetilde P_{M\nu}=\nu\frac{d}{d\nu}P_{<M\nu}$, and $\mathcal R_\nu$ is the sum of all the errors just estimated above.
It has output frequency $\lesssim_M\nu$, and
\eqref{a-tail-remainder-norms} controls the term involving $\tilde{P}_{M\nu}(a_{<\nu}-1)$. Therefore, setting
$a_{r,\mu}=\int_0^\infty P_\mu\mathcal R_\nu\,\frac{d\nu}{\nu}$ and arguing as in Proposition \ref{VBunbil}, we obtain
\begin{equation}\label{a-remainder-bounds}
\begin{split}
\|a_{r,\mu}\|_{L_t^\infty L_x^2}
+\|a_{r,\mu}\|_{L_t^4L_x^{\frac{16}{7}}}
&\lesssim_K\epsilon^2c_\mu\mu^{1-s_c},
\\
\|a_{r,\mu}\|_{L_t^2L_x^\infty}
&\lesssim_K\epsilon^2c_\mu\mu^{\frac74-s_c}.
\end{split}
\end{equation}
Indeed, only $\nu\gtrsim_M\mu$ contributes, and both powers
$1-s_c$ and $7/4-s_c$ are negative for $s_c\geq2$.

The principal integrand $-(P_{<\delta\nu}a_{<\nu})Op(\mathfrak l^l_{<\delta\nu})\eta_\nu$ is localized at output frequency $\nu$. Using \eqref{a-auxiliary-bootstrap}, and taking $\epsilon$ small enough relative to $K$, we have
\begin{equation*}
\begin{split}
\|(P_{<\delta\nu}a_{<\nu})Op(\mathfrak l^l_{<\delta\nu})\eta_\nu\|_{L_t^\infty L_x^2}
&\lesssim\epsilon c_\nu\nu^{1-s_c},
\\
\|(P_{<\delta\nu}a_{<\nu})Op(\mathfrak l^l_{<\delta\nu})\eta_\nu\|_{L_t^4L_x^p}
&\lesssim_p\epsilon c_\nu\nu^{\frac{15}{8}-\frac2p-s_c},
\quad 8<p\leq\infty.
\end{split}
\end{equation*}
We note importantly that the implicit constants here are independent of $K$. Consequently \eqref{a-integral-principal},
\eqref{a-remainder-bounds}, and Bernstein from $L_x^{16/7}$
prove \eqref{aL2} and \eqref{a48} with $K$ independent constants.

We next justify that the above bounds also hold uniformly in the geometric truncation parameter and also close
\eqref{a-auxiliary-bootstrap}. The truncated variant of the above identity is
\begin{equation}\label{a-truncated-integral}
\begin{split}
a_{<\alpha}-1
=\int_0^\alpha\left(-(P_{<\delta\nu}a_{<\nu})Op(\mathfrak l^l_{<\delta\nu})\eta_\nu+\mathcal R_\nu\right)\frac{d\nu}{\nu}+P_{\geq M\alpha}(a_{<\alpha}-1).
\end{split}
\end{equation}
The same estimates apply to the shortened integral. The last
term is controlled by \eqref{a-tail-remainder-norms}. Finally, to close the bootstrap \eqref{a-auxiliary-bootstrap}, use
the summability of the $\dot{B}_{2,1}^2$ frequency envelope, as in the proof of
\eqref{B1rsummed}. More precisely, summing the $L_t^{\infty}L_x^2$ estimate with a factor of
$\mu$ and the $L_t^4L_x^{16}$ estimate with a factor of
$\mu^{1/4}$ gives, uniformly in $\alpha$,
\begin{equation*}
\|a_{<\alpha}-1\|_Y\lesssim\epsilon+C_K\epsilon^2.
\end{equation*}
This improves the bootstrap for $K$ sufficiently large and
then $\epsilon$ sufficiently small. 
\\

Finally, we prove the bilinear bounds. The kernel argument for the principal term used in Corollary \ref{Abilbounds}, now applied
to the $(\eta,u)$ bilinear bound, gives
\begin{equation*}
\begin{split}
\|(P_{<\delta\nu}a_{<\nu})Op(\mathfrak l^l_{<\delta\nu})\eta_\nu\,u_\lambda^{x_0}\|_{L_{t,x}^2}\lesssim
\epsilon^2c_\nu c_\lambda
\nu^{1-s_c}\lambda^{-s_c}\min(\nu,\lambda)^{\frac34}
\bigl(1+\min(\nu,\lambda)|x_0|\bigr)^{\frac12}
\end{split}
\end{equation*}
when $\nu$ and $\lambda$ are separated. Here, we placed the coefficient $P_{<\delta\nu}a_{<\nu}$ in $L_t^{\infty}L_x^{\infty}$ which is bounded by \eqref{a-auxiliary-bootstrap}. Since
$P_\mu\left((P_{<\delta\nu}a_{<\nu})Op(\mathfrak l^l_{<\delta\nu})\eta_\nu\right)=0$ unless $\nu\approx\mu$, the same
estimate after applying $P_\mu$ and integrating gives the
desired principal contribution in both frequency regimes. For the remainder terms, \eqref{a-remainder-bounds} gives
\begin{equation*}
\begin{split}
\|a_{r,\mu}u_\lambda^{x_0}\|_{L_{t,x}^2}
&\lesssim
\|a_{r,\mu}\|_{L_t^2L_x^\infty}
\|u_\lambda\|_{L_t^\infty L_x^2}
\lesssim\epsilon^3c_\mu c_\lambda
\mu^{\frac74-s_c}\lambda^{-s_c},
\qquad \mu\ll\lambda,
\\
\|a_{r,\mu}u_\lambda^{x_0}\|_{L_{t,x}^2}
&\lesssim
\|a_{r,\mu}\|_{L_t^4L_x^{\frac{16}{7}}}
\|u_\lambda\|_{L_t^4L_x^{16}}
\lesssim\epsilon^3c_\mu c_\lambda
\mu^{1-s_c}\lambda^{\frac34-s_c},
\qquad \mu\gg\lambda.
\end{split}
\end{equation*}
This proves \eqref{aulh} and \eqref{auhl}. The same argument applied to \eqref{a-truncated-integral}
proves the estimates for $P_\mu a_{<\alpha}$ uniformly in
$\alpha$. Finally, we prove \eqref{aulhtrans}. This is very similar to the corresponding estimate proved in Proposition \ref{VBunbil}. Repeating the remainder
estimates above with two $L_t^4L_x^{2p}$ factors, where $p>4$,
gives for any $\mu>0$,
\begin{equation}\label{arLpbound}
\|a_{r,\mu}\|_{L_t^2L_x^p}
\lesssim_p
\epsilon^2c_\mu\mu^{\frac74-\frac2p-s_c}.
\end{equation}
Choose $4<p<\frac{2}{1-\sigma}$. Proposition
\ref{basicMproperties} then gives
\begin{equation}
\begin{split}
M_{\mu,2}^{\sigma}(u_\lambda^{x_0},a_{r,\mu})\lesssim_{\sigma,p}
\mu^{\frac2p}
\|u_\lambda\|_{L_t^\infty L_x^2}
\|a_{r,\mu}\|_{L_t^2L_x^p}\lesssim_\sigma
\epsilon^3c_\lambda c_\mu\lambda^{-s_c}\mu^{\frac74-s_c}.
\end{split}
\end{equation}
For the principal contribution in \eqref{a-integral-principal},
we use the translation-averaged bound for $\eta_\nu$ in Section \ref{sec:bootstrap} along with Lemma \ref{annularcontrol} and the
same order one kernel argument as in Proposition \ref{VBunbil}.
Indeed, $P_{<\delta\nu}a_{<\nu}$ is uniformly bounded, and only
$\nu\approx\mu$ contributes after applying $P_\mu$. The kernel
therefore costs one factor of $\mu$, giving the right-hand
side of \eqref{aulhtrans}. This proves the estimate.
The same argument applies to the geometrically truncated $a_{<\alpha}$ in \eqref{a-truncated-integral} by applying Lemma \ref{hightailabound}, Bernstein, and the same estimate as used with the remainder terms $a_{r,\mu}$ above. This completes the proof.
\end{proof}
\subsection{Error bounds for \texorpdfstring{$a$}{}} To conclude this section, we need the following error estimate that will allow us to peel off the principal part of $a$ in our refined paradifferential expansion of the water-waves system later. More specifically, we will isolate the precise form of the linear component $a^{[1]}$ (relative to suitably frequency localized inputs) which was first derived formally in Lemma \ref{l:a-exp} and then rigorously establish the two key errors estimates that we will used later. 

\begin{proposition}[Difference bound for $a$]\label{diffabound}
Assume the standing bootstrap bounds in Section \ref{sec:bootstrap}.
We have the following estimates.
\begin{enumerate}
\item Let $L$ be a translation-trilinear form with rapidly decaying kernel, normalized at the corresponding frequency scales. Let $\lambda,\mu>0$ and $M\gg1$ with $\lambda\ll M\lambda\ll\mu$. Define $a_{<\lambda}^{[1]}\coloneq -T_{\mathfrak{l}}T_{a^l}\eta_{<\lambda}$. Then there holds
\begin{equation}\label{adiffmulti}
\mu^{-\frac{1}{2}}\|L(u^{x_0}_{\mu},u_{\mu},\nabla P_{<M\lambda}(a_{<\lambda}-a_{<\lambda}^{[1]}))\|_{L_{t,x}^1}\lesssim \epsilon^4\mu^{-2s_c}c_{\mu}^2(1+\mu|x_0|)^{\frac12}.
\end{equation}
\item For every $\lambda>0$, set
$a_\lambda\coloneq\lambda\frac{d}{d\lambda}a_{<\lambda}$ and
$a_\lambda^{[1]}\coloneq-T_{a^l}T_{\mathfrak l}\eta_\lambda$.
Then there holds
\begin{equation}\label{aprincipalL2Linf}
\|a_\lambda-a_\lambda^{[1]}\|_{L_t^2L_x^\infty}
\lesssim_\delta\epsilon^2c_\lambda\lambda^{-\frac14}.
\end{equation}
\end{enumerate}
\end{proposition}
\begin{proof}
We begin with the first, and much more difficult multilinear bound \eqref{adiffmulti}. The second bound \eqref{aprincipalL2Linf} will be a simple corollary of the estimates used to obtain the first bound. We will remark on this at the end. Our starting point is to use the fundamental theorem of calculus to write
\begin{equation}\label{aaverage}
\begin{split}
a_{<\lambda}^l-1=P_{<M\lambda}(a_{<\lambda}-1)
&=\int_{0}^{\lambda}P_{<M\nu}\nu\frac{d}{d\nu}a_{<\nu}\frac{d\nu}{\nu}+\int_{0}^{\lambda}\tilde{P}_{M\nu}(a_{<\nu}-1)\frac{d\nu}{\nu}
\end{split}
\end{equation}
where $\tilde{P}_{M\nu}\coloneq \nu\frac{d}{d\nu}P_{<M\nu}$. By Lemma \ref{hightailabound}, with $M$ sufficiently large, the latter term satisfies
\begin{equation*}
\int_{0}^{\lambda}\|\nabla\tilde{P}_{M\nu}a_{<\nu}\|_{L_t^1L_x^{\infty}}\frac{d\nu}{\nu}\lesssim \epsilon^5\int_{0}^{\lambda}\nu^{\frac{1}{2}}\frac{d\nu}{\nu}\lesssim \epsilon^5\lambda^{\frac{1}{2}}.
\end{equation*}
This is sufficient to estimate the contribution of the last term in \eqref{aaverage} when inserted into \eqref{adiffmulti}. Indeed, we can use $L_t^{\infty}L_x^2$ estimates for $u_{\mu}$ and $u_{\mu}^{x_0}$. Now, we handle the main part of \eqref{aaverage}. By definition, we have
\begin{equation*}
T_{\mathfrak l}T_{a^l}\eta_{<\lambda}
=\int_0^\lambda T_{\mathfrak l}T_{a^l}\eta_\nu\frac{d\nu}{\nu}.
\end{equation*}
The integrand is supported at frequency $\approx\nu$, so $P_{<M\lambda}$ is the identity on it for $\nu\leq\lambda$. Therefore, it suffices to show
\begin{equation}\label{dlambdaaest}
\Big\|L\left(u_{\mu}^{x_0},u_{\mu},\nabla\left(P_{<M\nu}\nu\frac{d}{d\nu}a_{<\nu}+T_{\mathfrak l}T_{a^l}\eta_\nu\right)\right)\Big\|_{L_{t,x}^1}\lesssim \epsilon^4c_\mu^2\mu^{-2s_c}\nu^{\frac12}(1+\mu|x_0|)^{\frac12}
\end{equation}
for $\nu\leq\lambda$. Throughout the proof, all constants in the frequency cutoffs are fixed, $M\delta\ll1$, and $M\nu\ll\mu$. To obtain \eqref{dlambdaaest}, we begin by expanding $\nu\frac{d}{d\nu}a_{<\nu}$ more precisely using the linearization formula \eqref{linearizedaformula}. First, we remove $L_t^1L_x^{\infty}$ errors coming from high-frequency tails. More precisely, we have
\begin{equation}\label{datruncated}
\begin{split}
\nu\frac{d}{d\nu}a_{<\nu}
&=-L^l_{<\nu}(a^l_{<\nu}\eta_\nu)+\frac{\Gamma^l_{<\nu}}{1+|\nabla\eta_{<\nu}|^2}\eta_\nu
\\
&\quad+(L^l_{<\nu}\Phi^{l,0}_{<\nu})L^l_{<\nu}B^l_{<\nu}-[L_{<\nu},B^l_{<\nu}](L_{<\nu}\Phi^0_{<\nu})
\\
&\quad+(M^l_{<\nu}\Phi^{l,0}_{<\nu})\cdot L^l_{<\nu}V^l_{<\nu}-[L_{<\nu},V^l_{<\nu}]\cdot(M_{<\nu}\Phi^0_{<\nu})+O_{L_t^1L_x^{\infty}}(\epsilon^2\nu^{-\frac12}),
\end{split}
\end{equation}
where $\Phi^{l,0}_{<\nu}=\psi_\nu-B^l_{<\nu}\eta_\nu$ and
\begin{equation*}
\begin{split}
\Gamma^l_{<\nu}
&=|M^l_{<\nu}V^l_{<\nu}|^2+|M^l_{<\nu}B^l_{<\nu}|^2+|L^l_{<\nu}V^l_{<\nu}|^2+|L^l_{<\nu}B^l_{<\nu}|^2
\\
&\quad+a^l_{<\nu}\Delta\eta_{<\nu}+2\nabla\eta_{<\nu}\cdot\nabla a^l_{<\nu}.
\end{split}
\end{equation*}
The output projection defining $a^l_{<\nu}$ is the one in \eqref{aaverage}. For convenience, the other output projections defining the superscript $l$ may be fattened by fixed factors. In particular, the projection on an operator output is taken sufficiently larger than the frequencies on its inputs, so that the errors can be estimated using as high-frequency tails. To justify \eqref{datruncated}, Proposition \ref{Bernsteinelliptic}, Proposition \ref{basicopbounds} and Lemma \ref{hightailabound} give
\begin{equation*}
\begin{split}
\|L_{<\nu}(a_{<\nu}\eta_\nu)-L^l_{<\nu}(a^l_{<\nu}\eta_\nu)\|_{L_t^1L_x^\infty}
&\lesssim\epsilon^4\nu^{-\frac12},\hspace{5mm}\|A_{<\nu}v_{<\nu}-A^l_{<\nu}v^l_{<\nu}\|_{L_t^1L_x^\infty}\lesssim\epsilon^4.
\end{split}
\end{equation*}
where above, and for the remainder of the proof, we write $A\in\{L,M\}$ and $v\in\{B,V\}$. For the first estimate, we place three factors of $\nabla^2\eta_{<\nu}$ in $L_t^4\dot B_{\infty,1}^0$ and the input in $L_t^4\dot{B}_{\infty,1}^1$, and use Lemma \ref{atailbound} for the difference of the inputs. The second estimate follows in the same way using also Proposition \ref{VBunbil}. Combining it with
\begin{equation*}
\|A_{<\nu}v_{<\nu}\|_{L_t^\infty L_x^\infty}+\|A^l_{<\nu}v^l_{<\nu}\|_{L_t^\infty L_x^\infty}\lesssim\epsilon\nu^{\frac12}
\end{equation*}
gives
\begin{equation*}
\|\Gamma_{<\nu}-\Gamma^l_{<\nu}\|_{L_t^1L_x^\infty}\lesssim\epsilon^2\nu^{\frac12}.
\end{equation*}
For the terms containing $a$ in $\Gamma_{<\nu}$, we use Lemma \ref{hightailabound} again. Multiplication by $\frac{\eta_\nu}{1+|\nabla\eta_{<\nu}|^2}$ then gives the required $O_{L_t^1L_x^\infty}(\epsilon^2\nu^{-1/2})$ error.
\\

The same tail argument with two factors of $\nabla^2\eta_{<\nu}$ also gives
\begin{equation}\label{a-difference-Phi-tail}
\|A_{<\nu}\Phi^0_{<\nu}-A^l_{<\nu}\Phi^{l,0}_{<\nu}\|_{L_t^{\frac43}\dot B_{\infty,1}^0}\lesssim\epsilon^3\nu^{-\frac78}.
\end{equation}
Here $\Phi^{l,0}_{<\nu}-\Phi^0_{<\nu}=(B_{<\nu}-B^l_{<\nu})\eta_\nu$; its contribution is estimated using Proposition \ref{Bernsteinelliptic}, Corollary \ref{Bestcor}, Proposition \ref{VBunbil} and $\|\eta_\nu\|_{L_t^\infty L_x^\infty}\lesssim\epsilon\nu^{-1}$. Pairing \eqref{a-difference-Phi-tail} with $\|A_{<\nu}v_{<\nu}\|_{L_t^4\dot B_{\infty,1}^0}\lesssim\epsilon\nu^{3/8}$ controls the errors in the two product terms in the second and third line of \eqref{datruncated}. Similarly, expanding the commutators gives
\begin{equation*}
[L_{<\nu},v_{<\nu}](A_{<\nu}\Phi^0_{<\nu})
=[L_{<\nu},v^l_{<\nu}](A_{<\nu}\Phi^0_{<\nu})+O_{L_t^1L_x^\infty}(\epsilon^2\nu^{-\frac12}).
\end{equation*}
Combining the above estimates proves \eqref{datruncated}. After applying $\nabla P_{<M\nu}$ and Bernstein, the contribution of all of these errors in \eqref{dlambdaaest} are perturbative by using two $L_t^{\infty}L_x^2$ bounds for $u_\mu$ and $u_\mu^{x_0}$. Now, we estimate the contribution of the remaining terms in \eqref{datruncated}, one by one.
\subsubsection{The principal error term} Here, we aim to show that
\begin{equation}\label{principalaerr}
\|L(u_\mu^{x_0},u_\mu,\nabla P_{<M\nu}
(L^l_{<\nu}(a^l_{<\nu}\eta_\nu)
-T_{\mathfrak l}T_{a^l}\eta_\nu))\|_{L_{t,x}^1}\lesssim \epsilon^4 c_{\mu}^2\mu^{-2s_c}\nu^{\frac{1}{2}}(1+\mu|x_0|)^{\frac12}.
\end{equation}
We first record the following estimate that we will need
for every $x_0\in\mathbb R^2$:
\begin{equation}\label{transaveua}
\left\|
(a^l_{<\nu}-a^l_{<\delta\nu})u_\mu^{x_0}
\right\|_{L_{t,x}^2}
\lesssim_\delta
\epsilon^2c_\mu\mu^{-s_c}\nu^{-\frac14}
(1+\nu|x_0|)^{\frac12}.
\end{equation}
By the fundamental theorem of calculus, the linearization formula \eqref{linearizedaformula} and arguing as in the estimates following \eqref{a-truncated-integral} we obtain
\begin{equation}\label{anulin}
a^l_{<\nu}-a^l_{<\delta\nu}
=
-\int_{\delta\nu}^{\nu}
(P_{<\delta\rho}a_{<\rho})
Op(\mathfrak l^l_{<\delta\rho})\eta_\rho
\frac{d\rho}{\rho}
+\mathcal E_a,
\end{equation}
where
\begin{equation*}
\|\mathcal E_a\|_{L_t^2L_x^p}
\lesssim_{\delta,p}
\epsilon^2\nu^{-\frac14-\frac2p},
\qquad p>4.
\end{equation*}
Indeed, using the notation from \eqref{a-truncated-integral} and using Lemma \ref{atailbound} and Bernstein to estimate any high frequency tail errors, we have
\begin{equation*}
\mathcal E_a=\int_{\delta\nu}^{\nu}P_{<M\nu}\mathcal R_\rho
\,\frac{d\rho}{\rho}+O_{L_{t}^2L_x^p}(\epsilon^2\nu^{-\frac{1}{4}-\frac{2}{p}})
\end{equation*}
and its bound follows by repeating the
remainder estimates used to prove \eqref{arLpbound} using $s_c=2$ and then integrating over
$\delta\nu\leq\rho\leq\nu$.
We pair the remaining principal contribution in \eqref{anulin},  with $u_\mu^{x_0}$ using
\eqref{uab-bi-unbal} and the kernel bound for
$Op(\mathfrak l^l_{<\delta\rho})$. Since $\rho\approx_\delta\nu$, the $L_{t,x}^2$ norm of this pairing is bounded by the right hand side of \eqref{transaveua}.
For the remainder, Bernstein, and the estimate for $\mathcal{E}_a$ above gives
\begin{equation*}
\|\mathcal E_a\|_{L_t^2L_x^\infty}
\lesssim_{\delta,p}
\nu^{\frac2p}\|\mathcal E_a\|_{L_t^2L_x^p}
\lesssim_{\delta,p}\epsilon^2\nu^{-\frac14}.
\end{equation*}
We then place $u_\mu^{x_0}$ in $L_t^\infty L_x^2$.
Combining these estimates gives \eqref{transaveua}.
\\

We now turn to establishing \eqref{principalaerr}. We first write
\begin{equation}\label{aprincipal3terms}
\begin{split}
L^l_{<\nu}(a^l_{<\nu}\eta_\nu)
-T_{\mathfrak l}T_{a^l}\eta_\nu
&=L^l_{<\nu}
\bigl((a^l_{<\nu}-a^l_{<\delta\nu})\eta_\nu\bigr)+(L^l_{<\nu}-T_{\mathfrak l})
(a^l_{<\delta\nu}\eta_\nu)
\\
&+T_{\mathfrak l}
\bigl(a^l_{<\delta\nu}\eta_\nu-T_{a^l}\eta_\nu\bigr).
\end{split}
\end{equation}
We estimate these three terms separately. For the remainder of this estimate we will write $\mathcal{Q}_{\nu}$ to generically denote an error term (which will be quartic or higher order) that satisfies the $L_t^1L_x^{\infty}$ bound
\begin{equation}\label{QnuL1Linf}
\|\mathcal Q_\nu\|_{L_t^1L_x^\infty}
\lesssim_{\delta,M}\epsilon^4\nu^{-\frac12}.
\end{equation}
Their contribution in \eqref{principalaerr} is allowed, since Bernstein and two $L_t^{\infty}L_x^2$ bounds give
\begin{equation*}
\begin{split}
\|L(u_\mu^{x_0},u_\mu,\nabla P_{<M\nu}\mathcal Q_\nu)\|_{L_{t,x}^1}\lesssim\nu\|u_\mu\|_{L_t^\infty L_x^2}^2
\|\mathcal Q_\nu\|_{L_t^1L_x^\infty}\lesssim
\epsilon^6c_\mu^2\mu^{-2s_c}\nu^{\frac12}.
\end{split}
\end{equation*}
For the first term, we set
\begin{equation*}
\begin{split}
f_\beta=S_\beta\bigl((a^l_{<\nu}-a^l_{<\delta\nu})\eta_\nu\bigr),\hspace{5mm}f_\beta^{\mathrm{pr}}=S_\beta\bigl((a^l_{<\nu}-a^l_{<\delta\nu}
-\mathcal E_a)\eta_\nu\bigr),
\qquad \beta\lesssim_M\nu.
\end{split}
\end{equation*}
where $\mathcal{E}_{a}$ is as above. By \eqref{anulin} and the bootstrap bounds,
\begin{equation*}
\|a^l_{<\nu}-a^l_{<\delta\nu}-\mathcal E_a\|_{L_t^4L_x^\infty}
\lesssim_\delta\epsilon\nu^{-\frac18},
\qquad
\|\eta_\nu\|_{L_t^4L_x^\infty}
\lesssim\epsilon\nu^{-\frac98}.
\end{equation*}
We first replace the outer elliptic operator by $T_L$.
With $f_\beta$ held fixed, the linearization formula gives
\begin{equation*}
\begin{split}
(L_{<\nu}-L_{<\delta\beta})f_\beta
=\int_{\delta\beta}^{\nu}
\left[
-L_{<\rho}(\eta_\rho L_{<\rho}f_\beta)
+\eta_\rho L_{<\rho}^2f_\beta
\right]\frac{d\rho}{\rho}.
\end{split}
\end{equation*}
We can easily replace the elliptic operators in this integral by
$T_{\mathfrak l}$ modulo a quartic error of the form $\mathcal{Q}_{\nu}$. More precisely, we obtain
\begin{equation*}
\begin{split}
P_{<M\nu}L^l_{<\nu}
\bigl((a^l_{<\nu}-a^l_{<\delta\nu})\eta_\nu\bigr)=
P_{<M\nu}T_L
\bigl((a^l_{<\nu}-a^l_{<\delta\nu})\eta_\nu\bigr)
+\mathcal C_\nu+\mathcal Q_\nu,
\end{split}
\end{equation*}
where
\begin{equation*}
\begin{split}
\mathcal C_\nu
=P_{<M\nu}
\int_{\beta\lesssim_M\nu}\int_{\delta\beta}^{\nu}
\left[
-T_{\mathfrak l}
(\eta_\rho T_{\mathfrak l}f_\beta^{\mathrm{pr}})
+\eta_\rho T_{\mathfrak l}^2f_\beta^{\mathrm{pr}}
\right]\frac{d\rho}{\rho}\frac{d\beta}{\beta}.
\end{split}
\end{equation*}
The error term $\mathcal Q_\nu$ comes from replacing
$f_\beta$ by $f_\beta^{\mathrm{pr}}$ in the linearization
remainder, replacing the elliptic operators there by
$T_{\mathfrak l}$, and high frequency tail terms.
These errors satisfy \eqref{QnuL1Linf} by using Propositions
\ref{basicopbounds}, \ref{paraexp} and
\ref{Bernsteinelliptic}. For instance, the operator errors are estimated
using four $L_t^4L_x^\infty$ factors. For the terms
containing $\mathcal E_a$, we instead use its
$L_t^2L_x^p$ bound, two $L_t^4L_x^\infty$ factors,
and Bernstein at the final output frequency
$\lesssim M\nu$. We omit the details, as these estimates
are similar to many of those carried out in previous sections. To estimate the contribution of $\mathcal C_\nu$ in \eqref{principalaerr}, we pair $\eta_\rho$ with
$u_\mu^{x_0}$, put the other $u_\mu$ in $L_t^\infty L_x^2$,
and put the two factors defining $f_\beta^{\mathrm{pr}}$ in
$L_t^4L_x^\infty$. We note carefully that one must retain the exterior derivative $\nabla$ when the product $\eta_{\rho}T_{\mathfrak{l}}f_{\beta}^{\mathrm{pr}}$ in the first term in the integrand above has small
output frequency $\alpha$. This is because the kernel of the outside $T_{\mathfrak{l}}$ introduces
a translation at scale $\alpha^{-1}$, giving a factor
$(1+\rho/\alpha)^{1/2}$ in the bilinear estimate. However, this is compensated by the fact that there are two "derivatives" applied at the output frequency (i.e. $\nabla$ and $T_{\mathfrak{l}}$), and we have
\begin{equation*}
\int_{\alpha\lesssim_\delta\rho}
\alpha^2\beta\left(1+\frac{\rho}{\alpha}\right)^{\frac12}
\frac{d\alpha}{\alpha}
\lesssim_\delta\rho^2\beta.
\end{equation*}
Consequently, we obtain
\begin{equation*}
\begin{split}
\|L(u_\mu^{x_0},u_\mu,\nabla\mathcal C_\nu)\|_{L_{t,x}^1}&\lesssim_{\delta,M}
\epsilon^5c_\mu^2\mu^{-2s_c}(1+\mu|x_0|)^{\frac12}
\nu^{-\frac54}
\int_{\beta\lesssim_M\nu}\int_{\delta\beta}^{\nu}
\beta\rho^{\frac34}
\frac{d\rho}{\rho}\frac{d\beta}{\beta}
\\
&\qquad\lesssim_{\delta,M}
\epsilon^5c_\mu^2\mu^{-2s_c}\nu^{\frac12}
(1+\mu|x_0|)^{\frac12}.
\end{split}
\end{equation*}
It remains to estimate the contribution of
$P_{<M\nu}T_L((a^l_{<\nu}-a^l_{<\delta\nu})\eta_\nu)$. For this, we apply Proposition \ref{paraexp} to $f_\beta$, at its
output frequency $\beta$. The principal term and the
zeroth-order term have kernels bounded by an order one
kernel at scale $\beta$, since
\begin{equation*}
\|\nabla^2\eta_{<\delta\beta}\|_{L_t^\infty\dot B^0_{\infty,1}}
\lesssim\epsilon\beta.
\end{equation*}

Here, we use two unbalanced bilinear estimates where we pair
$a^l_{<\nu}-a^l_{<\delta\nu}$ with $u_\mu$ using
\eqref{transaveua}, and $\eta_\nu$ with $u_\mu^{x_0}$ using
\eqref{uab-bi-unbal}. Using $\beta$ to denote the output frequency of the product $(a^l_{<\nu}-a^l_{<\delta\nu})\eta_\nu$, the kernels in these bilinear estimates (for each fixed $\beta$) produce
two translation factors of size $(1+\frac{\nu}{\beta})^{\frac{1}{2}}$. This is compensated for by two factors of $\beta$; one coming from the gradient applied to $T_{L}((a^l_{<\nu}-a^l_{<\delta\nu})\eta_\nu)$, and another from the operator $T_L$.  Consequently, modulo errors of the form $\mathcal Q_\nu$, the kernel bounds in
Proposition \ref{paraexp} give
\begin{equation*}
\begin{split}
&\|L(u_\mu^{x_0},u_\mu,\nabla P_{<M\nu}
T_L((a^l_{<\nu}-a^l_{<\delta\nu})\eta_\nu))\|_{L_{t,x}^1}
\\
&\lesssim_{\delta,M}
\epsilon^4c_\mu^2\mu^{-2s_c}\nu^{-\frac32}
(1+\mu|x_0|)^{\frac12}
\int_{\beta\lesssim\nu}
\beta^2\left(1+\frac{\nu}{\beta}\right)
\frac{d\beta}{\beta}
\\
&\lesssim_{\delta,M}
\epsilon^4c_\mu^2\mu^{-2s_c}\nu^{\frac12}
(1+\mu|x_0|)^{\frac12}.
\end{split}
\end{equation*}

Here the higher-order remainder in Proposition
\ref{paraexp} can be easily seen to belong to  $\mathcal Q_\nu$. This completes the estimate for the first term in \eqref{aprincipal3terms}.
\\

For the second term, set
\begin{equation*}
h_\beta=S_\beta(a^l_{<\delta\nu}\eta_\nu),
\qquad \beta\approx\nu.
\end{equation*}
This is nonzero only when $\beta\approx\nu$ because the coefficient $a_{<\delta\nu}^l$ is localized to frequency $\ll\nu$. We can write the second term in \eqref{aprincipal3terms} in the form 
\begin{equation}\label{secondtermrep}
(L^l_{<\nu}-T_{\mathfrak l})
(a^l_{<\delta\nu}\eta_\nu)=P_{<M\nu}\int_{\beta\approx\nu} (L_{<\nu}-Op(\mathfrak{l}^l_{<\delta\beta}))h_{\beta}\frac{d\beta}{\beta}.   
\end{equation}
We again linearize $L_{<\nu}$ by replacing $L_{<\nu}$ by $L_{<\delta\beta}$ before applying Proposition \ref{paraexp}.
The principal quadratic term, after replacing the operators in the linearization by their principal parts, are
\begin{equation*}
P_{<M\nu}
\int_{\beta\approx\nu}\int_{\delta\beta}^{\nu}
\left[
-T_{\mathfrak l}(\eta_\rho T_{\mathfrak l}h_\beta)
+\eta_\rho T_{\mathfrak l}^2h_\beta
\right]\frac{d\rho}{\rho}\frac{d\beta}{\beta}.
\end{equation*}
Here $\rho\approx_\delta\nu$. We pair $\eta_\rho$ with
$u_\mu$ and $\eta_\nu$ with $u_\mu^{x_0}$, placing
$a^l_{<\delta\nu}$ in $L_{t,x}^\infty$. For output frequency $\alpha$ inside the first $T_{\mathfrak l}$ factor, the two bilinear pairings introduce the translation factors $(1+\rho/\alpha)^{1/2}$ and $(1+\nu/\alpha)^{1/2}$, respectively. Since $\rho\approx \nu$, their product is bounded by $(1+\nu/\alpha)$. From the outer gradient factor (which by Bernstein contributes one factor of $\alpha$) and the kernel bounds for operator $T_{\mathfrak{l}}$ (which contributes a second factor), the integral in the output frequency $\alpha$ is bounded by
\begin{equation*}
\nu\int_{\alpha\lesssim_\delta\nu}
\alpha^2\left(1+\frac{\nu}{\alpha}\right)
\frac{d\alpha}{\alpha}
\lesssim_\delta\nu^3.
\end{equation*}
This shows the principal quadratic part of the replacement error between $L_{<\nu}$ and $L_{<\delta\beta}$ in \eqref{secondtermrep} is perturbative. 
We also show that the same perturbative bound holds for the errors in replacing the elliptic
operators when extracting the principal part above. Up to terms containing
two operator errors (which can be placed in the $\mathcal{Q}_{\nu}$ box by arguing as above), these are the following expressions,
integrated in $\rho$ and $\beta$ as above:
\begin{equation*}
\begin{split}
&-T_{\mathfrak l}
\bigl(\eta_\rho(L_{<\rho}-T_{\mathfrak l})h_\beta\bigr),\hspace{5mm}\eta_\rho T_{\mathfrak l}
(L_{<\rho}-T_{\mathfrak l})h_\beta,\hspace{5mm}\eta_\rho(L_{<\rho}-T_{\mathfrak l})
T_{\mathfrak l}h_\beta,\hspace{5mm}-(L_{<\rho}-T_{\mathfrak l})
\bigl(\eta_\rho T_{\mathfrak l}h_\beta\bigr).
\end{split}
\end{equation*}
All of the above terms are cubic or higher order. For the first three, we can pair the  $\eta_{\rho}$ factor with $u_\mu^{x_0}$, place $u_\mu$ in $L_t^\infty L_x^2$ and put the remaining two factors in $L_t^4\dot{B}_{\infty,1}^0$ to obtain a perturbative contribution. For the last expression, one has to integrate the internal product
$\eta_\rho T_{\mathfrak l}h_\beta$ over output frequencies
$\alpha$ and then linearize $L$ again to compute the error generated by $L_{<\rho}-T_{\mathfrak{l}}$. The resulting terms are of the form $\mathcal{Q}_{\nu}$ as well as a term that can be estimated similarly to $\mathcal{C}_{\nu}$ above. We omit the simple modifications for this. To summarize, the full error in replacing $L_{<\nu}$ with $L_{<\delta\beta}$ in  \eqref{secondtermrep} is perturbative. The remaining part we have to estimate is then
\begin{equation*}
P_{<M\nu}\int_{\beta\approx\nu}(L_{<\delta\beta}-Op(\mathfrak{l}^l_{<\delta\beta}))h_{\beta}\frac{d\beta}{\beta}    
\end{equation*}
For this term, we apply Proposition \ref{paraexp}. The zeroth-order term
$L^0_{<\delta\beta}h_\beta$ is estimated using
\eqref{K0bound}. We can pair its
$\nabla^2\eta_{<\delta\beta}$ factor with $u_\mu$ and
$\eta_\nu$ with $u_\mu^{x_0}$. This again gives a perturbative contribution. The higher-order remainder term in \eqref{paraexp} contributes a term of the form $\mathcal{Q}_{\nu}$. Indeed,
\begin{equation*}
\begin{split}
\|\mathcal E_{h_\beta}\|_{L_t^1L_x^\infty}\lesssim
\nu^{-2}
\|\nabla^2\eta_{<\delta\beta}\|_{L_t^4\dot B^0_{\infty,1}}^3
\|a^l_{<\delta\nu}\|_{L_{t,x}^\infty}
\|\eta_\nu\|_{L_t^4L_x^\infty}\lesssim\epsilon^4\nu^{-\frac12}.
\end{split}
\end{equation*}
This completes the estimate for the second term in \eqref{aprincipal3terms}.
\\

Finally, for the last term in \eqref{aprincipal3terms}, using the definition of $T_{a^l}$,
\begin{equation*}
a^l_{<\delta\nu}\eta_\nu-T_{a^l}\eta_\nu
=
\int_{\kappa\approx\nu}
(a^l_{<\delta\nu}-a^l_{<\delta\kappa})
S_\kappa\eta_\nu\,\frac{d\kappa}{\kappa}.
\end{equation*}

Applying \eqref{transaveua} and \eqref{uab-bi-unbal}, we pair
$a^l_{<\delta\nu}-a^l_{<\delta\kappa}$ with $u_\mu$
and $S_\kappa\eta_\nu$ with $u_\mu^{x_0}$, respectively. Using this, and then integrating over $\kappa\approx\nu$ gives
\begin{equation*}
\begin{split}
&\|L(u_\mu^{x_0},u_\mu,\nabla P_{<M\nu}T_{\mathfrak l}
(a^l_{<\delta\nu}\eta_\nu-T_{a^l}\eta_\nu))\|_{L_{t,x}^1}
\\
&\lesssim_\delta
\nu^2\int_{\kappa\approx\nu}
\left(\epsilon^2c_\mu\mu^{-s_c}\nu^{-\frac14}\right)
\left(\epsilon^2c_\mu\mu^{-s_c}\nu^{-\frac54}\right)
(1+\mu|x_0|)^{\frac12}
\frac{d\kappa}{\kappa}
\\
&\lesssim_\delta
\epsilon^4c_\mu^2\mu^{-2s_c}\nu^{\frac12}
(1+\mu|x_0|)^{\frac12}.
\end{split}
\end{equation*}
Combining this with the estimates for the other two terms in \eqref{aprincipal3terms} finally proves \eqref{principalaerr}.
\subsubsection{The term $\Gamma^l_{<\nu}$}
Next, we show that the contribution of $\Gamma^l_{<\nu}$ is perturbative. First, by Bernstein and H\"older,
\begin{equation*}
\begin{split}
\left\|L\left(u^{x_0}_\mu,u_\mu,\nabla P_{<M\nu}\left(\frac{\eta_\nu|A^l_{<\nu}v^l_{<\nu}|^2}{1+|\nabla\eta_{<\nu}|^2}\right)\right)\right\|_{L_{t,x}^1}&\lesssim\nu\|u_\mu\|_{L_t^\infty L_x^2}\|L(u_\mu^{x_0},\eta_\nu)\|_{L_{t,x}^2}\|A^l_{<\nu}v^l_{<\nu}\|_{L_t^4\dot B_{\infty,1}^0}^2
\\
&\lesssim\epsilon^5c_\mu^2\mu^{-2s_c}\nu^{\frac12}(1+\mu|x_0|)^{\frac12}.
\end{split}
\end{equation*}
This handles the four square terms in the definition of $\Gamma_{<\nu}^l$. For the remaining two terms, use two unbalanced bilinear estimates, pairing $\eta_\nu$ with $u_\mu$ and $u_\mu^{x_0}$ with either $\Delta\eta_{<\nu}$ or $\nabla a^l_{<\nu}$. Proposition \ref{abilinearbounds} give the factors $\nu^{-5/4}$ and $\nu^{3/4}$ (along with a factor of $(1+\mu|x_0|)^{\frac12}$), respectively. Hence, the full contribution of $\Gamma_{<\nu}$ in \eqref{dlambdaaest} is perturbative.

\subsubsection{The two product terms}
Bernstein and H\"older give
\begin{equation*}
\|L(u_\mu^{x_0},u_\mu,\nabla P_{<M\nu}F)\|_{L_{t,x}^1}
\lesssim\nu\|L(u_\mu^{x_0},L^l_{<\nu}\Phi^{l,0}_{<\nu})\|_{L_{t,x}^2}\|L(u_\mu,L^l_{<\nu}B^l_{<\nu})\|_{L_{t,x}^2}.
\end{equation*}
Lemma \ref{Phiest}, Proposition \ref{basicopbounds}, and \eqref{a-difference-Phi-tail} give
\begin{equation*}
A^l_{<\nu}\Phi^{l,0}_{<\nu}
=T_{\mathfrak a}\varphi_\nu+O_{L_t^2L_x^\infty}(\epsilon^2\nu^{-\frac34})
\end{equation*}
Hence, we have
\begin{equation*}
\|L(u_\mu^{x_0},A^l_{<\nu}\Phi^{l,0}_{<\nu})\|_{L_{t,x}^2}\lesssim\epsilon^2c_\mu\mu^{-s_c}\nu^{-\frac34}(1+\mu|x_0|)^{\frac12}.
\end{equation*}
On the other hand, we also have the unbalanced estimate
\begin{equation*}
\|L(u_\mu,L^l_{<\nu}v^l_{<\nu})\|_{L_{t,x}^2}\lesssim\epsilon^2c_\mu\mu^{-s_c}\nu^{\frac14}.
\end{equation*}
Combining these estimates gives
\begin{equation*}
\|L(u_\mu^{x_0},u_\mu,\nabla P_{<M\nu}(A^l_{<\nu}\Phi^{l,0}_{<\nu}L^l_{<\nu}v^l_{<\nu}))\|_{L_{t,x}^1}\lesssim\epsilon^4c_\mu^2\mu^{-2s_c}\nu^{\frac12}(1+\mu|x_0|)^{\frac12}.
\end{equation*}
This handles the two product terms in \eqref{datruncated}.

\subsubsection{Commutator errors}
Now we study the contribution of the two commutators in \eqref{datruncated}. We show the details for $[L_{<\nu},B^l_{<\nu}](L_{<\nu}\Phi^0_{<\nu})$, since the other term is similar. Fix $K$ sufficiently large relative to $M$ and $\delta^{-1}$. Relative to $P_{<K\nu}$, we decompose the commutator crudely as follows
\begin{equation*}
\begin{split}
[L_{<\nu},B^l_{<\nu}](L_{<\nu}\Phi^0_{<\nu})
&=[L_{<\nu},B^l_{<\nu}]P_{\geq K\nu}(L_{<\nu}\Phi^0_{<\nu})
\\
&\quad+\sum_{\beta\lesssim K\nu}[L_{<\nu},P_{\geq\delta\beta}B^l_{<\nu}]P_\beta P_{<K\nu}(L_{<\nu}\Phi^0_{<\nu})
\\
&\quad+\sum_{\beta\lesssim K\nu}[L_{<\nu},P_{<\delta\beta}B^l_{<\nu}]P_\beta P_{<K\nu}(L_{<\nu}\Phi^0_{<\nu})
\\
&\coloneq\mathcal C_1+\mathcal C_2+\mathcal C_3.
\end{split}
\end{equation*}
Proposition \ref{highfreqcommutator} and Proposition \ref{highfreqtailcomm} applies to $\mathcal C_1$ because of the choice of $K$. Combining this with \eqref{a-difference-Phi-tail} gives
\begin{equation*}
\begin{split}
\|\nabla P_{<M\nu}\mathcal C_1\|_{L_t^1L_x^\infty}
&\lesssim\nu\|\nabla B^l_{<\nu}\|_{L_t^4\dot B_{\infty,1}^0}\|P_{\geq K\nu}L_{<\nu}\Phi^0_{<\nu}\|_{L_t^{\frac43}\dot B_{\infty,1}^0}
\lesssim\epsilon^4\nu^{\frac12}.
\end{split}
\end{equation*}
Therefore, the contribution of $\mathcal{C}_1$ is perturbative by placing $u_\mu$ and $u_\mu^{x_0}$ in $L_t^{\infty}L_x^2$. For the remaining terms, set
\begin{equation*}
\mathcal E_\nu=L_{<\nu}\Phi^0_{<\nu}-T_{\mathfrak l}\varphi_\nu.
\end{equation*}
Lemma \ref{Phiest}, Proposition \ref{Bernsteinelliptic} and Proposition \ref{basicopbounds} gives for $p>4$
\begin{equation}\label{a-difference-Phi-error}
\|P_{<K\nu}\mathcal E_\nu\|_{L_t^2L_x^p}\lesssim_p\epsilon^2\nu^{-\frac34-\frac2p},
\qquad
\|P_{<K\nu}\mathcal E_\nu\|_{L_t^4L_x^{\frac{16}{7}}}\lesssim\epsilon^2\nu^{-\frac32}.
\end{equation}
For $\mathcal{C}_2$, we expand the commutator in each summand $\mathcal{C}_2^{\beta}$
\begin{equation*}
\mathcal{C}_2^{\beta}=L_{<\lambda}(P_{\geq\delta\beta}B^l_{<\lambda}P_{\beta}(L_{<\lambda}\Phi_{<\lambda}))-P_{\geq\delta\beta}B^l_{<\lambda}L_{<\lambda}P_{\beta}(L_{<\lambda}\Phi_{<\lambda})
\end{equation*}
To estimates this term, one peels off the principal part by replacing each $L$ by $T_{\mathfrak{l}}$ and $\Phi_{<\lambda}$ with $\varphi_{\lambda}$. For the principal part, we use the unbalanced pairings
$(u_\mu,P_{\geq\delta\beta} B^l_{<\lambda})$ and
$(u_\mu^{x_0},\varphi_\lambda)$ alongside Proposition \ref{VBunbil}. All of the subprincipal errors are handled as in the proof of \eqref{principalaerr}. We omit the details since they are similar.
\\

It remains to estimate $\mathcal C_3$. We apply Proposition \ref{weighted-commutator-bound} with
\begin{equation*}
f=P_{<\delta\beta}B^l_{<\nu},\qquad g_\beta=P_\beta P_{<K\nu}\mathcal E_\nu,
\end{equation*}
and $h_1=u_{\mu}$ and $h_2=u_{\mu}^{x_0}$ to obtain
\begin{equation}\label{firstC3comm}
\begin{split}
&\Big\|L\left(u_{\mu}^{x_0},u_{\mu},\nabla P_{<M\nu}[L_{<\nu},f]g_{\beta}\right)\Big\|_{L_{t,x}^1}\lesssim_{\sigma,p,r,M,K}
\nu
M_{\nu,2}^{\theta}\left(u_\mu,\nabla_xf\right)
M_{\nu,2}^{\theta}\left(u_\mu^{x_0},g_\beta\right)
\\
&\hspace{20mm}+
\left(\frac{\nu}{\beta}\right)^{1+\frac2r}
M_{\nu,2}^{\theta}\left(u_\mu,\nabla_xf\right)
\|u_\mu^{x_0}\|_{L_t^\infty L_x^2}
\|\nabla_x^2\eta_{<\nu}\|_{L_t^4\dot B_{\infty,1}^0}
\||D|^{\frac2p}g_\beta\|_{L_t^4L_x^p}
\\
&\hspace{20mm}+
M_{\nu,2}^{\theta}\left(u_\mu,\nabla_x^2\eta_{<\nu}\right)
\|u_\mu^{x_0}\|_{L_t^\infty L_x^2}
\|\nabla_xf\|_{L_t^4\dot B_{\infty,1}^0}
\||D|^{\frac2p}g_\beta\|_{L_t^4L_x^p}.
\end{split}
\end{equation}
where $\theta=\frac{\sigma}{2}$ and $\sigma$ and $r$ are as in the statement of Proposition \ref{weighted-commutator-bound} with their specific values to be chosen below. Next, by Proposition \ref{basicMproperties} and \eqref{a-difference-Phi-error} we have
\begin{equation*}
\begin{split}
M_{\nu,2}^\theta(u_\mu,g_\beta)\lesssim\epsilon^3c_\mu\mu^{-s_c}\nu^{-\frac34},\hspace{4mm}\||D|^{\frac2p}g_\beta\|_{L_t^4L_x^p}
\lesssim\epsilon^2\beta^{\frac78}\nu^{-\frac32}
\end{split}
\end{equation*}
Moreover, by Proposition \ref{VBunbil} and Lemma \ref{annularcontrol}, we have
\begin{equation*}
\begin{split}
M_{\nu,2}^\theta(u_\mu,\nabla f)\lesssim\epsilon^2c_\mu\mu^{-s_c}
\nu^{1-\theta}\beta^{\theta-\frac34},\hspace{5mm}
M_{\nu,2}^\theta(u_\mu,\nabla^2\eta_{<\nu})\lesssim\epsilon^2c_\mu\mu^{-s_c}\nu^{\frac34}.
\end{split}
\end{equation*}
The bound above for $M_{\nu,2}^\theta(u_\mu,g_\beta)$ also
holds with $u_\mu$ replaced by $u_\mu^{x_0}$, since it uses
only the energy bound
$\|u_\mu^{x_0}\|_{L_t^\infty L_x^2}
\lesssim\epsilon c_\mu\mu^{-s_c}$.

For the first term in \eqref{firstC3comm}, we combine the
pairing with $\nabla f$ with the bound for
$M_{\nu,2}^\theta(u_\mu^{x_0},g_\beta)$.
For the second term, we use the same pairing with $\nabla f$,
the energy bound, and the $L_t^4\dot{B}_{\infty,1}^0$ bound for
$\nabla^2\eta_{<\nu}$ and the above $L_t^4L_x^p$ bound for $|D|^{2/p}g_\beta$, retaining
the factor $(\nu/\beta)^{1+2/r}$.
For the third term, we instead use the pairing with
$\nabla^2\eta_{<\nu}$, $L_t^{\infty}L_x^2$ bound  for $u_{\mu}^{x_0}$, and the $L_t^4$
bounds for $\nabla f$ and $|D|^{2/p}g_\beta$ (the former coming from Proposition \ref{VBunbil}).
Thus the three contributions, sum to
\begin{equation*}
\begin{split}
&\epsilon^5c_\mu^2\mu^{-2s_c}\nu^{\frac54-\theta}
\sum_{\beta\lesssim K\nu}\beta^{\theta-\frac34}+\epsilon^6c_\mu^2\mu^{-2s_c}\nu^{\frac{11}{8}-\theta+\frac2r}
\sum_{\beta\lesssim K\nu}\beta^{\theta-\frac78-\frac2r}+\epsilon^6c_\mu^2\mu^{-2s_c}\nu^{-\frac34}
\sum_{\beta\lesssim K\nu}\beta^{\frac54}
\\
&\qquad\lesssim\epsilon^5c_\mu^2\mu^{-2s_c}\nu^{\frac12}.
\end{split}
\end{equation*}
where above, we choose $1<\sigma<2$ and $4<p<r<4/(2-\sigma)$ (so that our application of Proposition \ref{weighted-commutator-bound} is valid) and also $r,p$ and $\sigma$ so that
\begin{equation*}
\theta=\frac\sigma2>\frac78+\frac2r
\end{equation*}
to ensure the sums converge. For example, $\sigma=63/32$, $p=32$ and $r=64$ satisfy these conditions. Finally, the principal part of $\mathcal C_3$ is
\begin{equation*}
\mathcal C_4=\sum_{\beta\approx\nu}[L_{<\nu},P_{<\delta\beta}B^l_{<\nu}]P_\beta T_{\mathfrak l}\varphi_\nu.
\end{equation*}
There are only finitely many input frequencies $\beta\approx\nu$. For $g_\beta=P_\beta T_{\mathfrak l}\varphi_\nu$, the kernel bounds for $T_{\mathfrak{l}}\varphi$, Lemma \ref{annularcontrol} and Proposition \ref{phiu-bootstrap} give
\begin{equation*}
M_{\nu,2}^\theta(u_\mu^{x_0},g_\beta)\lesssim\epsilon^2c_\mu\mu^{-s_c}\nu^{-\frac34}(1+\mu|x_0|)^{\frac12},
\qquad
\||D|^{\frac2p}g_\beta\|_{L_t^4L_x^p}\lesssim\epsilon\nu^{-\frac58}.
\end{equation*}
for $p>8$. Applying Proposition \ref{weighted-commutator-bound} once more,
now with $\beta\approx\nu$ and the same choices of $p,r,\sigma$
as above, the line in \eqref{firstC3comm} is bounded by
\begin{equation*}
\nu\left(\epsilon^2c_\mu\mu^{-s_c}\nu^{\frac14}\right)
\left(\epsilon^2c_\mu\mu^{-s_c}\nu^{-\frac34}(1+\mu|x_0|)^{\frac12}\right)
=\epsilon^4c_\mu^2\mu^{-2s_c}\nu^{\frac12}(1+\mu|x_0|)^{\frac12}.
\end{equation*}
For the second line in \eqref{firstC3comm}, the factor
$(\nu/\beta)^{1+2/r}$ is bounded. Using the pairing with
$\nabla f$, and $L_t^4$ bounds for
$\nabla^2\eta_{<\nu}$ and $|D|^{2/p}g_\beta$, respectively,
gives
\begin{equation*}
\left(\epsilon^2c_\mu\mu^{-s_c}\nu^{\frac14}\right)
\left(\epsilon c_\mu\mu^{-s_c}\right)
\left(\epsilon\nu^{\frac78}\right)
\left(\epsilon\nu^{-\frac58}\right)
=\epsilon^5c_\mu^2\mu^{-2s_c}\nu^{\frac12}.
\end{equation*}
For the third term in \eqref{firstC3comm}, we instead use
the pairing with $\nabla^2\eta_{<\nu}$, and the $L_t^4$ bounds for $\nabla f$ and
$|D|^{2/p}g_\beta$, respectively, obtaining
\begin{equation*}
\left(\epsilon^2c_\mu\mu^{-s_c}\nu^{\frac34}\right)
\left(\epsilon c_\mu\mu^{-s_c}\right)
\left(\epsilon\nu^{\frac38}\right)
\left(\epsilon\nu^{-\frac58}\right)
=\epsilon^5c_\mu^2\mu^{-2s_c}\nu^{\frac12}.
\end{equation*}
Combining everything above, we have finally proved \eqref{dlambdaaest} for the contribution of the first commutator in \eqref{datruncated}. Replacing $B,L\Phi^0$ by $V,M\Phi^0$ gives the other commutator in the same way. Combining all the preceding estimates, gives
\begin{equation*}
\begin{split}
&\mu^{-\frac12}\|L(u_\mu^{x_0},u_\mu,\nabla P_{<M\lambda}(a_{<\lambda}+T_{\mathfrak l}T_{a^l}\eta_{<\lambda}))\|_{L_{t,x}^1}
\\
&\qquad\lesssim\epsilon^4c_\mu^2\mu^{-2s_c}(1+\mu|x_0|)^{\frac12}
\mu^{-\frac12}\int_0^\lambda\nu^{\frac12}\frac{d\nu}{\nu}
\lesssim\epsilon^4c_\mu^2\mu^{-2s_c}(1+\mu|x_0|)^{\frac12}.
\end{split}
\end{equation*}
This finally proves \eqref{adiffmulti}. 

Now, we prove the second error bound \eqref{aprincipalL2Linf}, which is essentially a corollary of the estimates used above. First, the estimates leading to \eqref{a-integral-principal} gives
\begin{equation*}
P_{<M\lambda}a_\lambda+
(P_{<\delta\lambda}a_{<\lambda})Op(\mathfrak l^l_{<\delta\lambda})\eta_\lambda=O_{L_t^2L_x^{\infty}}(\epsilon^2c_\lambda\lambda^{-\frac14}).
\end{equation*}
We can then replace $P_{<M\lambda}a_{\lambda}$ in this bound by $a_{\lambda}$. Indeed, the error term $P_{>M\lambda}a_\lambda$ consists of high frequency tails. For these, we can use the decomposition leading
to \eqref{datruncated}, and estimate each error in $L_t^2L_x^{\infty}$ instead of $L_t^1L_x^{\infty}$ (using for instance, two $L_t^4\dot{B}_{\infty,1}^0$ bounds). Hence, we have
\begin{equation}\label{aL2Linfintermediate}
a_\lambda+
(P_{<\delta\lambda}a_{<\lambda})Op(\mathfrak l^l_{<\delta\lambda})\eta_\lambda=O_{L_t^2L_x^{\infty}}(\epsilon^2c_\lambda\lambda^{-\frac14}).
\end{equation}
Finally, subtracting \eqref{a-truncated-integral} from
\eqref{a-integral-principal} gives
$\|a-a_{<\lambda}\|_{L_t^4L_x^\infty}\lesssim_\delta
\epsilon\lambda^{-1/8}$ and thus, we also have $\|a_{<\alpha}-a_{<\lambda}\|_{L_t^4L_x^\infty}\lesssim_\delta
\epsilon\lambda^{-1/8}$ for $\alpha>0$ satisfying $\alpha\approx_{\delta}\lambda$. Hence, we have
\begin{equation*}
(P_{<\delta\lambda}a_{<\lambda})Op(\mathfrak l^l_{<\delta\lambda})\eta_\lambda=T_{a^l}Op(\mathfrak{l}^l_{<\delta\lambda})\eta_{\lambda}+O_{L_t^2L_x^{\infty}}(\epsilon^2c_\lambda\lambda^{-\frac14})    
\end{equation*}
Similarly, by Bernstein and the fundamental theorem of calculus, we also have
\begin{equation*}
T_{a^l}Op(\mathfrak{l}^l_{<\delta\lambda})\eta_{\lambda}=T_{a^l}T_{\mathfrak{l}}\eta_{\lambda}+O_{L_t^2L_x^{\infty}}(\epsilon^2c_\lambda\lambda^{-\frac14}).   
\end{equation*}
Combining this with the previous bound and \eqref{aL2Linfintermediate} gives the bound \eqref{aprincipalL2Linf} and completes the proof of the proposition.
\end{proof}
\section{Expansions for \texorpdfstring{$D_t a$}{}}\label{sec:Dtaexpansions}
In this section, our main goal will be to prove the necessary error estimates and expansions related to the (truncated) material derivative of the Taylor coefficient $D_t a$. 
\subsection{The equation for \texorpdfstring{$D_t^l a$}{}}
We begin with the following lemma, which identifies the truncated material derivative of $a_{<\lambda}$ in a sufficiently precise way for proving estimates later. The reason we track this truncated evolution rather than the full one is because the material derivative of $a$ will only appear as a low frequency paracoefficient in our analysis later, rather than in bulk.
\begin{lemma}[Evolution of the truncated Taylor coefficient] Fix $\lambda>0$. Define \begin{equation*}
D_t^l=\partial_t+V^l_{<\lambda}\cdot\nabla,\qquad
V^l_{<\lambda}=P_{<M\lambda}V_{<\lambda}.
\end{equation*}
As in subsecton \ref{subsec:initialphi}, define also the kinematic and Bernoulli errors:
\begin{equation*}
\begin{split}
\kappa_{<\lambda}&=(\partial_t+V_{<\lambda}\cdot\nabla)\eta_{<\lambda}-B_{<\lambda},
\\
\rho_{<\lambda} &=\left. \left( \partial_t\phi_{<\lambda} +\frac12|\nabla_{x,z}\phi_{<\lambda}|^2+z \right) \right|_{z=\eta_{<\lambda}},
\end{split}
\end{equation*} 
where $\phi_{<\lambda}$ denotes the harmonic extension of $\psi_{<\lambda}$ into the domain $\Omega_{<\lambda}\coloneq\{z<\eta_{<\lambda}(t,x)\}$. Define also \begin{equation*} \mathcal{P}F \coloneq \left. \partial_z\Delta_D^{-1}F \right|_{z=\eta_{<\lambda}}, \end{equation*} where $\Delta_D^{-1}$ denotes the inverse Dirichlet Laplacian for $\Omega_{<\lambda}$, and write $\mathcal{H}f$ for the harmonic extension of a function $f$ into the domain $\Omega_{<\lambda}$. Then there holds \begin{equation}\label{Dtaeqn}  
\begin{split} D_t^l a_{<\lambda} =&a_{<\lambda}\nabla\cdot V_{<\lambda} -L_{<\lambda}(\kappa_{<\lambda}a_{<\lambda}) -\kappa_{<\lambda} \left.\partial_z^2p_{<\lambda}\right|_{z=\eta_{<\lambda}} \\ &-4\mathcal P\left( \nabla_{x,z}^2\phi_{<\lambda}: \nabla_{x,z}^2\mathcal H\eta_{<\lambda} \right)-2\mathcal P\left( \operatorname{tr} \big((\nabla_{x,z}^2\phi_{<\lambda})^3\big) \right) 
\\ 
&+4\mathcal P\left( \nabla_{x,z}^2\phi_{<\lambda}: \nabla_{x,z}^2\Delta_D^{-1} |\nabla_{x,z}^2\phi_{<\lambda}|^2 \right)+2\mathcal P\left( \nabla_{x,z}^2\phi_{<\lambda}: \nabla_{x,z}^2\mathcal H\rho_{<\lambda} \right)
\\
&+\bigl(V^l_{<\lambda}-V_{<\lambda}\bigr)\cdot\nabla a_{<\lambda}.
\end{split} \end{equation} 
\end{lemma} 
\begin{proof} Temporarily denote the ``interior'' material derivative for the domain $\Omega_{<\lambda}$ by \begin{equation*} \mathbf D_t =\partial_t +\nabla_{x,z}\phi_{<\lambda}\cdot\nabla_{x,z}. 
\end{equation*} 
For a function $F$ on $\Omega_{<\lambda}$, the chain rule gives
\begin{equation*}
\begin{split}
D_t^l\left( \left.F\right|_{z=\eta_{<\lambda}} \right)
&= \left.\mathbf D_tF\right|_{z=\eta_{<\lambda}} +\kappa_{<\lambda} \left.\partial_zF\right|_{z=\eta_{<\lambda}}
\\
&+\bigl(V^l_{<\lambda}-V_{<\lambda}\bigr)\cdot
\nabla_x\left(\left.F\right|_{z=\eta_{<\lambda}}\right).
\end{split}
\end{equation*} Since $\left.p_{<\lambda}\right|_{z=\eta_{<\lambda}}=0$ and $a_{<\lambda} = -\left.\partial_zp_{<\lambda}\right|_{z=\eta_{<\lambda}}$, we get  \begin{equation*} \left.\mathbf D_tp_{<\lambda}\right|_{z=\eta_{<\lambda}} = \kappa_{<\lambda}a_{<\lambda} \end{equation*} and therefore, 
\begin{equation*} \begin{split} D_t^l a_{<\lambda} &= -\left. \partial_z\mathbf D_tp_{<\lambda} \right|_{z=\eta_{<\lambda}} +a_{<\lambda}\nabla\cdot V_{<\lambda}-\kappa_{<\lambda} \left. \partial_z^2p_{<\lambda} \right|_{z=\eta_{<\lambda}}
\\ 
&+\bigl(V^l_{<\lambda}-V_{<\lambda}\bigr)\cdot\nabla a_{<\lambda}. \end{split} \end{equation*} Here we used \begin{equation*} \left. \nabla_{x,z}p_{<\lambda} \right|_{z=\eta_{<\lambda}} = a_{<\lambda}(\nabla\eta_{<\lambda},-1) \end{equation*} and the harmonicity of \(\phi_{<\lambda}\). From the decomposition $1=\mathcal{H}+\Delta^{-1}_D\Delta$ and the definition of $L_{<\lambda}$, we have \begin{equation*} \left. \partial_z\mathbf D_tp_{<\lambda} \right|_{z=\eta_{<\lambda}} = L_{<\lambda} (\kappa_{<\lambda}a_{<\lambda}) +\mathcal P\left( \Delta_{x,z}\mathbf D_tp_{<\lambda} \right). \end{equation*} 
Now set \begin{equation*} Q_{<\lambda} = \partial_t\phi_{<\lambda} +\frac12|\nabla_{x,z}\phi_{<\lambda}|^2+z. \end{equation*} Since \begin{equation*} \Delta_{x,z}Q_{<\lambda} = |\nabla_{x,z}^2\phi_{<\lambda}|^2, \qquad \Delta_{x,z}p_{<\lambda} = -|\nabla_{x,z}^2\phi_{<\lambda}|^2,  \end{equation*} and $\left. (Q_{<\lambda}+p_{<\lambda}) \right|_{z=\eta_{<\lambda}} = \rho_{<\lambda}$, we have \begin{equation*} Q_{<\lambda}+p_{<\lambda} = \mathcal H\rho_{<\lambda}. \end{equation*} Consequently, \begin{equation*} \mathbf D_t\nabla_{x,z}^2\phi_{<\lambda} = -\nabla_{x,z}^2p_{<\lambda} +\nabla_{x,z}^2\mathcal H\rho_{<\lambda} -(\nabla_{x,z}^2\phi_{<\lambda})^2. \end{equation*} Therefore, we obtain
\begin{equation*} \begin{split} \Delta_{x,z}\mathbf D_tp_{<\lambda} =& \mathbf D_t\Delta_{x,z}p_{<\lambda} +2\nabla_{x,z}^2\phi_{<\lambda}: \nabla_{x,z}^2p_{<\lambda}
\\ 
=& 4\nabla_{x,z}^2\phi_{<\lambda}: \nabla_{x,z}^2p_{<\lambda}-2\nabla_{x,z}^2\phi_{<\lambda}: \nabla_{x,z}^2\mathcal H\rho_{<\lambda} +2\operatorname{tr} \big((\nabla_{x,z}^2\phi_{<\lambda})^3\big). \end{split} \end{equation*} Finally, \begin{equation*} p_{<\lambda} = \mathcal H\eta_{<\lambda}-z -\Delta_D^{-1} |\nabla_{x,z}^2\phi_{<\lambda}|^2. \end{equation*} Substitution completes the proof
\end{proof}

\subsection{Estimates for \texorpdfstring{$D_t^l a$}{}} In this subsection, we prove the main multi-linear expansion involving the term $D_t^l a$. This will specifically arise as part of the principal term coming from the commutator between $T_{\mathfrak{q}}$ and $\partial_t+T_{V^l}\cdot\nabla$ in the equation for $\varphi$ when we derive the full equation for the complex variable $u$ in the next section.
\begin{proposition}\label{Dtacommutatorestimates}
Assume the bootstrap hypotheses in Section \ref{sec:bootstrap}. Using the notation in the previous subsection, we define the principal linear part of $D_t^l a_{<\lambda}$,
\begin{equation*}
(D_t^l a_{<\lambda})^{[1]}\coloneq T_{a^l}\nabla\cdot V_{<\lambda}.
\end{equation*}
For every $\mu>0$ and $x_0\in\mathbb R^2$, there is a decomposition
\begin{equation}\label{Dtacommutatorfluxdecomp}
-\frac{i}{2}\overline u_\mu^{x_0}P_\mu
\left(T_{\frac{\mathfrak q}{a}(D_t^l a-(D_t^l a)^{[1]})}\varphi\right)
=\mathcal R_{a,\mu}+\nabla_x\cdot Q_{a,\mu},
\end{equation}
where, for every $\frac43<p\leq2$,
\begin{equation}\label{Dtacommutatorfluxbound}
\|\mathcal R_{a,\mu}\|_{L_{t,x}^1}
+\mu^{\frac2p-1}\|Q_{a,\mu}\|_{L_t^1L_x^p}
\lesssim_p
\epsilon^4c_\mu^2\mu^{-2s_c}(1+\mu|x_0|)^{\frac12}.
\end{equation}
\end{proposition}
We note that above, $T_{\frac{\mathfrak q}{a}(D_t^l a-(D_t^l a)^{[1]})}$ is defined using the same geometric regularization convention as in Subsection \ref{subsec:principal}. For instance, its coefficient at scale $\lambda$ is given by
\begin{equation*}
P_{<M\lambda}\left(\frac{\mathfrak q_{<\lambda}}{a_{<\lambda}}
\left(D_t^l a_{<\lambda}-(D_t^l a_{<\lambda})^{[1]}\right)\right).
\end{equation*}
The first main ingredient in proving this proposition is the following estimate for the directly perturbative parts in the expansion of $D_t^l a_{<\lambda}$, after removing the contribution $-L_{<\lambda}(a_{<\lambda}\kappa_{<\lambda})$.
\begin{proposition}\label{Dtaestimates}
With the notation above, let $L$ be a translation invariant $4$-linear form with rapidly decaying kernel. We have the following two bounds.
\begin{enumerate}
\item For each dyadic number $\mu$ and $\lambda\ll\mu$ there holds
\begin{equation*}
\big\|L\left(u_{\mu}^{x_0},u_{\mu},(D_t^l a_{<\lambda}-(D_t^l a_{<\lambda})^{[1]}+L_{<\lambda}(a_{<\lambda}\kappa_{<\lambda}))\right)\big\|_{L_{t,x}^1}\lesssim \epsilon^4c_{\mu}^2\mu^{-2s_c}.    
\end{equation*}
The same bound holds with $|D|^{\frac{1}{2}}\varphi$ in place of $u_{\mu}$.
\item For every $\lambda>0$, let $a_\lambda=\lambda\frac{d}{d\lambda}a_{<\lambda}$ and $V_{\lambda}\coloneq S_{\lambda}V$. Define $(D_t^l a_\lambda)^{[1]}=T_{a^l}\nabla\cdot V_\lambda$. Then for every $8\leq p\leq\infty$, there holds
\begin{equation}\label{DtaprincipalL2L8}
\|D_t^l a_\lambda-(D_t^l a_\lambda)^{[1]}\|_{L_t^2L_x^p}
\lesssim_\delta
\epsilon c_\lambda\lambda^{\frac14-\frac2p}.
\end{equation}
\end{enumerate}
\end{proposition}
Only the first multilinear bound above will be relevant for proving Proposition \ref{Dtacommutatorestimates}. The latter bound is milder, but will be needed later in the interaction Morawetz analysis to control some of the error terms. We now proceed by proving Proposition \ref{Dtaestimates}. For the multilinear bound, we will show the details for when the translation-invariant 4-linear form above contains the two factors of $u_{\mu}$ and $u^{x_0}_{\mu}$. The case with $|D|^{\frac{1}{2}}\varphi$ in place of $u_{\mu}$ is identical, except one uses Proposition \ref{phiu-bootstrap} in place of the bootstrap hypotheses in section \ref{sec:bootstrap} in the proof below. We move through each term in \eqref{Dtaeqn} one by one. In our analysis, we will need to rewrite many of the above expressions in the appropriate flattened variables. We abbreviate $\eta_{<\lambda}$ by $\eta$ for the remainder of this subsection. We hereon adopt the notation of subsection \ref{subsec:flattenedcoords} using the graph $\eta_{<\lambda}$. We also set as shorthand
\begin{equation*}
U=(V_{<\lambda},B_{<\lambda}),\hspace{5mm}\zeta\coloneq\nabla_{z,x}\mathcal{H}\eta_{<\lambda}.
\end{equation*}
Before estimating each term in the expansion of $D_t^l a_{<\lambda}$, we need the following Proposition and corollary which we will use frequently
\begin{proposition}[Weighted harmonic-extension bounds]
\label{weighted-harmonic-extension-bounds}
Assume the bootstrap hypotheses in Section \ref{sec:bootstrap}. Let
$\lambda,\mu,\nu,\gamma$ be dyadic frequencies with
$\mu\gg\lambda$ and let $1<\sigma<2$. Let $c_\gamma$ denote the
frequency envelope for $u_0$. Moreover, assume
$\eta=\eta_{<\lambda}$ and
$U=U_{<\lambda}=(V_{<\lambda},B_{<\lambda})$. Then
\begin{equation}\label{weighted-zeta-gradient}
\begin{split}
\sup_{z\leq0}
M_{\gamma,2}^{\frac{\sigma}{2}}
\left(
u_\mu^{x_0},
D_iP_\gamma w_\zeta(z)
\right)
&\lesssim_\sigma
\epsilon^2
\gamma^{\frac34}
\mu^{-s_c}c_\gamma c_\mu
\left(1+\mu|x_0|\right)^{\frac12},
\\
\sup_{z\leq0}
M_{\gamma,2}^{\frac{\sigma}{2}}
\left(
u_\mu^{x_0},
P_\gamma w_\zeta(z)
\right)
&\lesssim_\sigma
\epsilon^2
\gamma^{-\frac14}
\mu^{-s_c}c_\gamma c_\mu
\left(1+\mu|x_0|\right)^{\frac12}.
\end{split}
\end{equation}
Moreover, we have
\begin{equation}\label{weighted-U-one-frequency}
\begin{split}
\sup_{z\leq0}
M_{\gamma,2}^{\frac{\sigma}{2}}
\left(
u_\mu^{x_0},
D_iP_\gamma w_U(z)
\right)
&\lesssim_\sigma
\epsilon^2
\gamma^{\frac14}
\mu^{-s_c}c_\mu
\left(1+\mu|x_0|\right)^{\frac12},
\\
\sup_{z\leq0}
M_{\gamma,2}^{\frac{\sigma}{2}}
\left(
u_\mu^{x_0},
P_\gamma w_U(z)
\right)
&\lesssim_\sigma
\epsilon^2
\gamma^{-\frac34}
\mu^{-s_c}c_\mu
\left(1+\mu|x_0|\right)^{\frac12}.
\end{split}
\end{equation}
Both of the differentiated bounds hold with
$\nabla_{z,x}$ in place of $D_i$.
\end{proposition}

Dyadically summing the first estimatein \eqref{weighted-U-one-frequency}, we obtain the following
useful corollary that we will refer to several times in the sequel.

\begin{corollary}\label{summationbilinear}
If $\frac32<\sigma<2$, then
\begin{equation}\label{weighted-U-low-summed}
\sup_{z\leq0}
M_{\nu,2}^{\frac{\sigma}{2}}
\left(
u_\mu^{x_0},
D_iP_{\lesssim\nu}w_U(z)
\right)
\lesssim_\sigma
\epsilon^2
\nu^{\frac14}
\mu^{-s_c}c_\mu
\left(1+\mu|x_0|\right)^{\frac12},
\end{equation}
whereas, if $1<\sigma<\frac32$, then
\begin{equation}\label{weighted-U-high-summed}
\sup_{z\leq0}
M_{\nu,2}^{\frac{\sigma}{2}}
\left(
u_\mu^{x_0},
D_iP_{\gtrsim\nu}w_U(z)
\right)
\lesssim_\sigma
\epsilon^2
\nu^{\frac14}
\mu^{-s_c}c_\mu
\left(1+\mu|x_0|\right)^{\frac12}.
\end{equation}
\end{corollary}

Now, we prove the proposition.

\begin{proof}
We first prove the gradient estimate in \eqref{weighted-zeta-gradient} with $\nabla_{z,x}$ in place of $D_i$. Let
$\widetilde P_\gamma$ be a sufficiently fattened annular
projection. We begin with the preliminary estimate
\begin{equation}\label{eta-output-gradient-bound}
\begin{split}
\sup_{z\leq0}
M_{\gamma,2}^{\frac{\sigma}{2}}
\left(
u_\mu^{x_0},
\widetilde P_\gamma\nabla_{z,x}w_\eta(z)
\right)\lesssim_\sigma
\epsilon^2
\gamma^{-\frac14}
\mu^{-s_c}c_\mu c_\gamma
\left(1+\mu|x_0|\right)^{\frac12}.
\end{split}
\end{equation}
Indeed, using sufficiently fattened cutoffs, we write $\eta
=
\widetilde P_\gamma\eta
+
\left(
P_{\ll\gamma}+P_{\gg\gamma}
\right)\eta.$
For the "matched-frequency" term, Proposition \ref{bilinearext} give
\begin{equation*}
\begin{split}
\sup_{z\leq0}
M_{\gamma,2}^{\frac{\sigma}{2}}
\left(
u_\mu^{x_0},
\widetilde P_\gamma\nabla_{z,x}
w_{\widetilde P_\gamma\eta}(z)
\right)&\lesssim
\gamma
M_{\gamma,2}^{\frac{\sigma}{2}}
\left(
u_\mu^{x_0},
\widetilde P_\gamma\eta
\right)+\sup_{z\leq0}
M_{\gamma,2}^{\frac{\sigma}{2}}
\left(
u_\mu^{x_0},
P_{\gg\gamma}\nabla_{z,x}
w_{\widetilde P_\gamma\eta}(z)
\right)
\end{split}
\end{equation*}
The first term on the right hand side is controlled by the right hand side of \eqref{eta-output-gradient-bound} by directly using the unbalanced bilinear estimate.  For the high-frequency tail term, we use the change-of-scale estimate in Proposition
\ref{basicMproperties} and Proposition \ref{weightedoffdiag} with $p=\infty$ and the slowly varying property to estimate
\begin{equation*}
\begin{split}
\sup_{z\leq0}
M_{\gamma,2}^{\frac{\sigma}{2}}
\left(
u_\mu^{x_0},
P_{\gg\gamma}\nabla_{z,x}
w_{\widetilde P_\gamma\eta}(z)
\right)&\lesssim
\sum_{\alpha\gg\gamma}
\left(
\frac{\gamma}{\alpha}
\right)^{1-\frac{\sigma}{2}}
\sup_{z\leq0}
M_{\alpha,2}^{\frac{\sigma}{2}}
\left(
u_\mu^{x_0},
P_\alpha\nabla_{z,x}
w_{\widetilde P_\gamma\eta}(z)
\right)
\\
&\lesssim_\sigma
\epsilon^2
\mu^{-s_c}c_\mu
\left(1+\mu|x_0|\right)^{\frac12}
\|\nabla_x\widetilde P_\gamma\eta\|_
{L_t^\infty\dot B_{\infty,1}^0}
\sum_{\alpha\gg\gamma}
\left(
\frac{\gamma}{\alpha}
\right)^{1-\frac{\sigma}{2}}
\alpha^{-\frac14}c_\alpha
\\
&\lesssim_\sigma
\epsilon^3
\gamma^{-\frac14}
\mu^{-s_c}c_\mu c_\gamma
\left(1+\mu|x_0|\right)^{\frac12}.
\end{split}
\end{equation*}
where we used $\sigma<2$ above. For the remaining boundary datum, Proposition
\ref{weightedoffdiag}, applied with $p=\infty$, gives
\begin{equation*}
\begin{split}
\sup_{z\leq0}
M_{\gamma,2}^{\frac{\sigma}{2}}
\left(
u_\mu^{x_0},
\widetilde P_\gamma\nabla_{z,x}
w_{(P_{\ll\gamma}+P_{\gg\gamma})\eta}(z)
\right)&\lesssim_\sigma
\epsilon^2
\gamma^{-\frac14}
\mu^{-s_c}c_\mu c_\gamma
\left(1+\mu|x_0|\right)^{\frac12}
\left\|
\nabla_x\eta
\right\|_{L_t^\infty\dot B_{\infty,1}^0}
\\
&\lesssim_\sigma
\epsilon^3
\gamma^{-\frac14}
\mu^{-s_c}c_\mu c_\gamma
\left(1+\mu|x_0|\right)^{\frac12}.
\end{split}
\end{equation*}
This proves \eqref{eta-output-gradient-bound}. Recall that the boundary trace $\zeta_j$ is defined by $\zeta_j
=
\left.D_jw_\eta\right|_{z=0}$. We next show that
\begin{equation}\label{zeta-boundary-pairing}
\begin{split}
M_{\gamma,2}^{\frac{\sigma}{2}}
\left(
u_\mu^{x_0},
\widetilde P_\gamma\zeta_j
\right)
\lesssim_\sigma
\epsilon^2
\gamma^{-\frac14}
\mu^{-s_c}c_\mu c_\gamma
\left(1+\mu|x_0|\right)^{\frac12}.
\end{split}
\end{equation}
First, let us take $j=1,2$. Fixing a sufficiently small $\delta>0$, we can decompose using the definition of $\zeta_j$ and $D_j$,
\begin{equation*}
\begin{split}
\widetilde P_\gamma\zeta_j
={}&
\widetilde P_\gamma\partial_j\eta
-\widetilde P_\gamma\left(
P_{<\delta\gamma}\partial_j\eta\,
\left.\partial_zw_\eta\right|_{z=0}
\right)-
\widetilde P_\gamma\left(
P_{\gtrsim\delta\gamma}\partial_j\eta\,
\left.\partial_zw_\eta\right|_{z=0}
\right).
\end{split}
\end{equation*}
The first term is bounded by \eqref{eta-output-gradient-bound},
evaluated at $z=0$. For the low-frequency coefficient, the support
properties give, after fattening $\widetilde P_\gamma$ if necessary,
\begin{equation*}
\begin{split}
\widetilde P_\gamma\left(
P_{<\delta\gamma}\partial_j\eta\,
\left.\partial_zw_\eta\right|_{z=0}
\right)=
\widetilde P_\gamma\left(
P_{<\delta\gamma}\partial_j\eta\,
\left.\widetilde P_\gamma\partial_zw_\eta\right|_{z=0}
\right).
\end{split}
\end{equation*}
Consequently, boundedness of $\tilde{P}_{\gamma}$ with respect to $M_{\gamma,2}^{\frac{\sigma}{2}}$ and \eqref{eta-output-gradient-bound} implies
\begin{equation*}
\begin{split}
M_{\gamma,2}^{\frac{\sigma}{2}}
\left(
u_\mu^{x_0},
\widetilde P_\gamma\left(
P_{<\delta\gamma}\partial_j\eta\,
\left.\partial_zw_\eta\right|_{z=0}
\right)
\right)&\lesssim
\|\nabla\eta\|_{L_{t,x}^\infty}
\sup_{z\leq0}
M_{\gamma,2}^{\frac{\sigma}{2}}
\left(
u_\mu^{x_0},
\widetilde P_\gamma\partial_zw_\eta(z)
\right)
\\
&\lesssim_\sigma
\epsilon^3
\gamma^{-\frac14}
\mu^{-s_c}c_\mu c_\gamma
\left(1+\mu|x_0|\right)^{\frac12}.
\end{split}
\end{equation*}
For the high-frequency coefficient, the product estimate \eqref{weighted-output-product},
the change-of-scale estimate in Proposition
\ref{basicMproperties}, and the unbalanced bilinear estimate give
\begin{equation*}
\begin{split}
M_{\gamma,2}^{\frac{\sigma}{2}}
\left(
u_\mu^{x_0},
\widetilde P_\gamma\left(
P_{\gtrsim\delta\gamma}\partial_j\eta\,
\left.\partial_zw_\eta\right|_{z=0}
\right)
\right)&\lesssim
\|\partial_zw_\eta\|_{L_z^\infty L_{t,x}^\infty}
\sum_{\kappa\gtrsim\gamma}
\left(
\frac{\gamma}{\kappa}
\right)^{1-\frac{\sigma}{2}}
M_{\kappa,2}^{\frac{\sigma}{2}}
\left(
u_\mu^{x_0},
\nabla_xP_\kappa\eta
\right)
\\
&\quad\lesssim
\epsilon^3
\gamma^{-\frac14}
\mu^{-s_c}c_\mu
\left(1+\mu|x_0|\right)^{\frac12}
\sum_{\kappa\gtrsim\gamma}
\left(
\frac{\gamma}{\kappa}
\right)^{\frac54-\frac{\sigma}{2}}
c_\kappa
\\
&\quad\lesssim_\sigma
\epsilon^3
\gamma^{-\frac14}
\mu^{-s_c}c_\mu c_\gamma
\left(1+\mu|x_0|\right)^{\frac12}.
\end{split}
\end{equation*}
where above, we used the extension bound in Proposition \ref{flattenedbasicbounds}:
\begin{equation*}
\|\partial_zw_\eta\|_{L_z^\infty L_{t,x}^\infty}
\lesssim
\|\nabla_x\eta\|_{L_t^\infty\dot B_{\infty,1}^0}
\lesssim\epsilon,
\end{equation*}
and that the final dyadic sum converges because $\sigma<2$. This proves
\eqref{zeta-boundary-pairing} for $j=1,2$; for $j=3$, it follows
directly from \eqref{eta-output-gradient-bound} because $\zeta_3=\partial_zw_{\eta}|_{z=0}$. We now decompose the boundary datum itself as
\begin{equation*}
\zeta_j
=
\widetilde P_\gamma\zeta_j
+
\left(
P_{\ll\gamma}+P_{\gg\gamma}
\right)\zeta_j.
\end{equation*}
For the matched-frequency term, Proposition \ref{bilinearext} gives
\begin{equation*}
\begin{split}
\sup_{z\leq0}
M_{\gamma,2}^{\frac{\sigma}{2}}
\left(
u_\mu^{x_0},
\nabla_{z,x}P_\gamma
w_{\widetilde P_\gamma\zeta_j}(z)
\right)&\lesssim
\gamma
M_{\gamma,2}^{\frac{\sigma}{2}}
\left(
u_\mu^{x_0},
\widetilde P_\gamma\zeta_j
\right)+\sup_{z\leq0}
M_{\gamma,2}^{\frac{\sigma}{2}}
\left(
u_\mu^{x_0},
P_{\gg\gamma}\nabla_{z,x}
w_{\widetilde P_\gamma\zeta_j}(z)
\right).
\end{split}
\end{equation*}
The first term can be estimated by the right-hand side of the first estimate in \eqref{weighted-zeta-gradient} using \eqref{zeta-boundary-pairing}. To estimate the additional high-output term, we use the change-of-scale
estimate in Proposition \ref{basicMproperties}, sum over its output
frequencies, and apply Proposition \ref{weightedoffdiag} with
$p=\infty$:
\begin{equation*}
\begin{split}
\sup_{z\leq0}
M_{\gamma,2}^{\frac{\sigma}{2}}
\left(
u_\mu^{x_0},
P_{\gg\gamma}\nabla_{z,x}
w_{\widetilde P_\gamma\zeta_j}(z)
\right)&\lesssim
\sum_{\alpha\gg\gamma}
\left(
\frac{\gamma}{\alpha}
\right)^{1-\frac{\sigma}{2}}
\sup_{z\leq0}
M_{\alpha,2}^{\frac{\sigma}{2}}
\left(
u_\mu^{x_0},
P_\alpha\nabla_{z,x}
w_{\widetilde P_\gamma\zeta_j}(z)
\right)
\\
&\quad\lesssim_\sigma
\epsilon^2
\mu^{-s_c}c_\mu
\left(1+\mu|x_0|\right)^{\frac12}
\|\nabla_x\widetilde P_\gamma\zeta_j\|_
{L_t^\infty\dot B_{\infty,1}^0}
\sum_{\alpha\gg\gamma}
\left(
\frac{\gamma}{\alpha}
\right)^{1-\frac{\sigma}{2}}
\alpha^{-\frac14}c_\alpha
\\
&\quad\lesssim_\sigma
\epsilon^3
\gamma^{\frac34}
\mu^{-s_c}c_\mu c_\gamma
\left(1+\mu|x_0|\right)^{\frac12}.
\end{split}
\end{equation*}
Here we used Bernstein and the bootstrap bound to obtain
\begin{equation*}
\|\nabla_x\widetilde P_\gamma\zeta_j\|_
{L_t^\infty\dot B_{\infty,1}^0}
\lesssim
\gamma
\|\nabla_x\widetilde P_\gamma\zeta_j\|_
{L_t^\infty L_x^2}
\lesssim
\epsilon\gamma,
\end{equation*}
and used the same summation as in the proof of
\eqref{eta-output-gradient-bound}. The remaining boundary datum is supported away from frequency
$\gamma$. Applying Proposition \ref{weightedoffdiag} componentwise
with $p=2$, we obtain
\begin{equation*}
\begin{split}
\sup_{z\leq0}
M_{\gamma,2}^{\frac{\sigma}{2}}
\left(
u_\mu^{x_0},
\nabla_{z,x}P_\gamma
w_{(P_{\ll\gamma}+P_{\gg\gamma})\zeta_j}(z)
\right)&\lesssim_\sigma
\epsilon^2
\gamma^{-\frac14+\frac22}
\mu^{-s_c}c_\mu c_\gamma
\left(1+\mu|x_0|\right)^{\frac12}
\|\nabla_x\zeta|_{z=0}\|_{L_t^\infty\dot B_{2,1}^0}
\\
&\qquad\lesssim
\epsilon^3
\gamma^{\frac34}
\mu^{-s_c}c_\mu c_\gamma
\left(1+\mu|x_0|\right)^{\frac12}.
\end{split}
\end{equation*}
Above, we used the bound
\begin{equation}\label{zetabound}
\|\nabla_x\zeta|_{z=\eta(t,x)}\|_{L_t^\infty\dot B_{2,1}^0}
\lesssim\epsilon.
\end{equation}
This bound follows from the identity $\zeta(x,\eta(x))=Dw_{\eta}(x,0)$ and using Proposition \ref{besovparabounds} to obtain
\begin{equation*}
\|\nabla_{x}Dw_{\eta}(x,0)\|_{L_t^{\infty}\dot{B}_{2,1}^0}\lesssim \|\nabla^2_x\eta\|_{L_t^{\infty}\dot{B}_{2,1}^0}\|\nabla_{z,x}w_{\eta}\|_{L_z^{\infty}L_t^{\infty}\dot{B}_{\infty,1}^0}+\|\nabla_{z,x}^2w_{\eta}\|_{L_z^{\infty}L_t^{\infty}\dot{B}_{2,1}^0}    
\end{equation*}
and then using the bounds
\begin{equation*}
\|\nabla_{z,x}w_{\eta}\|_{L_z^{\infty}L_t^{\infty}\dot{B}_{\infty,1}^0}+\|\nabla_{z,x}^2w_{\eta}\|_{L_z^{\infty}L_t^{\infty}\dot{B}_{2,1}^0}\lesssim \|\nabla_x\eta\|_{L_t^{\infty}\dot{B}_{\infty,1}^0}+\|\nabla_x^2\eta\|_{L_t^{\infty}\dot{B}_{2,1}^0}\lesssim\epsilon.  
\end{equation*}
The $\dot{B}_{\infty,1}^0$ bound for $\nabla_{z,x}w_{\eta}$ above follows from Proposition \ref{flattenedbasicbounds}, while the $\dot{B}_{2,1}^0$ bound for $\nabla_{z,x}^2w_{\eta}$ follows from a slight modification of the same proposition (using exactly the same proof strategy).
\\

Now, combining the above matched and off-diagonal estimates proves
\begin{equation*}
\begin{split}
\sup_{z\leq0}
M_{\gamma,2}^{\frac{\sigma}{2}}
\left(
u_\mu^{x_0},
\nabla_{z,x}P_\gamma w_\zeta(z)
\right)&\lesssim_\sigma
\epsilon^2
\gamma^{\frac34}
\mu^{-s_c}c_\mu c_\gamma
\left(1+\mu|x_0|\right)^{\frac12}.
\end{split}
\end{equation*}

We next prove the gradient version of \eqref{weighted-U-one-frequency}. Using sufficiently
fattened cutoffs and the linearity of the harmonic extension in its
boundary datum, we write
\begin{equation*}
P_\gamma w_U
=
P_\gamma w_{\widetilde P_\gamma U}
+
P_\gamma w_{(P_{\ll\gamma}+P_{\gg\gamma})U}.
\end{equation*}
For the matched-frequency term, Proposition \ref{bilinearext} gives
\begin{equation}\label{matchedU}
\begin{split}
&\sup_{z\leq0}
M_{\gamma,2}^{\frac{\sigma}{2}}
\left(
u_\mu^{x_0},
\nabla_{z,x}P_\gamma w_{\widetilde P_\gamma U}(z)
\right)\lesssim
\gamma
M_{\gamma,2}^{\frac{\sigma}{2}}
\left(
u_\mu^{x_0},
\widetilde P_\gamma U
\right)+\sup_{z\leq0}
M_{\gamma,2}^{\frac{\sigma}{2}}
\left(
u_\mu^{x_0},
P_{\gg\gamma}\nabla_{z,x}
w_{\widetilde P_\gamma U}(z)
\right).
\end{split}
\end{equation}
If $\gamma\lesssim\lambda$, Proposition \ref{VBunbil} gives 
\begin{equation*}
\gamma
M_{\gamma,2}^{\frac{\sigma}{2}}
\left(
u_\mu^{x_0},
\widetilde P_\gamma U
\right)
\lesssim_\sigma
\epsilon^2
\gamma^{\frac14}
\mu^{-s_c}c_\mu
\left(1+\mu|x_0|\right)^{\frac12}.
\end{equation*}
Suppose instead that $\gamma\gg\lambda$. Choose $q$ so that $4<q<\frac{4}{2-\sigma}$,
which is possible because $\sigma>1$. Proposition
\ref{basicMproperties} and the flattened harmonic-extension bound Proposition \ref{flattenedbasicbounds}
give 
\begin{equation*}
\begin{split}
\sup_{z\leq0}
M_{\gamma,2}^{\frac{\sigma}{2}}
\left(
u_\mu^{x_0},
\nabla_{z,x}P_\gamma w_{\widetilde P_\gamma U}(z)
\right)\lesssim_{\sigma,q}
\gamma^{\frac2q}
\|u_\mu\|_{L_t^\infty L_x^2}
\|\nabla_xP_\gamma U\|_{L_t^2L_x^q}.
\end{split}
\end{equation*}
By Propositions \ref{Bernsteinelliptic} and \ref{VBunbil},
\begin{equation*}
\begin{split}
\gamma^{\frac2q}
\|\nabla_xP_\gamma U\|_{L_t^2L_x^q}
&\lesssim
\gamma^{\frac2q-\frac58}
\|\nabla^2\eta_{<\lambda}\|_{L_t^4\dot B_{q,1}^0}
\||D|^{\frac58}U\|_{L_t^4\dot B_{\infty,1}^0}
\\
&\lesssim
\epsilon^2
\gamma^{\frac2q-\frac58}
\lambda^{\frac78-\frac2q}c_\lambda
\\
&\lesssim
\epsilon^2
\gamma^{\frac14}c_\gamma.
\end{split}
\end{equation*}
Here the last inequality follows from the slowly varying property of
the frequency envelope and
$\frac78-\frac2q>0$. The additional high-output term in \eqref{matchedU} is estimated using the change-of-scale estimate in
Proposition \ref{basicMproperties} and then summing in its output frequency and then using Proposition \ref{weightedoffdiag} with $p=\infty$,
\begin{equation*}
\begin{split}
\sup_{z\leq0}
M_{\gamma,2}^{\frac{\sigma}{2}}
\left(
u_\mu^{x_0},
P_{\gg\gamma}\nabla_{z,x}
w_{\widetilde P_\gamma U}(z)
\right)&\lesssim
\sum_{\alpha\gg\gamma}
\left(
\frac{\gamma}{\alpha}
\right)^{1-\frac{\sigma}{2}}
\sup_{z\leq0}
M_{\alpha,2}^{\frac{\sigma}{2}}
\left(
u_\mu^{x_0},
P_\alpha\nabla_{z,x}
w_{\widetilde P_\gamma U}(z)
\right)
\\
&\quad\lesssim_\sigma
\epsilon^2
\mu^{-s_c}c_\mu
\left(1+\mu|x_0|\right)^{\frac12}
\|\nabla_x\widetilde P_\gamma U\|_
{L_t^\infty\dot B_{\infty,1}^0}
\sum_{\alpha\gg\gamma}
\left(
\frac{\gamma}{\alpha}
\right)^{1-\frac{\sigma}{2}}
\alpha^{-\frac14}c_\alpha
\\
&\quad\lesssim_\sigma
\epsilon^3
\gamma^{\frac14}
\mu^{-s_c}c_\mu c_\gamma
\left(1+\mu|x_0|\right)^{\frac12}
\\
&\quad\lesssim_\sigma
\epsilon^2
\gamma^{\frac14}
\mu^{-s_c}c_\mu
\left(1+\mu|x_0|\right)^{\frac12}.
\end{split}
\end{equation*}
In the third line, we used Bernstein and Proposition \ref{VBunbil}. Combining the above estimates, we then obtain
\begin{equation}\label{weighted-U-matched-part}
\begin{split}
&\sup_{z\leq0}
M_{\gamma,2}^{\frac{\sigma}{2}}
\left(
u_\mu^{x_0},
\nabla_{z,x}P_\gamma w_{\widetilde P_\gamma U}(z)
\right)\lesssim_\sigma
\epsilon^2
\gamma^{\frac14}
\mu^{-s_c}c_\mu
\left(1+\mu|x_0|\right)^{\frac12}.
\end{split}
\end{equation}
The remaining boundary datum is supported away from frequency
$\gamma$. Proposition \ref{weightedoffdiag}, applied 
with $p=4$, and Sobolev embedding along with Proposition \ref{VBunbil} give
\begin{equation*}
\begin{split}
\sup_{z\leq0}
M_{\gamma,2}^{\frac{\sigma}{2}}
\left(
u_\mu^{x_0},
\nabla_{z,x}P_\gamma
w_{(P_{\ll\gamma}+P_{\gg\gamma})U}(z)
\right)&\lesssim_\sigma
\epsilon^2
\gamma^{-\frac14+\frac24}
\mu^{-s_c}c_\mu c_\gamma
\left(1+\mu|x_0|\right)^{\frac12}
\left\|
\nabla_xU
\right\|_{L_t^\infty\dot B_{4,1}^0}
\\
&\qquad\lesssim
\epsilon^3
\gamma^{\frac14}
\mu^{-s_c}c_\mu c_\gamma
\left(1+\mu|x_0|\right)^{\frac12}
\\
&\qquad\lesssim
\epsilon^2
\gamma^{\frac14}
\mu^{-s_c}c_\mu
\left(1+\mu|x_0|\right)^{\frac12}.
\end{split}
\end{equation*}
This proves the differentiated $U$-estimate with $\nabla_{z,x}$. Since
\begin{equation*}
D_iP_\gamma w
=
\partial_iP_\gamma w
-(\partial_i\eta)\partial_zP_\gamma w
\end{equation*}
and $\|\nabla\eta\|_{L_{t,x}^\infty}\ll1$, both differentiated
estimates also hold with $D_i$ in place of $\nabla_{z,x}$. Finally,
boundedness of Littlewood-Paley projections in the bilinear norm gives, for
$w=w_U$ or $w=w_\zeta$,
\begin{equation*}
\begin{split}
M_{\gamma,2}^{\frac{\sigma}{2}}
\left(
u_\mu^{x_0},
P_\gamma w
\right)
&\lesssim
\gamma^{-1}
M_{\gamma,2}^{\frac{\sigma}{2}}
\left(
u_\mu^{x_0},
\nabla_xP_\gamma w
\right).
\end{split}
\end{equation*}
The two undifferentiated estimates now follow from their
differentiated counterparts. This completes the proof.
\end{proof}
\begin{proof}[Proof of Proposition \ref{Dtaestimates}]
Now, we proceed with Proposition \ref{Dtaestimates}. We go term by term and will first work to establish the more difficult multilinear bound. After proving this, we will turn to the $L_t^2L_x^8$.
\subsubsection{The term
$-4\mathcal{P}(\nabla_{z,x}^2\phi_{<\lambda}:
\nabla_{z,x}^2\mathcal{H}\eta_{<\lambda})$} This is the most delicate term since it is quadratic at the principal level, and we we will need to use two full unbalanced bilinear estimates. We begin by writing the source term $\nabla_{z,x}^2\phi_{<\lambda}:\nabla_{z,x}^2\mathcal{H}\eta_{<\lambda}$ in the flattened variables. The vector field
\(\nabla_{x,z}\phi_{<\lambda}\) is componentwise harmonic and has
boundary trace $U$. Consequently, by
uniqueness of the harmonic extension, its pullback is $w_U$. Hence, $(\nabla_{z,x}^2\phi_{<\lambda})_*=Dw_U$. We also have  $(\nabla_{z,x}^2\mathcal{H}\eta_{<\lambda})_*=Dw_{\zeta}$. Therefore, in flattened coordinates, we are estimate (by slight abuse of notation)
\begin{equation*}
-4\mathcal P\left(
D_iw_{U_j}D_iw_{\zeta_j}
\right).
\end{equation*}
Since $D_iD_iw_{U_j}=0$, we therefore have
\begin{equation*}
D_iw_{U_j}D_iP_\beta w_{\zeta_j}
=
D_i\left(
D_iw_{U_j}P_\beta w_{\zeta_j}
\right).
\end{equation*}
Using sufficiently fat cutoffs, we decompose the term as
\begin{equation}\label{eta-output-decomposition}
\begin{split}
\mathcal E_{\mathrm{lo}}
&=
-4\sum_{\beta}
\left.
\mathcal{P}P_{<\beta}D_i
\left(
D_iw_{U_j}P_\beta w_{\zeta_j}
\right)
\right|_{z=\eta_{<\lambda}},
\\
\mathcal E_{\mathrm{out}}
&=
-4\sum_{\beta}
\sum_{\alpha\gtrsim\beta}
\left.
\mathcal{P}P_\alpha
\left(
D_iw_{U_j}D_iP_\beta w_{\zeta_j}
\right)
\right|_{z=\eta_{<\lambda}}.
\end{split}
\end{equation}
Thus the derivative is pulled onto the output only in
$\mathcal E_{\mathrm{lo}}$. Now, we fix several parameters, whose motivation will be come clear momentarily, $0<\delta<\frac18$ and set
\begin{equation*}
\sigma_0^+
=
\frac34+\delta,
\qquad
\sigma_0^-
=
\frac34-\delta,
\qquad
\sigma_1
= \frac{3}{4}.
\end{equation*}
Then $1<\sigma_0^\pm+\sigma_1<2$. Applying Corollary \ref{summationbilinear} with parameters
$2\sigma_0^+=\frac32+2\delta$ and
$2\sigma_0^-=\frac32-2\delta$, respectively, and taking $x_0=0$,
we obtain
\begin{equation}\label{eta-two-sigma-U-bounds}
\begin{split}
\sup_{z\leq0}
M_{\nu,2}^{\sigma_0^+}
\left(
u_\mu,
D_iP_{\lesssim\nu}w_U(z)
\right)\lesssim
\epsilon^2\nu^{\frac14}\mu^{-s_c}c_\mu,\hspace{5mm}\sup_{z\leq0}
M_{\nu,2}^{\sigma_0^-}
\left(
u_\mu,
D_iP_{\gtrsim\nu}w_U(z)
\right)\lesssim
\epsilon^2\nu^{\frac14}\mu^{-s_c}c_\mu.
\end{split}
\end{equation}
On the other hand, Proposition
\ref{weighted-harmonic-extension-bounds}, applied with parameter
$\sigma=\frac32$ gives
\begin{equation}\label{eta-zeta-own-scale}
\begin{split}
\sup_{z\leq0}
M_{\beta,2}^{\sigma_1}
\left(
u_\mu^{x_0},
P_\beta w_{\zeta_j}(z)
\right)\lesssim
\epsilon^2\beta^{-\frac14}
\mu^{-s_c}c_\mu c_\beta
\left(1+\mu|x_0|\right)^{\frac12}.
\end{split}
\end{equation}
Now, we estimate $\mathcal E_{\mathrm{lo}}$. We begin by splitting each summand inside $D_i$ at frequency
$\beta$:
\begin{equation*}
D_iw_{U_j}
=
D_iP_{\lesssim\beta}w_{U_j}
+
D_iP_{\gtrsim\beta}w_{U_j}.
\end{equation*}
Using Proposition \ref{weightedbilinear},
\eqref{eta-two-sigma-U-bounds}, and
\eqref{eta-zeta-own-scale} give
\begin{equation}\label{eta-low-output-estimate}
\begin{split}
\left\|
L\left(
h_\mu,
\mathcal E_{\mathrm{lo}}
\right)
\right\|_{L_{t,x}^1}
&\lesssim
\sum_{\beta}
\sup_{z\leq0}
M_{\beta,2}^{\sigma_0^+}
\left(
u_\mu,
D_iP_{\lesssim\beta}w_U(z)
\right)
M_{\beta,2}^{\sigma_1}
\left(
u_\mu^{x_0},
P_\beta w_\zeta(z)
\right)
\\
&+
\sum_{\beta}
\sup_{z\leq0}
M_{\beta,2}^{\sigma_0^-}
\left(
u_\mu,
D_iP_{\gtrsim\beta}w_U(z)
\right)
M_{\beta,2}^{\sigma_1}
\left(
u_\mu^{x_0},
P_\beta w_\zeta(z)
\right)
\\
&\lesssim
\epsilon^4
\mu^{-2s_c}c_\mu^2
\left(1+\mu|x_0|\right)^{\frac12}
\end{split}
\end{equation}
where we used $\|c_{\beta}\|_{l^1_{\beta}}\lesssim 1$. We next estimate $\mathcal E_{\mathrm{out}}$. By writing $D_i=\partial_i+\partial_i\eta\partial_z$ and using $\|\nabla\eta\|_{L^\infty_{t,x,z}}\lesssim\epsilon$, we can discard the coefficients $\partial_i\eta$ from the weighted norm. For $\alpha\gtrsim\beta$, the change-of-scale
estimate in Proposition \ref{basicMproperties}, applied to the
flat, $\beta$-localized derivatives, gives
\begin{equation}\label{eta-zeta-off-scale-gradient}
\begin{split}
\sup_{z\leq0}
M_{\alpha,2}^{\sigma_1}
\left(
u_\mu^{x_0},
D_iP_\beta w_{\zeta_j}(z)
\right)
&\lesssim
\sup_{z\leq0}
M_{\alpha,2}^{\sigma_1}
\left(
u_\mu^{x_0},
\nabla_{z,x}P_\beta w_{\zeta_j}(z)
\right)
\\
&\lesssim
\alpha^{\frac{1}{4}}\beta^{-\frac{1}{4}}
\sup_{z\leq0}
M_{\beta,2}^{\sigma_1}
\left(
u_\mu^{x_0},
\nabla_{z,x}P_\beta w_{\zeta_j}(z)
\right)
\\
&\lesssim
\epsilon^2\alpha^{\frac{1}{4}}\beta^{\frac12}
\mu^{-s_c}c_\mu c_\beta
\left(1+\mu|x_0|\right)^{\frac12}.
\end{split}
\end{equation}
where we used Proposition
\ref{weighted-harmonic-extension-bounds} in the final line. At the output frequency $\alpha$, split inside $D_i$ as
\begin{equation*}
D_iw_{U_j}
=
D_iP_{\lesssim\alpha}w_{U_j}
+
D_iP_{\gtrsim\alpha}w_{U_j}.
\end{equation*}
For every $\alpha\gtrsim\beta$, Proposition \ref{weightedbilinear} supplies the factor $\alpha^{-1}$. Therefore,
Proposition \ref{weightedbilinear},
\eqref{eta-two-sigma-U-bounds}, and
\eqref{eta-zeta-off-scale-gradient} shows that $\left\|
L\left(
h_\mu,
\mathcal E_{\mathrm{out}}
\right)
\right\|_{L_{t,x}^1}$ can be estimated by
\begin{equation}\label{eta-high-output-estimate}
\begin{split}
&\sum_{\beta}
\sum_{\alpha\gtrsim\beta}
\alpha^{-1}
\sup_{z\leq0}
\Bigg[
M_{\alpha,2}^{\sigma_0^+}
\left(
u_\mu,
D_iP_{\lesssim\alpha}w_U(z)
\right)+
M_{\alpha,2}^{\sigma_0^-}
\left(
u_\mu,
D_iP_{\gtrsim\alpha}w_U(z)
\right)
\Bigg]M_{\alpha,2}^{\sigma_1}
\left(
u_\mu^{x_0},
D_iP_\beta w_\zeta(z)
\right)
\\
&\qquad\lesssim
\epsilon^4
\mu^{-2s_c}c_\mu^2
\left(1+\mu|x_0|\right)^{\frac12}
\sum_\beta c_\beta
\sum_{\alpha\gtrsim\beta}
\alpha^{-\frac{1}{2}}\beta^{\frac12}
\\
&\qquad\lesssim
\epsilon^4
\mu^{-2s_c}c_\mu^2
\left(1+\mu|x_0|\right)^{\frac12}.
\end{split}
\end{equation}
Combining \eqref{eta-low-output-estimate} and
\eqref{eta-high-output-estimate}, we conclude that
\begin{equation*}
\left\|
L\left(
h_\mu,
-4\mathcal P\left(
\nabla_{z,x}^2\phi_{<\lambda}:
\nabla_{z,x}^2\mathcal H\eta_{<\lambda}
\right)
\right)
\right\|_{L_{t,x}^1}
\lesssim
\epsilon^4
\mu^{-2s_c}c_\mu^2
\left(1+\mu|x_0|\right)^{\frac12}.
\end{equation*}
\subsubsection{The term
$-2\mathcal{P}(\mathrm{tr}((\nabla_{z,x}^2\phi_{<\lambda})^3))$}
We use the symmetrized frequency decomposition convention from
\cite[(7.5)]{ifrim2025sharp}, (specifically in terms of the ordering of the
indices rather than the explicit choice of regularization used there). We make
this decomposition directly in the flattened variables. Thus,
symmetrizing, we write
\begin{equation*}
D_jw_{U_k}D_kw_{U_i}D_iw_{U_j}
=
3\sum_{\beta}
D_jP_\beta w_{U_k}
D_kP_{<\beta}w_{U_i}
D_iP_{<\beta}w_{U_j},
\end{equation*}
where the sum is understood in the symmetrized sense of
\cite[(7.5)]{ifrim2025sharp}, including the comparable-frequency
interactions. We now essentially repeat the computation in
\cite[(7.9)]{ifrim2025sharp} in order to place at least one
derivative on the output, while ensuring that the lower-frequency
factors remain differentiated. Since the derivatives $D_i$
commute, we have
\begin{equation*}
\begin{split}
&D_jP_\beta w_{U_k}
D_kP_{<\beta}w_{U_i}
D_iP_{<\beta}w_{U_j}
\\
&=
D_j\left(
P_\beta w_{U_k}
D_kP_{<\beta}w_{U_i}
D_iP_{<\beta}w_{U_j}
\right)-
P_\beta w_{U_k}
D_kD_jP_{<\beta}w_{U_i}
D_iP_{<\beta}w_{U_j}
\\
&
-
P_\beta w_{U_k}
D_kP_{<\beta}w_{U_i}
D_i\left(
D_jP_{<\beta}w_{U_j}
\right)
\\
&=
D_j\left(
P_\beta w_{U_k}
D_kP_{<\beta}w_{U_i}
D_iP_{<\beta}w_{U_j}
\right)-
\frac12D_k\left(
P_\beta w_{U_k}
D_jP_{<\beta}w_{U_i}
D_iP_{<\beta}w_{U_j}
\right)
\\
&+
\frac12D_kP_\beta w_{U_k}
D_jP_{<\beta}w_{U_i}
D_iP_{<\beta}w_{U_j}-P_\beta w_{U_k}
D_kP_{<\beta}w_{U_i}
D_i\left(
D_jP_{<\beta}w_{U_j}
\right).
\end{split}
\end{equation*}
The first two terms have the desired form. Consequently, suppressing
pullbacks by a slight abuse of notation, we can write schematically
\begin{equation}\label{cubictraceid}
\begin{split}
-2\mathcal P
\left(
\operatorname{tr}
\left((\nabla_{z,x}^2\phi_{<\lambda})^3\right)
\right)
={}&
\sum_{\beta}
\mathcal P D
L\left(
P_\beta w_U,
DP_{<\beta}w_U,
DP_{<\beta}w_U
\right)
+
\mathcal E_4,
\end{split}
\end{equation}
where
\begin{equation*}
\begin{split}
\mathcal E_4
=&
-3\sum_{\beta}
\mathcal P\left(
D_kP_\beta w_{U_k}
D_jP_{<\beta}w_{U_i}
D_iP_{<\beta}w_{U_j}
\right)
\\
&+
6\sum_{\beta}
\mathcal P\left(
P_\beta w_{U_k}
D_kP_{<\beta}w_{U_i}
D_i\left(
D_jP_{<\beta}w_{U_j}
\right)
\right).
\end{split}
\end{equation*}

Our next goal is to show that $\mathcal E_4$ has the same
output-divergence structure. Define the "divergence defects"
\begin{equation*}
d_\beta
\coloneq
D_kP_\beta w_{U_k},
\qquad
d_{<\beta}
\coloneq
D_jP_{<\beta}w_{U_j}.
\end{equation*}
These arise because the flattened frequency projected pieces are not
individually divergence-free. We compute
\begin{equation}\label{quarticdivfree}
\begin{split}
&
-3d_\beta
D_jP_{<\beta}w_{U_i}
D_iP_{<\beta}w_{U_j}
+
6P_\beta w_{U_k}
D_kP_{<\beta}w_{U_i}
D_id_{<\beta}
\\
&\qquad=
D_i\left(
6P_\beta w_{U_k}
D_kP_{<\beta}w_{U_i}
d_{<\beta}
-
3P_\beta w_{U_i}d_{<\beta}^2
\right)
\\
&\qquad
-
3\bigg(
d_\beta
\left(
D_jP_{<\beta}w_{U_i}
D_iP_{<\beta}w_{U_j}
-
d_{<\beta}^2
\right)+
2d_{<\beta}
D_iP_\beta w_{U_k}
D_kP_{<\beta}w_{U_i}
\bigg),
\end{split}
\end{equation}
where we used
\begin{equation*}
D_iD_kP_{<\beta}w_{U_i}
=
D_kd_{<\beta}.
\end{equation*}
It remains to observe that the last line of
\eqref{quarticdivfree} vanishes after summation in $\beta$. Define
\begin{equation*}
\Psi(V)
\coloneq
\left(D_kV_k\right)
D_jV_iD_iV_j
-
\frac13\left(D_kV_k\right)^3.
\end{equation*}
Since $w_U$ is divergence-free in the sense that
$D_iw_{U_i}=0$, we have $\Psi(w_U)=0$. On the other hand, the
symmetrized trilinear convention of
\cite[(7.5)]{ifrim2025sharp} gives
\begin{equation*}
\begin{split}
\Psi(w_U)
=
\sum_{\beta}
\bigg(d_\beta
\left(
D_jP_{<\beta}w_{U_i}
D_iP_{<\beta}w_{U_j}
-
d_{<\beta}^2
\right)+
2d_{<\beta}
D_iP_\beta w_{U_k}
D_kP_{<\beta}w_{U_i}
\bigg).
\end{split}
\end{equation*}
Thus the sum of the last two lines of
\eqref{quarticdivfree} vanishes. Consequently,
\begin{equation}\label{quartic-output-divergence}
\mathcal E_4
=
\mathcal P D_i\mathcal Q_{4,i},
\end{equation}
where
\begin{equation}\label{quartic-flux}
\mathcal Q_{4,i}
\coloneq
\sum_{\beta}
\left(
6P_\beta w_{U_k}
D_kP_{<\beta}w_{U_i}
d_{<\beta}
-
3P_\beta w_{U_i}d_{<\beta}^2
\right).
\end{equation}
We now estimate the full expression. Since the identities above
depend only on the symmetrized ordering and not on the particular
choice of regularization, we apply them once more using a
Littlewood--Paley decomposition in the flattened variables. Thus,
below, $\beta$ denotes the largest actual flattened frequency and
ranges over all dyadic numbers. 
\\

For outputs below $\beta$, we use the divergence identities
\eqref{cubictraceid} and \eqref{quartic-output-divergence}. For
outputs comparable to or larger than $\beta$, we instead retain the
original cubic trace since there is no gain from the output. Accordingly, we write
\begin{equation*}
-2\mathcal P
\left(
\operatorname{tr}
\left((\nabla_{z,x}^2\phi_{<\lambda})^3\right)
\right)
=
\mathcal C_{3,\mathrm{lo}}
+
\mathcal E_{4,\mathrm{lo}}
+
\mathcal C_{\mathrm{hi}},
\end{equation*}
where
\begin{equation*}
\begin{split}
\mathcal C_{3,\mathrm{lo}}
&\coloneq
\sum_\beta
\mathcal P P_{<\beta}D
L\left(
P_\beta w_U,
DP_{<\beta}w_U,
DP_{<\beta}w_U
\right),
\\
\mathcal E_{4,\mathrm{lo}}
&\coloneq
\sum_\beta
\mathcal P P_{<\beta}D_i\mathcal Q_{4,\beta,i},
\\
\mathcal C_{\mathrm{hi}}
&\coloneq
\sum_\beta\sum_{\alpha\gtrsim\beta}
\mathcal P P_\alpha
L\left(
DP_\beta w_U,
DP_{<\beta}w_U,
DP_{<\beta}w_U
\right),
\end{split}
\end{equation*}
and $Q_{4,\beta,i}$ is the $\beta$ summand in $\mathcal{Q}_{4,i}$. The term
$\mathcal C_{\mathrm{hi}}$ is the combined high-output contribution
of $\mathcal C_3+\mathcal E_4$, since their sum agrees with the
original cubic trace. Now fix a finite $p>8$ and choose
\begin{equation*}
2-\frac8p<\sigma<2.
\end{equation*}
Proposition
\ref{weightedbilinear} gives
\begin{equation*}
\begin{split}
\left\|
L\left(
h_\mu,
\mathcal C_{3,\mathrm{lo}}
\right)
\right\|_{L_{t,x}^1}
&\lesssim
\sum_\beta
\sup_{z\leq0}
\Bigg[
M_{\beta,2}^{\frac{\sigma}{2}}
\left(
u_\mu^{x_0},
P_\beta w_U(z)
\right)
M_{\beta,2}^{\frac{\sigma}{2}}
\left(
u_\mu,
DP_{<\beta}w_U(z)
DP_{<\beta}w_U(z)
\right)
\Bigg].
\end{split}
\end{equation*}
Proposition
\ref{weighted-harmonic-extension-bounds} gives
\begin{equation*}
\sup_{z\leq0}
M_{\beta,2}^{\frac{\sigma}{2}}
\left(
u_\mu^{x_0},
P_\beta w_U(z)
\right)
\lesssim
\epsilon^2\beta^{-\frac34}
\mu^{-s_c}c_\mu
\left(1+\mu|x_0|\right)^{\frac12},
\end{equation*}
while Proposition
\ref{basicMproperties} gives
\begin{equation*}
\begin{split}
&
\sup_{z\leq0}
M_{\beta,2}^{\frac{\sigma}{2}}
\left(
u_\mu,
DP_{<\beta}w_U(z)
DP_{<\beta}w_U(z)
\right)\lesssim_p
\beta^{\frac4p}
\|u_\mu\|_{L_t^\infty L_x^2}
\left\|
DP_{<\beta}w_U
\right\|_{L_z^\infty L_t^4L_x^p}^2.
\end{split}
\end{equation*}
The latter estimate is valid because
$\sigma>2-\frac8p$. A slight adaptation of the harmonic-extension estimate Proposition \ref{flattenedbasicbounds}, Bernstein,
and Proposition \ref{VBunbil} give
\begin{equation}\label{cubic-trace-low-factor}
\left\|
DP_{<\beta}w_U
\right\|_{L_z^\infty L_t^4L_x^p}^2
\lesssim
\epsilon^2
\beta^{\frac34-\frac4p}
(c_\beta)^2.
\end{equation}
Using also $\|u_\mu\|_{L_t^\infty L_x^2}
\lesssim
\epsilon\mu^{-s_c}c_\mu$ and summing the $\beta$ frequency envelope, we obtain
\begin{equation}\label{cubic-trace-low-output}
\left\|
L\left(
h_\mu,
\mathcal C_{3,\mathrm{lo}}
\right)
\right\|_{L_{t,x}^1}
\lesssim_p
\epsilon^5
\mu^{-2s_c}c_\mu^2
\left(1+\mu|x_0|\right)^{\frac12}.
\end{equation}
The quartic flux has the same low-output estimate. Indeed,
\begin{equation*}
\begin{split}
\left\|
d_{<\beta}
\right\|_{L_z^\infty L_t^4L_x^p}
&\lesssim
\left\|
DP_{<\beta}w_U
\right\|_{L_z^\infty L_t^4L_x^p},
\\
\left\|
6DP_{<\beta}w_U\,d_{<\beta}
-
3d_{<\beta}^2
\right\|_{L_z^\infty L_t^2L_x^{p/2}}
&\lesssim
\left\|
DP_{<\beta}w_U
\right\|_{L_z^\infty L_t^4L_x^p}^2.
\end{split}
\end{equation*}
Hence, Proposition
\ref{weightedbilinear}, Proposition
\ref{weighted-harmonic-extension-bounds}, and
\eqref{cubic-trace-low-factor} give
\begin{equation}\label{quartic-trace-low-output}
\left\|
L\left(
h_\mu,
\mathcal E_{4,\mathrm{lo}}
\right)
\right\|_{L_{t,x}^1}
\lesssim_p
\epsilon^5
\mu^{-2s_c}c_\mu^2
\left(1+\mu|x_0|\right)^{\frac12}.
\end{equation}
It remains to estimate the high-output contribution. Proposition
\ref{weightedbilinear} and Proposition \ref{basicMproperties} gives
\begin{equation*}
\begin{split}
\left\|
L\left(
h_\mu,
\mathcal C_{\mathrm{hi}}
\right)
\right\|_{L_{t,x}^1}
&\lesssim
\sum_\beta\sum_{\alpha\gtrsim\beta}
\alpha^{-1}
\sup_{z\leq0}
\Bigg[
M_{\alpha,2}^{\frac{\sigma}{2}}
\left(
u_\mu^{x_0},
DP_\beta w_U(z)
\right)M_{\alpha,2}^{\frac{\sigma}{2}}
\left(
u_\mu,
DP_{<\beta}w_U(z)
DP_{<\beta}w_U(z)
\right)
\Bigg].
\end{split}
\end{equation*}
As in the preceding subsubsection, we first bound the coefficients in the full variable coefficient derivative $D$ in $L^\infty$ and then apply the
change-of-scale estimate in Proposition \ref{basicMproperties} only to the flat, $\beta$-localized
derivatives. Since $\sigma/2<1$, Proposition
\ref{basicMproperties} and Proposition
\ref{weighted-harmonic-extension-bounds} give
\begin{equation}\label{cubic-trace-high-scale-change}
\begin{split}
\sup_{z\leq0}
M_{\alpha,2}^{\frac{\sigma}{2}}
\left(
u_\mu^{x_0},
DP_\beta w_U(z)
\right)
&\lesssim
\left(
\frac{\alpha}{\beta}
\right)^{1-\frac{\sigma}{2}}
\sup_{z\leq0}
M_{\beta,2}^{\frac{\sigma}{2}}
\left(
u_\mu^{x_0},
\nabla_{z,x}P_\beta w_U(z)
\right)
\\
&\lesssim
\epsilon^2
\left(
\frac{\alpha}{\beta}
\right)^{1-\frac{\sigma}{2}}
\beta^{\frac14}
\mu^{-s_c}c_\mu
\left(1+\mu|x_0|\right)^{\frac12}.
\end{split}
\end{equation}
For the second factor, Proposition \ref{basicMproperties} and \eqref{cubic-trace-low-factor} give
\begin{equation*}
\begin{split}
M_{\alpha,2}^{\frac{\sigma}{2}}
\left(
u_\mu,
DP_{<\beta}w_U
DP_{<\beta}w_U
\right)&\lesssim_p
\alpha^{\frac4p}
\|u_\mu\|_{L_t^\infty L_x^2}
\left\|
DP_{<\beta}w_U
\right\|_{L_z^\infty L_t^4L_x^p}^2
\\
&\lesssim_p
\epsilon^3
\alpha^{\frac4p}
\beta^{\frac34-\frac4p}
\mu^{-s_c}c_\mu
(c_\beta)^2.
\end{split}
\end{equation*}
From $\frac{\sigma}{2}
>
1-\frac4p
>
\frac4p$, it follows that
\begin{equation}\label{cubic-trace-high-output}
\begin{split}
\left\|
L\left(
h_\mu,
\mathcal C_{\mathrm{hi}}
\right)
\right\|_{L_{t,x}^1}\lesssim_p
\epsilon^5
\mu^{-2s_c}c_\mu^2
\left(1+\mu|x_0|\right)^{\frac12}.
\end{split}
\end{equation}
Combining \eqref{cubic-trace-low-output},
\eqref{quartic-trace-low-output}, and
\eqref{cubic-trace-high-output}, we conclude that
\begin{equation}\label{full-cubic-trace-estimate}
\left\|
L\left(
h_\mu,
-2\mathcal P
\left(
\operatorname{tr}
\left((\nabla_{z,x}^2\phi_{<\lambda})^3\right)
\right)
\right)
\right\|_{L_{t,x}^1}
\lesssim_p
\epsilon^5
\mu^{-2s_c}c_\mu^2
\left(1+\mu|x_0|\right)^{\frac12}.
\end{equation}
\subsubsection{The term $4\mathcal{P}(\nabla_{z,x}^2\phi_{<\lambda}:\nabla_{z,x}^2\Delta_D^{-1}|\nabla_{z,x}^2\phi_{<\lambda}|^2)$} To simplify notation, we set $\mathcal C_F
\coloneq
\mathcal P\left(
D_iw_U\cdot D_iF
\right)$ and also define
\begin{equation*}
F
\coloneq
\nabla_{z,x}\Delta_D^{-1}
\left(
D_iw_U\cdot D_iw_U
\right),
\qquad
F_\beta
\coloneq
\nabla_{z,x}\Delta_D^{-1}
P_\beta\left(
D_iw_U\cdot D_iw_U
\right),
\end{equation*}
so that $F=\sum_\beta F_\beta$. Let $4<p<\infty$ and $8<q<2p$. We will always choose $\sigma$ so that $2-\frac{4}{p}<\sigma<2$. For this sub-subsection, we will additionally always choose $p$ so that $2-\frac{4}{p}<\frac{3}{2}$. This will be relevant momentarily because will either impose $\sigma>\frac{3}{2}$ or $\sigma<\frac{3}{2}$ depending on the term we are estimating. We first record the estimates for $F$ that will be used below.
Using the symmetrized frequency decomposition convention of
\cite[(7.5)]{ifrim2025sharp}, together with $D_iD_iw_{U_{\kappa}}=0$, we have
\begin{equation}\label{quadratic-source-divergence}
\begin{split}
D_iw_U\cdot D_iw_U
=
2\sum_\kappa
D_iw_{U_\kappa}\cdot D_iw_{U_{<\kappa}}=
2\sum_\kappa
D_i\left(
w_{U_\kappa}\cdot D_iw_{U_{<\kappa}}
\right).
\end{split}
\end{equation}
For the matched-frequency pieces, where
$\gamma\lesssim\kappa$, 
H\"older and Bernstein, together with Proposition
\ref{flattenedbasicbounds} and \eqref{ellipticw}
(the latter controlling the undifferentiated extension
by the fundamental theorem of calculus in $z$), give
\begin{equation}\label{quadratic-source-bounds}
\begin{split}
\left\|
P_\beta\left(
w_{U_\kappa}\cdot D_iw_{U_\gamma}
\right)
\right\|_{L_z^\infty L_t^2L_x^p}
&\lesssim
\epsilon^2\beta^{\frac4q-\frac2p}
\kappa^{-\frac58-\frac2q}
\gamma^{\frac38-\frac2q}
c_{\kappa}c_{\gamma},
\\
\left\|
P_\beta\left(
D_iw_{U_\kappa}\cdot D_iw_{U_\gamma}
\right)
\right\|_{L_z^\infty L_t^2L_x^p}
&\lesssim
\epsilon^2\beta^{\frac4q-\frac2p}
\kappa^{\frac38-\frac2q}
\gamma^{\frac38-\frac2q}
c_{\kappa}c_{\gamma}
\end{split}
\end{equation}
where we used also Proposition \ref{VBunbil} to estimate the contributions of $U_{\kappa}$ and $U_{\gamma}$. Here the factor $\beta^{\frac4q-\frac2p}$ follows from Bernstein
from $L_x^{q/2}$ to $L_x^p$. Observe importantly that $\frac38-\frac2q>0$ and 
$\frac4q-\frac2p>0.$ We will need this to retain summability. Next, for $\gamma\lesssim\kappa$, let
\begin{equation*}
F_{\beta,\kappa,\gamma}
\coloneq
\nabla_{z,x}\Delta_D^{-1}
P_\beta\left(
D_iw_{U_\kappa}\cdot D_iw_{U_\gamma}
\right).
\end{equation*}
Consequently, with the sum understood in the symmetrized sense of
\cite[(7.5)]{ifrim2025sharp}, we have
\begin{equation*}
F_\beta
=
2\sum_{\gamma\lesssim\kappa}
F_{\beta,\kappa,\gamma}.
\end{equation*}
Using the divergence form in \eqref{quadratic-source-divergence}
for the first estimate below, and the nondivergence form for the
second, the elliptic estimate for $\Delta_D^{-1}$ and \eqref{quadratic-source-bounds} give
\begin{equation}\label{F-geometric-bounds}
\begin{split}
\left\|
F_{\beta,\kappa,\gamma}
\right\|_{L_z^\infty L_t^2L_x^p}
&\lesssim
\epsilon^2\beta^{\frac4q-\frac2p}
\kappa^{-\frac58-\frac2q}
\gamma^{\frac38-\frac2q}
c_\kappa c_\gamma,\hspace{5mm}\left\|
DF_{\beta,\kappa,\gamma}
\right\|_{L_z^\infty L_t^2L_x^p}\lesssim
\epsilon^2\beta^{\frac4q-\frac2p}
\kappa^{\frac38-\frac2q}
\gamma^{\frac38-\frac2q}
c_\kappa c_\gamma.
\end{split}
\end{equation}
The first estimate is useful in the region $\beta\lesssim\kappa$. The latter estimate, we will need later. When $\beta\gtrsim\kappa$, we  also have from the non-divergence form of \eqref{quadratic-source-divergence},
\begin{equation}\label{F-geometric-bounds2}
\|F_{\beta,\kappa,\gamma}\|_{L_z^{\infty}L_t^2L_x^p}\lesssim \epsilon^2\beta^{-1+\frac{4}{q}-\frac{2}{p}}\kappa^{\frac{3}{8}-\frac{2}{q}}\gamma^{\frac{3}{8}-\frac{2}{q}}c_{\kappa}c_{\gamma}    
\end{equation}
In particular, by summation in $\gamma$ and $\kappa$ and $\beta$, we obtain
\begin{equation}\label{F-summed-zero-derivative-bound}
\sum_{\beta}
\beta^{\frac14+\frac2p}
\|F_\beta\|_{L_z^\infty L_t^2L_x^p}
\lesssim
\epsilon^2.
\end{equation}
With this preliminary estimate in hand, we are ready to estimate the full $\mathcal C_F$. We decompose $C_F$ based on whether the output frequency of the product $D_iw_U\cdot D_iF$ is larger or smaller than the frequency $\beta$ in each $F_{\beta}$. That is, we write $\mathcal C_F
=
\mathcal C_{F,\mathrm{hi}}
+
\mathcal C_{F,\mathrm{lo}}$, where
\begin{equation*}
\mathcal C_{F,\mathrm{hi}}
\coloneq
\sum_\beta\sum_{\alpha\gtrsim\beta}
\left.
\partial_z\Delta_D^{-1}
P_\alpha\left(
D_iw_U\cdot D_iF_\beta
\right)
\right|_{z=\eta_{<\lambda}},
\end{equation*}
and
\begin{equation*}
\mathcal C_{F,\mathrm{lo}}
\coloneq
\sum_\beta
\left.
\partial_z\Delta_D^{-1}
P_{<\beta}\left(
D_iw_U\cdot D_iF_\beta
\right)
\right|_{z=\eta_{<\lambda}}.
\end{equation*}
Before proceeding further, we record the translation-free estimate
\begin{equation}\label{translation-free-F-bound}
M_{\nu,2}^{\frac{\sigma}{2}}
\left(
u_\mu,G
\right)
\lesssim_p
\nu^{\frac2p}
\|u_\mu\|_{L_t^\infty L_x^2}
\|G\|_{L_t^2L_x^p},
\end{equation}
which follows from Proposition \ref{basicMproperties}. This is valid because $\sigma>2-4/p$. Now we estimate $\mathcal C_{F,\mathrm{hi}}$, which by definition
corresponds to the region $\alpha\gtrsim\beta$. Within this region,
we retain the decomposition of $F_\beta$ into the pieces
$F_{\beta,\kappa,\gamma}$ and split further according to whether
$\kappa\lesssim\alpha$ or $\alpha\ll\kappa$. The complementary
region $\alpha\ll\beta$ belongs to $\mathcal C_{F,\mathrm{lo}}$ and
will be treated separately below. Now, we choose parameters $\sigma_-$ and $\sigma_+$ satisfying
\begin{equation*}
2-\frac4p<\sigma_-<\frac32<\sigma_+<2.
\end{equation*}
Suppose first that $\kappa\lesssim\alpha$. Since we are considering
$\mathcal C_{F,\mathrm{hi}}$, we also have $\beta\lesssim\alpha$. We leave the derivative on $F_{\beta,\kappa,\gamma}$ and use
complementary fat cutoffs to decompose inside the geometric
derivative:
\begin{equation*}
D_iw_U
=
D_iP_{\lesssim\alpha}w_U
+
D_iP_{\gg\alpha}w_U.
\end{equation*}

For the contribution of $D_iP_{\lesssim\alpha}w_U$, we use
$\sigma=\sigma_+>\frac32$. Proposition
\ref{weightedbilinear}, followed by Corollary \ref{summationbilinear},
\eqref{translation-free-F-bound}, and the second estimate in
\eqref{F-geometric-bounds}, gives
\begin{equation*}
\begin{split}
&\sum_{\alpha\gtrsim\kappa}
\sum_{\beta\lesssim\alpha}
\left\|
L\left(
h_\mu,
\left.
\partial_z\Delta_D^{-1}
P_\alpha\left(
D_iP_{\lesssim\alpha}w_U\cdot
D_iF_{\beta,\kappa,\gamma}
\right)
\right|_{z=\eta_{<\lambda}}
\right)
\right\|_{L_{t,x}^1}
\\
&\qquad\lesssim
\sum_{\alpha\gtrsim\kappa}
\sum_{\beta\lesssim\alpha}
\alpha^{-1}
\sup_{z\leq0}
\Bigg[
M_{\alpha,2}^{\frac{\sigma_+}{2}}
\left(
u_\mu^{x_0},
D_iP_{\lesssim\alpha}w_U(z)
\right)
M_{\alpha,2}^{\frac{\sigma_+}{2}}
\left(
u_\mu,
D_iF_{\beta,\kappa,\gamma}(z)
\right)
\Bigg]
\\
&\qquad\lesssim_p
\epsilon^5
\mu^{-2s_c}c_\mu^2
\left(1+\mu|x_0|\right)^{\frac12}
\kappa^{\frac38-\frac2q}
\gamma^{\frac38-\frac2q}
c_\kappa c_\gamma
\sum_{\alpha\gtrsim\kappa}
\alpha^{-\frac34+\frac2p}
\sum_{\beta\lesssim\alpha}
\beta^{\frac4q-\frac2p}
\\
&\qquad\lesssim_p
\epsilon^5
\mu^{-2s_c}c_\mu^2
\left(1+\mu|x_0|\right)^{\frac12}
\kappa^{\frac38-\frac2q}
\gamma^{\frac38-\frac2q}
c_\kappa c_\gamma
\sum_{\alpha\gtrsim\kappa}
\alpha^{-\frac34+\frac4q}
\\
&\qquad\lesssim_p
\epsilon^5
\mu^{-2s_c}c_\mu^2
\left(1+\mu|x_0|\right)^{\frac12}
\left(\frac{\gamma}{\kappa}\right)^{
\frac38-\frac2q}
c_\kappa c_\gamma.
\end{split}
\end{equation*}
Here the $\beta$-sum converges at low frequencies because
$q<2p$, while the $\alpha$-sum converges because $q>8$. For the contribution of $D_iP_{\gg\alpha}w_U$, we instead use
$\sigma=\sigma_-<\frac32$ and
Corollary \ref{summationbilinear}. Since
$\sigma_->2-\frac4p$, the translation-free estimate
\eqref{translation-free-F-bound} remains available. Repeating the
preceding calculation gives
\begin{equation*}
\begin{split}
&\sum_{\alpha\gtrsim\kappa}
\sum_{\beta\lesssim\alpha}
\left\|
L\left(
h_\mu,
\left.
\partial_z\Delta_D^{-1}
P_\alpha\left(
D_iP_{\gg\alpha}w_U\cdot
D_iF_{\beta,\kappa,\gamma}
\right)
\right|_{z=\eta_{<\lambda}}
\right)
\right\|_{L_{t,x}^1}
\\
&\qquad\lesssim_p
\epsilon^5
\mu^{-2s_c}c_\mu^2
\left(1+\mu|x_0|\right)^{\frac12}
\left(\frac{\gamma}{\kappa}\right)^{
\frac38-\frac2q}
c_\kappa c_\gamma.
\end{split}
\end{equation*}
Summing over $\gamma\lesssim\kappa$ and then $\kappa\in \mathfrak{c}^{\mathbb{Z}}$ handles the region with both $\alpha\gtrsim\beta$ and   $\kappa\lesssim\alpha.$ We now move to the remaining region in $\mathcal C_{F,\mathrm{hi}}$ which is  $\beta\lesssim\alpha\ll\kappa.$ Before performing any frequency decomposition, we first use the fact that $D_iD_iw_U=0$ to write
\begin{equation*}
D_iw_U\cdot D_iF_{\beta,\kappa,\gamma}
=
D_i\left(
D_iw_U\cdot F_{\beta,\kappa,\gamma}
\right).
\end{equation*}
Now, we decompose inside $D_iw_U$ as before:
\begin{equation*}
\begin{split}
D_i\left(
D_iw_U\cdot F_{\beta,\kappa,\gamma}
\right)
=
D_i\left(
D_iP_{\lesssim\alpha}w_U\cdot
F_{\beta,\kappa,\gamma}
\right)+
D_i\left(
D_iP_{\gg\alpha}w_U\cdot
F_{\beta,\kappa,\gamma}
\right).
\end{split}
\end{equation*}

For the first contribution, we take again $\sigma=\sigma_+>\frac32$, use Proposition
\ref{weightedbilinear},
 Corollary \ref{summationbilinear},
\eqref{translation-free-F-bound}, and the first estimate in
\eqref{F-geometric-bounds} give
\begin{equation*}
\begin{split}
&
\sum_{\alpha\ll\kappa}\sum_{\beta\lesssim\alpha}\left\|
L\left(
h_\mu,
\left.
\partial_z\Delta_D^{-1}
P_\alpha D_i\left(
D_iP_{\lesssim\alpha}w_U\cdot
F_{\beta,\kappa,\gamma}
\right)
\right|_{z=\eta_{<\lambda}}
\right)
\right\|_{L_{t,x}^1}
\\
&\qquad\lesssim
\sum_{\alpha\ll\kappa}\sum_{\beta\lesssim\alpha}\sup_{z\leq0}
\Bigg[
M_{\alpha,2}^{\frac{\sigma_+}{2}}
\left(
u_\mu^{x_0},
D_iP_{\lesssim\alpha}w_U(z)
\right)
M_{\alpha,2}^{\frac{\sigma_+}{2}}
\left(
u_\mu,
F_{\beta,\kappa,\gamma}(z)
\right)
\Bigg]
\\
&\qquad\lesssim_p
\sum_{\alpha\ll\kappa}\sum_{\beta\lesssim\alpha}\epsilon^5
\mu^{-2s_c}c_\mu^2
\left(1+\mu|x_0|\right)^{\frac12}
\alpha^{\frac14+\frac2p}
\beta^{\frac4q-\frac2p}
\kappa^{-\frac58-\frac2q}
\gamma^{\frac38-\frac2q}
c_\kappa c_\gamma
\\
&\qquad\lesssim \epsilon^5
\mu^{-2s_c}c_\mu^2
\left(1+\mu|x_0|\right)^{\frac12}
\left(\frac{\gamma}{\kappa}\right)^{
\frac38-\frac2q}
c_\kappa c_\gamma.
\end{split}
\end{equation*}

For the second contribution, we take
$\sigma=\sigma_-<\frac32$, and obtain (analogously to the complemntary bound in the previous region)
\begin{equation*}
\begin{split}
&
\sum_{\alpha\ll\kappa}\sum_{\beta\lesssim\alpha}\left\|
L\left(
h_\mu,
\left.
\partial_z\Delta_D^{-1}
P_\alpha D_i\left(
D_iP_{\gg\alpha}w_U\cdot
F_{\beta,\kappa,\gamma}
\right)
\right|_{z=\eta_{<\lambda}}
\right)
\right\|_{L_{t,x}^1}
\\
&\qquad\lesssim_p
\sum_{\alpha\ll\kappa}\sum_{\beta\lesssim\alpha}\epsilon^5
\mu^{-2s_c}c_\mu^2
\left(1+\mu|x_0|\right)^{\frac12}
\alpha^{\frac14+\frac2p}
\beta^{\frac4q-\frac2p}
\kappa^{-\frac58-\frac2q}
\gamma^{\frac38-\frac2q}
c_\kappa c_\gamma
\\
&\lesssim \epsilon^5
\mu^{-2s_c}c_\mu^2
\left(1+\mu|x_0|\right)^{\frac12}
\left(\frac{\gamma}{\kappa}\right)^{
\frac38-\frac2q}
c_\kappa c_\gamma
\end{split}
\end{equation*}
Summing over $\gamma\lesssim\kappa$ and $\kappa\in \mathfrak{c}^{\mathbb{Z}}$, and using the estimate in the other region from before, we obtain finally
\begin{equation}\label{second-cubic-high-output}
\left\|
L\left(
h_\mu,
\mathcal C_{F,\mathrm{hi}}
\right)
\right\|_{L_{t,x}^1}
\lesssim_p
\epsilon^5
\mu^{-2s_c}c_\mu^2
\left(1+\mu|x_0|\right)^{\frac12}.
\end{equation}
It remains to estimate the low-output contribution. We can write
\begin{equation*}
\mathcal C_{F,\mathrm{lo}}
=
\sum_\beta\sum_{\alpha\ll\beta}
\left.
\partial_z\Delta_D^{-1}
P_\alpha
D_i\left(D_iw_U\cdot F_\beta
\right)
\right|_{z=\eta_{<\lambda}}.
\end{equation*}
As before, we split
\begin{equation*}
D_iw_U
=
D_iP_{\lesssim\alpha}w_U
+
D_iP_{\gg\alpha}w_U.
\end{equation*}

For the first contribution, we take
$\sigma=\sigma_+>\frac32$. Proposition \ref{weightedbilinear} and then Corollary \ref{summationbilinear} and
\eqref{translation-free-F-bound} gives
\begin{equation*}
\begin{split}
&\sum_\beta\sum_{\alpha\ll\beta}
\left\|
L\left(
h_\mu,
\left.
\partial_z\Delta_D^{-1}
P_\alpha D_i\left(
D_iP_{\lesssim\alpha}w_U\cdot F_\beta
\right)
\right|_{z=\eta_{<\lambda}}
\right)
\right\|_{L_{t,x}^1}
\\
&\qquad\lesssim
\sum_\beta\sum_{\alpha\ll\beta}
\sup_{z\leq0}
\Bigg[
M_{\alpha,2}^{\frac{\sigma_+}{2}}
\left(
u_\mu^{x_0},
D_iP_{\lesssim\alpha}w_U(z)
\right)
M_{\alpha,2}^{\frac{\sigma_+}{2}}
\left(
u_\mu,
F_\beta(z)
\right)
\Bigg]
\\
&\qquad\lesssim_p
\epsilon^3
\mu^{-2s_c}c_\mu^2
\left(1+\mu|x_0|\right)^{\frac12}
\sum_\beta\sum_{\alpha\ll\beta}
\alpha^{\frac14+\frac2p}
\|F_\beta\|_{L_z^\infty L_t^2L_x^p}
\end{split}
\end{equation*}
For the second contribution, we take
$\sigma=\sigma_-<\frac32$. Using Proposition
\ref{weightedbilinear}, Corollary \ref{summationbilinear},
and \eqref{translation-free-F-bound}, we get also
\begin{equation*}
\begin{split}
&\sum_\beta\sum_{\alpha\ll\beta}
\left\|
L\left(
h_\mu,
\left.
\partial_z\Delta_D^{-1}
P_\alpha D_i\left(
D_iP_{\gg\alpha}w_U\cdot F_\beta
\right)
\right|_{z=\eta_{<\lambda}}
\right)
\right\|_{L_{t,x}^1}
\\
&\qquad\lesssim_p
\epsilon^3
\mu^{-2s_c}c_\mu^2
\left(1+\mu|x_0|\right)^{\frac12}
\sum_\beta\sum_{\alpha\ll\beta}
\alpha^{\frac14+\frac2p}
\|F_\beta\|_{L_z^\infty L_t^2L_x^p}.
\end{split}
\end{equation*}
The two estimates above, together with
\eqref{F-summed-zero-derivative-bound}, yield
\begin{equation}\label{second-cubic-low-output}
\begin{split}
\left\|
L\left(
h_\mu,
\mathcal C_{F,\mathrm{lo}}
\right)
\right\|_{L_{t,x}^1}
&\lesssim_p
\epsilon^3
\mu^{-2s_c}c_\mu^2
\left(1+\mu|x_0|\right)^{\frac12}
\sum_\beta
\beta^{\frac14+\frac2p}
\|F_\beta\|_{L_z^\infty L_t^2L_x^p}
\\
&\lesssim_p
\epsilon^5
\mu^{-2s_c}c_\mu^2
\left(1+\mu|x_0|\right)^{\frac12}.
\end{split}
\end{equation}
Combining \eqref{second-cubic-high-output} and
\eqref{second-cubic-low-output}, we conclude that
\begin{equation*}
\left\|
L\left(
h_\mu,
\mathcal C_F
\right)
\right\|_{L_{t,x}^1}
\lesssim_p
\epsilon^5
\mu^{-2s_c}c_\mu^2
\left(1+\mu|x_0|\right)^{\frac12}.
\end{equation*}
\subsubsection{The term
$2\mathcal{P}(\nabla_{z,x}^2\phi_{<\lambda}:
\nabla_{z,x}^2\mathcal{H}\rho_{<\lambda})$}

In flattened coordinates, we are dealing with
\begin{equation*}
2\mathcal P\left(
D_iw_{U_j}D_iD_jw_{\rho_{<\lambda}}
\right).
\end{equation*}
Since $D_iD_iw_{U_j}=0$, we can as usual, write
\begin{equation}\label{rho-cubic-divergence}
D_iw_{U_j}D_iD_jw_{\rho_{<\lambda}}
=
D_i\left(
D_iw_{U_j}D_jw_{\rho_{<\lambda}}
\right).
\end{equation}
Now, fix $0<\delta<\frac18$ and set
\begin{equation*}
\sigma_0^+
=
\frac34+\delta,
\qquad
\sigma_0^-
=
\frac34-\delta,
\qquad
\sigma_1
=
1+\delta,
\qquad
\sigma^\pm
=
\sigma_0^\pm+\sigma_1.
\end{equation*}
Then $\sigma^-=\frac74$
and $\sigma^+=\frac74+2\delta<2.$ From Proposition \ref{basicMproperties} we have the general estimate
\begin{equation}\label{rho-asymmetric-M-bounds}
\begin{split}
M_{\lambda,1}^{\sigma^\pm}
\left(
h_\mu,
f_\pm g
\right)
\lesssim
M_{\lambda,2}^{\sigma_0^\pm}
\left(
u_\mu^{x_0},
f_\pm
\right)
M_{\lambda,2}^{\sigma_1}
\left(
u_\mu,
g
\right).
\end{split}
\end{equation}
Since $\sigma_1>1$, the weight in the second factor is
integrable. As in the previous term, we now split
\begin{equation*}
D_iw_{U_j}
=
D_iP_{\lesssim\lambda}w_{U_j}
+
D_iP_{\gg\lambda}w_{U_j}.
\end{equation*}
For the first term we use $\sigma_0^+>\frac34$. By Proposition \ref{basicMproperties} and Corollary \ref{summationbilinear}
\begin{equation}\label{rho-low-U-bound}
\begin{split}
\sup_{z\leq0}
M_{\lambda,2}^{\sigma_0^+}
\left(
u_\mu^{x_0},
D_iP_{\lesssim\lambda}w_{U_j}(z)
\right)&\lesssim
\sum_{\gamma\lesssim\lambda}
\left(
\frac{\lambda}{\gamma}
\right)^{1-\sigma_0^+}
\sup_{z\leq0}
M_{\gamma,2}^{\sigma_0^+}
\left(
u_\mu^{x_0},
D_iP_\gamma w_{U_j}(z)
\right)
\\
&\lesssim
\epsilon^2
\lambda^{\frac14}
\mu^{-s_c}c_\mu
\left(1+\mu|x_0|\right)^{\frac12}.
\end{split}
\end{equation}
For the second term we instead use
$\sigma_0^-<\frac34$. The corresponding high-frequency estimate and the same argument gives
\begin{equation}\label{rho-high-U-bound}
\begin{split}
\sup_{z\leq0}
M_{\lambda,2}^{\sigma_0^-}
\left(
u_\mu^{x_0},
D_iP_{\gg\lambda}w_{U_j}(z)
\right)\lesssim
\epsilon^2
\lambda^{\frac14}
\mu^{-s_c}c_\mu
\left(1+\mu|x_0|\right)^{\frac12}.
\end{split}
\end{equation}
On the other hand, since $\sigma_1>1$, Proposition \ref{basicMproperties} gives
\begin{equation*}
\begin{split}
\sup_{z\leq0}
M_{\lambda,2}^{\sigma_1}
\left(
u_\mu,
D_jw_{P_{<M\lambda}\rho_{<\lambda}}(z)
\right)\lesssim
\|u_\mu\|_{L_t^\infty L_x^2}
\left\|
Dw_{\rho_{<\lambda}}
\right\|_{L_z^\infty L_t^2L_x^\infty}
\lesssim
\epsilon^3
\lambda^{-\frac14}
\mu^{-s_c}c_\mu.
\end{split}
\end{equation*}
Here we used the assumed estimate
\begin{equation*}
\left\|
Dw_{\rho_{<\lambda}}
\right\|_{L_z^\infty L_t^2L_x^\infty}
\lesssim\|\nabla_x\rho_{<\lambda}\|_{L_t^2\dot{B}_{\infty,1}^0}\lesssim 
\epsilon^2\lambda^{-\frac14}
\end{equation*}
where in the first inequality, we used Proposition \ref{flattenedbasicbounds}, and in the second, we used the estimate \eqref{rhobasicbound}. Now, applying Proposition \ref{weightedbilinear}, and \eqref{rho-asymmetric-M-bounds}
separately to the two pieces of $D_iw_{U}$ and using \eqref{rho-low-U-bound} and \eqref{rho-high-U-bound}, we obtain
\begin{equation*}
\begin{split}
\left\|
L\left(
h_\mu,
2\mathcal P D_i\left(
D_iP_{\lesssim\lambda}w_{U_j}
D_jw_{\rho_{<\lambda}}
\right)
\right)
\right\|_{L_{t,x}^1}\lesssim
\epsilon^5
\mu^{-2s_c}c_\mu^2
\left(1+\mu|x_0|\right)^{\frac12}.
\end{split}
\end{equation*}
\subsubsection{The term $-\kappa_{<\lambda} \left.\partial_z^2p_{<\lambda}\right|_{z=\eta_{<\lambda}}$} We begin by rewriting this term in a slightly easier to handle form. Differentiating the equation
\begin{equation*}
\Delta_{z,x}p_{<\lambda}
=
-
\left|
\nabla_{z,x}^2\phi_{<\lambda}
\right|^2
\end{equation*}
in $z$ and using $\left.
\partial_zp_{<\lambda}
\right|_{z=\eta_{<\lambda}}
=
-a_{<\lambda},$ we obtain
\begin{equation*}
\partial_z^2p_{<\lambda}|_{z=\eta_{<\lambda}}
=
-L_{<\lambda}a_{<\lambda}
-
\mathcal P\left(
\partial_z
\left|
\nabla_{z,x}^2\phi_{<\lambda}
\right|^2
\right)
\end{equation*}
Consequently,
\begin{equation*}
\begin{split}
-\kappa_{<\lambda}
\left.
\partial_z^2p_{<\lambda}
\right|_{z=\eta_{<\lambda}}
=
\kappa_{<\lambda}L_{<\lambda}a_{<\lambda}+
\kappa_{<\lambda}\mathcal P\left(
\partial_z
\left|
\nabla_{z,x}^2\phi_{<\lambda}
\right|^2
\right).
\end{split}
\end{equation*}
Using the $L_t^2L_x^\infty$ bound for $\kappa_{<\lambda}$
established in the proof of \eqref{kappaubil}, together with
the energy bound for $u_\mu$, Proposition \ref{abilinearbounds},
and Proposition \ref{basicopbounds}, we obtain
\begin{equation*}
\begin{split}
\|L(u_{\mu},u_{\mu}^{x_0},-\kappa_{<\lambda}\partial_z^2p_{<\lambda}|_{z=\eta_{<\lambda}})\|_{L_{t,x}^1}&\lesssim \|L(u_{\mu}^{x_0},\kappa_{<\lambda})\|_{L_{t,x}^2}\|L(u_{\mu},L_{<\lambda}a_{<\lambda})\|_{L_{t,x}^2}
\\
&+\|L(u_{\mu}^{x_0},\kappa_{<\lambda})\|_{L_{t,x}^2}\|u_{\mu}\|_{L_{t}^{\infty}L_x^2}\|\mathcal{P}\left(
\partial_z
\left|
\nabla_{z,x}^2\phi_{<\lambda}
\right|^2
\right)\|_{L_{t}^2L_x^{\infty}}  
\\
&\lesssim +\epsilon^3\mu^{-2s_c}\lambda^{-\frac{3}{4}}c^2_{\mu}(1+\mu|x_0|)^{\frac12}\|\mathcal{P}\left(
\partial_z
\left|
\nabla_{z,x}^2\phi_{<\lambda}
\right|^2
\right)\|_{L_{t}^2L_x^{\infty}}
\\
&+\epsilon^4\mu^{-2s_c}c_{\mu}^2(1+\mu|x_0|)^{\frac12}
\end{split}
\end{equation*}
Using \eqref{flattened-weighted-elliptic-equation} (or more precisely, its proof), we bound 
\begin{equation*}
\begin{split}
\|\mathcal{P}\left(
\partial_z
\left|
\nabla_{z,x}^2\phi_{<\lambda}
\right|^2
\right)\|_{L_{t}^2L_x^{\infty}}\lesssim \|D_iw_U\cdot D_iw_U\|_{L_{z}^{\infty}L_t^2\dot{B}_{\infty,1}^0}.
\end{split}
\end{equation*}
 Fix $p>4$. Proposition \ref{flattenedbasicbounds} and Proposition \ref{VBunbil} and the
bootstrap bound give
\begin{equation*}
\left\|
Dw_{U_\nu}
\right\|_{L_z^\infty L_t^4L_x^{2p}}
\lesssim
\epsilon
\nu^{\frac38-\frac1p}
c_\nu.
\end{equation*}
For $\gamma\lesssim\kappa\lesssim\lambda$ and
$\alpha\lesssim\kappa$, H\"older and Bernstein give
\begin{equation*}
\begin{split}
\left\|
P_\alpha\left(
Dw_{U_\kappa}\cdot Dw_{U_\gamma}
\right)
\right\|_{L_z^\infty L_t^2L_x^\infty}
&\lesssim
\alpha^{\frac2p}
\left\|
Dw_{U_\kappa}
\right\|_{L_z^\infty L_t^4L_x^{2p}}
\left\|
Dw_{U_\gamma}
\right\|_{L_z^\infty L_t^4L_x^{2p}}
\\
&\lesssim
\epsilon^2
\alpha^{\frac2p}
\kappa^{\frac38-\frac1p}
\gamma^{\frac38-\frac1p}
c_\kappa c_\gamma.
\end{split}
\end{equation*}
Consequently,
\begin{equation*}
\begin{split}
\left\|
Dw_U\cdot Dw_U
\right\|_{L_z^\infty L_t^2\dot B_{\infty,1}^0}
&\lesssim
\epsilon^2
\sum_{\gamma\lesssim\kappa\lesssim\lambda}
\sum_{\alpha\lesssim\kappa}
\alpha^{\frac2p}
\kappa^{\frac38-\frac1p}
\gamma^{\frac38-\frac1p}
c_\kappa c_\gamma
\\
&\lesssim
\epsilon^2
\sum_{\gamma\lesssim\kappa\lesssim\lambda}
\kappa^{\frac34}
\left(
\frac{\gamma}{\kappa}
\right)^{\frac38-\frac1p}
c_\kappa c_\gamma
\\
&\lesssim
\epsilon^2\lambda^{\frac34}.
\end{split}
\end{equation*}

\subsubsection{The term $(V^l_{<\lambda}-V_{<\lambda})\cdot\nabla a_{<\lambda}$}
This is simply a high-frequency tail error and is straightforward to estimate.  Indeed, since $V^l_{<\lambda}-V_{<\lambda}=-P_{\geq M\lambda}V_{<\lambda}$,
Bernstein and Proposition \ref{VBunbil} and Proposition \ref{Bernsteinelliptic} give, for $M$ sufficiently large,
\begin{equation*}
\|V^l_{<\lambda}-V_{<\lambda}\|_{L_t^1L_x^\infty}
\lesssim
\||D|P_{>M\lambda}V_{<\lambda}\|_{L_t^1\dot B^0_{2,1}}
\lesssim \epsilon^4\lambda^{-1}.
\end{equation*}
By \eqref{a-coefficient-bounds}, we also have $\|\nabla a_{<\lambda}\|_{L_{t,x}^\infty}\lesssim\epsilon\lambda$. Hence, by H\"older and placing both $u_\mu$ factors in $L_t^\infty L_x^2$, we obtain
\begin{equation*}
\left\|L\left(u_\mu^{x_0},u_\mu,
(V^l_{<\lambda}-V_{<\lambda})\cdot\nabla a_{<\lambda}\right)\right\|_{L_{t,x}^1}
\lesssim \epsilon^7 c_\mu^2\mu^{-2s_c}.
\end{equation*}
Combining the preceding estimates finally proves the first bound in Proposition \ref{Dtaestimates}.

We now move on to proving the second estimate \eqref{DtaprincipalL2L8} for $8\leq p\leq\infty$. In this part of the proof, to keep notation less cumbersome, we use a dot to denote $\lambda\frac{d}{d\lambda}$. We begin by differentiating the identity \eqref{Dtaeqn}. We note here, that the the term $-L_{<\lambda}(a_{<\lambda}\kappa_{<\lambda})$, which was removed in the first bound, must now be included. Writing $d_V=V_{<\lambda}-P_{<\lambda}V$ as in \eqref{dBVdecomp}, a simple computation gives the identity
\begin{equation*}
\begin{split}
D_t^l a_\lambda-T_{a^l}\nabla\cdot V_\lambda &=\lambda\frac{d}{d\lambda}\left(
(\partial_t+V_{<\lambda}\cdot\nabla)a_{<\lambda}
-a_{<\lambda}\nabla\cdot V_{<\lambda}\right)
\\
&+a_\lambda\nabla\cdot V_{<\lambda}
-\dot V_{<\lambda}\cdot\nabla a_{<\lambda}
+(V^l_{<\lambda}-V_{<\lambda})\cdot\nabla a_\lambda
\\
&+(a_{<\lambda}-T_{a^l})\nabla\cdot V_\lambda
+a_{<\lambda}\nabla\dot d_V.
\end{split}
\end{equation*}
Here $V_\lambda=S_{\lambda}V$, so that $\dot V_{<\lambda}=V_\lambda+\dot d_V$. 

For the second line we use H\"older with two $L_t^4L_x^{2p}$ factors. Indeed, by \eqref{linearizedaformula}, \eqref{a-coefficient-bounds}, the linearization formula for $V$, Proposition \ref{basicopbounds} and Proposition \ref{VBunbil}, we have
\begin{equation}\label{usefulboundsaV}
\begin{split}
\|a_\lambda\|_{L_t^4L_x^{2p}}+\lambda^{-1}\|\nabla a_\lambda\|_{L_t^4L_x^{2p}}
&\lesssim_\delta\epsilon c_\lambda\lambda^{-\frac18-\frac1p},
\hspace{5mm}
\|\nabla V_{<\lambda}\|_{L_t^4L_x^{2p}}\lesssim\epsilon\lambda^{\frac38-\frac1p},
\\
\|\dot V_{<\lambda}\|_{L_t^4L_x^{2p}}
&\lesssim\epsilon c_\lambda\lambda^{-\frac58-\frac1p},
\hspace{5mm}
\|\nabla a_{<\lambda}\|_{L_t^4L_x^{2p}}\lesssim\epsilon\lambda^{\frac78-\frac1p}.
\end{split}
\end{equation}
Moreover, by Proposition \ref{Bernsteinelliptic} and Proposition \ref{VBunbil}, we have the tail bound $\|V^l_{<\lambda}-V_{<\lambda}\|_{L_t^4L_x^{2p}}\lesssim\epsilon\lambda^{-\frac58-\frac1p}$. For the paraproduct in the last line, we can use the fundamental theorem of calculus in \eqref{a-truncated-integral} and \eqref{a-coefficient-bounds} to obtain
\begin{equation*}
\|a_{<\lambda}-P_{<M\delta\sigma}a_{<\delta\sigma}\|_{L_t^4L_x^{2p}}
\lesssim_\delta\epsilon\lambda^{-\frac18-\frac1p},
\qquad
\|\nabla V_\lambda\|_{L_t^4L_x^{2p}}\lesssim\epsilon c_\lambda\lambda^{\frac38-\frac1p}.
\end{equation*}
whenever $\sigma\approx\lambda$. Thus $(a_{<\lambda}-T_{a^;})\nabla\cdot V_{\lambda}$ is bounded by $O_{L_t^2L_x^p}(\epsilon^2c_\lambda\lambda^{\frac14-\frac2p})$.

We next estimate the terms containing $\dot d_V$ or
$\kappa_{<\lambda}$. We begin with the corresponding estimate for
$\dot d_B$, which we will also need. Since
$d_B=B_{<\lambda}-P_{<\lambda}B$, the linearization formula for $B$
gives
\begin{equation*}
\dot d_B
=L_{<\lambda}\Phi_{<\lambda}(\psi_\lambda,\eta_\lambda)
+\eta_\lambda L_{<\lambda}B_{<\lambda}-S_\lambda B.
\end{equation*}
By Lemma \ref{Phiest} and the three estimates following
\eqref{Br4threeerrors}, we have
\begin{equation*}
L_{<\lambda}\Phi_{<\lambda}(\psi_\lambda,\eta_\lambda)
=T_{\mathfrak l}\varphi_\lambda
+O_{L_t^2L_x^p}(\epsilon^2c_\lambda\lambda^{-\frac34-\frac2p}).
\end{equation*}
Similarly, by Corollary \ref{principalVBexp} and
Proposition \ref{Leibnizcommutator} applied to
$[S_\lambda,T_{\mathfrak l}]\varphi$, we have
\begin{equation*}
S_\lambda B
=T_{\mathfrak l}\varphi_\lambda
+O_{L_t^2L_x^p}(\epsilon^2c_\lambda\lambda^{-\frac34-\frac2p}).
\end{equation*}
Thus the linear principal terms in $\dot d_B$ cancel. For the remaining product term, we can use Proposition \ref{basicopbounds} and Proposition \ref{VBunbil} to directly estimate
\begin{equation*}
\|\eta_\lambda L_{<\lambda}B_{<\lambda}\|_{L_t^2L_x^p}
\lesssim
\|\eta_\lambda\|_{L_t^4L_x^{2p}}
\|L_{<\lambda}B_{<\lambda}\|_{L_t^4L_x^{2p}}\lesssim\epsilon^2c_\lambda\lambda^{-\frac34-\frac2p}.
\end{equation*}
Using the above bounds, along with Bernstein and Proposition
\ref{Bernsteinelliptic} for the high frequency tails,
we conclude that
\begin{equation}\label{dotdB}
\|\nabla^j\dot d_B\|_{L_t^2L_x^p}
\lesssim_\delta\epsilon^2c_\lambda\lambda^{j-\frac34-\frac2p},
\hspace{5mm} j=0,1,2.
\end{equation}
The same argument for $\dot d_V$, with $L_{<\lambda}$ and
$T_{\mathfrak l}$ replaced by $M_{<\lambda}$ and $T_{\mathfrak m}$,
gives also
\begin{equation}\label{dotdV}
\|\nabla^j\dot d_V\|_{L_t^2L_x^p}
\lesssim_\delta\epsilon^2c_\lambda\lambda^{j-\frac34-\frac2p},
\hspace{5mm} j=0,1,2.
\end{equation}
Finally, we can obtain a similar bound for $\dot\kappa_{<\lambda}$ by differentiating the identity preceding
\eqref{kappaform} to write
\begin{equation*}
\dot\kappa_{<\lambda}
=-\dot d_B+\dot d_V\cdot\nabla\eta_{<\lambda}
+d_V\cdot\nabla\eta_\lambda
+\lambda\frac{d}{d\lambda}\mathcal B_{<\lambda}(V,\nabla\eta).
\end{equation*}
For the third term, we use \eqref{dvdbL2Linf} and
$\|\nabla\eta_\lambda\|_{L_t^\infty L_x^p}
\lesssim\epsilon c_\lambda\lambda^{-\frac2p}$.
For the last term, differentiating the definition of
$\mathcal B_{<\lambda}$ gives
\begin{equation}\label{calBderivative}
\begin{split}
\lambda\frac{d}{d\lambda}\mathcal B_{<\lambda}(V,\nabla\eta)
&=S_{\lambda}V\cdot \nabla\eta_{<\lambda}
+P_{<\lambda}V\cdot S_{\lambda}\eta
-S_\lambda(V\cdot\nabla\eta).
\end{split}
\end{equation}
Expanding each factor of $V$ and $\nabla\eta$ above as a sum over dyadic pieces $V_{\beta}$ and $\eta_{\alpha}$, we see that $\mathcal{B}_{<\lambda}(V_{\beta},\nabla\eta_{\alpha})$ can be written as $-[S_{\lambda},V_{\beta}]\nabla\eta_{\alpha}$ or $-[S_{\lambda},\nabla\eta_{\alpha}]V_{\beta}$ when $\beta\ll\alpha\approx\lambda$ and $\alpha\ll\beta\approx\lambda$ respectively. In both regimes, Proposition \ref{Leibnizcommutator} gives
\begin{equation*}
\left\|\lambda\frac{d}{d\lambda}
\mathcal B_{<\lambda}(V_\beta,\nabla\eta_\alpha)\right\|_{L_t^2L_x^p}
\lesssim
\frac{\min\{\alpha,\beta\}}{\lambda}
\|V_\beta\|_{L_t^4L_x^{2p}}
\|\nabla\eta_\alpha\|_{L_t^4L_x^{2p}},
\end{equation*}
which sums in the lower frequency to give $O_{L_t^2L_x^p}(\epsilon^2c_{\lambda}\lambda^{-\frac34-\frac2p})$. By estimating each term individually, we have the same estimate when summing over the balanced frequencies $\alpha\approx\beta\gtrsim\lambda$. It is easy to see that for any other pair of $\alpha$ and $\beta$, \eqref{calBderivative} vanishes. Combining these estimates with the bound for $\dot d_B$ and $\dot d_V$ above and Bernstein, it is easy to show the bound
\begin{equation*}
\|\nabla^j\dot\kappa_{<\lambda}\|_{L_t^2L_x^p}
\lesssim_\delta\epsilon^2c_\lambda\lambda^{j-\frac34-\frac2p},
\hspace{5mm} j=0,1,2.
\end{equation*}
Integrating the $j=0$ estimate from $\lambda$ to infinity and then
applying Bernstein, with Proposition \ref{Bernsteinelliptic} for the high frequency tails, gives also
\begin{equation*}
\|\nabla^j\kappa_{<\lambda}\|_{L_t^2L_x^p}
\lesssim_\delta\epsilon^2c_\lambda\lambda^{j-\frac34-\frac2p},
\qquad j=0,1,2.
\end{equation*}
Finally, applying the linearization formula for $L_{<\lambda}$
to $a_{<\lambda}\kappa_{<\lambda}$ gives
\begin{equation*}
\begin{split}
\lambda\frac{d}{d\lambda}
\bigl(L_{<\lambda}(a_{<\lambda}\kappa_{<\lambda})\bigr)&=L_{<\lambda}\bigl(a_\lambda\kappa_{<\lambda}
+a_{<\lambda}\dot\kappa_{<\lambda}\bigr)
\\
&-L_{<\lambda}\bigl(\eta_\lambda
L_{<\lambda}(a_{<\lambda}\kappa_{<\lambda})\bigr)
+\eta_\lambda L_{<\lambda}^2(a_{<\lambda}\kappa_{<\lambda}).
\end{split}
\end{equation*}
By the bounds above, Proposition \ref{besovparabounds} and the \eqref{usefulboundsaV}, we have
\begin{equation*}
\|a_\lambda\kappa_{<\lambda}
+a_{<\lambda}\dot\kappa_{<\lambda}\|_{L_t^2\dot B^1_{p,1}}
\lesssim\epsilon^2c_\lambda\lambda^{\frac14-\frac2p},
\qquad
\|a_{<\lambda}\kappa_{<\lambda}\|_{L_t^2\dot B^1_{p,1}}
\lesssim\epsilon^2\lambda^{\frac14-\frac2p}.
\end{equation*}
We apply Proposition \ref{basicopbounds} to the first line and
estimate the two terms in the second line together using
\eqref{a-zero-order-commutator} and
$\|\nabla\eta_\lambda\|_{L_t^\infty\dot B^0_{\infty,1}}
\lesssim\epsilon c_\lambda$. This gives
\begin{equation*}
\left\|\lambda\frac{d}{d\lambda}
\bigl(L_{<\lambda}(a_{<\lambda}\kappa_{<\lambda})\bigr)
\right\|_{L_t^2L_x^p}
\lesssim_\delta\epsilon^2c_\lambda\lambda^{\frac14-\frac2p}.
\end{equation*}
We next estimate the $\lambda$ derivative of the other boundary term $-\kappa_{<\lambda}\left.\partial_z^2p_{<\lambda}\right|_{z=\eta_{<\lambda}}.$ Differentiating $p_{<\lambda}|_{z=\eta_{<\lambda}}=0$ twice in space
and using the Laplace equation for $p_{<\lambda}$ gives
\begin{equation}\label{pressuretraceterm}
(1+|\nabla\eta_{<\lambda}|^2)\left.\partial_z^2p_{<\lambda}\right|_{z=\eta_{<\lambda}}
=-\left.|\nabla_{x,z}^2\phi_{<\lambda}|^2\right|_{z=\eta_{<\lambda}}
-2\nabla\eta_{<\lambda}\cdot\nabla a_{<\lambda}
-a_{<\lambda}\Delta\eta_{<\lambda}.
\end{equation}
We can rewrite the first term on the right hand side in terms of $V$ and $B$ by recalling that $\partial_{x_j}\phi_{<\lambda}$ and $\partial_z\phi_{<\lambda}$
are the harmonic extensions of $V_{<\lambda}^j$ and $B_{<\lambda}$,
respectively. By the definitions of $M_{<\lambda}$ and $L_{<\lambda}$, we therefore have
\begin{equation*}
-\left.|\nabla_{x,z}^2\phi_{<\lambda}|^2\right|_{z=\eta_{<\lambda}}
=-\sum_{j=1}^2|M_{<\lambda}V_{<\lambda}^j|^2
-2|M_{<\lambda}B_{<\lambda}|^2
-|L_{<\lambda}B_{<\lambda}|^2.
\end{equation*} 
The $L_t^{\infty}L_x^2$ bounds for $V$, $a$ and $\eta$, and Bernstein give an $L_{t,x}^\infty$ bound of
$O(\epsilon\lambda)$ for \eqref{pressuretraceterm} and a bound of $O_{L_{t,x}^{\infty}}(\epsilon c_\lambda\lambda)$
for its $\lambda$ derivative. Thus, by the product rule and the $j=0$
bounds for $\kappa_{<\lambda}$ and $\dot\kappa_{<\lambda}$, we obtain
\begin{equation*}
\begin{split}
\left\|\lambda\frac{d}{d\lambda}
\left(-\kappa_{<\lambda}\left.\partial_z^2p_{<\lambda}\right|_{z=\eta_{<\lambda}}\right)
\right\|_{L_t^2L_x^p}
&\lesssim_\delta
\epsilon\lambda\|\dot\kappa_{<\lambda}\|_{L_t^2L_x^p}
+\epsilon c_\lambda\lambda\|\kappa_{<\lambda}\|_{L_t^2L_x^p}
\\
&\lesssim_\delta\epsilon^3c_\lambda\lambda^{\frac14-\frac2p}.
\end{split}
\end{equation*}
It remains to estimate the other four terms
\begin{equation}\label{remainingfour}
\begin{split}
&-2\mathcal P\bigl(
\mathrm{tr}((\nabla_{x,z}^2\phi_{<\lambda})^3)
\bigr),
\hspace{5mm}
-4\mathcal P\bigl(
\nabla_{x,z}^2\phi_{<\lambda}:
\nabla_{x,z}^2\mathcal H\eta_{<\lambda}
\bigr),
\\
&2\mathcal P\bigl(
\nabla_{x,z}^2\phi_{<\lambda}:
\nabla_{x,z}^2\mathcal H\rho_{<\lambda}
\bigr),
\hspace{5mm}
4\mathcal P\bigl(
\nabla_{x,z}^2\phi_{<\lambda}:
\nabla_{x,z}^2\Delta_D^{-1}
|\nabla_{x,z}^2\phi_{<\lambda}|^2
\bigr)
\end{split}
\end{equation}
in \eqref{Dtaeqn}. It is not so convenient to differentiate these with respect to $\lambda$ directly because we would then have to obtain linearization formulas for $\mathcal{P},\mathcal{H}$ and $\Delta_D^{-1}$ and also $\phi_{<\lambda}$. It is better to first rewrite these terms into a more manageable form. Our starting point for this is to observe, using the definition of $L$ and $\mathcal{P}$, the general identity
\begin{equation*}
\mathcal P\bigl(
\nabla_{x,z}\mathcal H f\cdot\nabla_{x,z}\mathcal H g
\bigr)
=\frac12\bigl(fL_{<\lambda}g+gL_{<\lambda}f
-L_{<\lambda}(fg)\bigr)\coloneq Q_a(f,g).
\end{equation*}
Indeed, $\mathcal H f\,\mathcal H g-\mathcal H(fg)$
vanishes on the boundary, and its Laplacian is $2\nabla_{z,x}\mathcal{H}f\cdot\nabla_{z,x}\mathcal{H}g$. We also write 
\begin{equation*}
\begin{split}
J=(-a_{<\lambda}\nabla\eta_{<\lambda},a_{<\lambda}-1),
\hspace{5mm}
Y=\left.\nabla_{x,z}\mathcal H\rho_{<\lambda}
\right|_{z=\eta_{<\lambda}},\hspace{5mm}K_j=\Delta_D^{-1}\partial_j
|\nabla_{x,z}^2\phi_{<\lambda}|^2
\end{split}
\end{equation*}
where $j=1,2,3$. Since
$p_{<\lambda}=\mathcal H\eta_{<\lambda}-z
-\Delta_D^{-1}|\nabla_{x,z}^2\phi_{<\lambda}|^2$,
we can then write
\begin{equation*}
\nabla_{x,z}\Delta_D^{-1}
|\nabla_{x,z}^2\phi_{<\lambda}|^2
-\nabla_{x,z}\mathcal H\eta_{<\lambda}
=\mathcal HJ+K.
\end{equation*}
To see this, we observe that the left hand side has trace $J$ on the boundary
and its Laplacian (in $x$ and $z$) is
$\nabla_{x,z}|\nabla_{x,z}^2\phi_{<\lambda}|^2$,
while each component of $K$ has zero boundary data.
Consequently, the above two identities allow us to rewrite the contribution of the terms \eqref{remainingfour} in the equation \eqref{Dtaeqn} as
\begin{equation}\label{fourtermsexp}
4\sum_{j=1}^3Q_a(U_j,J_j)
+2\sum_{j=1}^3Q_a(U_j,Y_j)
-2\mathcal P\bigl(
\mathrm{tr}((\nabla_{x,z}^2\phi_{<\lambda})^3)
\bigr)
+4\mathcal P\bigl(
\nabla_{x,z}^2\phi_{<\lambda}:\nabla_{x,z}K
\bigr).
\end{equation}
Before applying $\lambda\frac{d}{d\lambda}$, we first collect the corresponding undifferentiated estimates for the expressions above. We will then show how these can be transferred easily to the differentiated bounds. Set $\theta=\frac14-\frac2p$, so that
$0\leq\theta\leq\frac14$. By symmetry, and the estimate \eqref{a-zero-order-commutator} we have
\begin{equation}\label{Qaestimate}
\|Q_a(f_\alpha,g_\gamma)\|_{L_x^p}
\lesssim
\min\{\alpha,\gamma\}
\max\{\alpha,\gamma\}^{\frac1{12}+\theta}
\|f_\alpha\|_{L_x^{p_1}}
\|g_\gamma\|_{L_x^{p_2}}.
\end{equation}
Here we used the commutator structure in $Q_a$ to ensure the derivative is shifted onto
the lower frequency factor (which by Bernstein produces the factor of $\min\{\alpha,\gamma\}$), and then applied Bernstein, which gives a factor of $\max\{\alpha,\gamma\}^{\frac{1}{12}+\theta}$. 

For the contribution of $Q_a(U_j,J_j)$, we apply this estimate with $p_1=p_2=12$, and then use Proposition \ref{VBunbil}, Proposition \ref{abilinearbounds} and the bootstrap bounds for the $\nabla\eta$ factors appearing in $J$. We first consider the contribution in the region $\alpha,\gamma\lesssim\lambda$, obtaining
\begin{equation}\label{QUJ}
\begin{split}
\sum_{\alpha,\gamma\lesssim\lambda}
\|Q_a(P_{\alpha}U_j,P_{\gamma}J_j)\|_{L_t^2L_x^p}
&\lesssim_\delta
\epsilon^2\sum_{\alpha,\gamma\lesssim\lambda}
c_\alpha c_\gamma
\min\{\alpha,\gamma\}
\max\{\alpha,\gamma\}^{\frac1{12}+\theta}
\alpha^{-\frac{19}{24}}\gamma^{-\frac7{24}}
\\
&=\epsilon^2\sum_{\alpha\leq\gamma\lesssim\lambda}
c_\alpha c_\gamma\gamma^\theta
\left(\frac{\alpha}{\gamma}\right)^{\frac5{24}}
+\epsilon^2\sum_{\gamma<\alpha\lesssim\lambda}
c_\alpha c_\gamma\alpha^\theta
\left(\frac{\gamma}{\alpha}\right)^{\frac{17}{24}}
\\
&\lesssim\epsilon^2\lambda^\theta.
\end{split}
\end{equation}
Here we used $\theta\geq0$ (and additionally the $\ell^1$ summability of the frequency envelope when $\theta=0$) to bound the larger frequency factor by $\lambda^\theta$. The contribution coming from when one of $\alpha$ or $\gamma$ is much larger than $\lambda$ consists of high frequency tails. Here, we apply the same bound \eqref{Qaestimate} and then use Proposition \ref{Bernsteinelliptic} or Lemma \ref{hightailabound} depending on which factor ($U$ or $J$) is at high frequency, and then Bernstein. This gives the same bound as above. Hence,
\begin{equation*}
\sum_{\alpha,\gamma>0}
\|Q_a(P_{\alpha}U_j,P_{\gamma}J_j)\|_{L_t^2L_x^p}
\lesssim
\epsilon^2\lambda^\theta.
\end{equation*}
Next, we estimate $Q_a(U_j,Y_j)$. First, using the definition of $M_{<\lambda}$ and $L_{<\lambda}$, we can write
\begin{equation*}
Y=\bigl(M_{<\lambda}\rho_{<\lambda},
L_{<\lambda}\rho_{<\lambda}\bigr).
\end{equation*}
As a consequence, we have the following estimate, which follows from Proposition \ref{basicopbounds},
\begin{equation}\label{BesovY14}
\begin{split}
\||D|^{\frac{1}{4}+\theta}Y\|_{L_t^2\dot B^{0}_{\infty,1}}
\lesssim_\delta
\|\nabla\rho_{<\lambda}\|_{L_t^2\dot B^{\frac14+\theta}_{\infty,1}}
+\lambda^{\frac14+\theta}
\|\nabla\rho_{<\lambda}\|_{L_t^2\dot B^0_{\infty,1}}\lesssim_\delta\epsilon^2\lambda^\theta.
\end{split}
\end{equation}
where in the last inequality, we interpolated and used \eqref{rhobasicbound} for $j=0$ and $j=2$.

By the $L_t^\infty L_x^2$ bounds for $B$ and $V$ in Proposition \ref{VBunbil} and Bernstein, we also have
$\|U_\alpha\|_{L_t^\infty L_x^6}
\lesssim\epsilon c_\alpha\alpha^{-\frac56}$.
Now, applying \eqref{Qaestimate} with $(p_1,p_2)=(6,\infty)$ and taking $\gamma$ as the frequency of $Y$,
we can first sum in $\alpha$, since
\begin{equation*}
\begin{split}
\sum_{\alpha>0}
\min\{\alpha,\gamma\}
\max\{\alpha,\gamma\}^{\frac1{12}+\theta}
\alpha^{-\frac56}
&=\gamma^{\frac14+\theta}
\left(
\sum_{\alpha\leq\gamma}
\left(\frac{\alpha}{\gamma}\right)^{\frac16}
+\sum_{\alpha>\gamma}
\left(\frac{\gamma}{\alpha}\right)^{\frac34-\theta}
\right)
\\
&\lesssim\gamma^{\frac14+\theta}
\end{split}
\end{equation*}
where we used $0\leq\theta\leq\frac14$.
The remaining sum in $\gamma$ is controlled by the Besov estimate for $Y$ in \eqref{BesovY14}. Consequently,
\begin{equation}\label{QUY}
\left\|\sum_{j=1}^3Q_a(U_j,Y_j)\right\|_{L_t^2L_x^p}
\lesssim
\epsilon\|Y\|_{L_t^2\dot B^{\frac14+\theta}_{\infty,1}}
\lesssim\epsilon^3\lambda^\theta
\end{equation}

For the two cubic terms in \eqref{fourtermsexp}, we first
consider the term containing
$\mathrm{tr}((\nabla_{x,z}^2\phi_{<\lambda})^3)$.
We decompose the boundary data $U$ using a Littlewood-Paley decomposition
and then apply harmonic extension. By Proposition
\ref{flattenedbasicbounds} and then Proposition \ref{VBunbil}, we have
\begin{equation*}
\|Dw_{P_\alpha U}\|_{L_z^\infty L_x^q}
\lesssim_q
\|\nabla_xP_\alpha U\|_{\dot B^0_{q,1}}
\lesssim_q\alpha\|P_\alpha U\|_{L_x^q}.
\end{equation*}
In particular, we have the three bounds
\begin{equation*}
\begin{split}
\|Dw_{P_\alpha U}\|_{L_z^\infty L_t^\infty L_x^2}
\lesssim\epsilon c_\alpha\alpha^{-\frac12},
\hspace{5mm}
\|Dw_{P_\alpha U}\|_{L_t^4L_z^\infty L_x^{12}}
\lesssim\epsilon c_\alpha\alpha^{\frac5{24}},\hspace{5mm}\|Dw_{P_\alpha U}\|_{L_z^\infty L_t^\infty L_x^\infty}
\lesssim\epsilon c_\alpha\alpha^{\frac12}.
\end{split}
\end{equation*}
By the elliptic estimate for $\mathcal P$, Bernstein,
and dyadic summation, we have for each $0<\alpha\lesssim\lambda$,
\begin{equation*}
\begin{split}
&\left\|
\mathcal P\mathrm{tr}\bigl(
Dw_{P_\alpha U}(Dw_{P_{<\alpha}U})^2
\bigr)\right\|_{L_t^2L_x^p}
\\
&\qquad\lesssim
\sum_{p>0}\rho^{-1}
\left\|P_\rho\mathrm{tr}\bigl(
Dw_{P_\alpha U}(Dw_{P_{<\alpha}U})^2
\bigr)\right\|_{L_t^2L_z^\infty L_x^p}
\\
&\qquad\lesssim
\left(
\sum_{\rho\leq\alpha}\rho^{\frac1{12}+\theta}
\|Dw_{P_\alpha U}\|_{L_t^\infty L_z^\infty L_x^2}
+\sum_{\rho>\alpha}\rho^{-\frac{11}{12}+\theta}
\|Dw_{P_\alpha U}\|_{L_{t,z,x}^\infty}
\right)
\|Dw_{P_{<\alpha}U}\|_{L_t^4L_z^\infty L_x^{12}}^2
\\
&\qquad\lesssim
\left(
\alpha^{\frac1{12}+\theta}
\|Dw_{P_\alpha U}\|_{L_t^\infty L_z^\infty L_x^2}
+\alpha^{-\frac{11}{12}+\theta}
\|Dw_{P_\alpha U}\|_{L_{t,z,x}^\infty}
\right)
\|Dw_{P_{<\alpha}U}\|_{L_t^4L_z^\infty L_x^{12}}^2
\\
&\qquad\lesssim
\epsilon^3c_\alpha^3\alpha^\theta.
\end{split}
\end{equation*}
Here, for $\rho\leq\alpha$, we used Bernstein to go from $L_x^p$ to $L^{\frac{3}{2}}$ and then applied H\"older.
For $\rho>\alpha$, we instead use Bernstein to go from $L_x^p$ to $L_x^6$. The remaining sum over $\alpha\lesssim\lambda$ converges
by summability of $c_\alpha$ and
$\alpha^\theta\lesssim\lambda^\theta$. Doing a similar analysis for the high frequency tails which appear in the range $\alpha\gg \lambda$, using Proposition \ref{Bernsteinelliptic} gives the same estimate after summing over $\alpha$. Symmetrizing and applying the above estimates, we obtain
\begin{equation*}
\left\|
\mathcal P\bigl(
\mathrm{tr}((\nabla_{x,z}^2\phi_{<\lambda})^3)
\bigr)\right\|_{L_t^2L_x^p}
\lesssim_\delta\epsilon^3\lambda^\theta.
\end{equation*}
One can carry out a similar analysis for the other cubic term containing $K$ to also get 
\begin{equation*}
\left\|
\mathcal P\bigl(
\nabla_{x,z}^2\phi_{<\lambda}:\nabla_{x,z}K
\bigr)\right\|_{L_t^2L_x^p}
\lesssim_\delta\epsilon^3\lambda^\theta.
\end{equation*}

Now, we explain why the same estimates hold (with an additional factor of $c_{\lambda}$) after applying $\lambda\frac{d}{d\lambda}$ to each term in \eqref{fourtermsexp}.
For the terms $Q_a(\dot U_j,J_j)$ and $Q_a(U_j,\dot J_j)$,
where the derivative falls on one of the arguments in $Q_a$, we can use the linearization formulas for $V$ and $B$, and the bounds for $a_\lambda$ above to obtain the bounds
\begin{equation*}
\|\dot U\|_{L_t^4\dot B^{\frac{19}{24}+\theta}_{12,1}}
+\|\dot J\|_{L_t^4\dot B^{\frac7{24}+\theta}_{12,1}}
\lesssim_\delta
\epsilon c_\lambda\lambda^\theta.
\end{equation*}
By dyadic summation and the same norm allocations as in
\eqref{QUJ} (with $\dot U$ in place of $U$ or $\dot J$
in place of $J$, depending on the term being estimated),
we obtain
\begin{equation*}
\sum_{j=1}^3\left(
\|Q_a(\dot U_j,J_j)\|_{L_t^2L_x^p}
+\|Q_a(U_j,\dot J_j)\|_{L_t^2L_x^p}
\right)
\lesssim_\delta
\epsilon^2c_\lambda\lambda^\theta.
\end{equation*}
For $Q_a(\dot U_j,Y_j)$, we instead use
\begin{equation*}
\|\dot U\|_{L_t^\infty\dot B^{\frac56}_{6,1}}
\lesssim_\delta\epsilon c_\lambda
\end{equation*}
which gives the analogue of the estimate \eqref{QUY} for $Q_a(\dot{U}_j,Y)$. For the contribution of $\dot{Y}$, we have to differentiate $\rho_{<\lambda}$. To do this, we apply $\lambda\frac{d}{d\lambda}$ to
\eqref{eq:rho-good-form} to obtain
\begin{equation*}
\|\nabla^j\dot\rho_{<\lambda}\|_{L_t^2L_x^\infty}
\lesssim_\delta
\epsilon^2c_\lambda\lambda^{j-\frac54},
\hspace{5mm} j=0,1,2.
\end{equation*}
This follows by a similar argument to the bounds for $\dot{\kappa}_{<\lambda}$ above, so we leave the details to the reader. Combining this with \eqref{rhobasicbound}, the linearization
formulas for $L$ and $M$ in \eqref{Blin} and \eqref{Vlin},
and the operator and commutator estimates in Proposition
\ref{basicopbounds} and \eqref{a-zero-order-commutator},
we obtain
\begin{equation*}
\|\dot Y\|_{L_t^2\dot B^{\frac14+\theta}_{\infty,1}}
\lesssim_\delta
\epsilon^2c_\lambda\lambda^\theta.
\end{equation*}
Using the same norm allocations and dyadic summations
as in \eqref{QUY}, this gives
\begin{equation*}
\sum_{j=1}^3
\|Q_a(U_j,\dot Y_j)\|_{L_t^2L_x^p}
\lesssim_\delta
\epsilon^3c_\lambda\lambda^\theta.
\end{equation*}
It remains to estimate the contributions where the $\lambda\frac{d}{d\lambda}$ derivative
falls on an elliptic operator, and the derivatives of the two
cubic terms. 

For the former case, we first describe a simple estimate that will be useful for treating these terms efficiently. To this end, it is convenient to work in the
flattened coordinate frame. We fix a function $f$ (that does not depend on $\lambda$) on the boundary and write as usual $w_f$, for its flattened harmonic extension (with respect to $\eta_{<\lambda}$). Since $f$ does not depend on $\lambda$, it follows that $\dot{w}_f$ has zero trace on the boundary. Therefore, differentiating \eqref{wgduhamel} gives the representation formula
\begin{equation}\label{dotwrep}
\nabla_{z,x}\dot w_f
=\mathcal R_0(A_\eta\nabla_{z,x}\dot w_f)
+\mathcal R_0(\dot A_\eta\nabla_{z,x}w_f).
\end{equation}
The entries of $\dot A_\eta$ contain at least one factor of $\nabla\eta_\lambda$,
and thus by Proposition \ref{besovparabounds} we have
\begin{equation*}
\|\dot A_\eta\|_{L_t^\infty\dot B^1_{2,1}}
\lesssim\epsilon c_\lambda.
\end{equation*}
We then absorb the first term on the right hand side of \eqref{dotwrep} as in the proof of
Proposition \ref{flattenedbasicbounds}. The same product estimates used there then give the linearized extension bound
\begin{equation}\label{linextbound}
\|\nabla_{z,x}\dot w_f\|_{L_z^\infty\dot B^0_{q,1}}
\lesssim_q
\epsilon c_\lambda\|\nabla_xf\|_{\dot B^0_{q,1}},
\qquad 1<q\leq\infty.
\end{equation}
Arguing similarly to the proof of \eqref{a-zero-order-commutator}
and differentiating again the flattened equation as above, we also obtain
the commutator bound
\begin{equation}\label{linearizedcombound}
\|[\dot L_{<\lambda},f]g\|_{L_x^q}
\lesssim_{q,q_1,q_2}
\epsilon c_\lambda
\|\nabla_xf\|_{\dot B^0_{q_1,1}}
\|g\|_{\dot B^0_{q_2,1}},
\hspace{5mm}
\frac1q=\frac1{q_1}+\frac1{q_2},
\hspace{5mm} 1<q,q_1,q_2\leq\infty.
\end{equation}
The same estimate as in \eqref{linextbound} applies when
differentiating $Dw_f$, since the differentiated coefficients
of $D$ also contain $\nabla\eta_\lambda$. Now, using that $(\dot L_{<\lambda})f=\partial_z\dot w_f|_{z=0}$, we apply these bounds to first estimate the contribution of $\dot Q_a(U_j,J_j)$ and $\dot Q_a(U_j,Y_j)$ above.

For these terms, we keep the commutator structure from above. For fixed $f$ and $g$, we have
\begin{equation*}
\begin{split}
\dot Q_a(f,g)
&=\frac12\bigl(
f\dot L_{<\lambda}g+g\dot L_{<\lambda}f
-\dot L_{<\lambda}(fg)
\bigr)
\\
&=\frac12\bigl(
g\dot L_{<\lambda}f-[\dot L_{<\lambda},f]g
\bigr).
\end{split}
\end{equation*}
Then as in \eqref{Qaestimate}, by using \eqref{linextbound} and \eqref{linearizedcombound}, we obtain 
\begin{equation*}
\|\dot Q_a(f_\alpha,g_\gamma)\|_{L_x^p}
\lesssim_\delta
\epsilon c_\lambda
\min\{\alpha,\gamma\}
\max\{\alpha,\gamma\}^{\frac1{12}+\theta}
\|f_\alpha\|_{L_x^{p_1}}
\|g_\gamma\|_{L_x^{p_2}}.
\end{equation*}
with the same choices of $p_1,p_2$ as in \eqref{Qaestimate}. The same norm allocations and summations as in
\eqref{QUJ} and \eqref{QUY} therefore give
\begin{equation*}
\sum_{j=1}^3\left(
\|\dot Q_a(U_j,J_j)\|_{L_t^2L_x^p}
+\|\dot Q_a(U_j,Y_j)\|_{L_t^2L_x^p}
\right)
\lesssim_\delta\epsilon^3c_\lambda\lambda^\theta
\end{equation*}
as required.

Finally, to estimate the two cubic terms, we use the same estimates as above.
When the derivative falls on one of the boundary functions $U$, the
linearization formulas for $B$ and $V$ and the high-freuqncy tail bounds used in
\eqref{dotdB} and \eqref{dotdV} give
\begin{equation*}
\|Dw_{\dot U}\|_{L_t^\infty L_z^\infty L_x^2}
\lesssim\epsilon c_\lambda\lambda^{-\frac12},
\qquad
\|Dw_{\dot U}\|_{L_{t,z,x}^\infty}
\lesssim\epsilon c_\lambda\lambda^{\frac12}.
\end{equation*}
Then, for example, for the in term in \eqref{fourtermsexp} that is cubic in $\nabla_{z,x}^2\phi_{<\lambda}$,
we put this $\lambda$ differentiated factor (which plays the same role as the high-frequency term $Dw_{P_{\alpha}U}$ in the undifferentiated estimate above) depending on its output frequency in one of two norms above and
the other two factors in $L_t^4L_z^\infty L_x^{12}$. Similarly to the undifferentiated cubic bound above, this gives
\begin{equation*}
\begin{split}
\left\|
\mathcal P\mathrm{tr}\bigl(Dw_{\dot U}(Dw_U)^2\bigr)
\right\|_{L_t^2L_x^p}\lesssim_\delta\epsilon^3c_\lambda
\left(
\lambda^{\frac1{12}+\theta}\lambda^{-\frac12}
+\lambda^{-\frac{11}{12}+\theta}\lambda^{\frac12}
\right)\lambda^{\frac5{12}}
\lesssim\epsilon^3c_\lambda\lambda^\theta.
\end{split}
\end{equation*}

When the $\lambda$ derivative falls on an elliptic operator in this cubic term, we
differentiate its flattened equation as above. For
$\mathcal P$, the same argument shows that the operator $\dot{\mathcal P}$ satisfies the Besov bound
\begin{equation*}
\|\dot{\mathcal{P}}F\|_{L_t^2\dot{B}_{p,1}^0}\lesssim \epsilon c_{\lambda}\sum_{\rho>0}\rho^{-1}\|P_{\rho}F\|_{L_t^2L_z^{\infty}L_x^p},    
\end{equation*}
which suffices to give
\begin{equation*}
\begin{split}
\left\|
\dot{\mathcal P}\bigl(
\mathrm{tr}((\nabla_{x,z}^2\phi_{<\lambda})^3)
\bigr)\right\|_{L_t^2L_x^p}\lesssim_\delta
\epsilon^4c_\lambda\lambda^{\frac14-\frac2p}.
\end{split}
\end{equation*}
The cubic term containing $K$ is estimated similarly, using in addition to above, the analogous frequency localized bounds for the operator
$\nabla_{x,z}\Delta_D^{-1}\partial_j$ that appears in the definition of $\nabla K$.

Combining all of the above estimates then finally gives
\begin{equation*}
\begin{split}
&4\left\|\lambda\frac{d}{d\lambda}
\sum_{j=1}^3Q_a(U_j,J_j)\right\|_{L_t^2L_x^p}
+2\left\|\lambda\frac{d}{d\lambda}
\sum_{j=1}^3Q_a(U_j,Y_j)\right\|_{L_t^2L_x^p}
\\
+&2\left\|\lambda\frac{d}{d\lambda}
\mathcal P\bigl(
\mathrm{tr}((\nabla_{x,z}^2\phi_{<\lambda})^3)
\bigr)\right\|_{L_t^2L_x^p}
+4\left\|\lambda\frac{d}{d\lambda}
\mathcal P\bigl(
\nabla_{x,z}^2\phi_{<\lambda}:\nabla_{x,z}K
\bigr)\right\|_{L_t^2L_x^p}
\\
&\lesssim_\delta
\epsilon^2c_\lambda\lambda^{\frac14-\frac2p}.
\end{split}
\end{equation*}
Together with the linearized estimates for the other boundary and transport terms above,
this finally proves \eqref{DtaprincipalL2L8}.
\end{proof}

We now prove Proposition \ref{Dtacommutatorestimates}. At this point, we now have to treat the contribution of the remaining term $L_{<\lambda}(a_{<\lambda}\kappa_{<\lambda})$. Unfortunately, this will not be directly perturbative. Rather, at leading order, it will exhibit a null structure similar to that observed in the initial paradifferential expansion of the $\varphi$ equation. We begin with the following "principal" decomposition of $L_{<\lambda}(a_{<\lambda}\kappa_{<\lambda})$, which isolates the bulk of this null term. 
\begin{lemma}\label{Lakappaprincipal}
Let $\delta_0>0$ be sufficiently small. For each $\lambda\ll\delta_0\mu\ll\mu$, we have the decomposition
\begin{equation}\label{Lakappaprincipalexp}
\begin{split}
L_{<\lambda}(a_{<\lambda}\kappa_{<\lambda})
&=T_{\mathfrak l}T_{a_{<\lambda}}\Bigg[
P_{<\lambda}
\int_{\delta_0\mu}^{\infty}
\int_{\delta\alpha}^{\alpha}
\Bigl(
\mathcal T^{B,\mathrm{main}}_{<\delta\beta}
-\nabla\eta_\beta\cdot T_{\mathfrak m}\varphi_\alpha
-\nabla\eta_\alpha\cdot T_{\mathfrak m}\varphi_\beta
\Bigr)
\frac{d\beta}{\beta}\frac{d\alpha}{\alpha}
\\
&-\nabla\eta_{<\lambda}\cdot P_{<\lambda}
\int_{\delta_0\mu}^{\infty}
\int_{\delta\alpha}^{\alpha}
\mathcal T^{V,\mathrm{main}}_{<\delta\beta}
\frac{d\beta}{\beta}\frac{d\alpha}{\alpha}
+O_{L_{t,x}^2}\left(\epsilon^3\mu^{-\frac74}\right)
\Bigg]+\mathcal R.
\end{split}
\end{equation}
Here we use the principal quadratic expressions from
\eqref{TBmain2ord} and \eqref{Vquadraticprincipal}, namely
\begin{equation*}
\begin{split}
\mathcal T^{B,\mathrm{main}}_{<\delta\beta}
&=-T_{\mathfrak l}
\left(
\eta_\beta T_{\mathfrak l}\varphi_\alpha
+\eta_\alpha T_{\mathfrak l}\varphi_\beta
\right)
+\eta_\beta T_{\mathfrak l}^2\varphi_\alpha
+\eta_\alpha T_{\mathfrak l}^2\varphi_\beta,
\\
\mathcal T^{V,\mathrm{main}}_{<\delta\beta}
&=-T_{\mathfrak m}
\left(
\eta_\beta T_{\mathfrak l}\varphi_\alpha
+\eta_\alpha T_{\mathfrak l}\varphi_\beta
\right)
+\eta_\beta T_{\mathfrak m}T_{\mathfrak l}\varphi_\alpha
+\eta_\alpha T_{\mathfrak m}T_{\mathfrak l}\varphi_\beta.
\end{split}
\end{equation*}
The remainder satisfies
\begin{equation}\label{Lakappaprincipalerror}
\|L(u_\mu^{x_0},u_\mu,\mathcal R)\|_{L_{t,x}^1}
\lesssim_{\delta,\delta_0}
\epsilon^4c_\mu^2\mu^{-2s_c}(1+\mu|x_0|)^{\frac12}.
\end{equation}
\end{lemma}

\begin{proof}
We first replace $L_{<\lambda}(a_{<\lambda}\kappa_{<\lambda})$
by the paraproduct $T_{\mathfrak l}T_{a_{<\lambda}}\kappa_{<\lambda}$, and show that the error is of the form $\mathcal{R}$ above. We write
\begin{equation*}
\begin{split}
L_{<\lambda}(a_{<\lambda}\kappa_{<\lambda})
-T_{\mathfrak l}T_{a_{<\lambda}}\kappa_{<\lambda}=
(L_{<\lambda}-T_{\mathfrak l})
(a_{<\lambda}\kappa_{<\lambda})
+T_{\mathfrak l}
\left(
a_{<\lambda}\kappa_{<\lambda}
-T_{a_{<\lambda}}\kappa_{<\lambda}
\right).
\end{split}
\end{equation*}
For the first term, we slightly adapt the estimates following
\eqref{rhoathreeterms} but replace every instance of
$\rho_{<\delta\lambda}$ with
$a_{<\lambda}\kappa_{<\lambda}$. We place $\kappa_{<\lambda}$
in $L_t^2L_x^\infty$, using the estimate preceding
\eqref{kappaubil}, and $a_{<\lambda}$ in $L_{t,x}^{\infty}$.
We pair the additional $\eta$ factor from the linearization
of $L_{<\lambda}$ with one of the $u_\mu$ factors and use one unbalanced bilinear estimate and estimate
the other $u_\mu$ factor in $L_t^\infty L_x^2$.
The subprincipal errors are treated using Proposition
\ref{paraexp}, while the high-frequency tails are handled
using Proposition \ref{Bernsteinelliptic} and
Lemma \ref{hightailabound} for the tails of $a_{<\lambda}$.
We omit the detailed calculations since they are routine
modifications of the bounds following \eqref{rhoathreeterms}.
\\

For the second term, we expand the paraproduct
$a_{<\lambda}-T_{a_{<\lambda}}$ as
\begin{equation*}
(a_{<\lambda}-T_{a_{<\lambda}})\kappa_{<\lambda}
=
\int_0^\infty
S_\beta a_{<\lambda}\,
P_{<(M\delta)^{-1}\beta}\kappa_{<\lambda}
\frac{d\beta}{\beta}.
\end{equation*}
Note here, that $T_{a_{<\lambda}}$ is a standard paraproduct with the scale $\lambda$ fixed. For $\beta\lesssim M\lambda$, we use Proposition \ref{abilinearbounds} 
to pair the $S_{\beta}a_{<\lambda}$ factor with $\overline u_\mu^{x_0}$.
We place the $\kappa$ factor in $L_t^2L_x^\infty$, using
the estimate preceding \eqref{kappaubil}, and estimate
$u_\mu$ in $L_t^\infty L_x^2$. Integrating in $\beta$, we find that
\begin{equation*}
\begin{split}
\int_{\beta\lesssim M\lambda}\|L(\overline{u}_{\mu}^{x_0},u_{\mu},T_{\mathfrak l}\bigl(S_{\beta}a_{<\lambda}P_{<(M\delta)^{-1}\beta}\kappa_{<\lambda}\bigr))\|_{L_{t,x}^1}\frac{d\beta}{\beta}\lesssim_{M,\delta,\delta_0}
\epsilon^5c_\mu^2\mu^{-2s_c}(1+\mu|x_0|)^{\frac12}.
\end{split}
\end{equation*}

For the contribution of $P_{>M\lambda}a_{<\lambda}$,
we use \eqref{atailbound} with $s=2$, $p=2$, and
$\sigma=3$ and Bernstein to obtain
\begin{equation*}
\begin{split}
\int_{\beta\gtrsim M\lambda}
\beta\|P_\beta P_{>M\lambda}a_{<\lambda}\|_{L_t^2L_x^\infty}
\frac{d\beta}{\beta}\lesssim
\lambda^{-1}
\||D|^3P_{>M\lambda}a_{<\lambda}\|_{L_{t,x}^2}\lesssim\epsilon^5\lambda^{\frac{3}{4}}.
\end{split}
\end{equation*}
Combining this with
$\|\kappa_{<\lambda}\|_{L_t^2L_x^\infty}
\lesssim\epsilon^2\lambda^{-\frac34}$
bounds this part of the paraproduct error in
$L_t^1L_x^\infty$ by $\epsilon^7$.
Estimating both $u_\mu$ factors in $L_t^{\infty}L_x^2$ and combining this with the $\beta\lesssim M\lambda$ bound above, we get
\begin{equation*}
\begin{split}
\left\|L\left(
\overline u_\mu^{x_0},u_\mu,
T_{\mathfrak l}\left(
a_{<\lambda}\kappa_{<\lambda}
-T_{a_{<\lambda}}\kappa_{<\lambda}
\right)
\right)\right\|_{L_{t,x}^1}\lesssim_{\delta,\delta_0}
\epsilon^5c_\mu^2\mu^{-2s_c}(1+\mu|x_0|)^{\frac12},
\end{split}
\end{equation*}
which satisfies \eqref{Lakappaprincipalerror}.
\\

We next expand the remaining principal component $T_{\mathfrak{l}}T_{a_{<\lambda}}\kappa_{<\lambda}$. For this, we recall the expansion of $\kappa$ used in the proof of
\eqref{kappaform}. Specifically, with the notation
\begin{equation*}
d_B=B_{<\lambda}-P_{<\lambda}B,
\qquad
d_V=V_{<\lambda}-P_{<\lambda}V,
\qquad
\mathcal B_{<\lambda}(f,g)
=(P_{<\lambda}f)(P_{<\lambda}g)-P_{<\lambda}(fg),
\end{equation*}
we have the identity
\begin{equation*}
\kappa_{<\lambda}
=-d_B+\nabla\eta_{<\lambda}\cdot d_V
+\mathcal B_{<\lambda}(V,\nabla\eta).
\end{equation*}
We begin by expanding the contribution of $d_B$. As in
\eqref{dBVdecomp}, we first write
\begin{equation*}
\begin{split}
d_B
={}&
\bigl(B_{<\lambda}-T_{\mathfrak l}\varphi_{<\lambda}\bigr)
-P_{<\lambda}\bigl(B-T_{\mathfrak l}\varphi\bigr)
+[T_{\mathfrak l},P_{<\lambda}]\varphi.
\end{split}
\end{equation*}
To expand the remaining differences, we use the fundamental
theorem of calculus to obtain
\begin{equation*}
\begin{split}
d_B
&=P_{\geq\lambda}\int_0^\lambda
\left(
\alpha\frac{d}{d\alpha}B_{<\alpha}
-T_{\mathfrak l}\varphi_\alpha
\right)\frac{d\alpha}{\alpha}-P_{<\lambda}\int_\lambda^\infty
\left(
\alpha\frac{d}{d\alpha}B_{<\alpha}
-T_{\mathfrak l}\varphi_\alpha
\right)\frac{d\alpha}{\alpha}
\\
&+[T_{\mathfrak l},P_{<\lambda}]\varphi.
\end{split}
\end{equation*}
We first estimate the commutator
$[T_{\mathfrak l},P_{<\lambda}]\varphi$.
We use Proposition \ref{Leibnizcommutator} and two unbalanced
bilinear estimates, with the pairings
\begin{equation*}
(u_\mu^{x_0},\nabla^2\eta_\beta),
\hspace{5mm}
(u_\mu,T_{\mathfrak l}\varphi_\alpha).
\end{equation*}
Using the $L^{\infty}_{t,x}$ bound for $a$ and the kernel bounds for $T_{\mathfrak{l}}$, we then obtain
\begin{equation*}
\begin{split}
\left\|L\left(
u_\mu^{x_0},u_\mu,
T_{\mathfrak l}T_{a_{<\lambda}}
[T_{\mathfrak l},P_{<\lambda}]\varphi
\right)\right\|_{L_{t,x}^1}&\lesssim_{\delta,\delta_0}
\epsilon^4c_\mu^2\mu^{-2s_c}(1+\mu|x_0|)^{\frac12}
\lambda
\int_{\alpha\approx\lambda}
\int_{\beta\lesssim\alpha}
\alpha^{-1}\alpha^{-\frac34}\beta^{\frac34}
\frac{d\beta}{\beta}\frac{d\alpha}{\alpha}
\\
&\lesssim_{\delta,\delta_0}
\epsilon^4c_\mu^2\mu^{-2s_c}(1+\mu|x_0|)^{\frac12}.
\end{split}
\end{equation*}
Thus this term satisfies \eqref{Lakappaprincipalerror}.

We now estimate the two integrals in the expression for $d_B$.
We apply the $B$ analogue of the successive linearizations
\eqref{VtailFTCstart} and \eqref{TVlinearize}, as in
the argument preceding \eqref{FBendpoint}.
More precisely, these give
\begin{equation}\label{BlinDtaexp}
\begin{split}
\alpha\frac{d}{d\alpha}B_{<\alpha}
&=
L^{l,\alpha}_{<\delta\alpha}
\Phi^l_{<\delta\alpha}(\psi_\alpha,\eta_\alpha)
+\eta_\alpha L^{l,\alpha}_{<\delta\alpha}
B^l_{<\delta\alpha}
\\
&+\int_{\delta\alpha}^{\alpha}
\beta\frac{d}{d\beta}
\left(
L^{l,\alpha}_{<\beta}
\Phi^l_{<\beta}(\psi_\alpha,\eta_\alpha)
+\eta_\alpha L^{l,\alpha}_{<\beta}B^l_{<\beta}
\right)\frac{d\beta}{\beta}
+O_{L_t^1L_x^\infty}(\epsilon^4\alpha^{-1}),
\end{split}
\end{equation}
where
$L^{l,\alpha}_{<\gamma}\coloneq
P_{<M\alpha}L_{<\gamma}$.
The error is the high-frequency tail term $B_\alpha^{\mathrm{err}}$ from the decomposition
preceding \eqref{principalpartBpara}. Its $L_t^1L_x^{\infty}$ bound follows from
Proposition \ref{Bernsteinelliptic}. 

For $\alpha\ll_M\lambda$, the remaining terms have
frequency below $\lambda$ and are annihilated by
$P_{\geq\lambda}$. Thus, after estimating these tails,
the first integral in the expression for $d_B$ only contributes
at scales $\alpha\approx_M\lambda$.
We next estimate all contributions (of both integrals) in the range $\lambda\lesssim \alpha<\delta_0\mu$ directly.

For the first term on the right hand side of \eqref{BlinDtaexp}, the $B$ analogue of
\eqref{survivingprincipal} gives
\begin{equation*}
\begin{split}
L^{l,\alpha}_{<\delta\alpha}
\Phi^l_{<\delta\alpha}(\psi_\alpha,\eta_\alpha)
-T_{\mathfrak l}\varphi_\alpha=
L^{l,\alpha}_{<\delta\alpha}
\bigl(
\Phi^l_{<\delta\alpha}(\psi_\alpha,\eta_\alpha)
-\varphi_\alpha
\bigr)
+\bigl(L^{l,\alpha}_{<\delta\alpha}
-T_{\mathfrak l}\bigr)\varphi_\alpha.
\end{split}
\end{equation*}
The last term also contains $\varphi_\alpha$ and an additional
$\nabla^2\eta_\beta$ factor, by the application of
Proposition \ref{paraexp}. Therefore, as with the commutator above, its contribution is perturbative and satisfies \eqref{Lakappaprincipalerror}. 
For the first term, by using 
\eqref{Phidifferror} and slightly adapting the estimates that follow, we see that $L^{l,\alpha}_{<\delta\alpha}
\bigl(
\Phi^l_{<\delta\alpha}(\psi_\alpha,\eta_\alpha)
-\varphi_\alpha
\bigr)$ also provides two unbalanced bilinear pairings, and thus contributes a perturbative error of the form \eqref{Lakappaprincipalerror}. Similarly, the term $\eta_\alpha L^{l,\alpha}_{<\delta\alpha}B^l_{<\delta\alpha}$ also provides two unbalanced bilinear pairings, $(u_{\mu},\eta_{\alpha})$ and $(\overline{u}_{\mu}^{x_0},L^{l,\alpha}_{<\delta\alpha}B^l_{<\delta\alpha})$, the latter of which can be estimated using Proposition \ref{VBunbil}. The quadratic terms in the differentiated integral in \eqref{BlinDtaexp}
satisfy the same estimates for $\lambda\lesssim\alpha<\delta_0\mu$: their two input frequencies satisfy
$\beta\approx_{\delta}\alpha$, which after applying the relevant linearization formulas for $L_{<\beta}^{l,\alpha}$, $\Phi_{<\beta}^l$ and $B_{<\beta}^l$ again provides two unbalanced bilinear pairs, and thus a contribution satisfying \eqref{Lakappaprincipalerror}. The cubic terms satisfy the same bound after placing
the additional $\nabla\eta$ factor in $L_{t,x}^\infty$.
Indeed, since $\beta\approx_\delta\alpha$, this absorbs
one derivative and leaves the same two unbalanced pairings.
\\

For $\alpha\geq\delta_0\mu$, the first line in
\eqref{BlinDtaexp} has output frequency comparable to
$\alpha$. Since $\lambda\ll\delta_0\mu$, its contribution and that
of $T_{\mathfrak l}\varphi_\alpha$ are annihilated by
$P_{<\lambda}$. It remains to consider
\begin{equation*}
-P_{<\lambda}
\int_{\delta_0\mu}^{\infty}
\int_{\delta\alpha}^{\alpha}
\beta\frac{d}{d\beta}
\left(
L^{l,\alpha}_{<\beta}
\Phi^l_{<\beta}(\psi_\alpha,\eta_\alpha)
+\eta_\alpha L^{l,\alpha}_{<\beta}B^l_{<\beta}
\right)
\frac{d\beta}{\beta}\frac{d\alpha}{\alpha}.
\end{equation*}
Applying the $B$ analogue of \eqref{TVsecondFTC} and following the error analysis leading to \eqref{TBmain2ord}, we obtain
\begin{equation*}
\begin{split}
\beta\frac{d}{d\beta}
\left(
L^{l,\alpha}_{<\beta}
\Phi^l_{<\beta}(\psi_\alpha,\eta_\alpha)
+\eta_\alpha L^{l,\alpha}_{<\beta}B^l_{<\beta}
\right)=
\mathcal T^{B,\mathrm{main}}_{<\delta\beta}
+O_{L_{t,x}^2}(\epsilon^3\alpha^{-\frac74}).
\end{split}
\end{equation*}
Here $\mathcal T^{B,\mathrm{main}}_{<\delta\beta}$ is the
expression defined by \eqref{TBmain2ord}.

The error is at least cubic. It consists of the last two
terms in the $B$ analogue of \eqref{TVsecondlinearization},
the integrated term in \eqref{TVsecondFTC}, and
the errors from replacing $\Phi$ and the elliptic operators
by their principal parts. We estimate these in the usual way using two
$L_t^4L_x^\infty$ bounds and one $L_t^\infty L_x^2$ bound. These are the same as the estimates used following
\eqref{TVsecondFTC} and \eqref{TV3maincont}, together with
\eqref{L2phierr}. It is important here that $\alpha\geq\delta_0\mu$, since
integrating the $L_{t,x}^2$ error gives
\begin{equation*}
\epsilon^3
\int_{\alpha\geq\delta_0\mu}
\int_{\delta\alpha}^{\alpha}
\alpha^{-\frac74}
\frac{d\beta}{\beta}\frac{d\alpha}{\alpha}
\lesssim_{\delta,\delta_0}
\epsilon^3\mu^{-\frac74}.
\end{equation*}
This is precisely the reason for isolating the $\alpha\geq\delta_0\mu$ region above; we could not estimate directly the low frequency errors directly in $L_{t,x}^2$ for this reason. In summary, combining everything above, we have shown that
\begin{equation*}
\begin{split}
T_{\mathfrak l}T_{a_{<\lambda}}d_B
&=T_{\mathfrak l}T_{a_{<\lambda}}\Bigg[
-P_{<\lambda}
\int_{\delta_0\mu}^{\infty}
\int_{\delta\alpha}^{\alpha}
\mathcal T^{B,\mathrm{main}}_{<\delta\beta}
\frac{d\beta}{\beta}\frac{d\alpha}{\alpha}
+O_{L_{t,x}^2}\left(\epsilon^3\mu^{-\frac74}\right)
\Bigg]+\mathcal R,
\end{split}
\end{equation*}
where $\mathcal R$ satisfies \eqref{Lakappaprincipalerror}.
\\

The analogous expansion for $d_V$ is identical, using
\eqref{VtailFTCstart}, \eqref{TVlinearize},
\eqref{TVsecondFTC}, and \eqref{Vquadraticprincipal}.
Its principal quadratic contribution is
\begin{equation*}
-P_{<\lambda}
\int_{\delta_0\mu}^{\infty}
\int_{\delta\alpha}^{\alpha}
\mathcal T^{V,\mathrm{main}}_{<\delta\beta}
\frac{d\beta}{\beta}\frac{d\alpha}{\alpha}.
\end{equation*}
The preceding estimates apply almost verbatim, except we keep the additional $\nabla\eta_{<\lambda}$ para-coefficient outside the
projection and estimate it in $L_{t,x}^\infty$.

It remains to consider
$\mathcal B_{<\lambda}(V,\nabla\eta)$.
We use the representation \eqref{bilformrep}. The part with input frequencies
$\lesssim_M\lambda$ contributes to \eqref{Lakappaprincipalerror} by two unbalanced
bilinear estimates and Proposition \ref{VBunbil}. For the remaining high-high contribution, we replace (where the input frequency of $V$ is taken to be $\alpha$) $V_\alpha$ by $T_{\mathfrak m}\varphi_\alpha$, using
\eqref{Fauxfactorbounds} and Proposition
\ref{Leibnizcommutator}. If $\alpha<\delta_0\mu$, we put
the replacement error in $L_t^2L_x^\infty$, pair the
$\nabla\eta$ factor with one $u_\mu$ factor, and put
the other $u_\mu$ factor in $L_t^{\infty}L_x^2$. This again gives a perturbative contribution satisfying \eqref{Lakappaprincipalerror}.
For $\alpha\geq\delta_0\mu$, we instead put the
$\nabla\eta$ factor in $L_t^\infty L_x^2$, which contributes the direct $L_{t,x}^2$ bound
\begin{equation*}
\int_{\alpha\geq\delta_0\mu}
\left(\epsilon^2\alpha^{-\frac34}\right)
\left(\epsilon\alpha^{-1}\right)
\frac{d\alpha}{\alpha}
\lesssim_{\delta,\delta_0}
\epsilon^3\mu^{-\frac74}.
\end{equation*}
It remains to estimate the principal quadratic contribution in $\mathcal{B}_{<\lambda}(T_{\mathfrak{m}}\varphi,\nabla\eta)$ in the region where one of the input frequencies $\alpha$ satisfies (say) $\alpha\geq M\lambda$. First, for $M\lambda\leq\alpha<\delta_0\mu$, we can use two unbalanced bilinear estimates by pairing the
$T_{\mathfrak m}\varphi_{\alpha}$ and $\nabla\eta$ factor (which has comparable frequency)
with the two $u_\mu$ factors. Hence, these terms are also included in $\mathcal R$. The remaining quadratic contribution is given by
\begin{equation*}
-P_{<\lambda}
\int_{\delta_0\mu}^{\infty}
\int_{\delta\alpha}^{\alpha}
\left(
\nabla\eta_\beta\cdot T_{\mathfrak m}\varphi_\alpha
+\nabla\eta_\alpha\cdot T_{\mathfrak m}\varphi_\beta
\right)
\frac{d\beta}{\beta}\frac{d\alpha}{\alpha}.
\end{equation*}
Combining the three contributions according to
$\kappa_{<\lambda}
=-d_B+\nabla\eta_{<\lambda}\cdot d_V
+\mathcal B_{<\lambda}(V,\nabla\eta)$
proves \eqref{Lakappaprincipalexp} and
\eqref{Lakappaprincipalerror}.
\end{proof}
Now, we have all of the ingredients we need to prove Proposition \ref{Dtacommutatorestimates}. \begin{proof}[Proof of Proposition \ref{Dtacommutatorestimates}]
By Proposition \ref{Dtaestimates}, it remains to treat the term
\begin{equation*}
\frac{i}{2}\overline u_\mu^{x_0}P_\mu
\left(T_{\frac{\mathfrak q}{a}L(a\kappa)}\varphi\right).
\end{equation*}
We begin by applying Lemma \ref{Lakappaprincipal}, placing the error term $\mathcal{R}$ in
\eqref{Lakappaprincipalerror} in $\mathcal R_{a,\mu}$. We next write the remaining principal contribution as
\begin{equation}\label{Dtaprincipalrho}
\begin{split}
\int_0^\infty \mathcal I_{a,\mu,\rho}\frac{d\rho}{\rho}
\coloneq
\frac{i}{2}\overline u_\mu^{x_0}P_\mu
\int_0^\infty\int_{\nu\approx\mu}
Op\left(
P_{<M\delta\nu}\left(
\frac{\mfq_{<\delta\nu}}{a_{<\delta\nu}}F_{\rho,\delta\nu}
\right)\right)
\varphi_\nu\frac{d\nu}{\nu}\frac{d\rho}{\rho},
\end{split}
\end{equation}
where we define $F_{\rho,\lambda}$ to be the expression
\begin{equation*}
\begin{split}
T_{\mathfrak l}T_{a_{<\lambda}}S_\rho
\int_{\delta_0\mu}^{\infty}\int_{\delta\alpha}^{\alpha}
\Bigl[
P_{<\lambda}\Bigl(
\mathcal T^{B,\mathrm{main}}_{<\delta\beta}
-\nabla\eta_\beta\cdot T_{\mathfrak m}\varphi_\alpha
-\nabla\eta_\alpha\cdot T_{\mathfrak m}\varphi_\beta
\Bigr)-\nabla\eta_{<\lambda}\cdot
P_{<\lambda}\mathcal T^{V,\mathrm{main}}_{<\delta\beta}
\Bigr]\frac{d\beta}{\beta}\frac{d\alpha}{\alpha}
\end{split}
\end{equation*}
where above, we used the identity $\int_0^\infty S_\rho\,d\rho/\rho=1$.
For each $\rho$, we first replace all of the paradifferential
operators $T_{\mathfrak a}$ inside the square brackets and the
coefficient $\nabla\eta_{<\lambda}$ with
$Op(\mathfrak a_{<\delta\rho})$ and
$\nabla\eta_{<\delta\rho}$, respectively. We claim all of the errors contribute either a perturbative of the form $\mathcal{R}_{a,\mu}$ in  \eqref{Dtacommutatorfluxbound} or they can be estimated by $O_{L_{t,x}^2}(\epsilon^3\rho\mu^{-\frac{7}{4}})$ in which case it contributes to the flux correction $Q_{a,\mu}$.

Indeed, applying the fundamental theorem of calculus in the frequency parameter to the differences of these expressions introduces
an additional $\nabla\eta_\gamma$ factor, with
$\delta\rho\lesssim\gamma\lesssim\alpha$. To estimate this, we consider separately the regions where $\gamma<\delta_0\mu\ll\mu$ or
$\gamma\geq\delta_0\mu\ll\mu$, where $\delta_0>0$ is some fixed, small constant. In the first range, we pair
$\nabla\eta_\gamma$ with one frequency $\mu$ factor ($u$ or $|D|^{\frac{1}{2}}\varphi$) and use one unbalanced bilinear estimate. We can then, as usual  place the remaining factors
$\eta_\beta$ and $\varphi_\alpha$ in $L_t^4L_x^\infty$,
and the other high frequency factor in
$L_t^\infty L_x^2$. Denoting the resulting contribution to
$\mathcal I_{a,\mu,\rho}$ by
$\mathcal I_{a,\mu,\rho}^{\mathrm{err}}$, we obtain
\begin{equation*}
\begin{split}
\|\mathcal I_{a,\mu,\rho}^{\mathrm{err}}\|_{L_{t,x}^1}
&\lesssim
\epsilon^5c_\mu^2\mu^{-2s_c}(1+\mu|x_0|)^{\frac12}
\rho\int_{\delta_0\mu}^{\infty}\alpha^{-\frac34}
\int_{\delta\rho\lesssim\gamma<\delta_0\mu}
\gamma^{-\frac14}\frac{d\gamma}{\gamma}\frac{d\alpha}{\alpha}
\\
&\lesssim
\epsilon^5c_\mu^2\mu^{-2s_c}(1+\mu|x_0|)^{\frac12}
\left(\frac{\rho}{\mu}\right)^{\frac34}.
\end{split}
\end{equation*}
Here the factor $\rho$ comes from the outer $T_{\mathfrak l}$,
and the integration in $\beta\approx_\delta\alpha$ and the boundedness of $a_{<\lambda}$ has been
absorbed into the implicit constant. Commuting
$\nabla\eta_{<\delta\rho}$ through $S_\rho P_{<\lambda}$
also produces a directly perturbative error by using Proposition \ref{Leibnizcommutator} and the same bilinear pairing, and $L_t^4L_x^{\infty}$ and $L_t^{\infty}L_x^2$ estimate as above. In summary, all the above errors contribute to $\mathcal{R}_{a,\mu}$ in \eqref{Dtacommutatorfluxbound}.

On the other hand, for the range $\gamma\geq\delta_0\mu$, we instead place
$\nabla\eta_\gamma$ direclty in $L_t^\infty L_x^2$ and the other
two factors at frequency $\alpha$ and $\beta$ in $L_t^4L_x^\infty$. Thus the corresponding
coefficient error $\mathcal E_{\rho,\lambda}$, including
the outer operators, satisfies
\begin{equation*}
\begin{split}
\|\mathcal E_{\rho,\lambda}\|_{L_{t,x}^2}
&\lesssim_{\delta,\delta_0}
\epsilon^3\rho
\int_{\delta_0\mu}^{\infty}\alpha^{-\frac34}
\int_{\delta_0\mu\lesssim\gamma\lesssim\alpha}
\gamma^{-1}\frac{d\gamma}{\gamma}\frac{d\alpha}{\alpha}
\\
&\lesssim_{\delta,\delta_0}
\epsilon^3\rho\mu^{-\frac74}.
\end{split}
\end{equation*}
We also include in $\mathcal E_{\rho,\lambda}$ the errors from
replacing compositions of the for $Op(\mathfrak{a})Op(\mathfrak{b})$ by $Op(\mathfrak{a}\mathfrak{b}$) for $\mathfrak{a},\mathfrak{b}\in\{\mathfrak{m},\mathfrak{l},\mathfrak{g}\}$ (see below).
These contain an additional $\nabla^2\eta$ factor and one
less derivative on the high-frequency inputs. By placing
$\varphi_\alpha$ in $L_t^\infty L_x^2$ and the two $\eta$
factors in $L_t^4L_x^\infty$, their contribution satisfies the same estimate
\begin{equation*}
\|\mathcal E_{\rho,\lambda}^{\mathrm{comp}}\|_{L_{t,x}^2}
\lesssim_{\delta,\delta_0}
\epsilon^3\rho
\int_{\delta_0\mu}^{\infty}\alpha^{-\frac{21}{8}}
\int_{\gamma\lesssim\alpha}\gamma^{\frac78}
\frac{d\gamma}{\gamma}\frac{d\alpha}{\alpha}
\lesssim_{\delta,\delta_0}
\epsilon^3\rho\mu^{-\frac74}.
\end{equation*}
We now study the contribution of the remaining principal term after the replacement errors estimated above. Here, we obtain a $L_{t,x}^2$ bound similar to the one above. Here, we have to be more careful because the remaining term is quadratic. The key observation is that this term exhibits a null structure similar to the one which was used in the estimate \eqref{nullL2bound}. To see this, we recall from \eqref{glmidentities}, the identities
\begin{equation}\label{glmidentities2}
\mathfrak l_{<\delta\rho}-\nabla\eta_{<\delta\rho}\cdot\mathfrak m_{<\delta\rho}=\mathfrak g_{<\delta\rho},
\qquad
\mathfrak g_{<\delta\rho}\mathfrak l_{<\delta\rho}=-i\xi\cdot\mathfrak m_{<\delta\rho}.
\end{equation}
The first identity combines (modulo the errors above) the outer operators in
$\mathcal T^{B,\mathrm{main}}$ and
$\mathcal T^{V,\mathrm{main}}$ into $-Op(\mathfrak g)$.
For the remaining terms, the second identity gives
\begin{equation*}
\begin{split}
\eta_\beta Op(\mathfrak g\mathfrak l)\varphi_\alpha
-\nabla\eta_\beta\cdot Op(\mathfrak m)\varphi_\alpha
&=-\nabla\cdot
\bigl(\eta_\beta Op(\mathfrak m)\varphi_\alpha\bigr)
+\eta_\beta Op(\nabla_x\cdot\mathfrak m)\varphi_\alpha.
\end{split}
\end{equation*}
The last term contains $\nabla^2\eta_{<\delta\rho}$ and is
included in $\mathcal E_{\rho,\lambda}$ by the composition
error estimate above. The same calculation applies with
$\alpha$ and $\beta$ interchanged. Consequently, we have a decomposition of the form
\begin{equation}\label{Flambdarhodecomp}
F_{\rho,\lambda}
=F_{\rho,\lambda}^{\mathrm{null}}
+\mathcal E_{\rho,\lambda}+\mathcal R_{\rho,\lambda},
\end{equation}
where $\mathcal R_{\rho,\lambda}$ produces the directly
perturbative contribution $\mathcal I_{a,\mu,\rho}^{\mathrm{err}}$,
and the remaining integral is
\begin{equation*}
\begin{split}
F_{\rho,\lambda}^{\mathrm{null}}
&=-T_{\mathfrak l}T_{a_{<\lambda}}S_\rho P_{<\lambda}
\int_{\delta_0\mu}^{\infty}\int_{\delta\alpha}^{\alpha}
Op\Bigl(
\mathfrak g_{<\delta\rho}(x,\xi+\zeta)
\mathfrak l_{<\delta\rho}(x,\zeta)
+i(\xi+\zeta)\cdot\mathfrak m_{<\delta\rho}(x,\zeta)
\Bigr)(\eta_\beta,\varphi_\alpha)
\frac{d\beta}{\beta}\frac{d\alpha}{\alpha}
\\
&-T_{\mathfrak l}T_{a_{<\lambda}}S_\rho P_{<\lambda}
\int_{\delta_0\mu}^{\infty}\int_{\delta\alpha}^{\alpha}
Op\Bigl(
\mathfrak g_{<\delta\rho}(x,\xi+\zeta)
\mathfrak l_{<\delta\rho}(x,\zeta)
+i(\xi+\zeta)\cdot\mathfrak m_{<\delta\rho}(x,\zeta)
\Bigr)(\eta_\alpha,\varphi_\beta)
\frac{d\beta}{\beta}\frac{d\alpha}{\alpha}.
\end{split}
\end{equation*}
Here $Op$ denotes the bilinear operator, with $\xi$ and $\zeta$
the frequencies of the $\eta$ and $\varphi$ inputs, respectively.
By \eqref{glmidentities2} and homogeneity, the quantized symbols above vanish when $\xi=-\zeta$. Thus after factoring out $\xi+\zeta$ in the symbol, we can expand $\eta$ and $\varphi$ in terms of $u$ and $\overline{u}$ and then apply the balanced estimate Lemma \ref{balancedrestrictedoutput}
and the bound $\|T_{\mathfrak{l}}T_{a_{<\lambda}}S_{\rho}\|_{L_{t,x}^2\to L_{t,x}^2}\lesssim \rho$ to obtain
\begin{equation}\label{Fnull} 
\|F_{\rho,\lambda}^{\mathrm{null}}\|_{L_{t,x}^2}
\lesssim
\epsilon^2\rho^{\frac32}\mu^{-\frac94}.
\end{equation}
With the above estimates, we are now ready to prove the error-flux decomposition \eqref{Dtacommutatorfluxdecomp}. We write (suppressing $\lambda$ for convenience) $F_\rho$ to denote the sum of
$F_{\rho,\lambda}^{\mathrm{null}}$, $\mathcal E_{\rho,\lambda}$,
and the $L_{t,x}^2$ error in \eqref{Lakappaprincipalexp} after
applying $T_{\mathfrak l}T_{a_{<\lambda}}S_\rho$. The preceding
estimates give
\begin{equation*}
\|F_\rho\|_{L_{t,x}^2}
\lesssim_{\delta,\delta_0}
\epsilon^2\rho^{\frac32}\mu^{-\frac94}
+\epsilon^3\rho\mu^{-\frac74}
\lesssim_{\delta,\delta_0}
\epsilon^2\rho\mu^{-\frac74},
\end{equation*}
where we used $\rho\lesssim\lambda\ll\mu$. We now argue as in the flux decomposition in the proof of \eqref{QPhibound}. We first write
\begin{equation*}
\begin{split}
T_{\frac{\mfq}{a}F_\rho}\varphi=T_{\nabla_x\cdot\left(\frac{\mfq}{a}
 |\xi|^{-\frac{1}{2}}\nabla\Delta^{-1}F_\rho\right)}|D|^{\frac{1}{2}}\varphi
-T_{\nabla_x\left(|\xi|^{-\frac{1}{2}}\frac{\mfq}{a}\right)\cdot
\nabla\Delta^{-1}F_\rho}|D|^{\frac{1}{2}}\varphi.
\end{split}
\end{equation*}
The second term supplies (from the derivative on the symbol) either an additional factor of $\nabla^2\eta_{<\lambda}$ or $\nabla a_{<\lambda}$. This gives a directly perturbative contribution. Indeed, by pairing
the differentiated coefficient with $\overline u_\mu^{x_0}$
using Proposition \ref{abilinearbounds} or a direct unbalanced
estimate for $\nabla^2\eta$ and combining this with the
$L_{t,x}^2$ bound for $F_\rho$ by applying Bernstein to $\nabla\Delta^{-1}F_{\rho}$ in $L_t^2L_x^{\infty}$, and placing  $|D|^{\frac12}\varphi$ in $L_t^{\infty}L_x^2$ (which is at frequency $\approx\mu$), we get
\begin{equation*}
\begin{split}
\left\|
P_\mu\left(
T_{\nabla_x\left(\frac{\mfq}{a}|\xi|^{-\frac{1}{2}}\right)\cdot
\nabla\Delta^{-1}F_\rho}|D|^{\frac{1}{2}}\varphi
\right)\overline u_\mu^{x_0}
\right\|_{L_{t,x}^1}&\lesssim_{\delta,\delta_0}
\epsilon^3c_\mu^2\mu^{-2s_c}(1+\mu|x_0|)^{\frac12}
\mu^{\frac34}\|F_\rho\|_{L_{t,x}^2}
\\
&\lesssim_{\delta,\delta_0}
\epsilon^5c_\mu^2\mu^{-2s_c}(1+\mu|x_0|)^{\frac12}
\frac{\rho}{\mu}
\end{split}
\end{equation*}
which is integrable over $\rho\lesssim\lambda\ll\mu$, and gives the required $L_{t,x}^1$ perturbative bound. It remains to consider the first term. As in the proofs of
\eqref{Netaauxtrilinear} and \eqref{QPhibound}, we can use the
divergence on the low frequency paracoefficient, together
with the estimates above, to obtain a decomposition
\begin{equation*}
\overline u_\mu^{x_0}P_\mu\left(
T_{\nabla_x\cdot\left(\frac{\mfq}{a}
\nabla\Delta^{-1}F_\rho\right)}\varphi
\right)
=\mathcal R_{a,\mu,\rho}+\nabla\cdot Q_{a,\mu,\rho},
\end{equation*}
where, for $\frac43<p\leq2$,
\begin{equation*}
\begin{split}
\|\mathcal R_{a,\mu,\rho}\|_{L_{t,x}^1}
+\mu^{\frac2p-1}\|Q_{a,\mu,\rho}\|_{L_t^1L_x^p}
\lesssim_{\delta,\delta_0,p}
\epsilon^4c_\mu^2\mu^{-2s_c}(1+\mu|x_0|)^{\frac12}
\left(\frac{\rho}{\mu}\right)^{\frac32-\frac2p}.
\end{split}
\end{equation*}
Since $\frac32-\frac2p>0$, integrating in
$\rho\lesssim\lambda\ll\mu$ and including the directly
perturbative errors above gives
\eqref{Dtacommutatorfluxdecomp} and
\eqref{Dtacommutatorfluxbound}. We leave the detailed verification of this last part to the reader, as it is similar to before.
\end{proof}
We end this section by recording two basic bilinear estimates satisfied by $D_t^la_{<\delta\lambda}$.
\begin{corollary}[Bilinear bound for $D_t^l a$]
\label{Dtabilinearcorollary}
Assume the bootstrap bounds in Section
\ref{sec:bootstrap}. Write $D_t^l=(\partial_t+V^l_{<\delta\lambda}\cdot\nabla)$. For every
$0<\lambda\lesssim\mu$ and $x_0\in\mathbb R^2$, there holds
\begin{equation}\label{Dtabilinearbound}
\left\|
\left(D_t^l a_{<\delta\lambda}\right)
u_\mu^{x_0}
\right\|_{L_{t,x}^2}
\lesssim_\delta
\epsilon^2c_\mu\mu^{-s_c}\lambda^{\frac14}(1+\lambda|x_0|)^{\frac12},
\hspace{5mm}
V_{[\lambda]}=P_{<M\delta\lambda}V_{<\delta\lambda}.
\end{equation}
The same estimate holds with
$u_\mu^{x_0}$ replaced by $(|D|^{\frac12}\varphi_\mu)^{x_0}$. We also have the following linearized variant 
\begin{equation}\label{Dtabilinearlin}
\begin{split}
\|\bigl(D_t^l a_{<\delta\lambda}\bigr)w_\mu^{x_0}\|_{L_{t,x}^2}\lesssim_\delta
\epsilon d_\mu\lambda^{\frac14}(1+\mu|x_0|)^{\frac12}.
\end{split}
\end{equation}
The same estimate holds with $w_{\mu}^{x_0}$ replaced by $(|D|^{\frac{1}{2}}r_{\mu})^{x_0}$
\end{corollary}
\begin{proof}
This is a simple corollary of \eqref{DtaprincipalL2L8}. For each
$\sigma>0$, taking $D_t^l=\partial_t+V^l_{<\sigma}\cdot\nabla$, we can write
\begin{equation*}
\sigma\frac{d}{d\sigma}(D_t^l a_{<\sigma})
=D_t^l a_\sigma
+(\sigma\partial_\sigma V^l_{<\sigma})\cdot\nabla a_{<\sigma}.
\end{equation*}
Taking $p=\infty$ in \eqref{DtaprincipalL2L8} gives
\begin{equation*}
\|D_t^l a_\sigma-T_{a^l}\nabla\cdot V_\sigma\|_{L_t^2L_x^\infty}
\lesssim_\delta\epsilon c_\sigma\sigma^{\frac14}.
\end{equation*}
For the transport term, two $L_t^4L_x^\infty$ bounds also give
\begin{equation*}
\|(\sigma\partial_\sigma V^l_{<\sigma})\cdot\nabla a_{<\sigma}
\|_{L_t^2L_x^\infty}
\lesssim_\delta
(\epsilon\sigma^{-\frac58})(\epsilon\sigma^{\frac78})
=\epsilon^2\sigma^{\frac14}.
\end{equation*}
For these two error terms, we place the remaining factor of $u_{\mu}^{x_0}$ or $w_{\mu}^{x_0}$ in $L_t^{\infty}L_x^2$ while for the contribution of $T_{a^l}\nabla\cdot V_{\sigma}$, we apply one
the unbalanced bilinear estimate using Proposition \ref{VBunbil}. Integrating in $0<\sigma<\delta\lambda$ gives
\eqref{Dtabilinearbound} and \eqref{Dtabilinearlin}. The same argument applies to
$(|D|^{\frac12}\varphi_\mu)^{x_0}$ and
$(|D|^{\frac12}r_\mu)^{x_0}$ using their corresponding
bilinear bounds. 
\end{proof}

\section{The equation for the complex variable \texorpdfstring{$u$}{}}
\label{s:u}
Here, we finally derive the refined paradifferential equation satisfied by the complex good variable $u=\eta+iT_{\mfq}\varphi$ from Definition \ref{complexvariable}. Our starting point is the paradifferential system for $(\eta,\varphi)$ in Theorems \ref{Netafinalcorollary} and \ref{Nphifinalcorollary},
\begin{equation}\label{Z-para-re}
\begin{cases}
\partial_t\eta+\nabla\cdot(T_{V^l}\eta)-T_G\varphi=N_{\eta},
\\
\left(\partial_t+T_{V^l}\cdot\nabla\right)\varphi+T_{a^l}\eta
=N_{\varphi}+N_{\varphi}^{\mathrm{flux}}+\mathcal R.
\end{cases}
\end{equation}
We will use $T_{\mfq}$ to symmetrize the principal terms. We will want to keep the principal contribution of the commutator of $T_{\mfq}$ with the (para) material derivative on the left-hand side. For this purpose, we define the auxiliary symbol $\bb$. At input frequency $\lambda$, writing $V_{[\lambda]}=P_{<M\delta\lambda}V_{<\delta\lambda}$ for the coefficient in $T_{V^l}$, we set\footnote{This is related 
to the notion of paramaterial derivative 
introduced later in Section~\ref{sec:para-functions}.}
\begin{equation}\label{ublogsymbol}
\bb_{<\delta\lambda}
\coloneq
\frac12\left(
\partial_t+V_{[\lambda]}\cdot\nabla_x
-\sum_{j,k=1}^2
\partial_{x_j}V_{[\lambda]}^k\,\xi_k\partial_{\xi_j}
\right)\ln\mfq_{<\delta\lambda}.
\end{equation}
Before deriving the equation for $u$, we first isolate the subprincipal error in the aforementioned commutator. This will still be quadratic at leading order, and not perturbative, but will have a more favorable balanced of derivatives (for normal form corrections) compared the principal part. 
\begin{lemma}\label{ulogcommutator}
Define
\begin{equation*}
\{V\cdot\xi,\mfq\}_{<\delta\lambda}\coloneq V_{[\lambda]}\cdot\nabla_x\mfq_{<\delta\lambda}
-\sum_{j,k}\partial_{x_j}V_{[\lambda]}^k\,\xi_k
\partial_{\xi_j}\mfq_{<\delta\lambda}.
\end{equation*}
Write $\mathcal{C} \coloneq i[\partial_t+T_{V^l}\cdot\nabla,T_{\mfq}]\varphi-T_{\bb}u$. Then the subprincipal commutator error $\mathcal{C}$ admits a decomposition
\begin{equation}\label{ulogcommdecomp}
\begin{split}
\mathcal C=
-T_{\bb}\overline u
+i\left([T_{V^l}\cdot\nabla,T_{\mfq}]
-T_{\{V\cdot\xi,\mfq\}}\right)\varphi
+\mathcal R_{\mathrm{comm}}
\end{split}
\end{equation}
where for every $\mu>0$ and $x_0\in\mathbb R^2$,
\begin{equation}\label{ulogcommrem}
\|(\mathcal R_{\mathrm{comm}})_\mu\overline u_\mu^{x_0}\|_{L_{t,x}^1}
\lesssim_\delta\epsilon^4c_\mu^2\mu^{-2s_c}(1+\mu|x_0|)^{\frac12}.
\end{equation}
We also have the following auxiliary bounds:
\begin{enumerate}
\item For every $\lambda>0$, there holds
\begin{equation}\label{ulogcommL2Linf}
\left\|
\mathcal{C}_\lambda
\right\|_{L_t^2L_x^\infty}
\lesssim_\delta\epsilon^2\lambda^{-\frac34},\hspace{5mm}\left\|
\mathcal{C}_{\lambda}
\right\|_{L_{t,x}^2}
\lesssim_\delta\epsilon^2c_\lambda\lambda^{\frac14-s_c}.
\end{equation}
\item For every $0<\lambda\ll\mu$ and $x_0\in\mathbb R^2$,
there holds
\begin{equation}\label{ulogcommauxtrilinear}
\left\|
\mathcal{C}_{\lambda}
\overline u_\mu^{x_0}u_\mu
\right\|_{L_{t,x}^1}
\lesssim_\delta
\epsilon^4c_\mu^2\lambda^{-1}\mu^{-2s_c}(1+\mu|x_0|)^{\frac12}.
\end{equation}
\end{enumerate}
\end{lemma}

\begin{proof}
We have
$[\partial_t,T_{\mfq}]=T_{\partial_t\mfq}$. Moreover, a simple computation using the definition of $\bb_{<\delta\lambda}$ \eqref{ublogsymbol} gives
\begin{equation*}
\partial_t\mfq_{<\delta\lambda}+\{V\cdot\xi,\mfq\}_{<\delta\lambda}
=2\bb_{<\delta\lambda}\mfq_{<\delta\lambda}.
\end{equation*}
Hence, using $u-\overline u=2iT_{\mfq}\varphi$, we obtain
\eqref{ulogcommdecomp} with
\begin{equation*}
\mathcal R_{\mathrm{comm}}
=2i(T_{\bb\mfq}-T_{\bb}T_{\mfq})\varphi.
\end{equation*}
Now, we prove the estimates for $\mathcal{R}_{\mathrm{comm}}$ and $\mathcal{C}$ above. For this, it will be slightly convenient to rescale $\mathfrak{q}$ to be order zero by defining $m_{<\delta\lambda}=|\xi|^{-1/2}\mfq_{<\delta\lambda}$. We also define $f=|D|^{1/2}\varphi$. This gives $\mathcal R_{\mathrm{comm}}=2i(T_{\bb m}-T_{\bb}T_m)f$. Now, fix the output frequency $\mu>0$ in \eqref{ulogcommrem}.

We begin by recording some simple bilinear estimates for the coefficients $\bb$ and $m$. For $|\xi|\approx\lambda$, $\sigma\approx\lambda$, and $\lambda\lesssim\mu$, we have
\begin{equation}\label{ulogcommcoefficients}
\begin{split}
\|\bb_{<\delta\lambda}(\cdot,\xi)u_\mu^z\|_{L_{t,x}^2}
&\lesssim_\delta \epsilon^2c_\mu\mu^{-s_c}\lambda^{\frac14}(1+\lambda|z|)^{\frac12},
\\
\|(\nabla_xm_{<\delta\lambda})(\cdot,\xi)f_\sigma^z\|_{L_{t,x}^2}
&\lesssim_\delta \epsilon^2c_\sigma\sigma^{-s_c}\lambda^{\frac34}(1+\lambda|z|)^{\frac12},
\\
\|(m_{<\delta\lambda}-m_{<\delta\sigma})(\cdot,\xi)f_\sigma^z\|_{L_{t,x}^2}
&\lesssim_\delta \epsilon^2c_\sigma\sigma^{-s_c}\lambda^{-\frac14}(1+\lambda|z|)^{\frac12}.
\end{split}
\end{equation}
Indeed, for the first estimate, $\bb_{<\delta\lambda}$ is a sum of bounded coefficients
multiplied by $\nabla V_{[\lambda]}$,
$(\partial_t+V_{[\lambda]}\cdot\nabla)a_{<\delta\lambda}$, and
$(\partial_t+V_{[\lambda]}\cdot\nabla)\nabla\eta_{<\delta\lambda}$.
For the $\nabla V_{[\lambda]}$ we use one unbalanced bilinear estimate in Proposition \ref{VBunbil}. For the material derivative
of $a_{<\delta\lambda}$ we use Corollary \ref{Dtabilinearcorollary}. For the material derivative of $\nabla\eta_{<\delta\lambda}$, we have the simple expansion
\begin{equation*}
\begin{split}
(\partial_t+V_{[\lambda]}\cdot\nabla)\nabla\eta_{<\delta\lambda}
&=\nabla B_{<\delta\lambda}+\nabla\kappa_{<\delta\lambda}
-\sum_k(\nabla V_{<\delta\lambda}^k)\partial_k\eta_{<\delta\lambda}
\\
&+(V_{[\lambda]}-V_{<\delta\lambda})\cdot\nabla\nabla\eta_{<\delta\lambda}.
\end{split}
\end{equation*}
We recall that  $\kappa_{<\delta\lambda}$ is defined by \eqref{kappadef} and measures the regularization error in the kinematic boundary condition. The left-hand side of this equation has frequency $\lesssim M\delta\lambda$, so we may insert
$P_{<CM\delta\lambda}$ on the right-hand side for some constant $C>1$. Proposition \ref{VBunbil}
then applies to $\nabla B_{<\delta\lambda}$ and $\nabla V_{<\delta\lambda}$, noting that in the latter
product $\partial_k\eta_{<\delta\lambda}$ has frequency $\lesssim\delta\lambda$. Moreover, the estimate preceding \eqref{kappaubil} gives
\begin{equation*}
\|\nabla P_{<CM\delta\lambda}\kappa_{<\delta\lambda}\|_{L_t^2L_x^\infty}
\lesssim_M \delta\lambda\|\kappa_{<\delta\lambda}\|_{L_t^2L_x^\infty}
\lesssim_\delta\epsilon^2\lambda^{\frac{1}{4}}.
\end{equation*}
The last term $(V_{[\lambda]}-V_{<\delta\lambda})\cdot\nabla\nabla\eta_{<\delta\lambda}$ satisfies the same
$L_t^2L_x^\infty$ bound by placing $V_{[\lambda]}-V_{<\delta\lambda}$ and $\nabla^2\eta_{<\delta\lambda}$ in $L_t^4L_x^\infty$.
For these errors we place $u_\mu$ in $L_t^\infty L_x^2$. This proves the
first estimate in \eqref{ulogcommcoefficients}.

For the second estimate, $m_{<\delta\lambda}$ is a smooth function of
$a_{<\delta\lambda}$, $\nabla\eta_{<\delta\lambda}$, and $\xi/|\xi|$. Applying $\nabla_x$
therefore gives either a factor of $\nabla a_{<\delta\lambda}$ or $\nabla^2\eta_{<\delta\lambda}$. Hence, we can use an unbalanced bilinear estimate, using Lemma \ref{hightailabound} and Proposition \ref{abilinearbounds} to handle the contributions of $\nabla a$.

For the third estimate in \eqref{ulogcommcoefficients}, the mean value formula gives a factor of
$\nabla(\eta_{<\delta\lambda}-\eta_{<\delta\sigma})$ or
$a_{<\delta\lambda}-a_{<\delta\sigma}$.
The first is supported at frequency $\approx\delta\lambda$ and, by Bernstein is handled in the same way as the contribution of $\nabla^2\eta$ in the second bound in \eqref{ulogcommcoefficients}.
For the Taylor coefficient difference, one simply repeats the proof of the bilinear estimates in Proposition \ref{abilinearbounds} by applying the fundamental theorem of calculus in the regularization variable and then suitably adjusting the integration bounds in \eqref{a-integral-principal}, and the estimates that follow. We leave the simple verification of this to the reader.

We now apply \eqref{ulogcommcoefficients} to the composition error
$\mathcal R_{\mathrm{comm}}=2i(T_{\bb m}-T_{\bb}T_m)f$. By definition, we can expand
\begin{equation}\label{paracomexp}
\begin{split}
&(T_{\bb m}-T_{\bb}T_m)f
\\
&=\int_0^\infty\int_0^\infty
\Bigl(
\operatorname{Op}\bigl(P_{<M\delta\lambda}
(\bb_{<\delta\lambda}m_{<\delta\lambda})\bigr)
-\operatorname{Op}(P_{<M\delta\lambda}\bb_{<\delta\lambda})
 \operatorname{Op}(P_{<M\delta\sigma}m_{<\delta\sigma})
\Bigr)
S_\lambda S_\sigma f\frac{d\lambda}{\lambda}\frac{d\sigma}{\sigma}
\\
&+\int_0^\infty\int_0^\infty
\operatorname{Op}(P_{<M\delta\lambda}\bb_{<\delta\lambda})
\bigl[\operatorname{Op}(P_{<M\delta\sigma}m_{<\delta\sigma}),S_\lambda\bigr]
S_\sigma f\frac{d\lambda}{\lambda}\frac{d\sigma}{\sigma}.
\end{split}
\end{equation}
We note that when projecting the above expression to frequency $\mu$, the integration ranges are restricted to $\sigma\approx\lambda\approx\mu$. In this region, for the commutator term with $S_\lambda$ above, Proposition \ref{Leibnizcommutator} gives a factor $\lambda^{-1}$ and places a spatial derivative on $m_{<\delta\sigma}$. Pairing the $\bb$ factor with $\overline u_\mu^{x_0}$ and the differentiated $m$ factor with $f_\sigma$, the first two bounds in \eqref{ulogcommcoefficients} give the required $L_{t,x}^1$ bound \eqref{ulogcommrem}. In the other integral in \eqref{paracomexp}, we first replace the composition of the two operators $\operatorname{Op}(P_{<M\delta\lambda}\bb_{<\delta\lambda})
 \operatorname{Op}(P_{<M\delta\sigma}m_{<\delta\sigma})$ by $\operatorname{Op}(P_{<M\delta\lambda}\bb_{<\delta\lambda}P_{<M\delta\sigma}m_{<\delta\sigma})$. By standard pseudo-differential calculus, the error that is generated by this replacement also differentiates the symbol $m_{<\delta\sigma}$ in space with a factor of $\lambda^{-1}$, and so we can again use two bilinear estimates. We then expand the remaining symbol as
\begin{equation*}
\begin{split}
P_{<M\delta\lambda}(\bb_{<\delta\lambda}m_{<\delta\lambda})
-(P_{<M\delta\lambda}\bb_{<\delta\lambda})
(P_{<M\delta\sigma}m_{<\delta\sigma})&=(P_{<M\delta\lambda}\bb_{<\delta\lambda})
(m_{<\delta\lambda}-m_{<\delta\sigma})
\\
&+(P_{<M\delta\lambda}\bb_{<\delta\lambda})
(m_{<\delta\sigma}-P_{<M\delta\sigma}m_{<\delta\sigma})
\\
&+[P_{<M\delta\lambda},m_{<\delta\lambda}]\bb_{<\delta\lambda}.
\end{split}
\end{equation*}
Here the last commutator differentiates $m_{<\delta\lambda}$ by Proposition \ref{Leibnizcommutator}. Hence, we can apply the estimates in \eqref{ulogcommcoefficients} to the contributions of the three terms above to see that the full $\mathcal{R}_{\mathrm{comm}}$ satisfies the $L_{t,x}^1$ bound \eqref{ulogcommrem}. In each case we pair the $\bb$ coefficient with
$\overline u_\mu^{x_0}$ and the differentiated or "differenced"
$m$ coefficient with $f_\sigma$. Combining everything gives
\begin{equation*}
\begin{split}
\|(\mathcal R_{\mathrm{comm}})_\mu\overline u_\mu^{x_0}\|_{L_{t,x}^1}&\lesssim_\delta \epsilon^4c_\mu^2\mu^{-2s_c}(1+\mu|x_0|)^{\frac12}
\int_{\lambda\approx\mu}\int_{\sigma\approx\lambda}
\lambda^{\frac14}
\left(
\lambda^{-1}\lambda^{\frac34}
+\lambda^{-\frac14}
\right)
\frac{d\sigma}{\sigma}\frac{d\lambda}{\lambda}
\\
&\lesssim_\delta\epsilon^4c_\mu^2\mu^{-2s_c}(1+\mu|x_0|)^{\frac12}.
\end{split}
\end{equation*}
This proves \eqref{ulogcommrem}.

It remains to prove the auxiliary bounds for $\mathcal C_\lambda$. By Bernstein, it suffices to prove only the $L_{t,x}^2$ bound in \eqref{ulogcommL2Linf} and \eqref{ulogcommauxtrilinear}. Since $2iT_mf=u-\overline u$, we have
\begin{equation*}
-T_{\bb}\overline u+\mathcal R_{\mathrm{comm}}
=2iT_{\bb m}f-T_{\bb}u.
\end{equation*}
For these two terms on the right we use the first estimate in
\eqref{ulogcommcoefficients}, and the fact that $m$ is bounded. This gives immediately the $L_{t,x}^2$ bound in \eqref{ulogcommL2Linf}. In \eqref{ulogcommauxtrilinear}, it provides instead two unbalanced pairs; one with $\bb$ and $\overline{u}_{\mu}^{x_0}$, and the other with $u_{\mu}$ and the remaining factor of $u$ or $f$. For the remaining part of $\mathcal{C}_{\lambda}$, we write
\begin{equation*}
\begin{split}
P_\lambda \left([T_{V^l}\cdot\nabla,T_{\mfq}]
-T_{\{V\cdot\xi,\mfq\}}\right)\varphi
=\int_{\sigma\approx\lambda}
P_\lambda \left([T_{V^l}\cdot\nabla,T_{\mfq}]
-T_{\{V\cdot\xi,\mfq\}}\right)\varphi_\sigma
\frac{d\sigma}{\sigma}.
\end{split}
\end{equation*}
By standard pseudodifferential calculus, the coefficients in the
integrand contain either a spatial derivative of $V_{[\sigma]}$
or a difference $\sigma(V_{[\tau]}-V_{[\sigma]})$ with
$\tau\approx\sigma$. For the latter we apply the fundamental
theorem of calculus in the truncation parameter using the linearization formula for $V$, as in
\eqref{Fauxcoefficientbounds} and the estimates following it.
Along with Proposition \ref{VBunbil}, this gives the same bilinear bound as above, but with $\nabla V_{[\sigma]}$ or $\sigma(V_{[\tau]}-V_{[\sigma]})$ used in place of $\bb$.

Pairing these coefficients with $u_\sigma$ or $f_\sigma$ in the
terms above and integrating over $\sigma\approx\lambda$ gives
the $L_{t,x}^2$ bound in \eqref{ulogcommL2Linf}. For
\eqref{ulogcommauxtrilinear}, we instead pair the coefficient
with $u_\mu$ and pair $|D|^{\frac{1}{2}}\varphi_{\sigma}$ with
$\overline u_\mu^{x_0}$ using the unbalanced bilinear estimate.
The two pairings give factors of $\sigma^{\frac14}$ and
$\sigma^{-\frac54}$, respectively, which gives the required
$\lambda^{-1}$ after integration over $\sigma\approx\lambda$.
\end{proof}

Before deriving the equation for $u$, we also replace $T_G\varphi$
in the $\eta$ equation by its principal part $T_{\mathfrak g}\varphi$.
The difference $(T_G-T_{\mathfrak g})\varphi$ is not perturbative.
Hueristically, it has the form
$L_{\mathrm{lh}}(D^2\eta,|D|^{-\frac12}u)$, which suggests that
it can be eliminated by a single normal form correction. To make
this correction precise, we need a more explicit expression than $T_G$.
The following lemma will allow us to extract the "main" part of this subprincipal error.

\begin{lemma}\label{uGsubprincipallemma}
Under the bootstrap hypotheses in Section \ref{sec:bootstrap},
we have
\begin{equation}\label{uGsubprincipalexp}
(T_G-T_{\mathfrak g})\varphi
=-\int_0^\infty
\int_0^{\delta\lambda}
P_{<M\lambda}[T_{\mathfrak g_{<\alpha}},\eta_\alpha]^{2+}
T_{\mathfrak l_{<\alpha}}\varphi_\lambda
\frac{d\alpha}{\alpha}\frac{d\lambda}{\lambda}
+\mathcal R_G.
\end{equation}
Here $\varphi_\lambda=S_\lambda\varphi$ and we define
\begin{equation*}
[T_{\mathfrak g_{<\alpha}},\eta_\alpha]^{2+}
\coloneq
[T_{\mathfrak g_{<\alpha}},\eta_\alpha]
-\frac1iOp\bigl(
\partial_\xi p_\alpha\cdot\nabla\eta_\alpha\bigr),
\end{equation*}
where $p_{\alpha}(x,\xi)$ denotes the symbol of $T_{\mathfrak{g}_{<\alpha}}$. That is, $[T_{\mathfrak g_{<\alpha}},\eta_\alpha]^{2+}$ is what remains after subtracting the principal Poisson bracket term from the commutator. Moreover, for every $\mu>0$ and $x_0\in\mathbb R^2$, the remainder satisfies
\begin{equation}\label{uGsubprincipalerror}
\|(\mathcal R_G)_\mu\overline u_\mu^{x_0}\|_{L_{t,x}^1}
\lesssim_\delta
\epsilon^4c_\mu^2\mu^{-2s_c}(1+\mu|x_0|)^{\frac12}.
\end{equation}
\end{lemma}

\begin{proof}
Fix $\lambda>0$. Applying the linearization formula \eqref{Glin} to
$G_{<\alpha}\varphi_\lambda$ and using the identity
$G_{<\alpha}L_{<\alpha}\varphi_\lambda
+\nabla\cdot M_{<\alpha}\varphi_\lambda=0$, we obtain
\begin{equation*}
\alpha\frac{d}{d\alpha}G_{<\alpha}\varphi_\lambda
=-[G_{<\alpha},\eta_\alpha]L_{<\alpha}\varphi_\lambda
-\nabla\eta_\alpha\cdot M_{<\alpha}\varphi_\lambda.
\end{equation*}
Integrating this identity from $0$ to $\delta\lambda$ gives
\begin{equation}\label{Gsubpint}
\begin{split}
G_{<\delta\lambda}\varphi_\lambda
=|D|\varphi_\lambda
-\int_0^{\delta\lambda}
\Bigl(
[G_{<\alpha},\eta_\alpha]L_{<\alpha}\varphi_\lambda
+\nabla\eta_\alpha\cdot M_{<\alpha}\varphi_\lambda
\Bigr)\frac{d\alpha}{\alpha}.
\end{split}
\end{equation}
Next, we compute the contribution of the principal parts of $G$, $L$,
and $M$ in the integrand above. With the minus sign included, this is given by
\begin{equation*}
-[T_{\mathfrak g_{<\alpha}},\eta_\alpha]
T_{\mathfrak l_{<\alpha}}\varphi_\lambda
-\nabla\eta_\alpha\cdot T_{\mathfrak m_{<\alpha}}\varphi_\lambda.
\end{equation*}
Writing $p_\alpha$ for the full symbol of $T_{\mathfrak g_{<\alpha}}$,
the principal term of $[T_{\mathfrak g_{<\alpha}},\eta_\alpha]$
has symbol $\frac1i\partial_\xi p_\alpha\cdot\nabla\eta_\alpha$.
Thus, by the definition of the $2+$ commutator, we have
\begin{equation*}
\begin{split}
&-[T_{\mathfrak g_{<\alpha}},\eta_\alpha]
T_{\mathfrak l_{<\alpha}}\varphi_\lambda
-\nabla\eta_\alpha\cdot T_{\mathfrak m_{<\alpha}}\varphi_\lambda
\\
&=-[T_{\mathfrak g_{<\alpha}},\eta_\alpha]^{2+}
T_{\mathfrak l_{<\alpha}}\varphi_\lambda
+iOp\bigl(\nabla\eta_\alpha\cdot\partial_\xi p_\alpha\bigr)
T_{\mathfrak l_{<\alpha}}\varphi_\lambda
-\nabla\eta_\alpha\cdot T_{\mathfrak m_{<\alpha}}\varphi_\lambda.
\end{split}
\end{equation*}
We have to group the last two terms together. We first replace the composition in the first of these terms by the operator $iOp((\nabla\eta_\alpha\cdot\partial_\xi p_\alpha)p_\alpha^{\mathfrak l})$. Here, as above we write $p_\alpha$ for the symbol of
$T_{\mathfrak g_{<\alpha}}$. We also write
$p_\alpha^{\mathfrak l}$ and $p_\alpha^{\mathfrak m}$ for the
symbols of $T_{\mathfrak l_{<\alpha}}$ and $T_{\mathfrak m_{<\alpha}}$. The remaining symbol is then
\begin{equation*}
\begin{split}
\nabla\eta_\alpha\cdot
\left(i p_\alpha^{\mathfrak l}\partial_\xi p_\alpha
-p_\alpha^{\mathfrak m}\right)
&=\nabla\eta_\alpha\cdot
\frac{|\xi|^2\nabla\eta_{<\alpha}
-(\nabla\eta_{<\alpha}\cdot\xi)\xi}
{\mathfrak g_{<\alpha}}
\\
&+\nabla\eta_\alpha\cdot\Bigl(
i(p_\alpha^{\mathfrak l}-\mathfrak l_{<\alpha})\partial_\xi p_\alpha
+i\mathfrak l_{<\alpha}\partial_\xi(p_\alpha-\mathfrak g_{<\alpha})
-(p_\alpha^{\mathfrak m}-\mathfrak m_{<\alpha})
\Bigr).
\end{split}
\end{equation*}
The first line on the right is precisely
$\alpha\frac{d}{d\alpha}\mathfrak g_{<\alpha}$.
The second line consists purely of high-frequency tail errors: more precisely, for $\mathfrak{a}=\mathfrak g,\mathfrak l,\mathfrak m$, we have
\begin{equation*}
p_\alpha^a(x,\xi)-a_{<\alpha}(x,\xi)
=-\int_0^\infty
\bigl(P_{\geq M\delta\sigma}a_{<\alpha}\bigr)(x,\xi)
S_\sigma(\xi)\frac{d\sigma}{\sigma}.
\end{equation*}
These contributions will be included in $\mathcal R_G$ below,
together with the composition error mentioned above.
Since $\mathfrak g(0,\xi)=|\xi|$, by combining the above identity with \eqref{Gsubpint} we then obtain
\begin{equation*}
(T_G-T_{\mathfrak g})\varphi
=-\int_0^\infty
\int_0^{\delta\lambda}
P_{<M\lambda}\left([T_{\mathfrak g_{<\alpha}},\eta_\alpha]^{2+}
T_{\mathfrak l_{<\alpha}}\varphi_\lambda
\right)\frac{d\alpha}{\alpha}\frac{d\lambda}{\lambda}
+\mathcal R_G.
\end{equation*}
Here $\mathcal R_G$ consists of the errors from replacing $G,L,M$
by their principal parts and the errors identified above. Explicitly, we have
\begin{equation}\label{RGdef}
\begin{split}
\mathcal R_G
&=-\int_0^\infty\int_0^{\delta\lambda}
P_{<M\lambda}\Bigl(
[G_{<\alpha}-T_{\mathfrak g_{<\alpha}},\eta_\alpha]L_{<\alpha}
+[T_{\mathfrak g_{<\alpha}},\eta_\alpha]
(L_{<\alpha}-T_{\mathfrak l_{<\alpha}})\Bigr)\varphi_{\lambda}\frac{d\alpha}{\alpha}\frac{d\lambda}{\lambda}
\\
&-\int_0^\infty \int_0^{\delta\lambda}
P_{<M\lambda}\bigl(\nabla\eta_\alpha\cdot(M_{<\alpha}-T_{\mathfrak m_{<\alpha}})\varphi_\lambda\bigr)
\frac{d\alpha}{\alpha}\frac{d\lambda}{\lambda}
\\
&+\int_0^\infty \int_0^{\delta\lambda}
P_{<M\lambda}\Bigl(
iOp\bigl(\nabla\eta_\alpha\cdot\partial_\xi p_\alpha\bigr)
T_{\mathfrak l_{<\alpha}}
-\nabla\eta_\alpha\cdot T_{\mathfrak m_{<\alpha}}
-Op\Bigl(\alpha\frac{d}{d\alpha}\mathfrak g_{<\alpha}\Bigr)
\Bigr)\varphi_\lambda
\frac{d\alpha}{\alpha}\frac{d\lambda}{\lambda}
\\
&+\int_0^\infty
\Bigl(
P_{<M\lambda}Op(\mathfrak g_{<\delta\lambda})
-Op(\mathfrak g^l_{<\delta\lambda})
\Bigr)\varphi_\lambda\frac{d\lambda}{\lambda}.
\end{split}
\end{equation}

We now estimate $\mathcal R_G$. The errors from replacing $G$, $L$,
and $M$ by their principal parts have integrand given by (omitting the sign and $P_{<M\lambda}$ projection)
\begin{equation}\label{operatorreplacements}
\begin{split}
[G_{<\alpha}-T_{\mathfrak g_{<\alpha}},\eta_\alpha]
L_{<\alpha}\varphi_\lambda
+[T_{\mathfrak g_{<\alpha}},\eta_\alpha]
(L_{<\alpha}-T_{\mathfrak l_{<\alpha}})\varphi_\lambda+\nabla\eta_\alpha\cdot
(M_{<\alpha}-T_{\mathfrak m_{<\alpha}})\varphi_\lambda.
\end{split}
\end{equation}
We show the details for the first term since it's the most difficult. The other terms are handled similarly. We first consider the contribution zeroth order part of $G_{<\alpha}-Op(\mathfrak g_{<\alpha})$, as given Proposition \ref{paraexp}. We let $v$ be the part of
$L_{<\alpha}\varphi_\lambda$ at frequency $\approx\lambda$.
Since $\alpha<\delta\lambda$, the product $\eta_{\alpha}v$ is at frequency $\approx\lambda$. Writing $G^0_{<\alpha}$ for this zeroth order part and $K^0_{\alpha,\lambda}$ for its kernel (restricted to the $|\xi|\approx\lambda$ annulus), we have
\begin{equation*}
[G^0_{<\alpha},\eta_\alpha]v(x)
=\int_{\mathbb R^2}K^0_{\alpha,\lambda}(x,y)
\bigl(\eta_\alpha(y)-\eta_\alpha(x)\bigr)v(y)dy.
\end{equation*}
By the fundamental theorem of calculus and the kernel bound \eqref{K0bound},
the kernel in this integral satisfies
\begin{equation*}
\begin{split}
&\left|K^0_{\alpha,\lambda}(x,y)
\bigl(\eta_\alpha(y)-\eta_\alpha(x)\bigr)\right|
\\
&\lesssim_N
\lambda^{-1}\lambda^2\langle\lambda(x-y)\rangle^{-N}
\int_0^1|\nabla\eta_\alpha(x+t(y-x))|\,dt
\int_{\mathbb R^2}\langle\theta\rangle^{-N}
|\nabla^2\eta_{<\alpha}(y+\lambda^{-1}\theta)|d\theta.
\end{split}
\end{equation*}
Here $\lambda^{-1}$ comes from the factor $|x-y|$ in the
difference of $\eta_\alpha$ and cancels the derivative from
$L_{<\alpha}\varphi_\lambda$. For 
$\mu\approx\lambda$, this bound allows us to use two unbalanced bilinear estimates. Indeed, we can pairr $\nabla\eta_\alpha$ with
$\varphi_\lambda$ and $\nabla^2\eta_{<\alpha}$ with
$\overline u_\mu^{x_0}$. This gives the bound
\begin{equation*}
\begin{split}
\|L(\nabla\eta_\alpha,\varphi_\lambda)\|_{L_{t,x}^2}
\|L(\nabla^2\eta_{<\alpha},\overline u_\mu^{x_0})\|_{L_{t,x}^2}
&\lesssim
\epsilon^4c_\mu^2\mu^{-2s_c}(1+\mu|x_0|)^{\frac12}
\lambda^{-\frac12}\alpha^{-\frac14}\alpha^{\frac34}
\\
&=
\epsilon^4c_\mu^2\mu^{-2s_c}(1+\mu|x_0|)^{\frac12}
\lambda^{-\frac{1}{2}}\alpha^{\frac{1}{2}}
\end{split}
\end{equation*}
which more than suffices upon integration over $\alpha<\delta\lambda$ and $\lambda\approx\mu$. The contribution with $\mu\ll\lambda$ is a low-output
frequency elliptic tail and is estimated using \eqref{bernelliptic2}, together with the higher order errors just below.

For the higher order error, write $\mathcal E^G_{\alpha,\lambda}$ 
for the remainder term in Proposition \ref{paraexp}. Here, we can expand the commmutator and estimate 
$\mathcal E^G_{\alpha,\lambda}(\eta_\alpha v)$ and
$\eta_\alpha\mathcal E^G_{\alpha,\lambda}v$ separately. Taking $L=4$ in this bound and applying H\"older and Bernstein gives
\begin{equation*}
\begin{split}
\|\mathcal E^G_{\alpha,\lambda}(\eta_\alpha v)\|_{L_t^1L_x^2}
+\|\eta_\alpha\mathcal E^G_{\alpha,\lambda}v\|_{L_t^1L_x^2}&\lesssim_\delta
\lambda^{-3}
\|\nabla^2\eta_{<\alpha}\|_{L_t^4\dot B^0_{\infty,1}}^4
\|\eta_\alpha\|_{L_{t,x}^\infty}
\lambda\|\varphi_\lambda\|_{L_t^\infty L_x^2}
\\
&\lesssim_\delta
\lambda^{-3}(\epsilon\alpha^{\frac78})^4
(\epsilon\alpha^{-1})
\lambda(\epsilon c_\lambda\lambda^{-s_c-\frac12})
\\
&=\epsilon^6c_\lambda\lambda^{-s_c}
\left(\frac{\alpha}{\lambda}\right)^{\frac52},
\end{split}
\end{equation*}
which again suffices. Finally, the part of
$L_{<\alpha}\varphi_\lambda$ given by $L_{<\alpha}\varphi_\lambda-v$ localized away from $|\xi|\approx\lambda$ is
handled by Proposition \ref{Bernsteinelliptic} after expanding the commutator as above. Combining everything above finally shows that the contribution of $[G_{<\alpha}-T_{\mathfrak{g}_{<\alpha}},\eta_\alpha]L_{<\alpha}\varphi_{\lambda}$ is perturbative. The other terms in \eqref{operatorreplacements} are handled similarly. Finally, the last two lines in \eqref{RGdef} come from subprincipal errors of the form $Op(ab)-Op(a)Op(b)$ or high-frequency tails. The high frequency tail erorrs are easily dispensed with as above using Proposition \ref{Bernsteinelliptic}. On the other hand, the composition error differentiates $p_\alpha^{\mathfrak l}$ in $x$, so we can instead use the same two bilinear pairings with $\nabla\eta_\alpha$ and $\nabla^2\eta_{<\alpha}$ as above. This finally gives
\eqref{uGsubprincipalerror} and completes the proof.
\end{proof}

We can now combine Lemmas \ref{ulogcommutator} and \ref{uGsubprincipallemma} with the expansions
for $N_\eta$ and $N_\varphi$ in the previous sections to obtain the equation for $u$.

\begin{theorem}\label{ufinaltheorem}
Under the bootstrap hypotheses in section \ref{sec:bootstrap}, the complex variable $u\coloneq \eta+iT_{\mfq}\varphi$ satisfies
\begin{equation}\label{u-para}
\left(
\partial_t+T_{V^l}\cdot\nabla
+\frac12T_{\nabla\cdot V^l}
-T_{\bb}
+iT_{\mathfrak\omega}
\right)u
=N(u),
\end{equation}
where the nonlinear term $N(u)$ is defined by
\begin{equation}\label{unonlinearitydef}
N(u)\coloneq
-\left(\frac12T_{\nabla\cdot V^l}+T_{\bb}\right)\overline u
+N_u^{[1]}+\sum_{k=2}^{4}N_u^{[k]}
+iT_{\mfq}N_{\varphi}^{\mathrm{flux}}
+\mathcal R_u.
\end{equation}
The subprincipal contribution is given by
\begin{equation}\label{u-para-err}
\begin{split}
N_u^{[1]}&=(T_{\mathfrak g}-T_{\mathfrak\omega}T_{\mfq})\varphi
-i(T_{\mfq}T_{a^l}-T_{\mathfrak\omega})\eta
\\
&+i\left([T_{V^l}\cdot\nabla,T_{\mfq}]
-T_{\{V\cdot\xi,\mfq\}}\right)\varphi
\\
&-\int_0^\infty
\int_0^{\delta\lambda}
P_{<M\lambda}[T_{\mathfrak g_{<\alpha}},\eta_\alpha]^{2+}
T_{\mathfrak l_{<\alpha}}\varphi_\lambda
\frac{d\alpha}{\alpha}\frac{d\lambda}{\lambda}.
\end{split}
\end{equation}
For the balanced terms, we define
\begin{equation}\label{u-multilinear-def}
N^{[k]}_u\coloneq N_{\eta}^{[k]}+iT_{\mfq}N_{\varphi}^{[k]},
\qquad 2\leq k\leq4,
\end{equation}
where the terms on the right are defined in Corollaries \ref{Netafinalcorollary}
and \ref{Nphifinalcorollary}, and $N_{\varphi}^{\mathrm{flux}}$ is as in
\eqref{varphifinalequation}.

We also have the following error bounds.
\begin{enumerate}
\item ($L_{t,x}^1$ bound for the remainder). For every $\mu>0$ and
$x_0\in\mathbb R^2$, there holds
\begin{equation}\label{u-remainder-def}
\|(\mathcal R_u)_\mu\overline u_\mu^{x_0}\|_{L_{t,x}^1}
\lesssim_\delta
\epsilon^4c_\mu^2\mu^{-2s_c}(1+\mu|x_0|)^{\frac12}.
\end{equation}
\item ($L_t^2L_x^{\infty}$ bound for $N(u)$). For every $\lambda>0$, we have
\begin{equation}\label{NuL2Linf}
\|N_{\lambda}(u)\|_{L_t^2L_x^\infty}
\lesssim \epsilon^2\lambda^{-\frac34}.
\end{equation}
\item (Multilinear decomposition for $N(u)$). For every $0<\lambda\ll\mu$
and $x_0\in\mathbb R^2$, there is a decomposition
\begin{equation}\label{Nuauxdecomp}
N_{\lambda}(u)\,\overline u_\mu^{x_0}u_\mu
=\mathcal R_u^{\mathrm{aux}}
+\nabla\cdot Q_u^{\mathrm{flux}},
\end{equation}
where, for every $\frac43<p\leq2$,
\begin{equation}\label{Nuauxtrilinear}
\|\mathcal R_u^{\mathrm{aux}}\|_{L_{t,x}^1}
+\lambda^{\frac2p-1}
\|Q_u^{\mathrm{flux}}\|_{L_t^1L_x^p}
\lesssim_p
\epsilon^4c_\mu^2\lambda^{-1}
\mu^{-2s_c}(1+\mu|x_0|)^{\frac12}.
\end{equation}
\item ($L_{t,x}^2$ bound after subtracting the balanced quadratic part).
For every $\lambda>0$, we have 
\begin{equation}\label{NuL2remainder}
\|(N(u)-N_u^{[2]})_\lambda\|_{L_{t,x}^2}
\lesssim_\delta
\epsilon^2c_\lambda\lambda^{\frac14-s_c}.
\end{equation}
\end{enumerate}
\end{theorem}
It certain cases in the remaining sections, it will be convenient to use the slightly strengthened version of \eqref{NuL2remainder} where we further isolate the so-called balanced part $N_{u}^{[2],\mathrm{bal}}$ of $N_u^{[2]}$ which consists of the part of $N_u^{[2]}$ where the two high input frequency factors in $N_u^{[2]}$ are at frequency comparable to the output. More precisely, we have the following.
\begin{corollary}\label{N2balremainder}
Under the hypotheses of Theorem \ref{ufinaltheorem}, let
$N_u^{[2],\mathrm{bal}}$ denote the part of $N_u^{[2]}$
where both input frequencies are comparable to the output
frequency. Then \eqref{NuL2remainder} can be strengthened to
\begin{equation}\label{NuL2balancedremainder}
\|(N(u)-N_u^{[2],\mathrm{bal}})_\lambda\|_{L_{t,x}^2}
\lesssim_\delta
\epsilon^2c_\lambda\lambda^{\frac14-s_c},
\qquad \lambda>0.
\end{equation}
\end{corollary}
\begin{proof}
This follows immediately from \eqref{NuL2remainder},
the definitions of $N_\eta^{[2]}$ and $N_\varphi^{[2]}$
in \eqref{Neta2def} and \eqref{Nphi2def}, and
Lemma \ref{balancedrestrictedoutput}.
\end{proof}
Now, we turn to the proof of the Theorem \ref{ufinaltheorem}, beginning with the derivation of the equation \eqref{u-para}.
\begin{proof}[Derivation of \eqref{u-para}]
Applying $\partial_t+T_{V^l}\cdot\nabla$ to the definition of $u$
and using \eqref{Z-para-re}, we obtain
\begin{equation*}
\begin{split}
\left(\partial_t+T_{V^l}\cdot\nabla\right)u
&=T_G\varphi-iT_{\mfq}T_{a^l}\eta
-T_{\nabla\cdot V^l}\eta
+i\left[\partial_t+T_{V^l}\cdot\nabla,T_{\mfq}\right]\varphi
\\
&+N_{\eta}+iT_{\mfq}N_{\varphi}
+iT_{\mfq}N_{\varphi}^{\mathrm{flux}}
+iT_{\mfq}\mathcal R.
\end{split}
\end{equation*}
For the commutator term, we apply Lemma \ref{ulogcommutator}.
Since $u+\overline u=2\eta$, the divergence term satisfies
\begin{equation*}
\frac12T_{\nabla\cdot V^l}u-T_{\nabla\cdot V^l}\eta
=-\frac12T_{\nabla\cdot V^l}\overline u.
\end{equation*}
Thus, after moving $T_{\bb}u$ to the left, the combined contribution
of the divergence term and the commutator on the right is
\begin{equation}\label{u-para-main}
-\left(\frac12T_{\nabla\cdot V^l}+T_{\bb}\right)\overline u
+i\left([T_{V^l}\cdot\nabla,T_{\mfq}]
-T_{\{V\cdot\xi,\mfq\}}\right)\varphi
+\mathcal R_{\mathrm{comm}}.
\end{equation}
On the other hand, adding $iT_{\mathfrak\omega}u$ to the left
and using the definition of $u$, we write
\begin{equation*}
\begin{split}
T_G\varphi-iT_{\mfq}T_{a^l}\eta+iT_{\mathfrak\omega}u=(T_G-T_{\mathfrak g})\varphi
+(T_{\mathfrak g}-T_{\mathfrak\omega}T_{\mfq})\varphi-i(T_{\mfq}T_{a^l}-T_{\mathfrak\omega})\eta.
\end{split}
\end{equation*}
Substituting \eqref{uGsubprincipalexp} for $(T_G-T_{\mathfrak g})\varphi$ gives \eqref{u-para-err}. Substituting
\eqref{Netafinalexpansion} and \eqref{Nphifinalexpansion} gives
\eqref{u-multilinear-def}. Finally, setting
\begin{equation*}
\mathcal R_u
=\mathcal R_\eta+iT_{\mfq}(\mathcal R_\varphi+\mathcal R)
+\mathcal R_{\mathrm{comm}}+\mathcal R_G
\end{equation*}
gives \eqref{u-para} and \eqref{unonlinearitydef}. The first two
terms are the remainders from the equations for $\eta$ and $\varphi$,
with the factor $\lambda^{1/2}$ from $T_{\mfq}$ in the latter.
Their bounds, together with \eqref{ulogcommrem} and \eqref{uGsubprincipalerror}, give
\eqref{u-remainder-def}.
\end{proof}
\begin{proof}[Proof of the auxiliary bounds]
Now we prove the remaining three auxiliary bounds \eqref{NuL2Linf}, \eqref{Nuauxtrilinear} and \eqref{NuL2remainder}, For these estimates it is convenient to expand the $\varphi$ nonlinearity using \eqref{Nphierrinitialdef},
$N_\varphi^{\mathrm{err}}=N_\varphi+N_\varphi^{\mathrm{flux}}+\mathcal R$. Instead of expanding $N_\eta$ and $N_\varphi$ in terms of their layered multi-linear components, we can then write $N(u)$ as
\begin{equation*}
\begin{split}
N(u)=N_\eta+iT_{\mfq}N_\varphi^{\mathrm{err}}
+\mathcal C-\frac12T_{\nabla\cdot V^l}\overline u+(T_G-T_{\mathfrak g})\varphi
+(T_{\mathfrak g}-T_{\mathfrak\omega}T_{\mfq})\varphi
-i(T_{\mfq}T_{a^l}-T_{\mathfrak\omega})\eta.
\end{split}
\end{equation*}
where $\mathcal C$ is the commutator error defined above. From the bounds \eqref{ulogcommL2Linf} and
\eqref{ulogcommauxtrilinear}, we can immediatly estimate the contribution of $\mathcal{C}$ in \eqref{NuL2Linf}, \eqref{Nuauxtrilinear} and \eqref{NuL2remainder}. The term
$T_{\nabla\cdot V^l}\overline u$ contributes the same estimates by using two unbalanced bilinear estimates and Proposition \ref{VBunbil} for its contribution to \eqref{Nuauxdecomp} or just Proposiiton \ref{VBunbil} and one bilinear unbalanced estimate pairing the $\nabla\cdot V$ coefficient with $\overline{u}$ for the $L_{t,x}^2$ bound, alongside Bernstein for the $L_t^2L_x^{\infty}$ bound. The contribution of $(T_G-T_{\mathfrak{g}})\varphi$ is handled similarly by using the kernel and error estimates in Proposition \ref{paraexp}. For the two composition errors, we use
$\mathfrak\omega\mfq=\mathfrak g$ and
$\mfq a=\mathfrak\omega$ and argue exactly as in \eqref{paracomexp} using the coefficient bounds
in \eqref{ulogcommcoefficients} to obtain again the same bounds.

It remains to assemble the contributions of $N_\eta$ and
$iT_{\mfq}N_\varphi^{\mathrm{err}}$. At output frequency
$\lambda$, the kernel of $P_{\lambda}T_{\mfq}$ has $L^1$ norm
$\lesssim_\delta\lambda^{1/2}$. Thus \eqref{NetaL2Linf},
\eqref{NphiL2Linf}, and \eqref{NphiinitialL2Linf} give
\eqref{NuL2Linf}. After subtracting
$N_u^{[2]}=N_\eta^{[2]}+iT_{\mfq}N_\varphi^{[2]}$, we can use
\eqref{NetaL2cubic}, \eqref{NphiL2cubic}, and
\eqref{NphiinitialL2} instead to obtain \eqref{NuL2remainder}.

Finally, for the flux-error decomposition, the contribution of $N_\eta$
is handled using \eqref{Netaauxdecomp}-\eqref{Netaauxtrilinear} while the contribution of $N_{\varphi}$ is handled using \eqref{Nphiauxtrilinear}. 
\end{proof}
\section{The equation for the linearized complex variable \texorpdfstring{$w$}{}}
\label{subsec:wfinal}

We now combine the expansions for $h$ and $r$ in Theorems
\ref{Netalinearized} and \ref{Nphilinearized} to obtain the
corresponding layered paradifferential equation for the linearized variable $w=h+iT_{\mfq}r$. The symbols
$\mfq$, $\mathfrak\omega$, and $\bb$ are the same as in the
equation for $u$.
\begin{theorem}\label{wfinaltheorem}
Let $0<T\leq1$. On $[0,T]$, assume the bootstrap bounds for the
background solution $u$ in Section \ref{sec:bootstrap} and for the
linearized variable $w$ in Subsection \ref{subsec:para-boot}.
Assume in addition that the frequency envelope for $u_0$ satisfies the additional regularity hypothesis $c_\lambda\lesssim\lambda^{-\frac1{12}}$ for every
$\lambda>0$. Then $w$ satisfies the equation
\begin{equation}\label{wfinalequation}
\left(
\partial_t+T_{V^l}\cdot\nabla
+\frac12T_{\nabla\cdot V^l}
-T_{\bb}
+iT_{\mathfrak\omega}
\right)w
=N^{\mathrm{lin}}(w),
\end{equation}
where
\begin{equation}\label{wnonlinearitydef}
N^{\mathrm{lin}}(w)
=-\left(\frac12T_{\nabla\cdot V^l}+T_{\bb}\right)\overline w
+\sum_{k=0}^{2}N_w^{\mathrm{lin},[k]}
+\mathcal R_w^{\mathrm{lin}}.
\end{equation}
With the $k=0$ terms defined in \eqref{Netalin0def} and
\eqref{Nphilin0def} and the $k=2$ terms defined in
\eqref{Netalin2def} and \eqref{Nphilin2def}, we set
\begin{equation*}
N_w^{\mathrm{lin},[k]}
\coloneq N_\eta^{\mathrm{lin},[k]}
+iT_{\mfq}N_\varphi^{\mathrm{lin},[k]},
\hspace{5mm} k=0,2.
\end{equation*}
On the other hand, the subprincipal low-high contribution is defined by
\begin{equation*}
\begin{split}
N_w^{\mathrm{lin},[1]}&=
-\int_0^\infty P_{<M\lambda}\int_0^{\delta\lambda}
[T_{\mathfrak g_{<\alpha}},\eta_\alpha]^{2+}
T_{\mathfrak l_{<\alpha}}r_\lambda
\frac{d\alpha}{\alpha}\frac{d\lambda}{\lambda}
\\
&+(T_{\mathfrak g}-T_{\mathfrak\omega}T_{\mfq})r
-i(T_{\mfq}T_{a^l}-T_{\mathfrak\omega})h
\\
&+i\left([T_{V^l}\cdot\nabla,T_{\mfq}]
-T_{\{V\cdot\xi,\mfq\}}\right)r.
\end{split}
\end{equation*}
Here $r_\lambda=S_\lambda r$ and $[T_{\mathfrak{g}_{<\alpha}},\eta_{\alpha}]^{2+}$ is defined as in Lemma \ref{uGsubprincipallemma}.
We also have the following bounds:
\begin{enumerate}
\item For every $\mu>0$ and $x_0\in\mathbb R^2$, there holds
\begin{equation}\label{wlinremainder}
\|(\mathcal R_w^{\mathrm{lin}})_\mu
\overline w_\mu^{x_0}\|_{L_{t,x}^1}
\lesssim_\delta
\epsilon^2d_\mu^2(1+\mu|x_0|)^{\frac12}.
\end{equation}
\item For every $\nu>0$, there holds 
\begin{equation}\label{wlinL2remainder}
\|(N^{\mathrm{lin}}(w)-N_w^{\mathrm{lin},[2]})_\nu
\|_{L_{t,x}^2}
\lesssim_\delta
\epsilon d_\nu\nu^{\frac14}.
\end{equation}
\item For every $0<\nu\ll\mu$ and $x_0\in\mathbb R^2$, there holds
\begin{equation}\label{wlinauxtrilinear}
\|(N^{\mathrm{lin}}(w))_\nu
\overline w_\mu^{x_0}u_\mu\|_{L_{t,x}^1}
\lesssim_\delta
\epsilon^2d_\nu d_\mu\nu^{\frac14}\mu^{-\frac54}
(1+\mu|x_0|)^{\frac12}.
\end{equation}
\item For every $0<\lambda\ll\mu$ and $x_0\in\mathbb R^2$,
the nonlinear term $N(u)$ in \eqref{unonlinearitydef} admits a decomposition
\begin{equation}\label{Nuwwauxdecomp}
N_\lambda(u)w_\mu\overline w_\mu^{x_0}
=\mathcal R_{u,w}^{\mathrm{aux}}
+\nabla\cdot Q_{u,w}^{\mathrm{flux}},
\end{equation}
where, for $\frac43<p\leq2$, we have
\begin{equation}\label{Nuwwauxtrilinear}
\|\mathcal R_{u,w}^{\mathrm{aux}}\|_{L_{t,x}^1}
+\lambda^{\frac2p-1}\|Q_{u,w}^{\mathrm{flux}}\|_{L_t^1L_x^p}
\lesssim_{\delta,p}
\epsilon^2d_\mu^2\lambda^{-1}(1+\mu|x_0|)^{\frac12}.
\end{equation}
\end{enumerate}
\end{theorem}
As with Corollary \ref{NuL2balancedremainder}, we also have its linearized analogue below.
\begin{corollary}\label{Nbalremainderlinearized}
Under the hypotheses of Theorem \ref{wfinaltheorem}, let
$N_w^{\mathrm{lin},[2],\mathrm{bal}}$ denote the part of
$N_w^{\mathrm{lin},[2]}$ where both input frequencies are
comparable to the output frequency. Then
\eqref{wlinL2remainder} can be strengthened to
\begin{equation}\label{wlinL2balancedremainder}
\|(N^{\mathrm{lin}}(w)
-N_w^{\mathrm{lin},[2],\mathrm{bal}})_\nu\|_{L_{t,x}^2}
\lesssim_\delta
\epsilon d_\nu\nu^{\frac14},
\qquad \nu>0.
\end{equation}
\end{corollary}
\begin{proof}
We use the definitions \eqref{Netalin2def} and
\eqref{Nphilin2def} and argue as in Lemma
\ref{balancedrestrictedoutput}, using the mixed bilinear
bound \eqref{wuab-bi-bal-boot}. For the $\varphi$ contribution,
we use the symbol identity
$\mathfrak g\mathfrak l=-i\xi\cdot\mathfrak m$,
which makes its quadratic symbol defining $N^{\mathrm{lin},[2]}$ vanish at opposite input
frequencies. We can then and argue as in the proof of \eqref{nullL2bound} but using the linearized bootstrap bounds instead of the nonlinear ones.
Summing over the input frequencies above $\nu$ and combining
with \eqref{wlinL2remainder} proves the bound.
\end{proof}
Now, we turn to the proof of Theorem \ref{wfinaltheorem}.
\begin{proof}
The derivation of \eqref{wfinalequation} and the proof of
\eqref{wlinremainder} follow the proof of Theorem
\ref{ufinaltheorem}, using \eqref{etalinearizedequation} and
\eqref{philinearizedequation} in place of \eqref{Z-para-re}
and the linearized bound \eqref{Dtabilinearlin} in the proof
of Lemma \ref{ulogcommutator}. For the expansion of $(T_G-T_{\mathfrak g})r$, one argues as in Lemma
\ref{uGsubprincipallemma}, using \eqref{hboot} and the corresponding linearized bilinear
estimates in place of the nonlinear bounds.

For \eqref{wlinL2remainder}, we use \eqref{NetalinL2remainder}
and \eqref{NphilinL2remainder}, noting the additional factor of $\nu^{\frac12}$
coming from $P_\nu T_{\mfq}$ in the latter. The remaining
paradifferential terms are estimated as in the proof of
\eqref{NuL2remainder}, using the linearized bootstrap bounds
and \eqref{Dtabilinearlin}.  Finally, for \eqref{wlinauxtrilinear}, this is handled slightly differently and involves no flux correction compared to the nonlinear case. This is thanks to the additional $\frac{1}{12}$ derivatives of regularity assumed on the background solution here. Indeed, for the contribution of $N^{\mathrm{lin}}(w)-N_w^{\mathrm{lin},[2]}$, we can apply \eqref{wlinL2remainder}
and H\"older, using
\begin{equation*}
T^{\frac16}\|w_\mu\|_{L_t^{6,\infty}L_x^4}\|u_\mu\|_{L_t^{6,\infty}L_x^4}
\lesssim\epsilon d_\mu\mu^{-\frac54}.
\end{equation*}
For $N_w^{\mathrm{lin},[2]}$, we use the four $L_t^{6,\infty}L_x^4$
bounds as in the proofs of \eqref{Netalinauxtrilinear} and
\eqref{Nphilinauxtrilinear}. In the $\varphi$ contribution we
insert $P_\nu T_{\mfq}$ before applying H\"older, so that
$\nu^{-\frac14}$ in \eqref{Nphilinauxtrilinear} becomes
$\nu^{\frac14}$. This proves \eqref{wlinauxtrilinear}.

The proof of \eqref{Nuwwauxdecomp} and \eqref{Nuwwauxtrilinear}
is the same as that of \eqref{Nuauxdecomp} and \eqref{Nuauxtrilinear}, respectively, by
using instead \eqref{Netawwauxtrilinear} and 
\eqref{Nphiwwauxtrilinear} as well as the
linearized $(u,w)$ bilinear estimates in place of the $(u,u)$ ones. 
\end{proof}

\section{Paracontrolled functions and paratransported functions and symbols}
\label{sec:para-functions}

The  paracoefficients in our layered paradifferential expansion of the water wave equations, completed in the previous section, are all functions of $\nabla \eta$, or, by extension, functions of $\nabla u$ and $\nabla \bu$. However, this class of function does not suffice for what follows, due to the fact that it is not closed with respect to paraproducts. To rectify this, in this section we introduce several classes of \emph{paracontrolled functions}. 

Another issue connected to functions of $u$ and $\bu$, and implicitly to our paracontrolled functions, is that they do not commute well with the (para)material derivative.  This is the case even with the classical Littlewood-Paley frequency localization, which has served us well so far. To address this 
difficulty, in this section we also introduce a suitable notion of \emph{ paratransported functions and symbols}.

Throughout this section we assume 
that $u$ is a complex valued function satisfying the water wave equation 
\eqref{u-para}, and also satisfying 
the bootstrap assumptions \eqref{uk-ee-boot}-\eqref{uL6L4-boot}. Our constructions will depend only on $u$,
but will also be tested against functions $w$ which will be assumed 
to satisfy in turn the bootstrap bounds \eqref{wk-ee-boot}-\eqref{wL6L4-boot}. Heuristically $w$  should be thought of as either a solution to the full equation, the paradifferential equation or the linearized equation, but the specific equation is not used, only the bootstrap bounds.

\subsection{A preliminary decomposition}

Our goal here is to consider functions $g$ of the form 
\[
g = G(\nabla u)
\]
with smooth $G$
and to discuss their paradifferential balanced/unbalanced expansion.
To begin with, the function $\nabla u$ is assumed to satisfy the uniform fixed time Besov bound 
\begin{equation}\label{small-Besov-fe}
\|P_\lambda \nabla u\|_{L^2} \lesssim C\epsilon c_\lambda \lambda^{-1}, \qquad c_\lambda \in \ell^1.
\end{equation}
It is easily seen that this class of functions is a Banach algebra. We do not assume $G(0) = 0$, so the identity is also included.

In a second step we also consider 
space time properties derived from our bootstrap assumptions \eqref{uk-ee-boot}-\eqref{uL6L4-boot}.

To produce a good decomposition for our function $g$, we recall that for $\nabla u$ we have the continuous Littlewood-Paley expansion
\[
u = \int_0^\infty u_\lambda \frac{d\lambda}{\lambda}.
\]
For $g$ above, this yields a similar expansion of the form 
\[
g = \int_0^\infty \nabla u_\lambda G'(\nabla u_{<\lambda}) \frac{d\lambda}{\lambda}:= \int_0^\infty g_\lambda \frac{d\lambda}{\lambda} .
\]
where the highest input frequency in $g_\lambda$ is comparable to $\lambda$.

As in our layered paradifferential expansions, we further decompose each $g_\lambda$ depending on the second highest frequency, into a balanced and an unbalanced component,
\[
\begin{aligned}
g_\lambda = & \ \nabla u_\lambda \int_0^\lambda \nabla u_\mu G''(\nabla u_{<\mu}) \frac{d\mu}{\mu}
\\
= & \ \nabla u_\lambda \int_0^{c\lambda} \nabla u_\mu  P_{<\lambda} G''(\nabla u_{<\mu}) \frac{d\mu}{\mu}
\\ & \ 
+ \ \nabla u_\lambda\left( \int_{c\lambda}^\lambda \nabla u_\mu G''(\nabla u_{<\mu}) \frac{d\mu}{\mu}
+ \int_0^{c\lambda} \nabla u_\mu  P_{\geq \lambda} G''(\nabla u_{<\mu}) \frac{d\mu}{\mu}\right)
\\
:= & \ g_\lambda^{ub} + g_\lambda^{bal}
\end{aligned}
\]
where $c$ is a small universal constant. 
Here the unbalanced component is localized at frequency $\lambda$, while the balanced component is localized at frequency at most $\lambda$.

To describe the  properties of these expressions we introduce appropriate reference functions, to be used akin to our frequency envelopes. Precisely, we set
\[
\mfu_\lambda(x) = \int_{\R^d} \int_{\lambda/2}^{\lambda} |\nabla u_\mu(x+y)| \frac{d\mu}{\mu} \frac{\lambda^{d} dy}{(1+\lambda |y|)^{N}}  
\]
where $N$ is a fixed universal large constant. 

We can describe the size of $g_\lambda^{ub}$ and $g_\lambda^{bal}$ in terms of these reference functions as follows:

\begin{lemma}\label{l:gub+bal}
Assuming \eqref{small-Besov-fe}, we have
\begin{equation}
 |g_\lambda^{ub}| + |g_\lambda^{bal}| \lesssim    \mfu_\lambda, \qquad 
 |\nabla g_\lambda^{ub}| + |\nabla g_\lambda^{bal}| \lesssim    \lambda \mfu_\lambda
\end{equation}
as well as 
the balanced bound
\begin{equation}
|g_\lambda^{bal}| \lesssim  \int_{c\lambda}^{\lambda} \mfu_\lambda \mfu_\mu \frac{d\mu}{\mu}.
\end{equation}
\end{lemma}
The proof is straightforward, and left for the reader.

We next consider paramaterial derivatives of $g_\lambda^{ub}$ and $g_\lambda^{bal}$:

\begin{lemma}\label{l:dtg-low}
Assuming a full set of bootstrap bounds, we have 
\begin{equation}
 |(\partial_t + V_{[\lambda]}\cdot \nabla) g_\lambda|   \lesssim \lambda^\frac12 \mfu_\lambda + h_\lambda
\end{equation}
where 
\begin{equation}
 \|h_\lambda\|_{L^2} \lesssim C^2 \epsilon^2 c_\lambda^2 \lambda^{-\frac34}.   
\end{equation}
The same bound holds for    $g_\lambda^{ub}$ and $g_\lambda^{bal}$. 
\end{lemma}
We remark that here we use the paramaterial derivative associated  to frequency $\lambda$, even though 
$g_\lambda^{bal}$ is only localized to frequencies $\lesssim \lambda$.
\begin{proof}
As $g_\lambda = \nabla u_\lambda G'(\nabla u_{<\lambda})$,  
we have 
\begin{equation}\label{Dtgl}
|(\partial_t + V_{[\lambda]}\cdot \nabla)g_\lambda|\lesssim 
 |(\partial_t + V_{[\lambda]}\cdot \nabla)\nabla u_\lambda|+  |\nabla u_\lambda| |(\partial_t + V_{[\lambda]}\cdot \nabla)\nabla u_{<\lambda}| 
\end{equation}
For the first term we use the paradifferential equation \eqref{u-para}, estimating directly the linear paradifferential terms 
by $\lambda^{\frac12} \mfu_\lambda$,
\begin{equation}
 | (\partial_t + V_{[\lambda]}\cdot \nabla) \nabla u_\lambda| \lesssim \lambda^{\frac12} \mfu_\lambda + 
 |\nabla V_{[\lambda]}| |\nabla u_\lambda|+
 | \nabla N_\lambda(u)|.
\end{equation}
In $N_\lambda(u)$ we separate the 
 balanced quadratic term $N_\lambda^{[2],bal}$, for which we directly have 
\[
|N_\lambda^{[2],bal}| \lesssim 
\lambda^\frac32\calL_\lambda(|u_\lambda|, |u_\lambda|) \lesssim 
\lambda^{-\frac12} \mfu_\lambda.
\]
Using also the uniform bound 
$|\nabla V_{[\lambda]}|\lesssim 
C \epsilon \lambda^{\frac12}$, 
this yields
\begin{equation}\label{Dtu-peel}
 | (\partial_t + V_{[\lambda]}\cdot \nabla) \nabla u_\lambda| \lesssim \lambda^{\frac12} \mfu_\lambda + | \nabla(N_\lambda(u)-N_\lambda^{[2],bal}(u))|.
\end{equation}
To conclude, for the last term we use the $L^2$ bound \eqref{NuL2remainder}
for the difference $N_\lambda - N_\lambda^{[2],bal}$.

For the second term in \eqref{Dtgl} we expand and argue in a similar fashion,
\[
\begin{aligned}
|\nabla (\partial_t + V_{[\lambda]}\cdot \nabla)u_{<\lambda}| \lesssim & \  
\int_0^\lambda |\nabla (\partial_t + V_{[\mu]}\cdot \nabla)u_{\mu}|
+ \mu|V_{[\lambda]} - V_{[\mu]}| |\nabla u_{\mu}| \frac{d \mu}{\mu}
\\
\lesssim & \ \int_0^\lambda \mu^{\frac12} \mfu_\mu + |\nabla (N_\mu(u)-N_\mu^{[2],bal}(u))|\,  d \mu
\end{aligned}
\]
Finally, we can apply a pointwise bound for $\mfu_\mu$ and the    $L^2$ bound \eqref{NuL2remainder} for the
second term, followed by a careful application of Bernstein's inequality. 
\end{proof}

In addition to the above property, we also need a counterpart to the above lemma which takes advantage 
of the quadratic nature of $g_\lambda^{bal}$:
\begin{lemma}\label{l:dtg-bal}
Assuming a full set of bootstrap bounds for $u$ as in Section~\ref{sec:bootstrap}, we have 
\begin{equation}
 |(\partial_t + V_{[\lambda]}\cdot \nabla) g_\lambda^{bal}|   \lesssim 
  \mfu_\lambda \int_{\mu < \lambda} \mu^\frac12 \mfu_\mu \frac{d\mu}{\mu} + \int_{\nu > c\lambda} \mfu_\nu h_{\lambda \nu} \frac{d\nu}{\nu}
\end{equation}
where 
\begin{equation}
 \|h_{\lambda \nu}\|_{L^2} \lesssim C^2 \epsilon^2 c_\lambda c_\nu \lambda^{\frac34} \nu^{-\frac32}.   
\end{equation}
\end{lemma}

\begin{proof}
For brevity we consider only the first term in $g^{bal}_\lambda$,
\[
g^{bal,1}_\lambda := \int_{c\lambda}^\lambda \nabla u_\mu G''(\nabla u_{<\mu}) \frac{d\mu}{\mu},
\]
as the second 
is better. We expand as in the previous proof,
\[
\begin{aligned}
|(\partial_t + V_{[\lambda]}\cdot \nabla) g_\lambda^{bal}| 
\lesssim & \ |(\partial_t + V_{[\lambda]}\cdot \nabla)\nabla u_\lambda| \int_{c\lambda}^\lambda \mfu_\mu  \frac{d\mu}{\mu} + \mfu_\lambda \int_{c\lambda}^\lambda 
|(\partial_t + V_{[\lambda]}\cdot \nabla) \nabla u_\mu| \frac{d\mu}{\mu}
\\ & \ + \mfu_\lambda \int_{c\lambda}^\lambda 
\mfu_\mu 
|(\partial_t + V_{[\lambda]}\cdot \nabla) u_{<\mu}|\frac{d\mu}{\mu}.
\end{aligned}
\]
The last two terms already have a 
$\mfu_\lambda$ prefactor, and it suffices to estimate them exactly as in the previous lemma. For the first term we can use \eqref{Dtu-peel}. The contribution of the first term in \eqref{Dtu-peel} is already acceptable, so we are left with the 
expression
\[
 | \nabla(N_\lambda(u)-N_\lambda^{[2],bal}(u))| \int_{c\lambda}^\lambda \mfu_\mu  \frac{d\mu}{\mu}
\]
where we can simply estimate the first factor in $L^2$ by \eqref{Nuauxdecomp}.

\end{proof}

\newcommand{\const}{\kappa}
\subsection{Paracontrolled function spaces}\
All paracoefficients in our layered paradifferential expansion are of the type described in the 
previous subsection. However, this class of functions is not closed with respect to paraproducts, so 
in this section we expand the class of admissible paracoefficients into an algebra which is closed with respect to paraproducts.
We begin with
\begin{definition}
$\coeff$ is the space of functions $g$ of the form 
\[
g = g_0 + \sum g_\lambda,
\]
where $g_0$ is a constant, $|g_0| \leq  \const$, and $g_\lambda$ are localized at frequency $\leq \lambda$ with the following properties:

\begin{itemize}
    \item Size:
    \begin{equation}\label{g-in-B-point}
|g_\lambda| \leq \const \sum_{\mu \leq \lambda} (\mu/\lambda) \mfu_\mu + g_\lambda^{err},
\qquad \lambda \|g_\lambda^{err}\|_{L^\infty L^2} + \lambda^\frac14 \|g_\lambda^{err}\|_{L^2 L^\infty} \lesssim 
\const C \epsilon c_\lambda
 \end{equation}  
 \item Paramaterial derivative:
 \begin{equation}\label{dtg-dec}
 |(\partial_t + V_{[\lambda]}\cdot \nabla) g_\lambda| \leq \const \sum_{\mu \leq \lambda} (\mu/\lambda) \lambda^{\frac12} \mfu_\mu + h_\lambda,
\end{equation}
 where 
 \begin{equation}
 \| h_\lambda\|_{L^2L^8} \leq \const C \epsilon c_\lambda.     
 \end{equation}
 \end{itemize}
\end{definition}
The norm in $\coeff$ is taken to be the best constant $\const$ for which all the estimates above hold.
We also use the notations $\coeff_\lambda$ and $\coeff_{<\lambda}$ with the obvious meaning. We note that the choice of the $L^2L^8$ norm above is somewhat robust, 
and does not have an intrinsic significance;
but having an $L^2$ norm in time is important, as it allows us to match this with an energy bound, to get an $L^2_{xt}$
bilinear bound.

\medskip

The initial class $\coeff$ is a stepping stone toward introducing the class of paracontrolled functions 
$\coeffu$, which we will refer to as \emph{unbalanced}:

\begin{definition}
$\coeffu$ is the linear space of functions $f$ of the form 
\begin{equation}\label{B-unbal}
f = \sum_{\lambda} g_{< c\lambda} \nabla u_\lambda  
\end{equation}
where $g_{<\mu} \in \coeff_{<\mu}$ for each $\mu$, with the induced norm.
\end{definition}

By definition, the functions in $\coeffu$ are paracontrolled by $\nabla u$. We will also use the notations $\coeffu_\lambda$ and $\coeffu_{<\lambda}$, with the obvious meaning. For $s > 0$ we also denote by $D^{s} \coeffu$ the output of $D^s$ applied to $\coeffu$, or equivalently 
the space of functions $f$ of the form 
\[
f = \sum_\lambda \lambda^s g_{< c\lambda} \nabla u_\lambda  
\]
with $g_{<c\lambda}$ as in the definition above.

\medskip

Complementing the unbalanced class, we also introduce its balanced counterpart:

\begin{definition}
$\coeffb$ is the linear space of functions $f$ of the form 
\[
f = \sum_{\lambda} f_{\lambda}   
\]
where, for each $\lambda$, $f_\lambda$ is localized at frequency $\lesssim \lambda$
and satisfies the properties:
\begin{itemize}
    \item Size: 
    \begin{equation}\label{b-bal-size}
    |f_\lambda| \leq \const \sum_{\mu \leq \lambda} (\mu/\lambda)^{\frac12} \mfu_\mu \mfu_\lambda
    \end{equation}
    \item Paramaterial derivative:
    \begin{equation}\label{Dt-bbal}
 |(\partial_t + V_{[\lambda]}\cdot \nabla) f_\lambda|   \leq
  \const \mfu_\lambda \sum_{\mu < \lambda} \mu^\frac12 \mfu_\mu +\calL_\lambda(\sum_{\nu > c\lambda} \mfu_\nu h_{\lambda \nu})
\end{equation}
where 
\begin{equation}
 \|h_{\lambda \nu}\|_{L^2 L^\frac85} \lesssim \const C^2 \epsilon^2 c_\lambda c_\nu \lambda^{\frac74} \nu^{-\frac94} .   
\end{equation}
\end{itemize}
with the norm given by the best constant $\const$ above.

\end{definition}
We will also use the notations $\coeffb_\lambda$ and $\coeffb_{<\lambda}$, with the obvious meaning. Also for $s > 0$ we will on occasion work with $D^s \coeffb$.

Adding the two spaces 
of functions, we define 
\begin{equation}
    \coeffs = \coeffu + \coeffb
\end{equation}

As a first example of functions in this space, we 
have 
\begin{lemma} 
Functions of the form 
$g = F(\nabla u)$ or $g_{<\lambda} = F(\nabla u_{<\lambda})$, where $F$ is smooth with $F(0)=0$, belong to $\coeffs$.
\end{lemma}

\begin{proof}
This follows directly by 
applying Lemmas~\ref{l:gub+bal},\ref{l:dtg-low},~\ref{l:dtg-bal} for the decomposition
$g = g^{bal} + g^{ub}$.
\end{proof}

A second example which plays a role in our analysis is

\begin{lemma} \label{l:a-in-B}
The function $a$ belongs to $\coeff$, while its principal part $a^{[1]}$ belongs to $\coeffs$.   
\end{lemma}

This is a consequence of the bounds for $a$ in Proposition~\ref{abilinearbounds}, the   error term bound \eqref{aprincipalL2Linf}, 
as well as the $D_t a$ bound in \eqref{DtaprincipalL2L8}.

\medskip

Next,  we have the key algebra properties for $\coeff$ and $\coeffs$:

\begin{proposition}\label{p:algebra}
a) $\coeff$ is an algebra.

b) $\coeffs \subset \coeff$.

c) $\coeffb, \coeffs$ are algebras.      
\end{proposition}

\begin{proof}
a) Consider a product $g_\lambda g_\mu$ with $g_\lambda \in \coeff_\lambda$, $g_\mu \in \coeff_\mu$ and $\mu \leq \lambda$.
Then the product is localized at frequency $\lesssim \lambda$, and the estimate \eqref{g-in-B-point} is easily carried over to the product. 
We now apply the paramaterial derivative,
\[
(\partial_t + V_{[\lambda]}\cdot \nabla) (g_\mu g_\lambda) = 
g_\mu (\partial_t + V_{[\lambda]}\cdot \nabla)  g_\lambda
+ g_\lambda (\partial_t + V_{[\mu]}\cdot \nabla) g_\mu 
+ g_\lambda (V_{[\lambda]}- V_{[\mu]})\cdot \nabla g_\mu. 
\]
For the first term we simply use the uniform bound on $g_\mu$, and \eqref{dtg-dec} for $g_\lambda$.
For the second term we use   \eqref{dtg-dec} for $g_\mu$, so we need to estimate
\[
g_\lambda( \kappa \sum_{\nu \leq \mu} (\nu/\mu) \mu^{\frac12} \mfu_\nu + h_\mu)
\]
The first term is directly bounded by $\lambda^\frac12 g_\lambda$, which suffices. The second term is estimated in $L^2 L^8$, also directly.

 Finally, for the last term we use the $\mu^{-\frac12}$ uniform bound on $V_{[\lambda]}- V_{[\mu]}$
 and the uniform bound on $g_\mu$.

b) The nontrivial inclusion is $\coeffu \subset \coeff$. By the algebra property for $\coeff$, it suffices to show that $\nabla u_\lambda \in \coeff_\lambda$. But this follows from Lemma~\ref{l:dtg-low}.

c) 
We first show that $\coeffb$ is an algebra.
Due to part (a), it suffices to consider a product bound of the form
\begin{equation}\label{BxBbal}
\coeff_{\mu} \cdot \coeffb_\lambda 
\to \coeffb_{\lambda}, \qquad 
\mu \leq  \lambda.
\end{equation}
Consider a product $f_\mu g_\lambda$ of this type.
Since for $f_\mu \in \coeff_\mu$ we have 
the pointwise bound 
$|f_\mu| \lesssim c_\mu$,
the size bound \eqref{b-bal-size} immediately follows.

It remains to consider the paramaterial derivative bound
\eqref{Dt-bbal}. This is again straightforward if the paramaterial derivative goes on $g_\lambda$, so it remains to consider the expression 
\[
(\partial_t + V_{[\lambda]}\cdot \nabla) f_\mu  g_\lambda = 
(\partial_t + V_{[\mu]}\cdot \nabla) f_\mu  g_\lambda
+ (V_{[\lambda]}-V_{[\mu]})\cdot \nabla f_\mu  g_\lambda
\]
For the first term we use the representation \eqref{dtg-dec}. The contribution of the first term in \eqref{dtg-dec}
is easily treated, while the contribution of the second term is bounded by
\[
h_\mu \mfu_\lambda \sum_{\nu \leq \lambda} (\nu/\lambda)^\frac12 \mfu_\nu ,
\]
where we use Bernstein to estimate 
\[
\| h_\mu  \sum_{\nu \leq \lambda} (\nu/\lambda)^\frac12 \mfu_\nu \|_{L^2 L^\frac83} \lesssim 
\|h_\mu\|_{L^2 L^8}  \sum_{\nu \leq \lambda} (\nu/\lambda)^\frac12 \mfu_\nu \|_{L^2 L^4} \lesssim C^2 \epsilon^2 c_\nu\lambda^{-\frac12}
\]
which suffices.

For the second term we use 
\[
\| V_{[\mu]} - V_{[\lambda]}\|_{L^\infty} 
\lesssim \epsilon \mu^{-\frac12},
\]
combined with the uniform bound for $g_\mu$.

\bigskip

We now show that $\coeffs$ is an algebra.
We consider the product of two dyadic pieces
$f_\mu g_\lambda$  with $\mu \leq \lambda$,
of which at least one is in $\coeffu$.
If $g_\lambda \in \coeffb_\lambda$, then \eqref{BxBbal} suffices. Suppose now that 
$g_\lambda \in \coeffu_\lambda$. If $\mu \ll \lambda$, then we treat $f_\mu$ as a $\coeff$
paracoefficient, and then the algebra property for $\coeff$ suffices. It remains to consider the case when $\mu \approx \lambda$. If $f_\mu \in \coeffb$ then \eqref{BxBbal} still applies. Else, it remains to show that 
\begin{equation}\label{BuxBu}
\coeffu_{\mu} \cdot \coeffu_\lambda 
\to \coeffb_{\lambda}, \qquad 
\mu \approx  \lambda.
\end{equation}
Using again \eqref{BxBbal} to discard $\coeff$ paracoefficients, it suffices 
to show that $\nabla u_\lambda \nabla u_\mu \in \coeffb$. But this is exactly as in 
Lemma~\ref{l:dtg-bal}. 

The proof of the $\coeffs$ algebra property is now complete.
\end{proof}

As defined, $\coeff$ is a unital algebra.
However $\coeffb,\coeffs$ are not. We denote by $\tcoeffb$ and $\tcoeffs$ their unital extensions, obtained by adding complex constant functions; these are still included in $\coeff$. Finally $\tcoeffu$ will carry a similar meaning, even though it is not an algebra.

\bigskip

   Since the Besov space $\dot B^{d/2}_{2,1} \subset L^\infty_x$ 
is an algebra, also with good Moser estimates,
it is fairly straightforward to see that 
for $u \in L^\infty_t \dot B^{d/2}_{2,1}$ the space 
$\coeff$ is also contained in $L^\infty_t \dot B^{d/2}_{2,1}$. The frequency envelope bounds in \eqref{uk-ee}
also carry over easily:

\begin{lemma}
Assume that $u$ satisfies the uniform energy bounds \eqref{uk-ee}. Then similar bounds hold for $g$ in the unit $\coeff$ ball,
\begin{equation}
\| P_\mu g\|_{L^2} \lesssim C\epsilon c_\mu \mu^{-d/2}.     
\end{equation}
\end{lemma}

The lemma is a direct consequence of \eqref{g-in-B-point}.
We can also easily transfer bilinear $L^2$ bounds to $\coeff$  functions. Indeed, a direct application of 
\eqref{g-in-B-point} yields:

\begin{lemma} \label{l:DB-w}
Assume that $w$ satisfies the bootstrap bounds \eqref{uk-ee-boot}, \eqref{uw-bi-unbal-boot}. Then similar bounds hold for $g_\mu \in \coeff_\mu$ of unit size, 
\begin{equation}
\| \nabla g_\mu w_\lambda\|_{L^2} \lesssim C^2\epsilon^2 c_\mu d_\lambda  \mu^{\frac34} ,
\qquad \mu \ll \lambda.
\end{equation}
\end{lemma}

We can also use \eqref{dtg-dec} to show that
the para-material derivative of $\coeff$  functions satisfies good 
bilinear $L^2$ bounds:

\begin{lemma} \label{l:Dtg}
Assume that $u$ satisfies the bootstrap bounds \eqref{uk-ee-boot}, \eqref{uab-bi-unbal-boot}, etc. Then  for $g_\mu \in \coeff_\mu$ of unit size we have
\begin{equation}
\| T_{D_t} g_\mu \, w_\lambda\|_{L^2} \lesssim C^2\epsilon^2 c_\mu c_\lambda  \mu^{\frac14}, \qquad
\mu \ll \lambda.
\end{equation}
\end{lemma}

Our next task is to show that balanced expressions can be inserted into quartic expressions to produce $L^1$ bounds:

\begin{lemma}\label{l:bal-quad}
a) Let $g_\mu \in D^\frac12 \coeffb_\mu$, of unit size. Then  we have
\begin{equation}\label{bal-quad1}
\| g_{\mu}\calL_\lambda(w_\lambda,w_\lambda)\|_{L^1} \lesssim C^2 \epsilon^2 c_\mu d_\lambda^2, \qquad \mu \ll \lambda.   
\end{equation}
b)  For $g_\mu \in \coeffb_\mu$, of unit size, we have 
\begin{equation}\label{bal-quad2}
 \| (\partial_t + V_{[\lambda]} \cdot \nabla)   g_{\mu}\calL_\lambda(w_\lambda,w_\lambda)\|_{L^1} \lesssim C^2 \epsilon^2 c_\mu d_\lambda^2,  \qquad \mu \ll \lambda.       
\end{equation}
\end{lemma}
\begin{proof} 
The bound \eqref{bal-quad1} is a direct consequence of \eqref{b-bal-size} via two bilinear $L^2$ bounds.
The same applies in \eqref{bal-quad2} with the contribution of the first term in \eqref{Dt-bbal}. It remains to bound the contribution of the second term in \eqref{Dt-bbal},
\[
L_{\mu} (\mfu_\nu h_{\mu\nu})  \calL_\lambda(w_\lambda,w_\lambda), \qquad \nu > \mu
\]
Here we seek to pair $\mfu_\nu$ and $w_\lambda$ in a bilinear $L^2$ bound. We consider several cases:

\medskip

i) $ \mu \leq \nu < \lambda$. Then we can use an unbalanced bilinear $L^2$ bound for the $(\nu,\lambda)$
pairing, loosing a factor of $(\nu/\mu)^\frac12$ due to the extended shift. Then we get
\[
\calL_{\mu} (\mfu_\nu h_{\mu\nu})   L_\lambda(w_\lambda,w_\lambda) \in \nu^{-\frac14} (\nu/\mu)^\frac12  L^2
\mu^\frac74 \nu^{-\frac94} L^2 L^\frac83_\nu L^\infty L^2
\]
where we can use Bernstein to place the middle factor in $L^2 L^\infty$, and the frequency powers combine favourably. 

\medskip

ii) $ \mu < \lambda \ll \nu$. Here we can argue in the same manner as above.

\medskip

iii) $\mu < \nu \approx \lambda$. In this case we we can no longer use an unbalanced  bilinear $L^2$ bound
for the $(\nu,\lambda)$ pairing. Instead, we directly combine $L^4 L^{16}$ bounds for $\mfu_\nu$ and $w_{\lambda}$ with the $L^2 L^\frac83_\nu$ bound to get $L^2L^1$.

\end{proof}

The above Lemma will allow us to treat balanced paracoefficients in a perturbative manner later in the proof of density flux relations. The same, however, cannot be said about unbalanced paracoefficients,
and in particular we do not have access to  bounds such as \eqref{bal-quad1} and \eqref{bal-quad2}. In particular, it will be important to provide a clean representation for 
the paramaterial derivative of unbalanced functions. The next lemma 
carries out the first step, where 
we show that  placing the paramaterial derivative on $g_{<\lambda}$ yields perturbative terms.

\begin{lemma}\label{l:B-unbal-peel}
Let $f_\lambda \in \coeffu_\lambda$,
of the form $f_\lambda = g_{<\lambda} \nabla u_\lambda$, with $g_{<\lambda} \in \coeff_{<\lambda}$, of unit size.
Then we have 
\begin{equation}
(\partial_t + V_{[\lambda]} \cdot \nabla) f_\lambda 
= g_{<\lambda} \nabla (\partial_t + V_{[\lambda]} \cdot \nabla) u_\lambda
+ r_\lambda+r_\lambda^0
\end{equation}
where $r_\lambda$ satisfies bounds as in \eqref{Dt-bbal}
and $r_\lambda^0$ satisfies
\begin{equation}
r_\lambda^0 \lesssim \mfu_\lambda h_\lambda, \qquad \|h_\lambda \|_{L^2 L^\infty} \lesssim C \epsilon \lambda^\frac14    
\end{equation}
\end{lemma}
\begin{proof}
For $g = g_\mu$ with $\mu \ll \lambda$ we expand $r_\lambda$ as 
\[
\begin{aligned}
r_\lambda = &\ (\partial_t + V_{[\mu]} \cdot \nabla) g_\mu \nabla u_\lambda + (V_{[\lambda]}- V_{[\mu]}) \cdot \nabla  g_\mu \nabla u_\lambda + \nabla \cdot V_{[\lambda]}  g_\mu \nabla u_\lambda
\\
= & \ r_\lambda^1 + r_\lambda^2 + r_\lambda^3
\end{aligned}
\]
For $r_\lambda^1$ it suffices to 
use the bounds in \eqref{dtg-dec},  where the $h_\mu$ factor localized to frequency $\leq \mu$ yields the $r_\lambda^0$ term. 

For $r_\lambda^2$ we consider the two terms in the $g_\mu$ bound provided by \eqref{g-in-B-point}. For the first term it suffices to use the uniform bound $|V_{[\lambda]}- V_{[\mu]}| \lesssim C \epsilon \mu^{-\frac12}$, while for the second we estimate
\[
\| (V_{[\lambda]}- V_{[\mu]}) \cdot \nabla  g^{err}_\mu\|_{L^2} \lesssim 
\| V_{[\lambda]}- V_{[\mu]}\|_{L^\infty L^2} \|g^{err}_\mu\|_{L^2 L^\infty}
\lesssim C \epsilon c_\mu \mu^{-\frac32} \cdot C \epsilon c_\mu \mu^{\frac54}
= C^2 \epsilon^2 c_\mu^2 \mu^{-\frac14}
\]
using Bernstein to get into $L^2 L^\frac83$.

Finally, for $r_\lambda^3$ we use the pointwise bound for $g_\mu$, and split 
\[
V_{[\lambda]} =  V_{[\lambda]}^{[1]} + (V_{[\lambda]} - V_{[\lambda]}^{[1]}),
\]
where for the first term we have
\[
\nabla V_{[<\lambda]}^{[1]} \lesssim \sum_{\mu< \lambda} \mu^{\frac12}
\mfu_\mu,
\]
while for the remainder we have, by Corollary~\ref{principalVBexp},
\[
\|V_{[\lambda]} - V_{[\lambda]}^{[1]}\|_{L^2 L^\infty} \lesssim  C \epsilon \lambda^{\frac34}.
\]

\end{proof}

\subsection{Symbols and (para)commutators}

We denote by $S^{s}_{\lambda}$ the class of symbols of bilinear operators $A$
whose symbols $a(\xi,\eta)$ are smooth 
on the $\lambda$ scale. We denote the classes of rough symbols with coefficients in 
the spaces $\coeff$, $\coeffu$ and $\coeffb$
by $\coeff_\lambda S^{s}_{\lambda}$, $\coeffu_\lambda S^{s}_{\lambda}, \coeffb_\lambda S^{s}_{\lambda}$.

For pseudodifferential operators we will use the left quantization.  
Some choice of quantization is needed 
in order to define the associated bilinear
forms. Our choice here is also akin to the left quantization. Precisely, for a symbol 
$a(x,\xi,\eta)$ we define the associated 
bilinear form
\begin{equation}
A(u,\bar v)(x):= \iint a(x,\xi,\eta) 
e^{ix(\xi-\eta)} \hat u(\xi) \ol{\hat{v}(\eta)}\,  d \xi d\eta.
\end{equation}

To evaluate commutators of spatial pseudodifferential operators and separate the 
principal part we recall the notion of 
para-Lie bracket from \cite{AIT-minimal},
\begin{equation}\label{para-Lie}
\{ p,q\}_{p} = T_{q_\xi} p_x - T_{p_\xi} q_x,
\end{equation}
where the paraproducts above only involve 
the spatial dependence.

At the space-time level, the only operator we
work with is the material derivative $D_t$. 
Its paradifferential counterpart is called the \emph{para-material derivative}, defined 
on functions as
\begin{equation}\label{para-material}
T_{D_t} g := (\partial_t + T_V \cdot \nabla) g.
\end{equation}
Applied to symbols $a(t,x,\xi)$, the natural extension is suggested by the above para-Lie bracket,
\begin{equation}\label{para-material+}
T_{D_t} a := \partial_t a  + \{ V \cdot \xi, a\}_{p}.
\end{equation}
By a slight abuse, we will call it the para-material derivative of $a$ as well.

Turning our attention to the class of paracontrolled functions $\coeffu$, we make two heuristic observations concerning their para-material derivative:

\begin{itemize}
\item The para-material derivative of $g_\mu$ plays  a purely perturbative role; this is already apparent in Lemma~\ref{l:B-unbal-peel}.

\item The para-material derivative of $u_\mu$
is given by $T_{D_t} u_\mu = \omega_{[\mu]}(x,D) u_\mu + N_\mu$ where the second component $N_\mu$ contains only (i) lower order $low-hi$ interactions  and (ii) $high-high$ interactions which will also play a perturbative role in many instances. In the next section we 
will see how to peel off  such contributions.
\end{itemize}

\subsection{ The restricted paramaterial derivative and paratransported symbols}
\label{sec:paratransport}


Our goal in this subsection  is to
 separate the leading linear  part of the para-material derivative of our paracontrolled $\coeffu$ functions.  Toward this goal, we introduce the \emph{restricted para-material derivative} as follows:

\begin{definition}
Let $f \in \coeffu$. Then we define its  restricted para-material derivative $\rT f$ by 
\begin{equation}
\rT f : =   i   T_{\omega} f 
\end{equation}
\end{definition}
This definition easily extends to symbols in $\coeffuu S^s_\lambda$ in the obvious manner, exactly as for the earlier para-material derivative.
 
To justify this definition, our goal
here will be to show that the difference 
\[
T_{D_t} f - \rT f
\]
is perturbative in a suitable sense.
Ideally, we would like to show that it satisfies an estimate akin to 
\eqref{bal-quad2} in Lemma~\ref{l:bal-quad}. This turns out to be almost true, modulo a flux type term which 
is also good enough in order to play only a perturbative role later on.

Our main result below asserts that indeed, the restricted paramaterial derivative captures the 
nonperturbative part of the paramaterial derivative of $\coeffu$ functions:

\begin{proposition}\label{p:Dt-Dtr}
Let $f \in \coeffu_{< \lambda}$, of unit size.
Then the error term 
\[
R^4_\lambda = (T_{D_t} f - \rT f)L_\lambda(w_\lambda,w_\lambda)
\]
admits a representation 
\begin{equation}
   R^4_\lambda = \partial Q^{err}_\lambda + R^{err}_\lambda ,
\end{equation}
with a flux correction $Q^{err}$ satisfying
\begin{equation}\label{Q-err}
 \| Q_\lambda^{err}\|_{L^1_{t}L^p_x} \lesssim C^2 \epsilon^2 d_\lambda^2 \lambda^{d(\frac12-\frac1p)-1}   
\end{equation}
for $p > \frac43$, and an $L^1$ error
\begin{equation}\label{R-err}
 \| R_\lambda^{err}\|_{L^1_{tx}} \lesssim C^2 \epsilon^2 d_\lambda^2    \end{equation}

\end{proposition}
\begin{proof}
Considering a dyadic decomposition of $f \in \coeffu_{<\lambda}$ we take a 
frequency $\mu \ll \lambda$ and  $f_\mu \in \coeffu_\mu$,
of the form $f_\mu = g_{<c\mu} \nabla u_\mu$, with $g_{<c\mu} \in \coeff_{<c\mu}$. We expand 
\[
T_{D_t} f - \rT f = 
(\partial_t + V_{[\mu]} \cdot \nabla) (g_{<c\mu} \nabla u_\mu) - i \omega_{[\mu]}(c,D) (g_{<c\mu} \nabla u_\mu) 
\]
By Lemma~\ref{l:B-unbal-peel} and by the proof of Lemma~\ref{l:bal-quad}, we can freely replace the expression 
$(\partial_t + V_{[\mu]} \cdot \nabla) (g_{<c\mu} \nabla u_\mu)$ by $g_{<c\lambda} \nabla (\partial_t + V_{[\lambda]} \cdot \nabla) u_\lambda$ in the definition of $R_\lambda^4$. We can also commute
$\omega_{[\mu](x,D)}$ inside $g_{<c\mu}$ modulo a commutator which can be bounded by 
\[
\mu^{-\frac12}\calL_\mu(\nabla g_{<c_\mu},\mfu_\mu), 
\]
and further by the right hand side in \eqref{Dt-bbal}. In addition, commuting $\omega_{[\mu](x,D)}$ with $\nabla$ also yields directly a balanced contribution as in \eqref{Dt-bbal}.

Thus we are left with the expression
\[
g_{<c_\mu} \nabla 
( \partial_t + V_{[\mu]} \cdot \nabla + i \omega_{[\mu](x,D)}) u_\mu
\]
Using the equation \eqref{u-para},
we have 
\[
 ( \partial_t + V_{[\mu]} \cdot \nabla + i \omega_{[\mu](x,D)}) u_\mu = (- \frac12 \nabla \cdot V_{[\mu]}+ \vartheta(x,D)_{[\mu]}) u_\mu + N_\mu(u)
\]
We next consider the contribution of the paradifferential term  $(- \frac12 \nabla \cdot V_{[\mu]}+ \vartheta(x,D)_{[\mu]}) u_\mu$, which we claim is directly balanced 
 and controlled as in \eqref{Dt-bbal}. To see this we recall that 
 \[
 \vartheta_{[\mu]}(x,\xi) = D_{t,[\mu]} \ln \mfq_{[\mu]}(x,\xi) = 
 \frac12 a_{[\mu]}^{-1} D_{t,[\mu]}  a_{[\mu]} + \dn^{-1}_{[\mu]}(x,\xi)
 D_{t,[\mu]} \dn_{[\mu]}
 \]
In the first term we can use to replace $\nabla V_{[\mu]}$ by $\nabla V^{[1]}_{[\mu]}$, then use Proposition~\ref{Dtacommutatorestimates}
to replace
$D_{t,[\mu]}  a_{[\mu]}$ 
with $(D_t a)^{[1]}_\mu$ and finally use Proposition~\ref{p:Dt-Dtr} to replace $D_{t,[\mu]} \dn^{-1}_{[\mu]}$ by $D^R_{t,[\mu]} \dn^{-1}_{[\mu]}$ modulo perturbative errors. The remaining symbol is in $D^\frac12 \coeffs S^9$, and a slight extension of Lemma~\ref{l:DB-w}
 applies.

  This leaves us with 
the expression $g_{<c_\mu} \nabla N_\mu$, whose contribution to $R_\lambda^4$ is
\[
g_{<c_\mu} \nabla N_\mu (u)\calL_\lambda(w_\lambda,w_\lambda)
\]
Neglecting $g_{<c_\mu}$ for the moment, the desired representation 
is provided by \eqref{Nuauxdecomp}. 
Here $g_{<\mu}$ is bounded and thus harmlessly added in the second term.
In the first term, adding $g_{<\mu}$ yields an extra source term
\[
\nabla g_{<\mu} Q^{err}_\lambda,
\]
which can be estimated in $L^1_{t,x}$ and placed into $R^{err}_\lambda$
using the $L^2_{t,x}$ bound $\|\nabla g_{<\mu}\|_{L^\infty L^2} \lesssim \epsilon$ and the $p=2$ case of \eqref{Q-err}.
\end{proof}

We remark that using the restricted 
paramaterial derivative on balanced functions is not so useful but harmless, as Lemma~\ref{l:bal-quad}
shows that both terms in the difference are perturbative.

Our final definition here introduces the notion of \emph{para-transported symbols}:
\begin{definition}
We say that a symbol $g \in \coeffs S^s$
is para-transported along the particle flow
if $\rT g$ is balanced, i.e. it belongs to $D^\frac12 \coeffb$,
in which case we write
\[
\rT g  \paraequal 0.
\]
\end{definition}
This definition is automatically satisfied for the balanced part of 
$\coeffs S^s$ symbols, so it only 
represents a requirement on the unbalanced part.

We also note that this definition can be
applied to scalar functions, but it is not very helpful there as it essentially singles out the class of balanced functions. However, 
it will be very useful for symbols, where we remark that smooth symbols are not para-transported. Instead, in the next section we construct several 
useful classes of paratransported symbols.

\subsection{ The good vector fields and Littlewood-Paley  multipliers.}

Classically one constructs the symbols for the Littlewood-Paley multipliers as
\[
p_\lambda(\xi) = \chi_\lambda(|\xi|^2) 
\]
with a smooth bump function $\chi_\lambda$ localizing to the $\lambda^2$ region.
However, if we use such multipliers to 
produce frequency localization in our setting, then we immediately run into a problem, 
which is that commuting the para-material derivative  $\partial_t + V_{[\lambda]} \partial_x$ with $P_\lambda$ yields unbalanced contributions, which are difficult to control or even correct later on. Alternatively, we could try to adapt 
our metric to the transport term. Using 
an exact transport is tempting but brings forth another set of difficulties, as 
the flow map seems difficult to control in our function classes.
Instead, we will seek an intermediate option, where we do not have exact transport, but instead a balanced para-material derivative. Even here, considerable care must be taken when defining a good notion of para-material
derivative. As it turns out, the right choice
is provided by the restricted para-material derivative defined in the previous subsection.
With  this notion in place, we can now state our main para-transport construction:

\begin{theorem}\label{t:good-g} Assume that $u$ is a solution to the full equation \eqref{u-para} satisfying the  (bootstrap) bounds \eqref{uk-ee-boot}-\eqref{uL6L4-boot}. Then

a) There exist vector fields $ \xis_j$
with $\tcoeffs$ coefficients 
and leading component $\xi_j$ which are para-transported along the particle flow, i.e.,
its restricted para-material derivative $\rT \xis_j$ has coefficients 
in $D^\frac12 \coeffb$. 

b) There exists a metric $\gs$ with $\tcoeffs$ coefficients which is para-transported along the particle flow, i.e.,
its restricted para-material derivative $\rT \gs_\lambda$ has coefficients  in $D^\frac12 \coeffb_{\lambda}$. 
\end{theorem}

\begin{proof}

It suffices to prove (a). Once this is done, we can define the 
metric $g$ as 
\[
g^s(\xi, \xi) :\paraequal |\xi^\sharp|^2.
\]

The condition we will impose is for the principal symbol of the commutator,
\[
\rT \gamma \paraequal 0, \qquad \gamma := \xis_j.
\]
Expanding 
\[
\xis_j = f^{k}_j \xi_k,
\]
with $f^{k}_j \in \coeffu$,
the above condition becomes 
\[
\rT f^{k}_j \paraequal   T_{ f^l_j} \partial_l V^k_{<\lambda} .
\]

To construct such coefficients $ f^l_j$, it is convenient to start with a vector field  $\calv$ at the level of $\eta$ so that
\begin{equation}\label{calv}
\rT \calv   \paraequal  V .
\end{equation}
Since $D_t \eta = B$ and $V = ML^{-1} B$, we can achieve this by setting
\begin{equation}\label{calv-def}
\calv = T_{ML^{-1}} \eta,
\end{equation}
which is easily verified to satisfy \eqref{calv}, and also has the property $\partial_x \calv \in \coeffu$.
Then it suffices to choose $f^k_j$ with the recursive property that 
\begin{equation}\label{fkj-recursion}
  f^k_j =   T_{ f^l_j} \partial_l \calv^k.
\end{equation}

We can solve this iteratively by
\[
f_j^{k,(m+1)} := f_j^{k,(m)} + T_{f_j^{l,(m)}} \partial_l \calv^k, \qquad f_j^{k,(0)} = \delta_j^k.
\]
Using the smallness bootstrap assumption one can easily show first convergence in $L^\infty_x$, 
then convergence in $\coeffu$ via Proposition~\ref{p:algebra}.
Then the relation \eqref{fkj-recursion} shows that $f^{k}_j \in\coeffu$.
\end{proof}

For later use, we write the paratransport equation for $\calv$:

\begin{lemma}\label{l:calv}
The vector field $\calv$ satisfies an equation of the form
\begin{equation}
T_{D_t} \calv = V + \calN(u),    
\end{equation}
where $\calN(u)$ satisfies bounds similar to $N(u)$. 
\end{lemma}
The proof is relatively straightforward and is omitted. We remark that modulo paradifferential commutator terms, the leading part of $\calN(u)$ is $T_{ML^{-1}} N(u)$.

\begin{remark}\
\begin{itemize}
    \item 
Having such a metric, it is then natural 
to  define the ``\emph{good  Littlewood-Paley projectors}'' $\LPs_{\lambda}$ with symbols
\begin{equation}\label{def:p-lambda}
\lps_\lambda = \chi(\lambda^{-2} \gs_\lambda (\xi, \xi)).
\end{equation}

\item We also remark that our construction allows us to arbitrarily choose the constant coefficient, leading part of $\xis_j$ and  $\gs$.

\end{itemize}
\end{remark}

\section{Density-flux identities}
\label{sec:density-flux}

Here we establish density-flux identities for a succession of three problems of increasing difficulty:

\begin{enumerate}[label=(\roman*)]
    \item For the linear, constant coefficient problem,
    \item for the linear paradifferential equation,
    \item for the full nonlinear flow, and
    \item for the linearized equation.
\end{enumerate}
The first case will serve as a model, while the next two cases will both be critical for the proof of our results.

In the last subsection, we use these density-flux identities in order to close the energy component of our bootstrap.

\subsection{The linear, constant coefficient case}\label{subsec:df-cc}
Here we consider the equation
\begin{equation}\label{eq:flat}
i \partial_t u - \omega_0(D) u = 0.
\end{equation}
For our problem  $\omega_0$ is the  constant coefficient dispersion relation is 
\[
\omega_0(\xi) = |\xi|^\frac12,
\]
though the exact form is not used in the computations below.

We just do the computations for the mass and the momentum, as frequency localization can be easily added in. To start with, we set
\[
\fM_0(u) = |u|^2
\]
to be the mass density and compute
\[
\partial_t \fM_0(u) = i B(u,\bu), \qquad b(\xi,\eta) = -\omega_0(\xi) + \omega_0(\eta),
\]
where $B(u, \bar{u})$ is a bilinear form, with $b(\xi, \eta)$ its associated symbol.
To express the right hand side in spatial divergence form we need to write
\begin{equation}\label{def-pj}
\omega_0(\xi) - \omega_0(\eta) = -(\xi-\eta)_j \fp^j_0 (\xi,\eta).
\end{equation}
The choice of the symbols $\fp_0^j$ with this 
properties is not uniquely determined. But 
we introduce a canonical way to choose them by using the fundamental theorem of calculus on the straight line between $\xi$ and $\eta$ in order to make the symmetric choice
\begin{equation}\label{def-pj0}
\fp^j_0(\xi,\eta) = -\int_0^1 \partial_j \omega_0(\xi+h(\eta-\xi))\, dh.
\end{equation}
Then \eqref{def-pj} holds, and we have 
\begin{equation}\label{df-m-cc}
\partial_t \fM_0(u) = \partial_j \fP^j_0(u).
\end{equation}

Similarly, we have a density-flux relation for the momentum,
\begin{equation}\label{df-p-cc}
\partial_t \fP^j_0(u)  = \partial_k \fE_0^{jk}(u) ,
\end{equation}
where $\fE_0^{jk}(u)$ is the energy density and $\fe^{jk}_0(\xi,\eta)$ is its symbol:
\[
\fe^{jk}_0(\xi,\eta) = \fp^j_0(\xi,\eta) \fp^k_0(\xi,\eta).
\]

We can also write frequency localized versions of these relations, with the localized symbols defined as
\begin{equation}
   \fm_\lambda(\xi,\eta) = p_\lambda(\xi) p_\lambda(\eta), \qquad \fp^j_\lambda(\xi,\eta) = \fp^j_0(\xi,\eta) p_\lambda(\xi) p_\lambda(\eta), \qquad \fe^{jk}_\lambda(\xi,\eta) = \fe^{jk}_0(\xi,\eta) p_\lambda(\xi) p_\lambda(\eta).
\end{equation}
Note that $p_\lambda (\cdot )$'s are the symbols of the Littlewood-Paley projectors. 
Our next task is to define counterparts 
of these quantities for the paradifferential equation.

\subsection{Density-flux relations in the paradifferential case}
Here we consider frequency localized solutions $w_\lambda$ for the inhomogeneous 
paradifferential equation
\begin{equation}\label{para-re}
   (\partial_t + V_{[\lambda]} \cdot \nabla + \frac12 (\nabla \cdot V_{[\lambda]})- \bb_{[\lambda]}(x,D) +
i \omega_{[\lambda]}(x,D)) w_\lambda  = f_\lambda,
\end{equation}
where both the vector field $V_{[\lambda]}$ and
the real scalar symbol $\omega_{[\lambda]}$ of order $\frac12$ depend on the lower frequencies of our water wave solution $u$,
\begin{equation}\label{para-coeff-reg}
\nabla V_{[\lambda]} \in D^\frac12 \coeffs_\lambda, 
\qquad \omega_{[\lambda]} \in \coeffs_\lambda S^\frac12_\lambda,
\end{equation}
and the real scalar symbol $\bb$ is given by 
\[
\bb_{[\lambda]}(x,\xi)  = \frac12 (T_{D_t} \ln \mfq_{[\lambda]})(x,D).
\]

Inspired by construction of paratransported symbols in the Section~\ref{sec:paratransport}, we would like to define 
symbols $\ms_\lambda(\xi,\eta)$, $\psj_\lambda(x,\xi,\eta)$ and their associated fluxes
so that the density-flux identities hold
exactly as in the constant coefficient case, with balanced errors. 
All the paracoefficients in our construction will belong to $\tcoeffs$. By a slight abuse of notation, we will assume they have at most unit size.

\medskip
As a preliminary step, we carry out a generic computation:

\begin{lemma}\label{l:df}
Consider a frequency localized density symbol $d_\lambda$ and its associated flux symbol $f^j_\lambda =  \fp^j d_\lambda$. 
Assume  that $d_{\lambda}$ is symmetric with $d_{\lambda} \in \tcoeffs S^0_\lambda$.  
Then we have the density-flux identity
\begin{equation}\label{df-para-lem}
\begin{aligned}
 (\partial_t + \nabla \cdot V_{[\lambda]}) D_\lambda(w_\lambda, \bw_\lambda) 
 = & \ \partial_j F_\lambda(w_\lambda, \bw_\lambda)  + D_\lambda(f_\lambda, \bw_\lambda)+  D_\lambda(w_\lambda, \bar{f}_\lambda)
 \\ & \ +
 R_\lambda^{(0)} +R_\lambda^{(\frac12)} + R_\lambda^{(0),bal},
 \end{aligned}
 \end{equation}
where the symbol of the main source term is given by
\begin{equation}\label{df-R0}
r_\lambda^{(0)} = T_{D_t} d_\lambda +  d_\lambda \bb_{[\lambda]},
\end{equation}
while the symbols of the secondary source terms satisfy
\[
r_\lambda^{(\frac12)} \in D \coeffs S^{-\frac12}_\lambda +\tcoeffs \nabla \vartheta S^{-1}_\lambda 
, \qquad r_\lambda^{(0),err} \in D^\frac12 \coeffb S^0_\lambda.
\]
\end{lemma}
Here by a slight abuse of notation we define the extension of $ \bb_{[\lambda]}$ as a function of two Fourier variables $(\xi,\eta)$ in a symmetric fashion as 
\begin{equation}
    \bb_{[\lambda]}(x,\xi,\eta) : = \frac12( \bb_{[\lambda]}(x,\xi) + \bb_{[\lambda]}(x,\eta)).
\end{equation}
One may view this lemma as a pure 
microlocal commutator formula, which does not make use of our bootstrap assumptions, other than through \eqref{para-coeff-reg}.
These will be used later on in order to estimate error terms such as  $R_\lambda^{(0),bal}$.
We make several remarks about the leading order source term $R_\lambda^{(0)}$.
\begin{remark}[On the choice of $r_\lambda^{(0)}$] \

\begin{enumerate}[label=(\roman*)]
\item If $d_\lambda$ does not depend on either $\xi$ or $\eta$, then the first term in $r_\lambda^{(0)}$ represents the (para)Poisson bracket of $D_{t,[\lambda]}$
and $d_\lambda$. So we could call the above expression a multilinear (para)Poisson bracket.

\item As we will see later, in the symbol of $r_\lambda^{(0)}$ 
the important part is its restriction to the diagonal $\xi = \eta$. Indeed, if $ r_\lambda^{(0)}$ were to vanish on the diagonal then we can rewrite the associated bilinear form as a divergence plus an $R_\lambda^{(\frac12)}$ term.  Furthermore, the symbol of $r_\lambda^{(0)}$  restricted to the diagonal $\xi = \eta$ depends only on the symbol of $d_\lambda$ on the diagonal.

\item Any potential balanced terms in $R_\lambda^{(0)}$ and also in $R_\lambda^{(\frac12)}$ can be harmlessly included into the error $R^{(0),bal}_\lambda$, see also the upcoming Lemma~\ref{l:DtB}.
\end{enumerate}
\end{remark}

\begin{proof}
Separating variables, we reduce to the 
multiplicative case 
\begin{equation}\label{d-sep}
d_\lambda (x,\xi,\eta) = g(x) h(\xi,\eta),
\end{equation}
where  $g \in \tcoeffs_{<\lambda}$
 and $h \in S^0_\lambda$, localized at frequency $\lambda$. Then we compute
\[
\begin{aligned}
 (\partial_t + \nabla \cdot V_{[\lambda]}) D_\lambda(w_\lambda, \bw_\lambda)    
= & \  (\nabla \cdot V_{[\lambda]})  D_\lambda(w_\lambda, \bw_\lambda) + 
(\partial_t + V_{[\lambda]}  \cdot \nabla) g
\, H(w_\lambda, \bw_\lambda) 
\\
& + g H((\partial_t + V_{[\lambda]}  \cdot \nabla) w_\lambda, \bar w_\lambda)+ g H(w_\lambda,(\partial_t + V_{[\lambda]}  \cdot \nabla)  \bar w_\lambda)
\\
& + g [ V_{[\lambda]} H(\nabla w_\lambda, \bar w_\lambda) - H(V_{[\lambda]} \nabla w_\lambda, \bar w_\lambda)] + symmetric
\\
= & \ E_1+ E_2 + E_3 + E_4 + E_5 + E_6.
\end{aligned}
\]    
Leaving the first  term $E_1$ untouched for now, in $E_2$ we consider the balance of frequencies in $V_{<\lambda}$ and $g$. Given a dyadic component $V_{\nu} \nabla g_\mu$ with $\nu,\mu \ll \lambda$,
if $\mu \lesssim \nu$, its contribution 
can be placed in $R_\lambda^{(0),bal}$.
Then we are left with
\[
E_2 = T_{D_t} g H(w_\lambda, \bw_\lambda) + R^{(0),bal}_\lambda.
\]
For $E_3$ and $E_4$ we use the $w_\lambda$ 
equation. This yields several contributions as follows:
\begin{enumerate}[label=(\alph*)]
\item  the contribution of the source term $f_\lambda$ in \eqref{df-para-lem}, 

\item the expression 
\[
- \frac12 g [H((\nabla \cdot V_{<\lambda})  w_\lambda, \bar w_\lambda)+  H(w_\lambda,  (\nabla \cdot V_{<\lambda}) \bar w_\lambda)],
\]
which cancels $E_1$ modulo a lower order commutator which can be placed into $R^{(\frac12)}_\lambda$, and whose balanced part can be included in $R^{(0),bal}_\lambda$,

\item the expression 
\[
- \frac12 g [H((\bb_{[\lambda]}(x,D)  w_\lambda, \bar w_\lambda)+  H(w_\lambda,  \bb_{[\lambda]}(x,D) \bar w_\lambda)],
\]
which yields the second term in $r_\lambda^{(0)}$ modulo a commutator term, which yields the second term in  $r_\lambda^{(\frac12)}$. 

\item the contribution of the dispersion relation
\begin{equation}\label{H-delta-omega}
- ig [H( \omega_{[\lambda]}(x,D) w_\lambda, \bar w_\lambda) -  H(w_\lambda,   \ol{\omega_{[\lambda]}(x,D)  w_\lambda})].
\end{equation}
This has principal symbol
\[
-i g(x) (\omega(x,\xi)-\omega(x,\eta)) h(\xi,\eta) ,
\]
plus a commutator term in $D\coeffs S^{-\frac12}_{\lambda}$, which can be placed 
in $R_\lambda^{\frac12}$. The principal symbol can be expanded as 
\[
-i g(\xi-\eta)_j p^j(x,\xi,\eta) h(\xi,\eta).
\]
Taking the derivative out, we arrive at
\[
\partial_j g F^j(w_\lambda,\bar w_\lambda) + 
R_\lambda^{\frac12}.
\]
\end{enumerate}
Lastly, we consider the commutator terms 
$E_5$ and $E_6$. For $E_5+E_6$ we write
the principal symbol 
\[
 g(\xi_j \partial_k V^j \partial_{\xi_k} h(\xi,\eta)   + \eta_j \partial_k V^j \partial_{\eta_k} h(\xi,\eta)) , 
\]
where the remainder can be placed 
again in $R_\lambda^{(\frac12)}$. Finally,
the above principal symbol is unbalanced only if the $V$ frequency is larger than the $g$ frequency.
\end{proof} 

We now examine the source terms in the density-flux relation of the last lemma.
As a warm up, we begin with the balanced source term:

\begin{lemma}
If $u$ and $w$ satisfy 
the bootstrap assumptions 
in \eqref{uk-ee-boot}-\eqref{uL6L4-boot}, respectively \eqref{wk-ee-boot}-\eqref{wL6L4-boot},
then the balanced source term $R^{(0),bal}_\lambda$ satisfies the  bound
\begin{equation}\label{R-bal-err}
     \| R_\lambda^{(0),bal}\|_{L^1_{t,x}} \lesssim C^2 \epsilon^2 d_\lambda^2.
\end{equation}
\end{lemma}

\begin{proof}
 Recalling that $r^{(0),bal} \in D^\frac12 \coeffb$, after Littlewood-Paley decomposition  and separation of variables, it suffices to consider expressions of the form 
\[
g_\mu\calL_\lambda(w_\lambda, \bar w_\lambda) \qquad \mu \ll \lambda,
\]
with $g_\mu \in D^\frac12 \coeffb$.
But then we can apply Lemma~\ref{l:bal-quad}.
 
\end{proof}

The above lemma will allow us to treat the contribution on $R_\lambda^{(0),bal}$ perturbatively later on. In the same vein, we more generally introduce
\begin{definition}
We denote by $R_\lambda^{err}$ any expression which satisfies the bound \eqref{R-bal-err}.
\end{definition}

We now examine more closely the leading 
source term symbol $r^{(0)}_\lambda$
in \eqref{df-R0}. One may be tempted to assert that $r^{(0)}_\lambda \in D^\frac12 \coeffs S^0_\lambda$. This might\footnote{ This has to do with the available bounds for $a$; we will return to this in the proof of Proposition~\ref{p:first-df}.}  be the case for the second term in it, but  not necessarily the case for the first term, due to the fact that its expression requires the full water wave equation, which includes contributions from higher frequencies. Our strategy will be to place such contributions 
in the error term.

\begin{lemma}
\label{l:DtB}
If $u$ and $w_\lambda$ satisfy 
the  bootstrap assumptions in \eqref{uk-ee-boot}-\eqref{uL6L4-boot}, and \eqref{wk-ee-boot}-\eqref{wL6L4-boot}, respectively, then
the  symbol $T_{D_t} d_\lambda$  admits the decomposition
\begin{equation}
   T_{D_t} d_\lambda = \rT d_\lambda + r_\lambda, 
\end{equation}
where  the associated density $R_\lambda$ admits a representation
\begin{equation}
   R_\lambda = \partial Q^{err}_\lambda + R^{err}_\lambda ,
\end{equation}
with a flux correction $Q^{err}$ satisfying
\begin{equation}\label{Q-err+}
 \| Q_\lambda^{err}\|_{L^1_{t}L^p_x} \lesssim C^2 \epsilon^2 d_\lambda^2 \lambda^{d(\frac12-\frac1p)-1}   
\end{equation}
for $p > \frac43$.
\end{lemma}

\begin{proof}
Using separation of variables, it suffices  to consider $d_\lambda$ of the form \eqref{d-sep}. Then we have 
\[
T_{D_t} d_\lambda - \rT d_\lambda = 
( T_{D_t} g_\lambda - \rT g_\lambda) h_\lambda,
\]
so it remains to consider the expression
\[
( T_{D_t} g_\lambda - \rT g_\lambda)H_\lambda(w_\lambda,\bw_\lambda). 
\]
But now we can directly use Proposition~\ref{p:Dt-Dtr}.

\end{proof}
\begin{remark} \label{r:df}
The above Lemma shows that we  can replace 
the para-material derivative $T_{D_t} d_\lambda$ in the $r^{(0)}$ expression in Lemma~\ref{l:df} with the corresponding restricted para-material 
derivative $\rT d_\lambda$
at the (small) price of allowing more general 
$R^{err}_\lambda$ sources and
a flux correction $Q^{err}_\lambda$ in \eqref{df-para-lem}. 
\end{remark}

Our main task now is to obtain the variable coefficient analogue of the density-flux relations \eqref{df-m-cc}, \eqref{df-p-cc}.
The result in the above lemma is not fully satisfactory, due to the leading error term $R^{(0)}$. 
On one hand we cannot estimate it directly using bilinear estimates, on the other hand, using its nonresonant structure to remove it using matching density  and flux corrections yields corrections 
which are at the same level as $D_\lambda$ and $F_\lambda$, and the problem repeats, leaving us with the unenviable task of having to do an infinite number of corrections.

We will still rely directly on this lemma for subprincipal density/flux relations, but for the principal part it is essential  to have good choices for the densities and the associated fluxes, in order to eliminate the $R^{(0)}$ type term. Precisely, we need symbols for our densities which paracommute
with the material derivative. In particular,
the classical Littlewood-Paley projectors do not have this property and need to be replaced.

We begin with the mass density,  but where two key changes are required relative to the constant coefficient choice in  Section~\ref{subsec:df-cc}:
\begin{enumerate}
    \item The classical Littlewood-Paley projectors  do not paracommute with the paramaterial derivative. To rectify this, we replace them with their paratransported versions $\LPs_\lambda$ defined using the good metric $\gs$ in Theorem~\ref{t:good-g}
as in \eqref{def:p-lambda}.

   \item The paracoefficient $\vartheta$ in   the paradifferential flow \eqref{para-re} 
is not balanced, and could potentially produce direct mass growth. To eliminate its effect, we use the fact that it is the logarithmic paramaterial derivative of $\mfq^{-\frac12}$.
This suggests that a mass density correction of $\mfq^{-\frac12}$ should rectify this growth. While this idea goes in the right direction, the symbol $\mfq$ contains an $a$ 
factor, which we do not know\footnote{This is not to say that $a$ is not in $\coeffs$;
likely it is, but proving it is not entirely straightforward.}
if it is in $\coeffs$. To address this, we replace\footnote{While this substitution is harmless here, it would not be at all acceptable in the diagonalization of the $(\eta,\varphi)$ system} $a$ 
by its principal part $a^{[1]} = T_a T_\mfl \eta$. Since $a$ has been showed to belong to $\coeff$, this modified weight is in $\coeffs$, and in effect in $\coeffu$. 
\end{enumerate}

Given the above discussion, it is natural to define the modified mass symbol as follows: 
\begin{equation}\label{def:ms}
 \ms_\lambda(x,\xi,\eta) := 
 \tps_{\lambda}(x,\xi) \tps_\lambda(x,\eta),
 \qquad \tps_\lambda(x,\xi): = \lps_{\lambda}(x,\xi)\lambda^\frac14 \tilde \mfq_\lambda^{-\frac12}(x,\xi).
\end{equation}
where the modified version of $\mfq_\lambda$ is given by
\begin{equation}
  \tilde \mfq_\lambda(x,\xi) =   (\dn/a^{[1]})^\frac12_{[\lambda]}.
\end{equation}
The normalization factor $\lambda^\frac14$ is simply to insure that the mass density symbol has order zero. With these definitions, we have the symbol regularity
\[
\ms \in \tcoeffs S^0.
\]

Once we have the mass density, following  Lemma~\ref{l:df}, we can define the corresponding flux symbols 
\begin{equation}\label{choose-mflux}
\fsj_{\lambda}(x,\xi,\eta):= \fp^j(x,\xi,\eta) \ms_\lambda(x,\xi,\eta) \in \tcoeffs S^{-\frac12}_\lambda.
\end{equation}

The next difficulty, compared with the constant coefficient case, is that we can no longer choose the above flux symbols as the next densities, as $\fp^j$ do not 
paracommute with the material derivative.
Instead, we take as our starting point the 
constant coefficient version of the momentum $\fp_0$ as in \eqref{def-pj0}, and  define the modified momentum density symbols
$\psj_\lambda$ by
\begin{equation}\label{choose-psj}
   \psj_\lambda(x,\xi,\eta) = 
   \tfp^j(x,\xi,\eta) \ms_\lambda(x,\xi,\eta)\in \tcoeffs S^{-\frac12}_\lambda., \qquad \tfp^j(x,\xi,\eta): = \fp^j_{0}(\xi^\sharp,\eta^\sharp),
   \end{equation}
and the modified momentum flux symbols
\begin{equation}\label{choose-pflux}
\esj_{\lambda}(x,\xi,\eta)=  \psj_\lambda(x,\xi,\eta) \fp^k(x,\xi,\eta)\in \tcoeffs S^{-1}_\lambda.
\end{equation}

\begin{remark}\label{lp-tight}
As a minor variation of the above setup,
it is sometimes convenient to work with a slightly tighter localization in frequency. 
Precisely, given a frequency $\xi_0$ of size
$\lambda$, we can redefine the symbols $\lps_\lambda$ to localize to a small neighbourhood of $\xi_0$, still on the dyadic scale,
\[
\lps_\lambda = \chi(\frac{\xis_{\lambda}-\xi_0}{c\lambda}), \qquad c \ll 1.
\]
This is akin to splitting the prior $\lps_\lambda$ choice into small angle sectors.
For the bulk of the analysis, this finer choice makes no difference, but in a few secondary cases it does help in order to simlify slightly the argument by eliminating certain resonant interactions. To keep 
the argument streamlined, we do not introduce a new notation to reflect this. 
Instead, we will refer to this remark whenever this tighter localization is used.

\end{remark}

Now we can state our first density-flux result for paratransported densities as follows:

\begin{proposition}\label{p:first-df}
Assume the bootstrap assumptions in \eqref{uk-ee-boot}-\eqref{uL6L4-boot}, and \eqref{wk-ee-boot}-\eqref{wL6L4-boot}, respectively, for both $u$ and $w$, including the mixed case.
Then we have the identities
\begin{equation}
\left\{
\begin{aligned}
&\partial_t \Ms_{\lambda}+\nabla\cdot (V_{[\lambda]}\Ms_{\lambda})=\partial_j(\Fsj_{\lambda}+Q_{\lambda,m}^{j,err})+R^{3}_{\lambda,m} + R^{4}_{\lambda,m}  + 2 \Re \Ms_\lambda(v_\lambda,\bar f_\lambda)\\
&\partial_t \Psj_{\lambda}+\nabla\cdot (V_{[\lambda]}\Psj_{\lambda})=\partial_k (\Esj_{\lambda}+Q_{\lambda,p}^{jk,err})+R^{3}_{\lambda,p} + R^{4}_{\lambda,p}  + 2 \Re \Psj_\lambda(v_\lambda,\bar f_\lambda),
\end{aligned}
\right.
\end{equation}  
where 
\begin{enumerate}[label=\roman*)]
    \item $R^3_{\lambda,m}$ and $R^3_{\lambda,p}$
 are bilinear forms with symbols in 
 $D\coeffs S^{-\frac12}_\lambda$, respectively $D\coeffs S^{-1}_\lambda$.

\item $Q^{err}_{\lambda,m}$ and $Q^{err}_{\lambda,p}$ are flux
error terms estimated in $L^1_t L^p_x$
for $p > \frac43$,
\begin{equation}\label{df-err-flux}
\| Q_{\lambda,m}^{err}\|_{L^1_{t}L^p_x} \lesssim C^2 \epsilon^2 d_\lambda^2 \lambda^{d(\frac12-\frac1p)-1},  
\qquad
\| Q_{\lambda,p}^{err}\|_{L^1_{t}L^p_x} \lesssim C^2 \epsilon^2 d_\lambda^2 \lambda^{d(\frac12-\frac1p)-\frac32}; 
\end{equation}

\item $R^4_{\lambda,m}$ and $R^4_{\lambda,p}$ are 
error terms estimated in $L^1$, 
\begin{equation}\label{df-err-source}
\|R^4_{\lambda,m}\|_{L^1} \lesssim C^4 \epsilon^4 d_\lambda^2 , 
\qquad
\|R^4_{\lambda,p}\|_{L^1} \lesssim C^4 \epsilon^4 d_\lambda^2 \lambda^{-\frac12}.
\end{equation}
\end{enumerate}
\end{proposition}

\begin{proof}
Both parts of the Proposition are obtained from Lemma~\ref{l:df}, taking into account several observations:
\begin{enumerate}[label= (\roman*)]
    \item By Remark~\ref{r:df}, 
 we can replace $T_{D_t} \Ms_\lambda$ with $\rT \Ms_\lambda$.
 \item Despite the above substitution,
 the Leibniz rule still applies modulo 
 balanced terms, as for $f,g \in \coeffs$ we have
 \[
 T_\omega(fg) \paraequal f T_\omega g + g T_\omega f, 
 \]
as all commutator type terms are automatically balanced.
 
 \item For the $a$ component we have  $\rT a^{[1]} \paraequal (D_t a)^{[1]}$.

 \item The paracoefficient difference 
$T_{D_t} a - (T_{D_t} a)^{[1]}$ is perturbative by \eqref{Dtacommutatorfluxdecomp}.
\end{enumerate}
Based these observations, it becomes straightforward to use \eqref{def:ms}, \eqref{choose-mflux}, \eqref{choose-psj} and \eqref{choose-pflux}
in order to verify that the symbols $\ms_\lambda$
and $\psj_\lambda$ satisfy the correct commutation relation with the 
para-material derivative,
\begin{equation}
\rT \ms_\lambda \paraequal -\ms_\lambda \rT \ln \tilde \mfq_{[\lambda]}^\frac12,
\qquad \rT \psj_\lambda \paraequal - \psj_\lambda \rT \ln \tilde \mfq_{[\lambda]}^\frac12.
\end{equation}

These properties show that the rectified\footnote{ as in Lemma~\ref{l:df}}  symbol $r^{(0)}_\lambda$ is balanced, $r^{(0)} \in D^\frac12\coeffb S^{0}_\lambda$, therefore it can be estimated by two bilinear $L^2$ bounds as in Lemma~\ref{l:bal-quad}. This leaves us with the $R^{(\frac12)}_\lambda$ terms, where we can again use Lemma~\ref{l:bal-quad} to bound in $L^1$  the balanced part, and include in $R^3_{\lambda,m}$ the unbalanced one. 

\end{proof}

The density-flux relations in the previous proposition are still not sufficiently accurate for our purposes due to the presence of the unbalanced $R^3$ source terms,  which are not perturbative. The redeeming feature is that these can be essentially viewed as cubic terms (albeit with a lower frequency paracoefficient), which are both lower order and nonresonant, and thus can be eliminated  with suitable density-flux corrections. This is the goal of the next theorem:

\begin{theorem}\label{t:second-df}
Assume the bootstrap assumptions for both $u$ and $w$: \eqref{uk-ee-boot}-\eqref{uL6L4-boot}, and \eqref{wk-ee-boot}-\eqref{wL6L4-boot}, respectively. Then there exist density corrections $ B^3_{\lambda,m}$ and $B^{3}_{\lambda,p}$ and flux corrections 
$ Q^3_{\lambda,m}$ and $Q^{3}_{\lambda,p}$
so that the associated modified mass 
and modified momentum
\begin{equation}
   \Mss_{\lambda} = \Ms_{\lambda}+ B^3_{\lambda,m} , 
   \qquad  \Pssj_{\lambda} = \Psj_{\lambda}+ B^{j,3}_{\lambda,p}
\end{equation}
satisfy the relations 
\begin{equation}
\left\{
\begin{aligned}
&\partial_t \Mss_{\lambda}+\nabla\cdot (V_{[\lambda]}\Mss_{\lambda})=\partial_j\cdot (\Fsj_{\lambda}+ Q_{\lambda,m}^{j,3}+ Q_{\lambda,m}^{j,err}
) +
2 \Re \Mss_{\lambda}(w_\lambda,\bar f_\lambda) + R^{err}_{\lambda,m}  
\\
&\partial_t \Pssj_{\lambda}+\nabla\cdot (V_{[\lambda]}\Pssj_{\lambda})=\partial_k\cdot (\Esj_{\lambda}+ Q_{\lambda,p}^{jk,3}+ Q_{\lambda,p}^{jk,err}) +
2 \Re\Pssj_{\lambda}(w_\lambda,\bar f_\lambda)+ R^{j,err}_{\lambda,p}  .
\end{aligned}
\right.
\end{equation}  
Here the corrections and the source terms are as follows: 

\begin{enumerate}[label=\roman*)]
    \item 
 The symbols for the density corrections   $B^{3}_{\lambda,m}$ and $B^{j,3}_{\lambda,p}$
 are in $D^\frac12\coeffu S^{-\frac12}_\lambda$,
 respectively $D^\frac12 \coeffu S^{-1}_\lambda$;

 \item The symbols for the flux corrections   $Q^{j,3}_{\lambda,m}$ and $Q^{jk,3}_{\lambda,p}$
 are in $D^\frac12 \coeffu S^{-1}_\lambda$,
 respectively $D^\frac12 \coeffu S^{-\frac32}_\lambda$, given by 
 \begin{equation}\label{b-q-relation}
  q^{j,3}_{\lambda,m} = b^{3}_{\lambda,m}  \fp^j, 
  \qquad q^{jk,3}_{\lambda,p} = b^{j,3}_{\lambda,p}  \fp^k;
 \end{equation}

\item The additional flux terms 
$Q_{\lambda,m}^{j,err}$, 
$Q_{\lambda,p}^{jk,err}$ satisfy the bounds \eqref{df-err-flux};

\item The errors $R^{err}_{\lambda,m}$,
$R^{j,err}_{\lambda,p}$ satisfy the bounds \eqref{df-err-source}.
\end{enumerate}

\end{theorem}

\begin{proof}
We provide the argument for the mass relation; the one for the momentum relation is absolutely identical. The starting point is provided by Proposition~\ref{p:first-df}, where $R^3_{\lambda,m}$ is the only term which needs to be corrected. We 
split $R^3_{\lambda,m}$ into a dyadic sum
\[
R^3_{\lambda,m} = \sum_{\mu  \ll \lambda} R^3_{\lambda,m,\mu}
\]
where we
will distinguish  two types of paracoefficients in the symbols $r^3_{\lambda,m,\mu} \in D\coeffs_{\mu} S^{-\frac12}_\lambda$, as follows:

\begin{itemize}
\item Balanced, i.e. $r^{3}_{\lambda,m} \in \mu \coeffb_\mu S^{-\frac12}_\lambda$.

\item Unbalanced, i.e. 
$r^{3}_{\lambda,m} \in \mu \coeffu_\mu S^{-\frac12}_\lambda$,
schematically of the form
\[
r^{3}_{\lambda,m} = \mu  \coeff_{\ll \mu} \nabla u_\mu S^{-\frac12}_\lambda.
\]

\end{itemize}

The balanced terms can be estimated directly 
using Lemma~\ref{l:bal-quad},
which is akin to using two unbalanced bilinear bounds,
and placed into the perturbative error $R^{err}$. So it suffices to correct the unbalanced terms. 

For this we revisit the computation in Lemma~\ref{l:df}. Assume, after a partial separation of variables, that the symbol of $R^3_{\lambda,m}$ has the form 
\[
r^3_{\lambda,m}(x,\xi,\eta) = \mu \lambda^{-\frac12} g_{\ll \mu} L_\mu(\nabla u_\mu,\xi,\eta),
\]
where $g_{\ll \mu} \in \coeff$ and $L_\mu$
is a linear form of $\calL_\mu$ type with a smooth dependence on $\xi$ and $\eta$ on the $\lambda$ scale. Here 
one should also consider the nearly identical case when $\nabla u_\mu$ is replaced by $\nabla \bu_\mu$.
Correspondingly, we seek a correction with symbol of the form\footnote{It will become a-posteriori clear that we want the same 
$g$.}
\[
b^3_{\lambda,m}(x,\xi,\eta) = \mu^{\frac12} \lambda^{-\frac12}
g_{\ll \mu} L^1_\mu(\nabla u_\mu,\xi,\eta).
\]
Following the same computation as in Lemma~\ref{l:df}, we
identify the unbalanced contributions
in $(\partial_t + \partial_x \cdot V_{[\lambda]}) B^3$ as follows:

\medskip

a) In $E_2$ we have the unbalanced term
\[
\mu^\frac12 \lambda^{-\frac12}(\partial_t + V_{[\mu]}\cdot \nabla ) (g_{\ll \mu} L^1_\mu(\nabla u_\mu,\xi,\eta)),
\]
which, as in Lemma~\ref{l:DtB},  can be replaced by 
\[
i\mu^\frac12 \lambda^{-\frac12} g_{\ll \mu} \omega_{[\mu]}(x,D) L^1_\mu(\nabla u_\mu,\xi,\eta)
\]
modulo perturbative $Q^{err}$ and $R^{err}$
flux and source term corrections. 
Here we note that the contribution
of the difference $(V_{[\lambda]}-V_{[\mu]})$ can be estimated perturbatively using two unbalanced bilinear $L^2$ bounds.

\medskip

b) In $E_3$ and $E_4$ we have 

(i) the $\nabla \cdot V$ commutator term, which is balanced even before any commuting

(ii) The lower order linear paradifferential term, again perturbative

(iii) The leading order term as in \eqref{H-delta-omega}, which yields the principal symbol
\[
i\mu^\frac12 \lambda^{-\frac12} g_{\ll \mu} L^1_\mu( \nabla u_\mu, \xi,\eta) (\xi_j-\eta_j)\fp^j (x,\xi,\eta) .
\]
Pulling out the $\partial_j$ derivative 
in the bilinear form yields the desired 
flux contribution,
\[
\partial_j Q^{3,j}(w,\bw), \qquad 
q^{3,j}(x,\xi,\eta)=\mu^\frac12 \lambda^{-\frac12} g_{\ll \mu} L^1_\mu( \nabla u_\mu, \xi,\eta) \fp^j (x,\xi,\eta), 
\]
which is consistent with the relation \eqref{b-q-relation}.
In addition, we obtain commutator terms where the $\partial_j$ 
derivative falls on coefficients.
The terms with a derivative on $g_{\ll \mu}$
or on $\fp^j$ are balanced and thus perturbative, leaving us with the expression
\[
 \mu^\frac12 \lambda^{-\frac12} g_{\ll \mu} \partial_j L^1_\mu( \nabla u_\mu, \xi,\eta) \fp^j (x,\xi,\eta) .
\]
Here only the $\fp^j$ coefficients at frequency $\ll \mu$ yield nonperturbative contributions,
so we can replace $\fp^j (x,\xi,\eta)$ by $\fp_{[\mu]}^j$.

c) In $E_5$ and $E_6$ we get expressions of the form
\[
\mu^\frac12 \lambda^{-\frac12} \tilde\calL_\mu( \nabla u_\mu, \nabla V_{[\lambda]}, w_\lambda,\bar w_\lambda),
\]
which is again balanced regardless of the size of the $V$ frequency.

\bigskip

Summarizing, we have arrived at the relation
\begin{equation}
 (\partial_t + \partial_x \cdot V_{[\lambda] }) B^3  - \partial_j Q^{j,3}= \mu^{\frac12} \lambda^{-\frac12}  g_{\ll \mu} ({i}\omega_{[\mu]}(x,D) +  \fp^j_{[\mu]}(x,\xi,\eta) \partial^j) L^1_\mu(\nabla  {u_\mu},  \xi,\eta), 
\end{equation}
modulo balanced contributions.
Matching our $R^3_{\lambda,m}$, it follows that we 
need to choose $L^1_\mu$ so that
\[
(i\omega_{[\mu]}(x,D) +  \fp^j_{[\mu]}(x,\xi,\eta) \partial^j) L^1_\mu(\nabla  u_\mu,  \xi,\eta)  = -\mu^{\frac12}\calL_\mu(\nabla   u_\mu,  \xi,\eta) + \text{balanced}.
\]
Since $\fp^j_{[\mu]}(x,\xi,\eta)$ has size $\lambda^{-\frac12}$ and $\mu \ll \lambda$, the pseudodifferential operator on the left restricted to frequency $\mu$ 
has an elliptic symbol of size $\mu^\frac12$.
Hence it suffices to set
\begin{equation}
L^1_\mu = -i\mu^{\frac12}   Op((
\omega_{[\mu]}(x,\zeta)+\fp^j_{[\mu]}(x,\xi,\eta) \zeta_j)^{-1})(x,D) L_\mu  ,
\end{equation}
where all subprincipal commutator terms are balanced and thus perturbative. This concludes the proof of the theorem.
\end{proof}

\subsection{The nonlinear case}
Here we consider frequency localized density-flux identities for the full   nonlinear equation, which we write in the form 
\begin{equation}
   (\partial_t + V_{[\lambda]} \cdot \nabla +\frac12 (\nabla \cdot V_{[\lambda]}) - \vartheta_{[\lambda]}(x,D) +
i\omega_{[\lambda]}) u_\lambda  = N_\lambda(u).
\end{equation}
Here we do not assume that $u$ is exactly localized at frequency $\lambda$, instead we 
simply ask that the above relation holds whenever a slightly larger dyadic frequency projector is applied to both sides of the above equation. This will allow the density and flux functionals to enforce the stricter localization via the sharp Littlewood-Paley symbols $\lps_\lambda$.

Our main density flux identities for $u$ are as follows:

\begin{theorem}\label{t:df-NL}
Assume the bootstrap assumptions \eqref{uk-ee-boot}-\eqref{uL6L4-boot}, for $u$ hold. Then there exist density corrections $ B^{NL}_{\lambda,m}$ and $B^{j,NL}_{\lambda,p}$ and flux corrections 
$ Q^{j,NL}_{\lambda,m}$ and $Q^{jk,NL}_{\lambda,p}$
so that the associated modified mass 
and modified momentum
\begin{equation}\label{MP-NL}
   \Msnl_{\lambda} = \Mss_{\lambda}+ B^{NL}_{\lambda,m} , 
   \qquad  \Psnlj_{\lambda} = \Pssj_{\lambda}+ B^{j,NL}_{\lambda,p}
\end{equation}
satisfy the relations 
\begin{equation}
\left\{
\begin{aligned}
& \ \partial_t \Msnl_{\lambda}+\nabla\cdot (V_{[\lambda]}\Msnl_{\lambda})=\partial_j (\Fsj_{\lambda}+ Q^{j,3}_{\lambda,m} + Q^{j,err}_{\lambda,m}
+  Q^{NL}_{\lambda,m})  + R^{err}_{\lambda,m}  \\
& \ \partial_t \Psnlj_{\lambda}+\nabla\cdot (V_{[\lambda]}\Psnlj_{\lambda})=\partial_k (\Esj_{\lambda}+ Q_{\lambda,p}^{jk,3}
+ Q_{\lambda,p}^{jk,err}
+  Q^{jk,NL}_{\lambda,m}) + R^{j,err}_{\lambda,p} . 
\end{aligned}
\right.
\end{equation}   
Here the corrections and the source terms are as follows:

\begin{enumerate}[label=\Alph*)]
    \item 
 The additional density corrections   $B^{NL}_{\lambda,m}$ and $B^{NL}_{\lambda,p}$
 are as follows:
 \begin{equation}
  B^{NL}_{\lambda,m} =  B^{NL,[1]}_{\lambda,m}+
  B^{NL,[2]}_{\lambda,m} + B^{NL,[3]}_{\lambda,m}+ B^{NL,[4]}_{\lambda,m}
 \end{equation}
 and similarly for $B^{NL}_{\lambda,p}$, where
 \begin{itemize}
     \item The paradifferential corrections  $B^{NL,[1]}_{\lambda,m}$ and  $B^{j,NL,[1]}_{\lambda,p}$  have the form
 \[
     B^{j,NL,[1]}_{\lambda,m} = Op(D^\frac12 \coeffu_{<\lambda} S^{-\frac12}_\lambda)(\bfu_\lambda,\bfu_\lambda), 
     \qquad B^{jk,NL,[1]}_{\lambda,p} = Op(
    D^\frac12 \coeffu_{<\lambda} S^{-1}_\lambda)(\bfu_\lambda,\bfu_\lambda);
    \]

\item The  quadratic $hh$ corrections
$B^{NL,[2]}_{\lambda,m}$ and $B^{j,NL,[2]}_{\lambda,p}$ have the form\footnote{ In other words, the leading entry in the paracoefficient has frequency exactly $\lambda$, possibly with further lower frequency 
paracoefficients.},
\[
B^{NL,[2]}_{\lambda,m} = \sum_{\mu \gtrsim \lambda} \lambda^{-\frac12} Op(\coeffu_\lambda S^{\frac12}_{\mu})(\bfu_\mu,\bfu_\mu)
\qquad
B^{j,NL,[2]}_{\lambda,p}= \sum_{\mu \gtrsim \lambda} \lambda^{-1} Op(\coeffu_\lambda S^{\frac12}_{\mu})(\bfu_\mu,\bfu_\mu);
\]

\item The balanced cubic nonresonant corrections
$B^{NL,[3]}_{\lambda,m}$ and $B^{j,NL,[3]}_{\lambda,p}$ have the form
\[
B^{NL,[3]}_{\lambda,m}=Op(\coeff_{<\lambda} S^{2}_\lambda)(\bfu_\lambda,\bfu_\lambda,\bfu_\lambda,\bfu_\lambda)
\qquad
B^{j,NL,[3]}_{\lambda,p}= Op(\coeff_{<\lambda} S^{\frac32}_\lambda)(\bfu_\lambda,\bfu_\lambda,\bfu_\lambda,\bfu_\lambda);
\]

\item The balanced quartic nonresonant corrections
$B^{NL,[4]}_{\lambda,m}$ and $B^{j,NL,[4]}_{\lambda,p}$  have the form
\[
B^{NL,[4]}_{\lambda,m}=Op(\coeff_{<\lambda} S^{3}_\lambda)(\bfu_\lambda,\bfu_\lambda,\bfu_\lambda,\bfu_\lambda,\bfu_\lambda)
\qquad
B^{j,NL,[4]}_{\lambda,p}= Op(\coeff_{<\lambda} S^{\frac52}_\lambda)(\bfu_\lambda,\bfu_\lambda,\bfu_\lambda,\bfu_\lambda,\bfu_\lambda);
\]
 \end{itemize}

  \item The additional flux corrections   $Q^{j,NL}_{\lambda,m}$ and  $Q^{jk,NL}_{\lambda,p}$
 are as follows:
 \begin{equation}
   Q^{j,NL}_{\lambda,m}  = Q^{j,NL,[1]}_{\lambda,m}+ Q^{j,NL,[2]}_{\lambda,m} + Q^{j,NL,[3]}_{\lambda,m}+ Q^{j,NL,[4]}_{\lambda,m},
 \end{equation}
 and similarly for  $Q^{jk,NL}_{\lambda,p}$, where
 \begin{itemize}
     \item The paradifferential corrections  $Q^{j,NL,[1]}_{\lambda,m}$ and $Q^{jk,NL,[1]}_{\lambda,p}$ have the form
     \[
     Q^{j,NL,[1]}_{\lambda,m} = Op(D^\frac12 \coeffu_{<\lambda} S^{-1}_\lambda)(\bfu_\lambda,\bfu_\lambda), 
     \qquad Q^{jk,NL,[1]}_{\lambda,p} = Op(
    D^\frac12 \coeffu_{<\lambda} S^{-\frac32}_\lambda)(\bfu_\lambda,\bfu_\lambda),
    \]
    where the components with phase rotation symmetry\footnote{i.e. with arguments $(u_\lambda,\bu_\lambda)$} satisfy  the relations 
    \begin{equation}\label{bq-relation-nl}
     q^{j,NL,[1]}_{\lambda,m} = \fp^j_\lambda
     b^{NL,[1]}_{\lambda,m}, \qquad
    q^{jk,NL,[1]}_{\lambda,p}  = \fp^k_\lambda
     b^{j,NL,[1]}_{\lambda,p};
    \end{equation}

\item The quadratic $hh$ corrections $Q^{j,NL,[2]}_{\lambda,m}$ and  $Q^{jk,NL,[2]}_{\lambda,p}$  have the form
\[
Q^{NL,[2]}_{\lambda,m} = \sum_{\mu \gtrsim \lambda} \lambda^{-\frac12} Op(\coeffu_\lambda S^{0}_{\mu})(\bfu_\mu,\bfu_\mu)
\qquad
Q^{j,NL,[2]}_{\lambda,p}= \sum_{\mu \gtrsim \lambda} \lambda^{-1} Op(\coeffu_\lambda S^{0}_{\mu})(\bfu_\mu,\bfu_\mu);
\]

\item The balanced cubic nonresonant corrections
 $Q^{j,NL,[3]}_{\lambda,m}$ and $Q^{jk,NL,[3]}_{\lambda,p}$ have the form
 \[
Q^{NL,[3]}_{\lambda,m}=Op(\coeff_{<\lambda} S^{\frac32}_\lambda)(\bfu_\lambda,\bfu_\lambda,\bfu_\lambda,\bfu_\lambda)
\qquad
Q^{j,NL,[3]}_{\lambda,p}= Op(\coeff_{<\lambda} S^{1}_\lambda)(\bfu_\lambda,\bfu_\lambda,\bfu_\lambda,\bfu_\lambda);
\]

\item The balanced quartic nonresonant corrections
 $Q^{j,NL,[4]}_{\lambda,m}$ and  $Q^{jk,NL,[4]}_{\lambda,p}$   have the form
\[
Q^{NL,[4]}_{\lambda,m}=Op(\coeff_{<\lambda} S^{\frac52}_\lambda)(\bfu_\lambda,\bfu_\lambda,\bfu_\lambda,\bfu_\lambda,\bfu_\lambda)
\qquad
Q^{j,NL,[4]}_{\lambda,p}= Op(\coeff_{<\lambda} S^{2}_\lambda)(\bfu_\lambda,\bfu_\lambda,\bfu_\lambda,\bfu_\lambda,\bfu_\lambda);
\]

 \end{itemize}

\item 
The additional flux terms  $Q^{j,err}_{\lambda,m}$ and $Q^{jk,err}_{\lambda,p}$ satisfy, with $p > \frac43$,
\begin{equation}\label{df-err-flux-NL}
\| Q_{\lambda,m}^{err}\|_{L^1_{t}L^p_x} \lesssim C^2 \epsilon^2 c_\lambda^2 \lambda^{-2s_c} \lambda^{d(\frac12-\frac1p)-1},  
\qquad
\| Q_{\lambda,p}^{err}\|_{L^1_{t}L^p_x} \lesssim C^2 \epsilon^2 c_\lambda^2 \lambda^{-2s_c} \lambda^{d(\frac12-\frac1p)-\frac32}; 
\end{equation}
The errors $R^{err}_{\lambda,m}$  $R^{j,err}_{\lambda,p}$ satisfy the perturbative bounds
\begin{equation}\label{df-err-nl-m}
\begin{aligned}
 \|R^{err}_{\lambda,m}\|_{L^1_{t,x}} \lesssim C^4 \epsilon^4 c_\lambda^2 \lambda^{-2s_c},
\qquad 
 \|R^{j,err}_{\lambda,p}\|_{L^1_{t,x}} \lesssim C^4 \epsilon^4 c_\lambda^2 \lambda^{-2s_c-\frac12}.
\end{aligned}
\end{equation}

\end{enumerate}
 \end{theorem}

\bigskip 

\begin{proof}
We consider the mass identities, as the 
momentum identities are similar.
The starting point of the proof is the result in Theorem~\ref{t:second-df}, applied with $w_\lambda = u_\lambda$ and $f_\lambda = N_\lambda(u)$, which 
yields the desired density-flux relation with the additional source term
\[
 \qquad R = 2\Re \Mss_\lambda(N_\lambda(u),\bar u_\lambda) .
\]
Expanding the nonlinearity as in Theorem~\ref{ufinaltheorem},
\[
N_\lambda(u) = N^{[1]}_\lambda(u)+
N^{[2]}_\lambda(u) + N^{[3]}_\lambda(u)+N^{[4]}_\lambda(u)
+ F_\lambda(u)
\]
the contribution of $F_\lambda$ is directly estimated perturbatively. It remains to compute corrections for the 
contributions  $R^{[k]} := 2\Re\Mss_\lambda(N^{[k]}_\lambda(u),\bar u_\lambda)$, $k=1,\dots,4$, of each of the first four terms. Before we do that, we briefly need to discuss
the role played by the Taylor coefficient $a$. 

\medskip

\textbf{0. The $a$ management.} 
The discussion here is necessary because of the somewhat limited information we have about $a$, compared to all other paracoefficients in our problem. Precisely, as a consequence of the results in Sections~\ref{sec:aexpansions} and 
\ref{sec:Dtaexpansions}, we know that $a \in \coeff$ (see Lemma~\ref{l:a-in-B}), but not\footnote{It is likely that this can be established with some substantial additional effort; if that were the case, this discussion would become unnecessary.} that $a \in \tcoeffs$. We remark that $a$ appears in all the paracoefficients in the $N^{[k]}$ terms
via the paradifferential operator $T_{\mfq}$ appearing in the diagonalization of the $(\eta,\varphi)$ system, in two ways:

\begin{itemize}
    \item From the substitution of $\varphi$ in terms of $u$ in the $N_{\eta,\lambda}$
    and in the $N_{\varphi, \lambda}$ expressions.
    \item From the definition of $N_{u,\lambda}$ in terms of $N_{\eta,\lambda}$
    and  $N_{\varphi, \lambda}$.
\end{itemize}

The way we deal with this depends on which $N^{[k]}_{u,\lambda}$ we consider, as follows:

\begin{enumerate}
\item In $N^{[1]}_{u,\lambda}$: Here the Taylor coefficient $a_{[\lambda]}$ appears only within commutator structures, interacting with other low frequency paracoefficients,
so for simplicity consider one such commutator term, $[a_{[\lambda]}g_{<\lambda},b(x,D)]$ where $g_{<\lambda} \in \tcoeffs_{<\lambda}$ and $b \in \tcoeffs_{<\lambda} S^0_\lambda$. Separating variables, it suffices to consider a commutator of the form 
\[
C = [a_{[\lambda]}g_{<\lambda},h_{<\lambda} H_\lambda(D)] = h_{<\lambda}[a_{[\lambda]}g_{<\lambda}, H_\lambda(D)], \qquad h_{<\lambda} \in \tcoeffs_{<\lambda}.
\]
By standard commutator formulas, this can be expressed in the form
\[
Cw_\lambda = \lambda^{-1} h_{<\lambda} L_{\lambda}(\nabla (a_{[\lambda]}g_{<\lambda}), w_\lambda).
\]
At this point we consider several cases depending on the frequencies  of the three paracoefficients $h_{<\lambda}$, $a_{[\lambda]}$, respectively $g_{<\lambda})$.
\begin{itemize}
    \item If the two highest frequencies are comparable, then the source term generated by the above expression is directly perturbative via two bilinear unbalanced bounds.
    The same remains true in the unbalanced case, if the gradient is not on the highest frequency factor.
    \item The gradient is on $a_{[\lambda]}$, which is also the highest frequency. 
    Then we can use Proposition~\ref{diffabound} to perturbatively replace $\nabla a_{[\lambda]}$ by $\nabla a_{\lambda}^{[1]}$ which is in $D\coeffu_{<\lambda}$.
    \item The gradient is on $g_{<\lambda}$, which is also the highest frequency. 
    Then we can use Lemma~\ref{l:bal-quad} to perturbatively discard the balanced part of 
    $g_{<\lambda}$, leaving us with $\nabla g_\lambda \in D\coeffu_{<\lambda}$.
    Here we can include $a$ as a paracoefficient, where the property $a \in \coeff$ suffices.
\end{itemize}
A version of this description applies to all terms in $N^{[1]}_{u,\lambda}$, so at the conclusion of this discussion, we have replaced $N^{[1]}_{u,\lambda}$ with an expression with unbalanced coefficients of the form
\[
N^{[1]}_{u,\lambda}= Op(D\coeffu_{<\lambda} S^{-\frac12}_\lambda) u_\lambda+
Op(D^\frac12 \coeffu_{<\lambda} S^{0}_\lambda) \bu_\lambda,
\]
which is what will be assumed in the future.

\item In $N^{[k]}_{u,\lambda}$ with $k \geq 2$. Here powers of $a$ appear multiplicatively in the paracoefficients undifferentiated, so no reduction as in the $N^{[1]}_{u, \lambda}$ case is available. However, the normal form type construction of the corrections below
is more robust, so we will simply use the $\coeff$ regularity of the paracoefficients.

\end{enumerate}
Finally, we remark that the same discussion applies adliteram in the study of the linearized equation later on.

\

\medskip

\textbf{1. The contribution of $N^{[1]}_\lambda(u)$.} By the discussion above we can assume  that $N^{[1]}_\lambda$ has two types of terms, namely 
\[
N^{[1]}_\lambda(u) = Op (D\coeffu_{<\lambda} S^{-\frac12}_\lambda) u_\lambda + Op (D^\frac12 \coeffu_{<\lambda} S^{0}_\lambda) \bu_\lambda
\]
The first term yields a source which is exactly like $R^3_{\lambda,m}$ in Theorem~\ref{t:second-df}. This yields 
a correction of the same type as there. 

The second term is larger, but contributes 
a bilinear form in $(\bu_\lambda,\bu_\lambda)$
which has a stronger nonresonance
(see the discussion in Section~\ref{subsec:resonant}),
which is why it yields a density correction of the same class; here the full phase $\omega(\xi)+\omega(\eta) \approx 2\lambda^{\frac12}$ is elliptic, so no flux correction is needed, and \eqref{bq-relation-nl} applies only to the $(u_\lambda,\bar u_\lambda)$ part.

There is an additional step here, which is to estimate the contribution of the paradifferential source term, namely the  nonlinearity $N_\lambda(u)$
to the time derivative of  $B^{NL,[1]}_{\lambda,m}$,
which is 
\[
Op(D^\frac12\coeffu_{<\lambda} S^{-\frac12}_\lambda)(u_\lambda,\overline{N_\lambda(u)})
\]
Here we can pair $D^\frac12\coeffu$ with $u_\lambda$ for an unbalanced bilinear $L^2_{t,x}$  bound of size $C^2 \epsilon^2 \lambda^{\frac14} c_\lambda^2$, therefore 
it would suffice to pair this with  an $L^2_{t,x}$ bound for $N_\lambda(u)$.
Such a bound holds for the bulk of  $N_\lambda(u)$, excluding only the balanced $N^{[2]}$ part, see \eqref{NuL2remainder}, 
\begin{equation}\label{N-lambda-2}
 \| N_\lambda(u)- N^{[2],bal}_\lambda\|_{L^2_{t,x}} \lesssim  C^2\epsilon^2 c_\lambda \lambda^{\frac14}.   
\end{equation}
The balanced component is $N_\lambda^{[2],bal}(u_\mu,u_\mu)$ corresponding to the case when $\mu \approx \lambda$ and the two input frequencies are very close. But in this case our narrower choice of the dyadic frequency localization, captured by the constant $\step$ in Section~\ref{sec:notations}, guarantees that 
the two input frequencies in $N^{[2],bal}$ are separate from the frequency of $u_\lambda$, therefore we can pair one of the $u_\mu$ factors with $u_\lambda$ in an unbalanced bilinear  $L^2_{t,x}$ bound, and the other $u_\mu$ factor with the $D^\frac12 \coeffs$ coefficient in a second unbalanced bilinear $L^2_{t,x}$ bound.

\medskip

\textbf{2. The good contribution of $N^{[2]}_\lambda(u)$.}  Here we recall that the 
$hh$ term $N^{[2]}_\lambda$ contains 
contributions of the form
\[
N^{[2]}_\lambda(u) = \sum_{\mu \gtrsim \lambda}
\lambda \mu^{\frac12} Op(\coeff_{<\lambda} S^0_\mu)(\bfu_\mu,\bar \bfu_\mu),
\]
therefore it generates a source term of the form 
\[
R^{NL,[2]}_{\lambda,m}  = \sum_{\mu > \lambda}
\lambda \mu^{\frac12} \Im \Mss_\lambda\big(Op(\coeff_{<\lambda}S^0_\mu)(\bfu_\mu,\bfu_\mu), \bar u_\lambda\big).
\]
As in the $N^{[1]}$ analysis, the case when 
the $\bfu_\mu$ entries are not $(u_\mu,\bu_\mu)$ is better, so we focus on the 
$(u_\mu,\bu_\mu)$ scenario.
This is also as in Theorem~\ref{t:second-df}, but with the role of the frequencies 
$\mu$ and $\lambda$ switched. Strictly speaking, switching roles reproduces the argument only for $\mu \gg \lambda$; for $\mu \approx \lambda$  a direct normal form computation applies, which we now describe.
The $\coeff_{<\lambda}$ paracoefficients as well as  the $\coeff_{<\lambda}$ paracoefficients in $\Mss$
play no role in this computation.
Indeed, as observed in Section~\ref{subsec:resonant},
for any choice of signs $\epsilon_j$ and balanced frequencies $\xi^j$  we have 
\[
\sum_{j= 1}^3 \epsilon_j \omega(x,\xi^j) \neq 0 
\qquad \text {on} \quad \sum_{j=1}^3 \epsilon_j \xi^j = 0.
\]
This directly yields a division relation of 
the form
\[
r^{NL,[2]}_{\lambda,m}(x,\xi^1,\xi^2,\xi^3) = i b^{NL,[2]}_{\lambda,m}(x,\xi^1,\xi^2,\xi^3) \sum_{j= 1}^3 \epsilon_j \omega(x,\xi^j)
+ i q^{NL,[3]}_{\lambda,m}(x,\xi^1,\xi^2,\xi^3) \sum_{j= 1}^3 \epsilon_j \xi^j, 
\]

Thus the corresponding normal form type 
density and flux corrections are of the type
\[
B^{NL,[2]}_{\lambda,m} =   \sum_{\mu\gtrsim \lambda} \lambda^{-\frac12} \mu^\frac12 Op( \coeff_{<\lambda} S^{0}_{\mu})(u_\mu,\bar u_\mu) u_\lambda, \qquad 
Q^{NL,[2]}_{\lambda,m} =  \sum_{\mu\gtrsim \lambda} \lambda^{-\frac12} Op( \coeff_{<\lambda} S^{0}_{\mu})(u_\mu,\bar u_\mu) u_\lambda, 
\]
canceling at leading order the $R^{NL,[2]}_{\lambda,m}$ term.
We now consider the remaining terms in the paramaterial derivative of 
$B^{NL,[2]}_{\lambda,m}$. The contribution of the paramaterial derivative of the $\coeff_{<\lambda}$ paracoefficients can be treated perturbatively 
using a bilinear $L^2$ bound as in Lemma~\ref{l:Dtg}, paired with a second unbalanced bilinear $L^2$ bound for the frequency pairing $(\lambda,\nu)$.
It remains to examine the contributions of $N_\mu(u)$, respectively $N_\lambda(u)$ to the time derivative of $B^{NL,[2]}_{\lambda,m}$,
which have the form 
\[
B^{NL,[2]}_{\lambda,m} (N_\mu(u),u_\mu)u_\lambda , \qquad B^{NL,[2]}_{\lambda,m} (u_\mu,u_\mu)  N_\lambda(u) 
\]
where $u_\mu$ may also be conjugated.
We separate this into two cases:

\medskip

a) $\mu \gg \lambda$. Here we can use again \eqref{N-lambda-2} (the part that works) combined with a bilinear $L^2_{t,x}$ bound (unbalanced in the first case, balanced in the second).  It remains to consider 
the bad part of $N^{[2]}$, namely $N^{[2],bal}$. In the first case we have three high frequencies interacting for a 
frequency $\lambda$ output, therefore two of them must be transversal, and then we can use an  unbalanced bilinear $L^2_{t,x}$ for a $(\mu,\mu)$  pairing and an unbalanced bilinear $L^2$ for a $(\lambda,\mu)$  pairing.
In the second case we can instead use two unbalanced bilinear $L^2$ for $(\lambda,\mu)$  pairings.

\medskip

b) $\mu \approx \lambda$. In this case the bound \eqref{N-lambda-2} no longer suffices even when it works, so we need to examine in greater detail the components of $N_\lambda$  (the $N_\mu$ contribution is identical since $\mu \approx \lambda$).

\smallskip

b1) The $N_\lambda^{[1]}$ term yields a cubic $u_\lambda$ form which is both lower order (an extra half derivative on the lower frequency para coefficient) and nonresonant and can be corrected. 

\smallskip

b2) The $N_{\lambda}^{[2]}$ term with a high frequency $\nu$ yields a quartic form with two $\lambda$ 
and two $\nu$ frequencies, which can be estimated with two bilinear $L^2_{t,x}$ bounds provided that  $\nu \gg \lambda$. \emph{We set aside the case when  $\nu \approx \lambda$, which may produce quartic doubly resonant interactions. This case will be discussed separately in part 5 of the proof below.} 

\smallskip

b3) The $N_{\lambda}^{[3]}$ term yields a quintic source term which is either transversal or nonresonant.
This is similar to the direct contribution  of $N_{\lambda}^{[4]}$,
which is discussed in Step 4 of the proof.

\smallskip

b4) The $N_\lambda^{[4]}$ term yields a six-linear balanced contribution, which can be directly estimated
using Strichartz bounds.

\smallskip

b5)  The $F_\lambda$ contribution. This is estimated directly using \eqref{u-remainder-def}.

\medskip

\textbf{3. The good contribution of $N^{[3]}_\lambda(u)$.}
The associated error has the form
\[
R^{NL,[3]}_{\lambda,m} =\lambda^{\frac52} Op(\coeff_{<\lambda} S^{0}_\lambda)(\bfu_\lambda,\bfu_\lambda,\bfu_\lambda,\bfu_\lambda).
\]
Here we distinguish two components,  a resonant one, with exactly two $\bar u_\lambda$, and a nonresonant one. We consider the nonresonant case for now, \emph{setting aside the doubly resonant case}.

For the nonresonant component, we claim that a direct normal form analysis yields a correction
\begin{equation}\label{B3-correction}
B^{NL,[3]}_{\lambda,m} =\lambda^{2} Op(\coeff_{<\lambda} S^{0}_\lambda)(\bfu_\lambda,\bfu_\lambda,\bfu_\lambda,\bfu_\lambda),
\quad
Q^{NL,[3]}_{\lambda,m} =\lambda^{\frac32} Op(\coeff_{<\lambda} S^{0}_\lambda)(\bfu_\lambda,\bfu_\lambda,\bfu_\lambda,\bfu_\lambda).
\end{equation}
To see this, we start from the observation in Section~\ref{subsec:resonant} that 
for any choice of signs $\epsilon_j$ without phase rotation symmetry we have 
\[
\sum_{j= 1}^4 \epsilon_j \omega(x,\xi^j) \neq 0 
\qquad \text {on} \quad \sum_{j=1}^4 \epsilon_j \xi^j = 0.
\]
This directly yields a division relation of 
the form
\[
r^{NL,[3]}_{\lambda,m}(x,\xi^1,\xi^2,\xi^3,\xi^4) = i b^{NL,[3]}_{\lambda,m}(x,\xi^1,\xi^2,\xi^3,\xi^4) \sum_{j= 1}^4 \epsilon_j \omega(x,\xi^j)
+ i q^{NL,[3]}_{\lambda,m}(x,\xi^1,\xi^2,\xi^3,\xi^4) \sum_{j= 1}^4 \epsilon_j \xi^j, \]
with symbols $r^{NL,[3]}_{\lambda,m}$ and $
q^{NL,[3]}_{\lambda,m}$ as in \eqref{B3-correction}.

Adding these density and flux corrections
directly cancels the $R^{NL,[3]}_{\lambda,m}$
source term, at the expense of additional source terms of the form
\[
\lambda^{2} Op(T_{D_t} \coeff_{<\lambda} S^{0}_\lambda)(\bfu_\lambda,\bfu_\lambda,\bfu_\lambda,\bfu_\lambda),
\]
as well as 
\[
\lambda^{2} Op(\coeff_{<\lambda} S^{0}_\lambda)(\bfu_\lambda,\bfu_\lambda,\bfu_\lambda,{\bf N}_\lambda(u))
\]
and 
\[
\lambda^{2} Op(\coeff_{<\lambda} S^{0}_\lambda)(\bfu_\lambda,\bfu_\lambda,\bfu_\lambda,(\nabla \cdot V_{[\lambda]} - \frac12 \vartheta_{[\lambda]}(x,D)) \bfu_\lambda).
\]
The first is dealt with using an unbalanced bilinear $L^2_{t,x}$ as in Lemma~\ref{l:Dtg} and three $L^6_{t,x}$ bounds. For the second we can use the $L^2_{t,x}$ bound \eqref{N-lambda-2} for $N_\lambda - N_{\lambda}^{[2],bal}$. The contribution of $N_{\lambda}^{[2],bal}$ yields a balanced quintic form, which is simply transferred to the next step.
Finally fot the third we can also use a bilinear unbalanced $L^2$ bound for the last argument, combined with three $L^6_{t,x}$ bounds.
 
\medskip

\textbf{4. The contribution of $N^{[4]}_\lambda(u)$.} 
The associated error has the form
\[
R^{NL,[4]}_{\lambda,m} =\lambda^{\frac72} Op(\coeff_{<\lambda} S^{0}_\lambda)(\bfu_\lambda,\bfu_\lambda,\bfu_\lambda,\bfu_\lambda,\bfu_\lambda).
\]
where any  $u_\lambda$ may be conjugated.
Here we either have a transversal interaction (i.e. two separated frequencies), in which case the above expression can be estimated perturbatively
using a bilinear $L^2$ bound and three 
Strichartz bounds, or else the output frequency has size $O(\lambda)$, in which case  a flux correction suffices, 
\[
Q^{NL,[4]}_{\lambda,m} = Op(\coeff_{<\lambda} S^{\frac{5}{2}}_{\lambda})( \bfu_\lambda,\bfu_\lambda,\bfu_\lambda,\bfu_\lambda, \bfu_\lambda).
\]

\medskip

\textbf{5. The doubly resonant contributions.} These are quartic doubly resonant expressions in $u_\lambda$
of the form 
\[
Op(\coeff_{<\lambda} S^\frac52_\lambda)(u_\lambda, \bu_\lambda,u_\lambda,\bu_\lambda),
\]
which for the principal part have phase rotation symmetry,
and cannot be directly estimated without additional structural properties. To estimate these terms
we will carry out two steps:
\begin{enumerate}[label=(\roman*)]
\item Exactly compute their principal symbol, and show that it satisfies the doubly null condition, 
i.e. it vanishes to second order on the diagonal $\xi^1 = \xi^2 = \xi^3 = \xi^4$ within the set $\Delta^4 \xi = 0$.

\item Estimate the quartic forms with double null condition using two balanced bilinear $L^2_{t,x}$ bounds,
modulo a flux term.
\end{enumerate}

We carefully remark that lower order terms, e.g., of the form $Op(D^\frac12 \coeff_{<\lambda} S^2_\lambda)(u_\lambda, \bu_\lambda,u_\lambda,\bu_\lambda)$  or 
$Op (T_{D_t} \coeff_{<\lambda} S^2_\lambda)(u_\lambda, \bu_\lambda,u_\lambda,\bu_\lambda)$ 
 can be estimated directly using a bilinear unbalanced bound and three $L^6_{t,x}$ estimates. This in particular indicates that the double null condition is only  needed at the level of the principal symbol, irrespective of the quantization. 
We now consider the steps (i),(ii) noted above.
\medskip

\emph{(i) Show that the double  null condition holds.}
We recall that doubly resonant contributions
can only arise in two cases:

\begin{enumerate}[label= (\alph*)]
    \item From the iterated balanced quadratic
    terms in $N^{[2]}_\lambda$.

    \item From the cubic terms in $N^{[3]}_\lambda$ with phase rotation symmetry, 
    i.e. with arguments $(u,\bu,u)$.
\end{enumerate}
We will show that, as heuristically indicated by the normal form computation in Section~\ref{sec:layered}, these two types of contributions satisfy the double null condition separately.
To begin with, we recall that 
\[
u =\eta + i T_\mfq \varphi, \qquad N_{u}^{[k]} = N_\eta^{[k]} + i T_\mfq N_\varphi^{[k]}.
\]
So at leading order we have 
\[
\eta = \frac12(u+\bu), \qquad \varphi = - \frac{i}2  T_{\mfq^{-1}}(u - \bu).
\]

\bigskip

\noindent \emph{(i)(a) The iterated $N^{[2]}$ contribution.}  

Here we first recall the expressions for $
N^{[2]}_{\eta,\lambda}$ and $N^{[2]}_{\varphi,\lambda}$: 
\begin{equation}
\left\{
\begin{aligned}
N^{[2]}_{\eta,\lambda} = & \  - \DN_{[\lambda]}(\eta_{\gtrsim \lambda}  \L_{[\lambda]} \varphi_{\gtrsim \lambda})
- \nabla \cdot ( \eta_{\gtrsim\lambda}  \M_{[\lambda]} \varphi_{\gtrsim\lambda})
\\
N^{[2]}_{\varphi,\lambda} = & \ \frac12 \L_{[\lambda]} \varphi_{\gtrsim \lambda} \DN_{[\lambda]} \varphi_{\gtrsim \lambda} - \frac12 \M_{[\lambda]} \varphi_{\gtrsim \lambda} \cdot \nabla \varphi_{\gtrsim \lambda}.
\end{aligned}
\right.
\end{equation}
We use these to compute the principal symbol for $N^{[2]}_u$, which we separate into 
\[
N^{[2]}_u = N^{[2],++}(u,u) +N^{[2],+-}(u,\bu) +N^{[2],--}(\bu,\bu).
\]
Here a direct computation yields the principal symbols
\[
\begin{aligned}
4n^{[2],++}(\xi,\eta) = & i \mfq^{-1}(\eta) (\dn(\xi+\eta) \mfl(\eta) 
+ i(\xi+\eta) \mfm(\eta)) + \frac{i}2 \mfq(\xi+\eta) \mfq^{-1}(\xi) \mfq^{-1}(\eta)(\mfl(\xi) \dn(\eta) - i\mfm(\xi) \cdot \eta) \\
4n^{[2],--}(\xi,\eta) = & i \mfq^{-1}(\eta) (\dn(\xi+\eta) \bar\mfl(\eta) 
- i(\xi+\eta) \bar\mfm(\eta)) + \frac{i}2 \mfq(\xi+\eta) \mfq^{-1}(\xi) \mfq^{-1}(\eta)(\bar \mfl(\xi) \dn(\eta) + i\bar \mfm(\xi) \cdot \eta) \\
 4n^{[2],+-}(\xi,\eta)= & \  i \mfq^{-1}(\xi) (\dn(\xi-\eta) \mfl(\xi) + i(\xi-\eta) \mfm(\xi))  - i \mfq^{-1}(\eta) (\dn(\xi-\eta) \bar \mfl(\eta) + i(\xi-\eta) \bar\mfm(\eta)) 
\\
& \ 
+ \frac{i}{2}\mfq(\xi-\eta) \mfq^{-1}(\xi)\mfq^{-1}(\eta)( \mfl(\xi) \dn(\eta) + \bar \mfl(\eta) \dn(\xi) +i \mfm(\xi) \cdot \eta - i \bar \mfm(\eta) \cdot \xi), 
\end{aligned}
\]
where the first two relations should also be symmetrized in $\xi$, $\eta$.

Next we simplify the above expressions modulo 
terms that vanish quadratically. We will need such expansions near $\xi = \eta$ and near $\xi = -2\eta$ for $n^{[2],++}$
and $n^{[2],--}$, and near $\xi = 2\eta$ for $n^{[2],+-}$.
We begin with the expansion near $\xi = \eta$.
Starting from \eqref{LvsM}, we obtain the relations
\[
\dn(\xi) \mfl(\xi) + i \xi \cdot \mfm(\xi) = 0, \qquad \dn(\xi) \bar \mfl(-\xi) - i \xi \cdot \bar \mfm(-\xi) = 0.
\]

As a slight extension of these, it follows that\footnote{This is equivalent to the relation $|\xi||\eta| \approx \pm \xi \cdot \eta$ near $\xi \in \R^{\pm} \eta$, relative to the euclidean metric on the free surface.} 
\[
\begin{aligned}
&\dn(\xi) \mfl(\eta) + i \xi \cdot \mfm(\eta) \approx 0 \qquad \text{ near }  \xi \in \R^+ \eta  \eta\\
&\dn(\xi) \bar \mfl(\eta) - i \xi \cdot \bar \mfm(\eta) \approx  0  \qquad \text{ near }  \xi \in \R^- \eta.
\end{aligned}
\]
This yields, after symmetrization,
\begin{equation}\label{n2ppmm}
\left\{ 
\begin{aligned}
 8n^{[2],++}(\xi,\eta) \approx & \   i \mfq(\xi+\eta) \mfq^{-1}(\xi) \mfq^{-1}(\eta)(\mfl(\xi) \dn(\eta) 
 + \mfl(\eta) \dn(\xi))
 \\
  8n^{[2],--}(\xi,\eta)  \approx & \   i\mfq(\xi+\eta) \mfq^{-1}(\xi) \mfq^{-1}(\eta)(\bar \mfl(\xi) \dn(\eta)+ \bar \mfl(\eta) \dn(\xi))
\end{aligned}
\right.
\qquad \text{near } \xi = \eta,
\end{equation}
and similarly
\begin{equation}
\left\{ 
\begin{aligned}
 4n^{[2],++}(\xi,\eta) \approx & \   i \mfq^{-1}(\eta) \dn(\xi+\eta) \mfl(\eta) 
 \\
 4n^{[2],--}(\xi,\eta)  \approx & \  i \mfq^{-1}(\eta) \dn(\xi+\eta) \bar\mfl(\eta) 
\end{aligned}
\right.
\qquad \text{near } \xi = -2\eta,
\end{equation}
with the mixed sign representation as follows
\begin{equation}
\begin{aligned}
4n^{[2],+-}(\xi,\eta)\approx & \ - i \mfq^{-1}(\eta) \dn(\xi-\eta) \bar \mfl(\eta)
\end{aligned}
\qquad \text{near } \xi = 2\eta.
\end{equation}

The corresponding source terms in the 
density-flux identities are $2 \Re R^{ab}$
where
\[
R^{ab} = \Ms(N^{[2],ab},\bar u) 
\]
with symbols whose values on the diagonal
$a \xi+ b \eta -\zeta = 0$ are
\[
r^{ab}(\xi,\eta,\zeta) = \Ms(a\xi+b\eta,\zeta) n^{[2],ab}(\xi,\eta).
\]
The associated density correction $B^{ab}$ 
with signature $(a,b,-)$ will then have 
the symbol
\[
b^{ab}(\xi,\eta,\zeta) = \frac{i r^{ab}(\xi,\eta,\zeta)}{a\omega(\xi) + b \omega(\eta) 
- \omega(\zeta)}.
\]
We are interested in this symbol on the diagonal $a \xi+ b \eta -\zeta = 0$ where 
all three frequencies are of comparable 
size $\lambda$. This insures that the denominator is uniformly nonzero, precisely of size $\lambda^\frac12$, and the symbol is equal to 
\[
b^{ab}(\xi,\eta,\zeta) = \frac{ i \ms(\zeta,\zeta) n^{[2],ab}(\xi,\eta,\zeta)}{a\omega(\xi) + b \omega(\eta) 
- \omega(\zeta)}
\]
Finally we consider the quartic terms arising
from the time derivative of $B^{ab}$ via the balanced $N^{[2]}$ component of the nonlinearity. There are three cases to consider, depending on whether we bifurcate 
the first, second or third entry.

The easiest case is when the third component is time differentiated, in which case we obtain  the quartic output $R^{3,abcd}$ with symbol
\[
r^{abcd} (\xi,\eta, \gamma, \sigma)  \; = b^{ab}(\xi,\eta,\zeta)\,\overline{n^{[2],cd}(\gamma,\sigma)}, \qquad \zeta = c \gamma + d \sigma .
\]
This has 
signature $(a,b,-c,-d)$, has
phase rotation symmetry when 
$a+b-c-d = 0$,  and we are interested in the case 
when all four frequencies are nearby.
But then the requirement that $\zeta$ stays away from zero forces matching signs 
$a=b=c=d$, which fixes the choice of $r$
on the diagonal to 
\[
r^{++++} (\xi,\eta, \gamma, \sigma)
= \frac{i \ms(\zeta,\zeta) n^{[2],++}(\xi,\eta) \overline{n^{[2],++}(\gamma,\sigma)}}
{\omega(\xi)+ \omega(\eta) - \omega(\zeta)},
\]
and the similar expression with four $-$ signs,
\[
r^{----} (\xi,\eta, \gamma, \sigma)
= \frac{  i \ms(\zeta,\zeta) n^{[2],--}(\xi,\eta) \overline{n^{[2],--}(\gamma,\sigma)}}
{-\omega(\xi) - \omega(\eta) - \omega(\zeta)}.
\]
In $r^{++++}$ the frequencies $\xi$ and $\eta$ correspond to $u$ factors and the frequencies $\gamma$ and $\sigma$ correspond to $\bar u$ factors. Relabeling these frequencies  $\xi^1,\xi^3$, respectively $\xi^2,\xi^4$ and substituting the formulas in \eqref{n2ppmm}, we arrive at   $\xi^1,\xi^3$, respectively $\xi^2,\xi^4$ and substituting the formulas in \eqref{n2ppmm}, we arrive at 
\[
r^{++++} (\xi^1,\xi^2, \xi^3, \xi^4)
= \frac{i\, \ms(\zeta,\zeta)\, \mfq^2(\zeta) \prod_{j=1}^{4} \mfq^{-1}(\xi^j)}
{16(\omega(\xi^1)+ \omega(\xi^3) - \omega(\zeta))}
\,\big(\mfl(\xi^1) \dn(\xi^3) + \mfl(\xi^3) \dn(\xi^1)\big)
\big(\bar \mfl(\xi^2) \dn(\xi^4) + \bar \mfl(\xi^4) \dn(\xi^2)\big),
\]
where on the diagonal $\Delta^4 \xi = 0$ we have $\zeta = \xi^1+\xi^3 = \xi^2+\xi^4$.
Symmetrizing between $(\xi^1,\xi^3)$ and $(\xi^2,\xi^4)$, the numerator simply changes sign due 
to the $i$ factor. The denominators, on the other hand agree up to quadratic factors, as 
\[
\omega(\xi^1)+ \omega(\xi^3) \approx \omega(\xi^2)+ \omega(\xi^4) \qquad \text{ on } \Delta^4 \xi = 0.
\]
Thus the double  null condition follows for the symmetrization of $r^{++++}$. The case of $r^{----}$ 
is similar.

\bigskip

Consider now the case when we bifurcate the frequency $\xi$ entry via $N^{[2]}$, 
$\xi = c \gamma + d \sigma$. Then we obtain the choice of $r$ on the diagonal as 
\[
r^{abcd} (\gamma, \sigma,\eta,\zeta)
= \frac{i\,\ms(\zeta,\zeta) n^{[2],cd}(\gamma,\sigma) n^{[2],ab}(\xi,\eta)}
{a\omega(\xi)+ b\omega(\eta) - \omega(\zeta)}.
\]
The signature is now $(ac,ad,b,-)$, and the phase rotation symmetry requires $ac+ad+b-1=0$. The fact that all four frequencies should be close but with $\xi$ 
away from zero forces $c = d$ and in turn $b = -$. We are left with two cases, 
$(a,b,c,d) = (+,-,+,+)$ and $(a,b,c,d) = (-,-,-,-)$. These cases are  not as immediate as the prior ones.

Consider the case $(a,b,c,d) = (+,-,+,+)$. Then 
\[
r^{+-++}(\eta,\zeta,\gamma,\sigma) = \frac{i\ms(\zeta,\zeta) n^{[2],++}(\gamma,\sigma) n^{[2],+-}(\xi,\eta)}
{\omega(\xi)-\omega(\eta) - \omega(\zeta)},
\qquad \xi = \gamma + \sigma = \eta+\zeta.
\]
This gives
\[
\begin{aligned}
r^{+-++}(\eta,\zeta,\gamma,\sigma) = & \ \frac{ i \ms(\zeta,\zeta)    \mfq(\gamma+\sigma) \mfq^{-1}(\gamma) \mfq^{-1}(\sigma)(\mfl(\gamma) \dn(\sigma) 
 + \mfl(\gamma) \dn(\sigma)) \mfq^{-1}(\eta) \dn(\zeta) \bar\mfl(\eta)}
{\omega(\xi)-\omega(\eta) - \omega(\zeta)}
\\
= & \ \frac{\mfq(\gamma+\sigma) \mfq^{-1}(\gamma) \mfq^{-1}(\sigma)\mfq^{-1}(\eta)
\mfq^{-1}(\zeta)}{\omega(\xi)-\omega(\eta) - \omega(\zeta)} 
i (\mfl(\gamma) \dn(\sigma) 
 + \mfl(\gamma) \dn(\sigma))  \ell(\zeta)\dn(\zeta) \bar\mfl(\eta),
\end{aligned}
\]
where $\ell$ denotes the real valued function
\[
\ell(\zeta) = \ms(\zeta,\zeta) \mfq(\zeta).
\]
Here the fraction is real and symmetric modulo quadratic terms, so it remains to symmetrize
the remaining expression. Symmetrization in 
$(\eta,\zeta)$ gives 
\[
\begin{aligned}
\ell(\zeta) \dn(\zeta) 
\mfl(\eta) + \ell(\eta) \dn(\eta) 
\mfl(\zeta) = & \ \frac12 ( \ell(\zeta)+  \ell(\eta))( \dn(\zeta) 
\bar \mfl(\eta) + \dn(\eta) 
\bar \mfl(\zeta))
 \\ & + \frac12 ( \ell(\zeta)-  \ell(\eta))( \dn(\zeta) 
\bar \mfl(\eta) - \dn(\eta) 
\bar \mfl(\zeta))
\\
\approx & \ \frac12 ( \ell(\zeta)+  \ell(\eta))( \dn(\zeta) 
\bar \mfl(\eta) + \dn(\eta) 
\bar \mfl(\zeta)).
\end{aligned}
\]
Inserting this into the earlier expression 
and symmetrizing with respect to the pairs 
$(\eta,\zeta)$ and $(\gamma,\sigma)$
finally leads to the expression
\[
i (( \ell(\zeta)+  \ell(\eta) -\ell(\gamma) -\ell(\sigma))( \dn(\zeta) 
\bar \mfl(\eta) + \dn(\eta) 
\bar \mfl(\zeta))(\mfl(\gamma) \dn(\sigma) 
 + \mfl(\gamma) \dn(\sigma)),
\]
where the symmetric difference in the first factor vanishes of second order on $\Delta^4 \xi = 0$.

Lastly, consider the case $(a,b,c,d) = (-,-,-,-)$. Then 
\[
r^{----}(\eta,\zeta,\gamma,\sigma) = \frac{i \ms(\zeta,\zeta) \overline{n^{[2],--}}(\gamma,\sigma) n^{[2],--}(\xi,\eta)}
{-\omega(\xi)-\omega(\eta) - \omega(\zeta)},
\qquad -\xi = \gamma + \sigma = \eta+\zeta.
\]
Substituting,
\[
\begin{aligned}
r^{----}(\eta,\zeta,\gamma,\sigma) = \frac{i \ms(\zeta,\zeta) \mfq(\gamma+\sigma) \mfq^{-1}(\gamma) \mfq^{-1}(\sigma)( \mfl(\gamma) \dn(\sigma)+  \mfl(\sigma) \dn(\gamma)) \mfq^{-1}(\eta) \dn(\zeta) \bar\mfl(\eta)}
{-\omega(\xi)-\omega(\eta) - \omega(\zeta)},
\end{aligned}
\]
and the symmetrization is now carried out exactly as in the previous case, to obtain the second order vanishing on the diagonal.
\bigskip

\noindent \emph{(i)(b) The  $N^{[3]}$ contribution.}  We recall that 
\begin{equation}
\left\{
\begin{aligned}
\ N^{[3]}_{\eta,\lambda} = & \ 
\L_{[\lambda]}^{*} (\eta_\lambda \DN_{[\lambda]}(\eta_\lambda \L_{[\lambda]}\varphi_\lambda))
- \frac12 \L_{[\lambda]}^{*}( \eta_\lambda^2 \DN_{[\lambda]}\L_{[\lambda]} \varphi_\lambda) - \frac12 \L_{[\lambda]}^{*} \DN_{[\lambda]} (\eta_\lambda^2 \L_{[\lambda]} \varphi_\lambda)
\\
\  N^{[3]}_{\varphi,\lambda} = & \ - \L_{[\lambda]}\varphi \DN_{[\lambda]}(\eta_\lambda \L_{[\lambda]}\varphi) + \L_{[\lambda]}\varphi \eta_\lambda \DN_{[\lambda]} \L_{[\lambda]}\varphi.
\end{aligned}
\right.
\end{equation}

We are now able to compute the principal symbol of the resonant component $N_u^{[3],res}$ of $N_u^{[3]}$, which contains all the terms with the $(u,\bu,u)$
entries, which we will use in this order.
For clarity we separate $N_u^{[3],res}$
into two components,
\[
N_u^{[3],res}= N_u^{[3],res,++}+ N_u^{[3],res,+-}
\]
depending on whether we select both $u$ entries in $\eta$ in the first equation and in $\varphi$ in the second, or the $(u,\bu)$
entries.

\medskip

\textbf{ The $++$ component of $N^{[3]}$.}   Denoting the $u$ frequencies by $\xi^1,\xi^3$
and the $\bar u$ frequency by $\xi^2$,
and
$\xi^4 = \xi^1-\xi^2+\xi^3$,  we have 
\[
\begin{aligned}
8 N_u^{[3],res,++} =&\  \mfq^{-1}(\xi^2) \bar{\mfl}(\xi^2) \bar\mfl(\xi^4)(\dn(\xi^3-\xi^2)
- 
\frac12\dn(\xi^2)
- \frac{1}2\dn(\xi^4)) 
\\ 
& + \mfq(\xi^4) \mfq^{-1}(\xi^1)\mfq^{-1}(\xi^3) \mfl(\xi^1)  \mfl(\xi^3) (- \dn(\xi^2-\xi^3)+ \dn(\xi^3)).
\end{aligned}
\]
Noting that $\mfq(\xi^2) \mfq(\xi^4) - \mfq(\xi^1) \mfq(\xi^3)$ vanishes to second order on the diagonal,
after some simple algebraic manipulations and symmetrization in $\xi^1,\xi^3$  we obtain, modulo good terms, 
\[
\begin{aligned}
16  N_u^{[3],res,++}(\xi^1,\xi^2,\xi^3) \approx &\ \mfq^{-1}(\xi^2) \left[ \bar \mfl(\xi^2) \bar\mfl(\xi^4) \left( \dn(\xi^3-\xi^2) +\dn(\xi^1-\xi^2) -
\dn(\xi^2) - \dn(\xi^4) \right) \right. \\
& \ 
\left. \qquad- \mfl(\xi^1)  \mfl(\xi^3) \left( \dn(\xi^3-\xi^2)+ \dn(\xi^1-\xi^2) -\dn(\xi^1) -\dn(\xi^3) \right) \right].
\end{aligned}
\]
The contribution of this expression to the multilinear form $2\Re \Ms_\lambda(N^{[3]}_\lambda,\bar u)$
will have the principal symbol on the diagonal $\Delta^4 \xi = 0$ equal to
\[
\ms_\lambda(\xi^4,\xi^4) N_u^{[3],res,++}(\xi^1,\xi^2,\xi^3).
\]
Here we note that on $\Delta^4 \xi = 0$ the functions $\dn(\xi^2) + \dn(\xi^4)$ and $\dn(\xi^1) + \dn(\xi^3)$ agree to second order, so we can identify them, arriving at the equivalent expression
\[
\ms_\lambda(\xi^4,\xi^4) \mfq^{-1}(\xi^2)   (\bar \mfl(\xi^2) \bar\mfl(\xi^4) -  \mfl(\xi^1)  \mfl(\xi^3)) ( \dn(\xi^3-\xi^2)+ \dn(\xi^1-\xi^2) -2\dn(\xi^{avg}) ).
\]
The last factor is already fully symmetric on $\Delta^4 \xi = 0$, while the middle factor has antisymmetry
with respect to the odd/even switch. Symmetrizing the full expression, we arrive at
\[
(\ell(\xi^2,\xi^4)-\ell(\xi^1,\xi^3))
(\bar \mfl(\xi^2) \bar\mfl(\xi^4) -  \mfl(\xi^1)  \mfl(\xi^3))( \dn(\xi^3-\xi^2)+ \dn(\xi^1-\xi^2) -2\dn(\xi^{avg}) )
\]
with $\ell$ defined as the real valued symmetric function
\[
\ell(\xi,\eta) = \ms_\lambda(\xi,\xi) \mfq^{-1}(\eta) + \ms_\lambda(\eta,\eta) \mfq^{-1}(\xi).
\]
Here the symmetry of $\ell$ guarantees that the first factor has the desired second order vanishing, as needed.

\medskip

\textbf{ The $+-$ component of $N^{[3]}$.}  
Using the same notations as in the first case, we have the principal symbol
\[
\begin{aligned}
8 N_u^{[3],res,+-} =&\  \mfq^{-1}(\xi^3) {\mfl}(\xi^3) \bar \mfl(\xi^4)(
(\dn(\xi^3-\xi^2)+\dn(\xi^1+\xi^3)  
- \dn(\xi^3) - \dn(\xi^4))
\\ 
& + \mfq(\xi^4) \mfq^{-1}(\xi^1)\mfq^{-1}(\xi^2) \mfl(\xi^1) \bar \mfl(\xi^2)
(- \dn(\xi^2-\xi^3)-\dn(\xi^1+\xi^3)+ \dn(\xi^1) + \dn(\xi^2)).
\end{aligned}
\]
Here we remark that, modulo quadratic terms in the difference
we have
\[
\dn(\xi^1+\xi^3) \approx \dn(\xi^1) + \dn(\xi^3),
\]
therefore
\[
\begin{aligned}
8 N_u^{[3],res,+-} \approx &\  \mfq^{-1}(\xi^3) {\mfl}(\xi^3) \bar \mfl(\xi^4)(
(\dn(\xi^3-\xi^2)+\dn(\xi^1)  
 - \dn(\xi^4))
\\ 
& + \mfq(\xi^4) \mfq^{-1}(\xi^1)\mfq^{-1}(\xi^2) \mfl(\xi^1) \bar \mfl(\xi^2)
(- \dn(\xi^2-\xi^3)-\dn(\xi^3) + \dn(\xi^2)).
\end{aligned}
\]
Now the last factors in both terms vanish 
when all frequencies are equal, therefore 
we can again identify both the $\mfq$'s and 
the $\mfl$'s to
write
\[
\begin{aligned}
8 N_u^{[3],res,+-} \approx &\  \mfq^{-1}(\xi^4) |{\mfl}(\xi^{avg})|^2 
(\dn(\xi^1)  
 - \dn(\xi^4)-\dn(\xi^3) + \dn(\xi^2)).
\end{aligned}
\]
Symmetrizing in $(\xi^1,\xi^3)$ immediately 
yields 
\[
\begin{aligned}
8 N_u^{[3],res,+-} \approx &\  \mfq^{-1}(\xi^4) |{\mfl}(\xi^{avg})|^2 
(  
 - \dn(\xi^4) + \dn(\xi^2)).
\end{aligned}
\]
Inserting this expression in the multilinear form $2\Re \Ms_\lambda(N^{[3]}_\lambda,\bar u)$
will have the principal symbol on the diagonal $\Delta^4 \xi = 0$ equal to
\[
8\ms_\lambda(\xi^4,\xi^4) N_u^{[3],res,+-}(\xi^1,\xi^2,\xi^3) \approx  \ell(\xi^4) |{\mfl}(\xi^{avg})|^2 
( \dn(\xi^2) - \dn(\xi^4)),
\] 
where $\ell(\xi) = \ms(\xi,\xi) \mfq^{-1}(\xi)$. Finally, symmetrization in $(\xi^2,\xi^4)$ yields
\[
\approx (\ell(\xi^4)-\ell(\xi^2) |{\mfl}(\xi^{avg})|^2 (
(  
 \dn(\xi^2) - \dn(\xi^4))
\]
with the desired second order vanishing on the diagonal.

\bigskip

\emph{(ii) $L^1$ bound for balanced quartic forms with double null condition.} 
We summarize the content of this step in the following
\begin{lemma}
Consider a symmetric quartic form $R^4_\lambda \in Op \coeff S^{\frac52}_\lambda$ with phase rotation symmetry whose symbol $r^4_\lambda(x,\xi^1,\xi^2,\xi^3,\xi^4)$ is localized at frequency $\lambda$ and vanishes 
to second order on the diagonal $\xi^1 = \xi^2 = \xi^3 = \xi^4$. Then, under our bootstrap assumptions on $u$ we can represent 
\begin{equation}
  R^4_\lambda(u_\lambda, \bu_\lambda,u_\lambda,\bu_\lambda)  
= R^{4,err}_\lambda (u_\lambda, \bu_\lambda,u_\lambda,\bu_\lambda)
+ \partial_x Q^{4,err}_\lambda (u_\lambda, \bu_\lambda,u_\lambda,\bu_\lambda)
\end{equation}
with $q^{4,err} \in \coeff S^{\frac52}_\lambda$  and $R^{4,err}_\lambda$ satisfying the bound 
\begin{equation}
    \| R^{4, err}_\lambda(u_\lambda, \bu_\lambda,u_\lambda,\bu_\lambda) \|_{L^1} \lesssim C^4 \epsilon^4 c_\lambda^4 \lambda^{-2s_c}.
\end{equation}

\end{lemma}
\begin{proof}
We claim that there exists a representation 
for $r^4_\lambda$ of the form 
\begin{equation}
r^{4}_\lambda = [(\xi^1- \xi^2)(\xi^3-\xi^4)+
(\xi^1-\xi^4)(\xi^3-\xi^2)]
p^{4,1}_\lambda + (\xi^1-\xi^3)(\xi^2-\xi^4)
p^{4,2}_\lambda + \Delta^4 \xi q^4_\lambda,
\end{equation}
where
\[
 p^{4,1}, p^{4,2} \in \coeff S^{\frac12}_\lambda, 
\quad q^4_\lambda \in  \coeff S^{\frac32}_\lambda.
\]
Then we can estimate  the contribution of the first three terms in $L^1_{x,t}$ using two bilinear balanced $L^2$ bounds. On the other hand for the last term we commute the derivative outside,
\[
Op( \Delta^4 \xi q^4_\lambda) = i \partial_x
Q^4_\lambda + Op (-i \partial_x q^4_\lambda).
\]
The first term is the desired flux term, whereas $\partial_x q^4_\lambda \in D\coeff S^{\frac32}_\lambda$, and its output can be estimated using a bilinear unbalanced $L^2$
bound and three $L^6$ bounds. 

It remains to prove the above claim. It suffices to establish the desired representation on $\Delta^4 \xi = 0$.
We consider a finite Taylor expansion of $r^4_\lambda$  around the diagonal. The terms which are at least cubic in the differences
can be easily represented in the desired form, so it remains to look at the quadratic terms.
Due to the symmetries of $r$, it suffices to consider the following terms:
\begin{itemize}
    \item $(\xi^1-\xi^3)^2$. After symmetrization this gives on $\Delta^4\xi = 0$
 \[
(\xi^1-\xi^3)^2+ (\xi^2-\xi^4)^2
= - 2 [(\xi^1- \xi^2)(\xi^3-\xi^4)+
(\xi^1-\xi^4)(\xi^3-\xi^2)]
 \]
as needed.

 \item $(\xi^1-\xi^2)(\xi^1-\xi^3)$. After symmetrization this gives on $\Delta^4\xi = 0$
 \[
(\xi^1-\xi^2)(\xi^1-\xi^3)+(\xi^1-\xi^4)(\xi^1-\xi^3)
=  (\xi^1- \xi^3)^2
 \]
 as considered before.
\end{itemize}
The proof of the lemma is now concluded.
\end{proof}
This concludes the proof of the theorem.
\end{proof}

\subsection{The linearized case} 
Here we consider density-flux identities for the linearized equation, for which we use 
the form \eqref{u-para-lin-re}, which we recall here:
\begin{equation}
\label{u-para-lin}
  (\partial_t + V_{[\lambda]} \cdot \nabla
  + \frac12 (\nabla \cdot V)-\bb(x,D)
   +
i \omega_{[\lambda]}) w_\lambda  = N^{lin}_\lambda w,
\end{equation}
where $w$ will denote the linearized variable throughout this subsection. 
For the source term on the right we use the expansion 
\begin{equation}
  N_\lambda^{lin}  =N_\lambda^{lin,[0]} + N_\lambda^{lin,[1]}+N_\lambda^{lin,[2]} + F_\lambda^{lin},
\end{equation}
where the source terms $N_\lambda^{lin,[j]}$
are described in  Section~\ref{subsec:wfinal}. We briefly recall 
the description of the three terms:

\begin{itemize}
    \item $N_\lambda^{lin,[1]}w$ is the linear ``lh" paradifferential term, of the form
\[
N_\lambda^{lin,[1]}w = Op D \coeffs_{<\lambda} S^{-\frac12}_\lambda w_\lambda
+ Op  D^\frac12 \coeffs_{<\lambda} S^{0}_\lambda \bw_\lambda,
\]
and is identical to the paradifferential operator in $N^{[1]}_{\lambda}$.

\item $N_\lambda^{lin,[2]}w$ is the ``hh" term,
which can be represented as a quadratic form with a paracoefficient,
\[
N_\lambda^{lin,[2]}w = \lambda  Op  \coeff_{<\lambda} S_{\lambda}^{\frac12} (\bfw_{\gtrsim \lambda}, \bfu_{\gtrsim \lambda}),
\]
and represents a (nonsymmetric) polarization of $N^{[2]}_\lambda$. Here we remark that the ``hh" contributions with a lower frequency on $v$ are directly perturbative.

\item $N_\lambda^{lin,[0]}w$ is the ``hl" term,
which is the only one where the  $v$ frequency is no longer the highest, $\lambda$ in this case, and is instead the second highest, call it $\nu$. It has the form
\[
N_\lambda^{lin,[0]}w = \sum_{\nu \ll \lambda} \lambda \nu^\frac12 \coeff_{< \nu}\calL_\nu(\bfw_\nu,\bfu_\lambda).
\]
\end{itemize}

We  will consider the linearized equation in $L^2$. Our main density-flux identities are as follows:

\begin{theorem}\label{t:lin-df}
Let $u$ be a solution to the full equation satisfying the bounds \eqref{uk-ee-boot}-\eqref{uL6L4-boot} with an $\ell^1$ frequency envelope $\epsilon c_\lambda$ which in addition satisfies 
$c_\lambda \lesssim \lambda^{-\frac1{12}}$.
Let $w$ be a solution to the linearized equation \eqref{u-para-lin} in the time interval $[0,1]$, with $L^2$ 
initial data $w(0)$ with frequency envelope $d_\lambda$. Assume that the bootstrap assumptions \eqref{wk-ee-boot}-\eqref{wL6L4-boot} for  the linearized variable $w$ hold in $[0,1]$. Then there exist density corrections $ B^{lin}_{\lambda,m}$ and $B^{lin}_{\lambda,p}$ and flux corrections 
$ Q^{lin}_{\lambda,m}$ and $Q^{lin}_{\lambda,p}$
so that the associated modified mass 
and modified momentum
\begin{equation}
   \Mslin_{\lambda} = \Mss_{\lambda}+ B^{lin}_{\lambda,m} , 
   \qquad  \Pslin_{\lambda} = \Pss_{\lambda}+ B^{lin}_{\lambda,p}
\end{equation}
satisfy the relations 

\newcommand{\Fsjlin}{\fF^{\sharp,j}}
\newcommand{\Esjlin}{\fE^{\sharp,jk}}

\begin{equation}
\left\{
\begin{aligned}
&\ \partial_t \Mslin_{\lambda}+\nabla\cdot (V_{[\lambda]}\Mslin_{\lambda})=\nabla\cdot (\Fsjlin_{\lambda}+ Q^{j,3}_{\lambda,m} + Q^{j,err}_{\lambda,m} +  Q^{j,lin}_{\lambda,m})  + R^{4}_{\lambda,m}  \\
&\ \partial_t \Pslin_{\lambda}+\nabla\cdot (V_{[\lambda]}\Pslin_{\lambda})=\nabla\cdot (\Esjlin_{\lambda}+ Q_{\lambda,p}^{jk,3}  + Q^{jk,err}_{\lambda,p}+ Q^{jk,lin}_{\lambda,p}) + R^{j,4}_{\lambda,p},  
\end{aligned}
\right.
\end{equation}  
where the corrections and the source terms are as follows: 

\begin{enumerate}[label=\Alph*)]
    \item 
 The density corrections   $B^{lin}_{\lambda,m}$ and $B^{lin}_{\lambda,p}$
 are as follows:
 \begin{equation}
  B^{lin}_{\lambda,m} = B^{lin,[0]}_{\lambda,m}+ B^{lin,[1]}_{\lambda,m}+
  B^{lin,[2]}_{\lambda,m},
  \qquad 
  B^{lin}_{\lambda,p} = B^{lin,[0]}_{\lambda,p}+ B^{lin,[1]}_{\lambda,p}+
  B^{lin,[2]}_{\lambda,p},
 \end{equation}
 where
 \begin{itemize}
  \item The ``hl" corrections  $B^{lin,[0]}_{\lambda,m}$ and $B^{lin,[0]}_{\lambda,p}$ 
have the form 
\[
B^{lin,[0]}_{\lambda,m} = \sum_{\nu \ll \lambda}\lambda \coeffu_{< \nu}\calL_\nu( \bfw_\nu,\bfu_\lambda,\bfw_\lambda), \qquad
B^{lin,[0]}_{\lambda,p} = \sum_{\nu \ll \lambda}\lambda^\frac12 \coeffu_{< \nu}\calL_\nu( \bfw_\nu,\bfu_\lambda,\bfw_\lambda);
\]
  
\item The paradifferential corrections  $B^{lin,[1]}_{\lambda,m}$ and $B^{lin,[1]}_{\lambda,p}$ have the form
\[
     B^{j,lin,[1]}_{\lambda,m} = Op(D^\frac12 \coeffu_{<\lambda} S^{-\frac12}_\lambda)(\bfw_\lambda,\bfw_\lambda), 
     \qquad B^{jk,lin,[1]}_{\lambda,p} = Op(
    D^\frac12 \coeffu_{<\lambda} S^{-1}_\lambda)(\bfw_\lambda,\bfw_\lambda);
    \]  

\item The quadratic ``hh" corrections
$B^{lin,[2]}_{\lambda,m}$ and $B^{lin,[2]}_{\lambda,p}$ have the form 
\[
B^{lin,[2]}_{\lambda,m}= \sum_{\mu \gtrsim \lambda} \lambda^\frac12 \mu^\frac12 \coeff_{<\lambda}\calL_\lambda(\bfw_\lambda,\bfu_\mu,\bfw_\mu), \qquad  B^{lin,[2]}_{\lambda,p}= \sum_{\mu \gtrsim \lambda}\mu^\frac12 \coeff_{<\lambda}\calL_\lambda(\bfw_\lambda,\bfu_\mu,\bfw_\mu).
  \]
 \end{itemize}

 \item The flux corrections   $Q^{lin}_{\lambda,m}$ and $Q^{lin}_{\lambda,p}$
 are as follows:
 \begin{equation}
  Q^{lin}_{\lambda,m} = Q^{lin,[0]}_{\lambda,m}+ Q^{lin,[1]}_{\lambda,m}+
  Q^{lin,[2]}_{\lambda,m},
  \qquad Q^{lin}_{\lambda,p} = Q^{lin,[0]}_{\lambda,p}+ Q^{lin,[1]}_{\lambda,p}+
  Q^{lin,[2]}_{\lambda,p},
 \end{equation}
 where
 \begin{itemize}
 \item The ``hl" corrections  $Q^{lin,[0]}_{\lambda,m}$ and $Q^{lin,[0]}_{\lambda,p}$ 
have the form
\[
Q^{lin,[0]}_{\lambda,m} = \sum_{\nu \ll \lambda}\lambda^\frac12 \coeffu_{< \nu}\calL_\nu( \bfv_\nu,\bfu_\lambda,\bfw_\lambda), \qquad
Q^{lin,[0]}_{\lambda,p} = \sum_{\nu \ll \lambda} \coeffu_{< \nu}\calL_\nu( \bfv_\nu,\bfu_\lambda,\bfw_\lambda),
\]
with the symbols for the $(u_\lambda,\bw_\lambda)$ and the $(\bu_\lambda,w_\lambda)$ components satisfying
\begin{equation}\label{bq-relation-lin}
     q^{j,lin,[0]}_{\lambda,m} = \fp^j_\lambda
     b^{lin,[0]}_{\lambda,m}, \qquad
    q^{jk,lin,[0]}_{\lambda,p}  = \fp^k_\lambda
     b^{j,lin,[0]}_{\lambda,p};
    \end{equation} 
     \item The paradifferential corrections  $Q^{lin,[1]}_{\lambda,m}$ and $Q^{lin,[1]}_{\lambda,p}$ 
     have the form
     \[
     Q^{j,lin,[1]}_{\lambda,m} = Op(D^\frac12 \coeffu_{<\lambda} S^{0}_\lambda)(\bfw_\lambda,\bfw_\lambda), 
     \qquad Q^{jk,lin,[1]}_{\lambda,p} = Op(
    D^\frac12 \coeffu_{<\lambda} S^{-\frac32}_\lambda)(\bfw_\lambda,\bfw_\lambda);
    \]  
with the symbol for the $(w_\lambda,\bw_\lambda)$ component    satisfying \eqref{bq-relation-nl}.

\item The quadratic ``hh" corrections
$Q^{lin,[2]}_{\lambda,m}$ and $Q^{lin,[2]}_{\lambda,p}$ have the form 
\[
Q^{lin,[2]}_{\lambda,m}= \sum_{\mu \gtrsim \lambda} \lambda^\frac12  \coeff_{<\lambda}\calL_\lambda(\bfw_\lambda,\bfu_\mu,\bar \bfw_\mu), \qquad  Q^{lin,[2]}_{\lambda,p}= \sum_{\mu \gtrsim \lambda} \coeff_{<\lambda}\calL_\lambda(\bfw_\lambda,\bfu_\mu,\bar \bfw_\mu).
  \]

 \end{itemize}

\item 
The additional flux terms $Q^{err}_{\lambda,m}$ and $Q^{err}_{\lambda,p}$ satisfy
\eqref{df-err-flux}.\\
The errors $R^{4}_{\lambda,m}$, 
$R^{j,4}_{\lambda,p}$ satisfy
\begin{equation}\label{df-err-lin-m}
\|R^4_{\lambda,m}\|_{L^1} \lesssim C\epsilon d_\lambda^2, 
\qquad \|R^{j,4}_{\lambda,p}\|_{L^1} \lesssim C \epsilon d_\lambda^2 \lambda^{-\frac12}.
\end{equation}

\end{enumerate}
 \end{theorem}

\begin{proof}
The starting point of the proof is the paradifferential result in Theorem~\ref{t:second-df}. The contribution of $F^{lin}_\lambda$ is directly perturbative,
so it suffices to obtain the corrections 
associated to the $N^{lin,[j]}_\lambda$ source terms. But this is easily achieved due to the following considerations:

\begin{itemize}
    \item $N^{lin,[1]}_\lambda$ is identical 
    to $N^{[1]}_\lambda$, so we already have the correction from Theorem~\ref{t:df-NL}.

\item $N^{lin,[2]}_\lambda$ is nearly identical  to $N^{[2]}_\lambda$, so the correction is also obtained exactly as in Theorem~\ref{t:df-NL}.

\item $N^{lin,[0]}_\lambda$ is similar to 
$N^{lin,[1]}_\lambda$, only with the $u$ and $w$  arguments interchanged. So, the construction of the correction is also similar; the repetitive details are omitted.
\end{itemize}

Once the corrections are constructed, the remaining step is to bound the contribution of $N^{lin} w$ in the time derivative  of the corrections $B^{lin}_{\lambda,m}$ and $B^{lin}_{\lambda,p}$. 
We discuss the mass correction; the momentum correction is similar. This is a bilinear form in $w$ with $u$ dependent coefficients,
so the bound we need for the $N^{lin} w$ contribution is 
\[
\|B^{lin}_{\lambda,m}(\bfw,{\bf N}^{lin} w)\|_{L^1}
\lesssim C \epsilon d_\lambda^2.
\]
We successively consider the three terms 
in $B_{\lambda,m}^{lin}$:

\medskip 

a) Corresponding to $B_{\lambda,m}^{lin,[0]}$,
we need to estimate the expressions 
\[
\lambda\calL_\nu(\bfw_\nu,\bfu_\lambda) {\bf N}^{lin}_\lambda w, \qquad \lambda\calL_\nu({\bf N}^{lin}_\nu w,\bfu_\lambda) \bfw_\lambda, \qquad \nu \ll \lambda.
\]
For the first term we have the bilinear $L^2$ bound 
\[
\|\calL_\nu(\bfw_\nu,\bfu_\lambda)\|_{L^2_{xt}}
\lesssim C \epsilon d_\nu c_\lambda \lambda^{-2} \nu^\frac34.
\]
This can be paired with the $L^2$ bound \eqref{wlinL2remainder} 
for the bulk of $N_\lambda^{lin} w$,
leaving us with the expression $\lambda\calL_\nu(\bfw_\nu,\bfu_\lambda) N^{lin,[2],bal_\lambda w}$.
Here we have three frequency $\lambda$ inputs, but two of them must be transversal, so we can again use two unbalanced bilinear $L^2_{tx}$ bounds.

The second term is estimated directly using the bound \eqref{wlinauxtrilinear}.

 \medskip 

b) Corresponding to $B_{\lambda,m}^{lin,[1]}$,
we need to estimate the expressions 
\[
\nu^\frac12 \lambda^{-\frac12}\calL_\nu(\coeffs_{\nu},\bfw_\lambda) N^{lin}_\lambda w, \qquad \nu \ll \lambda.
\]   
This is similar to the first term in part (a).

 \medskip 

c) Corresponding to $B_{\lambda,m}^{lin,[2]}$,
we need to estimate the expressions 
\[
\lambda^\frac12 \nu^\frac12 \bfw_\lambda\calL_\lambda(\bfu_\nu ,N^{lin}_\nu w), \qquad \lambda^\frac12 \nu^\frac12 N^{lin}_\lambda w\calL_\lambda(\bfu_\nu , \bfw_\nu), \quad \nu \gtrsim \lambda.
\]
We separate the analysis into two cases:
\medskip

c1) $\nu \gg \lambda$. In the first term 
we can use a bilinear $L^2_{t,x}$ bound for $\bfw_\lambda \bfu_\nu$, and then argue as for the first term in part (a), using the $L^2_{t,x}$ estimate for the bulk of $N^{lin}_{\nu} w$
and then treating directly the remaining term. The second term is estimated directly using \eqref{wlinauxtrilinear}.

\medskip

c2) $\nu \approx \lambda$. In this case 
we have four balanced frequencies interacting, in an expression of the form 
\[
\lambda^\frac52\calL_\lambda(\bfu_\lambda,\bfu_\lambda,\bfw_\lambda,\bfw_\lambda)
\]
with a negligible $\coeff$ paracoefficient.
This is exactly like the cubic balanced case in the equation, which in the linearized case is treated directly in a perturbative manner, using four $L^{6,\infty}_tL^4_x$ bounds and H\"older's inequality in time.

\end{proof}

\subsection{Closing the energy estimates}
\label{sec:energy-close}
In this subsection we use the density-flux identities in order to close the uniform energy bounds component of our bootstrap assumptions. To fix the notations we consider the linear paradifferential flow; the argument is similar for the full equation and for the linearized equation.

Integrating in space and time the density-flux relation for the mass 
for a solution $w$ to the linear paradifferential equation, 
we obtain the relation 
\begin{equation}
\int \Mss_\lambda(w)(T) \, dx = \int \Mss_\lambda(w)(0)\, dx +  \int_0^T \int R^4_{\lambda,m}(w)(t)\, dx dt  ,  
\end{equation}
where 
\[
\Mss_\lambda(w) = \Ms_\lambda(w) + B^3_{\lambda, m}(w).
\]

We estimate the leading part as follows:

\begin{proposition}\label{p:ee}
Assume that the bootstrap assumptions \eqref{uk-ee-boot}-\eqref{uL6L4-boot} hold for $u$, and that the bootstrap assumptions \eqref{wk-ee-boot}-\eqref{wL6L4-boot} hold for the linearized variable $w$. Then we have the $L^1_{tx}$ bounds
\begin{equation}\label{m<L2}
\begin{aligned}
 \|\Mss_\lambda(w)\|_{L^1} \lesssim & \  \| \LPs_\lambda w\|_{L^2}^2 + O( C \epsilon d_\lambda^2) ,
 \\
 \| \Pss_\lambda(w)\|_{L^1} \lesssim  & \ \lambda^{-\frac12} (\| \LPs_\lambda w\|_{L^2}^2 + O( C \epsilon d_\lambda^2)),
\end{aligned}
\end{equation}
as well as the equivalence
\begin{equation}\label{m=L2}
  \int \Mss_\lambda(w) \, dx \approx \| \LPs_\lambda w\|_{L^2}^2 + O( C \epsilon d_\lambda^2)  .
\end{equation}
Similar bounds hold for the full equation, respectively the linearized equation.
\end{proposition}

Using the above proposition, we can estimate 
the localized energy of $w$ by 
\[
\begin{aligned}
\| \LPs_\lambda w(T) \|_{L^2}^2 \lesssim 
& \int \Mss_\lambda(w)(T) \, dx+ O( C \epsilon d_\lambda^2)  
\\
\lesssim & \ \int \Mss_\lambda(w)(0) \, dx+
\|R^4_{\lambda,m}\|_{L^1}+
O( C \epsilon d_\lambda^2)
\\
\lesssim & \ \| \LPs_\lambda(w)(0)\|_{L^2}^2
+ O( C \epsilon d_\lambda^2)
\\
\lesssim & \ (1+ C \epsilon) d_\lambda^2,
\end{aligned}
\]
uniformly in $T$.

This proves the uniform bound \eqref{wk-ee}, 
and closes the energy part of the bootstrap for the linear paradifferential equation.
Nearly identical arguments apply for the full equation, respectively the linearized equation, proving the uniform bounds \eqref{uk-ee}, respectively \eqref{wk-ee}.

\begin{proof}[Proof of Proposition~\ref{p:ee}]
To begin with, we show that the bounds in the proposition hold for the leading part of the mass and momentum functional, namely $\Ms_\lambda(w)$ and $\Ps_\lambda(w)$.

We recall that symbol of $\Ms_\lambda$ is
\[
\ms_\lambda(x,\xi,\eta)= \lambda^{\frac12}
\mfq^{-\frac12}(x,\xi)  \mfq^{-\frac12}(x,\eta)
\lps(x,\xi) \lps(x,\eta).
\]
The symbol $\lambda^{\frac12}
\mfq^{-\frac12}(x,\xi)  \mfq^{-\frac12}(x,\eta)$ is real, positive and of size one 
at frequency $\lambda$, so the first 
estimate in \eqref{m<L2} as well as \eqref{m=L2} for $\Ms$ follow
after one standard use of the product rule
for pseudodifferential operators. The second bound in \eqref{m<L2} for the momentum is similar; there one only needs the $\lambda^{-\frac12}$
bound for the symbol.

The above step is common for the paradifferential equation, the full equation and the linearized equation. 

To complete the proof of the proposition it remains to prove perturbative fixed time bounds for the corrections. These depend on which flow we consider, and we collect them in the following lemma:

\begin{lemma}\label{l:corrections}
Under the same bootstrap assumptions as in Proposition~\ref{p:ee}, fixed time bounds  hold as follows:

a) For the paradifferential equation mass/momentum  corrections we have
\begin{equation}
\| B^3_{\lambda,m}\|_{L^1_{tx}} \lesssim C \epsilon  d_\lambda^2, \qquad \| B^3_{\lambda,p}\|_{L^1_{tx}} \lesssim C \epsilon  d_\lambda^2 \lambda^{-\frac12}; 
\end{equation}

b) For the nonlinear equation mass/momentum corrections we have
\begin{equation}
\| B^{NL}_{\lambda,m}\|_{L^1_{tx}} \lesssim C^3 \epsilon^3  c_\lambda^2 \lambda^{-2s_c},
\qquad \| B^{NL}_{\lambda,p}\|_{L^1_{tx}} \lesssim C^3 \epsilon^3  c_\lambda^2 \lambda^{-2s_c}\lambda^{-\frac12};
\end{equation}

c)  For the linearized equation mass/momentum corrections we have
\begin{equation}
\| B^{lin}_{\lambda,m}\|_{L^1_{tx}} \lesssim C \epsilon d_\lambda^2,
\qquad
\| B^{lin}_{\lambda,p}\|_{L^1_{tx}} \lesssim C \epsilon  d_\lambda^2 \lambda^{-\frac12}.
\end{equation}
\end{lemma}

\begin{proof}
We argue in the case of the mass associated densities, the momentum associated densities are similar.

a) Recalling that $B^3_{\lambda,m} =Op(D^\frac12\coeffu_{<\lambda} S^{-\frac12}_\lambda)(w_\lambda,\bw_\lambda)$, after separation of variables it suffices 
to estimate for $g \in \coeffu$
\[
\lambda^{-\frac12} \| D^\frac12 g_{<\lambda} \calL_\lambda(w_\lambda,\bw_\lambda)\|_{L^1_{t,x}}
\lesssim \lambda^{-\frac12} \| D^\frac12 g_{<\lambda}\|_{L^\infty} \|w_\lambda\|_{L^2}^2 \lesssim C \epsilon  d_\lambda^2.
\]
\medskip

b) The bound for $B^{[1],NL}_{\lambda,m}$
is similar to the one in part (a). The next term is $B^{[2],NL}_{\lambda,m}$, which is estimated by
\[
\begin{aligned}
\| B^{[2],NL}_{\lambda,m}\|_{L^1_{t,x}}
&\lesssim \sum_{\nu \gtrsim \lambda} \lambda^\frac12 \nu^\frac12 \|\calL(u_\lambda,u_\nu,u_\nu)\|_{L^1_{t,x}}\\
&\lesssim \sum_{\nu \gtrsim \lambda} \lambda^\frac12 \nu^\frac12 \|u_\lambda\|_{L^\infty} \|u_\nu\|_{L^2}^2\\
&\lesssim C^3 \epsilon^3 c_\lambda \lambda^{-2s_c}\sum_{\nu \gtrsim \lambda} c_\nu^2 (\lambda/\nu)^{7/2} ,
\end{aligned}
\]
 and conclude using the slowly varying property for $c_\lambda$.   The bounds for 
 $B^{[3],NL}_{\lambda,m}$ and $B^{[4],NL}_{\lambda,m}$ are similar but simpler.

\medskip

c) The bounds for $B^{[1],lin}_{\lambda,m}$ and $B^{[2],lin}_{\lambda,m}$ are as above
with a worse but still acceptable power, namely $1/2$, for $\lambda/\nu$ in the case of the latter. The new term here is $B^{[0],lin}_{\lambda,m}$, which we consider next. We have
\[
\| B^{[0],lin}_{\lambda,m}\|_{L^1}
\lesssim \sum_{\nu \ll \lambda}\lambda \|\calL_\nu( w_\nu,u_\lambda,w_\lambda)\|_{L^1}
\lesssim \sum_{\nu \ll \lambda}\lambda \|w_\nu\|_{L^\infty} \|u_\lambda\|_{L^2}\|w_\lambda\|_{L^2}
\lesssim C \epsilon \sum_{\nu \ll \lambda} (\nu/\lambda) c_\lambda d_\lambda d_\nu ,
\]
where we can again finish the last summation using the slowly varying property for the $d_\lambda$ envelope.
 
\end{proof}
The proof of the proposition is now concluded.
\end{proof}
\section{Bilinear estimates via interaction Morawetz}
\label{sec:Morawetz}

Our objective here is to use interaction Morawetz methods in order to prove 
our main bilinear estimates, and thus close our bootstrap assumptions, completing the proof of Theorems~\ref{t:local-fe3}, \ref{t:local-fe3-lin}, \ref{t:local-fe3-lin-para}. 
This needs to be done at three levels:
\begin{itemize}
    \item for the paradifferential equation;
    \item for the full equation;
    \item for the linearized equation.
\end{itemize}
The bounds for the first two equations are 
global in time, while those for the linearized equation are local in time. However, the latter fact is not due to any obstructions arising in this section, but 
rather to the need to treat balanced cubic terms as errors in the layered expansion of the linearized equation in earlier sections.
Precisely, this is due to the fact that the error bound \eqref{wlinremainder} in Theorem~\ref{wfinaltheorem} is local in time.

Of the three flows above, the case of the paradifferential equation contains the main difficulties, and will be considered first. For the remaining cases we will 
be able to largely reuse the analysis in the paradifferential case for the main part, 
supplemented with additional error estimates
in order to handle the contributions of the additional density-flux corrections.

As a side note, we remark that our arguments also establish that the bilinear estimates are valid for any pairings of solutions to the three flows, locally in time where the linearized flow is involved. 

\medskip

In our interaction Morawetz analysis, we will consider separately two cases:
\begin{itemize}
    \item unbalanced interactions, where we pair two waves with different frequencies $\mu \ll \lambda$.
    \item balanced interactions, where we pair two waves with comparable frequencies, denoted both by $\lambda$.
\end{itemize}

The first case is easier, though it does have the added difficulty that the powers of the frequency factors $\lambda$ and $\mu$ need to be carefully balanced. The second case is technically more difficult, as one needs to capture the null structure in the interaction Morawetz output. 

\bigskip

The starting point here is provided by the density-flux identities which, for
a solution $w$ to the paradifferential equation \eqref{para-re} localized at
frequency $\lambda$, take the form
\begin{equation}
\partial_t \Mss_\lambda  = \partial_j(\Fsj_\lambda - V^j_{[\lambda]} \Mss_\lambda+ Q^{j,3}_{\lambda,m}) + R_{\lambda,m}  ,  
\end{equation}
respectively
\begin{equation}
\partial_t \Pssj_\lambda  = \partial_k (\Esj_\lambda - V_{[\lambda]} \Pssj_\lambda + Q^{kj,3}_{\lambda,p}) + R^j_{\lambda,p},  
\end{equation}
with additional density and flux corrections in the case of the full and the linearized equations, 
where we need to describe the respective components:

\begin{itemize}
\item The sharp mass $\Ms_\lambda$ has a $\tcoeffs S^0_\lambda$ symbol localized to frequency $\lambda$ of the form
\[
\ms_\lambda(x,\xi,\eta) = \tps_\lambda(x,\xi) \tps_\lambda( x,\eta);
\]

\item The mass flux  $\Fsj_\lambda$ has symbol  in $\tcoeffs S^{-\frac12}_\lambda$,
\[
\fsj_\lambda (x,\xi,\eta) = \fp^j_\lambda(x,\xi,\eta) \ms_\lambda(x,\xi,\eta),
\quad 
 \fp^j_\lambda(x,\xi,\eta) = 
\int_0^1 \partial_{\xi_j} \omega_{<\lambda}(x, (1-h)\xi+h\eta) \, dh 
\]
again localized at frequency $\lambda$;

\item The sharp momentum has symbol  in $\tcoeffs S^{-\frac12}_\lambda$,
\begin{equation}
   \psj_\lambda(x,\xi,\eta) = 
   \tfp^j_\lambda(x,\xi,\eta)\ms_{\lambda}(x,\xi,\eta), \qquad \tfp^j_\lambda(x,\xi,\eta) = \fp^j_{0,\lambda}(\xi^\sharp,\eta^\sharp);
\end{equation}

\item  The sharp momentum flux  
has symbol  in $\tcoeffs S^{-1}_\lambda$
\begin{equation}
\esj_{\lambda}(x,\xi,\eta)=  \psj_\lambda(x,\xi,\eta) \fp_\lambda^k(x,\xi,\eta);
\end{equation}

\item The corrected mass  and momentum are 
\[
\Mss_\lambda = \Ms_\lambda + B^3_{\lambda,m}, \qquad \Pssj_\lambda = \Psj_\lambda + B^{j,3}_{\lambda,p}
\]
with normal form type corrections  $B^3_{\lambda,m}$, $B^{j,3}_{\lambda,p}$ which is cubic at leading order;

\item The density corrections $B^3_{\lambda,m}$, $B^{j,3}_{\lambda,p}$ are quadratic forms 
with low frequency paracoefficients with symbols of the form
\[
b_{\lambda,m}^{3} \in  D^\frac12\coeffu_{<\lambda}S_\lambda^{-\frac12}, 
\qquad b_{\lambda,p}^{j,3} \in  D^\frac12\coeffu_{<\lambda}S_\lambda^{-1};
\] 

\item The flux corrections $Q^3_{\lambda,m}$
$Q^3_{\lambda,p}$
are quadratic forms with low frequency paracoefficients with symbols of the form  
  \[
q_{\lambda,m}^{3} \in  D^\frac12\coeffu_{<\lambda} 
 S_\lambda^{-1}, \qquad q_{\lambda,p}^{3} \in  D^\frac12\coeffu_{<\lambda} 
 S_\lambda^{-\frac32}.
 \]  
\end{itemize}

If we now switch to the case  of the full equation or the linearized equation,
then we need to modify the densities and fluxes with  additional corrections 

\begin{itemize}
\item ``$hl$'' terms $B^{[0]}_{\lambda,m}$, $B^{[0]}_{\lambda,p}$ in the case of the linearized equation, 
of the form
\[
B^{[0]}_{\lambda,m} = \sum_{\mu \ll \lambda}
\lambda \coeff_{<\mu}\calL_\mu(\bfv_\mu, \bfv_\lambda,\bfu_\lambda), \qquad B^{[0]}_{\lambda,p} = \sum_{\mu \ll \lambda}
\lambda^\frac12 \coeff_{<\mu}\calL_\mu(\bfv_\mu, \bfv_\lambda,\bfu_\lambda),
\]
with corresponding flux terms
\[
Q^{[0]}_{\lambda,m} = \sum_{\mu \ll \lambda}
\lambda^\frac12 \coeff_{<\mu}\calL_\mu(\bfv_\mu, \bfv_\lambda,\bfu_\lambda), \qquad Q^{[0]}_{\lambda,p} = \sum_{\mu \ll \lambda}
 \coeff_{<\mu}\calL_\mu(\bfv_\mu, \bfv_\lambda,\bfu_\lambda);
\]

\item additional ``$lh$'' terms $B^{[1]}_{\lambda,m}$, $B^{[1]}_{\lambda,p}$,
and corresponding flux terms $Q^{[1]}_{\lambda,m}$, $Q^{[1]}_{\lambda,p}$,
in the case of both the 
full equation and the linearized equation, which are  in the same classes as  $B^3$,
$Q^3$ above, with the difference that 
terms without\footnote{e.g. with the pairings $(v_\lambda,u_\lambda)$ or $(\bv_\lambda,\bu_\lambda)$} phase rotation symmetry are also allowed.

\item ``$hh$'' terms $B^{[2]}_{\lambda,m}$, $B^{[2]}_{\lambda,p}$ in the case of both the  full equation and the linearized equation, of the form
\[
B_{\lambda,m}^{[2]} = \lambda^\frac12 \nu^\frac12\calL_\lambda(\bfu_{\nu},\bfv_\nu) \bfv_\lambda, \qquad B_{\lambda,p}^{[2]} =  \nu^\frac12\calL_\lambda(\bfu_{\nu},\bfv_\nu) \bfv_\lambda
\]
in the latter case,  with corresponding flux terms    
\[
Q_{\lambda,m}^{[2]} =  \lambda^\frac12 \calL_\lambda(\bfu_{\nu},\bfv_\nu) \bfv_\lambda, \qquad Q_{\lambda,p}^{[2]} =\calL_\lambda(\bfu_{\nu},\bfv_\nu) \bfv_\lambda
\]    
both with additional $\coeff_{<\lambda}$ 
paracoefficients;
    
\item Balanced quartic terms $B^{[3]}_{\lambda,m}$, $B^{[3]}_{\lambda,p}$ in the case of  the full equation,
\[
B^{[3]}_{\lambda,m} = OP(\coeff S^\frac52_\lambda)(\bfu_\lambda,\bfu_\lambda,\bfu_\lambda,\bfu_\lambda), \qquad
B^{[3]}_{\lambda,p} = OP(\coeff S^2_\lambda)(\bfu_\lambda,\bfu_\lambda,\bfu_\lambda,\bfu_\lambda)
\]
with corresponding flux terms 
\[
Q^{[3]}_{\lambda,m} = OP(\coeff S^2_\lambda)(\bfu_\lambda,\bfu_\lambda,\bfu_\lambda,\bfu_\lambda), \qquad
Q^{[3]}_{\lambda,p} = OP(\coeff S^\frac32_\lambda)(\bfu_\lambda,\bfu_\lambda,\bfu_\lambda,\bfu_\lambda),
\]
including complex conjugates but excluding phase rotation symmetric\footnote{i.e. with exactly two terms conjugated} terms;

\item Balanced quintic terms $B^{[4]}_{\lambda,m}$ in the case of  the  full equation,
\[
B^{[4]}_{\lambda,m} = OP(\coeff S^\frac72_\lambda)(\bfu_\lambda,\bfu_\lambda,\bfu_\lambda,\bfu_\lambda,\bfu_\lambda), \qquad
B^{[4]}_{\lambda,p} = OP(\coeff S^3_\lambda)(\bfu_\lambda,\bfu_\lambda,\bfu_\lambda,\bfu_\lambda,\bfu_\lambda)
\]
with corresponding flux terms 
\[
Q^{[4]}_{\lambda,m} = OP(\coeff S^3_\lambda)(\bfu_\lambda,\bfu_\lambda,\bfu_\lambda,\bfu_\lambda,\bfu_\lambda), \qquad
Q^{[4]}_{\lambda,p} = OP(\coeff S^\frac52_\lambda)(\bfu_\lambda,\bfu_\lambda,\bfu_\lambda,\bfu_\lambda,\bfu_\lambda),
\]
including any number of complex conjugates.

\end{itemize}

The interaction Morawetz functional corresponding to two solutions $v$ and $w$ 
for the paradifferential equation, localized at frequency $\lambda$, respectively $\mu$
is 
\begin{equation}\label{I}
\bI_{\lambda\mu}(v,w) :=   \iint a_j(x-y) (\Mss_\lambda(v)(x) \Pssj_{\mu}(w) (y) -  
\Pssj_\lambda(v)(x) \Mss_\mu(w) (y)) \, dx dy,
\end{equation}
where in order to fix the notations we will assume that $\mu \lesssim \lambda$. 
The functions $v$ and $w$ will be assumed to be two different solutions to the 
linear paradifferential equation, localized at frequency $\lambda$, respectively $\mu$,
placed under the envelope\footnote{Strictly speaking, in the paradifferential case we work with frequency localized functions, so 
there need be no slowly varying type connection between these values.}
type bounds $d_\lambda$, respectively $d_\nu$.
Alternatively we will work with solutions for the full equation, in which case we set $d_\lambda = c_\lambda \lambda^{-s_c}$,
and with solutions for the linearized equation, where $d_\lambda$ is the actual slowly varying $L^2$ frequency envelope for the solution.

For the weight $a$ we will consider variations of the  
homogeneous choice
\[
 a(x) = |x|
\]
associated to a localization scale $\nu$
such as 
\[
a_\nu(x) = \nu^{-1} \langle \nu x\rangle := \la x \ra_{\nu}.
\]
The scale $\nu$ will always be chosen with the property that $\nu \lesssim \mu$.
On occasion we will also use the alternative
notation
\[
a_r(x) = r \la x/r\ra, 
\]
in order to emphasize the lowest spatial scale $r$.

Computing the time derivative of the interaction functional and integrating by parts we get 
\begin{equation}\label{dtI}
\frac{d}{dt} \bI_{\lambda,\mu}(u,v) 
= \bJ^4_{\lambda\mu} + \bJ^5_{\lambda\mu} + \bJ^6_{\lambda\mu}+ \bR_{\lambda\mu},
\end{equation}
where, with the convention that all quadratic forms in $v$ are evaluated at $x$ and all quadratic forms
in $w$ are evaluated at $y$, we have:

\bigskip

\emph{A. The leading term} contains the quartic contributions,
\begin{equation}
\begin{aligned}
\bJ^4_{\lambda\mu} = \int a_{jk}(x-y)
& \ \left(  \Ms_\lambda(v) \Esj_\mu(w)
-  \Psj_\lambda(v) \Fsk_\mu(w) \right.
 \left.- \Fsk_\lambda (u) \Psj_\mu(v) 
+
\Esj_\lambda(v) \Ms_\mu(w)\right) \, dx dy;
\end{aligned}
\end{equation}

\medskip
\emph{B. The main error term} $\bJ^5_{\lambda\mu}$ contains quintic contributions, and in the paradifferential case has two components
\[
\bJ^5_{\lambda\mu} = \bJ^{5,1}_{\lambda \mu} + \bJ^{5,2}_{\lambda \mu} + \bJ^{5,3}_{\lambda \mu}
\]
as follows: 
\begin{enumerate}[label=\roman*)]
\item  A leading component containing the 
contribution of the transport terms:
\begin{equation}
\bJ^{5,main}_{\lambda \mu} = \int a_{jk}(x-y)
(V_{<\lambda}(x) - V_{<\mu}(y))(\Mss_\lambda(v) \Pssj_\mu(w)
- \Pss_\lambda(v) \Mss_\mu(w))
\, dx dy
\end{equation}
with corresponding modifications in the case of the full equation, respectively the paradifferential equation.
\bigskip

\item  A second term containing the paradifferential, ``lh'' corrections 
\begin{equation}
\begin{aligned}
\bJ^{5,[1]}_{\lambda \mu} =    \int & a_{jk}(x-y)
  \left(\Ms_\lambda(v) Q^{jk,[1]}_{\mu,p}(w) 
-  \Psj_\lambda(v) Q^{k,[1]}_{\mu,p}(w)
- Q^{j,[1]}_{\lambda,m}(v) \Psk_\mu(w)
+ Q^{jk}_{\lambda,p}(v) \Ms_\mu(w) \right.
\\
& \left.+
B^{[1]}_{\lambda,m}(v) \Esj_\mu(w)
-  \Fsj_\lambda(v) B^{k,[1]}_{\mu,p}(w)
-  B^{j,[1]}_{\lambda,p}(v) \Fsk_\mu(w)
+ 
\Esj_\lambda(v) B^{[1]}_{\mu,m}(w)
\right) \, dx dy .
\end{aligned}
\end{equation}
\end{enumerate}
Here we include both the paradifferential 
$B^3$ and $Q^3$ corrections, and the 
additional $B^{[1]}$ and $Q^{[1]}$ corrections arising both for the full equation and in the linearized case.

\medskip

\emph{ C. The secondary error $\bJ^6$}  contains the $BQ$ type terms,
\[
\bJ^{6}_{\lambda\mu} \! =\! \int \! a_{jk}(x-y) \left(  B_{\lambda,m}(v) Q_{\mu,p}(w)
\! -\!  B_{\lambda,p}(v) Q_{\mu,m}(w) - Q_{\lambda,m} (u) B_{\mu,p}(v) 
\!+\!
Q_{\lambda,p}(v) B_{\mu,m}(w)\right)   
dx dy
\]
for all our choices of density flux corrections.

\medskip

\emph{ D. The residual terms} of the form 
\[
\bJ^{err}_{\lambda\mu} \!=\! \int \! a_{j}(x-y) \left(  \Mss_\lambda(v) R^j_{\mu,p}(w)
-  \Pssj_\lambda(v) R_{\mu,m}(w) + R_{\lambda,m} (u) \Pssj_\mu(v) 
\!-\! R_{\lambda,p}(v) \Mss_\mu(w)\right)    
dx dy.
\]
\bigskip

When considering either solutions to the full equation or two the linearized equation,
we encounter all of the error terms above, with minor adjustments in the notations, but, more importantly,
we will have additional ``quintic" error terms in part B, which need to be considered separately.
Adding to the two item list above, these are as follows:

\begin{enumerate}[resume,label=\roman*)]
\item A second, ``hl" contribution in the case of the linearized equation,
\[
\begin{aligned}
\bJ^{5,[0]}_{\lambda \mu} = &\  \int a_{jk}(x-y)
\! \left(\Ms_\lambda(v) Q^{jk,[0]}_{\mu,p}(w) 
-  \Psj_\lambda(v) Q^{k,[0]}_{\mu,p}(w)
- Q^{j,[0]}_{\lambda,m}(v) P^{\sharp,k}_\mu(w)
+ Q^{jk}_{\lambda,p}(v) \Ms_\mu(w) \right.
\\
& \qquad \left.+
B^{[0]}_{\lambda,m}(v) E^{jk}_\mu(w)
-  P^j_\lambda(v) B^{k,[0]}_{\mu,p}(w)
-  B^{j,[1]}_{\lambda,p}(v) P^{k}_\mu(w)
+ 
E^{jk}_\lambda(v) B^{[0]}_{\mu,m}(w)
\right) \, dx dy.
\end{aligned}
\]

\item  A term containing the ``hh" corrections, arising both when we consider the full equation, and  in the case of the linearized equation,
\begin{equation}
\begin{aligned}
\bJ^{5,[2]} =    \int & a_{jk}(x-y)
\!\left(\Ms_\lambda(v) Q^{jk,[2]}_{\mu,p}(w) 
-  \Psj_\lambda(v) Q^{k,[2]}_{\mu,p}(w)
- Q^{j,[2]}_{\lambda,m}(v) P^{\sharp,k}_\mu(w)
+ Q^{jk,[2]}_{\lambda,p}(v) \Ms_\mu(w) \right.
\\
& \qquad \left.+
B^{[2]}_{\lambda,m}(v) E^{jk}_\mu(w)
-  P^j_\lambda(v) B^{k,[2]}_{\mu,p}(w)
-  B^{j,[2]}_{\lambda,p}(v) P^{k}_\mu(w)
+ 
E^{jk}_\lambda(v) B^{[2]}_\mu(w)
\right) \, dx dy 
\end{aligned}
\end{equation} 
with $(\nu,\nu)$ as the pair of high frequencies.

\item Two final terms containing the contributions of cubic and quartic terms in the equation,
arising only in the study of the full equation, with formulas similar to the ones above.

\end{enumerate}

\bigskip

Now we outline our strategy for the proof of the bilinear $L^2_{t,x}$ bounds. These will be obtained as follows, starting from the time integrated form of the relation \eqref{dtI}:

a) Fixed time bounds. Precisely, as a consequence of the estimates \eqref{m<L2} we have the energy type fixed time bound
\begin{equation}\label{I-bd}
|\bI_{\lambda\mu}| \lesssim \mu^{-\frac12} 
d_\lambda^2 d_\mu^2 .
\end{equation}

b) Space-time error bounds. Here we estimate 
the error terms  perturbatively, as 
\begin{equation}\label{bJ5-bd}
\left|\int_0^T \bJ^5_{\lambda\mu}\, dt \right| \lesssim 
C \epsilon \mu^{-\frac12} d_\lambda^2 d_\mu^2 ,
\end{equation}
\begin{equation}\label{bJ6-bd}
\left|\int_0^T \bJ^6_{\lambda\mu} dt \right| \lesssim 
C \epsilon \mu^{-\frac12} d_\lambda^2 d_\mu^2 ,
\end{equation}
\begin{equation}\label{bR}
\left|\int_0^T \bR_{\lambda\mu} dt \right| \lesssim 
C \epsilon \mu^{-\frac12} d_\lambda^2 d_\mu^2 
\end{equation}

c) At the conclusion of the first two steps, 
we have obtained the $\bJ^4_{\lambda\mu}$ 
bound
\begin{equation}
   \int_0^T \bJ^4_{\lambda\mu} dt \lesssim  
   \mu^{-\frac12} 
d_\lambda^2 d_\mu^2 + C \epsilon \mu^{-\frac12} d_\lambda^2 d_\mu^2,
\end{equation}
so the last task is to prove that the last time integral of $\bJ^4_{\lambda\mu}$ above
controls the desired bilinear $L^2_{t,x}$ norms, modulo additional perturbative  errors of the same magnitude $C \epsilon \mu^{-\frac12} d_\lambda^2 d_\mu^2$.

The above strategy works as described for the most part, with one small proviso, which is that in the 
$\bJ^4_{\lambda \mu}$ analysis we will encounter one error term whose treatment requires an 
auxiliary correction  to the interaction Morawetz functional $\bI_{\lambda \mu}$. 

\begin{remark}
Since for the density-flux identities we use the paratransported versions $\LPs_\lambda$  of the Littlewood-Paley projectors, the arguments in this section yield bilinear estimates for the functions $\LPs_\lambda u$
rather than the classical Littlewood-Paley pieces $P_\lambda u$. However, it is easily seen that at the level of bilinear estimates
the transition between one class of projectors and the other is relatively straightforward. The details are left for the reader.
\end{remark}

\medskip

As noted earlier, there are two cases,  (i) the \emph{unbalanced} case $\mu \ll \lambda$, and (ii)
the \emph{balanced} case $\lambda \approx \mu$, which for clarity we will discuss separately. Before we do that, we obtain some preliminary bounds which do not depend on the balance of $\lambda$ and $\mu$, but which allow us to dispense quickly of some of the easier terms.

\subsection{ Preliminary bounds}
We begin by observing that the bound \eqref{bR} is straightforward, simply by combining the error bound \eqref{df-err-source} with the energy type bounds in \eqref{m<L2}. This applies for all three equations.

We next consider the effect of the $B$ and $Q$
corrections, for which we prove a set  of  mild 
space-time bounds:

\begin{lemma}\label{l:Bxt}
Assume appropriate bootstrap bounds as in Section~\ref{sec:bootstrap}. Then:

a) For the paradifferential equation mass/momentum  corrections we have:
\begin{equation}\label{B-easy}
\begin{aligned}
\| B^3_{\lambda,m}\|_{L^2L^1
} \lesssim C^3 \epsilon \lambda^{-\frac14} d_\lambda^2, \qquad   \| Q^3_{\lambda,m}\|_{L^2L^1}\lesssim C^3 \epsilon \lambda^{-\frac34} d_\lambda^2,
\\
\| B^3_{\lambda,p}\|_{L^2L^1
} \lesssim C^3 \epsilon \lambda^{-\frac34} d_\lambda^2, \qquad   \| Q^3_{\lambda,p}\|_{L^2L^1}\lesssim C^3 \epsilon \lambda^{-\frac54} d_\lambda^2.
\end{aligned}
\end{equation}

b) For the nonlinear equation mass correction we have:
\begin{equation}\label{B-easy-NL}
\begin{aligned}
\| B^{NL}_{\lambda,m}\|_{L^2 L^1} \lesssim C^3 \epsilon^3 c_\lambda^2 \lambda^{-\frac14} \lambda^{-2s_c} ,
\qquad   \| Q^{NL}_{\lambda,m}\|_{L^2L^1}\lesssim  C^3 \epsilon^3 c_\lambda^2  \lambda^{-2s_c} \lambda^{-\frac34}, 
\\
\| B^{NL}_{\lambda,p}\|_{L^2 L^1} \lesssim C^3 \epsilon^3 c_\lambda^2 \lambda^{-\frac34} \lambda^{-2s_c} ,
\qquad   \| Q^{NL}_{\lambda,p}\|_{L^2L^1}\lesssim  C^3 \epsilon^3 c_\lambda^2  \lambda^{-2s_c} \lambda^{-\frac54}. 
\end{aligned}
\end{equation}

c)  For the linearized equation mass correction we have:
\begin{equation}\label{B-easy-lin}
\begin{aligned}
\| B^{lin}_{\lambda,m}\|_{L^2L^1} \lesssim C^3 \epsilon  d_\lambda^2  \lambda^{-\frac14}, \qquad   \| Q^{lin}_{\lambda,m}\|_{L^2L^1}\lesssim C^3 \epsilon  d_\lambda^2 \lambda^{-\frac34},
\\
\| B^{lin}_{\lambda,p}\|_{L^2L^1} \lesssim C^3 \epsilon  d_\lambda^2  \lambda^{-\frac34}, \qquad   \| Q^{lin}_{\lambda,p}\|_{L^2L^1}\lesssim C^3 \epsilon  d_\lambda^2 \lambda^{-\frac54}.
\end{aligned}
\end{equation}
\end{lemma}
\begin{proof}
We argue in the case of the $B^3_{\lambda,m}$ densities, the momentum associated densities and the flux bounds are similar.

a) Using the expression $B^3_{\lambda,m} =Op(D^\frac12\coeffu_{<\lambda} S^{-\frac12}_\lambda)(v_\lambda,\bv_\lambda)$, after separation of variables it suffices 
to use Lemma~\ref{l:DB-w} to estimate for $g_{<\lambda} \in \coeffu_{<\lambda}$
\[
\lambda^{-\frac12} \| D^\frac12 g_{<\lambda} \calL_\lambda(v_\lambda,\bv_\lambda)\|_{L^2 L^1}
\lesssim  \lambda^{-\frac12} \| D^\frac12 g_{<\lambda} v_\lambda\|_{L^2} \|v_\lambda\|_{L^\infty L^2}^2 \lesssim C \epsilon  d_\lambda^2 \lambda^{-\frac14}.
\]
\medskip

b) The bound for $B^{[1],NL}_{\lambda,m}$
is similar to the one in part (a). The next term is $B^{[2],NL}_{\lambda,m}$, which is estimated using an unbalanced bilinear $L^2$ bound by
\[
\begin{aligned}
\| B^{[2],NL}_{\lambda,m}\|_{L^2 L^1}
\lesssim &\  \sum_{\nu \gtrsim \lambda} \lambda^\frac12 \nu^\frac12 \|\calL_\lambda (\bfu_\lambda,\bfu_\nu,\bfu_\nu)\|_{L^1}
\lesssim \sum_{\nu \gtrsim \lambda} \lambda^\frac12 \nu^\frac12 \|u_\lambda u_\nu\|_{L^2} \|u_\nu\|_{L^2}^2
\\
\lesssim & \ C^3 \epsilon^3 c_\lambda \lambda^{-2s_c}\sum_{\nu \gtrsim \lambda} c_\nu^2 \lambda^{-\frac14} (\lambda/\nu)^{7/2} ,
\end{aligned}
\]
 and conclude using the slowly varying property for $c_\lambda$.   The bounds for 
 $B^{[3],NL}_{\lambda,m}$ and $B^{[4],NL}_{\lambda,m}$ are similar but simpler.

\medskip

c) The bounds for $B^{[1],lin}_{\lambda,m}$ and $B^{[2],lin}_{\lambda,m}$ are as above
with a worse but still acceptable power, namely $1/2$ for $\lambda/\nu$, in the case of the latter. The new term here is $B^{[0],lin}_{\lambda,m}$, which we consider next. We have
\[
\| B^{[0],lin}_{\lambda,m}\|_{L^1}
\lesssim \sum_{\nu \ll \lambda}\lambda \|\calL_\nu( \bfv_\nu,\bfu_\lambda,\bfv_\lambda)\|_{L^2 L^1}
\lesssim \sum_{\nu \ll \lambda}\lambda \| v_\nu u_\lambda\|_{L^2}\|v_\lambda\|_{L^\infty L^2}
\lesssim C \epsilon \lambda^{-\frac14} \sum_{\nu \ll \lambda} (\nu/\lambda)^\frac34 c_\lambda d_\lambda d_\nu ,
\]
where we can again finish the last summation using the slowly varying property for the $d_\lambda$ envelope.
 
\end{proof}
 These bounds  will suffice for $\bJ^6$, which we consider next. For clarity we consider a sample term,
 with the worst balance of the frequencies $\lambda$ and $\mu$. In the paradifferential case, using \eqref{B-easy}, we have
\[
\begin{aligned}
 \left|\iint \partial^2 a(x-y) B^3_{\lambda,m}(v)
 Q^3_{\mu,p}(w)\, dx dy dt \right| \lesssim & \ 
\| \partial^2 a\|_{L^\infty} \|B^3_{\lambda,m}(v)\|_{L^2 L^1_x}
 \|Q^3_{\mu,p}(w)\|_{L^2 L^1_y} 
 \\ \lesssim & \ \mu \cdot C \epsilon d_\lambda^2 \lambda^{-\frac14} \cdot 
 C \epsilon d_\mu^2 \mu^{-\frac54}
 \\
 = & \  C^2 \epsilon^2 d_\lambda^2 d_\mu^2 (\mu \lambda)^{-\frac14}.
\end{aligned}
\]
A similar reasoning applies in the fully nonlinear case and the linearized case, using \eqref{B-easy-NL} and \eqref{B-easy-lin} instead.
This suffices for the proof of \eqref{bJ6-bd}.

The bounds in the above lemma are not good enough to directly use for $\bJ^5$ type terms. However,
there are two $\bJ^5$ type terms, namely 
$\bJ^{5,[3]}$ and $\bJ^{5,[4]}$ which can be dealt with using a slightly better version of the above Lemma:

  \begin{lemma}\label{l:B34}
Assume appropriate bootstrap bounds as in Section~\ref{sec:bootstrap}. Then for the nonlinear flow we have
  \begin{equation}\label{vg-BQ3+}
  \begin{aligned}
\| B^{[3,4]}_{\lambda,m}\|_{L^1L^{p}}\lesssim C^3 \epsilon^3 c_\lambda^2
\lambda^{-2s_c} \lambda^{2(\frac34-\frac1p)}
, \qquad \| Q^{[3,4]}
_{\lambda,m}\|_{L^1L^{p}}\lesssim C^3 \epsilon^3 
\lambda^{-2s_c}\lambda^{2(\frac12-\frac1p)}, \qquad p > \frac43 
\\
\| B^{[3,4]}_{\lambda,p}\|_{L^1L^{p}}\lesssim C^3 \epsilon^3 c_\lambda^2
\lambda^{-2s_c} \lambda^{2(\frac12-\frac1p)}
, \qquad \| Q^{[3,4]}
_{\lambda,p}\|_{L^1L^{p}}\lesssim C^3 \epsilon^3 
\lambda^{-2s_c}\lambda^{2(\frac14-\frac1p)}, \qquad p > \frac43 .
\end{aligned}
\end{equation}
\end{lemma}

\begin{proof}
We consider the first term, as the others are similar. We recall that 
\[
B^{[3]}_{\lambda,m} = \lambda^{\frac52}\calL_\lambda(\bfu_\lambda,\bfu_\lambda,\bfu_\lambda,\bfu_\lambda),
\]
with a negligible bounded $\coeff$ paracoefficient, and at least two 
frequency separated entries. For these two entries we use a bilinear $L^2$ bound, 
and for the other two we use the $L^4 L^{8+}$
bound \eqref{L4Lp-boot}. The powers of $\lambda$ are easily determined by scaling.

Next we have 
\[
B^{[4]}_{\lambda,m} = \lambda^{\frac72}\calL_\lambda(\bfu_\lambda,\bfu_\lambda,\bfu_\lambda,\bfu_\lambda,\bfu_\lambda)
\]
also with a negligible bounded $\coeff$ paracoefficient.
This is even simpler, as we can directly use 
one energy bound $L^{\infty}_t L^2_x$ and four 
  $L^4_t L^{8+}_x$
bounds.  
\end{proof}

We now use H\"older's inequality to estimate 
the contribution of these terms to $\bJ^5_{\lambda\mu}$. We consider a typical term
\[
\left|\int_0^T \int \partial^2 a(x-y)\Mss_\lambda(v)(x)
Q^{[3]}_{\mu,p}(w)(y)\, dx dy \right|
\lesssim \| \Mss_\lambda(v)\|_{L^\infty_t L^1_x}
\|Q^{[3]}_{\mu,p}(w)\|_{L^1L^p} \|\partial^2 a \|_{L^{p'}},
\]
where we recall that $\partial^2 a(z) \lesssim \la z \ra_{\mu}^{-1}$, and $p$ is chosen so that $\frac43 < p < 2$. This produces no additional powers of $\lambda$, and the power of $\mu$ is correct by scaling considerations.

\medskip

Finally, we note that the same argument also allows us to estimate the contribution of  $Q^{4,err}$ to $\bJ^5$ (and belatedly justifies the choice of the function spaces for $Q^{4,err}$).

\bigskip

At this point, we have dealt with the terms $\bJ^6_{\lambda\mu}$, $R_{\lambda\mu}$ and 
two of the $\bJ^5$ terms, namely $\bJ^{5,[3]}_{\lambda\mu}$ and $\bJ^{5,[4]}_{\lambda\mu}$ as well as the contribution of $Q^{4,err}$.
From here on we separate the analysis into the unbalanced and the balanced case.
We continue with the unbalanced case $\mu \ll \lambda$.

\subsection{The unbalanced case \texorpdfstring{$ \mu \ll \lambda$}{}}
\label{sec:unbalanced-close}

In this subsection we use the interaction Morawetz functional introduced above, in order to prove the unbalanced bilinear $L^2$ bounds \eqref{uab-bi-unbal-trans-boot}, respectively  \eqref{uw-bi-unbal-trans-boot} and 
\eqref{ww-bi-unbal-trans-boot} and close  
this part of our bootstrap argument. The bulk of the analysis happens at the level of the paradifferential equation, and then the transition to the full equation and to the linearized equation only requires some milder error estimates.

\subsubsection{ The leading term \texorpdfstring{$\bJ^4$}{}}
The quartic form in the integrand has symbol
\begin{equation}
\begin{aligned}
j^{4,jk}_{\lambda \mu}(x,y,\xi^1,\xi^2,\xi^3,\xi^4) = & \ \ms_\lambda(x,\xi^1,\xi^2)  \ms_\mu(y,\xi^3,\xi^4)
(\fp^j_{\lambda}(x,\xi^1,\xi^2) - \fp^j_\mu(y,\xi^3,\xi^4))
\\& \ (\tfp^{k}_\lambda(x,\xi^1,\xi^2) - \tfp^{k}_\lambda(\xi^{3},\xi^4)) ,
\end{aligned}
\end{equation}
where one could also symmetrize in $j$ and $k$.
Here the symbols of the four entries above are smooth on the corresponding dyadic scales, with  paracoefficients which  belong  to $\tcoeffs_{<\lambda}$, respectively $\tcoeffs_{<\mu}$.  

We  contend that 
only the values of the above symbol on the diagonals $\xi^1= \xi^2$, $\xi^3 = \xi^4$
are important, as any remainder vanishing on 
any of the two diagonals will be 
shown to be perturbative.
Indeed, consider an entry 
with symbol of the form 
\[
J^{4,jk,err}_{\lambda \mu}(\xi^1,\xi^2,\xi^3,\xi^4)= 
(\xi^1-\xi^2)\tcoeffs_{<\lambda} S^{-\frac32}_\lambda(\xi^1,\xi^2) \tcoeffs_{<\mu} S^{-\frac12}_\mu (\xi^3,\xi^4) .
\]
The factor $\xi^1-\xi^2$ yields a derivative 
for the $(v_\lambda, \bv_\lambda)$ entry, which can be integrated by parts in $x$. Then the derivative will apply either to the $\tcoeffs_{<\lambda}$ paracoefficient or to the 
$D^2 a$ weight, yielding an error term of the form
\[
\begin{aligned}
\bJ^{4,jk,err}_{\lambda \mu} =  & \ \iint
\partial^2 a \,  Op( D\coeffs_{<\lambda} S^{-\frac32}_\lambda)(v_\lambda,\bv_\lambda)
Op(\tcoeffs_{<\mu} S^{-\frac12}_\mu) (w_\mu,\bw_\mu) 
\\
& \ \qquad + \partial^3 a \,  Op(\tcoeffs_{<\lambda} S^{-\frac32}_\lambda)(v_\lambda,\bv_\lambda)
Op(\tcoeffs_{<\mu} S^{-\frac12}_\mu) (w_\mu,\bw_\mu) \, dx dy dt
\\
: = & \ \bJ^{4,jk,err,1}_{\lambda \mu} + \bJ^{4,jk,err,2}_{\lambda \mu}.
\end{aligned}
\]
Here at separation scale $|x-y| \approx r \gg \mu^{-1}$, the derivatives of the weight $a$ satisfy $|\partial^2 a| \lesssim r^{-1}$, respectively $|\partial^3 a| \lesssim r^{-2}$, with all scales $\lesssim \mu^{-1}$ grouped together.

To estimate the first term we use one bilinear 
$L^2_{tx}$ bound as in Lemma~\eqref{l:DB-w} for $D\coeffs_{<\lambda} u_\lambda$ and one $L^2_{t,xy}$ bound \eqref{uab-bi-unbal-trans} for $u_\lambda(x) v_\mu(y)$, finally combined with 
an energy bound for the remaining $v_\mu$.
At scale $r$ this yields
\[
\| \bJ^{4,jk,err,1}_{\lambda \mu,r}\|_{L^1_t} \lesssim 
C \epsilon   d_\lambda^2 d_\mu^2 r^{-1}  \lambda^{-\frac34} r^\frac12 \mu^{\frac14} \mu^{-\frac12} ,
\]
with the $\epsilon$ factor arising from the differentiated paracoefficient. This suffices, also with favourable dyadic $r$ summation for $r \gtrsim \mu^{-1}$.

To estimate the second term we simply use two
$L^2_{t,xy}$ bounds \eqref{uab-bi-unbal-trans} for $u_\lambda(x) v_\mu(y)$, which at scale $r$
yields
\[
\| \bJ^{4,jk,err,2}_{\lambda \mu,r}\|_{L^1_t} \lesssim   d_\lambda^2 d_\mu^2
 r^{-2}  \lambda^{-\frac32} r  ,
\]
and which is again favorable, with good $r$ summation and smallness arising from $r \gg \mu^{-1}$.

We argue similarly for error terms 
with symbol of the form 
\[
J^{4,jk,err}_{\lambda \mu}(\xi^1,\xi^2,\xi^3,\xi^4)= 
(\xi^3-\xi^4)\tcoeffs_{<\lambda} S^{-\frac32}_\lambda(\xi^1,\xi^2) \coeffs_{<\mu} S^{-\frac12}_\mu (\xi^3,\xi^4).
\]
Now the factor $\xi^3-\xi^4$ yields a derivative 
for the $(v_\mu, \bv_\mu)$ entry, which can be integrated by parts in $y$. Then the derivative will apply either to the $\coeffs_{<\mu}$ paracoefficient or to the 
$D^2 a$ weight, yielding an error term of the form
\[
\begin{aligned}
\bJ^{4,jk,err}_{\lambda \mu} =  & \ \iint
\partial^2 a \,  Op( \tcoeffs_{<\lambda} S^{-\frac12}_\lambda)(v,\bv)
Op(D\coeffs_{<\mu} S^{-\frac32}_\mu) (w,\bw) 
\\
& \ \qquad + \partial^3 a \,  Op(\tcoeffs_{<\lambda} S^{-\frac12}_\lambda)(v,\bv)
Op(\tcoeffs_{<\mu} S^{-\frac32}_\mu) (w,\bw) \, dx dy dt
\\
: = & \ \bJ^{4,jk,err,1}_{\lambda \mu} + \bJ^{4,jk,err,2}_{\lambda \mu}.
\end{aligned}
\]
Estimating these two terms in the same manner as before we obtain
\[
\| \bJ^{4,jk,err,1}_{\lambda \mu,r}\|_{L^1_t} \lesssim 
C \epsilon   d_\lambda^2 d_\mu^2 r^{-1}  \lambda^{-\frac12} r^\frac12 \mu^{\frac14} \mu^{-\frac34},
\]
respectively
\[
\| \bJ^{4,jk,err,2}_{\lambda \mu,r}\|_{L^1_t} \lesssim 
 d_\lambda^2 d_\mu^2 r^{-2}  \lambda^{-\frac12} r \mu^{-1}, 
\]
both favorable and with good $r$ summation.

We now turn our attention to the symbol of $j^4$ restricted to the diagonal. Here we remark that
the symbols $\fp^j$ and $\tfp^{j}$ agree
at leading order (i.e. for the constant, $x$ independent  terms in the paracoefficients) but not in the nonconstant 
terms. This fact creates a major difficulty which  we need to understand how to defeat.

To keep the ideas streamlined, assume first, for the sake of the argument, that we do not have 
the above issue, i.e. that $\tfp^{j} = \fp^j$.
Then on the diagonal we have a $j^4_{\lambda\mu}$ factorization of the form\footnote{excluding the 
$m^\sharp(x,\xi^1,\xi^2)$
$m^\sharp(y,\xi^3,\xi^4)$ factors, which are unaffected by this discussion.}
\[
(\fp^j(x,\xi^1,\xi^2) - \fp^j(y,\xi^3,\xi^4))(\fp^{k}(x,\xi^1,\xi^2) - \fp^{k}(y,\xi^3,\xi^4)) = 
q^j(x,y,\xi^1,\xi^4) q^k(x,y,\xi^2,\xi^3),
\]
with symbols 
\[
q^j(x,y,\xi,\eta) = 
(\fp^j(x,\xi) - \fp^j(y,\eta))
\]
of size $O(\mu^{-\frac12})$, smooth on the corresponding dyadic scales, and elliptic in the vector valued sense. 

Substituting the $j^4_{\lambda\mu}$ symbols as above, modulo acceptable commutator errors our integrand becomes
\[
\bJ^4_{\lambda \mu} = \int a_{jk}(x-y) Q^j(v_\lambda,\bw_\mu)(x,y) \overline{Q^k(v_\lambda,\bw_\mu)}(x,y)\,  dx dy,
\]
where by $v_\lambda$ we have denoted the frequency localized portion of $v$, given by 
\begin{equation}
 v_\lambda = \tPs_\lambda(x,D) v,   
\end{equation}
and similarly for $w_\mu$.

Here we have the advantage that the integrand 
is pointwise nonnegative. Hence for $a  = a^r$ 
this allows us to  control
\begin{equation}\label{Q-unbal}
r^{-1} \| Q(v_\lambda,\bw_\mu) \|_{L^2_{txy}(B_r)}^2 \lesssim 
d_\lambda^2 d_\mu^2 \mu^{-\frac12}, \qquad r \gg \mu^{-1},
\end{equation}
where $B_r = \{|x-y| \lesssim r\}$.
Using the ellipticity of $Q$ as a vector valued symbol, the bound \eqref{Q-unbal} leads to 
\begin{equation}\label{pf-bi-unbal}
r^{-1}  \mu^{-1} \| v_\lambda(x) \bw_\mu(y) \|_{L^2_{txy}(B_r)}\lesssim 
d_\lambda d_\mu \mu^{-\frac12}, \qquad r \gg \mu^{-1}
\end{equation}
as needed in order to close our unbalanced 
bilinear $L^2_{t,x}$ bound \eqref{uab-bi-unbal-trans}.

\bigskip

We still have to deal with the two slightly different symbols $\fp$ and $\tfp$, for which we need another argument to handle the errors due the 
difference $\delta \fp = \fp - \tfp$. In our main expression for the $\bJ^4_{\lambda\mu}$ symbol, this yields an error 
\[
j^{4,jm,err} = \ms_\lambda(x,\xi^1,\xi^2) \ms_\mu(y,\xi^3,\xi^4) (\fp^j(x,\xi^1) - \fp^j(y,\xi^4))
(\delta \fp^m(x,\xi^2) - \delta \fp^m(y,\xi^3))
\]
with $\delta \fp \in \coeffs S^{-\frac12}$, without any constant coefficient term.
In order to attempt to treat this expression perturbatively,  we expand in terms of the highest frequency paracoefficient,
\[
\delta \fp = \sum_{\nu \ll \lambda} \delta \fp_\nu,
\]
where, separating variables and neglecting harmless lower frequency paracoefficients, 
we schematically write 
\[
\delta \fp_\nu(x,\xi)= g_\nu(x) h_\nu(\xi), \qquad g_\nu \in \coeffs_\nu, \quad h_\nu\in S^{-\frac12}_\nu.
\]

Discarding again perturbative commutator terms, its  contribution to $\bJ^4_{\lambda\mu}$ has the form 
\[
\bJ^{4,err}_{\lambda\mu,\nu}:=\int  a_{jk}(x-y) (g_\nu(x)H^j_\nu v_\lambda \bw_\mu- 
v_\lambda \ol{g_\nu(y) H^j_\nu w_\mu})(\ol {\fP^k v_\lambda} w_\mu - 
\bv_\lambda  \fP^k w_\mu) \, dx dy,
\]
where the coefficients $g_\nu(x)$ and $g_\nu(y)$ can be shown to be perturbatively interchangeable.

We estimate this at separation scale $r$ using a bilinear $L^2_{t,x}$ bound for the $(\nu,\lambda)$ pairing
and a bilinear $L^2_{t,xy}$ bound
for the $(\lambda,\mu)$ pairing, to obtain
\[
\|\bJ^{4,err}_{\lambda\mu,\nu}\|_{L^1_t} \lesssim C\epsilon d_\lambda^2 d_{\mu}^2 \mu^{-1} r^{-1} 
\nu^{-\frac14} r^\frac12 \mu^{\frac14} = C\epsilon d_\lambda^2 d_{\mu}^2
r^{-\frac12} (\nu\mu)^{-\frac14} \nu^{-\frac12}.
\]
We need to compare this bound with $\mu^{-\frac12}$,
so this error is perturbative only if 
\begin{equation}\label{easy-r}
r \gtrsim r_{\nu}:=(\nu \mu)^{-\frac12}.    
\end{equation}
This leaves us with the small $r$ range $\mu^{-1} < r < r_\nu$, where our error has the form
\[
\bJ^{4,err,bad}_{\lambda\mu,\nu}:=\int  a_{jk}(x-y)  (g_\nu(x)H^j_\nu v_\lambda \bw_\mu - 
g_\nu(y) v_\lambda \ol{H_\nu w_\mu})(\ol {\fP^k v_\lambda} w_\mu- 
\bv_\lambda \fP^k w_\mu) \chi_{r < r_\nu} \,dx dy.
\]
To eliminate this remaining term, our strategy will be 
to correct the interaction Morawetz functional.
We will  take the correction to be 
\[
\bI^{4,err}_{\lambda\mu,\nu} = -  \int  a_{j}(x-y)  (\delta \Psj_{\lambda,\nu} (v,\bv) \Ms_\mu(w,\bw)- 
\Ms_\lambda(v,\bv) \delta \Psj_{\mu,\nu}(w,\bw)) \chi_{r < r_\nu}  \,dx dy,
\]
where we no longer need to go to the double sharp setting. In the same way as in the case
of \eqref{I-bd}, we have the fixed time estimate
\[
|\bI^{4,err}_{\lambda\mu,\nu}| \lesssim C \epsilon c_\nu d_\lambda^2 d_\nu^2 \nu^{-\frac12}.
\]
Computing its time derivative,  is not difficult to see that the flux terms generated by the $x$ and $y$ components in this correction exactly cancel $\bJ^{4,err}_{\lambda\mu,\nu}$. This said, there are still several additional components of the time derivative which need to be estimated.

\begin{enumerate}[label=\roman*)]
\item terms arising from the mismatch of the frequency truncation for the velocity $V$ in the material derivative. These are similar to $J^{5,main}$
but with $\Psj$ replaced by $\delta \Psj$,
and are estimated exactly as in Section~\ref{sec:J5-main-unbal}.

\item Terms with a derivative falling on $\chi_{r < r_\nu}$ after the integration by parts.
This yields $r^{-1}$ contributions localized in $\{ r \approx r_\nu\}$, which are perturbative
according to the prior discussion.

\item Terms with a derivative falling on the leading coefficient $g_\nu$. 
Schematically these have the form
\[
\bJ^{4,err,bad,3}_{\nu}:=\int  a_{j}(x-y) \nabla g_\nu (H^j_\nu v_\lambda \ol{w_\mu}- 
 v_\lambda \ol{H^j_\nu w_\nu})(\ol {\fP^j v_\lambda} w_\mu- 
\ol{v_\lambda} \fP^j w_\mu) \chi_{r < r_\nu} \, dx dy.
\]
We again estimate this at scale $r < r_\nu$ using a bilinear $L^2_{tx}$ bound for the $(\nu,\lambda)$ pairing as in Lemma~\ref{l:DB-w} and a bilinear $L^2_{txy}$ bound for a $(\lambda,\mu)$ pairing, obtaining
\[
\|\bJ^{4,err,bad,3}_\nu\|_{L^1_t} \lesssim C\epsilon d_\lambda^2 d_\mu^2 \mu^{-1}  
\nu^{\frac34} r^\frac12 \mu^{\frac14} = 
C\epsilon d_\lambda^2 d_\mu^2 r^{\frac12} \nu^\frac34\mu^{-\frac14} \mu^{-\frac12} \leq C\epsilon d_\lambda^2 d_\mu^2 \nu^\frac12\mu^{-\frac12} \mu^{-\frac12},
\]
which is clearly perturbative with room to spare.

\item Terms arising from the para-material derivative of $g_\nu$. This is the main difference compared with our main interaction 
Morawetz computation, as our density $ g_\nu Q_\nu$ is no longer para-transported at leading order. Thus, in addition to terms considered before, here we also need to contend with 
\[
\bJ^{4,err,bad,4}_{\nu}:=\int  a_{j}(x-y) T_{D_t} g_\nu (Q^j_\nu v_\lambda \ol{ w_\nu}- 
v_\lambda \ol{Q^j_\nu w_\mu}) \ol{v_\lambda}  w_\mu \chi_{r < r_\nu} \, dx dy.
\]
We estimate this at scale $r < r_\nu$ using a bilinear $L^2_{tx}$ bound for the $(\nu,\lambda)$ pairing, provided by Lemma~\ref{l:Dtg}, and a bilinear $L^2_{txy}$ bound, in order to obtain
\[
\|\bJ^{4,err,bad,3}_\nu\|_{L^1_t} \lesssim C\epsilon d_\lambda^2 d_\mu^2 \mu^{-\frac12}  
\nu^{\frac14} r^\frac12 \mu^{\frac14} = C\epsilon d_\lambda^2 d_\mu^2
r^{\frac12} \nu^\frac14\mu^{\frac14} \mu^{-\frac12} \leq  C\epsilon d_\lambda^2 d_\mu^2 \mu^{-\frac12},
\quad r < r_\nu.
\]
This is still perturbative, though we remark that the $r_\nu$ threshold is tight in this analysis.

\end{enumerate}

\subsubsection{The main error term}
\label{sec:J5-main-unbal}
We recall that this term has the form 
\begin{equation}
\bJ^{5,main} = \int a_{jk}(x-y)
(V_{<\lambda}(x) - V_{<\mu}(y))(\Mss_\lambda(v) \Pssj_\mu(w)
- \Pssj_\lambda(v) \Mss_\mu(w))
\, dx dy,
\end{equation}  
where $a_{jk}$ satisfies 
\[
|a_{jk}(x-y)| \lesssim \la x-y \ra_r^{-1}, \qquad r \gg \mu^{-1}.
\]
Of the $MP$ terms, the first has symbol $\mu^{-\frac12}$ 
and the second has symbol $\lambda^{-\frac12}$ which is smaller.
There is no coupling between the two, so without loss of generality
we will consider the first one.
Also the contributions of the corrections 
$B^3_{\lambda}$ and $Q^3_\lambda$ are similar but with a better frequency balance (a $1/2$
derivative on a low frequency), so we omit 
them and focus on the leading terms involving $\Ms$ and $\Psj$.

We consider the error bound in the dyadic region $|x-y| \approx r$, with $r \gtrsim \mu^{-1}$.
Depending on the $V$ frequency, denoted by $\gamma$, and on the separation scale $r$, we have
two potential strategies, to treat the two $V$'s separately (better for small $r$) or to exploit the difference (better for large $r$).
We first separate into cases depending on whether $\gamma \gtrsim \mu$. 

\medskip

a) $\mu \lesssim \gamma \ll \lambda$. Here we only have the $V_{\gamma}(x)$ factor, so need to look at 

\begin{equation}
\begin{aligned}
\bJ^{5,main}_{a \gamma} = & \  \int_{|x-y| \approx r} a_{jk}(x-y)
V_{\gamma}(x) \Ms_\lambda(u)(x) \Psj_\mu(v)(y)
\, dx dy
\\
\lesssim  & \ \mu^{-\frac12} \int_{|x-y| \approx r} |a_{jk}(z)|
|V_{\gamma}(x)| |u_\lambda(x)|^2 |v_\mu(y)|^2
\, dx dy.
\end{aligned}
\end{equation}  
If  we apply two bilinear $L^2_{t,xy}$ bounds for $u_\lambda v_\mu$ as in \eqref{uab-bi-unbal-trans}
we get
\begin{equation}\label{strat1}
\|\bJ^{5,main}_{a \gamma}\|_{L^1_t} \lesssim C\epsilon 
\mu^{-\frac12} \gamma^{-\frac12}
\lambda^{-2s_c} \mu^{-2s_c} \mu^{\frac12} = C \epsilon \gamma^{-\frac12} \lambda^{-2s_c} \mu^{-2s_c},
\end{equation}
which gains in $\gamma/\mu$ but has a logarithmic divergence in $r$.
Alternatively we can apply a bilinear $L^2_{t,x}$ bound for the frequency pairing  $(\gamma, \lambda)$  and an $L^2_{t,xy}$ bilinear bound for the frequency pairing $(\mu, \lambda)$. 
This gives 
\begin{equation}\label{strat2}
\|\bJ^{5,main}_{a}\|_{L^1_t} \lesssim 
C\epsilon \mu^{-\frac12} r^{-1}  \gamma^{-\frac34} r^\frac12 \mu^\frac14 \lambda^{-2s_c} \mu^{-2s_c} = C\epsilon r^{-\frac12} \gamma^{-\frac34} \mu^{-\frac14} \lambda^{-2s_c} \mu^{-2s_c} 
\end{equation}
which suffices.  Also we note that the same applies if $V_\gamma(x)$ were replaced by $V_\gamma(y)$ (with $\gamma < \mu$ now), which is needed for (b1). 

\medskip

b) $\gamma \ll \mu$. Here we have 
both terms $V_\gamma(x)$ and $V_\gamma(y)$. We split the integral dyadically relative to the separation scale $|x-r| \approx r$, and divide the analysis into large and small $r$, depending on $\gamma$. The interesting threshold, given by the uncertainty principle,  is $r = \gamma^{-1}$.
For larger $r$ we can harmlessly separate the two $V$'s, while for smaller $r$ it is essential to take advantage of the difference structure.

\medskip

b1) Large $r$, $r \gtrsim \gamma^{-1}$. Separating the two $V$'s, the easiest strategy would be to use \eqref{strat2} for both of them. This suffices if $r > \gamma^{-\frac32} \mu^{\frac12}$, 
so we are left with the range
\begin{equation}\label{r-b1}
\gamma^{-1} \lesssim  r \ll \gamma^{-\frac32} \mu^\frac12.    
\end{equation}
In this range we need to bound the 
 expression 
 \[
J_{b1} =  \int_{|x-y| \approx r} a_{jk}(z)
V_{\gamma}(x) \Ms_\lambda(v) \Ps_\mu(w)
\, dx dy,
 \]
 and its counterpart with the $V_\gamma(y)$ factor.
 
In this case we take advantage of the nonresonant cubic interaction in $x$, and 
correct the interaction Morawetz functional
as if we were simply integrating $V_\gamma$ along the para-material flow, using the 
vector field $\calv_\gamma \in Op (\coeffu_{<\gamma} S^0_\gamma) \eta_\gamma$ defined in \eqref{calv}, which satisfies a paratransport equation  
of the form (see Lemma~\ref{l:calv})
\begin{equation}
(\partial_t + V_{<\gamma} \cdot \nabla) \calv_\gamma = V_\gamma + \calN_\gamma(u) ,
\end{equation}
where $\calN_\gamma(u)$ has properties similar to the nonlinearity $N_\gamma(u)$ in our paradifferential $u$ evolution. This suggests 
using a normal form correction in our interaction Morawetz functional of the form
\begin{equation}
\bI_{b1} =  \int_{|x-y| \approx r} a_{jk}(x-y)
\calv(x) \Ms_\lambda(v) \Ps_\mu(w)
\, dx dy.
\end{equation}
For consistency, we need to evaluate $I_{b1}$ at fixed time,
\[
| \bI_{b1}| \lesssim \mu^{-\frac12} r^{-1} \|\eta_{\gamma}\|_{L^\infty}  \lambda^{-2s_c} \mu^{-2s_c} \lesssim C \epsilon c_\gamma \mu^{-\frac12} r^{-1} \gamma^{-1}  \lambda^{-2s_c} \mu^{-2s_c},
\]
which suffices for $r > \gamma^{-1}$.

On the other hand, introducing the above correction produces new error terms of the form 
\[
\bJ_{b1} - \frac{d}{dt} \bI_{b1} = \bJ_{b1}^1 + \bJ_{b1}^2+ \bJ_{b1}^3+ \bJ_{b1}^4,
\]
where 
\[
\bJ_{b1}^1 =  \int_{|x-y| \approx r} a_{j}(x-y) \calv_{\gamma}(x) \partial_k \Fsk_\lambda(v)\Psj_\mu(w)
\, dx dy
\]
\[
\bJ_{b1}^2 =  \int_{|x-y| \approx r} a_{j}(x-y)
\calv_{\gamma}(x)  \Ms_\lambda(v) \partial_y \Esj_\mu(w)
\, dx dy,
\]
\[
\begin{aligned}
\bJ_{b1}^3 = &  \int_{|x-y| \approx r} a_{jk}(x-y)\Big(
V_{<\gamma}(x) \partial_x \calv_{\gamma}(x)  \Ms_\lambda(u)(x) \Ps_\mu(v)(y) +
\calv_{\gamma}(x) \partial_x( V_{<\lambda}(x)  \Ms_\lambda(v)) \Ps_\mu(w)
\\
& \qquad  + \calv_{\gamma}(x)   \Ms_\lambda(v) \partial_y(V_{<\mu}(y)  \Ps_\mu(w)) \Big)
\, dx dy,
\end{aligned}
\]
\[
\bJ_{b1}^4 =  \int_{|x-y| \approx r} a_{jk}(x-y) \calN_{\gamma}(x)  \Ms_\lambda(v)\Ps_\mu(w)
\, dx dy.
\]
We successively bound these error terms:

\medskip

\emph{(i) The bound for $\bJ_{b1}^1$.}
In $J_{b1}^1$ we integrate by parts in $x$, 
moving the derivative on either $a$ or on $\eta$, where the latter is worse, yielding a factor of $\gamma$. Then we estimate as in \eqref{strat2}, but with an extra $\gamma^{\frac12} \lambda^{-\frac12}$ factor.
This yields
\[
\|\bJ_{b1}^1\|_{L^1_t} \lesssim C \epsilon r^{-\frac12} \gamma^{-\frac14} \mu^{-\frac14} \lambda^{-\frac12} ,
\]
which is suitable in the full $r$ range \eqref{r-b1}.

\medskip

\emph{(ii) The bound for $\bJ_{b1}^2$.} In $\bJ_{b1}^2$ we integrate by parts in $y$, moving the derivative on $a$. This gains us a factor of $r^{-1}\gamma^{-\frac12} \mu^{-\frac12}$ compared to \eqref{strat2}. Thus we obtain 
\[
\|\bJ_{b1}^2\|_{L^1_t} \lesssim C \epsilon r^{-\frac32} \gamma^{-\frac54} \mu^{-\frac34},  
\]
which is even better than in the previous case.

\medskip

\emph{(iii) The bound for $\bJ_{b1}^3$.}
In $\bJ_{b1}^3$ we integrate by parts in both $x$ and $y$. The thresholds for the $V$ truncation are different in the three terms, so we consider two cases:

\smallskip

A. $V$ is at frequency $\sigma \gtrsim \gamma$. Then we harmlessly separate the $V$'s and we need to estimate an expression of the form
\[
\bJ_{b1A}^3 =  \  \int_{|x-y| \approx r} a_{jk}(y-x)
V_{\sigma}(x) \partial_x \calv_{\gamma}(x) \,  \Ms_\lambda(v) \Ps_\mu(w) \, dx dy,
\]
as well as its analogue with  $V_\sigma(y)$ and the better term with the derivative going to $a$. Here we can use two bilinear $L^2_{t,x}$ 
corresponding to the $(\sigma,\lambda)$ and $(\gamma,\lambda)$ pairings, 
obtaining
\begin{equation}\label{J3-b1A}
\|\bJ_{b1A}^3\|_{L^1_t} \lesssim C \epsilon \mu^{-\frac12} r^{-1} \gamma^{-\frac14} \sigma^{-\frac34},
\end{equation}
which suffices when $\sigma > \gamma$, and also for smaller $\sigma$ provided that $r > \sigma^{-1}$.

\smallskip

B. $V$ is at frequency $\sigma < \gamma$ and $\gamma^{-1} < r < \sigma^{-1}$. Then the better strategy is to exploit the $V$ difference, i.e. to consider the expression
\[
\begin{aligned}
\bJ_{b1B}^3 = &  \ \int_{|x-y| \approx r} \partial a_{jk}(y-x)
(V_{\sigma}(x)-V_\sigma(y)) \calv_{\gamma}(x)  \Ms_\lambda(v) \Ps_\mu(w)\, dx dy
\\
\approx  &  \ \mu^{-\frac12} r^{-1} \int_{|x-y| \approx r} 
\nabla V_\sigma \calv_\gamma  |v_\lambda(x)|^2 |w_\mu(y)|^2 \, dx dy,
\end{aligned}
\]
where we can use again two bilinear $L^2_x$ bounds for frequency pairings $(\sigma,\lambda)$ and $(\gamma,\lambda)$.

\bigskip

\emph{(iv) The bound for $\bJ_{b1}^4$.}
Here we use the expansion in \eqref{Nuauxdecomp}. Compared with that,
in our setting we have additional $\coeffs_{<\lambda}$ paracoefficients, but these are harmless because they are uniformly bounded, and we can estimate 
their gradient in $L^\infty_t L^2_x$. 
Hence, from \eqref{Nuauxdecomp} we obtain the representation
\begin{equation*}
\calN_{\nu}(u)\,\overline \Ms_\lambda(u)
=\mathcal R_u^{\mathrm{aux}}
+\nabla\cdot Q_u^{\mathrm{flux}},
\end{equation*}
where, for every $\frac43<p\leq2$,
\begin{equation*}\label{Nuauxtrilinear+}
\|\mathcal R_u^{\mathrm{aux}}\|_{L_{t,x}^1}
+\nu^{\frac2p-1}
\|Q_u^{\mathrm{flux}}\|_{L_t^1L_x^p}
\lesssim_p
\epsilon^4c_\lambda^2\nu^{-1}
\lambda^{-2s_c}.
\end{equation*}
The contribution of $\mathcal R_u^{\mathrm{aux}}$ to $\bJ_{b1}^4$ is 
easily estimated using $r > \nu^{-1}$.
For the contribution  of $Q_u^{\mathrm{flux}}$ we integrate by parts, placing the derivative on $a$.
This  yields the expression
\[
\iint \partial^3 a Q_u^{\mathrm{flux}}(x)
\Ps(w)(y)\, dx dy dt,
\]
which we can estimate by H\"older's inequality for $r >\nu^{-1}$.

\bigskip

b2) Small $r$, precisely $r \lesssim \gamma^{-1}$. In this case we have to use the difference of the two $V$'s and represent $\bJ^{5,main}$ in the form
\begin{equation}\label{bJ-51}
\begin{aligned}
\bJ^{5,main}_{b2,\gamma} = & \   \int_{|x-y| \approx r} a_{jk}(x-y) \chi_r(x-y)
(x-y) \cdot \nabla V_\gamma(x_h) \Ms_\lambda(v) \Psj_\mu(w)
\, dx dy,
\end{aligned}
\end{equation} 
with $V$ evaluated at some intermediate point $x_h = (1-h)x+ hy$ with $h \in [0,1]$.

 A direct approach in this case is to use a bilinear $L^2_{t,x}$ bound for the frequency pairing $(\gamma,\lambda)$ and a bilinear $L^2_{t,xy}$ bound for the frequency pairing $(\mu,\lambda)$ in order
to obtain
\begin{equation}\label{low-r}
\|   \bJ^{5,main}_{b2,\gamma} \|_{L^1_t} \lesssim  C \epsilon \mu^{-\frac12} \gamma^\frac14 r^{\frac12} \mu^\frac14\lambda^{-2s_c} \mu^{-2s_c},
\end{equation}
which suffices in the range $r < (\gamma \mu)^{-\frac12}$. 
We are then left with the intermediate range for $r$
\begin{equation}\label{r-b2}
(\gamma \mu)^{-\frac12} < r < \gamma^{-1}.
\end{equation}

\bigskip

 In this case both  the two $v$ frequencies and the two $w$ frequencies are $(\gamma \mu)^\frac12$ close, which makes this $5$-linear interaction nonresonant. As in the case (b1), this allows us to use a correction for the interaction Morawetz functional. However, unlike in the case (b1), we cannot close 
 the argument in a single step, and we will instead need to use an infinite\footnote{For fixed $\gamma$ 
 the number of corrections will be finite, comparable to $\log (\mu/\gamma)$, but this can be arbitrarily large so, for practical purposes, we can think of it as infinite.}
 sum of corrections.

To facilitate the summation of the infinite series of corrections we need some preparation. 
The first difficulty we face is that our corrections involve arbitrarily many derivatives of the weight
in $x-y$, which in turn satisfy good bounds only if the weights are analytic. This conflicts with our use 
of the dyadic cutoffs $\chi_r$. To eliminate this issue, we make several adjustments:

\begin{itemize} 
\item For reiteration purposes we 
group all $r$ scales in \eqref{r-b2}
together, and switch to the cutoff $\chi_{[(\gamma \mu)^{-\frac12},\gamma^{-1}]}(x-y)$.

\item Assuming we started with $a = a_\nu$ 
with $\nu > (\gamma \mu)^{\frac12}$, we 
can harmlessly replace it with $a_{(\gamma \mu)^{\frac12}}$. This is because the difference is small
\[
|a_\nu - a_{(\gamma \mu)^{\frac12}}|
\lesssim r^{-2} (\gamma \mu)^{-1}
\]
and its contribution can be estimated directly as in \eqref{low-r} in the range \eqref{r-b2}.  
Overall, we can harmlessly replace the weight $a$ with $a_{\nu}$ where 
$\gamma < \nu <  (\gamma \mu)^{\frac12}$.

\item After the above substitution, we can expand $\bJ_{b2}$ to the full ball 
$|x-y| \lesssim \gamma^{-1}$, and only insert a truncation at $|x-y| \approx \gamma^{-1}$. This insures we have a weight which is analytic inside the ball.

\end{itemize}

At the conclusion of these adjustments, we are left with the integral 
\begin{equation}
\begin{aligned}
\bJ^{5,main}_{b2,\gamma} = & \   \int \chi_{< \gamma^{-1}}(x-y)
a^{(1)}(x-y) \nabla V_\gamma(x_h) \Ms_\lambda(v) \Psj_\mu(w)
\, dx dy,
\end{aligned}
\end{equation} 
where the universally bounded function $a^{(1)}$ is given by
\[
a^{(1)} (z) : = a_{jk}(z)z, \qquad  a = a_\nu, \qquad \nu < (\gamma \mu)^{\frac12}.
\]
We still have  non-analytic contributions 
in the derivatives of the weight, but these  are supported in the region $|x-y|\approx \gamma^{-1}$, and this will allow us to treat them as  directly perturbative at the first occurrence.

\bigskip

 We correct the above expression exactly as  in case (b1), as if we were simply para-materially integrating $V_\gamma$.  The (first) correction has the form 
\[
\bI^{(1)}_{b2} =  \int\chi_{< \gamma^{-1}}(x-y)
a^{(1)}(x-y)  
\nabla \calv_\gamma(x_h)   \Ms_\lambda(v) \Ps_\mu(w)
\, dx dy.
\]
Here we evaluate $\bI^{(1)}_{b2}$ at fixed time as before,
\[
| \bI_{b2}^{(1)}| \lesssim \mu^{-\frac12} \|\nabla \calv_{\gamma}\|_{L^\infty}  \lambda^{-2s_c} \mu^{-2s_c} \lesssim C \epsilon c_\nu \mu^{-\frac12}  c_\gamma 
\]
which suffices.

On the other hand the new error we create is 
\[
\bJ_{b2} - \frac{d}{dt} \bI^{(1)}_{b2} = \bJ_{b2}^1 + \bJ_{b2}^2+ \bJ_{b2}^3+ \bJ_{b2}^4,
\]
where 
\[
\bJ_{b2}^1 =  \int \chi_{< \gamma^{-1}}(x-y)
a^{(1)}(x-y)
 \nabla \calv_{\gamma}(x_h) \partial_x P_\lambda(v) \Ps_\mu(w)
\, dx dz,
\]
\[
\bJ_{b2}^2 =  \int \chi_{< \gamma^{-1}}(x-y)
a^{(1)}(x-y)
\nabla \calv_{\gamma}(x_h)  \Ms_\lambda(v) \partial_y E_\mu(w)(y)
\, dx dy,
\]
\[
\begin{aligned}
\bJ_{b2}^3 = & \  \int \chi_{< \gamma^{-1}}(x-y)
a^{(1)}(x-y) \big(
V_{<\gamma}(x) \partial_x  \nabla \calv_{\gamma}(x)  \Ms_\lambda(v)\Ps_\mu(w)
+ \nabla \calv_\gamma(x_h)  \partial_x(V_{<\lambda} \Ms_\lambda(v)) \Ps_\mu(w)
\\
& \ +\nabla \calv_\gamma(x_h)   \Ms_\lambda(v) \partial_y(V_{<\lambda}(y) \Ps_\mu(w))\big)
\, dx dy,
\end{aligned}
\]
\[
\bJ_{b2}^4 =  \int \chi_{< \gamma^{-1}}(x-y)
a^{(1)}(x-y)
 \nabla \calN_{\gamma}(x_h) \Ms_\lambda(v) \Ps_\mu(w)
\, dx dy.
\]

We successively bound these error terms perturbatively, except for a single contribution arising out of $\bJ^3_{b2}$, which we will need to reiterate. For the perturbative contributions we can harmlessly reinsert the dyadic cutoffs $\chi_r(x-y)$, and we will have a favourable $r$ summation. It is only the 
single nonperturbative term where this is not acceptable.

\medskip

\emph{(i) The bound for $\bJ_{b2}^1$.}
In $\bJ_{b2}^1$ we integrate by parts in $x$,
moving the derivative on $\chi_r$ in the worst case. At scale $r$ this gains us a factor of $r^{-1} \gamma^{-\frac12} \lambda^{-\frac12}$ compared to \eqref{low-r}, leading to
\begin{equation}\label{low-r1}
\|\bJ_{b2}^1\|_{L^1_t} \lesssim C \epsilon  \mu^{-\frac12} \gamma^{-\frac14} r^{-\frac12} \lambda^{-\frac12}\mu^\frac14,
\end{equation}
which is good in the desired $r$ range.

\medskip

\emph{ (ii) The bound for $\bJ_{b2}^2$.}
In $\bJ_{b2}^2$ we integrate by parts in $y$, 
moving again the derivative on $\chi_r$.
This gains us a factor of $r^{-1}\gamma^{-\frac12} \mu^{-\frac12}$ compared to \eqref{low-r}, yielding
\begin{equation}\label{low-r2}
\|\bJ_{b2}^2\|_{L^1_t} \lesssim  C \epsilon \mu^{-\frac12} \gamma^{-\frac14} r^{-\frac12} \mu^{-\frac12}\mu^\frac14,
\end{equation}
which again suffices.

\medskip

\emph{(iii) The bound for $\bJ_{b2}^4$.}
This is similar to the case (b1).
For the component of $N_\gamma(u)$
which satisfies
\begin{equation}
\|  N_\gamma(u) |w_\lambda|^2\|_{L^1_{t,x}}  \lesssim C \epsilon^2 c_\gamma^2 \gamma^{-1}  
\end{equation}
the argument is direct, omitting the
$r$ decomposition. For the remainder
$N^{[2]}_{\gamma,\lambda}(u)$ we 
gain high frequencies so we are also 
allowed a bilinear $L^2_{xy}$ bound 
for the $(\lambda,\mu)$ frequency pairing.
\medskip

\emph{(iv) The bound for $\bJ_{b2}^3$.}
In $\bJ_{b2}^3$ we integrate by parts in both $x$ and $y$.  The thresholds for the $V$ truncation are different, so we again consider three cases:

\medskip

A. $V$ is at high frequency $\sigma \gtrsim \mu$. Then we separate the $V$ and, moving the derivative on $\chi$ in the worst case, we need to estimate
\[
\bJ_{b2}^3 =  \  \int_{|x-y| \approx r} r^{-1} \chi_r(y-x)
V_{\sigma}(x) \nabla \calv_{\gamma}(x)  \Ms_\lambda(u)(x) \Ps_\mu(v)(y) \, dx dy
\]
Here we can use two bilinear $L^2_x$ 
for the $(\sigma, \lambda)$ and $(\gamma, \lambda)$ pairings,
obtaining
\[
\|\bJ_{b2}^3\|_{L^1_t} \lesssim  C \epsilon  \mu^{-\frac12} r^{-1} \gamma^{-\frac14} \sigma^{-\frac34}
\]
which suffices. 
\medskip

B. $V$ is at low frequency  $\sigma \lesssim \gamma$. Then integrating by parts in both $x$ and $y$ we obtain  the $V$ difference
\[
\begin{aligned}
\bJ_{b2}^3 = &  \  \int_{|x-y| \approx r} r^{-1}
(V_{\sigma}(x)-V_\sigma(y)) \nabla \calv_{\gamma}(x_h) V_\sigma \Ms_\lambda(v)^2 \Ps_\mu(w)
\, dx dy
\\
\approx  &  \ \int_{|x-y| \approx r} 
\nabla V_\sigma \nabla \eta_{\gamma}(x) \Ms_\lambda(u)(x) \Ps_\mu(v)(y) \, dxdy,
\end{aligned}
\]
where we can use two bilinear $L^2$ bounds involving $\sigma$ and $\gamma$ 
to obtain
\[
\|\bJ_{b2}^3\|_{L^1_t} \lesssim C \epsilon\mu^{-\frac12} \sigma^\frac14 \gamma^{-\frac14},
\]
which suffices.

\smallskip

C. $V$ is at intermediate frequency $\gamma < \sigma < \mu$. Then integrating by parts
in both $x$ and $y$ yields two terms, one with the derivative on $\calv_\gamma$,
\begin{equation}\label{b2-C1}
\bJ_{b2}^{3,1} =  \  \int\chi_{< \gamma^{-1}}(x-y)
a^{(1)}(x-y)
V_{\sigma}(x) \partial_x^2 \calv_{\gamma}(x_h)  |\Ms_\lambda(v) \Ps_\mu(w) \, dx dy,
\end{equation}
and the other with the derivative on the weight  but with a $V$ difference
\[
\bJ_{b2}^{3,2} = \ \int \nabla(\chi_{< \gamma^{-1}}(x-y)
a^{(1)}(x-y))
(V_{\sigma}(x)-V_{\sigma}(y)) \nabla \calv_{\gamma}(x_h)  \Ms_\lambda(v) \Ps_\mu(w) \, dx dy.
\]
The first term is easy to estimate perturbatively,  using two bilinear $L^2$ bounds 
to get 
\[
\|\bJ_{b2}^{3,1}\|_{L^1_t} \lesssim C \epsilon
\mu^{-\frac12} \sigma^{-\frac34} \gamma^\frac34.
\]

The second term is more difficult. Here $\nabla \calv_\gamma$ is not so helpful, and instead plays merely the role of a bounded low frequency paracoefficient. For the rest of the expression 
we distinguish two cases:

\medskip

C1: $r >\sigma^{-1}$. This is identical to (b1)
and is treated in the same manner. Notably, this includes the case when the derivative applies to $\chi$.

\smallskip

C2: $r < \sigma^{-1}$. Here we have returned to the original situation in case (b2), with the additional factor  $\nabla \calv_\gamma$, so we need to reiterate the process. For clarity, we record the remaining integral, 
\[
\bJ_{b2}^{3,2} =  \  \int \chi_{< \sigma^{-1}}(x-y)
\nabla a^{(1)}(x-y)
(V_{\sigma}(x)-V_{\sigma}(y)) \nabla \calv_{\gamma}(x_h)  \Ms_\lambda(v) \Ps_\mu(w) \, dx dy.
\]

To facilitate reiteration, we make some simplifying observations:

\begin{itemize} 
\item We can harmlessly replace $\calv(x_h)$ with $\calv(x)$ modulo a directly perturbative 
contribution in $\bJ_{b2}^{3,2}$, of the same type as in \eqref{b2-C1}. 

\item We can replace the $V_{\sigma}$ difference in terms of its gradient.
\end{itemize}

Taking into account the above observations, we are left with the integral 
\[
\bJ_{b2}^{3,2} =  \  \int  \chi_{< \sigma^{-1}}(x-y)
 a^{(2)}(x-y)
\nabla V_{\sigma}(x_h) \nabla \calv_{\gamma}(x)  \Ms_\lambda(v) \Ps_\mu(w) \, dx dy,
\]
where 
\[
a^{(2)}(z) = \partial a^{(1)}(z) z.
\]
This requires a next correction of the form  
\[
\bI^{(2)}_{b2} =  \int \chi_{< \sigma^{-1}}(x-y)
a^{(2)}(x-y)  
\nabla \calv_\sigma(x_h)  \nabla \calv_\gamma \Ms_\lambda(v) \Ps_\mu(w)
\, dx dy,
\]
and the process repeats. We generate the sequence of successive corrections 
in the reiteration of (b2) which have the form 
\begin{equation}
\bI_{b2}^{(j)}= \int  
\chi_{< \gamma_j^{-1}}(x-y)
a^{(j)}(x-y) 
\calv_{\gamma_j}(x_h) \cdots  \calv_{\gamma_1}(x) \Ms_\lambda(v) \Ps_\mu(w)
\, dx dy
\end{equation}
indexed over $j$ and the increasing sequence of frequencies $[\gamma] = \{\gamma_1, \cdots \gamma_j\}$,
\[
\gamma = \gamma_1 < \cdots < \gamma_j < \mu,
\]
where the sequence of functions $a^{(j)}$ 
is defined recurrently by
\[
a^{(j+1)}(z) = z \partial a^{(j)}(z).
\]
By analyticity these functions satisfy the bounds 
\begin{equation}
|a^{(j)}| \lesssim C^j j!
\end{equation}
with similar symbol type bounds for their derivatives. 

This implies that the added contribution of the successive corrections yields 
the correct bound with 
a factor
\[
\sum_{j \geq 1} \sum_{[\gamma]} C^j j! \epsilon^j c_{\gamma_1} \cdots c_{\gamma_j}
\lesssim \sum_{j \geq 1} C^j \epsilon^j (\sum_{\gamma < \mu} c_\gamma)^j \lesssim C\epsilon,
\]
where the smallness of $\epsilon$ insures convergence. This concludes the proof of \eqref{bJ5-bd} for the main $\bJ^5$ term.

\subsubsection{The second error term}

This is the contribution of the paradifferential $lh$ density-flux correction, of  which the worst term has the form 
\[
 \bJ^{5,[1]}_{\lambda \mu} = \int   a_{jk}(x-y)
\nu^\frac32 \lambda^{-\frac12}
B^{[1]}_{\lambda,m}
\Esj_\mu(w)\,
 dx dy.
\]
Here we recall that $B^{[1]}_{\lambda,m}$
has symbol $b^{[1]}_{\lambda,m} \in D^\frac12 \coeffu_{<\lambda} S^{-\frac12}$, so schematically we 
can express it as a sum over frequencies $\nu \ll \lambda$ of expressions of the form
\[
 \bJ^{5,[1]}_{\lambda \mu, \nu} = \int a_{jk}(x-y)
\nu^\frac32 \lambda^{-\frac12}
u_\nu |v_\lambda|^2
|w_\mu|^2
 dx dy
\]
with additional harmless $\coeff_{<\nu}$ paracoefficients.

Given a dyadic spatial separation scale $r \gtrsim\mu^{-1}$, we use
a bilinear $L^2_{tx}$ estimate for the $(\nu,\lambda)$ pairing and an $L^2_{txy}$ estimate for the $(\lambda,\mu)$ pairing. This yields
\[
|\bJ^{5,[1]}_{\lambda \mu, \nu}| \lesssim r^{-1} \nu^{\frac14} r^\frac12 \mu^{\frac14} \mu^{-1} \lambda^{-\frac12} \lesssim C \epsilon \nu^{\frac14} \mu^{-\frac14} \lambda^{-\frac12},
\]
which is satisfactory with decay in $r$ 
for $r > \mu^{-1}$. 

At this point, we have completed the proof of the unbalanced bilinear translation averaged $L^2$ bounds in the case of the paradifferential equation. In what follows, we consider the additional $\bJ^5$ error terms which arise in the case of the full equation and the linearized equation.

\subsubsection{The third error term}

This is the $bJ^{5,[0]}$ term, occuring only in the case of the linearized equation,
of which the worst term has the form 
\[
 \bJ^{5,[1]}_{\lambda \mu} = \int a_{jk}(x-y)
\nu^\frac32 \lambda^{-\frac12}
B^{[0]}_{\lambda,m}
\Esj_\mu(w)\, 
 dx dy.
\]
Here we recall that $B^{[0]}_{\lambda,m}$
has the form
\[
B^{[0]}_{\lambda,m}(v) = \lambda \calL_\nu(\bfv_\nu, \bfv_\lambda, \bfu_\lambda)
\]
with summation over $\nu \ll \lambda$, and with additional harmless $\coeffs_{<\nu}$ paracoefficients. So we need to estimate expressions of the form
\[
 \bJ^{5,[0]}_{\lambda \mu, \nu} = \lambda \int a_{jk}(x-y)
v_\nu(x) v_\lambda(x) \bu_\lambda(x)
|w_\mu(y)|^2\, 
 dx dy.
\]
For this, given a spatial localization scale $r > \mu^{-1}$, we use
a bilinear $L^2_{tx}$ estimate for the $(\nu,\lambda)$ pairing and an $L^2_{t,xy}$ estimate for the $(\lambda,\mu)$ pairing. This yields
\[
| \bJ^{5,[0]}_{\lambda \mu, \nu}| 
\lesssim C \epsilon d_\nu d_\lambda d_\mu^2 \lambda r^{-1} \nu^\frac34 r^\frac12 \mu^\frac14 \lambda^{-2}
=  C \epsilon d_\nu d_\lambda d_\mu^2 \lambda (r\mu)^{-\frac12}  (\nu\mu \lambda^{-2})^\frac34 ,
\]
which suffices. 

\subsubsection{The fourth error term} This term arises from the contribution of $N^{[2]}$ to the 
density-flux corrections. It does not occur in the case of the paradifferential flow, 
but it appears both in the analysis of the full equation and in the analysis of the linearized equation.
The analysis is similar in the two cases, but with a worse balance of frequency factors in the case of the linearized equation,
due to the fact that this is considered at lower regularity, namely $L^2$. So we discuss the linearized case, leaving the case of the full equation as an exercise for the interested reader.

To fix the scales  we assume as before that $\mu < \lambda$,
and then we can have the ``hh" interaction associated to the frequency $\lambda$ factor 
or to the frequency $\mu$ factor. In the first case, comparing weights, we identify the worst case as
\[
\bJ^{5,2,1} = \lambda^{\frac12}\nu^\frac12  \mu^{-1}   \iint \partial^2 a (x-y)  \calL_\lambda(\bfv_\lambda,\bfv_\nu,\bfu_\nu)(x) |w_\mu(y)|^2 \, dx dy, \qquad \mu < \lambda < \nu  .
\]
In the second case, there are two possible  worst case scenarios depending on how $\nu$ and $\mu$ compare, which we display together as 
\[
\bJ^{5,2,2} = (1+\nu^\frac12 \mu^{-\frac12})   \iint \partial^2 a (x-y)  |v_\lambda(x)|^2 \calL_\lambda(w_\mu,w_\nu,u_\nu)(y) \, dx dy, \qquad \mu < \lambda , \nu  .
\]
We remark that, as opposed to the first case, the second case also allows 
for the balanced scenario $\nu = \lambda$, which 
will turn out to be the most difficult to control.

We begin with the first term, where the best strategy is to use a bilinear $L^2_{t,x}$ bound for the $(\lambda,\nu)$ pair and a bilinear $L^2_{t,xy}$ for the $(\mu,\lambda)$ pair.
At scale $r$ this gives 
\[
\|\bJ^{5,2,1}\|_{L^1_t} \lesssim  C \epsilon \lambda^{\frac12}\nu^\frac12 \mu^{-1} r^{-1} \lambda^{\frac34} r^\frac12 \mu^\frac14 \nu^{-2} d_\lambda d_\nu d_\mu^2,
\]
where the last factor accounts for the critical Sobolev regularity of $u$. Using $r > \mu^{-1}$
and summing with respect to $r$ we obtain
\[
\|\bJ^{5,2,1}\|_{L^1_t} \lesssim  C \epsilon \lambda^{-\frac14}\mu^{-\frac14} (\lambda /\nu)^\frac32 d_\lambda d_\nu d_\mu^2
\lesssim  C \epsilon \lambda^{-\frac14}\mu^{-\frac14} d_\lambda^2  d_\mu^2,
\]
which suffices.

The second term is estimated in a similar fashion, except for the case when $\nu= \lambda$,
which we consider next. In this case we want to estimate
\[
\bJ^{5,2,2} = 
\nu^\frac12\mu^{-\frac12}  \int \partial^2 a(x-y)  |v_\lambda|^2(x)  w_{\mu}(y) P_\mu(w_\nu \bu_\nu)(y) \, dx dy.
\]
We evaluate it in the region $|x-y| \approx r > \mu^{-1}$.
Then we use a $L^2_{t,xy}$ 
bound for $v_\lambda w_\mu$ and another $L^2_{t,y}$ bound for $P_\mu (w_\nu u_\nu)$.
This gives 
\[
\|\bJ^{5,2,2}_r\|_{L^1_t} \lesssim  C \epsilon \lambda^{\frac12} \nu^{-\frac12} r^{-1} r^\frac12\mu^\frac14 \mu^{-\frac12} 
\lambda^{\frac54} 
 \lambda^{-2} \lesssim   C \epsilon \lambda^{-\frac{1}4} \mu^{-\frac14} d_\lambda^3 d_\nu,
\]
which suffices by the envelope property.

\subsection{ The balanced case \texorpdfstring{$ \mu \approx \lambda$}{}}
\label{sec:balanced-close}

In this subsection we use the interaction Morawetz functional in order to prove the balanced bilinear $L^2$ bounds which are part of our bootstrap argument. Again, the key part the analysis happens at the level of the paradifferential equation. The
 transition to the full equation and to the linearized equation is quite similar to the unbalanced case, and is 
abbreviated here.

\subsubsection{ The leading term}

In the balanced case the leading term is 
\begin{equation}
\bJ^4_{\lambda} = \int a_{jk}(x-y)
J^{4,jm}_\lambda(v,\bv,w,\bw) \, dx dy,
\end{equation}
where we have\footnote{Strictly speaking we should work with two different but comparable frequencies $\mu \approx \lambda$, but this does not make any significant difference while cluttering the notations.}  
\begin{equation}
J^{4,jm}_\lambda(v,\bv,w,\bw) =    \Ms_\lambda(v) \Esj_\lambda(w)
-  \Fsk_\lambda(v) \Psj_\lambda(w)
-  \Psj_\lambda(v) \Fsk_\lambda(w)
+ \Esj_\lambda(v) \Ms_\lambda(w), 
\end{equation}
with the same notation convention that 
all $v$ dependent quantities are evaluated at $x$, and all $w$ dependent quantities are evaluated at $y$. 

The quartic form $J^{4,jm}$ above has symbol
\begin{equation}\label{j4-jm}
\begin{aligned}
j^{4,jk}_{\lambda}(x,y,\xi^1,\xi^2,\xi^3,\xi^4) = & \ \ms_\lambda(x,\xi^1,\xi^2) \ms_\lambda(y,\xi^3,\xi^4)
(\fp^j_{\lambda}(x,\xi^1,\xi^2) - \fp^j_\lambda(y,\xi^3,\xi^4))
\\ & \ (\tfp^{k}_\lambda(x,\xi^1,\xi^2) - \tfp^{k}_\lambda(y,\xi^3,\xi^4)),
\end{aligned}
\end{equation}
which vanishes to second order when $x=y$ and all frequencies are equal. We recall this is 
called the
\emph{double null condition}.

Just as in the unbalanced case, one difficulty we have here is that $\fp^j \neq \tfp$, although they do agree in the constant coefficient case.
Another is that the paracoefficients (which are functions of $\nabla u$) are evaluated at $x$ and then at $y$ and we need to account for this difference.
But if $x=y$, then the above symbol vanishes of second order when all frequencies are equal.

We begin by removing some of the difficulties.  We first consider the $x-y$ differences, recalling that the paracoefficients are in $\tcoeffs_{<\lambda}$.
Taking coefficient differences in only one of the factors in \eqref{j4-jm}
retains  a simple null condition in the second factor, and yields an error of the type
\[
R^{4,1}_{\nu} \approx \lambda^{-2} \int a_{jk}(x-y) (b_{\nu}(x) -  b_\nu(y))
v_\lambda(x) \bw_\lambda(y)   \mfQ(\bv_\lambda,w_\lambda) \, dx dy ,
\]
with $b_\nu \in \coeffs_\nu$, $\nu \ll \lambda$, also with different combinations of complex conjugates and null forms $\mfQ$.
On the other hand, taking 
coefficient differences in both of the factors in \eqref{j4-jm}
destroys the null condition and yields errors of the type
\[
R^{4,2}_{\nu\gamma} \approx \lambda^{-1} \int a_{jk}(x-y) (b_{\nu}(x) -  b_{\nu}(y))
(b_{\gamma}(x) - b_{\gamma}(y)) |v_\lambda(x)|^2  |w_\lambda(y)|^2 \, dx dy ,
\]
with $\nu,\gamma \ll \lambda$.
We can estimate both of these errors perturbatively. 

\medskip

For $R^{4,1}_\nu$ we localize to a dyadic spatial scale $|x-y| \approx r \gtrsim \lambda^{-1}$, and consider two cases depending on the balance of $r$ and $\nu^{-1}$:
\medskip

i) If $r > \nu^{-1}$, then we can separate the two $b_\nu$'s, use a bilinear $L^2_{t,x}$ ( or $L^2_{t,y}$) bound 
for the $(\nu,\lambda)$ frequency pairing, and a bilinear $L^2_{t,xy}$ bound for the null form. This yields 
\[
\| R^{4,1}_{\nu}\|_{L^1_t} \lesssim C \epsilon 
c_\nu \lambda^{-1} r^{-1} \nu^{-\frac14}
r^\frac12 \lambda^\frac14 = C \epsilon 
c_\nu (r \nu)^{-\frac12} (\nu/\lambda)^{\frac14} \lambda^{-\frac12},
\]
which suffices, with dyadic summation 
in $r \gtrsim \nu^{-1}$ and $\nu$.

\medskip

(ii) If $r < \nu^{-1}$, then we rewrite 
the above expression as 
\[
R^{4,1}_{\nu,r} \approx \lambda^{-2} \int_0^1 \int O(1) \nabla b_{\nu}(x_h)
v_\lambda(x) \bw_\lambda(y)   \mfQ(\bv_\lambda,w_\lambda) \, dx dy dh,
\]
and then repeat the above strategy for fixed $h$. This yields
\[
\| R^{4,1}_{\nu}\|_{L^1_t} \lesssim C \epsilon 
c_\nu \lambda^{-1}  \nu^{\frac34}
r^\frac12 \lambda^\frac14 = C \epsilon 
c_\nu (r \nu)^{\frac12} (\nu/\lambda)^{\frac14} \lambda^{-\frac12},
\]
which again suffices.

\bigskip

For $R^{4,2}_{\nu\gamma}$ we assume without loss of generality that $\gamma < \nu$.  At spatial scale $r$, and consider two cases depending on the balance of $r$ and $\nu^{-1}$, $\gamma^{-1}$:
\medskip

i) If $r > \gamma^{-1}$, then we can separate both the $b_\nu$'s and the $b_\gamma's$, and use two bilinear $L^2_{tx}$ ( or $L^2_{ty}$) bound 
for the $(\nu,\lambda)$ and the $(\gamma,\lambda)$ frequency pairings.
 This yields 
\[
\| R^{4,2}_{\nu\gamma,r}\|_{L^1_t} \lesssim C \epsilon 
c_\nu \lambda^{-1} r^{-1} \nu^{-\frac14}
\gamma^{-\frac14} = C \epsilon 
c_\nu c_\gamma (r \nu)^{-\frac12} (r \gamma)^{-\frac12}(\nu \gamma/\lambda^2)^{\frac14} \lambda^{-\frac12},
\]
which suffices, with appropriate dyadic summation 
in $r$, $\nu$ and $\gamma$.

\medskip

(ii) If $\nu^{-1} < r < \gamma^{-1}$, then we rewrite the above expression as 
\[
R^{4,2}_{\nu,r} \approx \lambda^{-1} \int_0^1 \int O(1) \nabla b_{\gamma}(x_h)
(b_{\nu}(x) -  b_{\nu}(y)) |v_\lambda(x)|^2 |w_\lambda(y)|^2\,      dx dy dh,
\]
and then repeat the above strategy for fixed $h$. This yields
\[
\| R^{4,2}_{\nu}\|_{L^1_t} \lesssim C \epsilon 
c_\nu \lambda^{-1}  \gamma^{\frac34} \nu^{-\frac14}  = C \epsilon 
c_\nu  c_\gamma(\gamma/ \nu)^{\frac34} (\nu/\lambda)^{\frac12} \lambda^{-\frac12},
\]
which has a logarithmic summation loss in $r$, countered by the $\gamma/\nu$ 
factor.

\medskip

(iii) If $r < \nu^{-1} < r $, then we rewrite the above expression as 
\[
R^{4,2}_{\nu,r} \approx \lambda^{-1} \int_0^1\int_0^1 \int O(r) \nabla b_{\gamma}(x_h)
\nabla b_{\nu}(x_{h_1}) |v_\lambda(x)|^2 |w_\lambda(y)|^2    \,  dx dy dh dh_1,
\]
and then estimate as above. This yields
\[
\| R^{4,1}_{\nu}\|_{L^1_t} \lesssim C \epsilon 
c_\nu \lambda^{-1} r \gamma^{\frac34} \nu^{\frac34}  = C \epsilon 
c_\nu c_\gamma (r \nu)^\frac12 (r \gamma)^\frac12
 (\gamma \nu/\lambda^2)^{\frac14} \lambda^{-\frac12},
\]
which suffices.

\bigskip

Our next claim is that only the values of the  symbols $j^{4,jm}_\lambda$ on the diagonal $\Delta^4 \xi = 0$ are important. To fix the notations, consider
a general symbol $j^4$, vanishing to second order when all frequencies are equal, of the form
\[
j^{4,err} = \lambda^{-3}\Delta^4 \xi (\xi^2-\xi^3) \tcoeffs_{<\lambda} S^0.
\]
Separating variables, it suffices to consider the case when 
\[
j^{4,err} = \lambda^{-3}\Delta^4 \xi \,  (\xi^2-\xi^3) b_{\nu}, \qquad b \in \tcoeffs_{\nu}, \quad \nu < \lambda.
\]
The factor of $\Delta^4 \xi$  corresponds to a $\partial_x + \partial_y$ derivative on the quadrilinear  $(v,\bv,w,\bw)$ form. Integrating by parts, this derivative vanishes on $a(x-y)$, so it can only fall onto the coefficient $b_\nu$, which yields an expression  of the form 
\[
R^{4,3}_{\nu} \approx \lambda^{-3} \int a_{jk}(x-y) \nabla b_{\nu}(x) v_\lambda(x) w_\lambda(y)   \mfQ(\bv_\lambda  w_\lambda) \,  dx dy   .
\]
This can be estimated using a bilinear $L^2_{t,x}$ ( or $L^2_{t,y}$) bound 
for the $(\nu,\lambda)$ frequency pairing, and a bilinear $L^2_{t,xy}$ bound for the null form, to obtain
\[
\| R^{4,3}_\nu\|_{L^1_t} \lesssim C \epsilon c_\nu \lambda^{-2} r^{-1} \nu^\frac34  r^\frac12 \lambda^\frac14 =
C \epsilon c_\nu (\lambda r)^{-\frac12}
(\nu/\lambda)^\frac34 \lambda^{-\frac12}.
\]
This is sufficient.

We now take advantage of the above discussion to take a closer look at the symbol $j^{4,jm}$ on $\Delta^4 \xi = 0$, where as discussed above we can safely assume that $x = y$ in the paracoefficients.

For the discussion that follows we 
will omit the spatial variable $x$, $y$ 
from notations, as they do not play a role.

Our initial goal is to separate the variables $\xi^1$, $\xi^2$, respectively $\xi^3$, $\xi^4$, without losing the structure of $ j^{4,jm}$.  Toward this goal, we define the diagonal symbols 
\[
\fp^{j}(x,\xi) :=\fp^j(x,\xi,\xi)  = \partial_{\xi_j} \omega(x,\xi).
\]
Since $\fp^j$ is symmetric, we can relate the original and diagonal symbols by 
\[
\fp^j(\xi^1,\xi^2) = \frac12 (\fp^j(\xi^1) + \fp^j(\xi^2)) + (\xi^1-\xi^2)^2 r(\xi^1,\xi^2) .
\]
We contend that we can discard the contribution of $r$ in $j^4_{jm}$.
On $\Delta^4 \xi = 0$ this equals
\[
(\fp^{\sharp,j}(\xi^1,\xi^2) -  \fp^{\sharp,j}(\xi^3,\xi^4)) (r^m(\xi^1,\xi^2) -  r^m(\xi^3,\xi^4))(\xi^1 - \xi^2)^2,
\]
and vanishes of order four when all frequencies are equal. Here we can 
integrate by parts corresponding to one of the $\xi^1- \xi^2$ factors, with the derivative either on the paracoefficient 
as in $R^{4,3}$, or on $a$, but while still retaining the double null condition, as in 
\[
R^{4,4} \approx \lambda^{-4} \int \partial a_{jk}(x-y) \mfQ( v_\lambda(x), \bw_\lambda(y))   \mfQ(\bv_\lambda(x)  w_\lambda(y))  \, dx dy  . 
\]
This we can estimate perturbatively with two $L^2$ bilinear translation averaged null form estimates using the fact that $r \gg \lambda^{-1}$.

The same argument applies to the corresponding $\fp^\sharp$ difference, so we have reduced the expression of the principal part to 
\[
4 j^{4,jm} = (\psj(\xi^1)
+ \psj(\xi^2)-  \psj(\xi^3)-  \psj(\xi^4)) (\fp^m(\xi^1)+ \fp^m(\xi^2) -  \fp^m(\xi^3) -\fp^m(\xi^4)).
\]
We rewrite this expression in the form
\begin{equation}
\begin{aligned}
   2 j^{4,jm} =&\  (\psj(\xi^1)
 -  \psj(\xi^4)) (\fp^m(\xi^2) -  \fp^m(\xi^3))
\\
& + (\psj(\xi^2)
-  \psj(\xi^3)) (\fp^m(\xi^1) -\fp^m(\xi^4))
\\ 
& - 2 \Delta^4 \psj(\xi) \Delta^4 \fp^m(\xi).
\end{aligned}    
\end{equation}
The contribution of the last term vanishes of order four when all frequencies are equal, and can be treated perturbatively. To see this,
we note that on $\Delta^4 \xi = 0$ 
we have a quadratic expansion\footnote{Here we recall that on $\Delta^4 \xi = 0$ we have $\xi^1-\xi^2 = \xi^4-\xi^3$ and $\xi^1-\xi^4 = \xi^2-\xi^3$ so the two pairs of factors can be used interchangeably. This explains the apparent lack of symmetry here.}
\[
\Delta^4 \fp^m(\xi) = O((\xi^1-\xi^2)(\xi^1-\xi^4))
\]
with a smooth division. Hence after one integration by parts, we obtain an error term akin to $R^{4,4}$.

After dispensing with the error, we consider 
the leading term. Bypassing for the moment the 
added difficulty that $\fp \neq \tfp$, we assume equality, which yields a nonnegative density for $\bJ^4_\lambda$, and provides a bound as follows:
\begin{equation}\label{balanced-bi-r}
\int a_{jk}(x-y) (\fP^j v_\lambda (x) \bar w_\lambda(y) - v_\lambda(x) \overline{\fP^j w}_\lambda(y))    (\overline{\fP^j v}_\lambda (x)  w_\lambda(y) - \bv(x) \fP^j  w_\lambda(y))  \, dx dy \lesssim \lambda^{-\frac12}.
\end{equation}

Now we can make the usual choice for $a$ so that 
\[
D^2 a \gtrsim r^{-1} \qquad |x| \approx r \gg \lambda^{-1}.
\]
This yields control over the quantity
\begin{equation}\label{main-bi-Pj}
\int_{|x-y| \lesssim \, r} \sum_j  |\fP^j v_{\lambda} (x) \bar w_\lambda(y) - v_\lambda(x) \overline{\fP^j w}_\lambda(y)|^2   \, dx dy \lesssim \lambda^{-\frac12} r, \qquad  r \gg \lambda^{-1},
\end{equation}
where the function of $x$ and $y$ that we bound above in $L^2$ can be expressed as
\[
\fP^j v_{\lambda} (x) \bar w_\lambda(y) - v_\lambda(x) \overline{\fP^j w}_\lambda(y)
= (\fP^j(x,D_x) - \fP^j(y,D_y)) [v_\lambda (x) \bar u_\lambda(y)].
\]
Here we recall that the paradifferential operators $\fP^j\in Op \tcoeffs_{<\lambda} S^{-\frac12}$ have order $-\frac12$ and their symbols equal the
gradient of the dispersion relation of $\omega$.  As the Hessian of $\omega$ is nondegenerate, it follows that the symbols $\fp^j$ span all directions, i.e. they are jointly elliptic.
Then at the symbol level we can write
\[
\fp^j(x,\xi) - \fp^j(y,\eta) = 
a(x,\xi,\eta) (\xi-\eta) + (\fp^j(x,\xi) - \fp^j(y,\xi)) ,
\]
where $a \in Op \tcoeffs_{<\lambda} S_\lambda^{-\frac32}$ is elliptic. We claim that 
the second term above has a perturbative contribution.  Separating variables, 
we need to show that the norm 
\[
Err= \lambda^{-1} \| (b_{\nu}(x)-b_\nu(y)) v_\lambda(x) \bw_\mu(y)\|^2_{L^2(\{|x-y| < r\}}
\]
is small compared to the right hand side 
in \eqref{main-bi-Pj}. For this we consider two cases:

\bigskip

(i) $r > \nu^{-1}$. Then we separate the
two $b_\nu$'s, and use either a bilinear $L^2_{tx}$ or a bilinear $L^2_{ty}$ bound
for the $(\nu,\lambda)$ frequency pairing, combined with an energy bound for the remaining factor. This yields
\[
Err \lesssim C^2 \epsilon^2 c_\nu^2 \lambda^{-1} \nu^{-\frac12}
\lesssim r \lambda^{-\frac12} (\nu/\lambda)^\frac12 (\nu r)^{-1}
\]
which suffices.

\bigskip

(ii) $r < \nu^{-1}$. Then we pull an $r$ factor out of the $b$ difference, and then estimate as in (i),
\[
\begin{aligned}
Err \lesssim & \  \int_0^1 \lambda^{-1} r^2 \| \nabla b_{\nu}(x_h) v_\lambda(x) \bw_\mu(y)\|^2_{L^2(\{|x-y| < r\}} \, dh
\\
& 
\lesssim C^2 \epsilon^2 c_\nu^2 \lambda^{-1} r^2 \nu^{\frac32}
\lesssim r \lambda^{-\frac12} (\nu/\lambda)^\frac12 (\nu r),
\end{aligned}
\]
which again suffices.

After removing the above error term from \eqref{main-bi-Pj}, we arrive at 
\begin{equation}\label{main-bi-Pj+}
\int_{|x-y| \lesssim \, r} \sum_j  |
A(x,D_x,D_y)(\nabla_x + \nabla_y)(v_{\lambda} (x) \bar w_\lambda(y))|^2   \, dx dy \lesssim \lambda^{-\frac12} r.
\end{equation}
Since $A$ is elliptic, we consider its approximate inverse $A^{-1}(x,D_x,D_y)
\in Op \tcoeffs_{<\lambda} S_\lambda^{-\frac32}$.
Standard pseudodifferential calculus shows that 
\[
A^{-1} A = I + Op D\coeffs S_\lambda^{-1},
\]
where the contribution of the error can be estimated perturbatively in $L^2$ exactly as in (ii) above.

Then, applying $A^{-1}$ to the function 
$A(x,D_x,D_y)(\nabla_x + \nabla_y)(v_{\lambda} (x) \bar w_\lambda(y))$
in \eqref{main-bi-Pj+}, we arrive at
our main balanced bilinear $L^2$ bound,
\begin{equation}\label{bi-u-bv}
\int_{|x-y| \lesssim r} \sum_j  |\nabla v_\lambda (x) \bar w_\lambda(y) + v_\lambda(x) \nabla \bar w_\lambda(y))|^2   \, dx dy \lesssim \lambda^{\frac52} r, \qquad r \gg \lambda^{-1},
\end{equation}
which closes our main balanced $L^2$ bootstrap bound \eqref{uab-bi-bal-trans}.

\bigskip
The remaining step is to show that we can treat the 
difference $\delta \fp = \fp - \fp^\sharp$ in $j^{4,jm}$ perturbatively. In our main expression this yields an error 
\[
\begin{aligned}
j^{4,jk,err}_{\lambda}(x,y,\xi^1,\xi^2,\xi^3,\xi^4) = & \ \ms_\lambda(x,\xi^1,\xi^2) \ms_\lambda(y,\xi^3,\xi^4)
(\fp^j_{\lambda}(x,\xi^1,\xi^2) - \fp^j_\lambda(y,\xi^3,\xi^4))
\\ & \ (\delta\fp^{k}_\lambda(x,\xi^1,\xi^2) - \delta \fp^{k}_\lambda(y,\xi^3,\xi^4))
\end{aligned}
\]
To attempt to treat this expression perturbatively,  we expand in terms of the highest frequency paracoefficient
in $\delta \fp$,
\[
\delta \fp = \sum_{\nu \ll \lambda} \delta \fp_\nu,
\]
where, separating variables, we schematically write
\[
\delta \fp_\nu = g_\nu h_\nu, \qquad h_\nu \in S^{-\frac12}_\lambda, \quad g_\nu \in \coeffs_{\nu}.
\]
After peeling off negligible commutator terms as in the previous computations, its contribution to $\bJ^4$ has the form
\[
\bJ^{4,err}_\nu:=\int  a_{jk}(x-y) g_\nu (H^j_\nu v_\lambda \bw_{\lambda}- 
v \ol{H^j_\nu w_\lambda})(\ol {\fP^k v_\lambda} w_{\lambda}- 
\bv_\lambda \fP^k w_\lambda) \, dx dy.
\]
 We estimate this at scale $r$ using a bilinear $L^2_{tx}$ bound for the $(\nu,\lambda)$ pairing
and a bilinear $L^2_{xy}$ bound, obtaining
\[
\|\bJ^{4,err}_\nu\|_{L^1_t} \lesssim C\epsilon c_\nu \lambda^{-1} r^{-1} 
\nu^{-\frac14} r^\frac12 \lambda^{\frac14} = C\epsilon c_\nu
r^{-\frac12} (\nu\lambda)^{-\frac14} \lambda^{-\frac12}.
\]
We need to compare this with $\lambda^{-\frac12}$,
so this error is perturbative only if 
\begin{equation} \label{small-r-j4-err}
r \gtrsim r_{\nu}:=(\nu \lambda)^{-\frac12}.    
\end{equation}
We are left with the small $r$ range $r < r_\nu$,
\[
\bJ^{4,err,bad}_{\nu}:=\int  a_{jk}(x-y) g_\nu (H^j_\nu u_\lambda \bv_\lambda- 
u_\lambda \ol{H^j_\nu v_\lambda})(\ol {\fP^k u_\lambda} v_\lambda- 
\bu_\lambda \fP^k v_\lambda) \chi_{r < r_\nu} \, dx dy.
\]
The case when $g_\nu$ is balanced, 
i.e. $g_{\nu} \in \coeffb$
can be easily dispensed with using 
Lemma~\ref{l:bal-quad}, which
gives
\[
\|\bJ^{4,err,bad}_{\nu}\|_{L^1_t} \lesssim 
C^2 \epsilon^2 c_\nu^2 r^{-1} \nu^{-\frac12} \lambda^{-1},
\]
and which suffices exactly in the range \eqref{small-r-j4-err}. It remains to consider the unbalanced part of $g_\nu$, i.e. the case when $g_\nu \in \coeffu_\nu$.

To eliminate $\bJ^{4,err,bad}_{\nu}$ in this last case, our strategy will be, as in the unbalanced case,  to correct the interaction Morawetz functional.
We will  take the correction to be 
\[
\bI^{4,err} = -  \int  a_{j}(x-y)  g_\nu(x) (\delta\Psj_\lambda (v,\bv) \Ms_\lambda(w,\bw)- 
M^\sharp_\lambda(v,\bv) \delta \Psj_\lambda(w,\bw)) \chi_{r < r_\nu} \, dx dy,
\]
where it suffices to work with $\Ms$ and $\Ps$, rather than using the more accurate double sharp correction.
Computing its time derivative,  is not difficult to see that the flux terms generated by the $x$ and $y$ components in this correction exactly cancel $\bJ^{4,err}_\nu$. This said, there are still several additional components of the time derivative of $\bI^{4,err}$ which need to be estimated.
Taking $\delta\Psj_\lambda$ of the form 
\[
\delta\Psj_\lambda = g_\nu Q_\nu, \qquad g_\nu \in \coeffs_\nu, \quad\nu \ll \lambda,
\]
these have the form:

\begin{enumerate}[label=\roman*)]
\item terms arising from the material derivative. These are similar to $J^{5,main}$
but with $\Psj$ replaced by $\delta \Psj$,
and are estimated exactly as in Section~\ref{sec:J5-main-bal}.

\item Terms with a derivative falling on $\chi_{r < r_\nu}$ after the integration by parts.
This yields $r^{-1}$ contributions localized in $\{ r \approx r_\nu\}$, which are perturbative
according to the prior discussion.

\item Terms with a derivative falling on $g_\nu$ when integrating by parts in the flux terms. These have the schematic form
\[
\bJ^{4,err,bad,2}_{\nu} \approx \lambda^{-2} \int  a_{j}(x-y) \nabla g_\nu v_\lambda \bw_\lambda \mfQ(\bv_\lambda, w_\lambda)) \chi_{r < r_\nu} \, dx dy
\]
where the important fact is that we retain a first order null condition.
We again estimate this at scale $r < r_\nu$ using a bilinear $L^2_{tx}$ bound for the $(\nu,\lambda)$ pairing as in Lemma~\ref{l:Dtg} and a bilinear $L^2_{txy}$ bound, obtaining
\[
\|\bJ^{4,err,bad,2}_\nu\|_{L^1_t} \lesssim C \epsilon c_\nu \lambda^{-1}  
\nu^{\frac34} r^\frac12 \lambda^{-\frac14} = 
C \epsilon c_\nu \lambda^{-1}  r^{\frac12} \nu^\frac34\lambda^{-\frac14} \lambda^{-\frac12} \lesssim C \epsilon c_\nu \lambda^{-1} \nu^\frac12\lambda^{-\frac12} \lambda^{-\frac12},
\]
which is clearly perturbative with room to spare.

iii) Terms arising from the para-material derivative of $g_\nu$. This is the main difference compared with our main interaction 
Morawetz computation, as the symbol of our density $\nabla g_\nu Q_\nu$ is no longer para-transported at leading order. Thus, in addition to terms considered before, here we also need to contend with 
\[
\bJ^{4,err,bad,3}_{\nu}:=\int  a_{j}(x-y) D_{t,<\nu} g_\nu (H^j_\nu (v,\bv) \Ms_\lambda(w,\bw)- 
\Ms_\lambda(v,\bv)  H^j_\nu(w,\bw)) \chi_{r < r_\nu} \, dx dy,
\]
which also retains a single null condition.
Estimating this at scale $r < r_\nu$ using a bilinear $L^2_{tx}$ bound for the $(\nu,\lambda)$ pairing, provided by Lemma~\ref{l:Dtg}, and a bilinear $L^2_{txy}$ bound, we obtain
\[
\|\bJ^{4,err,bad,3}_\nu\|_{L^1_t} \lesssim C \epsilon c_\nu \lambda^{-\frac12}  
\nu^{\frac14} r^\frac12 \lambda^{-\frac14} = 
C \epsilon c_\nu r^{\frac12} \nu^\frac14\lambda^{\frac14} \lambda^{-\frac12} \leq C \epsilon c_\nu \lambda^{-\frac12},
\qquad r < r_\nu,
\]
which is still perturbative, but the $r_\nu$ threshold is tight now. 
\end{enumerate}

\subsubsection{The main error term}
\label{sec:J5-main-bal}
This has the form
\begin{equation}
\bJ^{5,main} = \int a_{jk}(x-y)
(V_{<\lambda}(x) - V_{<\lambda}(y))(\Mss_\lambda(v)(x) \Pssj_\lambda(w)(y)
- \Pssj_\lambda(v)(x) \Mss_\lambda(w)(y))
\, dx dy.
\end{equation}  
Here we have the key feature that the quartic symbol in the last factor 
vanishes when  all four high frequencies are equal, in other words we have a single null condition for the quartic form. This allows us to use the same argument as in the unbalanced case. 

We carefully remark that in order to obtain the above cancellation, we assume $x = y$ in the symbol of the quartic form.  To reduce the problem to to this case there is an additional error term to consider,  with no cancellation but with two paracoefficients, of the form
\begin{equation}
\bJ^{5,main,err} \approx \lambda^{-\frac12}\int
\partial^2 a(x-y) (V_{<\lambda}(x) - V_{<\lambda}(y))(g_{<\lambda}(x) -  g_{< \lambda}(y)) |v_\lambda(x)|^2  |w_\lambda(y)|^2
\, dx dy
\end{equation}  
for a paracoefficient $g \in \coeffs$.

To to bound an expression of this type,
we denote the two low frequencies for $V$, respectively $g$ by $\nu$ and $\gamma$, where the worst case is $\gamma < \nu$. Assuming this and given a spatial distance scale $r$, 
we consider three cases:
\medskip

(i) $r > \gamma^{-1}$.  In this case we can neglect the differences and directly use two bilinear $L^2$ bounds for the frequency pairs
$(\nu,\lambda)$ respectively $(\mu,\lambda)$.
We obtain
\[
\|\bJ^{5,main,err}\|_{L^1_t} \lesssim C^2 \epsilon^2
c_{\nu} c_\gamma r^{-1} \nu^{-\frac34} \gamma^{-\frac14} \lambda^{-\frac12}.
\]

\medskip

(ii) $\nu^{-1}<r < \gamma^{-1}$. Here we 
take advantage of the $V$ difference to rewrite our error term as a superposition 
of expressions of the form
\[
\bJ^{5,main,err} \approx \lambda^{-\frac12}\int
 \nabla V_{\nu}(x_h)  g_{<\lambda}  |v_\lambda(x)|^2  |w_\lambda(y)|^2
\, dx dy.
\]
Using here two bilinear $L^2$ bounds as above 
yields
\[
\|\bJ^{5,main,err}\|_{L^1_t} \lesssim C^2 \epsilon^2
c_{\nu} c_\gamma  \nu^{\frac14} \gamma^{-\frac14} \lambda^{-\frac12},
\]
which again suffices, with the logarithmic 
summation loss on $r$ countered by the $\nu/\gamma$ ratio.

\medskip
(ii) $r < \nu^{-1}$. Here we 
take advantage of both differences to rewrite our error term as a superposition 
of expressions of the form
\[
\bJ^{5,main,err} \approx \lambda^{-\frac12} r \int
 \nabla V_{\nu}(x_h)  \nabla g_{\gamma}(x_{h_1})   |v_\lambda(x)|^2  |w_\lambda(y)|^2
\, dx dy.
\]
Now two bilinear $L^2$ bounds 
yield
\[
\|\bJ^{5,main,err}\|_{L^1_t} \lesssim C^2 \epsilon^2
c_{\nu} c_\gamma  r \nu^{\frac14} \gamma^{\frac44} \lambda^{-\frac12},
\]
which directly suffices.

Finally we estimate the 
errors arising from the $B$ correction 
terms. Since these are of many types,
it is convenient to place them under 
a common umbrella, where the null condition is no longer needed. This is provided by Lemma~\ref{l:Bxt}, whose use we discuss in two cases:

\medskip

a) The double $B$ term, i.e. an expression of the form
\begin{equation}
\bJ^{5,main,err} \approx \int
\partial^2 a(x-y) (V_{<\lambda}(x) - V_{<\lambda}(y)) B_{\lambda,m}(x)  B_{\lambda,p}(y)
\, dx dy.
\end{equation} 
We denote by $\nu \ll \lambda$ the frequency of $V$, and separate small and large spatial scales. For $r > \nu^{-1}$, we use the uniform bound for $V$,
\[
\begin{aligned}
\|\bJ^{5,main,err}\|_{L^1_t} \lesssim & C \epsilon c_\nu 
\nu^{-\frac12} \int_{|x-y|> \nu^{-1}} 
\frac{1}{|x-y|}  |B_{\lambda,m}(x) || B_{\lambda,p}(y)|
\, dx dy
\\
\lesssim & \ C \epsilon c_\nu 
\nu^{\frac12} \|B_{\lambda,m}\|_{L^2 L^1} \|B_{\lambda,p}\|_{L^2L^1}
\\
\lesssim & \ 
C \epsilon c_\nu 
\nu^{\frac12} \lambda^{-1},
\end{aligned} 
\]
which suffices.

For $r < \nu^{-1}$, we use the gradient of $V$
instead, still with a uniform bound. This yields
\[
\begin{aligned}
\|\bJ^{5,main,err}\|_{L^1_t} \lesssim & C \epsilon c_\nu 
\nu^{\frac12} \int_{|x-y| < \nu^{-1}} 
 |B_{\lambda,m}(x) || B_{\lambda,p}(y)|
\, dx dy
\\
\lesssim & \ C \epsilon c_\nu 
\nu^{\frac12} \|B_{\lambda,m}\|_{L^2 L^1} \|B_{\lambda,p}\|_{L^2L^1}
\\
\lesssim & \ 
C \epsilon c_\nu 
\nu^{\frac12} \lambda^{-1},
\end{aligned} 
\]
exactly as in the first case.

b) The single $B$ expression, for instance
i.e. an expression of the form
\begin{equation}
\bJ^{5,main,err} \approx \int
\partial^2 a(x-y) (V_{<\lambda}(x) - V_{<\lambda}(y)) B_{\lambda,m}(x)  \Ps_\lambda(w)
\, dx dy.
\end{equation} 
With the $V$ frequency still denoted by $\nu
\ll \lambda$, we distinguish the same cases based on the spatial separation scale $r$. 
If $r > \nu^{-1}$ then we keep the $V$'s separate. For the $V_\nu(y)$ expression 
we directly use an $L^2_{ty}$ bilinear bound 
to get
\[
\begin{aligned}
\|\bJ^{5,main,err}\|_{L^1_t} \lesssim & \ r^{-1} \| B_{\lambda,m}\|_{L^2_t L^1_x} 
\| V_{<\lambda}(y)   \Ps_\lambda(w)\|_{L^2_t L^1_y} 
\\ 
\lesssim & \ C \epsilon c_\nu r^{-1} \lambda^{-\frac14}  \nu^{-\frac32} \nu^{\frac34} \lambda^{-\frac12} 
\\
\lesssim & \ 
C \epsilon c_\nu (r\nu)^{-1}  (\nu/\lambda)^{\frac14} \lambda^{-\frac12}.
\end{aligned}
\]
For the $V_\nu(x)$ expression we instead 
 use an $L^2_{txy}$ bilinear bound followed by Bernstein in $x$ at frequency $\nu$ to get
 \[
\begin{aligned}
\|\bJ^{5,main,err}\|_{L^1_t} \lesssim & \ r^{-1} \| B_{\lambda,m}\|_{L^2_t L^1_x} 
\| \chi_r(x-y) V_{\nu}(x)   \Ps_\lambda(w)\|_{L^2_t L^\infty_x L^1_y} 
\\ 
\lesssim & \ r^{-1} \nu
\| B_{\lambda,m}\|_{L^2_t L^1_x} 
\| \chi_r(x-y) V_{\nu}(x)   \Ps_\lambda(w)\|_{L^2_t L^2_x L^1_y} 
\\
\lesssim & \ C \epsilon c_\nu r^{-1} \nu \lambda^{-\frac14}  \nu^{-\frac32} r^\frac12 \nu^{\frac14} \lambda^{-\frac12} 
\\
\lesssim & \ 
C \epsilon c_\nu (r\nu)^{-\frac12}  (\nu/\lambda)^{\frac14} \lambda^{-\frac12},
\end{aligned}
\]
same as before.

If $r < \nu^{-1}$ then we use instead the gradient of $V_\nu$ and estimate the expression 
\begin{equation}
\bJ^{5,main,err} \approx \int_0^1\int \chi_{r}(x-y)
 V_{<\lambda}(x_h) B_{\lambda,m}(x)  V_{<\lambda}(x_h) \Ps_\lambda(w)
\, dx dy dh
\end{equation} 
in the same manner as before:
\[
\begin{aligned}
\|\bJ^{5,main,err}\|_{L^1_t} \lesssim & \  \| B_{\lambda,m}\|_{L^2_t L^1_x} 
\| \chi_{<\nu^{-1}}(x-y) \nabla V_{\nu}(x_h)   \Ps_\lambda(w)(y)\|_{L^2_t L^\infty_x L^1_y} 
\\ 
\lesssim & \  \nu
\| B_{\lambda,m}\|_{L^2_t L^1_x} 
\| \chi_{<\nu^{-1}}(x-y) V_{\nu}(x_h)   \Ps_\lambda(w)(y)\|_{L^2_t L^2_x L^1_y} 
\\
\lesssim & \ C \epsilon c_\nu  \nu \lambda^{-\frac14}  \nu^{-\frac12} r^\frac12 \nu^{\frac14} \lambda^{-\frac12} 
\\
\lesssim & \ 
C \epsilon c_\nu (r\nu)^{\frac12}  (\nu/\lambda)^{\frac14} \lambda^{-\frac12}.
\end{aligned}
\]
This again suffices.

\subsubsection{The second error term}

This is the paradifferential $lh$ term. We group similar expressions, namely those arising from the $B_m$ corrections and those arising from the $B_p$ corrections.
In the case of the paradifferential equation 
all such corrections have phase rotation symmetry, so we begin with this case.
We consider one of the terms,
\[
\begin{aligned}
\bJ^{5,[1]} = \iint a_{jk}(x-y)&  \left(  
B^{[1],j}_{\lambda,p}(v) \Psk_\lambda(w)
-  Q^{[1],jk}_{\lambda,p}(v) \Ms_{\lambda}(w) \right.
\\ & \left.
-  \Ms_{\lambda}(v) Q^{[1],jk}_{\lambda,p}(w)
+ 
\Psk_\lambda(v) B^{[1],j}_{\lambda,m}(w)
 \right) \, dx dy ,
 \end{aligned}
\]
or the similar expression with $B^{[1]}$ and $Q^{[1]}$ replaced by $B^{[3]}$ and $Q^{[3]}$.
This is multilinear in $u_\nu,u_\lambda, u_\lambda, v_\lambda,v_\lambda$ with $\nu \ll \lambda$, and $\coeffuu$ paracoefficients. 
Furthermore, we know that by \eqref{b-q-relation}, for the corrections with phase rotation symmetry we have the symbol relation 
\begin{equation}\label{b-q[1]}
q^{[1],jk}_{\lambda,p} = \fp^j_\lambda b^{3,k}_{\lambda,p} ,
\end{equation}
therefore the above quadrilinear form satisfies the double null condition.

Hence, schematically, 
this expression has the form 
\[
 \lambda^{-2} \nu^\frac32 \iint \partial^2 a(x-y) 
v_\nu(x) v_\lambda(x) \bar w_\lambda(y) \mfQ( 
w_\lambda(y), \bu_\lambda(x))  \, dx dy.
\]
Consider a dyadic spatial separation scale $r \gtrsim\lambda^{-1}$. Then the above expression can be estimated at scale $r$ using a  bilinear $L^2_{tx}$ bound for the 
$(\nu,\lambda)$ frequency pair, and a 
balanced bilinear $L^2_{txy}$ bound.
\[
\|\bJ^{5,2}\|_{L^1_t} \lesssim C \epsilon c_\nu r^{-1} \nu^{\frac14} r^\frac12 \lambda^{\frac14} \lambda^{-1} \lambda^{-\frac12} \lesssim C \epsilon c_\nu (r\lambda)^{-1} (\nu/\lambda)^{\frac14}  \lambda^{-\frac12},
\]
which is satisfactory with decay in $r$ 
for $r > \lambda^{-1}$.

For the full equation and the linearized equation, we also have terms without phase rotation symmetry, e.g. $B^{[1]}(v_\lambda,v_\lambda)$ and $Q^{[1]}(v_\lambda,v_\lambda)$. The downside here is that we lose the symbol relation \eqref{b-q[1]}. However, by using tighter frequency localization cutoffs as in Remark~\ref{lp-tight}, we can insure that the two $v_\lambda$ entries have nearly opposite
frequency localizations, rather than nearly equal. This insures transversality in this case, and eliminates any need for a null condition.

At this point, we have completed the proof of the balanced bilinear translation averaged $L^2$ bounds \eqref{uab-bi-bal-trans} in the case of the paradifferential equation. In what follows, we consider the additional $\bJ^5$ error terms which arise in the case of the full equation and the linearized equation.

\subsubsection{The third error term}

This is the $bJ^{5,[0]}$ term, occurring only in the case of the linearized equation,
 which  has the form 
\[
\begin{aligned}
\bJ^{5,[0]} = \iint a_{jk}(x-y)&  \left(  
B^{[0],j}_{\lambda,p}(v) \Psk_\lambda(w)
-  Q^{[0],jk}_{\lambda,p}(v) \Ms_{\lambda}(w) \right.
\\ & \left.
-  \Ms_{\lambda}(v) Q^{[0],jk}_{\lambda,p}(w)
+ 
\Psk_\lambda(v) B^{[0],j}_{\lambda,m}(w)
 \right) \, dx dy .
 \end{aligned}
\]
This is similar to the previous case, using the fact that that the symbols for $B^{[0]}_\lambda$ and $Q^{[0]}_\lambda$ are related 
via \eqref{bq-relation-lin}, which guarantees that the above quadrilinear expression still has  a null structure.

\subsubsection{The fourth error term}
Here  we  have the ``hh" interaction associated to one of the frequency $\lambda$ factors 
\[
\bJ^{5,3} = \lambda^{-\frac12}\nu^\frac12    \iint \partial^2 a (x-y)  L(v_\lambda,u_\nu,v_\nu)(x) |w_\lambda(y)|^2 \, dx dy, \qquad \lambda < \nu  .
\]

Here the best strategy is to use a bilinear $L^2_{t,x}$ bound for the $(\lambda,\nu)$ pair, and a bilinear $L^2_{t,xy}$ bound for a second $(\lambda,\nu)$ pair.
This gives 
\[
\|\bJ^{5,3,1}\|_{L^1_t} \lesssim  \lambda^{\frac12}\nu^\frac12 \lambda^{-1} r^{-1} \lambda^{\frac34} r^\frac12 \lambda^\frac14 \lambda^2
\nu^{-4},
\]
where the last two factors account for the critical Sobolev regularity. Using $r > \nu^{-1}$
we obtain
\[
\|\bJ^{5,2,1}\|_{L^1_t} \lesssim  \lambda^{-\frac14}\nu^{-\frac14} (\lambda /\nu)^{7/2},
\]
which suffices. 

The strategy applies always when $\lambda \ll \nu$. However, there is also the case 
when $\nu \approx \lambda$. There 
all five frequencies are comparable, and 
it is possible to have four nearly equal frequencies, though not five. Still, in this case we have a single guaranteed transversality and we cannot use two bilinear
$L^2$ bounds. But then the five-linear interaction is nonresonant, so in particular 
the expression $L(v_\nu, u_\nu,v_\nu)$ 
is nonresonant and can be corrected. 
This computation is similar to the $B^{[2]}$
computation, and is omitted.

\section{Local and global well-posedness} \label{sec:rough}

\subsection{ Global solutions for  regular data}

Our starting point here is the classical local well-posedness in regular Sobolev spaces. The result in \cite{ABZ}, for instance, yields local solutions in $\H^s$ for $s > s_c+\frac12$.  

 Consider now initial data $(\eta_0,\psi_0)$ and associated good variable $\varphi_0$ with $\H^s$ regularity and satisfying also the smallness condition 
\eqref{small-data}. Our result in Theorem~\ref{t:local-fe3} guarantees 
that $\mathcal B^{s_c}$ frequency envelope bounds propagate on the time of existence of 
\cite{ABZ}. This gives uniform control of the $\H^s$ norm of the 
solution at high frequency, while the conserved energy gives control over the 
low frequency. Thus the $\H^s$ norm of the solution remains uniformly bounded
over the time of existence, which in turn implies that the solution is either global or extendable.
Thus a continuity  argument shows that the solution is global in time, with uniform $\H^s$ bounds, and satisfying all the  bounds in Theorem~\ref{t:local-fe3}
globally in time.

\subsection{ Global solutions for rough data}

Our aim here is to construct rough solutions as limits of smooth solutions.
We consider initial as follows:
\begin{enumerate}
    \item satisfies the smallness condition 
\eqref{small-data},
\item has $\dH^{s_c+\frac1{12}}$ regularity.
\end{enumerate}
We denote by $\left\{c_k\right\}_k$ a frequency envelope
for this initial data in $\dH^{s_c}$.

For this data we consider a smooth family of regularizations $(\eta^k_0,\psi^k_0)$ defined by 
\begin{equation}
  (\eta^k_0,\psi^k_0) :=   (P_{<k} \eta_0,P_{<k} \psi_0).
\end{equation}

\bigskip

\subsubsection{Regularized data}  We show that the regularized data sits under the same frequency envelope, with good difference 
bounds:

\begin{lemma}
For initial data $(\eta_0,\psi_0)$ as above, the associated regularized    
data satisfies uniform frequency envelope bounds:
\begin{equation}\label{data-k}
   \|  (P_j \eta^k_0, P_j \varphi^k_0)\|_{\dot H^{s_c}} \lesssim c_k 2^{-N(j-k)_+}.
\end{equation}
Furthermore, the initial good linearized variable
\begin{equation*}
(h^k_0,r^k_0) = (\partial_k \eta^k_0, 
\partial_k \psi^k_0- T_{B^{l}_0} \partial_k \eta^k_0)
\end{equation*}
satisfies
\begin{equation}\label{data-k-lin}
 \|  (h^k_0,r^k_0) \|_{\H^0} \lesssim 
 c_k 2^{-ks}.
\end{equation}
\end{lemma}
\begin{proof}
For \eqref{data-k} it suffices to consider the second component, for which we write
\[
\varphi_0^k = P_{<k} \psi_0 - T_{B^{k,l}_0} P_{<k} \eta_0=   P_{<k} \varphi_0 +
P_{<k} T_{B^l_0} \eta_0 -  T_{B^{l}_0} P_{<k} \eta_0 =  P_{<k} \varphi_0 +
[P_{<k}, T_{B^l_0}] \tilde P_{k} \eta_0.
\]
The commutator is localized at frequency $2^k$, and we can estimate it by
\[
\| [P_{<k}, T_{B^l_0}] \tilde P_{k} \eta_0\|_{\dot H^{s_c+\frac12}}
\lesssim 2^{-\frac{k}2} \| \nabla B^l_{0,<k}\|_{L^\infty} \|\tilde P_{k} \eta_0\|_{ \dot H^{s_c}} \lesssim \epsilon c_k.
\]

For \eqref{data-k-lin} it also suffices to consider the second component, for which we compute
\[
\begin{aligned}
r_0^k = & \ P_k \psi_0 - T_{B^l_0} P_k \eta_0 
\\
= & \  P_k \varphi_0 + P_k T_{B_0^l} \eta_0
- T_{B^l_0} P_k \eta_0 
\\
= & \  P_k \varphi_0 + [P_k, T_{B_0^l}] \eta_0.
\end{aligned}
\]
Now we estimate the commutator as before, 
and the proof is concluded.

\end{proof}
By the above Lemma, the regularized data $(\eta^k_0,\psi^k_0)$ yield global in time smooth solutions, denoted by 
$(\eta^k,\psi^k)$,
which satisfy uniformly the estimates in  Theorem~\ref{t:local-fe3}.

\subsubsection{Convergence of regularized solutions} 
Our goal here is to show that the solutions $(\eta^k,\psi^k)$ converge uniformly in the strong $\H^s$ topology
provided that $s \geq s_c+\frac1{12}$.

From Theorem~\ref{t:local-fe3} we have the uniform frequency envelope bounds 
\begin{equation}\label{uniform-uk}
\| P_j (\eta^k,\varphi^k) \|_{\dH^{s_c}}
\lesssim c_k 2^{-N(j-k)_+},
\end{equation}
globally in time, as well as a full set of linear and bilinear bounds.
On the other hand from Theorem~\ref{t:local-fe3-lin} we have 
the linearized bounds
\begin{equation}\label{uniform-diffk}
\| (h^k,r^k) \|_{\H^0}
\lesssim c_k 2^{-sk},
\end{equation}
locally in time. 

We seek to combine the two sets of bounds in order to show that 
\begin{equation}
 \| P_j (\eta^{k+1}-\eta^k,\varphi^{k+1} -\varphi^k) \|_{\dH^{s_c}}
\lesssim c_k 2^{-2|j-k|},   
\end{equation}
which in turn suffices to show convergence in $\dot{\mathcal B}^{s_c}$, 
and also in $\H^s$.

The high frequencies $j > k$ are controlled directly from \eqref{uniform-uk}. For the low frequencies we seek instead to use the bound \eqref{uniform-diffk}. For the first component, we have 
$\frac{d}{dk} \eta^k = \psi^k$, which suffices. We consider the second component, where we write
\[
\begin{aligned}
\frac{d}{dk} \varphi^k = &\  \psi^k - T_{B^{k,l}} 
\eta^k - T_{\frac{d}{dk} B^{k,l}} \eta^k
\\
= &\  r^k  - T_{\frac{d}{dk} B^{k,l}} \eta^k.
\end{aligned}
\]
The first term is controlled from 
\eqref{uniform-diffk}.  It remains to consider the second term, for which it suffices to 
show that we have the bound
\begin{equation}\label{dB-H12}
\| \frac{d}{dk} B^{k,l}_{[j]}\|_{\dot H^{-\frac12}} \lesssim
2^{-k s_c},
\end{equation}
uniformly in $j$,
where we recall that 
\[
B^{k,l}_{[j]} = B(P_{<j} \eta^k,
P_{<j} \psi^k) = L(P_{<j} \eta^k) P_{<j} \psi^k.
\]
Using the formula \eqref{dML} for the linearization of $L$ we compute
\[
\begin{aligned}
\frac{d}{dk} B^{k,l}_{[j]} = & \ - L(P_{<j}\eta^k) (\partial_k P_{<j} \eta^k
B^{k,l}_{[j]}) + \partial_k P_{<j} \eta_k L(P_{<j}\eta^k) B^{k,l}_{[j]}
+ L(P_{<j}\eta^k) \partial_k P_{<j} \psi^k
\\
=& \  L(P_{<j}\eta^k) ( 
P_{<j} (r^k + T_{B^{k,l}} h^k)
-P_{<j} h^k B^{k,l}_{[j]} ) 
+P_{<j} h^k L(P_{<j} \eta^k)B^{k,l}_{[j]}
\\
=& \  L(P_{<j}\eta^k)  
P_{<j} r^k + L(P_{<j}\eta^k) 
P_{<j} (
T_{B^{k,l}} h^k - h^k B^{k,l}_{[j]} ) 
+P_{<j} h^k L(P_{<j} \eta^k)B^{k,l}_{[j]}.
\end{aligned}
\]
Here $B^{k,l}_{[j]}$ are estimated uniformly in 
$\dot H^\frac32$ uniformly in $j$.
Since $L(P_{<j} \eta^k)$ is an order $1$ operator, the desired bound \eqref{dB-H12}
directly follows for the first and the last term. 
It remains to consider the middle term, where it 
would suffice to show that 
\[
\| P_{<j} (
T_{B^{k,l}} h^k - h^k B^{k,l}_{[j]} ) \|_{\dot H^\frac12}
\lesssim \| h^k\|_{L^2}.
\]
We can directly bound 
\[
\| h^k B^{k,l}_{[j]} - T_{B^{k,l}_{[j]}} h^k\|_{\dot H^\frac12} \lesssim \|h_k\|_{L^2} \|B^{k,l}_{[j]}\|_{\dot H^{-\frac32}}, 
\]
so it remains to show that 
\[
\| P_{<j} 
T_{B^{k,l}- B^{k,l}_{[j]}} h^k \|_{\dot H^\frac12}
\lesssim \| h^k\|_{L^2}.
\]
This in turn requires the uniform bound 
\[
\| B^{k,l}_{[j_1]}- B^{k,l}_{[j]}\|_{L^\infty} \lesssim  2^{-\frac{j_1}2}, \qquad j_1 < j, 
\]
which is a consequence 
of the bound 
\[
\|\frac{d}{dj} B^{k,l}_{[j]}\|_{L^\infty} \lesssim  \epsilon c_j 2^{-\frac{j}2}.
\]
Indeed, writing $L=L(P_{<j}\eta^k)$, we have
\[
\begin{aligned}
\frac{d}{dj} B^{k,l}_{[j]}
= & - L(\eta_j^k B^{k,l}_{[j]})
+\eta_j^k L B^{k,l}_{[j]}+L\psi_j^k
\\
= & L\varphi_j^k+\eta_j^k L B^{k,l}_{[j]}
-L\bigl((B^{k,l}_{[j]}-T_{B^{k,l}})\eta_j^k\bigr)
+L[P_j,T_{B^{k,l}}]\eta^k.
\end{aligned}
\]
The desired estimate follows by the same argument used to prove
\eqref{BBern} in Lemma \ref{Phiest}, using Bernstein's inequality
and Proposition \ref{Bernsteinelliptic} to control the high
frequency tails. The last commutator is estimated as in
\eqref{PhiCommGeneral}.

\subsubsection{ Weak Lipschitz dependence}

The objective here is to prove that, for smooth solutions as in Theorem~\ref{t:local-fe3}, we have

\begin{lemma}\label{l:weak-Lip}
Let $(\eta^1,\psi^1)$ and $(\eta^2,\psi^2)$ be solutions as in Theorem~\ref{t:global}, and $\varphi^1,\varphi^2$ be the associated good variables. Then we have
\begin{equation}\label{eq:weak-Lip}
 \|(\eta^2-\eta^1,\varphi^2-\varphi^1)\|_{L^{\infty}([0,T];\H)} \lesssim C(T) \|(\eta^2-\eta^1,\varphi^2-\varphi^1)(0)\|_{\H}    .
\end{equation}
    
\end{lemma}

As a corollary of this lemma, we conclude 
that our rough solutions are the unique limits of smooth solutions, which is our main uniqueness statement for rough solutions.

\begin{proof}
Since our rough solutions are uniformly approximated by regular solutions in the strong topology, it suffices to prove the bound \eqref{eq:weak-Lip} for regular solutions. 

We consider a linear path in the initial data space for $(\eta,\varphi)$, namely 
$\dot{\mathcal B}^{s_c} \cap \dH^{s_c+\frac1{12}}$, denoted by 
\[
(\eta^{\sigma}_0,\varphi^{\sigma}_0):= 
(2-\sigma)(\eta^1_0,\varphi^1_0)+(\sigma-1)(\eta^2_0,\varphi^2_0), \qquad \sigma \in [1,2],
\]
and denote by $(\eta^{\sigma},\varphi^{\sigma})$ the associated global solutions.

We also denote by $(h^\sigma,r^\sigma)$
the associated linearized variables.
We first claim that at the initial time we have
\begin{equation} \label{good-lin-data}
\| (h^\sigma_0,r^\sigma_0)\|_{L^{\infty}([0,T]:\H)} \lesssim 
 \|(\eta^2-\eta^1,\varphi^2-\varphi^1)(0)\|_{\H}.
\end{equation}
We only need to prove the bound for $r^\sigma_0$, which is represented as follows:
\[
\begin{aligned}
 r^\sigma_0 = & \ \frac{d}{d\sigma} \psi^\sigma_0 - T_{B^\sigma_0} h^\sigma_0   
 \\
 = & \ \frac{d}{d\sigma} \varphi^\sigma_0
 + T_{\frac{d}{d \sigma} B^\sigma_0} \eta_0^\sigma,
\end{aligned}
\]    
where computing as before we have
\[
\begin{aligned}
\frac{d}{d\sigma} B_0^\sigma 
= & \ L(\eta_{0}^\sigma) \partial_\sigma \varphi_{0} 
+ L(\eta_{0}^\sigma) (T_{B_{0}^\sigma}- B_0^\sigma) \partial_\sigma \eta_{0}
+ \partial_\sigma \eta_{0}^\sigma L(\eta_{0}^\sigma) B_0^\sigma .
\end{aligned}
\]
Using the $\dot H^{s_c-\frac12}$ bound for $B^\sigma_0$ and elliptic bounds for $L$, we obtain 
\begin{equation}
 \|    \frac{d}{d\sigma} B_0^\sigma\|_{\dot H^{-\frac12}} \lesssim \| (\partial_\sigma \eta_{0},\partial_\sigma \varphi_{0})\|_{\H}.
\end{equation}
This easily implies \eqref{good-lin-data}.

Propagating the bound \eqref{good-lin-data} in time via Theorem~\ref{t:local-fe3-lin}, it remains to show that at a time $t$ we can recover $(\partial_\sigma \eta,\partial_\sigma \varphi)$ from 
$(h^\sigma, r^\sigma)$ in the $\H$ topology. 

Reversing the previous computation, 
it suffices to establish the bound
\[
 \|    \frac{d}{d\sigma} B_0^\sigma\|_{\dot H^{-\frac12}} \lesssim \| (h^\sigma,q^\sigma)\|_{\H}.
\]
Here we use
\[
\begin{aligned}
\frac{d}{d\sigma} B^\sigma 
= & \ L(\eta^\sigma)r^\sigma 
+  h^\sigma L(\eta^\sigma) B^\sigma,
\end{aligned}
\]
and elliptic bounds for $L$.
\end{proof}

\subsubsection{Continuous dependence}
To show that our rough solutions have 
continuous dependence on the initial data
in the strong topology, it suffices to refer to the primer \cite{IT-primer} and the abstract result of \cite{ABITZ}, which allow us to conclude that continuous dependence follow from the 
higher energy bounds of Theorem~\ref{t:local-fe3} extended to rough solutions (tame estimates), and the weak-Lipschitz bounds of Lemma~\ref{l:weak-Lip}. 

\section{The higher dimensional case}
\label{sec:hid}

Our objective here is to provide a short comparative discussion of the proof of our result in higher space dimension $n\geq 4$. 
This can be seen as a strict subset of the 
$n=3$ case, so we will describe the
simplifications, distinguishing the cases 
$n = 4$ and $n \geq 5$. 

\subsection{ Changes in the bootstrap 
bounds} To begin with, the bootstrap bounds 
\eqref{uk-ee-boot}, \eqref{uab-bi-unbal-trans-boot}, \eqref{uab-bi-bal-trans-boot}
and \eqref{Taylorcoeffsmall} remain completely unchanged. The  bilinear $L^2_{t,x}$  bounds change as follows, in a manner consistent with scaling:

\begin{enumerate}[label=(\roman*)]
\item Unbalanced bilinear $L^2_{t,x}$ bounds: Let $u_{\mu}$ and $u_{\lambda}$ be unbalanced, with $\mu\ll \lambda$. Then 
\begin{equation} \label{uab-bi-unbal4}
\| u_\lambda \bu_\mu^{x_0}  \|_{L^2_{t,x}} \lesssim \epsilon^2 c_\lambda c_\mu \lambda^{-s_c} \mu^{-s_c+\frac{3}{4}+\frac{d-2}2}(1+\mu |x_0|)^\frac12;
\end{equation}

\item Balanced bilinear $L^2_{t,x}$ bound: 
\begin{equation} \label{uab-bi-bal4} 
\|P_{\alpha} (u_\lambda \bu_\mu^{x_0})  \|_{L^2_{t,x}} \lesssim \epsilon^2 c_\lambda c_\mu \lambda^{-2s_c+\frac{5}{4}}  \alpha^{\frac{d-3}2} (1+ \alpha |x_0|)^\frac12
, \qquad \alpha \lesssim \lambda \approx \mu .
\end{equation}
\end{enumerate}
These are obtained simply by making obvious dimension adjustments in the proof of Lemma~\ref{lem:boot}.

The mixed $L^p_{t,x}$ bounds also change depending on the dimension. For $n \geq 5$, the  estimate \eqref{uab-bi-bal4} directly implies
\begin{enumerate}[resume, label=(\roman*)]
\item
$L^4_{t,x}$ bound
\begin{equation}\label{uL45}
\|u_{\lambda}\|_{L_{t,x}^4}\lesssim \epsilon c_{\lambda}\lambda^{-s_c+ \frac{d}4-\frac18}, \qquad n \geq 5.
\end{equation}
\end{enumerate}
On the other hand for $n = 4$ the $L^4_{t,x}$ 
bound remains borderline inaccessible, and we  have instead
\begin{enumerate}[resume,label=(\roman*)']
\item $L_t^4L_x^p$ bound: For $4<p\leq\infty$, 
\begin{equation}\label{uL4Lp4}
\|u_{\lambda}\|_{L_t^4L_x^p}\lesssim_p \epsilon c_{\lambda}\lambda^{-s_c - \frac18 +d(\frac12-\frac{1}{p})}, \qquad n = 4, 
\end{equation}
as well as
\begin{equation}\label{uL6L44}
\|u_{\lambda}\|_{L_t^{p} L_x^4}\lesssim  \epsilon c_{\lambda}\lambda^{-s_c  \frac{1}{2p}- \frac{d}2}, \qquad n=4.
\end{equation}
\end{enumerate}

\subsection{ Changes in the layered paradifferential expansions}

The algebraic form of our expansions is independent of the dimension, so it does not change. What changes is the interpretation of the terms, as follows:

\begin{itemize}
    \item In dimension $n \geq 5$, the $L^4_{t,x}$
    bound \eqref{uL45} allows us to estimate 
    the $N^{[3]}$ and $N^{[4]}$ layers perturbatively, so they can be omitted
    from all expansions, and the proof of the error bounds in the expansions simplifies significantly; in particular
    the null condition for $N^{[3]}$ is no longer needed.

    \item In dimension $n = 4$, the near $L^4_{t,x}$  bound \eqref{uL45} still allows us to estimate the $N^{[4]}$ layer perturbatively, but $N^{[3]}$ is still needed. However, only a simple null condition for $N^{[3]}$ is necessary.
\end{itemize}

\subsection{ Changes in the local well-posedness results} These stem from the 
manner in which the cubic term $N^{[3]}$ 
is treated in the study of the 
linearized equation, as follows:

\medskip

a) In dimension $n = 4$, the $N^{[3]}$
type contribution is nearly perturbative
due to the bound \eqref{uL6L44}. Hence 
the condition $c_\lambda \lesssim \lambda^{-\frac{1}{12}}$ in Theorem~\ref{t:local-fe3-lin}
can be relaxed to 
\begin{equation}
c_\lambda \lesssim \lambda^{-\delta}, \delta > 0,    
\end{equation}
with $\delta$ arbitrarily small. This leads to the corresponding improvement in Theorem~\ref{t:local}, which now only requires $s > s_c$.

\medskip

b) In dimension $n \geq 5$, the $N^{[3]}$
type contribution is fully perturbative
due to the $L^4_{t,x}$ bound \eqref{uL45}. Hence 
the condition $c_\lambda \lesssim \lambda^{-\frac{1}{12}}$ in Theorem~\ref{t:local-fe3-lin}
can be fully eliminated, which leads to the corresponding improvement in Theorem~\ref{t:local}, which now only requires $s \geq  s_c$. 

Another significant 
improvement occurs in the study of the linearized equation. In lower dimension, 
the need to estimate the $N^{[3]}$ term 
perturbatively forced a bounded time interval. This is no longer the case in higher dimension $n \geq 5$, where the linearized bounds in Theorem~\ref{t:local-fe3-lin} are now global in time.

\bibliography{1d-global.bib}

\bibliographystyle{plain}

\end{document}